%% file: These.tex
\documentclass[11pt,a4paper]{book}

\input{Preambule}

\begin{document}

\frontmatter
\include{PageGarde}

\include{resumes}

\tableofcontents
\mainmatter

\include{intro-0}

\include{french-intro}

\include{moduli-0}
\include{analytic-0}
\include{tate-0}
\include{autom-0}
\include{cartan-0}

\include{euler-0}

\appendix
\include{algeom-0}
\include{reps-0}
\include{local-0}


\nocite{*}
\bibliographystyle{amsplain}
\bibliography{biblio}

{\footnotesize \textsc{Elie STUDNIA} --- \begin{otherlanguage}{french} Universit\'e Paris Cit\'e, Sorbonne Universit\'e, CNRS, IMJ-PRG, F-75013 Paris, France \end{otherlanguage} }

\end{document}

%% file: Preambule.tex
\usepackage[utf8]{inputenc}
\usepackage[OT2,T1]{fontenc}

\usepackage[a4paper]{geometry}
\usepackage[french,english]{babel}

\usepackage{fancyhdr}
\usepackage[final]{graphicx}
\usepackage{numprint} 

\usepackage[figurewithin=none]{caption}
\usepackage{zref}
\usepackage{dsfont}
\usepackage{enumitem}
\usepackage{ifthen}
\usepackage{verbatim}

\usepackage[nottoc]{tocbibind}

\usepackage{setspace}
\usepackage{amsmath,amssymb,makeidx,amsthm}
\usepackage{stmaryrd}
\usepackage{mathrsfs}
\usepackage[all]{xy} 
\usepackage{multirow} 
\usepackage{array}
\usepackage{mathtools}
\usepackage{tikz-cd}

\usepackage{url}
\usepackage[hypertexnames=false]{hyperref}

\newcommand{\R}{{\mathbb R}}
\newcommand{\N}{{\mathbb N}}
\newcommand{\Z}{{\mathbb Z}}
\newcommand{\Q}{{\mathbb Q}}
\renewcommand{\C}{{\mathbb C}}
\newcommand{\HH}{{\mathbb H}}
\newcommand{\F}{{\mathbb F}}
\DeclareSymbolFont{cyrletters}{OT2}{wncyr}{m}{n}
\DeclareMathSymbol{\Sha}{\mathalpha}{cyrletters}{"58}
\newcommand{\rar}{\rightarrow}
\newcommand{\lrar}{\longrightarrow}

\renewcommand{\P}{\mathscr{P}}
\newcommand{\Ell}{\ensuremath{\mathbf{Ell}}}
\newcommand{\MP}{\ensuremath{\mathscr{M}(\P)}}
\newcommand{\MPp}[1]{\ensuremath{\mathscr{M}(#1)}}
\newcommand{\CMP}{\ensuremath{\overline{\mathscr{M}}(\P)}}
\newcommand{\CMPp}[1]{\ensuremath{\overline{\mathscr{M}}(#1)}}
\newcommand{\CP}{\ensuremath{C_{\P}}}
\newcommand{\CPp}[1]{\ensuremath{C_{#1}}}

\renewcommand{\emptyset}{\varnothing}
\newcommand{\mrm}[1]{\mathrm{#1}}

\newcommand{\OO}{\mathcal{O}}
\newcommand{\Qbar}{\overline{\mathbb{Q}}}
\newcommand{\Del}[1]{\ensuremath{\mrm{D}\left(#1\right)}}
\newcommand{\Delp}[2]{\ensuremath{\mrm{D}_{#1}\left(#2\right	)}}
\newcommand{\Fr}{\mrm{Frob}}
\newcommand{\Tate}[2]{\ensuremath{\mathcal{T}_{#1}}#2}

\renewcommand{\tilde}[1]{\widetilde{#1}}
\newcommand{\Aut}[1]{\ensuremath{\operatorname{Aut}(#1)}}
\newcommand{\Auti}[2]{\ensuremath{\operatorname{Aut}_{#1}(#2)}}
\newcommand{\GL}[1]{\ensuremath{\operatorname{GL}_2(#1)}}
\newcommand{\GLn}[2]{\ensuremath{\operatorname{GL}_{#1}(#2)}}
\newcommand{\SL}[1]{\ensuremath{\operatorname{SL}_2(#1)}}
\newcommand{\PGL}[1]{\ensuremath{\operatorname{PGL}_2(#1)}}
\newcommand{\PSL}[1]{\ensuremath{\operatorname{PSL}_2(#1)}}

\DeclareMathOperator{\Sp}{Spec}

\DeclareMathOperator{\Real}{Re}

\DeclareMathOperator{\Tr}{Tr}
\DeclareMathOperator{\im}{im}

\usepackage{calc}

\makeatletter
\let\oldlabel\label
\renewcommand{\label}[1]{%
  \zref@labelbylist{#1}{special}
  \oldlabel{#1}
}

\zref@newlist{special}
\zref@newprop{section}{\arabic{section}.\arabic{subsection}}
\zref@addprop{special}{section}

\newcounter{propo}
\renewcommand{\thepropo}{\thesection.\arabic{propo}}
\makeatletter \@addtoreset{propo}{section}\makeatother

\makeatletter \@addtoreset{rema}{chapter}\makeatother

\newcounter{theop}
\renewcommand{\thetheop}{\thechapter.\Alph{theop}}
\makeatletter \@addtoreset{theop}{chapter}\makeatother

\newcounter{questionintro}
\renewcommand{\thequestionintro}{\arabic{questionintro}}

\newcounter{conjectureintro}
\counterwithout*{conjectureintro}{chapter}
\renewcommand{\theconjectureintro}{\arabic{conjectureintro}}

\newcounter{theoremintro}
\renewcommand{\thetheoremintro}{\Alph{theoremintro}}

\newcounter{propositionintro}
\renewcommand{\thepropositionintro}{\arabic{propositionintro}}

\newcounter{mytheo}
\renewcommand{\themytheo}{\arabic{mytheo}}

\newcounter{mycor}
\renewcommand{\themycor}{\arabic{mycor}}
\makeatletter \@addtoreset{mycor}{mytheo}\makeatother

\newcounter{frquestionintro}
\renewcommand{\thefrquestionintro}{\arabic{frquestionintro}}

\newcounter{frtheoremintro}
\renewcommand{\thefrtheoremintro}{\Alph{frtheoremintro}}

\newcounter{frpropositionintro}
\renewcommand{\thefrpropositionintro}{\arabic{frpropositionintro}}

\newcounter{frmytheo}
\renewcommand{\thefrmytheo}{\arabic{frmytheo}}

\newcounter{frmycor}
\renewcommand{\thefrmycor}{\arabic{frmycor}}
\makeatletter \@addtoreset{frmycor}{frmytheo}\makeatother

\newcommand{\lem}[2][None]{\refstepcounter{propo}\ifthenelse{\equal{#1}{None}}{}{\label{#1}} \smallskip \noindent \textbf{Lemma \thepropo\;} {\it #2}\smallskip}
\newcommand{\nott}[1]{\smallskip \noindent \textbf{Notation \;} {\it #1} \smallskip}
\newcommand{\theo}[2][None]{\refstepcounter{propo}\ifthenelse{\equal{#1}{None}}{}{\label{#1}} \smallskip \noindent \textbf{Theorem \thepropo\;} {\it #2} \smallskip}
\newcommand{\theoi}[2][None]{\refstepcounter{theop}\ifthenelse{\equal{#1}{None}}{}{\label{#1}}\smallskip \noindent \textbf{Theorem \thetheop\;} {\it #2} \smallskip}
\newcommand{\cori}[2][None]{\refstepcounter{theop}\ifthenelse{\equal{#1}{None}}{}{\label{#1}}\smallskip \noindent \textbf{Corollary \thetheop\;} {\it #2} \smallskip}

\newcommand{\prop}[2][None]{\refstepcounter{propo}\ifthenelse{\equal{#1}{None}}{}{\label{#1}} \smallskip \noindent \textbf{Proposition \thepropo\;} {\it #2} \smallskip}

\newcommand{\cor}[2][None]{\refstepcounter{propo}\ifthenelse{\equal{#1}{None}}{}{\label{#1}} \smallskip \noindent \textbf{Corollary \thepropo\;} {\it #2} \smallskip}

\newcommand{\defi}[2][None]{\refstepcounter{propo}\ifthenelse{\equal{#1}{None}}{}{\label{#1}} \smallskip \noindent {\bf Definition \thepropo\;} {#2} \smallskip}

\newcommand{\rem}[2][None]{\refstepcounter{propo}\ifthenelse{\equal{#1}{None}}{}{\label{#1}}\smallskip \noindent {\bf Remark\;} {#2} \smallskip}
\newcommand{\rems}[1]{\smallskip \noindent {\bf Remarks\;} {#1} \smallskip}

\newcommand{\demo}[1]{\noindent {\it Proof.\;--\;} #1\hfill$\Box$ \smallskip}

\newcommand{\qnn}[2][None]{\refstepcounter{propo}\ifthenelse{\equal{#1}{None}}{}{\label{#1}} \smallskip \noindent \textbf{Question \thepropo\;} {\it #2} \smallskip}

\newcommand{\questin}[2][None]{\refstepcounter{questionintro}\ifthenelse{\equal{#1}{None}}{}{\label{#1}} \smallskip \noindent \textbf{Question \thequestionintro\;} {\it #2}\smallskip}
\newcommand{\conjintro}[2][None]{\refstepcounter{conjectureintro}\ifthenelse{\equal{#1}{None}}{}{\label{#1}} \smallskip \noindent \textbf{Conjecture \theconjectureintro\;} {\it #2}\smallskip}
\newcommand{\theoin}[2][None]{\refstepcounter{theoremintro}\ifthenelse{\equal{#1}{None}}{}{\label{#1}} \smallskip \noindent \textbf{Theorem \thetheoremintro\;} {\it #2}\smallskip}
\newcommand{\defintro}[1]{\smallskip \noindent \textbf{Definition \;} {\it #1}\smallskip}
\newcommand{\propintro}[2][None]{\refstepcounter{propositionintro}\ifthenelse{\equal{#1}{None}}{}{\label{#1}} \smallskip \noindent \textbf{Proposition \thepropositionintro\;} {\it #2}\smallskip}

\newcommand{\mytheo}[2][None]{\refstepcounter{mytheo}\ifthenelse{\equal{#1}{None}}{}{\label{#1}} \smallskip \noindent \textbf{Theorem \themytheo\;} {\it #2}\smallskip}
\newcommand{\mycor}[2][None]{\refstepcounter{mycor}\ifthenelse{\equal{#1}{None}}{}{\label{#1}} \smallskip \noindent \textbf{Corollary \themycor\;} {\it #2}\smallskip}

\newcommand{\questinfrench}[2][None]{\refstepcounter{frquestionintro}\ifthenelse{\equal{#1}{None}}{}{\label{#1}} \smallskip \noindent \textbf{Question \thefrquestionintro\;} {\it #2}\smallskip}
\newcommand{\theoinfrench}[2][None]{\refstepcounter{frtheoremintro}\ifthenelse{\equal{#1}{None}}{}{\label{#1}} \smallskip \noindent \textbf{Th\'eor\`eme \thefrtheoremintro\;} {\it #2}\smallskip}
\newcommand{\defintrofrench}[1]{\smallskip \noindent \textbf{D\'efinition \;} {\it #1}\smallskip}
\newcommand{\propintrofrench}[2][None]{\refstepcounter{frpropositionintro}\ifthenelse{\equal{#1}{None}}{}{\label{#1}} \smallskip \noindent \textbf{Proposition \thefrpropositionintro\;} {\it #2}\smallskip}

\newcommand{\mytheofrench}[2][None]{\refstepcounter{frmytheo}\ifthenelse{\equal{#1}{None}}{}{\label{#1}} \smallskip \noindent \textbf{Th\'eor\`eme \thefrmytheo\;} {\it #2}\smallskip}
\newcommand{\mycorfrench}[2][None]{\refstepcounter{frmycor}\ifthenelse{\equal{#1}{None}}{}{\label{#1}} \smallskip \noindent \textbf{Corollaire \thefrmycor\;} {\it #2}\smallskip}

%% file: PageGarde.tex

\begin{center}
\begin{tabular}{c@{\hskip 7cm}c@{\hskip 1cm}}
\includegraphics[height=2.5cm]{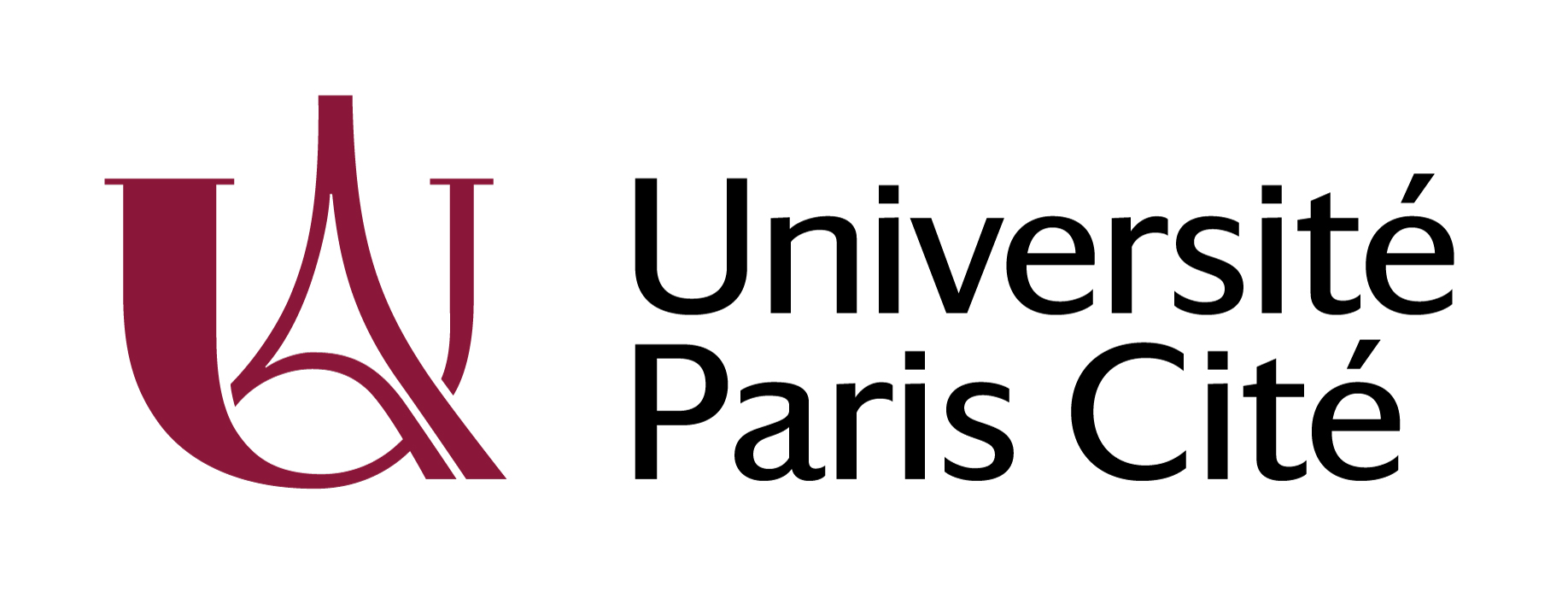} &
\end{tabular}
\end{center}
	
\begin{center}
		
\vspace*{.03\textheight}
\textsc{\LARGE Universit\'e Paris Cit\'e}\\[0.2cm] 
		\large \'Ecole doctorale Sciences Math\'ematiques de Paris-Centre (386)\\
		  IMJ-PRG \\ 
  
  			\vfill
 
	 		\rule{\textwidth}{0.8pt} \\ 
	 		\vspace{10pt}
	 		 { \LARGE \bfseries $L$-functions and rational points for Galois twists of modular curves} 
	 		 \vspace{10pt}
	 		 \rule{\textwidth}{0.8pt} \\ 
		\end{center}
		
		\vfill
		\begin{center}
			Par \textsc{\Large \'Elie Studnia}\\[1cm] 
			Th\`ese de doctorat de \textsc{\large Math\'ematiques}\\[1.2cm]
			Dirig\'ee par \textsc{\large Lo\"ic MEREL}\\[0.2cm]
			Présentée et soutenue publiquement le 06/11/2024 
		\end{center}
		
		\vspace{1cm}
Devant un jury compos\'e de : 
		\begin{center}

\begin{tabular}{lll}
M\textsuperscript{me} Anna \textsc{Cadoret}, Professeur des Universit\'es & Sorbonne Universit\'e  & examinatrice \\
M. Lo\"ic \textsc{Merel}, Professeur des Universit\'es & IMJ-PRG & directeur \\
M\textsuperscript{me} Giada \textsc{Grossi}, CR & Universit\'e Sorbonne Paris Nord & examinatrice \\
M. Vincent \textsc{Pilloni}, DR & Universit\'e Paris-Saclay & examinateur \\
M. Jan \textsc{Vonk}, Professor & Universiteit Leiden & examinateur\\
M. Pierre \textsc{Parent}, MCF-HDR & Universit\'e de Bordeaux & rapporteur \\
M. Henri \textsc{Darmon}, Professor & McGill University & rapporteur (excus\'e)
\end{tabular}

\end{center}


\newpage

\vspace*{\fill}

\noindent\begin{center}
\begin{minipage}[t]{(\textwidth-2cm)/2}
Institut de math\'ematiques de Jussieu-Paris Rive gauche. UMR 7586. \\
Bo\^ite courrier 247 \\
4 place Jussieu \\
\numprint{75252} Paris Cedex 05
\end{minipage}
\hspace{1.5cm}
\begin{minipage}[t]{(\textwidth-2cm)/2}
 Universit\'e Paris Cit\'e. \\
 \'Ecole doctorale de sciences\\
 math\'ematiques de Paris centre.\\
  Bo\^ite courrier 290 \\
 4 place Jussieu\\
 \numprint{75252} Paris Cedex 05
 \end{minipage}
\end{center}

%% file: resumes.tex
\chapter*{R\'esum\'es}

\doublespacing
\begin{center}{ \LARGE \bf
Fonctions L et points rationnels de tordues galoisiennes de courbes modulaires
}
\end{center}
\vspace*{2cm}

\onehalfspacing
\noindent { \Large \textbf{R\'esum\'e}}
\vspace*{0.75cm}

Dans cette th\`ese, on discute de l'existence de courbes elliptiques sur une base g\'en\'erale pour lesquelles la p-torsion est
un sch\'ema en groupes fix\'e (o\`u p est un nombre premier). On prouve tout d'abord que ce probl\`eme est repr\'esentable au
sens de Katz-Mazur par une courbe affine lisse, puis on d\'ecrit une autre manière de construire sa compactification sur un
corps, comme la tordue de X(p) par une repr\'esentation galoisienne bidimensionnelle modulo p. On \'etudie ensuite la
fonction L d'une telle courbe modulaire tordue sur le corps des rationnels et le signe de l'équation fonctionnelle
conjecturalement satisfaite par ses facteurs. Lorsque la représentation galoisienne poss\`ede une petite image, on affine
cette décomposition et on montre que les fonctions L de chacun des facteurs satisfont des \'equations fonctionnelles dont
on d\'etermine le signe. Gr\^ace au syst\`eme d'Euler de Beilinson-Flach, et en \'etendant des résultats de grande image
galoisienne d\^us \`a Loeffler, on montre la conjecture de Bloch-Kato en rang z\'ero pour chaque facteur sans multiplication
complexe, et on applique ce r\'esultat pour contr\^oler explicitement la mauvaise r\'eduction de courbes elliptiques pour
laquelle l'image de Galois est contenue dans le normalisateur d'un sous-groupe de Cartan non d\'eploy\'e, sous une
hypothèse suppl\'ementaire de non-annulation de certaines fonctions L de Rankin-Selberg, que l'on sait vérifier dans
certains cas.

\vspace*{1cm}

\subsection*{Mots-cl\'es}
\noindent points rationnels, courbes modulaires, formes modulaires, courbes elliptiques, repr\'esentations galoisiennes, fonctions
L, syst\`emes d'Euler, modularit\'e

\newpage


\doublespacing
\begin{center}{ \LARGE \bf
L-functions and rational points on Galois twists of modular curves
}
\end{center}
\vspace*{2cm}

\onehalfspacing
\noindent { \Large \textbf{Abstract}}
\vspace*{0.75cm}

In this thesis, we discuss the existence of elliptic curves over a general base scheme for which the p-torsion is a given
finite \'etale scheme, where p is a given prime. We first prove that this problem is representable in the sense of Katz and
Mazur by a smooth affine curve, then describe another way to construct the compactification of such a curve over a base
field, as a twist of X(p) by a two-dimensional Galois representation modulo p. We then study the L-function of such a
twisted modular curve over the field of rationals and compute the sign of the functional equation that its factors are
conjectured to satisfy. When the aforementioned Galois representation has small image, we show that this decomposition
can be refined further, and we prove that the associated L-functions satisfy a functional equation whose sign we also
determine. We then demonstrate the rank zero Bloch-Kato conjecture for every such factor without complex multiplication
using the Euler system of Beilinson-Flach elements and extending big Galois image results due to Loeffler. We finally
provide an application of such a result by giving an explicit bound on the bad reduction of elliptic curves such that the
image of the Galois action on its p-torsion is contained in the normalizer of a nonsplit Cartan subgroup, under the
additional hypothesis that certain Rankin-Selberg L-functions do not vanish, which is known to hold in certain cases.

\vspace*{1cm}

\subsection*{Keywords}
\noindent Rational points, modular curves, modular forms, Galois representations, Euler systems, elliptic curves, L-functions, modularity

\singlespacing
\newpage

%% file: intro-0.tex
\chapter*{Introduction}
\phantomsection
\addcontentsline{toc}{chapter}{Introduction}
\renewcommand{\thesection}{\arabic{section}}
\renewcommand{\thesubsection}{\thesection.\arabic{subsection}}
\section{Existing work}

\subsection{Background}

Given an elliptic curve $E/\Q$ and an integer $N \geq 1$, the subgroup $E[N]$ of $E(\overline{\Q})$ made with its $N$-torsion points is endowed with an action of $G_{\Q}$ by group automorphisms. By choosing a basis of the free $\Z/N\Z$-module $E[N]$ of rank two, this defines a continuous group homomorphism $\rho_{E,N}: G_{\Q} \rar \GL{\Z/N\Z}$ which encodes important information on the arithmetic of $E$. 

If $\varphi: E \rar F$ is an isogeny of elliptic curves over $\Q$ with degree prime to $N$, it induces an isomorphism $\rho_{E,N} \simeq \rho_{F,N}$. The Frey-Mazur conjecture originates from the following question asked by Mazur in \cite{FreyMazur}:

\questin[mazur-qn]{Does there exist an integer $N \geq 7$ and non-isogenous elliptic curves $E, F$ over $\Q$ such that $\rho_{E,N}$ and $\rho_{F,N}$ are conjugate? }

It is thus a statement about a \emph{comparison} of two different Galois representations. Such a problem was already considered by Serre \cite[Th\'eor\`emes 6", 7]{Serre-image-ouverte}, who showed that when $E, F$ are not isogenous over $\Qbar$\footnote{Serre's condition is slightly different, but he writes in \emph{loc.cit.} that he expects it to be equivalent, which was subsequently proved by Faltings \cite[Corollary 2 to Theorem 4]{Faltings86}.}, then the intersection of the fields $\Q(E_{tors})$ and $\Q(F_{tors})$ is a number field. This implies in particular that the set $S$ of integers $N \geq 1$ such that $\rho_{E,N}$ is conjugate to $\rho_{F,N}$ is finite. 

Serre also provides a explicit bound on $\max{S}$ in his paper about Cebotarev's theorem:

\theoin[serre-cebo]{(Serre \cite[Th\'eor\`eme 21]{Serre-Cebotarev}) Assume the generalized Riemann hypothesis. Let $E,F/\Q$ be non-isogenous elliptic curves and let $N_{E,F}$ denote the product of the primes of bad reduction for $E$ or $F$. There is a universal constant $c > 0$ such that one has  
\[\max{S} \leq c\log{N_{E,F}}(\log{\log{2N_{E,F}}})^6.\]
}

The key difference with Mazur's question is that Serre's result is \emph{non-uniform}: his bound depends on both elliptic curves $E,F$, while, in Mazur's question, the bound on $N$ is absolute, i.e. he asks for $N \geq 7$. 

A second motivation to discuss Mazur's question arises from the study of ternary Diophantine equations. The historical example is the study of the Fermat equation: 

\theoin[frey-on-fermat]{(Frey \cite[III.2]{Frey86}) Assume that modular representations satisfy level-lowering\footnote{This is Conjecture \underline{S} in \emph{loc.cit.}, which was proved shortly afterwards by Ribet \cite{Ribet90}.} and that semi-stable elliptic curves over $\Q$ are modular. Then, for any prime $p \geq 3$, the Fermat equation $x^p+y^p=z^p$ has no rational solution $(x,y,z)$ with $xyz \neq 0$.}

The main idea linking elliptic curves over $\Q$ to Diophantine equations comes from the so-called \emph{Frey-Hellegouarch curve}, first introduced by Hellegouarch \cite{HelleCrelle, THelle}. Let $a,b,c \in \Z$ be pairwise coprime integers such that $a+b+c=0$ and consider the elliptic curve $E/\Q$ with equation $y^2=x(x-a)(x+b)$ (with discriminant $(4abc)^2$). It can be checked (see for instance \cite[\S 4.1]{Conj-Serre}) that $E$ has good reduction at any prime $p \nmid 2abc$, semi-stable reduction at any odd prime $p \mid abc$. Let $p,\ell$ be distinct odd prime numbers such that the $\ell$-adic valuations of $a,b,c$ are all divisible by $p$ and not all zero. A local study of the elliptic curve $E/\Q_{\ell}$ \cite[IV, Appendix A.1]{SerreMcGill} shows that the representation $\rho_{E,p}$ is unramified at $\ell$. 

In particular, in the case of the Fermat equation where $(a,b,c)=(u^p,v^p,w^p)$, the representation $\rho_{E,p}$ is only ramified at $2$ and $p$; a more precise study of ramification at these primes (see for instance \cite[\S 4.2]{Conj-Serre}) shows that, in this situation, $\rho_{E,p}$ cannot be modular\footnote{The importance of level-lowering in this argument is that it states that if $\rho_{E,p}$ is modular, it is so \emph{in the smallest possible level}, which would be $2$.}.  

In a later work \cite{FreyFLT}, Frey establishes a similar connection between a generalized class of ternary Diophantine equations and isomorphisms between the $\rho_{E,N}$ for various elliptic curves $E/\Q$. 

\theoin[frey-on-mazur]{(Frey \cite[Proposition 3.5]{FreyFLT})
Let $S$ be a a finite set of prime numbers, and $\Z_S$ denote the set of $n \in \Z$ such that every prime divisor of $n$ is in $S$. Let $\mathcal{E}_S$ be the set of isomorphism classes of elliptic curves over $\Q$ such that their conductor is a square-free product of elements of $S$. Then the following statements are equivalent:
\begin{itemize}[noitemsep,label=$-$]
\item The set of $(x,y,z) \in \Z^3$ such that $(x,y,z)=1$ and there exists $a,b,c \in \Z_S, p \geq 5$ prime such that $ax^p+by^p+cz^p=0$ is finite.
\item There are finitely many elliptic curves $E/\Q$ (up to isomorphism) such that $E[2] \subset E(\Q)$ and $\rho_{E,p} \simeq \rho_{E_0,p}$ for some $E_0 \in \mathcal{E}_S$ and some prime $p \geq 5$.   
\end{itemize}} 

A theorem of Shafarevich \cite[IV, \S 1.4]{SerreMcGill} shows that, in the statement of Theorem \ref{frey-on-mazur}, $\mathcal{E}_S$ is finite. Moreover, there does not seem to be any conceptual reason for the triviality requirement on the Galois representation $\rho_{E,2}$: up to a small modification that we will discuss later, its presence or absence should not be meaningful. The same applies to the semi-stability condition on $E_0$. One possible statement of the Frey-Mazur conjecture is the following slight generalization of the second condition of Theorem \ref{frey-on-mazur}. 

\conjintro[fm-weak]{(See \cite[Conjecture 5]{FreyFLT}) Let $E/\Q$ be an elliptic curve. There are finitely many elliptic curves up to isomorphism $F/\Q$ such that the representations $\rho_{E,p}$ and $\rho_{F,p}$ are conjugate for some prime $p \geq 7$.}

It is in fact expected \cite[Conjecture 2.1]{Freitas-Kraus} that a uniform version of this conjecture holds. 

\conjintro[fm-unif]{There exists a constant $C_{FM}>0$ such that for any elliptic curves $E, F/\Q$ such that $\rho_{E,p}\simeq \rho_{F,p}$ for some prime $p > C_{FM}$, $E$ and $F$ are isogenous.} 

In \cite[\S 4.3]{Conj-Serre}, Serre discusses the following attempt towards solving the Diophantine equation $u^p+v^p+31w^p=0$, where $p \geq 5$ is a prime distinct from $31$, with \emph{the modular method}. If $(u,v,w)$ is a nontrivial solution, then, up to dividing by their greatest common divisor, we may assume that $u,v,w$ are pairwise coprime. Let $(A,B,C)$ be a permutation of $(u^5,v^5,31w^5)$ or $(-u^5,-v^5,-31w^5$) such that $B$ is even (hence divisible by $32$) and $A \equiv -1 \pmod{4}$. Let $E$ be the elliptic curve with equation $y^2=x(x-A)(x+B)$: by \cite[\S 4.1]{Conj-Serre}, it is semi-stable with full $2$-torsion and its primes of bad reduction are exactly those that divide $ABC$. Its minimal discriminant is $31^2\left(\frac{ABC}{31 \cdot 2^p}\right)^{2} 2^{2p-8}$. Therefore, by Serre's Conjecture \cite{KW1,KW2}, $\rho_{E,p}$ is the reduction of a newform $f \in \mathcal{S}_2(\Gamma_0(62))$. 

There are two Galois orbits of newforms in $\mathcal{S}_2(\Gamma_0(62))$, of which exactly one is rational, corresponding to the newform $f_{\Q}$\footnote{It is the newform $F$ in \cite[\S 4.3 Remarque (3)]{Conj-Serre}.}. A local argument at $\ell\in \{3,5\}$ as in \cite[\S 4.3 Lemme 1]{Conj-Serre}, shows, using the $q$-expansion of the newform with non-rational $q$-expansion given by the LMFDB \cite{lmfdb} that $f=f_{\Q}$, hence $E$ is congruent modulo $p$ to some elliptic curve of conductor $62$. Let $E_0$ be the only elliptic curve of conductor $62$ with full $2$-torsion by the LMFDB data \cite{lmfdb}: it is given by the minimal Weierstrass equation $y^2+xy+y=x^3-x^2-21x+41$ and (minimal) discriminant $62^2$. If Conjecture \ref{fm-weak} holds for $E_0$ then, when $p$ is large enough, the elliptic curve $E$ must be isogenous to $E_0$, hence isomorphic to $E_0$ since $E_0$ is the only elliptic curve of conductor $62$ with $E_0[2] \subset E_0(\Q)$. The minimal discriminants are equal up to a unit, thus $ABC=31 \cdot 32$, $p=5$, and $(u,v,w) \in \{(2,-1,-1),(-1,2,-1),(-2,1,1),(1,-2,1)\}$.\footnote{Since the submission of this thesis, Samir Siksek told me about the following results. First, as an application of the symplectic congruence criterion (Proposition \ref{ko-prop2}), Halberstadt and Kraus \cite[Proposition 2.2]{HK-Fermat} proved that the equation had no solution when $p \equiv 3\pmod{4}$. Using a different argument introduced by Kraus in \cite{Kraus-a3b3cp} to study the Diophantine equation $a^3+b^3=c^p$ (which works prime by prime and does not seem to have an obvious modular interpretation), Siksek proves in \cite[Theorem 12]{Siksek-PS} that the equation has no solution when $11 \leq p \leq 10^6$. These results are discussed in the PhD thesis of Heline Deconinck \cite[Chapter 4]{Deconinck-PhD}, which also rules out a further $1/12$-th of possible congruence classes for a prime $p$ modulo $3 \cdot 5 \cdot 7 \cdot 11$ such that the equation has a nontrivial solution.}

\subsection{Existing work about congruences of elliptic curves}

In \cite{FreyMazur}, Mazur suggests that Question \ref{mazur-qn} can be rephrased as a problem about finding rational points on certain twists of the modular curves $X(N)$. Let us make that observation more precise. Let $E, F$ be two elliptic curves over a field $K$, and $N \geq 1$ be an integer not divisible by the characteristic of $K$. Then $\rho_{E,N}$ and $\rho_{F,N}$ are conjugate if, and only if, there is an isomorphism $\varphi: E[N] \rar F[N]$ of finite \'etale group schemes over $K$. Both $E[N]$ and $F[N]$ are endowed with the bilinear perfect alternated Weil pairing $\langle -,\,-\rangle$, so there is a unique $\alpha \in (\Z/N\Z)^{\times}$ such that, for any $P,Q\in E[N](\overline{K})$, $\langle \varphi(P),\,\varphi(Q)\rangle = \alpha\langle P,\,Q\rangle$. We write $\alpha = \det{\varphi}$. For instance, if $\varphi$ comes from an isogeny $E \rar F$ of degree $d$ coprime to $N$, then $\det{\varphi}=d$. 

For $N \geq 3$ and $\alpha \in (\Z/N\Z)^{\times}$, there exists\footnote{Such a result seems to have been known to Mazur in \cite{FreyMazur}, but we are only aware of partial statements in the literature such as \cite[Proposition 1]{KO}, or \cite[I.7]{Mazur-Open} in the case of number fields. We prove a generalization of this statement in Chapter \ref{moduli-spaces}.} a smooth projective geometrically connected curve $X_E^{\alpha}(N)$ with a map $j: X_E^{\alpha}(N) \rar \mathbb{P}^1_K$ over $K$ such that $j^{-1}(\mathbb{A}^1_K)$ parametrizes the isomorphism classes of couples $(F/U,\iota)$, where $F \rar U$ is an elliptic curve over the $K$-scheme $U$ and $\iota: E[N]_U  \rar F[N]_U$ is an isomorphism of finite \'etale group schemes with determinant $\alpha$. In this parametrization, the morphism $j$ computes the $j$-invariant of the elliptic curve $F/U$. 

Let $Y_E^{\alpha}(N)$ be the open subscheme $j^{-1}(\mathbb{A}^1_K)$ of $X_E^{\alpha}(N)$. Its rational points are called the \emph{non-cuspidal} points of $X_E^{\alpha}(N)$. Conversely, the points of $X_E^{\alpha}(N) \backslash Y_E^{\alpha}(N)$ are called the \emph{cusps}. In particular, $X_E^1(N)$ always possesses a $K$-rational point attached to $(E/K,\mrm{id})$. Furthermore, an isogeny $\varphi: E \rar F$ of elliptic curves over $K$ with degree $d$ coprime to $N$ defines a non-cuspidal $K$-rational point in $X_E^d(N)$, as well as, for any $\alpha \in (\Z/N\Z)^{\times}$, an isomorphism $X_F^{\alpha}(N) \rar X_E^{d\alpha}(N)$.  

When $K=\Q$, as noted in \cite{FreyMazur}, the curve $X_{E_{\C}}^{\alpha}(N)$ is isomorphic to the classical curve $X(N)_{\C}$: indeed, because $\C$ is an algebraically closed field, the $N$-torsion subscheme of any elliptic curve over $\C$ is constant. It is in this sense that the $X_E^{\alpha}(N)$ are twists of $X(N)$. In particular, this lets us compute the genus of $X(N)$: it is zero when $N\leq 5$ (hence the requirement that $p \geq 7$ in Conjecture \ref{fm-weak}), one when $N=6$, and at least $3$ when $N \geq 7$. 

Hence, by Faltings's Theorem \cite[Theorem 7]{Faltings86}, every $X_E^{\alpha}(N)$ for $N \geq 7$ has finitely many rational points, so that Conjecture \ref{fm-weak} is equivalent to 

\conjintro[fm-weak-modular]{Let $E/\Q$ be an elliptic curve. There is a constant $C_E$ such that for any prime $p > C_E$ and any $\alpha \in \F_p^{\times}$, the non-cuspidal rational points of $X_E^{\alpha}(p)$ come from elliptic curves isogenous to $E$. }

Thus Conjecture \ref{fm-unif} is equivalent to 

\conjintro[fm-strong-modular]{There exists a constant $C > 0$ such that, for any elliptic curve $E/\Q$ and any $\alpha \in \F_p^{\times}$, for any prime $p > C$, the non-cuspidal rational points of $X_E^{\alpha}(p)$ come from elliptic curves isogenous to $E$.}

From now on, we assume, as in the statement of Conjecture \ref{fm-weak-modular}, that $N=p$ is a prime number. 

In order to relate Question \ref{mazur-qn} to Conjecture \ref{fm-weak-modular}, it is useful to be able to determine, given two elliptic curves $E, F/\Q$, 
\begin{enumerate}[noitemsep,label=(\roman*)]
\item\label{intro-wheniso} whether there is an isomorphism of group schemes $E[p] \simeq F[p]$, or in other words whether $\rho_{E,p} \simeq \rho_{F,p}$,
\item\label{intro-ifsymplectic} whenever it is the case, what the determinants of isomorphisms $\iota: E[p] \rar F[p]$ are.   
\end{enumerate}

\defintro{In the situation \ref{intro-wheniso}, we say that the elliptic curves $E,F$ are \emph{congruent} modulo $p$. }

The article \cite{KO} by Kraus and Oesterl\'e seems to be the first one discussing both of these questions and answering Question \ref{mazur-qn} for $p=7$. The article points out a practical answer to \emph{disprove} the existence of an isomorphism in many cases: it is sufficient (see Proposition 3 of \emph{loc.cit.}) to find a prime $\ell$ of common good reduction for $E$ and $F$ such that $a_{\ell}(E) \not\equiv a_{\ell}(F) \pmod{p}$. Indeed, one has $a_{\ell}(E) \equiv a_{\ell}(F) \pmod{p}$ for every prime $\ell$ of common good reduction for $E$ and $F$ if and only if the representations $\rho_{E,p}$ and $\rho_{F,p}$ have isomorphic semi-simplifications, which holds when\footnote{This is an equivalence if $\rho_{E,p}$ or $\rho_{F,p}$ is irreducible, which is the ``generic'' situation, as shown by Theorem \ref{mazur-y0}.} $E$ and $F$ are congruent modulo $p$. Under the assumptions that $E, F$ are modular\footnote{which is now known by works of Breuil, Conrad, Diamond, Taylor and Wiles \cite{wiles, moddiamond, mod3, fullmod}.}, the authors of \emph{loc.cit.} make this criterion effective. 

\propintro[ko-prop4]{(Kraus-Oesterl\'e \cite[Proposition 4]{KO}) Let $E, F$ be two modular elliptic curves over $\Q$. There is an effective bound $M = M(E,F)$ such that $\rho_{E,p}$ and $\rho_{F,p}$ have the same semi-simplification if, and only if, for any prime $\ell \leq M(E,F)$, the following two conditions are true:
\begin{itemize}[noitemsep,label=$-$]
\item if $\ell$ is a prime of common good reduction for $E,F$, one has $a_{\ell}(E) \equiv a_{\ell}(F) \pmod{p}$,   
\item if $\ell$ is a prime of good reduction for exactly one of the curves $E,F$, and of bad multiplicative reduction for the other one, then one has $a_{\ell}(E)a_{\ell}(F) \equiv \ell+1\pmod{p}$. 
\end{itemize}}

Replacing $\rho_{E,p}$ with its semi-simplification is usually not a problem: by work of Mazur \cite[Theorem 1]{FreyMazur}, the representation $\rho_{E,p}$ is irreducible when $p \geq 11$ and $p \neq 13$, except for a finite number of well-known cases. 

The authors of \emph{loc.cit.} also discuss \ref{intro-ifsymplectic}. Indeed, suppose that the group schemes $E[p]$ and $F[p]$ are isomorphic: we can multiply any isomorphism between them by $d \in \F_p^{\times}$, which multiplies its determinant by $d^2$. Thus, the set of determinants of isomorphisms $\iota: E[p] \rar F[p]$ is a subset of $\F_p^{\times}$ stable under multiplication by the subgroup $\F_p^{\times 2}$ of squares in $\F_p^{\times}$. 

\defintro{Let $E, F$ be two elliptic curves over $\Q$ which are congruent mod $p$. When there is an isomorphism $\iota: E[p] \rar F[p]$ of group schemes of determinant $1$, we say that $E[p]$ and $F[p]$ are \emph{symplectically isomorphic}, or that $E, F$ are \emph{symplectically congruent} modulo $p$. When there is an isomorphism $\iota: E[p] \rar F[p]$ of group schemes of determinant in $\F_p^{\times} \backslash \F_p^{\times 2}$, we say that $E[p]$ and $F[p]$ are \emph{anti-symplectically isomorphic}, or that $E, F$ are \emph{anti-symplectically} congruent modulo $p$.}

Note that it is possible for two elliptic curves $E, F$ over $\Q$ to be both symplectically congruent modulo $p$ and anti-symplectically congruent modulo $p$. By \cite[Proposition 1.1]{Cremona-Freitas} and \cite[\S 5.2 (iv)]{Serre-image-ouverte}, this is only possible if the image of $\rho_{E,p}$ is contained in a split Cartan subgroup. 

In order to find out whether two elliptic curves that are congruent modulo $p$ are symplectically or anti-symplectically congruent modulo $p$, the authors of \cite{KO} consider the invariants of the elliptic curve in \emph{bad reduction}, that is, attached to the restriction of $\rho_{E,p}$, $\rho_{F,p}$ to an inertia group at some prime $\ell$. 

\propintro[ko-prop2]{(Kraus-Oesterl\'e \cite[Proposition 2]{KO}) Let $E, F$ be two elliptic curves over $\Q$ with minimal discriminants $\Delta_E, \Delta_F$, and $p,\ell$ be two distinct prime numbers such that 
\begin{itemize}[noitemsep,label=$-$]
\item the elliptic curves $E, F$ both have bad multiplicative reduction at $\ell$, 
\item the group schemes $E[p]$ and $F[p]$ are isomorphic,
\item $v_{E,\ell}:= v_{\ell}(\Delta_E)$ is not divisible by $p$. 
\end{itemize}
Then $v_{F,\ell} := v_{\ell}(\Delta_F)$ is not divisible by $p$. Moreover, $E[p]$ and $F[p]$ are symplectically isomorphic if and only if $\left(\frac{v_{E,\ell}v_{F,\ell}}{p}\right)=1$.} 

Such local criteria have subsequently been developed in other settings, in particular in potential good reduction. We refer to \cite{Freitas-Kraus} for a list. %

Let $E, F$ be elliptic curves over $\Q$ that are congruent modulo $p$. It is possible by \cite[Proposition 1.1]{Cremona-Freitas} that one cannot determine, using solely local methods, whether the congruence is symplectic or anti-symplectic. A solution, suggested in \S 1.1 of \emph{loc.cit.}, is to use explicit models for the curves $X_E^{\alpha}(p)$. More precisely, for small values of $p$, one can construct part or all of the following data:
\begin{itemize}[noitemsep,label=$-$]
\item explicit equations for a family $X^{\alpha}(p)$ of proper curves over the ring $\Q[a,b]_{4a^3+27b^2}$ such that its specialization at the maximal ideal $(a-a_0,b-b_0)\subset \Q[a,b]_{4a^3+27b^2}$ for $(a_0,b_0) \in \Q^2$ yields a model of $X_E^{\alpha}(p)$ for the elliptic curve $E/\Q$ given by the Weierstrass equation $y^2=x^3+a_0x+b_0$.     
\item explicit equations for the $j$-invariant map $X^{\alpha}(p) \rar \mathbb{P}^1_{\Q[a,b]_{4a^3+27b^2}}$,
\item for any non-cuspidal point $P \in X^{\alpha}(p)(K)$ above the maximal ideal $(a-a_0,b-b_0)$ for $a_0, b_0 \in K$ for some number field $K$, corresponding to an isomorphism $E[p] \rar E'[p]$ of group schemes over $K$ with determinant $\alpha$, where $E$ denotes the elliptic curve $y^2=x^3+a_0x+b_0$ over $K$, a specification (for instance using the invariants $c_4, c_6$ of $E'$) of the $K$-isomorphism class of $E'$ \footnote{Indeed, when for instance $j(P)$ is neither $0$ or $1728$, the previous items only determine $E'$ up to quadratic twist.}.  
\end{itemize}

When $p \in \{3,5\}$ and $\alpha \in \F_p^{\times 2}$ (the \emph{symplectic case}), this was solved by Rubin and Silverberg \cite{Rubin-Silverberg}. An explicit model for $X_E^1(7)$ was first constructed by Halberstadt and Kraus \cite{HalbKraus}, and an explicit equation for $X_E^{-1}(7)$ was first described by Poonen, Schaefer and Stoll \cite{PSS}. These formulae were proved by a different method and completed by work of Fisher: he describes equations for $X_E^{\alpha}(p)$ for any $\alpha$ when $p \in \{3,5\}$ \cite{Fisher5}, when $p \in \{7,11\}$ \cite{Fisher-711}, when $p=13$ \cite{Fisher-13}.  

As Halberstadt and Kraus demonstrate in \cite{HalbKraus}, explicit models for the curve $X_E^{\alpha}(p)$ turn out to be convenient means for exhibiting congruences of elliptic curves modulo $p$, or even families of such congruences. Indeed, they describe an infinite family of $6$-uples of pairwise non-isogenous elliptic curves over $\Q$ which are pairwise symplectically congruent modulo $7$. In \cite{Fisher-711} and \cite{Fisher-13}, Fisher constructs infinite families of symplectic congruences for $p \in \{7, 11, 13\}$ and infinite families of anti-symplectic congruences for $p \in \{7,11,13\}$. 

In \cite[\S 3.7]{Cremona-Freitas}, the authors prove that there is only one pair of elliptic curves congruent modulo $17$ with both conductors smaller than $500 000$ (up to quadratic twist), given by the LMFDB labels \cite{lmfdb} $3675.e1$, $47775.q1$, and that the congruence is anti-symplectic. Fisher exhibits in \cite{Fisher-17} the first known example of symplectic congruence modulo $17$, with conductors $279809270 = 2\cdot 5\cdot 13\cdot 59\cdot 191^2$, and $3077901970 = 2\cdot 5\cdot 11 \cdot 13 \cdot  59 \cdot 191^2$, and conjectures in the same article that these are the only congruences between elliptic curves over $\Q$ modulo a prime $q \geq 17$. This conjecture is somewhat supported by \cite[Theorem 1.3]{Cremona-Freitas}, which states that this conjecture holds for elliptic curves of conductor at most $500 000$. 

To find the symplectic congruence modulo $17$ which we just discussed (and, in fact, the families of congruences modulo $13$ of \cite{Fisher-13}, as well as the family of anti-symplectic congruences modulo $11$, crucially relying on the work of \cite{Kumar}), Fisher does not in fact use the models for any of the $X_E^1(17)$. Instead, he directly considers a moduli space $Z(p,\alpha)$ for triples $(E,F,\iota)$ of elliptic curves such that $\iota: E[p] \rar F[p]$ is an isomorphism of determinant $\alpha$, known as a \emph{modular diagonal surface} first studied by Kani and Schanz \cite{MDQS}. The latter authors prove, among others, the following classification:

\theoin[kani-schanz-classif]{(Kani-Schanz \cite[Theorem 4]{MDQS})
\begin{itemize}[noitemsep,label=$-$]
\item $Z(p,\alpha)$ is a rational surface if and only if $p \leq 5$ or $p=7$ and $\alpha \in \F_p^{\times 2}$,
\item $Z(p,\alpha)$ is birational to an elliptic $K3$-surface if and only if $p =7$ and $\alpha \notin \F_p^{\times 2}$,
\item $Z(p,\alpha)$ is birational to an elliptic surface of Kodaira dimension $1$ if and only if $p=11$ and $\alpha \in \F_p^{\times 2}$,
\item $Z(p,\alpha)$ is a surface of general type in all other cases.  
\end{itemize}}

By studying the geometry of $Z(11,1)$, Kani and Rizzo \cite{KR} were able to prove that there were infinitely many families of symplectic congruences of elliptic curves (over $\Q$) modulo $11$.  

As explained in \cite{MDQS}, a consequence of this classification is that when $p \geq 13$, according to Lang's conjecture (see for instance \cite[Conjecture B]{CHM}), the rational points on $Z(p,\alpha)$ should be located on finitely many rational or elliptic curves (i.e. the curves with genus $g \leq 1$) up to finitely many exceptions. 

\theoin[bats-hit]{(Bakker-Tsimerman \cite[Theorem 29]{BaTs}) Let $N \geq 1$. There is an integer $M_N$ with the following property. Let $p > M_N$ be a prime number, $k$ be an algebraically closed field of characteristic zero. Let $B$ be a smooth projective curve over $k$ of gonality at most $N$ endowed with a nonconstant morphism $B \rar Z(p,\alpha)_k$. Then $B$ factors through a Hecke correspondence.}

Even assuming Lang's conjecture, to prove the Frey-Mazur conjecture in this case, one would have to prove that there are no ``exceptional'' points for large enough $p$. 

\subsection{Finding rational points on modular curves: methods and results}

There is a long-standing strategy about finding rational points of modular curves, pioneered by Mazur in \cite{MazurY1, FreyMazur} for the study of the modular curves $X_1(p)$ and $X_0(p)$. An advantage of this approach is that it enables one to find rational points for the curve of interest without need for an explicit model, which has been used to discuss the following question. 

\questin[serre-unif]{(Serre's uniformity question \cite[\S 4.3]{Serre-image-ouverte}) Does there exist a constant $C>0$ such that for any elliptic curve $E/\Q$ without complex multiplication and any prime $p > C$, the representation $\rho_{E,p}$ is surjective?}

As explained in the beginning of \cite{MazurY1}, local arguments by Serre \cite{Serre-image-ouverte} show that it is enough to rule out (when $p > 13$) the following three cases:
\begin{enumerate}[noitemsep,label=(\alph*)]
\item\label{suq-0} the representation $\rho_{E,p}$ is reducible,
\item\label{suq-s} the image of $\rho_{E,p}$ is contained in the normalizer of a split Cartan subgroup,
\item\label{suq-ns} the image of $\rho_{E,p}$ is contained in the normalizer of a nonsplit Cartan subgroup. 
\end{enumerate}

In all three situations, for a given $p$, there is a modular curve parametrizing elliptic curves with Galois image contained in the given subgroup. The focus of the articles \cite{MazurY1,FreyMazur} is \ref{suq-0}, (although they prove partial results regarding \ref{suq-s} as well).

\theoin[mazur-y0]{(Mazur \cite{FreyMazur}) Let $E/\Q$ be an elliptic curve such that $\rho_{E,p}$ is reducible for some prime $p > 37$. Then $E$ has complex multiplication by $\Q(\sqrt{-p})$ and in particular one has $p \leq 163$.}

The proof also yields more precise results on the torsion subgroups of elliptic curves. 

\theoin[mazur-y1]{(Mazur \cite{MazurY1}) Let $E/\Q$ be an elliptic curve. Then the torsion subgroup of $E(\Q)$ is isomorphic to one of the following groups: $\Z/n\Z$ for $n \leq 10$ or $n=12$, or $\Z/2\Z \times \Z/2n\Z$ for $n \leq 4$. In particular, if $p$ is a prime such that $E(\Q)$ has a $p$-torsion point, then one has $p \leq 7$.  }

Before coming back to Serre's uniformity question, let us mention that Theorem \ref{mazur-y1} has been generalized to arbitrary number fields by Merel \cite{Merel} and made explicit by Parent \cite{Parent-Torsion}. The proofs of these statements still follow Mazur's strategy. 

\theoin[torsion-nbfield]{(Merel \cite{Merel}, Oesterl\'e, Parent \cite{Parent-Torsion}) Let $d \geq 1$ be an integer. There exists an effectively computable bound $B(d)$ such that for any number field $K$ of degree at most $d$ and any elliptic curve $E/K$, the torsion subgroup of $E(K)$ has cardinality at most $B(d)$. The prime divisors of $B(d)$ are all at most $(1+3^{d/2})^2$. }

Case \ref{suq-s} in Serre's uniformity question was ruled out for large enough $p$ in \cite{BP}, and for all $p \geq 11$ except $p=13$ in \cite{BPR}. The arguments therein also follow Mazur's strategy. The recent articles \cite{LFL} (for primes $p \geq 1.4 \cdot 10^7$) and \cite{lombardo} (for all primes, using the calculations of \cite{BBM} for $p < 100$) use Mazur's strategy to show that in case \ref{suq-ns}, the image of $\rho_{E,p}$ is \emph{equal} to said normalizer. However, Mazur's approach is more difficult to adapt for this setting for reasons that we will get into. 

After describing these applications, we now review Mazur's strategy. 

Let $X$ be a compactified modular curve over $\Q$. Then $X$ is endowed with a natural morphism $j: X \rar \mathbb{P}^1_{\Q}$ mapping any noncuspidal point of $X$ to the $j$-invariant of the elliptic curve that it represents. The strategy consists in studying separately the points of $X(\Q)$ with non-integral $j$-invariant and the points of $X(\Q)$ with integral $j$-invariant. In terms of moduli spaces, this amounts to \emph{first} proving that elliptic curves corresponding to a rational point of $X$ have no (or few) primes of potentially multiplicative reduction, \emph{then} treating this case by other methods.   

The strategy to rule out the existence of rational points with non-integral $j$-invariant (or otherwise control them) is to use a \emph{formal immersion argument} following the spirit of \cite[Corollary 4.3]{FreyMazur}. More precisely, let $\ell$ be a prime number and $P \in X(\Q)$ be a point such that $j(P) \notin \Z_{(\ell)}$. Assume for the sake of simplicity that $\ell > 2$ and that $X$ has a smooth proper model $\mathcal{X}$ over $\Z_{(\ell)}$. Then $P$ reduces to a \emph{cuspidal point} $c \in \mathcal{X}(\F_{\ell})$. 

Let $A/\Q$ be a quotient of the Jacobian of $X$, so that $A$ is an Abelian variety with good reduction at $\ell$ with N\'eron model $\mathcal{A}$ over $\Z_{(\ell)}$, and let $f: X \rar A$ be a morphism. Then $f$ extends to a morphism $\mathcal{X} \rar \mathcal{A}$ still denoted by $f$, and we can usually assume that $f(c) = 0$. Then $f(P) \in \mathcal{A}(\Q)$ is a point whose reduction modulo $\ell$ vanishes. By a result of Raynaud \cite[Proposition 1.1]{FreyMazur}, this implies that $f(P)$ is either a non-torsion point, or zero. 

In practice, one often proves that $f(P)$ is a torsion point by showing that $A(\Q)$ is a finite group. In \cite{FreyMazur}, Mazur uses the \emph{Eisenstein quotient} $\tilde{J}$ of the Jacobian $J_0(p)$ of $X_0(p)$, whose construction and properties -- in particular the finiteness of $\tilde{J}(\Q)$ by descent -- are discussed in \cite{MazurY1}. Since this paper, a different method has emerged thanks to recent progress in the Bloch-Kato conjecture, in particular the work of Kolyvagin-Logachev \cite{KolLog} and Kato \cite{KatoBSD}. They prove instances of the \emph{Bloch-Kato conjecture}: that is, if the $L$-function of certain abelian varieties $A$ does not vanish at its central point, then the group $A(\Q)$ is finite. 

That is why, in \cite{Merel,Parent-Torsion,BP,BPR}, the \emph{Eisenstein quotient} is replaced with the \emph{winding quotient}, that is, the largest quotient $J_e$ of $J_0(p)$ such that the $L$-function of $J_e$ does not vanish at the central point. 

In contrast, the argument showing that $f(P)$ is nonzero is geometric and local. The cuspidal subscheme of $\mathcal{X}$ is well-understood, so it is usually feasible to choose $f$ such that $f(c_0)=0$, for some cuspidal point $c_0 \in X$ whose reduction modulo $\ell$ is also $c$. The \emph{formal immersion argument} consists in proving that the ring homomorphism $\widehat{f^{\sharp}}: \widehat{\OO_{A,0_{\F_{\ell}}}} \rar \widehat{\OO_{X,c}}$ is surjective. Indeed, this claim implies that $f^{\sharp}$ is an epimorphism in the category of Noetherian local rings, which implies that the morphism $f: \mathcal{X}(\Z_{\ell}) \rar \mathcal{A}(\Z_{\ell})$ is injective in the residue disk of $c$ (which contains $c_0$ and $P$): since $f(c_0)=0_{A_{\Q_{\ell}}}$, one has $f(P) \neq 0_{A_{\Q}}$. 

In applications, the formal immersion argument is convenient because the geometry of modular curves at their cusps is well-understood, including in finite characteristic, thanks to the $q$-expansion principle, which lets one describe the differential forms near a cusp. 

The second part of Mazur's strategy is to handle points $P \in X(\Q)$ corresponding to elliptic curves with an integral $j$-invariant. The method usually differs according to the problem. 

In the case of torsion on elliptic curves, it is enough to apply a specialization lemma: indeed, if an elliptic curve $E/K$ with potentially good reduction above a prime $\mathfrak{r}$ of residue characteristic $3$ has a $p$-torsion point (with $p > 7$), then the reduction of this point modulo $\mathfrak{r}$ (in the N\'eron model $\mathcal{E}$ of $E$) still has order $p$. Thus $p$ divides the cardinality of $\mathcal{E}(\OO_K/\mathfrak{r})$, which has an explicit bound in terms of the degree of the extension $K/\Q$ using Hasse's theorem (see for instance \cite[Proposition 1.3]{Parent-Torsion}). 

For \ref{suq-0}, Mazur studies the \emph{isogeny character} of the elliptic curve and proves that it is of the form $\alpha \omega_p^k$, where $\alpha$ is a character of $G_{\Q}$ of order dividing $12$, $\omega_p$ is the mod $p$ cyclotomic character of $G_{\Q}$ and $k \in \Z/(p-1)\Z$ can take at most $6$ possible values. By comparing the implied congruences with Hasse's theorem, he is then able to prove that either $p \in \{17,37\}$ or $\Q(\sqrt{-p})$ has class number one. 

For \ref{suq-s}, the argument in \cite{BP} is \emph{Runge's method}, which provides an upper bound on the $j$-invariant of integral points lying on the modular curve, which can be shown to contradict a lower bound coming from refinements (in this case, due to Pellarin \cite{Pellarin}) of Masser-W\"ustholz's \cite{MW} isogeny bound. This is also the idea used in \cite{BPR,LFL,lombardo}.   

Let us now explain why \ref{suq-ns} is more difficult to approach with Mazur's strategy. It is known by \cite[Theorem 1]{Chen-Merel-conj} that the $L$-function of every simple isogeny factor of the Jacobian $J_{ns}^+(p)$ of the associated modular curve $X_{ns}^+(p)$ vanishes with odd order. By the Birch and Swinnerton-Dyer conjecture, we thus expect that every isogeny factor of $J_{ns}^+(p)$ has infinitely many rational points. Moreover, the study of integral points of $X_{ns}^+(p)$ with Runge's method does not work, because this method relies on \emph{Runge's condition} \cite[\S 1]{BP} -- namely, that there are at least two Galois orbits of cusps -- which is not satisfied on $X_{ns}^+(p)$. 

A different approach has emerged in the past decade to compute modular points on rational curves, the so-called \emph{quadratic Chabauty} method, following ideas of Kim \cite{Kim05, KimTama08, Kim09}. 

Unlike Mazur's approach, it seems that no examples have been completely worked out without using an explicit model of the curve to effectively determine the rational points; however, when the curve of interest has a model, the quadratic Chabauty method applies in much greater generality than Mazur's method, since the Jacobian $J$ of the modular curve $X$ does not need to have a quotient with finitely many rational points any more. A somewhat typical condition (see \cite[Section 2]{Xns+17}) would be that $J$ satisfies
\[ \dim{J} = \mrm{rk}\,J(\Q),\, \mrm{rk}\, NS(J) > 1.  \]

 Using this method, the authors of \cite{Xns+13} were able to rule out $p=13$ in cases \ref{suq-s} (and \ref{suq-ns}, because it turns out that the associated modular curves are isomorphic when $p=13$), this completing the solution of \ref{suq-s}:

\theoin[uniformity-split]{(Bilu-Parent \cite{BP}, Bilu-Parent-Rebolledo \cite{BPR}, Balakrishnan-Dogra-M\"uller-Tuitman-Vonk \cite{Xns+13}) Let $p \geq 11$ be a prime number. Then any elliptic curve $E/\Q$ such that the image of $\rho_{E,p}$ is contained in the normalizer of a split Cartan subgroup (resp. of a nonsplit Cartan subgroup if $p=13$) has complex multiplication. } 

More recently, the authors of \cite{Xns+13} determine in \cite{Xns+17}, among other examples, the rational points of the curve $X_{ns}^+(17)$, but report that computational difficulties arose in the case of $X_{ns}^+(19)$. On the other hand, Dogra and Le Fourn \cite{LFD} proved that the quadratic Chabauty method could in principle be applied to any $X_{ns}^+(p)$ of genus at least two.

\section{Basic results of this thesis (Chapters 1-2)}

This thesis is a study of the arithmetic of the modular curves $X_E^{\alpha}(p)$ for primes $p \geq 7$ and elliptic curves $E/\Q$. 

\subsection{Chapter 1: Moduli problems}

In Chapter \ref{moduli-spaces}, we review the theory of moduli problems of elliptic curves as described in the work of Katz and Mazur \cite{KM}. The main goal is to give a complete proof of the existence of the modular curves $X_E^{\alpha}(p)$ over more general bases than fields (generalizing \cite[Proposition 2]{KO} and \cite[\S 1.7]{Mazur-Open}), as well as their basic geometric properties. In particular, we show that, like classical modular curves, these curves can be endowed with supplementary $\Gamma_0(m)$-structures yielding Hecke operators. Such a result seems to be folklore, but we are not aware of a fleshed-out proof. 

When applying Mazur's method to a modular curve $X$, as we saw, one constructs a morphism $f$ from the modular curve $X$ to a suitable quotients of its Jacobian. One then has to analyze the geometry of $f$ at the cusps, including in finite characteristic. That is why this chapter deals with moduli spaces over bases that are not fields. 

\mytheo[XG-exists]{Let $N \geq 3$ be an integer and $R$ be a regular excellent Noetherian $\Z[1/N]$-algebra. Let $G$ be a finite \'etale commutative group scheme over $R$ which is \'etale-locally isomorphic to $(\Z/N\Z)^{\oplus 2}$, endowed with a fibrewise perfect bilinear alternating pairing $(-,-)_G: G \times G \rar \mu_N$. Then there is a proper smooth $R$-scheme $X_G(N)$ of relative dimension one satisfying the following properties:
\begin{itemize}[noitemsep,label=$-$]
\item There is a finite locally free map $(j,\det): X_G(N) \rar \mathbb{P}^1_R \times (\Z/N\Z)^{\times}_R$, and the geometric fibres of $\det$ are connected. 
\item The reduced subscheme attached to the closed subspace $\{j=\infty\}$, called the \emph{cuspidal subscheme}, is finite \'etale over $R$ and defines an effective Cartier divisor. 
\item For any $R$-scheme $S$, $j^{-1}(\mathbb{A}^1_R)(S)$ is functorially isomorphic to the set of equivalence classes of $(E/S, \iota)$, where $E/S$ is a relative elliptic curve, $\iota: G_S \rar E[N]$ is an isomorphism of finite \'etale group schemes over $S$. Moreover, under this isomorphism, $j$ maps the couple $(E/S,\iota)$ to the $j$-invariant of $E/S$, and $\det$ maps the couple $(E/S,\iota)$ to the unique $\alpha \in (\Z/N\Z)^{\times}_R(S)$ such that the fibrewise perfect bilinear alternating pairing
\[G \times G \overset{\iota \times \iota}{\longrightarrow} E[N] \times E[N] \overset{\mrm{We}}{\rar} (\mu_N)_S\] is exactly $\alpha \cdot (-,-)_G$, where $\mrm{We}$ denotes the Weil pairing.  
\item The formation of $X_G(N), j,\det$, the cuspidal subscheme, and the above isomorphism of functors commutes with base change.  
\end{itemize}

For $\alpha \in (\Z/N\Z)^{\times}$, let $X_G^{\alpha}(N)$ denote the closed open subscheme $\det^{-1}(\alpha)$ of $X_G(N)$, and $J_G^{\alpha}(N)$ denote its relative Jacobian over $R$. For $n \in (\Z/N\Z)^{\times}$, the isomorphism of functors $(E/S,\iota,C) \longmapsto (E/S,n \cdot \iota,C)$ extends to an isomorphism $[n]: X_G^{\alpha}(N) \rar X_G^{n^2\alpha}(N)$.

For every prime number $\ell \nmid N$, for every $\alpha \in (\Z/N\Z)^{\times}$, there is a Hecke correspondence on $X_G^{\alpha}(N) \times X_G^{\ell\alpha}(N)$ inducing an endomorphism $T_{\ell}$ of $\prod_{\alpha \in (\Z/N\Z)^{\times}}{J_G^{\alpha}(N)}$. The various $T_{\ell}$ and $[n]$ as previously commute pairwise, and $T_{\ell}$ sends the component $J_G^{\alpha}(N)$ into the component $J_G^{\ell\alpha}(N)$. 
}

By choosing $N=p$ prime, and for $G$ the $p$-torsion of an elliptic curve endowed with its Weil pairing, we find:

\mycor[XEalpha-exists-hecke]{Let $R$ be an excellent Dedekind ring\footnote{This is for instance the case if $\mrm{Frac}(R)$ is of characteristic zero, if $R$ is a complete discrete valuation ring, or $R$ is essentially of finite type over a field.} and $E/R$ be an elliptic curve. Let $p \geq 3$ be a prime invertible in $R$ and $\alpha \in \F_p^{\times}$. There exists a smooth proper relative curve $X_E^{\alpha}(p) \rar \Sp{R}$ with geometrically connected fibres endowed with the following datum:
\begin{itemize}[noitemsep,label=$-$]
\item a finite locally free map $j: X_E^{\alpha}(p) \rar \mathbb{P}^1_{R}$ such that the reduced closed subscheme attached to the subspace $\{j=\infty\}$ is finite \'etale over $R$,
\item for every $R$-scheme $S$, a functorial isomorphism between $j^{-1}(\mathbb{A}^1_R)(S)$ and the equivalence classes of pairs $(F/S,\iota)$, where $F/S$ is an elliptic curve and $\iota: E[p]_S \rar F[p]$ is an isomorphism of group schemes such that the following two pairings are equal: 
\[E[p]_S \times E[p]_S \overset{\alpha \cdot \mrm{We}}{\lrar} (\mu_p)_S,\quad E[p]_S \times E[p]_S \overset{\iota \times \iota}{\longrightarrow} F[p]\times F[p] \overset{\mrm{We}}{\rar} (\mu_p)_S.\]
\item a relative Jacobian $J_E^{\alpha}(p)$,
\item for every prime $\ell \neq p$, Hecke operators $T_{\ell}: J_E^{\alpha}(p) \rar J_E^{\ell\alpha}(p)$ that pairwise commute,
\item for every $n \in \F_p^{\times}$, an isomorphism $[n]: X_E^{\alpha}(p) \rar X_E^{n^2\alpha}(p)$.  
\end{itemize}}

It remains to relate these curves $X_E^{\alpha}(p)$ to the classical modular curve $X(p)$ attached to the congruence subgroup $\Gamma(p)$. 

More precisely, let $X(p,p)$ be the compactification of the moduli scheme parametrizing elliptic curves over a $\Z[1/p]$-scheme along with a basis of their $p$-torsion subgroup: it is a smooth projective $\Z[1/p]$-scheme of relative dimension one, and is endowed with a left action of the group $\GL{\F_p}/\{\pm I_2\}$, a finite locally free morphism $j: X(p,p) \rar \mathbb{P}^1_{\Z[1/p]}$ (given by the $j$-invariant) and the Weil pairing morphism $X(p,p) \rar \Sp{\Z[1/p,\zeta_p]}$, which is smooth proper of relative dimension one with geometrically connected fibres. 

When $R$ is a field in which $p$ is invertible and $E/R$ is an elliptic curve, we realize the disjoint reunion of the curves $X_E^{\alpha}(p)$ as a Galois twist of (geometrically disconnected) modular curve $X(p,p)_R$ over $R$, and check that this identification is compatible with the structures identified above: the $j$-invariant and the Hecke operators. While this correspondence is implicitly used by e.g. Wiles \cite[p. 543]{wiles} for the purpose of the $3-5$-switch, and alluded to in the introduction of \cite{Virdol}, it has not, to the author's knowledge, been spelled out. 

There are several possible definitions of \emph{Galois twist}, and we mean here the following one (which is a special case of Proposition \ref{cocycle-twist}): 

\defintro{Let $X$ be a smooth quasi-projective scheme over a field $k$ with separable closure $k_s$. Let $\rho: \mrm{Gal}(k_s/k) \rar \mrm{Aut}(X)$ be a continuous group homomorphism, where $\mrm{Aut}(X)$ is endowed with the discrete topology. The \emph{twist} of $X$ by $\rho$ is the unique (up to unique isomorphism) smooth quasi-projective $k$-scheme $X'$ endowed with an isomorphism $j: X'_{k_s} \rar X_{k_s}$ of $k_s$-schemes such that, if $\beta: X'(k_s) \overset{\sim}{\rar} X'_{k_s}(k_s) \overset{j}{\rar} X_{k_s}(k_s) \overset{\sim}{\rar} X(k_s)$ and $g \in \mrm{Gal}(k_s/k)$, for any $P \in X'(k_s)$, one has $\beta(g \cdot P) = \rho(g)(g \cdot \beta(P))$. }

This construction is functorial in a natural sense (see Proposition \ref{twist-equiv}). 

\mytheo[XE-vs-Xrho]{Let $k$ be a field with separable closure $k_s$ and $p \geq 3$ be a prime number invertible in $k$. Let $E/k$ be an elliptic curve and $(P,Q) \in E(k_s)^2$ be a basis of its $p$-torsion. For each $\sigma \in \mrm{Gal}(k_s/k)$, let $\rho(\sigma) \in \GL{\F_p}$ be the matrix such that $\begin{pmatrix}\sigma(P)\\\sigma(Q)\end{pmatrix} = \rho(\sigma)\begin{pmatrix}P\\Q\end{pmatrix}$. Then $\rho: \mrm{Gal}(k_s/k) \rar \GL{\F_p}$ is a group anti-homomorphism. 

Let $X(p,p)_{\rho^{-1}}$ be the twist of $X(p,p)_k$ by the Galois representation $\rho^{-1}$. The twist by $\rho^{-1}$ of the morphism \[(j,\mrm{We}): X(p,p)_k  \rar \mathbb{P}^1_k \times_k \Sp{k \otimes \Z[\zeta_p]}\] is a morphism \[(j,\mrm{We}_{\rho^{-1}}): X(p,p)_{\rho^{-1}} \rar \mathbb{P}^1_k \times_k S,\] where $S$ denotes the constant $k$-scheme whose underlying set is the collection $\mu_p^{\times}(k_s)$ of primitive $p$-th roots of unity in $k_s$. For $\zeta \in \mu_p^{\times}(k_s)$, let $X_{\rho,\zeta}(p)$ denote the inverse image of the component $\zeta$ in $X(p,p)_{\rho^{-1}}$. 

Let $\xi = \langle P,\,Q\rangle_{E[p]} \in \mu_p^{\times}(k_s)$. For every $a \in \F_p^{\times}$, there is an isomorphism of $\mathbb{P}^1_k$-schemes $\iota_{E,P,Q}^a: X_E^a(p)_k \rar X_{\rho,\xi^a}(p)$ such that for any non-cuspidal point $(F/k_s,\iota) \in X_E^a(p)(k_s)$, one has \[\iota_{E,P,Q}^a(F/k_s,\iota) = (F/k_s,(\iota(P),\iota(Q))) \in X(p,p)(k_s)\supset X_{\rho,\xi^a}(p)(k_s).\] 

Moreover, the disjoint reunion of the $\iota_{E,P,Q}^a$ induces an isomorphism $\iota_{E,P,Q}^J$ between the direct product of the $J_E^a(p)$ and the Galois twist $J(p,p)_{\rho^{-1}}$ by $\rho^{-1}$ of the Jacobian\footnote{The notion of Jacobian of a non-geometrically connected curve does not seem entirely standard, but amounts to a minor modification of existing results on smooth proper schemes $C \rar S$ of relative dimension one with geometrically connected fibres. We discuss the necessary adaptations in Appendix \ref{jacobian-relative}.} of the curve $X(p,p)$ over $k$. The Abelian variety $J(p,p)_{\rho^{-1}}$ is endowed (through the twist) of Hecke operator $T_{\ell}$ for any prime $\ell \neq p$ and $[n]$ for any $n \in \F_p^{\times}$ (where $[n]$ maps a couple $(F/S,(A,B))$ -- where $F/S$ is an elliptic curve over the $k$-scheme $S$ and $(A,B)$ is a basis of $F[p](S)$ -- to $(F/S,(nA,nB))$). 

The isomorphism $\iota_{E,P,Q}^J$ preserves the Hecke operators and the action of $\F_p^{\times}$. 
}

For later use, we also state and prove in this chapter some well-known features of the arithmetic geometry of classical modular curves, namely: the Eichler-Shimura relation, a lemma about the action of the Hecke operator on the cuspidal subscheme, and a $q$-expansion principle. \\

\subsection{Chapter 2: Hecke and $\GL{\F_p}$-structure of $J(p,p)$}

In Chapter \ref{analytic}, we study the geometric properties of the modular curve $X(p,p)_{\Z[1/p]}$. Let $J(p,p)$ denote its relative Jacobian. Let $X_1(p,\mu_p)$ denote the base change to $\Z[1/p,\zeta_p]$ of the smooth proper modular curve over $\Z[1/p]$ whose noncuspidal points parametrize couples $(E/S, P)$ up to isomorphism, where $E/S$ is an elliptic curve, $P$ is a point of ``exact order $p$'' in the sense of \cite[(3.2)]{KM}, and let $J_1(p,\mu_p)$ denote the Jacobian of the relative curve $X_1(p,\mu_p) \rar \Sp{\Z[1/p]}$; it is endowed with modified Hecke operators $T_{\ell}'$ and modified diamond operators $\langle n\rangle$' (in order to take into account the action on $\Sp{\Z[1/p,\zeta_p]}$, see Definition \ref{modified-hecke-operator}). 

For each $t \in \F_p$ (resp. $t=\infty$), we define a morphism $u_t: J(p,p) \rar J_1(p,\mu_p)$ as follows by push-forward functoriality: we map a couple $(E,(P,Q))$ to $(E,tP+Q,\langle P,\,Q\rangle_E)$ (resp. to $(E,P,\langle P,\,Q\rangle_E)$), where $\langle P,\,Q\rangle_E$ is the Weil pairing on $E[p]$. 

Let $P(p)$ be the direct product of $p+1$ copies of $J_1(p,\mu_p)$ indexed by $\F_p \cup \{\infty\}$; we can construct an action of $\GL{\F_p}$ on $P(p)$ which commutes with every $T_{\ell}'$, and such that the modified diamond $\langle n\rangle'$ is exactly the action of $nI_2 \in \GL{\F_p}$. Let $u: J(p,p) \rar P(p)$ be the direct product of the $u_t$. Then $C(p):=\ker{u} = \cap_{t \in \F_p \cup \{\infty\}}{\ker{u_t}}$ is a proper group scheme over $\Z[1/p]$, and we let $C^0(p)$ be the (smooth proper over $\Z[1/p]$) N\'eron model of the connected component of $C(p)_{\Q}$ containing the unit point.

\mytheo[the-exact-sequence]{The morphism $u$ commutes to the action of $\GL{\F_p}$ and the Hecke operators, in the sense that $u \circ T_{\ell} = T'_{\ell}\circ u$. Moreover, there exists an explicit morphism $R: P(p) \rar P(p)$ equivariant for the $T'_{\ell}$ and the action of $\GL{\F_p}$ such that the sequence of Abelian schemes over $\Z[1/p]$ 
\[0 \rar C^0(p) \rar J(p,p) \overset{u}{\rar} P(p) \overset{R}{\rar} P(p),\]
is a complex whose fibres are exact up to isogeny. 
More precisely, for any ring homomorphism $\Z[1/p] \rar k$ with a field $k$, the groups of connected components of $(\ker{R})_k$ and $C(p)_k$ are annihilated by $2p(p+1)$. 
}

We also describe the holomorphic differentials on $X(p,p)_{\C}$ along with their action by the Hecke operators and $\GL{\F_p}$. Let us first introduce some notation for the rest of this introduction. 

\textbf{Notation for representation of $\GL{\F_p}$:}
\begin{itemize}[noitemsep,label=\tiny$\bullet$]
\item $\mathcal{D}$ denotes the set of characters $\F_p^{\times} \rar \C^{\times}$,
\item $\mrm{St}$ denotes the Steinberg representation of $\GL{\F_p}$: it is of dimension $p$ and can be realized over $\Z$.
\item For any $\alpha \in \mathcal{D}$, let $\mrm{St}_{\alpha}$ denote the twist of $\mrm{St}$ by the character $\alpha(\det)$ of $\GL{\F_p}$. 
\item For any $\alpha,\beta \in \mathcal{D}$, let $\pi(\alpha,\beta)$ denote the principal series representation of $\GL{\F_p}$ attached to the characters $\alpha,\beta$. 
\end{itemize}
For any left $\C[\GL{\F_p}]$-module $V$, let $V^{\vee}$ denote the \emph{right} $\C[\GL{\F_p}]$-module $\mrm{Hom}(V,\C)$.  \\

\textbf{Notation for spaces of modular forms:}
\begin{itemize}[noitemsep,label=\tiny$\bullet$]
\item $\mathscr{S}$ denotes the collection of $(f,\mathbf{1})$, where $\mathbf{1} \in \mathcal{D}$ is the trivial character, and $f$ is a newform of level $p$, weight $2$, and trivial character. 
\item $\mathscr{P}$ denotes the collection of $(f,\chi)$, where $f \in \mathcal{S}_2(\Gamma_1(p))$ is a normalized newform with nontrivial character $\chi$. The complex conjugation acts freely on $\mathscr{P}$ by exchanging $(f,\chi)$ with $(\overline{f},\chi^{-1})$, where $\overline{f} \in \mathcal{S}_2(\Gamma_1(p))$ the normalized newform whose $q$-expansion is $\sum_{n \geq 1}{\overline{a_n(f)}q^n}$.
\item $\mathscr{P}'$ denotes a set of representatives for $\mathscr{P}$ modulo the complex conjugation. 
\item $\mathscr{C}$ denotes the collection of $(f,\chi)$, where $f \in \mathcal{S}_2(\Gamma_1(p)\cap \Gamma_0(p^2))$ is a normalized newform with character $\chi \in \mathcal{D}$ and such that no twist of $f$ by a character in $\mathcal{D}$ is of level $p$.
\item $\mathscr{C}_M$ denotes the collection of the $(f,\mathbf{1}) \in \mathscr{C}$ such that $f$ has complex multiplication (necessarily by the field $\Q(\sqrt{-p})$, and $p \equiv 3\pmod{4}$). 
\item $\mathscr{C}_1$ denotes the collection of the $(f,\mathbf{1}) \in \mathscr{C} \backslash \mathscr{C}_M$, and $\mathscr{C}'_1$ is a set of representatives for $\mathscr{C}_1$ modulo twists by the unique quadratic character in $\mathcal{D}$. 
\item If $f \in \mathcal{S}_k(\Gamma_1(N))$ is a newform with character $\chi$, $L$ is a number field containing the coefficients of $f$ and $\lambda$ is a maximal ideal of $\OO_L$, $V_{f,\lambda}$ is the continuous irreducible representation $G_{\Q} \rar \GL{L_{\lambda}}$ such that, for any prime $q \nmid N$ coprime to the residue characteristic of $\lambda$, the characteristic polynomial of the arithmetic Frobenius at $q$ is $X^2-a_q(f)X+q^{k-1}\chi(q)$.   
\end{itemize}

\mytheo[differential-forms-xpp]{There is an identification $H^0(X(p,p)_{\C},\Omega^1) \simeq \bigoplus_{b \in \F_p^{\times}}{\mathcal{S}_2(\Gamma(p)) \otimes \Delta_{b,1}}$ compatible with the action of the $T_{\ell}$ and of $\GL{\F_p}$.  
Moreover, $H^0(X(p,p)_{\C}, \Omega^1)$ decomposes, as a right $\C[\GL{\F_p}]$-module endowed with the action of the $T_{\ell}$, as the following direct sum:
\begin{itemize}[noitemsep,label=$-$]
\item for every $(f,\mathbf{1}) \in \mathscr{S}$, for every $\alpha \in \mathcal{D}$, a copy of $\mrm{St}_{\alpha}^{\vee}$, 
\item for every $(f,\chi) \in \mathscr{P}$ modulo complex conjugation, for every $\alpha \in \mathcal{D}$, a copy of $\pi(\alpha,\alpha\chi)^{\vee}$,
\item for every $(f,\chi) \in \mathscr{C}$, a copy of $C_f^{\vee}$, for a certain irreducible cuspidal representation $C_f$ of $\GL{\F_p}$ attached to a pair $\{\phi_f,\,\phi_f^p\}$ of characters of $\F_{p^2}^{\times}$.  
\end{itemize}
For every factor as above (with $\alpha=\mathbf{1}$ in the third kind), for any prime number $\ell \neq p$ (resp. any $n \in \F_p^{\times}$), $T_{\ell}$ (resp. $nI_2$) acts on the factor by multiplication by $a_{\ell}(f)\alpha(\ell)$ (resp. $\alpha^2(n)\chi(n)$).   
}

\section{Results on the twisted modular Jacobians (Chapters 3-4)}

\subsection{Chapter 3: The Tate module of the $J_E^{\alpha}(p)$}

In Chapter \ref{tate-modules-twists}, we use the Eichler-Shimura relation and the results of Chapter \ref{analytic} to describe the Tate module of the Jacobian of $X(p,p)_{\Q}$, then of the $X_E^{\alpha}(p)$, where $E/\Q$ denotes an elliptic curve. 

In fact, the Tate module of Galois twists of $X(p,p)_{\Q}$ has already been discussed by Virdol \cite[Theorem 1.1]{Virdol} for the purposes of proving meromorphic continuation for its $L$-function in certain cases. 

\theoin[Virdolthm]{(Virdol \cite[Theorem 1.1]{Virdol}) Let $\rho: \mrm{Gal}(\Qbar/\Q) \rar \GL{\F_p}$ be a group homomorphism. One has an equality of $L$-functions 
\[L(s,X(p,p)_{\rho}) = \prod_{\pi}{L(s,\rho_{\pi,\ell}\otimes (\tilde{\varphi_{\pi}}\circ \rho))},\]
where
\begin{itemize}[noitemsep,label=$-$]
\item in the product, the $\pi = \pi_f \otimes \pi_{\infty}$ are cuspidal automorphic representations of $\GL{\mathbb{A}_{\Q}}$ of weight $2$, where $\pi_{\infty}$ is cohomological and $\pi_f$ a representation of $\GL{\mathbb{A}_{f,\Q}}$ such that $\pi_f^{\hat{\Gamma}(p)} \neq 0$, where $\hat{\Gamma}(p) \leq \GL{\mathbb{A}_{f,\Q}}$ is the principal adelic congruence subgroup of level $p$, \item $\rho_{\pi,\ell}$ is the $\ell$-adic Galois representation attached to $\pi$,
\item $\tilde{\varphi}_{\pi}$ is the representation of $\GL{\F_p}$ attached to $\pi_f^{\hat{\Gamma}(p)}$.
\end{itemize}}

This chapter differs from the work of Virdol on the following aspects. First, his statement and arguments rely on the adelic language of automorphic representations, which we try to avoid at this point in the text. Our classical approach makes the representations $\tilde{\varphi_{\pi}} \circ \rho$ explicit. Second, when the Galois twist comes from the $p$-torsion of an elliptic curve (as in Theorem \ref{XE-vs-Xrho}), Virdol computes the $L$-function (in fact, the Tate module) of the product of the $J_E^{\alpha}(p)$ over all $\alpha \in \F_p^{\times}$, while we are interested in the Tate modules of each one of the $X_E^{\alpha}(p)$, for an elliptic curve $E/\Q$.

In the rest of the introduction, $p \geq 7$ is a prime number and $E/\Q$ is an elliptic curve. Fix a basis $(P,Q)$ of the $\F_p$-vector space $E[p](\Qbar)$, and let $\xi=\langle P,\,Q\rangle_{E[p]} \in \mu_p^{\times}(\Qbar)$ (where $\mu_p^{\times}(\Qbar)$ is the set of roots of unity in $\Qbar$ with order $p$). Let $\rho: G_{\Q} \rar \GL{\F_p}$ be the homomorphism given by 
\[\forall \sigma \in G_{\Q},\, \rho(\sigma)\begin{pmatrix}\sigma(P)\\\sigma(Q)\end{pmatrix} = \begin{pmatrix}P\\Q\end{pmatrix},\]
so that $\rho'(\sigma)=(\rho(\sigma)^T)^{-1}$ is the matrix of $\sigma$ in the basis $(P,Q)$, and in particular $\det{\rho}$ is the inverse of the cyclotomic character modulo $p$. Note that in particular $\rho$ is the inverse of the homonymous anti-homomorphism in Theorem \ref{XE-vs-Xrho}. 

We choose a number field $F \subset \C$ which is Galois over $\Q$ and whose ring of integers $\OO_F$ contains the following algebraic numbers: 
\begin{itemize}[noitemsep,label=\tiny$\bullet$]
\item the $p(p-1)^2(p+1)$-th roots of unity, 
\item the Fourier coefficients of any $f$ for $(f,\chi) \in \mathscr{S}\cup\mathscr{P}\cup\mathscr{C}$, 
\item for every $(f,\mathbf{1}) \in \mathscr{C}_M$, for every prime $\ell$ split in $\Q(\sqrt{-p})$, the roots of the Hecke polynomial $X^2-a_{\ell}(f)X+\ell$ of $f$ at $\ell$. 
\end{itemize}

For every $(f,\mathbf{1}) \in \mathscr{C}_M$, we choose a compatible system of characters $(\psi_{f,\lambda})_{\lambda}$ (where $\lambda$ runs through the maximal ideals of $\OO_F$) with $\psi_{f,\lambda}: G_{\Q(\sqrt{-p})} \rar F_{\lambda}^{\times}$ attached to $f$. \\

When $p \equiv 3\pmod{4}$, the restriction to $\F_p^{\times}\SL{\F_p}$ of the cuspidal representation $V_4$ of $\GL{\F_p}$ attached to the character of order $4$ of $\F_{p^2}^{\times}$ is the sum of two irreducible subrepresentations $C^+$ and $C^-$, which are dual of each other and conjugate of each other under $\GL{\F_p}/\F_p^{\times}\SL{\F_p}$.

\mytheo[tate-module-XE]{Let $\ell$ be a prime number and $\lambda$ be a prime ideal of $\OO_F$ with residue characteristic $\ell$. Then, for any $\alpha \in \F_p^{\times}$, $\Tate{\ell}{J_E^{\alpha}(p)} \otimes_{\Z_{\ell}} F_{\lambda}$ decomposes as the direct sum of the following $F_{\lambda}[G_{\Q}]$-submodules: 
\begin{itemize}[noitemsep,label=$-$]
\item For every $(f,\mathbf{1}) \in \mathscr{S}$, a copy of $V_{f,\lambda} \otimes (\mrm{St}\circ \rho)$, 
\item For every $(f,\chi) \in \mathscr{P}'$, a copy of $V_{f,\lambda} \otimes (\pi(\mathbf{1},\chi)\circ \rho)$,
\item For every $(f,\mathbf{1}) \in \mathscr{C}'_1$, a copy of $V_{f,\lambda} \otimes (C_f \circ \rho)$, where $C_f$ the cuspidal representation of $\GL{\F_p}$ attached to $(f,\mathbf{1})$ by Theorem \ref{differential-forms-xpp},
\item For every $(f,\mathbf{1}) \in \mathscr{C}_M$, a copy of $V_{f,\lambda,E}^{\epsilon} := \mrm{Ind}_{\Q(\sqrt{-p})}^{\Q}{\psi_{f,\lambda} \otimes C^{\epsilon}(\rho_{|G_{\Q(\sqrt{-p})}})}$, where $\epsilon \in \{\pm\}$ only depends on $\xi^{\alpha\F_p^{\times 2}}$. 
\end{itemize}
The summands of this decomposition are stable under the Hecke operators $T_{\ell}$ with $\ell \equiv 1\pmod{p}$, and $T_{\ell}$ acts on the factor corresponding to $(f,\chi)$ by multiplication by $a_{\ell}(f)$. 
}

Note that the corresponding Galois representations do not depend on the choice of representative in $\mathscr{P}'$ or $\mathscr{C}'_1$ of a given element in $\mathscr{P}$ (or $\mathscr{C}'_1$). \\

\textbf{Notation: }If $(f,\chi) \in \mathscr{S} \cup \mathscr{P}' \cup \mathscr{C}'_1$, let $R_{f,E}$ denote the following Artin representation:
\begin{itemize}[noitemsep,label=\tiny$\bullet$]
\item $\mrm{St}\circ \rho$ if $(f,\chi) \in \mathscr{S}$,
\item $\pi(\mathbf{1},\chi)\circ \rho$ if $(f,\chi) \in \mathscr{P}'$,
\item $C_f \circ \rho$ if $(f,\chi) \in \mathscr{C}'_1$. 
\end{itemize}

\subsection{Chapter 4: Galois-theoretic root numbers for $(V_{f,\lambda} \otimes R_{f,E})_{\lambda}$ and $(V_{f,\lambda,E}^{\epsilon})_{\lambda}$}

In Chapter \ref{automorphy}, we investigate the $L$-function of the factors of the Tate module of $J_E^{\alpha}(p)$. After a reminder about the relevant aspects of the Langlands program, we compute the Galois-theoretic root numbers of the compatible systems $(V_{f,\lambda} \otimes R_{f,E})_{\lambda}$ (see above for the definition of $R_{f,E}$) and $(V_{f,\lambda,E}^{\epsilon})_{\lambda}$ (defined in Theorem \ref{tate-module-XE}). For the sake of simplicity, we report here the results when $E/\Q$ has good reduction at $p$, although the computations in Chapter \ref{automorphy} cover more general cases. The result depends on the following trichotomy about Galois representations of $G_{\Q_p}$ in $\GL{\F_p}$. 

\defintro{Let $\rho_p: G_{\Q_p} \rar \GL{\F_p}$ be a representation whose determinant is the inverse of the cyclotomic character modulo $p$. We say that $\rho_p$ is
\begin{itemize}[noitemsep,label=$-$]
\item \emph{wild} if it is wildly ramified, 
\item \emph{split} if it is the sum of two characters $G_{\Q_p} \rar \F_p^{\times}$,
\item \emph{Cartan} if the image of the inertia subgroup is contained in a nonsplit Cartan subgroup.
\end{itemize}
The representation $\rho_p$ is in exactly one of these cases. }

\mytheo[root-number-XE-short]{Assume that $E/\Q$ has good reduction at $p$. 
Define $U$ as the product, over every prime number $\ell$ inert in $\Q(\sqrt{-p})$, of the $(-1)^{m_{\ell}}$, where $m_{\ell} \in \Z$ is half the conductor exponent of the Artin representation $V_4 \circ \rho$ at $\ell$. Then there is a sign $\varepsilon_p$ depending only on the finite flat $\Z_p$-scheme $E[p]_{\Z_p}$, $\xi^{\F_p^{\times 2}}$, and the choice of $\psi_{f,\lambda}$ such that the products of the local Deligne-Langlands constants attached to the summands of $\Tate{\ell}{J_E^{\alpha}(p)} \otimes F_{\lambda}$ are given by Table \ref{signs-surjective-general-introduction}. }

\begin{table}[htb]
\centering
\begin{tabular}{|c|c|c|c|}
\hline
& \multicolumn{3}{c|}{Type of $\rho_{|G_{\Q_p}}$}\\
\hline
& Split & Wild & Cartan\\
\hline
$V_{f,\lambda} \otimes R_{f,E}, (f,\mathbf{1}) \in \mathscr{S}$ & $(-1)^{(p+1)/2}$ & $(-1)^{(p-1)/2}a_p(f)$ & $(-1)^{(p+1)/2}$\\
\hline
$V_{f,\lambda} \otimes R_{f,E}, (f,\chi) \in \mathscr{P}'$ & \multicolumn{3}{c|}{$(-1)^{(p-1)/2}$}\\
\hline
$V_{f,\lambda} \otimes R_{f,E}, (f,\mathbf{1}) \in \mathscr{C}'_1$ & $(-1)^{(p+1)/2}$ & $(-1)^{(p-1)/2}$ & $(-1)^{(p+1)/2}$\\
\hline 
$V_{f,\lambda,E}^{\epsilon}, (f,\mathbf{1}) \in \mathscr{C}_M$ & $(-1)^{(p+1)/4}U$ & $(-1)^{\frac{p+1}{4}}U\epsilon\varepsilon_p$ & $(-1)^{\frac{p+1}{4}}U$\\
\hline
\end{tabular}
\caption{Global root numbers for summands of the Tate module of $J_{E}^{\alpha}(p)$}
\label{signs-surjective-general-introduction}
\end{table}

\subsection{Chapter 4: Automorphic properties of the $(V_{f,\lambda} \otimes R_{f,E})_{\lambda}$ and $(V_{f,\lambda,E}^{\epsilon})_{\lambda}$}

Using the Rankin-Selberg method, we prove:

\mytheo[functional-equation-JE-generic]{If $(f,\mathbf{1}) \in \mathscr{S}$ (resp. $(f,\chi) \in \mathscr{P}'$, resp. $(f,\mathbf{1}) \in \mathscr{C}'_1$) and the Artin representation $R_{f,E}$ is automorphic, then the $L$-function of the self-dual (up to Tate twist) Galois representation $V_{f,\lambda} \otimes R_{f,E}$ has a holomorphic continuation and satisfies a functional equation whose sign is given (when $E$ has good reduction at $p$) in Table \ref{signs-surjective-general-introduction}. 

For every $(f,\mathbf{1}) \in \mathscr{C}_M$, the $L$-function of the self-dual (up to Tate twist) Galois representation $V_{f,\lambda,E}^{\epsilon}$ admits a meromorphic continuation and satisfies a functional equation whose sign is given in Table \ref{signs-surjective-general-introduction} when $E/\Q$ has good reduction at $p$.} 

As in Theorem \ref{root-number-XE-short}, the assumption that $E/\Q$ has good reduction at $p$ is made in this part of the introduction for convenience: it is not needed in Chapter \ref{automorphy}. 

When $\rho$ is surjective, the representations $R_{f,E}$ are irreducible Artin representations of large dimension with non-solvable image, and this does not seem to change after restricting the representations to the absolute Galois group of a totally real number field. As such, proving their automorphy does not seem tractable by current methods. 

Another approach to understand the $V_{f,\lambda} \otimes R_{f,E}$ might be to reduce them modulo a prime number. We prove, using results in symmetric power functoriality due to Newton and Thorne \cite{NewThor}\footnote{A weaker result that is sufficient for our purposes can be found in the published articles \cite{NewThorQ1,NewThorQ2}. } and in tensor product functoriality in a preprint by Arias-de-Reyna and Dieulefait \cite{AdRD}, that the $V_{f,\lambda}\otimes R_{f,E}$ are automorphic modulo $p$. 

Let us first reproduce, for reference, Arias-de-Reyna and Dieulefait's result.  

\theoin[adrd]{(Arias-de-Reyna -- Dieulefait \cite[Theorem 1.1]{AdRD}) Let $f \in \mathcal{S}_k(\Gamma(1))$ be a cuspidal modular form and let $\rho_{\bullet}(f)$ be the attached compatible
system of Galois representations. Let $n > 0$ be an integer and $\pi$ be a regular algebraic cuspidal polarized automorphic representation of $GL_n(\mathbb{A}_{\Q})$ of level coprime to $3$. Let $r_{\bullet}(\pi)$ be the attached compatible system of Galois representations. Assume that $\rho_{\bullet}(f) \otimes r_{\bullet}(\pi)$ is a regular and irreducible compatible system. Then $\rho_{\bullet}(f) \otimes r_{\bullet}(\pi)$ is automorphic, i.e. there exists a regular algebraic cuspidal polarized automorphic representation $f \otimes \pi$ of $GL_{2n}(\mathbb{A}_{\Q})$ such that the compatible system attached to $f \otimes \pi$ is isomorphic to $\rho_{\bullet}(f) \otimes r_{\bullet}(\pi)$.}

\mytheo[residual-automorphy-XE]{Assume that $E$ has good reduction at $3$ and that $\rho$ is onto. Assume furthermore Theorem \ref{adrd} holds. There exists a number field $F' \supset F$ and a prime ideal $\mathfrak{p}$ of $F'$ with residue characteristic $p$ such that for any $(f,\chi) \in \mathscr{S} \cup \mathscr{P}'_1 \cup \mathscr{C}'_1$, the semi-simplification of the reduction modulo $\mathfrak{p}$ of $V_{f,\mathfrak{p}} \otimes R_{f,E}$ is the direct sum of the semi-simplifications of the reductions modulo $\mathfrak{p}$ of automorphic Galois representations (with coefficients in $(F')_{\mathfrak{p}}$). }

While this result is certainly interesting, its practical utility for Diophantine questions seems unclear. 

\section{When the image of $\rho$ is contained in the normalizer of a nonsplit Cartan subgroup (Chapter 5)}

In Chapter \ref{small-image}, we assume that $\rho(G_{\Q})$ is contained in the normalizer $N$ of a nonsplit Cartan subgroup $C$. Let $K/\Q$ be the imaginary quadratic field corresponding to the character $G_{\Q} \rar N \rar N/C \simeq \{\pm 1\}$. We assume in this introduction that $p \geq 11$ for the sake of simplicity: then $K$ is inert at $p$ by Proposition \ref{almost-never-ramified}.  

\subsection{Decomposition of the Tate module and functional equation}

In such a situation, the decomposition we gave in Theorem \ref{tate-module-XE} is not optimal and can be split further. Let us fix, once and for all, a group isomorphism $\iota$ between $C$ and $\F_{p^2}^{\times}$ preserving the norm and the trace. This isomorphism is not quite canonical, but the pair $\{\iota,\iota^p\}$ is: as a consequence, there is a canonical identification between pairs $\{\psi,\psi^p\}$ of characters of $C$ and pairs $\{\phi,\phi^p\}$ of characters of $\F_{p^2}^{\times}$. 

\mytheo[tate-module-XE-cartan]{Let $\lambda$ be a prime ideal of $\OO_F$ with residue characteristic $\ell$. For any $\alpha \in \F_p^{\times}$, $\Tate{\ell}{J_E^{\alpha}(p)}\otimes_{\Z_{\ell}} F_{\lambda}$ is the direct sum of the following continuous representations of $G_{\Q}$ stable under the Hecke operators:
\begin{itemize}[noitemsep,label=\tiny$\bullet$]
\item for every $(f,\chi) \in \mathscr{S}\cup\mathscr{P}'\cup\mathscr{C}'_1\cup\mathscr{C}_M$, a copy of $V_{f,\lambda} \otimes \theta$, where $\theta$ is a certain explicit character of $\mrm{Gal}(K(\mu_p)/\Q)$ independent from $\rho$; there are two possible $\theta$ if $(f,\chi) \in \mathscr{P}'\cup\mathscr{C}'_1$ and one otherwise,
\item for every $(f,\chi) \in \mathscr{S}\cup\mathscr{P}'$, for every pair $\{\psi,\psi^p\}$ of distinct characters $C \rar F^{\times}$ such that $\psi(tI_2)=\chi(t)$ for $t \in \F_p^{\times}$, a copy of $V_{f,\lambda} \otimes \left[\left(\mrm{Ind}_C^{N}{\psi}\right)\circ \rho\right]$, 
\item for every $(f,\mathbf{1}) \in \mathscr{C}'_1$ (resp. $(f,\mathbf{1}) \in \mathscr{C}_M$), for every pair $\{\psi,\psi^p\} \neq \{\phi_f,\phi_f^p\}$ of distinct characters $C/\F_p^{\times}I_2 \rar F^{\times}$ (resp. up to the equivalence relation $\{\psi,\psi^p\} \sim \{\theta,\theta^p\}$ if and only if $\psi\theta^{-1}$ or $\psi^p\theta^{-1}$ has order at most two), a copy of $V_{f,\lambda} \otimes \left[\left(\mrm{Ind}_C^{N}{\psi}\right)\circ \rho\right]$.
\end{itemize}
On each summand, the Hecke operator $T_{\ell}$ (for any prime $\ell \equiv 1\pmod{p}$) acts by multiplication by $a_{\ell}(f)$.  
}

Compared to the decomposition of Theorem \ref{tate-module-XE}, such summands are much easier to understand, and their properties as Galois representations are better-known. The rest of Chapter \ref{small-image} applies Mazur's strategy to the question of rational points on $X_E^{\alpha}(p)$. The first step is to show that the $L$-functions of the factors described in the previous theorem extend to entire functions, satisfy functional equations, and to determine the sign of said functional equations.   

\mytheo[signs-functional-equation-cartan-intro]{Every summand described in Theorem \ref{tate-module-XE-cartan} is a self-dual Galois representation (up to Tate twist), and its $L$-function extends to an entire function and satisfies a functional equation whose sign is $-1$ except in the following cases: 
\begin{itemize}[noitemsep,label=$-$]
\item $V_{f,\lambda} \otimes \left[\left(\mrm{Ind}_C^{N}{\psi}\right)\circ \rho\right]$ with $(f,\mathbf{1}) \in \mathscr{S}$ and $\left(\mrm{Ind}_C^{N}{\psi}\right)\circ \rho$ unramified at $p$,
\item $V_{f,\lambda} \otimes \left[\left(\mrm{Ind}_C^{N}{\psi}\right)\circ \rho^{-1}\right]$ with $(f,\mathbf{1}) \in \mathscr{C}'_1 \cup \mathscr{C}_M$ if the following condition is satisfied. Let $M/\Q_p$ denote the unique quadratic unramified extension and $\mathfrak{a}: \mrm{Gal}(\overline{\Q_p}/M) \rar \widehat{M^{\times}}$ be the class field theory reciprocity homomorphism mapping an arithmetic Frobenius to a uniformizer; in particular, $\mathfrak{a}$ has a subquotient $\overline{\mathfrak{a}}$ mapping the inertia subgroup $I_p \triangleleft \mrm{Gal}(\overline{\Q_p}/\Q_p)$ to $\OO_M^{\times}/(1+p\OO_M) \simeq \F_{p^2}^{\times}$. The sign of the functional equation of the representation is one if and only if the pairs $\{\phi_f \circ \overline{\mathfrak{a}}, \phi_f^p \circ \overline{\mathfrak{a}}\}$ and $\{\psi\circ\rho_{|I_p}, \psi^p\circ\rho_{|I_p}\}$ of characters $I_p \rar F^{\times}$ are equal. 
\end{itemize}
In both cases, $E$ and its quadratic twists have bad reduction at $p$. In the first case, one has $\rho(I_p)=3C$, $p \equiv 2, 5 \pmod{9}$, and $\psi$ has order exactly $3$. 
}

\subsection{The Bloch-Kato conjecture in rank zero for the factors of $J_E^{\alpha}(p)$}

The next result of this chapter is the proof of the Bloch-Kato conjecture in analytic rank zero for factors of $J_E^{\alpha}(p)$. This is achieved using the Beilinson-Flach Euler system constructed by Kings, Loeffler and Zerbes \cite{KLZ15} and by weakening the conditions (see for instance \cite[Section 4.4]{bigimage}) under which Rankin-Selberg convolutions of modular forms satisfy the required \emph{large image assumptions} -- which is among the contents of Chapter \ref{obstructions-euler}, which we discuss later.

In the setting of Theorem \ref{tate-module-XE-cartan}, let $(f,\chi),(f',\chi') \in \mathscr{S} \cup \mathscr{P}' \cup \mathscr{C}'_1\cup\mathscr{C}_M$, and $\{\psi,\psi^p\}, \{\psi',(\psi')^p\}$ be two pairs of distinct characters $C \rar F^{\times}$. We say that the Galois representations \[V_{f,\lambda} \otimes \left[\left(\mrm{Ind}_C^{N}{\psi}\right)\circ \rho\right],\, V_{f',\lambda} \otimes \left[\left(\mrm{Ind}_C^{N}{\psi'}\right)\circ \rho\right] \] are \emph{in the same Galois orbit} if there is a $\sigma \in G_{\Q}$ such that for all but finitely many prime numbers $q$, one has 
\[\mrm{Tr}\left(\Fr_q\mid V_{f',\lambda} \otimes \left[\left(\mrm{Ind}_C^{N}{\psi'}\right)\circ \rho\right]\right)  = \sigma \left[\mrm{Tr}\left(\Fr_q \mid V_{f,\lambda} \otimes \left[\left(\mrm{Ind}_C^{N}{\psi}\right)\circ \rho\right]\right)\right].\]

\mytheo[bloch-kato-zero]{Let $V$ be one of the four-dimensional summands described in Theorem \ref{tate-module-XE-cartan}, attached to a triple $((f,\chi),\{\psi,\psi^p\})$. For any $\alpha \in \F_p^{\times}$, there is a natural quotient $A$ stable by the Hecke operators of $J_E^{\alpha}(p)$ such that $\Tate{\ell}{A} \otimes_{\Z_{\ell}} F_{\lambda}$ identifies with the direct sum of the summands described in Theorem \ref{tate-module-XE-cartan} that are in the same Galois orbit as $V$. Suppose furthermore that $f$ does not have complex multiplication and $L(V,1) \neq 0$, so that $(f,\chi) \in \mathscr{S} \cup \mathscr{C}'_1$. 

Then $A(\Q)$ is finite.}

To construct the quotient $A$, one uses the fact that the choice of basis $(P,Q)$ for $E[p]$ defines an isomorphism $\coprod_{\alpha \in \F_p^{\times}}{X_E^{\alpha}(p)_{\Qbar}} \cong X(p,p)_{\Qbar}$, and in particular endows $\coprod_{\alpha \in \F_p^{\times}}{X_E^{\alpha}(p)_{\Qbar}}$ with an action of $\GL{\F_p}$.

\subsection{The formal immersion theorem}

In the last part of the chapter, we assume that the representation $V$ of Theorem \ref{bloch-kato-zero} is of the first form in Theorem \ref{signs-functional-equation-cartan-intro}, which implies in particular that $p \equiv 2, 5\pmod{9}$. 

Fix any cusp $c \in X_E^{\alpha}(p)(\Qbar)$ and let $q \equiv 1\pmod{p}$ be a prime number. Then the morphism $x \longmapsto (T_q-q-1)(x-c): X_E^{\alpha}(p)_{\Qbar} \rar J_E^{\alpha}(p)_{\Qbar}$ is in fact defined over $\Q$ and is independent from the choice of $c$, because $T_q-q-1$ vanishes on the divisors supported on the cuspidal subscheme, and consider its composition $f$ with the projection $J_E^{\alpha}(p) \rar A$, where $A$ is as in Theorem \ref{bloch-kato-zero}. Using the description of the differentials over $\C$ from Chapter \ref{analytic}, we are able to run an analogue of Mazur's formal immersion argument (although our actual formulation is closer to the Chabauty-Coleman method, see \cite{PoonenMcCallum}) for the morphism $f: X_E^{\alpha}(p) \rar A$. 

\defintro{For every cube root of unity $j \in \F_{p^2}^{\times}$, let \[u_j=\sum_{\substack{b \in \F_p\\1+b\sqrt{-3} \in \F_{p^2}^{\times 3}j}}{e^{\frac{2i\pi b}{p}}}.\] We say that a prime number $\ell$ is \emph{exceptional} for $p$ if $\ell \equiv \pm 1 \pmod{p}$ and there is a prime ideal $\lambda$ of $\Z[e^{2i\pi/p}]$ with residue characteristic $\ell$ containing every $u_j$. } 

Direct computations with Magma \cite{magma} show that for any $p \leq 1500$, the three $u_j$ generate the unit ideal of $\Z[e^{\frac{2i\pi}{p}}]$. We thus tentatively expect that for any $p$, there are no exceptional primes\footnote{After the first submission of this thesis, further computations on the IMJ-PRG's computer were carried out with a slightly improved algorithm, showing in approximately $8$ hours and $30$ GB of memory that there were no exceptional primes for $11 \leq p \leq 5000$ with $p \equiv 2, 5 \pmod{9}$.}.  

\mytheo[formal-immersion-diophantine]{Assume that $\rho: G_{\Q} \rar N$ is surjective with Artin conductor $N_{\rho}$. Let $\psi: C \rar F^{\times}$ be a character of order $3$ and assume that there exists $(f,\mathbf{1}) \in \mathscr{S}$ such that the central value of the $L$-function of the Galois representation $V_{f,\lambda} \otimes \left[\left(\mrm{Ind}_C^N{\psi}\right)\circ \rho\right]$ is nonzero: this implies in particular that the image by $\rho$ of an inertia subgroup at $p$ is the proper subgroup of $C$ with index $3$. 

Let $E'/\Q$ be an elliptic curve congruent to $E$ modulo $p$. Then the conductor of $E'$ is equal to $p^2N_{\rho}L$, where $L$ is a square-free product of exceptional primes that do not divide $pN_{\rho}$. }

The known results on the effectivity of Serre's open image theorem, such as \cite[Theorems 1.5, 1.10]{Zywina-Surj}, do not seem to rule out the existence of a non-CM elliptic curve $E'$ as in Theorem \ref{formal-immersion-diophantine}, nor does this result seem to easily imply that $E'$ and $E$ are isogenous, even assuming that $E$ and $E'$ both have conductor $p^2N_{\rho}$. \footnote{Since the submission of this thesis, I realized that the main result of \cite{BBM} implies that, in such a situation, when $p \leq 100$, the elliptic curve $E'$ has complex multiplication.}  

\subsection{Application to congruences to the CM elliptic curve $y^2=x^3-p$}

In \cite{Kraus-Thesis}, Kraus computes the Serre weight of the $p$-torsion of an elliptic curve with additive reduction. In particular, when $p \equiv 5\pmod{9}$ and $E$ is the elliptic curve $E_p$ with Weierstrass equation $y^2=x^3-p$, we can check that:
\begin{itemize}[noitemsep,label=\tiny$\bullet$]
\item $\rho(G_{\Q})=N$,
\item the image by $\rho$ of an inertia subgroup at $p$ is properly contained in $C$,
\item $\left(\mrm{Ind}_C^N{\psi}\right) \circ \rho$ is attached to the unique newform $g \in \mathcal{S}_1(\Gamma_1(108))$ with rational coefficients and complex multiplication by $\Q(\sqrt{-3})$.
\end{itemize}

When $p$ is small, we can check numerically with MAGMA \cite{magma} that there are no exceptional primes for $p$ and with PARI/GP \cite{pari} that for some $(f,\mathbf{1}) \in \mathscr{S}$ that the Rankin-Selberg $L$-value $L(f \times g,1)$ appearing in Theorem \ref{formal-immersion-diophantine} is nonzero\footnote{At the time of submission of this thesis, this was computed for $p=23, 29$ for the newform $g$ of weight $1$, conductor $108$, rational coefficients, and with CM by $\Q(\sqrt{-3})$ which appears in Theorem \ref{asymptotic-conductor-bound-intro}. Since then, further computations were made on the IMJ-PRG's calculators with the help of P. Molin, showing that the same non-vanishing result held for all $p \leq 150$ with $p \equiv 2, 5\pmod{9}$. We already know that there are no exceptional primes in this range, this implies in particular that the Corollary immediately below holds for $p=41, 59$, and remains true if $p=113, 131, 149$, when the conclusion is replaced with ``$E'$ has the same conductor as $E_p$''.}. In conductor at most $500000$, Cremona and Freitas \cite[Theorem 1.3]{Cremona-Freitas} prove that elliptic curves that are congruent modulo some prime $q \geq 19$ are isogenous. This implies:

\mycor{Let $p=23$. If $E'$ is an elliptic curve congruent modulo $p$ to $E_{p}$, then $E'$ is isogenous to $E_{p}$.  }

For arbitrary primes $p$, one expects that there exists $(f,\mathbf{1})$ such that $L(f \times g, 1) \neq 0$, but this seems difficult to prove. The following result of Michel, proved by analytic number theory, shows that this is indeed the case when $p$ is large enough.    

\theoin[michel]{(Michel, \cite[Appendix, Theorem 0.1]{Michel}) Let $g \in \mathcal{S}_1(\Gamma_1(N))$ be a newform with character $\chi$. There exists a constant $C_g$ satisfying the following property: for any prime $p > C_g$ if $g$ does not have real coefficients or $\chi(p) = -1$, there exists $f \in \mathcal{S}_2(\Gamma_0(p))$ such that $L(f \times g, 1) \neq 0$.}

\mytheo[asymptotic-conductor-bound-intro]{There exists a constant $p_0 \geq 11$ satisfying the following property. Let $p \geq p_0$ be a prime number such that $p \equiv 5 \pmod{9}$ and $E$ be the elliptic curve with complex multiplication by $\Q(\sqrt{-3})$ given by the minimal Weierstrass equation $y^2=x^3-p$ (with conductor $108mp^2$, where $m=1$ if $p \equiv 3\pmod{4}$ and $m=4$ if $p \equiv 1\pmod{4}$). Then, for any elliptic curve $E'$ congruent to $E$ modulo $p$, the conductor of $E'$ is given by $108mp^2L$, where $L$ is a square-free product of primes $\ell \equiv \pm 1\pmod{p}$ that are exceptional for $p$. In particular, if there are no exceptional primes for $p$, $E'$ and $E$ have the same conductor and $j(E')$ is an integer. }

\section{Large image results for the Rankin-Selberg convolution of two modular forms (Chapter 6)}

Chapter \ref{obstructions-euler} is a preprint \cite{Studnia} by the author of this thesis. It arose from the use of theory of Euler systems in the proof of Theorem \ref{bloch-kato-zero}. 

Let $\OO_K$ be the ring of integers of a finite extension $K/\Q_p$ and $T$ be a free $\OO_K$-module of finite rank on which $G_{\Q}$ has a continuous $\OO_K$-linear action. An Euler system for $T$ can be seen as a collection of cohomology classes $c_n \in H^1(\Q(\mu_n),T)$ for certain integers $n \geq 1$ (including $1$) satisfying certain norm relations. The theory of Euler systems, as described in the standard references \cite{Rubin-ES, MR-ES}, aims to bound the size of Selmer groups attached to $\mrm{Hom}_{\OO_K}(T,K/\OO_K(1))$ given an Euler system $(c_n)_n$ for $T$ such that $c_1 \neq 0$. However, the theory requires certain assumptions about the Galois image of $G_{\Q}$ in $\mrm{Aut}(T)$, such as variants or strengthenings of the following condition, which we will call the \emph{adapted image condition}: there exists $\sigma \in G_{\Q(\mu_{p^{\infty}})}$ such that $T/(\sigma-1)T \otimes_{\OO_K} K$ is a one-dimensional vector space. 

In \cite{KLZ15}, Kings, Loeffler and Zerbes construct an Euler system for some of the Galois representations appearing in Theorem \ref{tate-module-XE-cartan}, and derive the bounds on Selmer groups that we use in the proof of Theorem \ref{bloch-kato-zero} assuming the adapted image condition. In \cite[Section 4.4]{bigimage}, Loeffler proves that the adapted image condition is often satisfied:

\theoin[loefflertheo]{(Loeffler \cite{bigimage}) Let $f \in \mathcal{S}_k(\Gamma_1(N_f)), g \in \mathcal{S}_1(\Gamma_1(N_g))$ be newforms with respective characters $\varepsilon_f, \varepsilon_g$ such that $k > 1$ and $f$ is not CM. Let $L$ be a number field containing the coefficients of $f$ and $g$ with ring of integers $\OO_L$.
\begin{itemize}[noitemsep,label=$-$]
\item If there is a maximal ideal $\lambda$ of $\OO_L$ such that $V_{f,\lambda} \otimes_{L_{\lambda}} V_{g,\lambda}$ satisfies the adapted image condition, then $\varepsilon_f\varepsilon_g \neq 1$.
\item Assume that $\varepsilon_f\varepsilon_g \neq 1$ and that $N_f$ and $N_g$ are coprime. Then, for all but finitely many maximal ideals $\lambda$ of $\OO_L$, $V_{f,\lambda}\otimes V_{g,\lambda}$ satisfies the adapted image condition.   
\end{itemize}}

In the same paper, the author suggests the following question. 

\questin[loeffler-question-intro]{(Loeffler, \cite[Remark 4.4.2]{bigimage}) Consider the notations of Theorem \ref{loefflertheo} and assume that $\varepsilon_f\varepsilon_g \neq 1$. Is it true that, for all but finitely many maximal ideals $\lambda$ of $\OO_L$, $V_{f,\lambda}\otimes V_{g,\lambda}$ satisfies the adapted image condition?}

Theorem \ref{loefflertheo} does not cover all the cases needed for the proof of Theorem \ref{bloch-kato-zero}: in the second case described in Theorem \ref{signs-functional-equation-cartan-intro}, the conductors of both modular forms are divisible by $p^2$, so Theorem \ref{loefflertheo} does not prove the adapted image condition. This is the motivation for our work on Question \ref{loeffler-question-intro}. First, we give weaker sufficient conditions under which the answer to Question \ref{loeffler-question-intro} is positive. 

\mytheo[better-sufficient-conditions]{Let us keep the notations of Theorem \ref{loefflertheo} and assume that $\varepsilon_f\varepsilon_g \neq 1$. The answer to Question \ref{loeffler-question-intro} is positive if any of the following conditions is satisfied:
\begin{itemize}[noitemsep,label=\tiny$\bullet$]
\item For any (primitive) Dirichlet character $\chi$ with conductor $C$ such that $f \otimes \chi$ is Galois-conjugate to $f$, $C$ is coprime to $N_g$. 
\item Any (primitive) Dirichlet character $\chi$ such that $f \otimes \chi$ is Galois-conjugate to $f$ is even.
\item $\varepsilon_f^2=\varepsilon_g^2=1$ and $f \otimes \varepsilon_g$ is not Galois-conjugate to $f$.
\item The group $\varepsilon_g\left(\bigcap_{\chi}{\ker{\chi}}\right)$ contains an element of order $2$ but no element of order $4$, where $\chi$ runs through (primitive) Dirichlet characters such that $f \otimes \chi$ is Galois-conjugate to $f$.
\end{itemize}}

This result implies a slight generalization of the twisted Birch and Swinnerton-Dyer conjecture proved in \cite[Theorem 11.7.4]{KLZ15} for twists of elliptic curves by two-dimensional odd irreducible Artin representations: we show in Corollary \ref{BSD-improved} that in the statement of \emph{loc.cit.}, we need no longer assume that the elliptic curve and the Artin representation have coprime conductors.  

We also show that the general answer to Question \ref{loeffler-question-intro} is negative: the following equivalence lets us describe infinite families of counterexamples.  

\mytheo[answer-is-negative]{Let $f \in \mathcal{S}_k(\Gamma_1(N))$ be a newform of even weight. The following are equivalent:
\begin{itemize}[noitemsep,label=$-$]
\item For any quadratic number field $K$, the answer to Question \ref{loeffler-question-intro} is positive for all but finitely couples $(f,g)$, where $g$ is a newform of weight one with real or complex multiplication by $K$.   
\item The answer to Question \ref{loeffler-question-intro} is positive for all couples $(f,g)$ where $g$ is a newform of weight one.
\item Any (primitive) Dirichlet character $\chi$ such that $f \otimes \chi$ is Galois-conjugate to $f$ is even.
\end{itemize}
}

We also describe explicit counter-examples not arising from the same obstruction as the one highlighted by the proof of Theorem \ref{answer-is-negative}. 

\mytheo[sporadic-counter-examples]{There exist newforms $f \in \mathcal{S}_2(\Gamma_0(63)), g \in \mathcal{S}_1(\Gamma_1(1452))$ with coefficients respectively in $\Q(\sqrt{3}),\Q(\sqrt{-3})$ such that, for all but finitely many prime ideals $\mathfrak{p}$ of $L=\Q(e^{2i\pi/12})$ with residue characteristic $p \equiv 5,7 \pmod{12}$, $V_{f,\mathfrak{p}} \otimes V_{g,\mathfrak{p}}$ does not satisfy the adapted image condition. }

\renewcommand{\thesection}{\thechapter.\arabic{section}}
\renewcommand{\thesubsection}{\thesection.\arabic{subsection}}

\textbf{Acknowledgements}\; I am deeply grateful for the help and support that my supervisor Lo\"ic Merel gave me, as well as his many suggestions, comments and further questions that helped shape this work. I would like to thank the referees Pierre Parent and Henri Darmon for their reviewing this manuscript and their comments. I am also grateful to Farrell Brumley, Anna Cadoret, Pierre-Henri Chaudouard, Laurent Clozel, Pierre Colmez, Fred Diamond, Marc Hindry, Qing Liu, David Loeffler, Barry Mazur, Jean-Fran\c cois Mestre and Pascal Molin for fruitful discussions regarding parts of this thesis, for their patience in answering my questions, and for pointing me towards the relevant literature. 
\newpage

%% file: french-intro.tex
\begin{otherlanguage}{french}
\chapter*{R\'esum\'e de la th\`ese}
\phantomsection
\addcontentsline{toc}{chapter}{R\'esum\'e de la th\`ese}
\setcounter{section}{0}
\renewcommand{\thesection}{\arabic{section}}
\renewcommand{\thesubsection}{\thesection.\arabic{subsection}}

Ce chapitre est une traduction en fran\c cais de la deuxi\`eme partie de l'introduction et r\'esume le contenu de cette th\`ese. Les r\'esultats de l'introduction en anglais ne sont pas cit\'es ici; les r\'ef\'erences sont internes, renvoient au corps de la th\`ese ou \`a la bibliographie.\\

\'Etant donn\'e une courbe elliptique $E/\Q$ et un nombre premier $p \geq 7$, il existe une courbe modulaire $X_E^{\alpha}(p)$ d\'efinie sur $\Q$ dont les points non cuspidaux param\'etrent les couples $(F,\iota)$, o\`u $F$ est une courbe elliptique sur un $\Q$-sch\'ema $S$ et $\iota: E[p]_S \rar F[p]$ est un isomorphisme de sch\'emas finis plats. Cette th\`ese vise \`a \'etudier l'arithm\'etique des courbes $X_E^{\alpha}(p)$. 

\section{R\'esultats de base (Chapitres 1-2)}

\subsection{Chapitre 1: les probl\`emes de modules}

Dans le Chapitre \ref{moduli-spaces}, nous rappelons la th\'eorie des espaces de modules de courbes elliptiques telle qu'elle est d\'ecrite dans le livre de Katz et Mazur \cite{KM}. L'objectif de ce chapitre est de donner une preuve compl\`ete de l'existence des courbes modulaires $X_E^{\alpha}(p)$ sur des bases plus g\'en\'erales que des corps, g\'en\'eralisant ainsi \cite[Proposition 2]{KO} et \cite[\S 1.7]{Mazur-Open}, et de d\'ecrire bri\`evement leur g\'eom\'etrie. En particulier, nous montrons que ces courbes peuvent \^etre dot\'ees de $\Gamma_0(m)$-structures suppl\'ementaires qui les munissent de correspondances de Hecke. Un tel r\'esultat semble bien connu, mais il semble relever du folklore et nous n'en avons pas trouv\'e une preuve d\'etaill\'ee dans la litt\'erature.   

L'int\'er\^et de travailler sur des bases g\'eom\'etriques au-del\`a de celles des corps appara\^it lorsque l'on tente d'appliquer la technique de l'immersion formelle introduite par Mazur dans \cite{MazurY1}. En effet, cette derni\`ere impose de comprendre la g\'eom\'etrie de morphismes de la courbe modulaire vers certains quotients de sa jacobienne, y compris en caract\'eristique finie.  

Le premier r\'esultat est le suivant. 

\mytheofrench[fr-XG-exists]{Soient $N \geq 3$ un entier et $R$ une $\Z[1/N]$-alg\`ebre noeth\'erienne r\'eguli\`ere et excellente. Soit $G$ un sch\'ema en groupes fini sur $R$, \'etale-localement isomorphe \`a $(\Z/N\Z)^{\oplus 2}$, muni d'un accouplement bilin\'eaire altern\'e parfait $(-,-)_G: G \times G \rar \mu_N$. Il existe un sch\'ema $X_G(N)$ propre et lisse sur $R$ de dimension relative $1$ satisfaisant les propri\'et\'es suivantes:
\begin{itemize}[noitemsep,label=$-$]
\item Il poss\`ede un morphisme fini plat $(j,\det): X_G(N) \rar \mathbb{P}^1_R \times (\Z/N\Z)^{\times}_R$ et les fibres de $\det$ sont g\'eom\'etriquement connexes. 
\item Le sous-sch\'ema ferm\'e r\'eduit de $X_G(N)$ associ\'e au sous-espace ferm\'e $\{j=\infty\}$, appel\'e le \emph{sous-sch\'ema cuspidal}, est fini \'etale sur $R$ et d\'efinit un diviseur de Cartier effectif. 
\item Le foncteur des points de $j^{-1}(\mathbb{A}^1_R)$ est isomorphe au foncteur qui \`a un $R$-sch\'ema $S$ associe l'ensemble des classes d'\'equivalence de couples $(E/S,\iota)$, o\`u $E/S$ est une courbe elliptique relative, et $\iota: G_S \rar E[N]$ est un isomorphisme de $S$-sch\'emas en groupes. \`A travers cet isomorphisme, $j$ envoie le couple $(E/S,\iota)$ sur $j(E)$ et $\det$ sur l'unique $\alpha \in (\Z/N\Z)^{\times}$ tel que l'accouplement bilin\'eaire altern\'e 
\[G \times G \overset{\iota \times \iota}{\longrightarrow} E[N] \times E[N] \overset{\mrm{We}}{\rar} (\mu_N)_S\] soit exactement $\alpha \cdot (-,-)_G$, o\`u $\mrm{We}$ d\'esigne l'accouplement de Weil.  
\item La construction de $X_G(N), j,\det$, du sous-sch\'ema cuspidal, et de l'isomorphisme ci-dessus commute au changement de base.  
\end{itemize}

Pour $\alpha \in (\Z/N\Z)^{\times}$, notons $X_G^{\alpha}(N)$ le sous-sch\'ema ouvert ferm\'e $\det^{-1}(\alpha)$ de $X_G(N)$, et $J_G^{\alpha}(N)$ la jacobienne de cette courbe relative sur $R$. Pour tout $n \in (\Z/N\Z)^{\times}$, l'isomorphisme de foncteurs $(E/S,\iota,C) \longmapsto (E/S,n \cdot \iota,C)$ s'\'etend en un isomorphisme $[n]: X_G^{\alpha}(N) \rar X_G^{n^2\alpha}(N)$.

Pour chaque nombre premier $\ell \nmid N$, pour tout $\alpha \in (\Z/N\Z)^{\times}$, on dispose d'une correspondance de Hecke sur $X_G^{\alpha}(N) \times X_G^{\ell\alpha}(N)$ induisant un endomorphisme $T_{\ell}$ du sch\'ema ab\'elien $\prod_{\alpha \in (\Z/N\Z)^{\times}}{J_G^{\alpha}(N)}$. Tous les $T_{\ell}$ et les $[n]$ commutent deux \`a deux, et $T_{\ell}$ envoie la composante $J_G^{\alpha}(N)$ dans la composante $J_G^{\ell\alpha}(N)$. 
}

Lorsque $N$ est un nombre premier $p \geq 5$ et $G$ d\'esigne la $p$-torsion d'une courbe elliptique munie de son accouplement de Weil, on obtient:  

\mycorfrench[fr-XEalpha-exists-hecke]{Soient $R$ un anneau de Dedekind excellent\footnote{C'est par exemple le cas si $\mrm{Frac}(R)$ est de caract\'eristique nulle, si $R$ est un anneau de valuation discret complet, ou si $R$ est essentiellement de type fini sur un corps.} et $E/R$ une courbe elliptique. Soit $p \geq 3$ un nombre premier inversible dans $R$ et $\alpha \in \F_p^{\times}$. Il existe un $R$-sch\'ema propre et lisse de dimension relative un $X_E^{\alpha}(p)$ aux fibres g\'eom\'etriquement connexes, muni \'egalement des donn\'ees suivantes:
\begin{itemize}[noitemsep,label=$-$]
\item un morphisme fini localement libre $j: X_E^{\alpha}(p) \rar \F_p^{\times}$ tel que le sous-sch\'ema ferm\'e r\'eduit associ\'e au sous-espace $\{j=\infty\}$ soit fini \'etale sur $R$,
\item un isomorphisme entre le foncteur des points de $j^{-1}(\mathbb{A}^1_R)$ et le foncteur qui \`a chaque $R$-sch\'ema $S$ associe les classes d'\'equivalence de couples $(F/S,\iota)$, o\`u $F/S$ est une courbe elliptique relative et $\iota: E[p]_S \rar F[p]$ est un isomorphisme de sch\'emas en groupes tel que les deux accouplements ci-dessous soient \'egaux: 
\[E[p]_S \times E[p]_S \overset{\alpha \cdot \mrm{We}}{\lrar} (\mu_p)_S,\quad E[p]_S \times E[p]_S \overset{\iota \times \iota}{\longrightarrow} F[p]\times F[p] \overset{\mrm{We}}{\rar} (\mu_p)_S,\]
\item une jacobienne relative $J_E^{\alpha}(p)$,
\item pour chaque nombre premier $\ell \neq p$, des op\'erateurs de Hecke $T_{\ell}: J_E^{\alpha}(p) \rar J_E^{\ell\alpha}(p)$ qui commutent deux \`a deux,
\item pour chaque $n \in \F_p^{\times}$, un isomorphisme $[n]: X_E^{\alpha}(p) \rar X_E^{n^2\alpha}(p)$. 
\end{itemize}}

Nous relions ensuite ces $X_E^{\alpha}(p)$ \`a la courbe modulaire classique $\Gamma(p)$. Plus pr\'ecis\'ement, soit $X(p,p)$ la courbe modulaire compacte sur $\Z[1/p]$ dont les sections non-cuspidales param\'etrent les courbes elliptiques munies d'une base de leur $p$-torsion. C'est un $\Z[1/p]$-sch\'ema projectif lisse de dimension relative $1$ muni d'une action \`a gauche du groupe $\GL{\F_p}/\{\pm I_2\}$, d'un morphisme fini localement libre $j: X(p,p) \rar \mathbb{P}^1_{\Z/[1/p]}$ (donn\'e par le $j$-invariant), et d'un morphisme d'accouplement de Weil $X(p,p) \rar \Sp{\Z[1/p,\zeta_p]}$, qui est propre, lisse de dimension relative $1$ et aux fibres g\'eom\'etriquement connexes.   
 
Si $R$ est un corps de caract\'eristique distincte de $p$ et $E/R$ est une courbe elliptique, nous pouvons \'egalement construire la r\'eunion disjointe des $X_E^{\alpha}(p)$ comme tordue galoisienne de la courbe modulaire (non g\'eom\'etriquement connexe) $X(p,p)_R$ sur $R$. Nous v\'erifions \'egalement que cette identification est compatible avec les op\'erateurs de Hecke et le $j$-invariant dont nous avons parl\'e pr\'ec\'edemment. L'\'equivalence des deux constructions est implicitement utilis\'ee dans la preuve de la modularit\'e des courbes elliptiques semi-stables par Wiles \cite[p. 543]{wiles} afin d'effectuer le \og$3-5$ switch\fg{}, et semble mentionn\'ee dans l'introduction de \cite{Virdol}, mais il ne semble pas qu'elle ait \'et\'e compl\`etement r\'edig\'ee. 

Il y a plusieurs d\'efinitions possibles pour la notion de \emph{tordue galoisienne}. Celle que nous utilisons est la suivante (qui est un cas particulier de la Proposition \ref{cocycle-twist}). 

\defintrofrench{Soit $X$ un sch\'ema lisse et quasi-projectif sur un corps $k$ de cl\^oture s\'eparable $k_s$. Soit $\rho: \mrm{Gal}(k_s/k) \rar \mrm{Aut}(X)$ un morphisme de groupes continu, o\`u $\mrm{Aut}(X)$ est muni de la topologie discr\`ere. Le \emph{tordu} de $X$ par $\rho$ est l'unique $k$-sch\'ema lisse et quasi-projectif $X'$ muni d'un isomorphisme $j: j: X'_{k_s} \rar X_{k_s}$ de $k_s$-sch\'emas tel que, si $\beta$ d\'esigne la composition $X'(k_s) \overset{\sim}{\rar} X'_{k_s}(k_s) \overset{j}{\rar} X_{k_s}(k_s) \overset{\sim}{\rar} X(k_s)$ et $g \in \mrm{Gal}(k_s/k)$, pour tout $P \in X'(k_s)$, l'on a $\beta(g \cdot P) = \rho(g)(g \cdot \beta(P))$. }

Cette construction est fonctorielle dans un sens naturel explicit\'e dans la Proposition \ref{twist-equiv}. On utilisera parfois \og tordue\fg{} (au f\'eminin) lorsque l'on consid\`ere $X$ comme une courbe ou une vari\'et\'e sur le corps $k$. 

\mytheofrench[fr-XE-vs-Xrho]{Soient $k$ un corps de cl\^oture s\'eparable $k_s$ et $p \geq 3$ un nombre premier inversible dans $k$. Soient $E/k$ une courbe elliptique et $(P,Q) \in E(k_s)^2$ une base de son sous-groupe de $p$-torsion. Pour chaque $\sigma \in \mrm{Gal}(k_s/k)$, soit $\rho(\sigma) \in \GL{\F_p}$ la matrice telle que $\begin{pmatrix}\sigma(P)\\\sigma(Q)\end{pmatrix} = \rho(\sigma)\begin{pmatrix}P\\Q\end{pmatrix}$. Alors $\rho: \mrm{Gal}(k_s/k) \rar \GL{\F_p}$ est un anti-homomorphisme de groupes. 

Soit $X(p,p)_{\rho^{-1}}$ le tordu de $X(p,p)_k$ par la repr\'esentation galoisienne $\rho^{-1}$. Le tordu par $\rho^{-1}$ du morphisme \[(j,\mrm{We}): X(p,p)_k  \rar \mathbb{P}^1_k \times_k \Sp{k \otimes \Z[\zeta_p]}\] est un morphisme \[(j,\mrm{We}_{\rho^{-1}}): X(p,p)_{\rho^{-1}} \rar \mathbb{P}^1_k \times_k S,\] o\`u $S$ est le $k$-sch\'ema constant d'ensemble sous-jacent l'ensemble $\mu_p^{\times}(k_s)$ des racines primitives $p$-i\`emes de l'unit\'e dans $k_s$. Pour chaque $\zeta \in \mu_p^{\times}(k_s)$, soit $X_{\rho,\zeta}(p)$ l'image r\'eciproque de la composante $\zeta$ de $S$ dans $X(p,p)_{\rho^{-1}}$. 

Soit $\xi = \langle P,\,Q\rangle_{E[p]} \in \mu_p^{\times}(k_s)$. Pour chaque $a \in \F_p^{\times}$, on dispose d'un isomorphisme de $\mathbb{P}^1_k$-sch\'emas $\iota_{E,P,Q}^a: X_E^a(p)_k \rar X_{\rho,\xi^a}(p)$ tel que pour tout point non-cuspidal $(F/k_s,\iota) \in X_E^a(p)(k_s)$, l'on a \[\iota_{E,P,Q}^a(F/k_s,\iota) = (F/k_s,(\iota(P),\iota(Q))) \in X(p,p)(k_s)\supset X_{\rho,\xi^a}(p)(k_s).\] 

Enfin, la r\'eunion disjointe des $\iota_{E,P,Q}^a$ induit un isomorphisme $\iota_{E,P,Q}^J$ entre le produit des $J_E^a(p)$ et le tordu galoisien $J(p,p)_{\rho^{-1}}$ par $\rho^{-1}$ de la jacobienne\footnote{La notion de jacobienne pour une courbe non g\'eom\'etriquement connexe ne semble pas compl\`etement standard \`a en juger \cite{MilJac} et \cite[Chapters 8-9]{BLR}, mais peut \^etre d\'egag\'ee au prix de modifications mineures des r\'esultats existants pour les sch\'emas $C \rar S$ propres lisses de dimension relative $1$ aux fibres g\'eom\'etriquement connexes. Nous d\'ecrivons dans l'Appendice \ref{jacobian-relative} les adaptations n\'ecessaires.} de la courbe $X(p,p)_k \rar \Sp{k}$. La vari\'et\'e ab\'elienne $J(p,p)_{\rho^{-1}}$ est munie (\`a travers sa structure de tordue galoisienne) de l'op\'erateur de Hecke $T_{\ell}$ pour chaque nombre premier $\ell \neq p$ et de l'automorphisme $[n]$ pour chaque $n \in \F_p^{\times}$ (envoyant le couple $(F/S,(A,B))$ -- o\`u $F/S$ est une courbe elliptique relative sur le $k$-sch\'ema $S$ et $(A,B)$ forme une base de $F[p](S)$ -- sur $(F/S,(nA,nB))$). 

L'isomorphisme $\iota_{E,P,Q}^J$ pr\'eserve les op\'erateurs de Hecke et l'action de $\F_p^{\times}$. 
}

Puisque nous nous en servons plus tard dans le manuscrit, nous montrons quelques r\'esultats classiques de la g\'eom\'etrie arithm\'etique des courbes modulaires, et notamment la relation d'Eichler-Shimura, un principe de $q$-d\'eveloppement, et une description de l'action des op\'erateurs de Hecke sur le sous-sch\'ema cuspidal. \\

\subsection{Chapitre 2: la structure de $J(p,p)$ comme module sur l'alg\`ebre de Hecke et $\GL{\F_p}$}

Dans le chapitre \ref{analytic}, nous \'etudions les propri\'et\'es g\'eom\'etriques de la courbe modulaire $X(p,p)_{\Z[1/p]}$, dont nous appelons la jacobienne relative $J(p,p)$. Soit $X_1(p,\mu_p)$ le changement de base \`a $\Z[1/p,\zeta_p]$ de la courbe modulaire propre et lisse sur $\Z[1/p]$ dont les points non-cuspidaux param\`etrent les classes d'isomorphisme de couples $(E/S,P)$, o\`u $E$ est une courbe elliptique relative sur le $\Z[1/p]$-sch\'ema $S$ et $P$ un point \og d'ordre exact $p$\fg{} au sens de \cite[(3.2)]{KM}, et soit $J_1(p,\mu_p)$ la jacobienne de cette courbe relative sur $\Z[1/p]$. Alors nous munissons $J_1(p,\mu_p)$ d'op\'erateurs de Hecke modifi\'es $T'_{\ell}$ et d'op\'erateurs diamants modifi\'es $\langle n\rangle'$ (voir la D\'efinition \ref{modified-hecke-operator}), afin de prendre en compte l'action des op\'erateurs usuels sur la composante $\Z[1/p,\zeta_p]$. 

Pour chaque $t \in \F_p$ (resp. $t=\infty$), on d\'efinit le morphisme $u_t: J(p,p) \rar J_1(p,\mu_p)$ induit par le morphisme $X(p,p) \rar X_1(p,\mu_p)$ suivant: on envoie le couple $(E,(P,Q))$ sur $(E,tP+Q,\langle P,\,Q\rangle_E)$ (resp. $(E,P,\langle P,\,Q \rangle_E)$), o\`u $\langle P,\,Q\rangle_E$ d\'esigne l'accouplement de Weil sur $E[p]$.   

Soit $P(p)$ le produit de $p+1$ copies de $J_1(p,\mu_p)$ index\'ees par $\F_p \cup \{\infty\}$. On peut construire une action de $\GL{\F_p}$ sur $P(p)$ qui commute \`a l'action de chaque $T_{\ell}'$, et telle que l'action de $\langle n\rangle'$ soit exactement l'action de $nI_2 \in \GL{\F_p}$. On appelle $u: J(p,p) \rar P(p)$ le produit des $u_t$ pour $t \in \F_p \cup \{\infty\}$. Alors $C(p) := \ker{u}=\bigcap_{t \in \F_p \cup \{\infty\}}{\ker{u_t}}$ est un sch\'ema en groupes propre sur $\Z[1/p]$, et on note $C^0(p)$ le mod\`ele de N\'eron de la composante connexe de l'unit\'e de $C(p)_{\Q}$, qui est un sch\'ema en groupes lisse et propre sur $\Z[1/p]$.

\mytheofrench[fr-the-exact-sequence]{Le morphisme $u$ commute \`a l'action de $\GL{\F_p}$ et des op\'erateurs de Hecke, au sens que $u \circ T_{\ell} = T'_{\ell}\circ u$ pour chaque nombre premier $\ell \neq p$. De plus, il existe un morphisme explicite $R: P(p) \rar P(p)$ commutant \`a l'action des $T'_{\ell}$ et de $\GL{\F_p}$ tel que la suite de sch\'emas ab\'eliens sur $\Z[1/p]$ 
\[0 \rar C^0(p) \rar J(p,p) \overset{u}{\rar} P(p) \overset{R}{\rar} P(p),\]
soit un complexe dont les fibres sont exactes \`a isog\'enie pr\`es. 
De plus, pour tout corps $k$ de caract\'eristique distincte de $p$, les groupes de composantes connexes de $(\ker{R})_k$ et de $C(p)_k$ sont d'exposant divisant $2p(p+1)$. 
}

On d\'ecrit aussi les diff\'erentielles holomorphes sur $X(p,p)_{\C}$ ainsi que l'action des op\'erateurs de Hecke et de $\GL{\F_p}$ dont elles sont munies. Pour le reste de l'introduction, d\'efinissons la notation suivante: 

\textbf{Notations pour les repr\'esentations de $\GL{\F_p}$:}
\begin{itemize}[noitemsep,label=\tiny$\bullet$]
\item $\mathcal{D}$ est l'ensemble des caract\`eres (\`a valeurs complexes) du groupe ab\'elien $\F_p^{\times}$,
\item $\mrm{St}$ est la repr\'esentation de Steinberg de $\GL{\F_p}$: elle est de dimension $p$ et peut \^etre r\'ealis\'ee sur $\Z$,
\item pour chaque $\alpha \in \mathcal{D}$, $\mrm{St}_{\alpha}$ est la tordue de la repr\'esentation $\mrm{St}$ par le caract\`ere $\alpha(\det)$ de $\GL{\F_p}$,
\item pour tous $\alpha,\beta \in \mathcal{D}$, $\pi(\alpha,\beta)$ d\'esigne la repr\'esentation de s\'erie principale de $\GL{\F_p}$ associ\'ee aux deux caract\`eres $\alpha$ et $\beta$.
\end{itemize}
Si $V$ est un $\C[\GL{\F_p}]$-module \`a gauche, on appelle $V^{\vee}$ le $\C[\GL{\F_p}]$-module \`a \emph{droite} $\mrm{Hom}(V,\C)$.  \\

\textbf{Notations pour les espaces de formes modulaires:}
\begin{itemize}[noitemsep,label=\tiny$\bullet$]
\item $\mathscr{S}$ est l'ensemble des $(f,\mathbf{1})$, o\`u $\mathbf{1} \in \mathcal{D}$ est le caract\`ere trivial et $f \in \mathcal{S}_2(\Gamma_0(p))$ est propre pour les op\'erateurs de Hecke et normalis\'ee. 
\item $\mathscr{P}$ est l'ensemble des $(f,\chi)$, o\`u $f \in \mathcal{S}_2(\Gamma_1(p))$ est une forme modulaire propre pour les op\'erateurs de Hecke, normalis\'ee et de caract\`ere $\chi \in \mathcal{D}$ non trivial. La conjugaison complexe agit sur $\mathscr{P}$ sans point fixe en \'echangeant $(f,\chi)$ et $(\overline{f},\chi^{-1})$, o\`u $\overline{f} \in \mathcal{S}_2(\Gamma_1(p))$ est la forme modulaire dont le $q$-d\'eveloppement est $\sum_{n \geq 1}{\overline{a_n(f)}q^n}$.
\item $\mathscr{P}'$ est un ensemble de repr\'esentants de $\mathscr{P}$ modulo la conjugaison complexe. 
\item $\mathscr{C}$ d\'esigne l'ensemble des $(f,\chi)$, o\`u $f \in \mathcal{S}_2(\Gamma_1(p)\cap \Gamma_0(p^2))$ est propre pour les op\'erateurs de Hecke, normalis\'ee, de caract\`ere $\chi \in \mathcal{D}$, et telle qu'aucune tordue de $f$ par un caract\`ere dans $\mathcal{D}$ ne soit de niveau $p$.
\item $\mathscr{C}_M$ est l'ensemble des $(f,\mathbf{1}) \in \mathscr{C}$ tels que $f$ ait multiplication complexe (le corps quadratique imaginaire \'etant alors n\'ecessairement $\Q(\sqrt{-p})$ avec $p \equiv 3\pmod{4}$). 
\item $\mathscr{C}_1$ est l'ensemble des $(f,\mathbf{1}) \in \mathscr{C} \backslash \mathscr{C}_M$, et $\mathscr{C}'_1$ est un ensemble de repr\'esentants pour $\mathscr{C}_1$ modulo la torsion par l'unique caract\`ere quadratique non trivial contenu dans $\mathcal{D}$. 
\item Si $f \in \mathcal{S}_k(\Gamma_1(N))$ est une forme propre pour les op\'erateurs de Hecke, primitive et normalis\'ee, si $L$ est un corps de nombres contenant les coefficients de $f$ et $\lambda$ un id\'eal maximal de $\OO_L$, $V_{f,\lambda}$ est la repr\'esentation irr\'eductible continue $L_{\lambda}$-lin\'eaire de dimension $2$ de $G_{\Q} $ telle que, pour chaque nombre premier $q \nmid N$ distinct de la caract\'eristique r\'esiduelle de $\lambda$, le polyn\^ome caract\'eristique d'un Frobenius arithm\'etique \`a $q$ soit $X^2-a_q(f)X+q^{k-1}\chi(q)$.   
\end{itemize}

\mytheofrench[fr-differential-forms-xpp]{Il existe un isomorphisme $H^0(X(p,p)_{\C},\Omega^1) \simeq \bigoplus_{b \in \F_p^{\times}}{\mathcal{S}_2(\Gamma(p)) \otimes \Delta_{b,1}}$ qui commute \`a l'action des $T_{\ell}$ et de $\GL{\F_p}$.  
De plus, comme $\C[\GL{\F_p}]$-module \`a droite muni de l'action des $T_{\ell}$, $H^0(X(p,p)_{\C}, \Omega^1)$ est isomorphe \`a la somme directe suivante:
\begin{itemize}[noitemsep,label=$-$]
\item pour chaque $(f,\mathbf{1}) \in \mathscr{S}$ et chaque $\alpha \in \mathcal{D}$, une copie de $\mrm{St}_{\alpha}^{\vee}$, 
\item pour chaque $(f,\chi) \in \mathscr{P}$ modulo la conjugaison complexe, pour chaque $\alpha \in \mathcal{D}$, une copie de $\pi(\alpha,\alpha\chi)^{\vee}$,
\item pour chaque $(f,\chi) \in \mathscr{C}$, un copie de $C_f^{\vee}$, pour une certaine repr\'esentation cuspidale irr\'eductible $C_f$ de $\GL{\F_p}$ associ\'ee \`a une paire $\{\phi_f,\,\phi_f^p\}$ de caract\`eres de $\F_{p^2}^{\times}$.  
\end{itemize}
Pour chacun des termes ci-dessus, (o\`u l'on pose $\alpha=\mathbf{1}$ dans le troisi\`eme cas), pour tout nombre premier $\ell \neq p$ (resp. pour tout $n \in \F_p^{\times}$), $T_{\ell}$ (resp. $nI_2$) agit par multiplication par $a_{\ell}(f)\alpha(\ell)$ (resp. $\alpha^2(n)\chi(n)$).   
}

\section{Les jacobiennes modulaires tordues (Chapitres 3-4)}

\subsection{Chapitre 3: le module de Tate de $J_E^{\alpha}(p)$}

Dans le chapitre \ref{tate-modules-twists}, on utilise la relation d'Eichler-Shimura et les r\'esultats du Chapitre \ref{analytic} afin de d\'ecrire le module de Tate de la jacobienne de $X(p,p)_{\Q}$, puis de celle de $X_E^{\alpha}(p)$, o\`u $E/\Q$ est une courbe elliptique. 

Pr\'ec\'edemment, Virdol \cite[Theorem 1.1]{Virdol} a obtenu une description du module de Tate de tordues de la courbe $X(p,p)_{\Q}$ par une repr\'esentation $\rho: G_{\Q} \rar \GL{\F_p}$, ce qui lui permet (\'egalement dans \emph{op.cit.}) de d\'emontrer que la fonction $L$ de cette courbe poss\`ede dans certains cas un prolongement m\'eromorphe. Voici son th\'eor\`eme.

\theoinfrench[fr-Virdolthm]{(Virdol \cite[Theorem 1.1]{Virdol}) Soit $\rho: \mrm{Gal}(\Qbar/\Q) \rar \GL{\F_p}$ un morphisme continu de groupes. On a l'\'egalit\'e suivante de fonctions $L$
\[L(s,X(p,p)_{\rho}) = \prod_{\pi}{L(s,\rho_{\pi,\ell}\otimes (\tilde{\varphi_{\pi}}\circ \rho))},\]
o\`u
\begin{itemize}[noitemsep,label=$-$]
\item dans le produit, les $\pi = \pi_f \otimes \pi_{\infty}$ sont des repr\'esentations cuspidales automorphes de $\GL{\mathbb{A}_{\Q}}$ de poids $2$, o\`u $\pi_{\infty}$ est cohomologique et $\pi_f$ une repr\'esentation de $\GL{\mathbb{A}_{f,\Q}}$ telle que $\pi_f^{\hat{\Gamma}(p)} \neq 0$, et $\hat{\Gamma}(p) \leq \GL{\mathbb{A}_{f,\Q}}$ d\'esigne le sous-groupe $\left[I_2+p\mathcal{M}_2(\Z_p)\right]\times \prod_{\ell \neq p}{\GL{\Z_{\ell}}}$,
\item $\rho_{\pi,\ell}$ est la repr\'esentation galoisienne $\ell$-adique associ\'ee \`a $\pi$,
\item $\tilde{\varphi}_{\pi}$ d\'esigne la repr\'esentation du groupe $\GL{\F_p}$ associ\'ee \`a $\pi_f^{\hat{\Gamma}(p)}$.
\end{itemize}}

Ce chapitre se distingue des r\'esultats de Virdol sous plusieurs aspects. D'abord, les \'enonc\'es de Virdol sont exprim\'es dans un langage ad\'elique de repr\'esentations automorphes, l\`a o\`u nous gardons pour l'instant un point de vue classique. Ceci permet en particulier de rendre explicites les repr\'esentations $\tilde{\varphi_{\pi}} \circ \rho$ en termes de formes modulaires. 

Deuxi\`emement, lorsque la repr\'esentation $\rho$ du Th\'eor\`eme \ref{fr-Virdolthm} provient de la $p$-torsion d'une courbe elliptique $E/\Q$ (comme dans le cas du Th\'eor\`eme \ref{fr-XE-vs-Xrho}), Virdol calcule en r\'ealit\'e la fonction $L$ (ou plut\^ot le module de Tate) du produit $\prod_{\alpha \in \F_p^{\times}}{J_E^{\alpha}(p)}$, alors que nous d\'ecrivons le module de Tate de chaque $X_E^{\alpha}(p)$ s\'epar\'ement. 

Dans le reste de ce r\'esum\'e, $p \geq 7$ est un nombre premier et $E/\Q$ est une courbe elliptique. Fixons une base $(P,Q)$ du $\F_p$-espace vectoriel $E[p](\Qbar)$ et soit $\xi=\langle P,\,Q\rangle_{E[p]} \in \mu_p^{\times}(\Qbar)$ (o\`u $\mu_p^{\times}(\Qbar)$ est l'ensemble des racines $p$-i\`emes primitives de l'unit\'e dans $\Qbar$). Soit $\rho: G_{\Q} \rar \GL{\F_p}$ le morphisme continu de groupes tel que 

\[\forall \sigma \in G_{\Q},\, \rho(\sigma)\begin{pmatrix}\sigma(P)\\\sigma(Q)\end{pmatrix} = \begin{pmatrix}P\\Q\end{pmatrix},\]
de sorte que $\rho'(\sigma)=(\rho(\sigma)^T)^{-1}$ soit la matrice de $\sigma$ dans la base $(P,Q)$. En particulier, le caract\`ere $\det{\rho}: G_{\Q} \rar \F_p^{\times}$ est l'inverse du caract\`ere cyclotomique, et $\rho$ est l'inverse de l'anti-morphisme de groupes appel\'e $\rho$ dans le Th\'eor\`eme \ref{fr-XE-vs-Xrho}. 

Soit $F \subset \C$ un corps de nombres galoisien sur $\Q$ et dont l'anneau des entiers $\OO_F$ contient les entiers alg\'ebriques suivants: 
\begin{itemize}[noitemsep,label=\tiny$\bullet$]
\item les racines $p(p-1)^2(p+1)$-i\`emes de l'unit\'e, 
\item les coefficients de Fourier de $f$, pour tout $(f,\chi) \in \mathscr{S}\cup\mathscr{P}\cup\mathscr{C}$, 
\item pour chaque $(f,\mathbf{1}) \in \mathscr{C}_M$, pour chaque nombre premier $\ell$ d\'ecompos\'e dans $\Q(\sqrt{-p})$, les racines du polyn\^ome de Hecke $X^2-a_{\ell}(f)X+\ell$ of $f$ at $\ell$. 
\end{itemize}

Pour chaque $(f,\mathbf{1}) \in \mathscr{C}_M$, on choisit un syst\`eme compatible de caract\`eres $(\psi_{f,\lambda})_{\lambda}$ (o\`u $\lambda$ parcourt les id\'eaux maximaux de $\OO_F$), o\`u $\psi_{f,\lambda}: G_{\Q(\sqrt{-p})} \rar F_{\lambda}^{\times}$ est associ\'e \`a $f$ par la th\'eorie des formes modulaires \`a multiplication complexe. \\

Lorsque $p \equiv 3\pmod{4}$, la restriction \`a $\F_p^{\times}\SL{\F_p}$ de la repr\'esentation cuspidale $V_4$ de $\GL{\F_p}$ associ\'ee \`a la paire des deux caract\`eres d'ordre $4$ de $\F_{p^2}^{\times}$ est la somme de deux repr\'esentations irr\'eductibles $C^+$ et $C^-$. Elles sont duales l'une de l'autre et \'echang\'ees apr\`es conjugaison par $\GL{\F_p}/\F_p^{\times}\SL{\F_p}$. \\

\mytheofrench[fr-tate-module-XE]{Soient $\ell$ un nombre premier et $\lambda$ un id\'eal premier de $\OO_F$ de caract\'eristique r\'esiduelle $\ell$. Alors, pour tout $\alpha \in \F_p^{\times}$, $\Tate{\ell}{J_E^{\alpha}(p)} \otimes_{\Z_{\ell}} F_{\lambda}$ se d\'ecompose comme la somme directe des sous-$F_{\lambda}[G_{\Q}]$-modules suivants: 
\begin{itemize}[noitemsep,label=$-$]
\item Pour chaque $(f,\mathbf{1}) \in \mathscr{S}$, une copie de $V_{f,\lambda} \otimes (\mrm{St}\circ \rho)$, 
\item Pour chaque $(f,\chi) \in \mathscr{P}'$, une copie de $V_{f,\lambda} \otimes (\pi(\mathbf{1},\chi)\circ \rho)$,
\item Pour chaque $(f,\mathbf{1}) \in \mathscr{C}'_1$, une copie de $V_{f,\lambda} \otimes (C_f \circ \rho)$, o\`u $C_f$ est la repr\'esentation cuspidale de $\GL{\F_p}$ associ\'ee \`a $(f,\mathbf{1})$ par le Th\'eor\`eme \ref{fr-differential-forms-xpp},
\item Pour chaque $(f,\mathbf{1}) \in \mathscr{C}_M$, une copie de $V_{f,\lambda,E}^{\epsilon} := \mrm{Ind}_{\Q(\sqrt{-p})}^{\Q}{\psi_{f,\lambda} \otimes C^{\epsilon}(\rho_{|G_{\Q(\sqrt{-p})}})}$, o\`u $\epsilon \in \{\pm\}$ ne d\'epend que de $\xi^{\alpha\F_p^{\times 2}}$. 
\end{itemize}
Les termes de cette d\'ecomposition sont stables par les op\'erateurs de Hecke $T_{\ell}$ tels que $\ell \equiv 1\pmod{p}$, et un tel op\'erateur $T_{\ell}$ agit sur le terme associ\'e \`a $(f,\chi)$ par multiplication par $a_{\ell}(f)$. 
}

Remarquons que la repr\'esentation galoisienne apparaissant dans le Th\'eor\`eme \ref{fr-tate-module-XE} ne d\'epend pas du choix de repr\'esentant dans $\mathscr{P}'$ ou $\mathscr{C}'_1$. \\

Lorsque $(f,\chi) \in \mathscr{S} \cup \mathscr{P}' \cup \mathscr{C}'_1$, soit $R_{f,E}$ la repr\'esentation d'Artin suivante:
\begin{itemize}[noitemsep,label=\tiny$\bullet$]
\item $\mrm{St}\circ \rho$ pour $(f,\chi) \in \mathscr{S}$,
\item $\pi(\mathbf{1},\chi)\circ \rho$ pour $(f,\chi) \in \mathscr{P}'$,
\item $C_f \circ \rho$ pour $(f,\chi) \in \mathscr{C}'_1$. 
\end{itemize}

\subsection{Chapitre 4: calcul galoisien des signes d'\'equation fonctionnelle pour $(V_{f,\lambda} \otimes R_{f,E})_{\lambda}$ et $(V_{f,\lambda,E}^{\epsilon})_{\lambda}$}

Dans le chapitre \ref{automorphy}, nous discutons les propri\'et\'es de les fonctions $L$ des facteurs du module de Tate de $J_E^{\alpha}(p)$ d\'ecrits dans le Th\'eor\`eme \ref{fr-tate-module-XE}. Apr\`es des rappels sur les repr\'esentations automorphes et certains aspects de la correspondance de Langlands, nous calculons les signes d'\'equation fonctionnelle galoisiens (c'est-\`a-dire les produits des facteurs $\varepsilon$ associ\'es aux repr\'esentations galoisiennes locales d\'efinis par Deligne et Langlands et d\'ecrits par exemple dans \cite[\S 3]{NTB}) des repr\'esentations galoisiennes $V_{f,\lambda} \otimes R_{f,E}$ et $V_{f,\lambda,E}^{\epsilon}$, o\`u $\lambda$ parcourt les id\'eaux maximaux de $\OO_F$. 

M\^eme si les r\'esultats du chapitre \ref{automorphy} sont plus g\'en\'eraux, nous donnons ici le r\'esultat lorsque $E/\Q$ a bonne r\'eduction en $p$. Ce signe d\'epend du comportement de la repr\'esentation $\rho_{|G_{\Q_p}}$, que nous divisons en trois possibilit\'es. 

\defintrofrench{Soit $\rho_p: G_{\Q_p} \rar \GL{\F_p}$ une repr\'esentation continue telle que $\det{\rho_p}$ soit l'inverse du caract\`ere cyclotomique modulo $p$. On dit que $\rho_p$ est 
\begin{itemize}[noitemsep,label=$-$]
\item \emph{sauvage} si l'inertie sauvage agit non trivialement, 
\item \emph{diagonale} si $\rho_p$ est la somme directe de deux caract\`eres de $G_{\Q_p}$ \`a valeurs dans $\F_p^{\times}$,
\item \emph{Cartan} si l'image de l'inertie est contenue dans un sous-groupe de Cartan non d\'eploy\'e de $\GL{\F_p}$. 
\end{itemize}
La repr\'esentation $\rho_p$ est dans exactement une des trois situations ci-dessus. }

\mytheofrench[fr-root-number-XE-short]{Supposons que $E/\Q$ ait bonne r\'eduction en $p$. 

Soit $U$ le produit, sur chaque nombre premier $\ell$ inerte dans $\Q(\sqrt{-p})$, des $(-1)^{m_{\ell}}$, o\`u $m_{\ell} \in \Z$ est la moiti\'e de la valuation $\ell$-adique du conducteur de la repr\'esentation d'Artin $V_4 \circ \rho$. Alors il existe un signe $\varepsilon_p$ d\'ependant seulement du sch\'ema en groupes fini plat $E[p]_{\Z_p}$, de $\xi^{\F_p^{\times 2}}$ et du choix de $\psi_{f,\lambda}$ tel que le signe galoisien de l'\'equation fonctionnelle pour les facteurs apparaissant dans la d\'ecomposition de $\Tate{\ell}{J_E^{\alpha}(p)} \otimes F_{\lambda}$ d\'ecrite dans le Th\'eor\`eme \ref{fr-tate-module-XE} soit donn\'e dans la table\ref{fr-signs-surjective-general-introduction}. }

\begin{table}[htb]
\centering
\begin{tabular}{|c|c|c|c|}
\hline
& \multicolumn{3}{c|}{Type of $\rho_{|G_{\Q_p}}$}\\
\hline
& Diagonale & Sauvage & Cartan\\
\hline
$V_{f,\lambda} \otimes R_{f,E}, (f,\mathbf{1}) \in \mathscr{S}$ & $(-1)^{(p+1)/2}$ & $(-1)^{(p-1)/2}a_p(f)$ & $(-1)^{(p+1)/2}$\\
\hline
$V_{f,\lambda} \otimes R_{f,E}, (f,\chi) \in \mathscr{P}'$ & \multicolumn{3}{c|}{$(-1)^{(p-1)/2}$}\\
\hline
$V_{f,\lambda} \otimes R_{f,E}, (f,\mathbf{1}) \in \mathscr{C}'_1$ & $(-1)^{(p+1)/2}$ & $(-1)^{(p-1)/2}$ & $(-1)^{(p+1)/2}$\\
\hline 
$V_{f,\lambda,E}^{\epsilon}, (f,\mathbf{1}) \in \mathscr{C}_M$ & $(-1)^{(p+1)/4}U$ & $(-1)^{\frac{p+1}{4}}U\epsilon\varepsilon_p$ & $(-1)^{\frac{p+1}{4}}U$\\
\hline
\end{tabular}
\caption{Signes galoisiens d'\'equation fonctionnelle pour les facteurs de la d\'ecomposition de $J_{E}^{\alpha}(p)$}
\label{fr-signs-surjective-general-introduction}
\end{table}

\subsection{Chapitre 4: automorphie de $(V_{f,\lambda} \otimes R_{f,E})_{\lambda}$ and $(V_{f,\lambda,E}^{\epsilon})_{\lambda}$}

Nous pouvons d\'emontrer gr\^ace \`a la m\'ethode de Rankin-Selberg que

\mytheofrench[fr-functional-equation-JE-generic]{Soit $(f,\mathbf{1}) \in \mathscr{S}$. (resp. $(f,\chi) \in \mathscr{P}'$, resp. $(f,\mathbf{1}) \in \mathscr{C}'_1$) Si la repr\'esentation d'Artin est automorphe, alors la fonction $L$ de la repr\'esentation galoisienne autoduale (\`a tordue \`a la Tate pr\`es) $V_{f,\lambda} \otimes R_{f,E}$ poss\`ede un prolongement holomorphe et v\'erifie une \'equation fonctionnelle dont le signe est donn\'e (lorsque $E$ a bonne r\'eduction en $p$) par la table \ref{fr-signs-surjective-general-introduction}. 

Pour tout $(f,\mathbf{1}) \in \mathscr{C}_M$, la fonction $L$ de la repr\'esentation galoisienne autoduale (\`a tordue \`a la Tate pr\`es) $V_{f,\lambda,E}^{\epsilon}$ poss\`ede un prolongement m\'eromorphe et v\'erifie une \'equation fonctionnelle dont le signe est donn\'e dans la table \ref{fr-signs-surjective-general-introduction} lorsque $E/\Q$ a bonne r\'eduction en $p$.} 

Comme dans le th\'eor\`eme \ref{fr-root-number-XE-short}, l'hypoth\`ese de bonne r\'eduction en $p$ pour $E/\Q$ est faite pour simplifier l'introduction, et n'est pas utilis\'ee dans le chapitre \ref{automorphy}. 

Lorsque $\rho$ est surjective, les repr\'esentations $R_{f,E}$ sont d'image finie, irr\'eductibles de grande dimension, et leur image n'est pas r\'esoluble. Ceci ne change pas lorsque l'on restreint $\rho$ au groupe de Galois d'un corps de nombres totalement r\'eel. Par cons\'equent, il semble que les techniques actuelles ne permettent pas de montrer que les repr\'esentations $R_{f,E}$ sont automorphes. 

Une autre possibilit\'e pour comprendre la repr\'esentation $V_{f,\lambda} \otimes R_{f,E}$ pourrait \^etre de la r\'eduire modulo un nombre premier. Nous d\'emontrons, gr\^ace aux deux r\'esultats de fonctorialit\'e de Langlands suivants,
\begin{itemize}[noitemsep,label=\tiny$\bullet$]
\item pour les puissances sym\'etriques pour le groupe $\operatorname{GL}(2)$ sur un corps totalement r\'eel d\'emontr\'eee par Newton et Thorne \cite{NewThor}\footnote{Le cas du corps de base $\Q$, qui suffit \`a nos besoins, est trait\'e dans les articles publi\'es \cite{NewThorQ1,NewThorQ2}. },
\item pour les produits tensoriels de la forme $\operatorname{GL}(2) \times \operatorname{GL}(n)$, montr\'e sous certaines hypoth\`eses dans une pr\'epublication d'Arias-de-Reyna et Dieulefait \cite{AdRD},
\end{itemize}
que $V_{f,\lambda}\otimes R_{f,E}$ sont automorphes modulo $p$ en un sens que nous pr\'ecisons. 

Pour r\'ef\'erence, nous donnons ici une traduction en fran\c cais de l'\'enonc\'e du r\'esultat d'Arias-de-Reyna et Dieulefait.  

\theoinfrench[fr-adrd]{(Arias-de-Reyna -- Dieulefait \cite[Theorem 1.1]{AdRD}) Soit $f \in \mathcal{S}_k(\Gamma(1))$ une forme modulaire parabolique propre pour les op\'erateurs de Hecke et soit $\rho_{\bullet}(f)$ le syst\`eme compatible de repr\'esentations associ\'e. Soit $n > 0$ un entier et $\pi$ une repr\'esentation automorphe r\'eguli\`ere alg\'ebrique cuspidale polaris\'ee de $GL_n(\mathbb{A}_{\Q})$ de niveau premier \`a $3$. Soit $r_{\bullet}(\pi)$ le syst\`eme compatible de repr\'esentations galoisiennes associ\'e \`a $\pi$. Supposons que $\rho_{\bullet}(f) \otimes r_{\bullet}(\pi)$ soit un syst\`eme compatible r\'egulier et irr\'eductible. Alors $\rho_{\bullet}(f) \otimes r_{\bullet}(\pi)$ est automorphe au sens suivant: il existe une repr\'esentation automorphe r\'eguli\`ere alg\'ebrique cuspidale polaris\'ee $f \otimes \pi$ de $GL_{2n}(\mathbb{A}_{\Q})$ tel que le syst\`eme compatible de repr\'esentations galoisiennes associ\'e \`a $f \otimes \pi$ soit isomorphe \`a $\rho_{\bullet}(f) \otimes r_{\bullet}(\pi)$.}

\mytheofrench[fr-residual-automorphy-XE]{Supposons que $\rho$ soit surjective et provienne d'une courbe elliptique $E/\Q$ ayant bonne r\'eduction en $3$. Supposons de plus que le Th\'eor\`eme \ref{fr-adrd} soit vrai. Il existe un corps de nombres $F' \supset F$ et un id\'eal premier $\mathfrak{p}$ de $F'$ de caract\'eristique r\'esiduelle $p$ telle que, pour tout $(f,\chi) \in \mathscr{S} \cup \mathscr{P}'_1 \cup \mathscr{C}'_1$, la semi-simplification de la r\'eduction modulo $\mathfrak{p}$ de $V_{f,\mathfrak{p}} \otimes R_{f,E}$ soit la somme directe des semi-simplifications de r\'eductions modulo $\mathfrak{p}$ de repr\'esentations galoisiennes \`a coefficients dans $(F')_{\mathfrak{p}}$ provenant de syst\`emes compatibles associ\'es \`a des repr\'esentations automorphes r\'eguli\`eres alg\'ebriques cuspidales et polaris\'ees. }

Bien que ce r\'esultat soit int\'eressant, il ne semble pas qu'il ait d'application claire \`a la d\'etermination des points rationnels de $X_{E}^{\alpha}(p)$.

\section{Cas o\`u l'image de $\rho$ est contenue dans le normalisateur d'un sous-groupe de Cartan non d\'eploy\'e (Chapitre 5)}

Dans le chapitre \ref{small-image}, nous supposons que $\rho(G_{\Q})$ est contenue dans le normalisateur $N$ d'un sous-groupe de Cartan non d\'eploy\'e $C$. Soit $K/\Q$ le corps quadratique imaginaire correspondant au caract\`ere $G_{\Q} \rar N \rar N/C \simeq \{\pm 1\}$. Pour simplifier, nous supposons dans cette introduction que $p \geq 11$. Dans ce cas, la Proposition \ref{almost-never-ramified} assure que $K$ est inerte en $p$.

\subsection{D\'ecomposition du module de Tate et \'equation fonctionnelle}

Dans une telle situation, la d\'ecomposition obtenue dans le th\'eor\`eme \ref{fr-tate-module-XE} n'est pas optimale et peut \^etre affin\'ee. Fixons une fois pour toutes un isomorphisme de groupes $\iota$ entre $C$ et $\F_{p^2}^{\times}$ pr\'eservant la norme et la trace. L'isomorphisme $\iota$ n'est pas canonique, mais la paire $\{\iota,\iota^p\}$ l'est: par cons\'equent, $\iota$ permet d'identifier canoniquement les paires $\{\psi,\psi^p\}$ de caract\`eres de $C$ et les paires $\{\phi,\phi^p\}$ de caract\`eres de $\F_{p^2}^{\times}$.

\mytheofrench[fr-tate-module-XE-cartan]{Soit $\lambda$ un id\'eal premier de $\OO_F$ de caract\'eristique r\'esiduelle $\ell$. Pour chaque $\alpha \in \F_p^{\times}$, $\Tate{\ell}{J_E^{\alpha}(p)}\otimes_{\Z_{\ell}} F_{\lambda}$ est la somme directe des sous-$G_{\Q}$-modules suivants qui sont stables par les op\'erateurs de Hecke:
\begin{itemize}[noitemsep,label=\tiny$\bullet$]
\item pour chaque $(f,\chi) \in \mathscr{S}\cup\mathscr{P}'\cup\mathscr{C}'_1\cup\mathscr{C}_M$, une copie de $V_{f,\lambda} \otimes \theta$, o\`u $\theta$ est un caract\`ere de $\mrm{Gal}(K(\mu_p)/\Q)$ ne d\'ependant pas de $\rho$ ou $\lambda$; on somme sur deux caract\`eres $\theta$ possibles si $(f,\chi) \in \mathscr{P}'\cup\mathscr{C}'_1$ et un seul sinon,
\item pour chaque $(f,\chi) \in \mathscr{S}\cup\mathscr{P}'$, pour chaque paire $\{\psi,\psi^p\}$ de caract\`eres distincts $C \rar F^{\times}$ telle que $\psi(tI_2)=\chi(t)$ pour tout $t \in \F_p^{\times}$, une copie de $V_{f,\lambda} \otimes \left[\left(\mrm{Ind}_C^{N}{\psi}\right)\circ \rho\right]$, 
\item pour chaque $(f,\mathbf{1}) \in \mathscr{C}'_1$ (resp. $(f,\mathbf{1}) \in \mathscr{C}_M$), pour chaque paire $\{\psi,\psi^p\} \neq \{\phi_f,\phi_f^p\}$ de caract\`eres distincts $C/\F_p^{\times}I_2 \rar F^{\times}$ (resp. modulo \'egalement la relation d'\'equivalence d\'efinie par: $\{\psi,\psi^p\} \sim \{\theta,\theta^p\}$ si et seulement si l'un des caract\`eres $\psi\theta^{-1}$ et $\psi^p\theta^{-1}$ est d'ordre au plus deux), une copie de $V_{f,\lambda} \otimes \left[\left(\mrm{Ind}_C^{N}{\psi}\right)\circ \rho\right]$.
\end{itemize}
L'op\'erateur de Hecke $T_{\ell}$ (pour $\ell \equiv 1 \pmod{p}$ premier) agit sur un facteur de l'une des formes ci-dessus par multiplication par $a_{\ell}(f)$.  
}

Lorsque l'on les compare aux facteurs apparaissant dans la d\'ecomposition d\'ecrite dans le th\'eor\`eme \ref{fr-tate-module-XE}, les facteurs d\'ecrits ci-dessus sont bien plus simples, et leurs propri\'et\'es en tant que repr\'esentations galoisiennes bien mieux comprises. Dans le reste du chapitre \ref{small-image}, on applique la strat\'egie de Mazur dans le but de d\'eterminer les points rationnels de $X_E^{\alpha}(p)$. Tout d'abord, on montre que les fonctions $L$ des facteurs apparaissant dans le th\'eor\`eme \ref{fr-tate-module-XE-cartan} poss\`edent des prolongements entiers, v\'erifient des \'equations fonctionnelles dont on calcule le signe. 

\mytheofrench[fr-signs-functional-equation-cartan-intro]{Chacune des repr\'esentations d\'ecrites dans le th\'eor\`eme \ref{fr-tate-module-XE-cartan} est autoduale (\`a une torsion \`a la Tate pr\`es), et sa fonction $L$ poss\`ede un prolongement analytique et satisfait une \'equation fonctionnelle de signe $-1$, sauf dans les cas suivants:
\begin{itemize}[noitemsep,label=$-$]
\item $V_{f,\lambda} \otimes \left[\left(\mrm{Ind}_C^{N}{\psi}\right)\circ \rho\right]$ pour $(f,\mathbf{1}) \in \mathscr{S}$, lorsque $\left(\mrm{Ind}_C^{N}{\psi}\right)\circ \rho$ est non ramifi\'ee en $p$,
\item $V_{f,\lambda} \otimes \left[\left(\mrm{Ind}_C^{N}{\psi}\right)\circ \rho\right]$ pour $(f,\mathbf{1}) \in \mathscr{C}'_1 \cup \mathscr{C}_M$ si la condition suivante est v\'erifi\'ee. Soit $M/\Q_p$ l'unique extension quadratique non ramifi\'ee et $\mathfrak{a}: \mrm{Gal}(\overline{\Q_p}/M) \rar \widehat{M^{\times}}$ l'isomorphisme du corps de classe (normalis\'e de sorte \`a envoyer un Frobenius arithm\'etique sur une uniformisante); en particulier, $\mathfrak{a}$ poss\`ede un sous-quotient $\overline{\mathfrak{a}}$ envoyant le groupe d'inertie $I_p \triangleleft \mrm{Gal}(\overline{\Q_p}/\Q_p)$ dans $\OO_M^{\times}/(1+p\OO_M) \simeq \F_{p^2}^{\times}$. Le signe de l'\'equation fonctionnelle est $1$ si et seulement si les paires $\{\phi_f \circ \overline{\mathfrak{a}}, \phi_f^p \circ \overline{\mathfrak{a}}\}$ et $\{\psi\circ\rho_{|I_p}, \psi^p\circ\rho_{|I_p}\}$ de caract\`eres $I_p \rar F^{\times}$ sont \'egales. 
\end{itemize}
Dans ces deux cas, $E$ ainsi que ses tordues quadratiques ont mauvaise r\'eduction en $p$. Dans le premier, on a de plus $\rho(I_p)=3C$, $p \equiv 2, 5\pmod{9}$ et $\psi$ est d'ordre exactement $3$. 
}

\subsection{La conjecture de Bloch-Kato en rang z\'ero pour les facteurs de $J_E^{\alpha}(p)$}

Le r\'esultat suivant de ce chapitre est la preuve de la conjecture de Bloch-Kato en rang z\'ero pour la plupart des facteurs quadri-dimensionnels du module de Tate de $J_E^{\alpha}(p)$. Pour ce faire, on utilise le syst\`eme d'Euler des \'el\'ements de Beilinson-Flach construit par Kings, Loeffler et Zerbes \cite{KLZ15}. Pour pouvoir utiliser les r\'esultats de \cite{KLZ15}, il est n\'ecessaire de montrer que ces repr\'esentations satisfont les \emph{hypoth\`eses de grande image} requises, ce que nous faisons en \'etendant les r\'esultats de \cite[Section 4.4]{bigimage}. Ceci est l'un des points principaux du chaputre \ref{obstructions-euler} sur lequel nous reviendrons.

Nous repla\c cant dans le cadre du th\'eor\`eme \ref{tate-module-XE-cartan}, soit $(f,\chi),(f',\chi') \in \mathscr{S} \cup \mathscr{P}' \cup \mathscr{C}'_1\cup\mathscr{C}_M$, et soient $\{\psi,\psi^p\}, \{\psi',(\psi')^p\}$ deux paires de caract\`eres distincts de $C$ \`a valeurs dans $F^{\times}$. On dit que les repr\'esentations galoisiennes \[V_{f,\lambda} \otimes \left[\left(\mrm{Ind}_C^{N}{\psi}\right)\circ \rho\right],\, V_{f',\lambda} \otimes \left[\left(\mrm{Ind}_C^{N}{\psi'}\right)\circ \rho\right] \] \emph{sont dans la m\^eme orbite sous l'action de Galois} s'il existe $\sigma \in G_{\Q}$ tel que pour tout nombre premier $q$ sauf un nombre fini, l'on a  
\[\mrm{Tr}\left(\Fr_{q}\mid V_{f',\lambda} \otimes \left[\left(\mrm{Ind}_C^{N}{\psi'}\right)\circ \rho\right]\right)  = \sigma \left[\mrm{Tr}\left(\Fr_{q} \mid V_{f,\lambda} \otimes \left[\left(\mrm{Ind}_C^{N}{\psi}\right)\circ \rho\right]\right)\right].\] Notons que cette notion ne d\'epend pas du choix de $\lambda$.

\mytheofrench[fr-bloch-kato-zero]{Soit $V$ l'un des facteurs quadridimensionnels apparaissant dans le th\'eor\`eme \ref{fr-tate-module-XE-cartan} et associ\'e au couple $((f,\chi),\{\psi,\psi^p\})$. Pour chaque $\alpha \in \F_p^{\times}$, $J_E^{\alpha}(p)$ poss\`ede un quotient $A$ stable par les op\'erateurs de Hecke tel que $\Tate{\ell}{A} \otimes_{\Z_{\ell}} F_{\lambda}$ s'identifie \`a la somme directe des facteurs apparaissant dans le th\'eor\`eme \ref{fr-tate-module-XE-cartan} qui sont dans la m\^eme orbite que $V$ sous l'action de Galois. Supposons de plus que $f$ n'ait pas de multiplication complexe et que $L(V,1) \neq 0$ (ainsi $(f,\chi) \in \mathscr{S} \cup \mathscr{C}'_1$). Alors le groupe $A(\Q)$ est fini.}

Pour construire le quotient $A$, on utilise le fait que le choix de base $(P,Q)$ pour $E[p]$ d\'efinit un isomorphisme $\coprod_{\alpha \in \F_p^{\times}}{X_E^{\alpha}(p)_{\Qbar}} \cong X(p,p)_{\Qbar}$, ce qui permet de munir $\coprod_{\alpha \in \F_p^{\times}}{X_E^{\alpha}(p)_{\Qbar}}$ d'une action de $\GL{\F_p}$.

\subsection{Le th\'eor\`eme d'immersion formelle}

Supposons maintenant que la representation $V$ du th\'eor\`eme \ref{fr-bloch-kato-zero} rel\`eve du premier cas du th\'eor\`eme \ref{fr-signs-functional-equation-cartan-intro}, ce qui implique notamment que $p \equiv 2, 5\pmod{9}$. 

Soient $c \in X_E^{\alpha}(p)(\Qbar)$ une pointe et $q \equiv 1\pmod{p}$ un nombre premier. Le morphisme $x \longmapsto (T_q-q-1)(x-c): X_E^{\alpha}(p)_{\Qbar} \rar J_E^{\alpha}(p)_{\Qbar}$ est en fait d\'efini sur $\Q$ et ne d\'epend pas du choix de $c$, parce que $T_q-q-1$ annule les diviseurs support\'es sur le sous-sch\'ema cuspidal. Consid\'erons sa composition $f$ avec la projection $J_E^{\alpha}(p) \rar A$, o\`u $A$ est le quotient d\'efini par le th\'eor\`eme \ref{fr-bloch-kato-zero}. Gr\^ace \`a notre description des diff\'erentielles sur $\C$ de la courbe obtenue dans le chapitre \ref{analytic}, nous pouvons appliquer un argument analogue \`a l'immersion formelle de Mazur (plut\^ot formul\'e dans le langage de la m\'ethode de Chabauty-Coleman d\'ecrite dans \cite{PoonenMcCallum}) au morphisme $f: X_E^{\alpha}(p) \rar A$. 

\defintrofrench{Pour chaque racine cubique de l'unit\'e $j \in \F_{p^2}^{\times}$, soit \[u_j=\sum_{\substack{b \in \F_p\\1+b\sqrt{-3} \in \F_{p^2}^{\times 3}j}}{e^{\frac{2i\pi b}{p}}}.\] Nous disons qu'un nombre premier $\ell$ est \emph{exceptionnel} pour $p$ si $\ell \equiv \pm 1 \pmod{p}$ et $\Z[e^{2i\pi/p}]$ poss\`ede un id\'eal premier $\lambda$ de caract\'eristique r\'esiduelle $\ell$ contenant chaque $u_j$. } 

Des calculs avec Magma \cite{magma} montrent que lorsque $p \leq 1500$, les trois $u_j$ engendrent l'id\'eal unit\'e de $\Z[e^{\frac{2i\pi}{p}}]$. Nous sommes donc tent\'es de conjecturer qu'il n'y a pas de nombres premiers exceptionnels quel que soit $p$.   

\mytheofrench[fr-formal-immersion-diophantine]{Supposons que $\rho: G_{\Q} \rar N$ soit surjective et soit $N_{\rho}$ son conducteur d'Artin. Soit $\psi: C \rar F^{\times}$ un caract\`ere d'ordre $3$ et supposons qu'il existe $(f,\mathbf{1}) \in \mathscr{S}$ tel que la valeur centrale de la fonction $L$ de la repr\'esentation $V_{f,\lambda} \otimes \left[\left(\mrm{Ind}_C^N{\psi}\right)\circ \rho\right]$ soit non nulle: en particulier, l'image par $\rho$ d'un sous-groupe d'inertie en $p$ est le sous-groupe d'indice $3$ de $C$. 

Soit $E'/\Q$ une courbe elliptique telle que les $G_{\Q}$-modules $E[p](\Qbar)$ et $E'[p](\Qbar)$ soient isomorphes. Alors $E'$ a pour conducteur $p^2N_{\rho}L$, o\`u $L$ est un produit de nombres premiers exceptionnels pour $p$ deux \`a deux distincts et ne divisant pas $pN_{\rho}$. }

Notons que les r\'esultats effectifs connus sur le th\'eor\`eme d'image ouverte de Serre (par exemple \cite[Theorems 1.5, 1.10]{Zywina-Surj}), ne permettent pas de montrer \`a partir de l\`a qu'une courbe $E'$ comme dans le th\'eor\`eme ci-dessus non CM ne pourrait exister. Ces r\'esultats ne semblent pas non plus suffisants pour montrer que $E$ et $E'$ sont isog\`enes, m\^eme \`a supposer qu'elles sont toutes deux de conducteur $p^2N_{\rho}$. 

\subsection{Applications aux congruences \`a la courbe elliptique CM d'\'equation $y^2=x^3-p$}

Dans \cite{Kraus-Thesis}, Kraus calcule le poids de Serre de la $p$-torsion d'une courbe elliptique \`a r\'eduction additive. Lorsque $p \equiv 5\pmod{9}$ et $E$ est la courbe elliptique $E_p$ d'\'equation de Weierstrass $y^2=x^3-p$, on peut v\'erifier que:
\begin{itemize}[noitemsep,label=\tiny$\bullet$]
\item $\rho(G_{\Q})=N$,
\item l'image par $\rho$ d'un sous-groupe d'inertie en $p$ est le sous-groupe $3C$,
\item $\left(\mrm{Ind}_C^N{\psi}\right) \circ \rho$ est associ\'ee \`a l'unique forme modulaire primitive $g \in \mathcal{S}_1(\Gamma_1(108))$ \`a multiplication complexe par $\Q(\sqrt{-3})$ et dont les coefficients de Fourier sont rationnels.
\end{itemize}

Lorsque $p$ n'est pas trop grand, on peut v\'erifier num\'eriquement gr\^ace \`a Magma \cite{magma} qu'il n'y a pas de nombres premiers exceptionnels pour $p$, et avec PARI/GP \cite{pari} que, pour au moins un $(f,\mathbf{1}) \in \mathscr{S}$, la valeur centrale de la fonction $L$ de Rankin-Selberg $L(f \times g,1)$ apparaissant dans le th\'eor\`eme \ref{fr-formal-immersion-diophantine} ne s'annule pas. Cremona and Freitas ont montr\'e \cite[Theorem 1.3]{Cremona-Freitas} que deux courbes elliptiques dont les conducteurs sont inf\'erieurs \`a $500000$ et dont les $G_{\Q}$-modules de $q$-torsion (sur $\Qbar$) sont isomorphes sont n\'ecessairement isog\`enes lorsque $q \geq 19$. De ces consid\'erations, l'on peut d\'eduire:

\mycorfrench{Soit $p = 23$ et soit $E'$ une courbe elliptique telle que $E'[p](\Qbar)$ et $E_p[p](\Qbar)$ sont isomorphes comme $G_{\Q}$-modules. Alors $E'$ et $E_p$ sont isog\`enes.   }

En g\'en\'eral, on s'attend \`a ce qu'il existe $(f,\mathbf{1}) \in \mathscr{S}$ tel que $L(f \times g, 1) \neq 0$, mais cela semble difficile \`a montrer en toute g\'en\'eralit\'e. Le r\'esultat suivant de Michel, d\'emontr\'e par des techniques de th\'eorie analytique des nombres, montre que tel est le cas lorsque $p$ est assez grand.    

\theoinfrench[fr-michel]{(Michel, \cite[Appendix, Theorem 0.1]{Michel}) Soit $g \in \mathcal{S}_1(\Gamma_1(N))$ une forme primitive normalis\'ee et propre pour les op\'erateurs de Hecke de caract\`ere $\chi$. Il existe une constante $C_g$ v\'erifiant la propri\'et\'e suivante: pour tout nombre premier $p > C_g$, si, pour tout $f \in \mathcal{S}_2(\Gamma_0(p))$, l'on a $L(f \times g,1)=0$, alors $g$ est \`a coefficients r\'eels et $\chi(p)=1$.}

\mytheofrench{Il existe une constante $p_0 \geq 11$ v\'erifiant la propri\'et\'e suivante. Soit $p \geq p_0$ un nombre premier congru \`a $5$ modulo $9$ et $E$ la courbe elliptique \`a multiplication complexe par $\Q(\sqrt{-3})$ donn\'ee par l'\'equation de Weierstrass minimale $y^2=x^3-p$ (de conducteur $108mp^2$, o\`u $m=1$ si $p \equiv 3\pmod{4}$ et $m=4$ si $p \equiv 1\pmod{4}$). Alors, pour toute courbe elliptique $E'/\Q$ telle que $E'[p](\Qbar)$ et $E[p](\Qbar)$ soient isomorphes comme $G_{\Q}$-modules, le conducteur de $E'$ est $108mp^2L$, o\`u $L$ est un produit de nombres premiers exceptionnels pour $p$, deux \`a deux distincts. En particulier, s'il n'y a pas de nombres premiers exceptionnels pour $p$, $E'$ est de conducteur $108mp^2$ et $j(E') \in \Z$. }

\section{R\'esultats de grande image pour le produit tensoriel de deux formes modulaires (Chapitre 6)}

Le chapitre \ref{obstructions-euler} reproduit presque exactement la pr\'e-publication \cite{Studnia} par l'auteur de cette th\`ese. Le projet provient de la tentative de l'auteur d'appliquer la th\'eorie des syst\`emes d'Euler dans la preuve du th\'eor\`eme \ref{fr-bloch-kato-zero}. 

Soit $\OO_K$ l'anneau de valuation d'une extension finie $K/\Q_p$ et soit $T$ un $\OO_K$-module libre de rang fini munie d'une action continue et $\OO_K$-lineaire de $G_{\Q}$. On peut voir un syst\`eme d'Euler pour $T$ comme une collection de classes de cohomologie $c_n \in H^1(\Q(\mu_n),T)$ pour certains entiers $n \geq 1$ (comprenant notamment $n=1$) satisfaisant certaines compatibilit\'es vis-\`a-vis de la corestriction. La th\'eorie des syst\`emes d'Euler, telle qu'elle est d\'ecrite dans les r\'ef\'erences classiques \cite{Rubin-ES, MR-ES}, permet de contr\^oler la taille des groupes de Selmer du $G_{\Q}$-module $\mrm{Hom}_{\OO_K}(T,K/\OO_K(1))$ lorsque l'on dispose d'un syst\`eme d'Euler $(c_n)_n$ pour $T$ tel que $c_1 \neq 0$. N\'eanmoins, pour pouvoir appliquer cette th\'eorie, l'image de $G_{\Q}$ dans $\mrm{Aut}(T)$ doit v\'erifier certaines hypoth\`eses, souvent des variantes ou renforcements de la condition suivante, que nous appelons \emph{condition d'image adapt\'ee}: il existe $\sigma \in G_{\Q(\mu_{p^{\infty}})}$ tel que $T/(\sigma-1)T \otimes_{\OO_K} K$ soit une droite sur $K$. 

Dans \cite{KLZ15}, Kings, Loeffler et Zerbes construisent un syst\`eme d'Euler pour pour certaines des repr\'esentations apparaissant dans le th\'eor\`eme \ref{fr-tate-module-XE-cartan}, et en d\'eduisent les bornes sur la taille du groupe de Selmer utilis\'ees dans la preuve du th\'eor\`eme \ref{bloch-kato-zero} en supposant que la condition d'image adapt\'ee est v\'erifi\'ee. Dans \cite[Section 4.4]{bigimage}, Loeffler montre que cette condition est en fait souvent v\'erifi\'ee:

\theoinfrench[fr-loefflertheo]{(Loeffler \cite{bigimage}) Soient $f \in \mathcal{S}_k(\Gamma_1(N_f)), g \in \mathcal{S}_1(\Gamma_1(N_g))$ deux formes primitives normalis\'ees et propres pour les op\'erateurs de Hecke de caract\`eres respectifs $\varepsilon_f, \varepsilon_g$ telles que $k > 1$ et $f$ ne soit pas \`a multiplication complexe. Soit $L$ un corps de nombres contenant les coefficients de $f$ et $g$, et d'anneau des entiers $\OO_L$.
\begin{itemize}[noitemsep,label=$-$]
\item S'il existe un id\'eal maximal $\lambda$ de $\OO_L$ tel que $V_{f,\lambda} \otimes_{L_{\lambda}} V_{g,\lambda}$ poss\`ede une image adapt\'ee, alors $\varepsilon_f\varepsilon_g \neq 1$.
\item Supposons que $\varepsilon_f\varepsilon_g \neq 1$ et que $N_f$ et $N_g$ soient premiers entre eux. Alors l'ensemble des id\'eaux maximaux $\lambda$ de $\OO_L$ tels que $V_{f,\lambda}\otimes V_{g,\lambda}$ n'ait pas d'image adapt\'ee est fini.   
\end{itemize}}

Dans le m\^eme article, l'auteur pose la question suivante.  

\questinfrench[fr-loeffler-question-intro]{(Loeffler, \cite[Remark 4.4.2]{bigimage}) Reprenons les notations du th\'eor\`eme \ref{loefflertheo} et supposons que $\varepsilon_f\varepsilon_g \neq 1$. Est-il vrai que, pour tout id\'eal maximal $\lambda$ de $\OO_L$ sauf \'eventuellement un nombre fini, $V_{f,\lambda}\otimes V_{g,\lambda}$ ait une image adapt\'ee?}

Le th\'eor\`eme \ref{fr-loefflertheo} ne permet pas \`a lui seul d'utiliser les r\'esultats sur les syst\`emes d'Euler requis pour la preuve de \ref{fr-bloch-kato-zero}. Dans la deuxi\`eme situation du th\'eor\`eme \ref{fr-signs-functional-equation-cartan-intro}, les conducteurs des deux formes modulaires sont divisibles par $p^2$, et le th\'eor\`eme \ref{fr-loefflertheo} ne peut donc pas montrer que la condition d'image adapt\'ee est v\'erifi\'ee pour tout id\'eal premier sauf un nombre fini. C'est l\`a la motivation pour notre \'etude de la question \ref{fr-loeffler-question-intro}. Nous donnons en premier lieu des conditions suffisantes plus g\'en\'erales sous lesquelles on peut r\'epondre \`a la question \ref{fr-loeffler-question-intro} par l'affirmative. 

\mytheofrench[fr-better-sufficient-conditions]{Gardons les notations du th\'eor\`eme \ref{fr-loefflertheo} et supposons que $\varepsilon_f\varepsilon_g \neq 1$. On peut r\'epondre par l'affirmative \`a la question \ref{fr-loeffler-question-intro} si l'une des conditions suivantes est v\'erifi\'ee:
\begin{itemize}[noitemsep,label=\tiny$\bullet$]
\item Pour tout caract\`ere de Dirichlet (primitif) $\chi$ de niveau $C$ tel que $f \otimes \chi$ soit conjugu\'ee \`a $f$ par l'action de $G_{\Q}$, $C$ et $N_g$ sont premiers entre eux. 
\item Tout caract\`ere de Dirichlet (primitif) $\chi$ tel que $f \otimes \chi$ soit conjugu\'e \`a $f$ par l'action de $G_{\Q}$ est pair.
\item $\varepsilon_f^2=\varepsilon_g^2=1$ et $f \otimes \varepsilon_g$ n'est pas conjugu\'e \`a $f$ par l'action de $G_{\Q}$.
\item Le sous-groupe fini $\varepsilon_g\left(\bigcap_{\chi}{\ker{\chi}}\right)$ de $\C^{\times}$ contient $-1$ mais pas $\pm i$, o\`u $\chi$ parcourt les caract\`eres de Dirichlet (primitifs) tels que $f \otimes \chi$ soit conjugu\'ee \`a $f$ sous l'action de $G_{\Q}$.
\end{itemize}}

Ce r\'esultat permet d'am\'eliorer la conjecture de Birch and Swinnerton-Dyer d\'emontr\'ee dans \cite[Theorem 11.7.4]{KLZ15} pour les tordues de courbes elliptiques par les repr\'esentations d'Artin bidimensionnelles, impaires et irr\'eductibles: nous montrons dans le corollaire \ref{BSD-improved} que, dans l'\'enonc\'e de \emph{loc.cit.}, nous n'avons plus besoin de supposer que la courbe elliptique et la repr\'esentation d'Artin aient des conducteurs premiers.

Nous montrons ensuite que la r\'eponse g\'en\'erale \`a la question \ref{fr-loeffler-question-intro} est \og non\fg{}, et d\'ecrivons des familles infinies de contre-exemples gr\^ace au r\'esultat suivant.  

\mytheofrench[fr-answer-is-negative]{Soit $f \in \mathcal{S}_k(\Gamma_1(N))$ une forme propre pour les op\'erateurs de Hecke, primitive et normalis\'ee telle que $k$ soit de poids pair. Les conditions suivantes sont \'equivalentes:
\begin{itemize}[noitemsep,label=$-$]
\item Pour tout corps quadratique $K$, on peut r\'epondre \og oui\fg{} \`a la question \ref{fr-loeffler-question-intro} pour tout couple $(f,g)$ sauf un nombre fini d'entre eux, o\`u $g$ est une forme modulaire parabolique primitive propre pour les op\'erateurs de Hecke et normalis\'ee poss\'edant multiplication r\'eelle ou complexe par $K$. 
\item La r\'eponse \`a la question \ref{fr-loeffler-question-intro} est \og oui\fg{} pour tout couple $(f,g)$, o\`u $g$ est une forme modulaire parabolique de poids un, primitive, propre pour les op\'erateurs de Hecke et normalis\'ee.
\item Tout caract\`ere de Dirichlet (primitif) $\chi$ tel que $f \otimes \chi$ soit conjugu\'e \`a $f$ sous l'action de $G_{\Q}$ est pair.
\end{itemize}
}

Nous d\'ecrivons \'egalement des contre-exemples explicites provenant d'une obstruction diff\'erente de celle apparaissant dans la preuve du th\'eor\`eme \ref{fr-answer-is-negative}. 

\mytheofrench[fr-sporadic-counter-examples]{Il existe des formes modulaires $f \in \mathcal{S}_2(\Gamma_0(63)), g \in \mathcal{S}_1(\Gamma_1(1452))$ propres pour les op\'erateurs de Hecke, primitives et normalis\'ees, \`a coefficients respectivement contenus dans $\Q(\sqrt{3}),\Q(\sqrt{-3})$ telles que, pour tout id\'eal maximal $\mathfrak{p}$ de $L=\Q(e^{2i\pi/12})$ de caract\'eristique r\'esiduelle $p \equiv 5,7 \pmod{12}$ sauf un nombre fini, $V_{f,\mathfrak{p}} \otimes V_{g,\mathfrak{p}}$ ne v\'erifie pas la condition d'image adapt\'ee. }

\renewcommand{\thesection}{\thechapter.\arabic{section}}
\renewcommand{\thesubsection}{\thesection.\arabic{subsection}}   

\end{otherlanguage}

%% file: moduli-0.tex
\chapter{Moduli spaces and twisted modular curves}
\label{moduli-spaces}

The material of this Chapter is not new or particularly difficult: most of these results seem to be folklore, but it is not clear whether they have been precisely spelled out somewhere. 

Its purpose is to prove\footnote{following discussions with Lo\"ic Merel and Barry Mazur} that two constructions for the curve of interest -- both as a moduli space (which was already known, and discussed in \cite[Proposition 1]{KO} or \cite[I.7-9]{Mazur-Open}) and a Galois twist (as discussed in \cite{Virdol}) -- are equivalent, and moreover that this equivalence respects all the natural constructions one wishes to do for modular curves: Atkin-Lehner automorphisms, adding isogeny structures and Hecke operators. While \cite{KO} or \cite{Mazur-Open} sketch the construction of the moduli space when the base ring is a field, we prove that the construction is valid over more general base rings (assumed to be regular excellent Noetherian to consider compactified moduli schemes).  

\theoi[XG-exists-pol-unpol]{Let $N \geq 3$ be an integer and $R$ be a regular excellent Noetherian $\Z[1/N]$-algebra. Let $m_0,m_1 \geq 1$ be two integers such that $m_0, m_1, N$ are pairwise coprime, $m_0$ is square-free and $m_1 \in R^{\times}$, and write $m=m_0m_1$. Let $G$ be a finite \'etale commutative group scheme over $R$ which is \'etale-locally isomorphic to $(\Z/N\Z)^{\oplus 2}$. Then there is a proper flat $R$-scheme $X_G(N,\Gamma_0(m))$ of relative dimension one with geometrically reduced fibres satisfying the following properties:
\begin{itemize}[noitemsep,label=$-$]
\item There is a finite locally free map $j: X_G(N,\Gamma_0(m)) \rar \mathbb{P}^1_R$ such that $X_G(N,\Gamma_0(m))$ is smooth over $\Sp{R}$ outside the inverse image under $j$ of the finitely many supersingular $j$-invariants in a fibre of residue characteristic dividing $m_0$. 
\item The reduced subscheme attached to the closed subspace $\{j=\infty\}$, called the \emph{cuspidal subscheme}, is finite \'etale over $R$ and defines an effective Cartier divisor. 
\item For any $R$-scheme $S$, $j^{-1}(\mathbb{A}^1_R)(S)$ is functorially isomorphic to the set of equivalence classes of $(E/S, \iota, C')$, where $E/S$ is a relative elliptic curve, $\iota: G_S \rar E[N]$ is an isomorphism of finite \'etale group schemes over $S$, and $C'$ is a cyclic subgroup scheme of $E$ of degree $m$ in the sense of \cite[(1.4), (3.4)]{KM}\footnote{This means that $C'$ is finite locally free of degree $m$ over $S$, and is fppf-locally over $S$ of the form $[P]+[2P]+\ldots+[mP]$ (as a Cartier divisor on $E$) for some $P \in E(S)$.}.   
\item The formation of $X_G(N,\Gamma_0(m)), j$, the cuspidal subscheme, and the above isomorphism of functors commutes with base change.  
\end{itemize}

Moreover, let $d,t \geq 1$ be integers such that $dt \mid m$, then the morphism of functors \[D_{d,t}: (E/S,\iota,C) \mapsto ((E/C[d])/S,\pi \circ \iota, C'),\] where $C[d]$ denotes the ``standard'' cyclic subgroup of $C$ of order $d$ in the sense of \cite[Theorem 6.7.2]{KM}, $\pi: E \rar E/C[d]$ denotes the natural isogeny and $C'$ denotes the subgroup of order $t$ of the cyclic subgroup $\pi(C)$ (see \cite[Theorems 6.7.2, 6.7.4]{KM}), extends to a finite locally free morphism of $R$-schemes $X_G(N,\Gamma_0(m)) \rar X_G(N,\Gamma_0(t))$ of constant rank mapping $X_G(N,\Gamma_0(m))$ to $X_G(N,\Gamma_0(t))$. 

For any prime $\ell \nmid N$, the correspondence $(D_{1,1},D_{\ell,1}): X_G(N,\Gamma_0(\ell)) \rar X_G(N) \times X_G(N)$ (where $X_G(N) := X_G(N,\Gamma_0(1))$) induces an endomorphism $T_{\ell}$ of the relative Jacobian $J_G(N)$ of $X_G(N)$, and the $T_{\ell}$ commute pairwise. 

Furthermore, for $n \in (\Z/N\Z)^{\times}$, there is an automorphism $[n]$ of $X_G(N,\Gamma_0(m))$ extending the isomorphism of functors $(E/S,\iota,C) \longmapsto (E/S,n\iota,C)$; $[n]$ commutes to any $D_{d,t}$. 

Assume finally that $G$ is endowed with a bilinear alternated pairing $\langle-,\,-\rangle_G: G \times G \rar \mu_N$ which is fibrewise perfect. Then there is a morphism $\det: X_G(N,\Gamma_0(m)) \rar (\Z/N\Z)^{\times}_R$ such that for any $(E/S,\iota,C) \in X_G(N,\Gamma_0(m))(S)$, $\det(F/S,\iota,C)$ is the unique $\alpha \in (\Z/N\Z)^{\times}_R(S)$ such that the pairings
\[G \times G \overset{\iota \times \iota}{\longrightarrow} F[N] \times F[N] \overset{\mrm{We}}{\rar} (\mu_N)_S, \quad G_S \times G_S \overset{\alpha \cdot \langle -,\,-\rangle_G}{\longrightarrow} (\mu_N)_S\] are equal, where $\mrm{We}$ denotes the Weil pairing.

For any $\alpha \in (\Z/N\Z)^{\times}$, $X_G^{\alpha}(N,\Gamma_0(m)) := \det^{-1}(\alpha)$ is a proper flat $R$-scheme with geometrically connected fibres. For any $\alpha \in (\Z/N\Z)^{\times}$ and $d,t \geq 1$ such that $dt \mid m$, $D_{d,t}$ maps $X_G^{\alpha}(N,\Gamma_0(m))$ to $X_G^{d\alpha}(N,\Gamma_0(t))$, so that $T_{\ell}$ maps $J_G^{\alpha}(N)$ into $J_G^{\ell\alpha}(N)$. For any $\alpha,n \in (\Z/N\Z)^{\times}$, $[n]$ is an isomorphism $X_G^{\alpha}(N,\Gamma_0(m)) \rar X_G^{\alpha n^2}(N,\Gamma_0(m))$.  
}

\rem{As noted in the introduction of this thesis, the notion of relative Jacobian for a proper smooth scheme $X \rar \Sp{R}$ of relative dimension one does not seem to be standard when the morphism does not have geometrically connected fibres, even when $R$ is a field. However, only minor modifications to the classical theory are needed to define and construct Jacobians in this setting, which we describe in Section \ref{degrees-jacobians}.}

\theoi[XG-same-as-twist]{Let us keep the notations of Theorem \ref{XG-exists-pol-unpol} and assume that $R$ is a field $k$ with separable closure $k_s$. 
Let $X(N,N)_k$ denote the curve attached by Theorem \ref{XG-exists-pol-unpol} to the constant group $(\Z/N\Z)^{\oplus 2}_k$, which is endowed with a left action of $\GL{\Z/N\Z}$ given for any noncuspidal point $(E/S,\iota) \in X(N,N)_k(S)$ by $M \cdot (E/S,\iota) = (E/S,\iota \circ M^T)$. 

Fix a basis $(P,Q)$ of $G(k_s)$, and let $\rho: \mrm{Gal}(k_s/k) \rar \GL{\Z/N\Z}$ be the function defined by \[\forall g \in \mrm{Gal}(k_s/k),\,\rho(g)\begin{pmatrix} P\\Q\end{pmatrix} = \begin{pmatrix}g(P)\\g(Q)\end{pmatrix}.\] Then $\rho$ is an anti-homomorphism, and let $X(N,N)_{\rho^{-1}}$ denote the Galois twist (see Proposition \ref{cocycle-twist}) of $X(N,N)_k$ by the homomorphism $\rho^{-1}: \mrm{Gal}(k_s/k) \rar \GL{\Z/N\Z}$. The curve $X(N,N)_{\rho^{-1}}$ is also endowed with a finite locally free morphism $j: X(N,N)_{\rho^{-1}} \rar \mathbb{P}^1_k$ representing the $j$-invariant, with Hecke correspondences and an action of $(\Z/N\Z)^{\times}$ coming from the twist by $\rho^{-1}$. 

Then there is a isomorphism of $\mathbb{P}^1_k$-schemes $\iota: X_G(N) \rar X(N,N)_{\rho^{-1}}$ such that, for any non-cuspidal $(E/k_s,\varphi) \in X_G(N)(k_s)$, one has \[\iota((E/k_s,\varphi)) = (E/k_s,(a,b) \mapsto aP+bQ) \in X(N,N)_k(k_s) \simeq X(N,N)_{\rho^{-1}}(k_s).\] After passing to the Jacobians, this isomorphism preserves the Hecke operators and the action of $[n]$. 

If $G$ is endowed with a perfect bilinear alternated pairing $\langle\cdot,\,\cdot\rangle_G: G \times G \rar \mu_N$, then $\det{\rho}$ is the mod $N$ cyclotomic character and the connected components of $X(N,N)_{\rho^{-1}}$ are smooth projective geometrically connected curves over $k$, naturally indexed by the primitive $N$-th roots of unity in $k_s$, and $\iota$ maps $X_G^{\alpha}(N)$ to the component of $X(N,N)_{\rho^{-1}}$ indexed by $\langle P,\,Q\rangle_G^{\alpha}$.    
}

As a preparation for the rest of the text, we also discuss some of the behavior of the curve $X(N)$ over $\Z[1/N]$: we prove an Eichler-Shimura relation (Corollary \ref{tate-eichler-shimura}) and describe in Section \ref{cuspidal-subscheme} its cuspidal subscheme and the $q$-expansion principle.  

This chapter is rooted in techniques from algebraic geometry, and in particular descent. Most of the relevant background, which we will use freely, is collected in Appendix \ref{alg-geom-prereq}. The main reference for this background remains the online textbook Stacks Project \cite{Stacks}. The main references for the modular aspects are \cite{DeRa} and \cite{KM}. For the proof of the Eichler-Shimura relation, we somewhat follow the argument in the proof of \cite[Appendix A, Theorem 5.16]{Ribet-Stein}, but (as noted in \emph{loc.cit.}), it is made significantly simpler by the theory of moduli schemes in characteristic $p$.

\nott{If $R$ is a ring and $S$ is a set, $S_R$ denotes the constant $R$-scheme with underlying set $S$. 
For any $N \geq 1$, 
\begin{itemize}[noitemsep,label=\tiny$\bullet$]
\item $\mu_N$ denotes the finite flat group scheme $\Sp{\Z[X]/(X^n-1)}$ over $\Z$, that is, the Cartier dual of the constant group scheme $(\Z/N\Z)_{\Z}$, 
\item $\mu_N^{\times}$ denotes the finite flat $\Z$-scheme $\Sp{\Z[X]/\Phi_N(X)}$, where $\Phi_N(X)$ is the $N$-th cyclotomic polynomial, which is \'etale above $\Z[1/N]$,
\item $\OO_N$ denotes the ring $\Z[1/N,\zeta_N]$, where $\zeta_N$ is a primitive $N$-th root of unity.
\end{itemize}
For every $d \mid N$, there is a natural closed immersion $\mu_d^{\times} \rar \mu_N$ (given on coordinate rings by $X \mapsto X$), and the induced map $\coprod_{d \mid N}{\mu_d^{\times}} \rar \mu_N$ is an isomorphism over $\Z[1/N]$.
}

\section{Moduli problems}
\subsection{Definitions}

We recall the framework introduced in \cite[Chapter 4]{KM}.

Given a ring $R$, we denote by $\Ell_R$ the category defined as follows. Its objects are relative elliptic curves $E \rar S$ (which we often denote as $E/S$), where $S$ is a $R$-scheme; given two relative elliptic curves $F \rar T$ and $E \rar S$, with $S,T$ being $R$-schemes, a morphism $F/T \rar E/S$ is a Cartesan diagram of $R$-schemes

\[
\begin{tikzcd}[ampersand replacement=\&]
F\arrow{r}{\alpha}\arrow{d} \& E\arrow{d}\\
T \arrow{r} \&  S
\end{tikzcd}
\]

where $\alpha$ preserves the zero section, so that $F \rar E \times_S T$ is an isomorphism \emph{of elliptic curves over $T$} (by \cite[Theorem 2.5.1]{KM}).

Note that if $R'$ is a $R$-algebra, we have a faithful forgetful functor $\Ell_{R'} \rar \Ell_R$.  

\defi{A \emph{moduli problem over $R$} is a contravariant functor $\P: \Ell_R \rar \mathbf{Set}$. 
If $R'$ is a $R$-algebra, the base change $\P_{R'}$ of $\P$ to $R$ is exactly the functor $\Ell_{R'} \rar \Ell_R \overset{\P}{\rar} \mathbf{Set}$. }

\defi{A moduli problem $\P$ over $R$ is \emph{rigid} if for any object $E/S$ of $(\mrm{Ell}/R)$, the action of $\mrm{Aut}(E/S)$ on $\P(E/S)$ (on the right) is free, that is: for any $\alpha \in \P(E/S)$, any $g \in \mrm{Aut}_R(E/S)$ with $g \neq \mrm{id}$, $\P(g)\alpha \neq \alpha$.}

\defi{A moduli problem $\P$ over $R$ is \emph{relatively representable} if for any object $E/S$ of $\Ell_R$, the functor $T \in \mathbf{Sch}_S \longmapsto \P(E_T/T) \in \mathbf{Set}$ is representable by a $S$-scheme $\P_{E/S}$.}

\defi{A moduli problem $\P$ over $R$ is \emph{representable} if there is an elliptic curve $\mathcal{E}/\MP$ (for some $R$-scheme $\MP$) and some $\alpha \in \P(\mathcal{E}/\MP)$ such that, for any object $E/S$ in $\Ell_R$, $f \in \mrm{Hom}_{\Ell_R}(E/S,\mathcal{E}/\MP) \longmapsto \P(f)\alpha \in \P(E/S)$ is an isomorphism. We call $\MP$ the moduli scheme, and we say that $(\mathcal{E}/\MP,\alpha)$ represents the moduli problem $\P$ (we will often omit $\alpha$ from the notation if it is not relevant).}

The following result is formal: 

\prop[base-change-representable]{Let $\P$ be a moduli problem over $R$ and let $R'$ be a $R$-algebra. Let us write $\P'=\P_{R'}$ for the base change of $\P$. 
\begin{itemize}[noitemsep,label=$-$]
\item If $\P$ is rigid, then so is $\P'$. 
\item If $\P$ is relatively representable, then any morphism $F/T \rar E/S$ in $\Ell_R$ induces a canonical isomorphism $\P_{F/T} \rar \P_{E/S} \times_S T$,
\item If $\P$ is relatively representable, then $\P'$ is relatively representable, and, for any object $E/S$ in $\Ell_{R'}$, we have an isomorphism $\P'_{E/S} \rar \P_{E/S}$ of $S$-schemes.
\item If $\P$ is representable by $(\mathcal{E}/\MP,\alpha)$, then $\P'$ is representable by \[(\mathcal{E}_{R'}/\MP_{R'},\P[\mathcal{E}_{R'}/\MP_{R'} \rar \mathcal{E}/\MP]\alpha).\] In particular, we have a canonical isomorphism $\MPp{\P'} \rar \MP_{R'}$ (we say that the formation of the moduli scheme commutes with base change).  
\end{itemize}}

Let $\mathscr{C}$ be a property of morphisms of schemes. We say that it is \emph{good} if it satisfies the following conditions:
\begin{itemize}[noitemsep,label=$-$]
\item it is stable under pre- and post-composition by isomorphisms, 
\item it is stable under arbitary base change, 
\item it satisfies fpqc descent in the following sense: if $f: X \rar Y$ is a morphism such that there exists a fpqc cover by maps $Y_i \rar Y$ such that $f_{Y_i}: X \times_Y Y_i \rar Y_i$ satisfies the property $\mathscr{C}$, then $f$ satisfies the property $\mathscr{C}$. 
\end{itemize}

All the properties given in Proposition \ref{fpqc-descent-prop} are good.

\defi{Let $\P$ be a relatively representable moduli problem over $R$ and $\mathscr{C}$ be a good property of morphisms of schemes. We say that $\P$ satisfies $\mathscr{C}$ over $\Ell_R$ if for every object $E/S$ in $\Ell_R$, then $\P_{E/S} \rar S$ has property $\mathscr{C}$. If $Z \subset \mathbb{A}^1_R$ is a locally closed subscheme of $\mathbb{A}^1_R$, we say that $\P$ satisfies $\mathscr{C}$ above $Z$ if, for every object $E/S$ in $\Ell_R$ such that $j(E/S) \in \mathbb{A}^1_R(S)$ factors through $Z$, then $\P_{E/S} \rar S$ satisfies property $\mathscr{C}$. 
}

\prop[adding-constants]{Let $R$ be a ring and $S_0$ be any $S$-scheme. The moduli problem $[S_0]$ on $\Ell_R$ defined by $E/S \longmapsto \mrm{Mor}_R(S,S_0)$ is relatively representable; for any relative elliptic curve $E/S$ ($S$ being a $R$-scheme), $[S_0]_{E/S}$ is naturally isomorphic to $S \times_R S_0$. In particular, if $\mathscr{C}$ is any good property of morphisms of schemes such that $S_0 \rar \Sp{R}$ satisfies $\mathscr{C}$, then $[S_0]$ satisfies $\mathscr{C}$.}

\demo{This is formal.}

\defi{Let $\P$ be a moduli problem over $R$, representable by the elliptic curve $\mathcal{E}/\MP$. There is a canonical $j$-invariant map $\MP \rar \mathbb{A}^1_R$, given by the $j$-invariant of $\mathcal{E}$.}

One checks formally that

\lem{The formation of the $j$-invariant map commutes with base change of representable moduli problems.}

\prop[base-change-good-property]{Let $\mathscr{C}$ be a good property of morphisms of schemes, $\P$ be a relatively representable moduli problem over $R$, and $\P'$ be the base change of $\P$ to an $R$-algebra $R'$. Then 
\begin{itemize}[noitemsep,label=$-$]
\item If $\P$ has property $\mathscr{C}$ over $\Ell_R$, then $\P'$ has property $\mathscr{C}$ over $\Ell_{R'}$,
\item If $R'$ is faithfully flat over $R$ and $\P'$ has property $\mathscr{C}$ over $\Ell_{R'}$, then $\P$ has property $\mathscr{C}$ over $\Ell_R$.
\end{itemize}}

We finally discuss the situation for morphisms of moduli problems. 

\prop{Let $T: \P \Rightarrow \P'$ be a natural transformation of relatively representable moduli problems over $\Ell_R$. Then, for every object $E/S$ in $\Ell_R$, $T$ induces a morphism of $S$-schemes $\P_{E/S} \rar \P'_{E/S}$ whose formation commutes with base change. 
If $\P,\P'$ are representable, then $T$ induces a map $\MP \rar \MPp{\P'}$ whose formation commutes with base change. }

We conclude this section with representability results.

\prop[product-rel-rep]{(\cite[(4.3.4)]{KM}) Let $\P, \P'$ be relatively representable moduli problems over $\Ell_R$. Then the moduli problem over $\Ell_R$
\[\P \times \P': E/S \longmapsto \P(E/S) \times \P'(E/S)\]
(with its natural action on the morphisms in $\Ell_R$) is relatively representable. Moreover, for every object $E/S$ of $\Ell_R$, one has a canonical isomorphism of $S$-schemes 
\[(\P \times \P')_{E/S} \rar \P_{E/S} \times_S \P'_{E/S}. \]

If furthermore $\P$ is representable by $(\mathcal{E}/\MP,\alpha)$, then $\P \times \P'$ is representable by \[\left((p')^{\ast}\mathcal{E}/\P'_{\mathcal{E}/\MP}, ((\P p')(\alpha),\beta)\right),\] where $p': \P'_{\mathcal{E}/\MP} \rar \MP$ is the map given by the relative representability of $\P'$, and 
\[\beta \in \P'((p')^{\ast}\mathcal{E}/\P'_{\mathcal{E}/\MP})\] corresponds to the identity under the identification \[\P'((p')^{\ast}\mathcal{E}/\P'_{\mathcal{E}/\MP}) \simeq \mrm{Hom}_{\MP}(\P'_{\mathcal{E}/\MP},\P'_{\mathcal{E}/\MP}).\]
}

\demo{This is \cite[(4.3.4)]{KM}, but since the proof (essentially an exercise in formalism) is not detailed in the book, we do it here. 
For the first claim, we simply note that if $E/S$ is an object of $\Ell_R$ and $T$ is any $S$-scheme, then we have, fonctorially in $T$, 
\begin{align*}
(\P \times \P')(E_T/T) &\simeq \P(E_T/T) \times \P'(E_T/T)\\
&\simeq \mrm{Hom}_S(T,\P_{E/S}) \times \mrm{Hom}_S(T,\P'_{E/S})\\
&\simeq \mrm{Hom}_S(T,\P_{E/S} \times_S \P'_{E/S}), 
\end{align*}
which concludes. 

Now let us consider the second claim. For the sake of easier notation, let $B=\MP$, $B'=\P'_{\mathcal{E}/\MP}$, and $\alpha'=(\P p')(\alpha)$. 
Given an object $E/S$ of $\Ell_R$, we have an isomorphism $\iota_{E/S}$ of functors on $\mathbf{Sch}_S$: $\iota_{E,S}: \P'(E_T/T) \rar \mrm{Mor}_S(T,\P'_{E/S})$. In particular, $\beta$ is the inverse image under $\iota_{\mathcal{E}/B}(B')$ of $\mrm{id}_{B'} \in \mrm{Mor}_B(B',B')$. 

Let $F/T$ be any object of $\Ell_R$: we need to prove that for any $a \in \P(F/T),\, b \in \P'(F/T)$, there is a unique morphism $f \in \Ell_R(F/T,(p')^{\ast}\mathcal{E}/B')$ such that $(\P f)\alpha'=a$, $(\P' f)(\beta)=b$. 

Let us pick $a \in \P(F/T)$ and $b \in \P'(F/T)$. There is a morphism $f_1 \in \Ell_R(F/T,\mathcal{E}/B)$ such that $a=(\P f_1)(\alpha)$. Let $g_1$ be the induced morphism $T \rar B$, so that $f_1$ determines an isomorphism $f_2: F/T \simeq g_1^{\ast}\mathcal{E}/T$. We thus can define $b'=(\P'(f_2^{-1}))(b) \in \P'(\mathcal{E} \times_{B} T/T)$. Let $(\iota_{\mathcal{E}/B}T)(b')=\tilde{f_3} \in \mrm{Mor}_{B}(T,B')$, and $f_3$ be the induced morphism $(\tilde{f_3})^{\ast}(p')^{\ast}\mathcal{E}/T \rar (p')^{\ast}\mathcal{E}/B'$ (noting that $\tilde{f_3}^{\ast}(p')^{\ast}=(p' \circ \tilde{f_3})^{\ast}=g_1^{\ast}$). 

By Yoneda's lemma, one has $\tilde{f_3}=(\iota_{\mathcal{E}/B}\tilde{f_3})((\iota_{\mathcal{E}/B} B')\beta)=(\iota_{\mathcal{E}/B}T)(\P'f_3)(\beta)$, and therefore $b'=(\P'f_3)\beta$, whence $b=(\P'(f_3 \circ f_2))(\beta)$. Moreover, one has \[\P(f_3 \circ f_2)(\alpha')=\P(p' \circ f_3 \circ f_2)(\alpha)=\P(f_1)(\alpha)=a.\] 

Now, we need to show injectivity. Let $\tilde{p}': (p')^{\ast}\mathcal{E}/B' \rar \mathcal{E}/B$ be the natural morphism. Let $f_1, f_2 \in \Ell_R(F/T, (p')^{\ast}\mathcal{E}/B')$ be such that $(\P f_i)(\alpha')=a$, $(\P'f_i)(\beta)=b$. Thus one has $a=(\P f_i)(\alpha')=(\P f_i) \circ (\P\tilde{p}')(\alpha)$, so that the two maps $\tilde{p}' \circ f_1=\tilde{p}' \circ f_2$. In particular, the two base maps $\tilde{f_i}: T \rar B'$ attached to the $f_i$ define the same structure of $B$-scheme on $T$, and the two isomorphisms $g_i \in \Ell_R(F/T, \tilde{f_i}^{\ast}(p')^{\ast}\mathcal{E}/T)$ induced by $f_i$ are the same, since $\tilde{f_i}^{\ast}(p')^{\ast}=q^{\ast}$.

Let $f'_i$ denote the morphism $q^{\ast}\mathcal{E}/T \rar (p')^{\ast}\mathcal{E}$ induced by $\tilde{f_i}$. We have $f_i=f'_i \circ g_i$, and since $g_1=g_2$ are isomorphisms, one knows that $(\P'f'_i)(\beta)$ is independent from $i$. Now, $(\P'f'_i)(\beta)=(\P'f'_i)((\iota_{\mathcal{E}/B}B')^{-1}\mrm{id}_{B'})=(\iota_{\mathcal{E}/B}T)^{-1}(\tilde{f_i})$. Therefore, $\tilde{f_1}=\tilde{f_2}$, hence $f'_1=f'_2$ and $f_1=f_2$. 

}

\prop[over-Ell-implies-over-R]{Let $\mathscr{C}$ be a good property of morphisms of schemes satisfying in addition the following assumptions: 
\begin{itemize}[noitemsep,label=$-$]
\item it is stable under post-composition by surjective smooth affine maps of relative dimension one (resp. by finite flat maps), 
\item if $f: X \rar Y$ is a finite \'etale surjective map and $g: Y \rar Z$ is such that $g \circ f$ satisfies $\mathscr{C}$, then $g$ satisfies $\mathscr{C}$.
\end{itemize}
Let $\P$ be a representable moduli problem over $\Ell_R$ having property $\mathscr{C}$ over $\Ell_R$ (resp. above some locally closed subscheme $Z$ of $\mathbb{A}^1_R$). Then $\MP \rar \Sp{R}$ has property $\mathscr{C}$ (resp. the $j$-invariant of $\MP \rar \mathbb{A}^1_R$ has property $\mathscr{C}$ above $Z$).
}

\demo{Since everything is Zariski-local with respect to $R$, we can assume that there is an odd prime $\ell$ invertible in $R$. The conclusion then follows from the properties of the following commuting diagram:
\[
\begin{tikzcd}[ampersand replacement=\&]
\MPp{\P,[\Gamma(\ell))]_R} \arrow{r}{f_1}\arrow{d}{\pi_1} \& \MPp{[\Gamma(\ell)]_R}\arrow{d}{j_{\ell}}\arrow{r}{\sigma} \& \Sp{R}\\
\MP \arrow{r}{j_{\P}} \&  \mathbb{A}^1_R \arrow{ur} \&
\end{tikzcd}
\]

Indeed, $\pi_1$ is finite \'etale by \cite[Theorem 3.7.1]{KM} while $f_1$ satisfies $\mathscr{C}$ (resp. satisfies $\mathscr{C}$ above $Z$), and $\sigma$ is smooth affine of relative dimension one by \cite[Corollary 4.7.2]{KM}, while $j_{\ell}$ is finite by \cite[Theorem 5.1.1, Proposition 8.2.2]{KM} and thus flat (over $\Z[1/\ell]$, this follows from miracle flatness \cite[Lemma 00R4]{Stacks}, since both schemes are smooth surjective $\Z[1/\ell]$-schemes where every nonempty open subset has dimension $2$, and we conclude by base change). 
}

\cor[over-Ell-concrete]{Let $\P$ be a representable moduli problem over $\Ell_R$. 
\begin{itemize}[noitemsep,label=$-$]
\item If $\P$ is affine (resp. quasi-compact) over $\Ell_R$, then $\MP$ is affine (resp. quasi-compact).
\item If $\P$ is surjective (resp. locally quasi-finite, locally of finite type, locally of finite presentation, finite with $R$ Noetherian, flat) over $\Ell_R$, its $j$-invariant $\MP \rar \mathbb{A}^1_R$ is surjective (resp. locally quasi-finite, locally of finite type, locally of finite presentation, finite, flat). In particular, if $\P$ is locally quasi-finite of finite presentation and flat over $\Ell_R$, then $\MP \rar \Sp{R}$ is Cohen-Macaulay \cite[Section 00N7]{Stacks} and its nonempty fibres are pure of relative dimension one.
\item If $\P$ is affine of finite presentation over $\Ell_R$ and \'etale above some open subscheme $U \subset \mathbb{A}^1_R$, then $j^{-1}(U) \subset \MP \rar \Sp{R}$ is smooth surjective of relative dimension one.
\end{itemize} 
}

\demo{The first point follows from the properties of affine and quasi-compact morphisms by fpqc descent and \cite[Lemma 01ZT]{Stacks}. The second point follows from Proposition \ref{fpqc-descent-prop} for fpqc descent, \cite[Section 036M, Section 036K]{Stacks} for the ``finite \'etale descent'' conditions, and the fact that finite flat morphisms are clearly Cohen-Macaulay. The dimension of the fibres follows from the results of \cite[Section 02NM]{Stacks}. For the last statement, we use the Proposition in the same way along with \cite[Lemma 036U]{Stacks} to check that the ``finite \'etale descent'' condition is satisfied. }

\prop[representability-descends]{Let $R \rar R'$ be a faithfully flat \'etale map, and let $\P$ be a moduli problem on $\Ell_R$. Assume that:
\begin{itemize}[noitemsep,label=$-$]
\item $\P_{R'}$ is representable and affine,
\item For all objects $E/S$ of $\Ell_R$, $T \in \mathbf{Sch}_S \longmapsto \P(E_T/T)$ is an fpqc sheaf. 
\end{itemize}
Then $\P$ is representable.
}

\demo{We need to show that $\P$ is relatively representable, affine over $\Ell_R$, and rigid, by \cite[(4.7.0)]{KM}\footnote{The statement in the book discusses the case $R=\Z$, but it applies to any case, since the moduli problems used to cover $\Ell_{\Z}$, the ``Legendre'' and ``na\"ive level three'' problems (over $\Z[1/2]$ and $\Z[1/3]$ are still representable by smooth affine curves over $R[1/2]$ and $R[1/3]$ for any ring $R$, by Proposition \ref{base-change-representable}), and the argument does not use any property of $\Z$ as a base ring.}. 

Let us first show that $\P$ is rigid. Let $T$ be a $R$-algebra, $E/T$ be an elliptic curve, $\alpha \in \P(E/T)$ and $u \in \mrm{Aut}(E/T)$ such that $(\P u)(\alpha)=\alpha$. Let $T'=T \otimes_R S$ and $\alpha',u'$ be the base changes of $\alpha,u$ to the elliptic curve $E_{T'}/T'$. Since $\P_S$ is representable, and since $\alpha' \in \P_S(E_{T'}/T')$, and $(\P_Su')(\alpha')=\alpha'$, one has $u'=\mrm{id}$. Thus the automorphisms $\mrm{id}$ and $u$ of the $T$-scheme $E$ become equal after a faithfully flat base change with respect to $T$: so $u=\mrm{id}$. 

Now, we show that $\P$ is relatively representable. Let $E/T$ be an elliptic curve where $T$ is a $R$-scheme; we want to show that the pre-sheaf on the category of $T$-schemes $U \longmapsto \P(E_U/U)$ is representable. The hypothesis states that this pre-sheaf is an fpqc sheaf, and its restriction to the subcategory of schemes over $T_{R'} := T \times_R \Sp{R'}$ ($T_{R'} \rar T$ being affine faithfully flat) is representable by a scheme $X'$ affine over $T_{R'}$. By Proposition \ref{desc-affine-schemes}, this sheaf is representable by an affine $T$-scheme $X$.}

\subsection{$\Gamma_0$-structures and degeneracy maps}
\label{gamma0m-structure}

Recall the following definition inspired from \cite[\S 4.14]{KM}. 

\defi{Let $N \geq 1$ be an integer and $R$ be a ring. Let $\mathbf{FGp}^N_R$ be the category defined as follows:
\begin{itemize}[noitemsep,label=$-$]
\item Its objects are finite flat locally free commutative group schemes $G/S$, where $S$ is a $R$-scheme, and $G$ is annihilated by the multiplication by $N$,
\item A morphism $G'/S' \rar G/S$ is a Cartesian diagram as below such that the induced map $f: G' \rar G \times_S S'$ is an isomorphism of $S'$-group schemes. \[\begin{tikzcd}[ampersand replacement=\&] G' \arrow{r}{f} \arrow{d} \& G\arrow{d} \\S' \arrow{r}\& S\end{tikzcd}\]
\end{itemize}

A \emph{level $N$ structure} on a moduli problem $\P$ on $\Ell_R$ is a functor $F: \mathbf{FGp}^N_R \rar \mathbf{Set}$ such that $\P$ is the composition of $F$ and $E/S \in \Ell_R \longmapsto E[N]/S \in \mathbf{FGp}^N_R$. Unless it is both relevant and not obvious, when we claim that a moduli problem $\P$ on $\Ell_R$ is of level $N$, we omit the functor $F$ from the notation. }

\defi{Let $\P,\P'$ be moduli problems on $\Ell_R$, endowed with level $N$ structures by the functors $F,F'$. A morphism of moduli problems of level $N$ is a natural transformation $f: F \Rightarrow F'$. In particular, for every object $E/S$ of $\Ell_R$ and every $\alpha \in \mathscr{P}(E/S)=F(E[N]/S)$, we can define $f(\alpha) \in F'(E[N]/S)=\P'(E/S)$; this construction defines in particular a morphism of functors $\P \Rightarrow \P'$.}

\rems{\begin{itemize}[noitemsep,label=$-$]
\item If $S_0$ is a $R$-scheme, the moduli problem $[S_0]$ on $\Ell_R$ has a natural level one structure. 
\item Suppose that $(\P,F)$ is a moduli problem of level $N$, and $f: E \rar E'$ is an isogeny of elliptic curves over a given $R$-scheme $S$ of degree prime to $N$. Then $f$ induces an isomorphism of $S$-group schemes $E[N] \rar E'[N]$, thus induces an isomorphism $\P(E/S) \rar \P(E'/S)$. When furthermore $E=E'$ and $f$ is the multiplication by some integer $d \equiv 1 \mod{N}$, $f$ acts trivially on $E[N]$: this means that $(\Z/N\Z)^{\times}$ acts\footnote{Note that this action depends on the choice of $F$.} on $\P$; we will denote this action by $[\cdot]$. Note that this action preserves the $j$-invariant.  
\item We can define notions of base change similarly to the previous section. 
\item If $N \geq 1$ is an integer and $\P$ is a moduli problem on $\Ell_R$ endowed with the level $N$ structure $F$, and $M$ is any multiple of $N$, then we can endow $\P$ with the natural level $M$ structure $F_{N \rar M}: \mathbf{FGp}^M_R \overset{G \mapsto G[N]}{\longrightarrow} \mathbf{FGp}^N_R \overset{F}{\rar} \mathbf{Set}$. 
\item One obviously has $(F_{N \rar M})_{M \rar M'}=F_{N \rar M'}$ if $M \mid M'$.
\item If $(\P,F)$ is a moduli problem of level $N$, $\hat{\Z}^{\times}$ acts on $\P$ by projecting to $(\Z/N\Z)^{\times}$ and applying $[\cdot]$. If $N \mid M$, then the induced action of $\hat{\Z}^{\times}$ on $\P$ given by $(\P,F_{N \rar M})$ is the same action. 
\item We thus define, for any integer $m \geq 1$, the category of moduli problems of level coprime to $m$: its objects are moduli problems endowed with a level $n$ structure for some $n \geq 1$ coprime to $m$, and a morphism $(\P,F) \rar (\P',F')$, when $F,F'$ are respectively level $n,n'$ structures on $\P,\P'$ (with $nn'$ coprime to $m$) is any natural transformation of functors $F_{n \rar d} \Rightarrow F_{n' \rar d}$ for some $d$ divisible by $n,n'$ and coprime to $m$. The composition is defined similarly, if necessary by raising $d$. The induced morphism $\P \Rightarrow \P'$ is equivariant for the actions of $\hat{\Z}^{\times}$ (and any $\Z_p^{\times}$ for $p \mid m$ acts trivially). 
\end{itemize}
}

\prop[product-basic]{Let $\P,\P'$ be moduli problems on $\Ell_R$ of levels $N,N'$. Then the moduli problem $\P \times \P'$ is endowed with a natural level $M$ structure, where $M$ is the least common multiple of $N$ and $N'$.}

\demo{This is formal.}

\defi[basic-moduli-definition]{Let $m \geq 1$ be an integer. Let us recall the definitions of the basic moduli problems $[\Gamma(m)]_R, [\Gamma_1(m)]_R,[\Gamma_0(m)]_R$ over a ring $R$, following \cite[Chapter 3]{KM}. In the case of $[\Gamma_1(m)]_R$ and $[\Gamma(m)]_R$, we will always assume that $m \in R^{\times}$. 

The moduli problem $[\Gamma(m)]_R$ (resp. $[\Gamma_1(m)]_R$) is defined as follows: 
\begin{itemize}[noitemsep,label=$-$]
\item Let $E$ be an elliptic curve over a $R$-scheme $S$. Then $[\Gamma(m)]_R(E/S)$ (resp. $[\Gamma_1(m)]_R(E/S)$) is the set of group homomorphisms $\phi: (\Z/m\Z)^{\oplus 2} \rar E[m]$ (resp. $\phi: \Z/m\Z \rar E[m](S)$) such that for any geometric point $\overline{s}$ of $S$, the induced map $\phi: (\Z/m\Z)^{\oplus 2} \rar E[m](\overline{s})$ (resp. $\phi: \Z/m\Z \rar E[m](\overline{s})$) is an isomorphism (resp. is injective). Equivalently, $[\Gamma(m)]_R(E/S)$ (resp. $[\Gamma_1(m)]_R(E/S)$) denotes the collection of couples $(P,Q) \in E[m](S)$ (resp. of $P \in E[m](S)$) such that for any geometric point $\overline{s}$ of $S$, the image of $(P,Q)$ in $E[m](\overline{s})$ is a $\Z/m\Z$-basis of this group (resp. the image of $P$ in $E[m](\overline{s})$ has order exactly $m$).  
\item Given a morphism $F/T \rar E/S$ in $\Ell_R$, the morphism $[\Gamma(m)]_R(E/S) \rar [\Gamma(m)]_R(F/T)$ (resp. $[\Gamma_1(m)]_R(E/S) \rar [\Gamma_1(m)]_R(F/T)$) is induced by the map $E(S) \rar F(T)$.  
\end{itemize}

The moduli problem $[\Gamma_0(m)]_R$ is defined as follows (as in \cite[(3.4)]{KM}):
\begin{itemize}[noitemsep,label=$-$]
\item Given an elliptic curve $E$ over some $R$-scheme $S$, $[\Gamma_0(m)]_R(E/S)$ denotes the collection of $S$ subgroup schemes $C$ of $E$ that are finite locally free of degree $m$ over $S$, and admit fppf-locally over $S$ a generator -- in the sense that fppf-locally on $S$, there is a section $P \in C(S)$ such that as Cartier divisors over $E$, $C$ is equal to $\sum_{i=1}^m{[iP]}$. 
\item Given a morphism $F/T \rar E/S$ in $\Ell_R$, the morphism $[\Gamma_0(m)]_R(E/S)\rar [\Gamma_0(m)]_R(F/T)$ is given by $C \longmapsto C_T$. 
\end{itemize}
}

\rem{\cite[(3.4)]{KM} also gives another definition for the moduli problem $[\Gamma_0(m)]_R$, as the collection of isomorphism classes of isogenies $E \rar E'$ whose kernel is finite locally free of degree $m$ over $S$ and fppf-locally cyclic. The two definitions are equivalent because, for any relative elliptic curve $E/S$ and any subgroup scheme $G$ of $E$ which is finite locally free over $S$, $G$ is the kernel of a finite locally free isogeny $\pi: E \rar E'$ by \cite[Lemma 07S7]{Stacks}\footnote{I am grateful to MathOverflow user SashaP for pointing out this result to me.}.  
}

\prop[basic-moduli]{On any ring $R$, and for any integer $m \geq 1$, the moduli problems $[\Gamma_0(m)]_R$, $[\Gamma_1(m)]_R$, $[\Gamma(m)]_R$ are relatively representable, finite locally free of constant rank on $\Ell_R$, \'etale above $\Z[1/m]$ (over $\Ell_R$). They are endowed with natural level $m$ structures. Moreover, the induced action of $\hat{\Z}^{\times}$ on $[\Gamma_0(m)]_R$ is trivial.}

\demo{The relative representability and the properties of the moduli problems follow from \cite[Lemma 4.12.1, Theorem 5.1.1]{KM} and Proposition \ref{base-change-good-property}). The definition of the level $m$ structure follows from the discussion in \cite[Chapter 3]{KM} (using for instance \cite[Proposition 1.10.6]{KM}). 
As for the action of $\hat{\Z}^{\times}$ on $[\Gamma_0(m)]_R$, note that, for any $d \geq 1$ coprime to $m$, if $E/S$ is an elliptic curve and $C \leq E[m]$ is a cyclic subgroup of order $m$ (in the sense of \cite[(3.4)]{KM}, that is: a closed subgroup scheme of $E[m]$, finite locally free of rank $m$ over $S$, admitting a generator fppf-locally on $S$), then $dC=C$.
}

\prop[rep-plus-gamma0]{Let $R$ be a ring and $\P$ be a representable moduli problem on $\Ell_R$ of level $N$. Let $m$ be an integer coprime to $N$. Then the moduli problem $\P \times [\Gamma_0(m)]_R$ is representable by a scheme $\MPp{\P,[\Gamma_0(m)]_R}$ endowed with a finite flat map of locally constant rank to $\MP$ (which is \'etale on the invertible locus of $m$). This moduli problem has a natural level $Nm$ structure. }

\demo{$[\Gamma_0(m)]_R$ is relatively representable by Proposition \ref{basic-moduli}, so $\P \times [\Gamma_0(m)]_R$ is representable by Proposition \ref{product-rel-rep}. The level structure is given by Proposition \ref{product-basic}.}

\prop[degeneracies]{Let $m$ be a positive integer and $d,t$ be integers such that $dt \mid m$. For any representable moduli problem $\P$ on $\Ell_R$ of level $N$ coprime to $m$, there exists a \emph{degeneracy map} $D_{d,t}: \MPp{\P,[\Gamma_0(m)]_R} \rar \MPp{\P,[\Gamma_0(t)]_R}$ defined by the following rule. Let $E/S$ be an elliptic curve endowed with a cyclic subgroup $C$ of order\footnote{As in \cite{KM}, this denotes a finite $S$-subgroup scheme of $E$, locally free of degree $m$ over $S$, and fppf locally cyclic.} $m$, and $\alpha \in \P(E/S)$. Let $C_d$ be the ``standard'' cyclic subgroup scheme of $C$ of order $d$ in the sense of \cite[Theorem 6.7.2]{KM}, and let $\pi: E \rar E'=E/C_d$ be the cyclic isogeny of degree $d$. Then $D_{d,t}$ maps $(E/S,\alpha,C)$ to $(E'/S,\pi(\alpha),C')$, where $C' \leq E'[t]$ is the finite flat subgroup scheme over $S$, cyclic of degree $t$ defined as the standard subgroup of the finite flat cyclic subgroup $\pi(C)$ of $E'$ (which has order $m/d$ by \cite[Theorem 6.7.4]{KM}).
}

\defi{The $D_{d,t}$ defined in Proposition \ref{degeneracies} are the \emph{degeneracy maps}.}

\cor[elem-degeneracies]{We keep the setting of Proposition \ref{degeneracies}.
\begin{itemize}[noitemsep,label=\tiny$\bullet$]
\item The formation of any $D_{d,t}$ is natural with respect to $\P$ (for morphisms of representable moduli problems of level $N$ coprime to $m$, such as the action of $(\Z/N\Z)^{\times}$) and commutes with base change. 
\item $D_{1,m}$ is the identity, and every $D_{1,q}$ (for $q \mid m$) is a map of $\MP$-schemes (where the structure map comes from Proposition \ref{rep-plus-gamma0}).
\item Suppose that $d,t$ are divisors of $m$ such that $dt \mid m$, and let $r,s$ be divisors of $t$ such that $rs \mid t$. Then $D_{r,s} \circ D_{d,t} = D_{dr,s}$.  
\end{itemize}
}

\demo{The first two statements are completely formal. The last assertion follows from \cite[Theorem 6.7.4]{KM}. }

\prop[degn-cartesian]{Let $F: \P \Rightarrow \P'$ be a morphism of representable moduli problems of level $N$ on $\Ell_R$. Let $m,d,t\geq 1$ be such that $dt \mid m$ and $m$ is coprime to $N$. Then the following diagram is Cartesian:
\[
\begin{tikzcd}[ampersand replacement=\&]
\MPp{\P,[\Gamma_0(m)]_R} \arrow{r}{D_{d,t}}\arrow{d}{F} \& \MPp{\P,[\Gamma_0(t)]_R}\arrow{d}{F}\\
\MPp{\P',[\Gamma_0(m)]_R} \arrow{r}{D_{d,t}} \&  \MPp{\P',[\Gamma_0(t)]_R}
\end{tikzcd}
\]
}

\demo{Let $\tau_1=(E'/S,\alpha',C') \in \MPp{\P,[\Gamma_0(t)]_R}(S), \tau_2=(E/S,\beta,C) \in \MPp{\P',[\Gamma_0(m)]_R}(S)$ be such that $F(\tau_1)=D_{d,t}(\tau_2) \in \MPp{\P',[\Gamma_0(t)]_R}(S)$. 
Then there exists a cyclic isogeny $\pi: E \rar E'$ of degree $d$ whose kernel is the standard subgroup $C_d$ of degree $d$ of $C$, such that $C'$ is the standard cyclic subgroup of degree $t$ of $\pi(C)$ and $F(\alpha')=\pi(\beta)$. 

Let $\alpha \in \P(E/S)$ be the inverse image of $\alpha'$ under $\pi$, then $(E/S,\alpha,C) \in \MPp{\P,[\Gamma_0(m)]_R}(S)$ lies above $\tau_1$ and $\tau_2$. 

Suppose that $(K/S,\gamma,D) \in \MPp{P,[\Gamma_0(m)]_R}(S)$ is a point lying over $(\tau_1,\tau_2)$. Then there is an isomorphism of elliptic curves $f: K \rar E$ mapping $F(\gamma)$ to $\beta$ and $D$ to $C$. Moroever, there is a cyclic isogeny $g: K \rar E'$ of degree $d$ such that $\alpha'=g(\gamma)$, $\ker{g}$ is the standard subgroup of $D$ of order $d$, and $C'$ is the standard cyclic subgroup of $g(D)$ of order $t$. 

Then $g, \pi \circ f$ are cyclic isogenies $K \rar E'$ with the same kernel. Thus there exists an automorphism $\psi$ of $E'$ such that $\psi \circ g= \pi \circ f$, whence \[\psi(F(\alpha'))=\psi(F(g(\gamma)))=F(\pi(f(\gamma)))=\pi(f(F(\gamma)))=\pi(\beta)=F(\alpha').\]

Since $\P'$ is representable, it is rigid, hence $\psi$ is the identity. Therefore, $\gamma=f^{-1}(\alpha)$, and $(F/S,\gamma,D)$ is isomorphic to $(E/S,\alpha,C)$, so we are done.
}

\lem[rel-coprime-elliptic-curves]{Let $E,F$ be relative elliptic curves over a basis $S$, endowed with cyclic subgroups $C,D$ of coprime degrees $c,d \geq 1$. Suppose we are given an isomorphism $\iota: E/C \rar F/D$ of elliptic curves (thus endowing $F$ with the structure of finite locally free $E/C$-scheme of degree $d$). Then $E \times_{E/C} F$ is an elliptic curve over $S$.}

\demo{Since $E, F, E/C$ are proper commutative group schemes over $S$, so is $G=E \times_{E/C} F$. Moreover, $G \rar S$ factors through $G \rar F \rar S$, with $G \rar F$ finite locally free and $F \rar S$ smooth of relative dimension one, hence $G$ is a proper flat commutative group scheme over $S$ of relative dimension one.  

Moreover, $G_{S[1/c]} \rar F_{S[1/c]}$ is the base change of the \'etale isogeny $E_{S[1/c]} \rar (E/C)_{S[1/c]}$, so it is \'etale and $G_{S[1/c]} \rar S$ is smooth of relative dimension one. Similarly, $G \rar S$ is smooth above $S[1/d]$ of relative dimension one. Since $c,d$ are coprime, $G \rar S$ is a proper smooth commutative group scheme of relative dimension one. 

We finally need to check that the geometric fibres of $G \rar S$ are connected, so we may assume that $S$ is the spectrum of an algebraically closed field. Let $p: E \rar E/C, q: F \rar F/D$ be the projections and $G^0$ be the connected component of $G$ containing the unit: $G^0$ is an elliptic curve over $S$. 

Suppose that $G^0$ contains the kernel $\{0_E\} \times_{E/C} D$ of the first projection $G \rar E$: then the isogeny of elliptic curves $G_0 \subset G \rar E$ is finite locally free of the same degree $d$ as $G \rar E$. Since $G^0$ is a closed open subscheme of $G$, this implies that $G^0=G$. The map $(p^{\vee}\iota^{-1}q,c\cdot\mrm{id}): F \rar E \times_{E/C} F$ factors through $G^0$ and its image contains $\{0_E\} \times_{E/C} D$ (which is the image of $D$), whence the conclusion. }

\prop[two-degn-cartesian]{Let $\P$ be a representable moduli problem of level $N$ on $\Ell_R$. Let $m \geq 1$ be an integer coprime to $N$. Let $s',t' \geq 1$ be two coprime divisors of $m$, and $s,t \geq 1$ be divisors of $s',t'$ respectively. Then the following diagram is Cartesian:
\[
\begin{tikzcd}[ampersand replacement=\&]
\MPp{\P,[\Gamma_0(m)]_R} \arrow{r}{D_{s,\frac{m}{s'}}}\arrow{d}{D_{t,\frac{m}{t'}}} \& \MPp{\P,[\Gamma_0\left(\frac{m}{s'}\right)]_R}\arrow{d}{D_{t,\frac{m}{s't'}}}\\
\MPp{\P,[\Gamma_0\left(\frac{m}{t'}\right)]_R} \arrow{r}{D_{s,\frac{m}{s't'}}} \&  \MPp{\P,[\Gamma_0\left(\frac{m}{t's'}\right)]_R}
\end{tikzcd}
\]
}

\demo{It is straightforward to check that the diagram commutes using \cite[Lemma 3.5.1, Theorem 6.7.4]{KM}. Let $s_{\infty},t_{\infty}$ be the greatest divisors of $m$ which have the same prime divisors as $s',t'$ respectively. Then $s_{\infty},t_{\infty},\frac{m}{s_{\infty}t_{\infty}}$ are pairwise coprime positive integers. 

If $C$ is a cyclic subgroup of some elliptic curve $E/S$, let $C[d]$ denote its standard cyclic subgroup of degree $d$. 

\emph{Step 1: Surjectivity}

Let \[x_1=(E_1/S,\alpha_1,C_1) \in \MPp{\P,[\Gamma_0(\frac{m}{t'})]_R}(S), x_2=(E_2/S,\alpha_2,C_2) \in \MPp{\P,[\Gamma_0(\frac{m}{s'})]_R}(S)\] be such that $D_{s,\frac{m}{s't'}}(x_1)=D_{t,\frac{m}{s't'}}(x_2)$. 

Thus, there exists an isomorphism $\iota: E_1/C_1[s] \rar E_2/C_2[t]$ identifying the standard cyclic subgroups of order $\frac{m}{s't'}$ of $C_1/C_1[s]$ and $C_2/C_2[t]$ such that \[(\P\iota)(\alpha_1 \mod{C_1[s]})=\alpha_2 \mod{C_2[t]}.\] 

Since $s,t$ are pairwise coprime, $E = E_1 \times_{E_1/C_1[s]} E_2$ is a relative elliptic curve over $S$ by Lemma \ref{rel-coprime-elliptic-curves}. Let $\pi_i: E \rar E_i$ be the natural projection: $\pi_1,\pi_2$ are isogenies of degrees $t,s$ respectively. Let $\alpha \in \P(E/S)$ be such that $(\P\pi_1)(\alpha)=\alpha_1$: then 
\begin{align*}
(\P\pi_2)(\alpha) \pmod {C_2[t]} &= (\P\iota)((\P\pi_1)(\alpha) \mod{C_1[s]})=(\P \iota)(\alpha_1 \mod{C_1[s]})\\
&=\alpha_2 \mod{C_2[t]},
\end{align*} whence $(\P\pi_2)(\alpha)=\alpha_2$.

Since $t$ is coprime to $\frac{m}{t_{\infty}}$ (resp. $s$ is coprime to $\frac{m}{s_{\infty}}$), there exists a unique cyclic subgroup $C'_1$ of order $\frac{m}{t_{\infty}}$ (resp. $C'_2$ of order $\frac{m}{s_{\infty}}$) of $E$ such that $\pi_1(C'_1)=C_1[m/t_{\infty}]$ (resp. $\pi_2(C'_2)=C_2[m/s_{\infty}]$).  

Since $s \mid \frac{m}{t_{\infty}}$, $t \mid \frac{m}{s_{\infty}}$ and $st$ is coprime to $\frac{m}{s_{\infty}t_{\infty}}$, one has 
\begin{align*}
\iota(\pi_1\left(C'_1\left[\frac{m}{s_{\infty}t_{\infty}}\right]\right) \mod{C_1[s]}) &= \iota(C_1\left[\frac{m}{s_{\infty}t_{\infty}}\right] \mod{C_1[s]})\\
&=\iota(C_1 \mod{C_1[s]})\left[\frac{m}{s_{\infty}t_{\infty}}\right]=(C_2 \mod{C_2[t]})\left[\frac{m}{s_{\infty}t_{\infty}}\right]\\
&= C_2\left[\frac{m}{s_{\infty}t_{\infty}}\right]\mod{C_2[t]}=\pi_2\left(C'_2\left[\frac{m}{s_{\infty}t_{\infty}}\right]\right), 
\end{align*}

so that $C'_1$ and $C'_2$ have the same $\frac{m}{s_{\infty}t_{\infty}}$-torsion and there exists by \cite[Lemma 3.5.1]{KM} a cyclic subgroup $C$ of $E$ with degree $m$ such that $C\left[\frac{m}{t_{\infty}}\right]=C'_1, C\left[\frac{m}{s_{\infty}}\right]=C'_2$. 

We claim that \[D_{s,\frac{m}{s'}}((E,\alpha,C))=x_2,\,D_{t,\frac{m}{t'}}((E,\alpha,C))=x_1.\] It is enough to prove that \[\ker{\pi_1}=C[t],\,\ker{\pi_2}=C[s],\,\pi_1(C[t_{\infty}])[\frac{t_{\infty}}{t'}]=C_1[\frac{t_{\infty}}{t'}],\,\pi_2(C[s_{\infty}])[\frac{s_{\infty}}{s'}]=C_2[\frac{s_{\infty}}{s'}].\]

Let $p_1: E_1 \rar E_1/C_1[s],\, p_2: E_2 \rar E_2/C_2[t]$ be the projections, then \[p_2\circ \pi_2(C[t])=p_2(\pi_2(C[t_{\infty}])[t])=p_2(C_2[t_{\infty}][t])=p_2(C_2[t])=0.\]

Similarly, the kernel of $p_2\circ \pi_2=\iota\circ p_1 \circ \pi_1$ contains $C[s]$. Thus $p_2 \circ \pi_2$ is an isogeny of degree $st$ whose kernel contains the cyclic subgroup $C[s]+C[t]=C[st]$ of degree $st$: hence $C[st] \rar \ker{p_2 \circ \pi_2}$ is a closed immersion of finite locally free $S$-schemes of same rank. Zariski-locally on $S$, it is a surjective morphism of two finite free $\OO(S)$-modules of same rank, hence an isomorphism, whence $C[st] \rar \ker{p_2 \circ \pi_2}$ is an isomorphism. 

Since $p_2,\pi_2$ are composable cyclic isogenies of coprime degrees, by \cite[Proposition 6.7.10]{KM}, $\ker{\pi_2}$ is the standard cyclic subgroup of $\ker{p_2 \circ \pi_2}=C[st]$ of degree $\deg{\pi_2}=s$, so it is $C[s]$ by \cite[Theorem 6.7.4]{KM}. Similarly, one has $\ker{\pi_1}=C[t]$. The subgroups $\pi_1(C[t_{\infty}])[\frac{t_{\infty}}{t'}],\,C_1[\frac{t_{\infty}}{t'}]$ are cyclic with degree dividing $t_{\infty}$. Since $p_1$ has degree coprime to $t_{\infty}$, they are equal if $p_1(\pi_1(C[t_{\infty}])[\frac{t_{\infty}}{t'}])=p_1(C_1[\frac{t_{\infty}}{t'}])$. Now, 
\begin{align*}
p_1(\pi_1(C[t_{\infty}])[\frac{t_{\infty}}{t'}])&=\iota^{-1}(p_2 \circ \pi_2(C[t_{\infty}]))[\frac{t_{\infty}}{t'}] = \iota^{-1}(p_2(C_2[t_{\infty}])[\frac{t_{\infty}}{t'}])\\
&=\iota^{-1}((C_2[t_{\infty}] \mod{C_2[t]})[\frac{t_{\infty}}{t'}])=\iota^{-1}(C_2 \mod{C_2[t]})[\frac{m'}{s't'}][\frac{t_{\infty}}{t'}]\\
&=(C_1 \mod{C_1[s]})[\frac{m'}{s't'}][\frac{t_{\infty}}{t'}]=p_1(C_1[\frac{t_{\infty}}{t'}]),
\end{align*}

and the same holds for $C_2$, which concludes the proof of the claim. \\

\emph{Step 2: Injectivity} 

Let $z=(F,\beta,D) \in \MPp{\P,[\Gamma_0(m)]_R}(S)$ be such that $D_{s,\frac{m}{s'}}(z)=x_2,D_{t,\frac{m}{t'}}(z)=x_1$. Let us show that $z$ is the same point as $(E,\alpha,C)$. 

Since $D_{t,\frac{m}{t'}}((F,\beta,D))=x_1$, there is a cyclic isogeny $j_1: F \rar E_1$ of degree $t$ with kernel $D[t]$ such that $(\P j_1)(\beta)=\alpha_1$, $C_1=j_1(D)[\frac{m}{t'}]$. Similarly, there is a cyclic isogeny $j_2: F \rar E_2$ of degree $s$ with kernel $C[s]$ and such that $(\P j_2)(\beta)=\alpha_2$ and $C_2=j_2(D)[\frac{m}{s'}]$. 

The couple $j:=(j_1,j_2): F \rar E=E_1 \times_{E_1/C_1[s]} E_2$ is a morphism of $S$-elliptic curves with $\pi_1 \circ j=j_1$, so $j$ is an isogeny of degree one, i.e. an isomorphism. In particular, \[(\P j)(\beta)=(\P \pi_1)^{-1}(\P j_1)(\beta)=(\P \pi_1)^{-1}\alpha_1=\alpha.\] 
The subgroup scheme $C' := j(D)$ of $E$ is cyclic of degree $m$. To prove that $C=C'$, it is thus enough to find two integers $a,b$ with least common multiple $m$ such that $C'[a]=C[a]$ and $C'[b]=C[b]$.  
Now, one has $\pi_1(C')[\frac{m}{t'}]=j_1(D)[\frac{m}{t'}]=C_1[\frac{m}{t'}]$, so $\pi_1(C'[\frac{m}{t_{\infty}}])=C_1[\frac{m}{t_{\infty}}]=\pi_1(C[\frac{m}{t_{\infty}}])$. Since $\pi_1$ has degree prime to $\frac{m}{t_{\infty}}$, one has $C'[\frac{m}{t_{\infty}}]=C[\frac{m}{t_{\infty}}]$. Similarly, $C'[\frac{m}{s_{\infty}}]=C[\frac{m}{s_{\infty}}]$, whence $C=C'$, which concludes the proof. }

\prop[AL-auto]{Let $m$ be a positive integer and $d \mid m$ be positive and prime to $\frac{m}{d}$. For any representable moduli problem $\P$ on $\Ell_R$ of level coprime to $m$, there is an \emph{Atkin-Lehner automorphism} $w_d$ of $\MPp{\P,[\Gamma_0(m)]_R}$ defined by the following rule. Let $E/S$ be an elliptic curve endowed with a cyclic subgroup $C$ of order $m$, and $\alpha \in \P(E/S)$. Let $C_d, C_{\frac{m}{d}}$ be the cyclic subgroups of $C$ of respective degrees $d,\frac{m}{d}$, $\pi: E \rar E'$ the cyclic isogeny of kernel $C_d$. Let $C'$ the finite flat $S$-subgroup scheme of $E'[m]$ of degree $m$ defined as follows: its $d$-torsion is the kernel of $\pi^{\vee}$, and its $\frac{m}{d}$-torsion is $\pi(C)=\pi(C_{\frac{m}{d}})$. Then $w_d$ maps $(E/S,\alpha,C)$ to $(E'/S,\pi(\alpha),C')$. 

Furthermore,
\begin{itemize}[noitemsep,label=$-$]
\item The formation of any $w_d$ is natural with respect to $\P$ (for morphisms of representable moduli problems of level $N$ coprime to $m$, such as the action of $(\Z/N\Z)^{\times}$) and commutes with base change,
\item If $d,d' \mid m$ are positive integers such that $d$ is coprime to $\frac{m}{d}$ and $d'$ is coprime to $\frac{m}{d'}$, then one has $w_d \circ w_{d'} = [\delta]\circ w_{\frac{dd'}{\delta^2}}$, where $\delta$ is the greatest common divisor of $d,d'$.
\item If $d,t \mid m$ are positive integers such that $dt \mid m$ and $d$ coprime to $\frac{m}{d}$, then $D_{d,t}=D_{1,t} \circ w_d$. 
\end{itemize}}

\demo{That the definition constructs a morphism of $R$-schemes of $\MPp{\P,[\Gamma_0(m)]_R}$ into itself follows from \cite[Lemma 3.5.1, Theorem 6.7.4]{KM}. In fact, since $d$ and $\frac{m}{d}$ are coprime, we do not need to involve in this case the notion of standard subgroup; we will use freely the combination of these results throughout the proof. The naturality and commutation to base change is immediate. As in the previous proof, if $D$ is a cyclic subgroup of order $k$ of some elliptic curve $F/T$, and $\ell \geq 1$ is a divisor of $k$, then $D[\ell]$ denotes the standard cyclic subgroup of $D$ of order $\ell$. 

Note that the last identity is a direct computation. 

 To check that the $w_d$ are automorphisms, it is enough to prove that $w_d \circ w_{d'} = [\delta] \circ w_{\frac{dd'}{\delta^2}}$ when $d,d'$ are prime powers, which we now assume.
	
When $d=d'$, $w_d$ maps a triple $(E/S,\alpha,C)$ to $(E'/S,\pi(\alpha),C')$, where $\pi: E \rar E'=E/C[d]$ is the cyclic isogeny of degree $d$ whose kernel is exactly the $d$-torsion subgroup of $C$, and $C'$ is the subgroup $\pi(C)+\ker{\pi^{\vee}}$ (by \cite[Lemma 3.5.1, Theorem 6.7.4]{KM}). Thus 
\begin{align*}
w_d^2(E/S,\alpha,C)&=(E'/C'[d],(\pi^{\vee}\circ \pi)(\alpha),\pi^{\vee}(C'))=(E'/\ker{\pi^{\vee}},[d]\alpha,\pi^{\vee}(\pi(C))+\ker{(\pi^{\vee})^{\vee}})\\
&=(E/S,[d]\alpha,dC+C[d])=[d](E/S,\alpha,C).\end{align*}

When $d \neq d'$ (so that $d,d',\frac{m}{dd'}$ are pairwise coprime), $w_{d'}$ maps a triple $(E/S,\alpha,C)$ to $(E_1/S,\pi_1(\alpha),C_1)$, where $\pi_1: E \rar E_1$ is the isogeny whose kernel is the $d'$-torsion of $C$, and $C_1=\ker{\pi_1^{\vee}}+\pi_1(C)$; then $w_d(E_1/S,\pi_1(\alpha),C_1)=(E_2/S,\pi_2(\pi_1(\alpha)),C_2)$, where $\pi_2: E_1 \rar E_2$ is the isogeny whose kernel is the $d$-torsion of $C_1$, and $C_2=\ker{\pi_2^{\vee}}+\pi_2(C_1)$. 

Since $\pi_2 \circ \pi_1$ is an isogeny of degree $dd'$, we need to prove that $\ker{\pi_2\circ\pi_1}=C[dd']$, that $C_2\left[\frac{m}{dd'}\right]=\pi_2(\pi_1(C))$, and that $C_2[dd'] = \ker{\pi_1^{\vee}\circ \pi_2^{\vee}}$. 

One has $C_2\left[\frac{m}{dd'}\right] = C_2\left[\frac{m}{d}\right]\left[\frac{m}{dd'}\right]=\pi_2(C_1)\left[\frac{m}{dd'}\right]$. Since $\pi_1,\pi_2$ have degree coprime to $\frac{m}{dd'}$, one thus has \[C_2\left[\frac{m}{dd'}\right]=\pi_2(C_1\left[\frac{m}{dd'}\right])=\pi_2(\pi_1(C)\left[\frac{m}{dd'}\right])=(\pi_2\circ \pi_1)(C)\left[\frac{m}{dd'}\right].\]  

Since $\pi_1,\pi_2$ have coprime degrees and are cyclic, $\pi_2 \circ \pi_1$ is cyclic by \cite[Proposition 6.7.10]{KM}, and its kernel $K=\pi_1^{-1}(C_1[d])=\pi_1^{-1}(\pi_1(C)[d])$ is such that $C[d]$ is its (standard) cyclic subgroup of degree $d$. Thus $K=\pi_1^{-1}(\pi_1(C))[d]=\ker{\pi_1}+C[d']=C[dd']$. 

Furthermore, $\ker{\pi_1^{\vee}\circ \pi_2^{\vee}}$ is a cyclic subgroup of $E_2$ of degree $dd'$. Since $\pi_1$ has degree prime to $d$, $(\ker{\pi_1^{\vee}\circ \pi_2^{\vee}})[d] = (\ker{\pi_2^{\vee}})[d] =C_2[d]$. Since $\pi_2$ has degree prime to $d'$, $(\ker{\pi_1^{\vee}\circ \pi_2^{\vee}})[d']$ is equal to the pre-image in $E_2[d']$ under $\pi_2^{\vee}$ of $\ker{\pi_1^{\vee}}=C_1[d']$. Since \[(\pi_2^{\vee})^{-1}(C_1[d'])=(\pi_2^{\vee})^{-1}(\pi_2^{-1}(\pi_2(C_1[d'])))=(d)^{-1}C_2[d']=C_2[d']\] (where all pre-images are being taken inside the $d'$-torsion), it follows that $\ker{\pi_1^{\vee}\circ \pi_2^{\vee}}=C_2[dd']$, which concludes. 
}

\prop[deg-AL-compat]{Let $\P$ be a representable moduli problem on $\Ell_R$ of level $N$ and $m,m' \geq 1$ be such that $N,m,m'$ are pairwise coprime. Let $\P'=\P \times [\Gamma_0(m)]_R$. There is a natural (with respect to $\P$ and morphisms of moduli problems of level coprime to $mm'$) isomorphism of moduli problems (of level $Nmm'$, and over $\P$) $\P' \times [\Gamma_0(m')]_R \simeq \P \times [\Gamma_0(mm')]_R$. 

Let $d,r,s \geq 1$ be divisors of $m'$ such that $d$ and $\frac{m'}{d}$ are coprime and $rs \mid m'$. Then, under the above isomorphism, the following diagrams commute:
\[
\begin{tikzcd}[ampersand replacement=\&]
\MPp{\P',[\Gamma_0(m')]_R} \arrow{d}{\sim} \& \& \&\MPp{\P',[\Gamma_0(m')]_R} \arrow{d}{\sim} \arrow{lll}{w_d(\P',[\Gamma_0(m')]_R)} \arrow{rrr}{D_{r,s}(\P',[\Gamma_0(m')]_R)} \&\&\& \MPp{\P',[\Gamma_0(s)]} \arrow{d}{\sim}\\
\MPp{\P,[\Gamma_0(mm')]_R} \&\& \&\MPp{\P,[\Gamma_0(mm')]_R} \arrow{lll}{w_d(\P,[\Gamma_0(mm')]_R)} \arrow{rrr}{D_{r,sm}(\P,[\Gamma_0(mm')]_R)}\& \&\&\MPp{\P,[\Gamma_0(sm)]_R}
\end{tikzcd}
\]}

\demo{This is completely formal, using \cite[Lemma 3.5.1, Theorem 6.7.4]{KM} for the relevant results on cyclic subgroups.}

\cor[flat-deg]{Let $\P$ be a representable moduli problem on $\Ell_R$ of level $N$, let $m \geq 1$ be coprime to $N$, and let $d,t \mid m$ be positive divisors such that $td \mid m$. Then $D_{d,t}$ is finite locally free, and \'etale above $\Z[1/m]$.}

\demo{When $d=1$, finiteness (resp. \'etaleness above $\Z[1/m]$) follows from the fact that $D_{1,t}$ is a morphism of finite locally free $\MP$-schemes (resp. finite \'etale $\MP$-schemes above $\Z[1/m]$) by \cite[Theorem 5.1.1]{KM}. For arbitrary $d$, it is formal to check that $D_{1,1} \circ D_{d,t}=D_{1,1} \circ w_d \circ D_{1,d}$, so that the composition $D_{1,1} \circ D_{d,t}$ is finite (resp. \'etale above $\Z[1/m]$), hence $D_{d,t}$ is finite (resp. \'etale above $\Z[1/m]$). To prove flatness and finite presentation, first assume that $\P=[\Gamma(\ell)]_{\Z[1/\ell]}$ for some prime $\ell \nmid 2m$.

In this case, finite presentation is automatic, the moduli schemes are regular affine with every nonempty open subset of dimension $2$ (because $\P$ is a modular family and the moduli problems $[\Gamma_0(s)]_{\Z}$ are regular by \cite[Theorem 5.1.1]{KM}), and $D_{d,t}$ is finite, so the miracle flatness theorem \cite[Lemma 00R4]{Stacks} applies and we are done. 

Let us treat the general case. The claim is Zariski-local with respect to $R$, so we may assume that there exists a prime $\ell \nmid 2Nm$ invertible in $R$. Since $[\Gamma(\ell)]_{\Z[1/\ell]}$ is representable finite \'etale locally of constant rank over $\Ell_{\Z[1/\ell]}$, for any divisor $s \geq 1$ of $m$, the forgetful map $\MPp{\P\times[\Gamma(\ell)]_R,[\Gamma_0(s)]_R} \rar \MPp{\P,[\Gamma_0(s)]_R}$ is finite \'etale surjective. By Proposition \ref{degn-cartesian}, the two squares in the following diagram are Cartesian.  
\[
\begin{tikzcd}[ampersand replacement=\&]
\MPp{\P,[\Gamma_0(m)]_R} \arrow{d}{D_{d,t}} \& \MPp{\P\times[\Gamma(\ell)],[\Gamma_0(m)]_R} \arrow{d}{D_{d,t}}\arrow{l}\arrow{r} \& \MPp{[\Gamma(\ell)]_R,[\Gamma_0(m)]_R}\arrow{d}{D_{d,t}}\\
\MPp{\P,[\Gamma_0(t)]_R} \& \MPp{\P\times[\Gamma(\ell)]_R,[\Gamma_0(t)]_R} \arrow{l}\arrow{r} \& \MPp{[\Gamma(\ell)]_R,[\Gamma_0(t)]_R}
\end{tikzcd}
\]

Therefore, the special case above shows that the rightmost vertical arrow is finite locally free, thus the same holds for the middle vertical arrow by base change, and thus for the leftmost vertical arrow as well by fpqc descent.
}

\subsection{Degeneracies for $\Gamma_1(N)$-structures}
\label{gamma1m-structure}

\defi{Let $N \geq 1$ be an integer, and $r \mid N$ be a positive divisor. Let $R$ be a $\Z[1/N]$-algebra. The moduli problem $\P=[\Gamma_1'(N,r)]_R$ on $\Ell_{R}$ is defined as follows: given a relative elliptic curve $E/S$, $\P(E/S)$ is the collection of couples $(P,C)$, where $P \in E(S)$ has exact order $N$, $C$ is a cyclic subgroup of $E/S$ of order $r$ such that the image of $P$ in $(E/C)(S)$ has exact order $N$.}

\prop[basic-with-bad]{Let $N,r \geq 1$ with $r \mid N$ and $R$ be a $\Z[1/N]$-algebra. Then $\P=[\Gamma_1'(N,r)]_R$ is relatively representable and finite \'etale of constant rank. It is endowed with a natural level $N$ structure. It is representable if $N \geq 4$.}

\demo{By Proposition \ref{base-change-good-property}, we may assume that $R=\Z[1/N]$. If $N \geq 4$, since $\P$ is a subfunctor of the representable finite \'etale moduli problem $\P'=[\Gamma_1(N)]_R \times [\Gamma_0(r)]_R$ (by \cite[Corollary 2.7.4, Theorem 5.1.1]{KM} and Proposition \ref{product-rel-rep}), $\P$ is rigid.

Next, we show that for any relative elliptic curve $E/S$, the subsheaf $T \longmapsto \P(E_T/T)$ of $T \longmapsto \P'(E_T/T)$ is representable by a closed open subscheme of $\P'_{E/S}$. We define a map from $\P'_{E/S}$ to a disjoint union of copies $S_d$ of $S$ indexed by divisors $d$ of $N$ as follows: for any $(P,C) \in \P'(E_T/T)$, the image of $P$ in the $S$-elliptic curve $E_T/C$ is a section of the finite \'etale $T$-group scheme $(E_T/C)[N]$ of constant rank. Therefore, its fibrewise order $\omega$ is a locally constant function on $T$, its values divide $N$, and the formation of $\omega$ commutes with respect to base change. Then we map $(P,C) \in \P'(E/S)$ to the structure map $\omega^{-1}(d) \subset T \rar S \simeq S_d$ for each $d \mid N$. Thus $\P(E/S)$ corresponds to the inverse image of $S_N$, whence the conclusion. }

\prop[bad-degeneracy]{Let $N,r \geq 1$ be integers with $r \mid N$ and $R$ be a $\Z[1/N]$-algebra. Let $\P$ be a relatively representable moduli problem on $\Ell_R$ of level prime to $N$. Assume that $\P$ is representable or that $N \geq 4$. Let $s,t \geq 1$ be such that $st \mid r$. 

There exists a \emph{degeneracy map} $D'_{s,t}: \MPp{\P,[\Gamma_1'(N,r)]_R} \rar \MPp{\P,[\Gamma_1'(N,t)]_R}$ defined by the following rule. Given an elliptic curve $E/S$ and a couple $(\alpha,(P,C)) \in (\P \times [\Gamma_1'(N,r)]_R)(E/S)$, let $C'$ be the standard subgroup of $C$ of order $s$, and $\pi: E \rar E':=E/C'$ the natural (\'etale) isogeny. Then $D'_{s,t}$ maps the point $(E/S,\alpha,(P,C))$ to $(E'/S,(\P\pi)(\alpha),(\pi(P),\pi(C)[t]))$, where $\pi(C)[t]$ is the standard cyclic subgroup of order $t$ of the cyclic subgroup $\pi(C)$ of order $r/s$.}

\demo{One still needs to check that $\pi(P)$ has exact order $N$ and that the image of $\pi(P)$ in $E'/\pi(C)[t]$ has exact order $N$. Since $N$ is invertible in $R$, everything can be checked fibrewise. Note that we have a successive composition $\pi_{\infty}: E \rar E' \rar E'/\pi(C)[t] \rar E/C$ of \'etale isogenies, and that $P,\pi_{\infty}(P)$ both have exact order $N$, whence the conclusion.  }

\lem[bad-degn-composition]{Let $N,r \geq 1$, $R$ and $\P$ are as in Proposition \ref{bad-degeneracy}. Let $s,t,u,v \geq 1$ be integers with $st \mid r$, and $uv \mid t$. Then the composition $D'_{u,v} \circ D'_{s,t}$ is exactly $D'_{su,t}$.}

\demo{This is a direct verification.}

\prop[bad-degn-cartesian]{Let $N,r,s,t \geq 1$ and $R$ be as in Proposition \ref{bad-degeneracy}, and let $F: \P \Rightarrow \P'$ be a morphism of relatively representable (representable if $N \leq 3$) moduli problems of level coprime to $N$ over $\Ell_R$. Then the following diagram is Cartesian:
\[
\begin{tikzcd}[ampersand replacement=\&]
\MPp{\P,[\Gamma_1'(N,r)]_R} \arrow{r}{D'_{s,t}} \arrow{d}{F}\& \MPp{\P,[\Gamma'_1(N,t)]_R} \arrow{d}{F} \\ 
\MPp{\P',[\Gamma_1'(N,r)]_R} \arrow{r}{D'_{s,t}} \& \MPp{\P',[\Gamma'_1(N,t)]_R}
\end{tikzcd}
\]
}

\demo{Let $E$ be an elliptic curve over a $R$-scheme $S$ and $\alpha \in \P(E/S)$. 

Let $(P,C) \in [\Gamma_1'(N,r)]_R(E/S)$, then 
\begin{align*}
F(D'_{s,t}((E,\alpha,(P,C)))) &=  F((E/C[s],\alpha \mod{C[s]},(P \mod{C[s]},C[st]/C[s]))) \\
&= (E/C[s],F(\alpha \mod{C[s]}),(P \mod{C[s]},C[st]/C[s])),
\end{align*} where $[n]$ means taking the standard cyclic subgroup of order $n$. 

On the other hand, since $F$ has level prime to $N$, it is functorial with respect of isogenies of degree dividing $N$, so that

\begin{align*}
D'_{s,t}(F((E,\alpha,(P,C)))) &= D'_{s,t}((E,F(\alpha),(P,C))) \\
&= (E/C[s],F(\alpha) \mod{C[s]},(P\mod{C[s]},C[st]/C[s]))\\
& = F(D'_{s,t}((E,\alpha,(P,C)))).
\end{align*}

Thus the diagram commutes. We then show that for any $R$-scheme $S$, the following diagram is Cartesian.  
\[
\begin{tikzcd}[ampersand replacement=\&]
\MPp{\P,[\Gamma_1'(N,r)]_R}(S) \arrow{r}{D'_{s,t}} \arrow{d}{F}\& \MPp{\P,[\Gamma'_1(N,t)]_R}(S) \arrow{d}{F} \\ 
\MPp{\P',[\Gamma_1'(N,r)]_R}(S) \arrow{r}{D'_{s,t}} \& \MPp{\P',[\Gamma'_1(N,t)]_R}(S)
\end{tikzcd}
\]

Indeed, let \[x_1=(E_1,\alpha',(P_1,C_1)) \in \MPp{\P',[\Gamma_1'(N,r)]_R}(S),\, x_2=(E_2,\beta,(P_2,C_2)) \in \MPp{\P,[\Gamma_1'(N,t)]_R}(S)\] be two points with $F(x_2)=D'_{s,t}(x_1)$. Thus, there exists an isogeny $\pi: E_1 \rar E_2$ such that $\ker{\pi}=C_1[s]$, $(\P'\pi)(\alpha')=F(\beta)$, $\pi(P_1)=P_2$ and $C_2=\pi(C_1[st])$. In particular, we can define $\alpha=(\P \pi)^{-1}(\beta)$, so that $(\P'\pi)F(\alpha)=F((\P\pi)(\alpha))=F(\beta)=(\P'\pi)(\alpha')$, hence $F(\alpha)=\alpha'$. Thus $x_0:=(E_1,\alpha,(P_1,C_1)) \in \MPp{\P,[\Gamma_1'(N,r)]_R}(S)$ is such that $F(x_0)=x_1$. But one immediately checks by construction that $D'_{s,t}(x_0)=x_2$. 

Let $y_0=(E_0,\alpha_0,(P_0,C_0)) \in \MPp{\P,[\Gamma_1'(N,r)]_R}(S)$ such that $F(y_0)=x_1$ and $D'_{s,t}(y_0)=x_2$. This means that there is an isomorphism $\iota: E_0 \rar E_1$ such that $F((\P\iota)(\alpha_0))=\alpha'$, $\iota(P_0)=P_1$, $\iota(C_0)=C_1$, and an isogeny $g: E_0 \rar E_2$ such that $\ker{g}=C_0[u]$, $\beta=(\P g)(\alpha_0)$, $C_2=\pi(C_0[uv])$, $P_2=\pi(P_0)$. Since $\pi \circ \iota$ and $g$ are isogenies with the same kernel and same target, there is an automorphism $\psi$ of $E_2$ such that $\psi \circ g=\pi \circ \iota$. 

Let us prove that $\psi$ is the identity, since $\P' \times [\Gamma_1(N)]_R$ is rigid, it is enough to show that $(\P' g)(F(\alpha_0))=(\P' \pi)\circ (\P'\iota)(F(\alpha_0))$, and that $g(P_0)=\pi(\iota(P_0))$. Indeed, one has 

\begin{align*}
(\P' \pi)\circ ( \P'\iota)(F(\alpha_0))&=(\P'\pi)\circ F((\P\iota)(\alpha_0))=(\P'\pi)(\alpha')=F(\beta)=(\P' g)(F(\alpha_0)),\\
 \pi(\iota(P_0))&=\pi(P_1)=P_2=g(P_0).
 \end{align*}

Thus $\pi \circ \iota=g$. Thus $(\P \iota)(\alpha_0)=(\P \pi)^{-1}(\P g)(\alpha_0)=(\P \pi)^{-1}(\beta)=\alpha$, and $\iota: y_0 \rar x_0$ is an isomorphism. }

\cor[bad-degn-fin-etale]{In the situation of Proposition \ref{bad-degeneracy}, $D'_{s,t}$ is finite \'etale.}

\demo{Suppose first that $k \geq 3$ is coprime to $N$, $R=\Z[1/(Nk)]$ and that $\P=[\Gamma(k)]_R$. Consider the following diagram, where $f,g$ are forgetful maps, finite \'etale (they are compositions of closed open immersions by the proof of Proposition \ref{basic-with-bad} and forgetting the $[\Gamma_1(N)]_R$ component, which is finite \'etale):

\[
\begin{tikzcd}[ampersand replacement=\&]
\MPp{[\Gamma(k)]_R,[\Gamma_1'(N,r)]_R} \arrow{r}{f} \arrow{d}{D'_{s,t}} \& \MPp{[\Gamma(k)]_R,[\Gamma_0(r)]_R}\arrow{d}{D_{s,t}} \\
\MPp{[\Gamma(k)]_R,[\Gamma_1'(N,t)]_R} \arrow{r}{g} \& \MPp{[\Gamma(k)]_R,[\Gamma_0(t)]_R} 
\end{tikzcd}
\] 

Since $D_{s,t}$ is finite \'etale by Corollary \ref{flat-deg}, $D'_{s,t}$ is finite \'etale.

Let us come back to the general case. The statement is Zariski-local with respect to $R$, so we may assume that some integer $k \geq 3$ coprime to $N$ is invertible in $R$. Then both squares in the following diagram are Cartesian, with all vertical maps being $D'_{s,t}$:

\[
\begin{tikzcd}[ampersand replacement=\&]
\MPp{[\Gamma(k)]_R,[\Gamma_1'(N,r)]_R} \arrow{d}{\alpha} \& \MPp{[\Gamma(k)]_R,[\Gamma_1'(N,r)]_R,\P} \arrow{l}\arrow{d}{\beta} \arrow{r}{\delta_r} \& \MPp{\P,[\Gamma_1'(N,r)]_R} \arrow{d}{\gamma}\\
\MPp{[\Gamma(k)]_R,[\Gamma_1'(N,t)]_R}  \& \MPp{[\Gamma(k)]_R,[\Gamma_1'(N,t)]_R,\P} \arrow{l} \arrow{r}{\delta_t} \& \MPp{\P,[\Gamma_1'(N,t)]_R} 
\end{tikzcd}
\]

By what we already saw, $\alpha$ is \'etale, thus so is $\beta$. Since $[\Gamma(k)]_R$ is finite \'etale over $\Ell_{\Z[1/k]}$, $\delta_r,\delta_t$ are finite \'etale. Thus, by \'etale descent, $\gamma$ is finite \'etale and we are done.}

\prop[mixed-degn-cartesian]{Let $N,r,s,t \geq 1$ be as in Definition \ref{bad-degeneracy} and $\P$ be a relatively representable moduli problem of level $L$ coprime to $N$, representable if $N \leq 3$. Let $m\geq 1$ be coprime to $NL$. Then, for any $u,v \geq 1$ such that $uv\mid m$, the following diagram is Cartesian:

\[
\begin{tikzcd}[ampersand replacement=\&]
\MPp{\P,[\Gamma_1'(N,r)]_R,[\Gamma_0(m)]_R} \arrow{r}{D'_{s,t}} \arrow{d}{D_{u,v}}\& \MPp{\P,[\Gamma'_1(N,t)]_R,[\Gamma_0(m)]_R} \arrow{d}{D_{u,v}} \\ 
\MPp{\P,[\Gamma_1'(N,r)]_R,[\Gamma_0(v)]_R} \arrow{r}{D'_{s,t}} \& \MPp{\P,[\Gamma'_1(N,t)]_R,[\Gamma_0(v)]_R}
\end{tikzcd}
\]
}

\demo{\emph{Step 1: The diagram commutes.}

Let $S$ be a $R$-scheme, and $x:=(E/S,\alpha,(P,C),D) \in \MPp{\P,[\Gamma_1'(N,r)]_R,[\Gamma_0(m)]_R}(S)$. Let us prove that $D_{u,v}(D'_{s,t}(x))=D'_{s,t}(D_{u,v}(x))$. 

Let $\pi_s: E \rar E/C[s]$ be the natural isogeny of degree $s$. One has \[D'_{s,t}(x)=(E/C[s],(\P\pi_s)(\alpha),(\pi_s(P),\pi_s(C)[t]),\pi_s(D)),\] where, as above, $[\cdot]$ denotes the standard cyclic subgroup of given degree. Note that $C[s]+D$ is a cyclic subgroup of $E$ of degree $sm$ by the factorization lemma \cite[(3.5)]{KM}, so $\pi_s(D)=\pi_s(D+C_s)$ is indeed a cyclic subgroup of degree $d$ by \cite[Theorem 6.7.4]{KM}. Let $\pi_{s,u}$ denote the quotient isogeny $E/C[s] \rar (E/(C[s]))/\pi_s(D)[u]$. Then, with $\pi^1=\pi_{s,u} \circ \pi_s: E \rar (E/(C[s]))/(\pi_s(D)[u])$ the cyclic isogeny (by \cite[Proposition 6.7.10]{KM}) of kernel $C[s]+D[u]$, one has:

\begin{align*}
D_{u,v}(D'_{s,t}(x)) &= \Big((E/C[s])/\pi_{s}(D)[u],&&(\P\pi_{s,u})(\P\pi_s)(\alpha),((\pi_{s,u}(\pi_s(P)),\pi_{s,u}(\pi_s(C)[t]))),\\
& &&\pi_{s,u}(\pi_s(D))[v]\Big)\\
&= \Big((E/C[s])/\pi_s(D)[u],&&(\P\pi^1)(\alpha),(\pi^1(P),\pi^1(C+D[u])[t]),\\
& && \pi^1(C[s]+D)[t]\Big).
\end{align*}

Similarly, if $\pi_u: E \rar E/C[u]$ is the natural isogeny and $\pi^2: E \rar (E/D[u])/\pi_u(C)[s]$ is the composition of $\pi_u$ and the quotient by $\pi_U(C)[s]$, one sees that 
\[D'_{s,t}(D_{u,v}(x)) = ((E/D[u])/\pi_u(C)[s],(\P\pi^2)(\alpha),(\pi^2(P),\pi^2(C+D[u])[t]),\pi^2(C[s]+D)[t])\]

Since $\pi^1,\pi^2$ are isogenies from $E$ with the same kernel, there is an isomorphism \[\iota: (E/(C[s]))/(\pi_s(D)[u])\rar (E/D[u])/\pi_u(C)[s]\] such that $\iota \circ \pi^1=\pi^2$, so the diagram commutes.

\emph{Step 2: Injectivity}

Suppose we have points \[x=(E,\alpha,(P,C),D),x'= (E',\alpha',(P',C'),D')\in \MPp{\P,[\Gamma_1'(N,r)]_R,[\Gamma_0(m)]_R}(S)\] such that $D'_{s,t}(x)=D'_{s,t}(y)$ and $D_{u,v}(x)=D_{u,v}(x')$. Then there are isomorphisms 
\begin{align*}
j_s&: E/C[s] \rar E/C'[s],\quad (\alpha,P,C,D) \mod{C[s]} \mapsto (\alpha',P',C',D') \mod{C'[s]},\\
j_u&: E/D[u] \rar E'/D'[u],\quad (\alpha,P,C,D) \mod{D[u]} \mapsto (\alpha',P',C',D') \mod{D'[u]}.
\end{align*}
The assumptions imply that the two compositions 
\begin{align*}
c_s&: E \rar E/C[s] \overset{j_s}{\rar} E'/C'[s] \rar E'/(C'[s]+D'[u]),\\
 c_u&: E \rar E/D[u] \overset{j_u}{\rar}E'/D'[u] \rar E'/(C'[s]+D'[u])
 \end{align*}
have the same kernel $C[s]+D[u]$, so they differ by an automorphism $\lambda$ of $E'/(C'[s]+D'[u])$. By construction, \[(\alpha',P') \mod{C'[s]+D'[u]} \in (\P \times [\Gamma_1(N)]_R)(E'/(C'[s]+D'[u]))\] is fixed by the automorphism $\lambda$: since $\P \times [\Gamma_1(N)]_R$ is rigid, $\lambda$ is the identity and $c_s=c_u$. 

Since $E' \rar E'/C'[s] \times_{E'/(C'[s]+D'[u])} E'/D'[u]$ is an isomorphism (it is a morphism of elliptic curves, and its kernel is contained in $C'[s]$ and in $D'[u]$, so it is a monomorphism thus an isomorphism), the couple $(E \rar E/C[s] \overset{j_s}{\rar} E'/C'[s],E \rar E/D[u] \overset{j_u}{\rar} E'/D'[u])$ defines an morphism $j: E \rar E'$ of elliptic curves. The composition of $j$ with the cyclic isogeny $E' \rar E'/C'[s]$ of degree $s$ is exactly $E \rar E/C[s] \overset{j_s}{\rar} E'/C'[s]$, which is an isogeny of degree $s$. So $j$ is not zero on any fibre, hence it is an isogeny of degree $\frac{s}{s}=1$, i.e. an isomorphism. 

Now, $(\P j)(\alpha) \mod{C'[s]}=(\P j_s)(\alpha \mod{C[s]}) = \alpha' \mod{C'[s]}$, and since $\P$ has level prime to $s$, $(\P j)(\alpha)=\alpha'$. The same reasoning shows that $(\P j)(D)=D'$. To prove that $(\P j)(P)=P'$ and $(\P j)(C)=C'$, one quotients instead by $D'[u]$, since $[\Gamma_1'(N,r)]_R$ has level prime to $u$. 

\emph{Step 3: Surjectivity}

Let finally 
\begin{align*}
x_1&=(E_1,\alpha_1,(P_1,C_1),D_1) \in \MPp{\P,[\Gamma_1'(N,r)]_R,[\Gamma_0(v)]_R}(S),\\ 
x_2&=(E_2,\alpha_2,(P_2,C_2),D_2) \in \MPp{\P,[\Gamma_1'(N,t)]_R,[\Gamma_0(m)]_R}(S)
\end{align*} be two points with the same image in $\MPp{\P,[\Gamma_1'(N,t)]_R,[\Gamma_0(v)]_R}(S)$. We want to prove that there exists $x_0 \in \MPp{\P,[\Gamma_1'(N,r)]_R,[\Gamma_0(m)]_R}(S)$ with $D_{u,v}(x_0)=x_1$, $D'_{s,t}(x_0)=x_2$. 

The assumption means that there exists an isomorphism $\iota: E_1/C_1[s] \rar E_2/D_2[u]$ mapping $\alpha_1 \mod{C_1[s]}$ to $\alpha_2 \mod{D_2[u]}$, $P_1 \mod{C_1[s]}$ to $P_2 \mod{D_2[u]}$, $(C_1 \mod{C_1[s]})[t] = C_2 \mod{D_2[u]}$, $D_1 \mod{C_1[s]} = (D_2 \mod{D_2[u]})[v]$. 

The scheme $E_0=E_1 \times_{E_1/C_1[s]} E_2$ is a relative elliptic curve over $S$; moreover, the projections $\pi_1: E_0 \rar E_1, \pi_2: E_0 \rar E_2$ are cyclic isogenies of degrees $u,s$ respectively. The assumptions imply that $(\P\pi_1)^{-1}(\alpha_1)=(\P \pi_2)^{-1}(\alpha_2)$ (because they have the same image in $E_1/C_1[s] \simeq E_2/D_2[u]$, through an isogeny of degree $su$ prime to the level of $\P$), and let $\alpha_0 \in \P(E_0/S)$ denote this element. Let $P_0 \in E_0(S)$ be the point $(P_1,P_2)$, then $P_0$ has exact order $N$ because so have $P_1$ and $P_2$. Because $u$ is coprime to $r$ (resp. $s$ is coprime to $m$), there is a unique cyclic subgroup $C_0$ (resp. $D_0$) of $E_0$ of order $r$ (resp. $m$) such that $\pi_1(C_0)=C_1$. (resp. $\pi_2(D_0)=D_2$).

We claim that $x_0=(E_0,\alpha_0,(P_0,C_0),D_0) \in \MPp{\P,[\Gamma_1(N,r)]_R,[\Gamma_0(m)]_R}(S)$ is such that $D'_{s,t}(x_0)=x_2$ and $D_{u,v}(x_0)=x_1$. 

To check that $x_0 \in \MPp{\P,[\Gamma_1(N,r)]_R,[\Gamma_0(m)]_R}(S)$, the only remaining assertion to check is that $P_0 \mod{C_0}$ has exact order $N$. Suppose that (in any geometric fibre) $dP_0 \mod{C_0} = 0$: then $dP_0 \in C_0$. Then, after applying $\pi_1$, we find that $dP_1 \in C_1$, hence $dP_1 \mod{C_1}=0$, so $N\mid d$, and we are done. 

By construction, $\pi_1: E_0 \rar E_1$ (resp. $\pi_2: E_0 \rar E_2$) is an isogeny of degree $u$ (resp. $s$) mapping $P_0$ to $P_1$ (resp. $P_2$), $\alpha_0$ to $\alpha_1$ (resp. $\alpha_2$) and $C_0$ to $C_1$ (resp. $D_0$ to $D_2$). So all we need to check is that $\ker{\pi_1}=D_0[u]$ and $D_1=\pi_1(D_0)[v]$, $\ker{\pi_2}=C_0[s]$ and $C_2=\pi_2(C_0)[t]$. 

The map $q_1: E_0 \overset{\pi_1}{\rar} E_1 \rar E_1/C_1[s]$ is an isogeny of degree $su$; its kernel contains $C_0[s]$ because $\pi_1(C_0[s])=\pi_1(C_0)[s]=C_1[s]$. Similarly, the map $q_2: E_0 \overset{\pi_2}{\rar} E_2 \rar E_2/D_2[u]$ is an isogeny of degree $su$ and its kernel contains $D_0[u]$. Since $q_2=\iota \circ q_1$, $q_1$ and $q_2$ have the same kernel $Q$, which is a subgroup scheme of $E$ (flat over $S$) of degree $su$ and contains the cyclic subgroup $C_0[s]+D_0[u]$ of degree $su$. In particular, $C_0[s]+D_0[u] \rar Q$ is a closed immersion of finite locally free $S$-schemes of degree $su$: it is an isomorphism (by checking the condition after each base change $\Sp{\OO_{S,s}} \rar S$, since a surjection of free modules of same rank is an isomorphism). In particular, the $q_i$ are cyclic, and by \cite[Proposition 6.7.10]{KM}, the kernel of $\pi_1$ (resp. $\pi_2$) is the standard cyclic subgroup of $C_0[s]+D_0[u]$ of degree $u$ (resp. $s$), i.e. $D_0[u]$ (resp. $C_0[s]$).   

Finally, one has \[q_2(C_0)[t]=\iota(q_1(C_0)[t])=\iota((\pi_1(C_0) \mod{C_1[s]})[t])=\iota((C_1 \mod{C_1[s]})[t])=C_2\mod{D_2[u]}.\] Now, $\pi_2(C_0)[t]$ and $C_2$ are two cyclic subgroups of $E_2$ of order $t$, and their image through the cyclic isogeny $E_2 \rar E_2/D_2[u]$ of degree $u$ coprime to $t$ are the same: hence $\pi_2(C_0[t])=C_2$. One shows similarly that $\pi_1(D_0)=D_1[v]$. 
}

\prop[two-bad-degn-cartesian]{Let $N,r \geq 1$, $R$ and $\P$ be as in Definition \ref{bad-degeneracy}. Let $s',t' \geq 1$ be coprime divisors of $r$, and $s,t \geq 1$ be divisors of $s',t'$ respectively. Then the following diagram is Cartesian:
\[
\begin{tikzcd}[ampersand replacement=\&]
\MPp{\P,[\Gamma_1'(N,r)]_R} \arrow{r}{D'_{s,\frac{r}{s'}}} \arrow{d}{D'_{t,\frac{r}{t'}}}\& \MPp{\P,[\Gamma'_1\left(N,\frac{r}{s'}\right)]_R} \arrow{d}{D_{t,\frac{r}{s't'}}} \\ 
\MPp{\P,[\Gamma_1'\left(N,\frac{r}{t'}\right)]_R} \arrow{r}{D'_{s,\frac{r}{s't'}}} \& \MPp{\P,[\Gamma'_1\left(N,\frac{r}{s't'}\right)]_R}
\end{tikzcd}
\]
}

\demo{This is similar to the proof of Propositions \ref{mixed-degn-cartesian} or \ref{two-degn-cartesian}.}

\defi[bad-AL]{Let $N\geq 1$ be an integer and $R$ be a $\Z[1/N]$-algebra. Let $\P$ be a relatively representable moduli problem over $\Ell_R$ (representable if $N \leq 3$) of level prime to $N$. Let $d \mid N$ be a positive integer coprime to $N/d$. We define as follows the morphism \[w'_d: \MPp{\P,[\Gamma_1(N)]_R} \times_{\Z[1/N]} (\mu_d^{\times})_{\Z[1/N]} \rar \MPp{\P,[\Gamma_1(N)]} \times_{\Z[1/N]} (\mu_d^{\times})_{\Z[1/N]}\] as the following morphism. Given an elliptic curve $E/S$, $\alpha \in \P(E/S)$, $P \in E(S)$ of exact order $N$, and $\zeta \in \mu_d^{\times}(S)$, let $\pi: E \rar E':=E/\langle N/d \cdot P\rangle$ be the natural isogeny of degree $d$. Then one has 
\[w'_d((E,\alpha,P,\zeta)) = (E',(\P\pi)(\alpha),Q,\zeta),\]
where $Q \in E'(S)[N]$ is the unique point such that $dQ=\pi(dP)$ and $\langle N/d \cdot P,\, N/d \cdot Q\rangle_{\pi}=\zeta$ (by \cite[(2.8.2)]{KM}, $N/d \cdot Q$ is well-defined, hence $Q$ is well-defined since $d$ and $N/d$ are coprime).
}

\rems{\begin{itemize}[noitemsep,label=$-$]
\item In this construction, $w'_d$ is obviously a morphism of $R \otimes_{\Z[1/N]} \Z[1/N,\mu_d]$-schemes, but the $(\mu_d^{\times})_{\Z[1/N]}$ cannot be ``factored into $\P$'' as a moduli problem of level one.
\item $w'_d$ has an obvious base change if we replace $\mu_d$ with $\mu_N$: in this case, the pairing condition becomes $\langle N/d \cdot P,\, N/d \cdot Q\rangle_{\pi}=\zeta^{N/d}$.
\item If one specifies instead some $\zeta_0 \in \mu_d^{\times}(\Sp{R})$ (that is, if $R$ is a $\Z[1/N,\zeta_d]$-algebra), then we can define a map $w'_d: \MPp{\P,[\Gamma_1(N)]_R} \rar \MPp{\P,[\Gamma_1(N)]_R}$ of $R$-schemes as the base change of Definition \ref{bad-AL} by $\zeta_0: \Sp{R} \rar \mu_d^{\times}$, and therefore satisfies the properties detailed in Lemmas \ref{bad-AL-natural}, \ref{bad-AL-ZnZ}, \ref{bad-AL-composition}.
\end{itemize}
}

\lem[bad-AL-natural]{The formation of $w'_d$ commutes with base change and morphisms of moduli problems of level coprime to $N$.}

\demo{This is formal.}

\lem[bad-AL-double]{With the notations of Definition \ref{bad-AL}, let $e=(e_p)_p \in \hat{\Z}^{\times}=\prod_p{\Z_p^{\times}}$ be such that $e_p=-1$ if $p \mid d$, and $e_p=d$ if $p \nmid d$. Then one has $w'_d \circ w'_d= [e]$. In particular, $w'_d$ is an isomorphism. }

\demo{Let $E/S$ be a relative elliptic curve, $\alpha \in \P(E/S)$, $P \in E(S)$ of exact order $N$, and $\zeta \in \mu_d^{\times}(S)$. 
Let $\pi: E \rar E':=E/\langle \frac{N}{d} \cdot P \rangle$ be the isogeny, and $Q$ be the point as in Definition \ref{bad-AL}, so that $(E'/S,(\P \pi)(\alpha),Q,\zeta) = w'_d(E/S,\alpha,P,\zeta)$. Then $\pi^{\vee}: E' \rar E$ is the degree $d$ isogeny whose kernel is generated by $\frac{N}{d} \cdot Q$. Moreover, one has $\pi^{\vee}(dQ)=\pi^{\vee}(\pi(dP))=d \cdot (dP)$. 

Finally, $\langle \frac{N}{d} \cdot Q, -\frac{N}{d} \cdot P\rangle_{\ker{\pi^{\vee}}} = \langle \frac{N}{d} \cdot P,\,\frac{N}{d} \cdot Q\rangle_{\ker{\pi}}=\zeta$ by \cite[(2.8.3)]{KM}. Therefore, if $P' \in E(S)$ is the unique point of exact order $N$ such that $dP'=d \cdot dP$, $\frac{N}{d} \cdot P'=-\frac{N}{d} \cdot P$, one has  
\begin{align*}
w'_d(w'_d((E/S,\alpha,P,\zeta)))&=w'_d((E'/S,(\P \pi)(\alpha),Q,\zeta))=(E/S,(\P \pi^{\vee})(\P \pi)(\alpha),P')\\
&=(E/S,[d]\alpha,e_1P), 
\end{align*}
where $e_1 \in (\Z/N\Z)^{\times}$ is congruent to $d$ mod $\frac{N}{d}$ and to $-1$ mod $d$. The conclusion follows.
}

\lem[bad-AL-ZnZ]{With the notation of Definition \ref{bad-AL}, let $e=(e_p)_p \in \hat{\Z}^{\times}$. Let $e'=(e'_p)_p \in \hat{\Z}^{\times}$ be given by $e'_p=e_p$ if $p \nmid d$, and $e'_p=e_p^{-1}$ otherwise. Then $w'_d \circ [e] = [e'] \circ w'_d$.}

\demo{This is a straightforward computation.}

\lem[bad-AL-cyclo]{With the notation of Definition \ref{bad-AL}, let $u \in (\Z/d\Z)^{\times}$, and $\underline{u}$ the automorphism of $\mu_d^{\times}$ whose cyclotomic character is $u$ (it does not act on the component $\MPp{\P,[\Gamma_1(N)]_P}$). Then one has $\underline{u}\circ [u] \circ w'_d = w'_d \circ \underline{u}$. }

\demo{This is also a direct computation.}

\lem[bad-AL-composition]{With the notations of Definition \ref{bad-AL}, after extending coefficients to $\OO_N$, if $d' \geq 1$ is another divisor of $N$ coprime to $\frac{dN}{d'}$, then $w'_{d'} \circ w'_d=[e] \circ w_{dd'}$, where $e \in (\Z/N\Z)^{\times}$ is congruent to $1$ mod $\frac{N}{d'}$ and to $d^{-1}$ mod $d'$.   }

\demo{Let $E/S$ be a relative elliptic curve, $\alpha \in \P(E/S)$, $\zeta \in \mu_N^{\times}(S)$ and $P \in E(S)$ be a point of exact order $N$. Let, for every $t \mid N$, $\zeta_t = \zeta^{N/t} \in \mu_N(S)$, which factors through $\mu_t^{\times}(S)$. 

Then $w'_d((E/S,\alpha,P,\zeta))=(E'/S,(\P\pi_1)(\alpha),P',\zeta)$, where $\pi_1: E \rar E' := E/\langle \frac{N}{d} \cdot P\rangle$ is the natural cyclic isogeny of degree $d$, $\langle \frac{N}{d} \cdot P,\,\frac{N}{d} \cdot P'\rangle_{\pi_1}=\zeta_d$, and $\pi_1(dP)=dP'$. 

Let $\pi_2: E' \rar E''=E'/\langle \frac{N}{d'} \cdot P'\rangle$ denote the natural cyclic isogeny of degree $d'$, and let $Q \in E''[N](S)$ be such that $d'Q=\pi_2(d'P')$ and $\langle \frac{N}{d'} \cdot P',\, \frac{N}{d'} \cdot Q\rangle=\zeta_{d'}$. Therefore, one has \[w'_{d'}(w_d((E/S,\alpha,P,\zeta)))=(E''/S,(\P\pi_2)(\P\pi_1)(\alpha),Q,\zeta).\] 

Since $\frac{N}{d'} \cdot P'=\frac{N}{dd'} \cdot d \cdot P'=\frac{N}{dd'} \cdot \pi_1(dP)=\pi_1(\frac{N}{d'} \cdot P)$, $\pi=\pi_2 \circ \pi_1$ is the cyclic isogeny (by \cite[Proposition 6.7.10]{KM}) of degree $dd'$ and of kernel generated by the point $\frac{N}{dd'} \cdot P$ of exact order $dd'$. In particular, we have $w'_{dd'}((E/S,\alpha,P,\zeta))=(E''/S,(\P\pi)(\alpha),R,\zeta)$ where $R \in E''(S)$ is another point of exact order $N$ such that $dd'R=\pi(dd'P)$, and $\langle \frac{N}{dd'} \cdot P,\, \frac{N}{dd'} \cdot R\rangle_{\pi}=\zeta_{dd'}$. 

To conclude, we have the following equalities using \cite[(2.8.3),(2.8.4.1)]{KM}:
\begin{align*}
dd'R=\pi(dd'P)&=\pi_2(\pi_1(dd'P))=\pi_2(d'\pi_1(dP))=\pi_2(d'dP')=d\pi_2(d'P')=d(d'Q)=dd'Q,\\
\langle \frac{N}{dd'}P,\,\frac{N}{dd'}Q\rangle_{\pi}^{(d')^2}&= \langle \frac{N}{d}P,\, \frac{N}{dd'}(d'Q)\rangle_{\pi_2\circ \pi_1} = \langle \frac{N}{d}P,\, \pi_2^{\vee}\left(\pi_2\left(\frac{N}{d}P'\right)\right)\rangle_{\pi_1}\\
&= \langle \frac{N}{d}P,\, d'\frac{N}{d}P'\rangle_{\pi_1}=\zeta_d^{d'}=\zeta_{d}^{ed'}=\zeta_{dd'}^{e(d')^2},\\
\langle \frac{N}{dd'}P,\,\frac{N}{dd'}Q\rangle_{\pi}^{d^2}&= \langle \frac{N}{d'}Q,\, \frac{N}{dd'}(dP)\rangle_{\pi_1^{\vee}\circ \pi_2^{\vee}}^{-1} = -\langle \frac{N}{d'}Q,\, \frac{N}{dd'} \cdot \pi_1\left(dP\right)\rangle_{\pi_2^{\vee}}^{-1}\\
&= \langle \frac{N}{d'}Q,\, \frac{dN}{dd'}P'\rangle_{\pi_2^{\vee}}^{-1}=\langle \frac{N}{d'}P',\,\frac{N}{d'}Q\rangle_{\pi_2} = \zeta_{d'}=\zeta_{d'}^{de}=\zeta_{dd'}^{d^2e}.
\end{align*}
}

\prop[bad-AL-good-deg]{Let $N \geq 1$ be an integer, $R$ be a $\Z[1/N]$-algebra and $\P$ be a relatively representable moduli problem over $\Ell_R$ (representable if $N \leq 3$) of level $L$ coprime to $N$. Let $d \geq 1$ be a divisor of $N$ coprime to $N/d$, $m \geq 1$ be an integer coprime to $NL$, and $u,v \geq 1$ such that $uv \mid m$. 
Let $u' \in (\Z/N\Z)^{\times}$ be congruent to $u$ mod $d$ and to $1$ mod $\frac{N}{d}$. 
Then the following diagram is Cartesian: 

\[
\begin{tikzcd}[ampersand replacement=\&]
\MPp{\P,[\Gamma_1(N)]_R,[\Gamma_0(m)]_R} \times_{\Z} \mu_d^{\times}\arrow{r}{w'_d} \arrow{d}{D_{u,v}} \&  \MPp{\P,[\Gamma_1(N)]_R,[\Gamma_0(m)]_R} \times_{\Z} \mu_d^{\times} \arrow{d}{D_{u,v}}\\
\MPp{\P,[\Gamma_1(N)]_R,[\Gamma_0(v)]_R} \times_{\Z} \mu_d^{\times}\arrow{r}{[u'] \circ w'_d}  \&  \MPp{\P,[\Gamma_1(N)]_R,[\Gamma_0(v)]_R} \times_{\Z} \mu_d^{\times}
\end{tikzcd}
\] 

}

\demo{Since the $w'_d$ are isomorphisms, it is enough to prove that the diagram commutes. Thus, let $E/S$ be a relative elliptic curve and $(E/S,\alpha,P,C,\zeta) \in (\MPp{\P,[\Gamma_1(N)]_R,[\Gamma_0(m)]_R} \times \mu_d^{\times})(S)$. 

Let $E'=E/\langle \frac{N}{d} \cdot P\rangle$, $\pi: E \rar E'$ the natural isogeny, and $Q \in E'(S)$ the unique point of exact order $N$ such that $dQ=\pi(dP)$ and $\langle \frac{N}{d} \cdot P,\,\frac{N}{d} \cdot Q\rangle_{\pi}=\zeta$. Let $P'=P \mod{C[u]}$, $\pi_C: E/C[u] \rar E'' := (E/C[u])/(\langle \frac{N}{d} \cdot P'\rangle)$ be the natural isogeny, and $Q' \in E''(S)$ the unique point of exact order $N$ such that $dQ'=\pi_C(P')$, and $\langle \frac{N}{d} \cdot P' ,\, \frac{N}{d} \cdot Q'\rangle_{\pi_C}=\zeta$. Note that $\pi,\pi_C$ have degree $d$ prime to $m$. 

It is clear that the only point to prove is that $Q'=u'Q \mod{\pi(C)[u]}$, i.e. that $dQ'=dQ \mod{\pi(C)[u]}$ and $\langle \frac{N}{d} \cdot P',\, \frac{N}{d} \cdot (Q \mod{\pi(C)[u]})\rangle_{\pi_C}=\zeta^{u}$. The first equality is clear because $dQ'=\pi_C(dP')$ and $dQ=\pi(dP)$. 

For the second equality, let $p: E \rar E/C[u]$ and $p': E' \rar E'/\pi(C)[u]$ be the natural isogenies of degree $u$, and let $Q''=p'(Q)$. By a direct computation, $\pi_C \circ (\pi_C^{\vee} \circ p'-p\circ \pi^{\vee})=0$, thus $\pi_C^{\vee} \circ p'=p\circ \pi^{\vee}$. Using \cite[(2.8.3),(2.8.4)]{KM}, we check that

\begin{align*}
\langle \frac{N}{d}P',\,\frac{N}{d}Q''\rangle_{\pi_C} &= \langle \frac{N}{d}P',\,\frac{N}{d}Q\rangle_{(p')^{\vee}\circ \pi_C} = \langle \frac{N}{d}Q,\,\frac{N}{d}P'\rangle_{(p')^{\vee}\circ \pi_C}^{-1}\\
&= \langle \frac{N}{d}Q,\,\frac{N}{d}P'\rangle_{p\circ \pi^{\vee}}^{-1} = \langle \frac{N}{d}Q,\,\frac{N}{d}p^{\vee}(P')\rangle_{\pi^{\vee}}^{-1}=\langle u\frac{N}{d}P,\,\frac{N}{d}Q\rangle_{\pi}=\zeta^u,
\end{align*}
whence the conclusion. 
}

\section{Compactified moduli schemes}

We construct compactified moduli schemes attached to some moduli problems, using the material of Appendix \ref{situation-at-infty}. We follow \cite[Chapter 8]{KM} for most of the discussion. As in this source, all base rings are assumed to be regular excellent Noetherian. We will always assume that our moduli problems over $R$ are representable and finite over $\Ell_R$ (by Proposition \ref{base-change-good-property}, these conditions are stable under base change). In particular, by \cite[Proposition 8.2.2]{KM}, the $j$-invariant map of the moduli problems will always be finite (we could also have used Corollary \ref{over-Ell-concrete}).

\defi[normal-near-infinity-mp]{A moduli problem $\P$ over $\Ell_R$ is \emph{normal near infinity} if $\MP$ is normal near infinity in the sense of Definition \ref{normal-near-infinity-mp}. Lemma \ref{equiv-normal-near-infinity} tells us that, under our assumptions, this is equivalent to the definition in \cite[(8.6.2)]{KM}. 

In this case, we define the \emph{compactified moduli scheme} $\CMP$ as the compactification of the $j$-invariant morphism $\MP \rar \mathbb{A}^1_R$ in the sense of Proposition \ref{compactification-functor}. It is endowed with a finite morphism $\overline{j}: \CMP \rar \mathbb{P}^1_R$ extending the usual $j$-invariant. Its \emph{cuspidal subscheme} $\CP$ is the reduced scheme attached to the closed subspace $\{\overline{j}=\infty\}$ of $\CMP$. The proof of Proposition \ref{compactification-functor} shows that this is equivalent to the construction of \cite[(8.6.3)]{KM}, using \cite[Proposition 8.2.2]{KM}. }

\prop[finite-maps-extend]{Let $\P,\P'$ be representable moduli problems over $R$, normal near infinity. Let $f: \MP \rar \MPp{\P'}$ be any finite map of $R$-schemes. Then $f$ extends to a finite map $\CMP \rar \CMPp{\P'}$.}

\defi[normal-near-infinity-relative]{If $R'$ is a $R$-algebra and $\P'$ is a moduli problem over $\Ell_{R'}$, we say that $\P'$ is \emph{normal near infinity over $R$} if there exists a closed subscheme $Z$ of $\mathbb{A}^1_R$, finite over $R$, such that $\MPp{\P'}$ is normal outside $j^{-1}(Z_{R'})$.}

\cor[base-changed-morphisms]{Suppose that $R'$ is a $R$-algebra and that we have moduli problems $\P,\P'$ over $\Ell_R,\Ell_{R'}$ respectively. Suppose that $\P$ is normal near infinity, and that $\P'$ is normal near infinity over $R$ and that we have a morphism of moduli problems $\P'\Rightarrow \P_{R'}$. Then the base change map of moduli schemes $\MPp{\P'} \rar \MP_{R'}$ extends to a morphism of compatified moduli schemes $\CMPp{\P'} \rar \CMP_{R'}$ over $\mathbb{P}^1_{R'}$, and this maps sends $\CPp{\P'}$ to $ (\CP)_{R'}$. }

\defi[formation-commutes-base-change]{Suppose that we can take $\P'=\P_{R'}$ in the Corollary \ref{base-changed-morphisms}. The morphisms $\overline{m}: \CMPp{\P'} \rar \CMP_{R'}$ and $c: \CPp{\P'} \rar (\CP)_{R'}$ constructed by this Corollary are called the \emph{base change morphisms} for $\P$ and $R \rar R'$. We then say that the formation of the compactified moduli scheme (resp. the cuspidal subscheme) for $\P$ \emph{commutes with the base change $R \rar R'$} if $\overline{m}$ (resp. $c$) is an isomorphism.}

\rem{In the situation of Definition \ref{formation-commutes-base-change}, $\overline{m}$ is finite and an isomorphism above $\overline{j}^{-1}(\mathbb{A}^1_{R'})$.}

\cor[finite-extension-mp]{In the situation of the Definition \ref{formation-commutes-base-change}, let $Y$ be a finite $R$-scheme. Any morphism $\MP \rar Y$ of $R$-schemes extends to $\CMP \rar Y$, and the extension $\CMPp{\P'} \rar Y_{R'}$ is the morphism $\CMPp{\P'} \overset{\overline{m}}{\rar}\CMP_{R'} \rar Y_{R'}$.}

\demo{The existence of the extensions is Corollary \ref{extension-to-finite}. The two morphisms agree on the open subscheme $\MPp{\P'}$, $Y_{R'}$ is separated, so they are equal. }

\prop[normal-near-infinity-ascends]{(\cite[Proposition 8.6.7]{KM}) Let $R \rar R'$ be a smooth morphism of $R_0$-algebras and $\P$ be a moduli problem over $\Ell_R$ normal near infinity over $R_0$. Then $\P_{R'}$ is normal near infinity over $R_0$. Moreover, the formation of the compactified moduli scheme and the cuspidal subscheme commutes with base change.}

\demo{Let $Z \subset \mathbb{A}^1_{R_0}$ be a closed subscheme, finite over $R_0$, such that $\MP \backslash j^{-1}(Z_R)$ is normal. Since $\MPp{\P_{R'}} \rar \MP$ is smooth, $\MPp{\P_{R'}}$ is normal outside $j^{-1}(Z_{R'})$, hence $\P_{R'}$ is normal near infinity over $R$ and we have well-defined base change morphisms. Since normalization commutes with smooth base change \cite[Lemma 03GV]{Stacks}, the formation of the compactified moduli scheme commutes with the base change $R \rar R'$. 

Let us now discuss the formation of the cuspidal subscheme. The isomorphism of $\mathbb{P}^1_{R'}$-schemes $\CMPp{\P_{R'}} \rar \CMP_{R'}$ induces a surjective closed immersion $\iota: \CPp{\P_{R'}} \rar (\CP)_{R'}$. This map will be an isomorphism as long as $(\CP)_{R'}$ is reduced. Now, $(\CP)_{R'} \rar \CP$ is a smooth map to a reduced scheme, hence $(\CP)_{R'}$ is reduced, thus $\iota$ is an isomorphism.}

\lem[with-extra-functions]{Let $R \rar R'$ be a finite \'etale morphism of regular excellent Noetherian rings, and $\P$ be a finite representable moduli problem over $\Ell_R$ normal at infinity. Then the finite representable moduli problem $\P'=\P \times [\Sp{R'}]$ is normal at infinity, and we have natural isomorphisms $\CMPp{\P'} \rar \CMP \times_R \Sp{R'}$, $\CPp{\P'} \rar \CP \times_R \Sp{R'}$.}

\demo{First, note that by Proposition \ref{product-rel-rep}, $\P'$ is finite representable with moduli scheme $\MP \times_R \Sp{R'} \rar \Sp{R}$. If $Z \subset \MP$ is a closed subscheme finite over $\Sp{R}$ such that $\MP$ is normal outside $Z$, then $\MPp{\P'}$ is normal outside the closed subscheme $Z_{R'}$ which is finite over $R$. Hence $\P'$ is normal at infinity.
The isomorphism $\MPp{\P'} \rar \MP \times_R \Sp{R'}$ preserves the $j$-invariants, so it is an isomorphism of $\mathbb{P}^1_R$-schemes. Since relative normalization commutes with smooth base change \cite[Lemma 03GV]{Stacks}, the relative normalization of $\MP_{R'}$ in $\mathbb{P}^1_{R'}$ is exactly $\CMP_{R'}$. 

From the commutative diagram

\[
\begin{tikzcd}[ampersand replacement=\&]
\MPp{\P} \times_R \Sp{R'} \arrow{d}{j} \arrow{r}\& \MPp{\P'}\arrow{d}{j}\\
\mathbb{P}^1_{R'} \arrow{r} \& \mathbb{P}^1_R
\end{tikzcd}
\]
 we get a natural morphism $f: \CMP \times_R \Sp{R'} \rar \CMPp{\P'}$, which is an isomorphism outside the cuspidal subschemes. Since $f$ is a morphism of finite $\mathbb{P}^1_R$-schemes, $f$ is finite. Furthermore, $f$ is an isomorphism above $\MPp{\P'}$, which is dense and whose inverse image $\MP \times_R \Sp{R'}$, which is dense. Moreover, all the points of the target (resp. the source) of $f$ in the cuspidal subscheme (resp. in $f^{-1}(\CPp{\P'})$, i.e. in $\CP \times_R \Sp{R'}$) are normal: by definition for the target, by \'etale base change for the source. By Lemma \ref{kinda-zmt}, $f$ is an isomorphism.   

}

\prop[normal-near-infinity-descends]{Let $R \rar R'$ be a faithfully flat (resp. finite faithfully flat) morphism and $\P$ a moduli problem over $\Ell_R$. If $\P_{R'}$ is normal near infinity over some ring $R_0$, where $R$ is a $R_0$-algebra, (resp. normal near infinity), $\P$ is normal near infinity over $R_0$ (resp. $\P$ is normal near infinity). Moreover, if $R'$ is \'etale over $R$, then the formation of the compactified moduli scheme and the cuspidal subscheme for $\P$ commutes with the base change $R \rar R'$. }

\demo{Suppose that $\P_{R'}$ is normal near infinity over $R_0$: there is a finite subscheme $Z \subset \mathbb{A}^1_{R_0}$ such that $\MP \backslash j^{-1}(Z_{R'})$ is normal. Then the map $\MPp{\P_{R'}} \backslash j^{-1}(Z_{R'}) \rar \MP \backslash j^{-1}(Z_R)$ is faithfully flat from a normal scheme, hence its codomain is a normal scheme%
. In particular, $\P$ is normal at infinity over $R_0$. 

Suppose now that $\P_{R'}$ is normal near infinity and $R'$ is finite faithfully flat over $R$. There is a closed subscheme $Z' \subset \mathbb{A}^1_{R'}$, finite over $R'$ (thus affine), such that $\MPp{\P_{R'}} \backslash j^{-1}(Z')$ is normal. Because $R$ is Noetherian, the image $Z$ of $Z'$ in $\mathbb{A}^1_R$ is a closed subscheme, finite over $R$; by definition, $\MPp{\P_{R'}} \backslash j^{-1}(Z_{R'})$ is contained in $\MPp{\P_{R'}} \backslash j^{-1}(Z')$, hence is normal. Thus $\P_{R'}$ is normal at infinity over $R$. The previous point shows that $\P$ is normal at infinity over $R$. 

The last point follows from Proposition \ref{normal-near-infinity-ascends}.   
}

\defi{Let $R'$ be a $R$-algebra. A moduli problem $\P$ on $\Ell_{R'}$ is \emph{smooth at infinity over $R$} if the following properties hold:
\begin{itemize}[noitemsep,label=$-$]
\item $\P$ is normal near infinity.
\item There exists a closed subscheme $Z \subset \mathbb{A}^1_R$, finite over $R$, such that $\CMP \rar \Sp{R'}$ is smooth outside $j^{-1}(Z_{R'})$.  
\item $\CP$ is finite \'etale over $R$.
\end{itemize}
When $R=R'$, we say that $\P$ is smooth at infinity. This is then exactly the notion discussed in \cite[Definition 8.6.5]{KM}. 
}

\rem{Since the $j$-invariant map $\MP \rar \mathbb{A}^1_R$ is finite, this definition for $R'=R$ is equivalent to asking that $\MP \rar \Sp{R}$ is smooth at infinity in the sense of Definition \ref{smooth-at-infty-gnl}. }

\lem[smooth-at-infinity-implies-normal]{Let $R'$ be an $R$-algebra. If $\P$ is a moduli problem on $\Ell_{R'}$ which is smooth at infinity over $R$, then $\P$ is normal near infinity over $R$.}

\demo{This is because $\MP \backslash j^{-1}(Z_{R'}) \rar \Sp{R'}$ is smooth over a regular ring, hence $\MP \backslash j^{-1}(Z_{R'})$ is a regular scheme by \cite[Lemma 07NF]{Stacks}.}

\prop[smooth-at-infinity-base-change]{(\cite[Proposition 8.6.6]{KM}) Let $R_0 \rar R \rar R'$ be morphisms of rings. If $\P$ is a moduli problem on $\Ell_{R}$ which is smooth at infinity over $R_0$, then $\P_{R'}$ is smooth at infinity over $R_0$, and the formation of the compactified moduli scheme and the cuspidal subscheme for $\P$ commutes with the base change $R \rar R'$. Moreover, the cuspidal subscheme of $\P$ is an effective Cartier divisor in $\CMP$. }

\demo{This is an immediate consequence of Proposition \ref{smooth-at-infty-base-change}.}

\prop[smooth-at-infinity-descends]{Let $R_0 \rar R \rar R'$ be morphisms of rings with $R \rar R'$ faithfully flat smooth. Let $\P$ be a moduli problem on $\Ell_{R}$ such that $\P_{R'}$ is smooth at infinity over $R_0$. Then $\P$ is smooth at infinity over $R_0$ and the formation of the compactified moduli scheme and the cuspidal subscheme for $\P$ commutes with the base change $R \rar R'$.}

\demo{Let $Z \subset \mathbb{A}^1_{R_0}$ be a closed subscheme, finite over $R_0$, such that $\MPp{\P_{R'}}$ is smooth over $\Sp{R'}$ outside the $j$-invariants in $Z_{R'}$. Now, $\P$ is normal near infinity over $R_0$ by Proposition \ref{normal-near-infinity-descends}, and the formation of the compactified moduli scheme and the cuspidal subscheme commutes with base change by Proposition \ref{normal-near-infinity-ascends}. In particular, after making the fpqc base change $R \rar R'$, $\CMP_{R} \rar \Sp{R}$ is smooth of relative dimension one outside the $j$-invariants in $Z_R$, therefore $\CMP \rar \Sp{R}$ is smooth of relative dimension one outside the $j$-invariants in $Z_R$. 

We previously said that the isomorphism of compactified moduli schemes induces an isomorphism of cuspidal subschemes $\iota: C_{\P_{R'}} \rar (C_{\P})_{R'}$. Hence $C_{\P} \rar \Sp{R}$ becomes finite \'etale after the faithfully flat base change $R \rar R'$, hence $C_{\P} \rar \Sp{R}$ is finite \'etale.}

\prop[normal-near-inf-level]{Let $R$ be a regular excellent Noetherian ring and let $\P$ be a finite locally free representable moduli problem on $\Ell_R$ endowed with a level $N$ structure, for some $N \geq 1$, and assume that $\P$ is normal near infinity. Then the action of $(\Z/N\Z)^{\times}$ on $\MP$ extends to an action of $(\Z/N\Z)^{\times}$ on $\CMP$. This action is natural with respect to $\P$ (for morphisms of finite locally free representable moduli problems on $\Ell_R$ normal near infinity and endowed with a level $N$ structure).}  

\demo{This is a direct application of Proposition \ref{compactification-functor}.}

\prop[degn-AL-infty]{Let $N,m \geq 1$ be coprime integers, $R$ be a regular excellent Noetherian ring, and $\P$ be a finite locally free representable moduli problem on $\Ell_R$ endowed with a level $N$ structure, such that $\P \times [\Gamma_0(m)]_R$ is normal near infinity. Then, for every divisor $d \geq 1$ of $m$, the moduli problem $\P \times [\Gamma_0(t)]_R$ is normal near infinity. Moreover, 
\begin{itemize}[noitemsep,label=$-$]
\item for any $d,t \geq 1$ such that $dt \mid m$, there is a finite map $\CMPp{\P,[\Gamma_0(m)]_R} \rar \CMPp{\P,[\Gamma_0(t)]}$ of $R$-schemes extending the morphism $D_{d,t}$ of Definition \ref{degeneracies}. We denote this extension by $D_{d,t}$.  
\item for any $d \geq 1$ dividing $m$ and coprime to $m/d$, there is an automorphism of $R$-schemes of $\CMPp{\P,[\Gamma_0(m)]_R}$ extending $w_d$. We denote this extension by $w_d$.
\item the formation of $D_{d,t}$ and $w_d$ is natural with respect to $\P$ (for morphisms of representable moduli problems of level $N$ such that their product with $[\Gamma_0(m)]_R$ is normal near infinity). 
\item the relations between the various $w_d$ (for $d \geq 1$ dividing $m$ and coprime to $m/d$) and $D_{s,t}$ (for various $s,t \geq 1$ such that $st \mid m$) stated in Corollary \ref{elem-degeneracies} and Proposition \ref{AL-auto} still hold. 
\end{itemize}
}

\demo{For any $t \geq 1$ dividing $m$, $D_{1,t}: \MPp{\P,[\Gamma_0(m)]_R} \rar \MPp{\P,[\Gamma_0(t)]_R}$ is finite flat surjective. The rest is a direct application of the construction of Proposition \ref{compactification-functor}, since these maps are defined on the (non-compactified) moduli schemes.}

\prop[flatness-at-cusps]{In the situation of Proposition \ref{finite-maps-extend}, assume furthermore that $f$ is flat, that $\P,\P'$ are flat near infinity\footnote{in the sense that there exists an open subscheme $U \subset \mathbb{A}^1_R$ whose complement is finite over $\Sp{R}$ and such that $\P,\P'$ are flat over $\Ell_R$ above $U$} over $\Ell_R$, and that $\CMP$, $\CMPp{\P'}$ are respectively Cohen-Macaulay and regular at every point of the cuspidal subscheme. Then the compactifications $\overline{f}$ of $f$ and of the $j$-invariants are finite locally free. These hypotheses are satisfied if $f$ is flat and $\P,\P'$ are smooth at infinity. 
}

\demo{After taking a direct factor of $R$ (since it is Noetherian), we can assume that $R$ has no nontrivial idempotents. Since it is regular, $R$ is a domain. 

Let $x \in \CMP$ be contained in the cuspidal subscheme, $y$ be its image in $\mathbb{P}^1_R$ and $z$ be its image in $\Sp{R}$. We show that $\dim{\OO_{\CMP,x}}=\dim{\OO_{\mathbb{P}^1_R,y}}=\dim{\OO_{\Sp{R},z}}+1$. This implies in particular that $\dim{\OO_{\CMP,x}}=\dim{\OO_{\CMPp{\P'},f(x)}}$. The morphisms $\overline{f}^{\sharp}: \OO_{\mathbb{P}^1_R,y} \rar \OO_{\CMP,x}$ or $\overline{j}^{\sharp}: \OO_{\CMPp{\P'},f(x)} \rar \OO_{\CMP,x}$ are local, and they are localizations at a maximal ideal of finite morphisms (by Proposition \ref{compactification-functor}); and since their base rings are regular, while the target ring is Cohen-Macaulay, the morphisms are flat by miracle flatness \cite[Lemma 00R4]{Stacks}. 

Let us prove that $\dim{\OO_{\mathbb{P}^1_R,y}}=\dim{\OO_{\Sp{R},z}}+1$. By \cite[Lemma 00ON]{Stacks}, we are reduced to the case where $R$ is a field. Since $y$ is above the infinity section, it is a closed point of $\mathbb{P}^1$ and the conclusion follows.

Next, we need to prove that $\dim{\OO_{\mathbb{P}^1_R,y}} = \dim{\OO_{\CMP,x}}$. Let $U \subset \mathbb{P}^1_R$ an affine open subscheme containing the cuspidal subscheme such that $j$ is flat above $U \cap \mathbb{A}^1_R$ and $\OO(j^{-1}(U))$ is normal. Then $\OO(U)$ is a normal domain, so, since $\OO(U) \rar \OO(j^{-1}(U))$ is finite, it is enough to show by \cite[Lemma 00ON]{Stacks} that this map satisfies going down. Since $\OO(j^{-1}(U))$ is normal Noetherian, it is a finite product of domains, and, by \cite[Proposition 00H8]{Stacks}, the going down property holds as soon as $\OO(U)$ injects in each direct factor of $\OO(j^{-1}(U))$. Now, if $C \subset j^{-1}(U)$ is a closed open subscheme, then $j: C \rar U$ is still finite flat, and $U$ is irreducible, so $j(C)=U$. Hence $\OO(C)$ is faithfully flat over $\OO(U)$, and $\OO(U) \rar \OO(C)$ is injective. 
}

\cor[flat-deg-infty]{The $D_{d,t}$ defined in Proposition \ref{degn-AL-infty} are finite locally free if for every point $x \in X:=\CMPp{\P,[\Gamma_0(k)]_R}$ (for any divisor $k \geq 1$ of $m$) lying in the cuspidal subscheme, the local ring $\OO_{X,x}$ is regular, and if $X \rar \mathbb{P}^1_R$ is flat above some open subset $U \subset \mathbb{A}^1_R$ with $\mathbb{A}^1_R \backslash U$ finite over $R$. These conditions are satisfied if every $\P \times [\Gamma_0(k)]_R$ is smooth at infinity. }

To apply Proposition \ref{degn-AL-infty}, we thus need concrete criteria in order to decide when the moduli problems $\P \times [\Gamma_0(m)]_R$ are normal near infinity.

\prop[struct-gamma0m-smooth]{Let $\P$ be a finite representable moduli problem over $\Ell_R$, where $R$ is a Noetherian ring. Let $m_0,m_1 \geq 1$ be coprime integers such that $m_0$ is square-free and $m_1 \in R^{\times}$, with $m=m_0m_1$. Let $U \subset \mathbb{A}^1_R$ be an open subset such that its complement $Z$ is finite over $\Sp{R}$ and $\P_{|U}$ is \'etale over $U$. Then:
\begin{itemize}[noitemsep,label=$-$]
\item $\MPp{\P,[\Gamma_0(m)]_R]} \rar \Sp{R}$ is Cohen-Macaulay of relative dimension one with geometrically reduced fibres at every point of $j^{-1}(U)$, and smooth at every point which is not supersingular in residue characteristic $p \mid m_0$.
\item Assume furthermore that $R$ is regular excellent Noetherian of generic characteristic zero. Then $\P \times_R [\Gamma_0(m)]_R$ is smooth at infinity, the same statement as above holds for its compactified moduli scheme, and $\CMPp{\P,[\Gamma_0(m)]_R} \rar \CMP$ is finite locally free of constant rank. 
\item If furthermore $R'$ is any regular excellent Noetherian $R$-algebra, the conclusions of the previous point hold over $R'$. 
\end{itemize}}

\demo{The last point follows from the second one using Proposition \ref{smooth-at-infinity-base-change}. 

We know that the map $\MPp{\P,[\Gamma_0(m)]_R} \rar \MP$ is finite locally free of constant rank (\'etale over $\Z[1/m]$) and the same argument as Corollary \ref{over-Ell-concrete} shows that $j^{-1}(U) \subset \MPp{\P} \rar \Sp{R}$ is smooth of relative dimension one. Thus, $j^{-1}(U) \subset \MPp{\P,[\Gamma_0(m)]_R} \rar \Sp{R}$ is Cohen-Macaulay of relative dimension one and smooth above $R[1/m]$. 

To prove the claimed smoothness (and reducedness of the geometric fibres), we may therefore assume that $R$ is an algebraically closed field $k$ of characteristic $p \mid m_0$. After replacing $\P$ with $\P \times_R [\Gamma_0(m_0m_1/p)]_R$, we may assume that $(m_0,m_1)=(p,1)$. 

Fix some prime $\ell \nmid 2m$, then $\MPp{[\Gamma(\ell)]_k,[\Gamma_0(p)]_k,\P} \rar \MPp{\P,[\Gamma_0(p)]_k}$ is finite \'etale of constant positive rank, hence it is enough to prove the statement for the moduli problem $[\Gamma(\ell)]_k\times [\Gamma_0(p)] \times \P$. The map $\MPp{[\Gamma(\ell)]_k,[\Gamma_0(p)]_k,\P} \rar \MPp{[\Gamma(\ell)]_k,[\Gamma_0(p)]_k}$ is finite \'etale above $j^{-1}(U)$, so we are reduced to the case $\P=[\Gamma(\ell)]_k$. This then follows from \cite[Theorem 10.10.3 (5)]{KM} (with the same $p$, $N=\ell$, $\Gamma \leq \GL{\Z/p\ell\Z}$ being the subgroup of upper-triangular matrices that are trivial mod $\ell$), whose assumptions ((a), (b) and balancedness\footnote{This is the reason why $m_0$ had to be square-free.}) can be directly checked, and we come back from what \emph{loc. cit.} calls the ``canonical problems'' to our setting using \cite[Proposition 9.1.7]{KM}.

Note that, to use this statement, we need to base change from $R_0=\Z[1/\ell]$ (we are dealing with a representable moduli problem, so there is no issue) and check that the moduli problem $[\Gamma_0(p)]_{R_0}$ is indeed the quotient of $[\Gamma(p)]_{R_0}$ by the subgroup of \emph{lower}-triangular matrices (see Remark \ref{left-vs-right-actions}), which is \cite[Theorem 7.4.2]{KM}.

Now, we discuss the second point. Assume that $R$ is regular excellent Noetherian, and that $\P \times [\Gamma_0(m)]_R$ is smooth at infinity. Then the required smoothness for the compactified moduli scheme follows directly from the first point. The morphism $f: \CMPp{\P,[\Gamma_0(m)]_R} \rar \CMP$ is finite by assumption, flat away from the cusps because so is the moduli problem $[\Gamma_0(p)]_R$. Away from some finite $R$-subscheme containing $j^{-1}(Z)$ and the supersingular points in residue characteristic $p \mid m_0$ (using \cite[Lemma 12.5.4]{KM}), $f$ is a finite map of smooth $R$-schemes of relative dimension one, so it is flat by miracle flatness. 

Thus only the smoothness at infinity remains to prove. We assume that $R$ is regular excellent Noetherian, and (after possibly shrinking it), that $U$ is principal defined by a monic polynomial. If $m_0=1$, then the statement follows from \cite[Theorem 8.6.8]{KM}, so we assume $m_0 > 1$. Let $Y_0(m)^R$ be the coarse moduli scheme attached to $[\Gamma_0(m)]_R$, $Y_0(m):=Y_0(m)^{\Z[1/m_1]}$. By \cite[Theorem 5.1.1, Lemma 8.1.2]{KM}, $Y_0(m)$ is normal. By \cite[Proposition 8.1.6]{KM}, we have a canonical morphism of $\mathbb{A}^1_R$-schemes $B: Y_0(m)^R \rar Y_0(m) \times_{\Z[1/m_1]} \Sp{R}$. Let $V$ be a principal open subscheme of $\mathbb{A}^1_{\Z[1/m_1]}$ defined by a monic polynomial which is a multiple of $j(j-1728)$, and such that $V$ does not contain the supersingular $j$-invariants in residue characteristic $p \mid m_0$ (see \cite[Lemma 12.5.4]{KM}). By \cite[(8.5.4)]{KM}, $B$ is an isomorphism above $j^{-1}(V)$. 

The open subspace $V'=U \cap V$ of $\mathbb{A}^1_R$ is defined by a monic polynomial, hence it is the complement of a closed subscheme $Z'$ finite over $\Sp{R}$. Then $\pi: \MPp{\P,[\Gamma_0(m)]_R} \rar Y_0(m)^R$ is finite, and a morphism of smooth $R$-schemes of relative dimension one above $j^{-1}(V')$, hence by miracle flatness it is flat. 

Let us prove that $Y_0(m)$ is smooth at infinity in the sense of Definition \ref{smooth-at-infty-gnl}, and let $X_0(m)$ be its compactification. We apply \cite[Theorem 10.10.3 (5)]{KM} to each prime $p \mid m_0$, with $N=m_0m_1/p$, $\Gamma_1 \leq \GL{\Z/N\Z}$, $\Gamma_2 \leq \GL{\Z/N\Z}$ the subgroups consisting of upper-triangular matrices: the assumptions (a) and (b) are trivially satisfied and one checks directly that $\Gamma_1$ is balanced. By Proposition \ref{smooth-at-infty-base-change}, for any regular excellent Noetherian $\Z[m_1^{-1}]$-algebra $S$, that $Y_0(m)^S$ is smooth at infinity, and $X_0(m)^S := \overline{Y_0(m)^S} \rar X_0(m)_S$ is an isomorphism away from the (inverse images under finite maps) of the finite $R$-scheme $j^{-1}(\mathbb{A}^1_R \backslash V_R)$. 

Let us temporarily assume that $\pi$ is \'etale above $j^{-1}(V)$. Since $R$ has generic residue characteristic zero, we can then apply Proposition \ref{abhyankar-compactification} to $\MPp{\P,[\Gamma_0(m)]_R} \rar X_0(m)^R$, and the conclusion would follow. 

Thus, all we need to do is prove that $\pi$ is \'etale above $j^{-1}(V')$. Let $k$ be an algebraically closed field and a $R$-algebra. It is enough to show that the fibres of $\pi \times_R \Sp{k}$ above $j^{-1}(V')_k$ are \'etale. Let $BC: Y_0(m)^k \rar Y_0(m)^R_k$ be the canonical base change morphism, which is an isomorphism above $j^{-1}(V')$. 

Since, above $j^{-1}(V)$, the formation of the coarse moduli scheme for $[\Gamma_0(m)]$ commutes with base change, and since $\pi$ is flat above $j^{-1}(V')$, we may thus assume that $R$ is an algebraically closed field $k$ of characteristic $c$. Let $U' \subset \mathbb{A}^1_k$ be the affine open subset on which $j(j-1728)$ is invertible. If $c \neq 2$, consider the moduli problem $[U',1]$ of \cite[Theorem 8.4.9]{KM}, which yields the following commutative diagram of $U'$-schemes:

\[
\begin{tikzcd}[ampersand replacement=\&]
\MPp{\P,[\Gamma_0(m)]_k,[U',1]} \arrow{r}{\alpha}\arrow{d}{\beta}\& \MPp{\P,[\Gamma_0(m)]_k} \times_{\mathbb{A}^1_k} U'\arrow{d}{\pi'}\\
\MPp{[\Gamma_0(m)]_k,[U',1]} \arrow{r}{\delta} \& Y_0(m)^k\times_{\mathbb{A}^1_k} U'
\end{tikzcd}
\]

By \cite[Corollary 8.4.10]{KM}, $\alpha,\delta$ are finite \'etale. By assumption, $\beta$ is finite \'etale above $j^{-1}(U \cap U')$. Thus $\pi'$ is \'etale above $j^{-1}(U \cap U')$: since $\pi \otimes_R \Sp{k}=BC \circ \pi'$, and $BC$ is an isomorphism above $j^{-1}(V') \subset j^{-1}(U \cap U')$, $\pi \times_R \Sp{k}$ is \'etale above $j^{-1}(V')_k$. 

When $c = 2$, the same argument applies by replacing the moduli problem $[U',1]$ with the moduli problem $[\mu_4]$ of \cite[Theorem 8.4.11]{KM}, and replacing the references to \cite[Theorem 8.4.9, Corollary 8.4.10]{KM} with references to \cite[Theorem 8.4.11, Corollary 8.4.12]{KM}. 
 }

\cor[smooth-at-infinity-gamma0m]{Let $\P$ be a finite \'etale representable moduli problem over $\Ell_R$, for some regular excellent Noetherian ring $R$ of generic characteristic zero. Then, for all coprime $m_0,m_1 \geq 1$ with $m_0$ square-free and $m_1 \in R^{\times}$, $\P \times_R [\Gamma_0(m_0m_1)]_R$ is smooth at infinity. }

\prop[bad-degn-AL-infty]{Let $N \geq 1$ be an integer, $r \mid N$, and $R$ be a regular excellent Noetherian $\Z[1/N]$-algebra, and $\P$ be a finite locally free relatively representable (representable if $N \leq 3$) moduli problem on $\Ell_R$ with level prime to $N$, such that $\P \times [\Gamma_1'(N,r)]_R$ is normal near infinity. Then, for every divisor $d \geq 1$ of $r$, the moduli problem $\P \times [\Gamma_1'(N,d)]_R$ is normal near infinity. Moreover, 
\begin{itemize}[noitemsep,label=$-$]
\item for any $s,t \geq 1$ such that $st \mid r$, the morphism $D'_{s,t}$ of Definition \ref{degeneracies} extends to a finite map $D'_{s,t} : \CMPp{\P,[\Gamma_1'(N,r)]_R} \rar \CMPp{\P,[\Gamma_1'(N,t)]}$ of $R$-schemes.  
\item $D'_{s,t}$ is finite locally free if for every $d \mid r$ and every $x \in X:=\CMPp{\P,[\Gamma_1'(N,d)]_R}$ contained in the cuspidal subscheme, $\OO_{X,x}$ is regular. This is in particular the case if every $[\P] \times [\Gamma_1'(N,d)]_R$ is smooth near infinity. 
\item for any $d \geq 1$ dividing $N$ and coprime to $N/d$, there is an automorphism of $R$-schemes of $\CMPp{\P,[\Gamma_1(N)]_R}$ extending $w'_d$. We denote this extension by $w'_d$.
\item the formation of $D'_{s,t}$ and $w'_d$ is natural with respect to $\P$ (for morphisms of representable moduli problems of level $N$ such that their product with $[\Gamma_1'(N,r)]_R$ or $[\Gamma_1(N)]_R$ is normal near infinity). 
\item the relations between the various $w'_d$ and $D'_{s,t}$ stated in Lemmas \ref{bad-degn-composition}, \ref{bad-AL-ZnZ}, \ref{bad-AL-cyclo} (with the suitable extension of scalars) and \ref{bad-AL-composition} still hold. 
\end{itemize}
}

\demo{The proof is analogous to the proof of Proposition \ref{degn-AL-infty} (and Corollary \ref{flat-deg-infty}) using the results of Section \ref{gamma1m-structure}. }

\section{Application to our moduli problems}
\label{application-moduli}

\defi[torsion-group]{Let $N \geq 3$ be an integer, $S$ be a $\Z[1/N]$-scheme. A \emph{$N$-torsion group over $S$} is a finite \'etale group $G$ over $S$ which is \'etale-locally isomorphic to $(\Z/N\Z)^{\oplus 2}$.}

\defi[polarized-torsion-group]{A \emph{Weil pairing} for a $N$-torsion group $G$ is a bilinear alternating pairing $G \times_R G \rar \mu_N$ which is fibre-wise perfect. By Cartier duality \cite[\S 4]{Shatz-Groups}, this alternating pairing induces an morphism $w$ from $G$ to its Cartier dual. Since both groups are finite \'etale over $R$, $w$ is finite \'etale by \cite[Lemma 02GW]{Stacks}. Since $w$ is a fibrewise isomorphism, $w$ is an isomorphism. We call the datum of a $N$-torsion group over $R$ endowed with a Weil pairing a \emph{polarized $N$-torsion group}.  }

The moduli problems that we will consider are the following ones:

\defi{Let $G$ be a $N$-torsion group over a $\Z[1/N]$-algebra $R$. The moduli problem $\P_G$ on $\Ell_R$ is defined as follows: 
\begin{itemize}[noitemsep,label=$-$]
\item let $S$ be a $R$-scheme and $E/S$ be an elliptic curve, then $\P_G(E/S)$ is the set of isomorphisms $\iota: G \times_R S \rar E[N]$ of $S$-group schemes, 
\item for a morphism $F/T \rar E/S$ in the category $\Ell_R$, the morphism $\P_G(E/S) \rar \P_G(F/T)$ is induced by the isomorphism $F[N] \rar E[N] \times_S T$. 
\end{itemize}
It is endowed with a natural structure of moduli problem of level $N$. In particular, for any $m \geq 1$ coprime to $N$, the moduli problem $\P_G(\Gamma_0(m)) := \P_G \times [\Gamma_0(m)]_R$ is endowed with a natural level structure of level $Nm$. 
}

\rem{The moduli problem $[\Gamma(N)]_R$ over a $\Z[1/N]$-algebra $R$ is exactly the moduli problem $\P_G$ attached to the constant group scheme $(\Z/N\Z)^{\oplus 2}_R$.}

\defi{Let $G$ be a polarized $N$-torsion group with Weil pairing $W$. Given an isomorphism of $S$-group schemes $\varphi: G_S \rar E[N]$, the map $G_S \times G_S \overset{\varphi\times\varphi}{\rar} E[N] \times E[N] \overset{\mrm{We}}{\rar} \mu_N$ is an alternate pairing; therefore it can be written as $\alpha \circ W$, where $\alpha$ is an endomorphism of the $S$-group scheme $\mu_N$. In particular, $\alpha$ is a locally constant function $S \rar (\Z/N\Z)^{\times}$, which we call the determinant of $\varphi$ and denote $\det{\varphi}$.

For $u \in (\Z/N\Z)^{\times}$ (resp. and $m \geq 1$ coprime to $N$), let $\P_{G,u}$ (resp. $\P_{G,u}(\Gamma_0(m))$) denote the subfunctor of $\P_G$ (resp. $\P_G(\Gamma_0(m))$) parametrizing isomorphisms $\iota$ (resp. couples $(\iota,C)$) with $\det{\iota}=u$. }

\prop[PG-elementary]{Let $G$ be a $N$-torsion group over a $\Z[1/N]$-algebra $R$. Then the moduli problem $\P_G$ satisfies the following properties: 
\begin{enumerate}[noitemsep,label=$(\roman*)$]
\item\label{PG-elem-1} It is representable and finite \'etale of constant rank $r_N \geq 1$ over $\Ell_R$, by an elliptic curve $\mathcal{E}_G(N)/Y_G(N)$, with $Y_G(N) \rar \Sp{R}$ smooth affine of relative dimension one,
\item\label{PG-elem-2} Its $j$-invariant is finite locally free of constant rank $r'_N \geq 1$, \'etale on the locus where $j(j-1728)$ is invertible,
\item\label{PG-elem-3} If $R$ is regular excellent Noetherian, then $\P_G$ is smooth at infinity over $\Z[1/N]$: the formation of its compactified moduli scheme and cuspidal subschemes commutes with any base change,
\item\label{PG-elem-4} If $R$ is regular excellent Noetherian, the compactified moduli scheme $X_G(N)$ attached to $\P_G$ is smooth proper over $\Sp{R}$ of relative dimension one,
\item\label{PG-elem-5} If $R$ is regular excellent Noetherian, $X_G(N)$ is a regular scheme and the $j$-invariant map is finite locally free of rank $r'_N$, 
\item\label{PG-elem-6} The level $N$ structure on $\P_G$ defines a natural action of $(\Z/N\Z)^{\times}$ on $Y_G(N)$ which preserves the $j$-invariant. When $R$ is regular excellent Noetherian, it extends to an action of $(\Z/N\Z)^{\times}$ on $X_G(N)$ (still denoted $[\cdot]$) whose formation commutes with base change.
\end{enumerate} }

\demo{When $G$ is constant and $R=\Z[1/N]$, \ref{PG-elem-1}-\ref{PG-elem-4} (apart from the statement about base change) follow from Corollary \ref{over-Ell-concrete} and \cite[Theorem 5.1.1, Proposition 8.6.8]{KM}. Since $X_G(N)$ is smooth proper over $\Z[1/N]$ of relative dimension one, it is regular (with every nonempty open subset of dimension two), hence Cohen-Macaulay, and the $j$-invariant $X_G(N) \rar \mathbb{P}^1_{\Z[1/N]}$ is finite from a Cohen-Macaulay scheme to a regular scheme (every nonempty open subscheme of either having dimension $2$), hence is flat by miracle flatness. Finally, \ref{PG-elem-6} (apart from the commutation with base change) follows from the results on level structures and Proposition \ref{compactification-functor}. 

When $G$ is constant over a different base, the results (including the base change statements) follow from Proposition \ref{base-change-good-property}, Corollary \ref{over-Ell-concrete}, and Proposition \ref{smooth-at-infinity-base-change}. 

For general $G$, there is a faithfully flat \'etale map $R \rar R'$ such that $G_{R'}$ is constant. Since $\P(N)_{R'}$ is representable and affine over $\Ell_{R'}$, and since $T \longmapsto \P_G(E_T/T)$ is a fpqc sheaf on $\mathbf{Sch}_S$ for all objects $E/S$ of $\Ell_R$, Proposition \ref{representability-descends} proves that $\P_G$ is representable. By Proposition \ref{base-change-good-property}, $\P_G$ is finite \'etale of rank $r_N$ over $\Ell_R$. Thus the moduli scheme is smooth of relative dimension one over $R$, with finite flat $j$-invariant of rank $r'_N$ by Corollary \ref{over-Ell-concrete} (and base changing to $R'$ for the value of the rank). The $j$-invariant map $Y_G(N) \rar \mathbb{A}^1_R$ is then \'etale above the invertible locus of $j(j-1728)$ by \cite[Corollary 8.4.5]{KM}.

We use Proposition \ref{smooth-at-infinity-descends} (with $R_0=\Z[1/N]$ and the base change $R \rar R'$) to show that $\P_G$ is smooth at infinity over $\Z[1/N]$ and that the formation of its compactified moduli scheme and cuspidal subscheme commutes with the base change $R \rar R'$ (and thus with any base change by Proposition \ref{smooth-at-infinity-base-change}). Thus the following diagram is Cartesian:
\[
\begin{tikzcd}[ampersand replacement=\&]
X_{G_{R'}}(N) \arrow{r}{j'}\arrow{d}\& \mathbb{P}^1_{R'}\arrow{r}{\sigma'}\arrow{d}\& \Sp{R'}\arrow{d}\\
X_G(N) \arrow{r}{j}\& \mathbb{P}^1_R \arrow{r}{\sigma}\&\Sp{R}
\end{tikzcd}
\]

Since $j'$ is finite flat and \'etale away from the cusps and the vanishing locus of $j'(j'-1728)$, and $\sigma' \circ j'$ is smooth proper of relative dimension one, by fpqc descent %
 $j$ is finite flat and \'etale away from the cusps and the vanishing locus of $j(j-1728)$, and $\sigma\circ j$ is proper smooth of relative dimension one. 

\ref{PG-elem-6} follows from the fact that $Y_G(N) \rar X_G(N)$ is an epimorphism (for separated $R$-schemes), the fact that $\P$ has a natural level $N$ structure, and the fact that the formation of the compactified moduli scheme for $\P_G$, as well as the cuspidal subscheme, commutes with base change. }

\cor{If $\P$ is an \'etale representable moduli problem over $\Ell_R$ (where $R$ is any ring), then its $j$-invariant $\MP \rar \mathbb{A}^1_R$ is \'etale on the locus where $j(j-1728)$ is invertible.}

\demo{This statement is Zariski-local over $\Sp{R}$, so we may assume that some odd prime $\ell$ is invertible in $R$. Consider the same commuting diagram as in Proposition \ref{over-Ell-implies-over-R}:
\[
\begin{tikzcd}[ampersand replacement=\&]
\MPp{\P,[\Gamma(\ell))]_R} \arrow{r}{f_1}\arrow{d}{\pi_1} \& \MPp{[\Gamma(\ell)]_R}\arrow{d}{j_{\ell}}\\
\MP \arrow{r}{j_{\P}} \&  \mathbb{A}^1_R\&
\end{tikzcd}
\]
We just proved that $j_{\ell}$ was \'etale above the locus where $j(j-1728)$ is invertible, and $f_1$ is \'etale because $\P$ is. Since $\pi_1$ is finite \'etale, $j_{\P}$ is \'etale-locally \'etale above the locus where $j(j-1728)$ is invertible, hence it is \'etale above this locus by \cite[Lemma 036W]{Stacks}.}

\nott{We denote by $\mathcal{E}(N)/Y(N)$ the elliptic curve $\mathcal{E}_{(\Z/N\Z)^{\oplus 2}}(N)/Y_{(\Z/N\Z)^{\oplus 2}}(N)$. The action of $\GL{\Z/N\Z}$ on the right on $(\Z/N\Z)^{\oplus 2}$ by $\begin{pmatrix} a \\b \end{pmatrix} \mid M = M^T \begin{pmatrix} a \\b\end{pmatrix}$ induces a left action of $\GL{\Z/N\Z}$ by automorphisms on $\P(N)$, hence on $Y(N)$. 

Let $S$ be any $\Z[1/N]$-scheme and $G$ be a $N$-torsion group over $S$. The datum of an isomorphism $f: (\Z/N\Z)^{\oplus 2} \rar G$ of group schemes over $S$ is the same (by considering the image under $f$ of the two basis sections) as the datum of two sections $A,B \in G(S)$ such that $(a,b) \in (\Z/N\Z)^{\oplus 2} \longmapsto aA+bB \in G(S)$ is injective (thus an isomorphism). Under this identification, for any $S$-scheme, $Y(N)(S)$ fonctorially identifies with the set of $S$-isomorphism classes of triples $(E,P,Q)$ where $E/S$ is an elliptic curve, and $(P,Q)$ is a $\Z/N\Z$-basis of $E[N](S)$. Under this identification, a matrix $\begin{pmatrix}a & b\\c & d\end{pmatrix} \in \GL{\Z/N\Z}$ maps a $S$-point $(E,P,Q)$ of $Y(N)$ to $(E,aP+bQ,cP+dQ)$.    

In particular, this procedure identifies $[a]$ with the action of the matrix $aI_2 \in \GL{\Z/N\Z}$. }

\rem[left-vs-right-actions]{In \cite{KM} (for instance (10.3)), the authors let $\GL{\Z/N\Z} \cong \mrm{Aut}((\Z/N\Z)^{\oplus 2})$ act on $X(N)$ on the \emph{right}. Our definition of the action (which is on the \emph{left}) of $g \in \GL{\Z/N\Z}$ is exactly the action of the transpose $g^T$ of $g$ in \emph{loc.cit.}. This will become important when we will write down an explicit description of the cuspidal subscheme in Section \ref{cuspidal-subscheme}. 
} 

\prop[XN-Weil]{Let $(P^{univ},Q^{univ})$ be the basis of $\mathcal{E}(N)[N]$ as a finite group scheme over $Y(N)$. Then $\mrm{We} = \langle P^{univ},\,Q^{univ}\rangle_{\mathcal{E}(N)}$ defines a map $\mrm{We}: Y(N) \rar \Sp{\OO_N}$. It extends to a smooth proper map $\mrm{We}: X(N) \rar \Sp{\OO_N}$ with geometrically connected fibres. In particular, $Y(N),X(N)$ are connected, and the geometric fibres of $X(N) \rar \Sp{\Z[1/N]}$ (as well as those of its restriction to $Y(N)$) have exactly $\varphi(N)$ connected components. 

Furthermore, the natural left action of $\GL{\Z/N\Z}$ on $Y(N)$ extends uniquely to the compactified moduli scheme $X(N)$, preserves the $j$-invariant, and the following diagram commutes, for any $M \in \GL{\Z/N\Z}$:
 \[
\begin{tikzcd}[ampersand replacement=\&]
X(N)  \arrow{r}{M}\arrow{d}{\mrm{We}}\& X(N)\arrow{d}{\mrm{We}}\\
\Sp{\OO_N} \arrow{r}{\underline{\det{M}}} \& \Sp{\OO_N}  
\end{tikzcd}
\]
}

\demo{First, $\Sp{\OO_N} \rar \Sp{\Z[1/N]}$ is finite \'etale. So the extension of $\mrm{We}$ follows from Corollary \ref{extension-to-finite}. The extension of the action of $\GL{\Z/N\Z}$ follows from the previous results and Proposition \ref{compactification-functor}. Since $X(N) \rar \Sp{\Z[1/N]}$ is smooth proper of relative dimension one, and $\OO_N$ is finite \'etale over $\Z[1/N]$, $\mrm{We}: X(N) \rar \Sp{\OO_N}$ is smooth proper of relative dimension one. 

To check that the extension of the action of $\GL{\Z/N\Z}$ preserves the $j$-invariant, it is enough to check it on $Y(N)$, where it is obvious. Checking that the diagram commutes can also be done after restricting to $Y(N)$ (by Proposition \ref{compactification-functor}), where it is an easy computation. 

To discuss the geometric fibres of $\mrm{We}: X(N) \rar \Sp{\OO_N}$, we see that $\mrm{We}: Y(N) \rar \Sp{\OO_N}$ is the scheme representing the moduli problem $[\Gamma(N)]^{\Z[\zeta_N]-can}_{\OO_N}$  (see \cite[Propositions 9.3.1, (9.4.3.1)]{KM}), and apply \cite[Corollary 10.9.2]{KM}\footnote{More directly, we can directly follow the argument given in \cite{DeRa}: by considering the complex uniformization, the base change to $\C$ of $\mrm{We}: Y(N) \rar \Sp{\OO_N}$ is connected, whence $X(N)$ is connected, and it follows that $\OO(X(N)) \otimes_{\OO_N} \C \simeq \C$. Since $\OO(X(N))$ is a finite, normal, integral $\OO_N$-algebra, it is exactly $\OO_N$, and $\OO_{\Sp{\OO_N}} \rar \mrm{We}_{\ast}\OO_{X(N)}$ is an isomorphism. It follows by cohomology and base change that the formation of $\mrm{We}_{\ast}\OO_{X(N)}$ commutes with arbitrary base change, QED.}.
}

\lem[open-comp]{Let $X$ be a scheme of finite type over a field $k$, Cohen-Macaulay of relative dimension one. Let $Z \subset X$ be a closed subscheme consisting of finitely many regular points of $X$. Then $C \in \pi_0(X) \longmapsto C \backslash Z \in \pi_0(X \backslash Z)$ is a bijection.}

\demo{Suppose that $Y$ is any quasi-compact Noetherian scheme. Then we know that the connected components of $Y$ identify with the connected components of the graph $\mathcal{G}_Y$, whose vertices are the irreducible components of $Y$, and where two vertices (corresponding to irreducible components $Z \neq Z'$) are linked if and only if they have a nonempty intersection. 

Thus, it is enough to prove the exact same statement when ``connected components'' are replaced with ``irreducible components'', and show that, given two irreducible components $Y \neq Y'$ of $X$, $Y$ and $Y'$ meet if and only if $Y \backslash Z$ and $Y' \backslash Z$ meet. Note that since $X$ and $X \backslash Z$ are pure of dimension one, any subspace of either of them is an irreducible component if and only if it is closed, irreducible, and of Krull dimension one. 

In particular, given an irreducible component $Y \subset X$, $Y \backslash Z$ is closed in $X \backslash Z$ and is a nonempty open subset of $Y$ (otherwise $Y$ would be contained in the zero-dimensional closed subspace $Z$, a contradiction), hence is irreducible and of Krull dimension one (because $Y$ is of finite type over a field). Moreover, the reunion of the $Y \backslash Z$, where $Y$ runs through the finitely many irreducible components of $X$, is exactly $X \backslash Z$, hence the $Y \backslash Z$ are the irreducible components of $X \backslash Z$. 

Given two irreducible components $Y,Y' \subset X$, if $Y \backslash Z=Y' \backslash Z$, then $Y \cap Y'$ is a closed subpace of $Y$ (resp. $Y'$) containing the dense open subspace $Y \backslash Z$ (resp. $Y' \backslash Z$), hence $Y \cap Y'=Y=Y'$. 

Given two distinct irreducible components $Y,Y' \subset X$, a point $x \in Y \cap Y'$ has the property that $\OO_{X,x}$ contains two distinct minimal prime ideals, hence is not regular. Therefore $x \notin Z$. Thus $Y$ meets $Y'$ if and only if $Y \backslash Z$ meets $Y' \backslash Z$ and the conclusion follows.}

\cor[XG-comp]{If $G$ is a $N$-torsion group over any $\Z[1/N]$-algebra $R$, then the geometric fibres of $Y_G(N) \rar \Sp{R}$ (resp. $X_G(N) \rar \Sp{R}$ if $R$ is regular excellent Noetherian) have exactly $\varphi(N)$ connected components. }

\demo{ 
The situation when $G$ is constant and $R=\Z[1/N]$ is completely handled by Proposition \ref{XN-Weil} and Lemma \ref{open-comp}. Let us discuss the general case. Let $R'$ be a faithfully flat \'etale $R$-algebra such that $G_{R'}$ is constant. It is enough to prove the statements for $Y_{G_{R'}}(N) \rar \Sp{R'}$ (resp. $X_{G_{R'}}(N) \rar \Sp{R'}$ if $R$, and therefore $R'$, is excellent). But $Y_{G_{R'}}(N) \rar \Sp{R'}$ (resp. $X_{G_{R'}}(N) \rar \Sp{R'}$) is a base change of $Y(N) \rar \Sp{\Z[1/N]}$ (resp. $X(N) \rar \Sp{\Z[1/N]}$), so it satisfies the conclusion.}

\prop[XG-pol]{Suppose that $G$ is a polarized $N$-torsion group over some $\Z[1/N]$-algebra $R$. Then the rule $(E/S,\iota) \longmapsto \det{\iota} \in (\Z/N\Z)^{\times}_S$ defines a morphism $\det:  Y_G(N) \rar (\Z/N\Z)^{\times}$, which extends to $\det: X_G(N) \rar (\Z/N\Z)^{\times}$ if $R$ is regular excellent Noetherian. In either case, one has $\det \circ [u] = u^2 \cdot \det$ for $u \in (\Z/N\Z)^{\times}$. 
Let, for each $u \in (\Z/N\Z)^{\times}$, $Y_{G,u}(N)$ (resp. $X_{G,u}(N)$ if $R$ is regular excellent Noetherian) the inverse image of $u$ under $\det: Y_G(N) \rar (\Z/N\Z)^{\times}$ (resp. $\det: X_G(N) \rar (\Z/N\Z)^{\times}$), and $\mathcal{E}_{G,u}(N)=\mathcal{E}_G(N) \times_{Y_G(N)} Y_{G,u}(N)$. Then $\mathcal{E}_{G,u}(N)/Y_{G,u}(N)$ represents $\P_{G,u}$, and $\P_{G,u}$ is a finite \'etale surjective moduli problem over $\Ell_R$.  
Moreover, if $R$ is regular excellent Noetherian, $\P_{G,u}$ is normal near infinity and $X_{G,u}$ identifies with the compactified moduli scheme $\CMPp{\P_{G,u}}$; under this identification, the cuspidal subscheme of $X_{G,u}$ identifies with the inverse image of the cuspidal subscheme of $X_{G}(N)$ in $X_{G,u}(N)$.  
}

\demo{It is clear that $\det$ is well-defined. By Corollary \ref{extension-to-finite}, it extends to $X_G(N)$ if $R$ is regular excellent Noetherian. 
If $\alpha \in \P_G$ is such that $(\mathcal{E}_G(N)/Y_G(N),\alpha)$ represents $\P_G$, it is a tautology that $(\mathscr{E}_{G,u}(N),Y_{G,u}(N), \alpha_{|\det^{-1}(u)})$ represents $\P_{G,u}$.

We need to show that $\P_{G,u}$ is finite \'etale over $\Ell_R$. To do that, consider an elliptic curve $E/S$ where $S$ is a $R$-scheme. Then there is a morphism of $S$ schemes $\det: (\P_G)_{E/S} \rar (\Z/N\Z)_S^{\times}$ given by $(\iota: G_T \rar E_T[N]) \longmapsto \det{\iota} \in (\Z/N\Z)^{\times}_S(T)$. It is clear that $(\P_{G,u})_{E/S} = (\P_G) \times_{(\Z/N\Z)^{\times}_S} u$, where $u$ is the closed open immersion $S \rar (\Z/N\Z)_S^{\times}$ onto the component corresponding to $u$. Thus $(\P_{G,u})_{E/S} \rar (\P_G)_{E/S}$ is a closed open immersion of $S$-schemes into a finite \'etale $S$-scheme, whence the conclusion.  

To show that $\P_{G,u}$ is surjective, it is enough to show that for any elliptic curve $E/k$ over an algebraically closed field where $N$ is invertible, then $(\P_{G,u})_{E/k}$ is not the empty scheme, that is, that $\P_{G,u}(E/k)$ is not the empty set. Since $G_{k}$ is constant, this follows from the fact that the Weil pairing $E[N](k) \times E[N](k) \rar \mu_N(k)$ is surjective. 

A consequence of the previous results is that, when $R$ is regular excellent Noetherian, $Y_{G,u}(N)$ is a normal scheme, hence $\P_{G,u}$ is normal near infinity.

The closed open immersion $Y_{G,u} \rar Y_G(N)$ extends to a finite map $\CMPp{\P_{G,u}} \rar X_G(N)$. It comes from a morphism of moduli problems, so it preserves the $j$-invariant. Since the image of $Y_{G,u}(N)$ is contained in $X_{G,u}(N)$ and $Y_{G,u}(N)$ is dense in $\CMPp{\P_{G,u}}$, the finite map $\CMPp{\P_{G,u}} \rar X_G(N)$ factors through the closed open subscheme $X_{G,u}(N)$. This finite map $f: \CMPp{\P_{G,u}} \rar X_{G,u}(N)$ is thus between normal schemes, and is an isomorphism above the dense open subscheme $Y_{G,u}(N)$ whose pre-image is dense. By Lemma \ref{kinda-zmt}, $f$ is an isomorphism. Since $f$ clearly preserves the $j$-invariant, $f$ also preserves the cuspidal subschemes. 

Since $\det \circ [v]=v^2 \cdot \det$, $[v]$ induces an isomorphism $Y_{G,u}(N) \rar Y_{G,uv^2}(N)$ of $\mathbb{A}^1_R$-schemes, which extends by Proposition \ref{compactification-functor} to an isomorphism $X_{G,u}(N) \rar X_{G,uv^2}(N)$ of $\mathbb{P}^1_R$-schemes where $R$ is regular excellent Noetherian. }

\cor[XG-pol-prop]{Let $G$ be a polarized $N$-torsion group over some $\Z[1/N]$-algebra $R$. Then, for any $u \in (\Z/N\Z)^{\times}$:
\begin{enumerate}[noitemsep,label=$(\roman*)$]
\item\label{PGu-1} $Y_{G,u}(N) \rar \Sp{R}$ (resp. $X_{G,u}(N) \rar \Sp{R}$ if $R$ is regular excellent Noetherian) is smooth of relative dimension one with connected geometric fibres (resp. proper smooth of relative dimension one with connected geometric fibres). 
\item\label{PGu-2} The $j$-invariant $Y_{G,u}(N) \rar \mathbb{A}^1_R$ (resp. $X_{G,u}(N) \rar \mathbb{P}^1_R$ if $R$ is regular excellent Noetherian) is finite locally free, \'etale on the locus where $j$ is defined and $j(j-1728)$ is invertible.
\item\label{PGu-3} If $R$ is regular excellent Noetherian, then $\P_{G,u}$ is smooth at infinity over $\Z[1/N]$ and the formation of its compactified moduli scheme and cuspidal subscheme commute with base change. 
\end{enumerate}
}

\demo{\ref{PGu-1}, \ref{PGu-2}, \ref{PGu-3} all follow from Propositions \ref{PG-elementary} and \ref{XG-pol}, apart from the statement on geometric fibres. For each geometric point $x$ of $\Sp{R}$, the associated geometric fibres of $(Y_{G,u})_x$ (or $(X_{G,u})_x$) are non-empty, and their disjoint reunion has $\varphi(N)=|(\Z/N\Z)^{\times}|$ connected components by Corollary \ref{XG-comp}. It follows that each $(Y_{G,u})_x$ (or $(X_{G,u})_x$) is connected. 
}

\prop[XG-Gamma0m]{Let $m_0,m_1,N \geq 1$ be pairwise coprime integers with $N \geq 3$ and $m_0$ square-free. Let $m=m_0m_1$, $G$ be a $N$-torsion group over some $\Z[(Nm_1)^{-1}]$-algebra $R$. Then the following properties hold:
\begin{enumerate}[noitemsep,label=$(\roman*)$]
\item\label{PGm-1} $\P_G(\Gamma_0(m))$ is representable by a scheme $Y_G(N,\Gamma_0(m))$, it is finite locally free of constant rank $r_{N,m} \geq 1$ over $\Ell_R$, \'etale above $\Z[1/m]$,
\item\label{PGm-2} $Y_G(N,\Gamma_0(m)) \rar \Sp{R}$ is affine, Cohen-Macaulay of relative dimension one with geometrically reduced fibres, smooth outside the supersingular $j$-invariants in characteristic $p \mid m_0$, its geometric fibres have $\varphi(N)$ connected components, and the $j$-invariant $Y_G(N,\Gamma_0(m)) \rar \mathbb{A}^1_R$ is finite locally free of rank $r'_{N,m} \geq 1$, \'etale on the locus where $m_0j(j-1728)$ is invertible. 
\item\label{PGm-3} If $R$ is regular excellent Noetherian, $\P_G(\Gamma_0(m))$ is smooth at infinity over $\Z[(Nm_1)^{-1}]$, 
\item\label{PGm-4} If $R$ is regular excellent Noetherian, the compactified moduli scheme $X_G(N,\Gamma_0(m))$ of $\P_G(\Gamma_0(m))$ is Cohen-Macaulay of relative dimension one over $\Sp{R}$, smooth outside the supersingular $j$-invariants in characteristic $p \mid m_0$; its geometric fibres are reduced and possess $\varphi(N)$ connected components, 
\item\label{PGm-5} If $R$ is regular excellent Noetherian, the $j$-invariant $X_G(N,\Gamma_0(m)) \rar \mathbb{P}^1_R$ is finite locally free of constant rank $r'_{N,m}$,
\item\label{PGm-6} The natural action of $(\Z/Nm\Z)^{\times}$ on $\P_G(\Gamma_0(m))$ (as a moduli problem of level $Nm$) factors through the quotient to $(\Z/N\Z)^{\times}$ and is equivariant for the natural forgetful map $\P_G(\Gamma_0(m)) \rar \P_G$. If $R$ is regular excellent Noetherian, it extends as an equivariant action of $(\Z/N\Z)^{\times}$ on the morphism $X_G(N,\Gamma_0(m)) \rar X_G(N)$. 
\end{enumerate}}

\rem{When $G$ is constant, we write $Y_G(N,\Gamma_0(m))=Y(\P_G,\Gamma_0(m))$ for short, and similarly for the compactified moduli scheme.}

\demo{The representability of the moduli problem follows from Proposition \ref{product-rel-rep}, using \cite[Theorem 5.1.1]{KM} and Proposition \ref{PG-elementary}, since we are taking the Cartesian product of $\P_G$ (representable) and $[\Gamma_0(m)]_R$ (relatively representable). It is finite locally free of constant positive rank (only depending on $N,m$) and \'etale above $\Z[1/m]$ over $\Ell_R$, because $\P_G$ and $[\Gamma_0(m)]_R$ are. The results of \ref{PGm-2} (apart from the number of connected components of the geometric fibres and the ``constant rank $r'_{N,m} \geq 1$'' part of the $j$-invariant) are the consequence of Corollary \ref{over-Ell-concrete} and the first part of Proposition \ref{struct-gamma0m-smooth}. In particular, $\P_G(N,\Gamma_0(m))$ is normal near infinity.

When $R=\Z[(Nm_1)^{-1}]$ and $G$ is constant, since $\mathcal{E}(N)/Y(N)$ is a modular family in the sense of \cite[\S 4.11]{KM}, \cite[Theorem 5.1.1]{KM} shows that $Y(N,\Gamma_0(m))$ is a regular scheme. Thus, the $j$-invariant map $X(N,\Gamma_0(m)) \rar \mathbb{P}^1_{\Z[(Nm_1)^{-1}]}$ is finite between two regular Noetherian schemes (every nonempty open subset of either having dimension two), hence it is flat by the miracle flatness theorem \cite[Lemma 00R4]{Stacks}. It is thus finite locally free of constant rank $r'_{N,m} \geq 1$, since $\mathbb{P}^1_{\Z[(Nm_1)^{-1}]}$ is connected. 

In particular, $X(N,\Gamma_0(m)) \rar \Sp{\Z[(Nm_1)^{-1}]}$ is proper flat. Moreover, $Y(N,\Gamma_0(m)) \rar \mathbb{A}^1_{\Z}$ is finite locally free $r'_{N,m}$, and by taking arbitrary base change we find that the $j$-invariant map $Y_G(N,\Gamma_0(m)) \rar \mathbb{A}^1_R$ is always of constant rank $r'_{N,m}$ when $G$ is constant. Since $G$ is always \'etale-locally constant, it follows that $j: Y_G(N,\Gamma_0(m)) \rar \mathbb{A}^1_{\Z[(Nm_1)^{-1}]}$ has rank $r'_{N,m}$ for any $N$-torsion group $G$. 

When $R$ is regular excellent Noetherian of generic characteristic zero or $G$ is constant, \ref{PGm-3},\ref{PGm-4} (apart from the statement on the number of connected components of the geometric fibres) follow from Proposition \ref{struct-gamma0m-smooth} and Corollary \ref{smooth-at-infinity-gamma0m}. 

If $G$ be a $N$-torsion group over any regular excellent Noetherian $\Z[(Nm_1)^{-1}]$-algebra $R$, then there is a faithfully flat \'etale $R$-algebra $R'$ such that $G_{R'}$ is constant. Thus, by Proposition \ref{smooth-at-infinity-descends} $\P_G(\Gamma_0(m))$ is smooth at infinity over $\Z[(Nm)^{-1}]$, and the formation of its compactified moduli scheme and cuspidal subscheme commutes with base change, thus \ref{PGm-3}, \ref{PGm-4} (without, once again, the statement on the connected components of the geometric fibres) hold. 

The compactification of the $j$-invariant is a finite map $\overline{j}: X_G(N,\Gamma_0(m)) \rar \mathbb{P}^1_R$ of $R$-schemes, where the source is Cohen-Macaulay of relative dimension one over $R$, and the target is smooth over $R$ of relative dimension one (with $R$ regular). By miracle flatness, this implies that $\overline{j}$ is flat, hence locally free. Since it has rank $r'_{N,m}$ above the dense open subscheme $\mathbb{A}^1_R$, it is of constant rank $r'_{N,m}$.

The assertion \ref{PGm-6} is formal, using Proposition \ref{compactification-functor} and the definition of the natural level $Nm$ structure on $\P_G \times [\Gamma_0(m)]_R$, since it is clear that $(\Z/m\Z)^{\times}$ acts trivially on the moduli problem $[\Gamma_0(m)]_R$. 

To conclude, it is enough to prove the statement on the number of connected components of geometric fibres of $X(N,\Gamma_0(m)) \rar \Sp{\Z[(Nm_1)^{-1}]}$ (by base change, this will hold whenever $G$ is constant, and any $G$ is \'etale locally constant). By \cite[Lemma 0E0N]{Stacks}, we only need to prove it when $G$ is constant over $R=\Z[(Nm)^{-1}]$. It is thus enough to prove in this situation that $X_G(N,\Gamma_0(m)) \rar X(N)_{\Z[(Nm)^{-1}]} \overset{\mrm{We}}{\rar} \Sp{\OO_N[m^{-1}]}$ has connected geometric fibres (since this map is proper smooth). This is a consequence of \cite[Corollary 10.9.5]{KM}, by setting $N_1=N,K=Nm,N_2=m$ and $\Gamma \leq \GL{\Z/Nm\Z}$ the subgroup of upper-triangular matrices $M \equiv I_2 \mod{N}$.
}

\cor[XG-Gamma0m-regl]{Let $m_0,m_1,m,N,R,G$ be as in Proposition \ref{XG-Gamma0m}, and assume that $R$ is Noetherian. Then the following statements hold:
\begin{enumerate}[noitemsep,label=$(\roman*)$]
\item\label{PGm-7} If $R$ is $(R_k)$, $m_0 \in R$ is not a zero divisor, and [$k \leq 1$ or ($R/pR$ is $(R_k)$ for all $p \mid m_0$)], then $Y_G(N,\Gamma_0(m))$ is $(R_k)$,
\item\label{PGm-7b} If $R$ is $(S_k)$, then so is $Y_G(N,\Gamma_0(m))$,  
\item\label{PGm-8} If $R$ is reduced (resp. normal) and $m_0 \in R$ is not a zero divisor, $Y_G(N,\Gamma_0(m))$ is reduced (resp. normal),
\item\label{PGm-9} If $R$ is regular excellent Noetherian and $m_0 \in R$ is not a zero divisor, then $X_G(N,\Gamma_0(m))$ is normal,
\item\label{PGm-9b} If $R$ is regular excellent Noetherian and $R/pR$ is $(R_k)$ for all $p \mid m_0$, then $X_G(N,\Gamma_0(m))$ is $(R_k)$.
\end{enumerate}
}

\demo{Regularity ascends along smooth maps \cite[Lemma 07NF]{Stacks}, so $X_G(N,\Gamma_0(m))$ is regular at every point of its cuspidal subscheme. Hence \ref{PGm-9b} (resp. \ref{PGm-9}) follows from \ref{PGm-7} (resp. \ref{PGm-8}). The reunion of \ref{PGm-7} and \ref{PGm-7b} implies \ref{PGm-8} by, for instance, Serre's criterion for normality (and the analog for reducedness) \cite[Lemmas 031R, 031S]{Stacks}. \ref{PGm-7b} is a consequence of the fact that $Y_G(N,\Gamma_0(m))$ is Cohen-Macaulay of relative dimension one over $\Sp{R}$.

Let us finally discuss \ref{PGm-7}: we assume that $m_0 \in R$ is not a zero divisor. After a faithfully flat \'etale base change, we can assume that $G$ is constant. Since the kernel of $R \rar R \otimes \Q$ is finitely generated, it is contained in the $d$-torsion of $R$ for some integer $d \geq 1$ prime to $Nm$. Thus, after taking the Zariski-cover $\Sp{R[1/m]}$, $\Sp{R[1/d]}$, we can assume that either $m \in R^{\times}$ or $R$ is flat over $\Z$. The first case follows from ascending properties of smooth maps, so from now on we assume that $R$ is flat over $\Z$.

Let us first discuss the property $(R_0)$ (resp. $(R_1)$). Let $U \subset Y(N,\Gamma_0(m))$ the open subset excluding the supersingular $j$-invariants in characteristic $p \mid m$. Then the complement of $U$ is finite, hence by the going-down theorem \cite[Lemma 00HS]{Stacks} (applied to the flat map $Y_G(N,\Gamma_0(m)) \rar Y(N,\Gamma_0(m))$) the complement of $U_R$ in $Y(N,\Gamma_0(m))_R$ cannot contain points of codimension zero or one. Since $U \rar \Sp{\Z}$ is smooth, $U_R$ is $(R_0)$ (resp. $(R_1)$) and we are done.   

Now, we fix some $k \geq 0$ and assume that $R$ and each $R/pR$ (for $p \mid m_0$) is $(R_k)$. The map $Y_G(N,\Gamma_0(m))_R \rar Y(N,\Gamma_0(m))$ is affine and flat to a regular base scheme, and the rings of fibres are of the form $R \otimes_{\Z} F$, where $F$ is a finitely generated field extension of a field $k$ which is either finite or $\Q$. We need to prove that these $R \otimes_{\Z} F$ satisfy the property $(R_k)$. By considering the map $R \otimes_{\Z} k \rar R \otimes_{\Z} F$, all we need to show is that the fibre rings of $R \otimes_{\Z} k \rar R \otimes_{\Z} F$ are $(R_k)$. These fibres are regular by \cite[Lemma 0381]{Stacks} since $k$ is perfect and $F$ is finitely generated. 
}

\defi[XGm-weil-det]{If $G$ is constant over $R$, the following map (denoted $\mrm{We}$) is the \emph{Weil pairing}: \[X_G(N,\Gamma_0(m)) \rar X_G(N) \simeq X(N)_R \overset{\mrm{We}}{\rar} \Sp{\OO_N} \times_{\Z} \Sp{R}.\]
If $G$ is polarized, we call $\det$ (the determinant) the map $X_G(N,\Gamma_0(m)) \rar X_G(N) \overset{\det}{\rar} (\Z/N\Z)^{\times}_R$. }

\cor[XGm-conn]{The maps $\mrm{We}$ and $\det$ from Definition \ref{XGm-weil-det} are flat surjective with connected geometric fibres. In particular, $X(N,\Gamma_0(m))$ and $Y(N,\Gamma_0(m))$ are irreducible. }

\demo{In the case of the Weil pairing map, it is enough to prove the result when $R=\Z[1/N]$. Then $X(N,\Gamma_0(m)) \rar \Sp{\OO_N} \rar \Sp{\Z[1/N]}$ is surjective, and its geometric fibres have $\varphi(N)$ connected components. Since $\Sp{\OO_N}$ is finite \'etale over $\Sp{\Z[1/N]}$, $X(N,\Gamma_0(m))$ is proper flat over $\Sp{\OO_N}$ and in particular surjective. If $k$ is an algebraically closed field of characteristic not dividing $N$, it is easy to see that the geometric fibre of $X(N,\Gamma_0(m))$ at $\Z[1/N] \rar k$ is exactly the disjoint reunion of the geometric fibres of $X(N,\Gamma_0(m)) \rar \Sp{\OO_N}$ for all $\varphi(N)$ morphisms $\OO_N \rar k$. All of these $\varphi(N)$ geometric fibres are nonempty, but their disjoint reunion has exactly $\varphi(N)$ connected components: this shows that they must each be connected.

The map $\mrm{We}: X(N,\Gamma_0(m)) \rar \Sp{\OO_N[1/m_1]}$ is smooth with geometrically connected fibres. Therefore $X(N,\Gamma_0(m))_{\Q}$ is a regular connected scheme (hence irreducible), and it is dense in $X(N,\Gamma_0(m))$, whence $X(N,\Gamma_0(m))$ (and thus $Y(N,\Gamma_0(m))$) is irreducible. 

For the geometric fibres of $\det$, a similar argument reduces the problem to showing that for any $u \in (\Z/N\Z)^{\times}$ and any algebraically closed field $k$ of residue characteristic not dividing $N$, there is an element in $X(N,\Gamma_0(m))(k)$. Let $E/k$ be any elliptic curve (ordinary if $m$ vanishes in $k$, for instance with $j$-invariant taken outside $\F_{p^2}$ by \cite[Lemma 12.5.4]{KM}), there are isomorphisms of arbitrary determinant $G_k \rar E[N](k)$ since $G, E[N]$ are constant \'etale, and there is an arbitrary finite \'etale subgroup scheme of $E[m]$, cyclic of order $m$ (if needed by separating the $p$ part, where $p \mid m$ is the characteristic of $k$). }

\cor[flat-deg-constant]{In Corollary \ref{flat-deg}, the map $D_{d,t}$ is finite locally free of constant rank, which is independent from the moduli problem.}

\demo{Using the same reduction than the proof of Corollary \ref{flat-deg} (whose notations we keep), it is enough to prove the statement when $\P=[\Gamma(\ell)]_{\Z[1/\ell]}$ with $\ell \nmid 2m$ some prime. Assume more generally that $\P=[\Gamma(k)]_{\Z[1/k]}$ where $k \geq 3$ is any integer prime to $m$. Then $D_{d,t}$ is a finite locally free morphism $Y(k,\Gamma_0(m)) \rar Y(k,\Gamma_0(t))$ between regular irreducible $\Z[1/k]$-schemes, so that it is locally free of constant rank $r_k \geq 1$. By Proposition \ref{degn-cartesian}, for any primes $\ell,\ell' \nmid 2m$, we then have $r_{\ell}=r_{\ell\ell'}=r_{\ell'}$, and the conclusion follows.}

\cor{Let $m_0,m_1,m,N,R,G$ be as in Proposition \ref{XG-Gamma0m} and assume that $G$ is polarized. Let $u \in (\Z/N\Z)^{\times}$. Then the moduli problem $\P_{G,u}(\Gamma_0(m)) := \P_{G,u} \times_{\P_G} \P_G(\Gamma_0(m))$ over $\Ell_R$ is representable by the closed open subscheme $\det^{-1}(u)$ of $Y_G(N,\Gamma_0(m))$. This moduli problem is finite locally free surjective over $(\Ell_R)$. 
Moreover, if $R$ is regular excellent Noetherian, $\P_{G,u}(\Gamma_0(m))$ is normal near infinity and the closed open subscheme $X_{G,u}(N,\Gamma_0(m))=\det^{-1}(u)$ of $X_G(N,\Gamma_0(m))$ identifies with the compactified moduli scheme $\CMPp{\P_{G,u}}$; under this identification, the cuspidal subscheme $\CPp{\P_{G,u}}$ identifies with the inverse image of the cuspidal subscheme of $X_{G}(N,\Gamma_0(m))$ in $X_{G,u}(N,\Gamma_0(m))$.  
The $Y_{G,u}(N,\Gamma_0(m))$ (resp. the $X_{G,u}(N,\Gamma_0(m))$) seen as $R$-schemes, have geometrically connected fibres. 
}

\demo{The first part is completely formal, apart from the surjectivity. It is enough to show that for any elliptic curve $E$ over an algebraically closed field $k$ where $Nm_1$ is invertible, $\P_{G,u}(\Gamma_0(m))(E/k)$ is nonempty. This is clear since $G$ is \'etale, hence constant, and the Weil pairing $E[N] \times E[N] \rar \mu_N$ is onto, while, if some prime power $q \mid m$ is invertible in $k$, we can choose as $\Gamma_0(q)$-structure any cyclic subgroup generated by a point of order $q$, and if there is a prime $p$ that vanishes in $k$ and divides $m$, then $p$ prime to $m/p$ and the kernel of the Frobenius (see Section \ref{char-p}) is a $\Gamma_0(p)$-structure. 

The moduli problem $\P_{G,u}$ is normal near infinity because $Y_{G,u}(N,\Gamma_0(m)) \rar Y_G(N,\Gamma_0(m))$ is a closed open immersion. Its compactification in the sense of Proposition \ref{compactification-functor} is a finite map $f_u: \CMPp{\P_{G,u}(\Gamma_0(m))} \rar X_G(N,\Gamma_0(m))$ whose image is contained in the closed open subscheme $\det^{-1}(u)$. In other words, $f_u: \CMPp{\P_{G,u}(\Gamma_0(m))} \rar X_{G,u}(N,\Gamma_0(m))$ is a finite morphism and is an isomorphism above the dense open subscheme $Y_{G,u}(N,\Gamma_0(m))$ with dense inverse image. Both schemes are normal on the inverse image of the closed subspace $\{1/j = 0\}$. By Lemma \ref{kinda-zmt}, $f_u$ is an isomorphism, and in particular it induces an isomorphism from $\CPp{\P_{G,u}(\Gamma_0(m))}$ to the restriction to $\CPp{\P_G(\Gamma_0(m))}$ to $\det^{-1}(u)$ (since the induced map is an isomorphism on the underlying topological space, while the two subschemes are reduced).

For the last part, by Proposition \ref{XG-Gamma0m}, the disjoint reunion of the connected components of the geometric fibres of $Y_{G,u}(N,\Gamma_0(m))$ (resp. $X_{G,u}(N,\Gamma_0(m))$) over a geometric point $\overline{x}$ of $R$ (resp. assuming that $R$ is regular excellent Noetherian) over all $u \in (\Z/N\Z)^{\times}$ has exactly $\varphi(N)=|(\Z/N\Z)^{\times}|$ elements. Since every $Y_{G,u}(N,\Gamma_0(m)) \rar \Sp{R}$ is surjective, its geometric fibres have to be connected.}

\prop[ext-degn]{Let $N \geq 3$, $R$ be a $\Z[1/N]$-algebra (resp. $R$ is a regular excellent Noetherian $\Z[1/N]$-algebra) and $G$ be a $N$-torsion group over $R$. The collection of $Y_G(N,\Gamma_0(m))$ (resp. $X_G(N,\Gamma_0(m))$), for integers $m=m_0m_1$ coprime to $N$ with $m_0,m_1 \geq 1$, $m_0$ square-free and $m_1 \in R^{\times}$, is endowed with degeneracy maps $D_{d,t}$ and Atkin-Lehner automorphisms as in Propositions \ref{degeneracies} and \ref{AL-auto} (resp. Proposition \ref{degn-AL-infty}). They are finite locally free of constant rank, and commute with the action of $(\Z/N\Z)^{\times}$. The relations stated in Corollary \ref{elem-degeneracies} and Proposition \ref{AL-auto} still hold.}

\demo{This is a direct translation of what we previously did in this more concrete context. The specific statements are Propositions \ref{degeneracies}, \ref{AL-auto}, \ref{degn-AL-infty}, as well as Corollaries \ref{elem-degeneracies}, \ref{flat-deg-infty}, \ref{flat-deg-constant}.}

\prop[degen-AL-GL2]{Let $m_0,m_1,m,N,R,G$ be as in Proposition \ref{XG-Gamma0m}. Assume that $G$ is constant. 
\begin{itemize}[noitemsep,label=$-$]
\item Then $Y(N,\Gamma_0(m))$ has a natural action of $\GL{\Z/N\Z}$ which commutes to every degeneracy map and Atkin-Lehner automorphism. A matrix $\begin{pmatrix}a & b\\c & d\end{pmatrix}$ acts by \[(E/S,(P,Q),C) \longmapsto (E/S,(aP+bQ,cP+dQ),C).\] Moreover, the action of $aI_2$ (for $a \in (\Z/N\Z)^{\times}$) is exactly the action of $[a']$ for any $a' \in (\Z/Nm\Z)^{\times}$ congruent to $a$ modulo $N$.
\item If moreover $R$ is regular excellent Noetherian, this action extends to $X_G(N,\Gamma_0(m))$ and commutes to the extensions of degeneracy maps and Atkin-Lehner automorphisms; moreover, the action of $aI_2$ is given by any $[a']$ for any $a' \in (\Z/Nm\Z)^{\times}$ congruent to $a$ modulo $N$.
\end{itemize}}

\demo{This follows from the fact that $\P(N)$ has a natural level structure of level $N$, and that the natural action of $\GL{\Z/N\Z}$ on $\P(N)$ acts by morphisms of moduli problems of level $N$, as well as Corollary \ref{elem-degeneracies} and Proposition \ref{AL-auto}. The action of scalar matrices is an explicit computation. The second point follows from Proposition \ref{compactification-functor}. }

\prop[degen-AL-Weil]{Let $m_0,m_1,m,N,R,G$ be as in Proposition \ref{XG-Gamma0m}.
\begin{itemize}[noitemsep,label=$-$]
\item Suppose that $G$ is constant and that $d,t \geq 1$ are such that $dt \mid m$ (resp. $d \geq 1$ divides $m$ and is coprime to $m/d$). Then $\mrm{We} \circ D_{d,t} = \underline{d} \circ \mrm{We}$ (resp. $\mrm{We} \circ w_d = \underline{d} \circ \mrm{We}$).
\item Suppose that $G$ is polarized and that $d,t \geq 1$ are such that $dt \mid m$ (resp. $d \geq 1$ divides $m$ and is coprime to $m/d$). Then $\det \circ D_{d,t} = d \cdot \det$ (resp. $\det \circ w_d = d \cdot \det$).
\item When $R$ is regular excellent Noetherian, both identities extend to $X_G(N,\Gamma_0(m))$.
\end{itemize}}

\demo{The first two points are direct computations. The third point follows from the previous two and the fact that $S(X_G(N,\Gamma_0(m))) \rar S(Y_G(N,\Gamma_0(m)))$ is injective for any separated scheme $S$ by Proposition \ref{compactification-functor}. }

\bigskip

The methods of this section can be applied, \emph{mutans mutandis}, to prove the following statement:

\prop[X1N-gamma0]{Let $N \geq 5$ be an integer and $m=m_0m_1$ be coprime to $N$, with $m_0$ square-free and coprime to $m_1$. 
\begin{itemize}[noitemsep,label=\tiny$\bullet$]
\item For any $\Z[1/(Nm_1)]$-algebra $R$, the moduli problem $[\Gamma_1(N)]_R \times [\Gamma_0(m)]_R$ is representable, finite flat of constant rank over $\Ell_R$, \'etale above $R[m_0^{-1}]$. Write $Y_1(N,\Gamma_0(m))_R$ for the moduli scheme. The natural action of $(\Z/Nm\Z)^{\times}$ on $Y_1(N,\Gamma_0(m))$ factors through $(\Z/N\Z)^{\times}$.
\item $Y_1(N,\Gamma_0(m))_R \rar \Sp{R}$ is affine, Cohen-Macaulay of relative dimension one, smooth outside the supersingular $j$-invariants in characteristic $p \mid m_0$, with geometrically connected fibres.
\item When $R$ is regular excellent Noetherian, this moduli problem is smooth at infinity over $\Z[(Nm_1)^{-1}]$, its compactified moduli scheme $X_1(N,\Gamma_0(m))$ is finite locally free of constant rank over $\mathbb{P}^1_R$, proper Cohen-Macaulay of relative dimension one over $\Sp{R}$ with geometrically connected fibres, and smooth outside the supersingular $j$-invariants in characteristic $p \mid m_0$. The natural action of $(\Z/Nm\Z)^{\times}$ factors through $(\Z/N\Z)^{\times}$. 
\item The formation of $X_1(N,\Gamma_0(m))$ and its cuspidal subscheme commutes with arbitrary base change of regular excellent Noetherian $\Z[(Nm_1)^{-1}]$-algebras.
\end{itemize}
The conclusions of Proposition \ref{ext-degn} apply here after replacing $Y_G(N,\Gamma_0(m))$ with $Y_1(N,\Gamma_0(m))$ (resp. $X_G(N,\Gamma_0(m))$ with $X_1(N,\Gamma_0(m))$ if $R$ is regular excellent Noetherian).}

\section{Properties of Galois twists}

In this section, we fix an integer $N \geq 3$ and a field $k$ of characteristic not dividing $N$ with separable closure $k_s$. We fix a left representation $\rho: \mrm{Gal}(k_s/k) \rar \GL{\Z/N\Z}$.

\subsection{$Y_{\rho}(N)$ and its rational points}

\defi{The group $\GL{\Z/N\Z}$ acts on the following diagram of quasi-projective $k$-schemes (and trivially on the second row)
\[
\begin{tikzcd}[ampersand replacement=\&]
Y(N)_k \arrow{r}{\mrm{in}}\arrow{d}{j}\& X(N)_k \arrow{d}{\overline{j}}\& C_N \arrow{l}{\mrm{cl}} \arrow{d}\\
\mathbb{A}^1_k \arrow{r}{\mrm{in}} \& \mathbb{P}^1_k \& \Sp{k} \arrow{l}{\infty}
\end{tikzcd}
\] Twisting by $\rho$ thus gives a diagram of quasi-projective $k$-schemes:
\[
\begin{tikzcd}[ampersand replacement=\&]
Y_{\rho}(N)_k \arrow{r}{\mrm{in}}\arrow{d}{j}\& X_{\rho}(N) \arrow{d}{\overline{j}}\& C_{N,\rho} \arrow{l}{\mrm{cl}} \arrow{d}\\
\mathbb{A}^1_k \arrow{r}{\mrm{in}} \& \mathbb{P}^1_k \& \Sp{k} \arrow{l}{\infty}
\end{tikzcd}
\]
Moreover, the first cell is a Cartesian diagram, $j,\overline{j}$ are finite flat of constant rank, and $C_{N,\rho}$ is finite \'etale over $k$. Moreover, $X_{\rho}(N)$ is smooth over $\Sp{k}$.}

\demo{This follows from Propositions \ref{cocycle-twist} and \ref{twist-equiv}, as well as the properties of $X(N)$ discussed in the previous section, and the fact that finiteness, flatness, \'etaleness, smoothness descend along fpqc covers.}

\lem[right-representation]{Let $E/k$ be an elliptic curve and $P,Q \in E[N](k_s)$ be such that $E[N](k_s)$ is free over $\Z/N\Z$ with basis $(P,Q)$. For every $\sigma \in \mrm{Gal}(k_s/k)$, define $D(\sigma) \in \GL{\Z/N\Z}$ such that $D(\sigma)\begin{pmatrix}P \\ Q\end{pmatrix}=\begin{pmatrix}\sigma(P)\\\sigma(Q)\end{pmatrix}$. Then $D(\sigma\sigma')=D(\sigma')D(\sigma)$ for all $\sigma,\sigma' \in \mrm{Gal}(k_s/k)$.}

\demo{Simple computation.}

From this result, we obtain a description of the rational points of $Y_{\rho}(N)$. 

\prop[rational-Yrho]{For any subextension $L$ of $k_s/k$, $Y_{\rho}(N)(L)$ is in natural bijection with the $L$-isomorphism classes of triples $(E,P,Q)$ satisfying the following properties: 
\begin{itemize}[noitemsep,label=\tiny$\bullet$]
\item $E$ is an elliptic curve defined over $L$,
\item $P,Q \in E[N](k_s)$ and $E[N](k_s)$ is free over $\Z/N\Z$ with basis $(P,Q)$,
\item For all $\sigma \in \mrm{Gal}(k_s/L)$, $\begin{pmatrix}\sigma(P)\\\sigma(Q)\end{pmatrix} = \rho^{-1}(\sigma)\begin{pmatrix}P \\ Q\end{pmatrix}$.
\end{itemize}
If $x \in Y_{\rho}(N)(L)$ is given by the triple $(E,P,Q)$, then $j(x) \in \mathbb{P}^1(L)$ is the $j$-invariant of $E$.}

\demo{Let $k \subset L \subset k_s$ be a subextension. Then, by \cite[(4.3.1)]{KM}, $Y_{\rho}(N)(L)$ is in bijection with the triples $(E,P,Q)$, where $E$ is an elliptic curve over $k_s$, and $(P,Q)$ is a basis of $E[N](k_s)$ such that, for every $\sigma \in \mrm{Gal}(k_s/L)$, one has $(E,P,Q) \simeq_{\overline{k_s}} (\sigma(E),a\sigma(P)+b\sigma(Q),c\sigma(P)+d\sigma(Q))$, where $\rho(\sigma)=\begin{pmatrix}a & b\\c & d\end{pmatrix}$. In particular, $j(E) \in k_s$ is fixed by $G_L$, hence we can assume that $E$ is defined over $L$, and that, for all $\sigma \in G_L$ with $\rho(\sigma)=\begin{pmatrix}a & b\\c & d\end{pmatrix}$, there is a $k_s$-automorphism $t_{\sigma}$ of $E_{k_s}$ mapping $(P,Q)$ to $(a\sigma(P)+b\sigma(Q),c\sigma(P)+d\sigma(Q))$.  

We claim that $\sigma \longmapsto t_{\sigma}^{-1}$ satisfies the cocycle property of Proposition \ref{cocycle-twist}. Since the identity involves two $k_s$-automorphisms of $E_{k_s}$, it is enough, by rigidity, to check it on $E[N](k_s)$, and hence on $P$ and $Q$. Let $\sigma,\sigma' \in \mrm{Gal}(k_s/L)$ and $\rho(\sigma)=\begin{pmatrix} a &b\\c & d\end{pmatrix}, \rho(\sigma')=\begin{pmatrix}a'&b'\\c'&d'\end{pmatrix}$. Then 
\begin{align*}
t_{\sigma\sigma'}(P) &= (aa'+bc')\sigma(\sigma'(P))+(ab'+bd')\sigma(\sigma'(Q)) \\
&= a\sigma(a'\sigma'(P)+b'\sigma'(Q))+b\sigma(c'\sigma'(P)+d'\sigma'(Q))= a\sigma(t_{\sigma'}(P))+b\sigma(t_{\sigma'}(Q))\\
& = \sigma\circ t_{\sigma'}\circ \sigma^{-1}(a\sigma(P)+b\sigma(Q)) = \sigma t_{\sigma'}\sigma^{-1}t_{\sigma}(P),
\end{align*}
and similarly for $Q$, hence the cocycle property holds. 

Let $E'=E_{t^{-1}}$ be the twisted elliptic curve over $L$: we still have $P,Q \in E'[N](k_s)$, and $(E',P,Q)$ corresponds to the same point of $Y(N)(k_s)$ as $(E,P,Q)$. Let now $\sigma \in \mrm{Gal}(k_s/L)$, and write $\rho(\sigma)=\begin{pmatrix}a & b\\c & d\end{pmatrix}, \rho(\sigma)^{-1}=\begin{pmatrix}a' & b'\\c' & d'\end{pmatrix}$. Then one has 
\[\sigma_{E'}(P)=t_{\sigma}^{-1}\sigma_E(P) = t_{\sigma}^{-1}(a'(a\sigma(P)+b\sigma(Q))+b'(c\sigma(P)+d\sigma(Q)))= a'P+b'Q,\]

hence $\begin{pmatrix}\sigma_{E'}(P)\\\sigma_{E'}(Q)\end{pmatrix}=\rho(\sigma)^{-1}\begin{pmatrix}P \\ Q\end{pmatrix}$.

In other words, we proved that every point in $Y_{\rho}(N)(L)$ came from a triple satisfying the required conditions. Conversely, it is clear that every such triple produces a $k_s$-rational point on $Y_{\rho}(N)$ stable under $\mrm{Gal}(k_s/L)$, hence a $L$-rational point. 

To conclude, all we need to do is show that any two triples $(E,P,Q),(E',P',Q')$ satisfying the requested conditions and representing the same point in $Y(N)(k_s)$ are $L$-isomorphic. Indeed, these triples are isomorphic over $k_s$ by \cite[(4.3.1)]{KM}, so let $f: E_{K_s} \rar E'_{K_s}$ be a $k_s$-isomorphism mapping $(P,Q)$ to $(P',Q')$. One easily checks that the induced map $f[N]: E[N](k_s) \rar E'[N](k_s)$ commutes to $\mrm{Gal}(k_s/L)$, hence by rigidity the same holds for $f$, so that $f$ is defined over $L$ by Corollary \ref{inf-galois-desc}, and $f$ is by fpqc descent an isomorphism. 

The last statement is clear from the construction and the definition of $Y_{\rho}$. 
}

\prop[quadratic-twist-equiv]{Let $\chi: \mrm{Gal}(k_s/k) \rar \{\pm 1\}$ be a character, and let \[\rho'=\chi \otimes \rho: \mrm{Gal}(k_s/k) \rar \GL{\Z/N\Z}.\] The natural $k_s$-isomorphism between the following two diagrams
\[
\begin{tikzcd}[ampersand replacement=\&]
Y_{\rho}(N)_k \arrow{r}{\mrm{in}}\arrow{d}{j}\& X_{\rho}(N)_k \arrow{d}{\overline{j}}\& (C_N)_{\rho} \arrow{l}{\mrm{cl}} \arrow{d} \& Y_{\rho'}(N)_k \arrow{r}{\mrm{in}}\arrow{d}{j}\& X_{\rho'}(N)_k \arrow{d}{\overline{j}}\& (C_N)_{\rho'} \arrow{l}{\mrm{cl}}\arrow{d} \\ 
\mathbb{A}^1_k \arrow{r}{\mrm{in}} \& \mathbb{P}^1_k \& \Sp{k} \arrow{l}{\infty} \& \mathbb{A}^1_k \arrow{r}{\mrm{in}} \& \mathbb{P}^1_k \& \Sp{k} \arrow{l}{\infty}
\end{tikzcd}
\]
coming from the definition of a Galois twist is defined over $k$.

Moreover, for any finite subextension $L \subset k_s$ of $k$, the induced map $Y_{\rho}(N)(L) \rar Y_{\rho'}(N)(L)$ is described, through the description of Proposition \ref{rational-Yrho}, by the rule $(E,P,Q) \longmapsto (E \otimes \chi, P, Q)$, where $E \otimes \chi$ is the quadratic twist of $E$ by $\chi_{|\mrm{Gal}(k_s/L)}$.}

\demo{By Corollary \ref{inf-galois-desc}, since it involves quasi-compact $k$-schemes, the isomorphism is defined over $k$ as long as it commutes to $\mrm{Gal}(k_s/k)$. Since $X_{\rho}(N),X_{\rho'}(N)$ are proper smooth over $k$, it is enough to check that the induced bijection $Y_{\rho}(N)(k_s) \rar Y_{\rho'}(N)(k_s)$ commutes with the action of $\mrm{Gal}(k_s/k)$. 

Now, given $\sigma \in \mrm{Gal}(k_s/k)$ and an elliptic curve $E/k_s$ with a basis $(P,Q)$ of $E[N]$, let $x \in Y_{\rho}(N)(k_s), x' \in Y_{\rho}(N)(k_s)$ be the geometric points attached to the triple $(E,P,Q)$. If $\rho(\sigma)=\begin{pmatrix}a & b\\c & d\end{pmatrix}$, then one computes that
\begin{align*}
\sigma(x)&=(\sigma(E),a\sigma(P)+b\sigma(Q),c\sigma(P)+d\sigma(Q)),\\
\sigma(x')&=(\sigma(E),\chi(\sigma)(a\sigma(P)+b\sigma(Q)), \chi(\sigma)(c\sigma(P)+d\sigma(Q))).
\end{align*}
The two triples are isomorphic (by the automorphism $\chi(\sigma)$ of $\sigma(E)$) and therefore equal as geometric points of $Y_{\rho}(N)$.   
 
Let now $L/k$ be a subextension of $k_s/k$ and $(E,P,Q)$ be a triple as in Proposition \ref{rational-Yrho} parametrizing some $L$-point of $Y_{\rho}(N)$. Then its image in $Y_{\rho'}(N)(L)$ is parametrized by a triple $(E',P',Q')$ which satisfies the conditions of Proposition \ref{rational-Yrho} (but for the left representation $\rho'$) and is isomorphic to $(E,P,Q)$ over $k_s$. This Proposition shows that $(E',P',Q')$ is uniquely determined up to $L$-isomorphism. It is enough to check that $(E \otimes \chi, P, Q)$ satisfies the requested conditions. 
}

\subsection{Comparison with the $X_G$}

Recall the following classical result. 

\prop{Let $k$ be a field with separable closure $k_s$. The functor $G \longmapsto G(k_s)$, from finite \'etale commutative group schemes over $k$ to finite discrete $\mrm{Gal}(k_s/k)$-modules, is an equivalence of categories.}

\cor[group-basis]{Let $k$ be a field with separable closure $k_s$ and $N$ an integer invertible in $k$. Let $G$ be a finite \'etale group scheme over $k$, \'etale-locally isomorphic to $(\Z/N\Z)^{\oplus 2}$. Fix a basis $(A,B)$ of $G(k_s)$ over $\Z/N\Z$, and define, for every $\sigma \in \mrm{Gal}(k_s/k)$, $D_G(\sigma) \in \GL{\Z/N\Z}$ by $D_G(\sigma)\begin{pmatrix}A\\B\end{pmatrix}=\begin{pmatrix}\sigma(A)\\\sigma(B)\end{pmatrix}$, so that $\sigma \longmapsto D_G(\sigma)$ is a group anti-homomorphism. 

Let $E$ be an elliptic curve over $k$. Let $\mathscr{B}_{D_G}(E)$ denote the collection of bases $(P,Q)$ of $E[N](k_s)$ such that for every $\sigma \in \mrm{Gal}(k_s/k)$, one has $D_G(\sigma)\begin{pmatrix}P\\Q\end{pmatrix}=\begin{pmatrix}\sigma(P)\\\sigma(Q)\end{pmatrix}$. Then \[\varphi \in \mrm{Iso}_{k-Gp}(G,E[N]) \longmapsto (\varphi(A),\varphi(B)) \in \mathscr{B}_{D_G}(E)\] is a bijection.}

\prop[xrho-elementary]{Let $N \geq 3, m_0,m_1 \geq 1$ be integers such that $m_0,m_1,N$ are pairwise coprime, $m_0$ is square-free and $m_1N$ is invertible in $k$, and let $m=m_0m_1$. Let $X_{\rho}(N,\Gamma_0(m))$ denote the twist of $X(N,\Gamma_0(m))_k$ by $\rho$ (since $\GL{\Z/N\Z}$ acts on $X(N,\Gamma_0(m))$).

Then, the $j$-invariant, the Atkin-Lehner automorphisms and the degeneracy maps (from Proposition \ref{ext-degn}) can be twisted to finite locally free morphisms between the $X_{\rho}(N,\Gamma_0(m))$ (or with target $\mathbb{P}^1_k$) satisfying the same relations as Proposition \ref{ext-degn}. In particular, the $j$-invariant $X_{\rho}(N,\Gamma_0(m)) \rar \mathbb{P}^1_k$ is well-defined and $X_{\rho}(N,\Gamma_0(m))$ is smooth over $k$ at every point, except, when $m_0$ vanishes in $k$, the points whose $j$-invariant is supersingular. 

For every subextension $L/k$ of $k_s/k$, the points $x \in X_{\rho}(N,\Gamma_0(m))(L)$ with $j(x) \in \mathbb{A}^1_k(L)$ identify with the $L$-isomorphism classes of quadruples $(E,P,Q,C)$ such that:
\begin{itemize}[noitemsep,label=$-$]
\item $E$ is an elliptic curve over $L$ and $P,Q \in E(k_s)$ form a basis of $E[N](k_s)$,
\item $C$ is a cyclic subgroup of $E$ of degree $m$,
\item For all $\sigma \in \mrm{Gal}(k_s/L)$, $\begin{pmatrix}\sigma(P)\\\sigma(Q)\end{pmatrix} = \rho^{-1}(\sigma)\begin{pmatrix}P \\ Q\end{pmatrix}$.
\end{itemize}

Under this identification, the descriptions of the degeneracy maps (resp. the Atkin-Lehner automorphisms) in Proposition \ref{degeneracies} (resp. Proposition \ref{AL-auto}) still apply. }

\demo{The $X(N,\Gamma_0(m))_k$ are proper $k$-schemes of dimension one hence projective \cite[Lemma 0A26]{Stacks}, and they have an action of $\GL{\Z/N\Z}$ for which the degeneracy maps and Atkin-Lehner involutions are equivariant: the existence and regularity of $X_{\rho}(N,\Gamma_0(m))$ and twisted degeneracy maps (as well as their relations) is a consequence of Propositions \ref{cocycle-twist}, \ref{twist-equiv}, and the previous results on the $X(N,\Gamma_0(m))_k$. 

To prove that we identified correctly the rational points of $X_{\rho}(N,\Gamma_0(m))$ (as well as the claims about the behavior of the degeneracy maps and the Atkin-Lehner automorphisms), it is enough to show that every $(E,P,Q,C) \in X_{\rho}(N,\Gamma_0(m))(k_s)$ stable under $\mrm{Gal}(k_s/L)$ is represented by a unique $L$-isomorphism class of $(E',P',Q',C')$ as in the statement of the Proposition. To this end, we argue as in Proposition \ref{rational-Yrho}.

Let $L/k$ be a subextension of $k_s/k$. Let us call $(E,P,Q,C)$ a good quadruple over $L$ if the following conditions are satisfied: 
\begin{itemize}[noitemsep,label=$-$]
\item $E$ is an elliptic curve over $L$ and $P,Q \in E(k_s)$ form a basis of $E[N](k_s)$,
\item $C$ is a cyclic subgroup of $E$ of degree $m$,
\item For all $\sigma \in \mrm{Gal}(k_s/L)$, $\begin{pmatrix}\sigma(P)\\\sigma(Q)\end{pmatrix} = \rho^{-1}(\sigma)\begin{pmatrix}P \\ Q\end{pmatrix}$.
\end{itemize} 

In order to conclude the proof of the Proposition, we only need to show that every non-cuspidal point $(E,P,Q,C) \in X_{\rho}(N,\Gamma_0(m))(k_s)$ stable under $\mrm{Gal}(k_s/L)$ is represented by a unique $L$-isomorphism class of quadruple $(E',P',Q',C')$ which is good for $L$. The argument mimics the proof of Proposition \ref{rational-Yrho}.

Let $(E_1,P_1,Q_1,C_1), (E_2,P_2,Q_2,C_2)$ be two good quadruples over $L$ that are $k_s$-isomorphic. Then there is a $k_s$-isomorphism $\iota: (E_1)_{k_s} \rar (E_2)_{k_s}$ mapping $(P_1,Q_1,C_1)$ to $(P_2,Q_2,C_2)$. For any $\sigma \in \mrm{Gal}(k_s/L)$, the $k_s$-automorphism $\iota^{-1}\circ (\mrm{id} \times \Sp{\sigma^{-1}}) \circ \iota \circ (\mrm{id} \times \Sp{\sigma})$ of $(E_1)_{k_s}$ agrees with the identity on $E_1[N]$ (by the assumptions on $(P_1,Q_1)$ and $(P_2,Q_2)$), hence, by rigidity \cite[Corollary 2.7.4]{KM}, $\iota$ commutes with $\mrm{Gal}(k_s/L)$, hence $\iota$ is defined over $L$ by Corollary \ref{inf-galois-desc}. 

By Proposition \ref{rational-Yrho}, a non-cuspidal $L$-rational point with finite $j$-invariant on $X_{\rho}(N,\Gamma_0(m))$ is given by a $k_s$-isomorphism class of $(E,P,Q,C)$, where $E/L$ is an elliptic curve, $P,Q \in E[N](k_s)$ form a basis of $E[N](k_s)$, one has $\begin{pmatrix}\sigma(P)\\\sigma(Q)\end{pmatrix} = \rho^{-1}(\sigma)\begin{pmatrix}P \\ Q\end{pmatrix}$ for all $\sigma \in \mrm{Gal}(k_s/L)$, and $C \leq E_{k_s}$ is a cyclic subgroup of degree $m$. 

Choose $\sigma \in \mrm{Gal}(L/k_s)$ and let $\rho(\sigma)=\begin{pmatrix}a & b\\c & d\end{pmatrix}$. Then 
\begin{align*}[(E,P,Q,C)] &= \sigma([(E,P,Q,C)]) \\
&= [(E,a\sigma(P)+b\sigma(Q),c\sigma(P)+d\sigma(Q),\sigma(C))] = [(E,P,Q,\sigma(C))],\end{align*} so there is a $k_s$-automorphism $\iota$ of $E_{k_s}$ preserving $(P,Q)$ and mapping $C$ to $\sigma(C)$. By rigidity, $\iota$ is the identity, so $C=\sigma(C)$.

Hence, as a Weil divisor on $E_{k_s}$, $C$ is invariant by $\mrm{Gal}(k_s/L)$, so the Weil divisor attached to $C$ comes from a Weil divisor on $E_L$ (with degree $m$). Since $E_L$ is a regular scheme, $C$ comes from a Cartier divisor $C_0$ on $E_L$ of degree $m$. Since $C=(C_0)_{k_s}$ is a subgroup scheme of $E_{k_s}$, Corollary \ref{inf-galois-desc} implies that $C_0$ is a subgroup scheme of $E$ of degree $L$. Since $C$ is fppf-locally cyclic, by approximation, $C_0$ is fppf-locally cyclic, and $(E,P,Q,C_0)$ is a good quadruple for $L$ that represents $(E,P,Q,C)$. 
}

\prop[group-is-twist-unpolarized]{Let $G$ be a finite \'etale group scheme over $k$, \'etale-locally isomorphic to $(\Z/N\Z)^{\oplus 2}$. Choose a basis $(A,B)$ of $G(k_s)$ and define $D_G$ as in Corollary \ref{group-basis}, so that $\sigma \in \mrm{Gal}(k_s/k) \longmapsto D_G(\sigma)^{-1} \in \GL{\Z/N\Z}$ is a group homomorphism. For each integer $m \geq 1$ coprime to $N$ such that $m=m_0m_1$ with $m_0,m_1 \geq 1$ pairwise coprime, $m_0$ square-free, $m_1$ invertible in $k$, there is a natural isomorphism $\iota_m: X_G(N,\Gamma_0(m)) \rar X_{D_G^{-1}}(N,\Gamma_0(m))$ of projective $k$-schemes. Moreover, 
\begin{itemize}[label=$-$, noitemsep]
\item This isomorphism preserves the $j$-invariant,
\item This family of isomorphims commutes with the Atkin-Lehner automorphisms and all the degeneracy maps,
\item Under this isomorphism, if $L/k$ is a subextension of $k_s/k$ and $E$ is an elliptic curve over $L$, endowed with a cyclic subgroup $C$ of order $m$ and an isomorphism $\varphi: G_L \rar E[N]$, the image under $\iota_m$ of the $\overline{k}$-point $(E,\varphi,C)$ is the $\overline{k}$-point $(E,(\varphi(A),\varphi(B)),C)$ of $Y_{D_G^{-1}}(N,\Gamma_0(m))$. 
\end{itemize} 
}

\demo{Since $(A,B)$ is a $\Z/N\Z$-basis of $G(k_s)$, the morphism of $\mathbb{A}^1_{k_s}$-schemes (where the $j$-invariant provides the structure map) $\tau_m: Y_G(N,\Gamma_0(m))_{k_s} \rar Y(N,\Gamma_0(m))_{k_s}$ defined by \[(E/S,\varphi,C) \longmapsto (E/S,\varphi(A),\varphi(B),C)\] for any $k_s$-scheme $S$ is an isomorphism. Moreover, $\tau_m$ is compatible with the Atkin-Lehner automorphisms and the degeneracy maps. By Proposition \ref{compactification-functor}, it extends to an isomorphism $\overline{\tau}_m: X_G(N,\Gamma_0(m))_{k_s} \rar X(N,\Gamma_0(m))_{k_s}$ which is compatible with the Atkin-Lehner automorphisms and the degeneracy maps. 
 
Let us prove that, for any $\sigma \in \mrm{Gal}(k_s/k)$, the following diagram commutes:
\[
\begin{tikzcd}[ampersand replacement=\&]
X_G(N,\Gamma_0(m)) \times_k \Sp{k_s} \arrow{r}{\overline{\tau}_m}\arrow{d}{(\mrm{id},\underline{\sigma}^{-1})} \& X(N,\Gamma_0(m)) \times_k \Sp{k_s} \arrow{d}{D_G^{-1}(\sigma)\circ (\mrm{id},\underline{\sigma}^{-1})}\\
X_G(N,\Gamma_0(m)) \times_k \Sp{k_s} \arrow{r}{\overline{\tau}_m} \&  X(N,\Gamma_0(m)) \times_k \Sp{k_s}
\end{tikzcd}
\]
Then the uniqueness of the Galois twist will yield the conclusion. It is enough to check the commutation when restricted to $Y_G(N,\Gamma_0(m))$, because of the properties of compactified moduli schemes (see Proposition \ref{compactification-functor}). 

In other words, given two $k$-morphisms $f: S \rar Y_G(N,\Gamma_0(m))$, $g: S \rar \Sp{k_s}$, we need to prove that $D_G^{-1}(\sigma) \circ \mrm{pr}_1 \circ \tau_m(f,g) = \tau_m(f,\underline{\sigma}^{-1} \circ g)$. Let us write $\tau=\tau_m$ for the sake of brevity. 

Suppose that $f$ corresponds to a triple $(E,\varphi,C)$, where $E$ is an elliptic curve over $S$, $C \leq E[m]$ a cyclic locally free subgroup scheme of $E$ of degree $m$, and $\varphi: G_S \rar E[N]$ is an isomorphism. Then the first projection $p_1$ of $\tau(f,g) \in Y(N,\Gamma_0(m))(S)$ corresponds to the triple $(E/S,(\varphi(g^{\ast} A),\varphi(g^{\ast} B)),C)$. Thus $D_G^{-1}(\sigma) \circ p_1 \in Y(N,\Gamma_0(m))(S)$ is the triple 
\[\left(E/S,\varphi\left[g^{\ast}\begin{pmatrix}\sigma^{-1}(A),\sigma^{-1}(B)\end{pmatrix}\right],C\right)=(E/S,\left[\varphi(A \circ \underline{\sigma}^{-1} \circ g),\varphi(B \circ \underline{\sigma}^{-1} \circ g)\right],C),\]

which is exactly the triple given by $\tau(f,\underline{\sigma}^{-1} \circ g)$, and we are done.}

 \lem[inv-cycl-is-constant]{Let $k$ be a field with characteristic not dividing $N$. Let $\omega_k: \mrm{Gal}(k_s/k) \rar (\Z/N\Z)^{\times}$ be the mod $N$ cyclotomic character. Let $\OO_{N,k}$ be the finite \'etale $k$-algebra $k \otimes \Sp{\OO_N}$ (which is a finite \'etale $k$-scheme), and define $\rho: \sigma \in \mrm{Gal}(k_s/k) \longmapsto \mrm{id} \times \underline{\omega_k(\sigma)^{-1}}\in \mrm{Aut}_k(\Sp{\OO_{N,k}})$,  
Then $(\Sp{\OO_{N,k}})_{\rho}$ is naturally isomorphic to the disjoint union of $\varphi(N)$ copies of $k$, naturally indexed by $\mrm{Hom}(\OO_N,k_s)$.
}

\demo{Let $k'$ be the extension of $k$ generated by a primitive $N$-th root of unity $\zeta$. By the construction in the proof of Proposition \ref{cocycle-twist}, $(\Sp{\OO_{N,k}})_{\rho}$ is the spectrum of the sub-algebra of $A = k' \otimes_{\Z} \OO_N$ fixed by the automorphisms $\sigma \otimes \underline{\omega_k(\sigma)}$ for all $\sigma \in \mrm{Gal}(k'/k)$. 
Choose a morphism $\pi:\OO_N \rar k'$, and let, for every $a \in (\Z/N\Z)^{\times}$, $\mu_a$ be the following homomorphism of $k'$-algebras: \[\mu_a: \alpha \otimes \beta \in k' \otimes_{\Z} \OO_N \mapsto \alpha \cdot \pi(\underline{a}(\beta)) \in k'.\] 

Since $A$ is a finite $k'$-algebra, the kernel $\mathfrak{m}_a$ of $\mu_a$ is a maximal ideal of $A$. The $\mu_a$ are linear forms, so, for any $a,a' \in (\Z/N\Z)^{\times}$, either $\mathfrak{m}_a \neq \mathfrak{m}_{a'}$ or $\mu_a$ and $\mu_{a'}$ are collinear. Since the vectors $\begin{pmatrix}\mu_a(1)\\\mu_a(1 \otimes \zeta)\end{pmatrix}=\begin{pmatrix}1 \\\pi(\zeta)^a\end{pmatrix}$ (for $a \in (\Z/N\Z)^{\times}$) are pairwise linearly independent, the $\mathfrak{m}_a$ are pairwise distinct. 

By the Chinese reminder theorem \cite[Lemma 00DT]{Stacks}, the product \[\Pi=\prod_{a \in (\Z/N\Z)^{\times}}{\mu_a}: A \rar \prod_{a \in (\Z/N\Z)^{\times}}{k'}\] is surjective between $k'$-algebras of the same finite dimension, hence is an isomorphism of $k'$-algebras. 

For the action of $\mrm{Gal}(k'/k)$ previously defined on $A$, every $\mu_a$ is easily seen to be equivariant, so $\Pi$ induces an isomorphism of $k$-algebras $A^{\mrm{Gal}(k'/k)} \rar \prod_{a \in (\Z/N\Z)^{\times}}{k}$. 

To conclude, note that the $\pi \circ \underline{a}$ are exactly the ring homomorphisms $\OO_N \rar k_s$.
}

\prop[group-is-twist-polarized]{Consider the situation and notations of Proposition \ref{group-is-twist-unpolarized}, but assume furthermore that $G$ is polarized. Let $\zeta_k: \OO_N \rar k_s$ be the image of the polarization of $G$ at $(A,B)$. Then
\begin{enumerate}[noitemsep,label=(\roman*)]
\item\label{gtp-1} $\det{D_G}$ is the cyclomotic character mod $N$,
\item\label{gtp-2} if we let $D_G^{-1}$ act on $\OO_{N,k}$ through its determinant as in Lemma \ref{inv-cycl-is-constant}, the twist by $D_G^{-1}$ of the map $\mrm{We}: X(N,\Gamma_0(m))_k \rar \Sp{\OO_{N,k}}$ is \[\mrm{We}_{D_G^{-1}}: (X(N,\Gamma_0(m))_k)_{D_G^{-1}} \rar \coprod_{f: \OO_N \rar k_s}{\Sp{k}},\]
\item\label{gtp-3} For $(E,P,Q,C) \in Y(N,\Gamma_0(m))_k(k_s)$, $\mrm{We}_{D_G^{-1}}(E,P,Q,C)$ belongs to the component \[\langle P,\,Q\rangle_{E[N]} \in \mu_N(k_s) \simeq \mrm{Hom}(\OO_N,k_s),\] 
\item\label{gtp-4} the following diagram commutes, where the horizontal arrows are isomorphisms and the bottom line is made with constant $k$-schemes,
\[
\begin{tikzcd}[ampersand replacement=\&]
X_G(N,\Gamma_0(m)) \arrow{r}{\iota_m}\arrow{d}{\det} \& (X(N,\Gamma_0(m))_k)_{D_G^{-1}} \arrow{d}{(\mrm{We})_{D_G^{-1}}}\\
(\Z/N\Z)^{\times} \arrow{r}{u \longmapsto \underline{u}(\zeta_k)} \&  \mrm{Hom}(\OO_N,k_s)
\end{tikzcd}
\]
\item\label{gtp-5} let $d,t$ be divisors of $m$ and $D_{d,t}$ be the associated degeneracy map from Proposition \ref{ext-degn}. The following diagram commutes:
\[
\begin{tikzcd}[ampersand replacement=\&]
X_G(N,\Gamma_0(t)) \arrow{ddd}{\det} \arrow{rrr}{\iota_t} \& \& \& (X(N,\Gamma_0(t))_k)_{D_G^{-1}} \arrow{ddd}{(\mrm{We})_{D_G^{-1}}}\\
\& X_G(N,\Gamma_0(m)) \arrow{ul}{D_{d,t}} \arrow{r}{\iota_m}\arrow{d}{\det} \& (X(N,\Gamma_0(m))_k)_{D_G^{-1}} \arrow{d}{(\mrm{We})_{D_G^{-1}}}\arrow{ur}{(D_{d,t})_{D_G^{-1}}}\&\\
\& (\Z/N\Z)^{\times} \arrow{ld}{\cdot d} \arrow{r}{u \longmapsto \underline{u}(\zeta_k)} \&  \mrm{Hom}(\OO_N,k_s)\arrow{rd}{\circ \underline{d}}\& \\
(\Z/N\Z)^{\times} \arrow{rrr}{u \longmapsto \underline{u}(\zeta_k)} \& \&\& \mrm{Hom}(\OO_N,k_s)
\end{tikzcd}
\]

\end{enumerate}
}

\demo{\ref{gtp-1} follows from the theory of the Weil pairing, \ref{gtp-2} directly from Lemma \ref{inv-cycl-is-constant}, and \ref{gtp-3} is a tautology. For \ref{gtp-4}, the target scheme is a finite disjoint union of copies of $\Sp{k}$, so it is enough to check commutativity on a Zariski-dense subset of $Y_G(N,\Gamma_0(m))(\overline{k})$ (since the bottom arrow is easily seen to be an isomorphism). The explicit description of $\iota_m$ in Proposition \ref{group-is-twist-polarized} shows that the diagram commutes for $\overline{k}$-points of $Y_G(N,\Gamma_0(m))$ whose projection on $Y_G(N)$ is defined over $k_s$. To complete the proof of \ref{gtp-4}, we just note that $Y_G(N,\Gamma_0(m))_k \rar Y_G(N)_k$ is surjective and that $Y_G(N)(k_s)$ is infinite.   
For \ref{gtp-5}, the top cell commutes by Proposition \ref{group-is-twist-unpolarized}, the bottom cell commutes by a direct verification, the outer and inner cells commute by \ref{gtp-4}, the left cell commutes by Proposition \ref{degen-AL-Weil}, which concludes. 
}

\prop[polarized-twist-base-change]{Consider the situation and notations of Proposition \ref{group-is-twist-polarized}, and write $\rho=D_G^{-1}$. Let $K$ be a field extension of $k$ and $K_s$ be a separable closure of $K$. Fix an embedding $e: k_s \rar K_s$. Then $\iota$ defines a homomorphism $R_{e}: \mrm{Gal}(K_s/K) \rar \mrm{Gal}(k_s/k)$. Then, for some canonical morphism $\beta$, the following diagram commutes, and $X_{\rho\circ R_e}(N,\Gamma_0(m)) \overset{\beta}{\rar} X_{\rho}(N,\Gamma_0(m))_K$ is an isomorphism.  
\[
	\begin{tikzcd}[ampersand replacement=\&]
X_{G_K}(N,\Gamma_0(m)) \arrow{rd}{\det}\arrow{rrr}{\iota_{K,m}} \arrow{ddd} \& \& \& X_{\rho\circ R_e}(N,\Gamma_0(m))\arrow{ddd}{\beta}\arrow{ld}{\mrm{We}_{\rho}}\\
\& (\Z/N\Z)^{\times}_K \arrow{d}{=} \arrow[bend left=30]{r}{a \longmapsto \underline{a}(e \circ \zeta_k)} \& \mrm{Hom}(\OO_N,K_s)_K \arrow{d}{e^{-1}} \&\\
\& (\Z/N\Z)^{\times}_k \arrow{r}{a \longmapsto \underline{a}\zeta_k} \& \mrm{Hom}(\OO_N,k_s)_k \& \\
X_G(N,\Gamma_0(m)) \arrow{ur}{\det} \arrow{rrr}{\iota_{k,m}} \& \& \& X_{\rho}(N,\Gamma_0(m)) \arrow{ul}{\mrm{We}_{\rho}}
\end{tikzcd}
\]
Moreover, this diagram commutes to the degeneracy maps as in Proposition \ref{group-is-twist-polarized}. 

}

\demo{Note that all horizontal arrows in the diagram are isomorphisms and that both short vertical arrows are. The bottom and top cells commute by Proposition \ref{group-is-twist-polarized} (because $e(P),e(Q)$ define a basis of the polarized $p$-torsion group $G_K[p]$ over $K$, and that the attached representation is $\rho \circ R_e$). The left cell and the center cell commute by definition. Since the horizontal arrows are isomorphisms, there is a unique $\beta$ such that the outer cell commutes. For this choice of $\beta$, the rightmost cell commutes as well. Moreover, since the outer cell commutes, the horizontal maps are isomorphisms, and the leftmost arrow is the base change $k \rar K$ map, $X_{\rho\circ R_e}(N,\Gamma_0(m)) \overset{\beta}{\rar} X_{\rho}(N,\Gamma_0(m))_K$ is an isomorphism.}

\section{Hecke operators}
\label{hecke-operators}
In this section, we are mainly interested in the proper modular curves. Therefore our base rings will always be regular excellent Noetherian.

\subsection{$\Gamma_0(p)$-degeneracies in characteristic $p$}
\label{char-p}

Let $N,m,p \geq 1$ be pairwise coprime integers, where $p$ is prime and $N \geq 3$. Let $k$ be a perfect field of characteristic $p$ with separable (hence algebraic) closure $k_s$ and $\varphi$ be its absolute Frobenius. Let $G$ be a $N$-torsion group over $k$ and let $\P=\P_G(\Gamma_0(m))$ for the sake of brevity. It is a finite \'etale representable moduli problem over $\Ell_k$ endowed with a natural structure of level $Nm$.

By our previous results, the moduli problem $\P \times [\Gamma_0(p)]_k$ is representable by a finite flat $\MP$-scheme $\mathscr{M}(\P,\Gamma_0(p))=Y_G(N,\Gamma_0(mp))$ (the structure map to $\MP=Y_G(N,\Gamma_0(m))$ being $D_{1,m}$). Moreover, $\P \times [\Gamma_0(p)]_k$ is smooth at infinity and the compactified moduli scheme $\CMPp{\P,[\Gamma_0(p)]_k}$ is smooth over $k$ outside the supersingular $j$-invariants. By \cite[Theorem 13.4.7]{KM}, $Y_G(N,\Gamma_0(mp))$ is the disjoint union of two $\MP$-schemes $\MP_{(1,0)}$ and $\MP_{(0,1)}$ with crossings at the supersingular points. 

First, let us recall the usual definitions of Frobenius morphisms in characteristic $p$.

\defi{Any $\F_p$-scheme $X$ has an \emph{absolute Frobenius} $F_{abs}: X \rar X$, which induces the identity on points and $a \longmapsto a^p$ on rings of functions. 
Given a $S$-scheme $X$ with structure map $\sigma$, where $S$ is a $\F_p$-scheme, we denote by $X^{\varphi}$ the $S$-scheme defined by $X \times_S S$, where the map $S \rar S$ is exactly $F_{abs}$.
The \emph{relative Frobenius} of $X/S$ is then the map $(F_{abs,X},\sigma): X \rar X^{\varphi}$ of $S$-schemes. If $X$ is a relative elliptic curve, it is an isogeny of degree $p$ and the dual isogeny is the \emph{Verschiebung} (this is \cite[Lemma 12.2.1]{KM}). 
Suppose that $F_{abs}: S \rar S$ is an isomorphism, then we can define the $S$-scheme $X^{\varphi^{-1}} = X \times_S S$, where the map $S \rar S$ is exactly $F_{abs}^{-1}$. Then $(X^{\varphi^{-1}})^{\varphi}$ naturally identifies with $X$ (as a $S$-scheme). 
In particular, if in addition $X$ is a relative elliptic curve over $S$, there is an isogeny $X \rar X^{\varphi^{-1}}$ defined as the dual of the Frobenius $X^{\varphi^{-1}} \rar (X^{\varphi^{-1}})^{\varphi} \simeq X$. We will also call it the Verschiebung. 
} 

In particular, we have an Atkin-Lehner automorphism $w_p$ on $\mathscr{M}(\P,\Gamma_0(p))$ and the degeneracy map $D_{p,m}=w_p \circ D_{1,m}: \mathscr{M}(\P,\Gamma_0(p)) \rar \MP$. 

\prop{The automorphism $w_p$ preserves the supersingular locus of $\mathscr{M}(\P,\Gamma_0(p))$ and exchanges the ordinary loci of $\MP_{(1,0)}$ and $\MP_{(0,1)}$. }

\demo{The fact that the supersingular locus is preserved comes from \cite[Corollary 12.3.7]{KM} and the definition of $w_p$. Thus $w_p$ induces an automorphism of the ordinary locus of $\mathscr{M}(\P,\Gamma_0(p))$, which is the disjoint union of the ordinary loci of $\MP_{(0,1)}$ and $\MP_{(1,0)}$. Both are purely inseparable base changes of the smooth projective geometrically connected curve $\CMP_k \rar \Sp{k}$ with finitely many points removed, thus are irreducible schemes. So either $w_p$ preserves the ordinary loci of $\MP_{(0,1)}$ and $\MP_{(1,0)}$, or it exchanges them.

Let $E$ be an ordinary elliptic curve over $k_s$ and $\alpha \in \P(E/k_s)$ be a level structure. Consider the triple $(E/k_s,\alpha,\ker{F})$, where $F: E \rar E^{\varphi}$ is the Frobenius isogeny and $V$ is its dual. Then $x=(E/k_s,\alpha,\ker{F})$ defines a point in the ordinary locus of $\MP_{(1,0)}(k_s)$ by \cite[Proposition 13.4.4, Theorem 13.4.7]{KM}. But $w_p(x)=(E^{\varphi}/k_s, \alpha', \ker{V})$ for some $\alpha' \in \P(E^{\varphi}/k_s)$, and it is easy to see from \cite[(13.4.1-3), Theorem 13.4.7]{KM} that $w_p(x)$ is a point in the ordinary locus of $\MP_{(0,1)}(k_s)$, whence the conclusion. }

\prop[es-divisors]{Let $x \in \MP(k)$ be an ordinary point corresponding to the couple $(E/k,\alpha)$. Then the following equality of effective Cartier divisors holds: 
\[(D_{p,m})_{\ast}(D_{1,m})^{\ast}x = \varphi(x)+p \cdot \varphi'(x),\]
where $\varphi(x) \in \MP(k)$ is the couple $(E^{\varphi}/k, F(\alpha))$ with $F: E \rar E^{\varphi}$ the Frobenius isogeny, and $\varphi'(x) \in \MP(k)$ is the couple $(E^{\varphi^{-1}}/k,V(\alpha))$ with $V: E \rar E^{\varphi^{-1}}$ is the Verschiebung (attached to the elliptic curve $E^{\varphi^{-1}}/k$).}

\demo{By Lemma \ref{pushforward-is-weil}, pushing forward line bundles attached to Cartier divisors supported at normal points is the same as pushing forward said normal points.

Since $x$ is ordinary and the $D_{i,m}$ are finite flat, $(D_{1,m})^{\ast}x$ is a Cartier divisor on $\MPp{\P,[\Gamma_0(p)]_k}$ which is the sum of the inverse images $x_{(1,0)}$ and $x_{(0,1)}$ of $x$ on $\MP_{(1,0)}$ and $\MP_{(0,1)}$. By \cite[(13.5.6)]{KM}, $x_{(1,0)}$ is the Cartier divisor corresponding to a single (smooth) rational point of $\MP_{(1,0)}$; by definition of the $(1,0)$-cyclic moduli problem, one has \[x_{(1,0)}=(E/k,\alpha,\ker{F}) \in \MP_{(1,0)}(k).\] Thus $D_{p,m}(x_{(1,0)})$ is exactly the Cartier divisor corresponding to the $k$-rational point \[\phi(x) = (E^{\varphi},F(\alpha)) \in \MP(k).\]  

On the other hand, by \cite[(13.5.6)]{KM}, $\MP_{(0,1)} \rar \MP$ is a morphism of degree $p$ and induces a bijection on points with values in any perfect $k$-algebra. Since $x$ is an ordinary point, $x_{(0,1)}$ is an effective Cartier divisor of degree $p$, supported on the only point in $y \in \MP_{(0,1)}(k)$ such that $D_{1,m}(y)=x$, so that $x_{(0,1)}=p \cdot y$ as effective Cartier divisors. By definition of the $(0,1)$-cyclic moduli problem, it is easy to see that $(E/k, \alpha,\ker{V}) \in \MP_{(0,1)}(k)$ and $D_{1,m}((E/k, \alpha,\ker{V}))=x$, so that $y=(E/k, \alpha,\ker{V})$. To conclude, it is enough to note that $D_{p,m}(y)=\varphi'(x)$. 
}

\subsection{Hecke operators and the Eichler-Shimura relation}

\defi[has-good-hecke]{A finite flat representable moduli problem $\P$ over $\Ell_R$ ($R$ being a regular excellent Noetherian ring) of level $N$ \emph{has good Hecke operators} if the following conditions are satisfied for every couple $(m_0,m_1)$ of coprime positive integers that are coprime to $N$ and such that $m_0$ is square-free, $m_1 \in R^{\times}$:
\begin{itemize}[noitemsep,label=$-$]
\item $\P \times_R [\Gamma_0(m_0m_1)]_R$ is smooth at infinity. 
\item $\CMPp{\P,[\Gamma_0(m_0m_1)]_R} \rar \Sp{R}$ is proper, surjective, Cohen-Macaulay of relative dimension one with geometrically reduced fibres. 
\end{itemize}}

\lem[formal-has-good-hecke]{\begin{enumerate}[noitemsep,label=$(\roman*)$]
\item\label{fhgh-1} If a finite flat representable moduli problem $\P$ of level $N$ over $\Ell_R$ has good Hecke operators, then $\CMP \rar \Sp{R}$ is a relative curve, and any (compactified) degeneracy map $D_{d,t}$ is finite locally free of constant rank (independent from $\P$ and $R$) by Corollary \ref{flat-deg-infty}.  
\item\label{fhgh-2} If a finite flat representable moduli problem $\P$ of level $N$ over $\Ell_R$ has good Hecke operators and $R \rar R'$ is finite \'etale faithfully flat, then $\P \times [\Sp{R'}]$ has good Hecke operators (for the natural level $N$ structure).  
\item\label{fhgh-3} If a finite flat representable moduli problem $(\P,F)$ of level $N$ over $\Ell_R$ has good Hecke operators and $M$ is a multiple of $N$, then any base change $R \rar R'$ of $(\P,F_{N \rar M})$ has good Hecke operators as long as $R,R'$ have the same invertible integers coprime to $M$.
\item\label{fhgh-4} If a finite flat representable moduli problem $\P$ of level $N$ over $\Ell_R$ has good Hecke operators, so does $\P \times [\Gamma_0(m)]_R$ for any $m=m_0m_1$ with $m_0,m_1 \geq 1$, $m_0$ square-free and $m_1 \in R^{\times}$ for its natural structure of level $Nm_0m_1$ moduli problem.
\item\label{fhgh-5} If $\P$ is a finite \'etale representable moduli problem of level $N$ over $\Ell_R$ and $R$ has generic characteristic zero, then the two conditions are automatically verified for any base change of $R$.  
\item\label{fhgh-6} Let $\P$ be a finite flat representable moduli problem of level $N$ over $\Ell_R$ and $R'$ be a smooth faithfully flat $R$-algebra (hence regular excellent Noetherian) such that $\P_{R'}$ has good Hecke operators. If moreover $\P_{R'}$ is smooth at infinity over $R$, then $\P$ has good Hecke operators. 
\end{enumerate}}

\demo{\ref{fhgh-1} follows from Lemma \ref{good-base-curve} and Corollary \ref{flat-deg-infty}, \ref{fhgh-2} from Lemma \ref{with-extra-functions}, and \ref{fhgh-5} from (mostly) Proposition \ref{struct-gamma0m-smooth} and Corollary \ref{smooth-at-infinity-gamma0m}. \ref{fhgh-3} and \ref{fhgh-4} are formal verifications. Finally, \ref{fhgh-6} follows from Proposition \ref{smooth-at-infinity-descends} and Proposition \ref{struct-gamma0m-smooth}.
}

\cor[basic-has-good-hecke]{Let $R$ be a regular excellent Noetherian $\Z[1/N]$-algebra with $N \geq 3$ and $G$ be a $N$-torsion group over $R$. Let $m_0,m_1 \geq 1$ be coprime integers with product $m$ coprime to $N$ such that $m_0$ is square-free and $m_1 \in R^{\times}$. Then $\P_G(\Gamma_0(m))$ (endowed with its natural level $Nm$ structure) and $[\Gamma_1(N)]_R \times_R [\Gamma_0(m)]_R$ (endowed with its natural level $Nm$ structure, assuming that $N \geq 4$) have good Hecke operators. }

\defi[not-jacobian]{In the above situation, when $m_0 \in R^{\times}$, the relative curve is smooth. We denote its relative Jacobian\footnote{As noted in the introduction of this chapter, the notion of relative Jacobian for a proper smooth scheme $X \rar \Sp{R}$ of relative dimension one does not seem to be standard when the morphism does not have geometrically connected fibres, even when $R$ is a field. However, only minor modifications to the classical theory are needed to define and construct Jacobians in this setting, which we describe in Section \ref{degrees-jacobians}.} by $J_G(N,\Gamma_0(m))$ (resp. $J_1(N,\Gamma_0(m))$). }

\defi[picard-functors-definition-recopied]{We copy for convenience Definition \ref{picard-functors-definition}. If $X$ is a proper finitely presented $S$-scheme, $P_{X/S}$ denotes the functor on $\mathbf{Sch}_S$ given by $T \longmapsto \mrm{Pic}(X \times_S T)$. The \emph{relative Picard functor} of $X/S$ is the sheafification $\operatorname{Pic}_{X/S}$ of $P_{X/S}$ in the fppf topology (see Definition \ref{fppf-topology-definition}). If moreover the fibres of $X \rar S$ are pure of dimension one, we define $P^0_{X/S}$ as the sub-presheaf of $P_{X/S}$ such that $P_{X/S}(T)$ is made with the classes of those line bundles $\mathcal{L}$ on $X \times_S T$ such that for every $t \in T$, the line bundle $t^{\ast}\mathcal{L}$ on $X_{\kappa(t)}$ has degree zero on every connected component. The fppf-sheafification of $P^0_{X/S}$ is the subsheaf $\operatorname{Pic}^0_{X/S}$ of $\operatorname{Pic}_{X/S}$. 
}

\defi[hecke-exists]{Let $\P$ be a finite flat representable moduli problem of level $N$ over $\Ell_R$ which has good Hecke operators and $p \nmid N$ be a prime. The Hecke operator $T_p$ denotes the endomorphism $(D_{p,1})_{\ast}(D_{1,1})^{\ast}$ of any of the functors $P_{\CMP/\Sp{R}}$, $P^0_{\CMP/\Sp{R}}$, $\operatorname{Pic}_{\CMP/\Sp{R}}$, $\operatorname{Pic}^0_{\CMP/\Sp{R}}$. }

\lem[hecke-natural]{With the same notations as before:
\begin{enumerate}[noitemsep,label=$(\roman*)$]
\item The formation of $T_p$ commutes with arbitrary base change (of regular excellent Noetherian algebras). 
\item If the level $N$ structure on $\P$ is given by the functor $F$ and $M$ is any multiple of $N$ coprime to $p$, then the operator $T_p$ for the level $M$ structure $F_{N \rar M}$ on $\P$ is the same.
\item If $m=m_0m_1$ is coprime to $Np$, with $m_0$ square-free and coprime to $m_1 \in R^{\times}$, then, for any $d,t \geq 1$ with $dt \mid m$ (resp. any $d \geq 1$ with $d\mid m$ coprime to $m/d$), the formation of $T_p$ commutes to the push-forwards and pull-backs by $D_{d,t}: \CMPp{\P,[\Gamma_0(m)]_R} \rar \CMPp{\P,[\Gamma_0(t)]_R}$ (resp. to the push-forwards and pull-backs by the Atkin-Lehner automorphisms $w_d$ of $\CMPp{\P,[\Gamma_0(m)]_R}$). 
\item Let $C: \P \Rightarrow \P'$ be a flat natural transformation of finite flat representable moduli problems over $\Ell_R$ of level $N$ (for instance, $C$ can be $[a]$ for some $a \in (\Z/N\Z)^{\times}$). Then $C: \CMP \rar \CMPp{\P'}$ is finite locally free, and the formation of the Hecke operator $T_p$ commutes with the pull-back $C^{\ast}$ and the push-forward $[C_{\ast}]$. 
\end{enumerate}}

\demo{The first two points are clear. The last three points are applications of Proposition \ref{relative-picard-functoriality} and \ref{relative-picard-mixed-functoriality}. The Cartesian diagram away from the cusps is given by Propositions \ref{degn-cartesian} and \ref{two-degn-cartesian}, the required flatness follows from Proposition \ref{flatness-at-cusps} for $C$ and Corollary \ref{flat-deg-infty} for the degeneracy maps.

So we need to show that for any $R$-scheme $T$, $\iota: \MPp{\P,[\Gamma_0(m)]_R}_T \rar \CMPp{\P,[\Gamma_0(m)]_R}_T$ has scheme-theoretically dense image. The second assertion is true because $\iota$ is a finite flat base change of the scheme-theoretically dense inclusion $\mathbb{A}^1_T \rar \mathbb{P}^1_T$. 
}

\cor[hecke-weil]{Let $\P=\P_G(\Gamma_0(m))$, where $N,m_0,m_1 \geq 1$ are coprime integers, $N \geq 3$, $m_0$ is square-free, $R$ is a regular excellent Noetherian $\Z[(Nm_1)^{-1}]$-algebra, $m=m_0m_1$, and $G$ is the constant $N$-torsion group over $R$. Let $p \nmid Nm$ be a prime. 
Then $T_p$ commutes with the action of $\GL{\Z/N\Z}$. }
 
 \demo{$\GL{\Z/N\Z}$ acts by morphisms of moduli problems of level $Nm$. }

\prop[hecke-det]{Let $N,m_0,m_1 \geq 1$ be pairwise coprime integers, with $N \geq 3$ and $m_0$ square-free. Let $R$ be a regular excellent Noetherian $\Z[(Nm_1)^{-1}]$-algebra and $G$ be a polarized $N$-torsion group over $R$. Let $m=m_0m_1$ and $\P=\P_G(N,\Gamma_0(m))$. Let, for each $u \in (\Z/N\Z)^{\times}$, $\P_u=\P_{G,u}(N,\Gamma_0(m))$ and $X_{\P,u}$ be the associated compactified moduli scheme. Then there is an isomorphism of functors \[\iota: P_{\CMP/\Sp{R}}\simeq \prod_{u \in (\Z/N\Z)^{\times}}{P_{X_{\P,u}/\Sp{R}}}\] (preserving the $P^0$). Through this isomorphism, the operator $T_p$ decomposes as a direct product over $u \in (\Z/N\Z)^{\times}$ of the morphisms $(D_{p,1})_{\ast}(D_{1,1})^{\ast}: P_{X_{\P,u}/\Sp{R}} \rar P_{X_{\P,pu}/\Sp{R}}$, also denoted by $T_p$, and preserving the $P^0$. 

In particular, $\iota$ sheafifies to an isomorphism \[\iota: \operatorname{Pic}_{\CMP/\Sp{R}}\simeq \prod_{u \in (\Z/N\Z)^{\times}}{\operatorname{Pic}_{X_{\P,u}/\Sp{R}}}\] preserving the $\operatorname{Pic}^0$, and, through this isomorphism, $T_p$ is the direct sum over $u \in (\Z/N\Z)^{\times}$ of the morphisms $(D_{p,1})_{\ast}(D_{1,1})^{\ast}: \operatorname{Pic}_{X_{\P,u}/\Sp{R}} \rar \operatorname{Pic}_{X_{\P,pu}/\Sp{R}}$ .}

\demo{$D_{1,m}$ preserves the map to the constant scheme $(\Z/N\Z)^{\times}$, while $D_{p,m}$ multiplies it by $p$. We can then mimic the construction of the Hecke operators by considering instead, for each $u \in (\Z/N\Z)^{\times}$, the maps \[D^u_{1,m}: X_{G,u}(N,\Gamma_0(mp)) \rar X_{G,u}(N,\Gamma_0(m)),\,D^u_{p,m}: X_{G,u}(N,\Gamma_0(pm)) \rar X_{G,pu}(N,\Gamma_0(m)).\]}

\prop[hecke-commute]{Let $\P$ be a finite flat representable moduli problem over $\Ell_R$ of level $N$ having good Hecke operators. Let $p,q$ be distinct primes not dividing $N$. Then $T_p$ and $T_q$ commute. }

\demo{By Proposition \ref{hecke-natural}, the following diagram of pre-sheaves commutes (with all applications preserving $P^0$):
\[
\begin{tikzcd}[ampersand replacement=\&]
P_{\CMPp{\P,[\Gamma_0(p)]_R}/\Sp{R}} \arrow{r}{T_q}\arrow{d}{[u_{\ast}]} \& P_{\CMPp{\P,[\Gamma_0(p)]_R}/\Sp{R}} \arrow{d}{[u_{\ast}]}\\
P_{\CMP/\Sp{R}} \arrow{r}{T_q}\arrow{d}{u^{\ast}} \& P_{\CMP/\Sp{R}} \arrow{d}{u^{\ast}}\\
P_{\CMPp{\P,[\Gamma_0(p)]_R}/\Sp{R}} \arrow{r}{T_q} \& P_{\CMPp{\P,[\Gamma_0(p)]_R}/\Sp{R}} 
\end{tikzcd}
\]
where $u$ is either $D_{p,m}$ or $D_{1,m}$, whence the conclusion. 
}

\prop[hecke-bad-AL]{Let $\P$ be a finite flat relatively representable moduli problem over $\Ell_R$ of level $L$ prime to $N$, where $R$ is a regular excellent Noetherian $\Z[1/N]$-algebra, and assume that $\P$ is representable if $N \leq 3$ and that $\P'=\P \times_R [\Gamma_1(N)]_R$ has good Hecke operators. Let $d \mid N$ be a positive integer coprime to $N/d$, and $e \mid N/d$ be an integer. Let $p \nmid LN$ be a prime number, and let $p' \in (\Z/de\Z)^{\times}$ be congruent to $p$ mod $d$ and to $1$ mod $e$. Then, on $P_{\CMPp{\P,[\Gamma_1(N)]_R} \times \mu_{de}^{\times}/\Sp{R}}$, one has $([p'])_{\ast}(w'_d)_{\ast}T_p=T_p(w'_d)_{\ast}$.}

\demo{This is a consequence of Proposition \ref{relative-picard-mixed-functoriality}. The Cartesian diagrams are given in Proposition \ref{bad-AL-good-deg}, and the maps are finite locally free by Corollary \ref{flat-deg-infty}.}

\prop[eichler-shimura]{Suppose that $R=\F_p$ and $G$ is a constant $N$-torsion group over $R$. Let $m \geq 1$ be an integer coprime to $Np$. Then the Jacobian $J_G(N,\Gamma_0(m))$ of $X_G(N,\Gamma_0(m))$ is endowed with a Frobenius isogeny $F$ (its absolute Frobenius) of degree $p$ with dual isogeny $F'$, and the following identity holds on $J_G(N,\Gamma_0(m))$: 
\[T_p=F+[p] \cdot F'.\]}

\demo{First, recall that $\mrm{We}: X_G(N,\Gamma_0(m))_{\F_p} \rar \Sp{\OO_N \otimes \F_p}:= T$ is the Stein factorization of $X_G(N,\Gamma_0(m)) \rar \Sp{\F_p}$. Let $x,x' \in Y_G(N,\Gamma_0(m))(\overline{\F_p})$ be two ordinary points with $\mrm{We}(x)=\mrm{We}(x')$ and let us show that on the divisor $D=x-x'$ of $X_G(N,\Gamma_0(m))_{\overline{\F_p}}$, one has $T_p(D)=(F+[p]\cdot F')(D)$.

By Proposition \ref{es-divisors}, we know that $T_p(D)=\varphi(x)-\varphi(x')+p\cdot \varphi'(x)-p\cdot \varphi'(x')$ as Cartier divisors, where we extend $\varphi$ and $\varphi'$ by linearity. The goal is thus to show that (in the notation of this Proposition) $F(x-x')$ represents the Cartier divisor $\varphi(x)-\varphi(x')$ and that $[p]\cdot F'(x-x')$ represents the Cartier divisor $p\cdot(\varphi'(x)-\varphi'(x'))$. For the first part, note that $\varphi(x)$ (resp. $F(x-x')$) is exactly the pre-composition of the geometric point $x$ (resp. $x-x'$) by the absolute Frobenius of $\overline{\F_p}$. 

For the second part, we can write $x=\varphi(x_1),x'=\varphi(x'_1)$ for some ordinary geometric points $x_1,x'_1 \in Y_G(N,\Gamma_0(m))(\overline{\F_p})$, then one checks directly that $\varphi'(x)=[p] \cdot x_1, \varphi'(x')=[p] \cdot x'_1$, and $F'(x-x')=F'(F(x_1-x'_1))=p \cdot (x_1-x'_1)$, whence the conclusion.  

Now, $T_p$ and $F+[p] \cdot F'$ are two endomorphisms of $J_G(N,\Gamma_0(m))$, so there is a greatest closed subscheme $Z$ of $J_G(N,\Gamma_0(m))$ on which they coincide, and $Z$ is a subgroup scheme of $J_G(N,\Gamma_0(m))$. Since $J_G(N,\Gamma_0(m))$ is reduced, it is enough to show that $Z$ contains a Zariski-dense set of $\overline{F_p}$-points of $J_G(N,\Gamma_0(m))$. 

By Proposition \ref{difference-map-generates-jacobian}, there exists $r \geq 1$, such that the sum of $r$ difference maps \[\delta^r: (X_G(N,\Gamma_0(m))\times_T X_G(N,\Gamma_0(m)))^{\times_{\F_p} r} \rar J_G(N,\Gamma_0(m))\] is surjective. Let $U \subset X_G(N,\Gamma_0(m))_{\F_p}$ be the ordinary locus: then $V := (U \times_T U)^{\times_{\F_p} r}$ is dense in $(X_G(N,\Gamma_0(m))\times_T X_G(N,\Gamma_0(m)))^{\times_{\F_p} r}$, hence the image in $J_G(N,\Gamma_0(m))$ of $V(\overline{\F_p})$ (which is contained in $Z$) is Zariski-dense, whence the conclusion.  
}

\cor[tate-eichler-shimura]{Assume that $R=\Q$ and $G$ is a constant $N$-torsion group over $R$. Let $\ell$ be a prime and $p \nmid N\ell$ be another prime. Then, on the $\ell$-adic Tate module of $J_G(N)$, one has
\[\Fr_p^2-T_p \Fr_p + p[p]=0.\]}

\demo{Let $G_0$ be the constant $N$-torsion group on $\Z[1/N]$. Then $J_G(N)$ is the generic fibre of the Abelian scheme $J_{G_0}(N)$ (to which $T_p$ extends), so the isomorphism $\Tate{\ell}{J_G(N)} \rar \Tate{\ell}{J_{G_0}(N)_{\F_p}}$ given by the reduction modulo $p$ preserves $T_p$ and is equivariant for the (unramified) Galois action of a decomposition group $D_p \leq G_{\Q}$ above $p$. 
The conclusion follows from Proposition \ref{eichler-shimura}.}

\subsection{Link with the point of view of Galois twists}

Using the previous results from this section, the results on Jacobians (Propositions \ref{jacobian-relative-curve-exists}, \ref{jacobian-twist}, \ref{jacobian-twist-functoriality}) as well as the previous results about twists (Propositions \ref{group-is-twist-unpolarized}, \ref{group-is-twist-polarized}), the following statements are formal.

\prop[jac+hecke-twist]{Let $N \geq 3$ and $k$ be a field of characteristic not dividing $N$ with separable closure $k_s$. Let $G$ be a $N$-torsion group over $k$. Choose a basis $(A,B)$ of $G(k_s)$ and define $D_G$ as in Corollary \ref{group-basis}, so that $\sigma \in \mrm{Gal}(k_s/k) \longmapsto D_G(\sigma)^{-1} \in \GL{\Z/N\Z}$ is a group homomorphism. \smallskip

Let $m_0,m_1 \geq 1$ such that $m_0,m_1,N$ are pairwise coprime to $N$, $m_1$ is invertible in $k$ and $m_0$ is square-free, and $m=m_0m_1$. The isomorphism $\iota_m: X_G(N,\Gamma_0(m)) \rar (X(N,\Gamma_0(m))_k)_{D_G^{-1}}$ of Proposition \ref{group-is-twist-unpolarized} extends by push-forward functoriality to an isomorphism of the Jacobians $[(\iota_m)_{\ast}]: J_G(N,\Gamma_0(m)) \rar (J(N,\Gamma_0(m))_k)_{D_G^{-1}}$, where $(J(N,\Gamma_0(m))_k)_{D_G^{-1}}$ is the Galois twist of the Jacobian $J(N,\Gamma_0(m))_k$ where $\mrm{Gal}(k_s/k)$ acts by the push-forward action of $D_G^{-1}$. \smallskip

The collection of the $[(\iota_m)_{\ast}]$ (over all possible $m$) commutes with the push-forwards and the pull-backs of the twists by $\rho$ of the degeneracy maps and the Atkin-Lehner automorphisms. Moreover, for any prime $p \nmid Nm$, the finite flat maps $D_{1,m}, D_{p,m}: X_{\rho}(N,\Gamma_0(mp)) \rar X_{\rho}(N,\Gamma_0(m))$ induce an endomorphism $T_p=(D_{p,m})_{\ast}(D_{1,m})^{\ast}$ of $J_{\rho}(N,\Gamma_0(m))$, and $[(\iota_m)_{\ast}]$ commutes to the action of $T_p$.}

\prop[jac+hecke-twist-polarized]{Let us keep the notations of Proposition \ref{jac+hecke-twist}. Assume furthermore that $G$ is polarized, and let $\zeta_k: \OO_N \rar k_s$ be the image of the polarization of $G$ at $(A,B)$. For each $\xi \in \mu_N^{\times}(k_s)$, let $X_{\xi}(N,\Gamma_0(m))$ be the inverse image of $\xi$ under the morphism \[(\mrm{We})_{D_G^{-1}}: (X(N,\Gamma_0(m))_k)_{D_G^{-1}} \rar \coprod_{f: \OO_N \rar k_s}{\Sp{k}}\] defined in Proposition \ref{group-is-twist-polarized}. 

Then, for each $u \in (\Z/N\Z)^{\times}$, the isomorphism $\iota_{m,u}: X_{G,u}(N,\Gamma_0(m)) \rar X_{\underline{u}\zeta_k}(N,\Gamma_0(m))$ induced by $\iota_m$ extends to an isomorphism $[(\iota_{m,u})_{\ast}]$ of their Jacobians, such that, for all primes $p \nmid Nm$ invertible in $k$, the following diagrams commute:

\[
\begin{tikzcd}[ampersand replacement=\&]
J_G(N,\Gamma_0(m)) \arrow{d}{=} \arrow{r}{[(\iota_m)_{\ast}]} \& J_{\rho}(N,\Gamma_0(m)) \arrow{d}{=}\\
\prod_{u \in (\Z/N\Z)^{\times}}{J_{G,u}(N,\Gamma_0(m))} \arrow{r}{[(\iota_{m,u})_{\ast}]}\& \prod_{u \in (\Z/N\Z)^{\times}}{J_{\underline{u}\zeta_k}(N,\Gamma_0(m))} \\
J_{G,u}(N,\Gamma_0(m)) \arrow{r}{[(\iota_{m,u})_{\ast}]} \arrow{d}{T_p=[(D_{p,m})_{\ast}]D_{1,m}^{\ast}}\& J_{\underline{u}\zeta_k}(N,\Gamma_0(m)) \arrow{d}{T_p=[(D_{p,m})_{\ast}](D_{1,m})^{\ast}}\\
J_{G,pu}(N,\Gamma_0(m)) \arrow{r}{[(\iota_{m,pu})_{\ast}]}\& J_{\underline{up}\zeta_k}(N,\Gamma_0(m))
\end{tikzcd}
\]  
}

\subsection{The bad Hecke operators for $X_1(N,\Gamma_0(m))$}
\label{bad-hecke}

\prop[bad-has-good-hecke]{Let $N,r \geq 1$ with $r \mid N$ and $N \geq 4$ and $R$ be a regular excellent Noetherian $\Z[1/N]$-algebra. Then $\P=[\Gamma_1'(N,r)]_R$ with its level $N$ structure admits good Hecke operators.}

\demo{By base change, we may assume $R$ is a localization of $\Z[1/N]$. Since this moduli problem is finite \'etale by Proposition \ref{basic-with-bad}, the conditions are verified by Proposition \ref{struct-gamma0m-smooth}.}

\prop[bad-pb-with-geom-fibres]{If $N,r \geq 1$ are such that $r \mid N$, $N \geq 4$ and $R$ is a regular excellent $\Z[1/N]$-algebra, then, for any coprime integers $m_0,m_1 \geq 1$ coprime to $N$ such that $m_1 \in R^{\times}$ and $m_0$ square-free, the fibres of $\CMPp{[\Gamma'_1(N,r)] \times [\Gamma_0(m)]} \rar \Sp{R}$ are geometrically connected.}

\demo{It is enough to prove that $\CMPp{[\Gamma_1'(N,r)]_{R},[\Gamma_0(m)]_{R}} \rar \Sp{R}$ has connected geometric fibres when $R$ is a localization of $\Z[1/N]$ (and for any suitable $m$). By \cite[Lemma 0E0N]{Stacks}, we may assume that $R=\Q$, and it is then enough to show that for any $m \geq 1$ coprime to $N$, $\CMPp{[\Gamma_1'(N,r)]_{\Q},[\Gamma_0(m)]_{\Q}} \rar \Sp{\Q}$ is geometrically connected. 

There is a natural map $Y(N,\Gamma_0(m))_{\Q} \rar \MPp{[\Gamma_1'(N,r)]_{\Q},[\Gamma_0(m)]_{\Q}} \times_{\Q} \Sp{\Q(\mu_N)}$ (mapping a basis $(P,Q)$ of the $N$-torsion to the triple $(P,\langle N/r \cdot Q\rangle,\langle P,\,Q\rangle)$, the latter coordinate being the Weil pairing, and preserving the $\Gamma_0(m)$-structure) and we can easily check (by passing to $\Qbar$) that it is surjective. Thus, the compactification map $X(N,\Gamma_0(m))_{\Q} \rar \CMPp{[\Gamma_1'(N,r)]_{\Q},[\Gamma_0(m)]_{\Q}}_{\Q(\mu_N)}$ of $\Q(\mu_N)$-schemes is surjective. Its source is geometrically connected (as a $\Q(\mu_N)$-scheme), hence $\CMPp{[\Gamma_1'(N,r)]_{\Q},[\Gamma_0(m)]_{\Q}}_{\Q(\mu_N)} \times_{\Q(\mu_N)} \Sp{\Qbar}$ is connected, hence $\CMPp{[\Gamma_1'(N,r)]_{\Q},[\Gamma_0(m)]_{\Q}}$ is geometrically connected as a $\Q$-scheme, whence the conclusion. 
 }

\prop[bad-deg-cst]{In the situation of Proposition \ref{bad-degeneracy}, $D'_{s,t}$ has positive constant rank, independent from $\P$. }

\demo{Suppose that whenever $k \geq 3$ is coprime to $N$ and invertible in $R$, $D'_{s,t}$ for the moduli problem $[\Gamma(k)]_R$ is finite \'etale of constant rank $r_k > 0$. By Proposition \ref{bad-degn-cartesian}, the same holds for any moduli problem on $\Ell_R$ with a morphism to $[\Gamma(k)]_R$. In particular, if $k,k' \geq 3$ are invertible in $R$ and coprime to $N$, $r_{k}=r_{kk'}=r_{k'}$. By considering the case where $R=\Q$, we then see that all the $r_{k}$ are equal to some $\rho > 0$. 

Let $\P$ be a relatively representable moduli problem on $\Ell_R$ (representable if $N \leq 3$) of level coprime to $N$. Choose two coprime integers $k,k' \geq 3$ coprime to $N$, then $D'_{s,t}$ (for $\P$) has an \'etale cover by the $D'_{s,t}$ for the two representable moduli problems $\P \times [\Gamma(k)]_{R[1/k]}$ and $\P \times [\Gamma(k')]_{R[1/k']}$. These instances of $D'_{s,t}$ are, by Proposition \ref{bad-degn-cartesian}, base changes of instances of $D'_{s,t}$ for $[\Gamma(k)]_{\Z[1/k]},[\Gamma(k')]_{\Z[1/k]}$, so they are both finite \'etale of rank $\rho$ by the above. By descent, it follows that $D'_{s,t}$ for $\P$ is finite \'etale of rank $\rho$. 

 So all we have to do is show that $\MPp{[\Gamma(k)]_R,[\Gamma_1'(N,r)]_R}$ is connected. The morphism \[(E/S,(P,Q)) \in Y(Nk)_{\Z[1/Nk]}  \longmapsto (E/S,(kP,kQ),(NP,\langle \frac{kN}{r} \cdot Q\rangle)) \in \MPp{[\Gamma(k)]_R,[\Gamma'_1(N,r)]_R}\] is surjective (it is enough to prove it when $R$ is a field), so $\MPp{[\Gamma(k)]_R,[\Gamma_1'(N,r)]_R}$ is connected by Corollary \ref{XN-Weil}.}

\defi[Up-definition]{Let $N \geq 4, r,m_0,m_1 \geq 1$ be integers, with $r \mid N$, $N,m_0,m_1$ pairwise coprime and $m_0$ square-free. Let $R$ be a regular excellent Noetherian $\Z[(Nm_1)^{-1}]$-algebra and $p$ be a prime number be such that $pr \mid N$. For any divisor $i \geq 1$ of $N$, let $Y_1'(i)$ (resp. $X'_1(i)_R$) be the moduli scheme (resp. the compactified moduli scheme) attached to the finite \'etale representable moduli problem $[\Gamma_1'(N,i)]_R \times_R [\Gamma_0(m)]_R$. For $i \in \{1,p\}$, the degeneracy map $D'_{i,r}: Y'_1(pr)_R \rar Y'_1(r)_R$ extends to a map $D'_{i,r}: X'_1(pr)_R \rar X'_1(r)_R$. 

The operator $U_p$ is the endomorphism $(D'_{p,r})_{\ast}(D'_{1,r})^{\ast}$ of the functor $P_{X_1'(r)_R/\Sp{R}}$. In particular, $U_p$ preserves $P^0$ and $U_p$ extends to an endomorphism $U_p$ of the sheaf $\operatorname{Pic}_{X_1'(r)_R/\Sp{R}}$ preserving its subsheaf $\operatorname{Pic}^0_{X_1'(r)_R/\Sp{R}}$. }

\demo{By Proposition \ref{relative-picard-functoriality}, we only need to show that $D'_{1,r}, D'_{p,r}$ are finite locally free morphisms of relative curves over $R$. 
We know that the $X'_1(N,u,\Gamma_0(m))$ are relative curves over $\Sp{R}$ with geometrically connected fibres by Proposition \ref{bad-has-good-hecke}. 
Hence, it is enough to show that the $D'_{\ast,r}$ are finite locally free. 
They are the compactifications of finite locally free morphisms of constant rank by Corollaries \ref{bad-degn-fin-etale}, \ref{bad-deg-cst}, so they are finite, and the flatness follows from Proposition \ref{flatness-at-cusps} since the moduli problems are smooth at infinity.}

\prop[bad-hecke-natural]{With the notations of Definition \ref{Up-definition}:
\begin{itemize}[noitemsep,label=$-$]
\item The formation of $U_p$ commutes with base change. 
\item If $p,q$ are distinct primes dividing $N/r$, then $U_p$ and $U_q$ commute.
\item $U_p$ commutes with any push-forward or pull-back by $D_{u,v}$ (where $uv \mid m$). 
\item In particular, $U_p$ commutes with any $T_q$ for any $q$ coprime to $Nm$. 
\end{itemize}
} 

\demo{The first point is formal. The rest follows from Propositions \ref{relative-picard-functoriality}, \ref{relative-picard-mixed-functoriality}. To apply these results, we first need to check that all the maps involved indeed extend and are finite flat of constant rank, which is the case by Propositions \ref{degn-AL-infty}, \ref{bad-degn-AL-infty}, \ref{bad-deg-cst} and Corollary \ref{flat-deg-infty}. We also need some diagrams to be Cartesian above dense open subschemes, which is proved in the Propositions \ref{mixed-degn-cartesian} and \ref{two-bad-degn-cartesian}.}

\section{The cuspidal subscheme and the $q$-expansion principle}
\label{cuspidal-subscheme}

\nott{The ring $\Z((q))$ is the ring of Laurent series with coefficients in $\Z$, or, in other words, $(\Z[[q]])[q^{-1}]$. 
Let, for any even $k \geq 4$, \[E_k = 1-\frac{2k}{B_k}\sum_{n \geq 1}{\left(\sum_{d \mid n}{d^{k-1}}\right)q^n} \in \Z[[q]] \otimes \Q\] denote the normalized Eisenstein series coming from the classical theory of modular forms \cite[Section 1.1, Exercise 1.1.7]{DS}. One has $\frac{-2 \cdot 4}{B_4} = 240$, $\frac{-2 \cdot 6}{B_6} = -504$, so $E_4, E_6 \in \Z[[q]]$. }

\defi{The \emph{Tate curve} $\mrm{Tate}(q)$ is the elliptic curve over $\Z((q))$ with Weierstrass equation 
\[y^2 +xy= x^3+\left(-5\frac{E_4-1}{240}\right)x - \frac{1}{12}\left(5\frac{E_4-1}{240}+7\frac{E_6-1}{-504}\right)\] as in \cite[(A1.2.3)]{Katz}, endowed with the nowhere-vanishing differential $\omega = \frac{dx}{2y+x}=\frac{dy}{3x^2-\frac{5}{12}\frac{E_4-1}{240}+y}$.  
It satisfies the following properties: 
\begin{itemize}[noitemsep,label=$-$]
\item Its $j$-invariant is given by the series $j(q)=\frac{1}{q}+744+\cdots$ from classical modular forms (see for instance \cite[p. 7]{DS}). 
\item The invariant $c_4$ (resp. $c_6,\Delta$) attached to this Weierstrass equation is $E_4$ (resp. $-E_6$, resp. the classical series $\Delta(q)=q\prod_{n \geq 1}{(1-q^n)^{24}}$).  
\item For every $N \geq 1$, there is a short exact sequence of finite flat commutative group schemes over $\Z((q))$ 
\[0 \rar \mu_N \overset{\alpha}{\rar} \mrm{Tate}(q)[N] \overset{\beta}{\rar} \Z/N\Z \rar 0\] such that, for any $\Z((q))$-algebra $R$, $\zeta \in \mu_N(R)$, $P \in \mrm{Tate}(q)[N](R)$, one has \[\langle\alpha(\zeta),\,P\rangle_{\mrm{Tate}(q)[N]} = \zeta^{\beta(P)}.\] 
\end{itemize}
This is the same elliptic curve which is called $\mrm{Tate}(q)$ in \cite[(8.8)]{KM}.}

\rem{The computation of $c_6$ does not match the one given in \cite[(8.8)]{KM}: in \emph{loc.cit.}, the authors write that $c_6$-invariant of the Tate curve is $E_6$. This seems to be a mistake (rather than a different definition for the Tate curve): in particular, it implies by \cite[Theorem V.5.3]{AEC2} that the specialization of the Tate curve at any $q \in (3\Z_3)^3$ has \emph{nonsplit} multiplicative reduction. Thus \cite[Theorem V.3.1]{AEC2} contradicts the claim \cite[(8.8, T3)]{KM} about the torsion subscheme of the Tate curve. Given that \emph{loc.cit.} refers to \cite{DeRa, Katz} and that these both references use the same exact sequence for the torsion of the Tate curve (for instance \cite[Construction 1.15]{DeRa}), it seems reasonable that the correction to make is to clarify in \cite[(8.8)]{KM} that the $c_6$-invariant of the Tate curve is $-E_6$.}

\rem{Let $R$ be a regular excellent Noetherian ring and $\P$ be a finite representable moduli problem over $\Ell_R$ normal near infinity. Then the formal completion of $\CMP$ along its cuspidal subscheme is, by \cite[Lemma 8.11.2]{KM}, the normalization of $R[[1/j]]$ in the normal finite $R((1/j))$-scheme $\MP \times_{j, \mathbb{A}^1}\Sp{R((1/j))}$. Note that there exists a canonical continuous isomorphism $R[[1/j]] \simeq R[[q]]$ defined by $j = j(\mrm{Tate}(q))$. Thus, the formal completion of $\CMP$ along its cusps is naturally a $R[[q]]$-scheme. }

\prop[tate-curve-and-cusps]{Let $R$ be a regular excellent Noetherian ring and $\P$ be a finite representable moduli problem over $\Ell_R$ normal near infinity. Then the formal completion of $\CMP$ along its cuspidal subscheme is isomorphic, as a $R[[q]]$-scheme, to the normalization of $R[[q]]$ in the finite \'etale $R((q))$-scheme $\P_{\mrm{Tate}(q)/R((q))}/\{\pm 1\}$, where $\{\pm 1\}$ is the automorphism group of $\mrm{Tate}(q)/R((q))$. Moreover, this isomorphism is natural with respect to morphisms of moduli problems. }

\demo{This is \cite[Theorem 8.11.10]{KM}. }

\medskip 

The following notation is identical to the one introduced in \cite[(10.2)]{KM}. 

\nott{Let $N \geq 3$ and $\mrm{Surj}_N$ denote the set of surjective linear $f: (\Z/N\Z)^{\oplus 2} \rar \Z/N\Z$. Given $\Lambda \in \mrm{Surj}_N$, there is a unique vector $k_{\Lambda}$ such that $\Lambda(x) = \det(k_{\Lambda},x)$ for every $x \in (\Z/N\Z)^{\oplus 2}$. Then $k_{\Lambda}$ is a generator of $\ker{\Lambda}$. We also \emph{choose} a vector $\ell_{\Lambda}$ such that $\Lambda(\ell_{\Lambda})=1$. 

There is a left action of $\GL{\Z/N\Z}$ on $\mrm{Surj}_N$ given by $(g \cdot \Lambda)(x) = \Lambda(g^Tx)$ for $x \in (\Z/N\Z)^{\oplus 2}$ (seen as a column vector).
In the case of $\Lambda: (x,y) \mapsto x$, one sees that $k_{\Lambda}=(0,-1)$, and we choose $\ell_{\Lambda}=(1,0)$. 
}

\prop[km-10-2-5]{(See \cite[Corollary 10.2.5, (10.3.5)]{KM}.) Let $N \geq 3$ be an integer, $S$ be a $\Z[1/N]$-scheme and $G$ a finite \'etale commutative group scheme over $S$, with $N \cdot G = 0$ and a short exact sequence of abelian fppf sheaves $0 \rar \mu_N \rar G \overset{\pi}{\rar} \Z/N\Z \rar 0$. Let $G_1$ denote the closed subscheme $\pi^{-1}(1)$ of $G$.
There is an isomorphism $[\Gamma(N)]_{G/S} \rar \coprod_{\Lambda \in \mrm{Surj}_N}{\mu_N^{\times} \times_{\Z} G_1}$, given for any connected $S$-scheme $T$ by 
\[\left[\phi: (\Z/N\Z)^{\oplus 2}_T \overset{\sim}{\rar} G_T\right] \longmapsto \left(\phi(k_{(\pi\circ\phi)(T)}),\phi(\ell_{(\pi\circ\phi)(T)})\right)_{(\pi\circ\phi)(T)},\]
where $\left[(\pi\circ\phi)(T): (\Z/N\Z)^{\oplus 2}  \rar \Z/N\Z\right] \in \mrm{Surj}_N$. 

Moreover, if $T$ is a connected $S$-scheme, $\Lambda \in \mrm{Surj}_N$, $\zeta \in \mu_N(T), X \in G_1(T)$, and $(\zeta,X)_{\Lambda}$ is the image of $\iota \in [\Gamma(N)]_{G/S}(T)$, then, for $g \in \GL{\Z/N\Z}$, the image of $g \cdot \iota$ is given by $(\zeta',X')_{\Lambda'}$, where 
\begin{itemize}[noitemsep,label=\tiny$\bullet$]
\item $\Lambda' = \Lambda(g^T \cdot)$, $k_{\Lambda'} = (\det{g})(g^T)^{-1}k_{\Lambda}$, 
\item $\ell_{\Lambda'} = (g^T)^{-1}\ell_{\Lambda}+n(g,\Lambda)(g^T)^{-1}k_{\Lambda}$ for some $n(g,\Lambda) \in \Z/N\Z$,
\item $\zeta' = \zeta^{\det{g}}$ and $X'=X \cdot \zeta^{n(g,\Lambda)}$. 
\end{itemize}
}

\demo{This is the same as \cite[(10.2), (10.3)]{KM}, except that the compatibility with the action of $\GL{\Z/N\Z}$ is a slightly different calculation.
}

\cor[cocycle-ell-lambda]{Let $N \geq 3$ be an integer, $g, g' \in \GL{\Z/N\Z}$ and $\Lambda \in \mrm{Surj}_N$. Then one has 
\[\begin{pmatrix}\det{g'g} & 0\\n(g'g,\Lambda) & 1\end{pmatrix} = \begin{pmatrix}\det{g'} & 0\\n(g',g\cdot\Lambda) & 1\end{pmatrix}\begin{pmatrix} \det{g} & 0\\n(g,\Lambda) & 1\end{pmatrix}.\]}

\demo{This is a direct calculation.}

\medskip

\prop[formal-cuspidal-subscheme-gammaN]{Let $N \geq 3$ be an integer and fix a system $(\Lambda_i)_{i \in I}$ of representatives for $\mrm{Surj}_N/\{\pm 1\}$. Then the formal completion of $\mrm{We}: X(N) \rar \Sp{\OO_N}$ along its cuspidal subscheme is isomorphic to 
\[\coprod_{i \in I}{\mrm{Spf}(\OO_N[[q^{1/N}]])}.\]
 Moreover, let $(i, g) \in I \times \GL{\Z/N\Z}$, there exists $(j,g') \in I \times \{\pm g\}$ such that $\Lambda_i \circ (g')^T = \Lambda_j$. Then, the action of $g$ is the same as the action of $g'$; it induces a morphism from the $\Lambda_i$ component to the $\Lambda_j$ component, coming from the ring homomorphism $(\zeta_N,q^{1/N}) \longmapsto (\zeta_N^{\det{g}},\zeta_N^{n(g',\Lambda_i)}q^{1/N})$. }

\demo{This is a consequence of Proposition \ref{km-10-2-5}, which we applied to the description of the $N$-torsion group of $\mrm{Tate}(q)$ over $\Z((q))[N^{-1}]$ \cite[(8.8)]{KM} (because $-I_2$ acts without any fixed point on $\mrm{Surj}_N$ since $N > 2$). The conclusion then follows from Proposition \ref{tate-curve-and-cusps}.
Note that the action of $\GL{\Z/N\Z}$ must preserve $q$ since $q$ can be written as a formal power series in $1/j$. }

\cor[cuspidal-subscheme-gammaN]{Let $N \geq 3$ be an integer. Then the cuspidal subscheme of $X(N)$ is a finite disjoint reunion of copies of $(\mu_N^{\times})_{\Z[1/N]}$ indexed by $\mrm{Surj}_N/\{\pm 1\}$.}

\medskip

\prop[manin-drinfeld]{Let $N \geq 3$ be an integer and $\ell \equiv 1 \mod{N}$ be a prime. Let $J(N)$ be the relative Jacobian of $X(N) \rar \Sp{\Z[1/N]}$ and $\Delta: X(N) \times_{\Sp{\OO_N}} X(N) \rar J(N)$ be the difference map (see Propositions \ref{XN-Weil} and \ref{difference-map-relative-curve}). Let $\iota: C \rar X(N)$ be the cuspidal subscheme. Then $(T_{\ell}-\ell-1)\circ \Delta \circ (\iota \times \iota) = 0$. }

\demo{First, denote by $f: C' \rar J(N)$ the morphism $(T_{\ell}-\ell-1) \circ \Delta \circ (\iota \times \iota)$, where $C' := C \times_{\Sp{\Z[1/N,\zeta_N]}} C$ is a finite \'etale $\Z[1/N]$-scheme. Let $D$ be the scheme-theoretic image of $f$. Then $D$ is a closed subscheme of $J(N)$, so it is of finite type over $\Z[1/N]$. Since $C'$ is finite over $\Z[1/N]$, it is affine, and the morphism $C' \rar D$ is finite surjective: by \cite[Lemma 01YQ]{Stacks}, $D$ is affine, and thus $D$ is finite over $\Z[1/N]$. Moreover, since $C' \rar \Z[1/N]$ is finite \'etale, $\OO(D)$ is torsion-free, so $D \rar \Z[1/N]$ is flat and unramified, hence $D$ is finite \'etale over $\Z[1/N]$. To show that $D \rar J(N)$ is the unit section, we only need to check that the fibre of $D$ at $\ell$ is the unit section of $J(N)_{\F_{\ell}}$. 

Since $C$ is a disjoint reunion of copies of $\OO_N$ and $\ell \equiv 1\mod{N}$, $C'_{\F_{\ell}}$ is a finite reunion of copies of $\Sp{\F_{\ell}}$. In other words, we want to prove that for any two cuspidal points $x,y \in X(N)(\F_{\ell})$ with the same image in $\Sp{\F_{\ell} \otimes \Z[\zeta_N]}$, one has $T_{\ell}(x-y)=(\ell+1)(x-y)$. The conclusion follows from Proposition \ref{eichler-shimura}, since $[\ell]$ acts trivially on $J(N)$, $x-y \in J(N)(\F_{\ell})$, and the Frobenius isogeny on $J(N)_{\F_{\ell}}$ has degree $\ell$. }

\bigskip

\prop[q-expansion-setup]{Let $N \geq 3$ be an integer. Let $\Lambda \in \mrm{Surj}_N$ and $M = \begin{pmatrix} \alpha & 0\\u & 1\end{pmatrix} \in \GL{\Z/N\Z}$. Let $U \subset X(N)$ be the affine inverse image of the open subset $\{j \neq 0, 1728\} \subset \mathbb{P}^1_{\Z[1/N]}$. Then the datum 
\[\left(\mrm{Tate}(q),(k_{\Lambda},\ell_{\Lambda}) \longmapsto (\zeta_N^{\alpha},q^{1/N}\zeta_N^u)\right)\] defines a morphism of $\OO_N \otimes \Z[[q]]$-schemes $C_{(M,\Lambda)}: \Sp{\OO_N \otimes \Z[[q^{1/N}]]} \rar U \times_{\Z[1/j]} \Sp{\Z[[q^{1/N}]]}$ such that the composition 
\[\Sp{\OO_N} \rar \Sp{\OO_N\otimes \Z[[q^{1/N}]]} \overset{C_{(M,\Lambda)}}{\longrightarrow} U \times_{\Z[1/j]} \Sp{\Z[[q^{1/N}]]} \overset{\mrm{We}}{\rar} \Sp{\OO_N}\] is given by $\zeta \longmapsto \zeta^{\alpha}$. 

Given a regular excellent Noetherian $\Z[1/N]$-algebra $R$ and $\omega \in H^0(X(N)_R,\Omega^1_{X(N)_R/R})$, there is a unique $f \in \OO_N \otimes R[[q^{1/N}]]$ such that, for every $t \geq 1$, the pull-back of $\omega$ under the morphism $\Sp{\OO_N \otimes R \otimes \Z[q^{1/N}]/(q^t)} \rar \Sp{\OO_N \otimes R \otimes \Z[[q^{1/N}]]} \overset{C_{(M,\Lambda)}}{\longrightarrow} U_R$ can be written as $f \cdot d(q^{1/N})$.}

\demo{%
By \cite[(8.7),(8.8)]{KM}, the group of sections of the $N$-torsion subscheme of $\mrm{Tate}(q)$ over $\OO_N \otimes \Z((q))$ identifies with the group $\zeta_N^{\Z/N\Z} \times q^{\left(\frac{1}{N}\Z\right)/\Z}$. 

Consider the morphism $C_{(M,\Lambda)}^0: \Sp{\OO_N \otimes \Z((q^{1/N}))} \rar Y(N)$ given by the elliptic curve $(\mrm{Tate}(q),\iota)$, where $\iota: (\Z/N\Z)^{\oplus 2} \rar \mrm{Tate}(q)[N](\OO_N \otimes \Z((q^{1/N})))$ is the isomorphism mapping the $(k_{\Lambda},\ell_{\Lambda})$ to $(\zeta_N^{\alpha},q^{1/N}\zeta_N^u)$. The post-composition of this morphism with the $j$-invariant (resp. the Weil pairing) is exactly given by the inclusion $\Z[j] \subset \Z((q))$ (resp. the projection $\Sp{\OO_N \otimes \Z((q))}\rar \Sp{\OO_N}$). Let $U \subset Y(N)$ denote the inverse image of open subscheme of $\mathbb{P}^1_{\Z[1/N]}$ where $1/j$ is defined and $1728/j-1$ is invertible, then $U$ is affine since $\mathbb{P}^1_{\Z[1/N]}$ is separated. By \cite[(8.11.5)]{KM}, one has $j(\mrm{Tate}(q))(j(\mrm{Tate}(q))-1728) \in \Z((q))^{\times}$, so $C_{(M,\Lambda)}^0$ induces a morphism \[C_{(M,\Lambda)}^1: \Sp{\OO_N \otimes \Z((q^{1/N}))} \rar U \times_{\Sp{\Z[1/j]}} \Sp{\Z[[q]]}.\]

In $\Z[[q]]$, $\frac{1728}{j}-1 \in -1+q\Z[[q]]$ is invertible, so $U\times_{\Sp{\Z[1/j]}} \Sp{\Z[[q]]}$ is a finite $\OO_N \otimes \Z[[q]]$-scheme, so the image of $C_{(M,\Lambda)}^1$ factors through the spectrum of the integral closure $C$ of $A := \OO_N \otimes \Z[[q]]$ in $B := \OO_N \otimes \Z((q^{1/N}))$. The normalization of $\Z[[q]]$ in $\Z((q^{1/N}))$ is clearly $\Z[[q^{1/N}]]$; since normalization commutes with smooth base change \cite[Lemma 03GG]{Stacks} and $\OO_N$ is smooth over $\Z$, one has $C = \OO_N\otimes \Z[[q^{1/N}]]$.   

Hence $C_{(M,\Lambda)}^1$ factors through a finite morphism \[C_{(M,\Lambda)}: \Sp{\OO_N\otimes \Z[[q^{1/N}]]} \rar U \times_{\Z[1/j]} \Sp{\Z[[q^{1/N}]]}.\]

Let $R$ be a regular excellent Noetherian $\Z[N^{-1}]$-algebra. For every integer $t \geq 1$, consider the map $C_{(M,\Lambda),t}^R: \Sp{\OO_N \otimes R \otimes \Z[[q^{1/N}]]/(q^t)} \rar U_R \times_{\Z[1/j]} \Sp{\Z[[q^{1/N}]]}$. Since \[\Omega^1_{\left(\OO_N\otimes R\otimes \Z[[q^{1/N}]]/(q^t)\right)/R} \simeq \frac{(\OO_N \otimes R)[[q^{1/N}]]}{(q^t,Ntq^{t-1/N})}\cdot dq^{1/N},\]
 there exists a unique $f_{(M,\Lambda),t}(\omega) \in (\OO_N \otimes R)[[q^{1/N}]]/(q^t,Ntq^{t-1/N})$ such that \[f_{(M,\Lambda),t}(\omega)d(q^{1/N}) = (C_{(M,\Lambda),t}^R)^{\ast}\omega.\]

To conclude, it is enough to note that $\underset{\longleftarrow}{\lim}\,\left[\frac{(\OO_N \otimes R)[[q^{1/N}]]}{(q^t,Ntq^{t-1/N})}\right] \simeq (R \otimes \OO_N)[[q^{1/N}]]$.}

\defi[qexp-definition]{In the setting of Proposition \ref{q-expansion-setup}, $f$ is the $q$-expansion of the differential $\omega$ at the cusp datum $(M,\Lambda)$. Let $C_{(M,\Lambda)}^{\ast}$ denote the morphism \[\omega \in H^0(X(N)_R,\Omega^1_{X(N)_R/R}) \longmapsto f \in \OO_N \otimes R[[q^{1/N}]].\]  
}

\lem[qexp-elem]{Consider the setting of Proposition \ref{q-expansion-setup}. The formation of the $R$-linear morphism $C_{(M,\Lambda)}^{\ast}: H^0(X(N)_R,\Omega^1_{X(N)_R/R}) \rar \OO_N \otimes R[[q^{1/N}]]$ of $q$-expansion at $(M,\Lambda)$ commutes with base change, and its image is contained in the subring $\OO_N \otimes R \otimes (\Z[1/N])[[q^{1/N}]]$. Moreover, if $g \in \GL{\Z/N\Z}$, then the $q$-expansion of $\omega$ at $\left(\begin{pmatrix}\det{g} & 0\\n(g,\Lambda) & 1\end{pmatrix}M, g \cdot \Lambda\right)$ (the second coordinate is exactly $\Lambda \circ g^T$) is the same as the $q$-expansion at $(M,\Lambda)$ of $g^{\ast}\omega$. }

\demo{The first part of this statement is clear by Proposition \ref{differentials-are-free}. To prove the second part of the claim, write 
\[M=\begin{pmatrix}\alpha & 0\\u & 1\end{pmatrix},\, M'=\begin{pmatrix}\det{g} & 0\\n(g,\Lambda) & 1\end{pmatrix}M=\begin{pmatrix}\alpha\det{g} & 0\\ \alpha n(g,\Lambda)+u & 1\end{pmatrix}.\]
We will show that $C_{(M',g \cdot \Lambda)} = g \circ C_{(M,\Lambda)}$. In the notation of the proof of Proposition \ref{q-expansion-setup}, it is enough to check that $C_{(M',\Lambda\circ g^T)}^0 = g \circ C_{(M,\Lambda)}^0$. %

The morphism $g \circ C_{(M,\Lambda)}^0 \in Y(N)(\OO_N \otimes \Z((q^{1/N})))$ corresponds to the couple $(\mrm{Tate}(q),\lambda \circ g^T)$, where $\lambda: (\Z/N\Z)^{\oplus 2} \rar \mrm{Tate}(q)[N](\OO_N \otimes \Z((q^{1/N})))$ maps $(k_{\Lambda},\ell_{\Lambda})$ to $(\zeta_N^{\alpha},q^{1/N}\zeta_N^u)$. Thus, by Proposition \ref{km-10-2-5}, 
\begin{align*}
(\lambda \circ g^T)(k_{\Lambda'}) &= \lambda(\det{g} k_{\Lambda}) = \zeta_N^{\alpha\det{g}},\\
(\lambda \circ g^T)(\ell_{\Lambda'}) &= \lambda(\ell_{\Lambda})\lambda(n(g,\Lambda)k_{\Lambda})=q^{1/N}\zeta_N^{u+\alpha n(g,\Lambda)},
\end{align*}
hence $(\mrm{Tate}(q),\lambda \circ g^T)$ is exactly the morphism attached to $C_{(M',g \cdot \Lambda)}$. 
}

\rem{By considering a model of $X(N)$ and its cuspidal subscheme over $\Z[\zeta_N]$ (as in, say, \cite[Theorem 10.9.1]{KM}), one can show that the image of this map is contained in the subring $\Z[\zeta_N] \otimes R \otimes \Z[[q^{1/N}]]$, that is, that the denominators of the $q$-expansion are bounded. We do not do this here, because it is not relevant to the rest of the text.}

\prop[qexp-vs-formal-immersion]{In the setting of Proposition \ref{q-expansion-setup}, let $\mathfrak{p}$ be a prime ideal of $\OO_N \otimes R$ and $x \in X(N)_R$ be the image of the prime ideal $(\mathfrak{p},q^{1/N})$ under $C_{(M,\Lambda)}$. Then the local homomorphism $(C_{(M,\Lambda)})^{\sharp}: \OO_{X(N)_R,x} \rar (\OO_N \otimes R)_{\mathfrak{p}}[[q^{1/N}]]$ becomes an isomorphism after completion. 

In particular, for any $t \geq 1$, 
\[\omega_x \in \mathfrak{m}_{X(N)_R,x}^t\Omega^1_{X(N)_R/R,x} \Longleftrightarrow C_{(M,\Lambda)}^{\ast}\omega\in (\mathfrak{p},q^{1/N})^t(\OO_N \otimes R)_{\mathfrak{p}}[[q^{1/N}]].\] }

\demo{We keep the notation of the proof of Proposition \ref{q-expansion-setup}. 

Let $y \in B := \Sp{(\OO_N \otimes R)_{\mathfrak{p}} \otimes (\Z[1/N])[[q^{1/N}]]}$ be the point corresponding to the prime ideal $(\mathfrak{p},q^{1/N})$. The scheme morphism $C^{(M,\Lambda)}$ induces a local homomorphism of $R[1/j]_{(1/j)}$-algebras \[\mu = C_{(M,\Lambda)}^{\sharp}:\, \OO_{U_R,x} \rar \OO_{B,y}= (\OO_N \otimes R)_{\mathfrak{p}}[[q^{1/N}]].\]

By Proposition \ref{smooth-at-infinity-base-change}, there exists $j_N \in \OO_{U_R,x}$ such that the cuspidal subscheme of $U_R$ is exactly, in a neighborhood of $x$, cut out by the section $j_N$. Thus $(j_N)=\sqrt{(1/j)}$ as ideals of $\OO_{U_R,x}$. Because $U_R$ is a smooth $R$-scheme, $\OO_{U_R,x}$ is a regular local ring, hence a unique factorization domain by \cite[Lemma 0AG0]{Stacks}. The local ring $\OO_{U_R,x}/(j_N)$ is the ring of stalks of the cuspidal subscheme at $x$, so it is a regular local ring. Hence $j_N \in \OO_{U_R,x}$ is prime. 

The closed subspace cut out by the kernel $K$ of the morphism $\OO_{U_R,x} \rar \OO_{B,y}/(q^{1/N})$ is the cuspidal subspace. Since $\OO_{B,y}/(q^{1/N})$ is reduced, $K=j_N\OO_{U_R,x}$, so $\mu_q: \OO_{U_R,x}/(j_N) \rar \OO_{B,y}/(q^{1/N})$ is a homomorphism of $R$-algebras such that, by Proposition \ref{q-expansion-setup}, the following diagram commutes, and its vertical maps are localizations at prime ideals with the same restriction to $R$ by Corollary \ref{cuspidal-subscheme-gammaN}:

\[\begin{tikzcd}[ampersand replacement=\&]
\OO_N \otimes R \arrow{r}{\underline{\det{M}} \otimes \mrm{id}}\arrow{d} \& \OO_N \otimes R\arrow{d}\\
\OO_{U_R,x}/(j_N) \arrow{r}{\mu_q} \& \OO_{B,y}/(q^{1/N}).
\end{tikzcd}\]

Hence $\mu_q$ is an isomorphism. In particular, it implies that $\OO_{U_R,x} \rar \OO_{B,y}/(\mathfrak{p},q^{1/N})$ is onto. Since $\OO_{U_R,x}/(j_N)$ is the ring of stalks at $x$ of the finite \'etale (over $R$) cuspidal subscheme of $X(N)_R$, its maximal ideal is generated by $\mathfrak{p}$: hence the maximal ideal of $\OO_{U_R,x}$ is generated by $(\mathfrak{p},j_N)$, and the maximal ideal of $\OO_{B,y}$ is generated by $\mu(\mathfrak{m}_{U_R,x})$. Hence $\OO_{U_R,x} \rar \OO_{B,y}/(\mathfrak{p},q^{1/N})^t$ is onto for every $t \geq 1$. 

Thus the completion morphism $\hat{\mu}: \widehat{\OO_{U_R,x}} \rar \widehat{(\OO_N \otimes R)_{\mathfrak{p}}}[[q^{1/N}]]$ is a surjective homomorphism of regular local rings by \cite[Lemma 07NY]{Stacks}. We show that the two regular local rings have the same dimension, which implies that $\hat{\mu}$ is an isomorphism. By \cite[Lemma 07NV]{Stacks}, it is enough to show that $\OO_{U_R,x}$ and $(\OO_N \otimes R)_{\mathfrak{p}}[x]_{(x)}$ have dimension $1+\dim{R_{\mathfrak{p}'}}$, where $\mathfrak{p}'$ is the inverse image in $R$ of $\mathfrak{p}$. In fact, both rings are the rings of stalks at a point $z$ of a smooth $R$-scheme of relative dimension one, $z$ being mapped to the prime ideal $\mathfrak{p}'$, and $z$ being closed in the fibre at $\mathfrak{p}'$. The conclusion then follows from usual dimension theory (for instance \cite[Theorem 4.3.12]{QL} and \cite[Lemma 02G1]{Stacks}). 

Therefore, for every $t \geq 1$, $\mu_t: \OO_{U_R,x}/\mathfrak{m}_{U_R,x}^t \rar \OO_{B,y}/(\mathfrak{p},q^{1/N})^t$ is an isomorphism. 

Since the local sections $1/j$ and $j_N$ cut out the same closed subspace of $U_R$ in a neighborhood of $x$, $1/j \in \OO_{U_R,x}^{\times} j_N^s$ for some $s \geq 1$. By Proposition \ref{formal-cuspidal-subscheme-gammaN}\footnote{The formal variable $q$ in Proposition \ref{formal-cuspidal-subscheme-gammaN} is actually the same as the one in this proof, because $1/j$ can be written as the same formal power series in either variable. However, the $q^{1/N}$ do not have to be the same.}, one has $s=N$. Then $\mu(j_N)^N \in  \OO_{B,y}^{\times} \cdot (1/j) = \OO_{B,y}^{\times}q$, therefore $\mu(j_N) \in (\OO_N \otimes R)^{\times}_{\mathfrak{p}}q^{1/N}+q^{2/N}\OO_{B,y}$.

Now, $\Omega^1_{U_R/R,x}$ is a free $\OO_{U_R/R,x}$-module of rank one with basis $dj_N$ by \cite[\S 2.2 Proposition 7]{BLR}. Write $\omega_x = gdj_N$: then one has $C_{(M,\Lambda)}^{\ast}\omega=\frac{d\mu(j_N)}{dq^{1/N}}\mu(g) \cdot dq^{1/N}$, with $\mu(g),\mu(j_N) \in (\OO_N \otimes R)_{\mathfrak{p}}[[q^{1/N}]]$ and the derivative $\frac{d\mu(j_N)}{dq^{1/N}}$ is computed formally in $(\OO_N \otimes R)_{\mathfrak{p}}[[q^{1/N}]]$.  

We saw above that $\frac{d\mu(j_N)}{dq^{1/N}} \in (\OO_N \otimes R)_{\mathfrak{p}}[[q^{1/N}]]^{\times}$ was invertible. Therefore, for $t \geq 1$, $\omega_x \in \mathfrak{m}_{U_R,x}^t\Omega^1_{U_R/R,x}$ if and only if $g \in \mathfrak{m}_{U_R,x}^t$, if and only if $\mu(g) \in (\mathfrak{p},q^{1/N})^t\OO_{B,y}$ if and only if $\frac{d\mu(j_N)}{dq^{1/N}}\mu(g) \in (\mathfrak{p},q^{1/N})^t\OO_{B,y}$, whence the conclusion.

}

\prop[qexp-coefficients]{In the notation of Lemma \ref{qexp-elem}, the morphism \[C_{(M,\Lambda)}^{\ast}: H^0(X(N)_R,\Omega^1_{X(N)_R/R}) \rar \OO_N \otimes R \otimes (\Z[1/N])[[q^{1/N}]]\] is injective and its cokernel is flat. 

Therefore, if $S$ is a regular excellent Noetherian $R$-algebra, with $R \rar S$ injective, the inverse image of $\OO_N \otimes R \otimes (\Z[1/N])[[q^{1/N}]]$ under \[(C_{(M,\Lambda)}^{\ast}): H^0(X(N)_S,\Omega^1_{X(N)_S/S}) \rar \OO_N \otimes S \otimes (\Z[1/N])[[q^{1/N}]]\] is exactly $H^0(X(N),\Omega^1_{X(N)_R/R})$.}

\demo{By Proposition \ref{differentials-are-free}, and since $\Z[\zeta_N,N^{-1}][[q^{1/N}]]$ is a torsion-free abelian group, $H^0(X(N)_R,\Omega^1_{X(N)_R/R})$ and $ \Z[\zeta_N] \otimes R \otimes (\Z[1/N])[[q^{1/N}]]$ are flat $R$-modules. The same Proposition implies that $(C_{(M,\Lambda)}^S)^{\ast} \cong (C_{(M,\Lambda)}^R)^{\ast} \otimes_R S$ if $S$ is a regular excellent Noetherian $R$-algebra.  

Assume that $(C_{(M,\Lambda)}^{\Z[1/N]})^{\ast}$ is injective and remains injective after tensoring with $\F_p$ for any prime $p \nmid N$. Let $V$ denote the cokernel of $C_{(M,\Lambda)}^{\Z[1/N]})^{\ast}$: then for any prime $p \nmid N$, one has $\mrm{Tor}_1^{\Z[1/N]}(\F_p,V)=0$. Hence $V$ is flat over $\Z[1/N]$, and the sequence \[0 \rar H^0(X(N),\Omega^1_{X(N)/\Z[1/N]}) \rar \Z[\zeta_N] \otimes (\Z[1/N])[[q^{1/N}]] \rar V \rar 0\] remains exact after tensoring with any $\Z[1/N]$-module, so we are done. 

Since $(C_{(M,\Lambda)}^{\Z[1/N]})^{\ast}$ is injective if and only if $(C_{(M,\Lambda)}^{\Q})^{\ast}$ is injective, it is enough to show that $(C_{(M,\Lambda)}^{R})^{\ast}$ is injective for any field $R$ of characteristic not dividing $N$. 

Thus, let $R$ be any field of characteristic not dividing $R$, and let $\omega \in H^0(X(N)_R,\Omega^1_{X(N)_R/R})$ be in the kernel of $C_{(M,\Lambda)}^{\ast}$. Let $x \in X(N)_R$ be the image under $(C_{(M,\Lambda)}\times R$ of any ideal of the form $(\mathfrak{p},q^{1/N})$ for some prime ideal $\mathfrak{p} \subset R \otimes \Z[\zeta_N]$. Then, by Proposition \ref{qexp-vs-formal-immersion}, one has $\omega_x \in \cap_{t \geq 1}{\mathfrak{m}_{X(N)_R,x}^t\Omega^1_{X(N)_R/R,x}}$. Since $\Omega^1_{X(N)_R/R,x}$ is a free $\OO_{X(N)_R,x}$-module of rank one, by Krull's theorem \cite[Corollary 1.3.13]{QL}, $\omega_x = 0$. 

Since $X(N)_R \rar \Sp{R}$ is smooth of relative dimension one, $\Omega^1_{X(N)_R/R}$ is an invertible $\OO_{X(N)_R}$-module, and its global section $\omega$ vanishes on some open subscheme $U \subset X(N)_R$ such that $U$ contains the cuspidal subscheme of $X(N)_R$. Then $U$ is dense in the regular scheme $X(N)_R$, so that, for any affine open subscheme $V$, the restriction $\Omega^1_{X(N)_R/R}(V) \rar \Omega^1_{X(N)_R/R}(V \cap U)$ is injective, hence $\omega=0$.} 

Let us finish this section with a lemma that will help us compute the full $q$-expansion of a differential. 

\medskip

\lem[with-delta-c]{Let $N \geq 3$ be an integer and $c \in (\Z/N\Z)^{\times}$, then \[C_{(I_2,(x,y)\mapsto x)} \circ \underline{c} = \begin{pmatrix}1 & 0\\0 & c\end{pmatrix}C_{(I_2,(x,y)\mapsto x)}.\] }

\demo{The $\Z[1/N]$-scheme $X(N)$ is proper smooth and $\OO_N \otimes \Z[[q^{1/N}]]$ is reduced, so it is enough to show that $C^0_{(I_2,(x,y)\mapsto x)} \circ \underline{c} = \begin{pmatrix}1 & 0\\0 & c\end{pmatrix}C^0_{(I_2,(x,y)\mapsto x)}$. 

Since $k_{(x,y) \mapsto x} = (0,-1)$ and $\ell_{(x,y)\mapsto x} = (1,0)$, $C^0_{(I_2,(x,y)\mapsto x)}$ is the point \[(\mrm{Tate}(q),q^{1/N},\zeta_N^{-1}) \in Y(N)(\OO_N \otimes \Z((q^{1/N}))).\] 
Thus, the morphisms $C^0_{(I_2,(x,y)\mapsto x)} \circ \underline{c},\,\begin{pmatrix}1 & 0\\0 & c\end{pmatrix}C^0_{(I_2,(x,y)\mapsto x)}$ are attached to the following elliptic curves with level structure over $\OO_N \otimes \Z((q^{1/N}))$:
\[(\mrm{Tate}(q),q^{1/N},\zeta_N^{-c}),\, (\mrm{Tate}(q),q^{1/N},\zeta_N^{-c}),\] which are equal, whence the conclusion. 

}

\lem[qexp-plus-gl2]{Let $N \geq 3$ be an integer, $(M_1,\Lambda_1)$ be a cusp datum and $g \in \SL{\Z/N\Z}$. Then there exists a cusp datum $(M_2,\Lambda_2)=(M_2,\Lambda_1 \circ g^T)$ such that $g \circ C_{(M_1,\Lambda_1)} = C_{(M_2,\Lambda_2)}$. }

\demo{Since $X(N)$ is a proper $\Z[1/N]$-scheme and $\OO_N \otimes \Z[[q^{1/N}]]$ is reduced, it is enough to show that for some cusp datum $(M_2,\Lambda_2)$, one has $g \circ C^0_{(M_1,\Lambda_1)} = C^0_{(M_2,\Lambda_2)}$.

Let $M_1 = \begin{pmatrix} \alpha_1 & 0\\u_1 & 1\end{pmatrix}$, then $C^0_{(M_1,\Lambda)}$ is the point
\[(\mrm{Tate}(q),(k_{\Lambda_1},\ell_{\Lambda_1}) \mapsto (\zeta_N^{\alpha_1},q^{1/N}\zeta_N^{u_1})) \in Y(N)(\OO_N \otimes \Z((q^{1/N}))).\]
Thus $g \circ C^0_{(M_1,\Lambda_1)}$ is  the point 
\[
(\mrm{Tate}(q),((g^T)^{-1}k_{\Lambda_1},(g^T)^{-1}\ell_{\Lambda_1}) \mapsto (\zeta_N^{\alpha_1},q^{1/N}\zeta_N^{u_1}))  \in Y(N)(\OO_N \otimes \Z((q^{1/N}))),\] which is equal by Proposition \ref{km-10-2-5} to 
\[ (\mrm{Tate}(q), (k_{\Lambda_2},\ell_{\Lambda_2}-n(g,\Lambda_1)k_{\Lambda_2}) \mapsto (\zeta_N^{\alpha_1},q^{1/N}\zeta_N^{u_1})) \in Y(N)(\OO_N \otimes \Z((q^{1/N}))),\] which is equal to 
\[ (\mrm{Tate}(q), (k_{\Lambda_2},\ell_{\Lambda_2}) \mapsto (\zeta_N^{\alpha_1},q^{1/N}\zeta_N^{u_1+n(g,\Lambda_1)\alpha_1})) \in Y(N)(\OO_N \otimes \Z((q^{1/N}))),\] which is exactly $C^0_{(M_2,\Lambda_2)}$ with $M_2 = \begin{pmatrix} \alpha_1 & 0\\u_1+n(g,\Lambda_1)\alpha_1 & 1\end{pmatrix}$.

}

\newpage

%% file: analytic-0.tex
\chapter{The geometric situation}
\label{analytic}

\section{Statement of results}

We introduce the objects that we will be working with in this chapter, and state the results that we aim to prove. We refer to the next subsection for a more careful discussion of why (or how) the various constructions sketched here are well-defined, and mostly conform to standard notation. 

Let $p \geq 7$ be a prime number. 

Let $X(p,p)$ be the compactified modular curve over $\Z[1/p]$ parametrizing the elliptic curves endowed with a full level $p$ structure. This curve possesses an action of $\F_p^{\times}$ (denoted by $[-]$) and Hecke correspondences $T_{\ell}$ for $\ell \neq p$; all of them act on its Jacobian $J(p,p)$, we denote by $\mathbb{T}$ the commutative subring of $\mrm{End}(J(p,p))$ generated by the $T_{\ell}$ and the $[n]$ for any $n \in \F_p^{\times}$. For every $n \geq 1$ coprime to $p$, we can define a Hecke operator $T_n \in \mathbb{T}$ by the usual formulas (see Definition \ref{hecke-tn}). Moreover, there is a natural action of $\GL{\F_p}$ on $X(p,p)$, which induces by push-forward functoriality an action of $\GL{\F_p}$ on $J(p,p)$ (such that $nI_2$ acts by $[n]$ for any $n \in \F_p^{\times}$), which commutes to $\mathbb{T}$. We wish to understand the structure of this action of $\mathbb{T}[\GL{\F_p}]$ on $J(p,p)$. This will enable us in later chapters to discuss the Tate module of $J(p,p)$, then of its various twists by Galois representations. 

Let $X_1(p)$ be the compactified modular curve over $\Z[1/p]$ parametrizing elliptic curves endowed with a point of order $p$, and $X_1(p,\mu_p)=X_1(p)_{\Z[1/p,\zeta_p]}$, seen as a smooth projective curve over $\Z[1/p]$. Its non-cuspidal points parametrize triples $(E,P,\alpha)$, where $E$ is a relative elliptic curve over some base $S$, $P \in E(S)$ is a point of order $p$, and $\alpha \in \OO(S)$ is a primitive $p$-th root of unity. The curve $X_1(p,\mu_p)$ is endowed with natural Hecke correspondences $T_{\ell}$ for primes $\ell \neq p$, $U_p$ and diamond operators $\langle n\rangle$ for $n \in \F_p^{\times}$ (all coming from $X_1(p)$), as well as an Atkin-Lehner involution $w_p$. All of them induce endomorphisms of its Jacobian $J_1(p,\mu_p)$: $U_p,T_{\ell}$ and the $\langle n\rangle$ pairwise commute, but they do not commute with $w_p$.  

We define modified Hecke operators $T'_{\ell}$ (for primes $\ell \neq p$) and modified diamond operators $\langle n\rangle'$ (for $n \in (\Z/p\Z)^{\times}$) on $J_1(p,\mu_p)$ (see Definition \ref{modified-hecke-operator}). They generate a commutative subring $\mathbb{T}'$ of of $\mrm{End}(J_1(p,\mu_p))$ that we call the modified Hecke algebra. As above (see Definition \ref{hecke-tnprime}), we can define elements $T'_n \in \mathbb{T}'$ for any $n \geq 1$ coprime to $p$. By construction, $\mathbb{T}'$ commutes with $U_p$ and $w_p$. Let $P(p)$ be the direct product of $p+1$ copies of $J_1(p,\mu_p)$ indexed by $\F_p \cup \{\infty\}$; we can construct actions of $\GL{\F_p}$, $U_p$ and $w_p$ on $P(p)$ which all commute with $\mathbb{T}'$. 

For each $t \in \F_p$ (resp. $t=\infty$), we define a morphism $u_t: J(p,p) \rar J_1(p,\mu_p)$ as follows by push-forward functoriality: we map a couple $(E,(P,Q))$ to $(E,tP+Q,\langle P,\,Q\rangle_E)$ (resp. to $(E,P,\langle P,\,Q\rangle_E)$), where $\langle P,\,Q\rangle_E$ is the Weil pairing on $E[p]$. The morphism $u: J(p,p) \rar P(p)$ is the direct product of the $u_t$. Then $C(p):=\ker{u}$ is a proper group scheme over $\Z[1/p]$, and we let $C^0(p)$ be the connected component of $C(p)$ containing the image of the unit section. 

We define the endomorphism $R$ of $P(p)$ as follows: given $s,t \in \F_p \cup \{\infty\}$, the morphism $J_1(p,\mu_p)_s \subset P(p) \overset{R}{\rar} P(p) \rar J_1(p,\mu_p)_t$ is:

\begin{itemize}[noitemsep,label=\tiny$\bullet$]
\item $-U_p$ is $s=t$,
\item $w_p$ if $s \neq t$ and $\infty \in \{s,t\}$,
\item $\langle s-t\rangle w_p$ if $s,t$ are distinct elements of $\F_p$.
\end{itemize}

Then $R$ commutes with $\mathbb{T}'$ and $\GL{\F_p}$. 

\theoi[analytic-1]{The following sequence of Abelian varieties over $\Q$ 

\[0 \rar C^0(p) \rar J(p,p) \overset{u}{\rar} P(p) \overset{R}{\rar} P(p)\]

is exact up to isogeny. It is equivariant for the Hecke operators (that is, $T_{n},\langle n\rangle$ for $J(p,p)$ and $T'_{n},\langle n\rangle'$ for $P(p)$) and for the actions of $\GL{\F_p}$. }

Thanks to this exact sequence, we can describe the space of holomorphic differentials on $X(p,p)_{\C}$ as a $(\mathbb{T} \otimes \C)[\GL{\F_p}]$-module. To do this, we need some information about the representations of $\GL{\F_p}$, which is recalled in Appendix \ref{reps-gl2}. 

\nott{\begin{itemize}[noitemsep,label=$-$]
\item The Steinberg representation is denoted by $\mrm{St}$.
\item If $\psi: \F_p^{\times} \rar \C^{\times}$ is a caracter, $\mrm{St}_{\psi}$ denotes the twist of $\mrm{St}$ by the character $\psi(\det)$. 
\item Given two Dirichlet characters $\alpha,\beta: \F_p^{\times} \rar \C^{\times}$, the principal series representation attached to the pair $\{\alpha,\beta\}$ is $\pi(\alpha,\beta)$.
\item If $V$ is a left complex representation of $\GL{\F_p}$, $V^{\vee} := \mrm{Hom}(V,\C)$ is a right $\C[\GL{\F_p}]$-module. 
\end{itemize}}

\theoi[analytic-2]{The right $(\mathbb{T} \otimes \C)[\GL{\F_p}]$-module $H^0(X(p,p)_{\C},\Omega^1)$ is isomorphic to the direct sum of the following modules:
\begin{itemize}[noitemsep,label=$-$]
\item For every Dirichlet character $\psi: \F_p^{\times} \rar \C^{\times}$ and every newform $f \in \mathcal{S}_2(\Gamma_1(p))$ with character $\chi$, modulo the relation $(f,\psi) \sim (\overline{f},\chi\psi)$, the representation $\pi(\psi,\chi\psi)^{\vee}$ if $\chi \neq \mathbf{1}$ and $\mrm{St}_{\psi}^{\vee}$ if $\chi=1$. The pull-back by the Hecke operator $T_n$ for an integer $n \geq 1$ coprime to $p$ (resp. by $m_n:=nI_2$ for $n \in \F_p^{\times}$) acts by multiplication by $a_{n}(f)\psi(n)$ (resp. $\chi(n)\psi(n)^2$). 
\item For every newform $f \in \mathcal{S}_2(\Gamma_0(p^2) \cap \Gamma_1(p))$ with character $\chi$, an irreducible cuspidal representation $C_f^{\vee}$ of $\GL{\F_p}$ where pulling back by $T_{n}$ (resp. by $m_n$) is multiplication by $a_{n}(f)$ (resp. $\chi(n)$).    
\end{itemize} }

\theoi[analytic-3]{For any field $k$ of characteristic distinct from $p$, the sequence \[J(p,p)_k \overset{u_k}{\rar} P(p)_k \overset{R_k}{\rar} P(p)_k\] is exact up to isogeny. Moreover, the groups $2p(p+1)\pi_0(\ker{u_k})$ and $2p(p+1)\pi_0(\ker{R_k})$ vanish. }

\section{Constructions}

Let $X(p,p)$ (resp. $Y(p,p)$) be the compactified moduli scheme (resp. moduli scheme) attached to the finite \'etale representable moduli problem $[\Gamma(p)]_{\Z[1/p]}$; it is projective and smooth of relative dimension one over $\Z[1/p]$. We let $J(p,p)$ denote its Jacobian as a relative curve\footnote{In usual sources such as \cite{BLR} or \cite{MilJac}, the Jacobian of a smooth projective curve over some base is usually defined only when the geometric fibres are connected. We provide an extended definition in Appendix \ref{degrees-jacobians} and check what properties still hold in this more general setting. Many of the constructions below are defined and studied in detail and greater generality in Chapter \ref{moduli-spaces}.} over $\Z[1/p]$. 

By Corollary \ref{basic-has-good-hecke}, there are Hecke operators $T_{\ell}$ acting on $J(p,p)$ for all primes $\ell \neq p$. Since $\GL{\F_p}$ acts on $[\Gamma(p)]_{\Z[1/p]}$ (by a morphism of moduli problems of level $p$), it acts on $X(p,p)$ and thus on $J(p,p)$ by push-forward functoriality. Moreover, for any $n \in \F_p^{\times}$, the action of $nI_2$ on $X(p,p)$ (hence on $J(p,p)$ by push-forward functoriality) is exactly the natural action of $n \in \F_p^{\times}$ as an automorphism of moduli problems of level $p$ (see the beginning of Section \ref{gamma0m-structure}). By Corollary \ref{hecke-weil}, the action of $\GL{\F_p}$ commutes to any $T_{\ell}$ for any prime $\ell \neq p$. Let $\mathbb{T}$ be the commutative sub-algebra of $\mrm{End}(J(p,p))$ generated by the $T_{\ell}$ (for $\ell \neq p$ prime) and the $nI_2$ for $n \in \F_p^{\times}$. 

\defi[hecke-tn]{We define the elements $T_n \in \mathbb{T}$ (for $n \geq 1$ coprime to $p$) in the following way:
\begin{itemize}[noitemsep,label=$-$]
\item $T_1$ is the identity,
\item for any prime $\ell \neq p$, for any $k \geq 2$, $T_{\ell^k}=T_{\ell}T_{\ell^{k-1}}-\ell\cdot (\ell I_2)T_{\ell^{k-2}}$.
\item for coprime integers $m,n \geq 1$ such that $p \nmid mn$, $T_{mn}=T_mT_n$.
\end{itemize}}

\lem[hecke-tn-group]{For any integers $n,m \geq 1$ coprime to $p$ and with greatest common divisor $\delta$, one has
\[T_nT_m=\sum_{d \mid \delta}{d \cdot (dI_2)T_{\frac{mn}{d^2}}}.\]
In particular, for all integers $n,m_1,\ldots,m_r \geq 1$ coprime to $p$, $(nI_2)T_{m_1}\ldots T_{m_r}$ is a $\Z$-linear combination of certain $(q'I_2)\cdot T_q$, where $q',q \geq 1$ are integers such that $q(q')^2 \equiv n^2m_1\ldots m_r \pmod{p}$. Thus $\mathbb{T}$ is generated as a group by the $T_n$. 
}

\demo{The first part of the statement implies the second part by an easy induction. By multiplicativity, we can assume that $n,m$ are powers of the same prime number $q \neq p$. We can assume that $n \leq m$, and we proceed by induction over $n$. When $n \leq q$, the identity is a reformulation of $T_qT_m=T_{qm}-q(qI_2)T_{m/q}$, which is the definition of $T_{mq}$. 
When $n > q$, one has
\begin{align*}
T_nT_m &= T_qT_{n/q}T_m - q(qI_2)T_{n/q^2}T_m = \sum_{d \mid n/q}{d \cdot (dI_2) \cdot T_qT_{\frac{mn}{qd^2}}} - \sum_{d \mid n/q^2}{dq \cdot (qdI_2) \cdot T_{\frac{mn}{q^2d^2}}} \\
&= \sum_{d \mid n/q}{d \cdot (dI_2) \cdot (T_{\frac{nm}{d^2}}+q(qI_2)T_{\frac{mn}{q^2d^2}})} - \sum_{q\leq d \mid n/q}{d \cdot (dI_2) \cdot T_{\frac{mn}{d^2}}}\\
&= n \cdot (nI_2)T_{\frac{m}{n}}+\sum_{d \mid n/q}{d \cdot (dI_2)T_{\frac{mn}{d^2}}},
\end{align*}
whence the conclusion.
}

\medskip\noindent

Write $\mrm{We}: X(p,p) \rar \Sp{\Z[1/p,\zeta_p]}$ for the Weil pairing as defined in Proposition \ref{XN-Weil}. 

We treat elements of $\F_p^{\oplus 2}$ as column vectors, and choose the following system $S$ of representatives for $\mathbb{P}^1(\F_p)$: $S$ is made with the $(t,1)$ for $t \in \F_p$ and with $(1,0)$. 

Let $X_1(p,\mu_p)$ (resp. $Y_1(p,\mu_p)$) be the compactified moduli scheme (resp. the moduli scheme) attached to the finite \'etale representable moduli problem $[\Gamma_1(p)]_{\Z[1/p]} \times [\Z[1/p,\mu_p]]$ (of level $p$): it is projective and smooth of relative dimension one over $\Z[1/p]$, so it has a relative Jacobian $J_1(p,\mu_p)$. By Lemma \ref{with-extra-functions}, $X_1(p,\mu_p)$ is exactly the fibre product $X_1(p) \times_{\Z[1/p]} \Sp{\Z[1/p,\mu_p]}$. By Corollary \ref{basic-has-good-hecke}, there are Hecke operators $T_{\ell}$ for $\ell \neq p$ acting on $J_1(p,\mu_p)$. Instead of the notation of Section \ref{gamma0m-structure}, we will denote the action of $(\Z/p\Z)^{\times}$ (given by the level $p$ structure) by $\langle \cdot\rangle$ rather than $[\cdot]$.

Given $n \in \F_p^{\times}$, let $\underline{n}$ be the automorphism of $\Sp{\Z[1/p,\mu_p]}$ with cyclotomic character $n$. In particular, by base change, $n \longmapsto \underline{n}$ induces another action of $\F_p^{\times}$ on the $\Z[1/p]$-scheme $X_1(p,\mu_p)$: on non-cuspidal points, it maps $(E/S,P,\alpha)$ to $(E/S,P,\alpha^n)$. 

\lem[bracket-is-diamond]{For any $n \in \F_p^{\times}$, $\langle n\rangle$ is the base change of the classical diamond operator on $X_1(p)$, mapping a non-cuspidal point $(E/S,P)$ to $(E/S,nP)$. Moreover, the actions of the Hecke operators $T_{\ell}$ (for $\ell \neq p$ prime), $\langle n \rangle$ and $\underline{n}$ (for $n \in \F_p^{\times}$) commute.}

\demo{The commutation with the Hecke operators follow from the fact that both $\langle \cdot \rangle$ and $\underline{\cdot}$ act as morphisms of moduli problems of level $p$. For the description of $\langle \cdot \rangle$ as well as its commutation with $\underline{\cdot}$, both claims are formal verifications on non-cuspidal points. The conclusion follows from the fact that the open immersion $Y_1(p,\mu_p) \rar X_1(p,\mu_p)$ is scheme-theoretically dense. }

By the results of Section \ref{bad-hecke}, we have a ``bad'' Hecke operator $U_p$ acting on $J_1(p,\mu_p)$ commuting to every $T_{\ell}$, $\langle \cdot \rangle$, $\underline{\cdot}$, there is an Atkin-Lehner automorphism $w_p$ (called $w'_p$ in Section \ref{bad-hecke}) of $X_1(p,\mu_p)$ acting on $J_1(p,\mu_p)$ by push-forward functoriality. 

The operator $w_p$ acts in the following way: it maps a non-cuspidal point $(E/S,P,\alpha)$ to $(E'/S,Q,\alpha)$, where $E'=E/\langle P\rangle$, $\pi: E \rar E'$ is the natural isogeny, and $Q \in \ker{\pi^{\vee}}(S) \subset E'[p](S)$ is the unique point such that $\langle P,\,Q\rangle_{\ker{\pi}} = \alpha$ (in the sense of \cite[(2.8)]{KM}). In particular, by \cite[(2.8.6.1)]{KM}, for any $Q_1 \in E[p](S')$ (with $S' \rar S$ \'etale) such that $\langle P,\, Q_1\rangle_{E[N]} = \alpha$, then $Q=\pi(Q_1)$.

\lem[Up-commutes]{$U_p$ commutes with the $T_{\ell}$ (for $\ell \neq p$ prime), as well as the $\langle n\rangle$ and the $\underline{n}$ for $n \in \F_p^{\times}$.}

\demo{For the Hecke operators and the diamonds, it follows from Proposition \ref{bad-hecke-natural}. For $\underline{\cdot}$, it follows from Proposition \ref{relative-picard-mixed-functoriality} and the fact that the following diagram is Cartesian (adopting the notation of Section \ref{bad-hecke}), for $n \in \F_p^{\times}$:
 \[
 \begin{tikzcd}[ampersand replacement=\&]
 \CMPp{[\Gamma_1(p)],[\mu_p]} \arrow{d}{\underline{n}}\& \CMPp{[\Gamma_1'(p,p)],[\mu_p]} \arrow{l}{D'_{1,1}} \arrow{r}{D'_{p,1}} \arrow{d}{\underline{n}}\& \CMPp{[\Gamma_1(p)],[\mu_p]} \arrow{d}{\underline{n}}\\
   \CMPp{[\Gamma_1(p)],[\mu_p]} \& \CMPp{[\Gamma_1'(p,p)],[\mu_p]} \arrow{l}{D'_{1,1}} \arrow{r}{D'_{p,1}} \& \CMPp{[\Gamma_1(p)],[\mu_p]}
 \end{tikzcd}
 \]}

\lem[wp-relations]{For any $n \in \F_p^{\times}$, one has $\langle n\rangle w_p\langle n\rangle = w_p$ and $w_p \circ \underline{n} = \underline{n}\circ \langle n\rangle \circ w_p$. Moreover, $w_p$ is an involution and, for any prime $\ell \neq p$, $\langle \ell\rangle w_p T_{\ell}=T_{\ell}w_p$.
}

\demo{The behavior with respect to $\langle \cdot \rangle$ has been proved in Section \ref{gamma1m-structure} for moduli schemes and this extends to the compactified moduli scheme by Proposition \ref{bad-degn-AL-infty}. Similarly, $w_p \circ w_p$ is the action of $-1 \in \hat{\Z}^{\times}$; since $-1$ is an automorphism of any elliptic curve, it acts trivially on moduli schemes, so that $w_p \circ w_p$ is the identity. The identity about Hecke operators is a consequence of Proposition \ref{hecke-bad-AL}. 

For the last equality, it is enough to prove that it holds for non-cuspidal points. If $E/S$ is an elliptic curve and $\alpha \in \mu_p(S)$, $P \in E(S)$ is a point of exact order $p$, let $\pi: E \rar E'$ be the isogeny whose kernel is generated by $P$. Then $w_p(\underline{n}(E/S,P,\alpha))=(E',Q,\alpha^n)$, where $Q \in \ker{\pi^{\vee}}(S)$ is the unique point such that $\langle -,\,Q\rangle_{\ker{\pi}}$ maps $P$ to $\alpha^n$. Therefore, $w_p(E/S,P,\alpha)=(E'/S,n^{-1}Q,\alpha)$. Hence $\underline{n}\langle n\rangle w_p(E/S,P,\alpha)=w_p(\underline{n}(E/S,P,\alpha))$, whence the conclusion. }

\medskip\noindent

\defi[modified-hecke-operator]{The \emph{modified} Hecke operator $T'_{\ell}$ on $J_1(p,\mu_p)$ is $T'_{\ell}=\underline{\ell}\circ T_{\ell}$, for any prime $\ell \neq p$. For any $n \in \F_p^{\times}$, the operator $\langle n\rangle'$ on $X_1(p,\mu_p)$ is given by $\langle n\rangle \circ \underline{n^2}$. The modified Hecke algebra $\mathbb{T}'$ is the subring of $\mrm{End}(J_1(p,\mu_p))$ generated by the $T'_{\ell}$ (for $\ell \neq p$ prime) and the $\langle n\rangle'$ for $n \in \F_p^{\times}$. }

\lem[modified-hecke-commutes]{Every element of $\mathbb{T}'$ commutes with $w_p$ and $U_p$.}

\demo{This follows from Lemmas \ref{Up-commutes} and \ref{wp-relations}.}

\defi[hecke-tnprime]{For any $n \geq 1$ coprime to $p$, we define an element $T'_n \in \mathbb{T}'$ as follow:
\begin{itemize}[noitemsep,label=$-$]
\item $T'_1$ is the identity,
\item For any prime $\ell \neq p$ and any $k \geq 2$, $T'_{\ell^k}=T'_{\ell}T'_{\ell^{k-1}}-\ell\langle\ell\rangle'T'_{\ell^{k-2}}$,
\item For any integers $m,n \geq 1$ such that $p,m,n$ are pairwise coprime, $T'_{mn}=T'_mT'_n$.
\end{itemize}}

\medskip\noindent

\defi{If $x=(t,1) \in S$ (resp. $x=(1,0) \in S$), the morphism $u_x: X(p,p) \rar X_1(p,\mu_p)$ is the compactification (in the sense of Proposition \ref{compactification-functor}) of the morphism \[(E,(P,Q)) \in Y(p,p)(S) \longmapsto (E,tP+Q,\langle P,\,Q\rangle_E) \in Y_1(p,\mu_p)(S).\] }

\lem{Every $u_x$ is finite locally free of constant rank and induces an isomorphism $\OO(X_1(p,\mu_p)) \rar \OO(X(p,p))$. }

\demo{Any $u_x$ is a morphism of finite $\mathbb{P}^1_{\Z[1/p]}$-schemes (the structure map being the $j$-invariant in each case), so it is finite. Since $X(p,p),X_1(p,\mu_p)$ are regular schemes where every nonempty open subset has dimension two, $u_x$ is flat by miracle flatness, so it is locally free. Since $X_1(p,\mu_p)$ is connected by Corollary \ref{bad-pb-with-geom-fibres}, $u_x$ is locally free of constant rank. The same Corollary implies that $\OO(X_1(p,\mu_p))=\Z[1/p,\zeta_p]$, and Proposition \ref{XN-Weil} shows that one also have $\OO(X(p,p)) \simeq \Z[1/p,\mu_p]$, whence the conclusion.}

\medskip\noindent

Let $X_1(p,\mu_p)^S$ be the disjoint reunion of $p+1$ copies of $X_1(p,\mu_p)$ indexed by $S$; by \cite[Lemma 03GO]{Stacks}, $X_1(p,\mu_p)^S$ can also be seen as the compactification of $p+1$ copies of $Y_1(p,\mu_p)$ indexed by $S$. Thus, we can also see $X_1(p,\mu_p)^S$ as the compactified moduli scheme attached to the following moduli problem over $\Ell_{\Z[1/p]}$, mapping a relative elliptic curve $E/Z$ to the set $[\Gamma_1(p)](E/Z) \times \mu_p^{\times}(Z) \times \{f: Z \rar S \mid f \text{ locally constant}\}$. 

We let $\GL{\F_p}$ act on $X_1(p,\mu_p)^S$ as follows: if $M \in \GL{\F_p}$, $x,x' \in S, \lambda \in \F_p^{\times}$ are such that $\lambda (M^T)^{-1}x=x'$, then $M$ acts by 
\[X_1(p,\mu_p)_x \simeq X_1(p,\mu_p) \overset{\underline{\det{M}} \circ \langle \lambda\rangle }{\longrightarrow} X_1(p,\mu_p) \simeq X_1(p,\mu_p)_{x'}.\]

A straightforward calculation (using the fact that $\langle -1\rangle$ is the identity) yields:

\lem[elem-calc-X1pS]{This defines a left action of $\GL{\F_p}$ on $X_1(p,\mu_p)^S$. Moreover,
\begin{itemize}[label=$-$,noitemsep]
\item the matrix $\begin{pmatrix} a & \ast\\0 & d\end{pmatrix} \in \GL{\F_p}$ preserves the component $X_1(p,\mu_p)_{(0,1)}$ and acts on it by $\underline{ad} \circ \langle d\rangle$,
\item the matrix $aI_2 \in \GL{\F_p}$ acts by $\langle a\rangle'$,
\item for $t,u \in \F_p$, the matrix $T^t=\begin{pmatrix}1 & 0\\t & 1\end{pmatrix}$ realizes the obvious isomorphism \[X_1(p,\mu_p)_{(u,1)} \simeq X_1(p,\mu_p) \simeq X_1(p,\mu_p)_{(u-t,1)},\]
\item the matrix $W^{-1} := \begin{pmatrix} 0 & 1\\-1 & 0\end{pmatrix}$ realizes the obvious isomorphism \[X_1(p,\mu_p)_{(0,1)} \simeq X_1(p,\mu_p) \simeq X_1(p,\mu_p)_{(1,0)}.\]
\end{itemize}}

Now, $X_1(p,\mu_p)^S$ is a disjoint reunion of smooth relative curves, so it is a smooth relative curve, and in particular it has a relative Jacobian $P(p)$ which is naturally the product of the $J_1(p,\mu_p)$ indexed by copies of $S$. In particular, the action of $\GL{\F_p}$ (by automorphisms) on $X_1(p,\mu_p)^S$ extends by push-forward functoriality (by Proposition \ref{relative-picard-functoriality}) to an action on $P(p)$. 

We also have a diagonal action of $\mrm{End}(J_1(p,\mu_p))$ on $P(p)$: in particular, $\langle \cdot \rangle$, $\underline{\cdot}$, and $\mathbb{T}'$ generate a commutative subring $\mathbb{T}$ of $P(p)$, which commutes to the action of $\GL{\F_p}$ (since this action only permutes the coordinates and applies certain $\langle n \rangle, \underline{m}$, and these operators commute with $\mathbb{T'}$).  

By Proposition \ref{relative-picard-functoriality}, the sum $u: J(p,p) \rar P(p)$ of the push-forwards $[(u_x)_{\ast}]$ of the $u_x$ (over $x \in S$) is well-defined.

\lem[u-equiv-gl2]{$u$ commutes to the action of $\GL{\F_p}$.}

\demo{Let $M \in \GL{\F_p}$, $\lambda \in \F_p^{\times}$, and $x,x' \in S$ such that $\lambda (M^T)^{-1}x=x'$. It is enough to show that $M \circ u_x=u_{x'} \circ M$ (as morphisms $X(p,p) \rar X_1(p,\mu_p)^S$) by Proposition \ref{relative-picard-functoriality}. It is enough to prove this away from the cuspidal points. Let $S$ be a $\Z[1/p]$-scheme and $(E,P,Q) \in Y(p,p)(S)$, then one has 

\begin{align*}
M\cdot u_x((E,P,Q))&= M \cdot \left(\left(E,x^T\begin{pmatrix}P\\Q\end{pmatrix}\right),\langle P,Q\rangle_E,x\right)\\
&= \left(\left(E,\lambda x^T\begin{pmatrix}P\\Q\end{pmatrix}\right),\langle P,\, Q\rangle_E^{\det{M}},x'\right)\\
&= \left(\left(E,(x')^TM\begin{pmatrix}P\\Q\end{pmatrix}\right),\langle P,\,Q\rangle_E^{\det{M}},x'\right)\\
&= u_{x'}(M\cdot (E,P,Q)).
\end{align*}
 }
 
\textbf{Notation for $2 \times 2$ matrices:} We adopt, in the rest of the text, the following notation for $2 \times 2$ matrices.
\begin{itemize}[noitemsep,label=$-$] 
\item For $a,b \in \F_p^{\times}$, $\Delta_{a,b} \in \GL{\F_p}$ denotes the matrix $\mrm{diag}(a,b)$.  
\item If $n \in \Z$ is coprime to $p$, $\delta_n \in SL_2(\Z)$ is a matrix congruent to $\Delta_{n^{-1},n}$ modulo $p$. Sometimes, we will require additional conditions on the choice of $\delta_n$, which we will make explicit.  
\item If $n \in \F_p^{\times}$, $m_n$ is the matrix $\Delta_{n,n}$. 
\item $U \in \GL{\F_p}$ is the matrix $\begin{pmatrix}1 & 1\\0 & 1\end{pmatrix}$, and $T=\begin{pmatrix}1 & 0\\1 & 1\end{pmatrix}$. 
\item $W \in \GL{\F_p}$ is the matrix $\begin{pmatrix}0 & -1\\1 & 0\end{pmatrix}$. 
\item $B \leq \GL{\F_p}$ is the subgroup made with upper-triangular matrices. 
\end{itemize}

\textbf{Notation for representations of $\GL{\F_p}$:} Throughout the chapter, we will use the well-known classification of complex representations of $\GL{\F_p})$ as described in Appendix \ref{reps-gl2}. The irreducible complex representations of $\GL{\F_p}$ split in four classes:
\begin{itemize}[noitemsep,label=$-$]
\item The abelian ones, which are all of the form $\psi(\det)$, where $\psi: \F_p^{\times} \rar \C^{\times}$ is a character. 
\item The principal series representations $\pi(\alpha,\beta)$ of dimension $p+1$, attached to a pair of distinct characters $\alpha,\beta: \F_p^{\times} \rar \C^{\times}$.
\item The twisted Steinberg representations $\mrm{St}_{\psi} = \mrm{St} \otimes \psi(\det)$ of dimension $p$, where $\psi: \F_p^{\times} \rar \C^{\times}$ is a character. 
\item The cuspidal representations of dimension $p-1$, which are exactly those such that $U$ has no nonzero invariant vector.  
\end{itemize}
For any left $\C[\GL{\F_p}]$-module $M$, we will denote $M^{\vee}=\mrm{Hom}_{\C}(M,\C)$, which is a right $\C[\GL{\F_p}]$-module.  

We denote by $\mathcal{D}$ the group of characters $\F_p^{\times} \rar \C^{\times}$. 

\prop[Pp-over-GL2]{Let $M=\begin{pmatrix}a & \ast\\0 & d\end{pmatrix} \in B$ act on $J_1(p,\mu_p)$ by $\underline{ad}\langle d\rangle$. Let $V$ be any $\Z[1/p]$-scheme. There is a unique isomorphism of $\Z[\GL{\F_p}]$-modules \[\Z[\GL{\F_p}] \otimes_{\Z[B]} J_1(p,\mu_p)(V) \rar P(p)(V)\] such that $[I_2] \times J_1(p,\mu_p)(V) \rar P(p)(V)$ is simply the inclusion $J_1(p,\mu_p)=J_1(p,\mu_p)_{(0,1)} \rar P(p)$ coming from the inclusion $X_1(p,\mu_p)_{(0,1)} \rar X_1(p,\mu_p)^S$.

This isomorphism is functorial with respect to $V$ and commutes with $\mathbb{T}$.}

\demo{The morphism $v_0: J_1(p,\mu_p)(V) = J_1(p,\mu_p)_{(0,1)}(V) \rar P(p)(V)$ is an morphism of $\Z[B]$-modules by Lemma \ref{elem-calc-X1pS}, which is functorial in $V$. Thus it induces a morphism of $\Z[\GL{\F_p}]$-modules $v: \Z[\GL{\F_p}] \otimes_{\Z[B]} J_1(p,\mu_p)(V) \rar P(p)(V)$ which is functorial in $V$. Now, \[\Z[\GL{\F_p}] \otimes_{\Z[B]} J_1(p,\mu_p)(V) \simeq \bigoplus_{M \in T^{\F_p} \cup \{W^{-1}\}}{[M] \otimes J_1(p,\mu_p)(V)},\] and, for any $f \in J_1(p,\mu_p)(V)$, $v([M] \otimes f)=M \cdot v_0(f)$.

Let, for each $x \in S$, $\iota_x$ denote the morphism $X_1(p,\mu_p) = X_1(p,\mu_p)_x \subset X_1(p,\mu_p)^S$. Then $v_0$ is exactly $[(\iota_{(0,1)})_{\ast}]$. 

If $M=T^t$ (resp. $M=W^{-1}$), then $M \circ [(\iota_{(0,1)})_{\ast}]=[(\iota_{(-t,1)})_{\ast}]$ (resp. it equals $[(\iota_{(1,0)})_{\ast}]$) by Lemma \ref{elem-calc-X1pS}. The conclusion follows from the fact that $\prod_{x \in S}{[(\iota_x)_{\ast}]}J_1(p,\mu_p)^{\times S} \rar P(p)$ is an isomorphism. } 
 
 \medskip\noindent

\prop[u-with-hecke]{For any prime $\ell \neq p$ and any $x \in S$, $[(u_x)_{\ast}]\circ T_{\ell}=T'_{\ell}\circ [(u_x)_{\ast}]$. In particular, $u \circ T_{\ell}=T'_{\ell} \circ u$. }

\rem{This does \emph{not} follow directly from Proposition \ref{hecke-natural}, because, while $u_x$ is the compactification of the morphism of moduli problems \[U_x: \begin{pmatrix}P\\Q\end{pmatrix} \in [\Gamma(p)]_{\Z[1/p]}(E/Z) \longmapsto \left(x^T\begin{pmatrix}P\\Q\end{pmatrix},\,\langle P,\,Q\rangle_E\right) \in  [\Gamma_1(p)]_{\Z[1/p]} \times [\mu_p^{\times}](E/Z),\] this morphism of moduli problems does not respect the level $p$ structure. }

\demo{Let $X_{\ell}(p,p)$ (resp. $X_{1,\ell}(p,\mu_p)$) be the compactified moduli scheme attached to the finite flat moduli problem $[\Gamma(p)]_{\Z[1/p]} \times [\Gamma_0(\ell)]$ (resp. $[\Gamma_1(p)]_{\Z[1/p]} \times [\mu_p^{\times}] \times [\Gamma_0(\ell)]$) over $\Ell_{\Z[1/p]}$. Both moduli problems are smooth at infinity by the results of Section \ref{hecke-operators} and both compactified moduli schemes are relative elliptic curves with global ring of functions $\Z[1/p,\mu_p]$. Consider the following diagram:

\[
\begin{tikzcd}[ampersand replacement=\&]
X(p,p) \arrow{r}{u_x = U_x} \& X_1(p,\mu_p)\\
X_{\ell}(p,p) \arrow{r}{U_x}\arrow{u}{D_{1,1}}\arrow{d}{D_{\ell,1}} \& X_{1,\ell}(p,\mu_p)\arrow{u}{D_{1,1}} \arrow{d}{\underline{\ell} \circ D_{\ell,1}}\\
X(p,p) \arrow{r}{u_x = U_x} \& X_1(p,\mu_p)
\end{tikzcd}
\]

Note that all morphisms are morphisms of finite flat $\mathbb{P}^1_{\Z[1/p]}$-schemes (because all moduli problems are finite flat, so this holds away from the cusps by Corollary \ref{over-Ell-concrete}; the finiteness above the cusps is by construction, and the flatness at the cusps follows from miracle flatness). The top square is Cartesian away from the cusps, essentially by Proposition \ref{product-rel-rep}, and in particular it commutes (since the inclusion $Y_1(p,\mu_p) \rar X_1(p,\mu_p)$ is scheme-theoretically dense). 

If we can prove that the bottom square commutes and is Cartesian away from the cusps, then the result follows from Proposition \ref{relative-picard-mixed-functoriality}. 

Consider, then, the following diagram:
\[
\begin{tikzcd}[ampersand replacement=\&]
X_{\ell}(p,p) \arrow{rrr}{U_x}\arrow{dd}{D_{\ell,1}}\arrow{dr}{w_{\ell}} \& \& \& X_{1,\ell}(p,\mu_p) \arrow{dd}{\underline{\ell} \circ D_{\ell,1}}\arrow{dl}{\underline{\ell} \circ w_{\ell}}\\
\& X_{\ell}(p,p) \arrow{r}{U_x}\arrow{dl}{D_{1,1}} \& X_{1,\ell}(p,\mu_p) \arrow{dr}{D_{1,1}} \& \\
X(p,p) \arrow{rrr}{u_x = U_x}\& \& \& X_1(p,\mu_p)
\end{tikzcd}
\]

All the schemes are relative curves over $\Z[1/p]$ with ring of global functions $\OO(\Z[1/p,\mu_p])$, and all the maps (except potentially the $U_x$) are finite flat by Corollary \ref{flat-deg-infty}. Since the moduli problems $[\Gamma(p)],[\Gamma_1(p)],[\mu_p^{\times}]$ are finite \'etale over $\Ell_{\Z[1/p]}$, $U_x$ is a flat morphism of moduli problems, so that every $U_x$ in the diagram is flat away from the cusps, and finite (since it is a morphism of $\mathbb{P}^1_{\Z[1/p]}$-schemes). By miracle flatness (since all the moduli problems are smooth at infinity), it follows that every $U_x$ is finite flat. 

Thus, it is enough to show that the diagram commutes and the larger rectangle is Cartesian away from the cusps. Both triangles commute by Proposition \ref{degn-AL-infty}. Moreover, the top trapezoid is easily checked to commute away from the cusps, hence it commutes (because we can apply Proposition \ref{compactification-functor}), and the bottom trapezoid is commutative and Cartesian away from the cusps, as we proved above. Moreover, the top trapezoid is automatically Cartesian since two ``opposite'' maps are isomorphisms, and the conclusion follows.}

\medskip

\defi{The morphism $R: P(p) \rar P(p)$ is defined as follows: for any $x,y \in S$, the map $J_1(p,\mu_p) \overset{[(\iota_x)_{\ast}]}{\longrightarrow} P(p) \overset{R}{\rar} P(p) \overset{\iota_y^{\ast}}{\rar} J_1(p,\mu_p)$ is given by 
\begin{itemize}[noitemsep,label=\tiny$\bullet$]
\item $-U_p$ if $x=y$,
\item $w_p$ if $x \neq y$ and $(1,0) \in \{x,y\}$,
\item $\langle s-t\rangle w_p$ if $x=(s,1),y=(t,1)$ with $s,t \in \F_p$ distinct. 
\end{itemize}}

\medskip

\lem[good-open-subset-hecke]{Let $\ell$ be a prime number. There is a dense open subscheme $U \subset \mathbb{A}^1_{\Q}$ such that for any elliptic curve $E$ over a field $K$ of characteristic $0$, if $j(E) \in U$, then the $j(E/C)$ are all pairwise distinct, where $C$ runs through the cyclic subgroups of order $\ell$ of $E_{\overline{K}}$.}

\demo{Consider the moduli scheme $Y_1'(N,N)$ with $N=q\ell \geq 4$ (for some integer $q \geq 1$) attached to the moduli problem $[\Gamma_1'(N,N)]_{\Q}$ from Section \ref{gamma1m-structure}. By Proposition \ref{basic-with-bad}, Lemma \ref{open-comp} and Corollary \ref{bad-pb-with-geom-fibres}, $Y_1'(N,N)$ is smooth of relative dimension one, affine and integral. Consider the two morphisms $Y_1'(N,N) \rar \mathbb{A}^1_{\Q}$ given by \[(E,P,C) \longmapsto j(E/\langle qP\rangle),\, (E,P,C) \longmapsto j(E/C[\ell]),\] and let $Z \subset Y'_1(N,N)$ be the closed subscheme on which these morphisms agree: then either $Z=Y'_1(N,N)$ or $Z$ is a discrete union of finitely many closed points -- so that the open subset $\mathbb{A}^1\backslash j(Z)$ works.  

All we need to do, then, is find a single elliptic curve $E$ over some algebraically closed field $K$ and independent points $P,Q \in E[\ell]$ such that $j(E/\langle P\rangle) \neq j(E/\langle Q\rangle)$. Let $\tau \in \C$ be a complex number with positive imaginary part, and consider the points $P=\frac{\tau}{\ell}, Q=\frac{1}{\ell}$ of the elliptic curve $E=\frac{\C}{\Z\tau\oplus \Z}$. Then the $j$-invariant of $E/\langle P\rangle$ is $j(\tau/\ell)$, and the $j$-invariant of $E/\langle Q\rangle$ is $j(\ell\tau)$. Using the $q$-expansion of $j$, it is clear that when $\tau=it$ with $t \rar +\infty$, then $j(\tau/\ell) \sim e^{2\pi t/\ell}$, $j(\ell\tau) \sim e^{2\pi \ell t}$, so $j(\tau/\ell) \neq j(\ell\tau)$.}

\prop[complex-of-R]{$R$ commutes with the action of $\mathbb{T}'[\GL{\F_p}]$ and $R \circ u=0$.}

\demo{Every entry of the matrix of $R$ is $\mathbb{T}'$-equivariant by Lemma \ref{modified-hecke-commutes}, so $R$ commutes to $\mathbb{T}'$. 

Next, we show that $R$ commutes to $T$. A direct computation shows that $T \circ \iota_{(t,1)} = \iota_{(t-1,1)}$ for $t \in \F_p$, while $T \circ \iota_{(1,0)} = \iota_{(1,0)}$. Let $s,t \in \F_p$, then 
\begin{align*}
\iota_{(t,1)}^{\ast}T \circ R \circ [(\iota_{(s,1)})_{\ast}]&=(T^{-1} \circ \iota_{(t,1)})^{\ast} \circ R \circ [(\iota_{(s,1)})_{\ast}]=\iota_{(t+1,1)}^{\ast}\circ R \circ [(\iota_{(s,1)})_{\ast}],\\
\iota_{(t,1)}^{\ast}R \circ T \circ [(\iota_{(s,1)})_{\ast}]&=\iota_{(t,1)}^{\ast} \circ R \circ [(T \circ \iota_{(s,1)})_{\ast}] = \iota_{(t,1)}^{\ast} \circ R \circ [(\iota_{(s-1,1)})_{\ast}],
\end{align*}

and since, for $s,t \in \F_p$, $\iota_{(t,1)}^{\ast}R[(\iota_{(s,1)})_{\ast}]$ depends only from $s-t \in \F_p/\{\pm 1\}$, we see that $\iota_{(t,1)}^{\ast}R \circ T \circ [(\iota_{(s,1)})_{\ast}]=\iota_{(t,1)}^{\ast}T \circ R \circ [(\iota_{(s,1)})_{\ast}]$. Since, for $x \in S$, $\iota_x^{\ast}R[(\iota_{(1,0)})_{\ast}]$ and $\iota_{(1,0)}^{\ast}R[(\iota_{x})_{\ast}]$ only depend on whether $x = (1,0)$, we also see that $\iota_x^{\ast}\circ (R \circ T) \circ [(\iota_{(1,0)})_{\ast}]=\iota_x^{\ast}\circ (T \circ R)\circ [(\iota_{(1,0)})_{\ast}]$. 

Thus $R$ commutes with $T$. 

Similar computations using Lemma \ref{wp-relations} show that $R$ also commutes with $\Delta_{a,1}$ for $a \in \F_p^{\times}$ and $W$, so that $R$ commutes with $\GL{\F_p}$. 

Therefore, to prove that $R \circ u=0$, it is enough to show that $\iota_{(1,0)}^{\ast} \circ R \circ u = 0$, that is, that $\sum_{t \in \F_p}{w_p \circ [(u_{(t,1)})_{\ast}]} = U_p \circ u_{(1,0)}$ (as morphisms $J(p,p) \rar J_1(p,\mu_p)$). The closed subscheme $Z$ on which these morphisms agree is a closed subgroup scheme of $J(p,p)$, so, since $J(p,p)$ is smooth, it is enough to show that $Z$ contains a Zariski-dense subset of $J(p,p)(\Qbar)$.  

By Lemma \ref{good-open-subset-hecke}, let $U \subset D(j(j-1728)) \subset \mathbb{A}^1_{\Q}$ be a nonempty open subscheme such that for any elliptic curve $E$ with $j$-invariant in $U$, the $j$-invariants of the $E/C$, where $C$ runs through cyclic subgroups of $E$ of order $p$, are pairwise distinct. For some large $r \geq 1$, the points of $J(p,p)(\Qbar)$ in the image of $(j^{-1}(U) \times_{(\mu_N^{\times})_{\Q}} j^{-1}(U))^{\times r}$ are Zariski-dense in $J(p,p)(\Qbar)$: it is thus enough to show that for any $(E,P,Q),(E',P',Q') \in j^{-1}(U)(\Qbar) \subset Y(p,p)(\Qbar)$ such that $\langle P,\,Q\rangle_E=\langle P',\,Q'\rangle_{E'}$, $\sum_{t \in \F_p}{w_p \circ [(u_{(t,1)})_{\ast}]}$ and $U_p \circ u_{(1,0)}$ coincide on $[(E,P,Q)]-[(E',P',Q')]$.

In this case, since $[\Gamma_1'(p,p)]_{\Q}, [\Gamma_1(p)]_{\Q}$ are representable finite \'etale, their moduli schemes are finite \'etale above $U$, by \cite[Proposition 8.6.8]{KM}. Therefore, the pull-back under $D'_{1,1}$ of any Cartier divisor $[(E,S,\alpha)]$ (resp. $[(E',S',\alpha)]$) (where $S,S'$ are points of $E,E'$ of exact order $p$) is exactly the sum over all cyclic subgroups $C$ of $E$ of degree $p$ and not containing $S$ (resp. $C$ of $E'$ not containing $S'$) of the quadruples $[(E,S,C,\alpha)]$ (resp. $[(E,S,C',\alpha)]$). 

Now, the cyclic subgroups of $E[p]$ (resp. $E'[p]$) not containing $P$ (resp. $P'$) are exactly the $tP+Q$ (resp. $tP'+Q'$) over all $t \in \F_p$. If $t \in \F_p$, then by \cite[(2.8.6)]{KM}, $\langle tP+Q,\, -P\pmod{tP+Q}\rangle_{\langle tP+Q\rangle} = \langle tP+Q,\,-P\rangle_{E[p]} = \alpha$. Thus, using Lemma \ref{pushforward-is-weil} (for each map), we have in $\mrm{Div}(Y_1(p,\mu_p)_{\Qbar})$:
\begin{align*}
&\sum_{t \in \F_p}{w_p \circ [(u_{(t,1)})_{\ast}]((E,P,Q)-(E',P',Q'))} = \sum_{t \in \F_p}{w_p((E,tP+Q,\alpha)-(E',tP+Q,\alpha))} \\
&= \sum_{t \in \F_p}{(E/\langle tP+Q\rangle,-P \pmod{tP+Q},\alpha)-(E'/\langle tP'+Q'\rangle,-P' \pmod{tP'+Q'},\alpha)}\\
&= \sum_{t \in \F_p}{(E/\langle tP+Q\rangle, P \pmod{tP+Q},\alpha)}-\sum_{t \in \F_p}{(E'/\langle tP'+Q'\rangle,P' \pmod{tP'+Q'},\alpha)}\\
&= (D'_{p,1})_{\ast}(D'_{1,1})^{\ast}([(E,P,\alpha)-(E',P,\alpha)])\\
&= U_p\circ u_{(1,0)}([(E,P,Q)]-[(E',P',Q')]),
\end{align*}
 whence the conclusion.

}

\section{Uniformizing $P(p)$}
In the rest of this section, $\HH$ denotes the complex upper half-plane, endowed with the usual action of $\GL{\R}^+$ by homographies: $\begin{pmatrix}a & b\\c & d\end{pmatrix}\cdot \tau=\frac{a\tau+b}{c\tau+d}$. 

Given an positive integer $N$, $\Gamma_0(N),\Gamma_1(N),\Gamma(N)$ denotes as usual the subgroup of matrices $M \in SL_2(\Z)$ which are respectively upper-triangular modulo $N$, congruent to $\begin{pmatrix}1 & \ast\\0 & 1\end{pmatrix}$ modulo $N$, or congruent to $I_2$ modulo $N$.

\prop[uniformize-x1pmu]{Let $X_1(p)^{an}$ be the classical compact Riemann surface whose non-cuspidal points are given by the quotient $\Gamma_1(p) \backslash \HH$ (see for example \cite[Chapter 2]{DS}). There exists a morphism of locally ringed spaces $X_1(p)^{an} \times \F_p^{\times} \times S \rar X_1(p,\mu_p)^S$ making $X_1(p)^{an} \times \F_p^{\times} \times S$ the analytification of $X_1(p,\mu_p)^S_{\C}$ in the sense of \cite[Exp. XII, (1.1)]{SGA1}. Under this morphism, a point $(\tau,b,x) \in \HH \times \F_p^{\times} \times S$ is mapped to $\left(\left(\frac{\C}{\tau\Z\oplus\Z},\frac{1}{p}\right),e^{-2i\pi b/p},x\right) \in X_1(p,\mu_p)^S(\C)$. In particular, this uniformization preserves the $j$-invariant (i.e. maps the algebraic $j$-invariant to the analytic $j(\tau)=\frac{1}{q}+744+196884q+\ldots$). 

Under this identification, a matrix $M \in \GL{\F_p}$ acts as follows: given $x,x' \in S$ and $\lambda \in \F_p^{\times}$ such that $\lambda (M^T)^{-1}x=x'$, then the composition 
\[X_1(p)^{an} \times \F_p^{\times} \simeq X_1(p)^{an}\times \F_p^{\times} \times \{x\} \overset{M}{\rar} X_1(p)^{an} \times \F_p^{\times} \times \{x'\} \simeq X_1(p)^{an} \times \F_p^{\times}\]

is exactly $\langle \lambda\rangle \times [\det{M}]$, where $[\det{M}]: \F_p^{\times} \rar \F_p^{\times}$ is the multiplication by $\det{M}$.
}

\rem{The reason for the minus sign comes from the uniformization that we will adopt for $X(p,p)$, and the sign convention for the Weil pairing in \cite[(2.8.5.2)]{KM}. }

Since we will in fact need several variants of this result, we prove a generalized version of this statement. It is not a new result by any means, but we include a proof since it is not spelled out in \cite{KM} or the otherwise very thorough \cite[Appendix A]{Ribet-Stein}. 

\lem[uniformize-general]{Let $a,b,c \geq 1$ be pairwise coprime integers with $a \geq 3$ or $b \geq 4$. Let $X^a$ (resp. $X_a$) be the compact connected Riemann surface attached (see for example \cite[Chapter 2]{DS}) to the congruence subgroup $\Gamma=\Gamma(a)\cap \Gamma_1(b) \cap \Gamma_0(c)$ (resp. the compactified moduli scheme attached to the finite \'etale moduli problem $[\Gamma(a)]_{\C} \times [\Gamma_1(b)]_{\C} \times [\Gamma_0(c)]_{\C}$). 

There exists a morphism of locally ringed spaces $\mrm{An}: X^a \times (\Z/a\Z)^{\times} \rar X_a$ such that, for any $(\Gamma\tau,v) \in \HH/\Gamma \times (\Z/a\Z)^{\times}$, one has \[\mrm{An}(\Gamma\tau,v)=\left[\left(\frac{\C}{\tau\Z\oplus\Z},\frac{v\tau}{a},\frac{1}{a},\frac{1}{b},\langle \frac{1}{c}\rangle\right)\right].\] Moreover, $\mrm{An}$ satisfies the universal property of an analytification in the sense of \cite[Exp. XII, (1.1)]{SGA1}, and $j\circ \mrm{An}$ is exactly the modular $j$-invariant. }

\demo{First, we know by Corollary \ref{smooth-at-infinity-gamma0m} that $X_a$ is a smooth projective scheme of dimension one over $\C$. Therefore, its analytification $X'$ is (by \cite[Exp. XII, (2.6), (3.1), (3.2)]{SGA1}) a (perhaps disconnected) compact Riemann surface. 

Next, we construct a morphism of locally ringed spaces $Y^a \rar Y_a$, where $Y^a=\HH/\Gamma$ and $Y_a$ is the moduli scheme attached to the finite \'etale moduli problem $[\Gamma(a)]_{\C} \times [\Gamma_1(b)]_{\C} \times [\Gamma_0(c)]_{\C}$.\\

\emph{Step 1: A ``versal'' elliptic curve over $\HH$:}

Let $R=\OO(\HH)$ be the ring of holomorphic functions on $\HH$. It is a Bezout domain by \cite[(26.5)]{Forster}. In particular, the image of the obvious morphism $\theta: \HH \rar \Sp{R}$ of locally ringed spaces is dense, and meets every nonempty closed subscheme of $\Sp{R}$ cut out by a finitely generated ideal.    

Let $G_k(\tau)=\sum_{(m,n) \in \Z^2\backslash \{0\}}{(m\tau+n)^{-k}}$ for each even $k > 2$. Consider the cubic $\mathcal{E}^{an}$ over $R$ of equation $y^2=x^3-15G_4x-35G_6$ (that is, $\mathcal{E}^{an}$ is the closed subscheme of $\mathbb{P}^2_R$ with equation $y^2z=x^3-15G_4xz-35G_6z^3$). Its discriminant $\delta \in R$ is $16(4 \cdot 15^3G_4^3-35^2G_6^2)$, so it is a modular form of weight $12$ and level $1$. The first terms in the $q$-expansion of $\delta$ are given by (using \cite[p. 5]{DS})

\begin{align*}
\delta(\tau)&=16 \cdot 4\left(\frac{2 \cdot 15\pi^4}{90}+\frac{30(2i\pi)^4}{6}q+O(q^2)\right)^3-16 \cdot 27\left(\frac{70\pi^6}{945}+\frac{70(2i\pi)^6}{120}q+O(q^2)\right)^2\\
&= \pi^{12}\left[64 \cdot \left(\frac{1}{3}+80q+O(q^2)\right)^3-16\cdot 27\left(\frac{2}{27}-\frac{7\cdot 16}{3}q+O(q^2)\right)^2\right]\\
&= \pi^{12}\left[\frac{64}{27}+\frac{64 \cdot 80}{3}q+O(q^2)-\frac{64}{27}+\frac{64 \cdot 7 \cdot 16}{3}+O(q^2)\right]\\
&= \frac{64\pi^{12} \cdot 16(5+7)}{3}q+O(q^2)=(2\pi)^{12}q+O(q^2),
\end{align*}

hence $\delta=(2\pi)^{12}\Delta$ does not vanish on $\HH$ (where $\Delta$ is the weight $12$ normalized cusp form of level one). Hence $\mathcal{E}^{an}$ is an elliptic curve over $R$. Its (algebraic) $j$-invariant is $-1728\frac{(-60G_4)^3}{\delta} \in R$, which is, by \cite[p.7]{DS}, the modular $j$-invariant. We know that, for any $\tau \in \HH$, by \cite[Proposition VI.3.6]{AEC1}, 
\begin{align*}
z \in \C/(\tau\Z\oplus\Z) \longmapsto [z]_{\tau} := [\wp_{\tau\Z\oplus\Z}(z):2\wp'_{\tau\Z\oplus\Z}(z):1] \in &\left(\mathcal{E}^{an}\times_{\Sp{R},\tau} \Sp{\C}\right)(\C) \\ &\subset (\mathbb{P}^2_R \times_{\Sp{R},\tau} \Sp{\C})(\C)
\end{align*} defines an isomorphism. 

Let $c,d \in \Z$ and $N \geq 1$ be coprime integers (that is, $c\Z+d\Z+N\Z=\Z$), and $P \in \mathcal{E}^{an}(R)$ be given by $[\wp_{\tau\Z\oplus\Z}\left(\frac{c\tau+d}{N}\right):2\wp'_{\tau\Z\oplus\Z}\left(\frac{c\tau+d}{N}\right):1]$ (it is a well-defined point, because this can be checked at each $\tau \in \HH$). For any $e \mid N$, there is a (universal) closed subscheme $Z_e$ of $\Sp{R}$ on which $eP=0$: it is a base change of the diagonal $\mathcal{E}^{an} \rar \mathcal{E}^{an}\times_R \mathcal{E}^{an}$, which is a finitely presented closed immersion by \cite[Lemma 0818]{Stacks}. Hence $Z_e \rar \Sp{R}$ is a finitely presented closed immersion and $Z_e$ is cut out by a finitely generated ideal. But $Z_e$ does not meet $\theta(\HH)$ unless $e=N$, in which case $Z_e \supset \theta(\HH)$. It follows that $P \in \mathcal{E}^{an}(R)$ is a point of exact order $N$. 

Similarly, given functions $f:\HH \rar \C$, $g: \HH \rar \HH$, we write $[f(\tau)]_{g(\tau)}$ or $[f]_g$ for the point $(\tau \longmapsto [\wp_{g(\tau)\Z\oplus\Z}(f(\tau)):2\wp'_{g(\tau)\Z\oplus\Z}(f(\tau)):1]) \in \operatorname{Mor}_{\C}(\Sp{R},\mathcal{E}^{an})$. When there is no index, it is assumed that $g$ is the identity. 

Let $R'$ be the product of copies of $R$ indexed by $(\Z/a\Z)^{\times}$. The enriched elliptic curve 
\[\epsilon=\left(\mathcal{E}^{an}_{R'}/R',\left(\left(\left[\frac{v\tau}{a}\right]\right)_{v \in (\Z/a\Z)^{\times}},\left[\frac{1}{a}\right]_{R'}\right),\left[\frac{1}{b}\right]_{R'},\langle\left[\frac{1}{c}\right]_{R'}\rangle\right) \]

 defines a morphism $\upsilon: \Sp{R'} \rar Y_a$. \\

\emph{Step 2: Passing to the quotient by $\Gamma$}

Now, there is a natural left action of $SL_2(\Z)$ on $\HH$, which induces a natural (componentwise) action of $SL_2(\Z)$ on $\Sp{R'}$. Let us show that $\upsilon$ factors through a morphism $\upsilon': \Sp{(R')^{\Gamma}} \rar Y_a$. Since $Y_a$ is affine by Proposition \ref{over-Ell-concrete}, it is enough to show that $\epsilon$ is invariant under $\Gamma$. 

Let $g=\begin{pmatrix}u & v\\w & t\end{pmatrix} \in\Gamma$. Then $g \cdot \epsilon$ is the elliptic curve $E$ over $R'$ given by the equation \[y^2=x^3-15G_4\left(\frac{u\cdot+v}{w\cdot+t}\right)x-35G_6\left(\frac{u\cdot+v}{w\cdot+t}\right) \Leftrightarrow y^2=x^3-15(w\tau+t)^4G_4x-35(w\tau+d)^6G_6, \] with attached structure
\begin{align*}
\left(\left[z\frac{g \cdot \tau}{a}\right]_{g \cdot \tau},\left[\frac{1}{a}\right]_{g \cdot \tau},\left[\frac{1}{b}\right]_{g \cdot \tau},\langle\left[\frac{1}{c}\right]_{g \cdot \tau}\rangle\right).
\end{align*}

Since one can straightforwardly check that, with $(\wp,\wp_g)=(\wp_{\tau\Z\oplus\Z},\wp_{g\cdot\tau\Z\oplus\Z})$,
\[\wp_{g}(z)=(w\tau+t)^2\wp((w\tau+t)z),\quad \wp'_{g}(z)=(w\tau+t)^3\wp((w\tau+t)z),\] this data is equal to
\begin{align*}
\left(\left[(w\tau+t)^2\wp\left(z \frac{u\tau+v}{a}\right):2(w\tau+z)^3\wp'\left(z\frac{u\tau+v}{a}\right):1\right]\right)_{z \in (\Z/a\Z)^{\times}},\\
\left(\left[(w\tau+t)^2\wp\left(\frac{w\tau+t}{a}\right):2(w\tau+z)^3\wp'\left(\frac{w\tau+t}{a}\right):1\right]\right)_{z \in (\Z/a\Z)^{\times}},\\
\left(\left[(w\tau+t)^2\wp\left(\frac{w\tau+t}{b}\right):2(w\tau+z)^3\wp'\left(\frac{w\tau+t}{b}\right):1\right]\right)_{z \in (\Z/a\Z)^{\times}},\\
\left\langle\left(\left[(w\tau+t)^2\wp\left(\frac{w\tau+t}{c}\right):2(w\tau+z)^3\wp'\left(\frac{w\tau+t}{c}\right):1\right]\right)_{z \in (\Z/a\Z)^{\times}}\right\rangle.
\end{align*}

Since $(u,v) \equiv (t,w) \equiv (1,0) \pmod{a}$, since $w \equiv 0 \pmod{bc}$ and $u \equiv t \equiv 1 \pmod{b}$, this data is equal to 
\begin{align*}
\left(\left[(w\tau+t)^2\wp\left(z \frac{\tau}{a}\right):2(w\tau+z)^3\wp'\left(z\frac{\tau}{a}\right):1\right]\right)_{z \in (\Z/a\Z)^{\times}},\\
\left(\left[(w\tau+t)^2\wp\left(\frac{1}{a}\right):2(w\tau+z)^3\wp'\left(\frac{1}{a}\right):1\right]\right)_{z \in (\Z/a\Z)^{\times}},\\
\left(\left[(w\tau+t)^2\wp\left(\frac{1}{b}\right):2(w\tau+z)^3\wp'\left(\frac{1}{b}\right):1\right]\right)_{z \in (\Z/a\Z)^{\times}},\\
\left\langle t \cdot \left(\left[(w\tau+t)^2\wp\left(\frac{1}{c}\right):2(w\tau+z)^3\wp'\left(\frac{1}{c}\right):1\right]\right)_{z \in (\Z/a\Z)^{\times}}\right\rangle.
\end{align*}

The map $[x:y:z] \in \mathcal{E}_{R'}^{an} \longmapsto [(w\tau+t)^2x:(w\tau+t)^3y:z] \in E'$ is an isomorphism and maps the additional data attached to $\epsilon$ to exactly the data above, so $\mathcal{E}^{an}_{R'}$ (with its level $\Gamma$ structure) is isomorphic to $g^{\ast}\mathcal{E}^{an}_{R'}$. Hence $\upsilon$ is invariant by the action of $\Gamma$, and $\upsilon$ factors through a morphism $\upsilon': \Sp{(R')^{\Gamma}} \rar Y_a$. Since the ring of functions of $Y^a \times (\Z/a\Z)^{\times} \simeq \Gamma \backslash (\HH \times (\Z/a\Z)^{\times})$ is exactly $(R')^{\Gamma}$, we get a morphism $\lambda: Y^a \times (\Z/a\Z)^{\times} \rar Y_a$ of locally ringed spaces.  

Given $(\tau,v) \in \HH \times (\Z/a\Z)^{\times}$, $\lambda(\tau,v)$ is by construction the following elliptic curve with level $\Gamma$ structure \[y^2=x^3-15G_4(\tau)x-35G_6(\tau),\left(\left[\frac{v\tau}{a}\right]_{\tau},\left[\frac{1}{a}\right]_{\tau},\left[\frac{1}{b}\right]_{\tau}, \left\langle\left[\frac{1}{c}\right]_{\tau}\right\rangle\right).\]

Thus, its uniformization is $\left(\frac{\C}{\tau\Z\oplus\Z},\frac{v\tau}{a},\frac{1}{a},\frac{1}{b},\langle \frac{1}{c}\rangle\right)$, and in particular the $j$-invariant of this elliptic curve is exactly $j(\tau)$, so that $j_{alg} \circ \lambda$ is exactly the modular $j$-invariant. \\

\emph{Step 3: $\lambda$ induces a bijection $Y^a \times (\Z/a\Z)^{\times} \rar Y_a(\C)$}

Let $(\tau,z),(\tau,z') \in \HH \times (\Z/a\Z)^{\times}$ be such that $\lambda(\tau,z)=\lambda(\tau',z')$. By analytification, there is an isomorphism of complex tori $f: \frac{\C}{\tau\Z\oplus\Z} \rar \frac{\C}{\tau'\Z\oplus\Z}$ mapping $\frac{z\tau}{a}$ (resp. $\frac{1}{a}$, $\frac{1}{b}$, $\langle \frac{1}{c}\rangle$) to $\frac{z'\tau'}{a}$ (resp. $\frac{1}{a}$, $\frac{1}{b}$, $\langle \frac{1}{c}\rangle)$. By \cite[Proposition 1.3.2]{DS}, $f$ is of the form $q \longmapsto mq+t$ where $m,t \in \C$ and $m(\tau\Z\oplus\Z)=\tau'\Z\oplus\Z$. Since $f(0)=0$, we can take $t=0$. Since $f$ preserves the Weil pairing, we must have $e^{-2i\pi z/a}=e^{-2i\pi z'/a}$ by \cite[(2.8.5.2)]{KM}, so $z=z'$. 

We know that there are $u,v,w,t \in \Z$ such that $\tau'=m(u\tau+v),1=m(w\tau+t)$. In particular, $ut-vw=\pm 1$. Since $\tau'=\frac{u\tau+v}{w\tau+t} \in \HH$, we know that $ut-vw=1$. Therefore, \begin{align*}
m\frac{z\tau}{a} &\equiv \frac{z\tau'}{a}\equiv m\frac{z(u\tau+v)}{a}  \pmod{m(\tau\Z\oplus\Z)},\\
m\frac{1}{a} &\equiv \frac{1}{a}\equiv m\frac{w\tau+t}{a} \pmod{m(\tau\Z\oplus\Z)},\\
m\frac{1}{b} &\equiv \frac{1}{b}\equiv m\frac{w\tau+t}{b} \pmod{m(\tau\Z\oplus\Z)},\\
m\langle \frac{1}{c}\rangle &\equiv \langle \frac{1}{c}\rangle \equiv m\langle \frac{w\tau+t}{c}\rangle \pmod{m(\tau\Z\oplus\Z)},
\end{align*}

hence $g=\begin{pmatrix} u & v\\w & t\end{pmatrix} \in \Gamma$ and $g\tau=\tau'$. So $\lambda: Y^a \times (\Z/a\Z)^{\times} \rar Y_a(\C)$ is injective. 

Let $E$ be an elliptic curve over $\C$ and $(P_a,Q_a)$ be a basis of its $a$-torsion, $P_b \in E(\C)$ have exact order $b$, and $C \leq E$ be a cyclic subgroup of order $c$. By \cite[Theorem VI.5.1]{AEC1} there is some $\tau \in \HH$ such that $E$ is isomorphic to the elliptic curve $F$ of equation $y^2=x^3-15G_4(\tau)x-35G_6(\tau)$, and we can find coefficients $u,v,w,t \in \Z/abc\Z$ such that $ut-wx=\delta \in (\Z/abc\Z)^{\times}$ is congruent to $1$ mod $bc$ and such that under the isomorphism $E \rar F$, $(P_a,Q_a,P_b,C)$ is mapped to 
\[\left(\left[\frac{u\tau+v}{a}\right]_{\tau},\left[\frac{w\tau+t}{a}\right]_{\tau},\left[\frac{w\tau+t}{b}\right]_{\tau},\langle\left[\frac{w\tau+t}{c}\right]_{\tau}\rangle\right),\]

where $[z]_{\tau}$ denotes the point $\left[\wp_{\tau\Z\oplus\Z}(z):2\wp'_{\tau\Z\oplus\Z}(z):1\right]$ for $z \in \C/(\tau\Z\oplus\Z)$. 

In particular, let $(u_1,v_1)=\delta^{-1}(u,v)$, and pick lifts $u',v',w',t' \in \Z$ of $u,v,w,t$ such that $u't'-w'v'=1$, and let $\tau'=\frac{u'\tau+v'}{w'\tau+t'}$. Let $F'$ be the complex elliptic curve with equation $y^2=x^3-15G_4(\tau')-35G_6(\tau')$, so that \[\alpha: [x:y:z] \in F \mapsto [(w'\tau+t')^2x:(w'\tau+t')^3y:z] \in F'\] is an isomorphism. 
Then one has 
\begin{align*}
&\alpha\left(\left[\frac{u\tau+v}{a}\right]_{\tau},\left[\frac{w\tau+t}{a}\right]_{\tau},\left[\frac{w\tau+t}{b}\right]_{\tau},\langle\left[\frac{w\tau+t}{c}\right]_{\tau}\rangle\right)\\
&= \left(\left[\frac{\delta\tau'}{a}\right]_{\tau'},\left[\frac{1}{a}\right]_{\tau'},\left[\frac{1}{b}\right]_{\tau'},\langle\left[\frac{1}{c}\right]_{\tau'}\rangle\right),
\end{align*}

so that $(E,P_a,Q_a,P_b,C) \in Y_a(\C)$ is $\lambda((\tau',\delta))$.  \\

\emph{Step 4: The analytification of $\lambda: Y^a \times (\Z/a\Z)^{\times} \rar Y_a \subset X_a$ extends holomorphically to $X^a \times (\Z/a\Z)^{\times}$}

Let $Y'$ be the analytification of $Y_a$, it is a (possibly disconnected) Riemann surface contained in $X'$ and such that $X' \backslash Y'$ is finite. The morphism $\lambda$ factors through a holomorphic morphism $\lambda': Y^a \times (\Z/a\Z)^{\times}\rar Y'$, which is bijective by Step 3 and \cite[Exp. XII, (1.1)]{SGA1} (hence a biholomorphism, but we will not use it). Let $\iota: Y^a \times \Z/a\Z)^{\times} \overset{\lambda'}{\rar} Y' \subset X'$. Our claim is that $\iota$ extends holomorphically to $X^a \times (\Z/a\Z)^{\times}$. 

Let $c \in (X^a \backslash Y^a)\times (\Z/a\Z)^{\times}$, and let $f$ be any meromorphic function on $X'$. Let us show that $f \circ \iota$ extends meromorphically to $c$. Let $\beta: D \rar X^a \times (\Z/a\Z)^{\times}$ be a biholomorphism to some open subspace of $X^a \times (\Z/a\Z)^{\times}$, where $D$ is the open unit disk, $\beta(D \backslash \{0\}) \subset Y^a \times (\Z/a\Z)^{\times}$ and $\beta(0)=c$. So it is enough to check that the singularity at $0$ of $f \circ \iota\circ \beta$ is not essential. If it were, by Casorati-Weierstrass (see for instance \cite[Theorem V.3.2]{LangComplex}), there would be a dense $S \subset \C$ such that for any $s \in S$, $(f \circ \iota\circ \beta)^{-1}(s)$ is infinite. Since $\iota,\beta$ are injective and $f$ has finite fibres, and we would get a contradiction.  

Hence, for any holomorphic $\phi: X' \rar \mathbb{P}^1(\C)$, $\phi \circ \iota: Y^a \times (\Z/a\Z)^{\times}$ extends holomorphically to $X^a \times (\Z/a\Z)^{\times}$. By \cite[14.13]{Forster}, this implies that for any sequence $x_n \in Y^a \times (\Z/a\Z)^{\times}$ converging to some point of $X^a \times (\Z/a\Z)^{\times}$, $\iota(x_n)$ has at most one limit point. Since $X'$ is compact, this implies that $\iota$ extends continuously. By picking charts and using Morera's theorem (see for instance \cite[Theorem III.7.9]{LangComplex}, although the book does not mention the relevant corollary: a continuous function on an open subset of $\C$ which is holomorphic outside finitely many points is holomorphic), it is straightforward to check that the continuous extension of $\iota$ to $X^a \times (\Z/a\Z)^{\times}$ is holomorphic. \\

\emph{Step 5: $\iota$ is isomorphic to the inclusion $Y^a \times (\Z/a\Z)^{\times} \subset X^a \times (\Z/a\Z)^{\times}$.}

Since $X^a \times (\Z/a\Z)^{\times}$ is compact, its image under $\iota$ is compact and contains a co-finite open subspace of $X'$, so $\iota$ is surjective. Since the restriction of $\iota$ to the co-finite open subset $Y^a \times (\Z/a\Z)^{\times}$ is injective (by Step 3 and \cite[Exp. XII, (1.1)]{SGA1}), the degree of $\iota$ with respect to any connected component of $X'$ (\cite[4.24]{Forster}) is one. Hence $\iota$ is an isomorphism.

}

\cor[differentials-x1pmu]{This identification induces an isomorphism of right $\C[\GL{\F_p}]$-modules \[\iota_{1,p}: H^0(X_1(p,\mu_p)^S_{\C},\Omega^1) \rar (\mathcal{S}_2(\Gamma_1(p)) \otimes \C[\F_p^{\times}]) \otimes_{\Z[B]} \Z[\GL{\F_p}]\] such that:
\begin{itemize}[noitemsep,label=\tiny$\bullet$]
\item $\begin{pmatrix}a & \ast\\0 & d\end{pmatrix} \in B$ acts on $\mathcal{S}_2(\Gamma_1(p)) \otimes \C[\F_p^{\times}]$ by $\langle d\rangle \otimes [ad]$, where $[ad]: \F_p^{\times} \rar \F_p^{\times}$ is the multiplication by $ad$. 
\item If $t \in \F_p$, $M=T^t$ and $x=(t,1)$ (resp. $M=W$ and $x=(1,0)$), then, for all $f \in \mathcal{S}_2(\Gamma_1(p)), b \in \F_p^{\times}$, $\iota_{1,p}^{-1}(f \otimes [b] \otimes M)$ is the holomorphic differential on $X_1(p)^{an} \times \F_p^{\times} \times S$ which is equal to $f(\tau)d\tau$ on the component $(b^{-1},x)$ and zero on the other components.     
\end{itemize} 
}

\demo{The uniformization of Proposition \ref{uniformize-x1pmu} yields by \cite[Exp. XII, Corollaire 4.3]{SGA1} an isomorphism of right $\C[\GL{\F_p}]$-modules \[\iota: H^0(X_1(p,\mu_p)^S_{\C},\Omega^1) \rar \bigoplus_{(b,x) \in \F_p^{\times} \times S}{\mathcal{S}_2(\Gamma_1(p)) \otimes [(b,x)]}.\] So all we need to do is construct a $\C[\GL{\F_p}]$-isomorphism \[j: \bigoplus_{(b,x) \in \F_p^{\times} \times S}{\mathcal{S}_2(\Gamma_1(p)) \otimes [(b,x)]} \rar (\mathcal{S}_2(\Gamma_1(p)) \otimes \C[\F_p^{\times}]) \otimes_{\Z[B]} \Z[\GL{\F_p}]\] with the correct compatibilities. 

Let $x,x' \in S$ and $M \in \GL{\F_p}$, then write $\lambda (M^T)^{-1}x=x'$ for some $\lambda \in \F_p^{\times}$. Then we know that, given a differential $\omega=f(\tau)d\tau$ with support on the component $(b,x')$ of $X_1(p) \times \F_p^{\times} \times S$, $M^{\ast}\omega = (\langle \lambda\rangle)^{\ast}(f(\tau)d\tau)$ and has support on the component $(b(\det{M})^{-1},x)$. In other words, we know that $M^{\ast}(f \cdot [(b,x')])= ((\langle \lambda \rangle)^{\ast}f) \cdot [(b(\det{M})^{-1},x)]$. 

Let us define $j$ by
\[ f \cdot (b,(t,1)) \longmapsto  f \otimes [b^{-1}] \otimes [T^t],\quad f \cdot (b,(1,0)) \longmapsto [b^{-1}] \otimes [W].\]
Then the second condition is satisfied for $\iota_{1,p} := j \circ \iota$ and $\iota_{1,p}$ is clearly a $\C$-linear isomorphism. It remains to show that the first condition is satisfied and that $\iota_{1,p}$ is a morphism of $\C[\GL{\F_p}]$-modules. 

Let $M=\begin{pmatrix} a &\ast\\0 & d\end{pmatrix} \in B$ and $\omega \in H^0(X_1(p,\mu_p)_{\C}^S,\Omega^1)$ be such that its support is contained in the component $(b^{-1},(0,1))$. Since $d(M^T)^{-1}\begin{pmatrix}0 \\1\end{pmatrix}=\begin{pmatrix}0 \\1\end{pmatrix}$, the support of the differential $M^{\ast}\omega=(\underline{ad})^{\ast}\langle d\rangle^{\ast}\omega$ is contained in the component $((abd)^{-1},(0,1))$. Let $\iota_{1,p}(\omega)=f \cdot [(b,(0,1))]$, then one has $\iota_{1,p}(M^{\ast}\omega) = (\langle \lambda\rangle)^{\ast}f \cdot [(abd,(0,1))]$.  \\

Finally, we need to show that $\iota_{1,p}$ is a morphism of right $\C[\GL{\F_p}]$-modules. 

Let $f \in \mathcal{S}_2(\Gamma_1(p)), b \in \F_p^{\times}$, $x=(t,1)$ and $M=T^t$ (resp. $x=(1,0)$ and $M=W$), and $N \in \GL{\F_p}$. Let $\omega$ be the differential equal to $f(\tau)d\tau$ on the component $(b^{-1},x)$. Then $N^{\ast}\omega$ is the differential $(\langle\lambda\rangle^{\ast}f)(\tau)d\tau$ on the component $(b^{-1}(\det{N})^{-1},x')$ where $\lambda (N^T)^{-1}x'=x$, that is, $\lambda x'=N^T x$.
Thus $\iota_{1,p}(\omega) = f \otimes [b] \otimes M$, while $\iota_{1,p}(N^{\ast}\omega)=(\langle \lambda\rangle^{\ast} f) \otimes [b\det{N}] \otimes N'$, where $N'$ is $T^s$ (resp. $W$) if $x'=(s,1)$ (resp. $(1,0)$). 

Since $\lambda L_2(N')=\lambda x'=N^Tx=C_2(N^TM^T)=L_2(MN)$, we can write $MN=\begin{pmatrix} a & \ast\\0  &\lambda\end{pmatrix}N'$ where $a\lambda=a\lambda\det{N'}=\det{M}\det{N}=\det{N}$. 

Thus 
\begin{align*}
\iota_{1,p}(\omega) \mid N &= (f \otimes [b]) \otimes [MN] = \left((f \otimes [b]) \mid \begin{pmatrix}\lambda^{-1}\det{N}&\ast\\0 & \lambda\end{pmatrix} \right)\otimes [N'] \\
&= \left(\langle \lambda\rangle^{\ast}f \otimes [\det{N} b]\right),
\end{align*}
 
whence the conclusion. 

}

\prop[hecke-x1pmu]{Let $\mathcal{S}_2(\Gamma_1(p),\F_p^{\times},S)$ denote $(\mathcal{S}_2(\Gamma_1(p)) \otimes \C[\F_p^{\times}]) \otimes_{\Z[B]} \Z[\GL{\F_p}]$ for the sake of easier notation.
Let $\ell \neq p$ be a prime number and $n \in \F_p^{\times}$. Then the following diagrams commute:
\small{\[
\begin{tikzcd}[ampersand replacement=\&]
\& H^0(X_1(p,\mu_p)^S_{\C},\Omega^1) \arrow{r}{\iota_{1,p}} \arrow{d}{T_{\ell}^{\ast}} \& \mathcal{S}_2(\Gamma_1(p),\F_p^{\times},S) \arrow{d}{(T_{\ell} \otimes \mrm{id}) \otimes \mrm{id}}\&\\
\& H^0(X_1(p,\mu_p)^S_{\C},\Omega^1) \arrow{r}{\iota_{1,p}}  \& \mathcal{S}_2(\Gamma_1(p),\F_p^{\times},S) \&\\
H^0(X_1(p,\mu_p)^S_{\C},\Omega^1) \arrow{r}{\iota_{1,p}} \arrow{d}{\langle n\rangle^{\ast}} \& \mathcal{S}_2(\Gamma_1(p),\F_p^{\times},S) \arrow{d}{(\langle n\rangle \otimes \mrm{id}) \otimes \mrm{id}}\& H^0(X_1(p,\mu_p)^S_{\C},\Omega^1) \arrow{r}{\iota_{1,p}} \arrow{d}{\underline{n}^{\ast}} \& \mathcal{S}_2(\Gamma_1(p),\F_p^{\times},S)\arrow[swap]{d}{(\mrm{id} \otimes \underline{n}) \otimes \mrm{id}}\\
H^0(X_1(p,\mu_p)^S_{\C},\Omega^1) \arrow{r}{\iota_{1,p}}  \& \mathcal{S}_2(\Gamma_1(p),\F_p^{\times},S) \& H^0(X_1(p,\mu_p)^S_{\C},\Omega^1) \arrow{r}{\iota_{1,p}}  \& \mathcal{S}_2(\Gamma_1(p),\F_p^{\times},S)
\end{tikzcd}
\]}
In each case, the Hecke operators on the right hand side are those defined analytically, with the convention of \cite[Chapter 5.2]{DS}. 
}

\demo{Since all arrows are morphisms of right $\C[\GL{\F_p}]$-modules, it is enough to prove the claim when replacing all the diagrams by
\[
\begin{tikzcd}[ampersand replacement=\&] 
H^0(X_1(p,\mu_p)_{\C},\Omega^1) \arrow{r}{\iota_{1,p}} \arrow{d} \&\mathcal{S}_2(\Gamma_1(p)) \otimes \C[\F_p^{\times}] \arrow{d}\\
H^0(X_1(p,\mu_p)_{\C},\Omega^1) \arrow{r}{\iota_{1,p}} \& \mathcal{S}_2(\Gamma_1(p)) \otimes \C[\F_p^{\times}]
\end{tikzcd}
\]
For $\underline{n}$, the statements are straightforward from the construction of $\iota_{1,p}$. 

For $\langle n\rangle$, the statement is a direct consequence of the fact that the following diagram commutes:
\[
\begin{tikzcd}[ampersand replacement=\&]
(\Gamma_1(p)\backslash \HH) \times \F_p^{\times} \arrow{d}\arrow{rrr}{\Gamma_1(p)\tau\mapsto\Delta_{n^{-1},n}\Gamma_1(p)\tau}\&\&\& (\Gamma_1(p) \backslash \HH) \times \F_p^{\times}\arrow{d}\\
X_1(p,\mu_p)\arrow{rrr}{\langle n\rangle} \&\&\& X_1(p,\mu_p)
\end{tikzcd}
\]

Now, we only discuss Hecke operators. Let $\ell \neq p$ be a prime and $\omega \in H^0(X_1(p,\mu_p)_{\C},\Omega^1)$ which is nonzero only on the component $b^{-1}$ (for $b \in \F_p^{\times}$), and is equal to $f(\tau)d\tau$ on this component, for some $f \in \mathcal{S}_2(\Gamma_1(p))$: thus $\iota_{1,p}(\omega)=f \otimes [b]$. 

Then $T_{\ell}^{\ast}\omega=(D_{1,1}^{\ast})^{\ast}((D_{\ell,1})_{\ast})^{\ast}\omega$. Now, by Lemma \ref{uniformize-general} (and Lemma \ref{with-extra-functions}), $X_1(p,\mu_p,\Gamma_0(\ell))$ can be uniformized as $[((\Gamma_1(p) \cap \Gamma_0(\ell)) \backslash \HH^{\ast}] \times \F_p^{\times}$, by mapping $(\tau,z)$ to $\left(\frac{\C}{\tau\Z\oplus\Z},\frac{1}{p},\langle\frac{1}{\ell}\rangle,e^{-2i\pi z/p}\right)$, 

and the following diagram commutes:
\[
\begin{tikzcd}[ampersand replacement=\&]
\HH \times \F_p^{\times} \arrow{d}\& \HH \times \F_p^{\times} \arrow{l}{M}\arrow{d}\arrow{r}{N}\& \HH \times \F_p^{\times} \arrow{d}\\
X_1(p,\mu_p) \& X_1(p,\mu_p,\Gamma_0(\ell)) \arrow{l}{D_{\ell,1}} \arrow{r}{D_{1,1}} \& X_1(p,\mu_p)
\end{tikzcd}
\]

where $M(\tau,z)=(\delta_{\ell}\cdot (\ell\tau),z)$, where $\delta_{\ell} \in \Gamma(\ell)$ is any matrix whose reduction mod $p$ is $\Delta_{\ell^{-1},\ell}$, and $N(\tau,z)=(\tau,z)$. Hence the pull-back of $\omega'=D_{\ell,1}^{\ast}\omega$ to (the middle) $\HH \times \F_p^{\times}$ is $\ell(f\mid \delta_{\ell})(\ell\tau)d\tau$ on the component $b^{-1}$ and $0$ on the other components.   

To compute $(D_{1,1}^{\ast})^{\ast}\omega'$, we consider the following commutative diagram, where $\upsilon,\upsilon'$ are uniformizations:
\[
\begin{tikzcd}[ampersand replacement=\&]
\HH \times \F_p^{\times} \arrow{dr}{\upsilon'}\arrow{r}{(g_v)_{v}} \&\coprod\limits_{v \in \mathbb{P}^1(\F_{\ell})}{\HH \times \F_p^{\times}} \arrow{r}{N}\arrow{d}{(\pi_v)_{v}} \&\HH\times \F_p^{\times}\arrow{d}{\upsilon}\\
\& X_1(p,\mu_p,\Gamma_0(\ell)) \arrow{r}{D_{1,1}} \& X_1(p,\mu_p)
\end{tikzcd}
\]

Given $v=[c:d] \in \mathbb{P}^1(\F_{\ell})$, $\pi_v$ is the map $\tau \longmapsto \left(\frac{\C}{\tau\Z\oplus\Z},\frac{1}{p},\langle \frac{c\tau+d}{\ell}\rangle\right)$. Pick, for each $v=[c:d] \in \mathbb{P}^1(\F_p)$, a matrix $M_v \in \Gamma_1(p)$ such that $M_v \equiv \begin{pmatrix} \ast&\ast\\c & d\end{pmatrix} \pmod{\ell}$. Then define $\pi_v,g_v$ as $M_v,M_v^{-1}$ respectively. By Lemma \ref{good-open-subset-hecke}, the above diagram is Cartesian above a certain co-finite open subscheme $U$ of $X_1(p,\mu_p)$. 

Let $D_{\ell} = \begin{pmatrix}\ell & 0\\0 & 1\end{pmatrix}, D'_{\ell} = \begin{pmatrix}1& 0\\0 & \ell\end{pmatrix} \in \mathcal{M}_2(\Z)$. Then, on the component $b^{-1}$ of $\upsilon^{-1}(U)$, one has
\begin{align*}
\upsilon^{\ast}(T_{\ell}^{\ast}\omega) &= \upsilon^{\ast}(D_{1,1}^{\ast})^{\ast}\omega' = \sum_{v \in \mathbb{P}^1(\F_{\ell})}{(g_v^{-1})^{\ast}(\ell f\mid\delta_{\ell})(\ell\tau)d\tau)}\\
&= \sum_{v \in \mathbb{P}^1(\F_{\ell})}{(f \mid \delta_{\ell}D_{\ell}M_v)(\tau)d\tau}=\sum_{M \in (\Gamma_0(\ell) \cap \Gamma_1(p)) \backslash \Gamma_1(p)}{(f \mid \delta_{\ell} D_{\ell})(\tau)d\tau}
\end{align*}

Here, we always consider the weight two action of $\GL{\R}^+$, that is, \[\left(f \mid\begin{pmatrix}a & b\\c & d\end{pmatrix}\right)(\tau) = \frac{ad-bc}{(c\tau+d)^2}f\left(\frac{a\tau+b}{c\tau+d}\right),\] so that scalar matrices act trivially. Now, let us show that \[M \in (\Gamma_0(\ell) \cap \Gamma_1(p))\backslash \Gamma_1(p) \longmapsto \delta_{\ell}D_{\ell}M \in \Gamma_1(p)\backslash \Gamma_1(p)D'_{\ell}\Gamma_1(p)\] is a bijection. 

Let $M \in \Gamma_1(p)$ and $M'=\delta_{\ell}D_{\ell}M$. Then, by \cite[Lemma 3.29]{shimura-hecke}, $M' \in \Gamma_1(p)D'_{\ell}\Gamma_1(p)$ as long as the following conditions are satisfied: the reductions in $\mathcal{M}_2(\F_p)$ of $D'_{\ell}\Gamma_1(p)$ and $\Gamma_1(p)D'_{\ell}$ are equal and contain $M' \pmod{p}$, and $M' \in \Gamma(1)D'_{\ell}\Gamma(1)$. 

Since $\Gamma(1)D'_{\ell}\Gamma(1)$ is the set of matrices with integral entries and determinant $\ell$, it contains $M'$. One easily checks that the reduction mod $p$ of $D'_{\ell}\Gamma_1(p)$ and that of $\Gamma_1(p)D'_{\ell}$ is exactly $\begin{pmatrix}1 & \ast\\0 & \ell\end{pmatrix}$, and that $M' \equiv \Delta_{1,\ell}M \equiv \begin{pmatrix}1 & \ast\\0 & \ell\end{pmatrix}\pmod{p}$. Thus $M' \in \Gamma_1(p)D'_{\ell}\Gamma_1(p)$.

Now, let $M,M' \in \Gamma_1(p)$. Then 
\begin{align*}
\delta_{\ell}D_{\ell}M \in \Gamma_1(p)\delta_{\ell}D_{\ell}M' &\Leftrightarrow M(M')^{-1} \in D_{\ell}^{-1}\delta_{\ell}^{-1}\Gamma_1(p)\delta_{\ell}D_{\ell}\\
&\Leftrightarrow M(M')^{-1} \in D_{\ell}^{-1}\Gamma_1(p)D_{\ell} \cap \Gamma_1(p)\\
&\Leftrightarrow M(M')^{-1} \in \Gamma_1(p) \cap \Gamma_0(\ell)\Leftrightarrow M \in (\Gamma_1(p) \cap \Gamma_0(\ell))M',
\end{align*}
and we are done. 

Now, $\upsilon^{\ast}(T_{\ell}^{\ast}\omega)$ and $(T_{\ell}f)(\tau)d\tau$ are holomorphic differentials on $\HH \times \F_p^{\times}$ and they agree except at points contained in finitely many (discrete) $SL_2(\Z)$-orbits. Hence they are agree, whence the conclusion. }

\cor[diagonalize-x1pmu-T]{Given a normalized newform $f \in \mathcal{S}_2(\Gamma_1(p))$ with character $\chi$, $\psi \in \mathcal{D}$ and $M \in \{T^t\mid t\in \F_p\} \cup \{W\}$, let 
\[u_{f,\psi,M} = \iota_{1,p}^{-1}\left(f \otimes \left(\sum_{t \in \F_p^{\times}}{\psi^{-1}(t)[t]}\right) \otimes M\right) \in H^0(X_1(p,\mu_p)^S_{\C},\Omega^1).\] 

Then the $u_{f,\psi,M}$ form a basis of $H^0(X_1(p,\mu_p)^S_{\C},\Omega^1)$. Moreover, every $u_{f,\psi,M}$ is an eigenvector for the action of $\mathbb{T}$. For a prime $\ell \neq p$, and $n \in \F_p^{\times}$, one has 
\[T_{\ell}^{\ast}u_{f,\psi,M} = a_{\ell}(f)f,\quad \langle n\rangle^{\ast}u_{f,\psi,M} = \chi(n)f,\quad \underline{n}^{\ast}u_{f,\psi,M} = \psi(n)f.\]
}

\demo{By Proposition \ref{hecke-x1pmu}, it is a direct verification that \[T_{\ell}^{\ast}u_{f,\psi,M}=a_{\ell}(f)f, \underline{n}^{\ast}u_{f,\psi,M}=\psi(n), \langle n\rangle^{\ast}u_{f,\psi,M}=\chi(n)u_{f,\psi,M}.\] We show that the $\iota_{1,p}(u_{f,\psi,M})$ form a basis of $(\mathcal{S}_2(\Gamma_1(p)) \otimes \C[\F_p^{\times}]) \otimes_{\Z[B]} \Z[\GL{\F_p}]$, which implies the conclusion. It is enough to prove that the $f \otimes \sum_{t \in \F_p^{\times}}{\psi^{-1}(t)[t]}$, where $f \in \mathcal{S}_2(\Gamma_1(p))$ runs through normalized newforms and $\psi$ through elements of $\mathcal{D}$, form a basis of $\mathcal{S}_2(\Gamma_1(p)) \otimes \C[\F_p^{\times}]$. 

The normalized newforms of $\mathcal{S}_2(\Gamma_1(p))$ form a basis by \cite[Theorem 5.8.2]{DS}, since the only cusp form of weight $2$ and level $1$ is null. Hence, we only need to check that the $e_{\psi}=\sum_{t \in \F_p^{\times}}{\psi^{-1}(t)[t]}$ form a basis of $\C[\F_p^{\times}]$. There are $p-1=\dim_{\C}{\C[\F_p^{\times}]}$ of them, so we only need to prove that the $e_{\psi}$ are linearly independent. 

Let $g \in \F_p^{\times}$ be a generator, then $\underline{g}(e_{\psi})=\psi(g)e_{\psi}$ is an eigenvector for the multiplication by $g$ with eigenvalue $\psi(g)$. Since $\psi$ is a generator, $\psi(g)$ determines $\psi$, hence the $e_{\psi}$ are linearly independent. 
}

\lem[gamma1p-onlytwists]{Let $f \in \mathcal{S}_2(\Gamma_1(p))$ be a normalized newform with character $\chi$. Then the twists of $f$ by primitive Dirichlet characters are pairwise distinct, and the only primitive Dirichlet characters $\psi$ such that $f \otimes \psi$ has conductor $p$ are $1,\chi^{-1}$.}

\demo{Let $\psi \neq 1$ be a primitive Dirichlet character such that $f \otimes \psi$ has conductor $p$. By Corollary \ref{level-newform-bigtwist}, $\psi$ must have conductor $p$, and $\psi\chi,\psi$ must have distinct conductors. This implies that $\chi \neq 1$ and $\psi=\chi^{-1}$. 
If furthermore $f \otimes \psi=f$, then by comparing characters we see that $\chi=\chi^{-1}$, so that $\chi$ is a quadratic character. Since $\chi(-1)=1$, $\chi$ is the character of a real quadratic field and $f$ has real multiplication, which contradicts \cite[Theorem 4.5]{Antwerp5-Ribet}. }

\defi[eigenspaces-x1pmu-Tprime]{Given a normalized newform $f \in \mathcal{S}_2(\Gamma_1(p))$ with character $\chi$ and $\psi\in \mathcal{D}$, $H^0(X_1(p,\mu_p)_{\C}^S,\Omega^1)[f,\psi]$ denotes the subspace of $H^0(X_1(p,\mu_p)_{\C}^S,\Omega^1)$ on which every $T'_{\ell}$ acts by $a_{\ell}(f)\psi(\ell)$, and every $\langle n\rangle'$ acts by $\chi(n)\psi(n)^2$ (where $\ell \neq p$ is prime and $n \in \F_p^{\times}$).  }

\lem[eigenspaces-x1mpu-Tprimen]{With the same notation, for any integer $n \geq 1$ prime to $p$, $T'_n-a_n(f)\psi(n)$ annihilates $H^0(X_1(p,\mu_p)_{\C}^S,\Omega^1)[f,\psi]$. }

\demo{This is a direct induction using the classical relations between the Fourier coefficients of a newform (see for instance \cite[Proposition 5.8.5]{DS}). }

\prop[diagonalize-x1pmu-Tprime]{Each $H^0(X_1(p,\mu_p)^S_{\C},\Omega^1)[f,\psi]$ as in Definition \ref{eigenspaces-x1pmu-Tprime} is stable under $\mathbb{T}$, $w_p^{\ast}$, $R^{\ast}$ and $\GL{\F_p}$. 
Moreover, $H^0(X_1(p,\mu_p)^S_{\C},\Omega^1)$ is the direct sum of the $H^0(X_1(p,\mu_p)^S_{\C},\Omega^1)[f,\psi]$, where we sum over the couples $(f,\psi)$ as in Definition \ref{eigenspaces-x1pmu-Tprime}, and only count exactly one couple among each pair $\{(f,\psi),(\overline{f},\chi\psi)\}$. 

Finally, we have the following isomorphisms of right $\C[\GL{\F_p}]$-modules. 
If $\chi=1$,
\[[I_2] \in \pi^{\vee}(\psi,\psi) \longmapsto \iota_{1,p}^{-1}\left(f \otimes \left(\sum_{t \in \F_p^{\times}}{\psi^{-1}(t)[t]}\right) \otimes [I_2]\right)\in H^0(X_1(p,\mu_p)^S_{\C},\Omega^1)[f,\psi]\]
If $\chi \neq 1$,
\begin{align*}
\pi^{\vee}(\psi,\chi\psi) \oplus \pi^{\vee}(\chi\psi,\psi) &\rightarrow H^0(X_1(p,\mu_p)^S_{\C},\Omega^1)[f,\psi],\\
([I_2],0) & \longmapsto \iota_{1,p}^{-1}\left(f \otimes \left(\sum_{t \in \F_p^{\times}}{\psi^{-1}(t)[t]}\right) \otimes [I_2]\right)\\
(0,[I_2]) & \longmapsto \iota_{1,p}^{-1}\left(\overline{f} \otimes \left(\sum_{t \in \F_p^{\times}}{(\chi\psi)^{-1}(t)[t]}\right) \otimes [I_2]\right).
\end{align*}
}

\demo{The stability properties come from the fact that $w_p,R$ and the actions of $\mathbb{T}$ and $\GL{\F_p}$ commute with $\mathbb{T}'$. Moreover, the $u_{f,\psi,M}$ from Corollary \ref{diagonalize-x1pmu-T} form a basis of $H^0(X_1(p,\mu_p)^S_{\C},\Omega^1)$ made with eigenvectors for the action of $\mathbb{T}' \subset \mathbb{T}$. So, as a $\C[\mathbb{T}']$-module, $H^0(X_1(p,\mu_p)^S_{\C},\Omega^1)$ is the direct sum of its eigenspaces. 

By Corollary \ref{diagonalize-x1pmu-T}, the eigenspaces are always of the following form: the eigenvalue of $T'_{\ell}$ is $a_{\ell}(f)\psi(\ell)$, and that of $\langle n\rangle'$ is $\chi(n)\psi(n)^2$, where $f \in \mathcal{S}_2(\Gamma_1(p))$ is a normalized newform of character $\chi$ and $\psi \in \mathcal{D}$. 

Let $f,f' \in \mathcal{S}_2(\Gamma_1(p))$ be normalized newforms with characters $\chi,\chi'$, and $\psi,\psi' \in \mathcal{D}$, such that the eigenspaces of $H^0(X_1(p,\mu_p)_S^{\C},\Omega^1)$ attached to $(f,\psi)$ and $(f',\psi')$ are the same. Then for every prime $\ell \neq p$, one has $a_{\ell}(f)\psi(\ell)=a_{\ell}(f')\psi'(\ell)$ and for every $n \in \F_p^{\times}$, one has $\chi(n)\psi(n)^2=\chi'(n)\psi'(n)^2$. By strong multiplicity one \cite[Theorem 4.6.19]{Miyake}, the newforms $f \otimes\psi$ and $f' \otimes \psi'$ (see Proposition \ref{twist-exists}) are thus equal, so $f'$ must be the newform $f \otimes \overline{\psi'}\psi$. By Lemma \ref{gamma1p-onlytwists}, this implies that $\psi'=\psi$ (thus $(f,\psi)=(f',\psi')$) or $\psi=\overline{\chi}\psi'$ (so that $(f',\psi')=(\overline{f},\chi\psi)$). 

Conversely, if $(f',\psi') \in \{(f,\psi),(\overline{f},\chi\psi)\}$ where $f \in \mathcal{S}_2(\Gamma_1(p))$ is a normalized newform and $\psi \in \mathcal{D}$, then \[H^0(X_1(p,\mu_p)_S^{\C},\Omega^1)[f,\psi]=H^0(X_1(p,\mu_p)_S^{\C},\Omega^1)[f',\psi'].\]  

Therefore, \[H^0(X_1(p,\mu_p)^S_{\C},\Omega^1) = \bigoplus_{(f,\psi)}{H^0(X_1(p,\mu_p)^S_{\C},\Omega^1)[f,\psi]},\] where $f$ runs through normalized newforms in $\mathcal{S}_2(\Gamma_1(p))$, $\psi$ runs through elements of $\mathcal{D}$ and we count exactly one couple among each pair $\{(f,\psi), (\overline{f},\chi\psi)\}$ when $f$ has character $\chi \neq 1$.

Let $(f,\psi)$ be as in Definition \ref{eigenspaces-x1pmu-Tprime}, and let $\chi$ be the character of $f$. By what we proved above and Corollary \ref{diagonalize-x1pmu-T}, a basis of $H^0(X_1(p,\mu_p)^S_{\C},\Omega^1)[f,\psi]$ is given by the $u_{f,\psi,M}$ if $\chi=1$, and the $u_{f,\psi,M},u_{\overline{f},\chi\psi,M}$ otherwise, where $M$ runs through $\{T^t,\,t \in \F_p\}\cup \{W\}$. Thus, when $\chi=1$, $\iota_{1,p}(H^0(X_1(p,\mu_p)^S_{\C},\Omega^1)) = \C f \otimes \sum_{t \in \F_p^{\times}}{\psi^{-1}(t)[t]} \otimes_{\C[B]} \C[\GL{\F_p}]$, where $\begin{pmatrix}a & \ast\\0 & b\end{pmatrix}$ acts (on the right) on $f \otimes \sum_{t \in \F_p^{\times}}{\psi^{-1}(t)[t]}$ by multiplication by $\psi(ab)=\psi(a)\psi(b)$. In this case, 

\[[I_2] \in \pi^{\vee}(\psi,\psi) \longmapsto \iota_{1,p}^{-1}\left(f \otimes \left(\sum_{t \in \F_p^{\times}}{\psi^{-1}(t)[t]}\right) \otimes [I_2]\right)\in H^0(X_1(p,\mu_p)^S_{\C},\Omega^1)[f,\psi]\] is therefore an isomorphism. 

When $\chi \neq 1$, $\iota_{1,p}(H^0(X_1(p,\mu_p)^S_{\C},\Omega^1))[f,\psi]$ is equal to  
\[\left(\C \left[f \otimes \sum_{t \in \F_p^{\times}}{\psi^{-1}(t)[t]}\right] \oplus \C\left[\overline{f} \otimes \sum_{t \in \F_p^{\times}}{(\chi\psi)^{-1}(t)[t]}\right]\right) \otimes_{\C[B]} \C[\GL{\F_p}],\]%
and, given $M=\begin{pmatrix}a & \ast\\0 & b\end{pmatrix}$, one has 
\begin{align*}
\left(f\otimes \sum_{t \in \F_p^{\times}}{\psi^{-1}(t)[t]}\right) \mid M &= \psi(a)(\chi\psi)(b)\left(f\otimes \sum_{t \in \F_p^{\times}}{\psi^{-1}(t)[t]}\right),\\
\left(\overline{f}\otimes \sum_{t \in \F_p^{\times}}{(\psi\chi)^{-1}(t)[t]}\right) \mid M &= (\chi\psi)(a)\psi(b)\left(\overline{f}\otimes \sum_{t \in \F_p^{\times}}{(\psi\chi)^{-1}(t)[t]}\right),
\end{align*}
whence the conclusion. 
}

\lem[diagonalize-x1pmu-Tprime-R]{Let $(f,\psi)$ be as in Definition \ref{eigenspaces-x1pmu-Tprime}. Then $H^0(X_1(p,\mu_p)^S_{\C},\Omega^1)[f,\psi]$ is not contained in the kernel of $R^{\ast}$. }

\demo{Let $\omega'=u_{f,\psi,I_2} \in H^0(X_1(p,\mu_p)^S_{\C},\Omega^1)[f,\psi] \simeq H^0(P(p)_{\C},\Omega^1)[f,\psi]$, and $\omega$ its (nonzero) restriction to the closed open subscheme $X_1(p,\mu_p)_{(0,1)}$. 
Then, if $\chi$ denotes the character of $f$, for any $t \in \F_p$, 
\[\iota_{(t,1)}^{\ast}R^{\ast}(\iota_{(0,1)})_{\ast}\omega = (\iota_{(0,1)})^{\ast}R[(\iota_{(t,1)})_{\ast}]\omega = (\langle t\rangle w_p)^{\ast}\omega = \chi(t)w_p^{\ast}\omega \neq 0,\]
hence $R^{\ast}\omega' \neq 0$, and the conclusion follows.}

\section{Uniformizing $X(p,p)$}

\prop[uniformize-xpp]{Let $X(p)$ be the classical compact connected Riemann surface whose non-cuspidal points are given by the quotient $\Gamma(p) \backslash \HH$ (see for example \cite[Chapter 2]{DS}). 
There exists a morphism of locally ringed spaces $X(p) \times \F_p^{\times} \rar X(p,p)_{\C}$, making it the analytification of $X(p,p)_{\C}$ in the sense of \cite[Exp. XII, (1.1)]{SGA1}. Under this morphism, a point $(\tau,b) \in \HH \times \F_p^{\times}$ is mapped to the enriched elliptic curve $\left(\frac{\C}{\Z\tau\oplus\Z},\frac{b\tau}{p},\frac{1}{p}\right) \in Y(p,p)(\C)$. 
 
This uniformization induces the following action of $\GL{\F_p}$ on $X(p) \times \F_p^{\times}$: for $M \in \GL{\F_p}$, $(x,b) \in X(p) \times \F_p^{\times}$, let $M'=\Delta_{b\det{M},1}^{-1}M\Delta_{b,1} \in \SL{\F_p}$, then $M \cdot (x,b)=(M' \cdot x, b(\det{M}) x)$.}

\demo{This is a direct consequence of Lemma \ref{uniformize-general}. }

\cor[differentials-xpp]{This identification induces an isomorphism \[\iota_{p,p}: H^0(X(p,p)_{\C},\Omega^1) \rar \bigoplus_{b \in \F_p^{\times}}{\mathcal{S}_2(\Gamma(p)) \otimes \Delta_{b,1}},\] where, given $b \in \F_p^{\times}$ and $f \in \mathcal{S}_2(\Gamma(p))$, $f(\tau)d\tau$ identifies to a holomorphic differential on $X(p)$ (which we will also denote $f(\tau)d\tau$ or $f$), and $f \otimes \Delta_{b,1}$ denotes the differential equal to $f(\tau)d\tau$ on $X(p) \times \{b^{-1}\}$, and vanishing on the other components.

In particular, for any $M \in \GL{\F_p}, b \in \F_p^{\times}, f\in \mathcal{S}_2(\Gamma(p))$, one has \[M^{\ast}(\iota_{p,p}^{-1}(f \otimes \Delta_{b,1})) = \iota_{p,p}^{-1}((f \mid_2 M') \otimes \Delta_{b\det{M},1}),\] where $M'=\Delta_{b,1}M\Delta_{b\det{M},1}^{-1} \in \SL{\F_p}$. 
}

\demo{The analytification of $M^{\ast}\iota_{p,p}^{-1}(f \otimes \Delta_{b,1})$ is the differential on $X(p) \times \F_p^{\times}$ given by $M^{\ast}(f(\tau)\mathbf{1}_{z=b^{-1}}d\tau)$ for $(\tau,z) \in \HH \times \F_p^{\times}$. Then  
\begin{align*}
M^{\ast}(f(\tau)\mathbf{1}_{z=b^{-1}}d\tau)(\tau,z)&=f(\Delta_{z\det{M},1}^{-1}M\Delta_{z,1}\tau)\mathbf{1}_{z\det{M}=b^{-1}}d(\Delta_{z\det{M},1}^{-1}M\Delta_{z,1}\tau) \\
&= (f \mid \Delta_{(b\det{M})^{-1}\det{M},1}^{-1}M\Delta_{(b\det{M})^{-1},1})(\tau)\mathbf{1}_{z=(b\det{M})^{-1}}d\tau\\
& = (f \mid \Delta_{b,1}M\Delta_{b\det{M},1}^{-1})(\tau)\mathbf{1}_{z=(b\det{M})^{-1}}d\tau,
\end{align*}
whence the conclusion.}

~\\

Now, we prove for later reference (this will be needed in Chapter \ref{formal-immersion-application}) that the algebraic $q$-expansion defined in Chapter \ref{cuspidal-subscheme} corresponds to the usual $q$-expansion after analytification. We first need the following lemma. 

\lem[riemann-surface-to-algebraic-factors-through-affine]{Let $Y$ be a smooth proper $\C$-scheme, pure of dimension one, and $X$ be a connected noncompact Riemann surface. Let $\varphi: X \rar Y$ be a morphism of locally ringed spaces. Then $\varphi$ factors through $\Sp{\OO(X)}$. }

\demo{We may assume that $Y$ is connected.

First, we show the following: for any morphism of locally ringed spaces $f: X \rar Y$ from a connected Riemann surface to a smooth connected $\C$-scheme, for any open subscheme $U \subset Y$, the image of $f^{\sharp}: \OO(U) \rar \OO(f^{-1}(U))$ is made with meromorphic functions on $X$. Indeed, when $U,Y,f$ are fixed, the statement is local with respect to $X$: since $f^{-1}(Y \backslash U)$ is discrete, we may assume that $f^{-1}(Y \backslash U)$ is a point $x \in X$. Let $\alpha \in \OO(U)$. There is an affine open subscheme $V \subset Y$ containing $f(x)$ and a nonzero function $\beta \in \OO(V)$ such that $\alpha\beta \in \OO(V \cap U)$ extends to a function $\gamma$ on the open subscheme $(V \cap U) \cup \{f(x)\}$. Then $f^{\sharp}(\alpha)=\frac{f^{\sharp}(\gamma)}{f^{\sharp}(\beta)}$ is a quotient of holomorphic functions near $x$, whence the conclusion.  

Let us come back to the original statement. Let $x,y \in Y$ be two distinct closed points, then $U=Y \backslash \{x\}, V = Y \backslash \{y\}$ are smooth connected affine $\C$-schemes of dimension one by \cite[Lemmas 0A24, 0A28]{Stacks}, and $X=U \cup V$. Let $A \subset \OO(X \backslash \varphi^{-1}(x)), B \subset \OO(X \backslash \varphi^{-1}(y))$ be the images of $\OO(U)$ and $\OO(V)$ respectively: $A, B$ are $\C$-algebras of meromorphic functions on $X$, so, by Mittag-Leffler's theorem \cite[(26.5)]{Forster}, there are functions $\alpha,\beta \in \OO(X)$ vanishing exactly on $f^{-1}(x)$ and $f^{-1}(y)$ respectively such that $A \subset \OO(X)_{\alpha}, B \subset \OO(X)_{\beta}$. Since $\alpha,\beta$ have no common zero, they generate the unit ideal of $\OO(X)$ (see \cite[Exercise 26.4]{Forster}). 

Hence, $\varphi: X \backslash f^{-1}(x) \rar U$ factors through $\Sp{\OO(X)_{\alpha}}$ and $\varphi: X \backslash f^{-1}(y) \rar V$ factors through $\Sp{\OO(X)_{\beta}}$, and the two morphisms $\Sp{\OO(X)_{\alpha}} \rar U$, $\Sp{\OO(X)_{\beta}} \rar V$ agree when restricted to $\Sp{\OO(X)_{\alpha\beta}}$ (because they both come from $\varphi$). Hence they glue to a morphism $g: \Sp{\OO(X)} \rar Y$ such that the composition $X \overset{\mrm{can}}{\rar} \Sp{\OO(X)} \overset{g}{\rar} Y$ is equal to $\varphi$ above $U$ and above $V$, hence this composition is $\varphi$. }

\prop[q-expansion-vs-analytic]{Let $\mu: z \otimes \zeta_p \in \C \otimes \Z[\zeta_p] \longmapsto ze^{\frac{-2i\pi}{p}} \in \C$. 

Let $\omega \in H^0(X(p,p)_{\C},\Omega^1)$ and let $\sum_{n \geq 0}{a_nq^{n/p}} \in \C[[q^{1/p}]]$ be the image under $\mu$ of its $q$-expansion at the cusp datum $(I_2,(x,y) \mapsto x)$ (see Definition \ref{qexp-definition}). Then the $\Delta_{1,1}$ component of $\iota_{p,p}(\omega)$ is $\frac{2i\pi}{p}\sum_{n \geq 0}{a_ne^{\frac{2i\pi (n+1)\tau}{p}}}$. }

\demo{\emph{Step 1: Descending the uniformization.}

Let $\mathcal{M}(\HH/p\Z)$ denote the ring of $p$-periodic holomorphic functions $f$ on $\HH$ that are meromorphic at infinity: in other words, if we write $f(\tau)=g(e^{2i\pi\tau})$, where $g: D(0,1)\backslash \{0\} \rar \C$ is a holomorphic function, then the singularity of $g$ at $0$ is either removable or a pole. 

Consider the following diagram:
\[
\begin{tikzcd}[ampersand replacement=\&]
\Sp{\OO(\HH \times \F_p^{\times})} \arrow{r}\arrow[bend left=12,near start]{rr}{\mrm{Vers}} \& \Sp{\mathcal{M}(\HH/p\Z) \otimes \C^{\F_p^{\times}}} \arrow[swap]{r}{\mrm{V}}\& Y(p,p)_{\C} \arrow{r} \& X(p,p)_{\C}\\	
\Sp{\C((q^{1/p}))} \arrow[swap]{ur}{Q} \arrow{r}{\mu} \arrow{d}\& \Sp{\C \otimes \Z[\zeta_p] \otimes \Z((q^{1/p}))} \arrow[swap,bend right=10]{ur}{C_{(I_2,(x,y) \mapsto x)}^0}\arrow{d} \& \&\\
\Sp{\C[[q^{1/p}]]} \arrow{r}{\mu}\& \Sp{\C \otimes \Z[\zeta_p] \otimes \Z[[q^{1/p}]]} \arrow[swap, bend right=20]{rruu}{C_{(I_2,(x,y) \mapsto x)}} \& \&
\end{tikzcd}
\]
In this diagram, the unmarked vertical and horizontal arrows are the natural ones, and the maps $C_{(I_2,(x,y) \mapsto x)}, C^0_{(I_2,(x,y) \mapsto x)}$ are defined in Proposition \ref{q-expansion-setup} and its proof. For the sake of brevity, since we will only consider the cusp datum $(I_2,(x,y)\mapsto x)$, we write $C, C^0$ instead of $C_{(I_2,(x,y) \mapsto x)}, C^0_{(I_2,(x,y) \mapsto x)}$. So we need to define $\mrm{Vers}, V, Q$, and check that they make the diagram commute. 

The application $\mrm{Vers}$ denotes the versal elliptic curve as defined in the proof of Lemma \ref{uniformize-general}: to be precise, it is given by the equation $y^2=x^3-15G_4 \cdot x - 35G_6$ with $G_4, G_6 \in \OO(\HH \times \F_p^{\times})$ do not depend on the $\F_p^{\times}$ coordinate, and the level structure is given by 
\[(\tau,b) \in \HH \times \F_p^{\times}\longmapsto \left(\left[\wp_{\tau\Z\oplus \Z}\left(\frac{b\tau}{p}\right):2\wp_{\tau\Z\oplus\Z}'\left(\frac{b\tau}{p}\right):1\right],\left[\wp_{\tau\Z\oplus \Z}\left(\frac{1}{p}\right):2\wp_{\tau\Z\oplus\Z}'\left(\frac{1}{p}\right):1\right]\right).\]

Now, $G_4, G_6 \in \mathcal{M}(\HH/p\Z)$, since they are modular forms for $SL_2(\Z)$. Moreover, let $c,d \in \Z$ be integers not both divisible by $p$. Then $\tau \longmapsto \wp_{\tau\Z\oplus\Z}'\left(\frac{c\tau+d}{p}\right)$ is a $p$-periodic holomorphic function on $\HH$, given as the sum of the series
\[-2\sum_{(u,v) \in \Z^2}{\left(\frac{c\tau+d}{p}-(u\tau+v)\right)^{-3}},\]
which converges uniformly on rectangles of the form $\{|\mrm{Re}(\tau)| \leq D, \mrm{Im}(\tau) \geq t\}$ for any $D,t > 0$. This implies that $\wp_{\tau\Z\oplus\Z}'\left(\frac{c\tau+d}{p}\right)$ converges to a constant as $\mrm{Im}(\tau) \rar \infty$, hence \[\left[\tau \longmapsto \wp_{\tau\Z\oplus\Z}'\left(\frac{c\tau+d}{p}\right) \right]\in \mathcal{M}(\HH/p\Z).\] The calculation at the start of \cite[Section 4.6]{DS} shows that $\left[\tau \longmapsto \wp_{\tau\Z\oplus\Z}\left(\frac{c\tau+d}{p}\right) \right]\in \mathcal{M}(\HH/p\Z)$, so there is a unique morphism $V: \Sp{\mathcal{M}(\HH/p\Z) \otimes \C^{\F_p^{\times}}} \rar Y(p,p)_{\C}$ making the second row of the diagram commute. 

For any $f \in \mathcal{M}(\HH/p\Z)\otimes \C^{\F_p^{\times}}$, there is a unique power series $g(q^{1/p}):= Q^{\sharp}(f) \in \C((q^{1/p}))$ such that for all $\tau \in \HH$, one has $f(\tau,1)=g(e^{\frac{2i\pi \tau}{p}})$. The morphism $Q$ in the diagram is the morphism of schemes coming from $Q^{\sharp}$.

To check that the diagram commutes, we need to write down an isomorphism preserving the level structure between the Tate curve, given over $\C((q^{1/p}))$ by 
\[y^2+xy = x^3-5x\frac{E_4-1}{240}+\frac{1}{12}\left(-5\frac{E_4-1}{240}-7\frac{E_6-1}{-504}\right),\] and the curve with Weierstrass equation \[v^2 = u^3-15G_4u-35G_6 \Leftrightarrow v^2 = u^3-\frac{\pi^4}{3}E_4u-\frac{2\pi^6}{27}E_6.\]

A direct computation with \cite[(A1.2.1)]{Katz} shows that the appropriate change of variables is  
\[u = (2i\pi)^2\left(x+\frac{1}{12}\right),\, v = \frac{(2i\pi)^3}{2}(x+2y).\]

Now, let us discuss the level structures. The level structure induced on the Tate curve by $C_{(I_2,(x,y)\mapsto x)}$ is an isomorphism $\phi: (\Z/p\Z)^2 \longmapsto q^{\left(\frac{1}{p}\Z\right)/\Z}\zeta_p^{\Z/p\Z}$ such that $\phi((0,-1))=\zeta_p$ and $\phi((1,0))=q^{1/p}$. In particular, the image of this level structure under $\mu$ is given by $(a,b) \longmapsto e^{\frac{2i\pi}{p}(a\tau+b)} \in \C^{\ast}/e^{2i\pi\tau\Z}$, which is isomorphic to $(a,b) \longmapsto \frac{a\tau+b}{p} \in \frac{\C}{\tau\Z\oplus\Z}$. 

This proves that the diagram that we drew at the beginning of the proof commutes. \\

\emph{Step 2: Extension of the diagram to infinity.}

Let $\HH_p^{\ast}$ be the Riemann surface $\HH/p\Z \cup \{i\infty\}$, where the holomorphic structure near $i\infty$ makes the bijection $\tau \in \HH_p^{\ast} \longmapsto e^{2i\pi \frac{z}{p}} \in D(0,1)$ a biholomorphism. Let $\HH_p^{\ast} \rar X(p,p)^{an}$ be the morphism $\tau \longmapsto (\tau,1)$ with the uniformization of Proposition \ref{uniformize-xpp} (see \cite[Chapter 2.4]{DS} for the description of the cusps of the analytification). By Lemma \ref{riemann-surface-to-algebraic-factors-through-affine}, the composition $\HH_p^{\ast} \rar X(p,p)^{an} \rar X(p,p)_{\C}$ factors through a morphism $V_{1,\infty}: \Sp{\OO(\HH_p^{\ast})} \rar X(p,p)_{\C}$. 

Let $\mrm{ev}_1$ be the morphism of schemes attached to the homomorphism of $\C$-algebras \[\mathcal{M}(\HH/p\Z) \otimes \C^{\F_p^{\times}} \simeq \mrm{Fun}(\F_p^{\times},\mathcal{M}(\HH/p\Z)) \overset{f \mapsto f(1)}{\longrightarrow} \mathcal{M}(\HH/p\Z).\] By Step 1, the following diagram commutes   
\[
\begin{tikzcd}[ampersand replacement=\&]
\Sp{\mathcal{M}(\HH/p\Z)} \arrow{dr}{\mrm{ev}_{1}} \arrow{r}\& \Sp{\OO(\HH_p^{\ast})} \arrow{dr}{V_{1,\infty}} \& \\
 \& \Sp{\mathcal{M}(\HH/p\Z) \otimes \C^{\F_p^{\times}}} \arrow{r}{V} \& X(p,p)_{\C}\\
\Sp{\C((q^{1/p}))} \arrow{uu}{Q_1} \arrow{ru}{Q} \arrow{rru}[swap]{C^0\circ \mu} \arrow{rr}\& \& \Sp{\C[[q^{1/p}]]} \arrow{u}{C \circ \mu}
\end{tikzcd}
\]
The morphism $Q_{1,\infty}: \Sp{\C((q^{1/p}))} \rar \Sp{\OO(\HH_p^{\ast})}$ actually factors through a morphism $Q_{\infty}: \Sp{\C[[q^{1/p}]]} \rar \Sp{\OO(\HH_p^{\ast})}$ corresponding to the usual $q$-expansion, mapping the function $\sum_{n \geq 0}{a_ne^{\frac{2i\pi n\tau}{p}}}$ to the series $\sum_{n \geq 0}{a_nq^{n/p}}$. 

The morphisms $V_{1,\infty} \circ Q_{\infty}, C \circ \mu: \Sp{\C[[q^{1/p}]]} \rar X(p,p)_{\C}$ are morphisms of reduced separated $\C$-schemes and agree on the dense open subscheme $\Sp{\C((q^{1/p}))}$ of the source, hence they are equal. In other words, the following diagram commutes:

\[\begin{tikzcd}[ampersand replacement=\&]
\HH_p^{\ast} \arrow{r}{\tau \mapsto (\tau,1)} \arrow{d}\&  X(p,p)^{an} \arrow{d}\\
\Sp{\OO(\HH_p^{\ast})} \arrow{r}{V_{1,\infty}}\& X(p,p)_{\C}\\
\Sp{\C[[q^{1/p}]]} \arrow{u}{Q_{\infty}} \arrow{ru}{C \circ \mu}\&
\end{tikzcd}\]

\emph{Step 3: Conclusion.}

Let $x \in X(p,p)(\C)$ be the image of the maximal ideal of $\Sp{\C[[q^{1/p}]]}$ under $C \circ \mu$: its inverse image $x^{an} \in X(p,p)_{\C}^{an}$ is the point $(i\infty,1)$.     

For $z \in \HH$ and $\gamma \in SL_2(\Z)$, if $z$ and $\gamma \cdot z$ have imaginary part greater than one, then $\gamma = \pm \begin{pmatrix}1 & \ast\\0 & 1\end{pmatrix}$. Therefore, the function $e^{\frac{2i\pi\tau}{p}}$ is well-defined as an element of $\OO_{X(p,p)_{\C}^{an},(i\infty,1)}$ and is a uniformizer of this discrete valuation ring by the construction of the Riemann surface attached to $\Gamma(p)$ \cite[Chapter 2.4]{DS}. By Step 2, the following diagram commutes.
\[
\begin{tikzcd}[ampersand replacement=\&]
\OO_{X(p,p)_{\C},x} \arrow{r} \arrow[bend right=20]{rr}{(C \circ \mu)^{\sharp}} \& \OO_{X(p,p)^{an},x^{an}} \arrow{r}{e^{\frac{2i\pi}{p}} \mapsto q^{1/p}} \& \C[[q^{1/p}]]
\end{tikzcd}
\]

Moreover, by \cite[Exp. XII, Th\'eor\`eme 1.1]{SGA1} and Proposition \ref{qexp-vs-formal-immersion}, all three rings are discrete valuation rings (since $X(p,p)_{\C}$ is smooth of relative dimension one over $\C$ and its analytification is a Riemann surface), all three morphisms are local and induce isomorphisms after completing.     

Let $\sum_{n \geq 1}{b_ne^{\frac{2i\pi n\tau}{p}}}$ denote the component of $\Delta_{1,1}$ in $\iota_{p,p}(\omega)$. Since $d\tau = \frac{p}{2i\pi}e^{\frac{-2i\pi}{p}}d(e^{\frac{2i\pi}{p}})$, the pull-back of $\omega$ as a holomorphic differential form on $X(p,p)^{an}_{\C}$ can be written locally around $x^{an}$ as $\frac{p}{2i\pi}\left(\sum_{n \geq 0}{b_{n+1}e^{\frac{2i\pi n\tau}{p}}}\right) \cdot d(e^{\frac{2i\pi\tau}{p}})$. The local diagram commutes, so the pull-back of $\omega$ by $C \circ \mu$ as a continuous differential on the analytic space $\Sp{\C[[q^{1/p}]]/(q^t)}$ is given by $\frac{p}{2i\pi}\left(\sum_{n = 0}^{pt}{b_{n+1}q^{n/p}}\right)d(q^{1/p})$, and we proved the result. }

\medskip

\cor[full-analytic-qexp-at-infty]{Let $\omega \in H^0(X(p,p)_{\C},\Omega^1)$ and let $\sum_{n \geq 0}{a_nq^{n/p}} \in (\Z[\zeta_p] \otimes\C)[[q^{1/p}]]$ be its $q$-expansion at the cusp datum $(I_2,(x,y) \mapsto x)$. Let, for every $a \in \F_p^{\times}$, $\mu_a: \Z[\zeta_p] \otimes \C \rar \C$ be the homomorphism of $\C$-algebras such that $\mu_a(\zeta_p)= e^{-\frac{2i\pi a}{p}}$. Then  
\[\iota_{p,p}(\omega) = \frac{2i\pi}{p}\sum_{a \in \F_p^{\times}}{\left(\sum_{n \geq 1}{\mu_a(a_{n-1})e^{\frac{2i\pi n\tau}{p}}}\right) \otimes \Delta_{a^{-1},1}}.\]
}

\demo{For any $t \geq 2$, let $C_t$ be the following composition: \[\Sp{(\Z[\zeta_p] \otimes \C) \otimes \Z[[q^{1/p}]]/(q^{t/p})} \rar \Sp{(\Z[\zeta_p] \otimes \C) \otimes \Z[[q^{1/p}]]} \overset{C_{(I_2,(x,y)\mapsto x)}}{\longrightarrow} X(p,p)_{\C}.\] Then, for any $a \in \F_p^{\times}$, one has by Lemma \ref{with-delta-c}
\begin{align*}
(\Sp{\mu_a})^{\ast}C_t^{\ast}\omega &= (\Sp{\mu_1 \circ \underline{a}})^{\ast}C_t^{\ast}\omega = (\Sp{\mu_1})^{\ast}(\underline{a})^{\ast}C_t^{\ast}\omega\\
&= (\Sp{\mu_1})^{\ast}C_t^{\ast}(\Delta_{a,1}^{\ast}\omega).
\end{align*}
Therefore, the power series in front of the $\Delta_{1,1}$ component of $\iota_{p,p}(\Delta_{a,1}^{\ast}\omega)$ (hence the $\Delta_{a^{-1},1}$ component of $\iota_{p,p}(\omega)$) is $\frac{2i\pi}{p}\sum_{n\geq 1}{\mu_a(a_{n-1})e^{\frac{2i\pi n\tau}{p}}}$. }

\cor[no-lines]{For every $\psi\in \mathcal{D}$, the subspace of $H^0(X(p,p)_{\C},\Omega^1)$ on which $\GL{\F_p}$ acts by $\psi(\det)$ is trivial.
}

\demo{Let $\omega \in H^0(X(p,p)_{\C},\Omega^1)$ be such that for all $g \in \GL{\F_p}$, $g^{\ast}\omega=\psi(\det{g})\omega$. Write $\iota_{p,p}(\omega)=\sum_{b \in \F_p^{\times}}{f_b \otimes \Delta_{b,1}}$. Then, if $M \in \SL{\F_p}$, 
\begin{align*}
\sum_{b \in \F_p^{\times}}{f_b \otimes \Delta_{b,1}} &= \iota_{p,p}(\omega)=\iota_{p,p}(M^{\ast}\omega) = M^{\ast}\sum_{b \in \F_p^{\times}}{f_b \otimes \Delta_{b,1}}\\
& = \sum_{b \in \F_p^{\times}}{(f_b \mid \Delta_{b,1}M\Delta_{b,1}^{-1})\mid \Delta_{b,1}}, 
\end{align*}
so every $f_b$ is in $\mathcal{S}_2(\Gamma(1)) = \{0\}$, so $\omega=0$. }

\medskip

\prop[hecke-xpp]{Let $\ell \neq p$ be a prime number, then the following diagram commutes, where $\ell$ maps $\Delta_{b,1}$ to $\Delta_{b\ell,1}$:
\[
\begin{tikzcd}[ampersand replacement=\&]
H^0(X(p,p)_{\C},\Omega^1) \arrow{r}{\iota_{p,p}} \arrow{d}{T_{\ell}^{\ast}} \& \bigoplus_{b \in \F_p^{\times}}{\mathcal{S}_2(\Gamma(p)) \otimes \Delta_{b,1}} \arrow{d}{[\Gamma(p)D_{1,\ell}\Gamma(p)] \otimes \ell}\\
H^0(X(p,p)_{\C},\Omega^1) \arrow{r}{\iota_{p,p}}  \& \bigoplus_{b \in \F_p^{\times}}{\mathcal{S}_2(\Gamma(p)) \otimes \Delta_{b,1}}
\end{tikzcd}
\]}

\demo{This is proved in a completely similar way to Proposition \ref{hecke-x1pmu}. Consider the following commutative diagram of locally ringed spaces, where $\upsilon$ is always the uniformization defined in Lemma \ref{uniformize-general}: 
\[
\begin{tikzcd}[ampersand replacement=\&]
\HH \times \F_p^{\times} \arrow{d}{\upsilon}\& \HH \times \F_p^{\times} \arrow{l}{h} \arrow{r}{(g_w)_{w}} \arrow{rd}{\upsilon} \& \coprod\limits_{w \in \mathbb{P}^1(\F_{\ell})}{\HH \times \F_p^{\times}} \arrow{d}{\pi} \arrow{r}{h'} \& \HH \times \F_p^{\times}\arrow{d}{\upsilon}\\
Y(p,p)_{\C} \& \& Y(p,p,\Gamma_0(\ell))_{\C} \arrow{ll}{D_{\ell,1}} \arrow{r}{D_{1,1}} \& Y(p,p)_{\C}
\end{tikzcd}
\]
Here, $h(\tau,z)=(\delta_{\ell}(\ell\tau),\ell z)$ where $\delta_{\ell} \in \Gamma(\ell)$ is congruent mod $p$ to $\Delta_{\ell^{-1},\ell}$, $h'(\tau,z)=(\tau,z)$, $\pi$ maps a point $(\tau,z)$ on the component $w=[c:d] \in \mathbb{P}^1(\F_{\ell})$ to the class of the elliptic curve $\left(\frac{\C}{\tau\Z\oplus\Z},\frac{z\tau}{p},\frac{1}{p},\langle\frac{c\tau+d}{\ell}\rangle\right)$, and $g_v: (\tau,z) \longmapsto (M_v^{-1}\tau,z)$ where $M_v \in \Gamma(p)$ is such that the image in $\mathbb{P}^1(\F_{\ell})$ of its second row is $v$. Note that $(\pi_v)_v$ is a morphism of locally ringed spaces because so is $g_v$ and so is $\upsilon$. 

As in Proposition \ref{hecke-x1pmu}, the rightmost $2 \times 2$ square is Cartesian above a cofinite open subscheme of $Y(p,p)$ (by Lemma \ref{good-open-subset-hecke}) and let \[D_{\ell} = \begin{pmatrix}\ell & 0\\0 & 1\end{pmatrix},\,D'_{\ell} = \begin{pmatrix}1 & 0\\0 & \ell\end{pmatrix}.\] As in the proof of Proposition \ref{hecke-x1pmu}, it is enough to show that \[M \in (\Gamma_0(\ell)\cap\Gamma(p)) \backslash \Gamma(p) \longmapsto \delta_{\ell}D_{\ell}M \in \Gamma(p) \backslash \Gamma(p)D'_{\ell}\Gamma(p),\] is bijective, which is proved in a similar way. 
}

\medskip\noindent

\nott{We denote by $\Gamma'$ the congruence group $\Gamma_1(p) \cap \Gamma_0(p^2)$. Let $\mathcal{N}_0$ be the collection of $(f,\chi)$, where $f \in \mathcal{S}_2(\Gamma_1(p))$ is a normalized newform with character $\chi$, and $\mathcal{N}$ be the reunion of $\mathcal{N}_0$ and the $(f,\chi)$, where $f$ is a normalized newform in $\mathcal{S}_2(\Gamma')$ with character $\chi$.  

By Proposition \ref{twist-exists}, for any $\psi \in \mathcal{D}$, and $(f,\chi) \in \mathcal{N}$, then $(f\otimes \psi,\chi\psi^2) \in \mathcal{N}$. }

\medskip\noindent

\prop[diagonalize-xpp-Tbig]{Let $j: \mathcal{S}_2(\Gamma(p)) \rar \mathcal{S}_2(\Gamma')$ be the isomorphism given by $f \longmapsto f\mid \begin{pmatrix}p & 0\\0 & 1\end{pmatrix}$: it extends as an isomorphism \[j': \bigoplus_{b \in \F_p^{\times}}{\mathcal{S}_2(\Gamma(p)) \otimes \Delta_{b,1}} \rar \bigoplus_{b \in \F_p^{\times}}{\mathcal{S}_2(\Gamma')\otimes\Delta_{b,1}}.\] The following collection of vectors is a basis of $H^0(X(p,p)_{\C},\Omega^1)$ made with eigenvectors for $\mathbb{T}$ and the $\Delta_{a,b}$ (for $a,b \in \F_p^{\times}$): 
\begin{itemize}[noitemsep,label=\tiny$\bullet$]
\item $u_{f,\psi} := (j' \circ \iota_{p,p})^{-1}\left(f(\tau) \otimes \sum_{t \in \F_p^{\times}}{\psi^{-1}(t)[\Delta_{t,1}]}\right)$, where $(f,\chi)$ runs through $\mathcal{N}$ and $\psi$ runs through elements of $\mathcal{D}$,
\item $u^p_{f,\psi} := (j' \circ \iota_{p,p})^{-1}\left(f(p\tau) \otimes \sum_{t \in \F_p^{\times}}{\psi^{-1}(t)[\Delta_{t,1}]}\right)$, where $(f,\chi)$ runs through $\mathcal{N}_0$ and $\psi$ runs through elements of $\mathcal{D}$.
\end{itemize}
In either case, the eigenvalue for the Hecke operator $T_{\ell}$, for $\ell \neq p$ prime (resp. $m_n$, resp. $\Delta_{n,1}$ for $n \in \F_p^{\times}$) is $a_n(f)\psi(n)$ (resp. $\psi(n)^2\chi(n)$, resp. $\psi(n)$). }

\demo{For any prime $\ell \neq p$, one can directly check that $j(f \mid [\Gamma(p)\mrm{diag}(1,\ell)\Gamma(p)]) = T_{\ell}(j(f))$. By similar verifications, the induced action by $j' \circ \iota_{p,p}$ of $\mathbb{T}[(\Delta_{a,b})_{a,b} \in \F_p^{\times}]$ on $\bigoplus_{b \in \F_p^{\times}}{\mathcal{S}_2(\Gamma') \otimes \Delta_{b,1}}$ is the following: pulling back by $T_{\ell}$ (resp. $m_n$, resp. $\Delta_{n,1}$) acts by $T_{\ell} \otimes [\ell]$ (resp. $\langle n\rangle \otimes [n^2]$, resp. $\mrm{id}\otimes[n]$) where $[n]$ is the application $f \otimes \Delta_{b,1} \longmapsto f \otimes \Delta_{bn,1}$. 

Now, as we saw in the proof of Corollary \ref{diagonalize-x1pmu-T}, the $\sum_{t \in \F_p^{\times}}{\psi^{-1}(t)[\Delta_{t,1}]}$ (where $\psi$ runs through characters in $\mathcal{D}$) form a basis of $\sum_{t \in \F_p^{\times}}{\C[\Delta_{t,1}]}$, and these vectors are eigenvectors of any $[n]$ ($n \in \F_p^{\times}$) with eigenvalue $\psi(n)$. Moreover, the $f(\tau)$ for $(f,\chi) \in \mathcal{N}$ and the $f(p\tau)$ for $(f,\chi) \in \mathcal{N}_0$ form a basis of $\mathcal{S}_2(\Gamma')$ by \cite[Theorem 5.8.3]{DS}, and they are eigenvectors for the Hecke and diamond operators (this follows for instance from the proof of \cite[Proposition 5.6.2]{DS}). This implies the conclusion.}

\cor[first-decomp-xpp]{The action of $\mathbb{T}$ on $H^0(X(p,p)_{\C},\Omega^1)$ is diagonalizable. Its eigenspaces, which are stable under the right action of $\GL{\F_p}$, are the following form:
\begin{itemize}[noitemsep,label=\tiny$\bullet$]
\item Given $(f,\chi) \in \mathcal{N}_0$ and $\psi \in \mathcal{D}$, the subspace on which $T_{\ell}$ (resp. $m_n$) acts with eigenvalue $a_{\ell}(f)\psi(\ell)$ (resp. $\chi(n)\psi(n)^2$) has dimension $p$ if $\chi=1$ and $p+1$ otherwise. We let $\Omega^1[f,\psi]$ be this eigenspace. 
\item Given $(f,\chi) \in \mathcal{N}$ such that all the twists of $f$ by Dirichlet characters modulo $p$ have conductor $p^2$ (we say that $f$ is \emph{twist-primitive}), the subspace on which $T_{\ell}$ (resp. $m_n$) acts with eigenvalue $a_{\ell}(f)$ (resp. $\chi(n)$) has dimension $p-1$. We let $\Omega^1[f]$ be this eigenspace.
\end{itemize} }

\demo{The action of $\mathbb{T}$ on $H^0(X(p,p)_{\C},\Omega^1)$ is diagonal in the basis of Proposition \ref{diagonalize-xpp-Tbig}. 

First, we need to show that, for any vector $v$ in this basis, there exist $(f,\chi) \in \mathcal{N}_0$ and $\psi \in \mathcal{D}$ such that for every $\ell \neq p$ and every $n \in \F_p^{\times}$, $T_{\ell}^{\ast}v=a_{\ell}(f)\psi(\ell)v$, $m_n^{\ast}v=\chi(n)\psi(n)^2$, or there exists $(f,\chi) \in \mathcal{N}$ such that $f$ is twist-primitive and for every $\ell \neq p$ and every $n \in \F_p^{\times}$, $T_{\ell}^{\ast}v=a_{\ell}(f)v$ and $m_n^{\ast}v=\chi(n)v$. 

First, one always has $u^p_{f,\psi}\in \Omega^1[f,\psi]$. Next, let $(f,\chi) \in \mathcal{N}$ and $\psi \in \mathcal{D}$. If $f$ is twist-primitive, then its twist $g=f\otimes \psi \in \mathcal{S}_2(\Gamma')$ is a twist-primitive newform, and one directly checks that $u_{f,\psi} \in \Omega^1[g]$. If $f$ is not twist-primitive, then $f=g \otimes \theta$ for some character $\theta \in \mathcal{D}$ and $(g,\chi\theta^{-2}) \in \mathcal{N}_0$, so that it is easy to check that $u_{f,\psi} \in \Omega^1[g,\psi\theta]$.

Let $(f,\chi) \in \mathcal{N}$ and $\psi \in \mathcal{D}$, assume that either $f$ is twist-primitive and $\psi=\mathbf{1}$, or $(f,\chi) \in \mathcal{N}_0$. Let us count how many vectors in the basis of Proposition \ref{diagonalize-xpp-Tbig} lie in $\Omega^1[f,\psi]$. If $\psi' \in \mathcal{D}$ and $(g,\chi') \in \mathcal{N}$ (resp. $(g,\chi') \in \mathcal{N}_0$), the vector $u_{g,\psi'}$ (resp. $u^p_{g,\psi'}$) lies in $\Omega^1[f,\psi]$ if and only if, for every prime $\ell \neq p$ and every $n \in \F_p^{\times}$, $a_{\ell}(g)\psi'(\ell)=a_{\ell}(f)\psi(\ell)$ and $\chi'(n)\psi'(n)^2=\chi(n)\psi(n)^2$. When this condition is satisfied, the newforms $g \otimes \psi'$ and $f \otimes \psi$ are equal (by strong multiplicity one \cite[Theorem 4.6.19]{Miyake}), so $(g,\psi')=(f \otimes \psi\overline{\psi'},\psi')$ -- and the converse is true. 

When $(f,\chi) \in \mathcal{N}_0$, every choice of $\psi'$ (resp. by Lemma \ref{gamma1p-onlytwists}, only $\psi'=\psi$ if $\chi=\mathbf{1}$, or only $\psi'\in \{\psi,\chi\psi\}$ if $\chi \neq \mathbf{1}$) determines a suitable couple $(g,\psi')$. Therefore, $\Omega^1[f,\psi]$ has dimension $p=(p-1)+1$ if $\chi=\mathbf{1}$ and $p+1=(p-1)+2$ otherwise. When $f$ is twist-primitive and $\psi=\mathbf{1}$, $g$ always has conductor $p^2$, so $\Omega^1[f]$ contains no vector of the form $u^p_{\cdot,\cdot}$. Since every choice of $\psi'$ determines a couple $(g,\psi')$, such that $u_{g,\psi'} \in \Omega^1[f]$, $\dim{\Omega^1[f]}=p-1$. }

\medskip

\defi{We define the sets $\mathscr{S},\mathscr{P},\mathscr{C}$ as follows:
\begin{itemize}[noitemsep,label=\tiny$\bullet$]
\item $\mathscr{S}$ is the set of $(f,\mathbf{1}) \in \mathcal{N}_0$ (alternatively, $f$ is a newform with trivial character and conductor $p$).
\item $\mathscr{P}$ is the set of $(f,\chi) \in \mathcal{N}_0$ with $\chi$ nontrivial (alternatively, $f$ is a newform with conductor $p$ and nontrivial character). 
\item $\mathscr{C}$ is the set of $(f,\chi) \in \mathcal{N}$ where $f$ is twist-primitive, that is, every twist of $f$ by a Dirichlet character of conductor $p$ has conductor $p^2$.  
\end{itemize} }

\cor[first-decomp-xpp-Tn]{For any $(f,\chi) \in \mathcal{N}_0$ and $\psi \in \mathcal{D}$ (resp. any $(f,\chi) \in \mathcal{N}$ such that $f$ is twist-primitive), for any $n \geq 1$ coprime to $p$, $T_n$ acts on $\Omega^1[f,\psi]$ (resp. $\Omega^1[f]$) by $a_n(f)\psi(n)$ (resp. $a_n(f)$).}

\demo{Exactly like in the proof of Lemma \ref{eigenspaces-x1mpu-Tprimen}, it is a direct induction on $n$ using the definition of $T_n$ and the relations between Fourier coefficients of newforms. }

\medskip

\cor[primitive-is-cuspidal]{Let $(f,\chi) \in \mathscr{C}$. Let $\Omega^1[f]$ be the subspace of $H^0(X(p,p),\Omega^1)$ which $T_{\ell}$ (resp. $m_n$) acts by $a_{\ell}(f)$ (resp. $\chi(n)$). Then $\Omega^1[f]$ is an irreducible cuspidal (right) representation of $\GL{\F_p}$ with central character $\chi$. Let, for each $\psi\in \mathcal{D}$, 
$e_{\psi}=u_{f_{\psi},\psi^{-1}}$ with the notation of Proposition \ref{diagonalize-xpp-Tbig} where $f_{\psi}$ is the twist of $f$ by $\psi$. 

Then the $e_{\psi}$ form a basis of $\Omega^1[f]$. Moreover, the following relations hold (where the notations $\mathfrak{g}$ and $\lambda_p$ are defined in Appendix \ref{appendix-local-const}):
\begin{itemize}[noitemsep,label=\tiny$\bullet$]
\item For $a \in \F_p^{\times}$, $e_{\psi} \mid U^a = \sum_{\theta}{\frac{\theta(a)}{p-1}\mathfrak{g}(\overline{\theta})e_{\psi\theta}}$, where $\theta$ runs through elements of $\mathcal{D}$, 
\item For any $a \in \F_p^{\times}$, $e_{\psi}\mid \Delta_{a,1} = \psi^{-1}(a)e_{\psi}$.
\item $e_{\psi} \mid W = \lambda_p(f_{\psi})e_{\overline{\chi}/\psi}$. 
\end{itemize}
}

\demo{The second relation is a direct consequence from Proposition \ref{diagonalize-xpp-Tbig}. Note that by \cite[Theorem 4.6.17]{Miyake}, for every $\psi \in \mathcal{D}$, $a_p(f_{\psi})=0$.

For any function $\alpha:\F_p^{\times} \rar \C$, one easily computes that \[\alpha = \frac{1}{p-1}\sum_{\chi \in \mathcal{D}}{\left(\sum_{x \in \F_p^{\times}}{\alpha(x)\chi^{-1}(x)}\right)\chi}.\] Thus, for every $m \in \F_p^{\times}$, one has $e^{\frac{2i\pi m}{p}} = \frac{1}{p-1}\sum_{\chi \in \mathcal{D}}{\chi(m)\mathfrak{g}(\overline{\chi})}$ (where $\mathfrak{g}$ is the Gauss sum, every character being regarded as a mod $p$ Dirichlet character). 

Let $\psi \in \mathcal{D}$ and $a \in \F_p^{\times}$. Let $V=\begin{pmatrix}1 & 1/p \\0 & 1\end{pmatrix} \in SL_2(\Q)$. Then one has, in $\bigoplus_{t \in \F_p^{\times}}{\C[[q]] \otimes [\Delta_{t,1}]}$,
 \begin{align*}
 j'(\iota_{p,p}(e_{\psi}\mid U^a)) &= \sum_{x \in \F_p^{\times}}{\psi(x)f_{\psi}\mid_2 V^{xa} \otimes [\Delta_{x,1}]}= \sum_{\substack{n \geq 1\\x \in \F_p^{\times}}}{\psi(x)a_n(f)\psi(n)e^{\frac{2i\pi}{p}nxa}e^{2i\pi n\tau} \otimes [\Delta_{x,1}]}\\
 &= \frac{1}{p-1}\sum_{\substack{n \geq 1\\(x,\theta) \in \F_p^{\times}\times \mathcal{D}}}{\psi(x)a_n(f)\psi(n)\mathfrak{g}(\overline{\theta})\theta(nxa)e^{2i\pi n \tau} \otimes [\Delta_{x,1}]}\\
 &= \frac{1}{p-1}\sum_{\theta \in \mathcal{D}}{\mathfrak{g}(\overline{\theta})\theta(a)\sum_{\substack{n \geq 1\\x \in \F_p^{\times}}}{(\psi\theta)(x)(\psi\theta)(n)a_n(f)e^{2i \pi n\tau} \otimes [\Delta_{x,1}]}}\\
 &= \sum_{\theta}{\frac{\mathfrak{g}(\overline{\theta})\theta(a)}{p-1}f_{\psi\theta} \otimes \sum_{x \in \F_p^{\times}}{(\psi\theta)(x)[\Delta_{x,1}]}}= \sum_{\theta}{\frac{\theta(a)}{p-1}\mathfrak{g}(\overline{\theta})(j' \circ \iota_{p,p})(e_{\psi\theta})}.
 \end{align*}

Now, one formally has $(j' \circ \iota_{p,p})([h \otimes \Delta_{b,1}] \mid_2 W) = h\mid_2\begin{pmatrix}p & 0\\0 &1\end{pmatrix} \delta_{b^{-1}}\begin{pmatrix}0 & -1/p\\p & 0\end{pmatrix} \otimes [\Delta_{b,1}]$, so that
\begin{align*}
j'(\iota_{p,p}(e_{\psi}\mid W)) &= \sum_{x \in \F_p^{\times}}{\psi(x)\overline{\psi^2\chi(x)}(f_{\psi} \mid \begin{pmatrix}0 & -1\\p^2 & 0\end{pmatrix}\otimes [\Delta_{x,1}]}= \lambda_p(f_{\psi})\sum_{x \in \F_p^{\times}}{\overline{\chi\psi(x)} \overline{f}_{\overline{\psi}}\otimes [\Delta_{x,1}]}\\
&= \lambda_p(f_{\psi})\sum_{x \in \F_p^{\times}}{\overline{\chi\psi(x)} f_{\overline{\psi\chi}}\otimes [\Delta_{x,1}]} = \lambda_p(f_{\psi}) e_{\overline{\psi\chi}},
\end{align*}
so we are done. 
}

\defi[left-cusp-rep-f]{Given $(f,\chi) \in \mathscr{C}$, we define $C_f = \mrm{Hom}_{\C}(\Omega^1[f],\C)$: it is an irreducible cuspidal left\footnote{In particular, note that a matrix $M \in \GL{\F_p}$ maps $u \in C_f$ to $u(-\mid M)$: in particular, this operation preserves the central character instead of inverting it.} representation of $\GL{\F_p}$. By Appendix \ref{cuspidal-reps-definition}, $C_f$ is attached to a pair $\{\phi,\phi^p\}$ of characters of $\F_{p^2}^{\times}$. We will often arbitrarily choose a character in this pair and call it \emph{the} character attached to $C_f$. Every time that we do this, the result we are interested in does not depend on which character we chose.}

\rem{Given $(f,\chi)$ as in Definition \ref{left-cusp-rep-f}, it turns out that the character of $\F_{p^2}^{\times}$ attached to $C_f$ characterizes the automorphic representation of $\GL{\Q_p}$ attached to $f$ and therefore, by the local Langlands correspondence, the restriction to a decomposition group at $p$ of the Galois representation attached to $f$. This will prove to be important in later computations of local constants. We will discuss this relationship in more detail later, so let us simply derive, for later reference, a few elementary consequences of the construction. }

\smallskip

\cor[twisting-Cf]{Let $(f,\chi) \in \mathscr{C}$ and $\psi \in \mathcal{D}$. Then $C_{f\otimes \psi} \simeq C_f \otimes \psi(\det)$. In particular, if $\phi$ is one of the characters of $\F_{p^2}^{\times}$ attached to $C_f$, then $\phi\cdot \psi(N_{\F_{p^2}/\F_p})$ is one of the characters attached to $C_{f \otimes \psi}$.}

\demo{The second part of the claim is a consequence of the first part and Corollary \ref{cuspidal-twist}. To prove the first part of the claim, let $(e_{\alpha})_{\alpha\in\mathcal{D}},(e'_{\alpha})_{\alpha\in \mathcal{D}}$ be the bases of $\Omega^1[f]$ and $\Omega^1[f\otimes \psi]$ constructed in Corollary \ref{primitive-is-cuspidal}. 

The $\C$-linear homomorphism $e_{\alpha} \longmapsto e'_{\alpha\psi^{-1}}$ is a linear isomorphism $\Omega^1[f] \rar \Omega^1[f \otimes \psi]$. By Corollary \ref{primitive-is-cuspidal}, the induced morphism $\Omega^1[f] \otimes \psi(\det) \rar \Omega^1[f \otimes \psi]$ is equivariant for the (right) actions of $U$, $\Delta_{a,1}$ and $W$, so it is a $\C[\GL{\F_p}]$-isomorphism, and we are done.
} 

\cor[various-computations-cuspidal]{Let $(f,\chi) \in \mathscr{C}$ and $\phi: \F_{p^2}^{\times} \rar \C^{\times}$ be a character attached to $C_f$. Then the following identities hold:
\begin{itemize}[noitemsep,label=\tiny$\bullet$]
\item $\phi_{|\F_p^{\times}}=\chi$. 
\item Let $z \in \F_{p^2}^{\times} \backslash \F_p^{\times}$ have trace and norm $a,N \in \F_p^{\times}$ respectively (this amounts to requiring that $z^2\notin \F_p^{\times}$). Then 
\[-\phi(z)-\phi^p(z)=\frac{\chi(a)}{p-1}\sum_{\psi \in \mathcal{D}}{\psi\left(\frac{a^2}{N}\right)\lambda_p(f_{\psi})\mathfrak{g}(\overline{\chi\psi^2})}.\] 
\item Let $z \in \F_{p^2}^{\times} \backslash \F_p$ with norm $N$ such that $z^2 \in \F_p^{\times}$. Then 
\[-\phi(z)=-\phi^p(z)=\frac{1}{2}\sum_{\psi^2\chi=1}{\psi^{-1}(N)\lambda_p(f_{\psi})}.\]
\item If $\psi \in \mathcal{D}$ satisfies $\psi^2\chi=1$ and $a \in \F_p^{\times} \backslash \F_p^{\times 2}$, then $\phi(\sqrt{a})=-\psi^{-1}(-a)\lambda_p(f_{\psi})$.  
\item If $\chi=\mathbf{1}$ and $\lambda\in \mathcal{D}$ is the nontrivial quadratic character, then $\lambda_p(f)=-\lambda(-1)\lambda_p(f_{\lambda})$. 
\end{itemize}
}

\demo{For any $a \in \F_p^{\times}$, $aI_2$ acts on $C_f$ by $\phi(a)$ by definition of $\phi$, and also acts by $\chi(a)$ by construction.

For the second statement, $-\phi(z)-\phi^p(z)$ is the trace of the action of $\begin{pmatrix}0 & -N\\1 & a\end{pmatrix} = \Delta_{N,1}WU^a$ on $\Omega^1[f,1]$ by construction. Let $\psi \in \mathcal{D}$, then 
\begin{align*}
e_{\psi}\mid \Delta_{N,1}WU^a &= \psi^{-1}(N)e_{\psi}\mid WU^a = \psi^{-1}(N)\lambda_p(f_{\psi})e_{\overline{\chi}/\psi}\mid U^a\\
&= \frac{\psi^{-1}(N)}{p-1}\lambda_p(f_{\psi})\sum_{\theta \in \mathcal{D}}{\theta(a)\mathfrak{g}(\overline{\theta})e_{\theta\overline{\chi}/\psi}}.
\end{align*} 
Since $\theta \overline{\chi}/\psi=\psi$ if and only if $\theta=\chi\psi^2$, the trace of the action of $\Delta_{N,1}WU^a$ on $\Omega^1[f]$ is \[\sum_{\psi}{\frac{\psi^{-1}(N)}{p-1}\lambda_p(f_{\psi})\mathfrak{g}(\overline{\chi\psi^2})\chi\psi^2(a)}=\frac{\chi(a)}{p-1}\sum_{\psi}{\psi\left(\frac{a^2}{N}\right)\lambda_p(f_{\psi})\mathfrak{g}(\overline{\chi\psi^2})}.\]

For the third statement, since $\phi(-1)=\chi(-1)=1$, $-\phi(z)=-\phi(z^p)=\frac{1}{2}\mrm{Tr}(D_{N,1}W \mid \Omega^1[f])$. As above, $e_{\psi}\mid \Delta_{N,1}W=\psi^{-1}(N)\lambda_p(f_{\psi})e_{\overline{\chi}/\psi}$, so the trace is \[\sum_{\psi=\overline{\psi}/\chi}{\psi^{-1}(N)\lambda_p(f_{\psi})}=\sum_{\psi^2\chi=1}{\psi^{-1}(N)\lambda_p(f_{\psi})}.\]

For the fourth statement, let $a=z^2$. Since $\chi(-1)=1$, $\chi$ has two square roots $\alpha,\beta \in \mathcal{D}$. We already proved that $-\alpha(-a)\lambda_p(f_{\overline{\alpha}})-\beta(-a)\lambda_p(f_{\overline{\beta}})=2\phi(z)$. Since $f_{\overline{\alpha}},f_{\overline{\beta}} \in \mathcal{S}_2(\Gamma_0(p^2))$, $\lambda_p(f_{\overline{\alpha}})^2=\lambda_p(f_{\overline{\beta}})^2=1$, and $(-\alpha(-a))^2=(-\beta(-a))^2=\chi(-a)=\chi(a)$. Thus $\alpha(-a)\lambda_p(f_{\overline{\alpha}})$ and $\beta(-a)\lambda_p(f_{\overline{\beta}})$ have the same square and their sum is not zero, so they are equal.

In particular, if $\chi=\mathbf{1}$, $\phi(z)=-\lambda(-a)\lambda_p(f_{\lambda}) = -\lambda_p(f)$, since $\lambda$ and $\mathbf{1}$ are the two elements of $\mathcal{D}$ with square $\mathbf{1}$. The last point follows from the fact that $\lambda(a)=-1$ by definition of $\lambda$ and $a$. 
}

\cor[old-is-principal]{\begin{itemize}[label=$-$,noitemsep]
\item Let $(f,\mathbf{1}) \in \mathscr{S}$ and $\psi \in \mathcal{D}$. Then $\Omega^1[f,\psi]$ is (as a right $\C[\GL{\F_p}]$-module) the dual (in the same sense as Definition \ref{left-cusp-rep-f}) of a twisted Steinberg representation of central character $\psi^2$. 
\item Let $(f,\chi) \in \mathscr{P}$ and $\psi \in \mathcal{D}$. Then, in the same sense as above, $\Omega^1[f,\psi]$ is the dual of an irreducible principal series representation of $\GL{\F_p}$ with central character $\chi\psi^2$.
\end{itemize}
}

\demo{The claim about the central character is true by construction. Assume that $\Omega^1[f,\psi]$ (in either case) is not irreducible as a $\C[\GL{\F_p}]$-module. Then it has an irreducible subrepresentation of dimension at most $\frac{p+1}{2}$. Since $\frac{p+1}{2} < p-1$, the description of irreducible representations of $\GL{\F_p}$ shows that this irreducible subrepresentation must be a line, which contradicts Corollary \ref{no-lines}. 

Thus $\Omega^1[f,\psi]$ is an irreducible right complex representation of $\GL{\F_p}$ of dimension $p$ or $p+1$: it is thus the dual of a Steinberg representation, or of an irreducible principal series representation. }

\section{Uniformizing $u$}

\prop[uniformize-u]{Let us endow $X_1(p,\mu_p)_{\C}^S$ and $X(p,p)_{\C}$ with the analytifications defined in Propositions \ref{uniformize-x1pmu} and \ref{uniformize-xpp}. Let $(\tau,b) \in \HH \times \F_p^{\times}$. Then, for $x \in S$, the analytification of $u_{x}$ maps $(\Gamma(p)\tau,b)$ to the following points:
\begin{align*}
x=(t,1),\, t \in \F_p &:\quad  u_x(\Gamma(p)\tau,b) = (\Gamma_1(p)T^{bt}\tau,b,x),\\
x=(1,0) &:\quad u_x(\Gamma(p)\tau,b) = (\Gamma_1(p)\delta_bW\tau,b,x).
\end{align*}  
}

\demo{Let $t_0,b_0 \in \Z$ lifting $t,b$. One has, using \cite[(2.8.5.2)]{KM}, with $\omega=e^{-\frac{2i\pi b}{p}}$,
\begin{align*}
&u_{(t,1)}\left(\frac{\C}{\tau\Z\oplus\Z},\frac{b\tau}{p},\frac{1}{p}\right) = \left(\omega, (t,1),\left(\frac{\C}{\tau\Z\oplus\Z},\frac{tb\tau+1}{p}\right)\right) \\
&\quad \quad= \left(\omega,(t,1),\left(\frac{\C}{\tau\Z\oplus(t_0b_0\tau+1)\Z},\frac{b_0t_0\tau+1}{p}\right)\right) \\
&\quad \quad= \left(\omega,(t,1),\left(\frac{\C}{\tau'\Z\oplus\Z},\frac{1}{p}\right)\right),\quad \text{ with }\tau'=\begin{pmatrix}1 & 0\\t_0b_0 & 1\end{pmatrix}\cdot \tau \in \HH,\\
&u_{(1,0)}\left(\frac{\C}{\tau\Z\oplus\Z},\frac{b\tau}{p},\frac{1}{p}\right) = \left(\omega,(1,0),\left(\frac{\C}{\tau\Z\oplus\Z},\frac{b\tau}{p}\right)\right)= \left(\omega,(1,0),\langle b\rangle\left(\frac{\C}{\tau\Z\oplus\Z}, \frac{\tau}{p}\right)\right)\\
&\quad \quad= \left(\omega,(1,0),\langle b\rangle\left(\frac{\C}{\Z\oplus\tau^{-1}\Z},\frac{1}{p}\right)\right)= \left(\omega,(1,0),\langle b\rangle\left(\frac{\C}{-\tau^{-1}\Z\oplus\Z},\frac{1}{p}\right)\right)\\
&\quad \quad= \left(\omega,(1,0),\left(\frac{\C}{\tau'\Z\oplus\Z},\frac{1}{p}\right)\right), \quad \text{ with } \tau'=\delta_b\begin{pmatrix}0 & -1\\1 & 0\end{pmatrix} \cdot \tau.
\end{align*} 
It is then straightforward to compare with the given analytifications.}

\cor[uniformize-omega1u]{Let $u^{\ast,an}$ be the morphism such that the following diagram commutes (where $\mathcal{S}_2(\Gamma_1(p,\mu_p,S))$ is as in Proposition \ref{hecke-x1pmu}):
\[
\begin{tikzcd}[ampersand replacement=\&]
H^0(X_1(p,\mu_p)^S_{\C},\Omega^1)\arrow{d}{\iota_{1,p}} \arrow{r}{\sim}\& H^0(P(p)_{\C},\Omega^1) \arrow{r}{u^{\ast}} \& H^0(J(p,p),\Omega^1) \arrow{r}{\sim}\& H^0(X(p,p),\Omega^1) \arrow{d}{\iota_{p,p}}\\
\mathcal{S}_2(\Gamma_1(p),\F_p^{\times},S)\arrow{rrr}{u^{\ast,an}} \& \& \& \bigoplus_{b \in \F_p^{\times}}{\mathcal{S}_2(\Gamma(p)) \otimes \Delta_{b,1}}
\end{tikzcd}
\]  
Then $u^{\ast,an}$ commutes with the given actions of $\GL{\F_p}$ and is given by the following formulas:
\begin{align*}
(f,b,t) \in \mathcal{S}_2(\Gamma_1(p)) \times \F_p^{\times} \times \F_p &: u^{\ast,an}(f \otimes [b] \otimes T^t) = \left(f\mid_2 T^{t/b}\right) \otimes \Delta_{b,1},\\
(f,b) \in \mathcal{S}_2(\Gamma_1(p)) \times \F_p^{\times} &: u^{\ast,an}(f \otimes [b] \otimes W) = (\langle b^{-1}\rangle f)\mid_2 W) \otimes \Delta_{b,1}.
\end{align*}
}

\demo{By construction, for any $x \in S$, if the second row of $M \in \{T^t\mid t \in \F_p\} \cup \{W\}$ is not $x$, $u^{\ast,an}$ vanishes on $\mathcal{S}_2(\Gamma_1(p)) \otimes \C[\F_p^{\times}] \otimes [M]$. So what we need to show instead is that, for any $b \in \F_p^{\times}$: 
\begin{align*}
(f,t) \in \mathcal{S}_2(\Gamma_1(p)) \times \F_p &: u^{\ast}_{(t,1)}(\iota_{1,p}^{-1}(f \otimes [b] \otimes T^t)) = \iota_{p,p}^{-1}\left(\left(f\mid_2 T^{t/b}\right) \otimes \Delta_{b,1}\right),\\
f \in \mathcal{S}_2(\Gamma_1(p))  &: u_{(1,0)}^{\ast}\left(f \otimes [b] \otimes W\right) = \iota_{p,p}^{-1}\left(((\langle b^{-1}\rangle f)\mid_2 W) \otimes \Delta_{b,1}\right).
\end{align*}
For both identities, the differentials vanish on all the connected components of $X(p)\times\F_p^{\times}$ except on the $b^{-1}$ component. Since the differentials are all holomorphic on $X(p)$, it is enough to show equality after pulling-back to $\HH$ by $\tau \in \HH \longmapsto (\Gamma(p)\tau,b^{-1})$.

Let $t \in \F_p$, and choose lifts $t_0,b_0 \in \Z$ of $t,b$ such that $b_0 \mid t_0$. We then have on $\HH \times \{b^{-1}\}$: 
\begin{align*}
u_{(t,1)}^{\ast}(f(\tau)d\tau) &= f(T^{t_0/b_0}\tau)d(T^{t_0/b_0}\tau) = (f\mid_2 T^{t_0/b_0})(\tau)d\tau,\\   
u_{(1,0)}^{\ast}(f(\tau)d\tau) &= f(\delta_{b^{-1}}W\tau)d(\delta_{b^{-1}}W\tau) = ((\langle b^{-1}\rangle f)\mid_2 W)(\tau)d\tau.
\end{align*}  
}

\cor[becomes-naive]{There is a commutative diagram of right $\C[\GL{\F_p}]$-modules,
\[
\begin{tikzcd}[ampersand replacement=\&]
H^0(X_1(p,\mu_p)^S_{\C},\Omega^1) \arrow{r}{u^{\ast}} \arrow{d} \& H^0(X(p,p)_{\C},\Omega^1)\arrow{d}\\
\mathcal{S}_2(\Gamma_1(p)) \otimes_{\Z[B \cap \SL{\F_p}]} \Z[\GL{\F_p}] \arrow{r}\& \mathcal{S}_2(\Gamma(p)) \otimes_{\Z[\SL{\F_p}]} \Z[\GL{\F_p}].
\end{tikzcd}
\]
where the vertical arrows are isomorphisms, and the bottom arrow is the natural one (given by the inclusion $\mathcal{S}_2(\Gamma_1(p)) \subset \mathcal{S}_2(\Gamma(p))$).   
}

\demo{The obvious identification $\bigoplus_{b \in \F_p^{\times}}{\mathcal{S}_2(\Gamma(p)) \otimes \Delta_{b,1}} \simeq \mathcal{S}_2(\Gamma(p)) \otimes_{\C[\SL{\F_p}]} \C[\GL{\F_p}]$ is an isomorphism of $\C[\GL{\F_p}]$-modules by definition of the action on the right-hand side, and we let the rightmost map in the diagram be precisely the composition of $\iota_{p,p}$ with this identification. 
On the other hand, we can identify $(\mathcal{S}_2(\Gamma_1(p)) \otimes \C[\F_p^{\times}]) \otimes_{\C[B]} \C[\GL{\F_p}]$ with $\mathcal{S}_2(\Gamma_1(p)) \otimes_{\C[B \cap \SL{\F_p}]} \C[\GL{\F_p}]$ as follows: we map $(f \otimes [b]) \otimes M$ to $f \otimes \Delta_{b,1}M$. This construction is well-defined, because if $M=\begin{pmatrix}a & \ast \\0 & d\end{pmatrix}M'$ (with $M,M' \in \GL{\F_p}$), then \begin{align*}
(f \otimes [b]) \otimes M &= ((f \otimes [b]) \mid \begin{pmatrix}a & \ast\\0 & d\end{pmatrix}) \otimes M'= ((\langle d\rangle f) \otimes [bad]) \otimes M'\\
f \otimes \Delta_{b,1}M &= f \otimes \begin{pmatrix}ba & \ast \\0 & d\end{pmatrix}M' = f \otimes \begin{pmatrix} d^{-1} & \ast\\0 & d\end{pmatrix}\Delta_{bad,1}M' = (\langle d\rangle f) \mid \Delta_{bad,1}M'.
\end{align*}
That this identification is $\C[\GL{\F_p}]$-equivariant is then clear, and we let the leftmost vertical map be the composition of $\iota_{1,p}$ with this identification. Now, all we need to do is check that the diagram commutes, which is a direct consequence from the formulae in Corollary \ref{uniformize-omega1u}. 
}

\cor[sequence-is-exact]{The sequence 
\[H^0(P(p)_{\C},\Omega^1) \overset{R^{\ast}}{\longrightarrow} H^0(P(p)_{\C},\Omega^1) \overset{u^{\ast}}{\longrightarrow} H^0(J(p,p)_{\C},\Omega^1)\]
is exact. The image of $u^{\ast}$ is exactly the sum of all the non-cuspidal irreducible $\C[\GL{\F_p}]$-submodules of $H^0(J(p,p)_{\C},\Omega^1)$. In other words, 
\[\im{u^{\ast}} \oplus \bigoplus_{(f,\chi) \in \mathscr{C}}{\Omega^1[f]} \simeq H^0(J(p,p)_{\C},\Omega^1).\]}

\demo{Let $M$ be a right $\C[\SL{\F_p}]$-module and $M'=M \otimes_{\C[\SL{\F_p}]} \C[\GL{\F_p}]$. For \[m = \sum_{b \in \F_p^{\times}}{m_b\otimes\Delta_{b,1}} \in M',\] one has $m \mid U = \sum_{b \in \F_p^{\times}}{(m_b \mid U^b)\otimes \Delta_{b,1}}$, so that \[(M')^U = \bigoplus_{b \in \F_p^{\times}}{M \otimes \Delta_{b,1}}=M^U \otimes_{\C[B \cap \SL{\F_p}]} \C[B].\] 

Therefore, by Corollary \ref{becomes-naive} and the discussion of irreducible representations of $\GL{\F_p}$ in Appendix \ref{reps-gl2}, the image of $u^{\ast}$ is exactly the sum of all irreducible $\C[\GL{\F_p}]$-submodules which are subquotients of principal series representations. By Corollaries \ref{first-decomp-xpp}, \ref{primitive-is-cuspidal}, \ref{old-is-principal}, the image of $u^{\ast}$ is the sum of the $\Omega^1[f,\psi]$ for $(f,\chi) \in \mathscr{S} \cup \mathscr{P}$ and $\psi\in \mathcal{D}$. 

Let $(f,\chi) \in \mathscr{S}$ (resp. $(f,\chi) \in \mathscr{P}$) and $\psi\in \mathcal{D}$. Let $V_{f,\psi}$ be the $\C[\GL{\F_p}]$-submodule of $H^0(P(p)_{\C},\Omega^1)$ where $T'_{\ell}$ (resp. $\langle n\rangle'$) acts by multiplication by $a_{\ell}(f)\psi(\ell)$ (resp. $\chi(n)\psi(n)^2$). Since $R^{\ast}$ commutes to $\mathbb{T}'$ and $T'_{\ell}u=uT_{\ell}$, one has:
\begin{align*}
\Omega^1[f,\psi] \cap u^{\ast}H^0(P(p)_{\C},\Omega^1) &= u^{\ast}V_{f,\psi},\\
V_{f,\psi} \cap R^{\ast}H^0(P(p)_{\C},\Omega^1) &= R^{\ast}V_{f,\psi}. 
\end{align*}

It is thus enough to prove that $V_{f,\psi} \overset{R^{\ast}}{\rar} V_{f,\psi} \overset{u^{\ast}}{\rar} \Omega^1[f,\psi]$ is exact. 

By Lemma \ref{complex-of-R}, we know that the sequence is a complex of $\C[\GL{\F_p}]$-modules. By Lemma \ref{diagonalize-x1pmu-Tprime-R}, $R^{\ast}$ is never trivial, and we saw that $u^{\ast}$ was surjective. By Corollary \ref{old-is-principal}, $\Omega^1[f,\psi]$ has length one as a $\C[\GL{\F_p}]$-module. By Proposition \ref{diagonalize-x1pmu-Tprime}, $V_{f,\psi}$ has length two as a $\C[\GL{\F_p}]$-module, and the conclusion follows. 

We already know that $u^{\ast}R^{\ast}=0$. To show exactness, we may reduce to subspaces where every $T_{\ell}$ and every $m_n$ (resp. $T'_{\ell},\langle n\rangle'$) acts by a scalar. Then the relevant submodule of $H^0(J(p,p),\Omega^1)$ has length one over $\C[\GL{\F_p}]$, and is in the image of $u^{\ast}$ by the above. But the image of $R^{\ast}$ also has length at least one, and we saw that the relevant submodule of $H^0(P(p)_{\C},\Omega^1)$ had length two, whence the desired exactness.  
}

\cor[decomposition-hecke-xpp]{$H^0(X(p,p)_{\C},\Omega^1)$ is the direct sum of the following right $\mathbb{T}[\GL{\F_p}]$-modules, on which every element of $\mathbb{T}$ acts as a scalar:
\begin{itemize}[label=\tiny$\bullet$, noitemsep]
\item For $(f,\mathbf{1}) \in \mathscr{S}$, and $\psi \in \mathcal{D}$, $\Omega^1[f,\psi]$, which is isomorphic to the twisted Steinberg $\mrm{St}_{\psi}^{\vee}$. 
\item For $(f,\chi) \in \mathscr{P}$ and $\psi \in \mathcal{D}$, modulo the relation $(f,\psi) \sim (\overline{f},\chi\psi)$, $\Omega^1[f,\psi]$, which is isomorphic to the irreducible principal series $\pi(\psi,\chi\psi)^{\vee}$, 
\item For $(f,\chi) \in \mathscr{C}$, $\Omega^1[f]$ which is the irreducible cuspidal representation $C_f^{\vee}$.
\end{itemize}
In all these cases (with $\chi=\mathbf{1}$ in the first one, and $\psi=\mathbf{1}$ in the last one), the operator $T_n$ acts on the submodule as the scalar $a_{n}(f)\psi(n)$, and $m_n$ acts by the scalar $\psi(n)^2\chi(n)$. 
}

\demo{By Corollary \ref{sequence-is-exact} and Proposition \ref{diagonalize-x1pmu-Tprime}, since the action of $\mathbb{T}$ is diagonalizable, $H^0(X(p,p)_{\C},\Omega^1)$ is the direct sum of the $\Omega^1[f]$ for $(f,\chi) \in \mathscr{C}$ and of the images under $u^{\ast}$ of the $H^0(X_1(p,\mu_p)^S_{\C},\Omega^1)[f,\psi]$ for $(f,\chi) \in \mathscr{S} \cup \mathscr{P}$ and $\psi \in \mathcal{D}$, modulo the relation $(f,\chi,\psi) \sim (\overline{f},\overline{\chi},\chi\psi)$. The action of $\mathbb{T}$ is given by Corollary \ref{first-decomp-xpp-Tn}. 

By Corollaries \ref{old-is-principal} and \ref{sequence-is-exact}, $\Omega^1[f,\psi]$ is an irreducible quotient of $H^0(X_1(p,\mu_p)^S_{\C},\Omega^1)[f,\psi]$, isomorphic to a twisted Steinberg if $\chi=\mathbf{1}$ and an irreducible principal series otherwise. By Proposition \ref{diagonalize-x1pmu-Tprime}, $\Omega^1[f,\psi]$ is isomorphic to $\mrm{St}^{\vee}_{\psi}$ if $\chi=\mathbf{1}$ and $\pi(\psi,\chi\psi)^{\vee}$ otherwise (since $\pi(\chi\psi,\psi)$ and $\pi(\psi,\chi\psi)$ are isomorphic $\C[\GL{\F_p}]$-modules). 
}

\section{Connected components of $C(p)$ and $\ker{R}$}
\label{connected-components-kernels}

\defi{By Proposition \ref{Pp-over-GL2}, there exists a unique morphism $v: P(p) \rar J(p,p)$ of abelian schemes which commutes to the action of $\GL{\F_p}$ and such that the morphism $v \circ [(\iota_{(0,1)})_{\ast}]: J_1(p,\mu_p) \rar J(p,p)$ is exactly $u_{(0,1)}^{\ast}$. We write $q=u \circ v$: it is a $\GL{\F_p}$-endomorphism of $P(p)$. 
We denote by $\Delta: P(p) \rar P(p)$ the sum of the $\langle n\rangle$. 
}

\lem[descr-q]{Let $x,y \in S$. The morphism \[J_1(p,\mu_p) \overset{[(\iota_x)_{\ast}]}{\longrightarrow} P(p) \overset{q}{\longrightarrow} P(p) \overset{\iota_y^{\ast}}{\longrightarrow} J_1(p,\mu_p)\] is given by:
\begin{itemize}[noitemsep,label=\tiny$\bullet$,topsep=0pt]
\item when $x=y$, it is multiplication by $p$,
\item when $x \neq y$ and $(1,0) \in \{x,y\}$, it is $w_pU_p$,
\item when $x=(s,1),y=(t,1)$ with $s \neq t$, it is $\langle s-t\rangle w_pU_p$.
\end{itemize}
In particular, $q$ commutes with $\Delta$.}

\demo{First assume that $x=(0,1)$. When $y=(0,1)$, the morphism is exactly $[u_{(0,1)}]_{\ast}u_{(0,1)}^{\ast}$, where $u_{(0,1)}: X(p,p) \rar X_1(p,\mu_p)$ is a morphism of smooth proper connected curves (over $\Z[1/p]$). By considering its complex analytification from Proposition \ref{uniformize-u}, $u_{(0,1)}$ has degree $p$, so that $[u_{(0,1)}]_{\ast}u_{(0,1)}^{\ast}$ is exactly the multiplication by $p$. 

When $y=(t,1)$ (resp. $y=(1,0)$ and $t=1$), the morphism is exactly $[(u_{y})_{\ast}]u_{(0,1)}^{\ast}$. Consider the following diagram:
\[
\begin{tikzcd}[ampersand replacement=\&]
X_1(p,\mu_p) \arrow{d}{\mrm{id}} \& X(p,p) \arrow{l}{u_{(0,1)}} \arrow{d}{\alpha} \arrow{r}{u_{y}} \& X_1(p,\mu_p)\arrow{d}{w_p\langle t\rangle^{-1}}\\
X_1(p,\mu_p) \& X(\Gamma_1'(p,p),\mu_p) \arrow{l}{D'_{1,1}} \arrow{r}{D'_{p,1}} \& X'_1(p,\mu_p).
\end{tikzcd}
\]
In this diagram, $\alpha$ is the compactification of the map \[(E,P,Q) \longmapsto (E,Q,\langle P+t^{-1}Q\rangle,\langle P,\,Q\rangle) \text{ (resp. }(E,P,Q) \longmapsto (E,Q,\langle P\rangle,\langle P,\,Q\rangle)\text{).}\] It is then direct (using Section \ref{gamma1m-structure}, \cite[(2.8.2),(2.8.6.1)]{KM}, and the fact that $\langle -1\rangle$ is trivial) that the diagram commutes away from the cusps with $\alpha$ being an isomorphism, hence, by Proposition \ref{compactification-functor}, the diagram commutes with $\alpha$ being an isomorphism. By Proposition \ref{relative-picard-mixed-functoriality}, $[(u_{y})_{\ast}]u_{(0,1)}^{\ast} = [((w_p\langle t\rangle^{-1})^{-1})_{\ast}]U_p=\langle t\rangle w_pU_p$.   

To describe the other coordinates, use the above and notice that the following diagrams commute (since $q$ is $\GL{\F_p}$-equivariant and $\langle -1\rangle$ is trivial), for $s,t \in \F_p$: 
\[
\begin{tikzcd}[ampersand replacement=\&]
J_1(p,\mu_p) \arrow{r}{[(\iota_{(t,1)})_{\ast}]} \arrow{dr}{\iota} \& P(p) \arrow{r}{q} \& P(p) \arrow{r}{\iota_{(s,1)}^{\ast}} \& J_1(p,\mu_p)\\
\& P(p) \arrow{u}{T^{-t}} \arrow{r}{q}\& P(p) \arrow{u}{T^{-t}} \arrow{ru}{\alpha}\&\\
J_1(p,\mu_p) \arrow{r}{[(\iota_{(1,0)})_{\ast}]} \arrow{dr}{\iota} \& P(p) \arrow{r}{q} \& P(p) \arrow{r}{\iota_{(s,1)}^{\ast}} \& J_1(p,\mu_p)\\
\& P(p) \arrow{u}{W} \arrow{r}{q}\& P(p) \arrow{u}{W} \arrow{ru}{\beta}\&\\
J_1(p,\mu_p) \arrow{r}{[(\iota_{(1,0)})_{\ast}]} \arrow{dr}{\iota} \& P(p) \arrow{r}{q} \& P(p) \arrow{r}{\iota_{(1,0)}^{\ast}} \& J_1(p,\mu_p)\\
\& P(p) \arrow{u}{W} \arrow{r}{q}\& P(p) \arrow{u}{W} \arrow{ru}{\gamma}\&
\end{tikzcd}
\]
where $\alpha=\iota_{(s,1)}^{\ast} \circ (T^t)^{\ast} = \iota_{(s-t,1)}^{\ast}$, $\beta=\iota_{(s,1)}^{\ast} \circ (W^{-1})^{\ast} = \langle s\rangle\iota_{(-1/s,1)}^{\ast}$, and $\gamma=\iota_{(1,0)}^{\ast} \circ (W^{-1})^{\ast} = \iota_{(0,1)}^{\ast}$ and $\iota=[(\iota_{(0,1)})_{\ast}]$.
}

\lem[unif-wp]{The uniformization of the automorphism $w_p$ of $X_1(p,\mu_p)$ is given by \[(\Gamma_1(p)\tau,b) \longmapsto (\delta_b\Gamma_1(p)\begin{pmatrix}0 & -1\\p & 0\end{pmatrix}\tau,b).\] }

\demo{Let $\tau \in \HH$. Then $(\tau,b)$ corresponds to the class of $\left(\frac{\C}{\tau\Z\oplus\Z},\frac{1}{p},e^{-2i\pi b/p}\right)$. Let $E=\frac{\C}{\tau\Z\oplus\Z}$ and $\pi: E \rar E'$ the isogeny whose kernel is $\langle \frac{1}{p}\rangle$. Thus $w_p(\tau,b)$ is the class of the enriched elliptic curve $\left(\frac{\C}{\tau\Z\oplus\frac{1}{p}\Z},\pi(Q),e^{-2i\pi b/p}\right)$, where $Q \in E(\C)$ is any $p$-torsion point such that $\langle \frac{1}{p},\,Q\rangle_{E[p]}=\langle \frac{1}{p},\, \pi(Q)\rangle_{\pi} = e^{-2i\pi b/p}$. Thus $Q=\frac{-b\tau}{p}$ works, so that $w_p(\tau,b)$ represents the enriched elliptic curve
\begin{align*}
\left(\frac{\C}{\tau\Z\oplus\frac{1}{p}\Z},\frac{-b\tau}{p},e^{-2i\pi b/p}\right) &\simeq \left(\frac{\C}{\Z\oplus\frac{1}{p\tau}\Z},\frac{-b}{p},e^{-2i\pi b/p}\right)\\
&\simeq \left(\frac{\C}{\frac{-1}{p\tau}\Z\oplus\Z},-b\frac{1}{p},e^{-2i\pi b/p}\right) \simeq \langle b\rangle \left(\frac{\C}{\frac{-1}{p\tau}\Z\oplus\Z},\frac{1}{p},e^{-2i\pi b/p}\right),
\end{align*}
whence the conclusion. 
}

\lem[unif-Up]{The pull-back by $U_p$ on $H^0(J_1(p,\mu_p)_{\C},\Omega^1) \overset{\iota_{1,p}}{\simeq} \mathcal{S}_2(\Gamma_1(p)) \otimes \C[\F_p^{\times}]$ is exactly the classical Hecke operator $U_p$.}

\demo{By the proof of Proposition \ref{descr-q}, $U_p=w_p[(u_{(1,0)})_{\ast}]u_{(0,1)}^{\ast}$  
We proceed as in the proof of Proposition \ref{hecke-xpp}. Consider the following diagram: 
\[
\begin{tikzcd}[ampersand replacement=\&]
\HH \times \F_p^{\times} \arrow{d}{\upsilon}\& \HH \times \F_p^{\times} \arrow{l}{h} \arrow{r}{(g_w)_{w}} \arrow{rd}{\upsilon} \& \coprod\limits_{0 \leq w < p}{\HH \times \F_p^{\times}} \arrow{d}{\pi} \arrow{r}{h'} \& \HH \times \F_p^{\times}\arrow{d}{\upsilon}\\
Y_1(p,\mu_p)_{\C} \& \& Y(p,p)_{\C} \arrow{ll}{w_pu_{1,0}} \arrow{r}{u_{(0,1)}} \& Y(p,p)_{\C}
\end{tikzcd}
\]
where $h(\tau,z)=(\frac{\tau}{p},z)$, $h'(\tau,z)=(\tau,z)$, $\pi$ maps $(\tau,z)$ on the copy $w$ of $\HH \times \F_p^{\times}$ to the class of the enriched elliptic curve $\left(\frac{\C}{\tau\Z\oplus\Z},\frac{z(\tau+w)}{p},\frac{1}{p}\right)$, and $g_w(\tau,z)=(\tau-w,z)$. 

By Proposition \ref{uniformize-u} and Lemma \ref{unif-wp}, and, since, for suitable choices of $M,M' \in SL_2(\Z)$ congruent to $\Delta_{b^{-1},b}$ mod $p$, $M\begin{pmatrix}0 & -1\\p& 0\end{pmatrix}M'\begin{pmatrix}0 & -1\\1 &0\end{pmatrix}=\begin{pmatrix}1 & 0\\0 &p\end{pmatrix}$, the diagram commutes. Moreover, this diagram is Cartesian above a suitable co-finite open subscheme of $Y_1(p,\mu_p)$. 

Therefore, for $(f,b) \in  \mathcal{S}_2(\Gamma_1(p)) \times \F_p^{\times}$, one has \[U_p^{\ast}(\iota_{1,p}(f \otimes [b])) = \iota_{1,p}^{-1}\left[\sum_{0 \leq w < p}{f\mid\begin{pmatrix}1 & 0\\0 & p\end{pmatrix}\begin{pmatrix}1 & w\\0 & 1\end{pmatrix}} \otimes [b]\right] = \iota_{1,p}^{-1}(U_p(f) \otimes [b]),\] whence the conclusion.}

\medskip

\cor{The endomorphism $(w_pU_p+1)\Delta$ of $P(p)$ vanishes. }

\demo{It is enough to show that $(w_pU_p+1)\Delta$ is the null endomorphism of $J_1(p,\mu_p)_{\Q}$ (since this morphism preserves the summands of $P(p)$), and this is equivalent to the fact that the pull-back by $\Delta(w_p+U_p)$ is the trivial endomorphism of $H^0(J_1(p,\mu_p)_{\C},\Omega^1) \overset{\iota_{1,p}}{\simeq} \mathcal{S}_2(\Gamma_1(p)) \otimes \C[\F_p^{\times}]$.  

By Proposition \ref{hecke-x1pmu}, $\Delta^{\ast}$ acts on $\mathcal{S}_2(\Gamma_1(p)) \otimes \C[\F_p^{\times}]$ as the sum of all the classical diamond operators of $\mathcal{S}_2(\Gamma_1(p))$: thus $\frac{\Delta^{\ast}}{p-1}$ is a projection on the subspace $\mathcal{S}_2(\Gamma_0(p)) \otimes \C[\F_p^{\times}]$. Therefore, by Lemma \ref{unif-wp} $w_p^{\ast}\Delta^{\ast}$ is exactly the composition of $\Delta^{\ast}$ with the action of $\begin{pmatrix}0 & -1\\p & 0\end{pmatrix}$. 

Thus, by Lemmas \ref{unif-wp}, \ref{unif-Up}, for every $f \in \mathcal{S}_2(\Gamma_0(p))$, one has \[(w_p^{\ast}+U_p^{\ast})(f \otimes [b]) = \left(\sum_{i=0}^{p-1}{f \mid \begin{pmatrix}1 & i\\0 & p\end{pmatrix}+\begin{pmatrix}0 & -1\\p & 0\end{pmatrix}}\right) \otimes [b] := g \otimes [b].\] 

We can check that $\begin{pmatrix}0  & -1\\p & 0\end{pmatrix}=\begin{pmatrix}1 & 0\\0 & p\end{pmatrix}\begin{pmatrix}0 & -1\\1 & 0\end{pmatrix}$ and the $\begin{pmatrix} 1 & i\\0 & p\end{pmatrix}$ (for $0 \leq i < p$) form representatives for $\Gamma_0(p) \backslash \Gamma_0(p)\begin{pmatrix}1 & 0\\0 & p\end{pmatrix}\Gamma(1)$, since $\Gamma_0(p)\mrm{diag}(1,p)\Gamma(1)=\begin{pmatrix}1 & 0\\0 & p\end{pmatrix}\Gamma(1)$. Thus $g \in \mathcal{S}_2(\Gamma(1))=\{0\}$, whence the conclusion. }

\prop[descr-r]{There exists an endomorphism $r$ of $P(p)$ commuting with $\Delta$ such that $q-(2p-\Delta) = r \circ R$. }

\demo{We define $r$ such that for any $x,y \in S$, the morphism \[J_1(p,\mu_p) \overset{[(\iota_x)_{\ast}]}{\rar} P(p) \overset{r}{\rar} P(p) \overset{\iota_y^{\ast}}{\rar} J_1(p,\mu_p)\] is equal to
\begin{itemize}[noitemsep,label=\tiny$\bullet$,topsep=5pt]
\item $w_p\Delta$ when $x=y$,
\item $-w_p$ when $x \neq y$ and one of them is $(1,0)$,
\item $-\langle s-t\rangle w_p$ when $x=(s,1),y=(t,1)$ with distinct $s,t \in \F_p$. 
\end{itemize}
It is clear that $r$ commutes with $\Delta$. 

Given distinct $s,t \in \F_p$, 
\begin{align*}
&\iota_{(t,1)}^{\ast} \circ r \circ R \circ [(\iota_{(s,1)})_{\ast}] \\
&= (-w_p)w_p+(-\langle s-t\rangle w_p)(-U_p)+(w_p\Delta)(\langle s-t\rangle w_p)+\sum_{x \neq s,t}{(-\langle t-x\rangle w_p) \langle x-s\rangle w_p}\\
&= \langle s-t\rangle w_pU_p+w_p\Delta w_p-1-\sum_{x \neq s,t}{\langle \frac{t-x}{x-s}\rangle}\\
&= \langle s-t\rangle w_pU_p+\Delta-1-(\Delta-\langle -1\rangle) = \langle s-t\rangle w_pU_p = \iota_{(t,1)}^{\ast}\circ q\circ [(\iota_{(s,1)})_{\ast}].  
\end{align*}
If $s \in \F_p$,
\begin{align*}
\iota_{(s,1)}^{\ast}\circ r \circ R \circ [(\iota_{(1,0)})_{\ast}] &= (-w_p)(-U_p)+(w_p\Delta)w_p+\sum_{t \neq s}{(-\langle t-s\rangle w_p)w_p} = w_pU_p+\Delta-\Delta\\
&=w_pU_p= \iota_{(s,1)}^{\ast}\circ q\circ [(\iota_{(1,0)})_{\ast}],\\
\iota_{(1,0)}^{\ast}\circ r \circ R \circ [(\iota_{(s,1)})_{\ast}] &= (w_p\Delta)w_p + (-w_p)(-U_p)+\sum_{t \neq s}{(-w_p)(\langle s-t\rangle w_p)}\\
&= \Delta + w_pU_p - \sum_{t \neq s}{\langle (t-s)^{-1}\rangle} = w_pU_p=\iota_{(1,0)}^{\ast}\circ q\circ [(\iota_{(s,1)})_{\ast}],\\
\iota_{(s,1)}^{\ast}\circ r \circ R \circ [(\iota_{(s,1)})_{\ast}] &= (-w_p)(w_p)+(w_p\Delta)(-U_p)+\sum_{t \neq s}{(-\langle t-s\rangle w_p)(\langle t-s \rangle w_p)}\\
&= -p-\Delta = \iota_{(s,1)}^{\ast}\circ q\circ [(\iota_{(s,1)})_{\ast}]-(2p-\Delta),\\
\iota_{(1,0)}^{\ast}\circ r \circ R \circ [(\iota_{(1,0)})_{\ast}] &= (w_p\Delta)(-U_p)+\sum_{t \in \F_p}{(-w_p)w_p} = \Delta-p \\
&= \iota_{(1,0)}^{\ast}\circ q\circ [(\iota_{(1,0)})_{\ast}] - (2p-\Delta).
\end{align*}
}

\cor[identities-q]{The following identities hold: $(2p-\Delta)q = q^2$, $q \circ r \circ R=r \circ R \circ q=0$. 
Moreover, for any field $k$ with characteristic not dividing $p$, $\im{q_k}$ has finite index in $\ker{R_k}$, and $(2p-\Delta)\pi_0(\ker{q_k})=(2p-\Delta)\pi_0(\ker{R_k})=0$.}

\demo{Since $R \circ q=R \circ u \circ v=0$, one has $\mrm{im}(q_k) \subset \ker{R_k}$ for any field $k$ with characteristic not dividing $p$. In particular, $(q-(2p-\Delta))q=r \circ R \circ q=0$. Moreover, if $y \in \ker{R}(\overline{k})$, then $q(y)-(2p-\Delta)y=r(R(y))=0$, hence $(2p-\Delta)y \in \im{q}(\overline{k})$. Thus $(2p-\Delta)\ker{R_k} \subset \im{q_k}$. Since $\Delta^2=(p-1)\Delta$, one has $(p+1+\Delta)(2p-\Delta)=2p(p+1)$ so $2p(p+1)\ker{R_k} \subset \im{q_k}$. Thus, $\im{q_k}$ is the connected component of unity in $\ker{R_k}$. Since $(2p-\Delta)\ker{R_k} \subset \im{q_k}$, we have $(2p-\Delta)\pi_0(\ker{R_k})=0$. 

$q$ commutes with $\Delta$, so $0=q^2-(2p-\Delta)q=q^2-q(2p-\Delta)=q\circ r \circ R$. By Proposition \ref{descr-r}, over any field $k$ with characteristic not dividing $p$, $(2p-\Delta)\ker{q_k} = (r \circ R)(\ker{q_k}) \subset \im{r_k \circ R_k}$. Since $2p(p+1)$ is a multiple of $2p-\Delta$ (so that $2p-\Delta: P(p)_k \rar P(p)_k$ is surjective), $\im{r_k \circ R_k}$ is the connected component of unity of $\ker{q_k}$, and we are done.}

\cor{Over any field $k$ whose characteristic is not $p$, the exponents of $\pi_0(\ker{u_k})$, $\pi_0(\ker{q_k})$, $\pi_0(\ker{R_k})$ divide $2p(p+1)$.  }

\demo{$\im{q}$ and $\im{u}$ are abelian subvarieties of $P(p)$, contained with finite index in $\ker{R}$, so they are equal. Thus, $\im{v}+\ker{u}=J(p,p)$. Since the neutral component $(\ker{u})^0$ of $\ker{u}$ has finite index in $\ker{u}$, $\im{v}+(\ker{u})^0$ contains a multiple of $\im{v}+\ker{u}$, so that $\im{v}+(\ker{u})^0=J(p,p)$. Consequently, $v: \pi_0(\ker{q}) \rar \pi_0(\ker{u})$ is surjective. It is thus enough, by the previous Corollary, to show that $2p(p+1) \in \mrm{End}(P(p))(2p-\Delta)$ (as endomorphisms of $P(p)$), which is true since $2p(p+1)=(p+1+\Delta)(2p-\Delta)$. }

\newpage

%% file: tate-0.tex
\chapter{Tate modules for the Jacobian of $X(p,p)$ and its twists}
\label{tate-modules-twists}

\section{Introduction}

The purpose of this section is to compute the Tate modules of $J(p,p)$ and $P(p)$, as well as that Jacobians of Galois twists of $X(p,p)$ or $X_1(p,\mu_p)^S$, as a $\mathbb{T}[G_{\Q}]$-module (or $\mathbb{T}'[G_{\Q}]$-module). To this effect, we first describe the Tate modules of $J(p,p)$ (resp. $P(p)$) as $\mathbb{T}[G_{\Q} \times \GL{\F_p}]$-modules (resp. $\mathbb{T}'[G_{\Q} \times \GL{\F_p}]$), and prove that the general case follows from Proposition \ref{jacobian-twist}.

We use the letter $\mathcal{T}$ to denote Tate modules, in order to distinguish them from Hecke operators. In this section and the following ones, $\omega_p$ will always denote the mod $p$ cyclotomic character. 

\nott{Let $f \in \mathcal{S}_k(\Gamma_1(N))$ be a newform with primitive character $\chi$. Let $F$ be a number field containing the coefficients of $f$ and $\mathfrak{l}$ be a prime ideal of $\OO_F$ with residue characteristic $\ell$. Let $T_{f,\mathfrak{l}}$ be a continuous representation $G_{\Q} \rar \GL{\OO_{F_{\mathfrak{l}}}}$ satisfying the following properties:
\begin{itemize}[noitemsep,label=\tiny$\bullet$]
\item $T_{f,\mathfrak{l}}$ is unramified at every prime number $q\nmid N\ell$,
\item $\det{T_{f,\mathfrak{l}}}=\chi(\omega_N)\omega^{k-1}$, where $\omega_N$ (resp. $\omega$) is the mod $N$ cyclotomic character (resp. the $\ell$-adic cyclotomic character),
\item If $q\nmid \ell N$ is a prime, the characteristic polynomial of the arithmetic Frobenius at $q$ is $X^2-a_{q}(f)X+q^{k-1}\chi(q)$,
\item $V_{f,\mathfrak{l}} := T_{f,\mathfrak{l}} \otimes_{\OO_{F_{\mathfrak{l}}}} F_{\mathfrak{l}}$ is absolutely irreducible (by \cite[Theorem 2.3]{Antwerp5-Ribet} when $k \geq 2$ and by \cite[Th\'eor\`eme 4.1]{Del-Ser} when $k=1$ -- see also \cite[Th\'eor\`eme 6.1, Remarque 6.5]{Del-Ser}).
\end{itemize}
In particular, while all possible choices for $T_{f,\mathfrak{l}}$ need not be isomorphic over $\OO_{F_{\mathfrak{l}}}$, it is the case for $V_{f,\mathfrak{l}}$. }

The determination of the Tate module of $P(p)$ can be easily reduced to that of $J_1(p)$, whose structure is well-known.
Recall the following notations from Section \ref{analytic}. 

\nott{\begin{itemize}[label=$-$,noitemsep]
\item $\mathcal{D}$ denotes the set of Dirichlet characters $(\Z/p\Z)^{\times} \rar \C^{\times}$, 
\item $\mathscr{S}$ denotes the set of $(f,\mathbf{1})$ where $f \in \mathcal{S}_2(\Gamma_0(p))$ is a newform, 
\item $\mathscr{P}$ denotes the set of $(f,\chi)$, where $f \in \mathcal{S}_2(\Gamma_1(p))$ is a newform with character $\chi \neq \mathbf{1}$,
\item $\mathscr{C}$ denotes the set of $(f,\chi)$, where $f \in \mathcal{S}_2(\Gamma_1(p) \cap \Gamma_0(p^2))$ is a newform with character $\chi \in \mathcal{D}$, and such that no twist of $f$ by a character in $\mathcal{D}$ has conductor $p$. 
\item For $\alpha \in \mathcal{D}$, $\mrm{St}_{\alpha}$ is the Steinberg representation of $\GL{\F_p}$ twisted by $\alpha(\det)$. 
\item For $\alpha,\beta \in \mathcal{D}$, $\pi(\alpha,\beta)$ is the principal series representation of $\GL{\F_p}$ attached to the pair of characters $\{\alpha,\beta\}$. 
\item By Definition \ref{left-cusp-rep-f}, we can attach to any $(f,\chi) \in \mathscr{C}$ a cuspidal representation of $\GL{\F_p}$ named $C_f$. 
\end{itemize}}

We fix, in this introduction, a number field $F \subset \C$ which is Galois over $\Q$ and contains:
\begin{itemize}[noitemsep,label=$-$] 
\item the $p(p^2-1)(p-1)$-th roots of unity (so that any representation of $\GL{\F_p}$ can be realized over $F$), 
\item the coefficients of every newform of weight $2$, level dividing $p^2$, and character of conductor dividing $p$,
\item for every newform $f \in \mathcal{S}_2(\Gamma_0(p^2))$ with complex multiplication (necessarily by the quadratic imaginary field $K := \Q(\sqrt{-p})$), for every prime number $\ell$ split in $K$, the roots of the polynomial $X^2-a_{\ell}(f)X+\ell$.   
\end{itemize}

\theoi[tate-1]{(Proposition \ref{tate-module-Pp-untwisted})
Let $\mathfrak{l} \subset \OO_F$ be a maximal ideal with residue characteristic $\ell$. There is an isomorphism of $(\mathbb{T} \otimes_{\Z} F_{\mathfrak{l}})[G_{\Q} \times \GL{\F_p}]$-modules:
\[\Tate{\ell}{P(p)} \otimes_{\Z_{\ell}} F_{\mathfrak{l}} \rar \bigoplus_{\substack{(f,\chi) \in \mathscr{S} \cup \mathscr{P}\\\psi \in \mathcal{D}}}{V_{f\otimes\psi,\mathfrak{l}} \otimes \pi(\psi,\psi\chi)}.\]
On the factor attached to $((f,\chi),\psi) \in (\mathscr{S}\cup\mathscr{P}) \times \mathcal{D}$, the Hecke operator $T_{n}$ (resp. $\langle n\rangle$) for $n \geq 1$ coprime to $p$ acts by multiplication by $a_{n}(f)$ (resp. $\chi(n)$). 
}

 The case of $J(p,p)$ (as an Abelian variety over $\Q$) also seems known, as a similar claim is made in \cite{Virdol}, but we keep to a classical language and our argument seems somewhat different. Since $J(p,p)(\C)$ is a complex torus, we can relate its singular cohomology to $H^0(J(p,p)(\C),\Omega^1)$ (for the left and right $\mathbb{T}[\GL{\F_p}]$-module structures). Using the results of the previous Chapter, we can determine $H_1(X(p,p)(\C),\C)$ as a $\mathbb{T}[\GL{\F_p}]$-module. Determining the structure of the Tate module of $J(p,p)$ then amounts to studying the action of $G_{\Q}$, which is given by the Eichler-Shimura relation (Corollary \ref{tate-eichler-shimura}). 

\theoi[tate-2]{(Corollary \ref{decomposition-xpp-tate}) Let $\mathfrak{l} \subset \OO_F$ be a maximal ideal with residue characteristic $\ell$. Then, as a $(\mathbb{T}\otimes F_{\mathfrak{l}})[G_{\Q} \times \GL{\F_p}]$-module, $\Tate{\ell}{J(p,p)} \otimes_{\Z_{\ell}} F_{\mathfrak{l}}$ is isomorphic to the direct sum of the following modules:
\begin{itemize}[noitemsep,label=\tiny$\bullet$]
\item for each $((f,\mathbf{1}),\psi) \in \mathscr{S} \times \mathcal{D}$, a copy of $V_{f \otimes \psi,\mathfrak{l}} \otimes \mrm{St}_{\psi}$, on which, for any $n \geq 1$ coprime to $p$, $T_{n}$ (resp. $nI_2$) acts by $a_{n}(f)\psi(n)$ (resp. $\psi(n)^2$), 
\item for each $((f,\chi),\psi) \in \mathscr{P} \times \mathcal{D}$ modulo the relation $((f,\chi),\psi) \sim ((\overline{f},\overline{\chi}),\chi\psi)$, a copy of $V_{f \otimes \psi,\mathfrak{l}} \otimes \pi(\psi,\psi\chi)$, on which on which, for any $n \geq 1$ coprime to $p$, $T_{n}$ (resp. $nI_2$) acts by $a_{n}(f)\psi(n)$ (resp. $\psi(n)^2\chi(n)$),  
\item for each $(f,\chi) \in \mathscr{C}$, a copy of $V_{f ,\mathfrak{l}} \otimes C_f$, on which, for any $n \geq 1$ coprime to $p$, $T_{n}$ (resp. $nI_2$) acts by $a_{n}(f)$ (resp. $\chi(n)$), 
\end{itemize}}

A consequence of this result is the description of the Tate module of any Galois twist of $J(p,p)$, also studied by Virdol (his result is reproduced in the general Introduction as Theorem \ref{Virdolthm}), but our classical statement seems more explicit. 

\cori[tate-3]{(Corollary \ref{decomposition-xrho}) Let $k$ be a field of characteristic distinct from $p$ and separable closure $k_s$. Let $F$ be a large enough number field, $\mathfrak{l} \subset \OO_L$ be a maximal ideal with residue characteristic $\ell$ invertible in $k$. Let $\rho: \mrm{Gal}(k_s/k) \rar \GL{\F_p}$ be a continuous group homomorphism. Then, as a $(\mathbb{T} \otimes_{\Z} F_{\mathfrak{l}})[G_{k}]$-module, $\Tate{\ell}{\mrm{Jac}(X(p,p)_{\rho})} \otimes_{\Z_{\ell}} F_{\mathfrak{l}}$ is isomorphic to the direct sum of
\begin{itemize}[label=\tiny$\bullet$,noitemsep]
\item for each $((f,\mathbf{1}),\psi) \in \mathscr{S} \times \mathcal{D}$, a copy of $\psi(\omega_p \det{\rho}) \otimes V_{f,\mathfrak{l}} \otimes \mrm{St}(\rho)$, on which, for any $n \geq 1$ coprime to $p$, $T_n$ (resp. $nI_2$) acts by $a_{n}(f)\psi(n)$ (resp. $\psi(n)^2$), 
\item for each $((f,\chi),\psi) \in \mathscr{P} \times \mathcal{D}$ modulo the relation $((f,\chi),\psi) \sim ((\overline{f},\overline{\chi}),\chi\psi)$, a copy of $\psi(\omega_p\det{\rho}) \otimes V_{f,\mathfrak{l}} \otimes (\pi(1,\chi)(\rho))$, on which, for any $n \geq 1$ coprime to $p$, $T_{n}$ (resp. $nI_2$) acts by $a_{n}(f)\psi(n)$ (resp. $\psi(n)^2\chi(n)$), 
\item for each $(f,\chi) \in \mathscr{C}$, a copy of $V_{f ,\mathfrak{l}} \otimes C_f(\rho)$, on which, for any $n \geq 1$ coprime to $p$, $T_{n}$ (resp. $nI_2$) acts by $a_{n}(f)$ (resp. $\chi(n)$).  
\end{itemize}
}

When we twist $X(p,p)$ by a Galois representation $\rho$ coming from the $p$-torsion from an elliptic curve $E$ over some field $k$ (with separable closure $k_s$ and of characteristic distinct from $p$), or more generally when $\det{\rho}$ is the inverse of the mod $p$ cyclotomic character, Proposition \ref{group-is-twist-polarized} implies that the twist $X(p,p)_{\rho}$ is the disjoint reunion of the $p-1$ smooth projective geometrically connected curves $X_{E[p],a}(p)$ over $k$ (where $a$ runs through $\F_p^{\times}$), or in general $X_{\rho,\xi}(p)$, when $\xi$ runs through $\mu_p^{\times}(k_s)$. The remainder of the section is dedicated to the computation of the Tate modules of the Jacobians of these curves. 

First, we have to restrict the Hecke algebra $\mathbb{T}$ acting on $J(p,p)$ (hence on $\mrm{Jac}(X(p,p)_{\rho})$) to a certain sub-algebra. Let $\mathbb{T}_1$ denote the subalgebra of $\mathbb{T}$ generated by the $(mI_2) \cdot T_n$ where $n,m \geq 1$ are integers such that $nm^2 \equiv 1\pmod{p}$. The endomorphism of $J(p,p)_{\rho}$ induced by $a \in \mathbb{T}_1$ factors as a direct sum of $a_{\xi} \in \mrm{End}(\mrm{Jac}(X_{\rho,\xi}(p)))$ (over $\xi \in \mu_p^{\times}(k_s)$).

\theoi[tate-4]{(Corollary \ref{decomposition-xrho-connected-general}) Let $\ell$ be a prime number and $\mathfrak{l}$ be a prime ideal of $\OO_F$ of residue characteristic $\ell$. Let $\rho: G_{\Q} \rar \GL{\F_p}$ such that $\det{\rho}=\omega_p^{-1}$. Let $\iota: \Qbar \rar k_s$ be a ring homomorphism, inducing a continuous group homomorphism $R_{\iota}: \mrm{Gal}(k_s/k) \rar \mrm{Gal}(\Qbar/\Q)$. 

For $\xi \in \mu_p^{\times}(k_s)$, $\Tate{\ell}{J_{\rho,\xi}(p)}\otimes_{\Z_{\ell}} F_{\mathfrak{l}}$ is isomorphic as a $(\mathbb{T}_1 \otimes F_{\mathfrak{l}})[G_{\Q}]$-module to the direct sum of the following factors:
\begin{itemize}[noitemsep,label=\tiny$\bullet$]
\item $V_{f,\mathfrak{l}} \otimes \mrm{St}(\rho)$ for every $(f,\chi) \in \mathscr{S}$, 
\item $V_{f,\mathfrak{l}} \otimes \pi(1,\chi)(\rho)$ for every $(f,\chi) \in \mathscr{P}$ modulo complex conjugation,
\item $V_{f,\mathfrak{l}} \otimes C_f(\rho)$ for every $(f,\chi) \in\mathscr{C}$ modulo twists by characters in $\mathcal{D}$ and such that $f$ does not have complex multiplication,
\item for every $(f,\chi) \in \mathscr{C}$ modulo twists by characters in $\mathcal{D}$ such that $f$ has complex multiplication, $\mrm{Ind}_{G_{\Q(\sqrt{-p})}}^{G_{\Q}}{\psi_{f,\mathfrak{l}} \otimes C^+(\rho_{|G_{\Q(\sqrt{-p})}})}$, where $\psi_{f,\mathfrak{l}}: G_{\Q(\sqrt{-p})} \rar F_{\mathfrak{l}}^{\times}$ is one of the CM characters attached to $f$ (which one is selected depends only on $\xi$) and $C^+$ is a cuspidal representation of $\F_p^{\times}\SL{\F_p}$ defined in Section \ref{splitting-sl2}. 
\end{itemize}
On every summand, the Hecke operator $(mI_2) \cdot T_n \in \mathbb{T}_1$ (for $n,m \geq 1$ such that $nm^2\equiv 1\pmod{p}$) acts by multiplication by $a_n(f)\chi(m)$. 
}

We stated the result for the base field $\Q$, but the body of this Chapter contains a similar version of this result which is valid for an arbitrary base field $k$. Let us briefly discuss its proof. 

When $\rho$ is surjective, this corresponds exactly to the decomposition into irreducible factors of the Tate module of $\mrm{Jac}(X(p,p)_{\rho})$ (by Corollary \ref{tate-3} and the results of Section \ref{irreducibility-factors}), so the proof is easier. In order to deal with general $\rho$, we need a stronger version of Theorem \ref{tate-4}, supplied by Corollary \ref{decomposition-xrho-connected-sl20}, which keeps track of the action of $\SL{\F_p}$. This action preserves the connected components, but is not defined over $\Q$ any more. 

The argument for any field $k$ of characteristic zero then works as follows: consider an elliptic curve $E_0/\Q$ with surjective attached Galois representation $\rho_0: G_{\Q} \rar \GL{\F_p}$ (given by a basis $(P,Q)$ of its $p$-torsion and $\xi_0=\langle P,\,Q\rangle_{E[p]}$). Then $X_{\rho,\xi}(p)$ is a twist of $X_{\rho_0,\xi_0}(p)_k$ by the Galois cocycle $\rho\rho_0^{-1}: \mrm{Gal}(k_s/k) \rar \SL{\F_p}$ (where a common choice of embedding $\Qbar \rar \overline{k}$ defines a homomorphism $\mrm{Gal}(\overline{k}/k) \rar G_{\Q}$ and identifies $\mu_p^{\times}(\Qbar)$ with $\mu_p^{\times}(\overline{k})$). By keeping track the action over $\Qbar$ of $\SL{\F_p}$ on all components on $X_{E_0[p],a}(p)$ (or equivalently on the $X_{\rho_0,\xi_0}(p)$), we are able to describe both the action of $G_{\Q}$ and that of $\SL{\F_p}$ on the Tate module of $X_{\rho,\xi}(p)$. 

For fields of positive characteristic, the same argument shows that it is enough to establish the stronger version of Theorem \ref{tate-4} for one representation $\rho_0$ over the prime field, and we choose one representation $\rho_0$ which lifts to characteristic zero, in which case the result is already known by specialization\footnote{I am grateful to Anna Cadoret for her answers to my questions regarding this theory.}. The existence of the curves $X_G(p)$ over rings of mixed characteristic is crucial to this part of the argument.

\section{The Tate module of $P(p)$}

\lem{There is an isomorphism of $\Z[\GL{\F_p} \times G_{\Q}]$-modules \[I: \Z[\GL{\F_p}] \otimes_{\Z[B]} (\Z[\F_p^{\times}] \otimes J_1(p)(\Qbar)) \rar P(p)(\Qbar),\] such that:
\begin{itemize}[noitemsep,label=$-$]
\item $M=\begin{pmatrix} a & \ast\\0 & d\end{pmatrix}$ acts on $\Z[\F_p^{\times}] \otimes J_1(p)(\Qbar)$ by $\underline{ad} \otimes \langle d\rangle$, where $\underline{ad} \in \mrm{End}(\Z[\F_p^{\times}])$ is the morphism $[b] \longmapsto [bad]$. 
\item $\sigma \in G_{\Q}$ acts on $\Z[\F_p^{\times}] \otimes J_1(p)(\Qbar)$ by the tensor product of $\underline{\omega_p(\sigma)}$ and its standard action. 
\item For any $M \in \GL{\F_p}$, $a \in \F_p^{\times}$, $f \in J_1(p)(\Qbar)$, $\ell \neq p$ prime, and $n \in \F_p^{\times}$, one has \[I(M \otimes [a] \otimes (T_{\ell}f)) = T_{\ell}(M \otimes [a] \otimes f),\quad I(M \otimes [a] \otimes (\langle n\rangle f)) = \langle n\rangle(M \otimes [a] \otimes f).\]
\end{itemize}}

\demo{By Proposition \ref{Pp-over-GL2}, it is enough to construct an isomorphism of $\Z[G_{\Q}]$-modules $\Z[\F_p^{\times}] \otimes J_1(p)(\Qbar) \rar J_1(p,\mu_p)(\Qbar)$ respecting the $T_{\ell}$, the $\langle n\rangle$, and mapping $\underline{a} \in \mrm{End}(\Z[\F_p^{\times}])$ to $\underline{a} \in \mrm{End}(\mu_p^{\times})$. 

By Proposition \ref{jacobian-relative-curve-exists}, we have an identification of functors $J_1(p,\mu_p)(-) \simeq J_1(p)(- \times_{\Z[1/p]} \mu_p^{\times})$ (valued in Abelian groups), where $\mu_p^{\times}$ is the affine scheme with coordinate ring $\Z[p^{-1},\zeta_p]$. This identification is natural, so it commutes with Hecke and diamond operators. Moreover, by Proposition \ref{automorphism-global-functions}, the following diagram commutes for any $\Q$-algebra $R$ and any $a \in \F_p^{\times}$:

\[
\begin{tikzcd}[ampersand replacement=\&]
J_1(p,\mu_p)(R) \arrow{d}{\underline{a}} \arrow{r}{\simeq} \& J_1(p)(R \otimes \Q(\mu_p))\arrow{d}{\underline{a}^{-1}}\\
J_1(p,\mu_p)(R) \arrow{r}{\simeq} \& J_1(p)(R \otimes \Q(\mu_p))
\end{tikzcd}
\]

Therefore, it is enough to find an isomorphism $\iota: \Z[\F_p^{\times}] \rar J_1(p)(\Qbar) \rar J_1(p)(\Qbar \otimes \Q(\mu_p))$ such that the following diagram commutes for each $a \in \F_p^{\times}$ and each $\sigma \in G_{\Q}$:

\[
\begin{tikzcd}[ampersand replacement=\&]
J_1(p)(\Qbar \otimes \Q(\mu_p)) \arrow{d}{\sigma \otimes \underline{a}^{-1}} \arrow{r}{\iota} \& \Z[\F_p^{\times}] \otimes J_1(p)(\Qbar)\arrow{d}{\underline{a\omega_p(\sigma)} \otimes \sigma}\\
J_1(p)(\Qbar \otimes \Q(\mu_p)) \arrow{r}{\iota} \& \Z[\F_p^{\times}] \otimes J_1(p)(\Qbar)
\end{tikzcd}
\]

To define this automorphism, we will work in slightly greater generality so as to make the argument clearer. Write $K=\Q(\mu_p)$, $\Gamma=\mrm{Gal}(K/\Q)$, $J=J_1(p)_{\Q}$ and fix an embedding $j: K \rar \Qbar$. This defines a group homomorphism of $q: G_{\Q}\rar \Gamma$ (when $K$ is a general finite Galois extension of $\Q$, $q$ is only defined up to conjugacy; here, it is absolute, but it is easier to pretend for a while that it isn't). Let $\pi: \Qbar \otimes K \rar \Qbar^\Gamma$ such that for each $g \in \Gamma$, $(x,y) \in \Qbar \times K$, $\pi(x \otimes y)(g)=x\cdot j(g(y))$. Then we have the following relations:

\[\pi(x \otimes g(y))(h)=\pi(x \otimes y)(hg),\, \pi(\sigma(x) \otimes y)(g) = \sigma(\pi(x \otimes y)(q(\sigma)^{-1}g)).\]

There is an obvious isomorphism of Abelian groups $\lambda: J(\Qbar^{\Gamma}) \rar \Z[\Gamma] \otimes J(\Qbar)$ mapping a point $x \in J(\Qbar^{\Gamma})$ to the sum over $g \in \Gamma$ of the $[g] \otimes p_g(x)$, where $p_g$ denotes the projection $\Qbar^{\Gamma} \rar \Qbar$. For any $\sigma \in G_{\Q}$ and $\gamma \in \Gamma$, for any $x \in J(\Qbar \otimes K)$, one has 

\begin{align*}
\lambda(\pi((\sigma \otimes \gamma)(x))) &= \sum_{h \in \Gamma}{[h] \otimes p_h(\pi(\sigma \otimes \gamma)(x))} = \sum_{h \in \Gamma}{[h] \otimes \sigma(p_{q(\sigma)^{-1}h\gamma}(\pi(x)))} \\
&= \sum_{h \in \Gamma}{[q(\sigma)h\gamma^{-1}] \otimes \sigma(p_h(\pi(x)))} = (\mathbf{L}_{q(\sigma)}\mathbf{R}_{\gamma^{-1}} \otimes \sigma)(\lambda(\pi(x))),
\end{align*}

where $\mathbf{L},\mathbf{R}$ denote the left and right translations on $\Z[\Gamma]$, which pairwise commute. Then $\iota = \lambda \circ \pi$ satisfies the conditions.
}

\prop[tate-module-Pp-untwisted]{Let $\ell$ be a prime number and $F$ be a finite extension of $\Q_{\ell}$ containing all the $(p-1)$-th roots of unity and the coefficients of every newform in $\mathcal{S}_2(\Gamma_1(p))$. Let $\mathfrak{l}$ be the maximal ideal of $\OO_F$. We then have an isomorphism of $(\mathbb{T} \otimes_{\Z} F)[G_{\Q} \times \GL{\F_p}]$-modules:
\[\Tate{\ell}{P(p)} \otimes_{\Z_{\ell}} F \rar \bigoplus_{\substack{(f,\chi) \in \mathscr{S} \cup \mathscr{P}\\\psi \in \mathcal{D}}}{V_{f\otimes\psi,\mathfrak{l}} \otimes \pi(\psi,\psi\chi)}.\]
On the factor attached to $((f,\chi),\psi)$, $T_{q}$ (resp. $\langle n\rangle$) acts by $a_{q}(f)$ (resp. $\chi(n)$). 
}

\demo{This follows from the previous Lemma and the description of the Tate module of $J_1(p)$ as given by \cite[Theorem 5.12]{Ribet-Stein} (note the discussion of \emph{loc.cit.} below Definition 5.6; \emph{loc.cit.} does not discuss basic facts about the action of Hecke operators on classical modular forms, which can be found in \cite[Chapter 5]{DS}).}

\section{The singular homology of $J(p,p)$}

\prop[from-omega1-to-singular]{Let $A$ be a (not necessarily commutative) ring acting on some complex torus $J$. Assume that $H^0(J,\Omega^1) \simeq \bigoplus_{i}{M_i}$, where each $M_i$ is a $\C$-vector space on which $A$ acts on the right. 
Then \[H_1(J,\C) \simeq \bigoplus_i{\mrm{Hom}_{\C}(M_i,\C) \oplus \mrm{Hom}_{\C}(\overline{M_i},\C)},\] 
where $\overline{M}:=M \otimes_{\C} \C$, $\C \rar \C$ being the complex conjugation. }

\demo{We can write $J=V/\Lambda$, where $V$ is a complex vector space and $\Lambda \subset V$ is a closed discrete subgroup with $V \simeq \Lambda \otimes_{\Z} \R$. Then $A$ acts on $\Lambda$ and $V$ by group homomorphisms and $\C$-linear maps, and one has, as $A \otimes \C$-modules, one has the following isomorphisms:
\begin{align*}
H_1(J,\C) &\simeq \Lambda \otimes_{\Z} \C \simeq (\Lambda \otimes_{\Z} \R)\otimes_{\R} \C \simeq V \otimes_{\R} \C \simeq V \oplus \overline{V}\\
 H^0(J,\Omega^1) &\simeq \mrm{Hom}_{\C}(V,\C),
\end{align*}
The conclusion is then formal. 
  }

\cor[decomposition-singular-xpp]{The $(\mathbb{T} \otimes \C)[\GL{\F_p}]$-module $H_1(X(p,p)(\C),\C)$ decomposes as the direct sum of two copies of the following $(\mathbb{T} \otimes \C)[\GL{\F_p}]$-modules:
\begin{itemize}[label=$-$,noitemsep]
\item for each $((f,\mathbf{1}),\psi) \in \mathscr{S}\times \mathcal{D}$, a copy of $\mrm{St}_{\psi}$, where $T_{\ell}$ acts by $\psi(\ell)a_{\ell}(f)$,
\item for each $((f,\chi),\psi) \in \mathscr{P} \times \mathcal{D}$ modulo the relation $(f,\chi,\psi) \sim (\overline{f},\overline{\chi},\psi\chi)$, a copy of $\pi(\psi,\psi\chi)$, where $T_{\ell}$ acts by $a_{\ell}(f)\psi(\ell)$, 
\item for each $(f,\chi) \in \mathscr{C}$, a module isomorphic to $C_f$, where $T_{\ell}$ acts by $a_{\ell}(f)$.
\end{itemize}
}

\demo{By Corollary \ref{decomposition-hecke-xpp}, for each $((f,\chi),\psi) \in (\mathscr{S} \cup \mathscr{P}) \times \mathcal{D}$ (resp. $(f,\chi) \in\mathscr{C}$), there is a ring homomorphism $\mu_{f,\psi}: \mathbb{T} \rar \C$ (resp. $\mu_f: \mathbb{T} \rar \C$) such that 
\begin{align*}
\mu_{f,\psi}(T_{\ell}) &= a_{\ell}(f)\psi(\ell) && \text{( resp. }\mu_f(T_{\ell})=a_{\ell}(f) \text{)}\\
\mu_{f,\psi}(nI_2)&=\chi(n)\psi^2(n) && \text{( resp. }\mu_f(nI_2)=\chi(n) \text{)}
\end{align*}

For each $((f,\chi),\psi) \in (\mathscr{S} \cup \mathscr{P}) \times \mathcal{D}$ (resp. $(f,\chi) \in \mathscr{C}$), let $M_{f,\psi}$ (resp. $M_f$) be the $(\mathbb{T}\otimes \C)[\GL{\F_p}]$ which is isomorphic over $\C[\GL{\F_p}]$ to $\mrm{St}_{\psi}$ if $\chi=\mathbf{1}$ and $\pi(\psi,\chi\psi)$ otherwise (resp. $C_f$) and on which $\mathbb{T}$ acts through the character $\mu_{f,\psi}$ (resp. $\mu_f$). In particular, $M_{f,\psi}$ and $M_{\overline{f},\chi\psi}$ are isomorphic, and the complex conjugate of $M_{f,\psi}$ (resp. $M_f$) is $M_{\overline{f},\overline{\psi}}$ (resp. $M_{\overline{f}}$). 

By Corollary \ref{decomposition-hecke-xpp}, $H^0(X(p,p)(\C),\Omega^1) \simeq H^0(J(p,p)(\C),\Omega^1)$ is isomorphic to \[\left(\bigoplus_{\psi \in \mathcal{D}}{\bigoplus^{'}_{(f,\chi) \in \mathscr{S} \cup \mathscr{P}}{\mrm{Hom}(M_{f,\psi},\C)}}\right) \oplus \bigoplus_{(f,\chi) \in \mathscr{C}}{\mrm{Hom}(M_f,\C)}\] as a right $(\mathbb{T} \otimes \C)[\GL{\F_p}]$-module, where the $'$ means that for each $(f,\chi) \in \mathscr{P}$, we choose exactly one element in $\{(f,\chi),(\overline{f},\overline{\chi})\}$ to count in the sum (this choice being consistent for all $\psi$). The conclusion follows from Proposition \ref{from-omega1-to-singular}. 
}

\section{The Tate module of twists of $J(p,p)$}

\lem[from-vanishing-to-characteristic]{Let $f \in \mathcal{S}_k(\Gamma_1(N))$ be a normalized newform with weight $k \geq 2$ and character $\chi$. Let $\Q(\chi) \subset F \subset \C$ be a number field containing the coefficients of $f$, $\ell$ a prime number and $\mathfrak{l} \subset \OO_F$ a prime ideal of residue characteristic $\ell$. Let $V$ be a continuous two-dimensional representation of $G_{\Q}$ with coefficients in $F_{\mathfrak{l}}$ such that 
\begin{itemize}[noitemsep,label=$-$]
\item $V$ is unramified at all but finitely many primes,
\item for all but finitely many primes $q$, the endomorphism $\Fr_{q}^2-a_{q}(f)\Fr_{q}+q^{k-1}\chi(q)$ of $V$ vanishes, where $\Fr_{q}$ is the arithmetic Frobenius.
\end{itemize}
 Then $V$ is isomorphic to $V_{f,\mathfrak{l}}$.}

\demo{For each $\sigma \in G_{\Q}$, denote by $P_{\sigma}$ (resp. $Q_{\sigma}$) be the characteristic polynomial of $\sigma \in \mrm{End}_{F_{\mathfrak{l}}}(V_{f,\mathfrak{l}})$ (resp. $\sigma \in \mrm{End}_{F_{\mathfrak{l}}}(V)$). Given $\sigma \in G_{\Q}$, we denote by $V(\sigma)$ the endomorphism of $V$ induced by $\sigma$. 

Now, $\sigma \longmapsto P_{\sigma}$ and $\sigma \longmapsto Q_{\sigma}$ are continuous maps. Whenever $\sigma$ is an arithmetic Frobenius, the assumption implies that $P_{\sigma}(V(\sigma))=0$: by Cebotarev, $P_{\sigma}(V(\sigma))=0$ for each $\sigma \in G_{\Q}$. 

Let $S \leq G_{\Q}$ be the subgroup of $\sigma$ such that $V(\sigma)$ is scalar, and $T \leq G_{\Q}$ be the subgroup of elements $\sigma$ such that $\det(\sigma \mid V_{f,\mathfrak{l}}) = \det(\sigma \mid V)$. If $\sigma \in G_{\Q} \backslash S$, the gcd $R$ of $P_{\sigma}$ and $Q_{\sigma}$ is such that $R(V(\sigma))=0$. Since $V(\sigma)$ is not scalar, $R$ must have degree two, so $P_{\sigma}=Q_{\sigma}$ and $\sigma \in T$. Hence $G_{\Q}$ is the reunion of the two closed subgroups $S$ and $T$, so $G_{\Q} \in \{S,T\}$. 

Let us prove that $G_{\Q} \neq S$. Assume, for the sake of contradiction, that $G_{\Q}=S$: in this case, there is a character $\psi: G_{\Q} \rar F_{\mathfrak{l}}^{\times}$ corresponding to $V$, so that, for all $\sigma \in G_{\Q}$, $\psi(\sigma)$ is a root of $P_{\sigma}$. Let $V':= V_{f,\mathfrak{l}} \otimes \psi^{-1}$: then $V'$ is irreducible and $1$ is an eigenvalue of every $V'(\sigma)$ for $\sigma \in G_{\Q}$. Let $\Q^{ab}$ denote the maximal abelian extension of $\Q$: every $\sigma \in G_{\Q^{ab}}$ acts on $V'$ with determinant $1$, so acts unipotently on $V'$.  

Let us prove that $V'$ has a line stable under $G_{\Q^{ab}}$. If it is not the case, there exists $\sigma,\tau \in G_{\Q^{ab}}$ such that $\ker(V'(\sigma)-\mrm{id}),\ker(V'(\tau)-\mrm{id})$ are distinct lines. Then, in a suitable basis of $V'$, the matrices of $\sigma,\tau$ are $\begin{pmatrix}1 & a\\0 & 1\end{pmatrix},\begin{pmatrix}1 & 0\\b & 1\end{pmatrix}$ for some $a,b \in F_{\mathfrak{l}}^{\times}$. Then the trace of $V'(\sigma\tau)$ is $2+ab$, so $\sigma\tau \in G_{\Q^{ab}}$ does not act unipotently on $V'$, whence a contradiction. 

So $V'$ contains a line stable under $G_{\Q^{ab}}$: since each $\sigma \in G_{\Q^{ab}}$ acts unipotently on $V'$, $(V')^{G_{\Q^{ab}}}$ is nonzero. Since $V'$ is irreducible, one has $(V')^{G_{\Q^{ab}}}=V'$, so the action of $G_{\Q}$ on $V'$ is abelian: this contradicts our original assumption, namely that $G_{\Q}=S$.  

Therefore, $G_{\Q}=T$. Let $\sigma \in G_{\Q}$, then $(P_{\sigma}-Q_{\sigma}) \in K\cdot X$ vanishes when evaluated at $V(\sigma)$, so $P_{\sigma}=Q_{\sigma}$. Therefore, since $V_{f,\mathfrak{l}}$ is absolutely irreducible, one has 
\[\dim_{F_{\mathfrak{l}}}\mrm{Hom}_{F_{\mathfrak{l}}}(V,V_{f,\mathfrak{l}}) = \dim_{\overline{F_{\mathfrak{l}}}}\mrm{Hom}_{\overline{F_{\mathfrak{l}}}}(V \otimes_{F_{\mathfrak{l}}} \overline{F_{\mathfrak{l}}},V_{f,\mathfrak{l}} \otimes_{F_{\mathfrak{l}}} \overline{F_{\mathfrak{l}}}) \geq 1,\]
so there is a nonzero morphism of $F_{\mathfrak{l}}[G_{\Q}]$-modules $f: V \rar V_{f,\mathfrak{l}}$. Since $V_{f,\mathfrak{l}}$ is irreducible, $f$ is an isomorphism and we are done. 
}

\medskip

\nott{From now on, $F\subset \C$ denotes a number field which is Galois over $\Q$ and contains the $p(p-1)^2(p+1)$-th roots of unity and the Fourier coefficients of the newforms in $\mathcal{S}_2(\Gamma_1(p))$ or $\mathcal{S}_2(\Gamma_1(p) \cap \Gamma_0(p^2))$. }

\medskip

\prop[decomposition-xpp-twokinds]{Let $\ell$ be a prime number and $\mathfrak{l}$ be a prime ideal of $\OO_F$ dividing $(\ell)$. Then $\Tate{\ell}{J(p,p)} \otimes_{\Z_{\ell}} F_{\mathfrak{l}}$ is the direct sum of its eigenspaces for $\mathbb{T}$, which are stable under $G_{\Q} \times \GL{\F_p}$, and are of the following form:
\begin{itemize}[noitemsep,label=\tiny$\bullet$]
\item For each $((f,\mathbf{1}),\psi) \in \mathscr{S} \times \mathcal{D}$, a $2p$-dimensional $F_{\mathfrak{l}}$-vector space on which any $T_{n}$ acts by $a_{n}(f)\psi(n)$ (with $n \geq 1$ prime to $p$), which is isomorphic over $F_{\mathfrak{l}}[\GL{\F_p}]$ to the direct sum of two copies of $\mrm{St}_{\psi}$, and over $F_{\mathfrak{l}}[G_{\Q}]$ to the sum of $p$ copies of $V_{f \otimes \psi,\mathfrak{l}}$. 
\item For each $((f,\chi),\psi) \in \mathscr{P} \times \mathcal{D}$ modulo the relation $(f,\chi,\psi) \sim (\overline{f},\overline{\chi},\chi\psi)$, a $2(p+1)$-dimensional $F_{\mathfrak{l}}$-vector space where any $T_{n}$ acts by $a_{n}(f)\psi(n)$ (with $n\geq 1$ coprime to $p$) which is isomorphic over $F_{\mathfrak{l}}[\GL{\F_p}]$ to the direct sum of two copies of $\pi(\psi,\psi\chi)$ and over $F_{\mathfrak{l}}[G_{\Q}]$ to the sum of $p+1$ copies of $V_{f \otimes \psi,\mathfrak{l}}$. 
\item For each $f \in \mathscr{C}$, a $2(p-1)$-dimensional vector space where any $T_{n}$ acts by $a_{n}(f)$ (with $n \geq 1$ coprime to $p$) which is isomorphic over $F_{\mathfrak{l}}[\GL{\F_p}]$ to the direct sum of two copies of $C_f$, and over $F_{\mathfrak{l}}[G_{\Q}]$ to the sum of $p-1$ copies of $V_{f,\mathfrak{l}}$. 
\end{itemize}
}

\demo{As a $(\mathbb{T} \otimes F_{\mathfrak{l}})[\GL{\F_p}]$-module, $\Tate{\ell}{J(p,p)} \otimes_{\Z_{\ell}} F_{\mathfrak{l}}$ is isomorphic to $H_1(J(p,p)(\C),F_{\mathfrak{l}})$. Now, the $\mu_{f,\psi}, \mu_f$ (resp. $M_{f,\psi},M_f$) defined in the proof of Corollary \ref{decomposition-singular-xpp} are in fact characters $\mathbb{T} \rar F$ (resp. can be defined as modules over $F[\GL{\F_p}]$). An easy induction shows that for any $n \geq 1$ coprime to $p$, $\mu_f(T_n)=a_n(f)$ (resp. $\mu_{f,\psi}(n)=a_n(f)\psi(n)$).  

Let us fix a system of representatives $\mathscr{P}'$ for $\mathscr{S} \cup \mathscr{P}$ modulo the complex conjugation. Then $H_1(J(p,p)(\C),F)$ and 
\[\left(\bigoplus_{((f,\chi),\psi) \in \mathscr{P}' \times \mathcal{D}}{M_{f,\psi}} \oplus \bigoplus_{(f,\chi) \in \mathscr{C}}{M_{f}}\right)^{\oplus 2}\] are two $(\mathbb{T}\otimes F)[\GL{\F_p}]$-modules that are isomorphic after the base change $F \rar \C$, so they are isomorphic. Therefore, one has the isomorphism of $(\mathbb{T}\otimes F_{\mathfrak{l}})[\GL{\F_p}]$-modules

\[\Tate{\ell}{J(p,p)} \otimes_{\Z_{\ell}} F_{\mathfrak{l}} \simeq \left(\bigoplus_{((f,\chi),\psi) \in \mathscr{P}' \times \mathcal{D}}{M_{f,\psi}^{\oplus 2}} \oplus \bigoplus_{(f,\chi) \in \mathscr{C}}{M_{f}^{\oplus 2}}\right)\otimes_F F_{\mathfrak{l}}.\] 

Since every factor is the eigenspace of some character $\mathbb{T} \rar F$, every summand in this decomposition is stable under $G_{\Q}$.  

Let $(f,\chi) \in \mathscr{C}$. The eigenspace $T[f]$ of $T := \Tate{\ell}{J(p,p)} \otimes_{\Z_{\ell}} F_{\mathfrak{l}}$ for the character $\mu_f$ has dimension $2(p-1)$ and is endowed with a continuous action of $\GL{\F_p} \times G_{\Q}$; it is isomorphic over $F_{\mathfrak{l}}[\GL{\F_p}]$ to two copies of $C_f$. By Proposition \ref{cusp-on-diag}, for every character $\psi \in \mathcal{D}$, the submodule $T[f,\psi_D]$ of $T[f]$ on which each $\Delta_{a,b}$ (for $a,b \in \F_p^{\times}$) acts by $\chi(a)\psi(a)\overline{\psi(b)}$ is two-dimensional, stable under $G_{\Q}$, and $T[f]$ is the direct sum of the $T[f,\psi_D]$. 

Now, $T[f,\psi_D]$ is a continuous two-dimensional representation of $G_{\Q}$ over $F_{\mathfrak{l}}$. By the Eichler-Shimura relation (Corollary \ref{tate-eichler-shimura}), for any prime $q \nmid p\ell$, the polynomial $X^2-a_{q}(f)X+q\chi(q)$ vanishes at the arithmetic Frobenius $\Fr_{q}$. By Lemma \ref{from-vanishing-to-characteristic}, 
$T[f,\psi_D]$ is isomorphic to $V_{f,\mathfrak{l}}$ as a $F_{\mathfrak{l}}[G_{\Q}]$-module. 

Let now $((f,\chi),\psi) \in \mathscr{P}' \times \mathcal{D}$. The eigenspace $T[f,\psi]$ of $T$ for the character $\mu_{f,\psi}$ has dimension $2p$ (or $2(p+1)$ if $\chi \neq \mathbf{1}$) and is endowed with a continuous action of $\GL{\F_p} \times G_{\Q}$; it is isomorphic over $F_{\mathfrak{l}}[\GL{\F_p}]$ to the direct sum of two copies of $\mrm{St}_{\psi}$ (or $\pi(\psi,\psi\chi)$ if $\chi \neq \mathbf{1}$). By the Eichler-Shimura relation (Corollary \ref{tate-eichler-shimura}), $\Fr_{q}^2-a_{q}(f)\psi(q)\Fr_{q}+q\chi(q)\psi(q)^2=0$ for each prime $q \nmid p\ell$, so, by Lemma \ref{from-vanishing-to-characteristic}, we are done if we can show that $T[f,\psi]$ is the direct sum of two-dimensional representations of $G_{\Q}$. 

In the case of $\mrm{St}_{\psi}$, by Lemma \ref{steinberg-on-borel}, we can take the eigenspaces of $T[f,\psi]$ for the action of $U=\begin{pmatrix}1 & 1\\0&1\end{pmatrix}$. 

In the case of $\pi(\psi,\psi\chi)$, let $g \in \GL{\F_p}$ be the generator of a Cartan subgroup. For any $1 \leq k \leq p$, $g^k$ is diagonalizable with eigenvalues in $\F_{p^2} \backslash \F_p$, so is not conjugate to any element of $B$: by \cite[\S 7.2 Proposition 20]{SerreLinReps}, it acts on $\pi(\psi,\psi\chi)$ with trace $0$. On the other hand, $g^{p+1}$ is a scalar matrix equal to $\lambda I_2$ for some $\lambda \in \F_p^{\times}$, hence the eigenvalues of the action of $g$ on $\pi(\psi,\psi\chi)$ are exactly the $(p+1)$-th roots of $(\chi\psi^2)(\lambda)$, each one with multiplicity one. Therefore, $T[f,\psi]$ is the direct sum of the eigenspaces for the action of $g$, which all have dimension two and are stable under $G_{\Q}$. 
}

\medskip

\lem[two-actions-to-product]{Let $V$ be a finite-dimensional vector space over some non-Archimedean local field $K$, endowed with continuous $K$-linear actions of two profinite groups $G,H$ such that 
\begin{itemize}[noitemsep,label=$-$]
\item the actions of $G$ and $H$ commute,
\item there exists finite-dimensional $K$-vector spaces $A,B$ endowed with $K$-linear continuous actions of $G$ and $H$ respectively, with $A$ (resp. $B$) being absolutely irreducible over $G$ (resp. $H$), such that $V_{|G}$ (resp. $V_{|H}$) is isomorphic to a direct sum of copies of $A$ (resp. $B$).
\end{itemize} 
Then there is an isomorphism of $K[G \times H]$-modules $V \simeq (A \otimes B)^{\oplus d}$ for some $d \geq 1$.}

\demo{Let $W \subset V$ be a sub-$K[H]$-module isomorphic to $B$. Then $U=\mrm{Hom}_{K[H]}(W,V)$ is a finite-dimensional $K$-vector space on which $G$ acts (through its action on $V$). We have a natural evaluation morphism $\epsilon: U \otimes_K W \rar V$ which is $K[G \times H]$-linear. Since $V_{|H}$ is a direct sum of copies of $W$, $\epsilon$ is surjective. 

Since $B$ is absolutely irreducible, it satisfies Schur's Lemma. Therefore, $\dim_K{U} = \frac{\dim{V}}{\dim{W}}$, so $\epsilon$ is a linear surjective morphism between finite dimensional $K$-vector spaces, so it is an isomorphism.

Now, it is enough to show that $U$ is isomorphic as a $K[G]$-module to a direct sum of copies of $A$. This is direct, since $U \otimes W_{|G} \simeq V_{|G}$, and $W_{|G}$ is trivial, while $V_{|G}$ is a direct sum of copies of the absolutely irreducible $K[G]$-module $A$. 
}

\medskip

\cor[decomposition-xpp-tate]{Let $\ell$ be a prime number and $\mathfrak{l}$ be a prime ideal of $\OO_F$ with residue characteristic $\ell$. Then $\Tate{\ell}{J(p,p)} \otimes_{\Z_{\ell}} F_{\mathfrak{l}}$ is isomorphic to the direct sum of the following $(\mathbb{T}\otimes F_{\mathfrak{l}})[G_{\Q} \times \GL{\F_p}]$-submodules:
\begin{itemize}[noitemsep,label=\tiny$\bullet$]
\item for each $((f,\mathbf{1}),\psi) \in \mathscr{S} \times \mathscr{D}$, a copy of $V_{f \otimes \psi,\mathfrak{l}} \otimes \mrm{St}_{\psi}$, on which, for any $n \geq 1$ coprime to $p$, $T_{n}$ (resp. $nI_2$) acts by $a_{n}(f)\psi(n)$ (resp. $\psi(n)^2$), 
\item for each $((f,\chi),\psi) \in \mathscr{P} \times \mathscr{D}$ modulo the relation $((f,\chi),\psi) \sim ((\overline{f},\overline{\chi}),\chi\psi)$, a copy of $V_{f \otimes \psi,\mathfrak{l}} \otimes \pi(\psi,\psi\chi)$, on which on which, for any $n \geq 1$ coprime to $p$, $T_{n}$ (resp. $nI_2$) acts by $a_{n}(f)\psi(n)$ (resp. $\psi(n)^2\chi(n)$),  
\item for each $(f,\chi) \in \mathscr{C}$, a copy of $V_{f ,\mathfrak{l}} \otimes C_f$, on which, for any $n \geq 1$ coprime to $p$, $T_{n}$ (resp. $nI_2$) acts by $a_{n}(f)$ (resp. $\chi(n)$), 
\end{itemize}}

\rem{It seems interesting to study the integral behavior of the decomposition when $\ell=p$. 

More precisely, one can consider the following situation. Let $\mathfrak{p}$ be a prime ideal of $\OO_F$ dividing $p$. Fix $((f,\chi),\psi) \in (\mathscr{S} \cup \mathscr{P}) \times \mathcal{D}$: the largest quotient $\mathcal{T}_{f,\psi}$ of $\Tate{p}{J(p,p)} \otimes_{\Z_p} \OO_{F_{\mathfrak{p}}}$ on which $\mathbb{T}$ acts by the character $\mu_{f,\psi}$ is a free $\OO_{F_{\mathfrak{p}}}$-module of rank $2p$ (or $2(p+1)$ when $\chi$ is nontrivial), which is (after tensoring with $\Q$) the direct sum of two copies of $\mrm{St}_{\psi}$ (or $\pi(\psi,\psi\chi)$). But what is its isomorphism class as a $\OO_{F_{\mathfrak{p}}}$-module? 

Fix a uniformizer $\varpi$ for $F_{\mathfrak{p}}$. When $\chi$ is trivial, $p$ is not an Eisenstein prime by \cite[II (9.7), (14.2)]{MazurY1}, so $T_{f,\mathfrak{p}}/\varpi$ is absolutely irreducible. Therefore, the isomorphism class of $T_{f,\mathfrak{p}}$ as a $\OO_{F_{\mathfrak{p}}}[G_{\Q}]$-module is well-defined. Similarly, $\mrm{St}_{\psi}/\varpi$ is absolutely irreducible, so there is a unique (up to isomorphism) $\OO_{F_{\mathfrak{p}}}[\GL{\F_p}]$-module which is free of rank $p$ over $\OO_{F_{\mathfrak{p}}}$ and isomorphic to $\mrm{St}_{\psi}$ over $F_{\mathfrak{p}}$. 

Using the absolute irreducibilities on both sides, one can check that $\mathcal{T}_{f,\psi}$ is, over $\OO_{F_{\mathfrak{p}}}[\GL{\F_p}]$, the direct sum of two copies of $\mrm{St}_{\psi}$. The argument of Lemma \ref{two-actions-to-product} can be adapted (using the residual irreducibility of the two representations) so that its conclusion still holds, that is, $\mathcal{T}_{f,\psi} \simeq \mrm{St}_{\psi} \otimes T_{f\otimes \psi,\mathfrak{p}}$. 

When $\chi$ is not trivial, the argument does not work any more, because the $\Z_p[\GL{\F_p}]$-module $\pi(\psi,\psi\chi)$ is not residually irreducible, nor is it isomorphic to $\pi(\psi\chi,\psi)$ (while the two modules are isomorphic over $\Q_p[G]$), and $T_{f,\mathfrak{p}}/\varpi$ need not be residually irreducible (or $T_{f,\psi}$ need not be uniquely defined up to $\OO_{F_{\mathfrak{p}}}[G_{\Q}]$-isomorphism). 

However, when $T_{f,\mathfrak{p}}/\varpi$ is absolutely irreducible (so that $T_{f,\mathfrak{p}}$ defines a unique isomorphism class of $\OO_{F_{\mathfrak{p}}}[G_{\Q}]$-module), we can see by decomposing according to the eigenspaces of a nonsplit Cartan subgroup $C$ (as in the proof of Proposition \ref{decomposition-xpp-twokinds}) that $\mathcal{T}_{f,\psi}$ is a direct sum of $p+1$ copies of $T_{f,\mathfrak{p}}$. This is then enough to prove (with an adaptation of the argument of Lemma \ref{two-actions-to-product}) that $\mathcal{T}_{f,\psi} \simeq T_{f,\psi} \otimes \pi'$, where $\pi'$ is a $\OO_{F_{\mathfrak{p}}}[\GL{\F_p}]$-module isomorphic to either $\pi(\psi,\psi\chi)$ or $\pi(\psi\chi,\psi)$. 

Since $\pi(\psi\chi,\psi)_{\F_p},\pi(\psi,\psi\chi)_{\F_p}$ have the same semi-simplification, but are defined by different (nonsplit) exact sequences (see Appendix \ref{reps-gl2-carac-p}), one way to decide between $\pi(\psi,\psi\chi)_{\OO_K}$ and $\pi(\psi\chi,\psi)_{\OO_K}$ would be to describe explicitly the image in singular homology of the sum of the $u_x$: $H_1(X(p,p)(\C),\F_p) \rar \bigoplus_{x \in S}{H_1(X_1(p)(\C) \times \F_p^{\times},\F_p)}$, but this does not seem straightforward. 

The case of cuspidal representations adds another difficulty (since $T_{f,\mathfrak{p}}$ also needs not be residually irreducible in this case, so this notation may refer to several non-isomorphic $\OO_{F_{\mathfrak{p}}}[\GL{\F_p}]$-modules): given $(f,\chi) \in \mathscr{C}$, it is not clear how many isomorphism classes of $\OO_{F_{\mathfrak{p}}}[\GL{\F_p}]$-modules $M$ which are free of finite rank over $\OO_{F_{\mathfrak{p}}}$ and such that $M \otimes_{\Z} \Q \simeq C_f$.   
}

\bigskip

\lem[from-qbar-to-ks]{Let $k$ be a field with characteristic $q$ and separable closure $k_s$. Let $\Q_{\{q\}'}$ be the maximal extension of $\Q$ which is unramified at $q$ (so $\Q_{\{0\}'}=\Qbar$), and $\OO_{\Q_{\{q\}'}}$ be the ring of algebraic integers of $\Q_{\{q\}'}$. 
\begin{itemize}[noitemsep,label=$-$]
\item There exists a ring homomorphism $\iota: \OO_{\Q_{\{q\}'}} \rar k_s$. 
\item Any ring homomorphism $\iota: \OO_{\Q_{\{q\}'}} \rar k_s$ defines a continuous group homomorphism $R_{\iota}: \mrm{Gal}(k_s/k) \rar \mrm{Gal}(\Q_{\{q\}'}/\Q)$. When $q>0$, $\ker{\iota}$ is a maximal ideal, and the image of $R_{\iota}$ is its decomposition subgroup. 
\item If $\iota': \OO_{\Q_{\{q\}'}} \rar k_s$ is another ring homomorphism, there exists a $\sigma \in \mrm{Aut}(\Q_{\{q\}'}/\Q)$ such that $\iota'=\iota \circ \sigma$, and $R_{\iota'}=\sigma^{-1}R_{\iota}\sigma$.
\end{itemize}}

\demo{When $q = 0$, the existence of an embedding $\Qbar \rar k_s$ follows from classical Galois theory. When $q>0$, there exists a maximal ideal $\mathfrak{q} \in \OO_{\Q_{\{q\}'}}$ with residue characteristic $q$. Since $\OO_{\Q_{\{q\}'}}$ contains all the roots of unity with order prime to $q$ and is integral over $\Z$, $\OO_{\Q_{\{q\}'}}/\mathfrak{q}$ is the algebraic closure of $\F_q$ , and there exists an embedding $\OO_{\Q_{\{q\}'}}/\mathfrak{q} \rar k_s$ by Galois theory again. 

Let us fix a morphism $\iota: \OO_{\Q_{\{q\}'}} \rar k_s$. 

Let us check that $\iota(\OO_{\Q_{\{q\}'}})$ is stable under $\mrm{Gal}(k_s/k)$. Any element in $\iota(\OO_{\Q_{\{q\}'}})$ is a root of some irreducible monic polynomial $f \in \Z[t]$. Conversely, let $f \in \Z[t]$ be a monic irreducible polynomial of degree $d$ and $x \in k_s$ be a root of $f$. If $f \pmod{q}$ is not separable (which implies necessarily that $q > 0$), then, by considering its decomposition in irreducible factors in $\F_q[t]$, there exists a monic polynomial $g \in \Z[t]$ such that $\sqrt{(f \pmod{q})}=(g\pmod{q})$ and $g \pmod{q}$ is separable. In particular, $x$ is still a root of $g$. After replacing $f$ with one of the irreducible monic factors of $g$ in $\Z[t]$, we can therefore assume that $f \pmod{q}$ is separable. Then the roots of $f$ are contained in $\OO_{\{q\}'}$, so we can write $f(t)=\prod_{i=1}^d{(t-y_i)}$ with $y_i \in \OO_{\{q'\}}$. Thus, in $k_s$, $(f \pmod{q})(t)=\prod_{i=1}^d{(t-\iota(y_i))}$, so $x$ is one of the $y_i$. 

Thus $\iota(\OO_{\Q_{\{q\}'}})$ is exactly the subset of elements of $k_s$ that are roots of irreducible monic $f \in \Z[t]$. In particular, it is stable under $\mrm{Gal}(k_s/k)$. 

Let $\mathfrak{q}=\ker{\iota}$, then, by definition, the group of automorphisms of $(\OO_{\{q\}'},\mathfrak{q})$ identifies with the automorphism group of $\OO_{\Q_{\{q\}'}}/\mathfrak{q}\simeq \iota(\OO_{\Q_{\{q\}'}})$. Thus, any $t \in \mrm{Gal}(k_s/k)$ induces an automorphism of $\iota(\OO_{\Q_{\{q\}'}})$, and there is therefore a unique automorphism $R_{\iota}(t) \in \mrm{Gal}(\Q_{\{q\}'}/\Q)$ preserving $(\OO_{\{q\}'},\mathfrak{q})$ such that for any $z \in \OO_{\{q\}'}$, $\iota(R_{\iota}(t)(z)) = t(\iota(z))$.

Let $\sigma \in \mrm{Gal}(\Q_{\{q\}'}/\Q)$ and $\iota'=\iota\circ\sigma$. One can check directly that $R_{\iota\circ\sigma} = \sigma^{-1}R_{\iota}\sigma$. 

Let $\iota': \OO_{\Q_{\{q\}'}} \rar k_s$ be another ring homomorphism with kernel $\mathfrak{q}$. The above argument shows that $\iota'$ and $\iota$ have the same image in $k_s$, thus $(\iota')^{-1}\iota$ is an automorphism of $\OO_{\Q_{\{q\}'}}/\mathfrak{q}$, which lifts to an automorphism $\sigma \in \mrm{Gal}(\Q_{\{q\}'}/\Q)$ preserving $(\OO_{\{q\}'},\mathfrak{q})$. Hence $\iota'=\iota\circ\sigma$. 

Let $\iota': \OO_{\Q_{\{q\}'}} \rar k_s$ be another ring homomorphism with kernel $\mathfrak{q}'$. There exists some $\sigma \in \mrm{Gal}(\Q_{\{q\}'}/\Q)$ such that $\sigma(\mathfrak{q})=\mathfrak{q}'$, so $\iota" := \iota'\circ\sigma: \OO_{\Q_{\{q\}'}} \rar k_s$ has kernel $\mathfrak{q}$: by the previous paragraph, we can write $\iota"=\iota \circ \sigma'$ for some $\sigma' \in \mrm{Gal}(\Q_{\{q'\}}/\Q)$, hence $\iota'=\iota \circ \sigma'\circ \sigma^{-1}$, and we are done.} 

\medskip

\cor[abelian-scheme-change-field]{Let us keep the notations of the previous lemma, and fix a ring homomorphism $\iota: \OO_{\Q_{\{q\}'}} \rar k_s$ with kernel $\mathfrak{q}$. Let $A$ be an abelian scheme over $\Z_\{(q)\}$. Let $\ell$ be a prime number distinct from $q$. Then, for any $n\geq 1$, the obvious homomorphisms 
\begin{align*}
&A[\ell^n](\OO_{\Q_{\{q\}'}} \otimes \Z_{(q)}) \rar A_{\Q}[\ell^n](\Q_{\{q\}'}),\quad \quad A_{\Q}[\ell^n](\Q_{\{q\}'}) \rar A_{\Q}[\ell^n](\Qbar),\\
&A[\ell^n](\OO_{\Q_{\{q\}'}} \otimes \Z_{(q)}) \rar A[\ell^n]\left((\OO_{\Q_{\{q\}'}})_{\mathfrak{q}}/\mathfrak{q}\right) \overset{\iota}{\rar} (A_k)[\ell^n](k_s), &&
\end{align*}
are isomorphisms.}

\demo{When $q=0$, this is simply the fact that $A[\ell^n](\Qbar) \overset{\iota}{\rar} (A_k)[\ell^n](k_s)$ is an isomorphism for any abelian variety $A$ over $\Q$, which is classical (see \cite[Theorem 8.2, Remark 8.4]{MilAb}). When $q > 0$, the map labelled $\iota$ is an isomorphism for the same reason (since $\iota: (\OO_{\Q_{\{q\}'}})_{\mathfrak{q}}/\mathfrak{q} \rar k_s$ is a morphism of separably closed fields), and the other maps are isomorphisms by the N\'eron-Ogg-Shafarevich criterion \cite[\S 1]{GoodRed}.}

\bigskip

\defi[notation-xrho]{Let $k$ be a field with characteristic prime to $p$ and separable closure $k_s$, and $\rho: \mrm{Gal}(k_s/k) \rar \GL{\F_p}$. We define the smooth projective curve (that is, pure one-dimensional scheme) $X_{\rho}(p)$ over $k$ as the twist of $X(p,p)_k$ by the action of $\rho$ (this is closer to the notation of Definition \ref{xrho-elementary}). Its Jacobian variety (in the sense of Proposition \ref{jacobian-relative-curve-exists}) is denoted by $J_{\rho}(p)$. }

\cor[decomposition-xrho]{Let $k$ be a field of characteristic distinct from $p$ and separable closure $k_s$. Let $F,\ell,\mathfrak{l}$ be as in Corollary \ref{decomposition-xpp-tate} such that $\ell \in k^{\times}$. Let $\rho: \mrm{Gal}(k_s/k) \rar \GL{\F_p}$ be a continuous group homomorphism. Then, as a $(\mathbb{T} \otimes_{\Z} F_{\mathfrak{l}})[G_{k}]$-module, $\Tate{\ell}{J_{\rho}(p)} \otimes_{\Z_{\ell}} F_{\mathfrak{l}}$ is isomorphic to the direct sum of
\begin{itemize}[label=\tiny$\bullet$,noitemsep]
\item for each $((f,\mathbf{1}),\psi) \in \mathscr{S} \times \mathscr{D}$, a copy of $\psi(\omega_p \det{\rho}) \otimes V_{f,\mathfrak{l}} \otimes \mrm{St}(\rho)$, on which, for any $n \geq 1$ coprime to $p$, $T_n$ (resp. $nI_2$) acts by $a_{n}(f)\psi(n)$ (resp. $\psi(n)^2$), 
\item for each $((f,\chi),\psi) \in \mathscr{P} \times \mathscr{D}$ modulo the relation $((f,\chi),\psi) \sim ((\overline{f},\overline{\chi}),\chi\psi)$, a copy of $\psi(\omega_p\det{\rho}) \otimes V_{f,\mathfrak{l}} \otimes (\pi(1,\chi)(\rho))$, on which, for any $n \geq 1$ coprime to $p$, $T_{n}$ (resp. $nI_2$) acts by $a_{n}(f)\psi(n)$ (resp. $\psi(n)^2\chi(n)$), 
\item for each $(f,\chi) \in \mathscr{C}$, a copy of $V_{f ,\mathfrak{l}} \otimes C_f(\rho)$, on which, for any $n \geq 1$ coprime to $p$, $T_{n}$ (resp. $nI_2$) acts by $a_{n}(f)$ (resp. $\chi(n)$).  
\end{itemize}
}

\demo{This follows essentially from Corollary \ref{decomposition-xpp-tate} using \cite[Proposition 2.2]{Virdol}, but let us spell the argument out. Fix some algebraic closure $\overline{k}$ of $k$. By Proposition \ref{jacobian-twist}, $J_{\rho}(p)$ identifies with the twist of $J(p,p)_k$ by $\rho$ (by the push-forward functoriality) and inherits its action by $\mathbb{T}$ by Proposition \ref{twist-equiv}. By Propositions \ref{cocycle-twist} and \ref{twist-equiv}, the following diagram commutes, with all the maps commuting to $\mathbb{T}$:
\[
\begin{tikzcd}[ampersand replacement=\&]
J_{\rho}(p)(\overline{k}) \arrow{r}{j_{J(p,p)_k}} \arrow{d}{\sigma}\& J(p,p)(\overline{k})\arrow{d}{\rho(\sigma) \circ \sigma}\\ 
J_{\rho}(p)(\overline{k}) \arrow{r}{j_{J(p,p)_k}} \& J(p,p)(\overline{k}),
\end{tikzcd}
\]
with $j_{J(p,p)_k}$ a morphism of abelian varieties by Proposition \ref{jacobian-twist}, and all the maps in the diagram commute with $\mathbb{T}$ by Proposition \ref{twist-equiv}. In particular, $j_{J(p,p)_k}$ identifies $\Tate{q}{J_{\rho}(p)}$ with $\Tate{q}{J(p,p)_k}$ as $\mathbb{T}$-module, and the $G_k$-actions coincide if we let $\sigma \in G_k$ act on $\Tate{q}{J(p,p)_k}$ by the composition of its natural action and $\rho(\sigma)$. 

By the previous remark and Corollary \ref{decomposition-xpp-tate}, the induced structure of $(\mathbb{T} \otimes F_{\mathfrak{l}})[\mrm{Gal}(k_s/k)]$-module on $\Tate{\ell}{J(p,p)_k} \otimes_{\Z_{\ell}} F_{\mathfrak{l}}$ is exactly the one described in the statement of this Corollary. }

\rem{Each of the factors in the decomposition is invariant if we replace $\rho$ with one of its quadratic twists. This fact is also a direct consequence of Proposition \ref{quadratic-twist-equiv} after compactifying. }

\section{Irreducibility of the factors of $\Tate{\ell}{J(p,p)_{\rho}}$}
\label{irreducibility-factors}

Let $\rho: \mrm{Gal}(k_s/k) \rar \GL{\F_p}$ be a continuous group homomorphism, where $k$ is a field of characteristic distinct from $p$ with separable closure $k_s$, and assume that $\det{\rho}=\omega_p^{-1}$ (we recall that $\omega_p$ is the mod $p$ cyclotomic character). By Proposition \ref{group-is-twist-polarized}, the twist of $X(p,p)_k$ by $\rho$ has exactly $p-1$ connected components $X_{\rho,\alpha}(p)$ (each one of them geometrically connected) for $\alpha \in \mu_p^{\times}(k_s)$, so the twist of $J(p,p)_k$ by $\rho$ has $p-1$ summands indexed by $\alpha \in \mu_p^{\times}(k_s)$, corresponding to the Jacobian varieties of the $X_{\rho,\alpha}(p)$. In particular, the direct sum of the Tate modules of the Jacobians of the $X_{\rho,\alpha}(p)$ is given by Corollary \ref{decomposition-xrho}. 

Our goal is to exploit the decomposition given by Corollary \ref{decomposition-xrho} to describe the Tate modules of the Jacobians of the $X_{\rho,\alpha}(p)$. 

\medskip

\lem[sl2zq-absorbs]{Let $G$ be a finite group and $q \nmid 6|G|$ be a prime number. Then any closed subgroup $H \leq \SL{\Z_q} \times G$ whose first projection is onto is equal to $\SL{\Z_q} \times G_1$ for some subgroup $G_1 \leq G$. }

\demo{First, we show that $U := \begin{pmatrix}1 & 1\\0 & 1\end{pmatrix}  \in \SL{\Z_q}$ satisfies $(U,e_G) \in H$. Let $D \subset \Z_q$ consisting of those $a \in \Z_q$ such that $\left(\begin{pmatrix} 1 & a\\0 & 1\end{pmatrix},e_G\right) \in H$. Then $D$ is a closed subgroup of $\Z_q$. The assumption implies that $(U,g) \in H$ for some $g \in G$. Then $(U^|G|,e_G)=(U,g)^{|G|} \in H$, so $|G| \in D$, so $\overline{|G|\Z} \subset D$. Since $q \nmid |G|$, $|G|\Z$ is dense in $\Z_q$, so $1 \in D$ and $(U,e_G) \in H$. 

For each $m \in \SL{\Z_q}$, the above argument shows that $(U,e_G) \in (m,1)^{-1}H(m,1)$, so $(mUm^{-1},e_G) \in H$. Therefore, the subset $S$ of those $M \in \SL{\Z_q}$ such that $(M,e_G) \in H$ is a closed subgroup containing every conjugate of $U$, so $S=\SL{\Z_q}$ and $\SL{\Z_q} \times \{e_G\} \subset H$, whence the conclusion.}

\lem[subreps-product]{Let $k$ be a field, $V,W$ be two $k$-vector spaces of finite dimension, and $f: G \rar H$ be a group homomorphism. Suppose that $H$ acts on the $k$-vector spaces $V$ and $W$ in such a way that $V_{|G}$ is absolutely irreducible and $f(G)$ acts trivially on $W$. Then any subrepresentation of $V \otimes_k W$ (over $k[H]$) is of the form $V \otimes_k W'$ for some subrepresentation $W' \leq W$. 
}

\demo{Under these assumptions, $V$ satisfies Schur's Lemma. Therefore, the morphism \[w \in W \longmapsto (-\otimes w) \in \mrm{Hom}_{f(G)}(V,V \otimes_k W)\] is an injective $k$-linear map between finite dimensional $k$-vector spaces of same dimension, hence is an isomorphism. 

Let $X \subset V \otimes_k W$ be any $k[H]$-submodule such that $X$ does not contain $V \otimes w$ for any $w \neq 0$. Then there is no nonzero $k[f(G)]$-morphism $V_{|f(G)} \rar X_{|f(G)}$. However, $X_{|f(G)}$ embeds in a sum of copies of the absolutely irreducible representation $V_{|f(G)}$, so any irreducible $k[f(G)]$-module of $X$ is isomorphic to $V_{|f(G)}$. So $X$ cannot have any irreducible $k[f(G)]$-submodule, so $X=0$. 

Let $X \subset V \otimes_k W$ be any $k[H]$-submodule. Then $W'=\{w \in W,\, X \supset V \otimes w\}$ is a $k[H]$-submodule of $W$, and $V \otimes_k W' \subset X$. Then $X'=X/(V \otimes_k W') \subset V \otimes_k (W/W')$ satisfies the assumptions of the above paragraph, so $X'=0$ and $X=V \otimes_k W'$.}

\bigskip

\prop[tensor-product-with-artin-irreducible]{Let $f \in \mathcal{S}_k(\Gamma_1(N))$ be a normalized newform without complex multiplication of weight $k \geq 2$, $F$ be a number field containing the coefficients of $f$, $\mathfrak{l}$ be a maximal ideal of $\OO_F$ with residue characteristic $\ell$, and $\lambda \in \OO_{F_{\mathfrak{l}}}$ be a uniformizer. Let $L$ be a number field and $\sigma: G_{L} \rar \GLn{n}{\OO_F}$ be an absolutely irreducible Artin representation. Then the representation $\sigma \otimes V_{f,\mathfrak{l}}$ is irreducible as a $F_{\mathfrak{l}}[G_L]$-module. Moreover, if $\ell$ is large enough (with respect to $f$, $L$ and $|\sigma(G_L)|$), then $(\sigma \otimes_{\OO_F} T_{f,\mathfrak{l}})/\lambda$ is an irreducible $(\OO_{F_{\mathfrak{l}}}/\lambda)[G_{L(\mu_{\ell^{\infty}})}]$-module. }

\demo{Let $M \supset L$ be the number field such that $\ker{\sigma}=G_M$. By \cite[Proposition 4.4]{Antwerp5-Ribet}, $(V_{f,\mathfrak{l}})_{|G_M}$ is irreducible. It is absolutely irreducible, because otherwise one of the following two impossible situations occurs:
\begin{itemize}[noitemsep,label=$-$]
\item $V_{f,\mathfrak{l}} \otimes \overline{F_{\mathfrak{l}}}$ has a unique $\overline{F_{\mathfrak{l}}}$-line stable under $G_M$, so this line is stable under $\mrm{Gal}(\overline{F_{\mathfrak{l}}}/F_{\mathfrak{l}})$, so this line is of the form $D \otimes \overline{F_{\mathfrak{l}}}$ for some line $D \subset V_{f,\mathfrak{l}}$ stable under $G_M$, which is impossible.
\item $V_{f,\mathfrak{l}} \otimes \overline{F_{\mathfrak{l}}}$ is the sum of two characters, so $V_{f,\mathfrak{l}}$ is an abelian representation of $G_M$, thus contradicting \cite[Proposition 4.4]{Antwerp5-Ribet}.
\end{itemize}

Let $W$ be a $F_{\mathfrak{l}}$-vector space on which $G_L$ acts by $\sigma$. By Lemma \ref{subreps-product} (for $G_M \rar G_L$), any $F_{\mathfrak{l}}[G_L]$-submodule of $V_{f,\mathfrak{l}} \otimes W$ is of the form $V_{f,\mathfrak{l}} \otimes W'$ for some subrepresentation $W'$ of $W$. Since $W$ is absolutely irreducible, $W'$ is either zero or $W$, so $V_{f,\mathfrak{l}} \otimes W'$ is irreducible. 

Let us assume that $M/\Q$ is unramified above $\ell$ and $\ell\nmid 30|\sigma(G_L)|$ (this holds for all but finitely many $\ell$). Let $\Gamma \subset \GL{\OO_{F_{\mathfrak{l}}}} \times \GLn{n}{\OO_K}$ be the image of $G_{L(\Q(\mu_{\ell^{\infty}}))}$ in $\mrm{Aut}(T_{f,\mathfrak{l}}) \times \mrm{Aut}(T_{\sigma})$, where $T_{\sigma}$ is some free $\OO_F$-module of finite rank on which $G_L$ acts by $\sigma$. 

Now, $\Gamma$ is a closed subgroup of $\GL{\OO_{F_{\mathfrak{l}}}} \times \mrm{im}(\sigma)$ (the second factor is a finite group), and, by \cite[Theorem 2.2.2]{bigimage}, its first projection contains a conjugate of $\SL{\Z_{\ell}}$ when $\ell$ is large enough. Since $\ell$ is unramified in $M/\Q$, the Galois extensions $L(\mu_{\ell^{\infty}})$ and $M$ are linearly disjoint over $L$, so the second projection of $\Gamma$ is exactly $\mrm{im}(\sigma)$. By Lemma \ref{sl2zq-absorbs}, $\Gamma$ contains a conjugate of $\SL{\Z_{\ell}} \times \mrm{im}(\sigma)$. 

To show the conclusion, it is enough to show that $(T_{f,\mathfrak{l}} \otimes_{\OO_F} T_{\sigma})/\mathfrak{l}$ is an absolutely irreducible $\Gamma$-module. Consider the space $R=\overline{\F_{\ell}}^{\oplus 2}$ endowed with its natural action of $\SL{\Z_{\ell}}$: it is then enough to show that $R \otimes_{\OO_F} (T_{\sigma}/\mathfrak{l})$ is an irreducible $\overline{\F_{\ell}}[\SL{\Z_{\ell}} \times \mrm{im}(\sigma)]$-module.  

The action of $\SL{\Z_{\ell}}$ on $R$ is (absolutely) irreducible, and, by our choice of $\ell$, the action of $\mrm{im}(\sigma)$ on $T_{\sigma}/\mathfrak{l}$ is absolutely irreducible. We can then apply Lemma \ref{subreps-product} to the morphism $\SL{\Z_{\ell}} \subset \SL{\Z_{\ell}} \times \mrm{im}(\sigma)$, whence the conclusion.   
}

\medskip

\cor[decomposition-xrho-irreducible]{In the situation of Corollary \ref{decomposition-xrho}, if $k$ is a number field and $\rho$ is onto, then all the given decomposition factors are absolutely irreducible, except when $p \equiv 3\pmod{4}$, $f$ is a CM modular form and one of the following two conditions is satisfied:
\begin{itemize}[noitemsep,label=$-$]
\item $k \supset \Q(\sqrt{-p})$,
\item $\omega_p\det{\rho}: G_k \rar \F_p^{\times}/\F_p^{\times 2}$ is trivial.
\end{itemize}. }

\demo{Since $\rho: G_{\Q} \rar \GL{\F_p}$ is onto, the representations $\mrm{St}(\rho),\pi(1,\chi)(\rho),C_f(\rho)$ are absolutely irreducible, so, when $f$ is not CM, we can apply Proposition \ref{tensor-product-with-artin-irreducible} and the conclusion follows. 

If $f$ is CM, then there is a nontrivial Dirichlet character $\epsilon$ such that $f \otimes \epsilon=f$. By considering the character of $f$, $\epsilon^2=1$, so $\epsilon$ is a quadratic character. By \cite[Theorem 4.5]{Antwerp5-Ribet}, $\epsilon$ is the character of an imaginary quadratic field, so $\epsilon(-1)=-1$. By Corollary \ref{level-newform-bigtwist}, since $f$ has conductor $p^2$ and its character has conductor at most $p$, the level of $\epsilon$ must be a power of $p$. Since $p > 2$, $\epsilon$ is the unique nontrivial quadratic character of $\F_p^{\times}$, that is, $\epsilon=\left(\frac{\cdot}{p}\right)$ is the character of the imaginary quadratic field $\Q(\sqrt{-p})$. Then $\left(\frac{-1}{p}\right)=\epsilon(-1)=-1$ so $p \equiv 3 \pmod{4}$. 

When none of the two conditions are verified, let $\psi,\psi': G_{\Q(\sqrt{-p})} \rar \overline{F_{\mathfrak{l}}}^{\times}$ be the two characters attached to $V_{f,\mathfrak{l}}$ (see \cite[Theorem 4.5]{Antwerp5-Ribet}). Let $k_1=k(\sqrt{-p})$: then $\rho(G_{k_1})$ is a normal subgroup of $\GL{\F_p}$ with index at most $2$, and is not $\F_p^{\times}\SL{\F_p}$ (otherwise $\det{\rho}: G_k \rar \F_p^{\times}/\F_p^{\times 2}$ is trivial or agrees with $\omega_p$), so $\rho(G_{k_1})=\GL{\F_p}$. 

Therefore, $(V_{f,\mathfrak{l}} \otimes C_f(\rho))_{|G_{k_1}} \simeq [\psi \otimes C_f(\rho_{|G_{k_1}})] \oplus [\psi' \otimes C_f(\rho_{|G_{k_1}})]$, and the two factors are irreducible. Assume that $V_{f,\mathfrak{l}} \otimes C_f(\rho)$ was reducible, and let $M$ be an irreducible subrepresentation. Up to exchanging $\psi,\psi'$, we can assume that $M_{|G_{k_1}} \simeq \psi \otimes C_f(\rho_{|G_{k_1}})$. After conjugating by $G_k$, we also have $M_{|G_{k_1}} \simeq \psi' \otimes C_f(\rho_{|G_{k_1}})$. By taking determinants, we see that $\left(\frac{\psi}{\psi'}\right)^{p-1}=1$. This is impossible since, by \cite[Theorem 4.3]{Antwerp5-Ribet}, $\frac{\psi}{\psi'}$ has infinite image.}

\bigskip

Before discussing in more detail the case where $(f,\chi) \in \mathscr{C}$ is a CM modular form, let us recall the properties of CM modular forms. This statement is well-known: for instance, it seems implied by \cite{Antwerp5-Ribet}. But since \emph{loc. cit.} (or any other text I was able to locate) does not explicitly provide or prove such an statement, we do it here. 

\prop[cm-yields-compatible-systems]{Let $f \in \mathcal{S}_k(\Gamma_1(N))$ be a normalized newform with complex multiplication by some imaginary quadratic field $K$, where $k \geq 2$. There exists a number field $F\subset \C$ containing the image of $K$ and the coefficients of $f$, and two locally constant characters $\psi_{\sigma}, \psi_{\tau}: \mathbb{A}_{K,f}^{\times} \rar F^{\times}$ satisfying the following properties (where $\mrm{Hom}(K,F)=\{\sigma,\tau\}$ and $\mathbb{A}_{K,f}^{\times}$ is the group of finite id\`eles of $K$):
\begin{itemize}[noitemsep,label=$-$]
\item for any $x \in K^{\times}$ and $i \in \{\sigma,\tau\}$, one has $\psi_i(x)=i(x)^{k-1}$,
\item both $\psi_i$ are trivial on $\prod_{v \nmid N}{\OO_{K_v}^{\times}}$,
\item pre-composing by the action of $\mrm{Aut}(K/\Q)$ exchanges the two $\psi_i$, 
\item for any prime ideal $\lambda$ of $\OO_F$ with residue characteristic $\ell$, for each $i \in \{\sigma,\tau\}$, the continuous character \[x =(x_{p})_p \in (\mathbb{A}_{\Q,f} \otimes K)^{\times} \simeq \mathbb{A}_{K,f}^{\times} \mapsto \psi_{i}(x)i(x_{\ell})^{1-k} \in F_{\lambda}\] is trivial on $K^{\times}$ and corresponds under the global class field theory reciprocity morphism\footnote{normalized so as to map a uniformizer to an arithmetic Frobenius} to a character $\psi_{i,\lambda}: G_K \rar F_{\lambda}^{\times}$ such that $V_{f,\lambda} \simeq \mrm{Ind}_K^{\Q}{\psi_{i,\lambda}}$. 
\end{itemize}
In other words, $(\psi_{\sigma,\lambda})_{\lambda}, (\psi_{\tau,\lambda})_{\lambda}$ form compatible systems over $F$ in the sense of Definition \ref{compatible-system-definition}. 

In particular, for any maximal ideal $\lambda$ of $\OO_F$, for any prime number $q$ split in $K$ and coprime to $N\lambda$, if $q\OO_K = \mathfrak{q}\mathfrak{q}'$, then $\psi_{\sigma}(\mathfrak{q}),\psi_{\sigma}(\mathfrak{q}')$ are the two roots of $X^2-a_q(f)X+q^{k-1}\chi(q)$, where $\chi$ is the character of $f$. 
  }

\demo{By \cite[Theorem 4.5]{Antwerp5-Ribet}, there is an integer $m \geq 1$ divisible by $N\Delta_{K/\Q}$, a character $\psi: J^{m} \rar \C^{\times}$ and an embedding $\sigma: K \rar \C$ such that:
\begin{itemize}[noitemsep,label=\tiny$\bullet$]
\item $J^{m}$ is the group of fractional ideals of $\OO_K$ coprime to $m$,
\item for any $x \in K^{\times}$ prime to $m$ and congruent to $1$ modulo $m$, $\psi((x)) = \sigma(x)^{k-1}$,
\item for any prime number $q$ split in $K$ and coprime to $m$, let $\mathfrak{q},\mathfrak{q}'$ denote the prime ideals of $\OO_K$ dividing $q$, then $\psi(\mathfrak{q})+\psi(\mathfrak{q}')=a_q(f)$, 
\item the character $\chi$ of $f$ is the product of the Dirichlet character $\varepsilon$ attached to $K/\Q$ and the character $a \in (\Z \backslash \cup_{p \mid m}{p\Z}) \longmapsto a^{1-k}\psi((a))$, which is a Dirichlet character modulo $m$.  
\end{itemize}

Now, the subgroup of $J^m$ made with principal ideals generated by elements $x \in K^{\times}$ congruent to $1$ modulo $m$ has finite index. Therefore, for any fractional ideal $I$ of $K$ coprime to $m$, one has $|\psi(I)|=(\mathbf{N}_{K/\Q}I)^{\frac{k-1}{2}}$. 

By \cite[Corollary VII.6.14]{Neukirch-ANT}, $\psi$ extends to a continuous character $\psi_0: \mathbb{A}_K^{\times}/K^{\times} \rar \C^{\times}$ such that, for any $z \in (K \otimes \R)^{\times}$, one has $\psi_0(z) = \sigma(z)^{1-k}$. Because $\OO_K^{\times}$ is finite, $K^{\times}$ is a discrete subgroup of $\mathbb{A}_{K,f}^{\times}$, so we can consider the restriction $\psi_{\sigma}$ of $\psi_0$ to $\mathbb{A}_{K,f}^{\times}$. Thus $\psi_{\sigma}: \mathbb{A}_{K,f}^{\times} \rar \C^{\times}$ satisfies the following properties:
\begin{itemize}[noitemsep,label=\tiny$\bullet$]
\item $\psi_{\sigma}$ is trivial on $\prod_{v \mid m}{(1+m\OO_{K_v})}\prod_{v \nmid m}{\OO_{K_v}^{\times}}$,
\item for any $z \in K^{\times}$, $\psi_{\sigma}(z) = \sigma(z)^{k-1}$,
\item if $q \nmid Nm$ is a prime number split in $K$ and $\mathfrak{q},\mathfrak{q}'$ are the two prime ideals above $q$, if $z\in K_{\mathfrak{q}},\, z'\in K_{\mathfrak{q}'}$ are uniformizers, then $\psi_{\sigma}(z)+\psi_{\sigma}(z')=a_q(f)$ and $\psi_{\sigma}(z)\psi_{\sigma}(z') = \chi(q)q^{k-1}$,
\item if $n \in \Z$ is an integer coprime to $m$, and $z$ is the id\`ele equal to $n$ at every place $v \nmid m$ and equal to $1$ at other places, then $\psi_{\sigma}(z) = \varepsilon(n)\chi(n)$. 
\end{itemize}

Let $\tau: K \rar \C$ be the other embedding and $\psi_{\tau}$ be the conjugate of $\psi_{\sigma}$ under the action of $\mrm{Gal}(K/\Q)$. 

Let $q \nmid m$ be a prime number split in $K$. Let $\mathfrak{q},\mathfrak{q}'$ denote the prime ideals of $\OO_K$ dividing $q$ and fix uniformizers $z \in K_{\mathfrak{q}}, z' \in K_{\mathfrak{q}'}$. The numbers $\psi_{\sigma}(z),\psi_{\sigma}(z')$ are roots of $X^2-a_q(f)X+q^{k-1}\chi(q)$, so they are algebraic. Let $F_0$ be a number field containing $K$, $a_q(f)$, $\chi(q)$ and $\psi_{\sigma}(z),\psi_{\sigma}(z')$. There exists an integer $n \geq 1$ such that $\mathfrak{q}^n,(\mathfrak{q}')^n$ are principal with generators $x,x' \in 1+m\OO_K$, so that 
\[\psi_{\sigma}(z)^n = \psi_{\sigma}(x) = \sigma(x)^{k-1},\quad \psi_{\sigma}(z')^n=\psi_{\sigma}(x')=\sigma(x')^{k-1},\] 
whence \[\psi_{\sigma}(z)\OO_{F_0} = \sigma(\mathfrak{q})^{k-1}\OO_{F_0},\, \psi_{\sigma}(z')\OO_{F_0} = \sigma(\mathfrak{q}')^{k-1}\OO_{F_0}.\]
In particular, $\frac{\psi_{\sigma}(z)}{\psi_{\sigma}(z')}$ is not an algebraic integer.

Let $\theta: \mathbb{A}_{K,f}^{\times} \rar \C^{\times}$ be a locally constant character satisfying the following properties:
\begin{itemize}[noitemsep,label=\tiny$\bullet$]
\item There are $\alpha,\beta \in \Z$ such that for any $x \in K^{\times}$, $\theta(x)=\sigma(x)^{\alpha}\tau(x)^{\beta}$,
\item For all but finitely many primes $q \nmid Nm$ split in $K$, if $\mathfrak{q},\mathfrak{q}'$ are the two ideals of $\OO_K$ dividing $q$ and $z \in K_{\mathfrak{q}},\,z' \in K_{\mathfrak{q}'}$ are uniformizers, then 
\[(X-\theta(z))(X-\theta(z')) = X^2-a_q(f)X+q^{k-1}\chi(q).\]
\end{itemize}
Let us prove that $\theta \in \{\psi_{\sigma},\psi_{\tau}\}$. 

Let $I \subset m\OO_K$ be a nonzero ideal such that for any finite place $v$ of $\OO_K$, $\theta(1+I\OO_{K_v})=\{1\}$. Moreover, $x \in \mathbb{A}_{K,f}^{\times} \mapsto |\theta(x)|\|x\|_{\mathbb{A}_{K,f}}^{\frac{\alpha+\beta}{2}} \in \R^{+\times}$ is equal to $1$ on the finite index subgroup $K^{\times}\prod_{v \mid I}{(1+\OO_{K_v})}\prod_{v \nmid I}{\OO_{K_v}^{\times}}$ of $\mathbb{A}_{K,f}^{\times}$. Since $\R^{+\times}$ is uniquely divisible, for any $x \in \mathbb{A}_{K,f}^{\times}$, one has $|\theta(x)|\|x\|_{\mathbb{A}_{K,f}}^{\frac{\alpha+\beta}{2}} = 1$. By taking for $x$ a uniformizer of a split prime ideal $\mathfrak{q} \subset \OO_K$ with residue caracteristic $q$ coprime to $Im$, we see that $\alpha+\beta = k-1$.

By the Dirichlet density theorem \cite[Theorem VII.13.2]{Neukirch-ANT}, we can find infinitely many primes $q \nmid m$ such that $q\OO_K=\mathfrak{q}\mathfrak{q}'$ with $\mathfrak{q},\mathfrak{q}'$ distinct ideals of $\OO_K$, and such that $\mathfrak{q}$ is a principal ideal generated by some $x \in 1+I\OO_K$. For all but finitely many such $q$, $\theta(x) = \sigma(x)^{\alpha}\tau(x)^{\beta}$ is a root of $X^2-a_q(f)X+q^{k-1}\chi(q)$, so $\theta(x) \in \{\psi_{\sigma}(x),\psi_{\tau}(x)\} = \{\sigma(x)^{k-1},\tau(x)^{k-1}\}$. Since $\frac{\sigma(x)}{\tau(x)}$ is not a root of unity, we see that $(\alpha,\beta)\in \{(0,k-1),\,(k-1,0)\}$. 

After pre-composing $\theta$ by the action of $\mrm{Gal}(K/\Q)$, we may assume that $(\alpha,\beta)=(k-1,0)$, so that $\psi_f\theta^{-1}: \frac{\mathbb{A}_{K,f}^{\times}}{K^{\times}\prod_{v \mid I}{(1+I\OO_{K_v})}\prod_{v \nmid I}{\OO_{K_v}^{\times}}} \rar \C^{\times}$ is a well-defined character from a finite group, so its image is contained in the roots of unity. 

Let $q \nmid m\mathbf{N}_{K/\Q}(I)$ be a prime number that splits in $K$ with $q\OO_K=\mathfrak{q}\mathfrak{q}'$. Let $z \in K_{\mathfrak{q}},\,z' \in K_{\mathfrak{q}'}$ be uniformizers and assume that $\theta(z),\theta(z')$ are the two roots of $X^2-a_q(f)X+q^{k-1}\chi(q)$. Then $\theta(z) \in \{\psi_{\sigma}(z),\psi_{\sigma}(z')\}$ is equal to $\psi(z)$ up to a root of unity. We saw that $\frac{\psi_{\sigma}(z)}{\psi_{\sigma}(z')}$ is not a root of unity, hence $\theta(z)=\psi_{\sigma}(z)$. By Dirichlet's density theorem \cite[Theorem VII.13.2]{Neukirch-ANT}, one then has $\theta=\psi_{\sigma}$. 

Let $F \subset \C$ be the number field generated by the image in $\C$ of $F$ and $K$. Applying the above argument to any $\mrm{Aut}(\C/F)$-conjugate of $\psi_{\sigma}$, we see that $\psi_{\sigma}$ (hence $\psi_{\tau}$) takes all its values in $F^{\times}$. 

Let now $\lambda$ be a prime ideal of $\OO_F$ with residue characteristic $\ell$. By \cite[Theorem 4.5]{Antwerp5-Ribet}, $V_{f,\lambda}$ is a semi-simple continuous abelian representation of $G_{K}$, so it decomposes as the sum of two continuous characters $\alpha, \beta: G_K \rar \overline{F_{\lambda}}$ that are unramified away from $N\ell$. Let, for $i \in \{\sigma,\tau\}$ (defining a morphism $K \otimes \Q_{\ell} \rar F_{\lambda}$), \[\psi_{i,\lambda}: x=(x_{\ell}) \in (K \otimes \mathbb{A}_{\Q,f})^{\times} \simeq \mathbb{A}_{K,f}^{\times} \mapsto \psi_{i}(x)i(x_{\ell})^{1-k} \in F_{\lambda}.\]

Then $\psi_{i,\lambda}$ is a continuous character $\mathbb{A}_K^{\times}/K^{\times}\prod_{v \nmid m\ell}{\OO_{K_v}^{\times}} \rar F_{\lambda}^{\times}$. 

Note that we have an exact sequence \[0 \rar (K \otimes \R)^{\times} \rar \mathbb{A}_K^{\times}/K^{\times} \rar \mathbb{A}_{K,f}^{\times}/K^{\times} \rar 1\] of locally compact Abelian groups, where $(K \otimes \R)^{\times} \simeq \C^{\times}$ is connected, and the profinite group $\prod_{v \mid 36m}{(1+36m\OO_{K_v})^{\times}}\prod_{v \nmid 36m}{\OO_{K_v}^{\times}}$ is an open subgroup of $\mathbb{A}_{K,f}^{\times}/K^{\times}$, so the latter is locally profinite hence totally disconnected. Hence $(K \otimes \R)^{\times}$ is the connected component of unity in $\mathbb{A}_K^{\times}/K^{\times}$. 

By \cite[Proposition 15.42]{Harari}, the class field theory reciprocity isomorphism thus identifies $\frac{\mathbb{A}_{K,f}^{\times}}{K^{\times}\prod_{v \nmid m\ell}{\OO_{K_v}^{\times}}}$ with the abelianization of the Galois group $G_{K,m\ell}$ of $K$ unramified outside $m\ell$ (by mapping a uniformizer at a place $v \nmid m\ell$ to an arithmetic Frobenius). Let $\psi'_{i,\lambda}: G_K \rar F_{\lambda}^{\times}$ be the character corresponding to $\psi_{i,\lambda}$ under this identification: $\psi'_{i,\lambda}$ is unramified outside $m\ell$. 

Thus, for any prime $q \nmid m\ell$ split in $K$, if $\mathfrak{q}$ denotes a prime ideal of $\OO_K$ above $q$, then $\psi'_{i,\lambda}(\Fr_{\mathfrak{q}})$ is the value of $\psi_{i}$ at a uniformizer of $\mathfrak{q}$ (in $K_{\mathfrak{q}}$), so it is a root of $X^2-a_q(f)X+q^{k-1}\chi(q)$, so it is either $\alpha(\Fr_{\mathfrak{q}})$ or $\beta(\Fr_{\mathfrak{q}})$. In other words, the function $(\psi'_{i,\lambda}-\alpha)(\psi'_{i,\lambda}-\beta): G_{K,m\ell} \rar \overline{F_{\lambda}}$ vanishes at the Frobenius of every prime ideal in $K$ which is prime to $m\ell$ and split in $K/\Q$. By Cebotarev's theorem \cite[Neukirch VII.13.4]{Neukirch-ANT}, $(\psi'_{i,\lambda}-\alpha)(\psi'_{i,\lambda}-\beta)=0$, hence $G_{K,m\ell}$ is the reunion of one of the two subgroups $\{\alpha=\psi'_{i,\lambda}\}$ and $\{\beta=\psi'_{i,\lambda}\}$. Note that $\psi_{\sigma,\lambda} \neq \psi_{\tau,\lambda}$ by considering the behavior at the place $\ell$, and we already know that $\alpha \neq \beta$. Therefore, $\{\alpha,\beta\} = \{\psi'_{\sigma,\lambda},\psi'_{\tau,\lambda}\}$, whence the conclusion. }

\nott{From now on (unless specified otherwise), we assume that the number field $F \subset \C$ satisfies furthermore the following property: for every $(f,\chi) \in \mathscr{C}$ with complex multiplication (necessarily by $\Q(\sqrt{-p})$ and $p \equiv 3\pmod{4}$), for every prime $q$ split in $\Q(\sqrt{-p})$, $F$ contains the roots of the Hecke polynomial $X^2-a_q(f)X+q\chi(q)$. 

To each such $(f,\chi) \in \mathscr{C}$, one can thus attach the two compatible systems $(\psi_{\sigma,\mathfrak{l}})_{\mathfrak{l} \in \mrm{Max}(\OO_F)}$ (for $\sigma \in \mrm{Hom}(K,F)$) of characters $G_{\Q(\sqrt{-p})} \rar F_{\mathfrak{l}}$. We mean this affirmation in the sense of Definition \ref{compatible-system-definition}, but for the remainder of the chapter, it is enough to note that for any prime ideal $\mathfrak{q}$ of $\Q(\sqrt{-p})$ with residue charateristic distinct from $p$, $\psi_{\sigma,\mathfrak{l}}(\Fr_{\mathfrak{q}})$ does not depend on the choice of a maximal ideal $\mathfrak{l}$ of $\OO_F$ with residue characteristic coprime to $\mathfrak{q}$. 

For the remainder of the text, we will in fact denote these two compatible systems by $(\psi_{f,\pm,\mathfrak{l}})_{\mathfrak{l}}$, where the assignment of $\pm$ is arbitrary. 
  }

\bigskip

\prop[decomposition-xrho-reducible-cm]{Let $(f,\chi) \in \mathscr{C}$ with $f$ CM, so that $p \equiv 3 \pmod{4}$ and $f$ has complex multiplication by $\Q(\sqrt{-p})$. Let $\chi_1 \in \mathcal{D}$ be such that $\chi_1^2=\chi$, and let $f_1 = f\otimes \chi_1^{-1}$. Then $(f_1,\mathbf{1}) \in \mathscr{C}$ is independent from the choice of $\chi_1$, and $C_{f_1}$ is the representation $V_4$ of Section \ref{splitting-sl2}. 

Let $\mathfrak{l}$ be a prime ideal of $\OO_F$ with residue characteristic $\ell$, and $\psi,\psi': G_{\Q(\sqrt{-p})} \rar F_{\mathfrak{l}}^{\times}$ be the two continuous characters (conjugate under $G_{\Q}/G_{\Q(\sqrt{-p})}$) such that $(V_{f_1,\mathfrak{l}})_{|G_{\Q(\sqrt{-p})}}$ is the direct sum of $\psi$ and $\psi'$. 

Let $k$ be a field with separable closure $k_s$ on which $p\ell$ is invertible. For any subextension $k'/k$ of $k_s/k$, let $G_{k'} = \mrm{Gal}(k_s/k')$. Note that there is a unique conjugacy class $e$ of morphisms from $G_k$ to the maximal quotient $G_{\Q,p\ell}$ of the absolute Galois group of $\Q$ which is unramified outside $p\ell$. Let $\alpha \in k_s$ be any root of the polynomial $X^2-X-\frac{p+1}{4}$. Let $\rho: G_k \rar \GL{\F_p}$ be any representation. Then the following statements hold:

\begin{enumerate}[label=(\roman*),noitemsep]
\item \label{notalldets-decomp} If $\det{\rho} \subset \F_p^{\times 2}$, then $C_f(\rho) \simeq \chi_1(\det{\rho}) \otimes \left(C^+(\rho) \oplus C^-(\rho)\right)$, where $C^+$ and $C^-$ are the representations $V^+$ and $V^-$ defined in Section \ref{splitting-sl2}. %

\item \label{numberfield-decomp} If $\alpha \in k$, then the image of $e$ is contained in $G_{\Q(\sqrt{-p}),p\ell}$ and \[(V_{f,\mathfrak{l}}\circ e) \otimes C_f(\rho) \simeq \chi_1(\omega_p\det{\rho})\otimes [(\psi(e) \otimes V_4(\rho)) \oplus (\psi'(e) \otimes V_4(\rho))].\] 

\item \label{numberfield-decomp-best} When $e$ has open image and $\rho(G_k)$ contains a conjugate of the Borel subgroup (i.e. the subgroup of upper-triangular matrices), both summands in the decomposition \ref{numberfield-decomp} are irreducible and they are not isomorphic. 

\item \label{cyclotomic-decomp} If $\alpha \notin k$ and $\omega_p\det{\rho}: G_k \rar \F_p^{\times}/\F_p^{\times 2}$ is trivial, let $e_1=e_{|G_{k(\alpha)}}$, then, for any $\epsilon \in \{\pm\}$, \[(V_{f,\mathfrak{l}} \circ e) \otimes C_f(\rho) \simeq \chi_1(\omega_p\det{\rho}) \otimes \bigoplus_{\theta \in \{\psi\circ e_1,\psi'\circ e_1\}}{\mrm{Ind}_{k(\alpha)}^k(\theta  \otimes C^{\epsilon}(\rho_{|G_{k(\alpha)}}))}.\]

\item \label{cyclotomic-decomp-best} In the situation of \ref{cyclotomic-decomp}, if $e$ has open image and $\rho(G_k)$ contains a conjugate of the Borel subgroup, the two summands in the decomposition of $V_{f,\mathfrak{l}} \otimes C_f(\rho)$ are irreducible. 
\end{enumerate}}

\demo{We already saw in the proof of Corollary \ref{decomposition-xrho-irreducible} why $p \equiv 3 \pmod{4}$ and $f \otimes \left(\frac{\cdot}{p}\right)=f$. Since the character $\chi_1$ is defined up to multiplication by $\left(\frac{\cdot}{p}\right)$, $f \otimes \chi_1^{-1}$ is independent from the choice of $\chi_1$. It is then clear using Proposition \ref{twist-exists} that $(f_1,\mathbf{1}) \in \mathscr{C}$. 

Let $\phi: \F_{p^2}^{\times}/\F_p^{\times}$ be one of the two characters attached to $C_{f_1}$. Since $f_1=f_1 \otimes \left(\frac{\cdot}{p}\right)$, by Corollaries \ref{twisting-Cf} and \ref{cuspidal-twist}, $\phi\cdot \left(\frac{N_{\F_{p^2}/\F_p}(\cdot)}{p}\right) \in \{\phi,\phi^p\}$. The only possibility is $\phi\cdot \left(\frac{N_{\F_{p^2}/\F_p}(\cdot)}{p}\right) = \phi^p$, so $\phi^2=(\phi^p)^2$, hence $\phi^4=\phi^{2(p+1)-2p+2}=1$. Since $\phi$ does not factor through the norm, $\phi^2 \neq 1$, so $\phi$ has order $4$ and $C_{f_1} \simeq V_4$. 

This implies by Corollary \ref{twisting-Cf} that $C_f \simeq \chi_1(\det) \otimes C_{f_1}$. \ref{notalldets-decomp} is then a direct consequence of Proposition \ref{decomp-cusp-sl}. 

The decomposition in \ref{numberfield-decomp} is a consequence of the theory of complex multiplication described in Proposition \ref{cm-yields-compatible-systems} (or see \cite[Proposition 4.4, Theorem 4.5]{Antwerp5-Ribet}). For \ref{numberfield-decomp-best}, the irreducibility follows from Proposition \ref{cusp-on-borel}, and the two summands are not isomorphic because $\frac{\psi}{\psi'}$, hence $\frac{\psi(e)}{\psi'(e)}$ has infinite image by \cite[Theorem 4.3]{Antwerp5-Ribet}. 

From now on, assume that $\alpha \notin k$ and that $(\omega_p\det{\rho})(G_k) \subset \F_p^{\times 2}$, and fix $\epsilon \in \{\pm\}$. Let $k_1=k(\alpha)$, $e_1=e_{|G_{k_1}}$. Then
\begin{align*}
(V_{f,\mathfrak{l}}\circ e) \otimes C_f(\rho) &\simeq \chi_1(\det{\rho}\omega_p)\otimes (V_{f_1,\mathfrak{l}} \circ e)\otimes C_{f_1}(\rho) \\
&\simeq \chi_1(\omega_p\det{\rho}) \otimes (V_{f_1,\mathfrak{l}} \circ e) \otimes \mrm{Ind}_{k_1}^k{C^{\epsilon}(\rho_{|G_{k_1}})}\\
&\simeq \chi_1(\omega_p\det{\rho}) \otimes \mrm{Ind}_{k_1}^k\left[(V_{f_1,\mathfrak{l}} \circ e_1) \otimes C^{\epsilon}(\rho_{|G_{k_1}}) \right]\\
&\simeq \chi_1(\omega_p\det{\rho}) \otimes \mrm{Ind}_{k_1}^k\left[((\psi \oplus \psi')\circ e_1) \otimes C^{\epsilon}(\rho_{|G_{k_1}})\right]\\
&\simeq \chi_1(\omega_p\det{\rho}) \otimes \bigoplus_{\theta \in \{\psi\circ e_1,\psi'\circ e_1\}}{\theta \otimes C^{\epsilon}(\rho_{|G_{k_1}})},
\end{align*}
which proves \ref{cyclotomic-decomp}. 

Assume that $\rho(G_k)$ contains a conjugate of the Borel subgroup. The claim is that for any $\theta \in \{\psi\circ e_{|G_{k_1}},\psi'\circ e_{|G_{k_1}}\}$ and $\epsilon \in \{\pm\}$, $W=\mrm{Ind}_{k_1}^k{\theta \otimes C^{\epsilon}(\rho_{|G_{k_1}})}$ is irreducible. Since $\psi\circ e_{|G_{k_1}},\psi'\circ e_{|G_{k_1}}$ are conjugate under $G_{k}/G_{k_1}$ and the same holds for $C^+(\rho_{|G_{k_1}}),C^-(\rho_{|G_{k_1}})$, we may assume $\theta=\psi\circ e_{|G_{k_1}}$, so that 
\[W_{|G_{k_1}} \simeq \left[(\psi\circ e_{|G_{k_1}}) \otimes C^{\epsilon}(\rho_{|G_{k_1}})\right] \oplus \left[(\psi'\circ e_{|G_{k_1}}) \otimes C^{-\epsilon}(\rho_{|G_{k_1}}).\right]\]
Both summands are irreducible by Proposition \ref{cusp-sl-borel}. 

So if $M \subset W$ is a proper nonzero irreducible $G_k$-submodule, then $M_{|G_{k_1}}$ is isomorphic to $(\psi\circ e_{|G_{k_1}}) \otimes C^{\epsilon}(\rho_{|G_{k_1}})$ or its conjugate under $G_k/G_{k_1}$. Since $M_{|G_{k_1}}$ is isomorphic to its conjugate under $G_k/G_{k_1}$, we find that \[(\psi\circ e_{|G_{k_1}}) \otimes C^{\epsilon}(\rho_{|G_{k_1}}) \simeq (\psi'\circ e_{|G_{k_1}}) \otimes C^{-\epsilon}(\rho_{|G_{k_1}}).\]

Thus $\frac{\psi}{\psi'}$ is trivial on the finite-index subgroup $e(G_{k_1} \cap \ker{\rho})$ of $G_{\Q(\sqrt{-p})}$, which contradicts \cite[Theorem 4.3]{Antwerp5-Ribet}.  
}

\medskip

\bigskip

\section{Tate modules of the Jacobians of the connected components of $X(p,p)_{\rho}$}
\label{tate-module-connected-xrho}

Let $E$ be an elliptic curve over a field $k$ with characteristic distinct from $p$. Let $k_s$ be a separable closure of $k$ and $(P,Q)$ be a basis of $E[p](k_s)$. In the terminology of Section \ref{application-moduli}, $E[p]$ with its Weil pairing is a polarized $p$-torsion group\footnote{that is, a finite \'etale group scheme over $\Sp{k}$, \'etale locally isomorphic to $\F_p^{\oplus 2}$, endowed with a (fibrewise) perfect alternated bilinear pairing $E[p] \times E[p] \rar \mu_p$.} over $\Sp{k}$.

Let, for each $\sigma \in \mrm{Gal}(k_s/k)$, $\rho_1(\sigma) \in \GL{\F_p}$ be the matrix such that $\rho_1(\sigma)\begin{pmatrix}P \\Q\end{pmatrix}=\begin{pmatrix}\sigma(P)\\\sigma(Q)\end{pmatrix}$. By Lemma \ref{right-representation}, $\rho_1: \mrm{Gal}(k_s/k) \rar \GL{\F_p}$ is an anti-homomorphism, hence $\rho=\rho_1^{-1}$ is a group homomorphism with determinant $\omega_p^{-1}$. 

By Proposition \ref{group-is-twist-polarized}, we have a commutative diagram of $k$-schemes, where the bottom schemes are constant, the horizontal arrows are isomorphisms, and the fibres of the vertical morphisms are smooth projective geometrically connected $k$-schemes of dimension one:

\[
\begin{tikzcd}[ampersand replacement=\&]
X_{E[p]}(p) \arrow{rr} \arrow{d}{\det} \&\& X(p,p)_{\rho} \arrow{d}{\mrm{We}_{\rho}}\\
\F_p^{\times} \arrow{rr}{a \mapsto \underline{a}(\langle P,\,Q\rangle)} \&\&\mrm{Hom}_k(k \otimes_{\Z} \Z[1/p,\zeta_p], k_s).
\end{tikzcd}
\]

In fact, Proposition \ref{group-is-twist-polarized} tells that we have a similar diagram for any $\rho$ such that $\det{\rho}=\omega_p^{-1}$, since such a representation encodes precisely a polarized $p$-torsion group over $k$. In particular, $X(p,p)_{\rho}$ is the reunion of $p-1$ smooth projective geometrically connected $k$-schemes of dimension one, indexed by the $\F_p^{\times}$-torsor $\mrm{Hom}_k(k \otimes_{\Z} \Z[1/p,\mu_p], k_s) = \mu_p^{\times}(k_s)$ (because the rightmost vertical map does not depend on the choice of basis of $G(k_s)$).   

\bigskip

\defi[notation-xrhoxi]{Recall from Definition \ref{notation-xrho} that if $k$ is a field of characteristic not dividing $p$ with separable closure $k_s$ and $\rho: \mrm{Gal}(k_s/k) \rar \GL{\F_p}$ is a group homomorphism with $\det{\rho}=\omega_p^{-1}$, we denote by $X_{\rho}(p)$ the smooth proper $k$-scheme of dimension one $(X(p,p)_k)_{\rho}$. It is endowed (by Proposition \ref{group-is-twist-polarized}) with a canonical morphism to the disjoint union of copies of $\Sp{k}$ indexed by $\mu_p^{\times}(k_s)$; given $\xi \in \mu_p^{\times}(k_s)$, we denote by $X_{\rho,\xi}(p)$ the inverse image under said morphism of $\xi$. 
}

\bigskip

\lem[agreement-implies-twist]{Let $f \in \mathcal{S}_k(\Gamma_1(N)), g \in \mathcal{S}_{k'}(\Gamma_1(N'))$ be normalized newforms and $M>1$ be an integer such that, for every prime number $\ell \equiv 1\pmod{M}$ and coprime to $NN'$, one has $a_{\ell}(f)=a_{\ell}(g)$. Then $k=k'$ and $f,g$ are twists one of the other.}

\demo{We may assume that $NN' \mid M$. Let us fix a prime number $p \mid M$, a number field $F \subset \C$ containing the coefficients of $f$ and $g$ as well as the $\varphi(M)$-th roots of unity, and a prime ideal $\mathfrak{p}$ of $\OO_F$ with residue characteristic $p$. We consider the $\mathfrak{p}$-adic representations $V_f, V_g$ of $G_{\Q,M}$ attached to $f$ and $g$ respectively, where the index $M$ indicates that the representations are unramified at any prime not dividing $M$. 

For every prime $\ell \equiv 1 \pmod{M}$, $\mrm{Tr}(\Fr_{\ell}\mid V_f)=a_{\ell}(f)=a_{\ell}(g)=\mrm{Tr}(\Fr_{\ell}\mid V_g)$. By Cebotarev, the $\Fr_{\ell}$ (for $\ell \equiv 1\pmod{M}$) are dense in $G_{\Q(\mu_M),M}$, so the representations $V_{f,M} = (V_f)_{|G_{\Q(\mu_M),M}},V_{g,M} = (V_g)_{|G_{\Q(\mu_M),M}}$ have the same traces, hence the same characteristic polynomials. By considering their determinants, we find that $k=k'$. 

Moreover, $V_{f,M}$ and $V_{g,M}$ are semi-simple (they have finite image when $k=1$, and it is \cite[Proposition 4.2]{Antwerp5-Ribet} when $k > 1$), so they are isomorphic, hence the $\mrm{Gal}(\Q(\mu_M)/\Q)$-module $H := \mrm{Hom}_{G_{\Q(\mu_M),M}}(V_f,V_g)$ is nonzero. 

Therefore, there exists a character $\psi: \mrm{Gal}(\Q(\mu_M)/\Q) \rar F^{\times}$ such that for some nonzero $u \in H$, $\sigma(u)=\psi(\sigma)u$ for all $\sigma \in \mrm{Gal}(\Q(\mu_M)/\Q)$. Thus $u$ corresponds to a nonzero $G_{\Q}$-homomorphism $V_{f,\mathfrak{p}} \rar V_{g \otimes \psi,\mathfrak{p}}$. Since both representations are irreducible, $u$ is an isomorphism and $f=g \otimes \psi$ by strong multiplicity one.}

\bigskip

\nott{We denote by $\lambda$ the character $\left(\frac{\cdot}{p}\right): \F_p^{\times} \rar \{\pm 1\}$, it is the unique quadratic Dirichlet character of conductor $p$.  

Finally, given a prime $\ell$, the ring $\OO_{\{\ell\}'}$ denotes the ring of integers of the maximal extension $\Q_{\{\ell\}'}$ of $\Q$ which is unramified at $\ell$. When $\ell=0$, $\OO_{\{\ell\}'}$ denotes the ring $\overline{\Z}$ of algebraic integers. }

\medskip

\defi[notation-T1]{The subalgebra $\mathbb{T}_1$ is the subalgebra of $\mathbb{T}$ generated by the $T_n\cdot (qI_2)$, where $n,q \geq 1$ are integers such that $nq^2 \equiv 1\pmod{p}$. In fact, the $T_n \cdot (qI_2)$ generate $\mathbb{T}_1$ as a group, by Lemma \ref{hecke-tn-group}. }

\medskip 

For later use, we also define larger endomorphism algebras of $\mathbb{T}[\GL{\F_p}]$ that will usually be non-commutative or act non-rationally on twists of $J(p,p)$. 

\defi[notation-T1Gamma]{Let $G,H \leq \GL{\F_p}$ be subgroups. We denote by $\mathbb{T}_{1,G}(H)$ the subring of $\mathbb{T}[\GL{\F_p}]$ generated by the $(T_n \cdot (mI_2))u$, where:
\begin{itemize}[noitemsep,label=$-$]
\item $n \geq 1$ is coprime to $p$,
\item $u \in \Z[G]$ is a formal linear combination of matrices of determinant $(nm^2)^{-1} \in \F_p^{\times}$,
\item $u$ commutes to $H$. 
\end{itemize}
Note that $\mathbb{T}_{1,G}(H)$ is generated as a group by the $T_nu$ as above, by Lemma \ref{hecke-tn-group}, and that its subring $\mathbb{T}_1$ is central. When $H$ contains only scalar matrices, we write $\mathbb{T}_{1,G}(H) := \mathbb{T}_{1,G}$. } 

\medskip

The interest of the $\mathbb{T}_{1,G}(H)$ is the following:

\lem[T1Gamma-action]{Let $k$ be a field such that $p$ is invertible in $k$, and let $k_s$ be its separable closure. Let $\rho: \mrm{Gal}(k_s/k) \rar \GL{\F_p}$ be a group homomorphism with $\det{\rho}=\omega_p^{-1}$. Then:
\begin{itemize}[noitemsep,label=$-$]
\item $\mathbb{T}[\GL{\F_p}]$ naturally embeds in $\mrm{End}_{k_s}\left(\prod_{\xi}{J_{\rho,\xi}(p)}\right)$, 
\item for any $M \in \GL{\F_p}$ and any $\xi \in \mu_p^{\times}$, $M$ acts by a morphism $J_{\rho,\xi}(p)_{k_s} \rar J_{\rho,\xi^{\det{M}}}(p)_{k_s}$ of abelian varieties over $k_s$,
\item the image of $\mathbb{T}_{1,\GL{\F_p}}$ in $\mrm{End}(J_{\rho}(p)_{k_s})$ is contained in $\prod_{\xi \in \mu_p^{\times}(k_s)}{\mrm{End}_{k_s}(J_{\rho,\xi}(p)_{k_s})}$,
\item the image of $\mathbb{T}_{1,\GL{\F_p}}(\im{\rho})$ in $\mrm{End}(J_{\rho}(p)_{k_s})$ is contained in $\prod_{\xi \in \mu_p^{\times}(k_s)}{\mrm{End}_{k}(J_{\rho,\xi}(p))}$. 
\end{itemize}}

\demo{By Proposition \ref{jacobian-twist}, $\prod_{\xi}{J_{\rho,\xi}(p)}$ is the twist of $J(p,p)_k$ by $\rho$, which means that $\mathbb{T}[\GL{\F_p}]$ tautologically embeds in $\mrm{End}_{k_s}\left(\prod_{\xi}{J_{\rho,\xi}(p)}\right)$. 

Moreover, any $z \in \mathbb{T}[\GL{\F_p}]$ commuting to $\im{\rho}$ has a twist $z_{\rho}$ as an endomorphism of $\prod_{\xi}{J_{\rho,\xi}(p)}$ by Proposition \ref{twist-equiv}. 

To conclude, it is enough to show by Proposition \ref{jac+hecke-twist-polarized} that for any $u \in \GL{\F_p}$, $u$ naturally induces a morphism of $k_s$-schemes: $X_{\rho,\xi}(p)_{k_s} \rar X_{\rho,\xi^{\det{u}}}(p)_{k_s}$. In other words, we have to show that the leftmost cell in the following diagram commutes:
\[
\begin{tikzcd}[ampersand replacement=\&]
\mu_p^{\times}(k_s) \arrow{rrr}{j} \arrow{ddd}{\det{u}}\& \& \& (\mu_p^{\times})_{k_s} \arrow{ddd}{\det{u}}\\
\& X_{\rho}(p) \arrow{ul}{\mrm{We}_{\rho}} \arrow{r}{j}\arrow{d}{u} \& X(p,p)_{k_s} \arrow{ru}{\mrm{We}} \arrow{d}{u} \&\\
\& X_{\rho}(p) \arrow{dl}{\mrm{We}_{\rho}} \arrow{r}{j} \& X(p,p)_{k_s} \arrow{rd}{\mrm{We}} \&\\
\mu_p^{\times}(k_s) \arrow{rrr}{j} \& \& \& (\mu_p^{\times})_{k_s} 
\end{tikzcd}
\]
Note that all vertical and horizontal arrows in this diagram are isomorphisms. The top and bottom cells commute because of the properties of twists (Proposition \ref{twist-equiv} and Proposition \ref{group-is-twist-polarized}), and the same holds for the outer cell. The rightmost cell commutes because of Proposition \ref{XN-Weil}. The inner cell commutes by definition, so we are done. 

}

\medskip

\defi[notation-T1rho]{Let $k$ be a field with characteristic distinct from $p$ and separable closure $k_s$. Let $\rho: \mrm{Gal}(k_s/k) \rar \GL{\F_p}$ be a continuous group homomorphism such that $\det{\rho}=\omega_p^{-1}$. Let $\Gamma \leq \GL{\F_p}$ be a subgroup such that $\im{\rho}$ is contained in the normalizer of $\Gamma$. We denote by $\mathbb{T}_{1,\Gamma,\rho}$ the subalgebra of $\mathbb{T}[\Gamma \rtimes_{\rho} \mrm{Gal}(k_s/k)]$ generated by the $(T_n \cdot (mI_2))(M,\sigma)$, for $n, m \geq 1, M \in \Gamma, \sigma \in \mrm{Gal}(k_s/k)$ such that $m^2n\det{M} \equiv 1 \pmod{p}$. It acts by group homomorphisms on $J_{\rho,\xi}(p)(k_s)$ for every $\xi \in \mu_p^{\times}(k_s)$. 

The given elements generate $\mathbb{T}_{1,\Gamma,\rho}$ as a group and $\mathbb{T}_1$ is a central subring of $\mathbb{T}_{1,\Gamma,\rho}$.}

\medskip

Since we have to work with each connected component separately, we will consider $\mathbb{T}_1$-module structures rather than $\mathbb{T}$-module structures.

\defi[smaller-scp]{For the sake of easier notation, we will write:
\begin{itemize}[noitemsep,label=$-$]
\item $\mathscr{P}'$ for a system of representatives of $\mathscr{P}$ modulo the complex conjugation,
\item $\mathscr{C}_M$ for the $(f,\mathbf{1}) \in \mathscr{C}$ with $f$ CM,
\item $\mathscr{C}'_1$ for a system of representatives of $(f,\mathbf{1}) \in \mathscr{C}$ with $f$ not CM and modulo twists by $\lambda$.
\end{itemize}

If $(f,\chi) \in \mathscr{S} \cup \mathscr{P} \cup \mathscr{C}$ and $M$ is a $(\mathbb{T}_1 \otimes F)$-module, its \emph{$f$-eigenspace} denotes the greatest submodule of $M$ on which every $(nI_2) \cdot T_m$ (for $m,n \geq 1$ coprime to $p$ such that $mn^2\equiv 1\pmod{p}$) acts by $a_m(f)\chi(n)$.} 

\medskip

\lem[eigenspaces-t1-modules]{Let $M$ be a $(\mathbb{T}_1 \otimes F)$-module.
\begin{itemize}[noitemsep,label=$-$]
\item Let $(f,\chi) \in \mathscr{S} \cup \mathscr{P} \cup \mathscr{C}$. The functor mapping a $(\mathbb{T} \otimes F)$-module $M$ to its $f$-eigenspace is exact.
\item If $(f,\chi) \in \mathscr{P}$ (resp. in $\mathscr{C}$), the $f$-eigenspace of $M$ is the same as the $f'$-eigenspace of $M$, where $(f',\chi') \in \mathscr{P}'$ is either $(f,\chi)$ or its complex conjugate (resp. $(f',\chi') \in \mathscr{C}'_1 \cup \mathscr{C}_M$ is the unique element such that $f'$ is a twist of $f$). 
\item $M$ is the direct sum of its $f$-eigenspaces over all $(f,\chi) \in \mathscr{S} \cup \mathscr{P}' \cup \mathscr{C}'_1 \cup \mathscr{C}_M$.
\end{itemize}}

\demo{Let $\ell$ be any prime number. Then $A := \mathbb{T}_1 \otimes_{\Z} F$ acts faithfully on \[T := \Tate{\ell}{J(p,p)} \otimes_{\Z} F \simeq \bigoplus_{\mathfrak{l} \mid \ell}{\Tate{\ell}{J(p,p)} \otimes_{\Z_{\ell}} F_{\mathfrak{l}}}.\] 

By Proposition \ref{decomposition-xpp-tate}, $T$ is the sum of the following $A$-submodules:
\begin{itemize}[noitemsep, label=\tiny$\bullet$]
\item for each $((f,\chi),\psi) \in \mathscr{S} \times \mathcal{D}$, a nonzero submodule where, for any $n,m \geq 1$ coprime to $p$ such that $n^2m\equiv 1\pmod{p}$, $(nI_2) \cdot T_m$ acts by $\psi^2(n)a_m(f)\psi(m)=\psi(mn^2)a_m(f)=a_m(f)$,
\item for each $((f,\chi),\psi) \in \mathscr{P} \times \mathcal{D}$ (modulo the relation $((f,\chi),\psi) \sim ((\overline{f},\overline{\chi}),\psi\chi)$), a nonzero submodule where, for any $m,n \geq 1$ coprime to $p$ such that $n^2m\equiv 1\pmod{p}$, $(nI_2) \cdot T_m$ acts by $\psi^2(n)\chi(n)a_m(f)\psi(m)=\chi(n)a_m(f)$,
\item for each $(f,\chi) \in \mathscr{C}$, a nonzero submodule where, for any $m,n \geq 1$ coprime to $p$ such that $n^2m\equiv 1\pmod{p}$, $(nI_2) \cdot T_m$ acts by $\chi(n)a_m(f)$. 
\end{itemize}

Let $m,n \geq 1$ be integers coprime to $p$ such that $n^2m\equiv 1\pmod{p}$. 

For any $(f,\chi) \in \mathscr{P}$, one has $\overline{\chi}(n)a_m(\overline{f})=\overline{\chi}(nm)a_m(f)=\chi(n)a_m(f)$. 

For any $((f,\chi),\psi) \in \mathscr{C} \times \mathcal{D}$, one has $(\chi\psi^2)(n)a_m(f \otimes \psi) = \psi(mn^2)\chi(n)a_m(f)=\chi(n)a_m(f)$. Therefore, $T$ is the sum of its $f$-eigenspaces over $(f,\chi) \in \Upsilon := \mathscr{S} \cup \mathscr{P}' \cup \mathscr{C}'_1 \cup \mathscr{C}_M$. 

By Lemma \ref{agreement-implies-twist}, the eigenspaces attached to any distinct two elements of $\Upsilon$ are distinct. So $T$ is the direct sum, over all $(f,\chi) \in \Upsilon$ of its $f$-eigenspaces, which are all nonzero. 

By construction, for every $(f,\chi) \in \Upsilon$, there is a character $\mu_f: \mathbb{T}_1 \rar F$ mapping $(nI_2) \cdot T_m$ (for $m,n \geq 1$ coprime to $p$ such that $n^2m\equiv 1\pmod{p}$ is nonzero) to $\chi(n)a_m(f)$, and $A$ acts on $T$ through the $F$-algebra homomorphism $A \rar \prod_{(f,\chi) \in \Upsilon}{A/\ker{\mu_f}}$. Since the $F$-algebra homomorphisms $\mu_f: A \rar F$  are pairwise distinct, they have pairwise distinct kernels (which are maximal ideals), so that $A \rar \prod_{(f,\chi) \in \Upsilon}{A/\ker{\mu_f}}$ is onto. Since $A \rar \mrm{End}_F(T)$ is injective, $A \rar \prod_{(f,\chi) \in \Upsilon}{A/\ker{\mu_f}}$ is an isomorphism. 

In particular, taking the $f$-eigenspace is exactly equivalent to tensoring with $A \rar A/\ker{\mu_f}$, which is a flat ring homomorphism, whence the conclusion.
}

\medskip

\cor[decomposition-xrho-cyclotomic]{Let $\mathfrak{l}$ be a maximal ideal of $\OO_F$ with residue characteristic $\ell$. Let $k$ be a field with separable closure $k_s$, where $p\ell$ is invertible and $\rho: \mrm{Gal}(k_s/k) \rar \GL{\F_p}$ be a continuous group homomorphism with $\det{\rho}=\omega_p^{-1}$. Let $\OO$ be the ring of integers of the maximal extension of $\Q$ unramified at the characteristic of $k$, whose Galois group over $\Q$ is $G_{\Q,k}$. Fix a ring homomorphism $\iota: \OO \rar k_s$ and let $e: \mrm{Gal}(k_s/k) \rar G_{\Q,k}$ be the attached group homomorphism by Lemma \ref{from-qbar-to-ks}. If $p \equiv 3 \pmod{4}$, we denote by $\alpha \in k_s$ a root of the polynomial $X^2-X-\frac{p+1}{4}$. 

Then the eigenspaces under $\mathbb{T}_1$ of $\Tate{\ell}{J(p,p)_{\rho}} \otimes_{\Z_{\ell}} F_{\mathfrak{l}}$ are the following $(\mathbb{T}_1 \otimes F_{\mathfrak{l}})[\mrm{Gal}(k_s/k)]$-modules:
\begin{itemize}[noitemsep,label=$-$]
\item for $(f,\mathbf{1}) \in \mathscr{S}$, the $f$-eigenspace is the sum of $p-1$ copies of $(V_{f,\mathfrak{l}}\circ e) \otimes \mrm{St}(\rho)$,
\item for $(f,\chi) \in \mathscr{P}'$, the $f$-eigenspace is the sum of $p-1$ copies of $(V_{f,\mathfrak{l}} \circ e) \otimes \pi(\mathbf{1},\chi)(\rho)$,
\item for $(f,\mathbf{1}) \in \mathscr{C}'_1$, the $f$-eigenspace is the sum of $p-1$ copies of $(V_{f,\mathfrak{l}} \circ e)\otimes C_f(\rho)$,
\item for $(f,\mathbf{1}) \in \mathscr{C}_M$, the $f$-eigenspace is the sum of $\frac{p-1}{2}$ copies of 
\begin{enumerate}[noitemsep,label=(\roman*)]
\item $\bigoplus_{r \in \{\pm\}}{\mrm{Ind}_{k(\alpha)}^k{\left[(\psi_{f,r,\mathfrak{l}}\circ e) \otimes C^{r}(\rho_{|G_{k(\sqrt{-p})}})\right]}}$ if $\omega_p(\mrm{Gal}(k_s/k)) \not\subset \F_p^{\times 2}$ or equivalently $\alpha \notin k$ (if $k$ does not have characteristic two, this is equivalent to asking that $\sqrt{-p} \notin k$),
\item $\bigoplus_{\theta \in \{\psi_{f,+,\mathfrak{l}}\circ e,\psi_{f,-,\mathfrak{l}}\circ e\}}{\mrm{Ind}_{k(\alpha)}^k{\left[\theta \otimes C^{r}(\rho_{|G_{k(\sqrt{-p})}})\right]}}$ otherwise (either choice of $r \in \{\pm\}$ leads to the same result). 
\end{enumerate}
\end{itemize}
}

\demo{By Lemma \ref{eigenspaces-t1-modules}, the direct sum of the $\mathbb{T}_1$-eigenspaces for $(f,\chi) \in \mathscr{S} \cup \mathscr{P}' \cup \mathscr{C}'_1 \cup \mathscr{C}_M$ is $\Tate{\ell}{J(p,p)_{\rho}} \otimes_{\Z_{\ell}} F_{\mathfrak{l}}$. 

Let $((f,\chi),\psi) \in (\mathscr{S} \cup \mathscr{P}')\times \mathcal{D}$. Let $V$ be a $(\mathbb{T} \otimes F)$-module such that for all primes $q \neq p$ (resp. all $n \in \F_p^{\times}$), $T_q$ (resp. $nI_2$) acts by $a_q(f)\psi(q)$ (resp. $\chi(n)\psi^2(n)$). Then, for all primes $q \neq p$ and all $n \in \F_p^{\times}$ with $n^2q\equiv 1\pmod{p}$, $(nI_2)\cdot T_q$ acts on $V$ by \[\chi(n)\psi^2(n)a_q(f)\psi(q)=\psi(qn^2)a_q(f)\chi(n)=a_q(f)\chi(n),\] so $V$ is its own $f$-eigenspace as a $(\mathbb{T}_1 \otimes F)$-module. 

Thus, for the eigenspaces attached to elements of $\mathscr{S} \cup \mathscr{P}'$, the claim follows from Corollary \ref{decomposition-xrho}. 

Let $(f,\chi) \in \mathscr{C}$, then $f$ is the twist by some $\psi \in \mathcal{D}$ of a certain $(g,\mathbf{1}) \in \mathscr{C}'_1 \cup \mathscr{C}_M$. Let $V$ be a $(\mathbb{T} \otimes F)$-module such that for all primes $q \neq p$ (resp. all $n \in \F_p^{\times}$), $T_q$ (resp. $nI_2$) acts by $a_q(f)$ (resp. $\chi(n)$). Then, for all primes $q \neq p$ and all $n \in \F_p^{\times}$ with $n^2q\equiv 1\pmod{p}$, $(nI_2)\cdot T_q$ acts on $V$ by $\chi(n)a_q(f)=\psi^2(n)a_q(g)\psi(g)=a_q(g)$, so $V$ is its own $f$-eigenspace as a $(\mathbb{T}_1 \otimes F)$-module. 

The application $((f,\mathbf{1}),\psi) \in \mathscr{C}'_1 \times \mathcal{D} \longmapsto (f\otimes \psi,\psi^2) \in \mathscr{C}$ is injective and its image is exactly the collection of non-CM modular forms. Moreover, if $(f,\mathbf{1}) \in \mathscr{C}'_1 \cup \mathscr{C}_M$ and $\psi \in \mathscr{D}$, then, by Corollaries \ref{twisting-Cf} and \ref{cuspidal-twist}, one has 
\begin{align*}
(V_{f \otimes \psi,\mathfrak{l}} \circ e) \otimes C_{f\otimes \psi}(\rho) &\simeq (V_{f,\mathfrak{l}}\circ e) \otimes (\psi(\omega_p) \circ e) \otimes \psi(\det{\rho}) \otimes C_f(\rho) \\
&\simeq \psi(\omega_p\det{\rho}) \otimes (V_{f,\mathfrak{l}} \otimes e) \otimes C_f(\rho) \simeq (V_{f,\mathfrak{l}}\circ e) \otimes C_f(\rho),
\end{align*}
so by Corollary \ref{decomposition-xrho} the description of the eigenspace is correct for all $f \in \mathscr{C}'_1$. 

Any $(f,\mathbf{1}) \in \mathscr{C}_M$ has exactly $\frac{p-1}{2}$ twists in $\mathscr{C}$, so the $f$-eigenspace is the direct sum of $\frac{p-1}{2}$ copies of $(V_{f,\mathfrak{l}} \circ e) \otimes C_f(\rho)$. When $\alpha \in k$, this representation decomposes as the stated eigenspace by Propositions \ref{decomposition-xrho-reducible-cm}\ref{numberfield-decomp} and \ref{decomp-cusp-sl}. When $\alpha\notin k$, $[(V_{f,\mathfrak{l}} \circ e) \otimes C_f(\rho)]^{\oplus \frac{p-1}{2}}$ decomposes as the stated eigenspace by Proposition \ref{decomposition-xrho-reducible-cm}\ref{cyclotomic-decomp}.}
 
\medskip

\defi[xrho-factors]{Let $\mathfrak{l} \subset \OO_F$ be a maximal ideal with residue characteristic $\ell$. Let $k$ be a field with separable closure $k_s$ and such that $p\ell \in k^{\times}$. Let $\iota: \OO_{\{q\}'} \rar k_s$ be a ring homomorphism, which defines a continuous group homomorphism $e: \mrm{Gal}(k_s/k) \rar \mrm{Gal}(\Q_{\{\ell\}'}/\Q)$. Let $\rho: \mrm{Gal}(k_s/k) \rar \GL{\F_p}$ be a continuous group homomorphism with $\det{\rho}=\omega_p^{-1}$. For $(f,\chi) \in \mathscr{S} \cup \mathscr{P} \cup \mathscr{C}$, let $R_{f,\mathfrak{l}}(\rho)$ denote the representation: 
\begin{itemize}[noitemsep,label=\tiny$\bullet$]
\item $(V_{f,\mathfrak{l}}\circ e) \otimes \mrm{St}(\rho)$ if $(f,\chi) \in \mathscr{S}$, 
\item $(V_{f,\mathfrak{l}}\circ e) \otimes \pi(1,\chi)(\rho)$ if $(f,\chi) \in \mathscr{P}$,
\item $(V_{f,\mathfrak{l}}\circ e) \otimes C_f(\rho)$ if $(f,\chi) \in \mathscr{C}$.
\end{itemize}
The above proof implies that $R_{f,\mathfrak{l}}(\rho)$ does not change if we twist $f$ by a character in $\mathcal{D}$.}

\rem{The homomorphism $e$ is defined up to conjugation, so the isomorphism class of $R_{f,\mathfrak{l}}(\rho)$ does not depend on the choice of $\iota$. However, as we will see, this does not make the choice of $\iota$ irrelevant to our description of the Tate modules.} 

\medskip

\prop[decomposition-xrho-connected-surjective]{Let $k$ be a number field and $\rho: G_k \rar \GL{\F_p}$ be a surjective continuous group homomorphism with $\det{\rho}=\omega_p^{-1}$. Let $\ell$ be a prime number and $\mathfrak{l} \subset \OO_F$ be a prime ideal of residue characteristic $\ell$. For each $\xi \in \mu_p^{\times}(\overline{k})$, $\Tate{\ell}{[\mrm{Jac}(X_{\rho,\xi})]} \otimes_{\Z_{\ell}} F_{\mathfrak{l}}$ is a $(\mathbb{T}_1 \otimes F_{\mathfrak{l}})[G_k]$-module, so its eigenspaces (for $\mathbb{T}_1\otimes F$) are the following irreducible $(\mathbb{T}_1 \otimes F_{\mathfrak{l}})[G_k]$-modules:
\begin{itemize}[noitemsep,label=$-$] 
\item for $(f,\chi) \in \mathscr{S} \cup \mathscr{P}' \cup \mathscr{C}'_1$, the $f$-eigenspace is $R_{f,\mathfrak{l}}(\rho)$,
\item for $(f,\mathbf{1}) \in \mathscr{C}_M$, the $f$-eigenspace is isomorphic to $\mrm{Ind}_{k(\sqrt{-p})}^k{\psi_{f,r_{f,\rho}(\xi),\mathfrak{l}} \otimes C^+(\rho_{|G_{k(\sqrt{-p})}})}$, where $r_{f,\rho}: \xi \in \mu_p^{\times}(\overline{k})/\F_p^{\times 2} \longmapsto r_{f,\rho}(\xi) \in \{\pm\}$ is a bijection that only depends on $f$ and $\rho$. 
\end{itemize}
}

\demo{For each $v \in \F_p^{\times}$, the automorphism $vI_2$ of $X(p,p)$ (which commutes with Hecke correspondences) induces an isomorphism $vI_2: X_{\rho,\xi} \rar X_{\rho,v^2\xi}$, so that the isomorphism class of $\Tate{\ell}{[\mrm{Jac}(X_{\rho,\xi})]} \otimes_{\Z_{\ell}} F_{\mathfrak{l}}$ as a $(\mathbb{T}_1 \otimes F_{\mathfrak{l}})[G_k]$-module only depends on the class of $\xi$ in $\mu_p^{\times}(\overline{k})/\F_p^{\times 2}$. 

For each $v \in \F_p^{\times}$, the automorphism $\Delta_{v,1}$ of $X(p,p)$ (which also commutes with Hecke correspondences) induces an isomorphism $\Delta_{v,1}: (X_{\rho,\xi})_{\overline{k}} \rar (X_{\rho,v\xi})_{\overline{k}}$ respecting the Hecke correspondences, so the isomorphism class of $\Tate{\ell}{[\mrm{Jac}(X_{\rho,\xi})]} \otimes_{\Z_{\ell}} F_{\mathfrak{l}}$ as a $(\mathbb{T}_1 \otimes F_{\mathfrak{l}})$-module does not depend on $\xi$.

By Corollary \ref{decomposition-xrho-cyclotomic}, the eigenspaces of the $(\mathbb{T}_1 \otimes F_{\mathfrak{l}})[G_k]$-module $\bigoplus_{\xi \in \mu_p^{\times}(\overline{k})}{\Tate{\ell}{[\mrm{Jac}(X_{\rho,\xi})]} \otimes_{\Z_{\ell}} F_{\mathfrak{l}}}$ are:
\begin{itemize}[label=\tiny$\bullet$,noitemsep]
\item for $(f,\chi) \in \mathscr{S}\cup \mathscr{P}'\cup\mathscr{C}'_1$, the direct sum of $p-1$ copies of $R_{f,\mathfrak{l}} := R_{f,\mathfrak{l}}(\rho)$,
\item for each $(f,\mathbf{1}) \in \mathscr{C}_M$, the direct sum of $\frac{p-1}{2}$ copies of $R_{f,\mathfrak{l}} := R_{f,+,\mathfrak{l}}\oplus R_{f,-,\mathfrak{l}}$, with $R_{f,\epsilon,\mathfrak{l}} := \mrm{Ind}_{k(\sqrt{-p})}^k{\psi_{f,\epsilon,\mathfrak{l}} \otimes C^+(\rho_{|G_{k(\sqrt{-p})}})}$.
\end{itemize}

The $R_{f,\mathfrak{l}}$ (resp. $R_{f,r,\mathfrak{l}}$) are defined over $F$ in the following sense: for any prime ideal $\mathfrak{q} \subset \OO_k$ such that $\rho$ is unramified at $\mathfrak{q}$ and the residue characteristic of $\mathfrak{q}$ is coprime to $p\mathfrak{l}$, $\mrm{Tr}(\Fr_{\mathfrak{q}} \mid R_{f,\mathfrak{l}})$ (resp. $\mrm{Tr}(\Fr_{\mathfrak{q}} \mid R_{f,r,\mathfrak{l}})$) is in $\OO_F$ and does not depend on $\mathfrak{l}$ (in other words, they form a compatible system with coefficients in $F$ in the sense of Definition \ref{compatible-system-definition}). 

The $R_{f,\mathfrak{l}}$ (for $(f,\chi) \in \mathscr{S} \cup \mathscr{P}' \cup \mathscr{C}'_1$) and the $R_{f,\pm,\mathfrak{l}}$ (for $(f,\chi) \in \mathscr{C}_M$) are irreducible $F_{\mathfrak{l}}[G_k]$-modules by Corollary \ref{decomposition-xrho-irreducible} and Proposition \ref{decomposition-xrho-reducible-cm}.

Let $(f,\chi) \in \mathscr{S} \cup \mathscr{P}' \cup \mathscr{C}'_1$. For each $\xi \in \mu_p^{\times}(\overline{k})$, let $T[f,\xi]$ denote the $f$-eigenspace of $\Tate{\ell}{[\mrm{Jac}(X_{\rho,\xi})]} \otimes_{\Z_{\ell}} F_{\mathfrak{l}}$. Then the $p-1$ different $T[f,\xi]$ have the same dimension over $F_{\mathfrak{l}}$, and their sum is isomorphic to $R_{f,\mathfrak{l}}^{\oplus (p-1)}$, which is a sum of $p-1$ irreducible $(\mathbb{T} \otimes F_{\mathfrak{l}})[G_{k}]$-modules, so we must have $T[f,\xi] \simeq R_{f,\mathfrak{l}}$ for each $\xi$.

Let $(f,\mathbf{1}) \in \mathscr{C}_M$. For each $\xi \in \mu_p^{\times}(\overline{k})$, let $T[f,\xi]$ be the $f$-eigenspace of $\Tate{\ell}{[\mrm{Jac}(X_{\rho,\xi})]} \otimes_{\Z_{\ell}} F_{\mathfrak{l}}$. We know that \[\bigoplus_{\xi}{T[f,\xi]} \simeq \left(R_{f,\mathfrak{l}}\oplus R_{f,-,\mathfrak{l}}\right)^{\oplus (p-1)/2},\] that all the summands on the right-hand side being irreducible, and all the $T[f,\xi]$ have the same dimension over $F_{\mathfrak{l}}$. Thus each $T[f,\xi]$ is isomorphic to $R_{f,r_{f,\mathfrak{l},\xi},\mathfrak{l}}$, for some sign $r_{f,\mathfrak{l},\xi}$. 

Let us fix some $\xi_0 \in \mu_p^{\times}(\overline{k})$, and let $r_{f,\mathfrak{l}}=r_{f,\mathfrak{l},\xi_0}$. Thus $T[f,\xi_0]$ is isomorphic to $R_{f,r_{f,\mathfrak{l}},\mathfrak{l}}$. Since the isomorphism class of $T[f,\xi]$ only depends on the class of $\xi$ in $\mu_p^{\times}(\overline{k})/\F_p^{\times 2}$, every $T[f,u^2\xi_0]$ is isomorphic to $R_{f,r_{\mathfrak{l}},\mathfrak{l}}$. Therefore every $T[f,u\xi_0]$ with $u \notin \F_p^{\times 2}$ is isomorphic to $R_{f,-r_{f,\mathfrak{l}},\mathfrak{l}}$. 

We now need to show that $(r_{f,\mathfrak{l}})_{f \in \mathscr{C}_M}$ does not depend on $\mathfrak{l}$. 

To prove this, let $(r_f)_{f \in \mathscr{C}_M}, (r'_f)_{f \in \mathscr{C}_M}$ be two distinct families of signs, and let us show that there are infinitely many prime ideals $\mathfrak{q} \subset \OO_{k}$ (unramified for $\rho$ and with residual characteristic coprime to $p\ell$) such that the traces of the action of $\Fr_{\mathfrak{q}}$ on $\bigoplus_f{R_{f,r_f,\mathfrak{l}}}$ and $\bigoplus_f{R_{f,r'_f,\mathfrak{l}}}$ has distinct traces. 

Let $M \in \F_p^{\times}\SL{\F_p}$ be a non-scalar matrix with characteristic polynomial $(X-a)^2$ for some $a \in\F_p^{\times}$. Write $\alpha=\mrm{Tr}(M \mid C^+)$ and $\beta=\mrm{Tr}(M \mid C^-)$, so that $\alpha\beta(\alpha-\beta) \in F^{\times}$. 

Let $\mathfrak{q}$ be any prime such that $\rho$ is unramified at $\mathfrak{q}$, $\rho(\Fr_{\mathfrak{q}})\in \F_p^{\times}M$, and the residue characteristic of $\mathfrak{q}$ is coprime to $p\ell$. We compute that 
\begin{align*}
\mrm{Tr}(\Fr_{\mathfrak{q}} \mid \bigoplus_{f}{R_{f,r_f,\mathfrak{l}}}) &= \alpha\sum_{f}{\psi_{f,r_f,\mathfrak{l}}(\Fr_{\mathfrak{q}})} + \beta \sum_f{\psi_{f,-r_f,\mathfrak{l}}(\Fr_{\mathfrak{q}})}\\
\mrm{Tr}(\Fr_{\mathfrak{q}} \mid \bigoplus_{f}{R_{f,r'_f,\mathfrak{l}}}) &= \alpha\sum_{f}{\psi_{f,r'_f,\mathfrak{l}}(\Fr_{\mathfrak{q}})} + \beta \sum_f{\psi_{f,-r'_f,\mathfrak{l}}(\Fr_{\mathfrak{q}})}.
\end{align*}

If, for all but finitely many $\mathfrak{q}$ as above, these two traces are equal, one then has 

\[\alpha\sum_f{\left[\psi_{f,r_f,\mathfrak{l}}(\Fr_{\mathfrak{q}})-\psi_{f,r'_f,\mathfrak{l}}(\Fr_{\mathfrak{q}})\right]} = \beta \sum_f{\left[\psi_{f,-r_f,\mathfrak{l}}(\Fr_{\mathfrak{q}})-\psi_{f,-r'_f,\mathfrak{l}}(\Fr_{\mathfrak{q}})\right]}\] 

for all but finitely many $\mathfrak{q}$ as above. Fix some $g \in G_k$ such that $\rho(g)=M$. Then, by Cebotarev's theorem, for any $h \in H := \ker{\rho: G_k \rar \PGL{\F_p}} \triangleleft G_{k(\sqrt{-p})}$,  

\[\alpha\sum_f{\left[\psi_{f,r_f,\mathfrak{l}}(gh)-\psi_{f,r'_f,\mathfrak{l}}(gh)\right]} = \beta \sum_f{\left[\psi_{f,-r_f,\mathfrak{l}}(gh)-\psi_{f,-r'_f,\mathfrak{l}}(gh)\right]}.\]

Since $G_{k(\sqrt{-p})}/H \simeq \PSL{\F_p}$ is a perfect group, $H \rar G_{k(\sqrt{-p})}^{ab}$ is surjective. Thus, the following identity holds on $G_{k(\sqrt{-p})}^{ab}$:

\[\alpha\sum_f{\left[\psi_{f,r_f,\mathfrak{l}}-\psi_{f,r'_f,\mathfrak{l}}\right]} = \beta \sum_f{\left[\psi_{f,-r_f,\mathfrak{l}}-\psi_{f,-r'_f,\mathfrak{l}}\right]}.\]

The $\psi_{f,r,\mathfrak{l}}$ (over all possible $f$ and $r$) are pairwise distinct characters (because, by \cite[Theorem 4.3]{Antwerp5-Ribet}, $\frac{\psi_{f,+,\mathfrak{l}}}{\psi_{f,-,\mathfrak{l}}}:G_{\Q(\sqrt{-p})} \rar F_{\mathfrak{l}}^{\times}$ has infinite image), and $\alpha\beta(\alpha-\beta) \in F^{\times}$, this implies that $r_f=r'_f$ for all $f$.

Now, let $\mathfrak{l},\mathfrak{l}' \subset \OO_F$ be two distinct maximal ideals with residue characteristics $\ell,\ell'$, $\xi \in \mu_p^{\times}(\overline{k})$. Let $\mathfrak{q} \subset \OO_F$ be any prime ideal such that $\rho$ is unramified at $\mathfrak{q}$ and $\mathfrak{q}$ is coprime to $p\ell\ell'$. We can then compute in $F$ that 

\begin{align*}
\mrm{Tr}(\Fr_{\mathfrak{q}} \mid \bigoplus_f{R_{f,r_{f,\mathfrak{l},\xi},\mathfrak{l}}}) &= \mrm{Tr}(\Fr_{\mathfrak{q}}\mid \Tate{\ell}{\mrm{Jac}(X_{\rho,\xi})}) - \sum_{(f,\chi) \in \mathscr{S} \cup \mathscr{P}' \cup \mathscr{C}'_1}{\mrm{Tr}(\Fr_{\mathfrak{q}} \mid R_{f,\mathfrak{l}})}\\
&= \mrm{Tr}(\Fr_{\mathfrak{q}}\mid \Tate{\ell'}{\mrm{Jac}(X_{\rho,\xi})}) - \sum_{(f,\chi) \in \mathscr{S} \cup \mathscr{P}' \cup \mathscr{C}'_1}{\mrm{Tr}(\Fr_{\mathfrak{q}} \mid R_{f,\mathfrak{l}'})}\\
&= \mrm{Tr}(\Fr_{\mathfrak{q}} \mid \bigoplus_f{R_{f,r_{f,\mathfrak{l}',\xi},\mathfrak{l}'}}) = \mrm{Tr}(\Fr_{\mathfrak{q}} \mid \bigoplus_f{R_{f,r_{f,\mathfrak{l}',\xi},\mathfrak{l}}}).
\end{align*}

By the above argument, this identity implies that \[(r_{f,\mathfrak{l},\xi})_{(f,\mathbf{1}) \in \mathscr{C}_M} = (r_{f,\mathfrak{l}',\xi})_{(f,\mathbf{1}) \in \mathscr{C}_M},\] whence the conclusion. 
}

\bigskip

To be able to handle a representation $\rho$ over an arbitrary field, we need to keep track of the action of $\GL{\F_p}$ compatibly with the action of $G_{\Q}$. To do that, we need ways to describe representations of profinite groups when they are only well-known on open subgroups. 

\lem[make-tensor-product-irreducible]{Let $G$ be a profinite group and $K \triangleleft G$ be an open subgroup. Let $V_1,V_2$ be two finite-dimensional continuous representations of $G$ over a $p$-adic field $F$. Assume that $(V_1)_{|K}$ is absolutely irreducible and $(V_2)_{|K}$ is a direct sum of copies of $(V_1)_{|K}$. Then there is a representation $W$ of $G/K$ over $F$ such that $V_2 \simeq V_1 \otimes W$. }

\demo{By the assumptions, the $F[G/K]$-module $W := \mrm{Hom}_{F[K]}(V_1,V_2)$ has dimension $\frac{\dim{V_2}}{\dim{V_1}}$ over $F$, and the obvious map $V_1 \otimes W \rar V_2$ (which is $G$-equivariant) is onto. Since it is a linear morphism between $F$-vector space of the same finite dimension, it is an isomorphism.  }

\medskip

\lem[make-tensor-product-induced]{Let $G$ be a profinite group, $G_2 \triangleleft G$ be an open subgroup of index two, and $K \triangleleft G$ an open subgroup contained in $G_2$. Let $F$ be a $p$-adic field and $\psi, \psi': G_2 \rar F^{\times}$ be continuous characters, conjugate under $G/G_2$, such that $\psi^{-1}\psi'$ has infinite image. Let $V$ be a finite-dimensional continuous representation of $G$ over $F$ such that $V_{|K}$ is isomorphic to a direct sum of copies of $\left[\mrm{Ind}_{G_2}^G{\psi}\right]_{|K}$. Then there exists a continuous finite-dimensional representation $W$ of $G_2/K$ over $F$ such that $V \simeq \mrm{Ind}_{G_2}^G{\left[\psi \otimes W\right]}$.}

\demo{Let $W=\mrm{Hom}_{F[K]}(\psi,V)$, then $W$ is a finite-dimensional $F[G_2/K]$-module and its dimension over $F$ is exactly $\frac{\dim{V}}{2}$. The canonical morphism $u: \psi \otimes W \rar V$ is $F[G_2]$-linear, and its image is exactly the $\psi$-isotypic component $V_2$ of $V_{|K}$. Then the induced morphism $f: \mrm{Ind}_{G_2}^G{\left[\psi \otimes W\right]} \rar V$ is $F[G]$-equivariant between finite-dimensional $F$-vector spaces of same dimension. Its image contains the sub-$G$-module generated by $V_2$, which is $V$, so $f$ is onto, so $f$ is an isomorphism.} 

\medskip

\lem[twist-mix-gl2-galois]{Let $k$ be a field with characteristic not $p$ and $\rho: \mrm{Gal}(k_s/k) \rar \GL{\F_p}$ be a continuous group homomorphism with $\det{\rho}=\omega_p^{-1}$. Then the actions of $\GL{\F_p}$ and $\mrm{Gal}(k_s/k)$ by $\mathbb{T}$-endomorphisms of $J_{\rho}(p)(k_s)$ factor through the group $\GL{\F_p} \rtimes_{\rho}\mrm{Gal}(k_s/k)$, whose product law is $(M,g)(M',g')=(M\rho(g)M'\rho(g)^{-1},gg')$.  }

\demo{For any $N \in \GL{\F_p}$ and any $\sigma \in \mrm{Gal}(k_s/k)$, the following diagrams are commutative by Proposition \ref{cocycle-twist}, and the actions of $\mrm{Gal}(k_s/k)$ and $\GL{\F_p}$ commute on $J(p,p)(k_s)$. 
\[
\begin{tikzcd}[ampersand replacement=\&]
J_{\rho}(p)(k_s)\arrow{d}{N} \arrow{r}{j} \& J(p,p)(k_s) \arrow{d}{N} \& J_{\rho}(p)(k_s) \arrow{r}{j} \arrow{d}{\sigma} \& J(p,p)(k_s) \arrow{d}{\rho(\sigma)\circ \sigma}\\
J_{\rho}(p)(k_s) \arrow{r}{j} \& J(p,p)(k_s) \& J_{\rho}(p)(k_s) \arrow{r}{j} \& J(p,p)(k_s).
\end{tikzcd}
\]

Thus, for any $M \in \GL{\F_p}$ and $\sigma \in \mrm{Gal}(k_s/k)$, the $\mathbb{T}$-endomorphism of $J_{\rho}(p)(k_s)$ induced by $\sigma \circ M$ is the same as the one induced by $\rho(\sigma)M\rho(\sigma)^{-1} \sigma$, and the conclusion follows.}

\medskip

\prop[decomposition-xrho-connected-surjective+sl2]{In the situation of Proposition \ref{decomposition-xrho-connected-surjective}, write $G_{k,\rho} := \SL{\F_p} \rtimes_{\rho} \mrm{Gal}(\overline{k}/k)$ for the sake of easier notation and $\tilde{\rho}$ be the homomorphism $(M,\sigma) \in G_{k,\rho} \mapsto M\rho(\sigma) \in \GL{\F_p}$. It is endowed with a natural morphism to $G_k$, which we omit from the notation. When $p \equiv 3 \pmod{4}$, we denote by $G_2$ the subgroup $\SL{\F_p} \rtimes_{\rho} \mrm{Gal}(\overline{k}/k(\sqrt{-p})) = \tilde{\rho}^{-1}(\F_p^{\times}\SL{\F_p})$. 

For each $\xi \in \mu_p^{\times}(\overline{k})$, $\Tate{\ell}{[\mrm{Jac}(X_{\rho,\xi})]} \otimes_{\Z_{\ell}} F_{\mathfrak{l}}$ is a $(\mathbb{T}_1 \otimes F_{\mathfrak{l}})[G_{k,\rho}]$-module, so its eigenspaces (for $\mathbb{T}_1\otimes F$) are the following irreducible $(\mathbb{T}_1 \otimes F_{\mathfrak{l}})[G_{k,\rho}]$-modules:
\begin{itemize}[noitemsep,label=$-$] 
\item for $(f,\chi) \in \mathscr{S} \cup \mathscr{P}' \cup \mathscr{C}'_1$, the $f$-eigenspace is $V_{f,\mathfrak{l}} \otimes \mrm{St}(\tilde{\rho})$, 
\item for $(f,\chi) \in \mathscr{P}'$, the $f$-eigenspace is $V_{f,\mathfrak{l}} \otimes \pi(\mathbf{1},\chi)(\tilde{\rho})$,
\item for $(f,\mathbf{1}) \in \mathscr{C}'_1$, the $f$-eigenspace is $V_{f,\mathfrak{l}} \otimes C_f(\tilde{\rho})$,
\item for $(f,\mathbf{1}) \in \mathscr{C}_M$, there exists a bijection $r_{f,\rho}: \xi \in \mu_p^{\times}(\overline{k})/\F_p^{\times 2} \longmapsto r_f(\xi) \in \{\pm\}$ (which is independent from $\mathfrak{l}$) such that the $f$-eigenspace is isomorphic to $\mrm{Ind}_{G_2}^G{\psi_{f,r,\mathfrak{l}} \otimes C^+(\tilde{\rho}_{|G_2})}$. 
\end{itemize}
}

\demo{For $(f,\chi) \in \mathscr{S}$ (resp. $(f,\chi) \in \mathscr{P}'$, resp. $(f,\chi) \in \mathscr{C}'_1$), let $E_{f,\mathfrak{l},\xi}$ be the $f$-eigenspace of $\Tate{\ell}{[\mrm{Jac}(X_{\rho,\xi})]} \otimes_{\Z_{\ell}} F_{\mathfrak{l}}$ as a $(\mathbb{T}_1 \otimes F_{\mathfrak{l}})[G_{k,\rho}]$-module. Let furthermore $R_{f,\mathfrak{l}}(\tilde{\rho})$ denote the $(\mathbb{T}_1 \otimes F_{\mathfrak{l}})[G_{k,\rho}]$-module $V_{f,\mathfrak{l}} \otimes \mrm{St}(\tilde{\rho})$ (resp. $V_{f,\mathfrak{l}} \otimes \pi(1,\chi)(\tilde{\rho})$, resp. $V_{f,\mathfrak{l}} \otimes C_f(\tilde{\rho})$). 

By Proposition \ref{decomposition-xrho-connected-surjective}, all the $E_{f,\mathfrak{l},\xi}$ are isomorphic, as $(\mathbb{T}_1 \otimes F_{\mathfrak{l}})[G_k]$-modules, to $R_{f,\mathfrak{l}}(\tilde{\rho})$, and both are irreducible as $F_{\mathfrak{l}}[G_k]$-modules. In particular, the $E_{f,\mathfrak{l},\xi}$ and the $R_{f,\mathfrak{l}}(\tilde{\rho})$ are all irreducible as $(\mathbb{T}_1 \otimes F_{\mathfrak{l}})[G_{k,\rho}]$-modules. 

The direct sum of the $E_{f,\mathfrak{l},\xi}$ is the $f$-eigenspace of $\Tate{\ell}{J_{\rho}(p)} \otimes_{\Z_{\ell}} F_{\mathfrak{l}}$, which is isomorphic over $(\mathbf{T}_1 \otimes_{\Z} F_{\mathfrak{l}})[G_{k,\rho}]$ by Corollaries \ref{decomposition-xpp-tate} and \ref{decomposition-xrho-cyclotomic} to $R_{f,\mathfrak{l}}(\tilde{\rho})^{\oplus (p-1)}$. By uniqueness of the decomposition into irreducible factors, for every $\xi$, one has $E_{f,\mathfrak{l},\xi} \simeq R_{f,\mathfrak{l}}(\tilde{\rho})$. 

Let now $(f,\mathbf{1}) \in \mathscr{C}_M$. For $r \in \{\pm\}$, let $R_{f,r,\mathfrak{l}}(\tilde{\rho}) = \mrm{Ind}_{G_2}^G{\psi_{f,r,\mathfrak{l}} \otimes C^+(\tilde{\rho}_{|G_2})}$ as a $(\mathbb{T}_1\otimes F_{\mathfrak{l}})[G_{k,\rho}]$-module (where $(mI_2) \cdot T_n$ acts by $a_n(f)$ for any $n,m \geq 1$ with $nm^2\equiv 1\pmod{p}$). Let $E_{f,\mathfrak{l},\xi}$ be, for any $\xi \in \mu_p^{\times}(\overline{k})$, the $f$-eigenspace of $\Tate{\ell}{J_{\rho,\xi}(p)}$ as a $(\mathbb{T}_1\otimes F_{\mathfrak{l}})[G_{k,\rho}]$-module. 

We know by Proposition \ref{decomposition-xrho-connected-surjective} that $E_{f,\mathfrak{l},\xi}$ and $R_{f,r_{f,\rho}(\xi),\mathfrak{l}}(\tilde{\rho})$ are isomorphic when restricted to $\{1\} \times G_k$, and that they are irreducible as $F_{\mathfrak{l}}[G_k]$-modules. Therefore, the $E_{f,\mathfrak{l},\xi}$ and $R_{f,r,\mathfrak{l}}(\tilde{\rho})$ are irreducible $F_{\mathfrak{l}}[G_{k,\rho}]$-modules. 

Moreover, $\bigoplus_{\xi \in \mu_p^{\times}(\overline{k})}{E_{f,\mathfrak{l},\xi}}$ is the $f$-eigenspace of $\Tate{\ell}{J_{\rho}(p)} \otimes_{\Z_{\ell}} F_{\mathfrak{l}}$, which is isomorphic by Corollaries \ref{decomposition-xpp-tate} and \ref{decomposition-xrho-cyclotomic} to 

\[V_{f,\mathfrak{l}} \otimes C_f(\tilde{\rho})^{\oplus (p-1)} \simeq \mrm{Ind}_{G_2}^{G_{k,\rho}}{(\psi_{f,+,\mathfrak{l}} \oplus \psi_{f,-,\mathfrak{l}}) \otimes C^+(\tilde{\rho})^{\oplus \frac{p-1}{2}}} \simeq R_{f,+,\mathfrak{l}}(\tilde{\rho})^{\oplus \frac{p-1}{2}}\oplus R_{f,-,\mathfrak{l}}(\tilde{\rho})^{\oplus \frac{p-1}{2}}.\]

Therefore, there exists a morphism $r: \mu_p^{\times}(k_s) \rar \{\pm\}$ such that all of its fibres have cardinality $\frac{p-1}{2}$ and every $E_{f,\mathfrak{l},\xi}$ is isomorphic to one $R_{f,r(\xi),\mathfrak{l}}(\tilde{\rho})$. For any $n \in \F_p^{\times}$, $nI_2: E_{f,\mathfrak{l},\xi} \rar E_{f,\mathfrak{l},\xi^{n^2}}$ is an isomorphism of $(\mathbb{T}_1 \otimes F_{\mathfrak{l}})[G_{k,\rho}]$-modules; therefore, we may choose $r$ so that it is a bijection $\mu_p^{\times}(k_s)/\F_p^{\times 2} \rar \{\pm 1\}$. 

}

\bigskip

\lem[twist-xrho-sl2]{Let $k$ be a field with characteristic distinct from $p$ and $k_s$ be a separable closure for $k$. Let $\rho,\rho_0: \mrm{Gal}(k_s/k) \rar \GL{\F_p}$ be two continuous group homomorphisms with $\det{\rho}=\det{\rho_0}=\omega_p^{-1}$. Then the function $\rho\rho_0^{-1}: \mrm{Gal}(k_s/k) \rar \SL{\F_p}$ is a cocycle in the sense of Proposition \ref{cocycle-twist}; for any $\xi \in \mu_p^{\times}(k_s)$, the twist of $X_{\rho_0,\xi}(p)$ is $X_{\rho,\xi}(p)$, and this identification is compatible with the action of $\mathbb{T}_{1,\{1\}}$ on the Jacobians. }

\demo{For any $\sigma \in \mrm{Gal}(k_s/k)$, consider the following diagram:
\[
\begin{tikzcd}[ampersand replacement=\&]
X_{\rho_0}(p)_{k_s} \arrow{rr}{j_{\rho_0}} \arrow[swap]{d}{\mrm{id}\times \Sp{\sigma^{-1}}} \& \& X(p,p)_{k_s} \arrow{r}{j_{\rho}} \arrow[swap]{dl}{\rho_0(\sigma)\times \Sp{\sigma^{-1}}} \arrow{dd}{\rho(\sigma) \times \Sp{\sigma^{-1}}}\& X_{\rho}(p)_{k_s} \arrow{dd}{\mrm{id} \times \Sp{\sigma^{-1}}}\\
X_{\rho_0}(p)_{k_s} \arrow{r}{j_{\rho_0}} \arrow[swap]{d}{\rho(\sigma)\rho_0(\sigma)^{-1}} \& X(p,p)_{k_s} \arrow[swap]{rd}{\rho(\sigma)\rho_0(\sigma)^{-1}}\& \& \\
X_{\rho_0}(p)_{k_s} \arrow[swap]{rr}{j_{\rho_0}} \& \& X(p,p)_{k_s} \arrow{r}{j_{\rho}} \& X_{\rho}(p)_{k_s}
\end{tikzcd}
\]

This diagram implies that $\sigma \in \mrm{Gal}(k_s/k) \longmapsto \rho(\sigma)\rho_0(\sigma)^{-1} \in \SL{\F_p} \rar \mrm{Aut}_{k_s}(X_{\rho_0}(p))$ satisfies the cocycle relation of Proposition \ref{cocycle-twist}, and that $X_{\rho}(p)$ is the twist of $X_{\rho_0}(p)$ by this cocycle. Moreover, by Lemma \ref{T1Gamma-action}, because it takes its values in $\SL{\F_p}$, this cocycle leaves invariant the morphism to the constant $k$-scheme with underlying set $\mu_p^{\times}(k_s)$ defined in Proposition \ref{group-is-twist-polarized}. Thus the twist of $X_{\rho_0}(p) \rar \mu_p^{\times}(k_s)$ by the cocycle $\rho\rho_0^{-1}$ is exactly $X_{\rho}(p) \rar \mu_p^{\times}(k_s)$. 

Therefore, for every $\xi \in \mu_p^{\times}(k_s)$, $X_{\rho,\xi}(p)$ is the twist of $X_{\rho_0,\xi}(p)$ by the cocycle \[\rho\rho_0^{-1}: \mrm{Gal}(k_s/k) \rar \SL{\F_p} \subset \mrm{Aut}_{k_s}(X_{\rho_0,\xi}(p)).\] 

This twist preserves (by definition) the geometric actions of $\GL{\F_p}$, and by Propositions \ref{jacobian-twist}, \ref{jacobian-twist-functoriality}, it also preserves the actions of $\mathbb{T}$.
}

\medskip

\cor[decomposition-xrho-connected-sl20]{Let $\mathfrak{l} \subset \OO_F$ be a maximal ideal with residue characteristic $\ell$. Let $k_s$ be a separable closure of $k$ and $\rho: \mrm{Gal}(k_s/k) \rar \GL{\F_p}$ such that $\det{\rho}=\omega_p^{-1}$. Let $k$ be a field of characteristic zero with separable closure $\iota: \Qbar \rar k_s$ be a ring homomorphism, inducing a continuous group homomorphism $R_{\iota}: \mrm{Gal}(k_s/k) \rar \mrm{Gal}(\Qbar/\Q)$. Let $G_{k,\rho}$ be the group $\SL{\F_p} \rtimes_{\rho} \mrm{Gal}(k_s/k)$, which is endowed with a group homomorphism $\tilde{\rho}: (M,\sigma) \in G_{k,\rho} \longmapsto M\rho(\sigma)$. 

Then, for $(f,\chi) \in \mathscr{S}\cup\mathscr{P}'\cup\mathscr{C}'_1\cup\mathscr{C}_M$ and $\xi \in \mu_p^{\times}(k_s)$, the $f$-eigenspace of $\Tate{\ell}{J_{\rho,\xi}(p)}\otimes_{\Z_{\ell}} F_{\mathfrak{l}}$ is isomorphic as a $F_{\mathfrak{l}}[G_{k,\rho}]$-module to the following:

\begin{itemize}[noitemsep,label=\tiny$\bullet$]
\item $V_{f,\mathfrak{l}}(R_{\iota}) \otimes \mrm{St}(\tilde{\rho})$ when $(f,\chi) \in \mathscr{S}$, 
\item $V_{f,\mathfrak{l}}(R_{\iota}) \otimes \pi(1,\chi)(\tilde{\rho})$ when $(f,\chi) \in \mathscr{P}'$,
\item $V_{f,\mathfrak{l}}(R_{\iota}) \otimes C_f(\tilde{\rho})$ when $(f,\chi) \in\mathscr{C}'_1$,
\item when $(f,\chi) \in \mathscr{C}_M$, let $k'=k(\sqrt{-p})\subset k_s$, and $G_{k',\rho}=\SL{\F_p} \rtimes_{\rho} \mrm{Gal}(k_s/k')$. There is a unique $r \in \{\pm\}$ depending only in $\iota^{-1}(\xi)$ such that the eigenspace is isomorphic to:
\begin{itemize}[noitemsep,label=$-$]
\item $\mrm{Ind}_{G_{k',\rho}}^{G_{k,\rho}}{\psi_{f,r,\mathfrak{l}}(R_{\iota}) \otimes C^+(\tilde{\rho}_{|G_{k',\rho}})}$ if $k' \neq k$,
\item $\bigoplus_{s \in \{\pm\}}{\psi_{f,rs,\mathfrak{l}}(R_{\iota}) \otimes C^s(\tilde{\rho})}$ if $k=k'$. 
\end{itemize}
\end{itemize}
}

\demo{If $\overline{L}$ is a separable closure of some field extension $L/k$ endowed with a $k$-embedding $\iota': k_s \rar \overline{L}$, there is a natural a group homomorphism $R_{\iota'}: \mrm{Gal}(\overline{L}/L) \rar \mrm{Gal}(k_s/k)$. For any $\zeta \in \mu_p^{\times}(k_s)$, it is then direct to see that $X_{\rho,\zeta}(p)_L$ is naturally isomorphic to $X_{\rho\circ R_{\iota'},\iota'(\zeta)}(p)_L$. In particular, if the description is correct for $\rho$, it is correct for $\rho\circ R_{\iota'}$.  

Therefore, by Lemma \ref{twist-xrho-sl2}, to prove the Proposition, it is enough to find a single representation $\rho: G_{\Q} \rar \GL{\F_p}$ with determinant $\omega_p^{-1}$ for which the claim is correct. Thanks to Proposition \ref{decomposition-xrho-connected-surjective+sl2}, it is enough to find a single surjective representation $\rho: \mrm{Gal}(\Qbar/\Q) \rar \GL{\F_p}$ with $\det{\rho}=\omega_p^{-1}$. 

By \cite[Cor. 1 \`a Prop. 21]{Serre-image-ouverte}, if $E/\Q$ is a semi-stable elliptic curve with good reduction at $2$ -- such as the elliptic curve with LMFDB label \cite{lmfdb} $11.a3$ and Weierstrass equation $y^2-y=x^3-x^2$ -- then, for any basis $(P,Q)$ of $E[p]$, the group homomorphism $\sigma \longmapsto \rho(\sigma)$ (defined for any $\sigma \in G_{\Q}$ by $\rho(\sigma)\begin{pmatrix}\sigma(P)\\\sigma(Q)\end{pmatrix} = \begin{pmatrix}P\\Q\end{pmatrix}$ -- it is a group homomorphism by Lemma \ref{right-representation}) satisfies the requirements.  
}

\rem[distinguishing-cm-factors]{Let $E/\Q$ be an elliptic curve and fix a basis $(P,Q)$ of its $p$-torsion module. Let $\zeta = \langle P,\,Q\rangle \in \mu_p^{\times}(\Qbar)$ and $\rho: G_{\Q} \rar \GL{\F_p}$ the group homomorphism given by 
\[\forall g \in G_{\Q},\, \begin{pmatrix}g(P)\\g(Q)\end{pmatrix} = \rho(g)^{-1}\begin{pmatrix}P\\Q\end{pmatrix}. \]
For any $\alpha \in \F_p^{\times}$, $X_E^{\alpha}(p)$ identifies with $X_{\rho,\zeta^{\alpha}}(p)$ in a way that respects the Hecke correspondences. Can one describe the CM eigenspaces of the Tate module of $\mrm{Jac}(X_E^{\alpha}(p)) \otimes_{\Z_{\ell}} F_{\mathfrak{l}}$ in a way which is independent from the choice of $\rho$? 

A possible approach could be to reparametrize the CM eigenspaces in the following more canonical way. First, the two compatible systems of characters $\psi_{f,\pm,\mathfrak{l}}$ can be parametrized in a choice-free way by the two embeddings $\Q(\sqrt{-p}) \rar F$. Second, a completely choice-free parametrization for the set of representations $\{C^+,C^-\}$ of $\F_p^{\times}SL(V)$ (where $V$ is a two-dimensional vector space over $\F_p$, with $p \equiv 3\pmod{4}$) with coefficients in $F$ is by the pairs $([U],\xi)$, where $[U]$ represents a conjugacy class of unipotent element in $SL(V)$ and $\xi \in \mu_p^{\times}(F)/\F_p^{\times 2}$ corresponds to the eigenvalues of the representation at any element in $[U]$. Note that the pair $([U],\xi)$ describes the same representation as $([U^{-1}],\xi^{-1})$. 

Let $\rho_E: G_{\Q} \rar GL(E[p])$ be the Galois representation attached to $E$. Fix any basis $(P,Q)$ of $E[p](\Qbar)$. It defines the conjugacy class $[U]$ of the unipotent automorphism $(P,Q) \mapsto (P,P+Q)$ of $E[p]$, and the root of unity $\xi = \langle P,\,Q\rangle\in \mu_p^{\times}(\Qbar)$. Then $\xi^{\F_p^{\times 2}}$ is defined over\footnote{Note how the following inclusion is necessary to even define $\psi_{\sigma}$.} $\Q(\sqrt{-p}) \subset \Qbar$, so, for any embedding $\sigma: \Q(\sqrt{-p}) \rar F$, it makes sense to define $\sigma(\xi)$ as an element of $\mu_p^{\times}(F)/\F_p^{\times 2}$. Thus, for $(f,\mathbf{1}) \in \mathscr{C}_M$, the $f$-eigenspaces of $\Tate{\ell}{\mrm{Jac}(X_E^{1}(p))} \otimes_{\Z_{\ell}} F_{\mathfrak{l}}$ and $\Tate{\ell}{\mrm{Jac}(X_E^{-1}(p))} \otimes_{\Z_{\ell}} F_{\mathfrak{l}}$ are given (up to re-ordering) by the couple 
\[\mathcal{E}_{f,E} := \Big(\mrm{Ind}_{\Q(\sqrt{-p})}^{\Q}{\left[\psi_{\sigma} \otimes C_{[U],\sigma(\xi)}((\rho_E)_{|G_{\Q(\sqrt{-p})}})\right]},\, \mrm{Ind}_{\Q(\sqrt{-p})}^{\Q}{\left[\psi_{\sigma} \otimes C_{[U],\sigma(\xi)^{-1}}((\rho_E)_{|G_{\Q(\sqrt{-p})}})\right]}\Big),\]
where $\sigma: \Q(\sqrt{-p}) \rar F$ is an embedding. Note that $\mathcal{E}_{f,E}$ is independent from the choices of $\sigma$ and of the basis $(P,Q)$. Hence, it seems logical (and can probably be proved in the same way as Corollary \ref{decomposition-xrho-connected-sl20}) that one of the following is true: 
\begin{itemize}[noitemsep,label=$-$]
\item either, for every elliptic curve $E/\Q$, $\mathcal{E}_{f,E}$ is equal to the couple made with the $f$-eigenspaces of $\Tate{\ell}{\mrm{Jac}(X_E^{1}(p))} \otimes_{\Z_{\ell}} F_{\mathfrak{l}}$ and $\Tate{\ell}{\mrm{Jac}(X_E^{-1}(p))} \otimes_{\Z_{\ell}} F_{\mathfrak{l}}$ in this order,
\item or, for every elliptic curve $E/\Q$, $\mathcal{E}_{f,E}$ is equal to the couple made with the $f$-eigenspaces of $\Tate{\ell}{\mrm{Jac}(X_E^{-1}(p))} \otimes_{\Z_{\ell}} F_{\mathfrak{l}}$ and $\Tate{\ell}{\mrm{Jac}(X_E^{1}(p))} \otimes_{\Z_{\ell}} F_{\mathfrak{l}}$ in this order.
\end{itemize}

It seems numerically possible (when $p$ is small enough for $X_{E}^{\alpha}(p)$ to have an explicit model, for a given elliptic curve $E$) to decide which alternative is correct, but we were not able to find any general argument to privilege one choice over the other in general. In particular, is there any reason to believe that the same alternative is correct, no matter what $(f,\mathbf{1}) \in \mathscr{C}_M$ or $p$ are? 
}

\bigskip

We now want to prove the equivalent of Corollary \ref{decomposition-xrho-connected-sl20} in positive characteristic. The same argument as the previous proof shows that it is enough to prove the result for a single representation $\rho: \mrm{Gal}(\overline{\F_q}/\F_q) \rar \GL{\F_p}$. Using the ideas of Proposition \ref{decomposition-xrho-connected-surjective+sl2} does not work for several reasons: because $\mrm{Gal}(\overline{\F_q}/\F_q)$ is pro-cyclic, $(V_{f,\mathfrak{l}})_{|\mrm{Gal}(\overline{\F_q}/\F_q)}$ is never absolutely irreducible, and $\rho$ cannot be surjective because $\GL{\F_p}$ is not solvable. 

The idea is to take for $\rho$ the restriction of a global representation to a decomposition subgroup, and this idea is better expressed using the relative curves $X_{G,u}(p)$ as in Section \ref{application-moduli} and fundamental groups. First, we carefully specify the relationship between Galois groups and fundamental groups in the following Lemma (which is a simple verification). 

\medskip

\defi[fundamental-group]{Let $X$ be a scheme and $\text{\bf{Fin\'Et}}_X$ denote the category of finite \'etale surjective $X$-schemes. Let $\overline{x}: \Sp{\Omega} \rar X$ be a geometric point (that is, $\Omega$ is an algebraically closed field), the \emph{fundamental group} $\pi_1(X,\overline{s})$ is the group of automorphisms of the functor $\text{\bf{Fin\'Et}}_X \rar \mathbf{FinSet}$ given by $U \longmapsto U_{\overline{x}}$. }

\lem[pi1-vs-galois]{Let $k$ be a field with separable closure $k_s$. For every finite \'etale $k$-scheme $Y$, let $Y_{k_s}$ denote the underlying set of the finite \'etale $k_s$-scheme $Y \times_k \,\Sp{k_s}$, and let $\iota_Y$ denote the bijection $Y(k_s) \rar Y \times_k \Sp{k_s}$ defined as follows: for $f: \Sp{k_s} \rar Y$, $\iota_Y(f)$ is the image of the morphism $\Sp{k_s} \overset{(f,\mrm{id})}{\longrightarrow} Y \times_k\,\Sp{k_s}$. Then the collection of the $\iota_Y$ is an isomorphism of functors $\text{\bf{Fin\'Et}}_k \rar \mathbf{FinSet}$. Moreover, for any finite \'etale $k$-scheme $Y$ and every $\sigma \in \mrm{Gal}(k_s/k)$, the following diagram commutes:
\[
\begin{tikzcd}[ampersand replacement=\&]
Y(k_s) \arrow{rr}{P \longmapsto (P \circ \Sp{\sigma}) } \arrow{d}{\iota_Y}\& \& Y(k_s) \arrow{d}{\iota_Y}\\
Y \times_k \Sp{k_s} \arrow{rr}{\mrm{id} \times \Sp{\sigma^{-1}}}\& \& Y \times_k \Sp{k_s}
\end{tikzcd}
\]

In particular, if $\alpha$ is the geometric point $\Sp{\Omega} \rar \Sp{k}$ for some algebraically closed extension $\Omega$ of $k_s$, this defines a morphism $\mrm{Gal}(k_s/k) \rar \pi_1(\Sp{k},\alpha)$: for every finite \'etale $k$-scheme $Y$, any $\sigma \in \mrm{Gal}(k_s/k)$ acts on the set of points $Y_{\alpha}$ (which is the same as the set of points of $Y_{k_s}$) by $\mrm{id} \times \Sp{\sigma^{-1}}$. }

\medskip

\rem{In the situation of Definition \ref{fundamental-group}, let $x \in X$ be the image of $\overline{x}$, $\kappa(x)$ denote its residue field and $\overline{\kappa(x)}$ be the separable closure of $\kappa(x)$ in $\Omega$. Lemma \ref{pi1-vs-galois} means that $g \in \mrm{Gal}(\overline{\kappa(x)}/\kappa(x)) \longmapsto (\mrm{id}\times \Sp{g^{-1}}) \in \pi_1(X,\overline{x})$ is a group homomorphism and that it is compatible with the expected Galois actions on $X(\overline{\kappa(x)})$. }

\medskip

The following Lemma describes how to apply the formalism of fundamental groups in our situation.

\lem[xgp-family-pi1-elem]{Let $\ell$ be a prime number, $\mathfrak{l}$ be a prime ideal of $\OO_F$ containing $\ell$. Let $S$ be an affine, connected, regular excellent Noetherian $\Z[(p\ell)^{-1}]$-scheme. Let $G$ be a polarized $p$-torsion group over $S$ and let, for each $u \in \F_p^{\times}$, $X_u := X_{G,u}(p)$ (see Section \ref{application-moduli}), which is a smooth proper $S$-scheme of relative dimension one with connected geometric fibres. 

Let $\OO_{p\ell}$ be the ring of integers of the maximal algebraic extension $\Q_{\{p,\ell\}}$ of $\Q$ which is unramified outside $p$ and $\ell$. Let $\overline{s}$ be a geometric point of $S$ above a topological point $s \in S$. Choose a morphism $\iota_{\overline{s}}: \OO_{p\ell} \rar \kappa(\overline{s})$: it defines a morphism $R_{\iota}: \pi_1(S,\overline{s}) \rar \mrm{Gal}(\Q_{\{p,\ell\}}/\Q)$ (whose conjugacy class does not depend on $\iota_{\overline{s}}$). 

Let $(P,Q)$ be a basis of $G(\kappa(\overline{s}))$ with $\zeta=\langle P,\,Q\rangle_G \in \mu_p^{\times}(\kappa(\overline{s}))$. For each $g \in \pi_1(S,\overline{s})$, let $D(g) \in \GL{\F_p}$ be the unique matrix such that $\begin{pmatrix} g(P)\\g(Q)\end{pmatrix} = D(g)\begin{pmatrix}P\\Q\end{pmatrix}$. Then \[\rho_{\overline{s}} := D^{-1}: \pi_1(S,\overline{s}) \rar \GL{\F_p}\] is a group homomorphism. 
 
The mod $p$ cyclotomic character $\omega_p$ is also a well-defined homomorphism $\pi_1(S,\overline{s}) \rar \F_p^{\times}$. Its restriction to $\mrm{Gal}(\overline{\kappa(s)}/\kappa(s))$ is the usual mod $p$ cyclotomic character; moreover, $\det{\rho_{\overline{s}}}=\omega_p^{-1}$. 

Let $\overline{\kappa(s)}$ denote the separable closure of $\kappa(s)$ in $\kappa(\overline{s})$, then $\rho_s := (\rho_{\overline{s}})_{|\mrm{Gal}(\overline{\kappa(s)}/\kappa(s))}$ identifies with the representation attached to the basis $(P,Q)$ of the polarized $p$-torsion group $G_s$ over $\kappa(s)$, and $(R_{\iota_{\overline{s}}})_{|\mrm{Gal}(\overline{\kappa(s)}/\kappa(s))}$ is exactly the group homomorphism $R_{\iota}$ defined in Lemma \ref{from-qbar-to-ks}. Furthermore, every $(X_u)_{s}$ is isomorphic to $X_{\rho_s,\zeta^u}(p)$, and the extension of this isomorphism to the Jacobians of these curves respects the Hecke correspondences. 

Let $S_0$ be a finite \'etale connected $S$-scheme endowed with a geometric point $\overline{s}_0$ lifting $\overline{s}$ such that $P,Q$ lift to $G(S_0)$. There is a natural action of $\SL{\F_p}$ on $X_{G,u}(p)\times_S S_0$ for any $u$ such that the isomorphism $(X_{G,u} \times_S S_0) \times_{S_0} \overline{s}_0 \simeq X_{G_{\overline{s}},u} \rar X_{\rho_s,\zeta^u}(p) \times_{\kappa(s)} \overline{s}$ of Proposition \ref{group-is-twist-polarized} is $\SL{\F_p}$-equivariant. 

Moreover, for any $u \in \F_p^{\times}$, the actions of $\mathbb{T}_1$, $\pi_1(S,\overline{s})$ and $\SL{\F_p}$ on $\Tate{\ell}{J_{\rho_s,\zeta^u}(p)}$ factor through the algebra $\mathbb{T}_1[\SL{\F_p} \rtimes_{\rho} \pi_1(S,\overline{s})]$. 

If we replace $(P,Q)$ with a different basis $\begin{pmatrix}P' \\Q'\end{pmatrix}=M\begin{pmatrix}P \\Q \end{pmatrix}$ (for $M \in \GL{\F_p}$), then the new representation is exactly $M\rho_{\overline{s}}M^{-1}$; for $N \in \SL{\F_p}$, the action of $N$ attached to the new basis $(P',Q')$ is $(MNM)_{(P,Q)}^{-1}$, where the subscript means that we consider the action of the matrix attached to the basis $(P,Q)$.   

If we replace $\iota_{\overline{s}}$ with another ring homomorphism $\iota=\iota_{\overline{s}} \circ \sigma$ for some $\sigma \in \mrm{Gal}(\Q_{\{p,\ell\}}/\Q)$, then $R_{\iota} = \sigma^{-1}R_{\iota_{\overline{s}}}\sigma$. 

If $\overline{s'}$ is another geometric point of $S$, there is a (non-canonical) isomorphism $I$ between the functors $\text{\bf{Fin\'Et}}_S \rar \mathbf{FinSet}$ given by $Y \longmapsto Y_{\overline{s}}$, $Y \longmapsto Y_{\overline{s'}}$. Let $(P',Q',\zeta',\overline{s'}_0)$ be the image 	of $(P,Q,\zeta,\overline{s_0})$ under $I$ (in particular, this uniquely defines an action of $\SL{\F_p}$ on $X_{G,u}(p) \times_S \kappa(\overline{s'})$). Moreover, $I(\iota_{\overline{s}})$ defines a ring homomorphism $\iota_{\overline{s'}}: \OO_{p\ell} \rar \kappa(\overline{s'})$ and also induces an isomorphism $\pi_1(I): \pi_1(S,\overline{s}) \rar \pi_1(S,\overline{s'})$, in such a way that $\rho_{\overline{s}}=\rho_{\overline{s}}\circ \pi_1(I)$. The following diagram commutes:
\[
\begin{tikzcd}[ampersand replacement=\&]
\pi_1(S,\overline{s}) \arrow{d}{\pi_1(I)} \arrow{r}{R_{\iota_{\overline{s}}}} \& \mrm{Gal}(\Q_{\{p,\ell\}}/\Q) \arrow{d}{\mrm{id}}\\
\pi_1(S,\overline{s'}) \arrow{r}{R_{\iota_{\overline{s'}}}} \& \mrm{Gal}(\Q_{\{p,\ell\}}/\Q)
\end{tikzcd}\]

For every $u \in \F_p^{\times}$, $I$ induces an isomorphism $\Tate{\ell}{\mrm{Jac}(X_u)_{\overline{s}}} \otimes_{\Z_{\ell}} F_{\mathfrak{l}} \rar \Tate{\ell}{\mrm{Jac}(X_u)_{\overline{s'}}} \otimes_{\Z_{\ell}} F_{\mathfrak{l}}$ of $(\mathbb{T}_1 \otimes F_{\mathfrak{l}})[\SL{\F_p} \rtimes_{\rho_{\overline{s}}} \pi_1(S,\overline{s})]$-module (whose action on $\mrm{Jac}(X_u)_{\overline{s'}}$ is given through $\pi_1(I)$). 
}

\demo{That $D$ is an anti-homomorphism (and thus that $\rho_{\overline{s}}$ is a group homomorphism) is a direct verification similar to Lemma \ref{right-representation}. 

Let $E \subset \mrm{Frac}(\OO_{p\ell})$ be a finite Galois extension of $\Q$ which is unramified outside $p\ell$. Then $\Sp{\OO_E[(p\ell)^{-1}]}$ is a finite \'etale $\Z[(p\ell)^{-1}]$-scheme, so $\pi_1(S,\overline{s})$ acts on the topological space $\Sp{\OO_{E}[(p\ell)^{-1}] \times_{\Z} S \times_S \overline{s}} \simeq (\Sp{\OO_{E}})_{\overline{s}}$. The choice of $\iota_{\overline{s}}$ amounts to specifying a point $i \in (\Sp{\OO_{E}})_{\overline{s}}$, with $i=\left[\iota_{\overline{s}} \otimes 1: \OO_E \otimes \kappa(\overline{s})\rar \kappa(\overline{s})\right]$. By definition, $(\Sp{\OO_{E}})_{\overline{s}}$ is endowed with actions of $\mrm{Gal}(E/\Q)$ (right) and $\pi_1(S,\overline{s})$ (left) that commute. Since $\mrm{Gal}(E/\Q)$ is the automorphism group of $(\Sp{\OO_{E}})_{\overline{s}}$, we can define a homomorphism $\pi_1(S,\overline{s}) \rar \mrm{Gal}(E/\Q)$ by mapping $g \in \pi_1(S,\overline{s})$ to the unique $\sigma \in \mrm{Gal}(E/\Q)$ such that $(\Sp{\sigma} \times \mrm{id})i = g(i)$. 

This morphism is clearly compatible with extensions $E' \supset E$, so this defines a homomorphism $R_{\iota_{\overline{s}}}: \pi_1(S,\overline{s}) \rar \mrm{Gal}(\Q_{p\ell}/\Q)$. Note that when $g=\Sp{\alpha^{-1}}$ with $\alpha \in \mrm{Gal}(\overline{\kappa(s)}/\kappa(s))$, then 
\begin{align*}
g(i) &= \ker[\iota_{\overline{s}} \circ (\mrm{id}\times\alpha^{-1})] = \ker[\alpha\circ\iota_{\overline{s}}\circ (\mrm{id}\times\alpha^{-1})]\\
&= \ker[\iota_{\overline{s}} \circ (R_{\iota_{\overline{s}}}(\alpha) \times\alpha)] \circ (\mrm{id}\times\alpha^{-1})]\\
&= \ker[\iota_{\overline{s}} \circ (R_{\iota_{\overline{s}}}(\alpha) \times \mrm{id})] = (\Sp{R_{\iota_{\overline{s}}}(\alpha)} \times \mrm{id})i,
\end{align*}
where the $R_{\iota_{\overline{s}}}(\alpha)$ is in the sense of Lemma \ref{from-qbar-to-ks}, so our notation makes sense.

For any $g \in \pi_1(S,\overline{s})$, $g$ acts on the finite discrete group $(\mu_p)_{k_s}$, so this automorphism is given by some $\underline{a}$ for $a \in \F_p^{\times}$. The mod $p$ cyclotomic character is simply $g \longmapsto a$. Since we came from the finite \'etale $\Z[p^{-1}]$-scheme $(\mu_p)_{\Z[1/p]}$, this character factors through $\pi_1(S,\overline{s}) \rar G_{\Q,p}$ (this map is only defined up to conjugacy, but since $\omega_p$ is a character, it does not depend from the representative of the conjugacy class). 

For any field $k$ with separable closure $k_s$ such that $p \in k^{\times}$, for any $\sigma \in \mrm{Gal}(k_s/k)$, for any $f \in \mu_p^{\times}(k_s)$, $f \circ \Sp{\sigma} = \sigma \circ f^{\sharp} = f^{\sharp}\circ \underline{\omega_p(\sigma)}$. This implies by Lemma \ref{pi1-vs-galois} that the mod $p$ cyclotomic character on $\pi_1(S,\overline{s})$ is the composition of $\pi_1(S,\overline{s}) \rar G_{\Q,p\ell} \overset{\omega_p}{\rar} \F_p^{\times}$, and that its restriction to any Galois group is the usual mod $p$ cyclotomic character. 
 
This discussion implies in particular that $\det{D}=\omega_p$, thus $\det{\rho_{\overline{s}}}=\omega_p^{-1}$. 

For every $n \geq 1$ and $u \in \F_p^{\times}$, $\mrm{Jac}(X_u)[\ell^n]$ is a finite \'etale $S$-group scheme with an action of $\mathbb{T}_1$. So $\pi_1(S,\overline{s})$ acts on the finite $\mathbb{T}_1$-module $\mrm{Jac}(X_u)[\ell^n]_{\overline{s}}$. We can then take the inverse limit over all $n \geq 1$ and tensor with $F_{\mathfrak{l}}$. 

Most of the claimed compatibilities are direct verifications, using Propositions \ref{group-is-twist-polarized} and \ref{jac+hecke-twist-polarized} as well as Lemmas \ref{T1Gamma-action}, \ref{pi1-vs-galois} and \ref{twist-mix-gl2-galois}. The existence of $I$ is \cite[Proposition 5.5.1]{Szamuely}. We elaborate on the construction of $\iota_{\overline{s'}}$ and the fact that $R_{\iota_{\overline{s}}}=R_{\iota_{\overline{s}}}\circ \pi_1(I)$. 

Let $E/\Q$ be a finite Galois extension unramified outside $p\ell$. Then \[i'=I(\ker\left[\iota_{\overline{s}}: \OO_E \otimes \kappa(\overline{s}) \rar \kappa(\overline{s})\right])\] is a point of $\Sp{\OO_E \otimes \kappa(\overline{s'})}$, so is the kernel of a unique morphism $\iota_{\overline{s'}}: \OO_E \otimes \kappa(\overline{s'}) \rar \kappa(\overline{s'})$ of $\kappa(\overline{s'})$-algebras. The following diagram commutes and commutes to the action of $\mrm{Gal}(E/\Q)$. 
\[
\begin{tikzcd}[ampersand replacement=\&]
\Sp{(\OO_E)_{\overline{s}}} \arrow{r}{g}\arrow{d}{I} \& \Sp{(\OO_E)_{\overline{s}}}\arrow{d}{I}\\
\Sp{(\OO_E)_{\overline{s'}}} \arrow{r}{g} \& \Sp{(\OO_E)_{\overline{s'}}}
\end{tikzcd}
\]
Hence $R_{\iota_{\overline{s}}}$ and $R_{\iota_{\overline{s'}}} \circ \pi_1(I)$ have the same projections in $\mrm{Gal}(E/\Q)$, whence $R_{\iota_{\overline{s}}} = R_{\iota_{\overline{s'}}} \circ \pi_1(I)$. 

 }

\prop[decomposition-xgp-sl20]{Let $\mathfrak{l}$ be a maximal ideal of $\OO_F$ with residue characteristic $\ell$, and $S$ be a connected regular excellent Noetherian affine $\Z[(p\ell)^{-1}]$-scheme of characteristic zero. Let $G$ be a polarized $p$-torsion group over $S$ and let, for each $u \in \F_p^{\times}$, $X_u := X_{G,u}(p)$ (see Section \ref{application-moduli}), which is a smooth proper $S$-scheme of relative dimension one with connected geometric fibres. 

Let $\OO_{p\ell}$ be the ring of integers of the maximal algebraic extension $\Q_{\{p,\ell\}}$ of $\Q$ unramified outside $p$ and $\ell$. Let $\overline{s}$ be a geometric point of $S$ above a topological point $s \in S$; fix a morphism $\iota_{\overline{s}}: \OO_{p\ell} \rar \kappa(\overline{s})$, which induces a morphism $R_{\iota}: \pi_1(S,\overline{s}) \rar \mrm{Gal}(\Q_{\{p,\ell\}}/\Q)$. 

Let $(P,Q)$ be a basis of $G(\kappa(\overline{s}))$ with $\zeta=\langle P,\,Q\rangle_G \in \mu_p^{\times}(\kappa(\overline{s}))$, and, for each $g \in \pi_1(S,\overline{s})$, $\rho_{\overline{s}}(g) \in \GL{\F_p}$ be the unique matrix such that $\begin{pmatrix} g(P)\\g(Q)\end{pmatrix} = \rho_{\overline{s}}^{-1}(g)\begin{pmatrix}P\\Q\end{pmatrix}$. 

Thus $\rho_{\overline{s}}: \pi_1(S,\overline{s}) \rar \GL{\F_p}$ is a group homomorphism with determinant $\omega_p^{-1}$. 

Let $S_0$ be a finite \'etale connected $S$-scheme endowed with a geometric point $\overline{s}_0$ lifting $\overline{s}$ such that $P,Q$ lift to $G(S_0)$. When $p \equiv 3 \pmod{4}$, let $S_2$ be the finite \'etale $S$-scheme corresponding to the $\pi_1(S,\overline{s})$-set $\{\pm\}$, where the action is by $\left(\frac{\omega_p}{p}\right)$. Then, for a suitable geometric point $\overline{s}_2$ of $S_2$, the map $(S_0,\overline{s}_0) \rar (S,\overline{s})$ factors through $(S_2,\overline{s_2})$. Note that $S_2=S \times_{\Z} \Sp{\Z\left[\frac{1+\sqrt{-p}}{2}\right]}$.    

Let $G_{\overline{s},\rho} = \SL{\F_p} \rtimes_{\rho_{\overline{s}}} \pi_1(S,\overline{s})$ and $\tilde{\rho_{\overline{s}}}: (M,g) \in G_{\overline{s},\rho} \longmapsto M\rho_{\overline{s}}(g) \in \GL{\F_p}$, and $G_2 = \tilde{\rho_{\overline{s}}}^{-1}(\F_p^{\times}\SL{\F_p}) \simeq \SL{\F_p} \rtimes_{\rho_{\overline{s}}} \pi_1(S_2,\overline{s}_2)$.  

Then, for each $u \in \F_p^{\times}$, the eigenspaces of the $(\mathbb{T}_1 \otimes F_{\mathfrak{l}})[G_{\overline{s},\rho}]$-module $\Tate{\ell}{\mrm{Jac}(X_u)_{\overline{s}}} \otimes_{\Z_{\ell}} F_{\mathfrak{l}}$ are:

\begin{itemize}[noitemsep,label=$-$]
\item for $(f,\mathbf{1}) \in \mathscr{S}$, a copy of $R_{f,\mathfrak{l}} := (V_{f,\mathfrak{l}}\circ R_{\iota_{\overline{s}}}) \otimes \mrm{St}(\tilde{\rho_{\overline{s}}})$,
\item for $(f,\chi) \in \mathscr{P}'$, a copy of $R_{f,\mathfrak{l}} := (V_{f,\mathfrak{l}}\circ R_{\iota_{\overline{s}}}) \otimes (\pi(1,\chi))(\tilde{\rho_{\overline{s}}})$,
\item for $(f,\mathbf{1}) \in \mathscr{C}'_1$, a copy of $R_{f,\mathfrak{l}} := (V_{f,\mathfrak{l}}\circ R_{\iota_{\overline{s}}}) \otimes C_f(\tilde{\rho_{\overline{s}}})$,
\item for $(f,\mathbf{1}) \in \mathscr{C}_M$, there is a sign $r$ depending only on $\iota^{-1}(\zeta)$ such that it is: 
\begin{itemize}[noitemsep,label=\tiny$\bullet$]
\item if $S_2 \rar S$ is an isomorphism, a copy of $\bigoplus_{v \in \{\pm\}}{(\psi_{f,rv,\mathfrak{l}}\circ R_{\iota_{\overline{s}}}) \otimes C^v(\tilde{\rho_{\overline{s}}})}$. 
\item if $S_2 \rar S$ is not an isomorphism, a copy of \[R_{f,r,\mathfrak{l}} := \mrm{Ind}_{\pi_1(S_2,\overline{s}_2)}^{\pi_1(S,\overline{s})}\left[(\psi_{f,r,\mathfrak{l}}\circ R_{\iota_{\overline{s}}}) \otimes C^+((\tilde{\rho_{\overline{s}}})_{|\pi_1(S_2,\overline{s}_2)})\right].\] 
\end{itemize}
\end{itemize}
}

\demo{Assume first that $s \in S$ is the generic point. The map $\mrm{Gal}(\overline{\kappa(s)}/\kappa(s)) \rar \pi_1(S,\overline{s})$ is surjective \cite[Remark 5.7.12]{Szamuely} and the kernel of this morphism acts trivially on every $\Tate{\ell}{\mrm{Jac}(X_u)}$ by definition: therefore, by Corollary \ref{decomposition-xrho-connected-sl20}, the result holds for any $(u,\mathfrak{l})$. 

The claim then holds for every possible choice of $(\overline{s},P,Q,\iota_{\overline{s}},S_0,S_2)$ and any $(u,\mathfrak{l})$ thanks to the discussion of Lemma \ref{xgp-family-pi1-elem}. }

\medskip

\cor[decomposition-xrho-connected-sl2]{The statement of Corollary \ref{decomposition-xrho-connected-sl20} holds if we only assume that the characteristic of $k$ is not $p$ or $\ell$ after the following modifications are made:
\begin{itemize}[noitemsep,label=\tiny$\bullet$]
\item if the characteristic of $k$ is $q \neq 0$, then let $\OO_{\{p,\ell\}}$ be the ring of integers of the maximal algebraic extension $\Q_{\{p,\ell\}}$ of $\Q$ which is unramified outside $p,\ell$, and $\iota$ is a ring homomorphism $\OO_{\{p,\ell\}} \rar k_s$ inducing a morphism $R_{\iota}: \mrm{Gal}(k_s/k) \rar \mrm{Gal}(\Q_{\{p,\ell\}}/\Q)$.  
\item if $p \equiv 3 \pmod{4}$ and $k$ has characteristic $2$, then $k'$ is the subextension of $k_s$ generated by the roots of $X^2-X-\frac{p+1}{4}$, so that $\mrm{Gal}(k_s/k')=\omega_p^{-1}(\F_p^{\times 2})$.
\end{itemize}}

\demo{When the characteristic of $k$ is zero, this was already proved, so we assume that $k$ has positive characteristic $q$. By the same argument as the proof of Corollary \ref{decomposition-xrho-connected-sl20}, it is enough to prove the result for \emph{one} representation of $\mrm{Gal}(\overline{\F}_q/\F_q)$. 

When $q \neq 11$, we apply Proposition \ref{decomposition-xgp-sl2} to $S=\Sp{\Z[(11p\ell)^{-1}]}$ and $G$ being the $p$-torsion group of the elliptic curve $E/S$ with Weierstrass equation given by $y^2-y=x^3-x^2$ (LMFDB \cite{lmfdb} label $11.a3$) at a geometric point in the special fibre at $q$, and then we restrict to its absolute Galois group. 

When $q=11$, we apply Proposition \ref{decomposition-xgp-sl2} to $S=\Sp{\Z[(14p\ell)^{-1}]}$ and $G$ being the $p$-torsion group of the elliptic curve $E/S$ with Weierstrass equation given by $y^2+xy+y=x^3-x$ (LMFDB label $14.a5$) at a geometric point in the special fibre at $q$, and then we restrict to the absolute Galois group. 
}

\cor[decomposition-xrho-connected-general]{Let $\mathfrak{l}$ be a maximal ideal of $\OO_F$ with residue characteristic $\ell$ and $k$ be a field with separable closure $k_s$ and residue characteristic distinct from $p,\ell$. Let $\OO$ be the ring of integers of the maximal algebraic extension $\Q_{\{p,\ell\}}$ of $\Q$ unramified outside $p,\ell$ and $\iota: \OO \rar k_s$ be a ring homomorphism, inducing a continuous group homomorphism $R_{\iota}: \mrm{Gal}(k_s/k) \rar \mrm{Gal}(\Q_{\{p,\ell\}}/\Q)$. 

Let $\rho: \mrm{Gal}(k_s/k) \rar \GL{\F_p}$ be a continuous group homomorphism with $\det{\rho}=\omega_p^{-1}$. Then, for any $\xi \in \mu_p^{\times}(k_s)$, as a $(\mathbb{T}_1 \otimes F_{\mathfrak{l}})[\mrm{Gal}(k_s/k)]$-module, $\Tate{\ell}{J_{\rho,\xi}(p)}\otimes_{\Z_{\ell}}F_{\mathfrak{l}}$ is isomorphic to the following direct sum:
\begin{itemize}[noitemsep,label=$-$]
\item for every $(f,\chi) \in \mathscr{S} \cup \mathscr{P}' \cup \mathscr{C}'_1$, a copy of $R_{f,\mathfrak{l}}(\rho)$, on which, for all $n,m \geq 1$ such that $nm^2\equiv 1\pmod{p}$, $T_n \cdot (mI_2)$ acts by $\chi(m)a_n(f)$, 
\item for $(f,\chi) \in \mathscr{C}_M$, let $k'$ be the subextension of $k_s/k$ generated by the roots of the polynomial $X^2-X-\frac{p+1}{4}$. There exists a sign $r$ depending only on $\iota^{-1}(\xi)/\F_p^{\times 2}$ such that the factor (on which, for any $n,m \geq 1$ with $nm^2\equiv 1\pmod{p}$) is:
\begin{itemize}[label=\tiny$\bullet$,noitemsep]
\item $\mrm{Ind}_{\mrm{Gal}(k_s/k')}^{\mrm{Gal}(k_s/k)}{(\psi_{f,r,\mathfrak{l}}\circ R_{\iota}) \otimes C^+(\rho_{|\mrm{Gal}(k_s/k')})}$, when $k'\neq k$,
\item $\bigoplus_{v \in \{\pm\}}{(\psi_{f,rv,\mathfrak{l}}\circ R_{\iota}) \otimes C^v(\rho_{|\mrm{Gal}(k_s/k')})}$ when $k'=k$. 
\end{itemize} 
\end{itemize}
}

\medskip

\prop[decomposition-xgp-sl2]{The statement of Proposition \ref{decomposition-xgp-sl20} holds \emph{as is} when $S$ is not assumed to have characteristic zero.}  

\demo{This is exactly similar to the proof of Proposition \ref{decomposition-xgp-sl20} after replacing Corollary \ref{decomposition-xrho-connected-sl20} with Corollary \ref{decomposition-xrho-connected-sl2}.}

\bigskip

We conclude this section by discussing the action of $\mathbb{T}_{1,\{1\}}$ on the Tate modules of the connected components. 

\prop[connected-action-T11]{Let $k$ be a field with characteristic $q \neq p$, $\mathfrak{l}$ be a maximal ideal of $\OO_F$ with residue characteristic $\ell \neq q$. Let $k_s$ be a separable closure of $k$ and $\rho: \mrm{Gal}(k_s/k) \rar \GL{\F_p}$ be group homomorphism such that $\det{\rho}=\omega_p^{-1}$. Let $\xi \in \mu_p^{\times}(k_s)$. 
Let $\OO$ be the ring of integers of the maximal algebraic extension $\Q_{\{p,\ell\}}/\Q$ unramified outside $p\ell$ and fix an embedding $\iota: \OO \rar k_s$. This determines a group homomorphism $R_{\iota}: \mrm{Gal}(k_s/k) \rar \mrm{Gal}(\Q_{\{p,\ell\}}/\Q)$. 

Let $n \geq 1, M \in \GL{\F_p},\sigma \in \mrm{Gal}(k_s/k)$ be such that $n\det{M} \equiv 1\pmod{p}$. Then:
\begin{itemize}[noitemsep,label=$-$]
\item For any $(f,\chi) \in \mathscr{S}$ (resp. $(f,\chi) \in \mathscr{P}'$, resp. $(f,\chi) \in \mathscr{C}'_1$), for any $r \geq 1$, the action of $(T_nM\sigma)^r$ on the $f$-eigenspace of $\Tate{\ell}{J_{\rho,\xi}(p)} \otimes_{\Z_{\ell}} F_{\mathfrak{l}}$ is $a_n(f)^r\mu_r\mrm{Tr}(R_{\iota}(\sigma^r) \mid V_{f,\mathfrak{l}})$, where $\mu_r$ is the trace of the action of $(M\rho(\sigma))^r$ on $\mrm{St}$ (resp. on $\pi(1,\chi)$, resp. on $C_f$), 
\item If $(f,\chi) \in \mathscr{C}_M$, $T_nM=0$ if $n \notin \F_p^{\times 2}$. Otherwise, for any $r \geq 1$, the action of $(T_nM\sigma)^r$ on the $f$-eigenspace of $\Tate{\ell}{J_{\rho,\xi}(p)} \otimes_{\Z_{\ell}} F_{\mathfrak{l}}$ is 
\begin{itemize}[noitemsep,label=\tiny$\bullet$]
\item if $\omega_p(\sigma)^r \notin \F_p^{\times 2}$, it is $0$, 
\item if $r$ is even and $\omega_p(\sigma) \notin \F_p^{\times 2}$, let $V_4$ be the representation of $\GL{\F_p}$ described in Section \ref{splitting-sl2}, the result is \[\frac{1}{2}a_n(f)^r\mrm{Tr}(\sigma^r \mid V_{f,\mathfrak{l}})\mrm{Tr}((M\rho(\sigma))^r \mid V_4),\]  
\item if $\omega_p(\sigma) \in \F_p^{\times 2}$ and $s$ is the sign attached to $\iota^{-1}(\xi)$ defined in Corollary \ref{decomposition-xrho-connected-sl2}; \[a_n(f)^r\left[\psi_{f,s,\mathfrak{l}}(R_{\iota}(\sigma^r))\mrm{Tr}((M\rho(\sigma))^r \mid C^+) + \psi_{f,-s,\mathfrak{l}}(R_{\iota}(\sigma^r))\mrm{Tr}((M\rho(\sigma))^r \mid C^-)\right].\]  
\end{itemize}
\end{itemize}
}

\demo{When $M \in \SL{\F_p}$, the claim is contained in Proposition \ref{decomposition-xrho-connected-sl2}. When $\det{M} \in \F_p^{\times 2}$, there exists $M_1 \in \SL{\F_p}$ and $m \in \F_p^{\times}$ such that $M=mM_1$: hence $nm^2 \equiv 1\pmod{p}$. Therefore $(T_nM\sigma)^r=(T_n(mI_2)M_1\sigma)^r=(T_n \cdot (mI_2))^r(M_1\sigma)^r$ acts on the $f$-eigenspace by $(a_n(f)\chi(m))^r(M_1\sigma)^r$, and we are thus reduced to the case of $\SL{\F_p}$ (because $\chi$ is the central character of the representation of $\GL{\F_p}$).   

Now, assume that $M \notin \F_p^{\times}\SL{\F_p}$. When $(f,\chi) \in \mathscr{C}_M$, $\Tate{\ell}{J_{\rho,\xi}(p)} \otimes_{\Z_{\ell}} F_{\mathfrak{l}}$ is a submodule of the $f$-eigenspace of $\Tate{\ell}{J_{\rho}(p)} \otimes_{\Z_{\ell}}F_{\mathfrak{l}}$, on which $T_n$ acts by $a_n(f)=0$ by Corollary \ref{decomposition-xpp-tate}. So we assume that $(f,\chi) \notin \mathscr{C}_M$. 

Let $\tau_{(M,\sigma),\xi,r}$ be the trace of the action of $(T_nM\sigma)^r$ on the $f$-eigenspace of $\Tate{\ell}{J_{\rho,\xi}(p)} \otimes_{\Z_{\ell}} F_{\mathfrak{l}}$. Clearly, $\tau_{(M,\sigma),\xi,r}$ only depends on the $\SL{\F_p} \rtimes_{\rho} \mrm{Gal}(k_s/k)$-conjugacy class of $(M,\sigma)$. Since $\det{M} \in \F_p^{\times 2}$, $\F_p[M]$ contains matrices of every possible determinant. When $N \in \GL{\F_p}$ commutes to $M$, $(N\rho(\sigma)^{-1},\sigma) \in \GL{\F_p} \rtimes_{\rho} \mrm{Gal}(k_s/k)$ commutes to $(M,\sigma)$. Therefore, the $\SL{\F_p} \rtimes_{\rho} \mrm{Gal}(k_s/k)$-conjugacy class of $(M,\sigma)$ is, in fact, its $\GL{\F_p} \rtimes_{\rho} \mrm{Gal}(k_s/k)$-conjugacy class. 

For any $\xi_1 \in \mu_p^{\times}(k_s)$, any $(N,\sigma') \in \GL{\F_p} \rtimes_{\rho} \mrm{Gal}(k_s,k)$ induces isomorphisms from the $f$-eigenspace of $\Tate{\ell}{J_{\rho,\xi_1}(p)} \otimes_{\Z_{\ell}} F_{\mathfrak{l}}$ to that of $\Tate{\ell}{J_{\rho,\xi_1^{\det{N}}}(p)} \otimes_{\Z_{\ell}} F_{\mathfrak{l}}$ which commute to $T_n$: therefore, $\tau_{(N,\sigma')(M,\sigma')^r(N,\sigma')^{-1},\xi^{\det{N}},r} = \tau_{(M,\sigma),\xi,r}$. Hence $\tau_{(M,\sigma),\xi,r}$ does not depend on $\xi$. 

Therefore, $(p-1)\tau_{M,\xi,r}$ is exactly the trace of $(T_nM\sigma)^r$ acting on the $f$-eigenspace of $\Tate{\ell}{J(p,p)_{\rho}}$, which is, by Proposition \ref{jacobian-twist}, the trace of $(T_nM\rho(\sigma)\sigma)^r$ acting on $\Tate{\ell}{J(p,p)_k}$. By Corollary \ref{decomposition-xpp-tate} (and Lemma \ref{agreement-implies-twist}, Corollary \ref{abelian-scheme-change-field}), this trace is 
\begin{align*}
&\sum_{\psi \in \mathcal{D}}{\mrm{Tr}(R_{\iota}(\sigma^r) \mid V_{f\otimes \psi,\mathfrak{l}})a_n(f \otimes \psi)^r\psi(\det{(M\rho(\sigma))^r})\mrm{Tr}(M^r \mid R_f)} \\
&= 2(p-1)a_n(f)^r\mrm{Tr}((M\rho(\sigma))^r \mid R_f),
\end{align*}
 where $R_f$ is $\mrm{St}$ if $(f,\chi) \in \mathscr{S}$, $\pi(1,\chi)$ if $(f,\chi) \in \mathscr{P}'$, and $C_f$ if $(f,\chi) \in \mathscr{C}'_1$. }

\medskip

\cor[connected-action-T11rho]{Let $\rho,k,k_s,q,\ell,\mathfrak{l},\iota,R_{\iota}$ be as Proposition \ref{connected-action-T11}. Let $(f,\chi) \in \mathscr{S}$ (resp. $(f,\chi) \in \mathscr{P}'$, resp. $(f,\chi) \in \mathscr{C}'_1$), and let $Q_f=\mrm{St}$ (resp. $Q_f=\pi(\mathbf{1},\chi)$, resp. $Q_f=C_f$). Assume that $V_{f,\mathfrak{l}}\circ R_{\iota}$ is semisimple, for instance if the image of $R_{\iota}$ is open. 

Let $G_{k,\rho} = \GL{\F_p} \rtimes_{\rho} \mrm{Gal}(k_s/k)$, and $\tilde{\rho}: (M,g) \in G_{k,\rho}\longmapsto M\rho(g)$. Let $\tilde{R}_{\iota}$ be the homomorphism $(M,g) \in G_{k,\rho} \longmapsto R_{\iota}(g)$. 

Then, for any $\xi \in \mu_p^{\times}(k_s)$, the $f$-eigenspace of $\Tate{\ell}{J_{\rho,\xi}(p)} \otimes_{\Z_{\ell}} F_{\mathfrak{l}}$ is a semi-simple $\mathbb{T}_{1,\GL{\F_p},\rho} \otimes F_{\mathfrak{l}}$-module, and it is isomorphic to the restriction of the $\mathbb{T}[G_{k,\rho}] \otimes F_{\mathfrak{l}}$-module $(V_{f,\mathfrak{l}} \circ \tilde{R}_{\iota}) \otimes Q_f(\tilde{\rho})$, where $T_n \in \mathbb{T}$ (resp. $nI_2$) for $n \geq 1$ coprime to $p$ acts by $a_n(f)$ (resp. $\chi(n)$). }

\demo{First, by Proposition \ref{decomposition-xpp-tate}, there is a ring homomorphism $\mathbb{T} \rar F$ mapping $T_n$ (resp. $nI_2$) to $a_n(f)$ (resp. $\chi(n)$). Thus the $(\mathbb{T} \otimes F_{\mathfrak{l}})[G_{k,\rho}]$-module $(V_{f,\mathfrak{l}} \circ \tilde{R}_{\iota}) \otimes Q_f(\tilde{\rho})$ described in the statement is well-defined. Let us prove that it is semi-simple: it is enough to show that it is semi-simple as a continuous representation of $G_{k,\rho}$ with coefficients in $F_{\mathfrak{l}}$. 

By \cite[\S 2.7 Lemma]{BH}, it is enough to show that $(V_{f,\mathfrak{l}} \circ \tilde{R}_{\iota}) \otimes Q_f(\tilde{\rho})$ is semi-simple as a $F_{\mathfrak{l}}[\mrm{Gal}(k_s/k)]$-module\footnote{As previously, the reference discusses only complex smooth representations, but their argument applies \emph{as is}.}. This is the case, because this module is the direct sum of a certain number of copies of $V_{f,\mathfrak{l}}\circ R_{\iota}$, which is semi-simple by assumption (if $R_{\iota}$ has an image of finite index, it is a consequence of \cite[Theorem 4.4]{Antwerp5-Ribet}). 

By Proposition \ref{connected-action-T11}, the $f$-eigenspace $E_{f,\mathfrak{l},\xi}$ of $\Tate{\ell}{J_{\rho,\xi}(p)} \otimes_{\Z_{\ell}} F_{\mathfrak{l}}$ is a $\mathbb{T}_{1,\GL{\F_p},\rho} \otimes F_{\mathfrak{l}}$-module with the same traces as $(V_{f,\mathfrak{l}}\circ \tilde{R}_{\iota}) \otimes Q_f(\tilde{\rho})$. So the semi-simplification of $E_{f,\mathfrak{l}}$ is exactly $(V_{f,\mathfrak{l}}\circ \tilde{R}_{\iota}) \otimes Q_f(\tilde{\rho})$. 

Now, $\bigoplus_{\xi}{E_{f,\mathfrak{l},\xi}}$ is the $f$-eigenspace of $\Tate{\ell}{J_{\rho,\xi}(p)} \otimes_{\Z_{\ell}} F_{\mathfrak{l}}$, which is, by Corollary \ref{decomposition-xpp-tate} and Proposition \ref{jacobian-twist}, isomorphic over $(\mathbb{T} \otimes F_{\mathfrak{l}})[G_{k\rho}]$ to $\bigoplus_{\psi \in\mathcal{D}}{(\psi(\omega_p\det{\tilde{\rho}}) \otimes (V_{f,\mathfrak{l}}\circ \tilde{R}_{\iota}) \otimes C_f(\tilde{\rho})}$, where $T_n$ (resp. $nI_2$) acts on the summand $\psi$ by $a_n(f)\psi(n)$ (resp. $\chi(n)\psi(n)^2$) for any $n \geq 1$ coprime to $p$. 

Therefore, one has the following isomorphism of $\mathbb{T}_{1,\GL{\F_p},\rho} \otimes F_{\mathfrak{l}}$-modules:
\[\bigoplus_{\xi}{E_{f,\mathfrak{l},\xi}} \simeq \left((V_{f,\mathfrak{l}}\circ \tilde{R}_{\iota}) \otimes Q_f(\tilde{\rho})\right)^{\oplus (p-1)},\] whence every $E_{f,\mathfrak{l},\xi}$ is semi-simple and the conclusion follows. }

\newpage

%% file: autom-0.tex
\chapter{Automorphy, functional equations and root numbers for $L(J_{\rho},s)$}
\label{automorphy}

\section{Introduction and statement of results}

Let $\rho: G_{\Q} \rar \GL{\F_p}$ be a continuous representation such that $\det{\rho}=\omega_p^{-1}$ is the inverse of the cyclotomic character modulo $p$. Then $X_{\rho}(p) := X(p,p)_{\rho}$ is a smooth projective $\Q$-scheme of relative dimension one, and is, by Proposition \ref{group-is-twist-polarized}, the reunion of smooth projective geometrically connected $\Q$-schemes $X_{\rho,\xi}(p)$ indexed by $\xi \in \mu_p^{\times}(\Qbar)$. The Jacobian of $X_{\rho,\xi}(p)$ is noted $J_{\rho,\xi}(p)$. 

We will use the following notation from Sections \ref{analytic} and \ref{tate-modules-twists}:

\nott{\begin{itemize}[label=$-$,noitemsep]
\item $\mathcal{D}$ denotes the set of Dirichlet characters $(\Z/p\Z)^{\times} \rar \C^{\times}$, 
\item $\mathscr{S}$ denotes the set of $(f,\mathbf{1})$ where $f \in \mathcal{S}_2(\Gamma_0(p))$ is a newform, 
\item $\mathscr{P}$ denotes the set of $(f,\chi)$, where $f \in \mathcal{S}_2(\Gamma_1(p))$ is a newform with character $\chi \neq \mathbf{1}$,
\item $\mathscr{P}'$ denotes a set of representatives in $\mathscr{P}$ for $\mathscr{P}$ modulo the complex conjugation, 
\item $\mathscr{C}$ denotes the set of $(f,\chi)$, where $f \in \mathcal{S}_2(\Gamma_1(p) \cap \Gamma_0(p^2))$ is a newform with character $\chi \in \mathcal{D}$, and such that no twist of $f$ by a character in $\mathcal{D}$ has conductor $p$,
\item $\mathscr{C}'_1$ denotes a set of representatives of the form $(f,\mathbf{1})$ for the $(f,\chi) \in \mathscr{C}$ modulo twists such that $f$ does not have complex multiplication.
\item $\mathscr{C}_M$ denotes the collection of $(f,\mathbf{1}) \in\mathscr{C}$ such that $f$ has complex multiplication. In this case, $p \equiv 3\pmod{4}$ and $f$ has complex multiplication by $\Q(\sqrt{-p})$. 
\item $\mathbb{T}_1$ denotes the subalgebra of $\bigoplus_{\xi \in \mu_p^{\times}(\Qbar)}{\mrm{End}(J_{\rho,\xi}(p))}$ generated by the $(mI_2) \cdot T_n$ such that $n,m \geq 1$ are integers such that $nm^2\equiv 1\pmod{p}$. 
\item For $\alpha \in \mathcal{D}$, $\mrm{St}_{\alpha}$ is the Steinberg representation of $\GL{\F_p}$ twisted by $\alpha(\det)$. 
\item For $\alpha,\beta \in \mathcal{D}$, $\pi(\alpha,\beta)$ is the principal series representation of $\GL{\F_p}$ attached to the characters $\alpha,\beta$. 
\item Given $(f,\chi) \in \mathscr{C}$, we can attach by Definition \ref{left-cusp-rep-f} a cuspidal representation of $\GL{\F_p}$ named $C_f$. 
\item $\omega_p$ denotes the cyclotomic character modulo $p$. 
\end{itemize}}

The varieties $J_{\rho,\xi}(p)$ are Abelian varieties over $\Q$: according to a conjecture akin to the Birch and Swinnerton-Dyer conjecture, the $L$-functions of the eigenspaces (for the Hecke algebra $\mathbb{T}_1$) of its Tate module are entire, and their vanishing orders at $s=1$ match the ranks of the eigenspaces of $J_{\rho,\xi}(p)(\Q) \otimes \Q$ for $\mathbb{T}_1$.

By Theorem \ref{tate-4}, after extending scalars to a sufficiently large extension $F_{\mathfrak{l}}$ of $\Q_{\ell}$ (where $\ell$ is a prime that may be equal to $p$), the Tate module of $J_{\rho,\xi}(p)$ splits as a direct sum of $(\mathbb{T}_1 \otimes F_{\mathfrak{l}})[G_{\Q}]$-submodules (which we call \emph{eigenspaces} for $\mathbb{T}_1$, because $\mathbb{T}_1$ acts on them through a ring homomorphism $\mathbb{T}_1 \rar F_{\mathfrak{l}}$) of the following form:

\begin{itemize}[noitemsep,label=\tiny$\bullet$]
\item $V_{f,\mathfrak{l}} \otimes \mrm{St}(\rho)$, where $(f,\mathbf{1}) \in \mathscr{S}$, 
\item $V_{f,\mathfrak{l}} \otimes \pi(1,\chi)(\rho)$, where $(f,\chi) \in \mathscr{P}'$,
\item $V_{f,\mathfrak{l}} \otimes C_f(\rho)$, where $(f,\chi) \in \mathscr{C}'_1$,
\item $\mrm{Ind}_{G_{\Q(\sqrt{-p})}}^{G_{\Q}}{\psi_{f,\mathfrak{l}} \otimes C^+(\rho_{|G_{\Q(\sqrt{-p})}})}$, for $(f,\mathbf{1}) \in \mathscr{C}_M$, where $\psi_{f,\mathfrak{l}}: G_{\Q(\sqrt{-p})} \rar F_{\mathfrak{l}}^{\times}$ is a character attached to $f$ and $C^+$ is a cuspidal representation of $\F_p^{\times}\SL{\F_p}$ described in Section \ref{splitting-sl2}. 
\end{itemize}

The summand attached to some $(f,\chi) \in \mathscr{S}\cup\mathscr{P}'\cup\mathscr{C}'_1 \cup\mathscr{C}_M$ is called the \emph{$f$-eigenspace} of the Tate module of $J_{\rho,\xi}(p)$. The $f$-eigenspaces, for $(f,\mathbf{1}) \in \mathscr{C}_M$, are referred to as the \emph{CM eigenspaces} of the Tate module of the $J_{\rho,\xi}(p)$. The other eigenspaces are the \emph{non-CM eigenspaces}.  

It is naturally expected that these $L$-functions satisfy functional equations of the form $s \longleftrightarrow 2-s$. If we can prove it and determine the constant (in fact, the \emph{sign}, by Lemma \ref{everything-self-dual}), we may deduce the partiy of the vanishing order at $s=1$ of said factor. It seems reasonable to assume that these vanishing orders are ``often'' $0$ or $1$, and finding out their parity amounts to determining whether the corresponding eigenspace of $J_{\rho,\xi}(p)(\Q) \otimes_{\Z} F$ vanishes or not. 

In general, proving that the $L$-function of a $\ell$-adic Galois representation is entire and satisfies a functional equation requires automorphy lifting theorems. In our setting, such an automorphy result seems to be beyond reach for several reasons. First, the presence of large Artin representations means that our automorphic representations will not be regular (see for example \cite[Definition 1.6]{Clozel-regulier}), in which case the usual patching methods will not work. To the author's knowledge, even the mechanism by which a large-dimensional Artin representation of $G_{\Q}$ can be attached to an automorphic form is completely unknown in general. This is because there seems to be no clear Shimura variety attached to this situation, and, again to the author's knowledge, this means that the recent progress in higher Hida theory \cite{HigherHida} cannot be applied to such a situation. 

In our situation, in the ``generic'' case (say, $\rho$ is surjective and comes from the $p$-torsion of an elliptic curve $E[p]$), the common technique of base change to work over a better-behaved number field does not seem to work. For $\operatorname{GL}(2)$ (that is, two-dimensional Galois representations) base change from $\Q$ to any totally real number field is known by \cite{BaseChangeGL2, LanglandsSym5}, and base change for a \emph{solvable} extension of totally real fields is well-known \cite[Proposition 4.25]{GeeMLT}, but neither seems to simplify the situation. Indeed, since $\rho(G_{\Q})$ contains the perfect group $\SL{\F_p}$, for any number field $K$ such that $\rho(G_K)$ contains $\SL{\F_p}$ and any solvable extension $L/K$, $\rho(G_L)$ contains the $r$-th derived subgroup of $\rho(G_K)$, which contains $\SL{\F_p}$. Furthermore, since $\rho$ maps the complex conjugation to a conjugate of $\Delta_{1,-1} \in \GL{\F_p}$, if $F$ is any totally real number field (which we may assume is Galois), $F$ is fixed under every complex conjugation, so $\rho(G_F)$ is a normal subgroup of $\GL{\F_p}$ containing $\Delta_{1,-1}$, so contains $\det^{-1}(\{\pm 1\})$, and in particular contains $\SL{\F_p}$. 

A related idea is to use Brauer's theorem on induced characters in order to do away with the Artin representations, but, for instance in the case of $V_{f,\mathfrak{p}} \otimes \mrm{St}(\rho)$, it is the virtual difference of the representations $\mrm{Ind}_K^{\Q}{V_{f,\mathfrak{p}}} - V_{f,\mathfrak{p}}$ with $G_K = \rho^{-1}(B)$. The previous paragraph shows that $K$ is not a solvable extension of a totally real field, so that the base changed representation $(V_{f,\mathfrak{p}})_{|K}$ is not known (to the best of the author's knowledge) to be automorphic or possess an entire $L$-function.  

Thus, we can only properly study the $L$-function of $J_{\rho}(p)$ in specific situations, where for instance $\rho$ is not surjective. This is carried out in the next Chapter when the image of $\rho$ is contained in the normalizer of a nonsplit Cartan subgroup.  

If we instead \emph{assume} that the involved Artin representations are automorphic (in a sufficiently strong sense, see Definition \ref{artin-is-automorphic}), we can show that the $L$-functions of the factors of $J_{\rho}(p)$ are entire and satisfy functional equations whose sign can be computed on the Galois side by the theory of Rankin-Selberg products \cite[Theorem 9.1]{Cogdell}. This is the goal of the first part of this section, after reminders about the Langlands correspondence and the dictionary between modular forms and automorphic representations. 

\theo[autom-1]{(Lemma \ref{everything-self-dual}, Propositions \ref{automorphic-times-modform} and \ref{automorphic-times-modform-cmtwist}) 
\begin{itemize}[noitemsep,label=$-$]
\item For every eigenspace $V$ in the Tate module of each $J_{\rho,\xi}(p)$, $V$ is isomorphic to $V^{\ast}(1)$.
\item The $L$-function of any CM eigenspace of the Tate module of $J_{\rho,\xi}(p)$ is well-defined, extends meromorphically to the complex plane and satisfies the expected functional equation, whose sign is given by the product of the local Deligne-Langlands constants.  
\item Let $(f,\chi) \in \mathscr{S}$ (resp. $(f,\chi) \in \mathscr{P}'$, $(f,\chi) \in \mathscr{C}'_1$, $(f,\chi) \in \mathscr{C}_M$). Assume that the Artin representation $\mrm{St}\circ \rho$ (resp. $\pi(1,\chi) \circ \rho$, resp. $C_f \circ \rho$, resp. $C^+\circ \rho_{|\Q(\sqrt{-p})}$) is automorphic in the sense of Definition \ref{artin-is-automorphic} (with the obvious adjustments over $\Q(\sqrt{-p})$). Then the $L$-function of the $f$-eigenspace of the Tate module of $J_{\rho,\xi}(p)$ extends holomorphically to the complex plane and satisfies the expected functional equation (whose sign is given by the product of the local Deligne-Langlands constants).
\end{itemize}}

To compute the (expected) signs of the functional equation satisfied by the $L$-functions of the eigenspaces, we need to introduce a distinction depending on the behavior of $\rho$ as a representation of $G_{\Q_p}$, which is expanded on in Proposition \ref{trichotomy-rho}. 

\defi{We say that $\rho$ is 
\begin{itemize}[noitemsep,label=\tiny$\bullet$]
\item \emph{split} if the image of a decomposition subgroup at $p$ is contained in a subgroup conjugate in $\GL{\F_p}$ to the group of diagonal matrices,
\item \emph{wild} if it is not tamely ramified at $p$,
\item $a$-\emph{Cartan} if the restriction to an inertia subgroup of $\rho$ is irreducible. In this case, the image by $\rho$ of an inertia subgroup at $p$ is the index $a$ subgroup of a nonsplit Cartan subgroup of $\GL{\F_p}$, and $a$ is an odd divisor of $\frac{p+1}{2}$. 
\end{itemize}}

\defi{Suppose that $\rho$ is $a$-Cartan, and let $C \leq \GL{\F_p}$ be a nonsplit Cartan subgroup such that the image of $\rho$ under an inertia subgroup $I_p$ at $p$ is $aC$. Let $E/\Q_p$ be the unique quadratic unramified extension. Choose one of the two isomorphisms $\iota: C \rar \F_{p^2}^{\times}$ respecting the trace and norm. The associated \emph{reciprocity index} (see Definition \ref{recip-indices-definition}) is the index of the largest subgroup of $I_p$ on which $\iota \circ \rho_{|I_p}$ agrees with the class field theory isomorphism $G_E \rar \widehat{E^{\times}}/(1+p\OO_E)$ (mapping an arithmetic Frobenius to a uniformizer).} 

\rem{By the link between class field theory and the fundamental characters \cite[Proposition 3]{Serre-image-ouverte}, if $\rho$ is $a$-Cartan, then $1$ is one of its reciprocity indices if and only if $\rho^{\ast}$ has Serre weight $2$ (in the sense of \cite[\S 2.2]{Conj-Serre}).  

In fact, if $\rho: G_{\Q} \rar \GL{\F_p}$ has determinant $\omega_p^{-1}$ and is $a$-Cartan, then $a$ is the greatest common divisor of $k-1$ and $p^2-1$, where $k$ is the Serre weight of $\rho$ (it is the same result if $k$ denotes the Serre weight of $\rho^{\ast}$) by \cite[\S 2.2]{Conj-Serre}. 

} 

The following result is a simplified version of Proposition \ref{root-numbers-surjective-computation}, which also describes the somewhat more complicated result for the CM eigenspaces.  

\theo[autom-2]{For $(f,\mathbf{1}) \in \mathscr{C}'_1$, let $\omega_f$ be the order of the character of $\F_{p^2}^{\times}$ attached to $f$. If $\rho$ is $a$-Cartan, let $n_{f,\rho} \in \{0,1\}$ be equal to $1$ if and only if $\frac{p+1}{\omega_f}$ is divisible by $a$ but $\frac{2(p+1)}{\omega_f}$ is not divisible by any of the two reciprocity indices of $\rho$.

Then product of the local constants attached to the non-CM eigenspaces of the Tate module of $J_{\rho,\xi}(p)$ (that is, the expected sign of the functional equation satisfied by their $L$-function) is given in Table \ref{root-numbers-easy}. 
}

\begin{table}[htb]
\centering
\begin{tabular}{|c|c|c|c|}
\hline
& \multicolumn{3}{c|}{Type of $\rho_{|G_{\Q_p}}$}\\
\hline
& Split & Wild & $a$-Cartan\\
\hline
$V_{f,\mathfrak{l}} \otimes \mrm{St}(\rho), (f,\mathbf{1}) \in \mathscr{S}$ & $(-1)^{(p+1)/2}$ & $(-1)^{(p-1)/2}a_p(f)$ & $(-1)^{(p+a)/2}$\\
\hline
$V_{f,\mathfrak{l}}\otimes \pi(\mathbf{1},\chi)(\rho), (f,\chi) \in \mathscr{P}'$ & \multicolumn{3}{c|}{$(-1)^{(p-1)/2}$}\\
\hline
$V_{f,\mathfrak{l}} \otimes C_f(\rho), (f,\mathbf{1}) \in \mathscr{C}'_1$ & $(-1)^{(p+1)/2}$ & $(-1)^{(p-1)/2}$ & $(-1)^{(p+1)/2+n_{f,\rho}}$\\
\hline 
\end{tabular}
\caption{Signs of the functional equation for eigenspaces of $J_{\rho,\xi}(p)$}
\label{root-numbers-easy}
\end{table}

In Section \ref{asymptotic-conseq-elliptic-curves}, we discuss the application of Theorem \ref{autom-2} when $\rho$ comes from an elliptic curve and $p \rar \infty$.  

We also study $J_{\rho}(p)[p]$ as a $\F_p[G_{\Q}]$-module, which has a very explicit semi-simplification. We show in Section \ref{residual-automorphy} that, thanks to recent progress in Langlands functoriality due to Arias-de-Reyna and Dieulefait in the preprint \cite{AdRD}, and to Newton and Thorne \cite{NewThor}, $J_{\rho}(p)[p] \otimes_{\F_p} \overline{\F_p}$ has the same semi-simplification as a direct sum of mod $p$ reductions of automorphic $p$-adic Galois representations. However, the practical utility of such a result is unclear. 

\theo[autom-3]{(Proposition \ref{jrho-mod-p-automorphic}) Assume that $\rho^{\ast}$ is surjective, or that it is irreducible with Serre weight in $\{2,p+1\}$, and that its Artin conductor $N$ is coprime to $3$. Let $V$ be a non-CM eigenspace of the $p$-adic Tate module of $J_{\rho,\xi}(p)$. Assume furthermore that \cite[Theorem 1.1]{AdRD} (that is, Theorem \ref{adrd}) holds. Then there exists regular algebraic cuspidal unitary automorphic representations $\pi_1,\ldots,\pi_r$ of $\operatorname{GL}(n_i)$ over $\Q$ and continuous representations $V_1,\ldots,V_r$ of $G_{\Q}$ with coefficients in $F_{\mathfrak{p}}$ (where $F \subset \C$ is a sufficiently large number field and $\mathfrak{p}$ is a maximal ideal of $\OO_F$ with residue characteristic $p$) such that the following two conditions are satisfied:
\begin{itemize}[noitemsep,label=\tiny$\bullet$]
\item For every $i$, $L(V_i,s)$ and $L(\pi_i,s)$ are well-defined formal products that agree up to shifting $s$ and perhaps finitely many Euler factors. 
\item The semi-simplifications of the mod $\mathfrak{p}$ reductions of $\bigoplus_{i}{V_i}$ and $V$ are isomorphic $\overline{\F_p}[G_{\Q}]$-modules. 
\end{itemize}
}

I am grateful to Pierre-Henri Chaudouard and Laurent Clozel for their answers to my questions about the Langlands Correspondence, Rankin-Selberg convolutions, and automorphy lifting theorems.  

\section{Relations between classical modular forms and automorphic forms}
\label{automorphic-dictionary-refresher}

The material in this Section is not new, but we spell it out carefully as a preparation to the computations of Section \ref{root-numbers-surjective}. Our references for the theory of automorphic forms are the books of Gelbart \cite{Gelbart} and Bump \cite{Bump}, and we use \cite{BH} for its clear description of local representations of $\operatorname{GL}(2)$. Some of what we discuss is also covered in, for instance, \cite[\S 2, \S 3]{LWnewforms}. 

Let $\psi: (\Z/M\Z)^{\times} \rar \C^{\times}$ be a primitive Dirichlet character. We can, for instance by \cite[Proposition 3.1.2]{Bump}, attach to $\psi$ a finite-order continuous character $\psi_{\mathbb{A}}: \mathbb{A}_{\Q}^{\times}/\Q^{\times}\prod_{q \nmid M}{\Z_q^{\times}} \rar \C^{\times}$ such that, for every prime $q \nmid M$, $\psi(q)$ is the image under $\psi_{\mathbb{A}}$ of the unique id\`ele $(z_v)_v$ such that $z_v=1$ if $v \neq q$ and $z_q=q$. In particular, the character $\prod_{q \mid M}{\Z_q^{\times}} \subset \mathbb{A}_{\Q}^{\times} \overset{\psi_{\mathbb{A}}}{\rar} \C^{\times}$ is exactly $\psi^{-1}$. Moreover, if $p \mid M$ is a prime, and $(z_v)_v$ is the id\`ele such that $z_v=1$ for $v \neq p$ and $z_p=p$, then $\psi_{\mathbb{A}}(z)=\psi'_p(p)$, where $\psi=\psi_p\psi'_p$, where $\psi_p, \psi'_p: (\Z/M\Z)^{\times} \rar \C^{\times}$ are Dirichlet characters whose conductors are a power of $p$ and prime to $p$ respectively. 

Let $f \in \mathcal{S}_k(\Gamma_1(N))$ be a normalized newform. We can attach to $f$ a cuspidal unitary automorphic representation $\pi_f$ of $\GL{\mathbb{A}_{\Q}}$. The representation $\pi_f$ is the restricted tensor product, over all places $v$ of $\Q$, of smooth admissible representations $\pi_{f,v}$ of $\GL{\Q_v}$. 

To describe the features of this correspondence, we need to introduce additive characters for all the local fields involved. When $p$ is a prime, we define the additive character 
\[\theta_p: \Q_p \rar \Q_p/\Z_p \subset \Q/\Z \overset{\exp(2i\pi \cdot)}{\longrightarrow}\C^{\times}. \]

When $F/\Q_p$ is a finite extension, we write $\theta_F = \theta_p \circ \mrm{Tr}_{F/\Q_p}$, and $\eta_F$ (or $\eta_p$ when $F=\Q_p$) for the character $\theta_F(\cdot/p)$, which has \emph{level one} in the sense of \cite[\S 1.7]{BH} when $F/\Q_p$ is tamely ramified (see for instance \cite[Theorem III.2.6]{Neukirch-ANT}), that is, $\eta_F$ vanishes on the maximal ideal of $\OO_F$, but not on $\OO_F$.    

The default character of $\R$ is $\theta_{\infty}: x \in \R \longmapsto \exp(-2i\pi x)$ and if $F$ is a finite extension of $\R$, we let $\theta_F=\theta_{\infty}\circ \mrm{Tr}_{F/\R}$. Thanks to this condition, for any number field $F$, $\prod_v{\theta_v}$ is a nontrivial additive character of $\mathbb{A}_F/F$. 

For any local field $F$, we denote by $\|\cdot\|_F$ the normalized absolute value.

\prop[langlands-elem]{Let $f \in \mathcal{S}_k(\Gamma_1(N))$ be a normalized newform with primitive character $\chi$. 
\begin{enumerate}[noitemsep,label=(\roman*)]
\item\label{le-1} For any place $v$ of $\Q$, $\pi_{f,v}$ has central character $(\chi_{\mathbb{A}})_{|\Q_v^{\times}}$. 
\item\label{le-2} Let $\psi$ be any Dirichlet character and $v$ be any place of $\Q$, then $\pi_{f \otimes \psi,v} \simeq \pi_{f,v} \otimes (\psi_{\mathbb{A}})_{|\Q_v^{\times}}(\det)$.
\item\label{le-3} For any $p \nmid N$, let $\alpha,\beta \in \C^{\times}$ be the two roots of $X^2-a_p(f)X+p^{k-1}\chi(p)$, then $\pi_{f,p}$ is the principal series representation attached to the two unramified unitary characters $p \longmapsto \alpha p^{(1-k)/2}$, $p \longmapsto \beta p^{(1-k)/2}$ of $\Q_p^{\times}$.  
\item\label{le-4} Let $L_{\infty}(f,s)=2(2\pi)^{-s}\Gamma(s)$, then, for any place $v$ of $\Q$, $L_v\left(f,s+\frac{k-1}{2}\right)=L(\pi_{f,v},s)$. 
\item\label{le-5} If $p$ is a prime number, $\varepsilon(\pi_{f,p},s,\theta_p) = p^{-v_p(N)\left(s-\frac{1}{2}\right)}\varepsilon_p(f)$, in the sense of Definition \ref{epsilon-for-modular-forms}. 
\item\label{le-6} $\varepsilon(\pi_{f,\infty},s,\theta_{\R}) = (-i)^k$.
\item\label{le-7} If $k \geq 2$, then $\pi_{f,\infty}$ is the \emph{discrete holomorphic series} of weight $k$, that is, the representation $\sigma(|x|^{(k-1)/2},|x|^{(1-k)/2}\mathfrak{s})$, where $\mathfrak{s}$ is the sign function. If $k=1$, then $\pi_{f,\infty}$ is the \emph{limit of discrete holomorphic series} $\pi(|x|^{1/2},|x|^{-1/2})$. 
\end{enumerate} }

\demo{\cite[Theorem 5.19 (a)]{Gelbart} attaches to $f$ a certain automorphic representation $\pi$ of $\GL{\mathbb{A}_{\Q}}$ which appears in $L^2_0(\GL{\Q}\backslash \GL{\mathbb{A}_{\Q}},\psi)$ for a certain continuous character $\psi: \mathbb{A}_{\Q}^{\times}/\Q^{\times} \rar \C^{\times}$. By \cite[Theorems 3.3.3, 3.3.4]{Bump}, this representation factors as a restricted tensor product of local representations. By \cite[Theorem 5.19]{Gelbart} and the discussion at the beginning of Section 5C of \emph{loc.cit.}, $\pi$ has conductor $N$: in particular, it is enough to test \ref{le-5} at $s=1/2$ by \cite[Remark 6.13 (ii)]{Gelbart}. We define $\pi_f$ as the \emph{contragredient} of this representation (which is the one constructed in, say, \cite[Section 3.6]{Bump}), which has the correct central character (and factorizes as a restricted tensor product), which concludes for \ref{le-1}. 

The discussion at the beginning of \cite[Section 5C]{Gelbart} shows that, when $p \nmid N$, $\pi_{f,p}$ is the principal series representation attached to two unramified characters $\mu_1,\mu_2$ of $\Q_p^{\times}$, with $\mu_i=\|\cdot\|^{-s_i}$, \[\mu_1\mu_2=(\chi_{\mathbb{A}})_{|\Q_p^{\times}},\quad \mu_1(p)+\mu_2(p)=p^{(1-k)/2}a_p(f).\]
Hence \[p^{s_1+\frac{k-1}{2}}+p^{s_2+\frac{k-1}{2}}=a_p(f),\quad p^{\left(s_1+\frac{k-1}{2}\right)+\left(s_2+\frac{k-1}{2}\right)}=p^{k-1}\chi(p),\] thus $\mu_1(p)p^{(k-1)/2},\mu_2(p)p^{(k-1)/2}$ are the two roots of $X^2-a_p(f)X+p^{k-1}\chi(p)$. This proves \ref{le-3} and \ref{le-4} for finite places at primes not dividing $N$. 

For \ref{le-7}, the situation when $k \geq 2$ is discussed at the beginning of \cite[Section 5C]{Gelbart}, and for $k=1$ it is \cite[Remark 2.5.2]{GelbartFLT}. This implies in particular \ref{le-4} at the infinite place. 

In particular, if $\psi$ is a Dirichlet character, the automorphic representations $\pi_f \otimes \psi_{\mathbb{A}}(\det)$ and $\pi_{f \otimes v}$ have the same components at infinity and the sane components at all but finitely many places, so they are isomorphic by strong multiplicity one \cite[Theorem 3.3.6]{Bump}. 

Now, the function $R(f,s)=\frac{L\left(\pi_f,s-\frac{k-1}{2}\right)}{L(f,s)}$ can be written, by \cite[Theorem 6.15]{Gelbart}, as a product $\prod_{p\mid N}{\frac{1-a_pp^{-s}}{1-b_pp^{-s}}}$, with $|a_p| \leq p^{(k-1)/2}$ by \cite[Theorem 4.6.17]{Miyake} and $|b_p| \leq p^{(k-1)/2}$ because $\pi_{f,p}$ is unitary (and the $L$-factors of ramified representations have degree at most $1$ in $p^{-s}$). By Proposition \ref{functional-eqn} and \cite[Theorem 6.18]{Gelbart}, we have a functional equation $R(f,s)=\alpha R(\overline{f},k-s)$, where $R(\overline{f},s)$ is defined similarly and is of the form $\prod_{p \mid N}{\frac{1-a'_pp^{-s}}{1-b'_pp^{-s}}}$ with $|a'_p|,|b'_p| \leq p^{(k-1)/2}$, and $\alpha \in \C^{\times}$. In particular, for each $p$, there is a constant $\alpha_p \in \C^{\times}$ such that  

\[(1-a_pp^{-s})(1-b'_pp^{s-k})=\alpha_p(1-b_pp^{-s})(1-a'_pp^{s-k}),\quad \prod_p{\alpha_p}=\alpha.\] 

By computing the coefficients in $p^s$ and $p^{-s}$, we see that $b'_p=\alpha_pa'_p$ and $a_p=\alpha_pb_p$, so that, by computing the constant coefficient, we find 

\[1+\alpha_p^2\frac{a'_pb_p}{p^k}=\alpha_p(1+\frac{a'_pb_p}{p^k}),\] so either $\alpha_p=1$ or $a'_pb_p \neq 0$ and $\alpha_p=\frac{p^k}{a'_pb_p}$, so that either $\alpha_p=1$ or $|\alpha_p| > 1$. 

Let us show that the product of the $\alpha_p$ has module one. Indeed, the local constants in the functional equation for modular forms have module one by Propositions \ref{local-constant-elementary} and \ref{functional-eqn}. On the automorphic side, this follows from \cite[Theorem 6.18]{Gelbart} and \cite[Proposition 2.21.2 (ii)]{JL} for cuspidal representations, and \cite[Propositions 3.5, 3.6]{JL} reduce the claim to the case of $\operatorname{GL}(1)$ for unitary characters by \cite[Theorem 4.27]{Gelbart}\footnote{The definitions of \cite{Gelbart} and \cite{JL} are compatible.} which follows from the explicit formula \cite[(3.6.3)]{NTB}. 

Therefore, every $\alpha_p$ is equal to one, so $a_p=b_p$ and we are done. This proves \ref{le-4} in full generality. 

To establish \ref{le-5}, we already saw that it was enough to prove it at $s=1/2$, which is claimed in \cite[(2.12)]{Li-RS}. All that remains to do is prove \ref{le-6}, which is a consequence of the automorphic functional equation \cite[Theorem 6.18]{Gelbart}, \ref{le-5}, and the modular functional equation from Proposition \ref{functional-eqn} and Corollary \ref{local-to-global-constant}. 
}

\medskip

The above properties are sufficient to determine the global dictionary, both from a global point of view (using the modular \cite[Theorem 4.6.19]{Miyake} and automorphic \cite[Theorem 3.3.6]{Bump} multiplicity one theorems) and from a local point of view (using the converse theorem \cite[\S 27.1]{BH} characterizing irreducible smooth representations by the $L$-factors and $\varepsilon$-factors  attached to their twists). But we are interested in describing the correspondence more precisely in certain specific cases.

\prop[langlands-elem-2]{Let $f \in \mathcal{S}_k(\Gamma_1(N))$ be a normalized newform, and $\chi$ be its primitive Dirichlet character with conductor $M$. 
\begin{itemize}[noitemsep,label=$-$]
\item Let $p \mid N$ be prime such that $p^2 \nmid NM$. Then $\pi_{f,p}$ is the twist of the Steinberg representation of $\GL{\Q_p}$ by the unramified character $p \longmapsto a_p(f)p^{1-k/2}$.
\item Let $p \mid N$ be prime and coprime to $N/M$. Then $\pi_{f,p}$ is the principal series representation of $\GL{\Q_p}$ attached to the characters $\alpha,\beta: \Q_p^{\times} \rar \C^{\times}$ such that $\alpha$ is unramified and maps $p$ to $a_p(f)p^{(1-k)/2}$, and $\alpha\beta=(\chi_{\mathbb{A}})_{|\Q_p^{\times}}$. 
\end{itemize}
}

\demo{In the first case, Proposition \ref{langlands-elem} implies that $L(\pi_{f,p},s)=1-a_p(f)p^{-(k-1)/2-s}$ with $a_p(f)^2=p^{k/2-1}\chi(p)$ (Proposition \ref{bad-L-factor}), so that $\pi_{f,p}$ is not cuspidal by \cite[\S 24.5]{BH}. However, for any non trivial Dirichlet character $\psi$ whose conductor is a $p$-th power, $a_p(f \otimes \psi)=0$ by Corollary \ref{level-newform-bigtwist}, so $L(\pi_{f,p} \otimes (\psi_{\mathbb{A}})_{|\Q_p^{\times}}(\det),s)=1$. Assume that $\pi_{f,p}$ is a principal series representation attached to two characters $\alpha,\beta$ of $\Q_p^{\times}$. Since $L(\pi_{f,p},s) \neq 1$, by \cite[\S 23.4, 26.1]{BH}, we may assume that $\alpha$ is unramified. But then, there is a nontrivial character $\gamma: \Q_p^{\times}/p^{\Z} \rar \C^{\times}$ such that $\beta\gamma$ is unramified, so that $L(\pi_{f,p} \otimes \gamma(\det),s)\neq 1$. By construction, we may write $\gamma=\psi_{|\mathbb{A}})_{|\Q_p^{\times}}$ for some nontrivial Dirichlet character $\psi$ of conductor some power of $p$: this is a contradiction. So $\pi_{f,p}$ is a twist of the Steinberg representation, and by considering its $L$-function (using \cite[\S 26.1]{BH}) the twisting character must map $p$ to $a_p(f)p^{1-k/2}$. 

In the second case, $L(\pi_{f,p},s)=1-a_p(f)p^{-s-(k-1)/2} \neq 1$ by Proposition \ref{bad-L-factor}, so as above $\pi_{f,p}$ is not cuspidal. Assume that $\pi_{f,p}$ is a twist of a Steinberg representation by the character $\phi$: since $L(\pi_{f,p},s) \neq 1$, $\phi$ is unramified, and the computation above shows that $\phi(p)=a_p(f)p^{1-k/2}$. One then has $L(\pi_{f \otimes \overline{\chi},p},s)=L(\pi_{f,p} \otimes (\chi_{\mathbb{A}})_{|\Q_p^{\times}}(\det),s)=L(\phi(\chi_{\mathbb{A}})_{|\Q_p^{\times}},s+1/2)=1$ because the character is ramified. Yet one has $L(\pi_{f \otimes \overline{\chi},p},s)=L_p(\overline{f},s+(k-1)/2)=1-\overline{a_p(f)}p^{-s-(k-1)/2} \neq 1$, so we have a contradiction. Thus $\pi_{f,p}$ is the principal series induced from two characters $\alpha,\beta$. By considering its $L$-function, we can assume that $\alpha$ is unramified and maps $p$ to $a_p(f)p^{(1-k)/2}$. Moreover, $\alpha\beta$ is the central character of $\pi_{f,p}$, which is exactly $(\chi_{\mathbb{A}})_{|\Q_p^{\times}}$.}

\bigskip

\prop[langlands-classify-c]{Let $(f,\chi) \in \mathscr{C}$. Let $E/\Q_p$ be the unique unramified quadratic extension. There is a character $\alpha: E^{\times}/(1+p\OO_E) \rar \C^{\times}$, unique up to $E/\Q_p$-conjugation, such that $\pi_{f,p}$ is the representation $\pi_{\alpha}$ from \cite[\S 19.1]{BH} (so $\alpha$ does not factor through the norm $E^{\times} \rar \Q_p^{\times}$). Moreover, for each $u \in \Z_p^{\times}$, $\alpha(u)=\overline{\chi}(u)$ and $\alpha(p)=1$. 
If $\psi \in \mathcal{D}$, then the character thus attached to $(f \otimes \psi,\chi\psi^2) \in \mathscr{C}$ is exactly $\alpha \cdot \overline{\psi(N_{E/\Q_p})}$, with the convention that $\psi(p)=1$. }

\rem{If $\chi \in \mathcal{D}$, the character $\Q_p^{\times} \rar \C^{\times}$ that maps $p$ to $1$ and agrees with $\overline{\chi}$ on $\Z_p^{\times}$ is exactly $(\chi_{\mathbb{A}})_{|\Q_p^{\times}}$. }

\demo{By \cite[\S 24.3]{BH} and Proposition \ref{langlands-elem}, one has $n(\pi_{f,p},\theta_p)=2$ so $n(\pi_{f,p},\eta_p)=0$. Let us now show that $\pi_{f,p}$ is cuspidal. By \cite[\S 27.2]{BH} and Proposition \ref{langlands-elem}, it is enough to show that for any Dirichlet character $\psi$ whose conductor is a power $P$ of $p$, one has $a_p(f \otimes \psi)=0$. When $\psi \neq \overline{\chi}$, $f \otimes \psi$ has conductor $P^2$ by Corollary \ref{level-newform-bigtwist}, so by Proposition \ref{bad-L-factor} one has $a_p(f \otimes \psi)=0$. Finally, $a_p(f \otimes \overline{\chi})=a_p(\overline{f})=\overline{a_p(f)}=0$.  

By \cite[(25.2.3)]{BH}, since $\eta_p$ has level one, $\pi_{f,p}$ has normalized level zero, and the claim is a consequence of \cite[\S 19.1 Proposition]{BH} (and the above remark for the central character and the behavior under twists). }

\medskip

\cor[cuspidal-alpha-phi]{Let $(f,\chi) \in \mathscr{C}$ and $E/\Q_p$ be the quadratic unramified extension. By Proposition \ref{langlands-classify-c}, let $\alpha: E^{\times} \rar \C^{\times}$ be the character attached to $\pi_{f,p}$. Then the cuspidal representation $C_f$ is attached to the character $\overline{\alpha}_{|\OO_E^{\times}}$. }

\demo{To ease notation, write $\tau=\mrm{Tr}_{\F_{p^2}/\F_p}$ and $\nu=N_{\F_{p^2}/\F_p}$. For any $\psi \in \mathcal{D}$, using (in this order), Proposition \ref{langlands-elem}, Proposition \ref{langlands-classify-c}, and \cite[\S 24.3, 25.2, 25.4, 23.7]{BH}, 

\begin{align*}
\lambda_p(f \otimes \psi)&=\varepsilon_p(f \otimes \psi) = \varepsilon(\pi_{f \otimes \psi,p},\frac{1}{2},\theta_p) = \varepsilon(\pi_{\alpha \cdot \overline{\psi}(N_{E/\Q_p})},\frac{1}{2},\eta_p)\\
&= \frac{-p}{(p^4)^{1/2}}\sum_{u \in (\OO_E/p)^{\times}}{\overline{\alpha}(u)\psi(\nu(u))\eta_p(\tau(u))}.
\end{align*}

Let $M \in \GL{\F_p}$ be a generator of some nonsplit Cartan subgroup, and let $a,N \in \F_p^{\times}$ be its trace and determinant, respectively. Let $g \in (\OO_E/p)^{\times}$ be one of the two elements with trace $a$ and norm $N$. One has $\overline{\alpha}(g)\overline{\alpha}(g^p)=\overline{\alpha}(\nu(g))=\chi(N)=\phi(g)\phi(g^p)$ and, by Corollary \ref{various-computations-cuspidal}, 
\begin{align*}
\mrm{Tr}(M \mid C_f) &= \frac{\chi(a)}{p-1}\sum_{\psi \in \mathcal{D}}{\psi\left(\frac{a^2}{N}\right)\lambda_p(f_{\psi})\mathfrak{g}(\overline{\chi\psi^2})}\\
&= \frac{-\chi(a)}{p(p-1)}\sum_{\psi \in \mathcal{D}}{\psi\left(\frac{a^2}{N}\right)\sum_{t \in (\OO_E/p)^{\times}}{\overline{\alpha}(t)\psi(\nu(t))e^{2i\pi\frac{\tau(t)}{p}}}\sum_{u \in \F_p^{\times}}{\overline{\chi\psi^2}(u)e^{2i\pi u/p}}}\\
&= \frac{-1}{p(p-1)}\sum_{u \in \F_p^{\times}}{\sum_{t \in (\OO_E/p)^{\times}}{\sum_{\psi \in \mathcal{D}}{\overline{\alpha}(\frac{at}{-u})e^{2i\pi \frac{u}{p}\left(1-\frac{1}{a}\tau(\frac{at}{-u})\right)}\psi\left(\frac{1}{N}\nu(\frac{at}{-u})\right)}}}\\
&= \frac{-1}{p}\sum_{\substack{t \in (\OO_E/p)^{\times}\\\nu(t)=N}}{\sum_{u \in \F_p^{\times}}{\overline{\alpha}(t)e^{2i\pi\frac{u}{p}\left(1-\frac{\tau(t)}{a}\right)}}}=\frac{-1}{p}\sum_{\substack{t \in (\OO_E/p)^{\times}\\\nu(t)=N}}{\sum_{u \in \F_p}{\overline{\alpha}(t)e^{2i\pi\frac{u}{p}\left(1-\frac{\tau(t)}{a}\right)}}}\\
&= -\sum_{\substack{t \in (\OO_E/p)^{\times}\\\tau(t)=a\\\nu(t)=N}}{\overline{\alpha}(t)} = -\alpha(g)-\alpha(g^p),
\end{align*} 

so that $\overline{\alpha}(g) \in \{\phi(g),\phi^p(g)\}$, whence $\{\overline{\alpha}_{|\OO_E^{\times}},\,\overline{\alpha}^p_{|\OO_E^{\times}}\} = \{\phi,\phi^p\}.$  
}

\section{Basics on the Langlands correspondence}
\label{langlands-correspondence-refresher}

The following discussion of Weil groups is somewhat terse. We refer the reader wishing for more detail to Tate's article \cite{NTB}. 

In this paragraph, all our local fields are assumed to be of characteristic zero. 

To a non-Archimedean local field $F$ one attaches a \emph{Weil group} $W_F$: it is a topological group endowed with a reciprocity isomorphism $\mathfrak{a}_F: W_F^{ab} \rar F^{\times}$ and a continuous homomorphism $W_F \rar G_F$. Moreover, $\mathfrak{a}_F$ is the restriction of the class field theory isomorphism $G_F \rar \hat{F^{\times}}$, which is normalized so that it maps a \emph{geometric} Frobenius to a uniformizer of $F$. The Weil group is constructed as follows in \cite{NTB}: let $k$ be the residue field of $F$, then $W_F$ is the subgroup of $G_F$ made with the $\sigma \in G_F$ such that the image of $\sigma$ in $\mrm{Gal}(\overline{k}/k)$ is an integer power of the Frobenius (arithmetic or geometric). The topology on $W_F$ is defined by the following properties:
\begin{itemize}[noitemsep,label=$-$]
\item $I_F \subset W_F$ is an open subgroup, 
\item the topologies induced on $I_F$ by $I_F \leq W_F$ and $I_F \leq G_F$ are the same, 
\item the image of the map $W_F \rar G_F \rar G_k$ is exactly $\langle F\rangle \simeq \Z$, where $F$ is the Frobenius. We require that the induced map $W_F \rar \langle F\rangle \simeq \Z$ be continuous, where $\Z$ is endowed with the discrete topology. 
\end{itemize}

Finally, we denote by $\|\cdot\|: W_F \rar \R^{+\times}$ the continuous group homomorphism given by the composition $W_F \overset{\mathfrak{a}_F}{\rar} F^{\times} \overset{\|\cdot\|}{\rar} \R^{+\times}$. In particular, it maps an \emph{arithmetic} Frobenius to the cardinality of the residue field of $\OO_F$. 

When $F$ is a local Archimedean field, there is also a definition of a topological group $W_F$ along with its morphism to $G_F$, a reciprocity isomorphism $\mathfrak{a}_F: W_F^{ab} \rar F^{\times}$ and a continuous group homomorphism $\|\cdot\|: W_F \rar \R^{+\times}$. When $F$ is complex, $W_F$ is exactly $F^{\times}$; when $F$ is real, $W_F = \overline{F}^{\times} \cup j\overline{F}^{\times}$, where $j^2=-1$ and conjugation by $j$ on $\overline{F}^{\times}$ is the nontrivial automorphism in $\mrm{Gal}(\overline{F}/F)$. 

Its definition is less explicit, but any number field $F$ also has a Weil group $W_F$. It is endowed with a reciprocity isomorphism $\mathfrak{a}_F: W_F^{ab} \rar \mathbb{A}_F^{\times}/F^{\times}$, and with a surjective homomorphism to $G_F$ such that $\mathbb{A}_F^{\times}/F^{\times} \rar W_F^{ab} \rar G_F^{ab}$ is exactly the class field theory surjection (normalized like $\mathfrak{a}$: that is, it maps a uniformizer to a geometric Frobenius. 

The Weil group is well-behaved with respect to extensions of local fields or number fields, as well as localizations at various places. We refer to \cite[\S 1]{NTB} for more information. 

\medskip

\defi{When $F$ is a local field (resp. a global field), a \emph{representation} of $W_F$ over a field $K$ of characteristic $0$ (resp. over $\C$) is a group homomorphism $W_F \rar \operatorname{GL}(V)$ for some finite-dimensional $K$-vector space $V$ with open kernel (resp. a continuous homomorphism $W_F \rar \operatorname{GL}(V)$, where $\operatorname{GL}(V)$ is endowed with its Archimedean topology).}

\defi{When $F$ is local non-Archimedean, a \emph{Deligne representation} of $W_F$ with coefficients in a field $K$ of characteristic zero is a couple $(\rho,N)$, where $\rho$ is a representation of $W_F$ on a finite-dimensional $K$-vector space $V$ and $N$ is a (necessarily nilpotent) endomorphism of $V$ such that, for any $w \in W_F$, $\rho(w)N\rho(w)^{-1}=\|w\|N$. We say that $(\rho,N)$ is a semi-simple Deligne representation if $\rho$ is a semi-simple representation of $W_F$. If needed, a representation of $W_F$ corresponds to the Deligne representation $(\rho,0)$ of $W_F$.   }

\medskip

The interest of Deligne representations is captured by the following result \cite[(4.2.1)]{NTB}: 

\prop[l-adic-to-weil-deligne]{Let $F$ be a finite extension of $\Q_p$ and $E$ be a finite extension of $\Q_{\ell}$ where $p,\ell$ are distinct primes. Choose a Frobenius $\phi \in W_F$ and a pro-$\ell$-tame inertia character $t_{\ell}: I_F \rar \Q_{\ell}$. To any Deligne representation $(\rho,N)$ of $W_F$ on a finite-dimensional $E$-vector space $V$, one can attach $\rho_E: \phi^n\sigma \in W_F \longmapsto \rho(\phi^n\sigma)\exp(t_{\ell}(\sigma)N) \in \GLn{E}{V}$ (where $\sigma \in I_F$ and $n \in \Z$), which is a continuous group homomorphism $W_F \rar \operatorname{GL}(V)$ (where $\operatorname{GL}(V)$ is endowed with the $\ell$-adic topology). 

This operation defines a bijection between the isomorphism classes of Deligne representations of $W_F$ over $E$ and the continuous representations of $W_F$ over $E$, and it is independent from the choices of $t_{\ell}$ and $\phi$.
}

\rem{In \cite[Definition 12.2.1]{Getz-Hahn}, complex Deligne representations (or Weil-Deligne representations) are continuous representations of $W_F \times \SL{\C}$ over a finite-dimensional $\C$-vector space such that their $\SL{\C}$-component is algebraic. These two formalisms are equivalent, as is described in \cite[\S 2.1]{Gross-Reeder}. }

Let $V$ be a continuous representation of a $p$-adic field $F$ with coefficients in a $\ell$-adic field $E$. We denote the Weil-Deligne representation attached to $V_{|W_F}$ by $\Del{V}$, if necessary by adding $F$ or $E$ as an index. Usually, the choice of $E$ is immaterial, $F$ is $\Q_p$ for some prime $p$, and $V$ comes from a continuous representation of $G_{\Q}$ over $E$, and we will shorten $\Del{V_{|G_{\Q_p}}}$ in $\Delp{p}{V}$.

To complex representations (resp. Deligne representations) of $W_F$ when $F$ is an Archimedean (resp. non-Archimedean) local field, Deligne and Langlands attach $L$- and $\varepsilon$-factors, which are meromorphic functions of a complex variable $s$. The definition of the local $\varepsilon$ also requires the datum of a nontrivial additive character $\psi$ of $(F,+)$. While we will state the definition of the $L$-factors, the definition of the $\varepsilon$ is more involved, and we simply recall some of their properties in Proposition \ref{local-constants-elem} for easier manipulation. We refer to \cite{NTB} and \cite[Chapter 7]{BH} for more information. We follow the normalization of the latter, called ``Langlands's local constants'' in \cite[(3.6)]{NTB}.

\defi{Let $F$ be a non-Archimedean local field with Weil group $W_F$, inertia group $I_F$, geometric Frobenius $\phi$, and let $q$ be the cardinality of the residue field. For $t \geq -1$, let $I_F^t$ denote the higher inertia subgroup in upper numbering (see \cite[Chapter IV]{Serre-local}, for instance). 

If $(\rho,N)$ is a Deligne representation of $W_F$ on some finite-dimensional complex vector space $V$, its $L$-function is the meromorphic function (see \cite[(4.1.5)]{NTB} or \cite[(31.3.1)]{BH}) \[L(V,s)=\det(1-q^{-s}\rho(\phi)\mid (\ker{N})^{I_F})^{-1}.\]  

If $\rho$ is a representation of $W_F$ in a characteristic zero field, then its conductor exponent $f(\rho)$ is given by the formula in \cite[(4.5)]{Deligne-csts} (an equivalent definition is in \cite[\S 4]{CF6}). It is not too difficult, using \cite[VI.2 Exercise 1a]{Serre-local} as an intermediate step, to deduce that \[f(\rho)=\int_{-1}^{\infty}{\dim{V/V^{I_F^t}}\,dt}.\] (this formula appears in \cite[Definition 3.1.27]{Wiese-Galois}; because upper numbering commutes with quotients, the formula does not have to involve the kernel of $\rho$). 

If $V$ is a Deligne representation of $W_F$, there is a modified formula for $a(V)$ in \cite[(4.1.5)]{NTB}.}

We can then do the following direct computation using Proposition \ref{l-adic-to-weil-deligne}. 

\lem[l-adic-to-weil-deligne-numerics]{Let $E,F$ be two non-Archimedean local fields with coprime residue characteristics. Let $V$ be a finite dimensional $E$-vector space endowed with a continuous action $\rho: G_F\rar \GLn{E}{V}$. Let $I_F$ denote the inertia group of $G_F$ and $I_F^t$ denote the higher ramification groups in upper numbering. Then:
\begin{itemize}[noitemsep,label=$-$]
\item Fix an embedding $E \rar \C$. Then $L(\Del{V^{\ast}} \otimes_E \C, s)$ is exactly the evaluation of the rational fraction $\det(1-t\rho(\Fr_F)\mid V_{I_F}) \in E(t) \subset \C(t)$ at $q^{-s}$, where $q$ is the cardinality of the residue field of $F$ and $V_{I_F}$ is the space of \emph{co-invariants} of $V$ under the inertia. 
\item $f(\Del{V})=f(\Del{V^{\ast}})=\int_{-1}^{\infty}{\dim_E(V/V^{I_F^t})\,dt}$.
\end{itemize}}

\medskip

\prop[local-constants-elem]{Let $F$ be a non-Archimedean local field with uniformizer $\varpi$, $s \in \C$ and $\psi: F \rar \C^{\times}$ be nontrivial additive character. Let $n_{\psi}$ be the largest integer $n$ such that $\psi_{|\varpi^{-n}\OO_F}$ is trivial (so if $F$ is unramified over $\Q_p$, then $n_{\theta_F}=0$). Let $V$ be a complex Deligne representation of $W_F$ with conductor exponent $f(V)$. Then:
\begin{itemize}[noitemsep,label=\tiny$\bullet$]
\item For any $a=\mathfrak{a}_F(\alpha) \in F^{\times}$, $\varepsilon(V,s,\psi(a \cdot))=\det{V}(\alpha)\|a\|^{s-\frac{1}{2}}\varepsilon(V,s,\psi)$. 
\item Write $-1=\mathfrak{a}_F(c)$, then $\varepsilon(V,s,\psi)\varepsilon(V^{\ast},1-s,\psi)=\det{\rho}(c)$.  
\item For any $s' \in \C$, $\varepsilon(V,s+s',\psi)=\varepsilon(V \otimes \|\cdot\|^{s'},s,\psi)$. 
\item Let $E/F$ be a finite extension and assume that $V=\mrm{Ind}_E^F{V_1}$. Then there is a constant $\lambda_{E/F}(\psi)$ depending only on $E/F$ and $\psi$ (not on $V_1$) such that \[\varepsilon(V,s,\psi)=\varepsilon(V_1,s,\psi\circ \mrm{Tr}_{E/F})\lambda_{E/F}(\psi)^{\dim{V_1}}.\] 
\item The map $W \longmapsto \varepsilon(V \otimes W,s,\psi)$ (where $W$ is a representation of $W_F$) is additive.\footnote{This is false if we assume only that $W$ is a Deligne representation.} 
\item For any unramified $\chi: F^{\times} \rar \C^{\times}$, $\varepsilon(V \otimes \chi(\mathfrak{a}_F),s,\psi)=\varepsilon(V,s,\psi)\chi(\varpi^{n_{\psi}\dim{V}+f(V)})$. 
\item Let $W$ be an unramified representation of $W_F$ and $\phi \in W_F$ be a geometric Frobenius, then \[\varepsilon(V \otimes W,s,\psi)=\varepsilon(V,s,\psi)^{\dim{W}}\det{W}(\phi^{n_{\psi}\dim{V}+f(V)}).\]
\item For any unramified $\chi: F^{\times} \rar \C^{\times}$, $\varepsilon(\chi(\mathfrak{a}_F),s,\psi)=\chi(\varpi^{n_{\psi}})\|\varpi\|^{n_{\psi}(s-1/2)}$.
\item With $V=\chi(\mathfrak{a}_F)$ for some $\chi: F^{\times} \rar \C^{\times}$, for any $y \in \varpi^{f(V)+n_{\psi}}\OO_F^{\times}$, 
\[\varepsilon(V,s,\psi)=\chi(y)\|\varpi\|^{(s-1/2)(f(V)+n_{\psi})}\frac{\sum_{u \in (\OO_F/\varpi^{f(V)})^{\times}}{\chi^{-1}(u)\psi\left(\frac{u}{y}\right)}}{\left|\sum_{u \in (\OO_F/\varpi^{f(V)})^{\times}}{\chi^{-1}(u)\psi\left(\frac{y}{u}\right)}\right|}.\]
\end{itemize}
The first five claims also hold when $F$ is Archimedean.}

\demo{These are easy computations using \cite[\S 3.6]{NTB} and \cite[\S 29.4]{NTB}. }

\rem{In particular, this implies that the conductor exponent of a Deligne representation $V$ of $W_F$ is exactly the $n(V,\psi)$ of \cite[\S 29.4 Proposition]{BH} with $\psi: F/\OO_F \rar \C^{\times}$ an additive character whose restriction to $\varpi_F^{-1}\OO_F/\OO_F$ ($\varpi_F$ being a uniformizer of $F$) is nontrivial. }

These relations imply the following reformulation of \cite[Theorem (3.5.3)]{NTB}:

\prop[functional-eqn-weil-gp]{Let $K$ be a number field with absolute discriminant $\Delta$ and $V$ be a complex, finite-dimensional representation of $W_K$ with absolute conductor $M$. Let $\Lambda(V,s),\Lambda(V^{\ast},s)$ be the formal products $(M|\Delta|)^{s/2}\prod_{v}{L(V_{|W_v},s)},(M|\Delta|)^{s/2}\prod_{v}{L(V^{\ast}_{|W_v},s)}$ respectively (over all places of $K$). These formal products are convergent in some right-half plane, define holomorphic functions on this half-plane which extend as meromorphic functions on $\C$, and satisfy, for any nontrivial additive character $\psi=\prod_v{\psi_v}: \mathbb{A}_K/K \rar \C^{\times}$, the functional equation \[\Lambda(V,s)=\prod_v{\varepsilon(V_{|W_v},\frac{1}{2},\psi_v)}\Lambda(V^{\ast},1-s).\]
}

\medskip

The \emph{local Langlands correspondence} is the following result: for any $n \geq 1$ and any local field $F$, there exists a bijection between the isomorphism classes of complex irreducible smooth admissible representations of $\GLn{n}{F}$ and the isomorphism classes of semi-simple $n$-dimensional Deligne representations of $W_F$. This correspondence preserves the $L$ and $\epsilon$ constants on either side, as well as those of pairs (corresponding to the Rankin-Selberg convolution on the automorphic side and to the tensor product for representations of $W_F$), it respects the determinants, twists by characters of $F^{\times}$ (via $\mathfrak{a}_F$) and taking contragredients. A formal statement can be found in the beginning of \cite{BCarayol} or in \cite[Chapter 12.4]{Getz-Hahn}. 

Let us first spell out the situation for Dirichlet characters, since we will need it. 

\lem{Let $\omega: G_{\Q} \rar \hat{\Z}^{\times}$ be the cyclotomic character and $\chi: (\Z/N\Z)^{\times} \rar \C^{\times}$ be a primitive Dirichlet character. Then, for any prime number $p$, $\Delp{p}{\chi(\omega)} = \overline{\chi}_{\mathbb{A}} \circ \mathfrak{a}_{\Q_p}$. }

\demo{If $\sigma \in \mrm{Gal}(\overline{\Q_p}/\Q_p)$ lies in the inertia group, it is known by \cite[\S 3, Theorem 2]{CF6} (note the different normalization for the class field theory isomorphism) that $\mathfrak{a}_{\Q_p}(\sigma)=\pi_p(\omega(\sigma))$, where $\pi_p: \hat{\Z}^{\times} \rar \Z_p^{\times}$ is the projection. Moreover, $\omega$ acts trivially on the maximal unramified extension of $\Q_p$, which contains all the roots of unity in $\overline{\Q_p}^{\times}$ of order prime to $p$. Since $(\overline{\chi}_{\mathbb{A}})_{|\Z_p^{\times}}$ is exactly the $p$-part of $\chi$ (that is, a character $\chi_p: \Z_p^{\times} \rar \C^{\times}$ such that the conductor of the Dirichlet character $\chi'_p=\chi\chi_p^{-1}$ is prime to $p$), the claimed equality holds for $\sigma$. 

Let $F \in \mrm{Gal}(\overline{\Q_p}/\Q_p)$ be an inverse image under $\mathfrak{a}_{\Q_p}$ of $\frac{1}{p}$ (so that $F$ is an \emph{arithmetic} Frobenius). By \cite[\S 3, Theorem 2]{CF6}, $\pi_p(\omega(F))=1$ and for any $q \neq p$, $\pi_q(\omega(F))=p$, so that $\chi(\omega(F))=\chi'_p(p)$. Now, we compute that 

\[(\overline{\chi}_{\mathbb{A}})_{|\Q_p^{\times}}\left(\frac{1}{p}\right)=\prod_{\ell \neq p}{(\overline{\chi}_{\mathbb{A}})_{|\Z_{\ell}^{\times}}(p)}=\prod_{\ell \neq p}{\chi_{\ell}(p)}=\chi'_p(p).\]
}

\bigskip

Let $f \in \mathcal{S}_k(\Gamma_1(N))$ be a newform and $p$ be a prime number. We can attach Deligne representations of $W_{\Q_p}$ to $f$ in two seemingly distinct manners. A first possibility is to fix a number field $F \subset \C$ containing the coefficients of $f$ and a prime ideal $\lambda \subset \OO_F$ prime to $p$, then consider $\Delp{p}{V_{f,\lambda}} \otimes_{F_{\lambda}} \C$ for some $F$-embedding of $F_{\lambda}$ into $\C$. A second possibility is to consider the image $\sigma_{f,p}$ of $\pi_{f,p}$ under the local Langlands correspondence. 

When $p \nmid N$, $\sigma_{f,p}$ and $(V_{f,\lambda})_{|W_{\Q_p}}$ are both sums of unramified characters of $W_{\Q_p}$, so relating them is not difficult. The same relation still holds at primes $p \mid N$ (but $p$ coprime to $\lambda$), but this is significantly more delicate, and is the \emph{local-global compatibility} proved by Carayol in \cite{localglobal}.

\prop[local-global-modular-forms]{Let $f \in \mathcal{S}_k(\Gamma_1(N))$ be a normalized newform. Let $E \subset \C$ be a number field containing the coefficients of $f$ and $\lambda \subset \OO_E$ a prime ideal. 
\begin{itemize}[noitemsep,label=$-$]
\item Let $p$ be a prime number prime to $\lambda$ and $\sigma: E_{\lambda} \rar \C$ be a $E$-homomorphism. The isomorphism class of $\Delp{p}{V_{f,\lambda}^{\ast}} \otimes_{E_{\lambda}} \C$ does not depend on the choice of $\lambda$ or $\sigma$. We denote it by $D_{f,p}$ or sometimes $(V_f^{\ast})_{|W_{\Q_p}}$. 
\item The complex Deligne representation $\sigma_{f,p}$ of $W_{\Q_p}$ attached by the local Langlands correspondence to $\pi_{f,p}$ is $D_{f,p} \otimes \|\cdot\|^{(k-1)/2}$. In particular, $\varepsilon_p(f)=\varepsilon(D_{f,p},\frac{k}{2},\theta_p)$. 
\item The complex representation $\sigma_{f,\infty}$ of $W_{\R}$ attached to $\pi_{f,\infty}$ is exactly $\|\cdot\|^{(k-1)/2} \otimes \mrm{Ind}_{W_{\C}}^{W_{\R}}{z^{1-k}}$. 
\item Let $\chi: (\Z/M\Z)^{\times} \rar E^{\times}$ be a primitive Dirichlet character and $p$ be a prime. Then $D_{f,p} \otimes \chi_{\mathfrak{A}}(\mathfrak{a}_{\Q_p}) = D_{f \otimes \chi,p}$.  
\end{itemize}}

\demo{Let us first discuss the place at infinity. By Proposition \ref{langlands-elem}, when $k \geq 2$, $\pi_{f,\infty}$ is the discrete holomorphic series of weight $k$ ($\sigma(|x|^{(k-1)/2},|x|^{(1-k)/2}\mathfrak{s})$), so that by \cite[\S 3.1]{GelbartFLT} $\sigma_{f,\infty}$ is exactly the induced to $W_{\R}$ of the character $z \longmapsto \left(\frac{z}{|z|}\right)^{k-1}$ of $W_{\C}$. When $k \geq 1$, the same reference tells us that $\sigma_{f,\infty}$ is exactly $\mathbf{1}\oplus\mathfrak{s}$, where $\mathfrak{s}$ denotes the morphism $\mathfrak{s}\circ \mathfrak{a}_{\R}: W_{\R} \rar \{\pm 1\}$. 

When $k=1$, the $V_{f,\lambda}^{\ast}$ all have finite image and come from the same odd representation $G_{\Q} \rar \GL{F}$. In particular, its representation to $G_{\R}$ is exactly $1 \oplus \epsilon$, with $\epsilon$ the unique nontrivial character of $G_{\R}$, and thus its restriction to $W_{\R}$ is the correct one. 

Now, we only discuss finite places. First, assume that $k \geq 2$ is even (resp. $k \geq 2$ is odd). Carayol's theorem \cite[Theorem A]{localglobal} states that there is a finite extension $E'$ of $E_{\lambda}$ and a continuous two-dimensional representation $\sigma': G_{\Q} \rar \GL{E'}$ such that, for any finite place $v \neq \ell$ (where $\ell$ is the residue characteristic of $\lambda$) of $\Q$, $\Delp{v}{\sigma'}$ is associated to $\pi_{f,v}$ (resp. $\pi_{f,v} \otimes \|\det\|^{-1/2}$, see the normalization in \cite[\S 0.2]{localglobal}) under the so-called \emph{Hecke correspondence} -- a different normalization of the Langlands correspondence. Let $w \in \{0,1\}$ be an integer with the same parity as $k$.

Fix an arbitrary $E$-embedding $E' \rar \C$. Let $q \nmid N\ell$ be a prime, $I_q \leq G_{\Q}$ be an inertia group at $q$ and $\Fr_q$ an arithmetic Frobenius at $q$. Then, for any $s \in \Z$, 
\begin{align*}
\det(1-q^{-s}\Fr_q \mid (\sigma')^{I_q}) &= L(\|\cdot\|^s\cdot \Delp{q}{(\sigma')^{\ast}})^{-1}=L\left(\pi_{f,q}\otimes\|\det\|^{-w/2},s+\frac{1}{2}\right)^{-1}\\
&= \left(1-\frac{a_q}{q^{k/2}}q^{-s+w/2}+\chi(q)q^{-2s-1+w}\right) \\
&= \det(1-q^{-s}\Fr_q \mid V_{f,\lambda}(-(k-w)/2)),
\end{align*}
so that $\sigma'$ is unramified outside $\ell N$ and, by Cebotarev, $\sigma' \simeq V_{f,\lambda}(-(k-w)/2) \otimes_{E_{\lambda}} E'$. This implies in particular that, whenever $p \neq \ell$, \[\sigma_{f,p} \simeq \|\cdot \|^{(w-1)/2} \otimes \Delp{p}{(\sigma')^{\ast}} \otimes_{E'} \C \simeq \left[\Delp{p}{V_{f,\lambda}^{\ast}} \otimes_{E_{\lambda}} \C\right]\left(\frac{(k-w)+(w-1)}{2}\right),\] whence the conclusion (where $V(z)$ denotes the twist of $V$ by $\|\cdot\|^z$).

Regardless of the value of $k$, let $\chi$ be a Dirichlet character and $p$ is a prime number distinct from the residue characteristic of $\lambda$. Write $\omega$ for the cyclotomic character. Then \[D_{f \otimes \chi,p} = \Delp{p}{V_{f \otimes \chi,\lambda}^{\ast}} = \Delp{p}{V_{f,\lambda}^{\ast} \otimes \overline{\chi}(\omega)} = \Delp{p}{V_{f,\lambda}^{\ast}} \otimes \Delp{p}{\overline{\chi}(\omega)} = D_{f,p} \otimes \chi_{\mathbb{A}}(\mathfrak{a}_{\Q_p}).\] 

When $k=1$, while the result is certainly known, it is not covered by Carayol's theorem. Since we will need it, we give an elementary proof. Every $V_{f,\lambda}$ has finite image and is realizable over $E$ by \cite[Remarque 6.5 and footnote]{Del-Ser}, so Cebotarev's theorem shows that these representations are all isomorphic to some $\rho: G_{\Q} \rar \GL{E} \subset \GL{\C}$. We show that for all primes $p$, one has $\varepsilon_p(f)=\varepsilon(\Delp{p}{\rho^{\ast}},\frac{1}{2},\theta_p)$: using the compatibility with twists that we have just proved, the equality of conductors between $\rho^{\ast}$ and $f$ \cite[Theorem 4.6]{Del-Ser} and the converse theorem \cite[\S 27.1]{BH}, we are done. 

Let $p \mid N$ be a prime. Let $S$ be the set of primes $q \mid N$ distinct from $p$. Let $\chi$ be the character of $f$, and, for each $q \mid N$, $\chi=\chi_q\chi'_q$, where the conductor of $\chi_q$ is a power of $q$ and the conductor of $\chi'_q$ is coprime to $q$. 

For each $q \in S'$, choose an odd integer $n_q > v_q(N)+2$ such that the property of \cite[\S 29.4 Proposition (4)]{BH} is satisfied for $\Delp{q}{\rho^{\ast}}$ and the additive character $\theta_q$. Let now $\psi$ be an even Dirichlet character of conductor $M=\prod_{q \in S'}{q^{n_q}}$, and write as above $\psi=\prod_q{\psi_q}$, $\psi'_q=\psi\overline{\psi_q}$, where the conductor of $\psi_q$ is $q^{n_q}$, and let us write $\sigma=\rho \otimes \chi$. 

If $\omega: G_{\Q} \rar \hat{\Z}^{\times}$ is the cyclotomic character, one can check directly that, for any primitive Dirichlet character $\varphi$, and any prime $q$, $\overline{\varphi}(\omega)_{|W_{\Q_q}} = \varphi_{\mathbb{A}}(\mathfrak{a}_{\Q_q})$. 

Let $q \in S'$ and $Q=q^{n_q}$. Then, for a certain $c_q \in \Q_q^{\times}$,
\begin{align*}
\varepsilon(\Delp{q}{\sigma^{\ast}},\frac{1}{2},\theta_q) &= \varepsilon(\psi_{\mathbb{A}}(\mathfrak{a}_{\Q_q}) \otimes \Delp{p}{\rho^{\ast}},\frac{1}{2},\theta_q)\\
&= \det{\Delp{q}{\rho^{\ast}}}(\mathfrak{a}_{\Q_q}^{-1}(c_q))^{-1}\varepsilon((\psi_{\mathbb{A}})_{|\Q_q^{\times}},\frac{1}{2},\theta_q)^2,
\end{align*}
where, for all $x \in 1+q^{(n_q+1)/2}\Z_q$, $\overline{\psi_q}(1+x)=(\psi_{\mathbb{A}})_{|\Q_q^{\times}}(1+x)=\theta_p(c_qx)$. In particular, $c_q=\frac{z_q}{Q}$ for some $z_q \in \Z_q^{\times}$. Since $\det{\rho}=\overline{\chi}_{\mathbb{A}}(\mathfrak{a}_{\Q_q})$, one has 
\[\det{\Delp{q}{\rho^{\ast}}}(\mathfrak{a}_{\Q_q}^{-1}(c_q))^{-1} = (\overline{\chi}_{\mathbb{A}})_{|\Q_q^{\times}}(\frac{z_q}{Q})=\chi_q(z_q)\overline{\chi'_q(Q)}.\]

Finally, by \cite[(3.6.3)]{NTB}
\[\varepsilon((\psi_{\mathbb{A}})_{|\Q_q^{\times}},\frac{1}{2},\theta_q) = (\psi_{\mathbb{A}})_{|\Q_q^{\times}}(Q)\frac{\mathfrak{g}(\psi_q)}{|\mathfrak{g}(\psi_q)|}= \psi'_q(Q)\frac{\mathfrak{g}(\psi_q)}{|\mathfrak{g}(\psi_q)|},\]
whence 
\[\varepsilon(\Delp{q}{\sigma^{\ast}},\frac{1}{2},\theta_q) = \chi_q(z_q)\overline{\chi'_Q}(Q)\psi'_q(Q)^2\frac{\mathfrak{g}(\psi_q)^2}{Q}.\]

Now, by Proposition \ref{local-constant-bigtwist} and the proof of Lemma \ref{gauss-deligne},
\[\varepsilon_q(f \otimes \psi_q) = \overline{\chi'_q}(Q)\psi_q(-1)\frac{\mathfrak{g}(\chi_q\psi_q)}{\mathfrak{g}(\overline{\psi_q})}=\overline{\chi'_q}(Q)\psi_q(-1)\chi_q(-(-z_q))\frac{\mathfrak{g}(\psi_q)}{\mathfrak{g}(\overline{\psi_q})},\]
thus by Proposition \ref{local-constant-coprime-twist} and Lemma \ref{gauss-module}, 
\begin{align*}
\varepsilon_q(f \otimes \psi)&=\psi'_q(Q)^2\varepsilon_q(f \otimes \psi_q)=\psi'_q(Q)^2\overline{\chi'_q}(Q)\psi_q(-1)\chi_q(z_q)\frac{\mathfrak{g}(\psi_q)}{\mathfrak{g}(\overline{\psi_q})}\\
&= \psi'_q(Q)^2\overline{\chi'_q}(Q)\chi_q(z_q)\frac{\mathfrak{g}(\psi_q)}{\overline{\mathfrak{g}(\psi_q)}}=\varepsilon(\Delp{q}{\sigma^{\ast}},\frac{1}{2},\theta_q).
\end{align*}

By Proposition \ref{functional-eqn}, the constant of $f \otimes \psi$ at infinity is $-i$. Since $\rho$ is an odd Artin representation, by \cite[\S 3.2]{NTB}, one has 
\[\varepsilon_{\infty}(\sigma^{\ast}_{\infty},\frac{1}{2},\theta_{\infty})=\varepsilon_{\R}(\mathbf{1},\frac{1}{2},\theta_{\infty})\varepsilon_{\R}(\mathfrak{s},\frac{1}{2},\theta_{\infty})=\mathbf{1}(-1)\cdot 1 \cdot \mathfrak{s}(-1) \cdot i=-i,\]

Since the global $L$-functions of $\sigma^{\ast}$ (in the convention of \cite{NTB}) and $f \otimes \psi$ are equal by \cite[Th\'eor\`eme 4.6]{Del-Ser}, we must have $\varepsilon_p(f \otimes \psi)=\varepsilon(\Delp{p}{\sigma^{\ast}},\frac{1}{2},\theta_p)$. 

By Proposition \ref{local-constant-coprime-twist}, $\varepsilon_p(f \otimes \psi)=\psi(p^{v_p(N)})\varepsilon_p(f)$. On the other hand, by \cite[(3.6.5)]{NTB},
\[\varepsilon(\Delp{p}{\sigma^{\ast}},\frac{1}{2},\theta_p) = \varepsilon(\Delp{\rho^{\ast}} \otimes \psi_{\mathbb{A}}(\mathfrak{a}_{\Q_p}),\frac{1}{2},\theta_p) = \psi(p^{v_p(N)})\varepsilon(\Delp{\rho^{\ast}},\frac{1}{2},\theta_p), \]
whence the desired equality. }

\rem{In the case $k=1$, our argument is essentially a mimic of Deligne's proof of the existence of local constants (as exposed in for instance \cite[Chapter 7]{BH}).}

\cor[local-epsilon-modular-galois]{Let $f \in \mathcal{S}_k(\Gamma_1(N))$ be a newform with primitive character $\chi$, $p \mid N$ be a prime and $\sigma \in \mrm{Gal}(\Qbar/\Q)$ whose image under the cyclotomic character is $c=(c_q)_q \in \hat{\Z}^{\times}$. Write $\chi=\chi_p\chi'_p$, where the conductor of $\chi_p$ is a power of $p$ and that of $\chi'_p$ is prime to $p$. Then $D_{\sigma(f),p} \simeq D_{f,p} \otimes_{\C,\sigma} \C$ (for any extension of $\sigma$ as an automorphism of $\C$) and one has
\[ \varepsilon_p(\sigma(f))= \sigma(\varepsilon_p(f))\sigma(\chi_p(c_p))\frac{\sigma(p^{kv_p(N)/2})}{p^{kv_p(N)/2}}.\]
}

\demo{We extend $\sigma$ to an automorphism of $\C$. Let $F/\Q$ be a finite Galois extension containing the coefficients of $f$ with $F \subset \C$. For any maximal ideal $\mathfrak{l} \subset \OO_F$ with residue characteristic $\ell \neq p$, for any $F$-embedding $F_{\mathfrak{l}} \rar \C$,
\begin{align*}
D_{\sigma(f),p} \otimes_{F_{\mathfrak{l}}} \C &\simeq \Delp{p}{V_{\sigma(f),\mathfrak{l}}^{\ast}} \otimes_{F_{\mathfrak{l}}} \C \simeq \Delp{p}{V_{f,\sigma^{-1}(\mathfrak{l})}^{\ast} \otimes_{F_{\sigma^{-1}(\mathfrak{l})},\sigma} F_{\mathfrak{l}}} \otimes F_{\mathfrak{l}} \C\\
&\simeq D_{f,p} \otimes_{\C,\sigma} \C,
\end{align*}
hence $D_{\sigma(f),p} \simeq \sigma(D_{f,p})$. 

Let $\eta$ be the unramified character of $W_{\Q_p}$ mapping a Frobenius to $\frac{\sigma(\sqrt{p})}{\sqrt{p}}$. By the discussion around \cite[(3.6.10)]{NTB} and in its notation, one has 
\begin{align*}
\varepsilon_p(\sigma(f)) &= \varepsilon(D_{\sigma(f),p},\frac{k}{2},\theta_p) = \varepsilon_1(\sigma(D_{f,p})\|\cdot\|^{k/2},\theta_p)\\
&=\sigma(\varepsilon_1(D_{f,p}\|\cdot\|^{k/2}\eta^{kv_p(N)},\sigma^{-1}\circ\theta_p))=\sigma(\varepsilon_1(D_{f,p}\|\cdot\|^{k/2},\eta^{k},\theta_p(c_p^{-1}\cdot)))\\
&= \sigma(\varepsilon(D_{f,p}\otimes \eta^{k},\frac{k}{2},\theta_p(c_p^{-1}\cdot))) = \sigma\left[\det{D_{f,p}}(\mathfrak{a}_{\Q_p}^{-1}(c_p^{-1}))\varepsilon(D_{f,p} \otimes \eta^{k},\frac{k}{2},\theta_p)\right]\\
&= \sigma\left[\det{D_{f,p}}(\mathfrak{a}_{\Q_p}^{-1}(c_p^{-1}))\left(\frac{\sigma(\sqrt{p})}{\sqrt{p}}\right)^{kv_p(N)}\varepsilon(D_{f,p},\frac{k}{2},\theta_p)\right]\\
&= \sigma(\varepsilon_p(f))\frac{\sigma(p^{kv_p(N)/2})}{p^{kv_p(N)/2}} \sigma(\det{D_{f,p}}(\mathfrak{a}_{\Q_p}^{-1}(c_p^{-1}))).
\end{align*}

Let $\omega_N$ denote the mod $N$ cyclotomic character. One has 
\[\det{D_{f,p}}=\Delp{p}{\det{V_{f,\mathfrak{l}}^{\ast}}} = \|\cdot\|^{1-k}\Delp{p}{\chi^{-1}(\omega_N)}=\|\cdot\|^{1-k} \chi_{\mathbb{A}} \circ \mathfrak{a}_{\Q_p},\] so that $\det{D_{f,p}}(\mathfrak{a}_{\Q_p}^{-1}(c_p^{-1}))=\chi_{\mathbb{A}}(c_p^{-1})=\overline{\chi_p}(c_p^{-1})=\chi_p(c_p)$. 
} 

\bigskip

\section{Preparation to the computation of the root numbers}
\label{preparatory-lemmas}

This section is dedicated to various lemmas about the representations involved in Section \ref{root-numbers-surjective}. They will help in the computation, so we factor them out for clarity.  

From now on, $\rho: G_{\Q} \rar \GL{\F_p}$ is a continuous representation with $\det{\rho}=\omega_p^{-1}$. Let $F$ be a number field as in Section \ref{tate-module-connected-xrho}, that is: it contains the $p(p-1)^2(p+1)$-th roots of unity, the coefficients of all the newforms in $\mathscr{S} \cup \mathscr{P} \cup \mathscr{C}$, and, for each prime $\ell$ such that $\ell \in \F_p^{\times 2}$ and every $(f,\mathbf{1}) \in \mathscr{C}_M$, $F$ contains the roots of $X^2-a_{\ell}(f)X+\ell$. We say that such a number field is \emph{large enough}. 

\lem[everything-self-dual]{For any maximal ideal $\lambda \subset \OO_F$,
 \begin{itemize}[noitemsep,label=$-$]
\item for any $(f,\mathbf{1}) \in \mathscr{S}$, $S:= V_{f,\lambda} \otimes \mrm{St}(\rho)$ satisfies $S^{\ast}(1) \simeq S$
\item for any $(f,\chi) \in \mathscr{P}$, $P := V_{f,\lambda} \otimes \pi(1,\chi)(\rho)$ satisfies $P^{\ast}(1) \simeq P$.
\item for any $(f,\chi) \in \mathscr{C}$, $C := V_{f,\lambda} \otimes C_f(\rho)$ satisfies $C^{\ast}(1)\simeq C$. 
\item for any $(f,\mathbf{1}) \in \mathscr{C}_M$, any $\epsilon \in \{\pm \}$, any of the two characters $\psi: G_{\Q(\sqrt{-p})} \rar F_{\mathfrak{l}}$ attached to $f$, $I := \mrm{Ind}_{G_{\Q(\sqrt{-p})}}^{G_{\Q}}{\psi \otimes C^{\epsilon}(\rho_{|G_{\Q(\sqrt{-p})}})}$ satisifes $I^{\ast}(1) \simeq I$. 
\end{itemize}}

\demo{Using the results of Appendix \ref{reps-gl2} (Corollary \ref{steinberg-dual} for $\mathscr{S}$, Lemma \ref{principal-series-dual} and Proposition \ref{twist-principal-series} for $\mathscr{P}$ and Corollaries \ref{cuspidal-twist}, \ref{cuspidal-duals} for $\mathscr{C}$), one has:
\begin{align*}
S^{\ast}(1) &\simeq V_{f,\lambda}^{\ast}(1) \otimes \mrm{St}^{\ast}(\rho) \simeq V_{f,\lambda} \otimes \mrm{St}(\rho), && (f,\mathbf{1}) \in \mathscr{S}\\
P^{\ast}(1) &\simeq V_{f,\lambda}^{\ast}(1) \otimes \pi(1,\chi)^{\ast}(\rho) \simeq V_{f,\lambda} \otimes \overline{\chi}(\omega_p) \otimes \pi(1,\chi^{-1})(\rho)&&\\
&\simeq  V_{f,\lambda} \otimes \chi(\det{\rho}) \otimes \pi(1,\chi^{-1})(\rho) \simeq V_{f,\lambda} \otimes \pi(\chi,1)(\rho), && (f,\chi) \in \mathscr{P}\\
C^{\ast}(1) &\simeq V_{f,\lambda}^{\ast}(1) \otimes C_f^{\ast}(\rho) \simeq V_{f,\lambda} \otimes \chi(\det{\rho}) \otimes C_{\overline{f}}(\rho) &&\\
&\simeq V_{f,\lambda} \otimes C_{\overline{f} \otimes \chi}(\rho)\simeq C, && (f,\chi) \in \mathscr{C}.
\end{align*}

Let now $(f,\mathbf{1}) \in \mathscr{C}_M$ and $\psi,\psi': G_{\Q(\sqrt{-p})} \rar F_{\mathfrak{l}}^{\times}$ be the two characters attached to $f$: thus $\psi\psi'=\det{V_{f,\lambda}}$ is the cyclotomic character. Thus 
\[(\psi \otimes C^{\epsilon}(\rho_{|G_{\Q(\sqrt{-p})}}))^{\ast}(1) \simeq \psi^{-1} \otimes \psi \psi' \otimes (C^{\epsilon})^{\ast}(\rho_{|G_{\Q(\sqrt{-p})}})\simeq \psi' \otimes C^{-\epsilon}(\rho_{|G_{\Q(\sqrt{-p})}}).\]
Since $G_{\Q(\sqrt{-p})}=\rho^{-1}(\F_p^{\times}\SL{\F_p})$, $\psi \otimes C^{\epsilon}(\rho_{|G_{\Q(\sqrt{-p})}}))^{\ast}(1)$ is by Proposition \ref{decomp-cusp-sl} the conjugate under $G_{\Q}/G_{\Q(\sqrt{-p})}$ of $\psi \otimes C^{\epsilon}(\rho_{|G_{\Q(\sqrt{-p})}})$, whence the conclusion.  
}

\medskip

\lem[local-global-steinberg]{Let $f \in \mathcal{S}_k(\Gamma_1(N))$ be a normalized newform with character $\chi$, and $p \mid N$ be a prime number not dividing $N/p$ or the conductor of $\chi$. Then $D_{f,p}$ is the couple $(\rho,N)$, where $\rho: W_{\Q_p} \rar \GL{\C}$ is the unramified representation such that the geometric Frobenius acts by $a_p(f)\begin{pmatrix}1 & 0\\0 &p\end{pmatrix}$ and $N$ acts by $\begin{pmatrix}0 & 1\\0 & 0\end{pmatrix}$.   }

\demo{By Proposition \ref{langlands-elem-2}, $\pi_{f,p}$ is the twist by the unramified character $p \longmapsto p^{1-k/2}a_p(f)$ of the Steinberg representation of $\GL{\Q_p}$. By \cite[\S 33.3]{BH}, $\sigma_{f,p}$ is the couple $(\sigma,N)$ such that $\sigma$ is the sum of the unramified characters $\phi \longmapsto \alpha p^{-1/2}$, $\phi \longmapsto \alpha p^{1/2}$, with $\alpha=p^{1-k/2}a_p(f)$. Since $\phi N \phi^{-1}=p^{-1}N$, after conjugating if necessary we may assume that \[\sigma(\phi)=\begin{pmatrix}\alpha p^{-1/2} & 0\\0 & \alpha p^{1/2}\end{pmatrix},\quad N=\begin{pmatrix}0 & 1\\0 & 0\end{pmatrix}.\] By Proposition \ref{local-global-modular-forms}, $D_{f,p}$ is the twist of $\sigma_{f,p}$ by the unramified character $\|\cdot\|^{(1-k)/2}$, which maps $\phi$ to $p^{(k-1)/2}$. 
}

\medskip

\lem[local-global-p]{Let $(f,\chi) \in \mathscr{P}$. Then $D_{f,p}$ is the direct sum of two characters of $W_{\Q_p}$.}

\demo{By Proposition \ref{langlands-elem-2}, $\pi_{f,p}$ is an irreducible principal series representation, so we just need to apply \cite[\S 33.3]{BH}. 
}

\medskip

\lem[local-global-c]{Let $(f,\chi) \in \mathscr{C}$. Let $E/\Q_p$ be the unique quadratic unramified extension. Let $\phi: E^{\times}/(1+p\OO_E) \rar \C^{\times}$ be the character such that $\{\phi_{|\OO_E^{\times}},\phi^p_{|\OO_E^{\times}}\}$ is attached to the cuspidal representation $C_f$ and $\phi(p)=-p$. Then $D_{f,p} \simeq \left(\mrm{Ind}_{W_E}^{W_{\Q_p}}{\overline{\phi} \circ \mathfrak{a}_E},0\right)$.  }

\demo{We saw that $\pi_{f,p}$ was the representation $\pi_{\alpha}$ for some character $\alpha:E^{\times}/(1+p\OO_E) \rar \C^{\times}$ such that $\overline{\alpha}_{|\OO_E^{\times}}$ was attached to $C_f$ and $\alpha(p)=1$. By the ``tame'' Langlands correspondence of \cite[\S 34.4]{BH}, $\sigma_{f,p} \simeq \left(\mrm{Ind}_{W_E}^{W_{\Q_p}}{(\alpha\eta) \circ \mathfrak{a}_E}, 0\right)$, with $\eta: E^{\times}/\OO_E^{\times} \rar \C^{\times}$ nontrivial and quadratic. Since 
\[(\|\cdot\|_{W_{\Q_p}})_{|W_E}=\|\cdot\|_{\Q_p} \circ \mathfrak{a}_{\Q_p} \circ (W_E \subset W_{\Q_p}) = \|\cdot\|_{\Q_p} \circ N_{E/\Q_p} \circ \mathfrak{a}_E = \|\cdot\|_E \circ \mathfrak{a}_E=\|\cdot\|_{|W_E},\]

one has \[D_{f,p} \simeq \sigma_{f,p} \otimes \|\cdot\|^{-1/2} \simeq (\mrm{Ind}_{W_E}^{W_{\Q_p}}{(\alpha\eta\|\cdot \|_E^{-1/2}) \circ \mathfrak{a}_E},0),\] with $\alpha(p)\eta(p)\|p\|_E^{-1/2}=1 \cdot (-1) \cdot (p^{-2})^{-1/2}=-p$, whence the conclusion.  
}

\bigskip

\defi{Given a newform $f \in \mathcal{S}_k(\Gamma_1(N))$, we let $D_{f,\infty}$ (or $\Delp{\infty}{V_{f,\lambda}^{\ast}}$, for any maximal ideal $\lambda$ of any $\OO_F$, where $F \subset \C$ is a number field containing the coefficients of $f$) be the complex representation of $W_{\R}$ given by $\mrm{Ind}_{W_{\C}}^{W_{\R}}{z^{1-k}}$. }

\prop[constants-with-modular-forms-elem]{Let $f \in \mathcal{S}_k(\Gamma_1(N))$ be a newform with primitive character $\chi$. Let $F$ be a number field containing the coefficients of $f$ and $\lambda \subset \OO_F$ a maximal ideal.
\begin{enumerate}[noitemsep,label=(\roman*)]
\item \label{cmfe-1} Let $\sigma: G_{\Q} \rar \GL{n}{\C}$ be an Artin representation and $c$ be the complex conjugation. Then one has \[\varepsilon(\Delp{\infty}{\sigma},\frac{1}{2},\theta_{\infty})\varepsilon(\Delp{\infty}{\sigma^{\ast}},\frac{1}{2},\theta)=\det{\sigma}(c),\quad \varepsilon(\Delp{\infty}{(V_{f,\lambda} \otimes \sigma)^{\ast}},\frac{k}{2},\theta_{\infty})=(-i)^{kn}.\]
\item \label{cmfe-2} Let $p$ be a prime number prime to $\lambda$ such that $D_{f,p}$ is the direct sum of two characters $\alpha,\beta$ of $\Q_p^{\times}$, and $\sigma: G_{\Q_p} \rar \GLn{n}{\C}$ be a continuous representation such that $\sigma^{\ast} \simeq \sigma \otimes \chi$. Then $\alpha(-1)^n=\beta(-1)^n$ and \[\varepsilon(\Delp{p}{(V_{f,\lambda} \otimes \sigma)^{\ast}},\frac{k}{2},\theta_p)=\alpha(-1)^n\varepsilon(\Delp{p}{\sigma},\frac{1}{2},\theta_p)\varepsilon(\Delp{p}{\sigma^{\ast}},\frac{1}{2},\theta_p).\]
\item \label{cmfe-3} Let $p \mid N$ be a prime number prime to $\lambda$ and the conductor of $\chi$, and such that $p^2 \nmid N$. Let $\sigma: G_{\Q_p} \rar \GLn{n}{\C}$ be a continuous representation such that $\sigma^{\ast} \simeq \sigma \otimes \chi$. Let $I_p \leq G_{\Q_p}$ be the inertia group and $\Fr_p$ be an arithmetic Frobenius, then
\[\varepsilon(\Delp{p}{(V_{f,\lambda} \otimes \sigma)^{\ast}},\frac{k}{2},\theta_p)=\varepsilon(\Delp{p}{\sigma},\frac{1}{2},\theta_p)\varepsilon(\Delp{p}{\sigma^{\ast}},\frac{1}{2},\theta_p)\det(-p^{1-k/2}a_p(f)\Fr_p \mid \sigma^{I_p}).\]
\end{enumerate}}

\demo{\ref{cmfe-1}: The first equality is the local functional equation. For the second equality, $\sigma$ is an Artin representation, so the representation $\Delp{\infty}{\sigma}$ of $W_{\R}$ factors through the quotient $W_{\R} \rar G_{\R}$ and is the sum of the characters $\mathbf{1}$ and $\mathfrak{s}$. Since $\varepsilon(D_{f,\infty},\frac{k}{2},\theta_{\infty})=(-i)^k$, all we need to do is prove that $D_{f,\infty} \otimes \mathfrak{s} \simeq D_{f,\infty}$. This is true because $D_{f,\infty}$ is induced from the subgroup $W_{\C}$ of index two and $\mathfrak{s}$ is a character of $W_{\R}/W_{\C}$. 

Now, we consider the situation of \ref{cmfe-2} let $p$ be a prime number coprime to $\lambda$ and assume that $D_{f,p}$ is the direct sum of the characters $\alpha,\beta$ of $\Q_p^{\times}$, and that $\sigma: G_{\Q_p} \rar \GLn{n}{\C}$ is a continuous representation such that $\sigma^{\ast} \otimes \sigma \otimes \chi$. In particular, after taking $\det{\Delp{p}{-}}$ and evaluating at $\mathfrak{a}_{\Q_p}^{-1}(-1)$, we see that when $n$ is odd, $(\overline{\chi}_{\mathbb{A}})_{|\Q_p^{\times}}(-1)=1$.  

One has $(\alpha\beta) \circ \mathfrak{a}_{\Q_p} = \det{D_{f,p}} = (\|\cdot\|^{k-1}\overline{\chi}_{\mathbb{A}} \circ \mathfrak{a}_{\Q_p})^{-1} = \|\cdot\|^{1-k}\chi_{\mathbb{A}}(\mathfrak{a}_{\Q_p})$, so $\alpha(-1)^n=\beta(-1)^n$. Write $\alpha=\alpha_0\|\cdot\|^{(1-k)/2}, \beta = \beta_0\|\cdot\|^{(1-k)/2}$, $\beta_0(\overline{\chi}_{\mathbb{A}})_{\Q_p^{\times}}=\alpha_0^{-1}$, then one has
\begin{align*}
\varepsilon(\Delp{p}{(V_{f,\lambda} \otimes \sigma)^{\ast}},\frac{k}{2},\theta_p) &= \varepsilon(\Delp{p}{\sigma^{\ast}} \otimes \alpha_0(\mathfrak{a}_{\Q_p}),\frac{1}{2},\theta_p)\varepsilon(\Delp{p}{\sigma^{\ast}} \otimes \beta_0(\mathfrak{a}_{\Q_p}),\frac{1}{2},\theta_p)\\
&= \varepsilon(\Delp{p}{\sigma^{\ast}} \otimes \alpha_0(\mathfrak{a}_{\Q_p}),\frac{1}{2},\theta_p)\varepsilon(\Delp{p}{\sigma} \otimes ((\overline{\chi}_{\mathbb{A}})\beta_0)(\mathfrak{a}_{\Q_p}),\frac{1}{2},\theta_p)\\
&= \varepsilon(\Delp{p}{\sigma^{\ast}} \otimes \alpha_0(\mathfrak{a}_{\Q_p}),\frac{1}{2},\theta_p)\varepsilon(\Delp{p}{\sigma} \otimes \alpha_0^{-1}(\mathfrak{a}_{\Q_p}),\frac{1}{2},\theta_p)\\
&= \det{\sigma^{\ast} \otimes \alpha_0(\mathfrak{a}_{\Q_p})}(\mathfrak{a}_{\Q_p}^{-1}(-1))=\alpha_0(-1)^n\det{\sigma}(\mathfrak{a}_{\Q_p}^{-1}(-1))\\
&= \alpha(-1)^n\varepsilon(\Delp{p}{\sigma},\frac{1}{2},\theta_p)\varepsilon(\Delp{p}{\sigma^{\ast}},\frac{1}{2},\theta_p).
\end{align*}

Now we prove \ref{cmfe-3}. Assume instead that $p$ does not divide the conductor of $\chi$ and $v_p(N)=1$. Let $\alpha: W_{\Q_p} \rar \C^{\times}$ be the unramified character mapping a geometric Frobenius $\phi$ to $a_p(f)$. Then $D_{f,p}=(\rho,N)$, with $N=\begin{pmatrix}0 & 1\\0 & 0\end{pmatrix}$ and $\rho=\begin{pmatrix}\alpha & 0\\0 & \alpha\|\cdot\|^{-1}\end{pmatrix}$. By \cite[\S 31.3]{BH}\footnote{This formula yields the same result as \cite[(4.1.5)]{NTB}, because for a Deligne representation $(\rho,N)$ on a $\C$-vector space $V$, $(z,f) \in ((V/\ker{N}) \otimes \|\cdot\|^{-1})  \times V^{\ast}/\ker{N^{\ast}} \longmapsto f(Nz) \in \C$ is a perfect duality of representations of $W_F$. }, 
\begin{align*}
&\varepsilon(\Delp{p}{(V_{f,\lambda} \otimes \sigma)^{\ast}},\frac{k}{2},\theta_p) = \varepsilon(D_{f,p} \otimes \Delp{p}{\sigma^{\ast}},\frac{k}{2},\theta_p)\\
&= \varepsilon((\alpha \oplus \alpha \|\cdot\|^{-1})\otimes \Delp{p}{\sigma^{\ast}},\frac{k}{2},\theta_p)\\
& \quad\quad \times\frac{L(\Delp{p}{\sigma} \otimes \alpha^{-1},\frac{2-k}{2})L(\Delp{p}{\sigma} \otimes \alpha^{-1}\|\cdot\|,\frac{2-k}{2})L(\Delp{p}{\sigma^{\ast}} \otimes \alpha,\frac{k}{2})}{L(\Delp{p}{\sigma^{\ast}} \otimes \alpha,\frac{k}{2})L(\Delp{p}{\sigma^{\ast}} \otimes \alpha\|\cdot\|^{-1},\frac{k}{2})L(\Delp{p}{\sigma}\otimes \alpha^{-1}\|\cdot\|,\frac{2-k}{2})}\\
&= \varepsilon((\alpha \oplus \alpha \|\cdot\|^{-1})\otimes \Delp{p}{\sigma^{\ast}},\frac{k}{2},\theta_p)\frac{L(\Delp{p}{\sigma} \otimes \alpha^{-1},1-k/2}{L(\Delp{p}{\sigma^{\ast}} \otimes \alpha\|\cdot\|^{-1},k/2)}.
\end{align*}

Because $\sigma$ has finite image, if $I \leq W_{\Q_p}$ denotes the inertia subgroup,
\begin{align*}
& \frac{L(\Delp{p}{\sigma} \otimes \alpha^{-1},1-\frac{k}{2})}{L(\Delp{p}{\sigma^{\ast}} \otimes \alpha\|\cdot\|^{-1},\frac{k}{2})} = \frac{\det(1-\sigma^{\ast}(\phi)\alpha(\phi)\|\phi\|^{-1}p^{-k/2}\mid (\sigma^{\ast})^I)}{\det(1-\sigma(\phi)\alpha^{-1}(\phi)p^{k/2-1}\mid \sigma^I)}\\
&= \frac{\det(1-\sigma(\phi^{-1})\alpha(\phi)p^{1-k/2}\mid \sigma^I)}{\det(1-\sigma(\phi)\alpha^{-1}(\phi)p^{k/2-1} \mid \sigma^I)}= \det(-p^{1-k/2}\sigma(\phi^{-1})\alpha(\phi)\mid \sigma^I)\\
&= \det(-p^{1-k/2}\sigma(\Fr_p)a_p(f)\mid \sigma^I).
\end{align*}

On the other hand, the properties of local constants imply that, for any two unramified characters $\beta,\gamma$ of $W_{\Q_p}$, $\varepsilon((\beta \oplus \gamma)\otimes \Delp{p}{\sigma^{\ast}},\frac{k}{2},\theta_p)$ only depends on $\beta\gamma$. Therefore, since $\alpha^2\|\cdot\|^{-1}$ is the unramified character $\chi_{\mathbb{A}}(\mathfrak{a}_{\Q_p}) \|\cdot \|^{1-k}: \phi \longmapsto p^{k-1}\chi(p)$,

\begin{align*}
&\varepsilon(\Delp{p}{(V_{f,\lambda} \otimes \sigma)^{\ast}},\frac{k}{2},\theta_p) \\
&= \varepsilon(\Delp{p}{\sigma^{\ast}} \otimes (\|\cdot\|^{(1-k)/2} \oplus \|\cdot\|^{(1-k)/2}\chi_{\mathbb{A}}(\mathfrak{a}_{\Q_p})),\frac{k}{2},\theta_p)\det(-p^{1-k/2}\sigma(\Fr_p)a_p(f)\mid \sigma^I)\\
&=\varepsilon(\Delp{p}{\sigma^{\ast}},\frac{1}{2},\theta_p) \varepsilon(\Delp{p}{\sigma^{\ast} \otimes \overline{\chi}},\frac{1}{2},\theta_p)\det(-p^{1-k/2}\sigma(\Fr_p)a_p(f)\mid \sigma^I),  
\end{align*}
whence the conclusion.
}

\bigskip

To study the local constants at the place $p$, we need first to study the shape of representations $G_{\Q_p} \rar \GL{\F_p}$. 

\prop[trichotomy-rho]{Let $p \geq 5$ be a prime. Write $G=G_{\Q_p}$ and let $I \triangleleft G$ the inertia subgroup. Let $\rho: G \rar \GL{\F_p}$ be a continuous group homomorphism with $\det{\rho}=\omega_p^{-1}$. Then exactly one of the following statements is true:
\begin{itemize}[noitemsep,label=$-$]
\item $\rho$ is not tamely ramified. There are characters $q, \beta: G \rar \F_p^{\times}$ and a cocycle $\nu: G \rar \F_p(q)$ that does not vanish on the wild inertia subgroup such that $\rho$ is conjugate to a representation of the form $\sigma \longmapsto \beta(\sigma)\begin{pmatrix}q(\sigma)&\nu(\sigma)\\0  & 1\end{pmatrix}$. 
\item $\rho$ is conjugate to a representation of the form $\sigma \longmapsto \begin{pmatrix}\alpha(\sigma)&0\\0 & \beta(\sigma)\end{pmatrix}$ for some characters $\alpha,\beta: G \rar \F_p^{\times}$, and $\det: \rho(I) \rar \F_p^{\times}$ is an isomorphism. 
\item There exists a nonsplit Cartan subgroup $C \leq \GL{\F_p}$ such that $\rho^{-1}(C)$ is exactly the Galois group of the unramified quadratic extension of $\Q_p$. Moreover, $\rho(I)$ is a subgroup of $C$ of index $a \mid \frac{p+1}{2}$ prime to $p-1$. 
\end{itemize}
}

\demo{First, let us check that all statements are pairwise incompatible. Clearly, the first condition is not compatible with any of the other two, in which $\rho$ is tamely ramified. In the second situation, $\rho(I)$ has order $p-1$, while, in the third situation, $\rho(I)$ has order $\frac{p^2-1}{a} \geq \frac{p^2-1}{(p+1)/2}=2(p-1)$, so they are incompatible too. 

Now, we show that one of the three statements holds.

Since the wild inertia subgroup $I^+ \triangleleft G$ is a pro-$p$-group, the image of $I^+$ is either trivial (then $\rho$ is tamely ramified) or conjugate to the subgroup $U^{\Z}=\begin{pmatrix} 1 &\ast\\0 & 1\end{pmatrix}$. 

Suppose we are in the latter case. Then (after conjugating so that $\rho(I^+)=U^{\Z}$), $\rho(G)$ is contained in the normalizer of $U^{\Z}$, which is the group of upper-triangular matrices, so we can write $\rho(\sigma)=\beta(\sigma)\begin{pmatrix}q(\sigma)&\nu(\sigma)\\0  & 1\end{pmatrix}$ for some characters $q,\beta: G \rar \F_p^{\times}$, and $\nu: G \rar \F_p$. It can be directly checked that $\nu: G \rar \F_p(q)$ is a cocycle, and since $\rho(I^+)$ is nontrivial, $\nu(I^+)=\F_p$. 

So assume from now on that $\rho$ is tamely ramified. Let $s$ be a pro-generator of the pro-cyclic group $I/I^+$, which has a pro-order prime to $p$, so that $\det{\rho(s)}$ generates $\det{\rho}(I)=\omega_p^{-1}(I)=\F_p^{\times}$. Therefore, $\rho(s) \in \GL{\F_p}$ is non scalar and its order is prime to $p$, hence it is either diagonalizable in $\GL{\F_p}$, or a non-scalar element of some nonsplit Cartan subgroup $C$. 

In the former situation, after conjugating, we can assume that $\rho(s)=\begin{pmatrix}\alpha & 0\\0 & \beta\end{pmatrix}$ with $\alpha \neq \beta$, so that $\rho(x^p)=\rho(x)$ for any $x \in I$. In particular, if $z \in G$ is an arithmetic Frobenius, $\rho(z)\rho(x)\rho(z)^{-1}=\rho(zxz^{-1})=\rho(x^p)=\rho(x)$, so $\rho(G)$ is contained in the centralizer subgroup of $\rho(s)$, which consists of the diagonal matrices and is abelian. Therefore, $\rho(I)$ is a quotient of the image of $I$ in $G^{ab}$, which is $\Z_p^{\times}$. Since $\rho(I)$ is cyclic with order prime to $p$, it follows that $\rho(I)$ is a quotient of $\F_p^{\times}$. Since $\det: \rho(I) \rar \F_p^{\times}$ is surjective by the assumption on $\rho$, it is an isomorphism. 

We now assume that $\rho(s)$ is a non-scalar element of some nonsplit Cartan group $C$. Then $\rho(I)$ is a subgroup of $C$, so it is cyclic of index $a \mid p^2-1$. Since $\det{C}=\F_p^{\times}=\det{\rho(I)}$, we have $(p+1)C=(p+1)\rho(I)$, so that the greatest common divisor $(p+1)a$ and $p^2-1$ is exactly $p+1$. Thus $a$ is coprime to $p-1$: it is odd, so coprime to $2(p-1)$ and $a \mid \frac{p+1}{2}$. An explicit computation shows that the normalizer of $\frac{p+1}{2}C$ in $\GL{\F_p}$ is exactly the normalizer $N$ of $C$, so that $\rho(G) \subset N$. 
 
If $\rho(G) \subset C$, then the image of $\rho$ is commutative, so that $\rho(I)$ is a quotient of the image of $I$ in $G^{ab}$, which is $\Z_p^{\times}$. Since $\rho(I)$ is cyclic with order prime to $p$, $\rho(I)$ is a quotient of $\F_p^{\times}$, which is impossible since $|\rho(I)|=\frac{p^2-1}{a} \geq \frac{p^2-1}{(p+1)/2}=2(p-1)$. So $\rho^{-1}(C)$ is the Galois group of some quadratic unramified extension of $\Q_p$, that is, of \emph{the} quadratic unramified extension of $\Q_p$.
}

\medskip

\defi{In the situation of Proposition \ref{trichotomy-rho}, we say that $\rho$ is 
\begin{itemize}[noitemsep,label=\tiny$\bullet$]
\item \emph{wild} if it is not tamely ramified,
\item \emph{split} if it is conjugate to a representation of the form $\sigma \longmapsto \begin{pmatrix}\alpha(\sigma)&0\\0 & \beta(\sigma)\end{pmatrix}$ for some characters $\alpha,\beta: G \rar \F_p^{\times}$,
\item \emph{$a$-Cartan} if $\rho(G)$ is contained in the normalizer $N$ of a Cartan subgroup $C$, with $\rho^{-1}(C)$ the Galois group of the quadratic unramified extension of $\Q_p$ and $\rho(I)$ is a cyclic subgroup of $C$ of index $a \mid \frac{p+1}{2}$ coprime to $p-1$.  
\end{itemize}
Note that being wild, split, or $a$-Cartan only depends on the conjugacy class of $\rho$.}

\medskip

\lem[expand-wild]{Let $q,\beta: G_{\Q_p} \rar \F_p^{\times}$ be characters and $\nu: G_{\Q_p} \rar \F_p(q)$ be a cocycle whose restriction to the wild inertia subgroup does not vanish. Let $\rho(\sigma)= \beta(\sigma)\begin{pmatrix}q(\sigma)&\nu(\sigma)\\0 & 1\end{pmatrix}$ for every $\sigma \in G_{\Q_p}$ and assume that $\det{\rho}=\omega_p^{-1}$. Let $G=G_{\Q_p}$, $I \triangleleft G$ be the inertia subgroup, $I^+\triangleleft G$ the wild inertia subgroup. Then:
\begin{itemize}[noitemsep,label=\tiny$\bullet$]
\item the image of $(q\beta,\beta): I/I^+ \rar (\F_p^{\times})^2$ is cyclic with order $p-1$,
\item the extension $L/\Q_p$ defined by $G_L=\ker{q}$ is tamely ramified of degree dividing $p-1$,
\item $f_{L/\Q_p}$ is odd and $v_2(e_{L/\Q_p})=v_2(p-1)$,
\item $\nu_{|G_L}$ is a wildly ramified character of conductor exponent $c$ such that $c-1$ is coprime to $e_{L/\Q_p}$,
\item $q_{|I}$ is the $(c-1)$-th power of the tame inertia character of $\mrm{Gal}(L/\Q_p)_0$, so that $q_{|I}$ is equal to the $\frac{(c-1)(p-1)}{e}$-th power of the mod $p$ cyclotomic character, 
\item $c \leq e+1+\frac{e}{p-1}$,
\item $L$ is linearly disjoint from the quadratic unramified extension $E$ of $\Q_p$ and the Langlands constant $\lambda_{LE/\Q_p}(\theta_p)$ is $(-1)^{(p+1)/2}$, 
\item Let $M \in \GL{\F_p}$ be such that $M\rho M^{-1}=\beta'\begin{pmatrix}q'& \nu'\\0 & 1\end{pmatrix}$, then one has \[(\beta',q')=(\beta,q),\quad \nu'_{|G_L} \in \left(\frac{\det{M}}{p}\right)\F_p^{\times 2}\nu_{|G_L}.\]	 
\end{itemize}
}

\demo{$\sigma \longmapsto \begin{pmatrix} (q\beta)(\sigma)&0\\0&1\end{pmatrix}$ is a split representation, so the first claim follows from Proposition \ref{trichotomy-rho}. Let now be $L/\Q_p$ the extension such that $\ker{q}=G_L$, and $E/\Q_p$ be the unique quadratic unramified extension. Since $H := \mrm{Gal}(L/\Q_p)$ identifies with the image of $q$, it is a cyclic group of order dividing $p-1$, so that $L/\Q_p$ is tamely ramified. 

Let $a: \Q_p^{\times}/(1+p\Z_p) \rar \F_p^{\times}$ be such that $a \circ \mathfrak{a}_{\Q_p}=q$, so that $\ker{a}=N_{L/\Q_p}L^{\times}$. Since $(q\beta^2)(I)=\det{\rho}(I)=\F_p^{\times}$, $q(I) \not\subset \F_p^{\times 2}$, thus $a(\Z_p^{\times}) \not\subset \F_p^{\times 2}$, so that $\F_p^{\times}/a(\Z_p^{\times})$ has odd cardinality. Since $|a(\Z_p^{\times})|=e_{L/\Q_p}$ (which we denote $e$ for short), $v_2(e)=v_2(p-1)$. In particular, $v_2(e) \leq v_2(|H|) \leq v_2(p-1)$, so $v_2(|H|)=v_2(e)=v_2(p-1)$ and $f_{L/\Q_p}$ is odd, therefore $E$ and $L$ are linearly disjoint. 

$\nu: G_L \rar \F_p$ is a morphism, and its restriction to the wild inertia subgroup is onto. Therefore, its conductor exponent $c$ is at least $2$. Write $\nu_{|G_L}=\nu_0 \circ \mathfrak{a}_L$, for some group homomorphism $\nu_0: L^{\times} \rar \F_p$. For $x \in G, y \in G_L$, we can see that $\nu(xyx^{-1})=q(x)\nu(y)$, so, for any $\sigma \in H$ and $x \in L^{\ast}$, $\nu_0(\sigma(x))=q(\sigma)\nu_0(x)$.  

Let $K \subset L$ be the maximal subextension which is unramified over $\Q_p$: there are uniformizers $\pi \in K, \varpi \in L$, with $\varpi^e=\pi$. Let $H_0=\mrm{Gal}(L/K)$ be the inertia subgroup of $H$, and, for $s \in H_0$, $\kappa(s)=\frac{s(\varpi)}{\varpi} \pmod{\varpi\OO_L}$: $\kappa: H_0 \rar k_L^{\times}$ is an injective homomorphism, where $k_L$ is the residue field of $L$. Since $e=|H_0| \mid p-1$, $\kappa$ is in fact an injective homomorphism $\kappa: H_0 \rar \F_p^{\times}$. 

By definition of $c\geq 2$, $\nu_0: \frac{1+\varpi^{c-1}\OO_L}{1+\varpi^c\OO_L} \rar \F_p(q)$ is a surjective character. But, as a $k_L[H_0]$-module, $\frac{1+\varpi^{c-1}\OO_L}{1+\varpi^c\OO_L}$ is isomorphic to $k_L(\kappa^{c-1})$, so $\kappa^{c-1}=q$. Since $q: H_0 \rar \F_p^{\times}$ is injective, $c-1$ is coprime to $e$, and in particular $c$ is odd. 

Suppose that $e < p-1$. Then by \cite[II.5.5]{Neukirch-ANT}, there is a group isomorphism $\varpi\OO_L \rar 1+\varpi\OO_L$ mapping $\varpi^r\OO_L$ to $1+\varpi^r\OO_L$ for any $r \geq 1$. Therefore, $1+\varpi^{e+1}\OO_L=(1+\varpi\OO_L)^{p} \subset \ker{\nu}$, so $c \leq e+1$. If $e=p-1$, then there is a group isomorphism $\varpi^2\OO_L \rar 1+\varpi^2\OO_L$ mapping $\varpi^r\OO_L$ to $1+\varpi^r\OO_L$ for any $r \geq 2$. We see as above that $(1+\varpi^{p+1}\OO_L)\subset (1+\varpi^2\OO_L)^p$, so $c \leq p+1$. In either case, $c \leq e+1+\frac{e}{p-1}$.  

The fact that $\kappa=(\omega_p)_{|I}^{\frac{p-1}{e}}$ is a direct consequence of \cite[Propositions 1, 3]{Serre-image-ouverte} (because $\kappa$ is the character $\theta_e$, while $\theta_{p-1}=(\omega_p)_{|I}$ because of the normalization of the class theory isomorphism). 

All that remains to do is compute the Langlands constant $\lambda_{LE/\Q_p}(\theta_p)$. Write $p^{\ast}=(-1)^{\frac{p-1}{2}}p$. Since $H$ is cyclic, it contains at most one character of order two. Moreover, $q^{-1}(\F_p^{\times 2})/G_L$ is an index two subgroup of $H$, corresponding to a certain quadratic extension $M$ of $\Q_p$ contained in $L$. On the other hand, $q^{-1}(\F_p^{\times 2})=(\det{\rho})^{-1}(\F_p^{\times 2})=G_{\Q(\sqrt{p^{\ast}})}$, so $M=\Q_p(\sqrt{p^{\ast}})$. Thus, $H$ has a unique character with order two, which we denote by $\chi_2$; it is the character of the quadratic extension $\Q_p(\sqrt{p^{\ast}})/\Q_p$, and is also $\omega_p^{(p-1)/2}$. 

Let $\psi: H \rar \C^{\times}$ be any character (in particular, $\psi$ is at most tamely ramified) and $\eta$ be the unique quadratic unramified character of $G_{\Q_p}$: then $\varepsilon(\psi\eta,\frac{1}{2},\theta_p)$ is $\varepsilon(\psi,\frac{1}{2},\theta_p)\eta_{\psi}$, where $\eta_{\psi}$ is $1$ if $\psi$ is unramified and $-1$ otherwise. Since $H$ has exactly $[L:\Q_p]-f_{L/\Q_p}=f_{L/\Q_p}(e-1) \notin 2\Z$ ramified characters, one has

\begin{align*}
\lambda_{LE/\Q_p}(\theta_p) &= \frac{\varepsilon(\mrm{Ind}_{LE}^{\Q_p}{\mathbf{1}},\frac{1}{2},\theta_p)}{\varepsilon(\mathbf{1},\frac{1}{2},\theta_{EL})} = \prod_{\chi: H \rar \C^{\times}}{\varepsilon(\chi,\frac{1}{2},\theta_p)\varepsilon(\chi\eta,\frac{1}{2},\theta_p)}\\
& = \prod_{\chi: H \rar \C^{\times}}{\varepsilon(\chi,\frac{1}{2},\theta_p)\varepsilon(\chi,\frac{1}{2},\theta_p)\eta_{\chi}} =-\prod_{\chi: H \rar \C^{\times}}{\varepsilon(\chi,\frac{1}{2},\theta_p)\varepsilon(\chi^{-1},\frac{1}{2},\theta_p)}\\
&=-\prod_{\chi: H \rar \C^{\times}}{\chi(\mathfrak{a}_{\Q_p}^{-1}(-1))}.
\end{align*}

Since $H$ is cyclic with even order, the product of its characters is exactly $\chi_2$. Then $\chi_2(\mathfrak{a}_{\Q_p}^{-1})$ is the Teichm\"uller lift of $\omega_p^{(p-1)/2}(\mathfrak{a}_{\Q_p}^{-1}(-1))=(-1)^{(p-1)/2}$, thus $\lambda_{LE/\Q_p}(\theta_p)=(-1)^{\frac{p+1}{2}}$. 

Let now $\beta',q':G_{\Q_p} \rar \F_p^{\times}$ be characters $\nu': G_{\Q_p}\rar \F_p(q')$ be a cocycle such that, for some $M \in \GL{\F_p}$, $M\rho M^{-1}=\beta'\begin{pmatrix}q' & \nu'\\0 & 1\end{pmatrix}$. Since the wild inertia subgroup $I^+$ is a pro-$p$-group, $\beta_{|I^+}=q_{|I^+}=1$, and $M\rho M^{-1}$ is not tamely ramified (because $\rho$ isn't), so $\nu': I^+ \rar \F_p$ is a nonzero group homomorphism. In particular, we have $M\begin{pmatrix}1 & \nu_{|I^+}\\0 & 1\end{pmatrix}M^{-1}=\begin{pmatrix}1 & \nu'_{|I^+}\\0 & 1\end{pmatrix}$: a direct computation shows that $M=\begin{pmatrix}u & v\\0 & w\end{pmatrix}$ for some $u,w \in \F_p^{\times}$ and $v \in \F_p$. This implies that $(\beta q,\beta)=(\beta'q',\beta')$, hence $(\beta,q)=(\beta',q')$, and another direct computation shows that $\frac{u}{w}\nu_{|G_L}=\nu'_{|G_L}$, hence $\nu'_{|G_L} \in (\det{M})\F_p^{\times 2}\nu_{|G_L}$. 
}

\medskip

\cor[expand-wild-cvalue]{With the same notation as in Lemma \ref{expand-wild}, let $1 \leq n \leq p-1$ the smallest positive integer such that $q_{|I}=(\omega_p^n)_{|I}$. Then, unless $L=\Q_p(\mu_p)$ and $n=1$, the space of possible cocycles $\nu$ (up to coboundaries) is one-dimensional and $c-1=\frac{ne}{p-1}$. When $n=1$ and $L=\Q_p(\mu_p)$, the space of possible cocycles $\nu$ is $2$-dimensional and $c \in \{2,p+1\}$, and the subspace of cocycles $\nu$ for which $c=2$ is a line. 
Moreover, in this situation, one has $c=2$ if and only if $\begin{pmatrix}\omega_p & \nu\\0 & 1\end{pmatrix}$ is \emph{peu ramifi\'ee} in the sense of \cite[\S 2.4]{Conj-Serre}, if and only if this representation comes from a finite flat subgroup over $\Z_p$.}

\demo{By local Tate duality and the computation of the local Euler-Poincar\'e characteristic \cite[(7.2.6),(7.3.1)]{CNB}, \[\dim{H^1(\Q_p,\F_p(q))} = \dim{H^0(\Q_p,\F_p(q))}+\dim{H^2(\Q_p,\F_p(q))}+1=1+\dim{H^0(\Q_p,\F_p(\omega_pq^{-1}))},\] so $\dim{H^1(\Q_p,\F_p(q))}$ is equal to one if $q\neq\omega_p$ and two otherwise.

By Lemma \ref{expand-wild}, $n \equiv (c-1)\frac{p-1}{e}\pmod{p-1}$, and $c-1 \leq e+\frac{e}{p-1}$. If $e < p-1$, then $c-1 \leq e$, so $n,(c-1)\frac{p-1}{e}$ are congruent modulo $p-1$ and in the interval $[1,p-1]$, so they are equal. If $e=p-1$, then $c \leq p+1$ and $n \equiv c-1 \pmod{p-1}$. So unless $c=p+1$, $1 \leq n,c-1 \leq p-1$ are congruent modulo $p-1$ hence equal. So the only possibility for $n \neq (c-1)\frac{p-1}{e}$ is if $e=p-1$ and $n=1$: in this case $c \in \{2,p+1\}$. 

When $e=p-1$ and $n=1$, assume that $L \neq \Q_p(\mu_p)$. Then consider the following diagram:
\[
\begin{tikzcd}[ampersand replacement=\&]
1 \arrow{r} \& 1+\varpi^2\OO_L \arrow{d}{\alpha: x \mapsto x^p} \arrow{r}{\subset} \& 1+\varpi\OO_L \arrow{d}{\beta: x \mapsto x^p} \arrow{r}\& \F_p \arrow{r}\& 1\\
 1 \arrow{r} \& 1+\varpi^{p+1}\OO_L \arrow{r}{\subset} \& 1+\varpi^p\OO_L  \arrow{r}\& \F_p \arrow{r}\& 1
\end{tikzcd}
\]

By the snake lemma, since $\beta$ is injective (as $\mu_p \not\subset L$) and $\alpha$ is an isomorphism, $\gamma$ is injective, hence an isomorphism. Thus $\beta$ is onto, and every element in $1+\varpi^p\OO_L$ is a $p$-th power. Therefore $c \leq p$ and therefore $c=2$. When $L=\Q_p(\mu_p)$ and $n=1$, one has $q=\omega_p$. The rest of the assertions is contained in \cite[\S 2.4, 2.8]{Conj-Serre}.  }

\bigskip

\lem[expand-cartan]{Let $\rho: G_{\Q_p} \rar \GL{\F_p}$ be a continuous representation with $\det{\rho}=\omega_p^{-1}$. Assume that $\rho$ is $a$-Cartan: let $C$ be a Cartan subgroup with normalizer $N$ such that $\rho(G_{\Q_p})\subset N$ and $\rho^{-1}(C)$ is the subgroup of $G_{\Q_p}$ corresponding to the unramified quadratic extension $E/\Q_p$. There is a unique character $\sigma: E^{\times}/(1+p\OO_E) \rar C$ such that $\rho_{|G_E}=\sigma \circ \mathfrak{a}_E$. One has $\sigma(p)=-I_2$ and $\sigma(x)=x^{-1}I_2$ for each $x \in \F_p^{\times}$. Hence $\rho(G_E)$ is the image by $\rho$ of the inertia subgroup. Let $i,j: (\OO_E/p)^{\times} \rar C$ be the two natural isomorphisms that respect the trace and the norm, and, for $k \in \{i,j\}$, let $b_k$ be the index in $C$ of the subgroup
$\{s \in (\OO_E/p)^{\times} \mid \sigma(s)=k(s)^{-1}\}$. Then each $b_k$ is divisible by $a$, each $\frac{p^2-1}{b_k}$ is prime to $a$, and the least common multiple of $b_i,b_j$ is $p+1$. 
}

\demo{The existence of $\sigma$ is a consequence of local class field theory, and it is clear that $a$ is the index of $\sigma(\OO_E^{\times})$ in $C$. For any $x=\mathfrak{a}_E(s) \in \OO_E^{\times}$, 
\[\sigma(N_{E/\Q_p}x)=\sigma(x^{p+1})=\rho(s)^{p+1}=\omega_p^{-1}(s)I_2=((N_{E/\Q_p}\mathfrak{a}_E)(s))^{-1}I_2=(N_{E/\Q_p}x)^{-1}I_2,\] so, for any $y \in \F_p^{\times}$, $\sigma(y)=y^{-1}I_2$. 
Let $F \in G_{\Q_p}$ be an arithmetic Frobenius mapped under $\mathfrak{a}_{\Q_p}$ to $p^{-1}$: then $\rho(F) \in N \backslash C$, so $\rho(F)^2=-\det{\rho}(F_p)I_2=-\omega_p^{-1}(F_p)I_2=-I_2$. Since $F^2$ is the image of $F$ under the transfer map $G_{\Q_p}^{ab} \rar G_E^{ab}$, one has $\mathfrak{a}_E(F^2)=\mathfrak{a}_{\Q_p}(F)=p^{-1}$, thus $\sigma(p)=-1$. 

Let $I \triangleleft G_{\Q_p}$ be the inertia subgroup. Since $\rho(I)$ has even cardinality, it contains the unique element of $C$ with order $2$, which is $-I_2$. Since $G_E$ is generated by $I$ and $F^2$, $\rho(G_E)=\rho(I)$, so that $\rho(I)$ is a subgroup of $\rho(G_{\Q_p})$ of index two. 

For $k \in \{i,j\}$, let $U_k = \{s \in (\OO_E/p)^{\times} \mid \sigma(s)=k(s)^{-1}\}$: it is a subgroup of $(\OO_E/p)^{\times}$ and let $b_k$ denote its index in $C$. We saw that $\F_p^{\times} \subset U_k$, hence $b_k \mid p+1$. Moreover, since $j=i^p$, $U_i \cap U_j \subset \{s \in (\OO_E/p)^{\times}\mid s^{p-1}=1\}=(p+1)C$, so $p+1$ is the greatest common multiple of $b_i$ and $b_j$. Moreover, 
\[ab_kC=b_k\sigma(\OO_E^{\times})=\sigma(U_k)=k(U_k)^{-1}=b_kC,\]
so that the greatest common divisor of $ab_k$ and $p^2-1$ is $b_k$, hence $a$ and $\frac{p^2-1}{b_k}$ are coprime. Since $a$ divides $b_k\frac{p^2-1}{b_k}$, $a$ divides each $b_k$.
}

\defi[recip-indices-definition]{The pair $\{b_i,b_j\}$ of integers attached to $\rho$ as in Lemma \ref{expand-cartan} is called the \emph{set of reciprocity indices} for $\rho$. This pair only depends on $\rho$ up to conjugacy. }

\prop[local-constants-p-induced]{Let $\rho: G_{\Q_p} \rar \GL{\F_p}$ be an $a$-Cartan representation. Let $E/\Q_p$ be the quadratic unramified extension, $C$ a nonsplit Cartan subgroup such that $\rho^{-1}(C)=G_E$. Let $\psi: C \rar \C^{\times}$ be a character and $\psi(\rho)=\mrm{Ind}_{G_E}^{G_{\Q_p}}{[\psi \circ \rho]}$. 
\begin{enumerate}[noitemsep,label=(\roman*)]
\item \label{lcpi-1} Assume that $\psi_{|\F_p^{\times}}=1$ and let $(f,\mathbf{1}) \in \mathscr{S}$. Then $\psi(\rho)$ is self-dual and 
\[\varepsilon(D_{f,p} \otimes \Delp{p}{\psi(\rho)^{\ast}},1,\theta_p)=r\varepsilon(\psi(\rho),\frac{1}{2},\theta_p)^2,\]
where $r \in \{\pm 1\}$ is equal to $-1$ if and only if $\psi^a=1$, that is, if and only if $\psi(\rho)$ is unramified. 
\item \label{lcpi-2} Let $(f,\chi) \in \mathscr{P}$ and assume that for all $a \in \F_p^{\times}$, $\psi(aI_2)=\chi(a)$. Then 
\[\varepsilon(D_{f,p} \otimes \Delp{p}{\psi(\rho)^{\ast}},1,\theta_p)=\varepsilon(\psi(\rho),\frac{1}{2},\theta_p)\varepsilon(\psi(\rho)^{\ast},\frac{1}{2},\theta_p)\]
\item \label{lcpi-3} Let $(f,\mathbf{1}) \in \mathscr{C}$ and assume that $\psi_{|\F_p^{\times}}=1$. Let $\phi: (\OO_E/p)^{\times} \rar \C^{\times}$ be any of the two characters attached to $C_f$. Then 
\[\varepsilon(D_{f,p} \otimes \Delp{p}{\psi(\rho)^{\ast}},1,\theta_p)=: r \in \{\pm 1\},\]
and $r=-1$ if and only if $\overline{\psi}\circ \rho_{|I_p} \in \{\phi \circ \mathfrak{a}_E,\, \phi^p \circ \mathfrak{a}_E\}$.
\end{enumerate} }

\demo{Write $I_p \triangleleft G_{\Q_p}$ for the inertia subgroup, and $F \in G_{\Q_p}$ be an arithmetic Frobenius such that $\mathfrak{a}_{\Q_p}(F_p)=p^{-1}$. Let $N$ be the normalizer of $C$, so that $\rho(F_p)=-I_2$ by Lemma \ref{expand-cartan}. Note that $\psi(\rho)^{\ast} \simeq \overline{\psi}(\rho)$, and that $\psi^p(\rho) \simeq \psi(\rho)$, since the action of $N/C$ on $C$ by conjugation is given by taking $p$-th powers.  

\ref{lcpi-1}: Since $\psi(\F_p^{\times})=1$, $\psi^{p+1}=1$, so $\psi(\rho)^{\ast} \simeq \overline{\psi}(\rho) \simeq \psi^p(\rho) \simeq \psi(\rho)$. By Proposition \ref{constants-with-modular-forms-elem} \ref{cmfe-3}, we need to check that $\delta=\det(-a_p(f)F \mid \psi(\rho)^I)$ is $1$ if $\psi(\rho)$ is ramified and $-1$ otherwise. Now, $\psi(\rho)_{|I_p} \simeq (\psi\oplus\psi^p)\circ\rho_{|I_p}$, so it has determinant one. Therefore, either $\psi(\rho)^{I_p}$ is zero, or it is $\psi(\rho)$. Since $F \notin G_E$, $\mrm{Tr}(F \mid \psi(\rho))=0$, so $\delta=\det(F \mid \psi(\rho)^I)$ is $-1$ if $\psi(\rho)$ is unramified and $1$ otherwise.

\ref{lcpi-2}: One has \[\psi(\rho) \otimes \chi \simeq \psi^p(\rho) \otimes \overline{\chi}(\det{\rho}) \simeq (\psi^p\overline{\chi}(\det))(\rho) \simeq (\psi^p\overline{\psi}^{p+1})(\rho) \simeq \overline{\psi}(\rho),\] and $D_{f,p}$ is the sum of two characters $\alpha,\beta$ by Lemma \ref{local-global-c}, so we can apply Proposition \ref{constants-with-modular-forms-elem} \ref{cmfe-2}, which yields the result since $\alpha(-1)^2=1$.

\ref{lcpi-3}: Let $\phi_0$ be the unique extension of $\phi$ to $E^{\times}$ such that $\phi_0(p)=-p$. The representation $\rho_{|G_E}$ can be uniquely written as $\sigma \circ \mathfrak{a}_E$ for some character $\sigma: E^{\times}/(1+p\OO_E) \rar C$. By Lemma \ref{expand-cartan}, $\sigma_{|\F_p^{\times}}$ is the inverse of the identity and $\sigma(p)=-I_2$. We saw that $\lambda_{E/\Q_p}(\theta_p)=1$, so that, by Lemma \ref{local-global-c}, one has  

\begin{align*}
\varepsilon(D_{f,p} \otimes \Delp{p}{\psi(\rho)^{\ast}},1,\theta_p) &= \varepsilon((\overline{\phi_0} \otimes (\overline{\psi}(\sigma)\oplus\overline{\psi}^p(\sigma))\circ \mathfrak{a}_E),1,\theta_E)\\
&= \prod_{\zeta \in \{\psi(\sigma),\psi^p(\sigma)\}}{\varepsilon(\overline{\phi_0\zeta},1,\theta_E)}=\prod_{\zeta \in \{\psi(\sigma),\psi^p(\sigma)\}}{\varepsilon(\overline{\phi_0\zeta}\|\cdot\|^{1/2},\frac{1}{2},\theta_E)}.
\end{align*}

If $(\zeta\phi_0)_{|\OO_E^{\times}}$ is nontrivial, then $(\zeta\phi_0)_{|\OO_E^{\times}}$ is a nontrivial character of $\OO_E^{\times}/\Z_p^{\times}(1+p\OO_E)$, so \[\varepsilon(\overline{\phi_0\zeta}\|\cdot\|^{1/2}, \frac{1}{2},\theta_E)=\phi_0(p)\|p\|_E^{1/2}\zeta(p)\frac{\sum_{u \in (\OO_E/p)^{\times}}{(\phi_0\zeta)(u)\theta_E\left(\frac{u}{p}\right)}}{\left|\sum_{u \in (\OO_E/p)^{\times}}{(\phi_0\zeta)(u)\theta_E\left(\frac{u}{p}\right)}\right|}=-\zeta(p)(\phi_0\zeta)(\sqrt{a}),\] for any $a \in \F_p^{\times} \backslash \F_p^{\times 2}$ by Lemma \ref{local-epsilon-fp2}. Since $(\phi_0)_{|\F_p^{\times}}=\psi_{|\F_p^{\times}}=1$, when $\psi \circ \rho_{|I_p} \notin \{\overline{\phi} \circ \mathfrak{a}_E,\overline{\phi}^p\circ \mathfrak{a}_E\}$ one has $\varepsilon(D_{f,p} \otimes \Delp{p}{\psi(\rho)^{\ast}},1,\theta_p)=1.$ 

If $(\psi(\sigma)\phi_0)_{|\OO_E^{\times}}$ is trivial, then $\varepsilon(\overline{\phi_0\psi(\sigma)}\|\cdot\|^{1/2},\frac{1}{2},\theta_E)=1$. Moreover, $\psi^p(\sigma)\phi_0=\psi_0^{-(p-1)}$ is not trivial on $\OO_E^{\times}$, otherwise it would factor through the norm. Then 
\begin{align*}
\varepsilon(D_{f,p} \otimes \Delp{p}{\psi(\rho)^{\ast}},1,\theta_p) &= -\psi^p(\sigma(p))(\phi_0\psi^p(\sigma))(\sqrt{a})=-\psi(-1)\psi^{p-1}(\sigma(\sqrt{a}))\\
&=-\psi^{(p-1)/2}(\sigma(a))=-1.
\end{align*}
}

\section{Root numbers for the factors of $T_{\ell}J_{\rho}$}
\label{root-numbers-surjective}

Let $\rho: G_{\Q} \rar \GL{\F_p}$ be a continuous representation with determinant $\omega_p^{-1}$. We are interested in the $L$-functions for the factors of $J_{\rho}(p)$.

\defi[artin-is-automorphic]{An \emph{irreducible} Artin representation $\sigma: G_{\Q} \rar \GLn{n}{\C}$ is \emph{automorphic} if there exists a unitary cuspidal smooth automorphic representation $\pi \simeq \bigotimes'_v{\pi_v}$ of $\GLn{n}{\mathbb{A}_{\Q}}$ such that, at every place $v$ of $\Q$, the local Langlands correspondence (say, as in \cite{BCarayol} for finite places) attaches $\pi_v$ to $\sigma^{\ast}_{W_{\Q_v}}$. 

We say that an Artin representation $\sigma: G_{\Q} \rar \GLn{n}{\C}$ is \emph{automorphic} if all its irreducible components are automorphic. }

\medskip

For example, in this sense, every odd two-dimensional representation is automorphic by Serre's conjecture \cite{KW1,KW2} and Proposition \ref{local-global-modular-forms}. 

To make some notation easier, let us use the following definitions. 

\defi[compatible-system-definition]{Let $E,F$ be number fields. A compatible system of Galois representations of $G_{E}$ over $F$ is the datum, for each maximal ideal $\lambda \subset \OO_F$, of a continuous semi-simple finite-dimensional representation $V_{\lambda}$ of $G_{E}$ over $F_{\lambda}$ such that the family $(F_{\lambda})_{\lambda}$ satisfies the following conditions:
\begin{itemize}[noitemsep,label=$-$]
\item There exists a finite set $S$ of finite places of $E$ such that, for any maximal ideal $\lambda$ of $\OO_F$ with residue characteristic $\ell$, $V_{\lambda}$ is unramified at every place of $E$ outside $S$ and not dividing $\ell$.
\item For every finite place $v$ of $E$, there exists a Deligne representation $\Delp{v}{V}$ of $W_{E_v}$ over $F$ such that for any maximal ideal $\lambda$ of $\OO_F$ with residue characteristic prime to $v$, one has $V_v \otimes_F F_{\lambda} \simeq \Delp{v}{V_{\lambda}}$.   
\end{itemize}
}

\rem{There are many possible notions of compatible systems in the literature, such as the ones given in \cite[2.30--2.32]{GeeMLT}. Our notion is closer to what \emph{loc.cit.} calls a strongly compatible system, but we do not require any $p$-adic Hodge-theoretic condition. Similarly, while it would be convenient to treat infinite places on the same footing as finite places, for instance by using motives as sketched in \cite[\S 4]{NTB}, we can do so by hand. }

\medskip

With this definition, compatible systems contain Artin representations, the cyclotomic character $(\Q_p(1))_p$, and the $V_f := (V_{f,\lambda})_{\lambda}$ attached to modular forms. They are stable under the usual operations of sum, tensor product, taking duals, restriction, induction, conjugation under automorphisms of $E$, and extensions of scalars.  

In particular, the $R_f(\rho):=(R_{f,\mathfrak{l}}(\rho))_{\mathfrak{l}}$ as in Definition \ref{xrho-factors} form a compatible system over $F$. In fact, when $(f,\chi) \in \mathscr{S} \cup \mathscr{P}'$ (resp. $(f,\chi) \in \mathscr{C}'_1$), $R_f(\rho)$ is defined with coefficients over the subfield of $\C$ generated by the coefficients of $f$ (resp. by the coefficients of $f$ and $\Q(i,\phi_f)^+$ by Proposition \ref{cusp-coeffs})\footnote{This is because $\pi(1,\chi)$ (resp. $\mrm{St}$) can be realized over $\Q(\chi)$ (resp. $\Q$).}. 

By Proposition \ref{cm-yields-compatible-systems}, if $(\psi_{f,\lambda})_{\lambda}$ denotes one of the compatible systems over $G_{\Q(\sqrt{-p})}$ attached to some $(f,\mathbf{1}) \in \mathscr{C}_M$, $\left(\mrm{Ind}_{\Q(\sqrt{-p})}^{\Q}\left[\psi_{f,\lambda} \otimes C^{\epsilon}(\rho_{|G_{\Q(\sqrt{-p})}})\right]\right)_{\lambda}$ also forms a compatible system. 

\medskip

\defi{Let $V := (V_{\lambda})_{\lambda}$ be a compatible system of representations of $G_E$ over a number field $F$. Suppose given an embedding $F \rar \C$. Then $V$ has a global $L$-function (at least as a formal Dirichlet series), a global conductor ideal $\mathfrak{f}_V \subset \OO_E$ and an absolute conductor $N_V$ defined as 
\[L(V,s) := \prod_{v}{L(\Delp{v}{V^{\ast}},s)},\quad \mathfrak{f}_V=\prod_v{\mathfrak{p}_v^{f(\Delp{v}{V^{\ast}})}},\quad N_V=N_{E/\Q}\mathfrak{f}_V.\] where the product runs over finite places of $E$, and, at a finite place $v$ with associated maximal ideal $\mathfrak{p}_v$, the $L$-factor or conductor exponent can be computed, by Lemma \ref{l-adic-to-weil-deligne-numerics}, either on the Deligne representation or on any $V_{\lambda}$ with $\lambda$ coprime to the residue characteristic of $v$. }

\bigskip

\prop[automorphic-times-modform]{Let $f \in \mathcal{S}_k(\Gamma_1(N))$ a normalized newform and $Y: G_{\Q} \rar \GLn{d}{\C}$ be an automorphic Artin representation. The function \[\Lambda(V_f \otimes Y,s)=M^{s/2}2^d(2\pi)^{-ds}\Gamma(s)^{d}L(V_f \otimes Y,s),\] where $M$ and $L(V_f \otimes Y,s)$ are the conductor and the $L$-function, respectively, of the compatible system $(V_{f,\lambda} \otimes Y)_{\lambda}$ (where $\lambda$ runs through maximal ideals of the ring of integers of a number field containing the coefficients of $f$ and on which $Y$ is realizable) extends meromorphically to $\C$ and satisfies the following functional equation:
\[\Lambda(V_f \otimes Y,s)=\prod_v{\varepsilon(D_{f,v} \otimes \Delp{v}{Y^{\ast}},\frac{k}{2},\theta_v)} \cdot \Lambda(Y^{\ast} \otimes V_{\overline{f}},k-s).\]
When $k > 1$, or $k=1$ and $Y$ has no irreducible component isomorphic to $V_f^{\ast}$, then $\Lambda(V_f \otimes Y,s)$ extends to an entire function. 
We call $\prod_v{\varepsilon(D_{f,v} \otimes \Delp{v}{Y^{\ast}},\frac{k}{2},\theta_v)}$ the \emph{Galois-theoretic root number} of the compatible system $V_f \otimes Y$.
}

\demo{We can assume that $Y$ is irreducible and attached to a cuspidal unitary smooth automorphic representation $\pi_Y$. Since $(\Delp{\infty}{Y^{\ast}})_{|W_{\C}}$ is a sum of copies of the trivial representation, while $D_{f,\infty}$ is induced from $W_{\C}$, it follows that $\Lambda(V_f \otimes Y,s)=L\left(\pi_f \times \pi_Y,s-\frac{k-1}{2}\right)$. To reach the conclusion, we can simply apply the theory of the Rankin-Selberg convolution \cite[Theorem 11.7.1]{Getz-Hahn}, while noting that by \cite{BCarayol} the local Langlands correspondence preserves the $L$-factors of products of pairs. }

While we do not know that the Artin representations appearing in $\Tate{\ell}{J_{\rho}(p)}$ are automorphic, we do know that the expected functional equation is, by Lemma \ref{everything-self-dual}, of the form $\Lambda(s) = \epsilon \Lambda(2-s)$, so that $\epsilon$ (the \emph{global root number}) is a sign which can be (in most cases) computed.      

There is one situation where we can be more specific, as is already pointed out in \cite[\S 3.1]{Virdol}: it is that of the CM forms.  

\prop[automorphic-times-modform-cmtwist]{Let $f \in \mathcal{S}_k(\Gamma_1(N))$ be a newform with complex multiplication by an imaginary quadratic field $K$ such that $k\geq 2$. Let $F\supset K$ be a number field such that there exists a compatible system $\psi=(\psi_v)_v$ of characters of $G_K$ over $F$ such that $\mrm{Ind}_K^{\Q}{\psi}=V_f$ (see Proposition \ref{cm-yields-compatible-systems}), and $Y: G_K \rar \GLn{d}{F}$ be an Artin representation. 
Let $\Delp{\infty}{\psi^{\ast}}$ denote the character $z \in W_{\C} \longmapsto z^{1-k}$, $I$ be the compatible system $\mrm{Ind}_K^{\Q}(Y \otimes \psi)$ and $M$ be its conductor. Then the function $\Lambda(I,s) := M^{s/2}2^{d}(2\pi)^{-ds}\Gamma(s)^dL(I,s)$ admits a meromorphic continuation to the complex plane and satisfies the following functional equation: 
\[\Lambda(I,s) = \prod_{v}{\varepsilon(\Delp{v}{\psi^{\ast}} \otimes \Delp{v}{Y^{\ast}},\frac{k}{2},\theta_{K_v})}\Lambda(I^{\ast}(k-1),k-s).\]

We call $\prod_{v}{\varepsilon(\psi(\mathfrak{a}_{K_v}) \otimes \Delp{v}{Y^{\ast}},\frac{k}{2},\theta_{K_v})}$ the \emph{global root number} of the compatible system $I$.
}

\demo{We fix an embedding of $F$ in $\C$. By \cite[(3.5), Theorem 4.5]{Antwerp5-Ribet}, there is a continuous character $\psi_0: \mathbb{A}_K^{\times}/K^{\times}\rar \C^{\times}$ such that for any prime ideal $\mathfrak{p}$ with residue characteristic coprime to $N$ and any place $v$ with residue characteristic coprime to $p$, $\psi_0$ maps any uniformizer at the place $\mathfrak{p}$ to $\psi_v(\Fr_{\mathfrak{p}})$ (which is an element of $F$), where $\Fr_{\mathfrak{p}}$ is an arithmetic Frobenius. Moreover, for $z \in \C^{\times}$, $(\psi_0)(z)=z^{1-k}$. 

Let $V$ be the complex representation of $W_{\Q}$ defined by $\psi_0(\mathfrak{a}) \otimes Y^{\ast}$. After checking $L$-factors at every place and computing conductor exponents (as in \cite[Chapter VII: 10.4, 11.2, 11.7, 12.1]{Neukirch-ANT}), we see that $\Lambda(I,s)$ is exactly $\Lambda(V,s)$ in the notation of Proposition \ref{functional-eqn-weil-gp}, and that $\Lambda(I^{\ast}(k-1),s)$ is exactly $\Lambda(V^{\ast},s-(k-1))$, so that the claim is a direct consequence of the functional equation for representations of $W_K$ \cite[(3.5.4)]{NTB}.
}

\bigskip

As announced, the goal of this section is the computation of the Galois-theoretic root numbers attached to the factors of $J_{\rho}(p)$.

\prop[root-numbers-surjective-computation]{Let $F$ be a large enough number field in the sense of Section \ref{preparatory-lemmas}. Let $\mathfrak{l}$ be a maximal ideal of $\OO_F$ with residue characteristic $\ell$. 

Let, for every $(f,\mathbf{1}) \in \mathscr{C}'_1$, $\phi_f$ be one character of $\F_{p^2}^{\times}$ attached to $f$ and $\omega_f$ be its order. Then, if $\rho$ is $a$-Cartan, let $n_{f,\rho} \in \{0,1\}$ be the integer equal to one if, and only if, $\frac{p+1}{\omega_f}$ is divisible by $a$ and $\frac{2(p+1)}{\omega_f}$ is not divisible by any of the two reciprocity indices of $\rho$.  

Moreover, for $(f,\mathbf{1}) \in \mathscr{C}_M$, let, for each prime $\ell$ inert in $\Q(\sqrt{-p})$, $m_{\ell}$ be one half of the conductor exponent of $V_4(\rho)$ (it is an integer for all $\ell$ and equal to zero for all but finitely many), and let $U=\prod_{\ell}{(-1)^{m_{\ell}}}$. 

The global root numbers for the eigenspaces of $T_{\ell}{J_{\rho}(p)} \otimes_{\Z_{\ell}} F_{\mathfrak{l}}$ are given in Table \ref{table-root-numbers}. 

The entry in the ``CM row'' of the table means the following. Let $K=\Q(\sqrt{-p})$ and let $\psi:\mathbb{A}_K^{\times}/K^{\times} \rar \C^{\times}$ one of the two continuous characters attached to $f$ by Proposition \ref{cm-yields-compatible-systems}. We compute the $L$-function of $\mrm{Ind}_{\Q(\sqrt{-p})}^{\Q}{\psi(\mathfrak{a}_K) \otimes (C^{\epsilon})^{\ast}(\rho_{|G_K})}$ in the normalization of \cite{NTB} (that is, we consider the situation of Proposition \ref{automorphic-times-modform-cmtwist} with $Y=C^{\epsilon}(\rho_{|G_K})$). 

Furthermore, when $\rho_{|G_{\Q_p}}$ is wild, it is $\SL{\F_p}$-conjugate to a representation of the form $\beta \otimes \begin{pmatrix}q & \nu\\0 & 1\end{pmatrix}$, where $\beta,q: G_{\Q_p} \rar \F_p^{\times}$ are characters and $\nu: G_{\Q_p} \rar \F_p(q)$ is a cocycle whose restriction to the wild inertia subgroup is nonzero. Let $G_L=\ker{q}$, $d=[L:\Q_p]$ and $\nu_0: G_L \rar \F_p$ be the character given by $\nu=\nu_0(\mathfrak{a}_L)$, whose conductor exponent $n$ is even. There is an element $\gamma \in L^{\times}$ with normalized valuation $1-n$ such that, for any $x \in \OO_L$ with valuation at least $\frac{n}{2}$, $e^{\frac{2i\pi}{p}\nu_0(1+x)}=\eta_L\left(\gamma x\right)$. Let $K_p$ denote the completion of $K$ at its unique place dividing $p$ and $\gamma_0=N_{L/K_p}\gamma \in K_p^{\times}$. Then $\varepsilon_{\mathfrak{p}}(\psi,\nu)=i\psi(\gamma_0)\|\gamma_0\|^{1/2} \in \{\pm 1\}$. 
}

\begin{table}[htb]
\centering
\begin{tabular}{|c|c|c|c|}
\hline
& \multicolumn{3}{c|}{Type of $\rho_{|G_{\Q_p}}$}\\
\hline
& Split & Wild & $a$-Cartan\\
\hline
$R_{f,\mathfrak{l}}(\rho), (f,\mathbf{1}) \in \mathscr{S}$ & $(-1)^{(p+1)/2}$ & $(-1)^{(p-1)/2}a_p(f)$ & $(-1)^{(p+a)/2}$\\
\hline
$R_{f,\mathfrak{l}}(\rho), (f,\chi) \in \mathscr{P}$ & \multicolumn{3}{c|}{$(-1)^{(p-1)/2}$}\\
\hline
$R_{f,\mathfrak{l}}(\rho), (f,\mathbf{1}) \in \mathscr{C}'_1$ & $(-1)^{(p+1)/2}$ & $(-1)^{(p-1)/2}$ & $(-1)^{(p+1)/2+n}$\\
\hline 
CM factors & $(-1)^{(p+1)/4}U$ & $(-1)^{\frac{d+2}{4}}U\epsilon\varepsilon_{\mathfrak{p}}(\psi,\nu)$ & $(-1)^{\frac{p+2a-1}{4}}U$\\
\hline
\end{tabular}
\caption{Global root numbers for factors of $J_{\rho}(p)$}
\label{table-root-numbers}
\end{table}

\demo{Throughout the proof, we choose a decomposition subgroup $G_{\Q_p}$ of $G_{\Q}$, $I_p \leq G_{\Q_p}$ is the inertia subgroup, $I_p^+ \leq I_p$ is the wild inertia subgroup, and $F_p \in G_{\Q_p}$ is an arithmetic Frobenius such that $\mathfrak{a}_{\Q_p}(F_p)=\frac{1}{p}$. The unique quadratic unramified extension of $\Q_p$ is denoted by $E$, and we call $c \in G_{\Q}$ the complex conjugation.

\textbf{For $(f,\chi) \in \mathscr{P}$:}

Write $P=\pi(1,\chi) \circ \rho$. Then 
\[P^{\ast} \simeq \pi(1,\chi)^{\ast} \circ \rho \simeq \pi(1,\chi^{-1}) \circ \rho \simeq \chi^{-1}(\det{\rho}) \otimes (\pi(\chi,1)\circ \rho) \simeq \chi \otimes P. \]

By Proposition \ref{constants-with-modular-forms-elem} \ref{cmfe-1}, $\varepsilon(D_{f,\infty} \otimes \Delp{\infty}{P^{\ast}},1,\theta_{\infty})=(-i)^{2(p+1)}=1$. 

For any finite place $v \nmid p$, $D_{f,v}$ is the sum of two characters (when $v \neq p$, $D_{f,v}$ is unramified, and when $v=p$ it is Lemma \ref{local-global-p}), so by Proposition \ref{constants-with-modular-forms-elem} \ref{cmfe-2}, one has \[\varepsilon(D_{f,v} \otimes \Delp{v}{P^{\ast}},1,\theta_v)=\varepsilon(\Delp{v}{P \oplus P^{\ast}},1,\theta_v).\] Therefore, 
\begin{align*}
\prod_{v}{\varepsilon(D_{f,v} \otimes \Delp{v}{P^{\ast}},1,\theta_v)}&= \prod_{q < \infty}{\varepsilon(\Delp{v}{P \oplus P^{\ast}},\frac{1}{2},\theta_v)}=\varepsilon(\Delp{\infty}{P \oplus P^{\ast}},\frac{1}{2},\theta_{\infty})^{-1}\\
&=\det(c \mid P)^{-1}\end{align*}
Since $\rho(c) \sim \Delta_{1,-1}$, $P(c)^2$ is the identity and $P(c)$ has trace $2$ by Proposition \ref{character-from-principal-series}, so $\det(c \mid P)=(-1)^{(p-1)/2}$. \\

\textbf{For $(f,\chi) \in \mathscr{S}$:}

Let $S=\mrm{St} \circ \rho$. Since $\mrm{St}$ is self-dual, $S$ is self-dual too. Therefore, by Proposition \ref{constants-with-modular-forms-elem} \ref{cmfe-2}, for any finite place $v \neq p$, $D_{f,v}$ is the sum of two unramified characters, so \[\varepsilon(D_{f,v} \otimes \Delp{v}{S^{\ast}},1,\theta_v)=\varepsilon(\Delp{v}{S \oplus S^{\ast}},\frac{1}{2},\theta_v).\] By Proposition \ref{constants-with-modular-forms-elem}, it follows that 
\begin{align*}
\prod_{v}{\varepsilon(D_{f,v} \otimes \Delp{v}{S^{\ast}},1,\theta_v)} &= (-i)^{2p}\varepsilon_p(D_{f,p} \otimes \Delp{p}{S \oplus S^{\ast}},1,\theta_p)\prod_{v \neq p,\infty}{\varepsilon(\Delp{v}{S\oplus S^{\ast}},\frac{1}{2},\theta_v)}\\
&=-\det(-a_p(f)F_p\mid S^{I_p})\prod_{v \neq \infty}{\varepsilon(\Delp{v}{S \oplus S^{\ast}},\frac{1}{2},\theta_v)}\\
&=-\det(-a_p(f)F_p\mid S^{I_p})\varepsilon(\Delp{\infty}{S \oplus S^{\ast}},\frac{1}{2},\theta_v)^{-1}\\
&=-\det(c \mid S)\det(-a_p(f)F_p \mid S^{I_p}).
\end{align*}
By Proposition \ref{character-from-principal-series}, $S(c)$ has trace $1$ and order $2$, so $\det(c \mid S)=(-1)^{(p-1)/2}$. Therefore,  
\[\prod_{v}{\varepsilon(D_{f,v} \otimes \Delp{v}{S^{\ast}},1,\theta_v)} = (-1)^{(p+1)/2}\det(-a_p(f)F_p \mid S^{I_p}).\]~\\

\emph{When $\rho_{|G_{\Q_p}}$ is wild:}

After conjugating $\rho$ (which does not cause any issue), we may assume that $\rho$ maps the wild inertia subgroup $I_p^+$ to $U^{\Z}=\begin{pmatrix}1 &\ast\\0 & 1\end{pmatrix}$. Then $\rho(G_{\Q_p})$ is contained in the normalizer of $U^{\Z}$, which is the group $B$ of upper-triangular matrices. Lemma \ref{steinberg-on-borel} therefore implies that $S^{I_p}=S^{G_{\Q_p}}$ is a line, so that $\det(-a_p(f)F_p \mid S^{I_p})=-a_p(f)$, whence the conclusion. \\

\emph{When $\rho_{|G_{\Q_p}}$ is split:}

After conjugating $\rho$, we may assume that for $s \in G_{\Q_p}$, $\rho(s)=\begin{pmatrix}\alpha(s) & 0\\0 & \beta(s)\end{pmatrix}$ for some characters $\alpha,\beta: G_{\Q_p} \rar \F_p^{\times}$. By Lemma \ref{steinberg-on-diag}, $S_{|G_{\Q_p}} \simeq \mathbf{1} \oplus \bigoplus_{\chi \in \mathcal{D}}{\chi\left(\frac{\alpha}{\beta}\right)}$. 

Let $a$ be the cardinality of $\rho(I_p) \pmod{\F_p^{\times}I_2}$. Then $S^{I_p}_{|G_{\Q_p}} \simeq \mathbf{1} \oplus  \bigoplus_{\chi \in \mathcal{D}, \chi^{(p-1)/a}=1}{\chi\left(\frac{\alpha}{\beta}\right)}$, so that 

\[\det(-a_p(f)F_p \mid S^{I_p}) = (-a_p(f))^{1+\frac{p-1}{a}} \prod_{\chi \in \mathcal{D}, \chi^{\frac{p-1}{a}}=1}{\chi\left(\frac{\alpha}{\beta}(F_p)\right)}.\]

To prove that this quantity is one, it is enough to show that $\frac{p-1}{a}$ is odd (because $a_p(f)^2=1$ and the factors with $\chi$ and $\chi^{-1}$ in the product cancel out). Now, $\frac{p-1}{a}$ is the cardinality of the group $\rho(I_p) \cap \F_p^{\times} I_2$. Since $\det: \rho(I_p) \rar \F_p^{\times}$ is an isomorphism by Proposition \ref{trichotomy-rho}, $-I_2 \notin \rho(I_p)$, so $\rho(I_p) \cap \F_p^{\times}I_2$ does not have an element of order $2$: thus $\rho(I_p) \cap \F_p^{\times}I_2$ has odd order, and we are done. \\

\emph{When $\rho_{|G_{\Q_p}}$ is $a$-Cartan:}

Let $C$ be a nonsplit Cartan subgroup of $\GL{\F_p}$ with normalizer $N$ such that $\rho(G_{\Q_p}) \subset N$ and $\rho^{-1}(C) = G_E$. Our assumption states that $\rho(aI_p)=\rho(G_E)=aC$ by Lemma \ref{expand-cartan}. By Proposition \ref{steinberg-on-norm-cartan}, since $\det{\rho}(I_p)=\F_p^{\times}$, it follows that 

\[S^{I_p}_{|G_{\Q_p}} \simeq \bigoplus'_{\psi^a=1}{\mrm{Ind}_{E}^{\Q_p}{\psi\circ \rho_{|G_E}}} \simeq \mathbf{1}^{\oplus \frac{a-1}{2}} \oplus \epsilon^{\oplus \frac{a-1}{2}},\]

where the first direct sum is over characters $\psi: C/\F_p^{\times} \rar \C^{\times}$ such that $\psi^2 \neq 1$, $\psi^a=1$, and we only count one character among $\psi$ and $\psi^{-1}$, so that there are $\frac{a-1}{2}$ terms in the sum, and $\epsilon$ is the character $\mrm{Gal}(E/\Q_p) \overset{\sim}{\rar} \{\pm 1\}$. 

Thus $\det(-a_p(f)F_p \mid S^{I_p}) = (-a_p(f))^{a-1} (-1)^{\frac{a-1}{2}}=(-1)^{\frac{a-1}{2}}$ and we are done. \\

\textbf{For $(f,\mathbf{1}) \in \mathscr{C}'_1$:}

Let $W=C_f(\rho)$, so that, as in the proof of Lemma \ref{everything-self-dual}, $W^{\ast} = W$. Let $\phi: E^{\times}/(1+p\OO_E) \rar \C^{\times}$ be attached to $f$ as in Lemma \ref{local-global-c}, and $\omega$ the order of $\phi_{|\OO_E^{\times}}: (\OO_E/p)^{\times}/\F_p^{\times} \rar \C^{\times}$.  

First, we compute that $\varepsilon(D_{f,\infty} \otimes \Delp{\infty}{W^{\ast}},1,\theta_{\infty})=(-i)^{2(p-1)}=1$ by Proposition \ref{constants-with-modular-forms-elem} \ref{cmfe-1}. For any finite place $v \neq p$, $D_{f,p}$ is the direct sum of two unramified characters so, by Proposition \ref{constants-with-modular-forms-elem} \ref{cmfe-2}, $\varepsilon(D_{f,v} \otimes \Delp{v}{W^{\ast}},1,\theta_v)=\varepsilon(\Delp{v}{W \oplus W^{\ast}},\frac{1}{2},\theta_v)$. Therefore, using Lemma \ref{det-cusp} as above,

\begin{align*}
\prod_{v}{\varepsilon(D_{f,v} \otimes \Delp{v}{W^{\ast}},1,\theta_v)} &= \frac{\varepsilon(D_{f,p} \otimes \Delp{p}{W^{\ast}},1,\theta_p)}{\varepsilon(\Delp{p}{W \oplus W^{\ast}},\frac{1}{2},\theta_p)}\varepsilon(\Delp{\infty}{W \oplus W^{\ast}},\frac{1}{2},\theta_{\infty})^{-1}\\
&= \frac{\varepsilon(D_{f,p} \otimes \Delp{p}{W^{\ast}},1,\theta_p)}{\det(c \mid W)\det{W}(\mathfrak{a}_{\Q_p}^{-1}(-1))}\\
&=\varepsilon(D_{f,p} \otimes \Delp{p}{W^{\ast}},1,\theta_p)\left(\frac{\det{\rho}(c)}{p}\right)\left(\frac{\det{\rho}(\mathfrak{a}_{\Q_p}^{-1}(-1))}{W}\right)\\
&= \varepsilon(D_{f,p} \otimes \Delp{p}{W^{\ast}},1,\theta_p)(-1)^2=\varepsilon(D_{f,p} \otimes \Delp{p}{W^{\ast}},1,\theta_p).
\end{align*}

\emph{$\rho_{|G_{\Q_p}}$ is split:}

After conjugating, we may assume that $\rho(s)=\begin{pmatrix}\alpha(s) & 0\\0 & (\omega_p^{-1}\alpha^{-1})(s)\end{pmatrix}$ for some character $\alpha: G_{\Q_p} \rar \F_p^{\times}$. Write $\alpha=\alpha_0\omega_p^t$, for some $t \in \Z/(p-1)\Z$ and some unramified character $\alpha_0: G_{\Q_p} \rar \F_p^{\times}$. By Proposition \ref{cusp-on-diag}, one has 

\[D_{f,p} \otimes \Delp{p}{W^{\ast}}\simeq \bigoplus_{\psi \in \mathcal{D}}{\overline{\psi}(\alpha_0^2) \otimes \Delp{p}{(V_{f,\lambda} \otimes \psi(\omega_p^{2t-1}))^{\ast}}}\simeq \bigoplus_{\psi \in \mathcal{D}}{\overline{\psi}(\alpha_0^2) \otimes D_{f \otimes \psi^{2t-1},p}},\]
hence, since $\prod_{\psi \in \mathcal{D}}{\psi^2}=\mathbf{1}$,
\begin{align*}
\varepsilon(D_{f,p} \otimes \Delp{W^{\ast}},1,\theta_p)&=\prod_{\psi \in \mathcal{D}}{\varepsilon(\overline{\psi}(\alpha_0^2) \otimes D_{f \otimes \psi^{2t-1},p},1,\theta_p)}= \prod_{\psi \in \mathcal{D}}{\overline{\psi}(\alpha_0^2(F_p^{-2}))\varepsilon(D_{f \otimes \psi^{2t-1},p},1,\theta_p)}\\
&= \prod_{\psi \in \mathcal{D}}{\psi(\alpha_0(F_p))^4} \varepsilon_p(f \otimes \psi^{2t-1})=\prod_{\psi \in \mathcal{D}}{\lambda_p(f \otimes \psi^{2t-1})}.
\end{align*}

For any $\psi \in \mathcal{D}$, by Proposition \ref{local-constant-elementary}, $\lambda_p(f \otimes \psi^{2t-1})\lambda_p(f \otimes \overline{\psi}^{2t-1})=1$ because $f$ has real coefficients. Therefore, if $\lambda \in \mathcal{D}$ is the nontrivial quadratic character, Corollary \ref{various-computations-cuspidal} implies that $\varepsilon(D_{f,p} \otimes \Delp{p}{V^{\ast}},1,\theta_p)=\lambda_p(f)\lambda_p(f\otimes \lambda)=(-1)^{\frac{p+1}{2}}$, and we are done.\\

\emph{$\rho_{|G_{\Q_p}}$ is $a$-Cartan:}

Let $C \leq \GL{\F_p}$ be a Cartan subgroup with normalizer $N$ such that $\rho(G_{\Q_p}) \subset N$ and $\rho^{-1}(C)=G_E$. Let $\lambda \in \mathcal{D}$ be the nontrivial quadratic character and $\epsilon: \mrm{Gal}(E/\Q_p) \overset{\sim}{\rar} \{\pm 1\}$. By Lemma \ref{expand-cartan}, let $\sigma: E^{\times}/(1+p\OO_E) \rar C$ be the group homomorphism such that $\rho_{|G_E} = \sigma \circ \mathfrak{a}_E$, with $\sigma(p)=-I_2$ and $\sigma(x)=x^{-1}I_2$ for $x \in \F_p^{\times}$. Fix an isomorphism $\iota: C \rar (\OO_E/p)^{\times}$ preserving the trace and the norm, and let $b_{\iota},b_{\iota^p}$ be the two reciprocity indices. 

 By Proposition \ref{cusp-on-norm-cartan}, in the notation of Proposition \ref{local-constants-p-induced},
\[W_{|G_{\Q_p}} \simeq \bigoplus_{\chi \in \{1,\lambda\}}{\chi(\omega_p^{-1})\epsilon^{r_{\chi}}} \oplus \bigoplus^{'}_{\psi}{\psi(\rho_{|G_{\Q_p}})},\] 
where $r_{\chi} \in \{0,1\}$ depends on $\chi,\phi$ and $\psi$ runs through the characters $C/\F_p^{\times} \rar \C^{\times}$, where $\psi \notin \{\phi,\phi^p\}$, $\psi^2 \neq \mathbf{1}$, and we count exactly one character in each pair $\{\psi,\psi^{-1}\}$. 

If $\chi \in \{1,\lambda\}$, 
\begin{align*}
\varepsilon(D_{f,p} \otimes \Delp{p}{(\chi(\omega_p^{-1}) \otimes \epsilon^{r_{\chi}})^{\ast}},1,\theta_p) &= \varepsilon(D_{f \otimes \chi^{-1},p} \otimes \epsilon^{r_{\chi}},1,\theta_p)=\epsilon^{r_{\chi}}(F_p^{-2})\varepsilon(D_{f \otimes \chi^{-1},p},1,\theta_p)\\
&=\varepsilon_p(f \otimes \chi^{-1})=\lambda_p(f \otimes \chi^{-1})=\lambda_p(f \otimes \chi).
\end{align*}

By Corollary \ref{various-computations-cuspidal} and Proposition \ref{local-constants-p-induced}, 
\[\varepsilon(D_{f,p} \otimes \Delp{W^{\ast}},1,\theta_p)=(-1)^{\frac{p+1}{2}+m},\] where $m \in \{0,1\}$ has the same parity as the number of $\psi: C/\F_p^{\times} \rar \C^{\times}$ such that $\psi \neq \phi,\phi^p$, $\psi^2 \neq 1$, and $\psi\circ \sigma\circ \mathfrak{a}_E$ coincides with either $\overline{\phi}\circ \mathfrak{a}_E$ or $\overline{\phi}^p\circ \mathfrak{a}_E$ on $I_p$, and we count one character in each pair $\{\phi,\psi^p\}$. 

 Thus, $m$ has the same parity as the number of $\psi: C/\F_p^{\times} \rar \C^{\times}$ with order at most $3$ and distinct from $\phi \circ \iota,\phi^p \circ \iota$ and such that $\psi \circ \sigma$ coincides with $\overline{\phi}$ on $(\OO_E/p)^{\times}$. If $\psi: C/\F_p^{\times} \rar \C^{\times}$ is a character with order at most $2$ such that $\psi \circ \sigma$ coincides with $\overline{\phi}$ on $(\OO_E/p)^{\times}$, then $\phi^2_{|\OO_E^{\times}}=1$. Then $\phi$ factors through the norm $E^{\times} \rar \Q_p^{\times}$, which is impossible. 
 
 Therefore, $m$ has the same parity as the number of $\psi: C/\F_p^{\times} \rar \C^{\times}$ distinct from $\phi \circ \iota,\phi^p \circ \iota$ and such that $\psi \circ \sigma$ coincides with $\overline{\phi}$ on $(\OO_E/p)^{\times}$. 
 
Since $\sigma(\OO_E^{\times})=aC$, there exists some $\psi: C/\F_p^{\times} \rar \C^{\times}$ such that $\overline{\phi}_{|\OO_E^{\times}}=\psi \circ \sigma_{|\OO_E^{\times}}$ if and only if $\phi$ is a $a$-th power, or equivalently if $\omega \mid \frac{p+1}{a}$. In this case, there are exactly $a$ characters $\psi$ satisfying the above condition, where $a$ is odd. Otherwise, $m=0$. 
 
Assume that $\psi=\phi \circ \iota$ and $\psi=\phi^p \circ \iota$ both satisfy $\psi \circ \sigma_{|\OO_E^{\times}}=\overline{\phi}_{|\OO_E^{\times}}$. Then $\overline{\phi}_{|\OO_E^{\times}}=\overline{\phi}^p_{|\OO_E^{\times}}$, so $\phi^{p-1}_{|\OO_E^{\times}}=\phi^{p+1}_{|\OO_E^{\times}}=1$, so $\phi^2_{|\OO_E^{\times}}=1$, which we already saw was impossible. 

For $i \in \{1,p\}$, $\psi := \phi^i \circ \iota$ satisfies $(\psi \circ \iota) \circ \sigma_{|\OO_E^{\times}} = \overline{\phi}_{|\OO_E^{\times}}$ if and only if $\phi$ vanishes on the subgroup $\{\iota^i(\sigma(x))x\mid x \in (\OO_E/p)^{\times}\}$ of $(\OO_E/p)^{\times}$ (with cardinality $b_{\iota^i}$), if and only if $\omega \mid \frac{p^2-1}{b_{\iota^i}}$. Since $p \equiv -1\pmod{\omega}$ and $b_{\iota^i} \mid p+1$, $\omega \mid  \frac{p^2-1}{b_{\iota^i}}$ if and only if $\omega \mid \frac{2(p+1)}{b_{\iota^i}}$. \\

\emph{When $\rho_{|G_{\Q_p}}$ is wild:}

By Lemma \ref{expand-wild}, after conjugating, we can assume that $\rho_{|G_{\Q_p}}=\beta\begin{pmatrix}q & \nu\\0 & 1\end{pmatrix}$, for some characters $\beta, q: G_{\Q_p} \rar \F_p^{\times}$ such that the extension $L/\Q_p$ cutting out $\ker{q}$ is tamely ramified with ramification index $e$ such that $v_2(e)=v_2(p-1)$, and for some cocycle $\nu: G_{\Q_p} \rar \F_p(q)$ is a cocycle whose restriction to $I_p^+$ is surjective, and such that the conductor exponent $n \geq 2$ of the wildly ramified character $\nu_{|G_L}: G_L \rar \F_p^{\times}$ is such that $n-1$ is coprime to $e$ (so $n$ is even). 

Write $\nu=\nu_0 \circ \mathfrak{a}_L$ for some character $\nu_0: L^{\times} \rar \F_p$. By the proof of Lemma \ref{expand-wild}, for any $\sigma \in \mrm{Gal}(L/\Q_p)$, $\nu_0(\sigma(x))=q(x)\nu_0(\sigma)$. A consequence is that the restriction of $\nu_0$ to any proper subfield of $L$ vanishes. 

Let $K \subset L$ be the largest subextension which is unramified over $\Q_p$, and let $\pi \in K, \varpi \in L$ be uniformizers such that $\varpi^e=\pi$. Thus $e\nu_0(\varpi)=\nu_0(\pi)=0$, so $\nu_0(\varpi)=0$, and $\nu_0(\Q_p^{\times})=\{0\}$. 

By Proposition \ref{cusp-on-borel}, if $B$ denotes the subgroup of upper-triangular matrices and $U=\begin{pmatrix} 1 & 1\\0 & 1\end{pmatrix}$, then 
\[(C_f)_{|B}=\mrm{Ind}_{\F_p^{\times}U^{\Z}}^{B}{\left[a\begin{pmatrix}1 & b\\0 & 1\end{pmatrix} \longmapsto e^{2i\pi b/p}\right]}.\]

Since $\begin{pmatrix}a & \ast\\0 & b\end{pmatrix} \in B/\F_p^{\times}U^{\Z} \longmapsto ab^{-1} \in \F_p^{\times}$ is an isomorphism, by Mackey's formula \cite[Proposition 22]{SerreLinReps} and Lemma \ref{local-global-c},

\begin{align*}
\Delp{p}{W^{\ast}} &\simeq \bigoplus_{t \in \F_p^{\times}/q(G_{\Q_p})}{\mrm{Ind}_L^{\Q_p}{e^{-2i\pi t\nu_0/p} \circ \mathfrak{a}_L}}\\
D_{f,p} \otimes \Delp{p}{W^{\ast}} &\simeq \bigoplus_{t \in \F_p^{\times}/q(G_{\Q_p})}{\mrm{Ind}_{LE}^{\Q_p}{\left[(e^{-2i\pi t\nu_0/p} \circ \mathfrak{a}_L) \otimes \overline{\phi}\circ \mathfrak{a}_E\right]}}\\
&\simeq \bigoplus_{t \in \F_p^{\times}/q(G_{\Q_p})}{\mrm{Ind}_{LE}^{\Q_p}{\left[(e^{-2i\pi t\nu_0(N_{LE/L})/p} \cdot \overline{\phi}(N_{LE/E}))\circ \mathfrak{a}_{LE}\right]}}.
\end{align*}

Since $v_2(e)=v_2(p-1)$, $[\F_p^{\times}:q(G_{\Q_p})]$ is odd, so 

\begin{align*}
&\varepsilon(D_{f,p} \otimes \Delp{p}{W^{\ast}},1,\theta_p) = \lambda_{LE/\Q_p}(\theta_p)^{[\F_p^{\times}:q(G_{\Q_p})]}\prod_{t \in \F_p^{\times}/q(G_{\Q_p})}{\varepsilon(e^{\frac{-2i\pi t}{p}\nu_0(N_{LE/L})} \cdot \overline{\phi}(N_{LE/E}),1,\theta_E)}\\
&= (-1)^{(p+1)/2}\prod_{t \in \F_p^{\times}/q(G_{\Q_p})}{\varepsilon(e^{\frac{-2i\pi t}{p}\nu_0(N_{LE/L})} \cdot \overline{\phi}(N_{LE/E}),1,\eta_{LE}(p\cdot))}\\
&= (-1)^{(p+1)/2}\prod_{t \in \F_p^{\times}/q(G_{\Q_p})}{e^{\frac{-2i\pi t}{p}\nu_0(p^2)}\overline{\phi}(p^{[L:\Q_p]})\|p\|_{LE}^{1/2}\varepsilon(e^{\frac{-2i\pi t}{p}\nu_0(N_{LE/L})} \cdot \overline{\phi}(N_{LE/E}),1,\eta_{LE})}\\
&= (-1)^{(p+1)/2}\prod_{t \in \F_p^{\times}/q(G_{\Q_p})}{\varepsilon(e^{\frac{-2i\pi t}{p}\nu_0(N_{LE/L})} \cdot \overline{\phi}(N_{LE/E}),1,\eta_{LE})}.
\end{align*}

Let us fix, from now on, some $t \in \F_p^{\times}$. There exists some $\gamma \in L^{\times}$ such that for every $x \in \varpi^{n/2}\OO_L$, $e^{2i\pi \nu_0(x)/p}=\eta_L(\gamma x)$. Then $\gamma$ is a well-defined element of $\varpi^{1-n}\OO_L/\varpi^{1-n/2}\OO_L$ and it has valuation exactly $1-n$ (since $L$ is tamely ramified, $\eta_L$ vanishes on $\varpi\OO_L$ but not on $\OO_L$). 

For any $\sigma \in \mrm{Gal}(L/\Q_p)$ and $q_0 \in \Z_p^{\times} \cap \N$ lifting $q(\sigma)$, for any $x \in \varpi^{n/2}\OO_L$, 
\begin{align*}
\eta_L(q_0\sigma(\gamma)x)&=\eta_L(q_0\gamma\sigma^{-1}(x))=e^{2i\pi \nu_0(1+q_0\sigma^{-1}(x))/p}=e^{2i\pi \nu_0((1+\sigma^{-1}(x))^{q_0})/p}\\
&=e^{2i\pi q(\sigma)\nu_0(\sigma^{-1}(1+x))/p}=e^{2i\pi \nu_0(1+x)/p}=\eta_L(\gamma x),
\end{align*}
so $\sigma(\gamma) \equiv q(\sigma)^{-1}\gamma \pmod{\varpi^{1-n/2}\OO_L}$. 

Write $\gamma=\varpi^{1-n}\zeta$ for some $\zeta \in (\OO_L/\varpi^{n/2})^{\times}$. In the computation above, if we take $\sigma \in I_{\Q_p}$, by Lemma \ref{expand-wild}, we see that $\zeta \in (\OO_L/\varpi^{n/2})^{\times}$ is fixed by $\mrm{Gal}(L/K)$. Since $\mrm{Gal}(L/K)$ has order prime to $p$ and $1+\varpi^{n/2}\OO_L$ is a pro-$p$-group, $H^1(\mrm{Gal}(L/K),1+\varpi^{n/2}\OO_L)=\{0\}$. Thus $H^0(\mrm{Gal}(L/K),\OO_L^{\times}) \rar H^0(\mrm{Gal}(L/K),(\OO_L/\varpi^{n/2})^{\times})$ is onto and we can take $\gamma=\varpi^{1-n}\zeta$ for $\zeta \in \OO_K^{\times}$. Since $1-n$ and $e=[L:K]$ are coprime, $\mrm{Tr}_{L/K}(\gamma)=\zeta\mrm{Tr}_{L/K}(\varpi^{1-n})=0$, so $\eta_L(\gamma)=1$. Moreover, we have $\nu_0(\gamma)=(1-n)\nu_0(\varpi)+\nu_0(\zeta)=0$. 

The extension $LE/L$ is unramified, so that $N_{LE/L}: \frac{1+\varpi^r\OO_{LE}}{1+\varpi^{r+1}\OO_{LE}} \rar \frac{1+\varpi^r\OO_L}{1+\varpi^{r+1}\OO_L}$ is onto. Hence the conductor exponent of the character $e^{-2i\pi t \nu_0(N_{LE/L})/p}$ is $n$. Moreover, for any $z \in 1+\varpi\OO_{LE}$, $N_{LE/E}(z) \in (1+\varpi\OO_{LE}) \cap \OO_E=1+p\OO_E$, hence $\overline{\phi}(N_{LE/E})$ has conductor exponent two, so the conductor exponent of its product with $e^{-2i\pi t \nu_0(N_{LE/L})/p}$ is $n$. 

Let finally $a \in \OO_{LE}$, $a' \in \OO_{LE}$ its unique conjugate under $\mrm{Gal}(LE/L)$, and $x=\varpi^{n/2}a$. Then
\begin{align*}
\overline{\phi}(N_{LE/E}(1+x))e^{-2i\pi t \nu_0(N_{LE/L}(1+x))/p} &= e^{-\frac{2i\pi t}{p}\nu_0(1+(a+a')\varpi^{n/2}+aa'\varpi^n)}=\left(e^{\frac{2i\pi}{p}\nu_0(1+\mrm{Tr}_{LE/L}(x))}\right)^{-t} \\
&= \eta_L(\gamma \mrm{Tr}_{EL/L}(x))^{-t}=\eta_{LE}(-t\gamma x).
\end{align*}

We can then apply \cite[\S 23.6 Proposition]{BH} (note that by their definition of level in \cite[\S 1.8]{BH}, the character has level $n-1$): 
\begin{align*}
&\varepsilon(e^{-2i\pi t\nu_0(N_{LE/L})/p} \cdot \overline{\phi}(N_{LE/E}),1,\eta_{LE})\\
&=p^{-\frac{2(n-1)f_{L/\Q_p}}{2}}e^{\frac{2i\pi t}{p}\nu_0(N_{LE/L}(-t\gamma))} \cdot \overline{\phi}^{-1}(N_{LE/E}(-t\gamma))\eta_{LE}(-t\gamma)\\
&= p^{-(n-1)f_{L/\Q_p}}e^{\frac{2i\pi t}{p} \nu_0(t^2\gamma^2)}\overline{\phi}^{-1}(N_{L/\Q_p}(-t\zeta \varpi^{1-n}))\\
&=p^{-(n-1)f_{L/\Q_p}}\phi(p)^{(n-1)f_{L/\Q_p}}=-1,
\end{align*}

and therefore $\varepsilon(D_{f,p} \otimes \Delp{p}{W^{\ast}},1,\theta_p)=-(-1)^{(p+1)/2}=(-1)^{(p-1)/2}$. \\

\textbf{When $(f,\mathbf{1}) \in \mathscr{C}_M$:}

Let $\psi$ be one of the two compatible systems of Galois representations of $K:=\Q(\sqrt{-p})$ with coefficients in (completions of) $F$ attached to $f$ and $\epsilon$ be a sign. Let $\psi^0: \mathbb{A}_K^{\times}/K^{\times} \rar \C^{\times}$ be the continuous adelic character of $K$ attached to $f$ as in \cite[(3.5)]{Antwerp5-Ribet}, normalized so that for any two finite places $v$ of $K$, $\lambda$ of $F$, where the residue characteristic of $v$ is coprime to $p\lambda$, $\psi^0$ maps a uniformizer at $v$ to the image of the arithmetic Frobenius by $\psi_{\lambda}$, and $\psi^0$ maps a complex number $z \in \C^{\times}$ to $\sigma(z)^{1-k}$ for some embedding $\sigma: K \rar \C$.

Since $f$ has real coefficients, $\overline{\psi^0}$ is the other complex character attached to $f$, and $\overline{\psi^0}=\psi^0\circ\tau$, where $\tau$ is the automorphism of $K/\Q$. Thus the character $|\psi^0|^2: \mathbb{A}_K^{\times}/K^{\times} \rar \R^{+\times}$ maps any uniformizer at a split prime $v$ of $K$ to its residue characteristic and any complex number $z$ to $\sigma(z)^{-1}\overline{\sigma(z)}^{-1}$, so that $\|\psi^2\|=\|\cdot\|_{\mathbb{A}}^{-1}$. 

Let $\mathfrak{p}$ denote the only prime ideal of $K$ above $p$. 
\begingroup
\newcommand{\pp}{\mathfrak{p}}
\newcommand{\Kpp}{K_{\mathfrak{p}}}

Because $f$ has conductor $p^2$, $\psi$ has conductor $\pp$ and $\psi^0$ is trivial on the compact open subgroup $(1+\pp\OO_{\Kpp})\prod_{v \nmid p}{\OO_{K_v}^{\times}}$. As we saw, the global root number for the compatible system $\mrm{Ind}_K^{\Q}[\psi \otimes C^{\epsilon}(\rho_{|G_K})]$ is $\prod_{v}{\varepsilon(\psi^0(\mathfrak{a}_{K_v}) \otimes \Delp{v}{(C^{\epsilon})^{\ast}(\rho_{|G_{K}})},1,\theta_{K_v})}$, where the product is over all places $v$ of $K$. \smallskip

\emph{The infinite place:}

The representation $\psi^0(\mathfrak{a}_{\C}) \otimes \Delp{\infty}{(C^{\epsilon})^{\ast}(\rho_{|G_{K}})}$ of $W_{\C}$ is the sum of $p-1$ copies of the character $z^{-1}$ (for some continuous automorphism $z$ of $\C$) of $W_{\C} \simeq \C^{\times}$, so, by \cite[(3.2.5), (3.6.1)]{NTB}, the local constant is $(-\varepsilon(z^{-1},1,e^{2i\pi\mrm{Tr}_{\C/\R}}))^{\frac{p-1}{2}}=(-i)^{\frac{p-1}{2}}=i(-i)^{\frac{p+1}{2}}=i(-1)^{\frac{p+1}{4}}$. \\

\emph{Split places:}

Let $\ell$ be a prime that splits in $K$ and let $v,v'$ be the two places of $K$ above $\ell$. For $w \in \{w,w'\}$, let $I_w$ be the conductor ideal of $C^{\epsilon}(\rho_{|G_K})$ at $w$. We view it, when relevant, as an ideal of $\OO_K$ which is a power of the maximal ideal defining $w$. The representations $C^{\epsilon}(\rho_{|G_{K_w}})$ and $C^{-\epsilon}(\rho_{|G_{K_w}})$ are dual of each other by Proposition \ref{decomp-cusp-sl}, so they have the same conductor ideal. But $C^{-\epsilon}(\rho_{|G_K})$ is the conjugate by $\mrm{Gal}(K/\Q)$ of $C^{\epsilon}(\rho_{|G_K})$, so that $I_{v'}$ is the image under $\mrm{Gal}(K/\Q)$ of the conductor ideal of $C^{-\epsilon}(\rho_{|G_{K_{v}}})$, which is $I_v$: so $I_{v'}=\tau(I_v)$, so the ideals have the same norm.  

Therefore, $\psi^0(I_v)\psi^0(I_{v'})=\psi^0(I_v\tau(I_v))=(\psi^0\overline{\psi^0})(I_v)=\mathbf{N}I_v$, and 
\begin{align*}
&\prod_{w \in \{v,v'\}}{\varepsilon(\psi^0(\mathfrak{a}_{K_w}) \otimes \Delp{w}{(C^{\epsilon})^{\ast}(\rho_{|G_{K}})},1,\theta_{K_w})}\\
&= \prod_{w \in \{v,v'\}}{\psi(I_w)(\mathbf{N}I_w)^{-1/2}\varepsilon(\Delp{w}{(C^{\epsilon})^{\ast}(\rho_{|G_{K}})},\frac{1}{2},\theta_{K_w})}\\
&= \psi(I_vI_{v'})(\mathbf{N}I_v)^{-1}\prod_{w \in \{v,v'\}}{\varepsilon(\Delp{w}{(C^{\epsilon})^{\ast}(\rho_{|G_{K}})},\frac{1}{2},\theta_{K_w})}=\prod_{w \in \{v,v'\}}{\varepsilon(\Delp{w}{(C^{\epsilon})^{\ast}(\rho_{|G_{K}})},\frac{1}{2},\theta_{K_w})}.
\end{align*}
~\\

\emph{Inert places:}

Let $\ell$ be a prime number which is inert in $K$ and let $\lambda=\ell\OO_K$ be the unique prime ideal above $\ell$. Let $\lambda^{n_{\ell}}$ be the conductor ideal of $C^{\epsilon}(\rho_{G_{K_{\lambda}}})$, which we also view as an id\`ele of $K$ (supported at $\lambda$) when useful. Using the link between conductors and local constants, and the behavior of the latter with respect to induced representations, the conductor exponent of $C(\rho)$ at $\ell$ is then $2n_{\ell}$. By \cite[(3.5)]{Antwerp5-Ribet}, since $f$ has trivial character, $(\psi^0)_{|K_{\lambda}}(\ell)=-\ell$, so
\begin{align*}
&\varepsilon(\psi^0(\mathfrak{a}_{K_{\lambda}}) \otimes \Delp{\lambda}{(C^{\epsilon})^{\ast}(\rho_{|G_{K}})},1,\theta_{K_{\lambda}}) = \psi^0(\lambda)(\ell^2)^{-\frac{n_{\ell}}{2}}\varepsilon(\Delp{\lambda}{(C^{\epsilon})^{\ast}(\rho_{|G_{K}})},\frac{1}{2},\theta_{K_{\lambda}})\\
& = (\psi^0)_{|K_{\lambda}}(\ell)^{n_{\ell}}\ell^{-n_{\ell}}\varepsilon(\Delp{\lambda}{(C^{\epsilon})^{\ast}(\rho_{|G_{K}})},\frac{1}{2},\theta_{K_{\lambda}})=(-1)^{n_{\ell}}\varepsilon(\Delp{\lambda}{(C^{\epsilon})^{\ast}(\rho_{|G_{K}})},\frac{1}{2},\theta_{K_{\lambda}}).
\end{align*}
~\\

\emph{Temporary conclusion:}

The root number is thus 

\[i(-1)^{\frac{p+1}{4}}\frac{\varepsilon(\psi^0(\mathfrak{a}_{\Kpp}) \otimes \Delp{\pp}{(C^{\epsilon})^{\ast}(\rho_{|G_{K}})},1,\theta_{\Kpp})}{\varepsilon(\Delp{\pp}{(C^{\epsilon})^{\ast}(\rho_{|G_{K}})},\frac{1}{2},\theta_{\Kpp})}\prod_{v}{\varepsilon(\Delp{v}{(C^{\epsilon})^{\ast}(\rho_{|G_{K}})},\frac{1}{2},\theta_{K_v})}\prod_{\left(\frac{\ell}{p}\right)=-1}{(-1)^{n_{\ell}}},\]

where the first product, which we denote by $\Pi$, is over the finite places $v$ of $K$. Since $C^{\epsilon}(\rho_{|G_K})$ is an Artin representation, its restriction to $W_{\C}$ is trivial, and thus $\Pi$ is the global root number of (the compatible system attached to) the Artin representation $C^{\epsilon}(\rho_{|G_K})$, hence of (the compatible system attached to) the Artin representation $\mrm{Ind}_K^{\Q}{C^{\epsilon}(\rho_{|G_K})}=V_4(\rho)$ in the notation of Section \ref{splitting-sl2}. Since $V_4$ can be realized over $\R$ by Proposition \ref{cusp-coeffs}, a theorem of Frohlich and Queyrut \cite{sign-ortho} states that $\Pi=1$. 

Hence our global root number is 
\[i(-1)^{\frac{p+1}{4}}\frac{\varepsilon(\psi^0_{\pp}(\mathfrak{a}_{\Kpp}) \otimes \Delp{\pp}{(C^{\epsilon})^{\ast}(\rho_{|G_{K}})},1,\theta_{\Kpp})}{\varepsilon(\Delp{\pp}{(C^{\epsilon})^{\ast}(\rho_{|G_{K}})},\frac{1}{2},\theta_{\Kpp})}\prod_{\left(\frac{\ell}{p}\right)=-1}{(-1)^{\frac{m_{\ell}}{2}}},\]

where $m_{\ell}$ is the conductor exponent at $\ell$ of $V_4(\rho)$. 

Before delving in the computation of the local constants at $\pp$, let us discuss $\psi^0_{\pp} := \psi^0_{|\Kpp}$. This character is trivial over $1+\pp\OO_{\Kpp}$ and is not constant over $\OO_{\Kpp}^{\times}$. Moreover, we know that for $x \in \Kpp^{\times}$, $\psi^0_{\pp}(N_{\Kpp/\Q_p}(x))=\|x\|_{\Kpp}^{-1}$. Since $\psi^0$ is trivial over $K^{\times} \prod_{v \neq \pp}{\OO_{K_v}^{\times}}$, we see that $\psi^0_{\pp}(x)=\sigma(x)$ for $x \in \pm\sqrt{-p}^{\Z}$. 

For $x \in \OO_{\Kpp}^{\times}$, $x \equiv \tau(x) \pmod{\pp}$, so $1=\|x\|_{\Kpp}^{-1}=\psi^0_{\pp}(N_{\Kpp/\Q_p}x)=\psi^0_{\pp}(x^2)$, so $\psi^0_{\pp}$ is a nontrivial quadratic character $(\OO_{\Kpp}/\pp)^{\times} \rar \C^{\times}$, and therefore, for $x \in \OO_{\Kpp}^{\times}$, $\psi^0_{\pp}(x)=\left(\frac{x \pmod{\pp}}{p}\right)$. \\

\emph{When $\rho_{|G_{\Q_p}}$ is split:}

After conjugating $\rho$ by $\SL{\F_p}$, we may assume that there are characters $\alpha,\beta: G_{\Q_p} \rar \F_p^{\times}$ such that $\rho(\sigma)=\begin{pmatrix}\alpha(\sigma) & 0\\0 & \beta(\sigma)\end{pmatrix}$ for every $\sigma \in G_{\Q_p}$. Then $\alpha\beta=\omega_p^{-1}$ and $G_{\Kpp}$ is the kernel of $\alpha\beta: G_{\Q_p} \rar \F_p^{\times}/\F_p^{\times 2}$. Write $\alpha=\alpha_0(\mathfrak{a}_{\Q_p})$ and $\beta=\beta_0(\mathfrak{a}_{\Q_p})$. By Proposition \ref{cusp-sl-diag}, 
\begin{align*}
C^{\epsilon}(\rho_{|G_{\Kpp}}) &\simeq \bigoplus_{\chi: \F_p^{\times 2} \rar \C^{\times}}{\chi \circ \left(\frac{\alpha}{\beta}\right)_{|G_{\Kpp}}}\\
\Delp{\pp}{(C^{\epsilon}(\rho_{|G_{\Kpp}}))^{\ast}} &\simeq \bigoplus_{\chi: \F_p^{\times 2} \rar \C^{\times}}{\chi \circ \left(\frac{\alpha_0}{\beta_0}\right) \circ N_{K/\Q_p}\circ \mathfrak{a}_{\Kpp}}.
\end{align*}
Since $p \equiv 3\pmod{4}$, all the characters in the above sum (except the trivial character) come in pairs $\{\chi,\chi^{-1}\}$. Since
\[\varepsilon(\chi \circ \left(\frac{\alpha_0}{\beta_0}\right) \circ N_{K/\Q_p}\circ \mathfrak{a}_{\Kpp},\frac{1}{2},\theta_{\Kpp})\varepsilon(\chi^{-1} \circ \left(\frac{\alpha_0}{\beta_0}\right) \circ N_{K/\Q_p}\circ \mathfrak{a}_{\Kpp},\frac{1}{2},\theta_{\Kpp})=\chi\circ \frac{\alpha_0}{\beta_0} \circ N_{K/\Q_p}(-1)=1,\] and $\varepsilon(\mathbf{1},\frac{1}{2},\theta_{\Kpp})=1$, we get

\[\varepsilon(\Delp{\pp}{(C^{\epsilon}(\rho_{|G_{\Kpp}}))^{\ast}},\frac{1}{2},\theta_{\Kpp})=1.\]

Let $\chi: \F_p^{\times 2} \rar \C^{\times}$ be a nontrivial character and $\gamma=\chi\left(\frac{\alpha_0}{\beta_0}\right): \Q_p^{\times} \rar \C^{\times}$. Recall that $\mrm{Gal}(\Kpp/\Q_p)=\{\mrm{id},\tau\}$ acts on $W_{\Kpp}$ by conjugation, and that this action preserves the local constants. Then one has  
\begin{align*}
&\varepsilon(\psi^0_{\pp}(\mathfrak{a}_{\Kpp}) \otimes \gamma(N_{K/\Q_p}\mathfrak{a}_{\Kpp}),1,\theta_{\Kpp})\varepsilon(\psi^0_{\pp}(\mathfrak{a}_{\Kpp}) \otimes \gamma^{-1}(N_{K/\Q_p}\mathfrak{a}_{\Kpp}),1,\theta_{\Kpp})\\
&= \varepsilon(\psi^0_{\pp}(\mathfrak{a}_{\Kpp}) \otimes \gamma(N_{K/\Q_p}\mathfrak{a}_{\Kpp}),1,\theta_{\Kpp})\varepsilon(\psi^0_{\pp}(\mathfrak{a}_{\Kpp})\|\cdot\|^1 \otimes \gamma^{-1}(N_{K/\Q_p}\mathfrak{a}_{\Kpp}),0,\theta_{\Kpp})\\
&= \varepsilon(\psi^0_{\pp}(\mathfrak{a}_{\Kpp}) \otimes \gamma(N_{K/\Q_p}\mathfrak{a}_{\Kpp}),1,\theta_{\Kpp})\varepsilon((\overline{\psi}^0_{\pp})^{-1}(\mathfrak{a}_{\Kpp}) \otimes \gamma^{-1}(N_{K/\Q_p}\mathfrak{a}_{\Kpp}),0,\theta_{\Kpp})\\
&= \varepsilon(\psi^0_{\pp}(\mathfrak{a}_{\Kpp}) \otimes \gamma(N_{K/\Q_p}\mathfrak{a}_{\Kpp}),1,\theta_{\Kpp})\varepsilon((\psi^0_{\pp})^{-1}(\mathfrak{a}_{\Kpp})\circ \tau \otimes \gamma^{-1}(N_{K/\Q_p}\mathfrak{a}_{\Kpp})\circ \tau,0,\theta_{\Kpp}\circ \tau)\\
&= \varepsilon(\psi^0_{\pp}(\mathfrak{a}_{\Kpp}) \otimes \gamma(N_{K/\Q_p}\mathfrak{a}_{\Kpp}),1,\theta_{\Kpp})\varepsilon((\psi^0_{\pp})^{-1}(\mathfrak{a}_{\Kpp}) \otimes \gamma^{-1}(N_{K/\Q_p}\mathfrak{a}_{\Kpp}),0,\theta)\\
&=\psi^0_{\pp}(-1)\gamma(N_{\Kpp/\Q_p}(-1))=-1.
\end{align*}

Since there are $\frac{1}{2}\left(\frac{p-1}{2}-1\right)=\frac{p-3}{4}$ pairs $\{\chi,\chi^{-1}\}$ of characters of $\F_p^{\times 2}$, it follows by the classical computation of the Gauss sum (Proposition \ref{gauss-sum}) 
\begin{align*}
&\varepsilon(\psi^0_{\pp}(\mathfrak{a}_{\Kpp}) \otimes \Delp{\pp}{(C^{\epsilon}(\rho_{|G_{\Kpp}}))^{\ast}},1,\theta_{\Kpp}) = (-1)^{\frac{p-3}{4}}\varepsilon(\psi^0_{\pp}(\mathfrak{a}_{\Kpp})\|\cdot\|^{1/2},\frac{1}{2},\theta_{\Kpp})\\
&= (-1)^{\frac{p-3}{4}}\times \sigma(p)(p^{-2})^{\frac{1}{2}}\frac{\sum_{x \in (\OO_{\Kpp}/\pp)^{\times}}{(\psi^0_{\pp})^{-1}(x)\theta_{\Kpp}\left(\frac{x}{p}\right)}}{\left|\sum_{x \in (\OO_{\Kpp}/\pp)^{\times}}{(\psi^0_{\pp})^{-1}(x)\theta_{\Kpp}\left(\frac{x}{p}\right)}\right|} = (-1)^{\frac{p-3}{4}}\frac{\sum_{x \in \F_p^{\times}}{\left(\frac{x}{p}\right)e^{\frac{4i\pi}{p}x}}}{\left|\sum_{x \in \F_p^{\times}}{\left(\frac{x}{p}\right)e^{\frac{4i\pi}{p}x}}\right|}\\
&= (-1)^{\frac{p-3}{4}}\left(\frac{2}{p}\right)i=-i,
\end{align*}
whence the conclusion. \\

\emph{When $\rho_{|G_{\Q_p}}$ is $a$-Cartan:}

Let $D \leq \GL{\F_p}$ be a nonsplit Cartan subgroup such that $\rho(I_p)=\rho(G_E)=aD$ (by Lemma \ref{expand-cartan}) and let $N$ be the normalizer of $D$. Let $D_2=D \cap \F_p^{\times}\SL{\F_p}$ and $N_2=N \cap \F_p^{\times}\SL{\F_p}$. Note that $a$ is an odd divisor of $\frac{p+1}{2}$, so $a\mid \frac{p+1}{4}$. 

Let $\eta: \mrm{Gal}(\Kpp E/\Kpp) \rar \{\pm 1\}$ be a character, trivial if and only if $p \equiv 3\pmod{8}$. Let $X$ be a set of characters $\chi: D_2/\F_p^{\times} \rar \C^{\times}$ with order at least $3$ such that for any $\chi: D_2/\F_p^{\times} \rar \C^{\times}$, $X$ contains exactly one of the characters $\chi$ and $\chi^p=\chi^{-1}$. By Proposition \ref{cusp-sl-norm-cartan}, 
\[C^{\epsilon}(\rho_{|G_{\Kpp}}) \simeq \eta \oplus \bigoplus_{\chi \in X}{\mrm{Ind}_{\Kpp E}^{\Kpp}{\left[\chi \circ \rho_{|G_{\Kpp E}}\right]}}.\]

Since $\theta_{\Kpp}$ is trivial on $\sqrt{-p}^{-1}\OO_{\Kpp}$ but not on $p^{-1}\OO_{\Kpp}$, and $\psi^0_{\pp}$ has conductor ideal $\mathfrak{p}$, one has 
\begin{align*}
\varepsilon(\overline{\eta},\frac{1}{2},\theta_{\Kpp})&=\overline{\eta}(\mathfrak{a}_{\Kpp}^{-1}(\sqrt{-p}))=(-1)^{\frac{p-3}{4}}\\
\varepsilon(\psi^0_{\pp}(\mathfrak{a}_{\Kpp}) \otimes \overline{\eta},\frac{1}{2},\theta_{\Kpp})&=\overline{\eta}(\mathfrak{a}_{\Kpp}^{-1}(\sqrt{-p}^{1+1}))\varepsilon(\psi^0_{\pp}(\mathfrak{a}_{\Kpp})\|\cdot\|^{1/2},\frac{1}{2},\theta_{\Kpp})\\
&= \psi^0_{\pp}(p)\|p\|^{1/2}\frac{\sum_{x \in (\OO_{\Kpp}/\pp)^{\times}}{(\psi^0_{\pp})^{-1}(u)\theta_{\Kpp}\left(\frac{u}{p}\right)}}{\left|\sum_{x \in (\OO_{\Kpp}/\pp)^{\times}}{(\psi^0_{\pp})^{-1}(u)\theta_{\Kpp}\left(\frac{u}{p}\right)}\right|}=i\left(\frac{2}{p}\right)=(-1)^{\frac{p+1}{4}}i.
\end{align*}

Write $\rho_{|W_E}=\sigma \circ \mathfrak{a}_E$ for some character $\sigma: E^{\times}/(1+p\OO_E) \rar D$ with image $aD$, so that one has $\sigma(p)=-I_2$ and, for each $x \in \F_p^{\times}$, $\sigma(x)=x^{-1}I_2$. 
Let $\chi \in X$ and define
\begin{align*}
\varepsilon_{\chi}&:=\frac{\varepsilon(\psi^0_{\pp}(\mathfrak{a}_{\Kpp}) \otimes \Delp{\pp}{\mrm{Ind}_{\Kpp E}^{\Kpp}{\left[\chi \circ \rho_{|G_{\Kpp E}}\right]}^{\ast}},1,\theta_{\Kpp})}{\varepsilon(\Delp{\pp}{\mrm{Ind}_{\Kpp E}^{\Kpp}{\left[\chi \circ \rho_{|G_{\Kpp E}}\right]}^{\ast}},1,\theta_{\Kpp})}\\
&= \frac{\varepsilon((\psi^0_{\pp} \circ N_{\Kpp E/\Kpp})(\mathfrak{a}_{\Kpp E}) \otimes (\overline{\chi}\circ \rho_{|W_{\Kpp E}}),1,\theta_{\Kpp E})}{\varepsilon(\overline{\chi}\circ \rho_{|W_{\Kpp E}},\frac{1}{2},\theta_{\Kpp E})}\\
&= \frac{\varepsilon(\psi^0_{\pp}(N_{\Kpp E/\Kpp}) \cdot \overline{\chi}(\sigma \circ N_{\Kpp E/E}),1,\theta_{\Kpp E})}{\varepsilon(\overline{\chi}(\sigma \circ N_{\Kpp E/E}),\frac{1}{2},\theta_{\Kpp E})}
\end{align*} 

Since $\Kpp E/\Kpp$ is unramified, $\psi^0_{\pp}(N_{\Kpp E/\Kpp})$ is the only character of order two of $(\OO_{\Kpp E}/\pp)^{\times}$. Since $N_{\Kpp E/E}((\Kpp E)^{\times})$ is the subgroup generated by $p$ and the squares of $\OO_E^{\times}$, one has \[(\sigma\circ N_{\Kpp E/E})(\OO_{\Kpp E}^{\times})=(\sigma \circ N_{\Kpp E/E})((\Kpp E)^{\times})=aD_2.\] 

We have to distinguish three cases: 
\begin{enumerate}[label=(\alph*),noitemsep]
\item\label{cmc1} $\overline{\chi}(\sigma \circ N_{\Kpp E/E})$ is trivial on the units, hence trivial. This is equivalent to $\chi(aD_2)=\{1\}$, or equivalently to $\chi^a=1$. Since $a$ is odd, there are $\frac{a-1}{2}$ such characters. 
\item\label{cmc2} $(\overline{\chi}\circ\sigma \circ N_{\Kpp E/E})(\OO_{\Kpp E}^{\times})=\{\pm 1\}$: in this case, the restrictions to $\OO_{\Kpp E}^{\times}$ of $\psi^0_{\pp}(N_{\Kpp E/\Kpp})$ and $(\overline{\chi}\circ\sigma \circ N_{\Kpp E/E})$ are two nontrivial characters $(\OO_{E \Kpp}/\pp)^{\times} \rar \{\pm 1\}$, so they are equal, and their product is the unramified character mapping the uniformizer $\sqrt{-p}$ to $-p$. This case occurs whenever $\chi^a \neq 1$ but $\chi^{2a}=1$. There are $\frac{2a-2}{2}-\frac{a-1}{2}=\frac{a-1}{2}$ such characters. 
\item\label{cmc3} In the remaining $\frac{p-3}{4}-(a-1)=\frac{p+1}{4}-a$ cases, $\psi^0_{\pp}(N_{\Kpp E/\Kpp}) \cdot \overline{\chi}(\sigma \circ N_{\Kpp E/E})$ and $\overline{\chi}(\sigma \circ N_{\Kpp E/E})$ are both characters $(\Kpp E)^{\times}/(1+\pp\OO_{\Kpp E}) \rar \C^{\times}$, and neither vanishes on $\OO_{\Kpp E}^{\times}$. 
\end{enumerate} 

We let $k=\OO_E/p$ and $T,\nu: k \rar \F_p$ denote the trace and norm, respectively. Let \[\alpha: (\Kpp E)^{\times}/\Z_p^{\times}(1+\pp\OO_{\Kpp E}) \rar \C^{\times}\] denote any ramified character. Since $-1$ is not a square in $\F_p^{\times}$, $\theta_{\Kpp E}$ is trivial on $\sqrt{-p}^{-1}\OO_{\Kpp E}$ but not on $p^{-1}\OO_{\Kpp E}$, we can apply Lemma \ref{local-epsilon-fp2}, so that 
\begin{align*}
\varepsilon(\alpha,s,\theta_{\Kpp E})&=\varepsilon(\alpha\|\cdot\|^{s-\frac{1}{2}},\frac{1}{2},\theta_{\Kpp E})=\alpha(p)\|p\|^{s-\frac{1}{2}}\frac{\sum_{x \in (\OO_{E\Kpp}/\pp)^{\times}}{\alpha^{-1}(x)\theta_{\Kpp E}(x/p)}}{\left|\sum_{x \in (\OO_{E\Kpp}/\pp)^{\times}}{\alpha^{-1}(x)\theta_{\Kpp E}(x/p)}\right|}\\
&=\alpha(p)p^{-4s+2}\frac{\sum_{x \in k^{\times}}{\alpha^{-1}(x)e^{\frac{4i\pi}{p}T(x)}}}{\left|\sum_{x \in k^{\times}}{\alpha^{-1}(x)e^{\frac{4i\pi}{p}T(x)}}\right|}\\
&=p^{-4s+2}\alpha(p)\frac{\sum_{x \in k^{\times}}{\alpha^{-1}(x)e^{\frac{2i\pi}{p}T(x)}}}{\left|\sum_{x \in k^{\times}}{\alpha^{-1}(x)e^{\frac{2i\pi}{p}T(x)}}\right|}\\
&= p^{-4s+2}\alpha(p)\alpha^{-1}(\sqrt{-1})=p^{-4s+2}\alpha(p\sqrt{-1}).
\end{align*}
When $\alpha$ is instead unramified, one has
\[\varepsilon(\alpha,s,\theta_{\Kpp E})=\varepsilon(\alpha\|\cdot\|^{s-\frac{1}{2}},\frac{1}{2},\theta_{\Kpp E})=\alpha(\sqrt{-p})p^{1-2s}.\]

Note that in all cases, the kernels of $\psi^0_{\pp}(N_{\Kpp E/\Kpp})$ and $(\overline{\chi}\circ\sigma \circ N_{\Kpp E/E})$ both contain $\sqrt{-1}^{\Z}\Z_p^{\times}(1+\pp\OO_{\Kpp E})$. Moreover, $(\overline{\chi}\circ\sigma \circ N_{\Kpp E/E})(p)=\overline{\chi}(\sigma(p^2))=1$.

Therefore, in all three cases, $\varepsilon_{\chi}=\varepsilon(\psi^0_{\pp}(N_{\Kpp E/\Kpp}) \cdot \overline{\chi}(\sigma \circ N_{\Kpp E/E}),1,\theta_{\Kpp E})$.   

\ref{cmc1}: in this case, $\varepsilon_{\chi}=\varepsilon(\psi^0_{\pp}(N_{\Kpp E/\Kpp}),1,\theta_{\Kpp E})=\psi^0_{\pp}(N_{\Kpp E/\Kpp}(p\sqrt{-1}))p^{-2}=1$.\\

\ref{cmc2}: in this case, $\varepsilon_{\chi} = \varepsilon(\sqrt{-p} \longmapsto -p,1,\theta_{\Kpp E})=-pp^{1-2}=-1$.\\

\ref{cmc3}: in this case, \[\varepsilon_{\chi}= \varepsilon(\psi^0_{\pp}(N_{\Kpp E/\Kpp}) \cdot \overline{\chi}(\sigma \circ N_{\Kpp E/E}),1,\theta_{\Kpp E})=p^{-4+2}\psi^0_{\pp}(N_{\Kpp E/\Kpp}(p))=p^{-2}\psi^0_{\pp}(p^2)=1.\] 

Finally, the global root number is $i(-1)^{\frac{p+1}{4}}U\frac{(-1)^{\frac{p+1}{4}}i}{(-1)^{\frac{p-3}{4}}}(-1)^{\frac{a-1}{2}}=(-1)^{\frac{p+2a-1}{4}}U$. \\

\emph{When $\rho_{|G_{\Q_p}}$ is wild:}

Let us discuss why the result is well-defined. There are three aspects: we can change the decomposition subgroup at $p$, we can conjugate $\rho_{|G_{\Q_p}}$ by some matrix such that it remains of the form $\beta\begin{pmatrix} q & \nu\\0 & 1\end{pmatrix}$, and we can change $\gamma$. We first explain why the second does not matter, then discuss the first change. The independence of the result with respect to the choice of $\gamma$ will be dealt with at the end of the proof. 
 
Indeed, let $M \in \GL{\F_p}$ be such that $M\left(\beta \begin{pmatrix}q & \nu\\0 & 1\end{pmatrix}\right)M^{-1}=\beta'\begin{pmatrix} q' & \nu' \\0 & 1\end{pmatrix}$. By Lemma \ref{expand-wild}, we see that $(\beta,q)=(\beta',q')$ and that $\nu'_{|G_L}=r\nu_{G_L}$, where $r=\frac{M_{1,1}}{M_{2,2}} \in \det{M}\F_p^{\times 2}$.  

Let $\gamma\in L^{\times},\gamma_0 \in \Kpp^{\times}$ be attached to $\nu_{|G_L}$. Fix some $r_0 \in \Z_p^{\times}$ lifting $r$, then $\gamma'=r_0\gamma$ works for $\nu'_{|G_L}$, and thus we may take $\gamma'_0=r_0^{[L:\Kpp]}\gamma_0$. Since $v_2([L:\Kpp])=v_2([L:\Q_p])-1=v_2(p-1)-1=0$, the class of $r_0^{[L:\Kpp]}$ in $\F_p^{\times}$ is a square if and only if $r$ is, if and only if $\det{M} \in \F_p^{\times 2}$. 

In particular, $\psi^0_{\pp}(\gamma'_0)\|\gamma'_0\|^{1/2}=\left(\frac{\det{M}}{p}\right)\psi^0_{\pp}(\gamma_0)\|\gamma_0\|^{1/2}$. This is what we wanted to prove, because $C^{\epsilon}(M\rho_{|G_K} M^{-1})=C^{\left(\frac{\det{M}}{p}\right)\epsilon}(\rho_{|G_K})$ by Proposition \ref{decomp-cusp-sl}.

Now, we discuss the behavior with respect to changing the embedding $\Qbar \rar \overline{\Q_p}$. Let $\sigma,\sigma': \Qbar \rar \overline{G_{\Q_p}}$ be two embeddings inducing homomorphisms $R_{\sigma}, R_{\sigma'}: G_{\Q_p} \rar G_{\Q}$. Then $\sigma'=\sigma \circ \alpha$ for some $\alpha \in G_{\Q}$, hence $R_{\sigma'}(t)=\alpha^{-1}R_{\sigma}(t)\alpha$ for every $t \in G_{\Q_p}$. The morphism $\rho \circ R_{\sigma}$ (resp. $\rho \circ R_{\sigma'}$) is $\SL{\F_p}$-conjugate to $\beta_{\sigma}\begin{pmatrix}q_{\sigma}&\nu_{\sigma}\\0 & 1\end{pmatrix}$ (resp. $\beta_{\sigma'}\begin{pmatrix}q_{\sigma'}&\nu_{\sigma'}\\0 & 1\end{pmatrix}$). 

Therefore, the morphisms $\beta_{\sigma}\begin{pmatrix}q_{\sigma}&\nu_{\sigma}\\0 & 1\end{pmatrix},\,\beta_{\sigma'}\begin{pmatrix}q_{\sigma'}&\nu_{\sigma'}\\0 & 1\end{pmatrix}: G_{\Q_p} \rar \GL{\F_p}$ are $\SL{\F_p}$-conjugate if $\det{\rho(\alpha)}\in \F_p^{\times 2}$ and conjugate by a matrix in $\GL{\F_p}\backslash \F_p^{\times 2}\SL{\F_p}$ otherwise. Note that $\det{\rho}(\alpha) \in \F_p^{\times 2}$ if and only if $\omega_p(\alpha) \in \F_p^{\times 2}$, if and only if $\alpha \in G_K$, that is, if and only if $\sigma_{|K}=\sigma'_{|K}$. 

By Lemma \ref{expand-wild}, $(\beta_{\sigma},q_{\sigma})=(\beta_{\sigma'},q_{\sigma'})$, the extension $L$ is the same in both cases, and $(\nu_{\sigma'})_{|G_L} = r(\nu_{\sigma})_{|G_L}$ for some $r \in \Z_p^{\times}$ which is a square if and only if $\omega_p(\alpha)$ is one.

Moreover, using the material of \cite[\S 1]{NTB}, one sees that there are morphisms $\tilde{\sigma},\tilde{\sigma'}: W_{\Kpp} \rar W_K$, and an automorphism $\tilde{\alpha}$ of $W_K$ such that the following diagram commutes, where the unmarked maps are the natural ones:
\[
\begin{tikzcd}[ampersand replacement=\&]
\mathbb{A}_K^{\times}/K^{\times} \arrow{dd}{\alpha}\& \& W_K^{ab} \arrow{dd}\arrow{ll}{\mathfrak{a}_K} \& \& W_K \arrow{rr}\arrow[near end]{dd}{\tilde{\alpha}}\arrow{ll} \& \& G_K \arrow[bend left=90,swap,looseness=2]{dd}{\alpha \ast \alpha^{-1}}\\
\& \Kpp^{\times}\arrow{ul}{\sigma'_{\pp}}\arrow{dl}{\sigma_{\pp}} \& \& W_{\Kpp}^{ab} \arrow[near end]{ll}{\mathfrak{a}_{\Kpp}} \arrow{ul}\arrow{dl}\& \& W_{\Kpp} \arrow{ll}\arrow{ul}{\tilde{\sigma'}}\arrow{dl}{\tilde{\sigma}} \arrow{r} \& G_{\Kpp} \arrow{u}{R_{\sigma'}}\arrow{d}{R_{\sigma}}\\
\mathbb{A}_K^{\times}/K^{\times} \&\& W_K^{ab} \arrow{ll}{\mathfrak{a}_K} \& \& W_K \arrow{ll}\arrow{rr}\& \& G_K 
\end{tikzcd}
\]

The claim we want to prove is thus $\psi^0(\sigma_{\pp}(\gamma_0))=\psi^0(\sigma'_{\pp}(\gamma'_0))$, where $\gamma \in L^{\times}, \gamma_0 \in \Kpp^{\times}$ are attached to $(\nu_{\sigma})_{|G_L}$ and $\gamma' = r\gamma, \gamma'_0=r^{[L:\Kpp]}\gamma_0 \in \Kpp$ are attached to $(\nu_{\sigma'})_{|G_L}$. As above, $[L:\Kpp]$ is odd, so $r^{[L:\Kpp]}$ is a square if and only if $\omega_p(\alpha) \in \F_p^{\times 2}$. Since $\gamma_0$ has odd valuation (as we will see later in the proof), 
\begin{align*}
\psi^0(\sigma_{\pp}(\gamma_0)) &= \psi^0\circ \alpha \circ \sigma'_{\pp}(r^{-[L:\Kpp]\gamma'_0})=\psi^0_{\pp}(r)^{-[L:\Kpp]}\psi^0_{\pp}(\sigma'_{\pp}(\gamma'_0))\frac{\overline{\psi^0_{\pp}}}{\psi^0_{\pp}}(\sigma'_{\pp}(\gamma'_0))\\
&=\left(\frac{\omega_p(\alpha)}{p}\right)\psi^0(\sigma'_{\pp}(\gamma'_0))\left(\frac{\omega_p(\alpha)}{p}\right),
\end{align*}
so we are done.

Now, we can compute the sign after fixing a decomposition subgroup and conjugating $\rho$ by a matrix in $\SL{\F_p}$ so that it has the form $\beta\begin{pmatrix} q & \nu\\0& 1\end{pmatrix}$.  

After conjugating by a matrix in $\SL{\F_p}$, we may assume that $\rho_{|G_{\Q_p}}=\beta \otimes \begin{pmatrix}q & \nu\\0 & 1\end{pmatrix}$ with $\beta,q: G_{\Q_p} \rar \F_p^{\times}$ characters and $\nu: G_{\Q_p} \rar \F_p(q)$ a cocycle. Define $L,\nu_0,n,\gamma,\zeta,e,\omega,\varpi$ as in the case ``$(f,\mathbf{1}) \in \mathscr{C}'_1$ and $\rho_{|G_{\Q_p}}$ wild''. The maximal unramified subextension of $K$ is denoted by $L_0$ to avoid a conflict of notation. Since $q\beta^2=\omega_p^{-1}$, one has $q^{-1}(\F_p^{\times 2})=\omega_p^{-1}(\F_p^{\times 2})=G_{\Kpp}$, so $G_{\Kpp}$ is contained in $L$ and $[\F_p^{\times 2}:q(G_{\Kpp})]=[\F_p^{\times}:q(G_{\Q_p})]$.  

By Proposition \ref{cusp-sl-borel} and Mackey's formula, 
\[\Delp{\pp}{C^{\epsilon}(\rho_{|G_{\Kpp}})^{\ast}} \simeq \bigoplus_{t \in \F_p^{\times 2}/q(G_{\Kpp})}{\mrm{Ind}_L^{\Kpp}\left[e^{\frac{-2i\pi t\epsilon}{p}\nu_0} \circ \mathfrak{a}_L\right]}.\]
Therefore, 
\begin{align*}
\varepsilon_q &:= \frac{\varepsilon(\psi^0_{\pp}(\mathfrak{a}_{\Kpp}) \otimes \Delp{\pp}{(C^{\epsilon})^{\ast}(\rho_{|G_{K}})},1,\theta_{\Kpp})}{\varepsilon(\Delp{\pp}{(C^{\epsilon})^{\ast}(\rho_{|G_{K}})},\frac{1}{2},\theta_{\Kpp})} = \prod_{t \in \F_p^{\times 2}/q(G_{\Kpp})}{\frac{\varepsilon((\psi^0_{\pp}(N_{L/\Kpp})e^{\frac{-2i\pi t\epsilon}{p}\nu_0})\circ \mathfrak{a}_L,1,\theta_L)}{\varepsilon(e^{\frac{-2i\pi t\epsilon}{p}\nu_0}\circ \mathfrak{a}_L,\frac{1}{2},\theta_L)}}\\
&= \prod_{t \in \F_p^{\times 2}/q(G_{\Kpp})}{\frac{\varepsilon((\psi^0_{\pp}(N_{L/\Kpp})e^{\frac{-2i\pi t\epsilon}{p}\nu_0})\|\cdot\|^{\frac{1}{2}},\frac{1}{2},\theta_L)}{\varepsilon(e^{\frac{-2i\pi t\epsilon}{p}\nu_0},\frac{1}{2},\theta_L)}}.
\end{align*}

As in the $\mathscr{C}'_1$ case, $\psi^0_{\pp}(N_{L/\Kpp})$ is trivial on $1+\varpi\OO_L$. Since $[L:\Kpp]$ divides $\frac{[L:\Q_p]}{2} \mid \frac{p-1}{2}$ which is odd, $N_{L/\Kpp}(\OO_L^{\times})$ is a subgroup of $\OO_K^{\times}$ with odd index, so it contains an element which is not a square, and $\psi^0_{\pp}(N_{L/\Kpp})$ does not vanish on $\OO_L^{\times}$. So $\psi^0_{\pp}(N_{L/\Kpp})$ has level zero in the sense of \cite[\S 1.8]{BH}. On the other hand, we saw in Lemma \ref{expand-wild} that $\nu_0$ had positive level $n-1$ in the sense of \cite[\S 1.8]{BH}.  

It is clear from the definition of $\gamma$ that whenever $x \in \varpi^{n/2}\OO_L$, $e^{\frac{-2i\pi t\epsilon}{p}\nu_0}(1+x)=\theta_L\left(\frac{-t\epsilon\gamma}{p}x\right)$, so, by the stability theorem \cite[\S 23.8]{BH}, one has 
\begin{align*}
\varepsilon_q^{-1} &= \prod_{t \in \F_p^{\times 2}/q(G_{\Kpp})}{\psi^0_{\pp}(N_{L/\Kpp}\frac{-t\epsilon\gamma}{p})\|\frac{-t\epsilon\gamma}{p}\|^{1/2}} =  \left(-\epsilon\psi^0_{\pp}(N_{L/\Kpp}\gamma)\|\gamma\|^{1/2}\right)^{[\F_p^{\times 2}:q(G_{\Kpp})]},\end{align*}

Let $\gamma_0=N_{L/\Kpp}(\gamma) \in \pp^{f_{L/\Kpp}(1-n)}\OO_{\Kpp}^{\times}$ with $\gamma_0 \in L^{\times}$. Since \[f_{L/\Kpp}(1-n)[\F_p^{\times 2}:q(G_{\Kpp})]=\frac{(p-1)(1-n)}{e}\] is an odd integer, it follows that \[\varepsilon_q^{-1}=-\epsilon (-1)^{\frac{[\F_p^{\times 2}:q(G_{\Kpp})]-1}{2}}[\psi^0_{\pp}(\gamma_0)\|\gamma_0\|^{1/2}] \in \{\pm i\}.\] Note that $[\F_p^{\times 2}:q(G_{\Kpp})]=[\F_p^{\times}:q(G_{\Q_p})]=\frac{p-1}{[L:\Q_p]}$ and by Lemma \ref{expand-wild}, $p-1 \equiv [L:\Q_p] \equiv 2\pmod{4}$. 

Let $d=[L:\Q_p]$. The global sign is thus \[\epsilon\left(\prod_{\left(\frac{\ell}{p}\right)=-1}{(-1)^{n_{\ell}}}\right)(-1)^{\frac{p+1}{4}+\frac{p-1-d}{2d}}\left[i\psi^0_{\pp}(\gamma_0)\|\gamma_0\|^{1/2}\right].\]

Now, since $d \equiv p-1 \equiv 2\pmod{4}$,
\begin{align*}
\frac{p+1}{4}+\frac{p-1-d}{2d}&=\frac{\frac{d}{2}(p+1)+p-1-d}{2d}=\frac{\left(\frac{d}{2}+1\right)(p-1)}{2d}\\
&=\frac{(p-1)\left(\frac{d}{2}+1\right)}{4\frac{d}{2}}\equiv \frac{(p-1)\left(\frac{d}{2}+1\right)}{4} \pmod{2}\\
&\equiv \frac{p-1}{2}\frac{d+2}{4} \equiv \frac{d+2}{4} \pmod{2}.  
\end{align*}

In the above, we saw that $\gamma_0$ had indeed odd valuation. Moreover, if we replace $\gamma$ with another allowed $\gamma' \in L^{\times}$, then $\gamma'\in (1+\varpi^{n/2}\OO_L)\gamma$, so $N_{L/\Kpp}\left(\frac{\gamma'}{\gamma}\right) \in 1+\pp\OO_{\Kpp} \subset \ker{\psi^0_{\pp}}$.

\endgroup
}

\rem{It is surprising and not obvious \emph{a priori}, that, except for the CM factors, the global root number for $R_f(\rho)$ only depends from $a_p(f)$ (which, by Proposition \ref{bad-L-factor}, encodes whether $(f,\chi)$ is in $\mathscr{S}$, $\mathscr{P}'$, or $\mathscr{C}$), on the restriction of $\rho$ to an inertia group at $p$, and on $p\pmod{4}$. While it seems natural \emph{a posteriori} that $a_p(f)$ and the restriction of $\rho$ to an inertia group at $p$ contribute to the sign, the relevance of $p\pmod{4}$ seems less easy to explain. }

\bigskip

\section{Consequences for elliptic curves}
\label{asymptotic-conseq-elliptic-curves}

In this Section, we fix an elliptic curve $E/\Q$, and consider, for any prime $p \geq 7$, the representation $\rho: G_{\Q} \rar \GL{\F_p}$ given by $\rho(\sigma)\begin{pmatrix}\sigma(P)\\\sigma(Q)\end{pmatrix}= \begin{pmatrix}P\\Q\end{pmatrix}$ for some basis $(P,Q)$ of $E[p]$ (as in Proposition \ref{group-is-twist-unpolarized}). By the results of \cite{Serre-image-ouverte}, we know that either $E$ has complex multiplication, or, for all but finitely many $p$, $\rho(G_{\Q})=\GL{\F_p}$. This is another reason why Proposition \ref{decomposition-xrho-connected-surjective} and the sign computation of the previous section are important: they describe, in some sense, the situation for all but finitely many $p$. 

The following result is a standard estimate of the cardinalities of $\mathscr{S},\mathscr{P},\mathscr{C}$:

\lem[spc-size]{
The following estimates fold, for any $p \geq 5$ and any $\epsilon\in \{\pm 1\}$:
\begin{align*}
&|\mathscr{S}| \in \frac{p+1}{12} + \left[-\frac{7}{6},0\right],\,&&|\mathscr{P}| \in \frac{(p-11)(p-3)}{24}+\left[0,\frac{7}{6}\right],\\
&\Bigg|\frac{1}{2}|\mathscr{S}|-\Big|\{(f,\chi) \in \mathscr{S}\mid a_p(f)=\epsilon\}\Big|\Bigg| \leq \frac{3\sqrt{p}\ln{11e^2 p}}{2\pi},\, &&|\mathscr{P}'| =\frac{1}{2}|\mathscr{P}|,\\
&|\mathscr{C}| \in \frac{(p-1)(p^2-2p-27)}{48}+\frac{p-1}{2}\left[0,\frac{7}{6}\right],\,&& |\mathscr{C}'_1|+\frac{1}{2}|\mathscr{C}_M| \in \frac{p^2-2p-27}{48}+\left[0,\frac{7}{6}\right],\\
&|\mathscr{C}_M| \in \{0,h(\Q(\sqrt{-p}))\},\,&& |\mathscr{C}_M| \leq \frac{\sqrt{p}\ln{4e^2p}}{\pi}.
\end{align*}}

\demo{$|\mathscr{S}|$ is the genus of $X_0(p)$, which is known to be in $\frac{p+1}{2}-\left[0,\frac{7}{6}\right]$ (for instance \cite[Exercise 3.1.4]{DS}). Furthermore, $|\mathscr{S}|+|\mathscr{P}|$ is the genus of $X_1(p)$, which is given in \cite[Figure 3.4]{DS}, and $|\mathscr{P}'|=\frac{1}{2}|\mathscr{P}|$. Moreover, by Corollaries \ref{differentials-xpp} and \ref{decomposition-hecke-xpp}, one has 
\begin{align*}
&(p-1)\left(|\mathscr{C}|+\frac{p-1}{2}\mathscr{S}+\frac{p+1}{2}\dim_{\C}{\mathcal{S}_2(\Gamma_1(p))}\right)\\
&=p(p-1)|\mathscr{S}+(p+1)(p-1)|\mathscr{P}'|+(p-1)|\mathscr{C}|=(p-1)\dim_{\C}{\mathcal{S}_2(\Gamma(p))},
\end{align*}
 which yields successive estimates for $\mathscr{P},\mathscr{P}',\mathscr{C}$. 

Next, by considering the possible twists, $|\mathscr{C}|=\frac{p-1}{2}|\mathscr{C}_M|+(p-1)|\mathscr{C}'_1|$. So all we have to do is compute $|\mathscr{C}_M|$. We know that it is zero if $p \equiv 1\pmod{4}$, so we assume $p \equiv 3 \pmod{4}$. By \cite[Theorem 4.5]{Antwerp5-Ribet}, $|\mathscr{C}_M|$ is exactly the number of continuous characters $\mathbb{A}_K^{\times}/K^{\times} \rar \C^{\times}$ (with $K:=\Q(\sqrt{-p})$) that are equal to $z \longmapsto \sigma(z)^{-1}$ for $z \in (K \otimes \R)^{\times}$ and some fixed isomorphism $\sigma: K \otimes \R \rar \C$, are trivial on $(1+\sqrt{-p}\OO_{K_p})\prod_{v \nmid \sqrt{-p}\infty}{\OO_{K_v}^{\times}}$, and are given on $\OO_{K_p}$ (as in the proof of Theorem \ref{root-numbers-surjective-computation}) by the unique nontrivial character $(\OO_{K_p}/\sqrt{-p})^{\times} \rar \{\pm 1\}$. 

If there exists such a character (that is, if $|\mathscr{C}_M| > 0$), then $|\mathscr{C}_M|$ is exactly the number of characters $\mathbb{A}_K^{\times}/\C^{\times}K^{\times}\prod_{v}{\OO_{K_v}^{\times}} \rar \C^{\times}$, therefore $|\mathscr{C}_M|$ is the class number $h(\Q(\sqrt{-p}))$ of $\Q(\sqrt{-p})$. By \cite[Lemma 4.3]{QC8+}, one has $h(\Q(\sqrt{-p})) \leq \frac{\sqrt{p}\ln{4e^2p}}{\pi}$. \\

Since $|\mathscr{S}|=\sum_{\epsilon}{|\{(f,\chi) \in \mathscr{S}\mid a_p(f)=\epsilon\}|}$, it is enough to conclude to prove that 
\[\frac{1}{2}|\mathscr{S}| - \frac{3\sqrt{p}\ln{11e^2p}}{2\pi} \leq |\{(f,\chi) \in \mathscr{S}\mid a_p(f)=-1\}| \leq \frac{1}{2}|\mathscr{S}|.\]

We mimic the argument of \cite[Section 4]{QC8+}. By Proposition \ref{local-constant-principal}, the cardinality of the set $\{(f,\chi) \in \mathscr{S}\mid a_p(f)=-1\}$ is equal to the dimension of the subspace of $\mathcal{S}_2(\Gamma_0(p))$ fixed under the Atkin-Lehner involution $w_p$, hence is equal to the genus of the quotient $X_0^+(p)$ of $X_0(p)$ by $w_p$. Note that if $p \geq 5$ is congruent to $1$ modulo $4$, $\Q(\sqrt{-p})$ never has class number one by genus theory (since $\Q(i,\sqrt{p})$ is a quadratic extension of $\Q(\sqrt{-p})$ which is everywhere unramified). By applying the Riemann-Hurwitz formula as in Proposition 4.4 of \emph{loc.cit.}, we find using Theorem 4.2 and Lemma 4.3 of \emph{loc.cit.} that
\[2|\mathscr{S}|-2 = 2(2|\mathscr{S}^-|-2)+|\mrm{Fix}(w_p)| \in 4|\mathscr{S}^-|-2+\left[0,\frac{\sqrt{p}\ln{4e^2p}}{\pi}+\frac{\sqrt{4p}\ln{16e^2p}}{\pi}\right].\]
The conclusion follows since 
\begin{align*}
\frac{\sqrt{p}\ln{4e^2p}}{\pi}+\frac{\sqrt{4p}\ln{16e^2p}}{\pi} &\leq \frac{\sqrt{p}}{\pi}\left[\ln{4e^2p}+\ln{16e^4p^2}\right] \leq \frac{\sqrt{p}}{\pi}\ln{2^{10}e^6p^3}\\
&\leq\frac{3\sqrt{p}\ln{11e^2p}}{\pi}.
\end{align*}

}

\lem[always-recip-one]{Suppose $E$ has good ordinary (resp. good supersingular) reduction at $p$. Then $\rho_{|G_{\Q_p}}$ is split or wild (resp. $1$-Cartan, and one of its reciprocity indices is one).}

\demo{This follows from the local study of Serre \cite[\S 1.11]{Serre-image-ouverte}. More precisely, the case ``split or wild'' for good ordinary reduction is exactly the discussion ``b) or c)'' in the Corollary to Proposition 11 of \emph{loc.cit}. For good supersingular reduction, Proposition 12 of \emph{loc.cit.} shows that the image of the inertia subgroup is a nonsplit Cartan subgroup, so that $\rho$ is $1$-Cartan. Moreover, the Corollary to Proposition 12 of \emph{loc.cit.} shows using Proposition 3 that $1$ is a reciprocity index for $\rho_{|G_{\Q_p}}$ (because Serre's normalization of class field theory isomorphism is the inverse of ours). }

\cor{As $p \rar \infty$, uniformly in $E,\ell,\xi \in \mu_p^{\times}(\Qbar)$ as long as $E$ has good reduction at $p$, the asymptotics for the number of eigenspaces for $(\mathbb{T}_1 \otimes \overline{\Q_{\ell}})$ of $\Tate{\ell}{J_{\rho,\xi}(p)} \otimes_{\Z_{\ell}} \overline{\Q_{\ell}}$ with a given sign are given in Table \ref{number-of-signs}. Recall that the eigenspaces corresponding to $(f,\chi) \in \mathscr{C}_M$ have dimension $p-1$ (and there are $O(\sqrt{p}\log{p})$ of them), and the ones corresponding to non-CM newforms have dimension in $\{2p-2,2p,2p+2\}$. }

\begin{table}[htb]
\centering
\begin{tabular}{|c|c|c|c|}
\hline
& \multicolumn{3}{c|}{Type of $\rho_{|G_{\Q_p}}$}\\
\hline
& Split & Wild & $a$-Cartan\\
\hline
$(-1)^{\frac{p-1}{2}}$ & $\frac{p^2}{48}+O(p)$ & $\frac{p^2}{24} + O(p)$ & $\frac{p^2}{48}+O(p)$\\
\hline
$(-1)^{\frac{p+1}{2}}$ & $\frac{p^2}{48}+O(p)$ & $\frac{p}{24} + O(\sqrt{p}\log{p})$ & $\frac{p^2}{48}+O(p)$\\
\hline
\end{tabular}
\caption{Asymptotic number of eigenspaces in the Tate module of $J_{\rho,\xi}(p)$ with a given global root number}
\label{number-of-signs}
\end{table}

\section{On residual automorphy for the $p$-torsion of $J_{\rho}(p)$}
\label{residual-automorphy}

Let $F$ be a number field as in Section \ref{preparatory-lemmas}, that is: it contains the $p(p-1)^2(p+1)$-th roots of unity, the coefficients of all the newforms in $\mathscr{S} \cup \mathscr{P} \cup \mathscr{C}$, and, for each prime $\ell$ split in $\Q(\sqrt{-p})$ and every $(f,\mathbf{1}) \in \mathscr{C}_M$, $F$ contains the roots of $X^2-a_{\ell}(f)X+\ell$. Fix a prime ideal $\mathfrak{p} \subset \OO_F$ with residue characteristic $p$, and an embedding $\OO_F/\mathfrak{p} \rar \overline{\F_p}$. 

For instance, any character $\chi: \F_p^{\times} \rar \C^{\times}$ takes its values in $\OO_F^{\times}$, so can be reduced modulo $\mathfrak{p}$ to a character $\chi_{\mathfrak{p}}: \F_p^{\times} \rar \overline{\F_p}^{\times}$. Since $\chi_{\mathfrak{p}}$ takes its values in the $(p-1)$-torsion of $\overline{\F_p}^{\times}$, that is, in $\F_p^{\times}$. Thus $\chi_{\mathfrak{p}}$ is an endomorphism of the group $\F_p^{\times}$, so it is of the form $x \longmapsto x^{a_{\chi}}$ for some integer $1 \leq a_{\chi} \leq p-1$.  

Similarly, let $\phi: \F_{p^2}^{\times} \rar \C^{\times}$ be a character. It takes its values in $\OO_F^{\times}$, so we can consider its reduction $\phi_{\mathfrak{p}}: \F_{p^2}^{\times} \rar (\OO_F/\mathfrak{p})^{\times} \subset \overline{\F_p}^{\times}$. Since $\phi_{\mathfrak{p}}$ takes its values in the $(p^2-1)$-torsion of $\overline{\F_p}^{\times}$, which is $\F_{p^2}^{\times}$. Thus, $\phi_{\mathfrak{p}}: \F_{p^2}^{\times} \rar \F_{p^2}^{\times}$ can be uniquely written as $x \longmapsto x^{a+bp}$, with $0 \leq a,b < p$ and $(a,b) \neq (p-1,p-1)$. For instance, $a \equiv b\pmod{2}$ if and only if $\phi$ is even. 

Replacing $\phi$ with $\phi^p$ amounts to exchanging $a$ and $b$, and $\phi$ factors through the norm $\F_{p^2}^{\times} \rar \F_p^{\times}$ if and only if it is the case for $\phi_{\mathfrak{p}}$, which is the case if and only if $a=b$. Assuming that $\phi$ does not factor through the norm, we set $b_{\phi}=\min(a,b), a_{\phi}=\max(a,b)$. Note that these quantities only depend on the pair $\{\phi,\phi^p\}$, so in particular, they are independent from the embeddings $\OO_F/\mathfrak{p} \subset \overline{\F_p}$, $\F_{p^2} \subset \overline{\F_p}$. 

Let $\rho: G_{\Q} \rar \GL{\F_p}$ be a representation such that $\det{\rho}=\omega_p^{-1}$, and let $\rho^{\ast}: G_{\Q} \rar \GL{\F_p}$ denote its contragredient: thus $\det{\rho^{\ast}}=\omega_p$. 

\medskip

\prop[decomposition-modp-sspl]{The $\overline{\F_p}[G_{\Q}]$-module $J_{\rho}(p)[p](\Qbar) \otimes_{\F_p} \overline{\F_p}$ has the same semi-simplification as the following direct sum of $\overline{\F_p}[G_{\Q}]$-modules:
\begin{itemize}[noitemsep,label=\tiny$\bullet$]
\item $\left[\bigoplus_{(f,\mathbf{1}) \in \mathscr{S}}{T_{f,\mathfrak{p}}/\mathfrak{p} \otimes \mrm{Sym}^{p-1}{\rho^{\ast}}}\right]^{\oplus (p-1)}$,
\item $\left[\bigoplus_{(f,\chi) \in \mathscr{P}'}{(T_{f,\mathfrak{p}}/\mathfrak{p}) \otimes \left[\mrm{Sym}^{p-1-a_{\chi}}{\rho^{\ast}} \oplus (\omega_p^{-a_{\chi}} \otimes \mrm{Sym}^{a_{\chi}}{\rho^{\ast}})\right]}\right]^{\oplus (p-1)}$,
\item $\bigoplus_{(f,\chi) \in \mathscr{C}}{T_{f,\mathfrak{p}}/\mathfrak{p} \otimes \left[(\omega_p^{1-a_f} \otimes \mrm{Sym}^{a_f-b_f-2}{\rho^{\ast}}) \oplus (\omega_p^{-b_f} \otimes \mrm{Sym}^{(p-1)-(a_f-b_f)}{\rho^{\ast}})\right]}$, where $b_f, a_f$ are the two integers attached to the character $\phi_f: \F_{p^2}^{\times} \rar \OO_F^{\times}$ (so that $0 \leq b_f < a_f \leq p$). 
\end{itemize}}

\demo{Let $\chi: \F_p^{\times} \rar \OO_F^{\times}$ be a nontrivial character, then, by Lemma \ref{central-character-and-contragredient} and Proposition \ref{twist-principal-series}, 
\[\pi(1,\chi)(\rho) \simeq \pi(1,\chi)^{\circ \ast}(\rho^{\ast}) \simeq \chi(\det)^{-1}\pi(1,\chi)(\rho^{\ast}) \simeq \pi(\chi^{-1},1)(\rho^{\ast}).\] 
By Corollary \ref{decomp-principal-series}, the semi-simplification of $\pi(\chi_{\mathfrak{p}}^{-1},1)(\rho^{\ast})$ is the direct sum of $\omega_p^{-a_{\chi}} \otimes \mrm{Sym}^{a_{\chi}}\rho^{\ast}$ and of $\mrm{Sym}^{p-1-a_{\chi}}\rho^{\ast}$. 

The $\F_p[\GL{\F_p}]$-module $\mrm{St}$ is isomorphic to $\mrm{Sym}^{p-1}$, so, by Lemma \ref{central-character-and-contragredient}, one has \[\mrm{St}(\rho)=\mrm{St}^{\circ \ast}(\rho^{\ast}) \simeq \mrm{St}(\rho^{\ast}) \simeq \mrm{Sym}^{p-1}\rho^{\ast}.\]

Now, if $(f,\chi) \in \mathscr{C}$, $C_f(\rho)=C_f^{\circ \ast}(\rho^{\ast})=\chi^{-1}(\det{\rho^{\ast}}) \otimes C_f(\rho^{\ast}) \simeq \chi^{-1}(\omega_p) \otimes C_{f}(\rho^{\ast})$. By \cite[Lemma 4.2]{Prasad-cusp}, the semi-simplification of $C_{\phi} \pmod{\mathfrak{p}}$ is the direct sum of $\omega_p^{b_{\phi}+1} \otimes \mrm{Sym}^{a_{\phi}-b_{\phi}-2}{R}$ and $\omega_p^{a_{\phi}-b_{\phi}} \otimes \mrm{Sym}^{(p-1)-(a_{\phi}-b_{\phi})}{R}$, where $R$ is the tautological representation of $\GL{\F_p}$. Since $\chi^{-1}(\omega_p)=\omega_p^{-a_f-b_f}$, we can then apply Corollary \ref{decomposition-xrho}. 
}

\bigskip

The shape of this decomposition suggests that we can use results in Langlands functoriality to prove that the components described in the previous propositions are mod $p$ reductions of Galois representations attached to automorphic representations. The ones we will use are \cite[Theorem A]{NewThor} for symmetric power functoriality, and the so far unpublished preprint \cite[Theorem 1.1]{AdRD} for $\operatorname{GL}_2 \times \operatorname{GL}_n$ functoriality.   

\medskip

\prop[jrho-mod-p-automorphic]{Assume that $\rho^{\ast}$ is onto, or is irreducible with Serre weight in $\{2,p+1\}$, and that its Artin conductor $N$ is coprime to $3$. Let $V$ be a $\overline{\F_p}[G_{\Q}]$-module of dimension $2r$ over $\overline{\F_p}$ which is one of the summands appearing in Proposition \ref{decomposition-modp-sspl}. Then, there exists a regular algebraic cuspidal polarized automorphic representation $\pi$ of $\GLn{2r}{\mathbb{A}_{\Q}}$ with coefficients in some number field $F'$ and a prime ideal $\mathfrak{p}_0$ of residue characteristic $p$ such that, if $V_{\pi,\mathfrak{p}_0}$ is the $\mathfrak{p}_0$-adic representation of $G_{\Q}$ attached to $\pi$, then the reduction mod $\mathfrak{p}_0$ of $V_{\pi,\mathfrak{p}_0}$ has the same semi-simplification as $V$ (for a suitable embedding $\OO_{F'}/\mathfrak{p}_0 \rar \overline{\F_p}$). }

\demo{By Proposition \ref{decomposition-modp-sspl}, $V \simeq T_{f,\mathfrak{p}}/\mathfrak{p} \otimes \mrm{Sym}^a\rho^{\ast} \otimes u(\omega_p)$ for some $(f,\chi) \in \mathscr{S} \cup \mathscr{P} \cup \mathscr{C}$, where $u: \F_p^{\times} \rar \F_p^{\times}$ is the reduction mod $\mathfrak{p}$ of a character $u_0: \F_p^{\times} \rar \OO_F^{\times}$. Clearly, one has $a=r-1$. 

Let us show that $\rho^{\ast}$ is attached to a non-CM modular form $g \in \mathcal{S}_k(\Gamma_1(N'))$ for some integer $N'$ coprime to $3p$. By Serre's conjecture (now a theorem \cite{KW1,KW2}), $\rho^{\ast}$ is attached to a normalized newform $g_1 \in \mathcal{S}_k(\Gamma_1(N))$ such that $k \equiv 2 \pmod{p-1}$. If $g_1$ is not CM, we can take $(g,N')=(g_1,N)$. Otherwise, $\rho^{\ast}$ is irreducible and its Serre weight $k$ lies in $\{2,p+1\}$. By Cebotarev's theorem, there is a prime $q > 3Np$ such that $\rho^{\ast}(\Fr_q)$ is conjugate to the image under $\rho^{\ast}$ of the complex conjugation. In particular, $q \equiv -1\pmod{p}$, so the characteristic polynomial of $\rho^{\ast}(\Fr_q)$ is $(X-1)(X-q)$. Hence, by Diamond's theorem\footnote{I am grateful to Fred Diamond for suggesting to approach this question under the lens of level-raising congruences and for directing me towards his theorem.} \cite[Theorem 5.1]{Ribet-report}, $\rho^{\ast}$ arises from an eigenform in $\mathcal{S}_k(\Gamma_1(N) \cap \Gamma_0(q))$ with level divisible by $q$. By \cite[Theorem 4.2]{Ribet-report}, $\rho^{\ast}$ arises from a newform $g \in \mathcal{S}_k(\Gamma_1(N)\cap \Gamma_0(q))$. By Proposition \ref{langlands-elem-2}, $\pi_{g,q}$ is the twist of a Steinberg representation. Therefore, by Proposition \ref{local-global-modular-forms} and the construction of the local Langlands correspondence in \cite[(33.3)]{BH}, $D_{g,q}$ is a special representation of $W_{\Q_q}$, so $g$ cannot have complex multiplication, and we let $(g,N')=(g,Nq)$. \smallskip

Let $F' \subset \C$ be a number field containing $F$ and the Fourier coefficients of $g$. Let $\mathfrak{p}_0$ be a prime ideal of $\OO_{F'}$ above $\mathfrak{p}$ such that for almost every prime $q$, $a_q(g) \pmod{\mathfrak{p}_0} = \mrm{Tr}(\rho^{\ast}(\Fr_q))$. 

The argument of \cite[\S 2]{Ribet-report} 
 shows that, after multiplying by a power of the Eisenstein series $E_{p-1} \in 1+\frac{p}{T}\Z[[q]]$ (where $T$ is an integer coprime to $p$) or enlarging $F'$ if necessary, there is a normalized eigenform $f_1 \in \mathcal{S}_{k'}(\Gamma(1))$ with coefficients in $F'$ such that, for all primes $\ell \neq p$, $a_{\ell}(f_1) \equiv a_{\ell}(f\otimes u_0) \mod{\mathfrak{p}_0}$, and such that $k' > (a+1)(k+1)$. In particular, $f_1$ does not have complex multiplication and $g$ is not in the Galois orbit of a twist of $f_1$. 

By this construction, $V$ and $[T_{f_1,\mathfrak{p}_0} \otimes \mrm{Sym}^a{T_{g,\mathfrak{p}_0}}]/\mathfrak{p}_0$ has the same semi-simplification as $V$. 

By \cite[Theorem A]{NewThor}, there is an regular algebraic cuspidal polarized automorphic representation $\pi_g$ of $\GL{r}{\mathbb{A}_{\Q}}$ whose attached $\lambda$-adic representation is $\mrm{Sym}^a{V_{g,\lambda}}$ for all maximal ideals $\lambda \subset \OO_{F'}$: in particular, $\pi_g$ is unramified away from the primes dividing $N'p$, so its level is coprime by $3$. The conclusion follows if we can apply \cite[Theorem 1.1]{AdRD}, so all we need to do is check its hypotheses of regularity and irreducibility. 

For any maximal ideal $\lambda \subset \OO_{F'}$, the Hodge-Tate weights of $V_{f_1,\lambda} \otimes \mrm{Sym}^a{V_{g,\lambda}}$ are pairwise distinct, so the compatible system $(V_{f_1,\lambda} \otimes \mrm{Sym}^a{V_{g,\lambda}})_{\lambda}$ is regular. 

By \cite[Theorem 3.2.2]{bigimage}, for any large enough prime $\ell$ and any prime ideal $\lambda \supset \OO_{F'}$ with residue characteristic $\ell$, the image of $G_{\Q(\mu_{\ell^{\infty}})}$ on $V_{f_1,\lambda} \oplus V_{g,\lambda}$ contains (for a suitable choice of bases) $\SL{\Z_{\ell}} \oplus \SL{\Z_{\ell}}$, so by Lemma \ref{subreps-product} $V_{f_1,\lambda} \otimes \mrm{Sym}^a{V_{g,\lambda}}$ is irreducible (because the action of $\SL{\Z_{\ell}}$ on $\mrm{Sym}^t(\overline{\Q_{\ell}}^{\oplus 2})$ is irreducible\footnote{This well-known fact can be shown as follows: the action of $U=\begin{pmatrix}1 & 1\\0 & 1\end{pmatrix}$ on an irreducible subrepresentation is unipotent, so this irreducible subrepresentation must contain some nonzero vector in $\ker(U-\mrm{id})$. Such a vector is, up to a scalar, $x=\begin{pmatrix}1 \\ 0\end{pmatrix}^k$, and the $U^T$-translates of $x$ generate $V:=\mrm{Sym}^t(\overline{\Q_{\ell}}^{\oplus 2})$, because $\dim{V}=t+1$ and $(U^T-\mrm{id})^tx \neq 0$, $(U^T-\mrm{id})^{t+1}=0$.}). 
}

\rems{\begin{itemize}[noitemsep,label=$-$]
\item \cite[Remark 1.4]{AdRD} suggests that the assumption $3 \nmid N$ is not necessary, but does not prove it. 
\item The purpose of the other assumptions on $\rho^{\ast}$ is to make sure that it is attached to a non-CM modular form in characteristic zero. It seems likely \cite{Diamond-personal} that the other assumptions can be weakened, but we have not been able to prove it. 
\end{itemize}}

\newpage

%% file: cartan-0.tex
\chapter{Representations valued in the normalizer of a nonsplit Cartan subgroup}
\label{small-image}

\section{Statement of results}

Throughout this chapter, $p\geq 7$ is a prime and $\rho: G_{\Q} \rar \GL{\F_p}$ is a continuous representation such that $\det{\rho}=\omega_p^{-1}$, and $\rho(G_{\Q})$ is contained in the normalizer $N$ of a nonsplit Cartan subgroup $C$. We call $K$ the imaginary quadratic field such that $G_K = \rho^{-1}(C)$ (see Proposition \ref{is-usually-cartan}). The goal of this section is to find out what more information we can tell about the various curves $X_{\rho,\xi}(p)$. 

The reason why we give this case a particular focus is that when $p \geq 17$ is a prime and $E/\Q$ is an elliptic curve, one of the following statements is true, by \cite{MazurY1,FreyMazur,BP,BPR}:
\begin{itemize}[noitemsep,label=$-$]
\item $E$ has complex multiplication,
\item $p=17$ (resp. $p=37$) and $E$ is a quadratic twist of an elliptic curve in the LMFDB isogeny class \cite{lmfdb} with label $14450.b$ (resp. $1225.b$),
\item The image of the group homomorphism $G_{\Q} \rar \mrm{Aut}(E[p])$ is contained in the normalizer of a nonsplit Cartan subgroup,
\item The group homomorphism $G_{\Q} \rar \mrm{Aut}(E[p])$ is onto.
\end{itemize}

Serre's uniformity question \cite[\S 4.3]{Serre-image-ouverte}, still currently open, asks whether the third case is redundant, at least for large enough primes $p$. This amounts to asking whether, for a large enough prime $p$, the rational points of a certain modular curve $X_{ns}^+(p)$ only correspond to cusps or CM points. An affirmative answer is suggested by the main theorem of \cite{Xns+13}, which uses a new method for the determination of rational points, the \emph{Chabauty-Kim} method, and shows that for elliptic curves $E/\Q$ without complex multiplication, the action of $G_{\Q}$ on $E[13]$ is either reducible, surjective, or has image contained in an exceptional subgroup. In \cite{Xns+17}, a similar method is used in \S 5.1 to describe all the elliptic curves corresponding to this last case, and in \S 5.5 the authors show that the third case is redundant for $p=17$. They run into computational issues for $p=19$, but \cite[Theorem 1, Remark 1]{LFD} proves that the method can in principle be made to work for all primes $p \geq 13$. 

In order to describe the results of this section, let us recall the following notation from the previous sections. 

\textbf{Notation for spaces of modular forms:}
\begin{itemize}[label=\tiny$\bullet$,noitemsep]
\item $\mathcal{D}$ denotes the set of Dirichlet characters $(\Z/p\Z)^{\times} \rar \C^{\times}$, 
\item $\mathscr{S}$ denotes the set of $(f,\mathbf{1})$ where $f \in \mathcal{S}_2(\Gamma_0(p))$ is a newform, 
\item $\mathscr{P}$ denotes the set of $(f,\chi)$, where $f \in \mathcal{S}_2(\Gamma_1(p))$ is a newform with character $\chi \neq \mathbf{1}$,
\item $\mathscr{P}'$ denotes a set of representatives for $\mathscr{P}$ modulo the involutive action of the complex conjugation, 
\item $\mathscr{C}$ denotes the set of $(f,\chi)$, where $f \in \mathcal{S}_2(\Gamma_1(p) \cap \Gamma_0(p^2))$ is a newform with character $\chi \in \mathcal{D}$, and such that no twist of $f$ by a character in $\mathcal{D}$ has conductor $p$,
\item $\mathscr{C}'_1$ denotes a set of representatives of the form $(f,\mathbf{1})$ for the $(f,\chi) \in \mathscr{C}$ modulo twists such that $f$ does not have complex multiplication.
\item $\mathscr{C}_M$ denotes the collection of $(f,\mathbf{1}) \in\mathscr{C}$ such that $f$ has complex multiplication. In this case, $p \equiv 3\pmod{4}$ and $f$ has complex multiplication by $\Q(\sqrt{-p})$.
\end{itemize}

\textbf{Notation for representations of $\GL{\F_p}$:} 
\begin{itemize}[noitemsep,label=\tiny$\bullet$]
\item $I_2 \in \GL{\F_p}$ is the identity matrix.
\item For $\alpha \in \mathcal{D}$, $\mrm{St}_{\alpha}$ is the Steinberg representation of $\GL{\F_p}$ twisted by $\alpha(\det)$. 
\item For $\alpha,\beta \in \mathcal{D}$, $\pi(\alpha,\beta)$ is the principal series representation of $\GL{\F_p}$ attached to the characters $\alpha,\beta$. 
\item Given $(f,\chi) \in \mathscr{C}$, we can attach by Definition \ref{left-cusp-rep-f} a cuspidal representation of $\GL{\F_p}$ named $C_f$. 
\end{itemize}    

\textbf{Notation for Hecke algebras:}
\begin{itemize}[noitemsep,label=\tiny$\bullet$]
\item $\mathbb{T}$ is the subalgebra of $\mrm{End}\left(\prod_{\xi \in \mu_p^{\times}(\Qbar)}{J_{\rho,\xi}(p)}\right)$ generated by the Hecke operators $T_n$ and the $(mI_2)$, where $n \geq 1$ is coprime to $p$ and $m \in \F_p^{\times}$. 
\item $\mathbb{T}_1$ is the subalgebra of $\mathbb{T}$ generated by the Hecke operators $(mI_2)\cdot T_n$, where $m,n \geq 1$ are integers such that $mn^2 \equiv 1\pmod{p}$. It is contained in $\bigoplus_{\xi \in \mu_p^{\times}(\Qbar)}{\mrm{End}(J_{\rho,\xi}(p))}$.
\end{itemize}

Let $f \in \mathcal{S}_k(\Gamma_1(N))$ be a normalized newform with character $\chi$ and with coefficients contained in the ring of integers $\OO_L$ of a number field $L$. Let $\lambda$ be a maximal ideal of $\OO_L$, then $V_{f,\lambda}$ denotes the continuous irreducible representation $G_{\Q} \rar \GL{L_{\lambda}}$ such that for any prime $q$ coprime to $N\lambda$, the characteristic polynomial of the arithmetic Frobenius $\Fr_q$ is $X^2-a_q(f)X+q^{k-1}\chi(q)$. 

We choose a number field $F \subset \C$ which is Galois over $\Q$ and whose ring of integers $\OO_F$ contains the following algebraic numbers: 
\begin{itemize}[noitemsep,label=\tiny$\bullet$]
\item the $p(p-1)^2(p+1)$-th roots of unity, 
\item the Fourier coefficients of any $f$ for $(f,\chi) \in \mathscr{S}\cup\mathscr{P}\cup\mathscr{C}$, 
\item for every $(f,\mathbf{1}) \in \mathscr{C}_M$, for every prime $\ell$ split in $\Q(\sqrt{-p})$, the roots of the Hecke polynomial $X^2-a_{\ell}(f)X+\ell$ of $f$ at $\ell$. 
\end{itemize}

Our first step is to describe the Tate module of $J_{\rho,\xi}(p)$ (for $\xi \in \mu_p^{\times}(\Qbar)$) as a Galois-module. 

\theoi[cartan-1]{(See Proposition \ref{tate-module-xrho-nsplit-cartan}) Let $\mathfrak{l}$ be a maximal ideal of $\OO_F$ with residue characteristic $\ell$. For any $\xi \in \mu_p^{\times}(\Qbar)$, $\Tate{\ell}{J_{\rho,\xi}(p)} \otimes_{\Z_{\ell}} F_{\mathfrak{l}}$ is isomorphic to the direct sum of $(\mathbb{T}_1 \otimes F_{\mathfrak{l}})[G_{\Q}]$-modules of the following form:
\begin{itemize}[noitemsep,label=\tiny$\bullet$]
\item $V_{f,\mathfrak{l}} \otimes \alpha(\rho)$, where $(f,\chi) \in \mathscr{S} \cup \mathscr{P}' \cup \mathscr{C}'_1 \cup \mathscr{C}_M$, and $\alpha: N \rar F^{\times}$ is a character such that $\alpha(tI_2)=\chi(t)$ for $t \in \F_p^{\times}$. 
\item $V_{f,\mathfrak{l}} \otimes \left[\mrm{Ind}_C^N{\psi}\right](\rho)$, $(f,\chi) \in \mathscr{S} \cup \mathscr{P}' \cup \mathscr{C}'_1 \cup \mathscr{C}_M$, where $\psi: C\rar F^{\times}$ is a character such that $\psi^{p-1} \neq \mathbf{1}$, $\psi(tI_2)=\chi(t)$ for every $t \in \F_p^{\times}$, and, if $(f,\chi) \in \mathscr{C}'_1 \cup\mathscr{C}_M$, $\psi$ is not one of the two characters attached to $f$. 
\end{itemize}
The Hecke operator $(mI_2) \cdot T_n$ acts on each summand by $\chi(m)a_n(f)$. Moreover, for every submodule $V$ as above, one has $V^{\ast}(1) \simeq V$ as $F_{\mathfrak{l}}[G_{\Q}]$-submodules. 
}

As a corollary, the $L$-function of $J_{\rho,\xi}(p)$ splits into factors according to the decomposition given in Theorem \ref{cartan-1}. By the classical Rankin-Selberg method, we can show that these factors extend as entire functions on the complex plane and satisfy a functional equation, whose signs we compute. 

We report here the computation of this sign when $K$ is inert at $p$ (which is the case for the application to elliptic curves when $p \geq 11$, see Lemma \ref{almost-never-ramified}). 

\theoi[cartan-2]{(See Corollary \ref{functional-equations-jrho-ncartan} and Section \ref{signs-functional-equation-cartan}) Let $F$ be a large enough number field and $\mathfrak{l} \subset \OO_F$ be a maximal ideal of residue characteristic $\ell$, and fix $\xi \in \mu_p^{\times}(\Qbar)$. Let $V$ be one of the summands of $\Tate{\ell}{J_{\rho,\xi}(p)} \otimes_{\Z_{\ell}} F_{\mathfrak{l}}$ (as a $F_{\mathfrak{l}}[G_{\Q}]$-module) described in Theorem \ref{cartan-1}. Then the $L$-function of $V$ extends to an entire function on the complex plane and satisfies the predicted functional equation. When $K$ is inert at $p$, the signs of these functional equations are the following:
\begin{itemize}[noitemsep,label=$-$]
\item If $V$ is two-dimensional, or four-dimensional with $(f,\chi) \in \mathscr{P}'$, the sign of the functional equation is $-1$. 
\item If $V$ is four-dimensional with $(f,\mathbf{1}) \in \mathscr{S}$, the sign is $1$ if, and only if, the representation $\left[\mrm{Ind}_C^N{\psi}\right](\rho)$ is unramified at $p$. 
\item If $V$ is four-dimensional with $(f,\mathbf{1}) \in \mathscr{C}'_1\cup\mathscr{C}_M$, let $\phi: \F_{p^2}^{\times} \rar F^{\times}$ be the character attached to $f$, $I_p$ be an inertia subgroup at $p$, and $\mathfrak{a}: I_p \rar \F_{p^2}^{\times}$ be given by the class field theory reciprocity homomorphism (for $E$ and normalized so that it maps an arithmetic Frobenius to a uniformizer). The sign of the functional equation is $1$ if, and only if\footnote{Since the first submission of this thesis, I realized that this criterion could be written in a clearer way: it is equivalent to ask $\psi^{k-1} \in \{\phi_f,\phi_f^p\}$, where $k$ is the Serre weight of the contragredient of $\rho$.}, $\phi\circ\mathfrak{a}$ is one of the characters $\psi(\rho_{|I_p}), \psi^p(\rho_{|I_p})$.
\end{itemize} 
}

\cori[cartan-2cor]{(Corollary \ref{good-reduction-all-minusone}) When $\rho$ comes from an elliptic curve with good reduction at $p$, then the signs of all the functional equations in Theorem \ref{cartan-2} are equal to $-1$. In particular, in this case, the central value of the $L$-function of all the factors described in Theorem \ref{cartan-1} vanishes. }


Nonetheless, in certain situations, the signs of some factors of the Tate module of $J_{\rho,\xi}(p)$ can be equal to one, and it is then possible for their central $L$-value to not vanish. In such a situation, thanks to the Euler systems constructed by Kato \cite{KatoBSD} and Kings, Loeffler and Zerbes \cite{KLZ15}, we can prove some cases of the Bloch-Kato conjecture in rank zero. In the case of two-dimensional factors, the result is a direct application of Kato's work; however, in the case of four-dimensional factors, applying the Euler system method requires that a certain \emph{large image} condition be satisfied. This verification is carried out in Chapter \ref{obstructions-euler}. 

\theoi[cartan-3]{(See Section \ref{bloch-kato-conjecture}) Let $F$ be a large enough number field and $\mathfrak{l} \subset \OO_F$ be a maximal ideal of residue characteristic $\ell$, and fix $\xi \in \mu_p^{\times}(\Qbar)$. Let $V$ be one of the summands of $\Tate{\ell}{J_{\rho,\xi}(p)} \otimes_{\Z_{\ell}} F_{\mathfrak{l}}$ (as a $F_{\mathfrak{l}}[G_{\Q}]$-module) described in Theorem \ref{cartan-1} such that $L(V,1) \neq 0$.
\begin{enumerate}[noitemsep,label=(\roman*)]
\item Suppose that $V$ is of the form $V_{f,\mathfrak{l}} \otimes \alpha(\rho)$, for some character $\alpha: N \rar F^{\times}$. Then, for all but finitely many $\mathfrak{l}$, the divisible $\OO_{F_{\mathfrak{l}}}$-linear Selmer group attached to $V$ is finite.  
\item Suppose that $V$ is of the form $V_{f,\mathfrak{l}} \otimes \left[\mrm{Ind}_C^N{\psi}\right](\rho)$ for some character $\psi: C \rar F^{\times}$ with $(f,\mathbf{1}) \notin \mathscr{C}_M$. Then, for all but finitely many $\mathfrak{l}$ with residue characteristic $\ell$ such that $f$ is ordinary at $\mathfrak{l}$ and that $\psi^{p-1}(\Fr_{\ell}) \neq 1$, the divisible $\OO_{F_{\mathfrak{l}}}$-linear Selmer group attached to $V$ is finite.  
\end{enumerate}
}

The finiteness of the Mordell-Weil group in a suitable quotient of the Jacobian is the main ingredient to apply Mazur's \emph{formal immersion argument}. It was first used in \cite[Corollary 4.3]{FreyMazur} and then became a standard technique in the determination of rational points on modular curves without reference to a model (see for instance \cite{Merel,BP,BPR,LFL}). We use here a more explicit version of this argument, inspired by Coleman's approach to the Chabauty method \cite{ColChab}, which was explained to me by Samuel Le Fourn. 

\defi{Let $p \equiv 2, 5\pmod{9}$ be a prime number. We say that a prime number $\ell$ is \emph{exceptional} for $p$ if $\ell \equiv \pm 1 \pmod{p}$ and there exists a prime ideal $\lambda$ of $\Z[e^{2i\pi/p}]$ with residue characteristic $\ell$ such that, for every cube root of unity $j \in \F_{p^2}^{\times}$, \[u_j:=\sum_{\substack{b \in \F_p\\1+b\sqrt{-3} \in \F_{p^2}^{\times 3}j}}{e^{\frac{2i\pi b}{p}}} \in \lambda.\] } 

Direct computations with Magma \cite{magma} show that for any $p \leq 1500$, the three $u_j$ generate the unit ideal of $\Z[e^{\frac{2i\pi}{p}}]$. We thus tentatively expect that for any $p$, there are no exceptional primes\footnote{Since the first submission of this thesis, this claim has been verified for $p \leq 5000$.}.

\theoi[cartan-4]{(Theorem \ref{formal-immersion-application}) Assume that 
\begin{itemize}[noitemsep,label=$-$]
\item $\rho: G_{\Q} \rar N$ is surjective of conductor $N_{\rho}$, 
\item the image under $\rho$ of an inertia group at $p$ is $3C$,
\item $X_{\rho,\xi}(p)$ has a non-cuspidal rational point corresponding to an elliptic curve $E$ over $\Q$ for some $\xi$,
\item there exists $(f,\mathbf{1}) \in \mathscr{S}$ such that, if $\psi: C \rar F^{\times}$ is a character with order $3$, and $g$ is the weight one holomorphic newform attached to the irreducible odd Artin representation $\left[\mrm{Ind}_C^N{\psi}\right](\rho)$, then $L(f \times g,1) \neq 0$. 
\end{itemize}
Then $p \equiv 2, 5 \pmod{9}$ and the conductor of $E$ is of the form $p^2N_{\rho}L$, where $L$ is a square-free product of exceptional primes for $p$ that do not divide $pN_{\rho}$. 
}

By our previous discussion about Serre's uniformity question, it is expected that the elliptic curve $E$ has complex multiplication when $p \geq 13$ (in particular, it has potentially good reduction everywhere and $L=1$), but this remains an open question. The recent works of Le Fourn and Lemos \cite{LFL} then Furio and Lombardo \cite{lombardo} show that unless $E$ has complex multiplication, one has $\rho(G_{\Q})=N$. By Lemma \ref{cartan-is-potentially-good}, the image of an inertia subgroup at $p$ is either $C$ or $3C$, and the latter case can only occur if $p \equiv 2, 5\pmod{9}$. \\

Using the computation of the Serre weight of the $p$-torsion of an elliptic curve by Kraus \cite{Kraus-Thesis}, we find that, for the elliptic curve $E_p$ with Weierstrass equation $y^2=x^3-p$ for $p \equiv 5\pmod{9}$, the Artin representation appearing in Theorem \ref{cartan-4} is the one attached to the unique newform $g \in \mathcal{S}_1(\Gamma_1(108))$ with rational coefficients and CM by $\Q(\sqrt{-3})$. For $p=23$, I checked numerically in PARI/GP with the help of Pascal Molin the non-vanishing of $L(f \times g,1)$ for some $(f,\mathbf{1}) \in \mathscr{S}$, which implies the following result. \footnote{Since the submission of this thesis, I was able to carry out further numerical computations, showing that the non-vanishing condition held for $23 \leq p \leq 150$ and $p \equiv 5 \pmod{9}$. Therefore, the proof of Corollary \ref{cartan-4cor} carries over to $p=41, 59$.} 

\cori[cartan-4cor]{Let $p=23$. Any elliptic curve $E/\Q$ congruent modulo $p$ to the elliptic curve $E_p$ is isogenous to $E_p$.}

Using a result of analytic number theory of Michel \cite[Appendix, Theorem 0.1]{Michel} pointed out to me by Farrell Brumley, the conclusion of Theorem \ref{cartan-4} holds for all but finitely many $E_p$, in the following sense.

\theoi[cartan-5]{(Proposition \ref{diophantine-consequence-asymptotic}) There exists a constant $p_0 \geq 11$ such that, for any prime $p \geq p_0$ congruent to $5$ modulo $9$, for any elliptic curve $E'/\Q$ such that $E'[p](\Qbar)$ is isomorphic, as a $G_{\Q}$-module, to the $p$-torsion of the elliptic curve $E_p$ with Weierstrass equation $y^2=x^3-p$. Then the conductor of $E'$ is of the form $108mp^2L$, where $L$ is a square-free product of exceptional primes for $p$, and $m=1$ if $p \equiv 3\pmod{4}$ and $m=4$ if $p \equiv 1\pmod{4}$. }

\section{The Tate module of the $J_{\rho,\xi}(p)$}

Before discussing the object of this section, we describe the representation $\rho$ in more detail. 

\prop[is-usually-cartan]{There exists an imaginary quadratic field $K/\Q$ such that $\rho^{-1}(C)=G_K$, and $K$ is inert or ramified at $p$. Moreover,
\begin{itemize}[noitemsep,label=$-$]
\item If $K/\Q$ is unramified at $p$, then $\rho$ is $a$-Cartan (see Proposition \ref{trichotomy-rho}) for some odd $a \mid \frac{p+1}{2}$. In particular, $[C:\rho(G_K)] \mid a$. 
\item If $p \equiv 1 \pmod{4}$ and $K/\Q$ is ramified at $p$, then $\rho$ is $\frac{p+1}{2}$-Cartan and $\rho(I_p \cap G_K)= \F_p^{\times}I_2$ for any inertia subgroup $I_p \leq G_{\Q}$ at $p$.
\item If $p \equiv 3 \pmod{4}$ and $K/\Q$ is ramified at $p$, then $\rho$ is split and $\rho(G_{\Q_p} \cap G_K)\in \{\F_p^{\times 2}I_2,\F_p^{\times}I_2\}$.     
\end{itemize}}

\demo{Since $\det{\rho}=\omega_p^{-1}$ is an odd character, the image of the complex conjugation under $\rho$ is conjugate to $\begin{pmatrix}-1 & 0\\0 & 1\end{pmatrix}$, so that $\rho(G_{\Q}) \not\subset C$, and $\rho^{-1}(C)$ is thus an open subgroup of $G_{\Q}$ with index two not containing the complex conjugation: it corresponds to an imaginary quadratic field $K$.

Fix an embedding $G_{\Q_p} \leq G_{\Q}$, let $I_p \leq G_{\Q_p}$ be the inertia subgroup. The group $N$ does not contain any element of order $p$, so $\rho$ is tamely ramified at $p$, so it cannot be wild (see Proposition \ref{trichotomy-rho}). 

If $p$ splits in $K$, then the restriction of $\rho$ to a decomposition subgroup at $p$ is contained in $C$, so Proposition \ref{trichotomy-rho} implies that $\rho$ is split: $\rho(G_{\Q_p})  \subset C$ is comprised of matrices which are diagonalizable over $\F_p$, $\rho(G_{\Q_p}) \subset \F_p^{\times}I_2$, hence $\F_p^{\times}=\omega_p^{-1}(G_{\Q_p}) \subset \det{\F_p^{\times}I_2}$, a contradiction. So $K$ is inert or ramified at $p$. 

Assume that $\rho_{|G_{\Q_p}}$ is split. Then $\rho(G_{\Q_p})$ is a subgroup of $N$ made with matrices with eigenvalues in $\F_p^{\times}$. Thus $\rho(G_{\Q_p}) \cap C=\rho(G_{\Q_p} \cap G_K)$ is scalar. Moreover, by Proposition \ref{trichotomy-rho}, $\rho(I_p)$ is a cyclic subgroup of order $p-1$ and has surjective determinant, so it is not scalar: hence $\rho(I_p) \not\subset C$, $\rho(I_p) \cap C=\rho(I_p)^2$, and $K$ is ramified at $p$. Thus $\rho(I_p) \cap C$ is a cyclic subgroup of order $\frac{p-1}{2}$ of $\F_p^{\times}I_2$, so $\rho(I_p) \cap C=\F_p^{\times 2}I_2$ and $\F_p^{\times 2}=\det{\rho(I_p)^2}=\det{\rho(I_p) \cap C}=\det{\F_p^{\times 2}I_2}=\F_p^{\times 4}$ so $p \equiv 3 \pmod{4}$.

Assume that $\rho$ is $a$-Cartan and $K/\Q$ is ramified at $p$. Let $E/\Q_p$ denote the quadratic unramified extension, there is a Cartan subgroup $D \leq \GL{\F_p}$ such that $\rho(I_p)=aD$, and $\rho(I_p \cap G_K)$ is an index two subgroup of $aD$ contained in $C$, so it is equal to $2aD \subset C$. Because $C$ is cyclic, this implies $2aD=2aC$. Moreover, any $z \in N \backslash C$, $\mrm{Tr}(N)=0$, so $N^2=-(\det{N})I_2$ and every element in $aD \backslash C $ has order dividing $2(p-1)$. So $aD$ is contained in the reunion of the two subgroups $2aD$ and its subgroup of $2(p-1)$-torsion. Since $a$ is odd, $aD \neq 2aD$, so $aD$ is its subgroup of $2(p-1)$-torsion, thus $a$ is a multiple of $\frac{p^2-1}{2(p-1)}=\frac{p+1}{2}$. Thus $a=\frac{p+1}{2}$; since $a$ is odd, one has $p \equiv 1 \pmod{4}$, and $\rho(I_p) \cap C=2aD=(p+1)D=\F_p^{\times}I_2$. 

The conclusion follows.}
\medskip

\nott{For the rest of the chapter, $\epsilon$ is the unique nontrivial character of $N/C \simeq \{\pm 1\}$: then $\epsilon_K = \epsilon\circ\rho$ is the Galois character attached to the imaginary quadratic field $K$, and we let $M$ be the conductor of the corresponding primitive Dirichlet character. The group homomorphism $G_K \rar C$ induced by $\rho$ is denoted $\rho_K$. }

\medskip

\defi{For $(f,\chi) \in \mathscr{S} \cup \mathscr{P} \cup \mathscr{C}$, we define $H_f$ as the following set of complex characters of $N$, and $H_f(\rho)$ as the collection of the $\alpha \circ \rho$ for $\alpha \in H_f$: 
\begin{itemize}[label=\tiny$\bullet$, noitemsep] 
\item If $(f,\chi) \in \mathscr{S}$, then $H_f$ consists of the single character $\left(\frac{\det}{p}\right)\epsilon_K^{\mathbf{1}(p \equiv 3\pmod{4})}$. Thus $H_f(\rho)$ contains the single character $\left(\frac{\omega_p}{p}\right)\epsilon_K^{\mathbf{1}(p \equiv 3\pmod{4})}$. 
\item If $(f,\chi) \in \mathscr{P}$, then $H_f$ contains the two characters $\psi(\det)\epsilon^{\mathbf{1}(\psi(-1)=-1)}$, for every $\psi \in \mathcal{D}$ such that $\psi^2=\chi$. Thus $H_f(\rho)$ consists of the two characters $\overline{\psi}(\omega_p)\epsilon_K^{\mathbf{1}(\psi(-1)=-1)}$ where $\psi$ is as in the previous sentence.
\item If $(f,\chi) \in \mathscr{C}$ and $f$ does not have complex multiplication, then $H_f$ contains the two characters $\psi(\det)\epsilon^r$, for every $\psi \in \mathcal{D}$ such that $\psi^2=\chi$, and $r \in \Z/2\Z$ is such that $(-1)^r\lambda_p(f \otimes \overline{\psi})=1$. Thus $H_f(\rho)$ consists of the two characters $\overline{\psi}(\omega_p)\epsilon_K^r$, for $(\psi,r)$ as in the previous sentence.
\item If $(f,\chi) \in \mathscr{C}$ and $f$ has complex multiplication, then $H_f$ (resp. $H_f(\rho)$) consists of the character $\psi(\det)\epsilon^{\frac{p-3}{4}}$ (resp. $\overline{\psi}(\omega_p)\epsilon_K^{\frac{p-3}{4}}$), where $\psi \in\mathcal{D}$ is the only \emph{even} character such that $\psi^2=\chi$. 
\end{itemize}

We also define the set $\mathcal{I}_f$ of characters of $C$, endowed with the equivalence relation $\sim$, as follows:
\begin{itemize}[label=\tiny$\bullet$, noitemsep]
\item For $(f,\chi) \in \mathscr{S} \cup \mathscr{P}$, $\mathcal{I}_f$ is the set of characters $C \rar \C^{\times}$ such that $\psi^{p} \neq \psi$ and $\psi(aI_2)=\chi(a)$ for all $a \in \F_p^{\times}$. We write $\alpha \sim \beta$ if $\alpha \in \{\beta,\beta^p\}$. 
\item For $(f,\chi) \in \mathscr{C}$ such that $f$ has no complex multiplication, $\mathcal{I}_f$ is the set of characters $\psi: C \rar \C^{\times}$ such that $\psi^p \neq \psi$, $\psi(aI_2)=\chi(a)$ for all $a \in \F_p^{\times}$, and $\psi$ is not one of the two characters of $C$ attached to $C_f$. We write $\alpha \sim \beta$ if $\alpha \in \{\beta,\beta^p\}$.
\item For $(f,\chi) \in \mathscr{C}_M$, $\mathcal{I}_f$ is the set of characters $\psi: C \rar \C^{\times}$ such that $\psi^{2(p-1)} \neq 1$ and $\psi(aI_2)=\chi(a)$ for every $a \in \F_p^{\times}$. We write $\alpha \sim \beta$ if $\alpha\beta^{-1}$ or $\alpha\beta^{-p}$ has order dividing $2$.
\end{itemize}}

\rems{\begin{itemize}[noitemsep,label=$-$]
\item For $\alpha: C \rar \F_p^{\times}$, $\mrm{Ind}_K^{\Q}{\alpha(\rho_K)}$ is reducible if and only if $\alpha(\rho_K)=\alpha^p(\rho_K)$, i.e. if and only if $\alpha^{p-1}(\rho_K)=1$. 
\item For any $(f,\chi) \in \mathscr{S}\cup\mathscr{P}\cup\mathscr{C}$, for any $\sigma \in G_{\Q}$, there are bijections $\alpha \in H_f \longmapsto \sigma\circ\alpha \in H_{\sigma}(f)$, $\psi \in \mathcal{I}_f \longmapsto \sigma\circ\psi \in \mathcal{I}_{\sigma}(f)$, by Corollaries \ref{cuspidal-alpha-phi} and \ref{local-epsilon-modular-galois}.
\item For any $(f,\chi) \in \mathscr{C}$ and $\alpha \in \mathcal{D}$, there are bijections $\beta \in H_f \rar \alpha\beta(\det)\in H_{f \otimes \alpha}$ and $\psi \in \mathcal{I}_f \longmapsto \beta(\det)\psi \in \mathcal{I}_{f \otimes \beta}$. This is a consequence of Corollaries \ref{twisting-Cf} and \ref{cuspidal-twist}. 
\end{itemize} }

\lem[everything-self-dual-again]{Let $(f,\chi) \in \mathscr{S} \cup \mathscr{P}' \cup \mathscr{C}'_1 \cup \mathscr{C}_M$. 
\begin{itemize}[noitemsep,label=$-$]
\item For any $\psi \in H_f(\rho)$, the newform $f \otimes \psi$ has trivial character\footnote{As always, we identify characters of $G_{\Q}$ and Dirichlet characters by means of the cyclotomic character $G_{\Q} \rar \hat{\Z}^{\times}$.}. 
\item For any $\psi \in \mathcal{I}_f$, one has $\left[\mrm{Ind}_K^{\Q}{\psi(\rho_K)}\right]^{\ast} \simeq \left[\mrm{Ind}_K^{\Q}{\psi(\rho_K)}\right] \otimes \chi(\omega_p)$.
\item For any $\psi \in \mathcal{I}_f$, one has $\det{\mrm{Ind}_K^{\Q}{\psi(\rho_K)}}=\epsilon_K\overline{\chi}(\omega_p)$. 
\end{itemize}
}

\demo{The first point follows directly from the description of $H_f(\rho)$: the square of any $\alpha \in H_f(\rho)$ is $\overline{\chi}(\omega_p)$. 
For the second point, since $G_K=\rho^{-1}(C)$ and $N/C$ acts on $C$ by $z \longmapsto z^p$, one has 
\begin{align*}
\left[\mrm{Ind}_K^{\Q}{\psi(\rho_K)}\right]^{\ast} &\simeq \left[\mrm{Ind}_K^{\Q}{\psi^p(\rho_K)}\right]^{\ast} \simeq \left[\mrm{Ind}_K^{\Q}{\overline{\psi}^p(\rho_K)}\right]\\
&\simeq \left[\mrm{Ind}_K^{\Q}{\psi(\rho_K) \cdot \overline{\psi}(\rho_K^{p+1})}\right] \simeq \left[\mrm{Ind}_K^{\Q}{\psi(\rho_K) \cdot \overline{\psi}(\omega_p^{-1}I_2)}\right]\\
&\simeq \psi(\omega_p I_2) \otimes \left[\mrm{Ind}_K^{\Q}{\psi(\rho_K)}\right] \simeq \chi(\omega_p) \otimes \left[\mrm{Ind}_K^{\Q}{\psi(\rho_K)}\right].
\end{align*}
For the third point, the restriction to $G_K$ of $\det{\mrm{Ind}_K^{\Q}{\psi(\rho_K)}}$ is the product of the conjugates under $G_{\Q}/G_K$ of $\psi(\rho_K)$, so it is $\psi(\rho_K)\cdot \psi^p(\rho_K)=\psi(\rho_K^{p+1})=\psi(\omega_p^{-1}I_2)=\overline{\chi}(\omega_p)$. Therefore, $\det{\mrm{Ind}_K^{\Q}{\psi(\rho_K)}} \in \{\overline{\chi}(\omega_p),\epsilon_K\overline{\chi}(\omega_p)\}$. Moreover, the image under $\det{\mrm{Ind}_K^{\Q}{\psi(\rho_K)}}$ of the complex conjugation $c$ is a symmetry with trace $0$ (since $K$ is imaginary), so the eigenvalues of the complex conjugation are $1$ and $-1$ (with multiplicity on each), and its determinant is $-1$. Since $\overline{\chi}(\omega_p(c))=1$ and $\epsilon_K(c)=-1$, the determinant has to be $\epsilon_K\overline{\chi}(\omega_p)$. 
}

We define the group homomorphism $\tilde{\rho}: (M,\sigma) \in N\rtimes_{\rho} G_{\Q} \longmapsto M\rho(\sigma) \in N$ and $\tilde{\rho_K}$ as the induced morphism $\tilde{\rho}^{-1}(C) \rar C$. 

Recall the following objects from Chapter \ref{tate-modules-twists}, introduced in Definitions \ref{notation-T1Gamma}, \ref{notation-T1rho}.

\medskip

\nott{We call $\mathbb{T}_{1,N}(N)$ the subalgebra of $\mrm{End}\left(\prod_{\xi \in \mu_p^{\times}(\Qbar)}{J_{\rho,\xi}(p)_{\Qbar}}\right)$ generated by the $T_nu$, where $n \geq 1$ is an integer coprime to $p$, $u \in \Z[N]$ commutes to $N$ and is a linear combination of matrices of determinant $n^{-1} \in \F_p^{\times}$. This subalgebra is contained in $\bigoplus_{\xi \in \mu_p^{\times}(\Qbar)}{\mrm{End}(J_{\rho,\xi}(p))}$.

We denote by $\mathbb{T}_{1,N,\rho}$ the subalgebra of $\mathbb{T}[N \rtimes_{\rho} G_{\Q}]$ generated by the $(T_n \cdot (mI_2))(M,\sigma)$, for $n, m \geq 1, M \in N, \sigma \in G_{\Q}$ such that $m^2n\det{M} \equiv 1 \pmod{p}$. It acts by group homomorphisms on $J_{\rho,\xi}(p)(\Qbar)$ for every $\xi \in \mu_p^{\times}(\Qbar)$. }

\bigskip

\prop[decomposition-ncartan-xrho-T11rho]{For any $\xi \in \mu_p^{\times}(k_s)$, for any $(f,\chi) \in \mathscr{S}\cup\mathscr{P}\cup\mathscr{C}'_1 \cup\mathscr{C}_M$, the $f$-eigenspace in $\Tate{\ell}{J_{\rho,\xi}(p)} \otimes_{\Z_{\ell}} F_{\mathfrak{l}}$ is isomorphic as a $\mathbb{T}_{1,N,\rho}\otimes_{\Z} F_{\mathfrak{l}}$-module to the restriction to $\mathbb{T}_{1,N,\rho}\otimes F_{\mathfrak{l}}$ of the following $(\mathbb{T}\otimes F_{\mathfrak{l}})[N \rtimes_{\rho} G_{\Q}]$-module: 
\[\left[\bigoplus_{\alpha \in H_f}{V_{f,\mathfrak{l}} \otimes \alpha(\tilde{\rho})}\right] \oplus \left[\bigoplus_{\psi \in \mathcal{I}_f/\sim}{V_{f,\mathfrak{l}} \otimes \mrm{Ind}_{\tilde{\rho}^{-1}(C)}^{N\rtimes_{\rho} G_{\Q}}{\psi(\tilde{\rho_K})}}\right],\] where $T_n \in \mathbb{T}$ (resp. $nI_2 \in\mathbb{T}$) for $n \geq 1$ coprime to $p$ acts by $a_n(f)$ (resp. $\chi(n)$). 
}

\bigskip

\demo{\emph{First case: $(f,\chi) \notin \mathscr{C}_M$.} 

By Corollary \ref{connected-action-T11rho}, the $f$-eigenspace of $\Tate{\ell}{J_{\rho,\xi}(p)} \otimes_{\Z_{\ell}} F_{\mathfrak{l}}$ is isomorphic as a $\mathbb{T}_{1,\GL{\F_p},\rho}\otimes F_{\mathfrak{l}}$-module to the restriction to $\mathbb{T}_{1,\GL{\F_p},\rho}\otimes F_{\mathfrak{l}}$ of $V_{f,\mathfrak{l}} \otimes R_f(\tilde{\rho})$, where $R_f$ is $\mrm{St}$ if $(f,\chi) \in \mathscr{S}$, $\pi(\mathbf{1},\chi)$ if $(f,\chi) \in \mathscr{P}'$, and $C_f$ if $(f,\chi) \in \mathscr{C}'_1$, and $T_n$ (resp. $nI_2$) for any $n \geq 1$ coprime to $p$ acts by $a_n(f)$. 

The goal is then simply to describe the decomposition of $(R_f)_{|N}$ into irreducible representations. It is enough to show that $\mrm{St}_{|N}$ (resp. $\pi(1,\chi)_{|N}$, resp. $C_f(\rho)_{|N}$) is isomorphic to $\bigoplus_{\psi \in H_f}{\psi} \oplus \bigoplus_{\psi \in \mathcal{I}_f/\sim}{\mrm{Ind}_C^{N}{\psi}}$ (the latter sum is well-defined because if $\psi,\psi' \in \mathcal{I}_f$ are such that $\psi \sim \psi'$, then $\psi$ is conjugate to $\psi'$ by $N/C$). 

When $(f,\chi) \in \mathscr{S}$, this is a consequence from Lemma \ref{steinberg-on-norm-cartan}. When $(f,\chi) \in \mathscr{P}'$, this also follows directly from Lemma \ref{principal-on-norm-cartan}. 

When $(f,\chi)=(f,\mathbf{1}) \in \mathscr{C}'_1$, we apply Proposition \ref{cusp-on-norm-cartan}. The summands made from induced representations correspond to each other, so what we need to prove is that the (abelian) characters match. Let $\phi_f: C \rar \C^{\times}$ be one of the two characters attached to $C_f$, then, by Proposition \ref{cusp-on-norm-cartan}, the two characters appearing in $(C_f)_{|N}$ are the $\psi\epsilon^r$, where $\psi \in \mathcal{D}$ runs through the characters such that $\psi^2(u)=\phi_f(uI_2)$ for each $u \in \F_p^{\times}$, and $r$ is such that $(-1)^r\psi(-1)$ is $-1$ if $\phi_f(u)=\psi(u^2)$ for every $u \in \sqrt{\F_p^{\times} I_2}$, and $1$ otherwise. Fix some $a \in \F_p^{\times} \backslash \F_p^{\times 2}$, then it is clear that $(-1)^r\psi(-1)=-\frac{\phi_f(\sqrt{a})}{\psi(a)}$, so $(-1)^r=-\frac{\phi_f(\sqrt{a})}{\psi(-a)}=\lambda_p(f \otimes \alpha) \in \{\pm 1\}$ by Corollary \ref{various-computations-cuspidal}. \\

\emph{Second case: $(f,\chi) \in \mathscr{C}_M$.}

By Corollary \ref{decomposition-xrho-connected-sl2}, as a $\mathbb{T}_1[\SL{\F_p} \rtimes_{\rho} G_{\Q}]$-module, the $f$-eigenspace of $\Tate{\ell}{J_{\rho,\xi}(p)} \otimes_{\Z_{\ell}} F_{\mathfrak{l}}$ is isomorphic, for some sign $r$ only depending on $\xi^{\F_p^{\times 2}}$, to
\[\mrm{Ind}_{\SL{\F_p} \rtimes_{\rho} G_{\Q(\sqrt{-p})}}^{\SL{\F_p} \rtimes_{\rho} G_{\Q}}{\psi_{f,r,\mathfrak{l}} \otimes C^+(\tilde{\rho_K})} \simeq \mrm{Ind}_{\SL{\F_p} \rtimes_{\rho} G_{\Q(\sqrt{-p})}}^{\SL{\F_p} \rtimes_{\rho} G_{\Q}}{\psi_{f,+,\mathfrak{l}} \otimes C^r(\tilde{\rho_K})},\]
where the extensions of $\tilde{\rho}$ (resp. $\tilde{\rho_K}$) to $\SL{\F_p} \rtimes_{\rho} G_{\Q}$ (resp. $\SL{\F_p} \rtimes_{\rho} G_{\Q(\sqrt{-p})}$) are given by $(M,g) \longmapsto M\rho(g) \in \GL{\F_p}$. 

For any $\psi,\psi' \in \mathcal{I}_f$, the character $\psi^2_{|C \cap \SL{\F_p}}$ is by definition nontrivial, and $\psi \sim \psi'$ if and only if the restrictions of $\psi,\psi'$ to $C \cap \F_p^{\times}\SL{\F_p}$ are equal or $p$-th powers one of the other, that is, if and only if the restrictions of $\psi,\psi'$ to $C \cap \F_p^{\times}\SL{\F_p}$ are conjugate under the action of $N$. Therefore, by Proposition \ref{cusp-sl-norm-cartan}, one has
\[C^r_{|N \cap \F_p^{\times}\SL{\F_p}} \simeq \left[\epsilon^{\frac{p-3}{4}} \oplus \bigoplus_{\chi \in \mathcal{I}_f/\sim}{\left(\mrm{Ind}_C^N{\chi}\right)}\right]_{|N \cap \F_p^{\times}\SL{\F_p}}, \] and the conclusion follows since $V_{f,\mathfrak{l}} = \mrm{Ind}_{\Q(\sqrt{-p})}^{\Q}{\psi_{f,r,\mathfrak{l}}}$. 
}

\prop[tate-module-xrho-nsplit-cartan]{Let $\mathfrak{l} \subset \OO_F$ be a maximal ideal with residue characteristic $\ell$. 
For any $\xi \in \mu_p^{\times}(\Qbar)$, and any $(f,\chi) \in \mathscr{S}\cup\mathscr{P}\cup\mathscr{C}'_1 \cup\mathscr{C}_M$, the $f$-eigenspace in $\Tate{\ell}{J_{\rho,\xi}(p)} \otimes_{\Z_{\ell}} F_{\mathfrak{l}}$ is isomorphic (as a $F_{\mathfrak{l}}[G_{\Q}]$-module) to 
\[\left[\bigoplus_{\alpha \in H_f(\rho)}{V_{f,\mathfrak{l}} \otimes \alpha}\right] \oplus \left[\bigoplus_{\psi \in \mathcal{I}_f/\sim}{V_{f,\mathfrak{l}} \otimes \mrm{Ind}_K^{\Q}{\psi(\rho_K)}}\right].\]
The summands of this decomposition are isomorphic to their Tate-twisted duals and do not depend on the choice of representatives for $\mathcal{I}_f/\sim$. Each summand is stable under the (geometric) action of $\mathbb{T}_{1,N} \otimes F_{\mathfrak{l}}$, is isomorphic to the direct sum of two copies of the same irreducible $\mathbb{T}_{1,N} \otimes F_{\mathfrak{l}}$-module, and two distinct summands have non-isomorphic attached $\mathbb{T}_{1,N} \otimes F_{\mathfrak{l}}$-modules. 
Moreover, every summand is irreducible as a $F_{\mathfrak{l}}[G_{\Q}]$-module, except in the following cases:
\begin{itemize}[noitemsep,label=$-$]
\item $\psi \in \mathcal{I}_f$ is such that $\psi(\rho_K)^{p-1}=1$. In this case, $\mrm{Ind}_K^{\Q}{\psi(\rho_K)}$ is the sum of two finite order characters $\alpha_{\psi},\beta_{\psi}: G_{\Q} \rar F^{\times}$. Moreover, $\alpha_{\psi}\beta_{\psi}=\epsilon_K\overline{\chi}(\omega_p)$ and $\alpha_{\psi}\beta_{\psi}^{-1}=\epsilon_K$. In this case, $V_{f,\mathfrak{l}} \otimes \alpha_{\psi}, V_{f,\mathfrak{l}} \otimes \beta_{\psi}$ are irreducible and isomorphic to their Tate-twisted duals. 
\item When $(f,\mathbf{1}) \in \mathscr{C}_M$, $\theta: G_K \rar F_{\mathfrak{l}}^{\times}$ one of the two characters attached to $f$, and $\psi \in \mathcal{I}_f$ such that $\psi(\rho_K)^{p-1} \neq 1$:
\begin{itemize}[noitemsep,label=\tiny$\bullet$] 
\item If $K \neq \Q(\sqrt{-p})$ and $\psi(\rho_K)^{4}=1$, $V_{f,\mathfrak{l}} \simeq \mrm{Ind}_K^{\Q}{\psi(\rho_K)}$ is the direct sum of two irreducible representations $\mrm{Ind}_{\Q(\sqrt{-p})}^{\Q}{\theta \otimes \alpha_{\psi}}$, $\mrm{Ind}_{\Q(\sqrt{-p})}^{\Q}{\theta \otimes \beta_{\psi}}$, where the characters $\alpha_{\psi},\beta_{\psi}$ are the two lifts to $G_{\Q(\sqrt{-p})}$ of the quadratic character $\psi(\rho_K)_{|G_{K(\sqrt{-p})}}$. Taking the Tate-twisted dual exchanges these two representations.  
\item When $(f,\mathbf{1}) \in \mathscr{C}_M$ (so that $p \equiv 3 \pmod{4}$), $K=\Q(\sqrt{-p})$ and $\psi \in \mathcal{I}_f$ is such that $\psi(\rho_K)^{p-1} \neq 1$, let $\theta: G_K \rar F_{\mathfrak{l}}^{\times}$ denote one of the characters attached to $f$. Then $V_{f,\mathfrak{l}}$ decomposes as the direct sum of the two following irreducible representations, both of which are isomorphic to their Tate-twiste duals: $\mrm{Ind}_K^{\Q}{\theta \otimes \psi(\rho_K)}$ and $\mrm{Ind}_K^{\Q}{\theta \otimes \psi^p(\rho_K)}$. 
\end{itemize}
\end{itemize}
}

\demo{Let $(f,\chi) \in \mathscr{S} \cup \mathscr{P}' \cup \mathscr{C}'_1 \cup \mathscr{C}_M$.

The group $N/C$ acts on $C$ by $z \mapsto z^p$, so, for any $\psi \in \mathcal{I}_f$, $G_{\Q}/G_K$ exchanges $\psi(\rho_K)$ and $\psi^p(\rho_K)$: therefore $\mrm{Ind}_K^{\Q}{\psi(\rho_K)} \simeq \mrm{Ind}_K^{\Q}{\psi^p(\rho_K)}$. 

Furthermore, let $\alpha: C \rar \{\pm 1\}$ be the unique character of order two. Then $\alpha(\rho_K)=\left(\frac{\det{\rho_K}}{p}\right)$ is the restriction to $G_K$ of the character $\left(\frac{\omega_p}{p}\right)$, and therefore $\mrm{Ind}_K^{\Q}{(\alpha\psi)(\rho_K)} \simeq \left(\frac{\omega_p}{p}\right) \otimes \mrm{Ind}_K^{\Q}{\psi(\rho_K)}$. Since $V_{f,\mathfrak{l}} \otimes \left(\frac{\omega_p}{p}\right) \simeq V_{f,\mathfrak{l}}$ if $(f,\mathbf{1}) \in \mathscr{C}_M$, it follows that $V_{f,\mathfrak{l}} \otimes \mrm{Ind}_K^{\Q}{(\alpha\psi)(\rho_K)} \simeq V_{f,\mathfrak{l}} \otimes \mrm{Ind}_K^{\Q}{\psi(\rho_K)}$. 

This proves that the given decomposition is independent from a choice of representatives for $\mathcal{I}_f/\sim$. By Lemma \ref{everything-self-dual-again}, all the summands are isomorphic to their Tate-twisted duals (since $V_{f,\mathfrak{l}} \simeq V_{f,\mathfrak{l}}^{\ast}(1) \otimes \chi(\omega_p)$). 

The decomposition itself is a specialization of Proposition \ref{decomposition-ncartan-xrho-T11rho}, which also implies the claimed stability under the geometric action of $\mathbb{T}_{1,N} \otimes F_{\mathfrak{l}}$. The Proposition also implies that the action of $\mathbb{T}_{1,N}$ on the summand $V_{f,\mathfrak{l}} \otimes \alpha(\rho)$ for $\alpha \in H_f$ (resp. $V_{f,\mathfrak{l}} \otimes \mrm{Ind}_C^N{\psi(\rho_K)}$ for $\psi \in \mathcal{I}_f$) is the restriction to $\mathbb{T}_{1,N}$ of the action of $\mathbb{T}[N]$ on $\alpha^{\oplus 2}$ (resp. $\left[\mrm{Ind}_C^N{\psi}\right]^{\oplus 2}$), where $T_n$ (resp. $mI_2$) acts by multiplication by $a_n(f)$. 

The claims that we need to prove are the following:
\begin{enumerate}[noitemsep,label=(\roman*)]
\item\label{tmxsc-1} For $\alpha \in H_f$, the datum of the $a_n(f)\chi(m)\alpha(M)$ for $n,m \geq 1$ and $M \in N$ such that $nm^2\det{M} \equiv 1\pmod{p}$ determines $\alpha$.  
\item\label{tmxsc-2} For any $\psi \in \mathcal{I}_f$, $\mrm{Ind}_C^N{\psi}$ is an irreducible $(\mathbb{T}_{1,N}\otimes F_{\mathfrak{l}})$-module, 
\item\label{tmxsc-3} For any $\psi,\psi' \in \mathcal{I}_f$, if $\mrm{Ind}_C^N{\psi}$ and $\mrm{Ind}_C^N{\psi'}$ are isomorphic $(\mathbb{T}_{1,N} \otimes F_{\mathfrak{l}})$-modules, then $\psi\sim\psi'$.  
\end{enumerate}

A result of Serre \cite[Corollaire 2 au Th\'eor\`eme 15]{Serre-Cebotarev} states that if $f$ is not CM, then, for any $u \in \F_p^{\times}$, there are infinitely many prime numbers $q \equiv u \pmod{p}$ such that $a_q(f) \neq 0$. This result implies directly \ref{tmxsc-1} when $(f,\chi) \notin \mathscr{C}_M$, and $H_f$ is a singleton when $f$ is CM, so \ref{tmxsc-1} is always true.

This result also shows that when $(f,\chi) \notin \mathscr{C}_M$ and under the assumptions of \ref{tmxsc-3}, the representations $\mrm{Ind}_C^N{\psi}$ and $\mrm{Ind}_C^N{\psi'}$ of $N$ over $F$ are isomorphic, which implies that $\psi' \in \{\psi,\psi^p\}$ thus $\psi,\psi'$. 

Let us prove \ref{tmxsc-3} when $(f,\chi) \in \mathscr{C}_M$. Indeed, under the assumptions of \ref{tmxsc-3}, the representations $\mrm{Ind}_C^N{\psi},\mrm{Ind}_C^N{\psi'}$ are isomorphic when reduced to $N\cap \SL{\F_p}$. This implies that the characters $\psi_{|C \cap \SL{\F_p}},\psi'_{|C \cap \SL{\F_p}}$ are conjugate by $N$. This implies that $\psi \sim \psi'$. 

Now, let us prove \ref{tmxsc-2}. Let $\psi \in \mathcal{I}_f$ be such that $\mrm{Ind}_C^N{\psi}$ is a reducible $(\mathbb{T}_{1,N} \otimes F_{\mathfrak{l}})$-module: its semi-simplification is the sum of characters $\alpha,\beta: \mathbb{T}_{1,N} \otimes F_{\mathfrak{l}} \rar F_{\mathfrak{l}}$; both of them map $T_n \cdot (mI_2)$ to $a_n(f)\chi(m)$ for any $n,m \geq 1$ such that $nm^2\equiv 1\pmod{p}$. Therefore, $\mrm{Ind}_C^N{\psi}$ is reducible as a $C \cap \SL{\F_p}$-module, which implies that $\psi^{p-1}_{|C \cap \SL{\F_p}}=\mathbf{1}$. Since, by construction, $\psi_{|\F_p^{\times}}^{p-1}=1$, we see that $\psi^{2(p-1)}=\mathbf{1}$. In particular, if $\chi=\mathbf{1}$, then $\psi^4=\mathbf{1}$, and therefore $(f,\chi) \notin \mathscr{C}_M$.  

Let $M \in N$; by Serre's result, there exists a prime number $n \geq 1$ such that $n\det{M} \equiv 1 \pmod{p}$ and $a_n(f) \neq 0$. We let $M$ act on $\mrm{Ind}_C^N{\psi}$ (recall its only structure is that of a $\mathbb{T}_{1,N}$-module) by $a_n(f)^{-1}\cdot T_nM$. Since $(T_nM)^{p^2-1}$ acts by $T_n^{p^2-1}=a_n(f)^{p^2-1} \neq 0$ (by Lemma \ref{hecke-tn-group}), $T_n,M$ are automorphisms. Let $m \geq 1$ be an integer such that $m\det{M} \equiv 1\pmod{p}$, let us show that $T_mM \in \mathbb{T}_{1,N}$ acts by $a_m(f)M$. 

Indeed, by Serre's result, there exists a prime $q \nmid nmp$ such that $a_q(f) \neq 0$ and $q \equiv \det{M} \pmod{p}$. Then $(T_nT_q)(T_mM)=T_{mq}(T_nM)$ acts by $a_{nq}(f)(T_mM)=a_n(f)a_q(f)(T_mM)$ and $a_{mq}(f)a_n(f)M=a_n(f)a_m(f)a_q(f)M$, so that $T_mM$ acts by $a_m(f)M$. 

It is straightforward to check that this construction defines a group action of $N$ on the $(\mathbb{T}_{1,N} \otimes F_{\mathfrak{l}})$-module $\mrm{Ind}_C^N{\psi}$, whose traces are exactly that of $N$ on $\mrm{Ind}_C^N{\psi}$; since this $(\mathbb{T}_{1,N}\otimes F_{\mathfrak{l}})$-module was reducible, it implies that $\mrm{Ind}_C^N{\psi}$ is reducible as a $F_{\mathfrak{l}}[N]$-module, which is impossible by definition of $\mathcal{I}_f$. \\

Now, we discuss the irreducibility conditions for the Galois representations. 

Let $\psi \in \mathcal{I}_f$. By Mackey's criterion \cite[Proposition 23]{SerreLinReps}, $\mrm{Ind}_K^{\Q}{\psi(\rho_K)}$ is irreducible if and only if $\psi(\rho_K)$ is disjoint from its conjugate under $G_{\Q}/G_K$, which is $\psi^p(\rho_K)$, as we saw. This is equivalent to $\psi(\rho_K) \neq \psi^p(\rho_K)$, or equivalently $\psi^{p-1}(\rho_K) \neq 1$. 

When $\psi^{p-1}(\rho_K)=1$, $\mrm{Ind}_K^{\Q}{\psi(\rho_K)}$ is the direct sum of two characters $\alpha_{\psi},\beta_{\psi}$ whose restriction to $G_K$ is $\psi(\rho_K)$. Therefore, $\alpha_{\psi}=\beta_{\psi}\epsilon_K$, and $\alpha_{\psi}\beta_{\psi}=\epsilon_K\overline{\chi}(\omega_p)$ (by the third claim in Lemma \ref{everything-self-dual-again}), so that $\alpha_{\psi}^2=\beta_{\psi}^2=\overline{\chi}(\omega_p)$. Therefore the representations $V_{f,\mathfrak{l}} \otimes \alpha_{\psi}, V_{f,\mathfrak{l}} \otimes \beta_{\psi}$ are irreducible and isomorphic to their Tate-twisted duals, and their sum is $V_{f,\mathfrak{l}} \otimes \mrm{Ind}_K^{\Q}{\psi(\rho_K)}$. 

When $\psi^{p-1}(\rho_K) \neq 1$ and $(f,\chi) \notin \mathscr{C}_M$, by Proposition \ref{tensor-product-with-artin-irreducible}, $V_{f,\mathfrak{l}} \otimes \mrm{Ind}_K^{\Q}{\psi(\rho_K)}$ is irreducible. 

Now, we assume that $\psi^{p-1}(\rho_K) \neq 1$ and $(f,\chi)=(f,\mathbf{1}) \in \mathscr{C}_M$. Let $\theta: G_{\Q(\sqrt{-p})} \rar F_{\mathfrak{l}}^{\times}$ be one of the characters attached to $f$ and $\theta'$ be its conjugate under $G_{\Q}/G_K$, so that \[V_{f,\mathfrak{l}} \otimes \mrm{Ind}_K^{\Q}{\psi(\rho_K)} \simeq \mrm{Ind}_{\Q(\sqrt{-p})}^{\Q}{\left\{\theta \otimes \left[\mrm{Ind}_{K}^{\Q}{\psi(\rho_K)}\right]_{|G_{\Q(\sqrt{-p})}}\right\}}.\]

 By Mackey's criterion, this representation is irreducible if and only if $\theta \otimes \left[\mrm{Ind}_{K}^{\Q}{\psi(\rho_K)}\right]_{|G_{\Q(\sqrt{-p})}}$ is irreducible and is not isomorphic to its conjugate under $G_{\Q}/G_{\Q(\sqrt{-p})}$, which is the representation $\theta' \otimes \left[\mrm{Ind}_{K}^{\Q}{\psi(\rho_K)}\right]_{|G_{\Q(\sqrt{-p})}}$. Over a finite index normal subgroup of $G_{\Q(\sqrt{-p})}$, these representations are isomorphic to $\theta^{\oplus 2},(\theta')^{\oplus 2}$ respectively. Since $\frac{\theta}{\theta'}$ has infinite image by \cite[Theorem 4.3]{Antwerp5-Ribet}, $\theta \otimes \left[\mrm{Ind}_{K}^{\Q}{\psi(\rho_K)}\right]_{|G_{\Q(\sqrt{-p})}}$ is not isomorphic to its conjugate under $G_{\Q}/G_{\Q(\sqrt{-p})}$. 

Thus $V_{f,\mathfrak{l}} \otimes \mrm{Ind}_K^{\Q}{\psi(\rho_K)}$ is irreducible if and only if $I := \left[\mrm{Ind}_{K}^{\Q}{\psi(\rho_K)}\right]_{|G_{\Q(\sqrt{-p}}}$ is irreducible. This condition is false when $K=\Q(\sqrt{-p})$, which we will treat separately, so we assume for the time being that $K \neq \Q(\sqrt{-p})$. Then $I \simeq \mrm{Ind}_{K(\sqrt{-p})}^{\Q(\sqrt{-p})}{\psi((\rho_K)_{|K(\sqrt{-p})})}$ is reducible, by Mackey's criterion, if and only if $\psi((\rho_K)_{|G_{K(\sqrt{-p})}})=\psi^p((\rho_K)_{|G_{K(\sqrt{-p})}})$, which is equivalent to $\psi^{p-1}((\rho_K)_{|G_{K(\sqrt{-p})}})=1$. Since $G_{K(\sqrt{-p})}=\rho^{-1}(C \cap \F_p^{\times}\SL{\F_p})$, this is equivalent to $\psi^{2(p-1)}(\rho_K)=1$. 

If $I$ is reducible, since $\psi^{p+1}=1$, one has $\psi^4(\rho_K)=1$. Let $\alpha,\beta$ be the two lifts of $\psi(\rho_{K})_{|G_{K(\sqrt{-p})}}$ to $G_{\Q(\sqrt{-p})}$, so that $I \simeq \alpha\oplus\beta$, $\alpha=\beta\epsilon_K$ and $\alpha\beta=\det{I}=\epsilon_K$. In particular, $\alpha,\beta$ are quadratic characters. Suppose that $\alpha$ is its own conjugate under $\mrm{Gal}(\Q(\sqrt{-p})/\Q)$: then it extends to a character $\alpha_0: G_{\Q} \rar \C^{\times}$ whose kernel contains $\rho^{-1}(4C)$. Thus $\alpha_0$ factors through $\rho$. Let $g \in G_{\Q} \backslash G_K$, then $\rho(g)^2$ is scalar so $g^2 \in G_{K(\sqrt{-p})}$ and $\alpha_0(g)^2=\alpha_0(g^2)=\psi(\rho_K(g)^2)=1$, so that $\alpha_0^2=1$. In particular, the character $g \in \rho(G_K) \longmapsto \psi(g)\alpha_0(g)^{-1}$ is trivial on the index-two subgroup $\rho(G_K) \cap \F_p^{\times}\SL{\F_p}$ but takes the values $\pm i$, whence a contradiction. Thus $\mrm{Gal}(\Q(\sqrt{-p})/\Q)$ exchanges $\alpha$ and $\beta$, and 
\[\left[\mrm{Ind}_{\Q(\sqrt{-p})}^{\Q}{\theta \otimes \alpha}\right]^{\ast}(1) \simeq \mrm{Ind}_{\Q(\sqrt{-p})}{\Q}{\theta^{-1}(1) \otimes \alpha^{-1}}\simeq \mrm{Ind}_{\Q(\sqrt{-p})}^{\Q}{\theta' \otimes \alpha} \simeq \mrm{Ind}_{\Q(\sqrt{-p})}^{\Q}{\theta \otimes \beta}.\] 

Now assume that $K=\Q(\sqrt{-p})$. Then 
\begin{align*}
V_{f,\mathfrak{l}} \otimes \mrm{Ind}_K^{\Q}{\psi(\rho_K)} &\simeq \mrm{Ind}_K^{\Q}{\theta \otimes \left[\psi(\rho_K) \oplus \psi^p(\rho_K)\right]}\\
& \simeq \left[\mrm{Ind}_K^{\Q}{\theta \otimes \psi(\rho_K)}\right] \oplus \left[\mrm{Ind}_K^{\Q}{\theta \otimes \psi^p(\rho_K)}\right].
\end{align*}

Since $(\theta\psi(\rho_K))^{\ast}(1) = \theta'\psi^{-1}(\rho_K)=\theta'\psi^p(\rho_K)$ is the conjugate under $\mrm{Gal}(K/\Q)$ of $\theta\psi(\rho_K)$, the summands above are the Tate twists of their duals. Since $\theta \otimes \psi(\rho_K)$ is isomorphic over some open subgroup of $G_K$ to $\theta$, and since $\frac{\theta}{\theta'}$ has infinite image, Macley's criterion implies that both summands are irreducible. 

}

\prop[functional-equations-jrho-ncartan]{Let $V$ be any of the following compatible systems:
\begin{enumerate}[noitemsep,label=(\alph*)]
\item\label{fencart-1} $(V_{f,\mathfrak{l}} \otimes \alpha)_{\mathfrak{l}}$ for $(f,\chi) \in \mathscr{E} := \mathscr{S} \cup \mathscr{P}'\cup\mathscr{C}'_1 \cup \mathscr{C}_M$ and $\alpha \in H_f(\rho)$,
\item\label{fencart-2} $(V_{f,\mathfrak{l}} \otimes \alpha)_{\mathfrak{l}}$ for $(f,\chi) \in \mathscr{E}$ and $\alpha: G_{\Q} \rar F^{\times}$ extending $\psi(\rho_K)$ for some $\psi \in \mathcal{I}_f$ such that $\psi(\rho_K)^{p-1}=1$,
\item\label{fencart-3} $\left(V_{f,\mathfrak{l}} \otimes \mrm{Ind}_K^{\Q}{\psi(\rho_K)}\right)_{\mathfrak{l}}$ for $(f,\chi) \in \mathscr{E}$ and $\psi \in \mathcal{I}_f$,
\item\label{fencart-4} if $K=\Q(\sqrt{-p})$ and $\theta$ is a compatible system of characters of $G_K$ attached to $(f,\chi) \in \mathscr{C}_M$, $\mrm{Ind}_K^{\Q}{\left[(\theta\psi)(\rho_K)\right]}$ for $\psi \in \mathcal{I}_f$.
\item\label{fencart-5} if $K \neq \Q(\sqrt{-p})$, when $\theta$ is a compatible system of characters of $G_{\Q(\sqrt{-p})}$ attached to $(f,\chi) \in \mathscr{C}_M$, $\psi \in \mathcal{I}_f$ is such that $\psi(\rho_K)^4=1$ and $\psi(\rho_K)^{p-1} \neq 1$, the compatible system $\left(\mrm{Ind}_{\Q(\sqrt{-p})}^{\Q}{\theta \otimes \alpha}\right)_{\mathfrak{l}}$, where $\alpha: G_{\Q(\sqrt{-p})} \rar \C^{\times}$ extends $\psi(\rho_K)_{|G_{K(\sqrt{-p})}}$.  
\end{enumerate}
Let $M$ be the conductor of $V$. Then \[\Lambda(V,s) := M^{s/2}[2(2\pi)^{-s}\Gamma(s)]^{\dim{V}/2}L(V,s)\] is defined on a right half-plane, extends to an entire function and satisfies the functional equation \[\Lambda(V,s)=\varepsilon(V)\Lambda(V^{\ast}(1),2-s),\] where $\varepsilon(V) \in \C^{\times}$. 
Moreover, $\varepsilon(V)=\prod_{v}{\varepsilon_v}$, where 
\begin{itemize}[noitemsep,label=\tiny$\bullet$]
\item In cases \ref{fencart-1} and \ref{fencart-2}, $v$ runs through the places of $\Q$ and $\varepsilon_v=\varepsilon_v(f \otimes \alpha)$ when $v$ is finite, and $\varepsilon_v(\infty)=(-i)^2=-1$. 
\item In case \ref{fencart-3}, $v$ runs through all places of $\Q$ and $\varepsilon_v=\varepsilon(D_{f,v} \otimes \Delp{v}{\left(\mrm{Ind}_K^{\Q}{\psi(\rho_K)}\right)^{\ast}},1,\theta_{\Q_v})$.
\item In case \ref{fencart-4}, $v$ runs through the places of $\Q(\sqrt{-p})$ and $\varepsilon_v=\varepsilon(\Delp{v}{(\theta\psi(\rho_K))^{\ast}},1,\theta_{K_v})$ for $v$ finite and $\varepsilon_{\infty}(z^{-1},1,\theta_{\C})$. 
\item In case \ref{fencart-5}, $v$ runs through the places of $\Q(\sqrt{-p})$ and $\varepsilon_v=\varepsilon(\Delp{v}{(\theta\alpha)^{\ast}},1,\theta_{K_v})$ for $v$ finite and $\varepsilon_{\infty}(z^{-1},1,\theta_{\C})$.
\end{itemize}
In all cases except \ref{fencart-5}, one has $V^{\ast}(1) \simeq V$ and $\varepsilon(V) \in \{\pm 1\}$.  
}

The definition of the $D_{f,v}$ as Weil-Deligne representations is recalled in Section \ref{langlands-correspondence-refresher}. 

\demo{In the first two cases, $V=V_{f \otimes \alpha,\mathfrak{l}}$ since $\alpha$ identifies with a Dirichlet character by the Kronecker-Weber theorem. The conclusion follows by Proposition \ref{functional-eqn} (for the correct functional equation), Corollary \ref{local-to-global-constant} (for the constant in the functional equation), and Proposition \ref{local-global-modular-forms} (for the correct Galois-theoretic interpretation).

In the third case, this follows from the Rankin-Selberg functional equation (Proposition \ref{automorphic-times-modform}), and the fact that $\mrm{Ind}_K^{\Q}{\psi(\rho_K)}$ is automorphic, because it is odd (and therefore either attached to a cuspidal weight one modular form, or the direct sum of two characters of $G_{\Q}$ of finite order). In the fourth and fifth cases, the result follows from Proposition \ref{automorphic-times-modform-cmtwist}. 

The claimed self-duality is proved in Proposition \ref{tate-module-xrho-nsplit-cartan}, and it implies that $\varepsilon(V)^2=1$. 

}

\section{Signs of the functional equation for factors of the Tate module of $J_{\rho,\xi}(p)$}
\label{signs-functional-equation-cartan}

\prop[signs-nsplit-Hf]{Let $(f,\chi) \in \mathscr{S} \cup \mathscr{P}' \cup \mathscr{C}'_1 \cup\mathscr{C}_M$, and $\alpha \in H_f(\rho)$. Let $V = \left(V_{f,\mathfrak{l}} \otimes \alpha\right)_{\mathfrak{l} \in \mrm{Max}(\OO_F)}$. The global root number of the compatible system $V$ (as defined in Proposition \ref{functional-equations-jrho-ncartan}) is:
\begin{itemize}[label=\tiny$\bullet$,noitemsep]
\item equal to $-\lambda_p(f \otimes \overline{\psi})$ (resp. $(-1)^{\frac{p+1}{4}}$) if $K$ is ramified at $p$, $(f,\chi) \in \mathscr{C}'_1$ (resp. $(f,\chi) \in \mathscr{C}_M$) and $\alpha \in \epsilon_K^{\Z}\overline{\psi}(\omega_p)$.
\item equal to $\eta(-p)a_p(f)$ if $K$ is ramified at $p$, $(f,\chi) \in \mathscr{S}$ and $p \equiv 3\pmod{4}$, where $\eta$ the quadratic Dirichlet character with conductor coprime to $p$ such that $\epsilon_K\eta$ has conductor $p$. 
\item equal to $-1$ otherwise. 
\end{itemize}
}

\demo{We can write $\alpha=\beta(\omega_p)\gamma(\omega)$, where $\beta \in \mathcal{D}$, $\gamma$ are Dirichlet characters with values in $F$, the conductor of $\gamma$ is coprime to $p$, and $\omega: G_{\Q} \rar \hat{\Z}^{\times}$ is the cyclotomic character. In all cases, $\gamma^2=\mathbf{1}$ and $\beta^2\chi=\mathbf{1}$. Let $P_1$ be the conductor of the newform $f \otimes \beta$ (with trivial character), so $P_1$ is a power of $p$. By Propositions \ref{functional-eqn} and Corollary \ref{global-constant-coprime-twist}, the global root number of $V$ is $(-i)^2\lambda_p(f \otimes \beta)\gamma(-P_1)$. 

Suppose that $\beta \neq \mathbf{1}$ and that $(f,\chi) \in \mathscr{S} \cup \mathscr{P}'$. Then, by Propositions \ref{functional-eqn} and \ref{local-constant-bigtwist}, $P_1=p^2$ and $\lambda_p(f \otimes \beta)=\beta(-1)$, so the global root number of $V$ is $-\alpha(\omega^{-1}(-1))=-\alpha(c)$, where $c$ is the complex conjugation.  

When $(f,\chi) \in \mathscr{S} \cup \mathscr{P}'$, after checking every case, we can see that $\beta = \mathbf{1}$ if and only if $K/\Q$ is ramified at $p$, $p \equiv 3 \pmod{4}$, and $(f,\chi) \in \mathscr{S}$. If one of these conditions is not satisfied, we can check case by case that $\alpha(c)=1$. When $K/\Q$ is ramified at $p$, $p \equiv 3\pmod{4}$ and $(f,\chi) \in \mathscr{S}$, $\beta=\mathbf{1}$, so $P_1=p$, $\lambda_p(f)=-a_p(f)$ by Propositions \ref{local-constant-principal} and \ref{bad-L-factor}, so the global root number of $V$ is $a_p(f)\gamma(p)$. 

When $(f,\chi) \in \mathscr{C}'_1 \cup \mathscr{C}_M$, then $(f \otimes \beta,\mathbf{1}) \in \mathscr{C}$, so $P_1=p^2$, and the global root number of $V$ is $-\gamma(-1)\lambda_p(f \otimes \beta)$.

If $(f,\chi) \in \mathscr{C}'_1$ and $K/\Q$ is inert at $p$, then $\beta=\overline{\psi}$ and $\gamma=\epsilon_K^r$ for some quadratic $\psi \in \mathcal{D}$, and $(-1)^r\lambda_p(f \otimes \overline{\psi})=1$, so that $\gamma(-1)=\lambda_p(f \otimes \beta)$, thus the global root number of $V$ is $-1$. 

If $(f,\chi) \in \mathscr{C}'_1$ and $K/\Q$ is ramified at $p$, $\beta=\overline{\psi}\lambda^r$ and $\gamma=\eta^r$ for some quadratic $\psi \in \mathcal{D}$, $(-1)^r\lambda_p(f \otimes \overline{\psi})=1$, and $\epsilon_K=\lambda\eta$, so, by Corollary \ref{various-computations-cuspidal}, \[-\gamma(-1)\lambda_p(f \otimes \beta)=-\eta(-1)^r(-\lambda(-1))^r\lambda_p(f \otimes \overline{\psi})=-\lambda_p(f \otimes \overline{\psi}).\] 

If $(f,\chi) \in \mathscr{C}_M$, $f=f \otimes \beta$, so we would like to compute $\lambda_p(f)$. Let $E/\Q_p$ be the unique quadratic unramified extension and $\phi: E^{\times}/(1+p\OO_E) \rar \C^{\times}$ be the character such that $\phi(p)=-p$ and $\phi_{|\OO_E^{\times}}$ has order $4$: then, by Lemma \ref{local-global-c}, Proposition \ref{decomposition-xrho-reducible-cm} and Lemma \ref{local-epsilon-fp2}, 
\begin{align*}
\lambda_p(f) &= \varepsilon_p(f) = \varepsilon(\mrm{Ind}_E^{\Q_p}{\phi(\mathfrak{a}_E)},1,\theta_p) = \lambda_{E/\Q_p}(\theta_p)\varepsilon(\phi(\mathfrak{a}_E)\|\cdot\|^{1/2},\frac{1}{2},\theta_E)\\
& = \phi_0(p)\|p\|^{1/2}\frac{\sum_{u \in (\OO_E/p)^{\times}}{\phi^{-1}(u)e^{\frac{2i\pi}{p}\mrm{Tr}_{E/\Q_p}(u)}}}{\left|\sum_{u \in (\OO_E/p)^{\times}}{\phi^{-1}(u)e^{\frac{2i\pi}{p}\mrm{Tr}_{E/\Q_p}(u)}}\right|}\\
&= -\phi(\sqrt{-1})=-(-1)^{\frac{p^2-1}{8}}=-(-1)^{\frac{p-1}{2}\frac{p+1}{4}}=-(-1)^{\frac{p+1}{4}}=(-1)^{\frac{p-3}{4}}.
\end{align*}

When $K/\Q$ is inert at $p$, $\gamma=\epsilon_K^{\frac{p-3}{4}}$, so the global root number of $V$ is $-(-1)^{\frac{p-3}{4}}(-1)^{\frac{p-3}{4}}=-1$. 

When $K/\Q$ is ramified at $p$, $\epsilon_K=\lambda\eta$, where the conductor of $\eta$ is coprime to $p$, so $\gamma=\eta^{\frac{p-3}{4}}$. Now, $\eta(-1)=\lambda(-1)(\lambda\eta)(-1)=1$, and the global root number of $V$ is $-(-1)^{\frac{p-3}{4}}=(-1)^{\frac{p+1}{4}}$. 
}

The results are gathered in Table \ref{root-numbers-cartan-Hf}:

\begin{table}[htb]
\centering
\begin{tabular}{|c|c|c|c|c|}
\hline
& $(f,\mathbf{1}) \in \mathscr{S}$ & $(f,\chi) \in \mathscr{P}'$ & $(f,\mathbf{1}) \in \mathscr{C}'_1$ & $(f,\mathbf{1}) \in \mathscr{C}_M$\\
\hline
$p \nmid M$ & \multicolumn{4}{c|}{$-1$}\\
\hline
$p \mid M$, $p \equiv 1\pmod{4}$ & \multicolumn{2}{c|}{$-1$} & $-\lambda_p(f \otimes \overline{\psi})$ & None\\
\hline
$p \mid M$, $p \equiv 3 \pmod{4}$ & $a_p(f)\eta(-p)$ & $-1$ & $-\lambda_p(f \otimes \overline{\psi})$ & $(-1)^{\frac{p+1}{4}}$\\
\hline
\end{tabular}
\caption{Functional equation signs for factors of $J_{\rho,\xi}(p)$ of the form $V_f \otimes \alpha$, where $\alpha=\overline{\psi}(\omega_p)\epsilon_K^r \in H_f(\rho)$ }
\label{root-numbers-cartan-Hf}
\end{table}

\lem[ramified-K-If-unram]{Let $(f,\mathbf{1}) \in \mathscr{S} \cup \mathscr{C}'_1 \cup\mathscr{C}_M$ and $\psi \in \mathcal{I}_f(\rho)$. Assume that $K/\Q$ is ramified at $p$.  
\begin{itemize}[noitemsep,label=$-$]
\item If $p \equiv 3\pmod{4}$, then $\psi(\rho_K)$ vanishes on a decomposition subgroup at $p$. 
\item If $p \equiv 1 \pmod{4}$, then $\psi(\rho_K)$ is an unramified character of order at most $2$, mapping an arithmetic Frobenius to $\psi(\sqrt{a})$ for any $a \in \F_p^{\times}\backslash \F_p^{\times 2}$. 
\end{itemize}}

\demo{We fix a decomposition group $G_{\Q_p} \subset G_{\Q}$, let $I_p \leq G_{\Q_p}$ be the inertia group. If $p \equiv 3\pmod{4}$, by Proposition \ref{is-usually-cartan}, $\rho(G_K \cap G_{\Q_p})$ is scalar, so is contained in the kernel of $\psi$. Thus we may assume that $p \equiv 1\pmod{4}$.

Therefore, $\rho$ is $\frac{p+1}{2}$-Cartan by Proposition \ref{is-usually-cartan}. By Lemma \ref{expand-cartan}, $\rho(G_{\Q_p})/\rho(I_p)$ has cardinality two, and so does $\rho(I_p)/\rho(I_p \cap G_K)$: since $\rho(I_p \cap G_K)=\F_p^{\times}I_2$, $|\rho(G_{\Q_p})|=4(p-1)$. By Lemma \ref{expand-cartan}, $\rho(G_{\Q_p})$ is not commutative, so $\rho(G_{\Q_p}) \not\subset C$ and $\rho(G_{\Q_p} \cap G_K)$ is an index two subgroup of $\rho(G_{\Q_p})$. Hence $\rho(G_{\Q_p} \cap G_K) \subset C$ has cardinality $2(p-1)$ and contains $\F_p^{\times}I_2=\rho(I_p \cap G_K)$. The conclusion follows.
}

\prop[signs-nsplit-If-not-irred]{Let $(f,\chi) \in \mathscr{S} \cup \mathscr{P}' \cup \mathscr{C}'_1 \cup\mathscr{C}_M$. Let $\psi \in \mathcal{I}_f(\rho)$ be such that $\psi(\rho_K)^{p-1}=1$, and let $\alpha, \beta: G_{\Q} \rar F^{\times}$ be the two characters extending $\psi(\rho_K)$, and $\alpha_0,\beta_0$ be the associated Dirichlet characters. The global root numbers of the compatible systems $V_f \otimes \alpha, V_f \otimes \beta$ are given by Table \ref{table-root-numbers-cartan-If-nonirred}. The pair of the two root numbers is given by Table \ref{table-root-numbers-cartan-If-nonirred-pair}.}

\begin{table}[htb]
\centering
\begin{tabular}{|c|c|c|c|c|}
\hline 
& $(f,\mathbf{1}) \in \mathscr{S}$ & $(f,\chi) \in \mathscr{P}'$ & $(f,\mathbf{1}) \in \mathscr{C}'_1$ & $(f,\mathbf{1}) \in \mathscr{C}_M$\\
\hline
$p \nmid G$ & $\gamma(-p)a_p(f)$ & None & $-\gamma(-1)\lambda_p(f)$ & $\gamma(-1)(-1)^{\frac{p+1}{4}}$\\
\hline
$p \mid G$ & \multicolumn{2}{c|}{$-\gamma(-1)$} & $\gamma(-1)\lambda_p(f)$ & $\gamma(-1)(-1)^{\frac{p+1}{4}}$\\
\hline
\end{tabular}
\caption{Functional equation signs for factors of $J_{\rho,\xi}(p)$ of the form $V_f \otimes \gamma$, where $\gamma$ is one of the characters of $G_{\Q}$ extending $\psi(\rho_K)$ for $\psi \in \mathcal{I}_f$ such that $\psi(\rho_K)^{p-1}=1$, and $\gamma$ has conductor $G$. }
\label{table-root-numbers-cartan-If-nonirred}
\end{table}

\begin{table}[htb]
\centering
\small{\begin{tabular}{|c|c|c|c|c|}
\hline 
& $(f,\mathbf{1}) \in \mathscr{S}$ & $(f,\chi)\in \mathscr{P}'$ & $(f,\mathbf{1}) \in \mathscr{C}'_1$ & $(f,\mathbf{1}) \in \mathscr{C}_M$\\
\hline
$p \nmid M, \psi(\rho_K)_{|I_p}=1$ & $\alpha_0(-p)a_p(f)$ & None & \multicolumn{2}{c|}{$\{\pm 1\}$}\\
\hline
$p \nmid M, \psi(\rho_K)_{|I_p} \neq 1$ & \multicolumn{4}{c|}{$\{\pm 1\}$}\\
\hline
$p \mid M$ & $\alpha_0(-1)\{1,\alpha_0(p)a_p(f)\}$ & $\{\pm 1\}$ & $\{-\alpha_0(-1)\lambda_p(f)\}$ & $\{\alpha_0(-1)(-1)^{\frac{p+1}{4}}\}$\\
\hline
\end{tabular}}
\caption{Given a factor of the Tate module of $J_{\rho,\xi}(p)$ of the form $V_f \otimes \mrm{Ind}_K^{\Q}{\psi(\rho_K)}$ with $\psi(\rho_K)^{p-1}=1$, the pair of signs of the functional equations of the two compatible systems $V_f \otimes \alpha, V_f \otimes \beta$, where $\alpha,\beta$ are the characters of $G_{\Q}$ extending $\psi(\rho_K)$. In this table, $\alpha_0$ denotes the primitive Dirichlet character attached to any character among $\alpha,\beta$ which is unramified at $p$.}  
\label{table-root-numbers-cartan-If-nonirred-pair}
\end{table}

\demo{By Proposition \ref{tate-module-xrho-nsplit-cartan}, one $\alpha\beta=\varepsilon_K\overline{\chi}(\omega_p)$ and $\alpha\beta^{-1}=\epsilon_K$, which implies $\alpha_0^2=\beta_0^2=\overline{\chi}$. Since $\psi(\rho_K)$ takes its values in a subgroup of $\mu_{p^2-1}$, it is tamely ramified above $p$. Since $K/\Q$ is tamely ramified at $p$, $\alpha,\beta$ are also tamely ramified at $p$, and the $p$-adic valuations of the conductors of $\alpha_0,\beta_0$ are at most one.

Let $\gamma \in \{\alpha_0,\beta_0\}$, we can write $\gamma=\delta\eta$, where $\eta$ is a quadratic character with conductor prime to $p$ and $\delta \in \mathcal{D}$ is such that $\delta^2\chi=\mathbf{1}$. 

Suppose that $\delta\neq \mathbf{1}$ and $(f,\chi) \in \mathscr{S} \cup \mathscr{P}'$. Then, by Proposition \ref{local-constant-bigtwist}, $f \otimes \delta \in \mathcal{S}_2(\Gamma_0(p^2))$ is a newform with functional equation sign $-\delta(-1)$ by Proposition \ref{functional-eqn}. By Corollary \ref{global-constant-coprime-twist}, the sign of the functional equation for $f \otimes \gamma = (f \otimes \delta) \otimes \eta$ is $-\delta(-1)\eta(-p^2)=-\gamma(-1)$. 

When $(f,\chi) \in \mathscr{S} \cup \mathscr{P}'$, $\delta=\mathbf{1}$ if and only if the conductor of $\gamma$ is coprime to $p$. Since $\delta^2\chi=\mathbf{1}$, this is only possible when $(f,\chi) \in \mathscr{S}$ and one of $\alpha,\beta$ is unramified at $p$. When $K/\Q$ is inert at $p$, $\alpha\beta^{-1}=\epsilon_K$, so $\alpha$ is unramified at $p$ if and only if $\beta$ is, if and only if $\mrm{Ind}_K^{\Q}{\psi(\rho_K)}$ is unramified at $p$. When $K/\Q$ is ramified at $p$ and $(f,\chi) \in \mathscr{S}$, since $\alpha^2=\beta^2=\mathbf{1}$ and $\alpha\beta=\epsilon_K$, exactly one among the characters $\alpha,\beta$ is unramified at $p$.

So assume that $(f,\chi) \in \mathscr{S}$ and $\delta=\mathbf{1}$. Then, by Proposition \ref{functional-eqn} and Corollary \ref{global-constant-coprime-twist}, the sign of the functional equation for $f \otimes \gamma$ is $-\eta(-p)\lambda_p(f)=\eta(-p)a_p(f)$ by Proposition \ref{bad-L-factor}. 

Suppose now that $(f,\chi)=(f,\mathbf{1}) \in \mathscr{C}'_1 \cup \mathscr{C}_M$. Then $(f \otimes \delta,\mathbf{1}) \in \mathscr{C}$, so the sign of the functional equation for $f \otimes \gamma=(f \otimes \delta) \otimes \eta$ is \[(-i)^2\eta(-p^2)\lambda_p(f\otimes \delta)=-\eta(-1)\lambda_p(f \otimes \delta)=-\gamma(-1)\delta(-1)\lambda_p(f \otimes \delta).\] By Corollary \ref{various-computations-cuspidal}, this is equal to $\gamma(-1)\lambda_p(f)$ if $\delta$ is nontrivial, and $-\gamma(-1)\lambda_p(f)$ if $\delta$ is trivial. As we saw in the proof of Proposition \ref{signs-nsplit-Hf}, when $f$ is CM, $\lambda_p(f)=(-1)^{\frac{p-3}{4}}$. 

This explains how Table \ref{table-root-numbers-cartan-If-nonirred} is filled. Now, we discuss the contents of Table \ref{table-root-numbers-cartan-If-nonirred-pair}. Since $\alpha_0^2=\beta_0^2=\overline{\chi}$, if $\chi$ is not trivial, then it is necessary for $\alpha_0,\beta_0$ to both be ramified at $p$. When $K/\Q$ is unramified at $p$, then the conductor of the quotient character $\alpha_0\beta_0^{-1}$ (which is the Dirichlet character attached to the quadratic extension $K/\Q$) is coprime to $p$, so $\alpha$ is ramified at $p$ if and only if $\beta$ is ramified at $p$. In this case, one has $\alpha_0(-1)\beta_0(-1)^{-1}=\epsilon_K(c)=-1$, where $c$ is the complex conjugation, which accounts for the first two rows of the table except when $(f,\chi) \in \mathscr{S}$ and $\psi(\rho_K)$ is unramified at $p$. In this case, note that $\frac{\alpha_0}{\beta_0}(-p) = \epsilon_K(c\Fr_p)=(-1)\cdot (-1)=1$. 

When $p \mid M$, by Proposition \ref{is-usually-cartan}, $\rho_K(I_p \cap G_K)$ is scalar. When $(f,\chi) \in \mathscr{P}'$, $\alpha$ and $\beta$ are both ramified at $p$, thus the two signs are $-\alpha_0(-1),-\beta_0(-1)$, whose quotient is $\epsilon_K(c)=-1$. 

Otherwise, $\chi=\mathbf{1}$, so $\psi(\rho_K(I_p \cap G_K))$ is trivial, so the restriction to an inertia subgroup $I_p$ at $p$ of $\alpha,\beta$ is quadratic. Since $\alpha=\beta\epsilon_K$, and the restriction to $I_p$ of $\epsilon_K$ is a nontrivial quadratic character, exactly one among the characters $\alpha,\beta$ is unramified at $p$, and we can assume that it is $\alpha$. The remaining entries in the table follow from noting that $-\beta_0(-1)=\alpha_0(-1)$, using the results previously written in Table \ref{table-root-numbers-cartan-If-nonirred}, and the fact that $\alpha_0(-p)=\psi(\sqrt{a})$ for any $a \in \F_p^{\times}\backslash \F_p^{\times 2}$. Moreover, $\alpha_0(p)=\psi(\rho_K(\Fr_K))$, so it is equal to $\psi(\sqrt{a})$ for any $a \in \F_p^{\times}\backslash \F_p^{\times 2}$ if $p \equiv 1 \pmod{4}$ and is equal to $1$ otherwise.  
}

\medskip

\prop[sign-nsplit-2dx2d]{Let $(f,\chi) \in \mathscr{S} \cup \mathscr{P}' \cup \mathscr{C}'_1 \cup\mathscr{C}_M$. Let $\psi \in \mathcal{I}_f(\rho)$ be such that $\psi(\rho_K)^{p-1}\neq 1$. Then the global root number of the compatible system $(V_{f,\mathfrak{l}} \otimes \mrm{Ind}_K^{\Q}{\psi(\rho_K)}$ is given by the following:
\begin{itemize}[noitemsep,label=$-$]
\item When $(f,\chi) \in \mathscr{S}$ and $K/\Q$ is unramified at $p$, the sign is $1$ if and only if $\psi(\rho_K)$ is unramified at $p$.
\item When $(f,\chi) \in \mathscr{S}$ and $K/\Q$ is ramified at $p$, the sign is $a_p(f)\psi(\sqrt{a})$ (where $a \in \F_p^{\times}\backslash \F_p^{\times 2}$) if $p \equiv 1\pmod{4}$, and $a_p(f)$ if $p \equiv 3\pmod{4}$. 
\item When $(f,\chi) \in \mathscr{P}'$, the sign is $-1$.
\item When $(f,\chi) \in \mathscr{C}'_1 \cup \mathscr{C}_M$ and $K/\Q$ is unramified at $p$, the sign is $1$ if and only if the restriction of $\psi(\rho_K)$ to an inertia subgroup at $p$ is the same as that of $\overline{\phi}(\mathfrak{a}_E)$ or $\overline{\phi}^p(\mathfrak{a}_E)$, where $E/\Q_p$ is the quadratic unramified extension. 
\item When $(f,\chi)\in \mathscr{C}'_1 \cup \mathscr{C}_M$ and $K/\Q$ is ramified at $p$, the sign is $1$. 
\end{itemize}
}

\demo{For easier notation, we write $\Psi := \mrm{Ind}_K^{\Q}{\psi(\rho_K)}$ and $\epsilon_K=\epsilon_{K,p}(\omega_p)\epsilon'_{K,p}(\omega)$, where $\omega: G_{\Q} \rar \hat{\Z}^{\times}$ is the cyclotomic character, $\epsilon_{K,p} \in \mathcal{D}$ and $\epsilon'_{K,p}$ has conductor coprime to $p$. We fix a decomposition subgroup at $p$, thereby inducing an injection $G_{\Q_p} \subset G_{\Q}$, and we denote by $I_p \triangleleft G_{\Q_p}$ the inertia subgroup and by $\Fr_p \in G_{\Q_p}$ an arithmetic Frobenius with $\mathfrak{a}_{\Q_p}(\Fr_p)=\frac{1}{p}$. 

We apply Proposition \ref{functional-equations-jrho-ncartan}. By Proposition \ref{constants-with-modular-forms-elem}\ref{cmfe-1}, the constant at infinity is $(-i)^{2 \cdot 2}=1$ in all cases. 
By Proposition \ref{constants-with-modular-forms-elem}\ref{cmfe-2} and Lemma \ref{everything-self-dual-again}, the constant at a prime $\ell \neq p$ is exactly $\varepsilon(\Delp{\ell}{\Psi \oplus \Psi^{\ast}},\frac{1}{2},\theta_{\ell})$. 

Therefore, the global root number is 
\begin{align*}
&\varepsilon(D_{f,p} \otimes \Delp{p}{\Psi^{\ast}},1,\theta_p)\prod_{\ell \neq p\infty}{\varepsilon(\Delp{\ell}{\Psi \oplus \Psi^{\ast}},\frac{1}{2},\theta_{\ell})}\\
&=\frac{\varepsilon(D_{f,p} \otimes \Delp{p}{\Psi^{\ast}},1,\theta_p)}{\varepsilon(\Delp{\infty}{\Psi \oplus \Psi^{\ast}},\frac{1}{2},\theta_{\infty})\varepsilon(\Delp{p}{\Psi \oplus \Psi^{\ast}},\frac{1}{2},\theta_p)}=\frac{\varepsilon(D_{f,p} \otimes \Delp{p}{\Psi^{\ast}},1,\theta_p)}{\det{\Psi(c)}\det{\Psi(\mathfrak{a}_{\Q_p}^{-1}(-1))}}.
\end{align*}

By Lemma \ref{everything-self-dual-again}, $\det{\Psi}=\epsilon_K\overline{\chi}(\omega_p)$, so $\det{\Psi}(c)=-1$, and $\det{\Psi}(\mathfrak{a}_{\Q_p}^{-1}(-1))=\epsilon_{K,p}(-1)$. Thus $\epsilon_{K,p}(-1) \in \{-1,1\}$, and it is equal to $-1$ if and only if $K/\Q$ is ramified at $p$ and $p \equiv 3\pmod{4}$. \\

\emph{When $(f,\chi) \in \mathscr{P}'$:} 

By Proposition \ref{constants-with-modular-forms-elem}\ref{cmfe-2} and Lemma \ref{everything-self-dual-again}, one has therefore
\[\varepsilon(D_{f,p} \otimes \Delp{p}{\Psi^{\ast}},1,\theta_p)=\varepsilon(\Delp{p}{\Psi \oplus \Psi^{\ast}},\frac{1}{2},\theta_p)=\det{\Psi}(\mathfrak{a}_{\Q_p}^{-1})=\epsilon_{K,p}(-1),\]
whence the conclusion. \\

\emph{When $(f,\chi) \in \mathscr{S}$:} 

By Proposition \ref{constants-with-modular-forms-elem}\ref{cmfe-3} and Lemma \ref{everything-self-dual-again}, one has 
\[\varepsilon(D_{f,p} \otimes \Delp{p}{\Psi^{\ast}},1,\theta_p)=\varepsilon(\Delp{p}{\Psi \oplus \Psi^{\ast}},\frac{1}{2},\theta_p)\det(-a_p(f)\Fr_p \mid \Psi^{I_p}),\] so the global root number is $-\det(-a_p(f)\Fr_p \mid \Psi^{I_p})$.  

When $K/\Q$ is unramified at $p$, $\det{\Psi}=\epsilon_K$ is trivial on $I_p$, so $\Psi^{I_p}$ is either the entire space of $\Psi$ or zero. In the latter case, the global root number is therefore $-1$. In the former case, the global root number is $-(-a_p(f))^2\det{\Psi}(\Fr_p)=-\epsilon_K(\Fr_p)=-(-1)=1$. 

When $K/\Q$ is ramified at $p$, $\Psi_{|G_{\Q_p}} \simeq \mrm{Ind}_{K\Q_p}^{\Q_p}{\psi(\rho_K)_{|G_{K\Q_p}}}$ and we can apply Lemma \ref{ramified-K-If-unram}. 

When $p \equiv 3\pmod{4}$, $\psi(\rho_K)_{|G_{K\Q_p}}$ is the trivial character, so $\Psi_{|G_{\Q_p}}$ is the direct sum of the trivial character and the (ramified) nontrivial character of $\mrm{Gal}(K\Q_p/\Q_p)$. Thus \[-\det(-a_p(f)\Fr_p\mid \Psi^{I_p})=-(-a_p(f))=a_p(f).\] 

When $p \equiv 1 \pmod{4}$, then, by Lemma \ref{ramified-K-If-unram}, $\psi(\rho_K)_{|G_{K\Q_p}}$ is the unramified character mapping an arithmetic Frobenius to $\psi(\sqrt{a})$ for any $a \in \F_p^{\times}\backslash \F_p^{\times 2}$. So $\Psi_{|G_{\Q_p}}$ is the direct sum of the two characters $\alpha,\beta$ of $G_{\Q_p}$ extending $\psi(\rho_K)_{|G_{K\Q_p}}$. Since $\alpha\beta^{-1}$ is the nontrivial character of $\mrm{Gal}(K\Q_p/\Q_p)$, at most one of these characters is unramified, and the unramified character of $G_{\Q_p}$ mapping an arithmetic Frobenius to $\psi(\sqrt{a})$ is one of $\alpha,\beta$. Thus \[-\det(-a_p(f)\Fr_p\mid \Psi^{I_p})=-(-a_p(f)\psi(\sqrt{a}))=a_p(f)\psi(\sqrt{a}).\]\\

\emph{When $(f,\chi) \in \mathscr{C}$:}

Let $E/\Q_p$ be the unique quadratic unramified extension and $\phi: E^{\times}/(1+p\OO_E) \rar \C^{\times}$ be a character such that $\phi_{|(\OO_E/p)^{\times}}$ is attached to $C_f$ and $\phi(p)=-p$. 
The case where $K/\Q$ is unramified at $p$ is a consequence of Proposition \ref{local-constants-p-induced}\ref{lcpi-3} (because $\epsilon_{K,p}(-1)=1$), so we may assume that $K/\Q$ is ramified at $p$. Let $\mathfrak{p}$ be the unique prime ideal of $K$ above $p$, and $K_p$ be the $\mathfrak{p}$-adic completion of $K$ and $\tau: \mrm{Gal}(K_p/\Q_p) \rar \{\pm 1\}$ the unique nontrivial character. 

We denote by $\lambda \in \mathcal{D}$ the unique nontrivial quadratic Dirichlet character of conductor $p$. 

The subgroup of norms $N_{K_p/\Q_p}K_p^{\times} \subset \Q_p^{\times}$ has valuation group $\Z$ and index $2$, so its intersection with $\Z_p^{\times}$ is $\Z_p^{\times 2}$. Hence the subgroup of norms $N_{EK_p/\Q_p}(EK_p)^{\times} \subset \Q_p^{\times}$ has index $4$, valuation group $2\Z$, contains $p^2\in N_{E/\Q_p}E^{\times} \cap N_{K_p/\Q_p}K_p^{\times}$, and contains $\Z_p^{\times 2}$: hence $N_{EK_p/\Q_p}(EK_p)^{\times}=p^{2\Z}\Z_p^{\times 2}$ is the subgroup of norms of $E(\sqrt{-p})$ and $EK_p=E(\sqrt{-p})$. Thus $\tau_{|G_E}=\left(\frac{\omega_p}{p}\right)$.

When $p \equiv 3\pmod{4}$, $\psi(\rho_K)_{|G_{K_p}}$ is trivial by Lemma \ref{ramified-K-If-unram}, so $\Psi_{|G_{\Q_p}}=\mathbf{1} \oplus \tau$. Since $D_{f,p}$ is induced from $E$, one has 
\[D_{f,p} \otimes \Delp{p}{\Psi^{\ast}} = D_{f,p} \otimes \left[\mathbf{1}\oplus \left(\frac{\omega_p}{p}\right)\right] = D_{f,p} \oplus D_{f \otimes \lambda,p},\]
 hence 
 \[\varepsilon(D_{f,p} \otimes \Delp{p}{\Psi^{\ast}},1,\theta_p)=\epsilon_p(f)\epsilon_p(f \otimes \lambda)=\lambda_p(f)\lambda_p(f \otimes \lambda)=-\lambda(-1)=(-1)^{\frac{p+1}{2}}=-\epsilon_{K,p}(-1)\] by Corollary \ref{various-computations-cuspidal}, whence the conclusion.

When $p \equiv 1 \pmod{4}$, $\psi(\rho_K)_{|G_{K_p}}$ is the unramified character $\Fr_K \longmapsto \psi(\sqrt{a})$ (for any $a \in \F_p^{\times} \backslash \F_p^{\times 2}$) with order dividing $2$ by Lemma \ref{ramified-K-If-unram}. Then $\Psi_{|G_{\Q_p}}$ is the direct sum of the unramified character $\eta: \Fr_p \in G_{\Q_p}\longmapsto \psi(\sqrt{a}) \in \{\pm 1\}$ and $\eta\tau$. Since $D_{f,p}$ is induced from $E$ and $\Psi_{|G_E}$ is isomorphic to the restriction to $G_E$ of $\mathbf{1} \oplus \left(\frac{\omega_p}{p}\right)$, one has as above

\begin{align*}
\varepsilon(D_{f,p} \otimes \Delp{p}{\Psi^{\ast}},1,\theta_p)&=\varepsilon(D_{f,p} \otimes (\mathbf{1} \oplus \left(\frac{\omega_p}{p}\right)^{-1}),1,\theta_p)=\varepsilon(D_{f,p} \oplus D_{f\otimes \lambda,p},1,\theta_p)\\
&= \epsilon_p(f)\epsilon_p(f \otimes \lambda)=\lambda_p(f)\lambda_p(f \otimes \lambda)=-\lambda(-1)=(-1)^{\frac{p+1}{2}}=-\epsilon_{K,p}(-1)
\end{align*} by Corollary \ref{various-computations-cuspidal}, whence the conclusion.
}

\cor[when-signs-are-1]{Assume that $K/\Q$ is unramified at $p$. Then, among the irreducible factors of $\Tate{\ell}{J_{\rho,\xi}(p)} \otimes_{\Z_{\ell}} F_{\mathfrak{l}}$ that are isomorphic to their Tate-twisted duals described in Proposition \ref{tate-module-xrho-nsplit-cartan}, the following ones have global root number $+1$:
\begin{itemize}[label=\tiny$\bullet$,noitemsep]
\item When $[C:\rho(G_K)] > 1$ and $(f,\chi) \in \mathscr{P}'$, there exists $\psi \in \mathcal{I}_f$ such that $\psi(\rho_K)^{p-1}=\{1\}$. Let $\alpha: G_{\Q} \rar F^{\times}$ be the (unique) odd character which extends $\psi(\rho_K)$, then $V_{f,\mathfrak{l}} \otimes \alpha$ works. 
\item When $\rho(G_K)=C$, $[C: \rho(I_p)] > 1$ and $(f,\mathbf{1}) \in\mathscr{S}$, there exists $\psi \in \mathcal{I}_f$ with $\psi(\rho_K)^{p-1}\neq \{1\}$ and $\psi(I_p)=\{1\}$. Then $V_{f,\mathfrak{l}} \otimes \mrm{Ind}_K^{\Q}{\psi(\rho_K)}$ works. 
\item When $\rho(I_p)=C$ and $b,b'$ are the two reciprocity indices of $\rho_{|G_{\Q_p}}$, for any $(f,\chi) \in \mathscr{C}'_1 \cup \mathscr{C}_M$ such that $b,b'$ do not divide $\frac{p+1}{\omega}$, where $\omega$ is the order of the character $\phi: C/\F_p^{\times} \rar \C^{\times}$ attached to $f$, then there exists $\psi \in \mathcal{I}_f$ such that $V_{f,\mathfrak{l}} \otimes \mrm{Ind}_K^{\Q}{\psi(\rho_K)}$ works. 
\end{itemize}
}

\demo{In the first case, we know that $\rho(G_K)=a'C$ with $a' \mid a$; by Proposition \ref{trichotomy-rho}, $a \mid \frac{p+1}{2}$ is coprime to $p-1$, thus $a' \mid \frac{p+1}{2}$ is odd. In particular, $\F_p^{\times}I_2 \subset \rho(G_K)$ and $[C:\rho(G_K)] > 2$. Let $\gamma \in \mathcal{D}$ with square $\chi$ and $\gamma_0=\gamma(\det): C \rar F^{\times}$. Then $(\gamma_0)(tI_2)=\chi(t)$ for any $t \in \F_p^{\times}$. 

The characters $\beta: C \rar F^{\times}$ such that $\beta(tI_2)=\chi(t)$ for every $t \in \F_p^{\times}$ are exactly the $\gamma_0\psi$, where $\psi$ runs through the characters $C/\F_p^{\times} \rar F^{\times}$. For any $\psi: C/\F_p^{\times} \rar F^{\times}$, one has 
\begin{align*}
\beta := \gamma_0\psi \in \mathcal{I}_f &\Leftrightarrow \beta^{p-1} \neq \mathbf{1} \Leftrightarrow \psi^{p-1} \neq \mathbf{1} \Leftrightarrow \psi^2\neq\mathbf{1},\\
\beta(\rho_K)^{p-1} = \{1\} &\Leftrightarrow \psi(\rho_K)^{p-1}=\{1\} \Leftrightarrow \psi^2(\rho_K)=\{1\}. 
\end{align*}

The group $C/(\F_p^{\times}\rho(G_K)) = C/\rho(G_K)$ is cyclic with cardinality at least $3$, so there is a character $\psi: C \rar F^{\times}$ such that $\psi(\F_p^{\times}\rho(G_K))=\{1\}$ and $\psi^2 \neq \mathbf{1}$. Then $\gamma_0\psi \in \mathcal{I}_f$ is such that $(\gamma_0\psi)(\rho_K)^{p-1}=\{1\}$, so that there are two characters $\alpha,\beta: G_{\Q} \rar \C^{\times}$ extending $(\gamma_0\psi)(\rho_K)$. Their ratio is the odd character $\epsilon_K$, so exactly one among $\alpha$ and $\beta$ is odd. Since $\psi(\rho_K)$ takes its values in $F^{\times}$ and $G_{\Q}$ is generated by $G_K$ and the complex conjugation, $\alpha, \beta$ also take their values in $F^{\times}$. The conclusion follows from Proposition \ref{signs-nsplit-If-not-irred}.

In the second case, $\rho(I_p)=aC$ for some odd $a \mid \frac{p+1}{2}$ which is greater than one. So $[C:\rho(I_p)] > 2$ and $\rho(I_p) \supset \F_p^{\times}$. Therefore, there exists a character $\psi: C/\rho(I_p) \rar F^{\times}$ with order at least $3$. Then $\psi(\rho_K)^{p-1} \neq \{1\}$, so $V_{f,\mathfrak{l}} \otimes \mrm{Ind}_K^{\Q}{\psi(\rho_K)}$ is irreducible and isomorphic to its Tate-twisted dual. Since $\psi(\rho(I_p))=\{1\}$, the conclusion follows from Proposition \ref{sign-nsplit-2dx2d}. 

In the third case, let $E/\Q_p$ be the quadratic unramified extension and fix an isomorphism $\iota: C \rar (\OO_E/p)^{\times}$ preserving the trace and norm. Let $\sigma: E/(1+p\OO_E) \rar C$ be such that $\rho_{|W_E} = \sigma(\mathfrak{a}_E)$. By Lemma \ref{expand-cartan}, $\sigma(x)=x^{-1}I_2$ for every $x \in \F_p^{\times}$. After exchanging $b,b'$ if needed, the following equalities hold:
\[\{x \in C\mid x\sigma(\iota(x))=I_2\}=bC,\quad \{x \in C\mid x\sigma(\iota^p(x))=I_2\}=b'C.\]
Let $\phi: C \rar F^{\times}$ be one of the characters attached to $C_f$ and $\omega$ be its order. It is enough to prove the following implication, by Proposition \ref{sign-nsplit-2dx2d}: if $b,b'$ do not divide $\frac{p+1}{\omega}$, then there exists $\psi: C/\F_p^{\times} \rar F^{\times}$ such that $\psi\circ\sigma\circ\iota =\overline{\phi}$ and $\psi \notin \{\phi,\phi^p\}$. 

Since $\rho(I_p)=C$, $\sigma\circ\iota$ is an automorphism of $C$, so we need to show that $\overline{\phi}\notin \{\phi\circ\sigma\circ\iota ,\phi\circ\sigma\circ\iota^p\}$. Let $(b_1,b_p)$ denote $(b,b')$, for any $j \in \{1,p\}$, 
\begin{align*}
\overline{\phi}=\phi\circ\sigma\circ\iota^j &\Leftrightarrow \phi\left(\{x\sigma(\iota^j(x))\mid x \in C\}\right)=\{1\} \Leftrightarrow \phi\left(\frac{p^2-1}{b_j}C\right)=\{1\}\\
&\Leftrightarrow \omega\mid \frac{p+1}{b_j}(p-1) \Leftrightarrow \omega\mid\frac{2(p+1)}{b_j},
\end{align*}
whence the conclusion.
}

\section{Application to elliptic curves}

The following result (which is independent from the setting defined at the beginning of the Chapter) comes essentially from \cite[Proposition 1.13]{Zywina}. In view of Corollary \ref{when-signs-are-1}, we have formulated it as a purely local result, because the image of the inertia group at $p$ also affects the sign computation. 

\lem[cartan-is-potentially-good]{Let $p \geq 7$ be a prime. Let $E/\Q_p$ be an elliptic curve and $\rho: G_{\Q_p} \rar \GL{\F_p}$ be a representation attached to the choice of a basis for $E[p]$ (in the sense of Proposition \ref{group-is-twist-polarized}, that is, $\det{\rho}=\omega_p^{-1}$). Assume that $\rho$ is $a$-Cartan. Then there exists a tamely ramified extension $L/\Q_p$ satisfying the following properties:
\begin{itemize}[noitemsep,label=$-$]
\item The ramification index $e_{L/\Q_p}$ properly divides $12$,
\item $E_L/L$ has good supersingular reduction,
\item One of the reciprocity indices of $\rho$ (see Lemma \ref{expand-cartan}) divides $e_{L/\Q_p}$,
\item $a \in \{1,3\}$,
\item If $a=3$, then $p \equiv 2,5\pmod{9}$ and the action of the inertia group $I_p$ on $E[2]$ is cyclic with order $3$. 
\end{itemize}
}

\demo{If $j(E) \notin \Z_p$, by \cite[Proposition 13]{Serre-image-ouverte}, a subgroup of index two of $\rho(I_p)$ is contained in a nonsplit Cartan subgroup and in a subgroup conjugate to $\begin{pmatrix}\ast & \ast\\0 & 1\end{pmatrix}$. Therefore $\{I_2\}$ has index at most two in $\rho(I_p)$, so $2 \geq |\rho(I_p)|\geq 2(p-1)$, whence a contradiction. So $j(E) \in \Z_p$ and $E$ has potentially good reduction at $p$ by \cite[Proposition VII.5.5]{AEC1}. By \cite[Theorem 2, Corollary 2]{GoodRed}, the image of $I_p$ in $\mrm{Aut}(E[m])$ is a quotient of the tame inertia group and does not depend on the choice of an integer $m \geq 3$ coprime to $p$. In particular, said image is cyclic with cardinality $e$ such that $\GL{\F_3}$ and $\GL{\Z/4\Z}$ contain cyclic subgroups of cardinality $e$. This implies that $e \in \{1,2,3,4,6\}$. 

Let $L/\Q_p$ be a finite extension such that $I_L$ is exactly the subgroup of $I_p$ of index $e$, which is the kernel of its action on $E[m]$ for every $m \geq 3$ coprime to $p$: so $L/\Q_p$ has ramification index $e$ and is tamely ramified. By the N\'eron-Ogg-Shafarevich criterion \cite[Theorem 1]{GoodRed}, $E_L/L$ has good reduction. By \cite[Proposition 11]{Serre-image-ouverte}\footnote{Note that our convention for $\rho$ is different than that of Serre.}, if $E_L/L$ has good ordinary reduction, then $\rho(I_L)$ is contained in a nonsplit Cartan subgroup and is conjugate to a subgroup of the form $\begin{pmatrix}1 & \ast\\0 & \ast\end{pmatrix}$, hence is trivial. Then $2(p-1) \leq |\rho(I_p)| \leq [I_p:I_L]=e \leq 6$, which is a contradiction. 

Hence $E_L/L$ has good supersingular reduction. We can then apply the discussion of \cite[\S 1.10, \S 1.11 (2)]{Serre-image-ouverte}. Since $\rho$ is tamely ramified, the denominators of the valuations of the elements of $E[p](\overline{\Q_p})$ are coprime to $p$; since $e \leq 6 < p$, the Newton polygon of the formal group attached to $E_L$ is the line $(1,e) \rar (p^2,0)$, so, by Proposition 10 of \emph{loc.cit.} (and taking into account the change of conventions with respect to Serre's) $\rho_{|I_L}$ is given by the $(-e)$-th power of a fundamental character of level two for $L$. 

Let $D \leq \GL{\F_p}$ be a nonsplit Cartan subgroup with $\rho(I_p)=aD=\rho(G_F)$, where $F/\Q_p$ is the quadratic unramified extension, by Lemma \ref{expand-cartan}. Let us fix an identification $\iota$ of $D$ with $(\OO_F/p)^{\times}$ preserving the trace and the norm. What we just wrote is that, after if necessary replacing $\iota$ with $\iota^p$, $\rho_{|I_p}(z)=\iota(\mathfrak{a}_F(z))^{-1}$ for any $z \in I_L$, so one of the reciprocity indices of $\rho$ divides $e$, so $a \in \{1,2,3,4,6\}$. Since $a$ is odd, $a \in \{1,3\}$. 

We proved in Lemma \ref{expand-cartan} that $a$ was coprime to $\frac{p^2-1}{b}$ for any reciprocity index $b$ of $\rho$, hence $a$ is coprime to $\frac{p^2-1}{e}$. Thus, if $a=3$, then $e \in \{3,6\}$, $a \mid \frac{p+1}{2}$ and $\frac{p^2-1}{e}$ is not divisible by $3$. Therefore, $v_3(p+1)=1$ so $p \equiv 2,5\pmod{9}$. 

Since $I_p/I_L$ is cyclic, the action of $I_p$ on $E[2]$ is a cyclic quotient of $\Z/e\Z$ by a $2$-group (since the group of automorphisms of $E[4]$ fixing $E[2]$ is a $2$-group) and it embeds in $\GL{\F_2}$, so the action of $I_p$ on $E[2]$ is cyclic with order $3$.
}

\medskip

Next, we show that in most situations, the extension $K/\Q$ is inert at $p$.   

\lem[almost-never-ramified]{Assume that $\rho_{|G_{\Q_p}}$ comes from an elliptic curve $E/\Q_p$ and that $K$ is ramified. Then $p=7$ and there exists a finite extension $L/\Q_p$ such that $E_L/L$ has good ordinary reduction and $e_{L/\Q_p} \in \{3,6\}$.}

\demo{If $p \equiv 1\pmod{4}$, by Lemma \ref{is-usually-cartan}, $\rho$ is $\frac{p+1}{2}$-Cartan. By Proposition \ref{cartan-is-potentially-good}, $\frac{p+1}{2} \leq a$, so $p \leq 7$. Since $p \equiv 1\pmod{4}$, one has $p=5$, which we excluded. So we assume that $p \equiv 3\pmod{4}$: thus $\rho$ is split and $\rho(G_{\Q_p} \cap G_K)$ is scalar. 

If $j(E) \notin \Z_p$, then, by \cite[Proposition 13]{Serre-image-ouverte}, $\rho(I_p)$ is conjugate to a subgroup of $\pm\begin{pmatrix}\ast & \ast\\0 & 1\end{pmatrix}$. Since $\rho(I_p \cap G_K)$ is scalar, we see that $\rho(I_p \cap G_K) \subset \{\pm I_2\}$, hence $4 \geq |\rho(I_p)|=p-1$ by Proposition \ref{trichotomy-rho}, a contradiction. As in the proof of Lemma \ref{cartan-is-potentially-good}, there is a finite extension $L/\Q_p$ with ramification index $e \in \{1,2,3,4,6\}$ such that $E_L/L$ has good reduction. 

If this reduction is supersingular, by the same argument as in the proof of Lemma \ref{cartan-is-potentially-good}, $\rho_{|I_L}$ agrees with $\left((\mathfrak{a}_F)_{|I_L}\right)^{-1}: I_L \rar (\OO_F/p)^{\times}$, where $F/\Q_p$ is the quadratic unramified extension. Thus $|\rho(I_L)|=\frac{p^2-1}{e}$, but $|\rho(I_L)| \leq |\rho(I_p)|=p-1$, so $e \geq p+1 \geq 8$, which is impossible. 

Hence $E_L/L$ has good ordinary reduction. By \cite[Proposition 11]{Serre-image-ouverte}, $\rho_{|I_L}$ is given by $\begin{pmatrix}1 & \ast\\0 & \omega_p^{-1}\end{pmatrix}$. Therefore $\omega_p^{-1}(I_L \cap G_K)=\{1\}$, hence $p-1 \mid 2e$. Therefore $p=7$ and $e \in \{3,6\}$. }

\medskip

\lem[what-about-cm]{Assume that $\rho$ comes from a elliptic curve $E/\Q$ with complex multiplication by $K$ and that $p \neq 7$. Then $K/\Q$ is inert at $p$, and $E$ acquires good reduction at $p$ after a finite extension of $\Q_p$ with ramification index dividing $4$ if $K=\Q(i)$, $6$ if $K=\Q(j)$, $2$ otherwise. In particular, unless $K=\Q(j)$, $\rho_{G_{\Q_p}}$ is $1$-Cartan.}

\demo{Let $I_p$ be an inertia group at $p$. Since $p \neq 7$, Lemma \ref{almost-never-ramified} implies that $K/\Q$ is inert at $p$, so by Lemma \ref{cartan-is-potentially-good} there is a finite extension $L/\Q_p$ with ramification index $e \in \{1,2,3,4,6\}$ such that $e$ is the cardinality of the image of $I_p$ on $E[m]$ for any $m \geq 3$ coprime to $p$, by \cite[Theorem 2]{GoodRed}. Our goal is to show that if $K \neq \Q(i)$ (resp. $K \neq \Q(j)$), then $e \mid 6$ (resp. $e \mid 4$). 

If $K \neq \Q(i)$ (resp. $K \neq \Q(j)$), there is an odd prime $\ell \equiv 3\pmod{4}$ (resp. $\ell \equiv 2\pmod{3}$) which splits in $K$. The image $J$ of $I_p$ in $\mrm{Aut}_{\Z_{\ell}}(\Tate{\ell}{E})$ is a cyclic group with cardinality $e$ and commutes to the image in $\mathcal{M}_2(\Z_{\ell})$ of the diagonalizable subgroup $\OO_K \otimes \Z_{\ell}$, so $J=g\begin{pmatrix}\alpha & 0\\0 & \beta\end{pmatrix}g^{-1}$ for some characters $\alpha,\beta: I_p \rar \Z_{\ell}^{\times}$ and $g \in \GL{\Z_{\ell}}$. Thus $J$ is cyclic of order $e$ and embeds in the torsion subgroup of $(\Z_{\ell}^{\times})^2$, which is isomorphic to $(\F_{\ell}^{\times})^{\oplus 2}$. By comparing exponents, we see that $e \mid \ell-1$, so $e \neq 4$ thus $e \mid 6$ (resp. $e$ prime to $3$ thus $e \mid 4$). }

\medskip

\prop[counting-If-signs]{Assume that $p>7$ and $\rho_{|G_{\Q_p}}$ comes from an elliptic curve $E/\Q_p$. Then $K/\Q$ is inert at $p$, $\rho_{|G_{\Q_p}}$ is $a$-Cartan for some $a \in \{1,3\}$, and let $b \in \{1,2,3,4,6\}$ be its smallest reciprocity index (which is divisible by $a$) and let $b'$ be the other reciprocity index. Let $(f,\mathbf{1}) \in \mathscr{C}'_1$ (resp. $(f,\mathbf{1}) \in \mathscr{C}_M$), let $\phi: C/\F_p^{\times}I_2 \rar \C^{\times}$ be a character attached to $C_f$, and let $\omega$ be the order of $\phi$. Let $I$ (resp. $I/2$) be the number of $\psi \in \mathcal{I}_f/\sim$ such that the global root number of the compatible system $V_f \otimes \mrm{Ind}_K^{\Q}{\psi(\rho_K)}$ is equal to $+1$. Then one has:
\begin{itemize}[label=\tiny$\bullet$,noitemsep]
\item $I \in \{0,\min(a,2)\}$,
\item $I > 0$ only if $b \geq 3$ and $a\omega \mid p+1$,
\item When $b=4$, $I=1$ if and only if $v_2(\omega)=v_2(p+1)$ and $(\omega,v_2(p+1)) \neq (4,2)$.
\item When $b=3$, $I=1$ if and only if $v_3(\omega)=v_3(p+1)$ and $(\omega,v_3(p+1)) \notin \{(3,1),(6,1)\}$. 
\end{itemize} 
}

\demo{$K/\Q_p$ is inert at $p$ by Lemma \ref{almost-never-ramified}, so all the claims before the study case by case follow from Lemma \ref{cartan-is-potentially-good} and Lemma \ref{expand-cartan}. Let $F/\Q_p$ be the quadratic unramified extension and $\sigma: F^{\times}/(1+p\OO_F) \rar C$ be such that $\rho_{|W_F}=\sigma(\mathfrak{a}_F)$. We fix an isomorphism $\iota: C \rar (\OO_F/p)^{\times}$ preserving the trace and the norm such that $\{x \in C\,x\sigma(\iota(x))=I_2\}=bC$. By definition, the image of $\sigma\circ\iota$ is $aC$ (thus it contains $\F_p^{\times}I_2=(p+1)C$), and, by Lemma \ref{expand-cartan}, $\sigma(x)=x^{-1}I_2$ for every $x \in \F_p^{\times}$, so $x \in C/\F_p^{\times}I_2 \longmapsto x\sigma(\iota(x)) \in C$ is well-defined. 

By Proposition \ref{sign-nsplit-2dx2d}, $I$ is the number of characters $\psi \in \mathcal{I}_f$ such that $\psi \circ \sigma \circ \iota = \overline{\phi}$, because $z \longmapsto z^p$ is an involution of $\mathcal{I}_f$ without any fixed point. Let thus $\Phi$ be the set of $\psi: C/\F_p^{\times}I_2 \rar \C^{\times}$ such that $\psi \circ \sigma \circ \iota = \overline{\phi}$. 

If $\psi\in \Phi$, then $\omega$ divides the order of $\psi$, so $\psi^2 \neq 1$. Thus $I$ is the difference of the cardinalities of $\Phi$ and of $\Phi \cap \{\phi,\phi^p\}$. 

Now, $\Phi$ is nonempty if, and only if, $\ker{\phi} \supset \ker{\sigma\circ\iota}$. By definition, $\ker{\phi}=\omega C$, while $[C:\ker{\sigma\circ\iota}]=|(\sigma\circ\iota)(C)|=\frac{p^2-1}{a}$, so $\Phi \neq \emptyset$ if and only if $\omega C \supset \frac{p^2-1}{a}C$, if and only if $a\omega \mid p^2-1$. 

Suppose that $a\omega \mid p^2-1$. One has $v_2(a\omega)=v_2(\omega)\leq v_2(p+1)$, $a \mid p+1$ is coprime to $p-1$ and $\omega \mid p+1$: write $\omega=2^z\omega'$ for some odd $\omega'$, then $a$, $\omega'$ are odd and divide $p+1$, so they are coprime to $p-1$: hence $a\omega' \mid (p+1)(p-1)$ is coprime to $p-1$, therefore $a\omega'\mid p+1$, hence $a\omega \mid p+1$. The converse is obvious, so $\Phi \neq\emptyset$ if, and only if, $a\omega \mid p+1$. 

When $a\omega \mid p+1$, $|\Phi|$ is the number of characters $\psi: C/\F_p^{\times}I_2 \rar \C^{\times}$ such that $\ker{\psi}$ contains $\sigma(\iota(C))=aC$. Therefore, $|\Phi|=a$. Since $\overline{\phi} \neq \overline{\phi}^p$, $\Phi$ does not contain both $\phi$ and $\phi^p$, so $I \in \{|\Phi|,|\Phi|-1\}$. 

Now, since $\omega \mid p+1$ and $b \mid p+1$,
\begin{align*}
\phi\in \Phi &\Leftrightarrow \phi\left(\{x \sigma(\iota(x))\mid x \in C\}\right)=\{1\} \Leftrightarrow \phi\left(\frac{p^2-1}{b}C\right)=\{1\} \\
&\Leftrightarrow \omega \mid \frac{p^2-1}{b} \Leftrightarrow \omega \mid \frac{2(p+1)}{b}.
\end{align*}

Let $b'$ be the other reciprocity index of $\rho_{|G_{\Q_p}}$, so that $a \mid b'$ and the least common multiple of $b$ and $b'$ is $p+1$. Then, similarly,
\begin{align*}
\phi^p\in \Phi &\Leftrightarrow \phi\left(\{x\sigma(\iota^p(x)) \mid x \in C\}\right)=\{1\} \Leftrightarrow \phi\left(\frac{p^2-1}{b'}C\right)=\{1\} \Leftrightarrow \omega \mid \frac{2(p+1)}{b'}.
\end{align*}

If $a=3$ and $I \neq 0$, then $|\Phi|=3$ so $a\omega \mid p+1$. By Lemma \ref{cartan-is-potentially-good}, this implies that $v_3(p+1)=1$ and $b \in \{3,6\}$, so $\omega$ is coprime to $3$ and $\omega$ divides $3\frac{2(p+1)}{b}=(p+1)\frac{6}{b}$, thus $\omega \mid\frac{2(p+1)}{b}$ and $\phi \in \Phi$. Hence $I=2$. 

So from now on we assume that $a=1$, so that $a\omega \mid p+1$, thus $|\Phi|=1$ and $I \in \{0,1\}$. When $b \leq 2$, $b\omega \mid p+1$, so $\phi\in \Phi$ and $I=0$.

There exists some integer $t$ such that $\sigma(\iota(x))=t\cdot x$ for all $x \in C$. Then $\frac{p^2-1}{b}$ (resp. $\frac{p^2-1}{b'}$) is the greatest common divisor of $p^2-1$ and $t+1$ (resp. $p^2-1$ and $pt+1$). Moreover, $a=1$, so $t$ is coprime to $p^2-1$. 
 
When $b=4$, note first that $p \equiv 3\pmod{4}$. Then $\phi \in \Phi$ if and only if $\omega \mid \frac{p+1}{2}$, which is equivalent to $v_2(\omega)<v_2(p+1)$. The least common multiple of $b,b'$ is $p+1$, so $\frac{p+1}{4} \mid b'$ and $b'=p+1$ if $v_2(p+1) > 2$. Thus, when $v_2(p+1) > 2$, $I=1$ if and only if $v_2(\omega) = v_2(p+1)$. When $v_2(p+1)=2$, the greatest common divisor of $p^2-1$ and $t+1$ is $\frac{p^2-1}{4}$, so $v_2(t+1)=1$, and $t \equiv 1\pmod{4}$. Thus $pt+1$ and $p^2-1$ are divisible by $4$, so their greatest common divisor $\frac{p^2-1}{b'}$ is also divisible by $4$, so $v_2(b') \leq 1$ and $b' \in \{\frac{p+1}{4},\frac{p+1}{2}\}$. Therefore, the greatest common divisor of $\frac{2(p+1)}{b'}$ and $p+1$ is $4$. 

Therefore, $I=1$ if and only if $v_2(\omega)=v_2(p+1)$ and $(\omega,v_2(p+1)) \neq (4,2)$.

When $b \in \{3,6\}$, note first that $p \equiv 2\pmod{3}$. Then $\phi \in \Phi$ if and only if $\omega \mid \frac{2(p+1)}{b}$, which is equivalent to $v_3(\omega) < v_3(p+1)$. Moreover, $b'$ is divisible by $\frac{p+1}{b}$, so, if $\phi^p \in \Phi$, one has $v_3(\omega)=v_3(p+1)$ and $\omega \mid \frac{2(p+1)}{b'}$, so $b'$ is coprime to $3$ and $v_3(p+1)=v_3(b)=1$. 

Assume now that $v_3(\omega)=v_3(p+1)=1$. Then $\frac{p^2-1}{b}$ (which is the greatest common divisor of $t+1$ and $p^2-1$) is coprime to $3$, so $t \not\equiv -1\pmod{3}$. Since $t$ is coprime to $p^2-1$, $t \equiv 1\pmod{3}$ and the greatest common divisor of $tp+1$ and $p^2-1$, which is $\frac{p^2-1}{b'}$, is divisible by $3$. Hence $b'$ is coprime to $3$ and $b' \in \{\frac{p+1}{3},\frac{p+1}{6}\}$. If $b'=\frac{p+1}{3}$, then $\frac{2(p+1)}{b'}=6$. If $b'=\frac{p+1}{6}$, since $p+1$ is the least common multiple of $b,b'$, $v_2(p+1)=v_2(b) \leq 1$, so $p \equiv 1\pmod{4}$, hence the greatest common divisor of $\frac{2(p+1)}{b'}=12$ and $p+1$ is $6$.  

Therefore, $I=1$ if and only if $v_3(\omega)=v_3(p+1)$ and $(\omega,v_3(p+1)) \notin \{(3,1),(6,1)\}$.

}

\cor[good-reduction-all-minusone]{Assume that $p > 7$ and that $\rho_{|G_{\Q_p}}$ comes from an elliptic curve $E/\Q_p$ with good reduction at $p$ up to quadratic twist. Then, in the decomposition of Proposition \ref{tate-module-xrho-nsplit-cartan}, all the given factors are irreducible and isomorphic to their Tate-twisted duals, and their global root numbers are all $-1$.}

\demo{By Lemma \ref{cartan-is-potentially-good}, $\rho_{|G_{\Q_p}}$ is $1$-Cartan and its smallest reciprocity index is at most $2$. By Lemma \ref{almost-never-ramified}, $K/\Q$ is inert at $p$, so, in the decomposition of Proposition \ref{tate-module-xrho-nsplit-cartan}, all the given factors are irreducible and isomorphic to their Tate-twisted duals. Indeed, the only possible exceptions are situations where $(f,\chi) \in \mathscr{S} \cup \mathscr{P}' \cup \mathscr{C}'_1 \cup \mathscr{C}_M$ (resp. $(f,\mathbf{1}) \in \mathscr{C}_M$) and $\psi$ is a character in $\mathcal{I}_f$ such that $\psi(\rho_K)^{p-1}=\{1\}$ (resp. $\psi(\rho_K)^2 \neq \{1\}$, $\psi(\rho_K)^4=\{1\}$). But since $\rho_{|G_{\Q_p}}$ is $1$-Cartan and $K/\Q$ is inert at $p$, $\rho(G_K)=C$, so this implies $\psi^{p-1}=\mathbf{1}$, which is impossible by definition of $\mathcal{I}_f$ (resp. $\psi^4=\mathbf{1}$, which means that $\psi$ is one of the characters attached to $f$, and this is impossible by Definition of $\mathcal{I}_f$).

 The fact that all the global root numbers are $-1$ can be read off Table \ref{root-numbers-cartan-Hf} for the factors of the form $V_f \otimes H_f(\rho)$, and off Proposition \ref{counting-If-signs} for the factors of the form $V_f \otimes \mrm{Ind}_K^{\Q}{\psi(\rho_K)}$. }

\bigskip

\rems{\begin{enumerate}[noitemsep,label=(\roman*)]
\item The main result of \cite{LFL} is that when $p > 1.4 \cdot 10^7$, and $\rho$ comes from an elliptic curve $E/\Q$ without complex multiplication, then $\rho(G_{\Q})=N$. The bound has been improved to $p > 37$ in the preprint \cite{lombardo}.
\item It is worth pointing out once again that if $\rho$ comes from an elliptic curve $E/\Q$, then $E$ defines a rational point on the modular curve $X_{ns}^+(p)$ (see for instance \cite{LFD}). Whether such points with non-CM $j$-invariants exist for any $p > 13$ is an open question.
\item Proposition \ref{counting-If-signs} suggests the following question: which divisors of $p+1$ arise as orders of $\phi_f$ for $(f,\mathbf{1}) \in \mathscr{C}'_1$? It seems possible to at least gather numerical evidence thanks to the work of Loeffler and Weinstein \cite{LWnewforms}. 
\end{enumerate}
}

\section{The Bloch-Kato conjecture for factors of the Tate module of $J_{\rho,\xi}(p)$}
\label{bloch-kato-conjecture}

The goal of this section is to prove that the factors of the Tate module of $J_{\rho,\xi}(p)$ attached to non-CM modular forms satisfy some version of the Bloch-Kato conjecture in rank zero.

\lem[limit-galois-action]{Let $(f,\chi) \in \mathscr{S} \cup \mathscr{C} \cup \mathscr{P}$, $\sigma \in G_{\Q}$ and $\chi$ be a primitive Dirichlet character. Assume that $(f',\chi') := (\sigma(f)\otimes \gamma,\sigma(\chi)\gamma^2) \in \mathscr{S} \cup \mathscr{P} \cup \mathscr{C}$. Then has $\gamma \in \mathcal{D}$ and the following holds: 
\begin{itemize}[noitemsep,label=$-$]
\item If $(f,\chi) \in \mathscr{S}$, then $\gamma=\mathbf{1}$ and $(f',\mathbf{1}) \in \mathscr{S}$. 
\item If $(f,\chi) \in \mathscr{P}$, then $\gamma \in \{\mathbf{1},\sigma(\chi^{-1})\}$ and $(f',\sigma(\chi)\gamma^2) \in \mathscr{P}$.
\item If $(f,\chi) \in \mathscr{C}$, then $(f',\gamma^2\sigma(\chi)) \in \mathscr{C}$ and $f'$ has complex multiplication if and only if $f$ does. 
\end{itemize}}

\demo{This follows from Lemma \ref{gamma1p-onlytwists} and Corollary \ref{level-newform-bigtwist}, as well as the definition of $\mathscr{C}$.}

\prop[representation-is-rational]{For every $(f,\chi) \in \mathscr{S}\cup\mathscr{P}\cup\mathscr{C}$, let $V_f$ denote the one-dimensional $F$-vector space on which $N$ acts trivially and $\mathbb{T}$ acts thusly: $T_n$ (resp. $nI_2$) acts by $a_n(f)$ (resp. $\chi(n)$). Let $(f,\chi) \in \mathscr{S}\cup\mathscr{P}\cup\mathscr{C}$. 

Let $\alpha \in H_f$. There is a prime ideal $I \subset \mathbb{T}_{1,N}(N)$ such that for any $(f',\chi') \in \mathscr{S} \cup \mathscr{P} \cup \mathscr{C}$:
\begin{itemize}[noitemsep,label=$-$]
\item for any $\alpha' \in H_{f'}$, $I \cdot (V_{f'} \otimes \alpha')$ is zero if, for some $\sigma \in G_{\Q}$ and $\gamma \in \mathcal{D}$, one has $(f',\alpha')=(\sigma(f) \otimes \gamma, \sigma(\alpha)\gamma(\det))$, and $V_{f'} \otimes \alpha'$ otherwise.
\item for any $\psi\in \mathcal{I}_{f'}$, $I \cdot (V_{f'} \otimes_F \mrm{Ind}_C^N{\psi}) = V_{f'} \otimes_F \mrm{Ind}_C^N{\psi}$.
\end{itemize}

Let $\psi \in \mathcal{I}_f$. There is a prime ideal $I \subset \mathbb{T}_{1,N}(N)$ such that for any $(f',\chi') \in \mathscr{S} \cup \mathscr{P} \cup \mathscr{C}$:
\begin{itemize}[noitemsep,label=$-$]
\item For any $\psi'\in \mathcal{I}_{f'}$, $I \cdot (V_{f'} \otimes_F \mrm{Ind}_C^N{\psi'}) = 0$ if, for some $\sigma \in G_{\Q}$ and $\gamma \in \mathcal{D}$, one has $f'=\sigma(f) \otimes \gamma$ and $\psi' \sim \sigma(\psi)\gamma(\det)$. Otherwise, $I \cdot (V_{f'} \otimes_F \mrm{Ind}_C^N{\psi'})= V_{f'} \otimes_F \mrm{Ind}_C^N{\psi'}$.
\item For any $\alpha \in H_{f'}$, $I \cdot (V_{f'} \otimes \alpha)= V_{f'} \otimes \alpha$.
\end{itemize}

}

\demo{
\emph{Step 1: Preliminaries.}
 
 The space $V_f$ is well-defined by (for instance) Corollary \ref{decomposition-hecke-xpp}. 

Let $(f',\chi') \in \mathscr{S} \cup \mathscr{P} \cup \mathscr{C}$. For $\alpha' \in H_{f'}$ (resp. $\psi' \in \mathcal{I}_{f'}$), let $M_{f',\alpha'}$ (resp. $M_{f',\psi'}$) denote the $(\mathbb{T} \otimes F)[N]$-module $V_{f'} \otimes \alpha'$ (resp. $V_{f'} \otimes \mrm{Ind}_C^N{\psi'}$). Note that if $\psi_1,\psi_2 \in \mathcal{I}_{f'}$ are such that $\psi_1 \sim \psi_2$, then $M_{f',\psi_1}$ and $M_{f',\psi_2}$ are isomorphic as $\mathbb{T}[N]$-modules. 

Let $\mathcal{M}$ be the collection of all the $M_{f',\alpha'}$ and $M_{f',\psi'}$ as above. Every $M \in \mathcal{M}$ is an absolutely irreducible $F[N]$-module. Therefore, by Schur's Lemma, $\mathbb{T}_{1,N}(N)$ (which is central in $\mathbb{T}[N]$) acts on $M$ by a $F$-valued character, which we call in the rest of this proof the \emph{central character}. Replacing $M$ with $M \otimes_{F,\sigma} F$ (for $\sigma \in G_{\Q}$) means post-composing the central character with $\sigma^{-1}$. 

Let $(f',\chi') \in \mathscr{S}\cup\mathscr{P}\cup\mathscr{C}$ and $\alpha' \in H_{f'}$ (resp. $\alpha' \in \mathcal{I}_{f'}$). It can be checked directly on generators that for every $\sigma \in G_{\Q}$ and every $g \in \mathbb{T}[N]$, $\mrm{Tr}(g \mid M_{\sigma(f'),\sigma\circ\alpha'}) = \sigma(\mrm{Tr}(g \mid M_{f',\alpha'}))$. Therefore, the central character of $M_{\sigma(f'),\sigma(\alpha')}$ is the post-composition of that of $M_{f',\alpha'}$ with $\sigma$. 

If in fact $(f',\chi') \in \mathscr{C}$ and $\gamma \in \mathcal{D}$, it can be checked on the generators that $M_{f'\otimes \gamma,\gamma(\det)\alpha'}$ and $M_{f',\alpha'}$ have the same traces for the actions of $\mathbb{T}_{1,N}$, so that they have the same central character. 

Let now $\alpha \in H_f$ (resp. $\alpha \in \mathcal{I}_f$). Let $I_{\alpha} \subset \mathbb{T}_{1,N}(N)$ be the kernel of the central character of $M_{f,\alpha}$ (which is a ring homomorphism $\mathbb{T}_{1,N}(N) \rar F$, so that $I_{\alpha}$ is a prime ideal). We will prove that, for any $M_{f',\alpha'} \in \mathcal{M}$, one of the following is true:
\begin{itemize}[noitemsep,label=$-$]
\item There exists $(\sigma,\gamma) \in G_{\Q} \times \mathcal{D}$ such that $f'=\sigma(f)\otimes\gamma$ and $\alpha'=\sigma(\alpha)\gamma(\det)$ (resp. $\alpha' \sim \sigma(\alpha)\gamma(\det)$) . Then $I_{\alpha}M_{f',\alpha'}=0$,
\item $IM_{f',\alpha'}=M_{f',\alpha'}$. 
\end{itemize}

If this claim is true, this clearly finishes the proof. Moreover, the implication in first clause is straightforward, since, if $(f',\alpha')=(\sigma(f)\otimes \gamma,\sigma(\alpha)\gamma(\det))$ (resp. $f'=\sigma(f)\otimes \gamma$ and $\alpha' \sim \sigma(\alpha)\gamma(\det)$), the central character of $M_{f',\alpha'}$ is the post-composition with $\sigma$ of the central character of $M_{f,\alpha}$. Since $IM_{f',\alpha'}$ is a $F[N]$-submodule of $M_{f',\alpha'}$, it is always $M_{f',\alpha'}$ or zero. Therefore, all that remains to prove is that if $I_{\alpha}M_{f',\alpha'}=0$ for $\alpha' \in H_{f'} \cup \mathcal{I}_{f'}$, then $\alpha' \in H_{f'}$ (resp. $\alpha' \in \mathcal{I}_{f'}$) and there exists $\sigma \in G_{\Q}$ and a Dirichlet character $\gamma \in \mathcal{D}$ such that $f'=\sigma(f) \otimes \gamma$ and $\alpha'=\sigma(\alpha)\gamma(\det)$ (resp. $\alpha' \sim \sigma(\alpha)\gamma(\det)$). 

Thus, from now on, we assume that $(f',\chi') \in \mathscr{S}\cup\mathscr{P}\cup\mathscr{C}$ and $\alpha'\in H_{f'}\cup\mathcal{I}_{f'}$ and that $I_{\alpha}M_{f',\alpha'}=0$. If we need another generic element of $\mathscr{S}\cup\mathscr{P}\mathscr{C}$, we will denote it by $(g,\theta)$, and let $\beta \in H_g \cup \mathcal{I}_g$. \\

\emph{Step 2: $\alpha' \in H_{f'}$ (resp. $\alpha' \in \mathcal{I}_{f'}$).}

Let $h \in C$ be a generator. Then $h^{p-1}$ generates the cyclic group $C \cap \SL{\F_p}$. If $\beta \in H_{g}$, then we can check that in each case, $\beta(h^{p-1})=1$. Thus the element $h^{p-1}+h^{1-p}-2 \in \mathbb{T}_{1,N}(N)$ annihilates $M_{g,\beta}$. Now, if $\beta \in \mathcal{I}_g$, $\beta(g^{p-1})$ is a $(p+1)$-th root of unity; it is distinct from $1$ because $\beta^{p-1}\neq \mathbf{1}$ by the definitions. Then the action of $h^{p-1}+h^{1-p}-2$ on $M_{g,\beta}$ is given by a nonzero algebraic integer. 

Hence, if $\alpha \in H_f$, $I_{\alpha}$ contains $h^{p-1}+h^{1-p}-2$, so $(h^{p-1}+h^{1-p}-2)M_{f',\alpha'}=0$, so $\alpha' \in H_{f'}$. 

If $\alpha \in \mathcal{I}_f$, there is a monic polynomial $\Pi \in \Z[t]$ with $\Pi(0) \neq 0$ such that $\Pi(h^{1-p}+h^{p-1}-2) \in I_{\alpha}$, hence $\Pi(h^{1-p}+h^{p-1}-2)M_{f',\alpha'}=0$. Therefore $h^{p-1}$ does not act trivially on $M_{f',\alpha'}$, thus $\alpha' \in \mathcal{I}_{f'}$. \\

\emph{Step 3: We may assume that $f'=f$ and that $M_{f,\alpha'}$ has the same central character as $M_{f,\alpha}$.}

$(\mathbb{T}_{1,N}(N)\otimes\Q)/I_{\alpha}$ is a $\Q$-algebra. It embeds in $F$ by two ring homomorphisms $\iota_{f,\alpha}$ and $\iota_{f',\alpha'}$ (given by the central characters of $M_{f,\alpha}$ and $M_{f',\alpha'}$): the images of $\iota_{f,\alpha}, \iota_{f',\alpha'}$ are sub-$\Q$-algebras of $F$, so they are fields. The morphism $\iota_{f',\alpha'}\circ \iota_{f,\alpha}^{-1}: \im{\iota_{f,\alpha}} \rar \im{\iota_{f',\alpha'}}$ is an isomorphism, so it extends as an automorphism $\sigma \in \mrm{Gal}(F/\Q)$. Thus the central character of $M_{f',\alpha'}$ is the post-composition with $\sigma$ of the central character of $M_{f,\alpha}$. In particular, $M_{\sigma^{-1}(f'), \sigma^{-1}\circ\alpha'}$ has the same central character as $M_{f,\alpha}$, so $I_{\alpha}M_{\sigma^{-1}(f'), \sigma^{-1}\circ\alpha} = 0$. 

By evaluating the central character, we see that for any $n \equiv 1\pmod{p}$, $a_n(\sigma^{-1}(f'))=a_n(f)$, so, by Lemma \ref{agreement-implies-twist}, there exists a character $\gamma \in \mathcal{D}$ such that $f=\sigma^{-1}(f') \otimes \gamma$. Thus $M_{f,\sigma^{-1}(\alpha')\gamma(\det)}$ has the same central character as $M_{f,\alpha}$. 

In this step, we assume that the problem is solved whenever $f=f'$ and $M_{f,\alpha}, M_{f',\alpha'}$ have the same central character. Therefore, there exists $\tau \in G_{\Q}$ such that $\tau(f)=f$ and $\sigma^{-1}(\alpha')\gamma(\det)=\tau(\alpha)$ (resp. $\sigma^{-1}(\alpha')\gamma(\det) \sim \tau(\alpha)$). Then one checks that $f'=\sigma(\tau(f))\otimes \sigma(\gamma^{-1})$ and $\alpha'=\sigma(\tau(\alpha))\sigma(\gamma^{-1}(\det))$ (resp. $\alpha'\sim\sigma(\tau(\alpha))\sigma(\gamma^{-1}(\det))$), whence the conclusion. \\

\emph{Step 4: One has $\alpha'=\alpha$ (resp. $\alpha' \sim \alpha$).}

We first assume that $f$ has no complex multiplication. 

It is enough to prove that for every $h \in N$, the trace of the action of $h$ on the $F$-vector spaces $M_{f,\alpha'}$ and $M_{f,\alpha}$ is the same. By checking the definitions, characters of $H_f$ are determined by their restrictions to $C$. Moreover, if $\psi \in \mathcal{I}_f$ and $h \in N$, then the action of $h$ on $M_{f,\psi}$ is zero. So it is enough to check the equality for $h \in C$.

In this case, let $q$ be a prime number such that $q\det{h}\equiv 1\pmod{p}$, then $T_q(h+h^p) \in \mathbb{T}_{1,N}(N)$. Since $M_{f,\alpha}$ and $M_{f,\alpha'}$ have the same central character, one has $a_q(f)\alpha(h)=a_q(f)\alpha'(h)$ (resp. $a_q(f)(\alpha(h)+\alpha(h)^p)=a_q(f)(\alpha'(h)+\alpha'(h)^p)$). By a result of Serre \cite[Cor. 2 au Th\'eor\`eme 15]{Serre-Cebotarev}, we can choose $q$ such that $a_q(f)\neq 0$, so $\mrm{Tr}(h \mid M_{f,\alpha})=\mrm{Tr}(h \mid M_{f,\alpha'})$, and we are done. 

When $f$ has complex multiplication, $H_f$ contains a single character, so we may assume that $\alpha,\alpha' \in \mathcal{I}_f$. By definition of $\sim$, it is enough to check that for any $h \in C$ with $\det{h} \in \F_p^{\times 2}$, $\mrm{Tr}(h \mid M_{f,\alpha})= \mrm{Tr}(h \mid M_{f,\alpha'})$. For such an $h$, there exists an integer $n \geq 1$ such that $n^2\det{h}=1$, so that $(nI_2)(h+h^p) \in \mathbb{T}_{1,N}(N)$, and therefore $\chi(n)\mrm{Tr}(h \mid M_{f,\alpha})=\chi(n)\mrm{Tr}(h \mid M_{f,\alpha'})$, whence the conclusion. 
}

\medskip 

\lem[make-abelian-variety]{Let $I \subset \mathbb{T}_{1,N}(\rho(G_{\Q}))$ (resp. $I \subset \mathbb{T}_{1,N}(N)$) be a subgroup. Then, for any $\xi \in \mu_p^{\times}(\Qbar)$, the abelian variety $A_I := J_{\rho,\xi}(p)/IJ_{\rho,\xi}(p)$ over $\Q$ is well-defined. Moreover, $(A_I)_{\Qbar}$ is endowed with an action of $\mathbb{T}_{1,\rho(G_{\Q})}$ (resp. $\mathbb{T}_{1,N}$) which extends to an action of $\mathbb{T}_{1,\rho(G_{\Q}),\rho}$ (resp. $\mathbb{T}_{1,N,\rho}$) on $A_I(\Qbar)$. Finally, $A_I$ only depends (up to $\Q$-isomorphism) on $I(\Tate{\ell}{J_{\rho,\xi}(p)} \otimes_{\Z_{\ell}} F_{\mathfrak{l}})$ for any maximal ideal $\mathfrak{l} \subset \OO_F$ with residue characteristic $\ell$.}

\demo{Since $\mathbb{T}$ is a finitely generated abelian group, so is $I$. By Lemma \ref{T1Gamma-action}, $\mathbb{T}_{1,N}(\rho(G_{\Q}))$ acts on $J_{\rho,\xi}(p)$. Hence $IJ_{\rho,\xi}(p)$ is an abelian subvariety of $J_{\rho,\xi}(p)$ by Proposition \ref{abelian-image-exists}, and the quotient $A_I$ is well-defined with an exact sequence $0 \rar I \cdot J_{\rho,\xi}(p)(\Qbar) \rar J_{\rho,\xi}(p)(\Qbar) \rar A_I(\Qbar) \rar 0$ by Proposition \ref{quotient-of-abelian-varieties}.  

The fact that $\mathbb{T}_{1,\rho(G_{\Q})}$ (resp. $\mathbb{T}_{1,N}$) acts on $(A_I)_{\Qbar}$ comes from Lemma \ref{T1Gamma-action} and the fact that $I$ commutes to $\mathbb{T}_{1,\rho(G_{\Q})}$ (resp. $\mathbb{T}_{1,N}$). The relationship of this action to that of $G_{\Q}$ is given by Lemma \ref{twist-mix-gl2-galois}.

If $I'$ is another subgroup such that $I(\Tate{\ell}{J_{\rho,\xi}(p)} \otimes_{\Z_{\ell}} F_{\mathfrak{l}})=I'(\Tate{\ell}{J_{\rho,\xi}(p)} \otimes_{\Z_{\ell}} F_{\mathfrak{l}})$, then, after replacing $I'$ with $I+I'$, we may assume that $I \subset I'$. Then $IJ_{\rho,\xi}(p) \subset I'J_{\rho,\xi}(p)$ are Abelian subvarieties of $J_{\rho,\xi}(p)$ such that 
\[\Tate{\ell}{IJ_{\rho,\xi}(p)} \otimes_{\Z_{\ell}} F_{\mathfrak{l}} = I(\Tate{\ell}{J_{\rho,\xi}(p)} \otimes_{\Z_{\ell}} F_{\mathfrak{l}})=I'(\Tate{\ell}{IJ_{\rho,\xi}(p)} \otimes_{\Z_{\ell}} F_{\mathfrak{l}})=\Tate{\ell}{I'J_{\rho,\xi}(p)} \otimes_{\Z_{\ell}} F_{\mathfrak{l}},\] so that $IJ_{\rho,\xi}(p)$ and $I'J_{\rho,\xi}(p)$ have the same dimension. They are therefore equal. 
}

\medskip

\cor[make-abelian-variety-for-Hf]{Let $\mathscr{E}=\mathscr{E}'=\mathscr{S}$ (resp. $(\mathscr{E},\mathscr{E}')=(\mathscr{P},\mathscr{P}')$, resp. $(\mathscr{E},\mathscr{E}')=(\mathscr{C} \backslash \mathscr{C}_M,\mathscr{C}'_1)$, resp. $\mathscr{E}=\mathscr{E}'=\mathscr{C}_M$). Let $\pi: \mathscr{E} \rar \mathscr{E}'$ be the identity (resp. picking the representative modulo complex conjugation, resp. picking the representative modulo twists by characters in $\mathcal{D}$, resp. the identity). Let $(f,\chi) \in \mathscr{E}$ and $\alpha \in H_f$. There exists a prime ideal $I \subset \mathbb{T}_{1,N}(N)$ satisfying the following properties: 
\begin{itemize}[noitemsep,label=$-$]
\item The Abelian variety $A_I = J_{\rho,\xi}(p)/IJ_{\rho,\xi}(p)$ over $\Q$ is well-defined, 
\item $A_I$ is endowed with an action of $\mathbb{T}_{1,N}$ over $\Qbar$ which extends to an action of $\mathbb{T}_{1,N,\rho}$ on $A_I(\Qbar)$,
\item For any maximal ideal $\mathfrak{l} \subset \OO_F$ with residue characteristic $\ell$, $\Tate{\ell}{A_I} \otimes_{\Z_{\ell}} F_{\mathfrak{l}}$ is isomorphic to the restriction to $\mathbb{T}_{1,N,\rho} \otimes F_{\mathfrak{l}}$ of the $(\mathbb{T} \otimes F_{\mathfrak{l}})[N \rtimes_{\rho} G_{\Q}]$-module 
\[\prod_{\substack{(f',\chi') \in \mathscr{E}'\\\alpha' \in H_{f'}}}{V_{f',\mathfrak{l}} \otimes \alpha'(\rho)},\] where $T_n$ (resp. $nI_2$, resp. $M \in N$) acts on the factor $V_{f',\mathfrak{l}} \otimes \alpha'(\rho)$ by $a_n(f')$ (resp. $\chi'(n)$, resp. $\alpha'(M)$), and we only take the product over triples $(f',\chi',\alpha')$ such that, for some $\sigma \in \mrm{Aut}(F/\Q)$ and some $\gamma \in \mathcal{D}$, one has $f'=\sigma(f) \otimes \gamma$ and $\alpha'=\sigma(\alpha)\gamma(\det)$. In particular, one has $(V_{f',\mathfrak{l}} \otimes \alpha'(\rho))_{\mathfrak{l}} = \sigma\left[(V_{f,\mathfrak{l}} \otimes \alpha(\rho))_{\mathfrak{l}}\right]$.
\item Let $\mathfrak{l} \subset \OO_F$ be a maximal ideal, and consider the character $\mathbb{T}_{1,N}(N) \rar F_{\mathfrak{l}}$ given by the action of $\mathbb{T}_{1,N}(N)$ on $V_{f,\mathfrak{l}} \otimes \alpha$. Then $\Tate{\ell}{A_I} \otimes_{\mathbb{T}_{1,N}(N)} F_{\mathfrak{l}} \simeq V_{f,\mathfrak{l}} \otimes \alpha(\rho)$. 
\end{itemize}
Moreover, the abelian subvariety $IJ_{\rho,\xi}(p)$ (and therefore $A_I$) does not depend on the choice of $I$. 
}

\demo{Proposition \ref{decomposition-ncartan-xrho-T11rho} yields the structure of $\Tate{\ell}{J_{\rho,\xi}(p)} \otimes_{\Z_{\ell}} F_{\mathfrak{l}}$ as a $\mathbb{T}_{1,N,\rho}$-module. The first two items are true for any subgroup $I \subset \mathbb{T}_{1,N}(N)$ by Lemma \ref{make-abelian-variety}; a prime ideal $I$ satisfying the third item is then given by Proposition \ref{representation-is-rational}.

Thus, only the fourth item needs to be checked. It amounts to checking the following assertion: let $\sigma \in G_{\Q}, \sigma \in \mathcal{D}$, $(f',\chi') \in \mathscr{E}'$ and $\alpha' \in H_{f'}$ be such that $(f',\alpha')=(\sigma(f)\otimes \gamma,\sigma(\alpha') \gamma(\det)$ and such that $\mathbb{T}_{1,N}(N)$ acts on $V_{f,\mathfrak{l}} \otimes \alpha$ and $V_{f',\mathfrak{l}} \otimes \alpha'$ by the same $F_{\mathfrak{l}}$-valued character, and let us show that $(f,\alpha)=(f',\alpha')$. 

Indeed, for any $n \equiv 1\pmod{p}$, $T_n$ acts on the first module by $a_n(f)$ and by $a_n(f')$ on the second one, so one has $a_n(f)=a_n(f')$. By Lemma \ref{agreement-implies-twist}, $f'$ is a twist of $f$. By construction of $\mathscr{E}'$, this implies that $f=f'$. Then $\alpha=\alpha'$ by Step 4 of the proof of Proposition \ref{representation-is-rational}.  
}

\medskip

\prop[BK-for-Hf]{In the situation of Proposition \ref{make-abelian-variety-for-Hf}, assume that $L(f \otimes \alpha(\rho),1) \neq 0$. 
Then for all $(f',\alpha')$ appearing in the product of the above Proposition, $L(f' \otimes \alpha'(\rho),1) \neq 0$. 
Moreover, for any prime $q$, the $q$-divisible Bloch-Kato Selmer group $\mrm{Sel}(A_I,A[q^{\infty}])$ is finite, and this Selmer group is zero for all but finitely many $q$. In particular, $A_I(\Q)$ is is finite. }

\demo{The non-vanishing comes from the fact that the newforms attached to the $f' \otimes \alpha'(\rho)$ are Galois-conjugate to the newform attached to $f \otimes \alpha(\rho)$ (and these constructions define newforms by Proposition \ref{twist-exists}), so we can apply \cite[Theorem 1]{Shimura-Periods}. The result is then a consequence of \cite[Theorem 14.2]{KatoBSD} and its CM analogue, proved in \S 15 of \emph{loc.cit.}. 
}

\medskip

\prop[make-abelian-variety-If]{
Define $\tilde{\rho}: (M,g) \in N \rtimes_{\rho} G_{\Q} \longmapsto N\rho(g)$ and $\tilde{\rho_K}: \tilde{\rho}^{-1}(C) \rar C$. Let $\pi,\mathscr{E},\mathscr{E}',(f,\chi)$ be as in Corollary \ref{make-abelian-variety-for-Hf}. Let $\psi \in \mathcal{I}_f$. There exists a prime ideal $I \subset \mathbb{T}_{1,N}(N)$ satisfying the following properties: 
\begin{itemize}[noitemsep,label=$-$]
\item The Abelian variety $A_I = J_{\rho,\xi}(p)/IJ_{\rho,\xi}(p)$ over $\Q$ is well-defined, 
\item $A_I$ is endowed with an action of $\mathbb{T}_{1,N}$ over $\Qbar$ which extends to an action of $\mathbb{T}_{1,N,\rho}$ on $A_I(\Qbar)$,
\item For any maximal ideal $\mathfrak{l} \subset \OO_F$ with residue characteristic $\ell$, $\Tate{\ell}{A_I} \otimes_{\Z_{\ell}} F_{\mathfrak{l}}$ is isomorphic to the restriction to $\mathbb{T}_{1,N,\rho} \otimes F_{\mathfrak{l}}$ of the $(\mathbb{T} \otimes F_{\mathfrak{l}})[N \rtimes_{\rho} G_{\Q}]$-module 
\[\prod_{(f',\chi') \in \mathscr{E}',\psi' \in \mathcal{I}_{f'}/\sim}{V_{f',\mathfrak{l}} \otimes \mrm{Ind}_K^{\Q}{\psi'(\tilde{\rho}_K)}},\] where:
\begin{itemize}[noitemsep,label=\tiny$\bullet$]
\item we sum over triples $(f',\chi',\psi')$ such that there exists $\sigma \in \mrm{Aut}(F/\Q)$ and $\gamma \in \mathcal{D}$ such that $(f',\psi')=(\sigma(f) \otimes \gamma,\sigma(\psi) \otimes \gamma(\det))$. In particular, in such a situation, for any maximal ideal $\mathfrak{l}$ of $\OO_F$, for all but finitely many prime numbers $q$, \[\mrm{Tr}\left(\Fr_q \mid V_{f',\mathfrak{l}} \otimes \mrm{Ind}_K^{\Q}{\psi'(\rho_K)}\right) = \sigma\left[\mrm{Tr}\left(\Fr_q \mid V_{f,\mathfrak{l}} \otimes \mrm{Ind}_K^{\Q}{\psi(\rho_K)}\right)\right].\] 
\item $T_n$ (resp. $nI_2$) acts on the factor $V_{f',\mathfrak{l}} \otimes \mrm{Ind}_K^{\Q}{\psi'(\tilde{\rho}_K)}$ by $a_n(f')$ (resp. $\chi'(n)$).    
\end{itemize}
\item Let $\mathfrak{l} \subset \OO_F$ be a maximal ideal, and consider the character $\mathbb{T}_{1,N}(N) \rar F_{\mathfrak{l}}$ given by the action of $\mathbb{T}_{1,N}(N)$ on $V_{f,\mathfrak{l}} \otimes \mrm{Ind}_C^N{\psi}$. Then $\Tate{\ell}{A_I} \otimes_{\mathbb{T}_{1,N}(N)} F_{\mathfrak{l}} \simeq V_{f,\mathfrak{l}} \otimes \mrm{Ind}_C^N{\psi(\rho_K)}$.
\end{itemize}
Moreover, the abelian subvariety $IJ_{\rho,\xi}(p)$ (and therefore $A_I$) does not depend on the choice of $I$. 
}

\demo{By mimicking the proof of Corollary \ref{make-abelian-variety-for-Hf}, the only two statements to prove are that $I$ is a prime ideal and the following: let $(f',\chi') \in \mathscr{E}'$, $\psi' \in \mathcal{I}_{f'}$, $\sigma \in G_{\Q}$ and $\gamma \in \mathcal{D}$, and assume that $(f',\psi')=(\sigma(f)\otimes \gamma,\sigma(\psi)\gamma(\det))$. Assume furthermore that $\mathbb{T}_{1,N}(N)$ acts by the same $F_{\mathfrak{l}}$-valued character on $V_{f,\mathfrak{l}} \otimes \mrm{Ind}_C^N{\psi}$ and $V_{f',\mathfrak{l}} \otimes \mrm{Ind}_C^N{\psi'}$. Then $f=f'$ and $\psi\sim\psi'$. 

As in the proof of Corollary \ref{make-abelian-variety-for-Hf}, by considering the action of $T_n$ for $n \equiv 1\pmod{p}$, we see that $a_n(f)=a_n(f')$, so by Lemma \ref{agreement-implies-twist}, $f$ is a twist of $f'$. By Step 4 in the proof of Proposition \ref{representation-is-rational}, one has $\psi\sim\psi'$.}

\prop[BK-for-If]{Let us consider the situation of Proposition \ref{make-abelian-variety-If}.
\begin{enumerate}[noitemsep,label=(\roman*)]
\item\label{BKIf1} Assume that $\psi(\rho_K)^{p-1} \neq \mathbf{1}$. Then the Artin representation $\mrm{Ind}_K^{\Q}{\psi(\rho_K)}$ is attached to a newform $g \in \mathcal{S}_1(\Gamma_1(D))$. 
\item\label{BKIf2} Assume moreover that $L(f \times g,1) \neq 0$. Let $(\sigma,\gamma) \in G_{\Q} \times \mathcal{D}$ be such that \[(f',\chi')=(\sigma(f)\otimes \gamma,\sigma(\chi)\gamma^2) \in \mathscr{E}'.\] Let $\psi'=\sigma(\psi)\gamma(\det)\in \mathcal{I}_{f'}$. Then the Artin representation $\mrm{Ind}_K^{\Q}{\psi'(\rho_K)}$ is attached to the newform $g'=\sigma(g) \otimes \gamma \in \mathcal{S}_1(\Gamma_1(D'))$, and one has $L(f' \otimes g',1) \neq 0$. 
\item\label{BKIf3} Assume furthermore that $f$ does not have complex multiplication and $K \neq \Q(\sqrt{-p})$. Let $\mathcal{Q}$ be a set of maximal ideals $\mathfrak{q} \subset \OO_F$ with residue characteristic $q$ satisfying the following conditions:
\begin{itemize}[label=\tiny$\bullet$]
\item $q \nmid 6pD\delta$, where $\delta$ is the cardinality of the image of $\mrm{Ind}_K^{\Q}{\psi(\rho_K)}(G_{\Q})$ in $\PGL{F}$. 
\item $\mrm{Ind}_K^{\Q}{\psi}(\Fr_q)$ is not scalar,
\item $f$ is ordinary at $\mathfrak{q}$.
\end{itemize}
Then, the set of prime numbers which are residue characteristics of some $\mathfrak{q} \in \mathcal{Q}$ has density $1-\delta^{-1}$. Moreover, for all but finitely many $\mathfrak{q} \in \mathcal{Q}$, if $\mathfrak{q}_0$ denotes the inverse image of $\mathfrak{q}$ under the character $\mathbb{T}_{1,N}(N) \rar F_{\mathfrak{q}}$ attached to $(f,\psi)$, the Bloch-Kato Selmer group $\mrm{Sel}(\Q,A_I[\mathfrak{q}_0^{\infty}])$ is finite. In particular, $A_I(\Q)$ is finite.  
\end{enumerate}
}

\rem{There are several possible normalizations for the Rankin-Selberg $L$-function attached to two newforms $f,g$ of weights $k \geq l$. 
\begin{itemize}[noitemsep,label=$-$]
\item The \emph{imprimitive Rankin-Selberg $L$-function} $L^{imp}(f,g,s)$ is the meromorphic continuation of the series $\sum_{n \geq 1}{a_n(f)a_n(g)n^{-s}}$; it is denoted by $D(s,f,g)$ in \cite{Shimura76}. 
\item We define \emph{Rankin-Selberg-Li $L$-function} by the construction of \cite[\S 2]{Li-RS}: more precisely, by Theorem 2.1 of \emph{loc.cit.}, there is a primitive Dirichlet character $\gamma$ such that the couple $(g \otimes \gamma, \overline{f}\otimes \gamma)$ satisfies the conditions (A), (B), (C) of \emph{loc.cit.}, and the Rankin-Selberg-Li $L$-function attached to $f,g$ is exactly the meromorphic continuation $L^{RSL}(f,g,s)$ of $s \longmapsto L_{g \otimes \gamma, \overline{f}\otimes \gamma}(s-\frac{k+l}{2}+1)\prod_{q \mid M"}{\theta_q(s-\frac{k+l}{2}+1)^{-1}}$ in the notation of \emph{loc.cit.} (the shift in $s$ is solely to ensure that the functional equations for $L^{imp}(f,g,s)$ and $L^{RSL}(f,g,s)$ are of the form $s \leftrightarrow k+l-1-s$). 
\item The \emph{Rankin-Selberg $L$-function} is the meromorphic continuation $L(f \times g,s)$ of the Euler product $\prod_p{L(D_{f,p}\otimes D_{g,p},p^{-s})}$. 
\end{itemize}
The following Lemma proves that, for the purposes of Proposition \ref{BK-for-If}, any of these $L$-functions works. 
}

\lem[multiple-rankin-selberg]{Let $f \in \mathcal{S}_k(\Gamma_1(A)), g \in \mathcal{S}_l(\Gamma_1(B))$ be two newforms with characters $\chi_f, \chi_g$ and weights $k \geq l$ such that $\chi_f\chi_g$ is nontrivial. Let $s$ be a complex number with $\Real{s} \geq \frac{k+l-1}{2}$. Then $L^{imp}(f,g,s), L^{RSL}(f,g,s), L(f\times g,s)$ are all defined, and they are either all zero, or all nonzero. }

\demo{We first show that for any prime $p$, the series $\sum_{n\geq 0}{a_{p^n}(f)a_{p^n}(g)p^{-ns}}$ is absolutely convergent and its sum is nonzero. Indeed, we know that the following formal identities hold: 
\[\sum_{n\geq 0}{a_{p^n}(f)p^{-ns}} = [(1-\alpha p^{-s})(1-\alpha' p^{-s})]^{-1},\quad, \sum_{n \geq 0}{a_{p^n}(g)p^{-ns}} = [(1-\beta p^{-s})(1-\beta'p^{-s})]^{-1}\] where $\alpha,\alpha',\beta,\beta'$ are algebraic integers such that $|\alpha|, |\alpha'| \leq p^{\frac{k-1}{2}}$ and $|\beta|,|\beta'|\leq p^{\frac{l-1}{2}}$ by the Ramanujan conjecture. Then, by \cite[Lemma 1]{Shimura76}, the following formal identity holds:
\[\sum_{n \geq 0}{a_{p^n}(f)a_{p^n}(g)p^{-ns}} = \frac{1-\alpha\alpha'\beta\beta' p^{-2s}}{(1-\alpha\beta p^{-s})(1-\alpha\beta' p^{-s})(1-\alpha'\beta p^{-s})(1-\alpha'\beta' p^{-s})}.\]

When $\Re{s} > \frac{k+l}{2}-1$, $u \in \{\alpha,\alpha'\}$, $v \in \{\beta,\beta'\}$, one has $|u v p^{-s}| < p^{\frac{k+l-2}{2}}p^{-\frac{k+l}{2}+1} < 1$, so, in the right half-plane $\{\Re{s} > \frac{k+l}{2}-1\}$, the previous formal identity holds in $\C$ and none of the factors in the right-hand side vanishes. 

As a consequence, for any primitive Dirichlet character $\gamma$, the claim for the couple $(f,g)$ and the couple $(f \otimes \gamma,g \otimes \overline{\gamma})$ are equivalent. Therefore, by \cite[Theorem 2.1]{Li-RS}, we may assume that $(g,\overline{f})$ satisfies the conditions A), B), C) of \emph{loc.cit.}. 

It is well-known that $\frac{1}{\Gamma}$ extends to an entire function \cite[Chapter XV, \S 2]{LangComplex}, and that, for any $s \in \C$ with $\Re{s} \geq 1$, $L(\chi_g\overline{\chi_f},s) \neq 0$ (for instance by combining \cite[Theorem 5.10]{IK} with the results of \cite[Section VII.2]{Neukirch-ANT}). By \cite[Theorem 2.2]{Li-RS} and the ensuing discussion, it is enough to show that, in the notation of \emph{loc.cit.}, for any $s \in \C$ with $\Re{s} \geq \frac{1}{2}$, one has $\theta_q(s) \neq 0$. This is done by inspecting all cases with the help of Proposition \ref{bad-L-factor}. }

Before coming back to Proposition \ref{BK-for-If}, we discuss the cohomology of Abelian varieties over local fields. This material seems standard, but we spell it out for the sake of completeness. 

\lem[galois-cohom-abelian]{Let $A$ be an Abelian variety over a non-Archimedian local field $K$ of characteristic zero. Let $A^{\vee}$ be the dual abelian variety.
\begin{enumerate}[noitemsep,label=(\roman*)]
\item\label{gca-1} Let $n \geq 1$. Local Tate duality and the Weil pairing induce a perfect pairing 
\[H^1(\mrm{Gal}(\overline{K}/K),A[n]) \times H^1(\mrm{Gal}(\overline{K}/K),A^{\vee}[n]) \rar \Q/\Z.\]  
\item\label{gca-2} Let $n \geq 1$ and $H^1_{\kappa}(K,A[n])$ be the image of the Kummer map $A(K)/nA(K) \rar H^1(K,A[n])$. Then $H^1_{\kappa}(K,A[n])$ and $H^1_{\kappa}(K,A^{\vee}[n])$ are orthogonal complements under local Tate duality. 
\item\label{gca-3} Assume that $q$ is a prime number distinct from the residue characteristic of $K$. Then, as $r\rar \infty$, the length of the cokernel of \[\mu_{A,r}: H^1(\mrm{Gal}(K^{nr}/K),(\Tate{q}{A})(K^{nr})) \rar H^1(\mrm{Gal}(K^{nr}/K),A[q^r](K^{nr}))\] is bounded. If $A$ has good reduction, this cokernel is trivial. 
\item\label{gca-4} Assume that $q$ is a prime number distinct from the residue characteristic of $K$. Then, one has
\[\underset{r \geq 1}{\sup}\,\left[\mrm{len}_{\Z_q}\left(\frac{H^1(\mrm{Gal}(K^{nr}/K),A[q^r](K^{nr}))+H^1_{\kappa}(K,A[q^r])}{H^1(\mrm{Gal}(K^{nr}/K),A[q^r](K^{nr}))\cap H^1_{\kappa}(K,A[q^r])}\right)\right] < \infty,\]
and this length is always zero when $A$ has good reduction. 
\end{enumerate}}

\demo{\ref{gca-1} is a consequence of \cite[\S 12, \S 16]{MilAb} and local Tate duality \cite[Theorem 7.2.15]{CNB}. 
For \ref{gca-2}, we show that $H^1_{\kappa}(K,A[n])$ and $H^1_{\kappa}(K,A^{\vee}[n])$ are orthogonal and that the product of their cardinalities is $|H^1(K,A[n])|$. Let $g$ be the dimension of $A$. By the local Euler-Poincar\'e characteristic calculation \cite[Theorem 7.3.1]{CNB}, the calculation of the degree of $A[n]$ \cite[Theorem 8.2]{MilAb}, and \cite[Lemma I.3.3]{MilneADT} 
\begin{align*}
|H^1(K,A[n])| &= |H^0(K,A[n])||H^2(K,A[n])| \|n^{2g}\|_K^{-1}= |A[n](K)||A^{\vee}[n]|\|n^g\|_K^{-1}\|n^g\|_K^{-1}\\
& = |A(K)/nA(K)||A^{\vee}(K)/nA^{\vee}(K)| = |H^1_{\kappa}(K,A[n])||H^1_{\kappa}(K,A^{\vee}[n])|
\end{align*}

To prove that $H^1_{\kappa}(K,A[n])$ and $H^1_{\kappa}(K,A^{\vee}[n])$ are orthogonal to one another, we use the fact that the pairing $H^r(K,A(\overline{K})) \times H^{1-r}(K,A^{\vee}(\overline{K})) \rar \Q/\Z$ (for $r=0,1$) of \cite[Corollary I.3.4]{MilneADT} is compatible with the local duality pairing. 

We deal with \ref{gca-3}. When $A$ has good reduction, the cokernel of $\mu_{A,r}$ is, by the N\'eron-Ogg-Shafarevich criterion \cite[Theorem 1]{GoodRed}, exactly equal to \[\mrm{coker}\left[H^1(\mrm{Gal}(K^{nr}/K),\Tate{q}{A}) \rar H^1(\mrm{Gal}(K^{nr}/K),A[q^r])\right] \subset H^2(\mrm{Gal}(K^{nr}/K),\Tate{q}{A}).\] Since $\mrm{Gal}(K^{nr}/K) \simeq \hat{\Z}$ has cohomological dimension one, $\mu_{A,r}$ is trivial. 

In general, let $M', M, M_2$ be the kernel, image, and cokernel of $(\Tate{q}{A})(K^{nr}) \rar A[q^r](K^{nr})$. We have the following exact sequences:
\begin{align*}
H^1(\mrm{Gal}(K^{nr}/K),(\Tate{q}{A})(K^{nr})) \rar H^1(\mrm{Gal}(K^{nr}/K),M) \rar H^2(\mrm{Gal}(K^{nr}/K),M'),\\
H^1(\mrm{Gal}(K^{nr}/K),M) \rar H^1(\mrm{Gal}(K^{nr}/K),A[q^r](K^{nr})) \rar H^1(\mrm{Gal}(K^{nr}/K),M_2)
\end{align*}
 which means that the length of the cokernel of $\mu_{A,r}$ is bounded by the sum of the lengths of $H^2(\mrm{Gal}(K^{nr}/K),M') = 0$ and $H^1(\mrm{Gal}(K^{nr}/K),M_2)$. By \cite[Lemma 1]{WashFLT}, the length of the cokernel of $\mu_{A,r}$ is at most the length of $M_2$. Now, by the cohomology long exact sequence for the inertia group of $K$, $M_2$ injects into $H^1(\mrm{Gal}(\overline{K}/K^{nr}),\Tate{q}{A})[q^r]$, so it is enough to prove that $H^1(\mrm{Gal}(\overline{K}/K^{nr}),\Tate{q}{A})$ is a finitely generated $\Z_q$-module. 
 
Let $I_q$ be the pro-$q$-quotient of the tame inertia quotient of $\mrm{Gal}(\overline{K}/K^{nr})$, corresponding to an extension $K^{q-tr}$. Then the kernel of $\mrm{Gal}(\overline{K}/K^{nr}) \rar I_q$ has supernatural order coprime to $q$, so, by the inflation-restriction exact sequence, $H^1(\mrm{Gal}(\overline{K}/K^{nr}),\Tate{q}{A}) \simeq H^1(I_q,(\Tate{q}{A})(K^{q-tr}))$. Since $I_q \simeq \Z_q$ is pro-cyclic, it is standard (the same argument as \cite[Lemma 1]{WashFLT} works) that $H^1(I_q,(\Tate{q}{A})(K^{q-tr}))$ is exactly the group of co-invariants of $(\Tate{q}{A})(K^{q-tr})$ under $I_q$. In particular, it is a subquotient of $\Tate{q}{A}$ so is finitely generated, and we are done. 

Now, we prove \ref{gca-4}. By \cite[Proposition I.3.8]{MilneADT}, there is an injection 
\begin{align*}
\frac{H^1(\mrm{Gal}(K^{nr}/K),A[q^r](K^{nr}))}{H^1(\mrm{Gal}(K^{nr}/K),A[q^r](K^{nr})) \cap H^1_{\kappa}(K,A[q^r])} &\rar H^1(\mrm{Gal}(K^{nr}/K),A(K^{nr})) \\&\simeq H^1(\mrm{Gal}(K^{nr}/K),\pi_0),\end{align*} where $\pi_0$ is the finite abelian group of geometric connected components of the special fibre of the N\'eron model of $A$. Thus 
\[\underset{r \geq 1}{\sup}\,\left[\mrm{len}_{\Z_q}\left(\frac{H^1(\mrm{Gal}(K^{nr}/K),A[q^r](K^{nr}))}{H^1_{\kappa}(K,A[q^r])\cap H^1(\mrm{Gal}(K^{nr}/K),A[q^r](K^{nr}))}\right)\right] < \infty\] and the same holds when $A$ is replaced with $A^{\vee}$. By local Tate duality, the following modules have the same length,  
\[\frac{H^1(\mrm{Gal}(K^{nr}/K),A^{\vee}[q^r](K^{nr}))}{H^1_{\kappa}(K,A^{\vee}[q^r])\cap H^1(\mrm{Gal}(K^{nr}/K),A^{\vee}[q^r](K^{nr}))},\,\, \frac{H^1(\mrm{Gal}(K^{nr}/K),A[q^r](K^{nr}))+H^1_{\kappa}(K,A[n])}{H^1(\mrm{Gal}(K^{nr}/K),A[q^r](K^{nr}))}, \] which proves the claim in arbitrary reduction. 

When $A$ has good reduction, the N\'eron model of $A$ is an Abelian scheme, so $\pi_0$ is the trivial group and $H^1(\mrm{Gal}(K^{nr}/K),A[q^r](K^{nr})) \subset H^1_{\kappa}(K,A[q^r])$. By \cite[Lemma 1]{WashFLT} and \cite[Lemma I.3.3]{MilneADT}, one has
\[H^1(\mrm{Gal}(K^{nr}/K),A[q^r](K^{nr}))| = |A[q^r](K)| = |A(K)/q^rA(K)| = |H^1_{\kappa}(K,A[q^r])|,\]
and $H^1(\mrm{Gal}(K^{nr}/K),A[q^r](K^{nr})) = H^1_{\kappa}(K,A[q^r])$.

}

We can now prove Proposition \ref{BK-for-If}. 

\demo{\ref{BKIf1}: By Lemma \ref{everything-self-dual-again}, the Artin representation $\psi_{\rho} := \mrm{Ind}_K^{\Q}{\psi(\rho_K)}$ is odd. By Proposition \ref{tate-module-xrho-nsplit-cartan}, $\psi_{\rho}$ is irreducible, so it is attached to a newform $g \in \mathcal{S}_1(\Gamma_1(D))$. \\
\ref{BKIf2}: The definition of $\psi'$ shows that the representation $\mrm{Ind}_K^{\Q}{\psi'(\rho_K)}$ is attached to the newform $g' := \sigma(g) \otimes \gamma \in \mathcal{S}_1(D')$. Since $\gamma \in \mathcal{D}$, $D' \mid Dp^2$ by Proposition \ref{twist-exists}.

The compatible system of representations of $G_{\Q}$ attached to $f \otimes g$ (resp. $f' \otimes g'$) is given by $\left(V_{f,\mathfrak{l}} \otimes \mrm{Ind}_K^{\Q}{\psi(\rho_K)}\right)_{\mathfrak{l}}$ (resp. $\sigma\left[\left(V_{f,\mathfrak{l}} \otimes \mrm{Ind}_K^{\Q}{\psi(\rho_K)}\right)_{\mathfrak{l}}\right]$). Moreover, by Lemma \ref{everything-self-dual-again}, the product of the central characters of $f$ and $g$ is $\epsilon_K$, so the imprimitive Rankin-Selberg $L$-series for $f \otimes g$ and $f' \otimes g'$ are, up to finitely many Euler factors (which are inconsequential by the proof of Lemma \ref{multiple-rankin-selberg}), conjugates by $G_{\Q}$. By Lemma \ref{multiple-rankin-selberg}, one has $L^{imp}(f \otimes g,1) \neq 0$, so, by \cite[Theorem 3]{Shimura76}, $L^{imp}(f' \otimes g',s)\neq 0$, and therefore $L(f' \otimes g',s) \neq 0$. \\

\ref{BKIf3}: The goal is to apply the properties of the Beilinson-Flach Euler system constructed in \cite{KLZ15}: we are essentially proving an analogue of Theorem 11.7.4, but the fact that we deal with a more general modular form rather than an elliptic curve makes the proof slightly more technical. \\

\emph{Step 1: Description of the Selmer conditions.}

First, for definiteness, we make the Selmer condition explicit. The algebra $\mathbb{T}' := \mathbb{T}_{1,N}(N)/(I)$ is a subring of the ring of endomorphisms of $A_I$; by Proposition \ref{make-abelian-variety-If}, $\mathbb{T}'$ is a domain and a finite free $\Z$-module. Thus, $\mathbb{T}' \otimes_{\Z} \Z_q$ is a finite reduced $\Z_q$-algebra without $\Z_q$-torsion. Since $\Z_q$ is Henselian, $\mathbb{T}' \otimes_{\Z} \Z_q$ is a finite product of local rings \cite[Lemmas 04GG, 04GM]{Stacks}, and the summands are exactly the $\widehat{\mathbb{T}'_{\mathfrak{q}'}}$, where $\mathfrak{q}'$ runs through ideals of $\mathbb{T}'$ dividing $q$. 

Thus, any $(\mathbb{T}'\otimes_{\Z} \Z_q)$-module $M$ is the direct sum of components that are $\widehat{\mathbb{T}'_{\mathfrak{q}'}}$-modules. We define $\Tate{\mathfrak{q}_0}{A_I}, \Tate{\mathfrak{q}_0}{A_I^{\vee}}$ as the $\widehat{\mathbb{T}'_{\mathfrak{q}_0}}$-components of $\Tate{q}{A_I}, \Tate{q}{A_I^{\vee}}$. If $M$ is in fact a torsion $\Z_q$-module, then the $\widehat{\mathbb{T}'_{\mathfrak{q}'}}$-component of $M$ is exactly $M[(\mathfrak{q}')^{\infty}]$. 

The Bloch-Kato condition at a place $v \neq q$ (resp. $v=q$) for $\Tate{\mathfrak{q}_0}{A_I}$ is the \emph{simple} (see \cite[Definition 11.2.3]{KLZ15}) local condition whose image in $H^1$ is given by the inverse image under the inclusion $H^1(\Q_v,\Tate{\mathfrak{q}_0}{A_I}) \rar H^1(\Q_v,\Tate{\mathfrak{q}_0}{A_I} \otimes_{\Z_q} \Q_q)$ of the unramified classes (resp. the crystalline classes $H^1_f(\Q_v,\Tate{\mathfrak{q}_0}{A_I} \otimes_{\Z_q} \Q_q)$ in the sense of \cite[(3.7.2)]{BKTama} when $v=q$). 

The Bloch-Kato condition at a place $v \neq q$ (resp. $v=q$) for $A_I[\mathfrak{q}_0^{\infty}]$ is the simple local condition whose image in $H^1$ is given by the image of the unramified (resp. crystalline) classes under the morphism $H^1(\Q_v,\Tate{\mathfrak{q}_0}{A_I} \otimes_{\Z_q} \Q_q) \rar H^1(\Q_v,\Tate{\mathfrak{q}_0}{A_I} \otimes \Q_q/\Z_q) \simeq H^1(\Q_v,A_I[\mathfrak{q}_0^{\infty}])$. 

We make the same definitions for the dual abelian variety $A_I^{\vee}$ of $A_I$, which is also endowed with an action of $\mathbb{T}'$. Now, since $A_I^{\vee}[\mathfrak{q}_0^{\infty}]$ is the (Tate-twisted) dual of $\Tate{\mathfrak{q}_0}{A_I}$ by \cite[\S 11, \S 16]{MilAb}, local duality \cite[Proposition 3.8]{BKTama}\footnote{The authors of \emph{loc.cit.} seem to omit this definition even though they are using it. This definition is from \cite[Remark I.3.6]{Rubin-ES}.} tells us that $H^1_f(\Q_v,A_I^{\vee}[\mathfrak{q}_0^{\infty}])$ and $H^1_f(\Q_v,\Tate{\mathfrak{q}_0}{A_I})$ are the orthogonal complements of each other. 

Note that $H^1_f(\Q_q,A_I[\mathfrak{q}_0^{\infty}])$ is exactly the kernel of $H^1(\Q_q,A_I[\mathfrak{q}_0^{\infty}]) \rar H^1(\Q_q,A_I(\overline{\Q_v}))$. Indeed, it is enough to prove the same result after replacing $\mathfrak{q}_0^{\infty}$ with $q^{\infty}$, and this is exactly \cite[Example 3.11]{BKTama}. 

After removing finitely many primes $q$, we may and do assume that $\mathbb{T}' \otimes \Z_{(q)}$ is a \emph{regular} semi-local ring, so the same holds for $\mathbb{T}' \otimes \Z_q$. This implies by Proposition \ref{make-abelian-variety-If} that $\Tate{q}{A_I}$ is a locally free $(\mathbb{T}' \otimes \Z_q)$-module of rank four, so by \cite[Lemma 02M9]{Stacks} it is free of rank four. \\

\emph{Step 2: The Beilinson-Flach Euler system applies.}

With the notation of \cite[\S 11.2]{KLZ15}, our first step is to prove that the Nekov\'a\v{r} cohomology group $\tilde{H}^2(\Z[(6pDq)^{-1}],\Tate{\mathfrak{q}_0}{A_I}; \Delta^{BK})$ is finite, where $\Delta^{BK}$ is the Bloch-Kato Selmer condition. 

It is enough to show that $\tilde{H}^2(\Z[(6pDq)^{-1}],\Tate{\mathfrak{q}_0}{A_I} \otimes_{\widehat{\mathbb{T}'_{\mathfrak{q}_0}}} \OO_{F_{\mathfrak{q}}};\Delta^{BK})$ is finite. By Proposition \ref{make-abelian-variety-If}, there is an injective $\OO_{F_{\mathfrak{q}}}[G_{\Q}]$-homomorphism $\Tate{\mathfrak{q}_0}{A_I} \otimes_{\widehat{\mathbb{T}'_{\mathfrak{q}_0}}} \OO_{F_{\mathfrak{q}}} \rar T_{f,\mathfrak{q}} \otimes_{\OO_{F_{\mathfrak{q}}}} T_{g,\mathfrak{q}}$ with finite cokernel.

Therefore, by the long exact sequence of Nekov\'a\v{r} cohomology, $\tilde{H}^2(\Z[(6pDq)^{-1}],\Tate{\mathfrak{q}_0}{A_I};\Delta^{BK})$ is finite if, and only if, $\tilde{H}^2(\Z[(6pDq)^{-1}],T_{f,\mathfrak{q}} \otimes_{\OO_{F_{\mathfrak{q}}}} T_{g,\mathfrak{q}};\Delta^{BK})$ is finite.

We are done if we can apply \cite[Theorem 11.7.3]{KLZ15} \footnote{It seems that the statement contains a slight typo: the statement should be about $M_{\OO}(f \otimes g)^{\ast}(\tau^{-1})$ rather than $M_{\OO}(f \otimes g)(\tau^{-1})$}, but there are several issues:
\begin{enumerate}[noitemsep,label=(\alph*)]
\item\label{1111} We need to check the assumptions of \cite[Hypothesis 11.1.1]{KLZ15}. 
\item\label{BI} We need to check that Hyp(BI) \cite[Hypothesis 11.1.2]{KLZ15} is satisfied. 
\item\label{NEZ} We need to check Hyp(NEZ) \cite[Hypothesis 11.1.4]{KLZ15}.  
\end{enumerate}

For \ref{1111}, $f$ and $g$ are ordinary at $\mathfrak{q}$ and have weights $2$ and $1$ respectively. Moreover, $q \geq 5$ is coprime to the levels of $f$ and $g$ by definition. The newform $f$ is $q$-distinguished by \cite[Remark 7.2.7]{KLZ15}, and, since it does not have complex multiplication, after removing finitely many $q$, it is non-Eisenstein modulo $\mathfrak{q}$ by \cite[Theorem 3.1]{Ribet2}. If we furthermore assume that $q$ does not divide $p-1$ times the cardinality of the image of $\mrm{Ind}_K^{\Q}{\psi(\rho_K)}$ (this only removes a finite number of primes), then $q$ does not divide the cardinality of the image of $G_{\Q}$ in $\mrm{End}(T_{g,\mathfrak{q}})$, and therefore $T_{g,\mathfrak{q}}/\mathfrak{q}$ is irreducible, so that $g$ is not Eisenstein at $\mathfrak{q}$.

Finally, suppose that $g$ is not $q$-distinguished: then the two eigenvalues of $\Fr_q$ acting on $T_{g,\mathfrak{q}}$ agree modulo $\mathfrak{q}$. Let $\zeta$ be the quotient of the two eigenvalues of $\Fr_q$. By definition of $\delta$, $\zeta$ is a $\delta$-th root of unity and we just assumed that $\zeta \equiv 1\pmod{\mathfrak{q}}$. Since $\delta$ and $q$ are coprime, $\zeta=1$ and $\Fr_q$ acts on $T_{g,\mathfrak{q}}$ by a scalar, which contradicts our definition of $\mathcal{Q}$. 

\ref{NEZ} is not an issue by the discussion below \cite[Hypothesis 11.1.4]{KLZ15}. That leaves \ref{BI}. Because $f$ has no complex multiplication, we can apply Proposition \ref{conditions-n-gri} (up to removing finitely many $\mathfrak{q}$): we only need to show that (after removing finitely many $\mathfrak{q}$, if necessary) there is a $\sigma \in G_{\Q}(\mu_{q^{\infty}})$ such that the co-invariants of $T_{f,\mathfrak{q}} \otimes_{\OO_{F_{\mathfrak{q}}}} T_{g,\mathfrak{q}}$ under $\sigma$ are free of rank one over $\OO_{F_{\mathfrak{l}}}$. We wish to apply Proposition \ref{special2}: because the character of $g$ is the quadratic Dirichlet character attached to $K$, in this case, we need to check that $K$ is not contained in the number field $H$ defined in Section \ref{sect-theoA}. By Lemma \ref{gamma1p-onlytwists}, either $(f,\mathbf{1}) \in \mathscr{C}'_1$, or $K=\Q$. In the latter case, if $\alpha$ is a Dirichlet character such that $f \otimes \alpha$ is Galois-conjugate to $f$, then, by Proposition \ref{level-newform-bigtwist}, $\alpha \in \mathcal{D}$ and $\alpha^2=1$. Hence $H \subset \Q(\sqrt{-p})$ does not contain $K$. In Proposition \ref{special2}, since the character of $g$ has square $\mathbf{1}$, the first possibility is not satisfied. Therefore the element $\sigma$ as stated exists.

Thus $\tilde{H}^2(\Z[(6pDq)^{-1}],\Tate{\mathfrak{q}_0}{A_I},\Delta^{BK})$ is finite. 

Therefore, by \cite[Proposition 11.2.9]{KLZ15}, $\mrm{Sel}(\Q,A_I^{\vee}[\mathfrak{q}_0^{\infty}])$ is finite. \\

\emph{Step 3: Finiteness of $\mrm{Sel}(\Q,A_I[\mathfrak{q}_0^{\infty}])$.}

Let $s \geq r \geq 0$, and consider the short exact sequence $0 \rar A_I^{\vee}[q^r] \rar A_I^{\vee}[q^s] \overset{q^r} A_I^{\vee}[q^{s-r}] \rar 0$ of $G_{\Q,qpD}$-modules. It induces a long exact sequence \[0 \rar A_I^{\vee}[q^{s-r}](\Q)/q^rA_I^{\vee}[q^s](\Q) \rar H^1(G_{\Q,qpD},A_I^{\vee}[q^r]) \rar H^1(G_{\Q,qpD},A_I^{\vee}[q^s]).\] By the Mordell-Weil theorem, if we take $r \geq r_q$ and $s \geq s+r_q$ for some $r_q > 0$, it implies that $\ker\left[H^1(G_{\Q,pqD},A_I^{\vee}[q^r]) \rar H^1(G_{\Q,pqD},A_I^{\vee}[q^{\infty}])\right]$ has bounded length as a $\Z_q$-module.

Therefore, by Lemma \ref{galois-cohom-abelian}, the Selmer group $H^1_{\mathcal{L}}(G_{\Q,pqD}, A_I^{\vee}[q^r][\mathfrak{q}_0^{\infty}])$ has bounded length as a $\Z_q$-module as $r \rar \infty$, where $\mathcal{L}$ denotes the following collection of local conditions: 
\begin{itemize}[noitemsep,label=$-$]
\item at places $v \neq q$, $\mathcal{L}_v$ is the set of unramified cohomology classes.  
\item at $v=q$, $\mathcal{L}_q$ is the projection to the $\mathfrak{q}_0$-component of the image of the Kummer map $A_I^{\vee}(\Q_q)/q^rA_I^{\vee}(\Q_q) \rar H^1(\Q_q,A_I^{\vee}[q^r])$. That this local condition yields the result is a consequence of \cite[Example 3.11]{BKTama}. 
\end{itemize}

We now apply the Greenberg-Wiles formula \cite[Theorem 8.7.9]{CNB} to the Galois module $A_I^{\vee}[q^r][\mathfrak{q}_0^{\infty}]$. By the Mordell-Weil theorem, \cite[Lemma 1]{WashFLT} and Lemma \ref{galois-cohom-abelian}, it follows that for large enough $r$,
\[\frac{|H^1_{\mathcal{L}}(\Q,A_I^{\vee}[q^r][\mathfrak{q}^{\infty}])|}{|H^1_{\mathcal{L}^{\vee}}(\Q,A_I[q^r][\mathfrak{q}^{\infty}])|} = \frac{|A_I^{\vee}(\Q)[\mathfrak{q}_0^{\infty}]||\left(A_I^{\vee}(\Q_q)/q^rA_I^{\vee}(\Q_q)\right)[\mathfrak{q}_0^{\infty}]|}{|A_I(\Q)[\mathfrak{q}_0^{\infty}]| |A_I^{\vee}(\R)[q^r][\mathfrak{q}_0^{\infty}]|},\]
where $\mathcal{L}^{\vee}$ is defined exactly in the same way as $\mathcal{L}$, but after replacing $A_I^{\vee}$ with $A_I$.  

We claim that
\[\ln{|\left(A_I^{\vee}(\Q_q)/q^rA_I^{\vee}(\Q_q)\right)[\mathfrak{q}_0^{\infty}]|} = 2er\ln{q} + O(1),\, \ln{|A_I^{\vee}(\R)[q^r][\mathfrak{q}_0^{\infty}]|}=2er\ln{q}+O(1),\]
where $e$ is the dimension of $\mathbb{T}'/\mathfrak{q}_0$ as a $\F_q$-vector space. 

By the theory of the $q$-adic exponential, there is a $\mathbb{T}' \otimes \Z_q$-submodule $M \subset A_I^{\vee}(\Q_q)$ such that $A_I^{\vee}(\Q_q)/M$ is finite and $M \otimes_{\Z_q} \Q_q \simeq T_0 (A_I^{\vee})_{\Q_q}$. Thus $M$ is free of rank two over $\mathbb{T}' \otimes \Z_q$, which implies that $\ln{|\left(A_I^{\vee}(\Q_q)/q^rA_I^{\vee}(\Q_q)\right)[\mathfrak{q}_0^{\infty}]|} = O(1)+\ln{|(M/q^r M)[\mathfrak{q}_0^{\infty}]|} = 2re + O(1)$. 

Similarly, by the theory of real Lie groups, $A_I^{\vee}(\R)$ is a compact Abelian Lie group, so there is an exact sequence $0 \rar \Lambda \overset{\exp}{\rar} A_I^{\vee}(\R) \rar C \rar 0$ of $\mathbb{T}'$-modules with $C$ finite and $\Lambda \otimes \Z_q \simeq \Tate{q}{A_I^{\vee}(\R)}$ a free $(\mathbb{T}'\otimes \Z_q)$-module of rank two. Thus 
\[\ln{|A_I^{\vee}(\R)[q^r][\mathfrak{q}_0^{\infty}]|} = O(1)+\ln{|(M/q^r M)[\mathfrak{q}_0^{\infty}]|} = 2re\ln{q}+O(1).\]

Therefore, $H^1_{\mathcal{L}^{\vee}}(\Q,A_I[q^r][\mathfrak{q}_0^{\infty}])$ has bounded length as $r \rar \infty$. By Lemma \ref{galois-cohom-abelian}, we can replace the system of local conditions $\mathcal{L}^{\vee}$ at the finitely many places $v \neq q$ of bad reduction by the local condition ``lies in the image of $A_I(\Q_v)$ under the local Kummer map''. Hence $\mrm{Sel}(\Q,A_I[q^r])[\mathfrak{q}_0^{\infty}]$ has bounded length as $r \rar \infty$, thus $\mrm{Sel}(\Q,A_I[q^{\infty}])[\mathfrak{q}_0^{\infty}]$ is finite. \smallskip

Recall that the Kummer exact sequence yields an injection $A_I(\Q) \otimes \Q_q/\Z_q \rar \mrm{Sel}(\Q,A_I[q^{\infty}])$, so that $\left(A_I(\Q) \otimes \Q_q/\Z_q\right)[\mathfrak{q}_0^{\infty}]$ is finite. Since this group is divisible, it is zero. Therefore the finitely generated $\widehat{\mathbb{T}'_{\mathfrak{q}_0}}$-module $A_I(\Q) \otimes_{\mathbb{T}'} \widehat{\mathbb{T}'_{\mathfrak{q}_0}}$ is $q$-divisible, so it is zero. Since $\widehat{\mathbb{T}'_{\mathfrak{q}_0}}$ is flat over $\mathbb{T}'$, there is no injective homomorphism $\mathbb{T}' \rar A_I(\Q)$. Thus, by the Mordell-Weil theorem, $\mrm{Ann}_{\mathbb{T}'}(A_I(\Q))$ is a nonzero ideal. Since $\mathbb{T}'$ is a domain, every proper quotient of $\mathbb{T}'$ is finite, hence $A_I(\Q)$ is finite. \\

\emph{Step 4: Density estimate.}

Since we considered all the primes in $\mathcal{Q}$ except for possibly finitely many, by Cebotarev's theorem, we need to prove that the set of residue characteristics of prime ideals at which $f$ is ordinary has density one. Now, let $q \neq p$ be a prime number which is unramified in $F$ and such that $f$ is not ordinary at any prime ideal $\mathfrak{q} \subset \OO_F$ with residue characteristic $q$. Then, for all $\mathfrak{q} \subset \OO_F$ with residue characteristic $q$, $a_q(f) \in \mathfrak{q}\OO_F$. Then $a_q(f) \in \cap_{\mathfrak{q} \mid q}{\mathfrak{q}\OO_F} = q\OO_F$. Every complex embedding of $a_q(f) \in \OO_F$ has absolute value at most $2\sqrt{q}$, which implies that $a_q(f)=0$. By Serre's result \cite[Cor. 2 au Th\'eor\`eme 15]{Serre-Cebotarev}, the set of prime numbers $q$ such that $a_q(f)=0$ has density zero, because $f$ has no complex multiplication. 
}

\medskip

\rem{Even if we do not assume any longer that $L(f \times g,1) \neq 0$, the results of Chapter \ref{obstructions-euler} are sufficient to prove that Hyp(BI) \cite[Hypothesis 11.1.2]{KLZ15} holds. Using \cite[Theorem 11.6.6]{KLZ15}, a similar argument to the proof of Proposition \ref{BK-for-If} shows that the rank of $A_I(\Q)$ over $\mathbb{T}'$ is bounded above by the vanishing order of the $q$-adic Rankin-Selberg $L$-function. 

Suppose now that $\rho$ comes from an elliptic curve with good reduction at $p$. One could make the following heuristic reasoning:
\begin{enumerate}[noitemsep,label=(\alph*)]
\item\label{heuristic-order} For ``many'' four-dimensional factors present in the decomposition of $\Tate{\ell}{J_{\rho,\xi}(p)} \otimes_{\Z_{\ell}} F_{\lambda}$, their complex $L$-function will vanish with order exactly one (since the sign of its functional equation is $-1$). 
\item\label{heuristic-cp} Hence, for ``many'' of these factors, the corresponding $q$-adic Rankin-Selberg $L$-function will vanish at the central value with order one. 
\item\label{heuristic-euler} Thus, a ``large'' quotient $A$ of $J_{\rho,\xi}(p)$ will satisfy that $\mrm{rk}\,A(\Q) \leq \frac{\dim{A}}{2}$.
\item\label{heuristic-chabauty} And therefore, one may apply the Chabauty-Coleman method (see \cite{PoonenMcCallum} for an account, or the following sections for an application in rank zero) to $X_{\rho,\xi}(p)$ for the quotient $A$, thus determining all its rational points.
\end{enumerate}

All of these steps, apart from \ref{heuristic-euler} which we just discussed, seem to be difficult. For example, in positive rank, one cannot rule out as easily with the Chabauty-Coleman method the existence of a rational point in the residue disk of a cusp, in the way that this works in rank zero (see the following sections).

}

\section{Reminder on the Chabauty-Coleman method}

The goal of this section is to state and prove the result that we will need for the next section. It is worth noting that this construction is somewhat \emph{ad hoc}. The so-called ``tiny'' integrals were generalized -- and, in particular, given meaning across residue disks -- by Coleman \cite{ColAb} in the abelian setting.  

\defi{Let $K$ be a non-Archimedean local field with residue field $k$ and $X$ be a proper $\OO_K$-scheme. By the valuative criterion for properness \cite[Lemma 0XB5]{Stacks}, for any finite extension $L/K$, the natural morphism $\mrm{Mor}_{\OO_K}(\Sp{\OO_L},X) \rar \mrm{Mor}_K(\Sp{L},X)$ is a bijection. 

Let $L/K$ be a finite extension and $\varpi \in \OO_L$ be a uniformizer. Given a point $x \in X(L)$, $x$ is the restriction of a point $x_0 \in X(\OO_L)$. The \emph{reduction} of $x$ is the point $\mrm{red}(x) \in X_k$ of the special fibre of $X$ which is the image under $x_0$ of the unique closed point of $\Sp{\OO_L}$. The point $\mrm{red}(x)$ is also endowed with a morphism $\Sp{\OO_L/\varpi} \rar \Sp{\kappa(\mrm{red}(x))}$; we will explicitly mention it when and how we need this additional data. 

Given a closed point $x \in X_k$, its \emph{residue disk} $]x\mkern1mu[$ is the set of closed points $y \in X_K$ (in particular, $y \in X_K(\kappa(y))$, with $\kappa(y)/K$ a finite extension) such that $\mrm{red}(y)=x$. If $L/K$ is a finite extension, $\varpi \in \OO_L$ is a uniformizer, and $x$ is endowed with a morphism $\Sp{\OO_L/\varpi} \rar \Sp{\kappa(x)}$ of $\OO_K$-schemes, we denote by $]x\mkern1mu[\,(L)$ the set of $L$-points of $X$ whose reduction (as $\Sp{\OO_L/\varpi}$-points, as defined in the above paragraph) is $x$.
}

\lem[syst-params]{Let $K$ be a finite extension of $\Q_p$ for some prime number $p$, and $\varpi \in \OO_K$ be a uniformizer: write $k=\OO_K/\varpi$. Let $X$ be a $\OO_K$-scheme of finite type and $x \in X_k$ be a closed point such that $X$ is smooth of relative dimension $d$ at $x$. Then $\OO_{X,x}$ is a regular ring of dimension $d+1$, and there exists $t_1,\ldots,t_d \in \OO_{X,x}$ such that $(\varpi,t_1,\ldots,t_d)$ generate its maximal ideal. Moreover, 
\begin{itemize}[label=$-$,noitemsep]
\item Let $K'$ be the finite unramified extension of $K$ with residue field $\kappa(x)$. There is a unique homomorphism of local $\OO_K$-algebras $\OO_{K'} \rar \widehat{\OO_{X,x}}$ inducing the identity on the residue fields, and it induces an isomorphism $\OO_{K'}[[T_1,\ldots,T_d]] \rar \widehat{\OO_{X,x}}$ (where $T_i$ is mapped to $t_i$). 
\item Let $L/K$ be a finite extension and $\varpi_L \in\OO_L$ be a uniformizer. Fix a homomorphism $\psi: \Sp{\OO_L/\varpi_L} \rar \Sp{\kappa(x)}$ of $\OO_K$-schemes: this defines a unique homomorphism of $K$-algebras $K' \rar L$. Moreover, let $H$ be the set of local homomorphisms of local $\OO_K$-algebras $\OO_{X,x} \rar \OO_L$ preserving $\psi$, then $f \in H \longmapsto (f(t_i))_{1 \leq i \leq d} \in (\varpi_L\OO_L)^d$ is a bijection.  
\end{itemize}
}

\demo{The statement is local: after replacing $X$ with an open subscheme containing $x$, we may assume that $X$ is a smooth affine $\OO_K$-scheme of relative dimension $d$. This implies by \cite[Lemma 07NF]{Stacks} that $\OO(X)$ is a regular ring, hence $\OO_{X,x}$ is a regular local ring. Moreover, by base change, if $k$ denotes the residue field of $\OO_K$, $X_k \rar \Sp{k}$ is smooth of relative dimension $d$, so $\OO(X) \otimes_{\OO_K} k$ is a regular ring, hence $\OO_{X,x}/(\varpi)$ is regular. Moreover, by \cite[Lemma 0AFF]{Stacks}, there exists a connected open subscheme $x \in U \subset X_k$ of dimension $d$; since $U$ is regular, it is an integral scheme. By \cite[Proposition 2.5.19]{QL}, the function field of $U$ has transcendence degree $d$. By the dimension formula \cite[Lemma 02JU, Lemma 02JB]{Stacks}, $\dim{\OO_{X,x}/(\varpi)}=d$, whence $\dim{\OO_{X,x}}=d+1$. 

As a consequence, $\OO_{X,x}$ is a regular local ring of dimension $d+1$, and by \cite[Corollary 4.2.15]{QL}, there are $t_1,\ldots,t_d \in \mathfrak{m}_{X,x}$ such that $(\varpi,t_1,\ldots,t_d)$ is a basis of $\mathfrak{m}_{X,x}/\mathfrak{m}_{X,x}^2$. By Nakayama, $(\varpi,t_1,\ldots,t_d)$ generate $\mathfrak{m}_{X,x}$. 

The residue field $\kappa(x)$ at $x$ is the residue field of the closed point $x \in X_k$, where $X_k$ is a finite type $k$-scheme, so that $\kappa(x)$ is a finite extension of $k$, whence $\kappa(x)$ is a finite field. Let $K'/K$ be the finite unramified extension such that $\kappa(x)$ identifies with the residue field of $K'$. Then $\OO_{K'}$ is an \'etale $\OO_K$-algebra, so it is formally \'etale by \cite[Lemma 04FF]{Stacks}. Therefore, for any $r \geq 1$, any homomorphism of $\OO_K$-algebras $\OO_{K'} \rar \OO_{X,x}/\mathfrak{m}_{X,x}^r$ lifts uniquely to a homomorphism of $\OO_K$-algebras $\OO_{K'} \rar \OO_{X,x}/\mathfrak{m}^{r+1}$. Hence there is a canonical homomorphism of $\OO_K$-algebras $\OO_{K'} \rar \widehat{\OO_{X,x}}$ which induces an isomorphism on the residue fields.

There is a natural homomorphism $f: \OO_{K'}[T_1,\ldots,T_d] \rar \widehat{\OO_{X,x}}$ mapping every $T_i$ to $t_i$. It maps any element outside $(\varpi,T_1,\ldots,T_d)$ to an element of $\OO_{K'}^{\times}+\mathfrak{m}_{X,x}\widehat{\OO_{X,x}} \subset \widehat{\OO_{X,x}}^{\times}$, so $f$ induces a morphism $f': \OO_{K'}[T_1,\ldots,T_d]_{(\varpi,T_1,\ldots,T_d)} \rar \widehat{\OO_{X,x}}$. Since $\mathfrak{m}_{X,x}$ is generated by $(\varpi,t_1,\ldots,t_d)$, $f'$ induces a surjection $\OO_{K'}[T_1,\ldots,T_d]/(\varpi,T_1,\ldots,T_d)^n \rar \OO_{X,x}/\mathfrak{m}_{X,x}^n$ for any $n \geq 1$. By \cite[Lemma 00NO]{Stacks}, this surjection is between $\Z_p$-modules with the same finite length, so it is an isomorphism. Therefore $f'$ is a completion. 

Let $L/K$ be a finite extension, $\varpi_L \in \OO_L$ a uniformizer, and endow $x$ with a morphism $\Sp{\OO_L/\varpi_L} \overset{\psi}{\rar} \Sp{\kappa(x)}$. The $\OO_K$-algebra $\OO_{K'}$ is formally \'etale, so the homomorphism $\kappa(x) \rar \OO_L/\varpi_L$ lifts uniquely to a $\OO_K$-homomorphism $\OO_{K'} \rar \OO_L$. By definition, $H$ is the set of commutative diagrams 
\[
\begin{tikzcd}[ampersand replacement=\&]
\OO_{X,x} \arrow{r}{y^{\sharp}}\arrow{d}\& \OO_L \arrow{d}\\
\kappa(x) \arrow{r}{\psi^{\sharp}} \& \OO_L/\varpi_L
\end{tikzcd}
\]  
where $y^{\sharp}$ is a local homomorphism and the two vertical arrows are the reductions in the quotient field. Thus $f \in H\longmapsto \widehat{f} \in \mrm{Hom}_{\OO_{K'}}(\widehat{\OO_{X,x}},\OO_L)$ is a bijection. Since $\widehat{\OO_{X,x}}$ is a ring of formal power series in the variables $t_1,\ldots,t_d$, $f \in \mrm{Hom}_{\OO_{K'}}(\widehat{\OO_{X,x}},\OO_L) \longmapsto (f(t_i))_{1 \leq  i \leq d} \in (\varpi_L\OO_L)^d$ is a bijection, whence the conclusion. 
} 

\medskip

\defi{In the situation of Lemma \ref{syst-params}, any $(t_1,\ldots,t_d)$ such that the maximal ideal of $\OO_{X,x}$ is generated by $(\varpi,t_1,\ldots,t_d)$ is called \emph{a system of parameters} for $X$ at $x$.}

\cor[smooth-residue-disks]{Let $K$ be a finite extension of $\Q_p$ with uniformizer $\varpi$ and $X$ be a proper $\OO_K$-scheme. Let $x\in X$ be a closed point (in the special fibre of $X$) such that $X \rar \Sp{\OO_K}$ is smooth at $x$. Let $(t_1,\ldots,t_d)$ be a system of parameters at $x$, and fix a finite extension $L/K$ with uniformizer $\varpi_L$ and a homomorphism $\Sp{\OO_L/\varpi_L} \rar \Sp{\kappa(x)}$. Then \[y \in \,\,]x\mkern1mu[\,(L) \longmapsto (y^{\sharp}(t_i))_{1 \leq i\leq d} \in (\varpi_L\OO_L)^d\] is a bijection, where $y^{\sharp}$ is the ring homomorphism $\OO_{X,x} \rar \OO_L$ given by the valuative criterion for properness and the definition of $]x\mkern1mu[\,(L)$. }

\demo{This is an application of Lemma \ref{syst-params}; the fact that $y\in\,\, ]x\mkern1mu[\,(L) \longmapsto y^{\sharp} \in H$ (with notation as in this Lemma) is a bijection is exactly the valuative criterion for properness. }

\medskip

\lem[params-differentials]{Let $K$ be a finite extension of $\Q_p$ for some prime $p$. Let $X$ be a $\OO_K$-scheme of finite type and $x \in X$ be a closed point in the special fibre at which $X$ is smooth of relative dimension $e$. Let $(t_1,\ldots,t_e)$ be a system of parameters at $x$. Then $\Omega^1_{X/\OO_K,x}$ is a free $\OO_{X,x}$-module of rank $e$ with basis $dt_1,\ldots,dt_e$. 
Moreover, there is a unique continuous derivation $\hat{d}: \widehat{\OO_{X,x}} \rar \Omega^1_{X/\OO_K,x} \otimes_{\OO_{X,x}} \widehat{\OO_{X,x}}$ extending $d: \OO_{X,x} \rar \Omega^1_{X/\OO_K,x}$.

 If $\varpi \in \OO_K$ is a uniformizer, for any $z_1,\ldots,z_e \in \widehat{\OO_{X,x}}$ such that $(\varpi,z_1,\ldots,z_e)$ generate the maximal ideal of $\widehat{\OO_{X,x}}$, the $\hat{d}z_1,\dots,\hat{d}z_e$ form a basis of $\Omega^1_{X/\OO_K,x} \otimes_{\OO_{X,x}} \widehat{\OO_{X,x}}$. }

\demo{There is a natural local homomorphism $\pi: \OO_K[T_1,\ldots,T_e]_{(\varpi,T_1,\ldots,T_e)} \rar \OO_{X,x}$; by construction, the maximal ideal of $\OO_{X,x}$ is generated by the image under $\pi$ of the maximal ideal of $\OO_K[T_1,\ldots,T_e]_{(\varpi,T_1,\ldots,T_e)}$. By Lemma \ref{syst-params}, the completion $\OO_K[T_1,\ldots,T_d]_{(\varpi,T_1,\ldots,T_d)} \rar \widehat{\OO_{X,x}}$ is finite free hence flat: since $\widehat{\OO_{X,x}}$ is faithfully flat over $\OO_{X,x}$ \cite[Lemma 00MC]{Stacks}, $\pi$ is flat. The residue fields of the source and range of $\pi$ are both finite, so the field extension induced by $\pi$ on the residue fields is finite separable. 

Let $U \subset X$ be a affine open subscheme containing $x$ such that every $t_i$ is contained in $\OO_X(U)$: there is a morphism $f: U \rar \mathbb{A}^e_{\OO_K}$ of $\OO_K$-schemes attached to $T_i \longmapsto t_i$. The morphism $f$ is between affine $\OO_K$-schemes of finite type, so it is of finite type, and it is \'etale at $x$ by \cite[Lemma 02GU (6)]{Stacks}: therefore, after reducing $U$ if necessary, we may assume that $f$ is \'etale. By \cite[Proposition 2.5]{BLR}, the map $\Omega^1_{\mathbb{A}^e_{\OO_K}/\OO_K,0_{\OO_K/\varpi}} \otimes_{\OO_K[T_1,\ldots,T_e]_{(\varpi,T_1,\ldots,T_e)}} \OO_{X,x} \rar \Omega^1_{X/\OO_K,x}$ is an isomorphism: hence $\Omega^1_{X/\OO_K,x}$ is free of rank $e$ over $\OO_{X,x}$ with basis $dt_1,\ldots,dt_e$.  

By the product rule, $d(\mathfrak{m}_{X,x}^r\OO_{X,x}) \subset \mathfrak{m}_{X,x}^{r-1}\Omega^1_{X/\OO_K,x}$. So $d$ extends to a continuous derivation $\hat{d}: \widehat{\OO_{X,x}}\rar \underset{\leftarrow}{\lim}\,\Omega^1_{X/\OO_K,x}/\mathfrak{m}_{X,x}^r\Omega^1_{X/\OO_K,x} \simeq \Omega^1_{X/\OO_K,x} \otimes_{\OO_{X,x}} \widehat{\OO_{X,x}}$. 

Now, let $z_1,\ldots,z_e \in \widehat{\OO_{X,x}}$ be such that $(\varpi,z_1,\ldots,z_e)$ generates the maximal ideal of $\widehat{\OO_{X,x}}$. To prove that $(\hat{d}z_1,\ldots,\hat{d}z_e)$ is a basis of the free $\widehat{\OO_{X,x}}$-module (of rank $e$) $\Omega^1_{X/\OO_K,x} \otimes_{\OO_{X,x}} \widehat{\OO_{X,x}}$, it is enough to prove that they form a basis of $\Omega^1_{X/\OO_K,x}/\mathfrak{m}_{X,x}$ as a $\kappa(x)$-vector space of dimension $e$. Indeed, if $z'_i \in \OO_{X,x}$ has the same image in $\OO_{X,x}/\mathfrak{m}_{X,x}^2$, then $dz'_i$ and $\hat{d}z_i$ have the same image in $\Omega^1_{X/\OO_K,x}/\mathfrak{m}_{X,x}$. Since $(\varpi,z_1,\ldots,z_e)$ is a basis of $\mathfrak{m}_{X,x}\widehat{\OO_{X,x}}/\mathfrak{m}_{X,x}^2\widehat{\OO_{X,x}}$, so is $(\varpi,z'_1,\ldots,z'_e)$; thus $(\varpi,z'_1,\ldots,z'_e)$ is a basis of $\mathfrak{m}_{X,x}/\mathfrak{m}_{X,x}^2$, whence, by Nakayama, the $z'_i$ form a system of parameters, and the $dz_i$ form a basis of $\Omega^1_{X/\OO_K,x}/\mathfrak{m}_{X,x}$, whence the conclusion.
}

\smallskip

\cor[formal-change-of-variables]{With the same notations as in Lemma \ref{params-differentials}, let $(t_1,\ldots,t_e)$ be a system of parameters at $x$ and $K'/K$ be the finite unramified extension such that the residue field of $K'$ is $\kappa(x)$. Let $P \in \OO_{K'}[[T_1,\ldots,T_e]]$, then $\hat{d}(P(t_1,\ldots,t_e))  = \sum_{i=1}^e{\frac{\partial P}{\partial T_i}(t_1,\ldots,t_e)dt_i}$. }

\demo{When $P$ is a scalar $\lambda \in \OO_{K'}$ which is a root of unity, there is an irreducible polynomial, then $\lambda^{p^u-1}=1$ for some $u \geq 1$, so $(p^u-1)\hat{d}\lambda=0$, so $\hat{d}\lambda=0$. The same holds when $P$ is a scalar $\lambda \in \OO_K$. Therefore the result is true when $P$ is constant. Moreover, the result is true when $P(T_1,\ldots,T_d)=T_i$ for some $i \geq 0$. Therefore, by the product rule, the result is correct whenever $P$ is a polynomial. By continuity, the result is correct. }

\medskip

\defi[ideal-first-coeffs]{In the situation of Lemma \ref{params-differentials}, assume that $e=1$, and let $t$ be a parameter of $X$ at $x$ (i.e. $\mathfrak{m}_{X,x}=(\varpi,t)$). Let $K'/K$ be the finite unramified extension such that the residue field of $K'$ identifies with $\kappa(x)$. Let $M \subset \Omega^1_{X/\OO_K,x}$ be a sub-$\OO_K$-module. Then there exists a sub-$\OO_K$-module $J \subset \OO_{K'}[[T]]$ such that the image of $M$ in $\Omega^1_{X/\OO_K,x} \otimes_{\OO_{X,x}} \widehat{\OO_{X,x}}$ is $\{P(t)dt,\, P \in J\}$. The \emph{ideal of first coefficients} $I \subset \OO_{K}$ of $M$ is the inverse image in $\OO_K$ of the ideal of $\OO_{K'}$ generated by the $P(0)$ for $P \in J$.    
}

\medskip

\lem[ideal-first-coeffs-well-defined]{In the situation of Definition \ref{ideal-first-coeffs}, the ideal of first coefficients does not depend on the choice of $t$. Moreover, if there is a finite extension $L/K$ such that $X \rar \Sp{\OO_K}$ factors through a morphism $f: X \rar \Sp{\OO_L}$ which is smooth of relative dimension one at $x$, then $L/K$ is unramified, so that $\Omega^1_{X/\OO_K,x} \rar \Omega^1_{X/\OO_L,x}$ is an isomorphism, and the ideal of first coefficients of $M$ for $X \rar \Sp{\OO_L}$ is generated by the ideal of first coefficients of $M$ for $X \rar \Sp{\OO_K}$.  }

\demo{If $t' \in \OO_{X,x}$ is another element such that $(\varpi,t')$ generate $\mathfrak{m}_{X,x}$, then there is a $P \in T\OO_{K'}[[T]]$ such that $P'(0) \in \OO_{K'}^{\times}$ and the identity $t'=P(t)$ holds in $\widehat{\OO_{X,x}}$. Let $J_t$ (resp. $J_{t'}$) be the sub-$\OO_K$-module of $\OO_{K'}[[T]]$ attached to $M$ and the parameter $t$ (resp. $t'$). For any $Q \in \OO_K[[T]]$, one has, in $\Omega^1_{X/\OO_K,x} \otimes_{\OO_{X,x}} \widehat{\OO_{X,x}}$, 
\[Q(t')dt'=Q(P(t))\hat{d}P(t) = (Q \circ P)(t) P'(t)dt = P'(t)(Q \circ P)(t)dt.\]
Thus $J_t = \{P'(T)Q(P(T))\mid Q \in J_{t'}\}$. The conclusion follows since $P(0)=0$ and $P'(0) \in \OO_{K'}^{\times}$.

Suppose that $X \rar \Sp{\OO_K}$ factors through a morphism $f: X \rar \Sp{\OO_L}$ for some finite extension $L/K$ such that $f$ is smooth at $x$. Then after replacing $X$ with an open subscheme, we may assume that $X \rar \Sp{\OO_K}$ and $X \rar \Sp{\OO_L}$ are smooth of relative dimension one. Since $x \in X$ lies in the special fibre for $X \rar \Sp{\OO_K}$, $p \notin \OO_{X,x}^{\times}$, so $f(x)$ is the closed point $s \in \Sp{\OO_L}$, thus $f$ is smooth and surjective. By \cite[Lemma 05B5]{Stacks}, $\Sp{\OO_L} \rar \Sp{\OO_K}$ is smooth. Since this morphism is \'etale at the generic point, it is smooth of relative dimension zero i.e. \'etale, and $L/K$ is unramified. 

We thus have an exact sequence $\Omega^1_{\OO_L/\OO_K,s} \otimes_{\OO_L} \OO_{X,x} \rar \Omega^1_{X/\OO_K,x} \rar \Omega^1_{X/\OO_L,x} \rar 0$; since $\OO_L$ is \'etale over $\OO_K$, the natural map $\Omega^1_{X/\OO_L,x} \rar \Omega^1_{X/\OO_K,x}$ is thus an isomorphism. In particular, $J_K$ and $J_L$ (with the same notation as in the beginning of the proof) are the same, and the ideal of first coefficients over $K$ (resp. over $L$) is generated by $\varpi^n$, where $n$ is the minimum of the valuations of the $P(0)$ for $P \in J_K$ (resp. $P \in J_L=J_K$), whence the conclusion. }

\medskip

\lem[ideal-first-coeffs-cover]{Let $f: X' \rar X$ be a morphism of $\OO_K$-schemes, and let $x' \in X'$ be a closed point of the closed fibre such that $X'$ and $X$ are smooth of relative dimension one over $\OO_K$ at $x'$ and $x=f(x)$ respectively. Assume that $f$ is smooth at $x'$. Let $M \subset \Omega^1_{X/\OO_K,x}$ be a sub-$\OO_K$-module. Let $f^{\ast}M \subset \Omega^1_{X'/\OO_K,x'}$ be the sub-$\OO_K$-module made with the $(f^{\ast}\omega)_{x'}$. Then $M$ and $f^{\ast}M$ have the same ideal of first coefficients. }

\demo{Let $t \in \OO_{X,x}$ be a parameter, it is enough to show that $f^{\sharp}(t) \in \OO_{X',x'}$ is a parameter. Now, $f$ is smooth at $x'$, so $f$ is smooth at $x'$ of relative dimension zero, so $f$ is \'etale at $x'$. Hence $\OO_{X,x} \rar \OO_{X',x'}$ is unramified, whence the conclusion.}

\medskip

Now, we turn to integration.

\lem[formal-poincare-lemma]{Let us keep the notations of Lemma \ref{params-differentials}, and let $\sum_{i=1}^e{P_i(t_1,\ldots,t_e)dt_i} \in \im{\hat{d}}$. For any $1 \leq i,j \leq e$, one has $\frac{\partial P_i}{\partial T_j}=\frac{\partial P_j}{\partial T_i}$. 

Furthermore, if $P_1,\ldots,P_e \in \OO_K[[T_1,\ldots,T_e]]$ are such that for any $1 \leq i, j \leq e$, $\frac{\partial P_i}{\partial T_j} = \frac{\partial P_j}{\partial T_i}$, then there exists a formal series $Q := \sum_{r \geq 1}{\frac{1}{r}Q_r(T_1,\ldots,T_e)} \in K[[T_1,\ldots,T_e]]$ such that:
\begin{itemize}[noitemsep,label=\tiny$\bullet$]
\item for any $1 \leq i \leq e$, one has $P_i = \frac{\partial Q}{\partial T_i}$, 
\item for every $j \geq 1$, $Q_j \in \OO_K[T_1,\ldots,T_e]$ is a homogeneous polynomial of degree $j$, 
\item for any $R \in K[[T_1,\ldots,T_e]]$ such that $\left(\frac{\partial R}{\partial T_i}\right)_{1 \leq i \leq e} = \left(P_i\right)_{1 \leq i \leq e}$, one has $R=Q+R(0)$.
\end{itemize}  }

\demo{That the image of $\hat{d}$ is contained in the $\OO_{K'}$-submodule of $\sum_{i=1}^e{P_i(t_1,\ldots,t_e)dt_i}$ as stated is a direct application of Corollary \ref{formal-change-of-variables}. 
Now, let $P_1,\ldots,P_e \in \OO_K[[T_1,\ldots,T_e]]$ such that for any $1 \leq i, j \leq e$, $\frac{\partial P_i}{\partial T_j} = \frac{\partial P_j}{\partial T_i}$. For every $r \geq 0$, let $P_i^r$ be the homogeneous component of $P_i$ of degree $r$. Since $\frac{\partial}{\partial T_i}$ maps homogeneous polynomials of degree $r$ to homogeneous polynomials of degree $r-1$, one has, for any $r \geq 0$ and any $1 \leq i, j \leq e$, $\frac{\partial P_i^r}{\partial T_j} = \frac{\partial P_j^r}{\partial T_i}$. Let, for every $r \geq 1$, $Q_r = \sum_{i=1}^e{T_iP_i^{r-1}(T_1,\ldots,T_e)}$: then $Q_r \in \OO_K[[T_1,\ldots,T_e]]$is homogeneous of degree $r$, and one has 
\[\frac{\partial Q_r}{\partial T_i} = P_i^{r-1}+\sum_{j=1}^e{T_j\frac{\partial P_j^{r-1}}{\partial T_i}} = P_i^{r-1}+\sum_{j=1}^e{T_j\frac{\partial P_i^{r-1}}{\partial T_j}} = P_i^{r-1}+(r-1)P_i^{r-1}=rP_i^{r-1}.\]

Thus $\frac{\partial Q}{\partial T_i}=P_i$ for every $1 \leq i \leq e$. Furthermore, if $R \in K[[T_1,\ldots,T_e]]$ is such that $\frac{\partial R}{\partial T_i}=P_i$ for every $1 \leq i \leq e$, then every derivative of the formal power series $R-Q \in K[[T_1,\ldots,T_e]]$ vanishes. Hence $R-Q$ is constant, whence the conclusion. 

}

\medskip

\prop[tiny-integral-well-defined]{Let $K/\Q_p$ be a finite extension and $X$ be a proper $\OO_K$-scheme. Let $x \in X$ be a closed point in the special fibre such that $X$ is smooth at $x$ of relative dimension $d$. Let $(t_1,\ldots,t_d)$ be a system of parameters at $x$. Let $K'/K$ be the finite unramified extension whose residue field is $\kappa(x)$ and consider a finite extension $L/K$ with uniformizer $\varpi_L$ endowed with a field homomorphism $\kappa(x) \rar \OO_L/\varpi_L$. Let $\omega \in \Omega^1_{X/\OO_K,x}$ be a \emph{closed} differential and let $A,B \in \,\,]x\mkern1mu[\,(L)$. 

We define the \emph{tiny integral} $\int_A^B{\omega}$ as follows: in $\Omega^1_{X/\OO_K,x} \otimes_{\OO_{X,x}} \widehat{\OO_{X,x}}$, there is a decomposition $\omega = \sum_{i=1}^e{P_i(t_1,\ldots,t_d)dt_i}$ for some $P_i \in \OO_{K'}[[T_1,\ldots,T_d]]$ such that $\frac{\partial P_i}{\partial T_j}=\frac{\partial P_j}{\partial T_i}$ for every $1 \leq i,j \leq d$. Let $Q \in K'[[T_1,\ldots,T_d]]$ be as in Lemma \ref{formal-poincare-lemma}, then \[\int_A^B{\omega} := Q(B^{\sharp}(t_1),\ldots,B^{\sharp}(t_d))-Q(A^{\sharp}(t_1),\ldots,A^{\sharp}(t_d)).\] 

This quantity is well-defined and does not depend on the choice of $(t_1,\ldots,t_d)$ (or $Q$). In fact, the formula remains correct for any $t_1,\ldots,t_d \in \mathfrak{m}_{X,x}\widehat{\OO_{X,x}}$ for which $\omega=\sum_{i=1}^d{F_i(t_1,\ldots,t_d)\hat{d}t_i}$ for some $F_i \in \OO_{K'}[[T_1,\ldots,T_d]]$ such that $\frac{\partial F_i}{\partial T_j} = \frac{\partial F_j}{\partial T_i}$. }

\demo{The differential $d: \Omega^1_{X/\OO_K,x} \rar \Omega^2_{X/\OO_K,x}$ maps $\mathfrak{m}_{X,x}^r\Omega^1_{X/\OO_K,x}$ into $\mathfrak{m}_{X,x}^{r-1}\Omega^2_{X/\OO_K,x}$, so it extends to a differential $\Omega^1_{X/\OO_K,x} \otimes_{\OO_{X,x}} \widehat{\OO_{X,x}} \rar \Omega^2_{X/\OO_K,x} \otimes_{\OO_{X,x}} \widehat{\OO_{X,x}} \simeq \Lambda^2_{\widehat{\OO_{X,x}}}\left[\Omega^1_{X/\OO_K,x} \otimes_{\OO_{X,x}} \widehat{\OO_{X,x}}\right]$; it also satisfies the product rule.  

We can write $\omega=\sum_{i=1}^d{a_idt_i}$ for all $1 \leq i \leq d$, where $a_i \in \OO_{X,x}$. Therefore, one has $\sum_{i=1}^r{da_i \wedge dt_i}=0$. We can write, in $\widehat{\OO_{X,x}}$, $a_i=P_i(t_1,\ldots,t_d)$ for some $P_i \in \OO_{K'}[[T_1,\ldots,T_d]]$, and the condition $\sum_{i=1}^r{da_i \wedge dt_i}=0$ means that for any $1 \leq i,j \leq d$, one has $\frac{\partial P_i}{\partial T_j} = \frac{\partial P_j}{\partial T_i}$. 

Let $Q=\sum_{r \geq 1}{\frac{1}{r}Q_r(T_1,\ldots,T_d)} \in K'[[T_1,\ldots,T_d]]$ be as in Lemma \ref{formal-poincare-lemma}. Let $L/K'$ be a finite extension with uniformizer $\varpi_L$, let us show that for any $b_1,\ldots,b_d \in \varpi_L\OO_L$, $Q(b_1,\ldots, b_d)$ is well-defined. Indeed, one has $\|Q_r(b_1,\ldots,b_d)\|_L \leq \|\varpi_L\|_L^r$, while $\|r\|_L \geq \frac{1}{r^{[L:\Q_p]}}$. Thus, for any $r \geq 1$, $\|r^{-1}Q_r(b_1,\ldots,b_d)\|_L \leq \frac{\|\varpi_L\|_L^r}{r^{[L:\Q_p]}}$, and the right-hand side goes to zero, so the series $\sum_{r \geq 1}{\frac{1}{r}Q_r(b_1,\ldots,b_d)}$ converges in $L$. 

Let $(t'_1,\ldots,t'_d)$ be another parameter system. There exists $V_1,\ldots,V_d \in \OO_{K'}[[T_1,\ldots,T_d]]$ such that $t_i=V_i(t'_1,\ldots,t'_d)$: thus \[\omega = \sum_{1 \leq i,j \leq d}{P_i(V_1,\ldots,V_d)(t'_1,\ldots,t'_d))\frac{\partial V_i}{\partial T_j}(t'_1,\ldots,t'_d)dt'_j}.\]
Let $P'_j = \sum_{1 \leq i \leq d}{P_i(V_1,\ldots,V_d)\frac{\partial V_i}{\partial T_j}} \in \OO_{K'}[[T_1,\ldots,T_d]]$, then 
\[\omega = \sum_{j=1}^d{P'_i(t'_1,\ldots,t'_d)dt'_j}.\] 
The formal series $Q' := Q(V_1,\ldots,V_d) \in K'[[T_1,\ldots,T_d]]$ satisfies $\frac{\partial Q'}{\partial T_i} = P'_i$, so all we need to prove that, for any $A \in \,\,]x\mkern1mu[\,(L)$,
\[Q'(A^{\sharp}(t'_1),\ldots,A^{\sharp}(t'_d)) = Q(A^{\sharp}(t_1),\ldots,A^{\sharp}(t_d)).\]
To prove this, we need to show the following:
\begin{enumerate}[noitemsep,label=(\roman*)]
\item\label{tiwd1} For any $1 \leq i \leq d$, $A^{\sharp}(t_i)=V_i(A^{\sharp}(t'_1),\ldots,A^{\sharp}(t'_d))$. 
\item\label{tiwd2} For any $x_1,\ldots,x_d \in \varpi_L\OO_L$, for any $S \in \OO_{K'}[T_1,\ldots,T_d]$, one has \[[S \circ (V_1,\ldots,V_d)](x_1,\ldots,x_d) = S(V_1((x_i)_{1 \leq i \leq d}),\ldots,V_d((x_i)_{1 \leq i \leq d})).\]
\item\label{tiwd3} For any $x_1,\ldots,x_d \in \varpi_L\OO_L$, one has \[\underset{r \rar \infty}{\lim}\,\left[\sum_{r \geq 1}{\frac{1}{r}[Q_r\circ (V_1,\ldots,V_d)](x_1,\ldots,x_d)}\right] = Q'(x_1,\ldots,x_d).\]
\item\label{tiwd4} For any $x_1,\ldots,x_d \in \varpi_L\OO_L$, one has \[\underset{r \rar \infty}{\lim}\,\left[\sum_{r \geq 1}{\frac{1}{r}[Q_r(A^{\sharp}(t_1),\ldots,A^{\sharp}(t_d))}\right] = Q(x_1,\ldots,x_d).\]  
\end{enumerate}

\ref{tiwd4} was proved earlier. For \ref{tiwd2}, the set of $S$ such that the identity holds true is a subring of $\OO_{K'}[T_1,\ldots,T_d]$ containing $\OO_{K'}$ and every $T_i$, so it is all of $\OO_{K'}[T_1,\ldots,T_d]$. \\

As to \ref{tiwd1}, fix some $r \geq 1$ and let $z \in \OO_{X,x}^{\times}$ be such that its image in $\kappa(x)^{\times}$ is a generator of this multiplicative group. There is a polynomial $\Pi \in \OO_{K}[Z,T_1,\ldots,T_d]$ such that $V_i - \Pi \in (\varpi_K,T_1,\ldots,T_d)^r$, with $\varpi_K \in \OO_K$ a uniformizer: thus $t_i -\Pi(z,t'_1,\ldots,t'_d) \in \mathfrak{m}_{X,x}^r$, hence $A^{\sharp}(t_i)-\Pi(A^{\sharp}(z),A^{\sharp}(t'_1),\ldots,A^{\sharp}(t'_d)) \in \varpi_L^r\OO_L$. The assumption on $V_i-\Pi$ implies that \[\Pi(A^{\sharp}(z),A^{\sharp}(t'_1),\ldots,A^{\sharp}(t'_d)) - V_i(A^{\sharp}(t'_1),\ldots,A^{\sharp}(t'_d)) \in \varpi_L^r\OO_L,\] hence $A^{\sharp}(t_i)-V_i(A^{\sharp}(t'_1),\ldots,A^{\sharp}(t'_d)) \in \varpi_L^r\OO_L$. Since $r$ was arbitrary, this completes the proof of \ref{tiwd1}. \\

Finally, let us prove \ref{tiwd3}. By construction, \[Q' = \sum_{r \geq 1}{\frac{1}{r}Q_r\circ (V_1,\ldots,V_d)},\quad Q_r \circ (V_1,\ldots,V_d)\in \frac{1}{r}(T_1,\ldots,T_d)^r\OO_{K'}[[T_1,\ldots,T_d]].\] Therefore, it is enough to show that \[\left\|\frac{1}{r}[Q_r\circ (V_1,\ldots,V_d)](x_1,\ldots,x_d)\right\|_L \leq \frac{\|\varpi_L\|_L^r}{r^{[L:\Q_p]}},\] because the right-hand side is summable as $r \rar \infty$. Let $x'_i=V_i((x_j)_{1 \leq j \leq d})\in \varpi_L\OO_L$ for every $1 \leq i \leq d$. The series $Q_r \in \OO_{K'}[[T_1,\ldots,T_d]]$ is a homogeneous polynomial of degree $r$, hence \[\|[Q_r\circ (V_1,\ldots,V_d)](x_1,\ldots,x_d) \| = \|Q_r(x'_1,\ldots,x'_d)\| \leq \|\varpi_L\|_L^r,\] and the bound follows. 

The proof of the claim that the formula remains correct if $(t_i)_i$ is not assumed to be a parameter system is exactly the same.  
}

\medskip
 
\prop[tiny-integral-properties]{Let $K/\Q_p$ be a finite extension and $X$ be a proper $\OO_K$-scheme. Let $x \in X$ be a closed point in the special fibre such that $X$ is smooth at $x$ of relative dimension $d$. Let $K'/K$ be the finite unramified extension whose residue field is $\kappa(x)$ and let $L/K$ be a finite extension with uniformizer $\varpi_L$ endowed with a field homomorphism $\kappa(x) \rar \OO_L/\varpi_L$. The following properties hold true:

\begin{itemize}[noitemsep,label=$-$]
\item Additivity: if $A,B,C \in \,\,]x\mkern1mu[\,(L)$ and $\omega \in \Omega^1_{X/\OO_K,x}$ is closed, then $\int_{A}^C{\omega} = \int_A^B{\omega}+\int_B^C{\omega}$. 
\item Linearity: for any $A,B \in \,\,]x\mkern1mu[\,(L)$, $\lambda,\lambda' \in \OO_K$, for any closed $1$-forms $\omega,\omega' \in \Omega^1_{X/\OO_K,x}$, one has $\int_A^B{\lambda\omega+\lambda'\omega'} = \lambda\int_A^B{\omega}+\lambda'\int_A^B{\omega'}$.
\item Continuity: if $A,B \in \,\,]x\mkern1mu[\,(L)$, $\omega \in \mathfrak{m}_{X,x}^r\Omega^1_{X/\OO_K,x}$ is closed, then $\|\int_{A}^B{\omega}\| \leq \|\varpi_L\|_L^{n(r)}$, where $n(r) \in \Z$ only depends on $L$ and $n(r) \underset{r \rar +\infty}{\longrightarrow} +\infty$. 
\item Fonctoriality in $L$: let $L'/K$ be another finite extension with uniformizer $\varpi_{L'}$ endowed with a homomorphism $\kappa(x) \rar \OO_{L'}/\varpi_{L'}$ and let $\sigma: L \rar L'$ be a field homomorphism respecting the given morphisms $\kappa(x) \rar \OO_L/\varpi_L, \kappa(x) \rar \OO_{L'}/\varpi_{L'}$. Then, for any $A,B \in\,\, ]x\mkern1mu[\,(L)$, $\sigma(A),\sigma(B) \in \,\,]x\mkern1mu[\,(L')$ and, for any closed $\omega \in \Omega^1_{X/\OO_K,x}$, one has $\sigma\left(\int_A^B{\omega}\right) = \int_{\sigma(A)}^{\sigma(B)}{\omega}$.  
\item Fonctoriality in $X$: let $f: X \rar Y$ be a morphism of proper $\OO_K$-schemes, and assume that $Y$ is smooth at the point $y=f(x)\in Y$ (which is automatically closed and contained in the closed fibre). Let $\omega \in \Omega^1_{Y,y}$ be a closed differential and $A,B \in \,\,]x\mkern1mu[\,(L)$. Then $f(A),f(B) \in \,\,]y\mkern1mu[\,(L)$ (where $\kappa(y) \rar \OO_L/\varpi_L$ factors through $f^{\sharp}: \kappa(y)\rar \kappa(x)$), $f^{\ast}\omega \in \Omega^1_{X,x}$ is closed, and $\int_A^B{f^{\ast}\omega} = \int_{f(A)}^{f(B)}{\omega}$.
\item Base change: let $X'=X \times_{\OO_K} \Sp{\OO_L}$ and $f: X' \rar X$ be the natural projection. The field homomorphism $\kappa(x) \rar \OO_L/\varpi_L$ defines a unique point $x'_0 \in X'(\OO_L)$, and let $x' \in f^{-1}(\{x\})$ be the image of the closed point: $x'$ is a closed point in the closed fibre of $X'$, and $X'$ is smooth over $\Sp{\OO_L}$ at $x$. Let $A,B \in\,\, ]x\mkern1mu[\,(L)$ and $\omega \in \Omega^1_{X/\OO_K,x}$ be closed, then $A,B$ have unique lifts $A',B' \in \,\,]x'\mkern1mu[\,(L)$, and $\int_A^B{\omega} = \int_{A'}^{B'}{f^{\ast}\omega}$.   
\end{itemize}}

\demo{Mirroring the construction and proof in Proposition \ref{tiny-integral-well-defined}, this reduces to statements on formal power series. }

\medskip

\cor[formal-immersion-tiny]{Consider the situation of Proposition \ref{tiny-integral-well-defined} with $d=1$ and $e < p-1$, where $e$ is the ramification index of $L/\Q_p$. Suppose that there exists distinct $A,B \in \,\,]x\mkern1mu[\,(L)$ such that $\int_A^B{\omega}=0$. Then the ideal of first coefficients of $\OO_K \cdot \omega$ is not $\OO_K$. }

\demo{Let $t$ be a parameter for $X$ at $x$: we can write $\omega = P(t)dt$ in $\Omega^1_{X/\OO_K,x} \otimes_{\OO_{X,x}} \widehat{\OO_{X,x}}$ for some $P=\sum_{n \geq 0}{c_nT^n} \in \OO_K[[t]]$. 
Let $a = A^{\sharp}(t), b = B^{\sharp}(t) \in \varpi_L\OO_L$. Then $a \neq b$ by Corollary \ref{smooth-residue-disks}, and by assumption \[c_0(a-b) = \sum_{n \geq 1}{\frac{c_{n+1}}{n+1}(b^{n+1}-a^{n+1})}.\]
Let us show for any $n \geq 2$, one has $\left\|\frac{a^n-b^n}{n}\right\|_L < \|a-b\|_L$. Since $p \neq 2$, the statement is true for $n=2$, so we may assume $n \geq 3$. We need to prove that $\left\|\sum_{i=0}^{n-1}{a^ib^{n-1-i}}\right\|_L < \|n\|_L$. It is enough to check that $v_L(n) < n-1$, where $v_L$ is the valuation on $L$ normalized so that $v_L(L^{\times})=\Z$. Let $q=p^r$ be the greatest power of $p$ which is not greater than $n$. Then $v_L(n) < (p-1)v_p(n) \leq (p-1)r \leq (1+(p-1)r)-1 \leq (1+(p-1))^r-1 \leq n-1$, whence the conclusion. 

In particular, $\|c_0(b-a)\|_L =\left\|\sum_{n \geq 1}{\frac{c_{n+1}}{n+1}(b^{n+1}-a^{n+1})}\right\|_L < \|b-a\|_L$, whence $v_L(c_0) > 0$.
}

\prop[integral-abelian-variety]{Let $K/\Q_p$ be a finite extension with residue field $k$ and $A$ be an Abelian scheme over $\OO_K$. Let $0_k$ (resp. $0_K$) be the unit point in the special fibre (resp. the generic fibre) of $A$. Let $\omega \in H^0(A,\Omega^1)$, and let $L/K$ be a finite extension. Then $\omega$ is closed, $]0_k\mkern1mu[\,(L)$ is a finite index subgroup of $A(L)$, and $P \in\,\, ]0_k\mkern1mu[\,(L) \longmapsto \int_{0_K}^P{\omega_{0_k}} \in L$ is a group homomorphism. }

\demo{Let $\varpi_L \in \OO_L$ be a group homomorphism. The reduction map is functorial with respect to proper $\OO_K$-schemes, so $\mrm{red}: A(L) \rar A(\OO_L/\varpi_L)$ is a group homomorphism. Since $]0_k\mkern1mu[\,(L)=\mrm{red}^{-1}(\{0_k\})$ by definition, $]0_k\mkern1mu[\,(L)$ is a subgroup of $A(L)$, and $A(L)/]0_k\mkern1mu[\,(L)$ injects in $A_k(\OO_L/\varpi_L)$. Since $\OO_L/\varpi_L$ is a finite extension of the finite field $k$, $A_k(\OO_L/\varpi_L)$ is finite, hence so is $A(L)/]0_k\mkern1mu[\,(L)$. 

Note that $H^0(A,\Omega^1_{A/\OO_K})$ embeds in $H^0(A_K,\Omega^1_{A_K/K})$ and is made of translation-invariant differentials by Proposition \ref{differentials-are-trivial-abelian}. Hence the differential $\omega$ is closed by Proposition \ref{differentials-are-closed}. By additivity, we have to show that for any $P,Q \in\,\, ]0_k\mkern1mu[\,(L)$, $\int_0^Q{\omega} = \int_{P}^{P+Q}{\omega}$. By base change, we may assume that $L=K$. Let $\tau: X \in A \longmapsto X+P \in A$ be the translation. Then $\int_0^Q{\omega} = \int_0^Q{\tau^{\ast}\omega} = \int_{\tau(0)}^{\tau(Q)}{\omega} = \int_P^{P+Q}{\omega}$, whence the conclusion. }

\medskip

The following result is the goal of the section. It is a version of Mazur's \emph{formal immersion argument}, first described in the proof of \cite[Corollary 4.3]{FreyMazur}, and then an important fixture in the determination of rational points on modular curves. Our formulation is perhaps less geometric, with its emphasis on differentials satisfying a \emph{rational} condition (as opposed to an \emph{integral} one), but it is convenient one for the purposes of the next section. 

\prop[formal-immersion-criterion]{Let $A$ be an Abelian scheme over $\Z_p$ for some odd prime $p$, $X$ be a smooth proper $\Z_p$-scheme of relative dimension one, and $f: X \rar A$ be a morphism of $\Z_p$-schemes. Let $x_0 \in X(\Z_p)$ be such that $f(x_0)=0_A$ and $x \in X$ be the image under $x_0$ of the closed point of $\Sp{\Z_p}$. Let $M \subset \Omega^1_{X/\Z_p,x}$ be the set of stalks of differentials $\omega \in H^0(X,\Omega^1)$ such that $\omega=f^{\ast}\omega_1$ for some $\omega_1 \in H^0(A_{\Q_p},\Omega^1)$. 

Assume that the ideal of first coefficients of $M$ is $\Z_p$. Then $f:\,\, ]x\mkern1mu[\,(\Q_p) \rar A(\Q_p) \otimes \Q$ is injective. In particular, $\{x_0\} = \,\,]x\mkern1mu[\,(\Q_p) \cap f^{-1}(A(\Q_p)_{tors})$. 
}

\demo{Let $y,y' \in\,\, ]x\mkern1mu[\,(\Q_p)$ and $d \geq 1$ be such that $d(f(y)-f(y'))=0$. Let $\omega \in M$ be such that the ideal of first coefficients of $\Z_p\cdot \omega$ is $\Z_p$: there exists $n \geq 1$ and $\omega_1 \in H^0(A,\Omega^1)$ such that $p^n\omega=(f^{\ast}\omega_1)_x$. Then 
\begin{align*}
\int_y^{y'}{\omega} &= p^{-n}\int_{f(y)}^{f(y')}{\omega_1} = p^{-n}\left(\int_0^{f(y')}{\omega_1}-\int_0^{f(y)}{\omega_1}\right)=p^{-n}\int_0^{f(y')-f(y)}{\omega_1}\\
&= \frac{1}{p^nd}\int_0^{d(f(y')-f(y))}{\omega_1} = \frac{1}{p^nd}\int_0^0{\omega_1}=0.
\end{align*}
This contradicts Corollary \ref{formal-immersion-tiny}. 
}

\medskip

\prop[formal-immersion-rational]{Let $A$ be an abelian variety over $\Q$ with good reduction at an odd prime $p$, $X$ be a smooth proper $\Z_{(p)}$-scheme of relative dimension one. Let $\mathcal{A}$ be the N\'eron model of $A$ over $\Z_{(p)}$. Let $f: X_{\Q} \rar A$ be a scheme morphism and $P \in X(\Q_p)$ be such that $f(P)=0_A$. Then $P$ lifts to some $P_0 \in X(\Z_p)$. Let $x \in X$ be the image under $P_0$ of the closed point of $\Sp{\Z_p}$. 

Suppose that $A(\Q)$ is finite and that some $Q \in X(\Q)$ distinct from $P$ specializes to $x$. Then, for every $\omega \in H^0(X,\Omega^1_{X/\Z_{(p)}}) \cap f^{\ast}H^0(A,\Omega^1_{A/\Q})$, $\omega_x \in \mathfrak{m}_{X,x}\Omega^1_{X/\Z_{(p)},x}$.}

\demo{That $P$ lifts to $P_0$ is exactly the valuative criterion for properness. Let $\mathcal{A}$ be the N\'eron model of $A$ over $\Z_{(p)}$; it is an Abelian scheme, and $f$ extends to a morphism of $\Z_{(p)}$-schemes $X \rar \mathcal{A}$, and let $\hat{f}: X_{\Z_p} \rar \mathcal{A}_{\Z_p}$ be its base change to $\Z_p$. 
Then $x$ identifies with a closed point $\hat{x} \in X_{\Z_p}$ of the special fibre, and $P,Q \in \,\,]\hat{x}\mkern1mu[\,(\Q_p)$ are distinct points. Moreover, $\hat{f}(P)=0_{A_{\Q_p}}$ and $\hat{f}(Q)=f(Q) \in A(\Q) \subset \mathcal{A}(\Q_p)$ is a torsion point. 

By Proposition \ref{formal-immersion-criterion}, for any $\omega \in H^0(X_{\Z_p},\Omega^1_{X/\Z_p})\cap f^{\ast}H^0(A_{\Q_p},\Omega^1_{A/\Q_p})$, one has \[\omega_{\hat{x}} \in \mathfrak{m}_{X_{\Z_p},\hat{x}}\Omega^1_{X_{\Z_p}/\Z_p,\hat{x}}.\] 

Let $\beta: X_{\Z_p} \rar X$ be the base change morphism and let $\omega \in H^0(X,\Omega^1_{X/\Z_{(p)}}) \cap f^{\ast}H^0(A,\Omega^1_{A/\Q})$. By the previous paragraph, $(\beta^{\ast}\omega)_{\hat{x}} \in \mathfrak{m}_{X_{\Z_p},\hat{x}}\Omega^1_{X_{\Z_p}/\Z_p,\hat{x}}$. The morphism \[\beta^{\ast}: \Omega^1_{X/\Z_{(p)},x} \otimes_{\OO_{X,x}} \OO_{X_{\Z_p},\hat{x}} \rar \Omega^1_{X_{\Z_p}/\Z_p,\hat{x}}\] is an isomorphism of free modules of rank one, so $\omega_x \in \mathfrak{m}_{X,x}\Omega^1_{X/\Z_{(p)},x}$. 

}

\section{An integrality result on $X_{\rho,\xi}(p)(\Q)$}

Let us recall the definition of exceptional prime given in the introduction of this Chapter. 

\defi[exceptional-primes]{Let $p \equiv 2, 5\pmod{9}$ be a prime number. We say that a prime number $\ell$ is \emph{exceptional} for $p$ if $\ell \equiv \pm 1 \pmod{p}$ and there exists a prime ideal $\lambda$ of $\Z[e^{2i\pi/p}]$ with residue characteristic $\ell$ such that, for every cube root of unity $j \in \F_{p^2}^{\times}$, \[u_j=\sum_{\substack{b \in \F_p\\1+b\sqrt{-3} \in \F_{p^2}^{\times 3}j}}{e^{\frac{2i\pi b}{p}}} \in \lambda.\] }

The goal of this section is to prove the following result.

\theo[formal-immersion-application]{Let $K/\Q$ be an imaginary quadratic field inert at $p$ and $N \leq \GL{\F_p}$ be the normalizer of a nonsplit Cartan subgroup $C$. Let $\rho: G_{\Q} \rar N$ be a surjective continuous group homomorphism such that $G_K=\rho^{-1}(C)$, $\det{\rho}=\omega_p^{-1}$, and $\rho$ is $3$-Cartan, so that $3 \mid p+1$. Let $\psi: C/\F_p^{\times} \rar F^{\times}$ be a character of order $3$ and $g \in \mathcal{S}_1(\Gamma_1(M_{3,\rho}))$ be the newform attached to the odd irreducible Artin representation $\mrm{Ind}_K^{\Q}{\psi(\rho_K)}$. 
Assume furthermore that:  
\begin{itemize}[noitemsep,label=$-$]
\item $X_{\rho,\xi}(\Q)$ contains a non-cuspidal point for some $\xi \in \mu_p^{\times}(\Qbar)$, corresponding to some elliptic curve $E/\Q$,
\item there exists $(f,\mathbf{1}) \in \mathscr{S}$ such that $L(f \times g,1) \neq 0$.
\end{itemize}
Then the conductor of $E$ is equal to $p^2M_{\rho}L$, where $M_{\rho}$ is the Artin conductor of $\rho$ (in the sense of Serre \cite{Conj-Serre}), and $L$ is a square-free product of exceptional primes that do not divide $M_{\rho}p$.  
  }
  
\rems{
\begin{itemize}[noitemsep,label=$-$]
\item The recent works of \cite{LFL,lombardo} show that the case where one has in fact $\rho(G_K)=3C$ is impossible, unless $E$ has CM. One could say that we replaced the global condition $\rho(G_K)\neq C$ with a local condition $\rho(I_p) \neq C$. 
\item By construction, one has $p \nmid M_{3,\rho}$, and, in this situation, the sign of the functional equation of $L(f \times g,s)$ is $+1$. 
\item It follows from \cite[Appendix, Thm 0.1]{Michel} that there is a bound $p_0$ depending only on $M_{3,\rho}$ such that if $p \geq p_0$, there exists $(f,\mathbf{1}) \in \mathscr{S}$ such that $L(f \times g,1) \neq 0$. We will use this result to give a more concrete application of Theorem \ref{formal-immersion-application} in Section \ref{diophantine-consequence-y2x3p}. 
\item It seems plausible \cite{Brumley} that the bound $p_0(M_{3,\rho})$ can be made effective. Can one expect a bound that does not depend on $M_{3,\rho}$? 
\item The curve $X_{\rho,\xi}(p)$ has genus at least two: by Faltings' theorem, it contains finitely many rational points. Theorem \ref{formal-immersion-application} thus holds for \emph{some} set of ``exceptional primes'': what is nontrivial is that we manage to make these primes explicit.  
\item Let $j \in \F_{p^2}^{\times}$ be a primitive third root of unity. Numerical computations with MAGMA \cite{magma} show that for $p \leq 1500$, the norms of $u_1, u_j, u_1-u_j, u_j-u_{j^2}$ generate the unit ideal of $\Z$, so there is no exceptional prime. However, the running time of the code grows quickly\footnote{The computation up to $p \leq 1200$ takes slightly under $9$ minutes; the computation up to $p \leq 1500$ takes slightly under $30$ minutes.} around $p \geq 1200$, making a direct computation of exceptional primes harder for larger primes $p$. We thus expect\footnote{As written elsewhere, further computations since this thesis was submitted established that the claim still held for $11 \leq p \leq 5000$.} that there are no exceptional primes for any prime $p \geq 11$ congruent to $2$ or $5$ modulo $9$. 
\item Suppose that $\rho$ comes from an elliptic curve $E_0 \neq E$: the theorem shows that the conductor of $E_0$ is also of the form $p^2M_{\rho}L'$, where $L'$ is a product of exceptional primes $\ell \equiv \pm 1\pmod{p}$ that are coprime to $M_{\rho}$. Even if $L=L'=1$, proving that $E$ and $E_0$ are isogenous still seems difficult. 
\end{itemize}}

Here is an outline of the proof. 

The first step is to show that the conductor of $E$ is of the form $p^2M_{\rho}L$, where $L$ is a square-free product of primes $\ell \equiv \pm 1\pmod{p}$ not dividing $pM_{\rho}$. This is accomplished by using Lemmas \ref{new-is-semistable} (for the primes distinct from $p$) and Lemma \ref{cartan-is-potentially-good} at $p$. The second step is to construct a suitable quotient $A_I$ of $J_{\rho,\xi}(p)$ with expected finite Mordell-Weil group: this is done using Corollary \ref{make-abelian-variety-If} and Proposition \ref{BK-for-If}. For the rest of the argument, we fix odd primes $q \equiv 1 \pmod{p}$, $\ell \mid L$ with $\ell \equiv \pm 1\pmod{p}$ not exceptional, an embedding  $\iota: \Qbar \rar \overline{\Q_{\ell}}$ inducing a group homomorphism $R_{\iota}: G_{\Q_{\ell}} \rar G_{\Q}$, and a polarized $p$-torsion group $G$ over $\Z_{(\ell)}$ such that $X_{G,1}(p)_{\Q} \simeq X_{\rho,\xi}(p)$.  

In particular, the rational point of $X_{\rho,\xi}(p)$ reduced modulo $\ell$ to a point $c_0 \in X_{G,1}(p)(\F_{\ell})$ lying in the cuspidal subscheme of $X_{G,1}(p)$. Since this cuspidal subscheme is \'etale, there is a unique section $c \in X_{G,1}(p)(\Z_{\ell})$ of the cuspidal subscheme lifting $c_0$. 

The morphism $f: P \in X_{G,1}(p)_{\Z_{(\ell)}} \mapsto (T_q-q-1)(P - c)\in (A_I)_{\Z_{(\ell)}}$ is well-defined over $\Q$, because, even if $c$ is not rational, $T_q-q-1$ annihilates the difference of any two cusps by Proposition \ref{manin-drinfeld}. So $f$ is well-defined over $\Z_{(\ell)}$. 

Finally, we prove that the ideal of first coefficients of the collection of stalks at $c_0$ of differentials $\omega \in H^0(X_{G,1}(p)_{\Z_{\ell}},\Omega^1)$ such that, for some $n \geq 1$, $\ell^n\omega_1 \in f^{\ast}H^0((A_I)_{\Q_{\ell}},\Omega^1)$ is all of $\Z_{\ell}$. Then $E,c \in \,\,]c_0\mkern1mu[\,(\Q_{\ell})$ are mapped by $f$ to torsion points, which contradicts Corollary \ref{formal-immersion-criterion}. 

Throughout this argument, we will use without mention basic facts about differentials on Abelian schemes and relative Jacobians proved in Appendix \ref{differentials-relative-jacobians}.

\medskip

\lem[new-is-semistable]{Let $\ell,p$ be distinct primes with $p > 3$. Let $E/\Q_{\ell}$ be an elliptic curve with conductor exponent $n_E$. Write $G_{\ell}$ (resp. $I_{\ell}$) for the absolute Galois group of $\Q_{\ell}$ (resp. its inertia subgroup). Then the conductor exponent $n_{E,p}$ of $E[p]$ (as a $G_{\ell}$-representation) is contained in $\{n_E-1,n_E\}$. If $n_{E,p} < n_E$, then $(n_{E,p}, n_E)=(0,1)$. In this case, the $\ell$-adic valuations of the $j$-invariant and the minimal discriminant of $E$ are divisible by $p$; moreover, one has $\mrm{Tr}(\Fr_{\ell} \mid E[p]) \equiv a_{\ell}(E)(\ell+1) \pmod{p}$ (with $a_{\ell}(E)=1$ if $E$ has split multiplicative reduction and $a_{\ell}(E)=-1$ if $E$ has nonsplit multiplicative reduction). 
}

\demo{Let $s_E$ (resp. $s_{E,p}$) be the Swan conductor exponents of $E[p]$ (resp. $\Tate{p}{E}$). Then one has by definition 
\[n_E=s_E+2-\mrm{rk}_{\Z_p}\,{H^0(I_{\ell},\Tate{p}{E})},\quad n_{E,p} = s_{E,p}+2-\dim_{\F_p}{H^0(I_{\ell},E[p])}.\]
Since $p \neq \ell$, one has $s_E=s_{E,p}$. Moreover, the homomorphism \[\tau: H^0(I_{\ell},\Tate{p}{E}) \otimes_{\Z_p} \F_p \rar H^0(I_{\ell},E[p])\] is injective, hence $n_E \geq n_{E,p}$, and equality holds if and only if $\tau$ is onto. 
 
Let us assume that $H^0(I_{\ell},E[p])$ is nonzero. Then, in a suitable basis of $\Tate{p}{E}$, the image of the action of $I_{\ell}$ is contained in the pro-$p$-group $\SL{\Z_p} \cap \begin{pmatrix} 1+p\Z_p & \ast \\ p\Z_p & 1+p\Z_p\end{pmatrix}$, so this action factors through the tame inertia quotient and $s_E=s_{E,p}=0$. Moreover, the image of the action of $I_{\ell}$ is a pro-cyclic subgroup of $\SL{\Z_p}\cap\begin{pmatrix} 1+p\Z_p &\ast\\p\Z_p & 1+p\Z_p\end{pmatrix}$: let $A \in \SL{\Z_p}$ be a pro-generator of this image. 

Since the action of $G_{\ell}$ by conjugation on the tame inertia quotient of $I_{\ell}$ is the multiplication by $\ell$, $A^{\ell}$ is similar to $A$. In particular, if $\alpha$ is an eigenvalue of $A$, then the characteristic polynomial of $A$ is $(X-\alpha)(X-\alpha^{-1})$, and $\alpha^{\ell}$ is also an eigenvalue of $A$: thus $\alpha^{\ell} \in \{\alpha,\alpha^{-1}\}$. Since $\alpha$ is contained in a quadratic extension of $\Q_p$ and $p > 3$, the multiplicative order of $\alpha$ is not divisible by $p$. Since $(X-\alpha)(X-\alpha^{-1}) \equiv (X-1)^2 \pmod{p}$ because $A \equiv \begin{pmatrix} 1 & \ast\\0 & 1\end{pmatrix}\pmod{p}$, the reduction of $\alpha$ modulo $p$ is $1$. Therefore, by Hensel's lemma, one has $\alpha=1$, so $A$ is similar in $\GL{\Z_p}$ to $\begin{pmatrix} 1 & a\\0 & 1\end{pmatrix}$ for some $a \in \Z_p$, and $n_E \leq 1$. 

On the other hand, if $H^0(I_{\ell},E[p])=0$, then $H^0(I_{\ell},\Tate{p}{E}) \otimes_{\Z_p} \F_p \rar H^0(I_{\ell},E[p])$ is an isomorphism so $n_E=n_{E,p}$. 

Thus, when $n_E \neq n_{E,p}$, $H^0(I_{\ell},E[p])$ and $n_{E,p} \leq n_E \leq 1$, hence $(n_{E,p},n_E)=(0,1)$ and $E[p]$ is invariant under $I_{\ell}$. Assume that this is the case: then $E$ has multiplicative reduction and $j(E) \notin \Z_{\ell}$. By \cite[Theorems V.5.2, V.5.3]{AEC2}, there exists $q_E \in \ell\Z_{\ell} \backslash \{0\}$ such that, for any extension $K/\Q_p$, if $E'$ denotes the elliptic curve $E$ if $E$ has split multiplicative reduction, and its unramified quadratic twist otherwise, $E'(K) \simeq K^{\times}/q_E^{\Z}$. If $\eta: G_{\ell} \rar \{\pm 1\}$ denotes the unramified character mapping $\Fr_{\ell}$ to $a_{\ell}(E)$, one has $E[p] \simeq \eta \otimes \begin{pmatrix}\omega_p & \ast\\0 & 1\end{pmatrix}$ as $G_{\ell}$-modules. 

In particular, since $E[p]$ is fixed under $I_{\ell}$, $q_E \in (\Q_{\ell}^{nr})^{\times p}$, hence $p \mid v_{\ell}(q_E)$. Moreover, one has $\mrm{Tr}(\Fr_{\ell}\mid E[p]) \equiv a_{\ell}(E)(\ell+1) \pmod{p}$. By \cite[Chapter IV.9, Theorem V.3.1]{AEC2}, $v_{\ell}(q_E)=v_{\ell}(\Delta_E)=-v_{\ell}(j(E))$, whence the conclusion. }

\medskip

\cor[step1]{Under the assumptions of Theorem \ref{formal-immersion-application}, the conductor of $E$ is of the form $p^2M_{\rho}L$, where $L$ is a square-free product of prime numbers $\ell \equiv \pm 1\pmod{p}$ not dividing $pM_{\rho}$. Moreover, for every $\ell \mid L$, $j(E) \notin \Z_{(\ell)}$. }

\demo{What is not given in Lemma \ref{new-is-semistable} is the conductor exponent at $p$ and the fact that $\ell \equiv \pm 1\pmod{p}$, if $\ell$ is a prime dividing the conductor of $E$ but not $pM_{\rho}$. By Lemma \ref{cartan-is-potentially-good}, $E$ has potentially good reduction at $p$ but bad reduction at $p$ (since $\rho$ is $3$-Cartan), thus the conductor exponent of $\rho$ at $p$ is $2$. If $E$ has bad reduction at the prime $\ell \nmid pM_{\rho}$, then, by Lemma \ref{new-is-semistable}, $E$ has multiplicative reduction at $\ell$. Therefore, $j(E) \notin \Z_{\ell}$ by \cite[Propositions VII.5.4, VII.5.5]{AEC1}, and, by \cite[Theorems V.5.2, V.5.3]{AEC2}, $\Fr_{\ell}$ acts on $E[p]$ by a matrix similar to $\pm \begin{pmatrix} \ell &\ast\\0 & 1\end{pmatrix}$. If this matrix is contained in $C$, then one has $\ell \equiv 1\pmod{p}$; if this matrix is in $N \backslash C$, its trace is zero and $\ell \equiv -1\pmod{p}$. }

\rem{The condition that $\ell \equiv \pm 1\pmod{p}$ implies in particular that $\ell \geq 2p-1 \geq 20$; it can be found in \cite[Theorem 1.5, Proposition 1.8]{Zywina-Surj}. We will mostly use this result to assume that $\ell \neq 2$, since the formal immersion argument at this prime is more intricate.  }

\bigskip

From now on, we assume that $p \equiv 2,5\pmod{9}$, so that $C \simeq C[3] \oplus 3C$. Let $(f,\mathbf{1}) \in \mathscr{S}$ and $\psi \in \mathcal{I}_f$ be a character of order three, and $V_f$ be the one-dimensional $F$-vector space on which $\mathbb{T}[N]$ acts as follows: $T_n$ (resp. $nI_2$, resp. $A \in N$) acts by $a_n(f)$ (resp. $1$, resp. $1$). Then the ring $\mathbb{T}_{1,N}(N)$ acts on $V_f \otimes_F \mrm{Ind}_C^N{\psi}$ by a character $\mathbb{T}_{1,N}(N) \rar F$, and we denote its kernel by $\mathfrak{I}$. We let $\mathfrak{I}_f$ denote the inverse image of $\mathfrak{I}$ in $\mathbb{T}_1$. 

Let $h$ be the weight one newform attached to the odd irreducible Artin representation $\left(\mrm{Ind}_C^N{\psi}\right)\circ \rho$. 

Let $\ell \equiv \pm 1\pmod{p}$ be a prime at which $\rho$ is unramified. There exists a polarized $p$-torsion group $\Gamma$ over $\Z_{(\ell)}$ endowed with a basis $(P,Q)$ of $\Gamma(\Qbar)$ such that for any $\sigma \in G_{\Q}$, $\rho(\sigma)\begin{pmatrix}\sigma(P)\\\sigma(Q)\end{pmatrix}=\begin{pmatrix}P\\Q\end{pmatrix}$. Let $\xi_0 =\langle P,\,Q\rangle_{\Gamma} \in \mu_p^{\times}(\Qbar)$. After replacing $\Qbar$ with a suitably large number field and using the valuative criterion for properness, we see that there are $P_0,Q_0 \in \Gamma(\overline{\Z}_{(\ell)})$ extending $P,Q$ (where $\overline{\Z}_{(\ell)} = \overline{\Z} \otimes \Z_{(\ell)}$). 

In particular, by Propositions \ref{group-is-twist-polarized} and \ref{jac+hecke-twist-polarized}, for any $t \in \F_p^{\times}$, $X_{\Gamma,t}(p)_{\Q}$ is naturally isomorphic to $X_{\rho,\underline{t}\xi_0}(p)$, and this identification respects the Hecke operators.

\lem[MDtrick]{Let $q \equiv 1\pmod{p}$ be a prime, $\Q \subset H \subset \Qbar$ be a number field, and $c \in X_{\rho,\xi}(H)$ be a point in the cuspidal subscheme for some $\xi \in \mu_p^{\times}(\Qbar)$. Then the morphism \[P \in X_{\rho,\xi}(p)_H \longmapsto (T_q-q-1)(P-c) \in J_{\rho,\xi}(p)_H\] descends to a morphism $X_{\rho,\xi}(p) \rar J_{\rho,\xi}(p)$, and this morphism is independent from the choice of $c$.}

\demo{Since $X_{\rho,\xi}(p), J_{\rho,\xi}(p)$ are smooth proper geometrically connected $\Q$-schemes, by Corollary \ref{inf-galois-desc}, it is enough to show that $P \in X_{\rho,\xi}(p)(\Qbar) \longmapsto (T_q-q-1)(P-c) \in J_{\rho,\xi}(\Qbar)$ commutes with the action of $G_{\Q}$. Since $T_q-(q+1) \in \mrm{End}(J_{\rho,\xi}(p))$, all we need to do is show that, for any two $c,c' \in X_{\rho,\xi}(p)(\Qbar)$ contained in the cuspidal subscheme, one has $(T_q-q-1)(c-c')=0$. This is exactly Proposition \ref{manin-drinfeld}.}

\rem{In fact, this morphism is independent from the choice of $\rho$. }

\medskip

\lem[quotient-vs-differentials]{Let $A$ be an Abelian variety over a field $k$ and $I$ be a subgroup of $\mrm{End}_k(A)$. Let $\pi: A \rar A/IA$ be the quotient homomorphism. Then $\pi^{\ast}H^0(A/IA,\Omega^1)$ is exactly the $I$-torsion subspace of $H^0(A,\Omega^1)$. }

\demo{Let $A'=A/IA$. For $V\in \{A,A'\}$, there is a perfect pairing $\mrm{Lie}(V) \times H^0(V,\Omega^1_{V/k}) \rar k$ (by Corollary \ref{differentials-are-infinitesimal}). The differentials in $\pi^{\ast}H^0(A',\Omega^1_{A'/k})$ are clearly in the $I$-torsion submodule of $H^0(A,\Omega^1_{A/k})$; their image is exactly the orthogonal of $\ker{\pi: \mrm{Lie}(A) \rar \mrm{Lie}(A')}$, which is the orthogonal of $I \cdot \mrm{Lie}(A)$, which is the $I$-torsion submodule of $H^0(A,\Omega^1_{A/k})$, whence the conclusion.}

\prop[step2]{Assume that $X_{\rho,\xi}(p)$ has a non-cuspidal rational point $P$ corresponding to an elliptic curve $E/\Q$ with conductor $p^2M_{\rho}L$ by Corollary \ref{step1} such that $\ell \mid L$. Assume furthermore that $L(f \times h, 1) \neq 0$. 

Then, there exists a point $c \in X_{\rho,\xi}(p)$ lying in the special fibre of the cuspidal subscheme such that for any $\omega \in H^0(X_{\Gamma,t}(p),\Omega^1_{X_{\Gamma,t}(p)/\Z_{(\ell)}})$, if $\mathfrak{I}\cdot \omega=0$, then $\omega_c \in \mathfrak{m}_{X_{\Gamma,t}(p),c}\,\Omega^1_{X_{\Gamma,t}(p)/\Z_{(\ell)}, c}$. }

\demo{Let $c_0 \in X_{\Gamma,t}(p)$ be the reduction of $E \in X_{\Gamma,t}(p)(\Q)$ modulo $\ell$. Since $j(E) \notin \Z_{(\ell)}$ by Corollary \ref{step1}, $c_0$ lies in the cuspidal subscheme of $X_{\Gamma,t}(p)$. The cuspidal subscheme of $X_{\Gamma,t}(p)$ is finite \'etale over $\Z_{(\ell)}$, so there is a $c \in X_{\Gamma,t}(p)(\Qbar)$ and a prime ideal $\mathfrak{l} \subset \overline{\Z}$ of residue characteristic $\ell$ such that $c_0$ is the reduction of $c$ modulo $\mathfrak{l}$. For the same reason, $c_{\mathfrak{l}} \in X_{\Gamma,t}(p)(\Qbar_{\mathfrak{l}})$ factors through some $c'_{\mathfrak{l}} \in X_{\Gamma,t}(p)(\Q_{\ell})$. 

Let $A = J_{\rho,\xi}(p)/\mathfrak{I}$: it is an Abelian variety over $\Q$ with good reduction at $\ell$ (since it is a quotient of $J_{\rho,\xi}(p)$ which has good reduction at $\ell$ by the N\'eron-Ogg-Shafarevich criterion \cite[Theorem 1]{GoodRed} and the description of its Tate module given in Proposition \ref{tate-module-xrho-nsplit-cartan}), and let $\mathcal{A}$ be its N\'eron model over $\Sp{\Z_{(\ell)}}$. 

Let $q \equiv 1\pmod{p}$ be a prime and $f: X_{\Gamma,t}(p) \rar \mathcal{A}$ be the morphism extending \[x \in X_{\Gamma,t}(p)_{\Q} \simeq X_{\rho,\xi}(p) \longmapsto ((T_q-q-1)(x-c) \pmod{\mathfrak{I}}) \in A_{\Q}\] as in Lemma \ref{MDtrick}. 

By Proposition \ref{BK-for-If}, $A(\Q)$ is a finite group, so we can apply Corollary \ref{formal-immersion-rational}: for any $\omega \in H^0(X_{\Gamma,t}(p),\Omega^1_{X_{\Gamma,t}(p)/\Z_{(\ell)}}) \cap f^{\ast}H^0(A,\Omega^1_{A/\Q})$, one has $\omega_{c_0} \in \mathfrak{m}_{X_{\Gamma,t}(p),c_0}\Omega^1_{X_{\Gamma,t}(p)/\Z_{(\ell)},c_0}$. 

The conclusion follows once we prove that $H^0(X_{\Gamma,t}(p),\Omega^1_{X_{\Gamma,t}(p)/\Z_{(\ell)}}) \cap f^{\ast}H^0(A,\Omega^1_{A/\Q})$ is exactly the $\mathfrak{I}$-torsion submodule of $H^0(X_{\Gamma,t}(p),\Omega^1_{X_{\Gamma,t}(p)/\Z_{(\ell)}})$. 

Since both submodules are saturated in $H^0(X_{\Gamma,t}(p),\Omega^1_{X_{\Gamma,t}(p)/\Z_{(\ell)}})$, it is enough to see that $f^{\ast}H^0(A,\Omega^1_{A/\Q})$ is the $\mathfrak{I}$-torsion submodule of $H^0(X_{\Gamma,t}(p)_{\Q},\Omega^1_{X_{\Gamma,t}(p)_{\Q}/\Q})$. 

By Proposition \ref{differentials-on-jacobian-curve}, it is enough to show that $T_q-q-1: J_{\rho,\xi}(p) \rar J_{\rho,\xi}(p)$ is an isogeny (note that it commutes to $\mathfrak{I}$) and that $\pi^{\ast}H^0(A,\Omega^1_{A/\Q})$ is exactly the $\mathfrak{I}$-torsion submodule of $H^0(J_{\rho,\xi}(p),\Omega^1_{J_{\rho,\xi}(p)/\Q})$. The second part is a consequence of Lemma \ref{quotient-vs-differentials}. 

By Corollary \ref{first-decomp-xpp}, the action of $T_q$ on $H^0(J_{\rho,\xi}(p),\Omega^1_{J_{\rho,\xi}(p)/\Q}) \otimes_{\Q} F$ is diagonalizable and all its complex eigenvalues are of the form $a_q(f)$ for $(f,\chi) \in \mathscr{S}\cup\mathscr{P}\cup\mathscr{C}$, so that $|a_q(f)| \leq 2\sqrt{q} < q+1$, thus $T_q-q-1: \mathbb{T}_1 \rar \mathbb{T}_1$ is injective. Hence $T_q-q-1$ divides (in $\mathbb{T}_1$) some positive integer $Q$. By \cite[\S 8]{MilAb}, $T_q-q-1$ is an isogeny, and we are done. 
}

\medskip

\lem[formal-immersion-is-geometric]{Let $X \rar \Sp{R}$ be a smooth proper morphism of relative dimension one. Let $c \in X$ and $V \subset H^0(X,\Omega^1_{X/R})$ be a subset. Let $R'$ be a $R$-algebra, $\beta: X_{R'} \rar X_R$ be the natural homomorphism and $c' \in X_{R'}$ be such that $\beta(c')=c$. Then the following are equivalent:
\begin{itemize}[noitemsep,label=$-$]
\item There exists $\omega \in V$ such that $\omega_c\notin \mathfrak{m}_{X,c}\,\Omega^1_{X/R,c}$. 
\item There exists $\omega \in V$ such that $(\beta^{\ast}\omega)_{c'} \notin \mathfrak{m}_{X_{R'},c'}\,\Omega^1_{X_{R'}/R',c'}$. 
\end{itemize}
}

\demo{Since $X \rar \Sp{R}$ is smooth of relative dimension one, we know that $\Omega^1_{X_R/R,c}$ is free of rank one over $\OO_{X_R/R,c}$ and that the natural map $\Omega^1_{X/R,c}\otimes_{\OO_{X,c}} \OO_{X_{R'},c'} \rar \Omega^1_{X_{R'}/R',c'}$ (given by $\beta^{\ast}$) is an isomorphism.  
Thus, for $\omega \in V$, $\omega_c \notin \mathfrak{m}_{X,c}\Omega^1_{X/R,c}$ if and only if $\omega_c$ is a basis of $\Omega^1_{X/R,c}$, if and only if $(\beta^{\ast}\omega)_{c'}$ is a basis of $\Omega^1_{X_{R'}/R',c'}$, if and only if $(\beta^{\ast}\omega)_{c'} \notin \mathfrak{m}_{X_{R'},c'}\Omega^1_{X_{R'}/R',c'}$.
}

\cor[step3]{Assume that $X_{\rho,\xi}(p)$ has a non-cuspidal rational point $P$ corresponding to an elliptic curve $E/\Q$ with conductor $p^2M_{\rho}L$ by Corollary \ref{step1} such that $\ell \mid L$ and that $L(f \times h, 1) \neq 0$. Then, for any flat local homomorphism $\Z_{(\ell)} \rar R$, there exists $c \in X(p,p)_R$ lying in the special fibre of the cuspidal subscheme such that for any $\omega \in H^0(X(p,p)_R,\Omega^1_{X(p,p)_R/R})$ such that $\mathfrak{I}\cdot \omega = 0$, one has $\omega_c \in \mathfrak{m}_{X(p,p)_R,c}\,\Omega^1_{X(p,p)_R/R,c}$.}

\demo{It is a consequence of Lemma \ref{formal-immersion-is-geometric} and Proposition \ref{step2}, as long as we can show that the $\mathfrak{I}$-torsion submodule of $H^0(X(p,p)_R,\Omega^1_{X(p,p)_R/R})$ is generated by the $\mathfrak{I}$-torsion submodule of $H^0(X(p,p)_{\Z_{(\ell)}},\Omega^1_{X(p,p)_{\Z_{(\ell)}}/\Z_{(\ell)}})$. This is true because $R$ is flat over $\Z_{(\ell)}$ (so in particular $H^0(X(p,p)_{\Z_{(\ell)}},\Omega^1_{X(p,p)_{\Z_{(\ell)}}/\Z_{(\ell)}}) \otimes_{\Z_{(\ell)}} R \rar H^0(X(p,p)_R,\Omega^1_{X(p,p)_R/R})$ is an isomorphism). } 

\bigskip

\lem[simplify-I]{Let $K/\Q$ be a field extension and $\omega \in H^0(X_{\rho,\xi}(p)_K,\Omega^1_{X_{\rho,\xi}(p)_K/K})$. Fix a generator $g \in C$ and let $\zeta_3 =g^{\frac{p^2-1}{3}}$. Then $\mathfrak{I}\cdot \omega=0$ if and only if the following conditions are satisfied:
\begin{itemize}[noitemsep,label=\tiny$\bullet$]
\item $\omega$ is the sum of its projections $\omega_{f'}$ to the $f'$-eigenspaces for $f' \in G_{\Q} \cdot f$ (one always has $(f',\mathbf{1}) \in \mathscr{S}$), 
\item $\omega \mid (1+\zeta_3+\zeta_3^2) = 0$,
\item $\omega \mid [g^{3(p-1)}+g^{3(1-p)}-2] = 0$, 
\item for any prime $q$ such that $q\det{g} \equiv 1\pmod{p}$ and $a_q(f) \neq 0$, choose integers $B,n_i,m_i$ with $B,n_i > 0$ and $n_i \equiv 1\pmod{p}$ such that $-Ba_q(f) = \sum_i{m_ia_{n_i}(f)}$. For every $f' \in G_{\Q} \cdot f$, one has $B\omega_{f'} \mid T_q(g+g^p) = \left(\sum_i{m_ia_{n_i}(f')}\right)\omega_{f'}$. 
\end{itemize}
Moreover, the $\mathfrak{I}$-torsion submodule of $H^0(X_{\rho,\xi}(p)_K,\Omega^1_{X_{\rho,\xi}(p)_K/K})$ has dimension $2|G_{\Q} \cdot f|$.  }

\demo{Because the action of $\mathbb{T}_{1,N}(N) \supset \mathfrak{I}$ is defined over $\Q$, it is enough to show the result when $K=\Q$. For the dimension, it is enough to check, by the proof of Proposition \ref{step2}, that $\dim{J_{\rho,\xi}(p)/\mathfrak{I}J_{\rho,\xi}(p)} = 2|G_{\Q} \cdot f|$. Since the only other Galois conjugate of $\psi: C/\F_p^{\times} \rar F^{\times}$ is $\psi^2=\psi^ p$, the conclusion follows from the description of the Tate module of $J_{\rho,\xi}(p)/\mathfrak{I}J_{\rho,\xi}(p)$ by Proposition \ref{make-abelian-variety-If}. 

The four conditions are necessary, because (by inspecting the Tate module of $A$ by Proposition \ref{make-abelian-variety-If}) $1+\zeta_3+\zeta_3^2, g^{3(p-1)}+g^{3(1-p)}-2, \sum_i{m_iT_{n_i}}-BT_q(g+g^p)$ are all contained in $\mathfrak{I}$, and the first condition is equivalent to $\mathfrak{I}_f\cdot \omega=0$. 

The claim is that the reunion of all four conditions is sufficient. In other words, it is enough to prove that $A$ is the quotient of $J_{\rho,\xi}(p)$ by $\mathfrak{I}_f, 1+\zeta_3+\zeta_3^2, g^{3(p-1)}+g^{3(1-p)}-2\mrm{id}$, and $\sum_i{m_iT_{n_i}}-BT_q(g+g^p)$ (for $q,m_i,n_i$ as given). 

It is enough to show that $\ker\left[\Tate{\ell}{J_{\rho,\xi}(p)} \otimes_{\Z_{\ell}} F_{\mathfrak{l}} \rar \Tate{\ell}{A} \otimes_{\Z_{\ell}} F_{\mathfrak{l}}\right]$ is exactly the submodule generated by the sum of the images of $\mathcal{I}_f, 1+\zeta_3+\zeta_3^2, g^{3(p-1)}+g^{3(1-p)}-2I_2$, and $\sum_i{m_iT_{n_i}}-BT_q(g+g^p)$. 

We use the description of Proposition \ref{decomposition-ncartan-xrho-T11rho} (and argue as in Proposition \ref{representation-is-rational}). The image of $\mathfrak{I}_f$ is made with those eigenspaces of $\Tate{\ell}{J_{\rho,\xi}(p)} \otimes_{\Z_{\ell}} F_{\mathfrak{l}}$ which are attached to newforms that are not Galois-conjugate to $f$. The image of $1+\zeta_3+\zeta_3^2$ is exactly made with the representations on which $\zeta_3$ acts as a matrix of order $3$ instead of the identity -- among others, it contains the spaces corresponding to $H_f$. Similarly, the image of $g^{3(p-1)}+g^{3(1-p)}-2I_2$ is exactly made with those subspaces where $g^{3(p-1)}$ does not act trivially. 

Thus, the summands of $\Tate{\ell}{J_{\rho,\xi}(p)} \otimes_{\Z_{\ell}} F_{\mathfrak{l}}$ which are not contained in (equivalently, do not meet) the image of $(\mathfrak{I}_f,1+\zeta_3+\zeta_3^2,g^{3(p-1)}+g^{3(1-p)}-2I_2)$ are those of the form $V_{f',\mathfrak{l}} \otimes \mrm{Ind}_{\tilde{\rho}^{-1}(C)}^{N\rtimes_{\rho} G_{\Q}}{\psi_1(\tilde{\rho_K})}$ with $(f',\mathbf{1}) \in \mathscr{S}$, $f' \in G_{\Q}\cdot f$, $\psi_1 \in\mathcal{I}_{f'}/\sim$, $\psi_1(\zeta_3) \neq 1$, $\psi_1(g)^{3(p-1)}=1$ hence $\psi_1(g)^6=1$.

Thus, in this situation, either $\psi_1 \in \{\psi,\psi^2\}$, or $\psi_1 \in \{\psi_6,\psi_6^p\}$ with $\psi_6 = \left(\frac{\det}{p}\right)\psi$. The goal is to show that for $(q,B,(m_i)_i,(n_i)_i)$ as above, $\sum_i{m_iT_{n_i}}-BT_q(g+g^p) \in \mathbb{T}_{1,N}(N)$ vanishes when applied to $V_{f',\mathfrak{l}} \otimes \mrm{Ind}_{\tilde{\rho}^{-1}(C)}^{N\rtimes_{\rho} G_{\Q}}{\psi(\tilde{\rho_K})}$ but is surjective when applied to $V_{f',\mathfrak{l}} \otimes \mrm{Ind}_{\tilde{\rho}^{-1}(C)}^{N\rtimes_{\rho} G_{\Q}}{\psi_6(\tilde{\rho_K})}$. 

Indeed, $\sum_i{m_iT_{n_i}}$ acts on $V_{f',\mathfrak{l}} \otimes \mrm{Ind}_{\tilde{\rho}^{-1}(C)}^{N\rtimes_{\rho} G_{\Q}}{\psi_1(\tilde{\rho_K})}$ by multiplication by \[\sum_{i}{m_ia_{n_i}(f')}=\sigma\left(\sum_{i}{m_ia_{n_i}(f)}\right) = \sigma(-Ba_q(f))=-Ba_q(f').\]

However, $BT_q(g+g^p)$ acts on $V_{f',\mathfrak{l}} \otimes \mrm{Ind}_{\tilde{\rho}^{-1}(C)}^{N\rtimes_{\rho} G_{\Q}}{\psi(\tilde{\rho_K})}$ by $Ba_q(f')(j+j^2)=-Ba_q(f')$, where $j$ is a primitive third root of unity. On the other hand, $BT_q(g+g^p)$ acts on $V_{f',\mathfrak{l}} \otimes \mrm{Ind}_{\tilde{\rho}^{-1}(C)}^{N\rtimes_{\rho} G_{\Q}}{\psi_6(\tilde{\rho_K})}$ acts by $-Ba_q(f')(j+j^2)=Ba_q(f')$. 

Therefore, we are done if there \emph{exists} one instance of $(B,q,(m_i)_i,(n_i)_i)$ as above. The existence of $q$ follows from Serre's theorem \cite[Cor. 2 au Th\'eor\`eme 15]{Serre-Cebotarev}. The $a_n(f)$ are algebraic integers: as in Proposition \ref{representation-is-rational}, the $B, (m_i)_i, (n_i)_i$ as above exist provided that $a_q(f)$ is in the number field generated by the $a_n(f)$ with $n \equiv 1\pmod{p}$. 

Let $\sigma \in \mrm{Gal}(F/\Q)$ fixing all the $a_n(f)$ for $n \equiv 1\pmod{p}$. Then $\sigma(f) \in \mathcal{S}_2(\Gamma_0(p))$ is a newform. By Lemmas \ref{agreement-implies-twist}, $\sigma(f)$ is a twist of $f$; by Lemma \ref{gamma1p-onlytwists}, $\sigma(f)=f$, hence $\sigma(a_q(f))=a_q(f)$, whence the conclusion. 
}

\bigskip

\prop[complex-differentials-vs-T1N]{Let $a \in \F_p^{\times}$, and write $X(p,p)_a$ for the connected component $X(p) \times \{a\}$ of $X(p,p)_{\C}$ (see Proposition \ref{uniformize-xpp}). By Corollary \ref{differentials-xpp}, $\iota_{p,p}$ identifies $H^0(X(p,p)_a,\Omega^1)$ with $\mathcal{S}_2(\Gamma(p)) \otimes \Delta_{a^{-1},1}$ as right $\C[\mathbb{T}_{1,\GL{\F_p}}]$-modules. 

Fix a generator $g \in C$ and let $\zeta_3=g^{\frac{p^2-1}{3}} \in C \cap \SL{\F_p}$ with order $3$. Then, a generating family for the $H^0(X(p,p)_a,\Omega^1)[\mathfrak{I}]$ is given by 
\[\mathcal{F}_a := \left(\left(f'\left(\frac{\tau}{p}\right) \otimes \Delta_{a^{-1},1}\right) \mid \sum_{r=0}^{\frac{p+1}{3}-1}{\Delta_{(\det{g})^{-3r},1}\,g^{3r}}\zeta_3^i\right)_{\substack{f' \in G_{\Q} \cdot f\\i \in \{0,1,2\}}},\] where we recall that $\Delta_{a,b}$ denotes the matrix $\begin{pmatrix} a & 0\\0 & b\end{pmatrix} \in \GL{\F_p}$.

The relations within $\mathcal{F}_a$ are generated by the collection, for every $f' \in G_{\Q}\cdot f$, of
\[\sum_{i=0}^2{\left[\left(f'\left(\frac{\tau}{p}\right) \otimes \Delta_{a^{-1},1}\right) \mid \sum_{r=0}^{\frac{p+1}{3}-1}{\Delta_{(\det{g})^{-3r},1}\,g^{3r}}\zeta_3^i\right]}=0.\] 
}

\demo{Fix a homomorphism $\iota: \Qbar \rar \C$ mapping $\xi_0$ to $e^{-\frac{2i\pi}{p}}$. By Proposition \ref{polarized-twist-base-change}, $X(p,p)_a$ identifies with $X_{\mathbf{id},\iota(\underline{a}\xi_0)} \simeq \left(X_{\rho,\underline{a}\xi_0}\right)_{\C}$.

We will prove that the only $\C$-linear relations within $\mathcal{F}_a$ are the ones given in the statement of the Proposition, and that the span of $\mathcal{F}_a$ (over $\C$) is contained in the $\mathfrak{I}$-torsion submodule of $H^0(X_{\rho,\underline{a}\xi_0}(p)_{\C},\Omega^1)$. Since the dimensions match by Lemma \ref{simplify-I}, this is enough to conclude. 

For any $f' \in G_{\Q} \cdot f$ and $a \in \F_p^{\times}$, let \[u_{f',a} = f'\left(\frac{\tau}{p}\right)\otimes\Delta_{a,1}.\]

\emph{Step 1:}
Let us show that the function $r \in \Z \longmapsto u_{f',a}^r := u_{f',a} \mid \Delta_{(\det{g})^{-r},1}\,g^r$ factors through $\Z/(p+1)\Z$. Indeed, let $r,s \in \Z$, then
\begin{align*}
u_{f',a} \mid \Delta_{(\det{g})^{-r-(p+1)s},1}g^{r+(p+1)s} &= u_{f',1} \mid \Delta_{a,1}\Delta_{(\det{g})^{-r-2s},1}\cdot ((\det{g})^sI_2)\cdot g^r\\
&= u_{f',1} \mid \Delta_{(\det{g})^{-s},(\det{g})^s}\Delta_{a,1}\Delta_{(\det{g})^{-r},1}\,g^r\\
&= u_{f',1}\mid \Delta_{a,1}\Delta_{(\det{g})^{-r},1}\,g^r= u_{f',a} \mid \Delta_{(\det{g})^{-r},1}\,g^r,
\end{align*}
whence the conclusion. \\

\emph{Step 2: The only relations between the $u_{f',a}^r$ are that $\sum_{r \in \Z/(p+1)\Z}{u_{f',a}^r}=0$.}

Let $q \equiv 1 \pmod{p}$ be a positive integer and $f' \in G_{\Q} \cdot f$. Then $T_q(u_{f',a}^r)=a_q(f')u_{f',a}^r$. Moreover, by Lemmas \ref{agreement-implies-twist} and \ref{gamma1p-onlytwists}, the collection $(a_q(f'))_{q \equiv 1\pmod{p}}$ determines $f'$. Therefore, let, for every $(f',a) \in (G_{\Q} \cdot f) \times \F_p^{\times}$, $U_{f',a}$ be the $\C[\SL{\F_p}]$-module generated by $u_{f',a}$, then the $U_{f',a}$ are in direct sum. Thus, it is enough to show that for any $f' \in G_{\Q} \cdot f$, the only relations in the family $(u_{f',a}^r)_{r,a}$ are that $\sum_{r \in \Z/(p+1)\Z}{u_{f',a}^r}=0$ for any $a$. Thus, fix $f' \in G_{\Q} \cdot f$. 

For any $\chi \in \mathcal{D}$, let 
\[\hat{u}_{f',\chi} = \sum_{a \in \F_p^{\times}}{\chi^{-1}(a)u_{f',a}} \in \Omega^1[f',\chi] \backslash \{0\}\] in the sense of Corollary \ref{decomposition-hecke-xpp}. For any $X = \begin{pmatrix} a & 0\\b & d\end{pmatrix} \in \GL{\F_p}$, a direct verification shows that $\hat{u}_{f',\chi}\mid X = \chi(\det{X})\hat{u}_{f',\chi}$. Therefore, by Corollary \ref{no-lines}, the $\C[\GL{\F_p}]$-module generated by $\hat{u}_{f',\chi}$ is isomorphic to $\mrm{St}_{\chi}$: thus, the $\C[\GL{\F_p}]$-modules generated by the $\hat{u}_{f',\chi}$ (where $\chi \in \mathcal{D}$ varies) form a direct sum. 

Moreover, for any $r \in \Z/(p+1)\Z$, one has
\[\sum_{a \in \F_p^{\times}}{\chi^{-1}(a)u_{f',a}^r} = \hat{u}_{f',\chi}\mid \Delta_{(\det{g})^{-r},1}\,g^r = \chi(\det{g})^{-r}\hat{u}_{f',\chi}\mid g^r.\]

By Lemma \ref{steinberg-vs-cartan}, the only (nontrivial) relation in the family $\left(\sum_{a \in \F_p^{\times}}{\chi^{-1}(a)u_{f',a}^r}\right)_{r \in \Z/(p+1)\Z}$ ($\chi$ being fixed) is that their sum vanishes. The preceding remarks show that the relations in the family $(\hat{u}_{f',\chi} \mid \Delta_{(\det{g})^{-r},1}\,g^r)_{\substack{r \in \Z/(p+1)\Z \\ \chi\in \mathcal{D}}}$ are generated by the fact that for any $\chi$,
\[\sum_{r \in \Z/(p+1)\Z}{\hat{u}_{f',\chi} \mid \Delta_{(\det{g})^{-r},1}\,g^r} = 0.\] 

Since $u_{f',a} = \frac{1}{p-1}\sum_{\chi \in \mathcal{D}}{\chi(a)\hat{u}_{f',\chi}}$, the family $(u_{f',a}^r)_{a,r}$ has rank $(p-1)(p+1)-|\mathcal{D}|=p(p-1)$, so it has exactly $p-1$ relations, and, for any $a \in \F_p^{\times}$, $\sum_{r \in \Z/(p+1)\Z}{u_{f',a}^r}=0$. Hence these are all the relations in the family $(u_{f',a}^r)_{(a,r) \in \F_p^{\times} \times \Z/(p+1)\Z}$, whence the conclusion. \\

\emph{Step 3: The span of $\mathcal{F}_a$ has dimension $2|G_{\Q} \cdot f|$ and the only relations between the vectors of $\mathcal{F}_a$ are those given in the statement.}

Fix $a \in \F_p^{\times}$. For $i \in \{0,1,2\}$ and $f' \in G_{\Q}\cdot f$, let $v_{f',a,i} = \sum_{r \in 3\Z/(p+1)\Z}{u_{f',a^{-1}}^r\mid \zeta_3^i}$, so that $\mathcal{F}_a = (v_{f',a,i})_{\substack{f' \in G_{\Q}\cdot f\\i \in \{0,1,2\}}}$.

Since $p-1 \mid \frac{p^2-1}{3}$, one has $(\det{g})^{\frac{p^2-1}{3}}=1$. Moreover, $\frac{p^2-1}{3}\equiv \frac{p+1}{3}\pmod{p+1}$, so, for $r \in \Z$, one has $u_{f',a^{-1}}^r \mid \zeta_3^i = u_{f',a^{-1}}^{r+i\frac{p+1}{3}}$ by Step 1. Hence 
\[v_{f',a,i} = \sum_{r \in 3\Z/(p+1)\Z}{u_{f',a^{-1}}^{r+i\frac{p+1}{3}}}.\]
Since $\frac{p+1}{3}$ is not divisible by $3$, the conclusion follows from Step 2. \\

\emph{Step 4: The vectors of $\mathcal{F}_a$ are contained in the $\mathfrak{I}$-torsion submodule.}
 
We use the criterion of Lemma \ref{simplify-I}. Let $f' \in G_{\Q}\cdot f$ and $i \in \{1,2\}$, write $f'=\sigma(f)$ for some $\sigma \in \mrm{Aut}(\C/\Q)$. For any $m,n \geq 1$ such that $m^2n\equiv 1\pmod{p}$, $(mI_2)\cdot T_n$ acts on $v_{f',a,i}$ by multiplication by $a_n(f')=\sigma(a_n(f))$. By definition of $\mathfrak{I}_f$, $\mathfrak{I}_f$ thus annihilates every $v_{f',a,i}$. 

Since $\frac{p+1}{3}$ is coprime to $3$, one has 
\[v_{f',i}\mid (1+\zeta_3+\zeta_3^2) = \sum_{r \in \Z/(p+1)\Z}{u_{f',a^{-1}}^r} = 0.\] Moreover, for any $r \in \Z$, $\Delta_{(\det{g})^{-r-3(p-1)},1}\,g^{r+3(p-1)} =\Delta_{(\det{g})^{-r},1}\,g^{r}g^{3(p-1)}$. Therefore, one has $u_{f,a^{-1}}^r \mid g^{3(p-1)} = u_{f,a^{-1}}^{r+3(p-1)}$, which implies that $v_{f',a,i} \mid g^{3(p-1)} = v_{f',a,i}$. 

Finally, write $-Ba_q(f) = \sum_{j}{w_ja_{n_j}(f)}$, where $q$ is a prime such that $a_q(f) \neq 0$ and $q\det{g} \equiv 1\pmod{p}$, the $w_j$ are integers, the $n_j$ are positive integers congruent to $1$ modulo $p$, and $B > 0$ is an integer. We need to show that $Bv_{f',a,i} \mid T_q(g+g^p) = \sum_j{w_ja_{n_j}(f')} \cdot v_{f',a,i}$ for each $(f',a,i) \in (G_{\Q} \cdot f) \times \F_p^{\times} \times \{0,1,2\}$. Indeed, let $f'=\sigma(f)$, then
\begin{align*}
&Bv_{f',a,i} \mid T_q(g+g^p) = B a_q(f') \sum_{r \in 3\Z/(p+1)\Z}{u_{f',a^{-1}q}\mid \Delta_{(\det{g})^{-r},1}\,g^r\zeta_3^i(g+g^p)}\\
&= Ba_q(f') \sum_{r \in 3\Z/(p+1)\Z}{u_{f',a^{-1}}\mid \left[\Delta_{(\det{g})^{-r-1},1}g^{r+1+i\frac{p+1}{3}}+\Delta_{(\det{g})^{-r-p},1}\,g^{r+p+i\frac{p+1}{3}}\right]}\\
&= Ba_q(f')\sum_{r \in 3\Z/(p+1)\Z}{u_{f',a^{-1}}^{r+i\frac{p+1}{3}+1}+u_{f',a^{-1}}^{r+i\frac{p+1}{3}-1}}\\
&= -Ba_q(f') v_{f',a,i} = \sigma(-Ba_q(f))v_{f',a,i} \\
&= \sigma\left(\sum_{j}{w_ja_{n_j}(f)}\right)v_{f',a,i} = \left(\sum_j{w_ja_{n_j}(f')}\right)v_{f',a,i},
\end{align*}
whence the conclusion.}

\lem[q-exp-s+nsplit]{Let $a \in \F_p^{\times}$, $\gamma=\begin{pmatrix} u & v \\ \ast & \ast\end{pmatrix} \in \SL{\F_p}$ and $(f',\mathbf{1}) \in \mathscr{S}$. Then
\begin{align*}
\left(f'\left(\frac{\tau}{p}\right) \otimes \Delta_{a,1}\right) \mid \gamma &= \sum_{n \geq 1}{a_n(f)e^{\frac{2in\pi}{p}\frac{av}{u}}}e^{\frac{2i\pi n\tau}{p}} \otimes \Delta_{a,1} &&\text{ if }u\neq 0,\\
&= -pa_p(f)\sum_{n \geq 1}{a_n(f)e^{2i\pi n\tau}} \otimes \Delta_{a,1} &&\text{otherwise.} 
\end{align*}
}

\demo{One has 
\[\left(f'\left(\frac{\tau}{p}\right) \otimes \Delta_{a,1}\right)\mid\gamma = \left[f'\left(\frac{\tau}{p}\right) \mid \begin{pmatrix} u & av\\a^{-1}w & x\end{pmatrix}\right] \otimes \Delta_{a,1}.\]
Since $f'\left(\frac{\tau}{p}\right)$ is invariant under the right action of the lower-triangular matrices in $\SL{\F_p}$, when $u \neq 0$, we may assume that $\begin{pmatrix} u & av \\a^{-1}w & x\end{pmatrix} = \begin{pmatrix} 1 & a\\0 & 1\end{pmatrix}$, in which case the conclusion is clear. 
Similarly, when $u = 0$, we may assume that $\begin{pmatrix} u & av\\a^{-1}w & x\end{pmatrix} = \begin{pmatrix} 0 & -1\\1 & 0\end{pmatrix}$. Then, by Proposition \ref{local-constant-principal},
\[f'\left(\frac{\tau}{p}\right)\mid \begin{pmatrix}0& -1\\1 & 0\end{pmatrix} = p f \mid \begin{pmatrix}1 & 0\\0 & p\end{pmatrix}\begin{pmatrix}0 & -1\\1 & 0\end{pmatrix} = p\lambda_p(f) f=-pa_p(f)f.\]}

Let us now fix an element $j \in \F_{p^2}^{\times}$ with order $3$ (hence $j \notin\F_p^{\times}$). There is a unique $\F_p$-basis $(\alpha,\beta)$ of $\F_{p^2}$ such that, for any $x=\begin{pmatrix}u & v\\\ast & \ast\end{pmatrix} \in C$, $u\alpha+v\beta$ has the same characteristic polynomial as $x$, and such that, if moreover $x=\zeta_3$, one has $u\alpha+v\beta=j$. We say that $(\alpha,\beta)$ is the \emph{adapted} basis of $\F_{p^2}$.

\defi[BK4S-sigma]{Given a $\F_p$-basis $\mathcal{B}=(u,v)$ of $\F_{p^2}$ and a third root of unity $j \in \Z/3\Z$, let
\[\Sigma_{\mathcal{B},t}=\sum_{\substack{au+bv \in \F_{p^2}^{\times 3}j^t/\F_p^{\times}\\a \neq 0}}{e^{\frac{2i\pi b/a}{p}}} \in \Z[\zeta_p].\]
}

\rem[BK4S-vs-exceptional]{Since $p \equiv -1\pmod{3}$, $-3$ is not a square in $\F_p^{\times}$, and one can check directly that $\Sigma_{(1,\sqrt{-3}),t} = u_{j^t}$. }

\lem[BK4S-elementary]{The following statements hold true:
\begin{itemize}[noitemsep,label=$-$]
\item When $\mathcal{B}$ is a fixed $\F_p$-basis of $\F_{p^2}$, the only $\Z$-linear relation between the $\Sigma_{\mathcal{B},t}$ is that their sum is zero. 
\item Let $(u,v)$ be a basis of $\F_{p^2}$, $\lambda \in \F_p$ and $t \in \Z/3\Z$, then $\Sigma_{(u+\lambda v,v),t}=e^{-2i\pi\lambda/p}\Sigma_{(u,v),t}$.
\item Let $(u,v)$ be a basis of $\F_{p^2}$, $s,t \in \Z/3\Z$ and $\alpha \in \F_{p^2}^{\times 3}j^s$. Then $\Sigma_{(\alpha u,\alpha v),t+s}=\Sigma_{(u,v),t}$.  
\item Let $(u,v)$ be a basis of $\F_{p^2}$, $t \in \Z/3\Z$ and $a,b \in \F_p^{\times}$. Then $\Sigma_{(au,bv),t}=\underline{ab^{-1}}(\Sigma_{(u,v),t})$.
\end{itemize}}

\demo{Let $\mathcal{B}=(u,v)$. For every $w\in \F_{p^2}^{\times}/\F_p^{\times}$, let $\omega(w)=e^{\frac{2i\pi}{p}ba^{-1}}$ if $w=(au+bv)\F_p^{\times}$ with $a \neq 0$, and $\omega(w)=0$ otherwise. We claim that the only $\Z$-linear relation between the $\omega(w)$ for $w \in \F_{p^2}^{\times}/\F_p^{\times}$ is that $\omega(\sqrt{-3})=0$ (since $p \equiv -1\pmod{3}$, $-3 \notin \F_p^{\times}$) and $\sum_w{\omega(w)}=0$. Indeed, it is direct to check that $\omega: \F_{p^2}^{\times}/\F_p^{\times} \rar \mu_p(\C) \cup \{0\}$ is a bijection. Now, note that $\Sigma_{\mathcal{B},t} = \sum_{w \in \F_p^{\times 2}j^t/\F_p^{\times}}{\omega(w)}$, whence the conclusion. 

The other claims are direct computations. 
}

\cor[exceptional-primes-criterion]{Let $\ell \neq p$ be a prime. The following are equivalent:
\begin{itemize}[label=$-$,noitemsep]
\item $\ell$ is exceptional.
\item For any $\F_p$-basis $\mathcal{B}$ of $\F_{p^2}$, there exists a prime ideal $\lambda \subset \Z[\zeta_p]$ dividing $\ell$ such that the $\Sigma_{\mathcal{B},t} \pmod{\lambda}$ are all zero. 
\item There exists a $\F_p$-basis $\mathcal{B}$ of $\F_{p^2}$ such that $\ell$ divides the norm of the (nonzero) ideal $(\Sigma_{\mathcal{B},0},\Sigma_{\mathcal{B},1},\Sigma_{\mathcal{B},2})$.  
\item For any $\F_p$-basis $\mathcal{B}$ of $\F_{p^2}$, $\ell$ divides the norm of the (nonzero) ideal $(\Sigma_{\mathcal{B},0},\Sigma_{\mathcal{B},1},\Sigma_{\mathcal{B},2})$.
\item For any primitive $p$-th root of unity $\zeta \in \overline{\F_{\ell}}$, there exists a third root of unity $j \in \F_{p^2}^{\times}$ such that 
\[ \sum_{\substack{t \in \F_p\\1+t\sqrt{-3} \in \F_{p^2}^{\times 3}j}}{\zeta^t}\neq 0.\]
\end{itemize}}

\demo{Let $\mathcal{B}$ be any $\F_p$-basis of $\F_{p^2}$. Let $\lambda \subset \Z[\zeta_p]$ be a prime ideal. Then the $\Sigma_{\mathcal{B},t} \pmod{\lambda}$ are all zero if and only if the ideal $I_{\mathcal{B}} := (\Sigma_{\mathcal{B},0}\Sigma_{\mathcal{B},1},\Sigma_{\mathcal{B},2})$ is contained in $\lambda$. Thus, given a prime number $\ell$, there exists a prime ideal $\lambda \subset \Z[\zeta_p]$ with residue characteristic $\ell$ such that the $\Sigma_{\mathcal{B},t} \pmod{\lambda}$ are all zero if and only if the norm of $I_{\mathcal{B}}$ is divisible by $\ell$. 

To conclude, thanks to Remark \ref{BK4S-vs-exceptional}, it is enough to prove that the Galois orbit of the ideal $I_{\mathcal{B}}$ does not depend on $\mathcal{B}$. Indeed, write $\mathcal{B}=(u,v)$ and let $\mathcal{B}'=(u',v')$ be another basis. Then, by Lemma \ref{BK4S-elementary}, $I_{\mathcal{B}}=I_{\frac{v'}{v}\mathcal{B}}=I_{\left(\frac{v'u}{v},v'\right)}$. We can write $\frac{v'u}{v} = au'+bv'$ for $(a,b) \in \F_p^{\times} \times \F_p$, so that, by Lemma \ref{BK4S-elementary}, $I_{\mathcal{B}'}=I_{(u',v')}=I_{(u'+ba^{-1}v',v')} = \underline{a}^{-1}I_{(au'+bv',v')} = \underline{a}^{-1}I_{\mathcal{B}}$. 
}

\medskip

\prop[BK4S-qexp-Itors-complex]{Let $(\alpha,\beta)$ be the adapted basis of $\F_{p^2}$. Let $a \in \F_p^{\times}$ and $M \in \SL{\F_p}$, write $\begin{pmatrix}\alpha'\\\beta'\end{pmatrix}=M^{-1}\begin{pmatrix}\alpha\\\beta\end{pmatrix}$. Then, for any $f' \in G_{\Q}\cdot f$, $a \in \F_p^{\times}$ and $i \in \{0,1,2\}$, one has 
\[(M\Delta_{a,1})^{\ast}\left[\left(f'\left(\frac{\tau}{p}\right)\otimes \Delta_{a^{-1},1}\right) \mid \sum_{r=0}^{\frac{p+1}{3}-1}{\Delta_{(\det{g})^{-3r},1}\,g^{3r}\zeta_3^i}\right] = \left(\sum_{n \geq 1}{a_n(f')E_ne^{\frac{2i\pi n\tau}{p}}}\right) \otimes\Delta_{1,1},\] where
\begin{itemize}[noitemsep,label=\tiny$\bullet$]
\item $E_n \in \Z[\zeta_p]$ for every $n \geq 1$,
\item If $p \nmid n$, then $E_n=\underline{a^{-1}n}(\Sigma_{(\alpha',\beta'),i})$,
\item If $p \mid n$, one has $E_n=-\frac{2(p+1)}{3}$ if $\beta' \in  \F_{p^2}^{\times 3}j^i$, and $E_n=\frac{p+1}{3}$ otherwise.
\end{itemize}
}

\demo{First, note that 
\begin{align*}
(M\Delta_{a,1})^{\ast}&\left[\left(f'\left(\frac{\tau}{p}\right)\otimes \Delta_{a^{-1},1}\right) \mid \sum_{r=0}^{\frac{p+1}{3}-1}{\Delta_{(\det{g})^{-3r},1}\,g^{3r}\zeta_3^i}\right] \\
&= \left(\sum_{r=0}^{\frac{p+1}{3}-1}{f'\left(\frac{\tau}{p}\right) \mid \Delta_{(a(\det{g})^{3r})^{-1},1}\,g^{3r}\zeta_3^iM\Delta_{a,1}}\right) \otimes \Delta_{1,1}.\end{align*}

For $0 \leq r < \frac{p+1}{3}$, the classes in $\mathbb{P}^1(\F_p)$ of the first rows of the $g^{3r}\zeta_3^i$ are the same as the classes of the $(u,v) \in \mathbb{P}^1(\F_p)$ such that $u\alpha+v\beta \in \F_{p^2}^{\times 3}j^i$. Now, if $g^{3r}\zeta_3^i = z\begin{pmatrix} u & v\\\ast &\ast\end{pmatrix}$, then the first row of $\Delta_{(a(\det{g})^{3r})^{-1},1}\,g^{3r}\zeta_3^sM\Delta_{a,1}$ is proportional to $(u',v') := (a(u[M]_{1,1}+v[M]_{2,1}),u[M]_{1,2}+v[M]_{2,2})$, which are exactly the coordinates of $u\alpha+v\beta$ in the basis $(a^{-1}\alpha',\beta')$.  

By Lemma \ref{q-exp-s+nsplit}, when $p \nmid n$, we obtain the expression $E_n=\underline{n}\Sigma_{(\alpha'/a,\beta),i}=\underline{na^{-1}}\Sigma_{(\alpha',\beta'),i}$ by Lemma \ref{BK4S-elementary}. When $p \mid n$, Lemma \ref{q-exp-s+nsplit} implies that the $n$-th coefficient in the $q$-expansion of $f'\left(\frac{\tau}{p}\right) \mid \Delta_{(a(\det{g})^{3r})^{-1},1}\,g^{3r}\zeta_3^iM\Delta_{a,1}$ is $a_n(f')$ if $u' \neq 0$, and $-pa_n(f')$ if $u'=0$. By the above, $u'=0$ occurs as a term for $g^{3r}\zeta_3^iM$ (for some $0 \leq r < \frac{p+1}{3}$) if and only if $\beta' \in \F_{p^2}^{\times 3}j^i$. The conclusion follows.}

\medskip

The following consequence of this computation is enough to finish the proof of Proposition \ref{formal-immersion-application}, by Corollary \ref{step3}. For this proof, we need the precise definition of the algebraic $q$-expansion given in Section \ref{cuspidal-subscheme}, as well as its link to the complex uniformization spelled out in Proposition \ref{q-expansion-vs-analytic}.

\medskip

\prop[formal-immersion-at-nonexceptional]{Let $\ell \equiv \pm 1\pmod{p}$ be a non-exceptional prime number and $\mathfrak{l} \subset \OO_F$ be a prime ideal with residue characteristic $\ell$. Then, for every $c \in X(p,p)_{\OO_{F,\mathfrak{l}}}$ lying in the special fibre of the cuspidal subscheme, there exists a differential $\omega \in H^0(X(p,p)_{\OO_{F,\mathfrak{l}}},\Omega^1_{X(p,p)_{\OO_{F,\mathfrak{l}}}/\OO_{F,\mathfrak{l}}})$ such that $\mathfrak{I}\cdot \omega=0$, and $\omega_c \notin \mathfrak{m}_{X(p,p)_{\OO_{F,\mathfrak{l}}},c}\,\Omega^1_{X(p,p)_{\OO_{F,\mathfrak{l}}}/\OO_{F,\mathfrak{l}},c}$. }

\demo{By Proposition \ref{qexp-vs-formal-immersion} and Corollary \ref{cuspidal-subscheme-gammaN}, it is enough to show the following: for any $\Lambda \in \mrm{Surj}_p$, there is a cusp datum $(M,\Lambda)$ such that, for any maximal ideal $\lambda$ of the finite \'etale $\OO_{F,\mathfrak{l}}$-algebra $R := \Z[1/p,\zeta_p] \otimes \OO_{F,\mathfrak{l}}$, there exists a differential $\omega \in H^0(X(p,p)_{\OO_{F,\mathfrak{l}}},\Omega^1_{X(p,p)_{\OO_{F,\mathfrak{l}}}/\OO_{F,\mathfrak{l}}})$ such that $\mathfrak{I} \cdot \omega = 0$ and the first coefficient of $C_{(M,\Lambda)}^{\ast}\omega \in R \otimes (\Z[1/p])[[q^{1/p}]]$ is not contained in $\lambda$.

Let $\Lambda \in \mrm{Surj}_p$. Let $\gamma \in \SL{\F_p}$ be such that $\Lambda = \Lambda_0 \circ \gamma^T$, where $\Lambda_0: (x,y) \in \F_p^{\oplus 2} \mapsto x \in \F_p$. Let $(M,\Lambda)$ be a cusp datum such that $C_{(M,\Lambda)} = \gamma\circ C_{(I_2,\Lambda_0)}$ by Lemma \ref{qexp-plus-gl2}. Thus, for any differential $\omega \in H^0(X(p,p)_{\C},\Omega^1_{X(p,p)_{\C}/\C})$, one has $C_{(M,\Lambda)}^{\ast}\omega = C_{(I_2,\Lambda_0)}^{\ast}(\gamma^{\ast}\omega)$. 

Since $\gamma^{\ast}$ is an automorphism of $H^0(X(p,p)_{\OO_{F,\mathfrak{l}}},\Omega^1_{X(p,p)_{\OO_{F,\mathfrak{l}}}/\OO_{F,\mathfrak{l}}})$, it is enough to prove the following statement by Proposition \ref{qexp-coefficients}: for any maximal ideal $\lambda$ of $\Z[\zeta_p] \otimes \OO_{F,\mathfrak{l}}$, there exists $\omega \in H^0(X(p,p)_{\C},\Omega^1_{X(p,p)_{\C}/\C})$ such that $\mathfrak{I} \cdot \omega = 0$, the $q$-expansion of $\gamma^{\ast}\omega$ is contained in $R \otimes (\Z[1/p])[[q^{1/p}]]$, and its first coefficient is not contained in $\lambda$. 

Since the automorphisms of $\Z[\zeta_p]$ act transitively on the set of maximal ideals of $R$, it is enough to prove the following statement by Lemma \ref{with-delta-c}: let $\lambda_0 \subset \Z[\zeta_p] \otimes \OO_{F,\mathfrak{l}}$ be the kernel of the morphism $\zeta_p \otimes z \in \Z[\zeta_p] \otimes \OO_{F,\mathfrak{l}} \mapsto ze^{-\frac{2i\pi}{p}} \in \C$ of $\OO_F$-algebras. Then, for every $a \in \F_p^{\times}$, there exists $\omega \in H^0(X(p,p)_{\C},\Omega^1_{X(p,p)_{\C}/\C})$ such that 
\[\mathfrak{I} \cdot \omega = 0,\quad C_{(I_2,\Lambda_0)}^{\ast}\left(\Delta_{1,a}^{\ast}\gamma^{\ast}\omega\right) = C_{(I_2,\Lambda_0)}^{\ast}\left[(\gamma\Delta_{1,a})^{\ast}\omega\right] \in (R\backslash \lambda_0) + q^{1/p}(R \otimes (\Z[1/p])[[q^{1/p}]]).\] 

We take \[\omega_i = \iota_{p,p}^{-1}\left(\left(f'\left(\frac{\tau}{p}\right)\otimes \Delta_{a^{-1},1}\right) \mid \sum_{r=0}^{\frac{p+1}{3}-1}{\Delta_{(\det{g})^{-3r},1}\,g^{3r}\zeta_3^i}\right)\] for each $i \in \{0,1,2\}$. 

The differential $\omega_i$ is annihilated by $\mathfrak{I}$ by Proposition \ref{complex-differentials-vs-T1N}. By Proposition \ref{BK4S-qexp-Itors-complex} and Corollary \ref{full-analytic-qexp-at-infty}, one has \[C_{(I_2,\Lambda_0)}^{\ast}\left[(\gamma\Delta_{1,a})^{\ast}\omega_i\right] = C_{(I_2,\Lambda_0)}^{\ast}\left[(\gamma\Delta_{a^{-1},a}\Delta_{a,1})^{\ast}\omega_i\right] \in R \otimes (\Z[1/p])[[q^{1/p}]].\] 

By Propositions \ref{BK4S-qexp-Itors-complex} and \ref{q-expansion-vs-analytic}, the image in $(\Z[\zeta_p] \otimes \OO_{F,\mathfrak{l}})_{\lambda_0}$ of the first coefficient of $C_{(I_2,\Lambda_0)}^{\ast}\left[(\gamma\Delta_{1,a})^{\ast}\omega_i\right]$ is $\Sigma_{(\alpha,\beta),i}$ for a certain $\F_p$-basis $(\alpha,\beta)$ of $\F_{p^2}$. Since $\ell$ is not exceptional, one of the $\Sigma_{(\alpha,\beta),i}$ is not contained in $\lambda_0$, so the corresponding $\omega_i$ suits the conditions. }

\bigskip

\section{Diophantine consequences for the elliptic curve $y^2=x^3-p$}
\label{diophantine-consequence-y2x3p}

In this section, $p \geq 11$ is a prime such that $p \equiv 5\pmod{9}$. 

\lem[eqn-minimal]{The Weierstrass equation $y^2=x^3-p$ over $\Z$ is minimal with discriminant $-2^4 \cdot 3^3\cdot p^2=-432p^2$ is minimal. It defines an elliptic curve $E_p/\Q$ with complex multiplication by $\Z[e^{2i\pi/3}]$ and conductor $108p^2$ if $p \equiv -1\pmod{4}$ and $432p^2$ if $p \equiv 1\pmod{4}$.}

\demo{By \cite[Chapter III.1]{AEC1}, the discriminant of the Weierstrass equation is \[\Delta = -\frac{(-216)^2(-4p)^2}{1728} = -2^{-6}\cdot 2^6 \cdot 3^{-3} \cdot 3^6 \cdot 2^4 \cdot p^2 = -2^4 \cdot 3^3 \cdot p^2=-432p^2.\]
Since it is not divisible by a twelfth power, the Weierstrass equation is minimal by \cite[Remark VII.1.1]{AEC1}. The elliptic curve $E_p$ has complex multiplication by $\Z[e^{2i\pi/3}]$, since $(x,y) \longmapsto (xe^{2i\pi/3},y)$ is an automorphism of $(E_p)_{\Q(e^{2i\pi/3})}$ of order $3$. We compute that $j(E_p)=0$, hence $E_p$ has potentially good reduction everywhere by \cite[Proposition VII.5.5]{AEC1}, thus multiplicative reduction nowhere by \cite[Proposition VII.5.4]{AEC1}, therefore $2,3,p$ are primes of additive reduction for $E_p$.  
In particular, $36p^2$ divides the conductor $N_p$ of $E_p$ and $v_p(N_p)=2$ by \cite[Theorem IV.10.2]{AEC2}. 

To compute the conductor exponents at $2$ and $3$, we apply Tate's algorithm \cite[Theorem IV.9.4]{AEC2}. At $3$, we stop of Step 3 of \emph{loc.cit.} and the reduction type is II with conductor exponent $v_3(\Delta)=3$. At $2$, we stop at Step 3 if $p \equiv 1\pmod{4}$ -- in which case the reduction is Type II and the conductor exponent is $v_2(\Delta)=4$ -- and at Step 5 if $p \equiv -1\pmod{4}$ -- in which case the reduction is of type IV and the conductor exponent is $v_2(\Delta)-2=2$.
}

\medskip

\prop[p-torsion-of-Ep]{The image of $G_{\Q}$ in $\mrm{Aut}(E_p[p])$ is exactly the normalizer $N$ of a nonsplit Cartan subgroup $C$. Moreover, the image of an inertia subgroup $I_p$ at $p$ is exactly $3C$. }

\demo{The image of $\Z[e^{2i\pi/3}]$ in $\mrm{End}(E_p[p]) \simeq \mathcal{M}_2(\F_p)$ is exactly $C \cup \{0\}$ for some nonsplit Cartan subgroup $C$, because $\F_p$ does not contain a primitive third root of unity. 

Let $f: G_{\Q} \rar \mrm{Aut}(E_p[p])$ be the morphism. For any prime number $q$, let $I_q \leq G_{\Q}$ be an inertia subgroup at $q$. Now, $f(G_{\Q(\mu_3)})$ commutes with the image of $\Z[e^{2i\pi/3}] = \mrm{End}_{\Q[e^{2i\pi/3}]}(E_p)$ in $\mrm{End}(E_p[p])$, so $f(G_{\Q(\mu_3)}) \subset C$. Hence $f(G_{\Q})$ is contained in the normalizer $N$ of $C$. 

By Lemma \ref{new-is-semistable}, the conductor exponent of the $\F_p[G_{\Q}]$-module $E_p[p]$ at $3$ is the same as the conductor of $E_p$, hence is $3$. In particular, the action of the inertia at $3$ on the $p$-torsion is wildly ramified, so $f(I_3)$ (hence $f(G_{\Q})$) contains an element of order $3$. By Lemma \ref{cartan-is-potentially-good}, $f(I_p) \supset 3C$. By Proposition \ref{is-usually-cartan}, $f(G_{\Q}) \not\subset C$. In particular, $f(G_{\Q})/3C$ is a subgroup of $N/3C \simeq \mathfrak{S}_3$ containing an element of order $3$ and not equal to the cyclic subgroup of order $3$. Hence $f(G_{\Q})=N$. 

Since $f(I_p) \subset C$, $f_{|I_p}$ is isomorphic to the representation $\begin{pmatrix}\chi^{k-1}& 0\\0 & \chi^{p(k-1)}\end{pmatrix}: I_p \rar \GL{\overline{\F_p}}$, where $k$ is the Serre weight of $f$ (see \cite[\S 2.2]{Conj-Serre}) and $\chi: I_p \rar \F_{p^2}^{\times}$ is the fundamental character of inertia of level $2$, which is a surjection $I_p \rar \F_{p^2}^{\times}$. By \cite[Th\'eor\`eme 1]{Kraus-Thesis}, the Serre weight of the representation $f$ is $\frac{p^2+6p+5}{6} = \frac{(p+1)(p+5)}{6}$. We can write $p=18t+5$, then the Serre weight of $f$ is $\frac{(18t+6)(18t+10)}{6}= (3t+1)(18t+10)\equiv 1 \pmod{3}$. Therefore, $\chi^{k-1}(I_p)$ is contained in the cubes of $\F_{p^2}^{\times}$, whence $f(I_p) \subset 3C$.  
} 

\medskip

\prop[weight-one-mf]{Let $C \leq \GL{\F_p}$ be a nonsplit Cartan subgroup with normalizer $N$. Let $(P,Q)$ be a basis of $E_p[p](\Qbar)$ such that the morphism $\rho: G_{\Q} \rar \GL{\F_p}$ satisfying \[\rho(\sigma)\begin{pmatrix}\sigma(P)\\\sigma(Q)\end{pmatrix} = \begin{pmatrix}P\\Q\end{pmatrix}\] has image $N$. Let $\psi: C \rar F^{\times}$ be a character of order $3$, and $g$ be the modular form attached to the odd, irreducible, Artin representation $\left(\mrm{Ind}_C^N{\psi}\right)\circ \rho$. Then $g$ is given by the LMFDB label \cite{lmfdb} $108.1.c.a$\footnote{Its $q$-expansion is $q-q^{7}-q^{13}-q^{19}+q^{25}+2q^{31}+\cdots$}. }

\demo{We already saw that $\left(\mrm{Ind}_C^N{\psi}\right)\circ \rho$ was attached to a weight one modular form $g$ with complex multiplication by $\Q(\sqrt{-3})$. The nebentypus of $g$ is exactly the character of $\Q(\sqrt{-3})$, since $\det\left[\mrm{Ind}_C^N{\psi}\right]$ is the unique nontrivial character of $N/C$. Since the traces of the representation $\mrm{Ind}_C^N{\psi}$ of $N$ are $2, 0, -1$, $g$ has rational coefficients by strong multiplicity one \cite[Theorem 4.6.19]{Miyake}. 

By Proposition \ref{p-torsion-of-Ep}, $\rho_g := \left(\mrm{Ind}_C^N{\psi}\right)\circ \rho$ is unramified at $p$, so by \cite[Th\'eor\`eme 4.6]{Del-Ser} the conductor of $g$ is exactly the conductor of the Galois representation $\rho_g$, which is only ramified at $2$ and $3$. Let us now determine the conductor exponents. 

We start at $3$. Since $\psi$ realizes an isomorphism of the $3$-torsion subgroup of $C$ on $\psi(C)$, the reduction modulo $p$ of $\rho_g$ has the same restriction to the wild inertia subgroup at $3$ as $E_p[p]$. Therefore, the Swan conductor exponent of $E_p[p]$ at $3$ and the Swan conductor exponent of $g$ at $3$ are equal. Moreover, $E_p$ has additive reduction at $3$, so the tame part of the conductor exponent of $E_p$ is $2$. Since the subgroup $C$ has no invariants in the representation $\mrm{Ind}_C^N{\psi}$, the tame part of the conductor exponent of $g$ at $3$ is also equal to $2$, so the conductor exponent of $g$ at $3$ is $3$. 

Now, we discuss the prime $2$. Because $\Q(\sqrt{-3})/\Q$ is unramified at $2$, the image under $\rho$ of an inertia subgroup $I_2$ at $2$ is contained in $C$. This implies that $\rho_g$ is tamely ramified at $2$, and that its conductor exponent is zero if $\rho(I_2)$ does not meet the $3$-torsion subgroup in $C$, and $2$ if $\rho(I_2)$ contains the $3$-torsion subgroup of $C$. Since $E_p$ has potentially good reduction at $2$ and $C$ is cyclic, by \cite[Corollary 1 to Theorem 2]{GoodRed}, the conductor exponent of $\rho_g$ is zero (resp. $2$) if $E_p$ acquires good reduction at $2$ on an extension of $\Q_2$ with ramification index coprime to $3$ (resp. divisible by $3$).  

If $L/\Q_2$ is the extension on which $E_p$ acquires good reduction, then the normalized valuation over $L$ of the minimal discriminant of $E$ over $\Q_2$ is divisible by $12$ by the formulas of \cite[Chapter III.1]{AEC1}, which implies that $v_L(2) \in 3\Z$, so the conductor of $g$ is $108$. 

By the LMFDB \cite{lmfdb}, there is exactly one newform with rational coefficients and complex multiplication by $\Q(\sqrt{-3})$ in $\mathcal{S}_1(\Gamma_1(108))$, whence the conclusion.

}

\prop[diophantine-consequence-asymptotic]{There exists a constant $p_0 \geq 11$ such that, for any prime $p \geq p_0$ congruent to $5$ modulo $9$, for any elliptic curve $E'$ such that $E'[p](\Qbar)$ is isomorphic, as a $G_{\Q}$-module, to the $p$-torsion of the elliptic curve $E_p$. Then the conductor of $E'$ is of the form $108mp^2L$, where $L$ is a square-free product of exceptional primes for $p$, and $m=1$ if $p \equiv 3 \pmod{4}$ and $m=4$ if $p \equiv 1\pmod{4}$. }

\demo{By Theorem \ref{formal-immersion-application}, the conclusion of the Proposition holds if there is $f \in \mathcal{S}_2(\Gamma_0(p))$ such that $L(f \times g,1) \neq 0$, where $g$ is the newform with LMFDB label \cite{lmfdb} $108.1.c.a$ as in Proposition \ref{weight-one-mf}. That there exists such an $f$ is given by a result of analytic number theory due to Michel \cite[Appendix, Theorem 0.1]{Michel}, since in this case the sign of the functional equation for $L(f \times g,s)$ is $+1$.  }

\medskip

\rem{For small enough $p$, one can also test the non-vanishing of the Rankin-Selberg $L$-value required by Theorem \ref{formal-immersion-application} numerically. Thanks to the help of Pascal Molin, I was able to check in PARI/GP that the conclusion of Theorem \ref{formal-immersion-application} held when $p=23$. In this case, we can also check directly that there are no exceptional primes, so the conductor of $E'$ is $108p^2=57132$. By \cite[Theorem 1.3]{Cremona-Freitas}, the only elliptic curves of conductor $108p^2$ congruent modulo $p$ to $E_{p}$ are isogenous to $E_{p}$. Therefore\footnote{Further computations after the submission of this thesis showed that there were no exceptional primes for $p$ when $p \equiv 2, 5\pmod{9}$ and $11 \leq p \leq 5000$, and that the conclusion of Theorem \ref{formal-immersion-application} held for $p$ and the elliptic curve $E_p$ when $23 \leq p \leq 150$ and $p \equiv 5\pmod{9}$. Thus, any elliptic curve congruent modulo $p \in \{23, 41, 59, 113, 131, 149\}$ to $E_p$ has the same conductor as $E_p$ and integral $j$-invariant. Finally, \cite[Theorem 1.3]{Cremona-Freitas} also applies for $p=59$ and $p=41$ (the latter after a suitable quadratic twist).}, for $p=23$, any elliptic curve congruent modulo $p$ to $E_{p}$ is isogenous to $E_{p}$. }

\newpage

%% file: euler-0.tex
\chapter{On obstructions to the Euler system method for Rankin-Selberg convolutions}
\label{obstructions-euler}

\section{Introduction}
Let $k,N \geq 1$ be integers. We denote by $\mathcal{S}_k(\Gamma_1(N))$ the space of cuspidal holomorphic modular forms with level $\Gamma_1(N)$. Let $f\in \mathcal{S}_k(\Gamma_1(N))$ be a cuspidal newform with character $\varepsilon_f$, and $L$ be a number field containing its coefficients. Given a prime ideal $\mathfrak{p} \subset \OO_L$, we can attach to $f$ a free $\OO_{L_{\mathfrak{p}}}$-module $T_f$ of rank two endowed with a continuous linear action of $G_{\Q}:=\mrm{Gal}(\Qbar/\Q)$. Works of Momose \cite{Momose} and Ribet \cite{Ribet2} show that, when $f$ has no complex multiplication and $k \geq 2$, the image of $G_{\Q}$ in $\Aut{T_f} \cong \GL{\OO_{L_{\mathfrak{p}}}}$ contains a conjugate of $\SL{\Z_p}$ for all but finitely many $\mathfrak{p}$. More recently, Loeffler \cite{bigimage} proved an \emph{adelic} open image theorem for such newforms, the analog of that of Serre \cite{Serre-image-ouverte} for elliptic curves without complex multiplication. 

Let $g \in \mathcal{S}_l(\Gamma_1(M))$ be another newform with character $\varepsilon_g$; after enlarging $L$ if necessary, we assume that its Fourier coefficients are contained in $L$. As above, we attach to $g$ a free $\OO_{L_{\mathfrak{p}}}$-module of rank two $T_g$ endowed with a continuous linear action of $G_{\Q}$. Let $T := T_f \otimes_{\OO_{L_{\mathfrak{p}}}} T_g$. The other goal of \cite{bigimage} was to study the image of $G_{\Q}$ in $\Aut{T}$. In particular, Loeffler proves that, when the following conditions are satisfied:

\begin{itemize}[noitemsep,label=\tiny$\bullet$]
\item $k,l \geq 2$,
\item neither $f$ nor $g$ has complex multiplication,
\item $f$ and $g$ are not twists one of another,
\item $\varepsilon_f\varepsilon_g \neq 1$,
\end{itemize}
 then, for all but finitely many $\mathfrak{p}$, there exists a $\sigma \in G_{\Q(\mu_{p^{\infty}})} := \mrm{Gal}(\Qbar/\Q(\mu_{p^{\infty}}))$ such that $T/(\sigma-1)T$ is free of rank one. Loeffler also discusses the case of $k \geq 2, l=1$ and asks the following question: 

\qnn[loefflerqn]{(Loeffler, see \cite[Remark 4.4.2]{bigimage}) Assume that $k \geq 2$, $l=1$, $f$ does not have complex multiplication, and $\varepsilon_f\varepsilon_g \neq 1$. Then does there exist, for all but finitely many $\mathfrak{p}$, $\sigma \in G_{\Q(\mu_{p^{\infty}})}$ such that $T/(\sigma-1)T$ is free of rank one?}

Loeffler proves in \cite[Proposition 4.4.1]{bigimage} that the answer is positive when $N$ and $M$ are coprime. 

Here is a few of the weaker sufficient conditions one can use to answer the question in the affirmative.

\theoi[weakerconds]{The answer to Question \ref{loefflerqn} is positive when any of the following conditions hold:
\begin{itemize}[noitemsep,label=$-$]
\item For any (primitive) Dirichlet character $\chi$ with conductor $C$ such that $f \otimes \chi$ is Galois-conjugate to $f$, $C$ is coprime to $M$ (in particular when $f$ has no inner twist). 
\item Any (primitive) Dirichlet character $\chi$ such that $f \otimes \chi$ is Galois-conjugate to $f$ is even.
\item $\varepsilon_f^2=\varepsilon_g^2=1$ and $f \otimes \varepsilon_g$ is not Galois-conjugate to $f$.
\item The group $\varepsilon_g\left(\bigcap_{\chi}{\ker{\chi}}\right)$ contains an element of order $2$ but no element of order $4$, where $\chi$ runs through (primitive) Dirichlet characters such that $f \otimes \chi$ is Galois-conjugate to $f$. 
\end{itemize}}

As a consequence, we obtain the following instance of the Bloch-Kato conjecture, which extends \cite[Theorem 11.7.4]{KLZ15} by removing some assumptions.

Recall the following definition:

\defi{Let $E$ be an elliptic curve over a number field $K$, and $p$ be a prime number. The $p^{\infty}$-Selmer group of $E$ is the subgroup 
\[\mrm{Sel}_{p^{\infty}}(E/K) = \ker\left[H^1(K,E[p^{\infty}]) \rar \prod_v{H^1(K_v,E(\overline{F_v}))}\right],\] where $v$ runs over the finite places of $K$ and the map is induced by localizing at the place $v$ and the inclusion $E[p^{\infty}] \subset E$.  }

\theoi[twistedbsd]{Let $E$ be an elliptic curve over $\Q$ without complex multiplication and $\rho$ be an odd two-dimensional irreducible representation of $G_{\Q}$, with splitting field $K$ and coefficients in some number field $L_0$ with ring of integers $\OO$. Let $\mathfrak{p}$ be a prime ideal of $\OO$ with residue characteristic $p$. Suppose that the following assumptions hold:
\begin{enumerate}[noitemsep]
\item $p$ is coprime to $30N_{\rho}N_E$, where $N_E$ (resp. $N_{\rho}$) is the conductor of $E$ (resp. $\rho$), 
\item if $\rho$ has dihedral projective image $D$, then the maximal cyclic subgroup of $D$ is not a $p$-group,  
\item the map $G_{\Q} \rar \Aut{T_pE}$ given by the Galois action on the $p$-adic Tate module of $E$ is onto,  
\item $E$ is ordinary at $p$ and $\rho(\Fr_p)$ has distinct eigenvalues mod $\mathfrak{p}$.
\end{enumerate} 
If $L(E,\rho,1) \neq 0$, then the group $\mrm{Hom}_{\Z_p[\mrm{Gal(F/\Q)}]}(\rho,\mrm{Sel}_{p^{\infty}}(E/F))$ is finite. 
}

We show, less constructively: 

\theoi[weakerqn]{If one replaces ``for all but finitely many $\mathfrak{p}$'' with ``for infinitely many $\mathfrak{p}$'' in Question \ref{loefflerqn}, then the answer is positive.}

We also show that the general answer to Question \ref{loefflerqn} is negative. Among other counter-examples, we find:

\theoi[INF-CEX]{(See Section \ref{sect-theoD}) Let $f \in \mathcal{S}_k(\Gamma_1(N))$ be a newform with character $\varepsilon_f$ without complex multiplication. Suppose that one of the following conditions holds:
\begin{itemize}[noitemsep]
\item $k$ is even and there exists a primitive Dirichlet character $\chi$ such that $f \otimes \chi$ is Galois-conjugate to $\chi$. 
\item $k \geq 3$ is odd, and there exist two Dirichlet characters $\chi_1,\chi_2$ such that each $f \otimes \chi_i$ is Galois-conjugate to $f$, and none of the restrictions to $\ker{\varepsilon_f}$ of $\chi_1,\chi_2,\chi_1\chi_2$ are trivial. 
\end{itemize} Then there exist infinitely many dihedral newforms $g \in \mathcal{S}_1(\Gamma_1(M))$ such that the answer to Question \ref{loefflerqn} for $(f,g)$ is negative.}

\theoi[counterex]{There exist newforms $f \in \mathcal{S}_2(\Gamma_0(63)), g \in \mathcal{S}_1(\Gamma_1(1452))$ with coefficients respectively in $\Q(\sqrt{3}),\Q(\sqrt{-3})$ satisfying the following property: for all but finitely many prime ideals $\mathfrak{p}$ of $L=\Q(e^{2i\pi/12})$ with residue characteristic $p \equiv 5,7 \pmod{12}$, there is no $\sigma \in G_{\Q(\mu_{p^{\infty}})}$ such that $\ker\left(\sigma-\mrm{id}\mid T_{f,\mathfrak{p}} \otimes_{\OO_{L_{\mathfrak{p}}}} T_{g,\mathfrak{p}}\right) \otimes_{\OO_{L_{\mathfrak{p}}}} L_{\mathfrak{p}}$ is a $L_{\mathfrak{p}}$-line.}

\rems{\begin{enumerate}[noitemsep,label=(\roman*)] 
\item The statement of Theorem \ref{twistedbsd} can be generalized by replacing elliptic curves with, for instance, abelian varieties attached to newforms of weight $2$, with only minimal changes to the proof. One would need to adapt the big image hypothesis, for instance by ensuring that the conclusion of Theorem \ref{weakerconds} holds. 
\item As explained in \cite[Remark 11.7.5(i)]{KLZ15}, applying Theorem \ref{twistedbsd} at a single prime shows that $\mrm{Hom}_{\OO[G_{\Q}]}(\rho,E(K) \otimes \OO)$ is zero. When generalizing Theorem \ref{twistedbsd} to more general modular forms, Theorem \ref{weakerqn} would imply that, provided that the special value is not zero, the relevant Mordell-Weil group is torsion. The point of Theorem \ref{twistedbsd} is to prove that, in addition, the $\rho$-part of $\Sha(E/K)[p^{\infty}]$ is finite for a positive density of such primes $p$.  
\item When $f \in \mathcal{S}_k(\Gamma_1(N))$, and $g \in \mathcal{S}_l(\Gamma_1(M))$ are newforms with characters $\varepsilon_f,\varepsilon_g$ respectively, then the Dirichlet character $\varepsilon_f\varepsilon_g$ has the same parity as $k+l$. In particular, when $k$ is even and $l=1$, $\varepsilon_f\varepsilon_g$ is odd, so cannot be trivial. The counter-examples in Theorem \ref{INF-CEX} and \ref{counterex} are thus non trivial. 
\item In the counter-example of Theorem \ref{counterex}, $g$ has complex multiplication, $f$ has even weight and an odd inner twist. While the answer to Question \ref{loefflerqn} is negative, this does not follow from the proof of Theorem \ref{INF-CEX}. Further counter-examples are given in Section \ref{sect-counterex}. 
\end{enumerate}}

The existence of a $\sigma$ as in Question \ref{loefflerqn} is important, because it seems necessary to apply the theory of Euler systems (see for instance \cite{MR-ES}), which aims to bound Selmer groups assuming the nonvanishing of certain $p$-adic $L$-functions. For the Galois-modules that we consider in this article, the Euler system would typically be the system of Beilinson-Flach elements constructed in \cite{KLZ15}. We do not see any reason why, for all our couples $(f,g)$ of counter-examples, the central archimedian Rankin-Selberg $L$-value $L(f \otimes g,1)$ should vanish. 

Hence, consider a couple $(f,g)$ that is a counter-example to Question \ref{loefflerqn} and assume that $L(f \otimes g,1) \neq 0$. Should one expect that, as in \cite[Theorem 11.7.3]{KLZ15}, for almost all primes $p$ ordinary for $f$ and $g$ and such that $f,g$ are $p$-distinguished\footnote{In the situation of Theorem \ref{counterex}, this last property holds for example when $p \equiv 5 \pmod{12}$, or when $p \equiv 7 \pmod{12}$ and $\left(\frac{p}{11}\right)=-1$.}, the Nekov\'a\v{r} cohomology group $\tilde{H}^2\left(\OO_L\left[\frac{1}{N_fN_gp}\right],T_f \otimes T_g;\Delta^{BK}\right)$ is finite\footnote{The statement of the theorem in \emph{loc. cit.} seems to contain a slight typo; one should (probably) read $M_{\OO}(f \otimes g)^{\ast}$ instead of $M_{\OO}(f \otimes g)$.}? However, how to go about proving such a statement without the theory of Euler systems does not seem clear at all.  \\

\textbf{Acknowledgements} \, I would like thank Lo\"ic Merel for encouraging me to work on this project. His careful reading and many comments helped me vastly improve the original version of this text. I would also like to thank David Loeffler for his help regarding the application of these results to the twisted Birch and Swinnerton-Dyer conjecture.

\section{The big image theorem for one modular form}
\label{sect-setup}

We first recall well-known results on twists and Galois conjugates of modular forms. This lets us define inner twists, who are crucial to stating the large image theorem for one modular form.  

The results on twists are mostly taken from \cite{AL78}, the lemmas on inner twists come from \cite{Momose}, and the presentation of the results on Galois images of this section comes almost entirely from \cite[Section 2]{bigimage}.  

First, we recall the definition and well-known properties (see \cite{AL78,Miyake}) of twists of newforms by Dirichlet characters. 

\prop[twist-exists]{Let $f \in \mathcal{S}_k(\Gamma_1(N))$ be a newform with character $\varepsilon_f$ of conductor $M$. Let $\chi$ be a primitive Dirichlet character of conductor $C$. There exists a unique new, cuspidal eigenform $f \otimes \chi \in \mathcal{S}_2(\Gamma_1(L))$ such that, for every prime number $p\nmid NC$, $a_p(f \otimes \psi)=a_p(f)\psi(p)$. We say that $f\otimes \chi$ is the \emph{twist} of $f$ by $\chi$. Moreover: 
\begin{enumerate}[label=$($\alph*$)$,noitemsep]
\item \label{twist-character} The character of $f \otimes \chi$ is the primitive Dirichlet character associated to $\varepsilon_f\chi^2$.
\item \label{level-twist-basic} $L$ divides the least common multiple $P$ of $N,CM,C^2$.
\item \label{level-twist-coprime} For every prime number $q \nmid C$, $v_q(L)=v_q(N)$.
\item \label{level-twist-local} Let $q$ be a prime number, and write $\chi=\chi_q\chi'$, where the conductor of $\chi_q$ is a power of $q$ and that of $\chi'$ is coprime to $q$. Then $v_q(L)$ only depends on $f$ and $\chi_q$.
\end{enumerate}}

\rem{If $\psi$ comes from the primitive Dirichlet character $\psi_0$, we will also denote by $f \otimes \psi$ the newform $f \otimes \psi_0$. }

\prop{Let $f \in \mathcal{S}_k(\Gamma_1(N))$ and write its $q$-expansion $f(\tau)=\sum_{n \geq 1}{a_n(f)q^n}$. Let $\sigma$ be any automorphism of $\C$. Then $\sigma(f)(\tau) = \sum_{n \geq 1}{\sigma(a_n(f))q^n} \in \mathcal{S}_k(\Gamma_1(N))$. Moreover:
\begin{enumerate}[label=$($\alph*$)$,noitemsep]
\item If $f$ has character $\varepsilon_f$, then $\sigma(f)$ has character $\sigma\circ \varepsilon_f$. 
\item If $f$ is a newform, $\sigma(f)$ is a newform. 
\item If $f$ is a newform, there is a number field containing all its coefficients. 
\end{enumerate} }

Let $f \in \mathcal{S}_k(\Gamma_1(N))$ be a newform with primitive character $\varepsilon_f$. Write its $q$-expansion as $f(\tau) = \sum_{n \geq 1}{a_n(f)q^n}$, and $L=\Q(a_n(f)\mid n \geq 1)$. 

Let now $\gamma: L \rar \C$ be a field homomorphism and $\chi$ be a primitive Dirichlet character. We say that $f$ has an \emph{inner twist} by $(\gamma,\chi)$ if the newforms $\gamma(f)$ and $f \otimes \chi$ are equal. 

\lem[momose15]{In such a situation, one can write $\chi=\lambda\varepsilon_f^t$, where $t \in \Z$ and $\lambda$ is quadratic. Moreover, the conductor of $\chi$ divides $N$ and $\gamma(L) = L$.}

\demo{Apart from the bound on the conductor of $\chi$, this is \cite[Lemma 1.5 (i)]{Momose}. 

Let $C$ be the conductor of $\chi$, and let $p$ be a prime such that $v_p(C) > v_p(N)$. Then by \cite[Theorem 4.1]{AL78} and Proposition \ref{twist-exists}\ref{level-twist-local}, the $p$-adic valuation of the conductor of $f \otimes \chi$ is $v_p(N)$, and is also $2v_p(C)$, and we get a contradiction.%
}

Let $\Gamma_f$ denote the set of $(\gamma,\chi)$ by which $f$ has an inner twist. 

\prop[abelian-dicho]{The formula $(\alpha,\chi) \cdot (\beta,\psi) = (\alpha\circ \beta, \chi\alpha(\psi))$ endows $\Gamma_f$ with an abelian group law. Moreover, the group homomorphism \[\pi: (\gamma,\chi) \in \Gamma_f \longmapsto \gamma \in \Aut{L}\] satisfies exactly one of the following properties:
\begin{itemize}[label=\tiny$\bullet$]
\item $\pi$ is injective with abelian image. Then the fixed field of $\mrm{im}(\pi)$ (which we identify to $\Gamma_f$) is a field $F$ such that $L/F$ is abelian with Galois group $\mrm{im}(\pi) \cong \Gamma_f$. 
\item There exists some $(\mrm{id},\eta) \in \Gamma_f$ where $\eta$ is the character of a quadratic field $K$ which is only ramified at primes dividing $N$. If $k \geq 2$, then $K$ is imaginary. 
\end{itemize}}

\demo{Checking that we have indeed defined a group law on $\Gamma_f$ such that $\pi$ is a homomorphism is formal. 

Let $(\sigma,\alpha),(\tau,\beta) \in \Gamma_f$. Then we can write $\alpha=\alpha_2\varepsilon^a$, $\beta=\beta_2\varepsilon^b$, with $a,b \in \Z$ and $\alpha_2,\beta_2$ quadratic. Then $\tau(\alpha)\alpha^{-1}=(\tau(\varepsilon)\varepsilon^{-1})^a=(\beta^2)^a=\varepsilon^{2ab}$: by symmetry $\tau(\alpha)\beta=\alpha\sigma(\beta)$. Because $\Aut{\C}$ acts through an abelian quotient on $\Q(\varepsilon)$, a straightforward computation shows that $\Gamma_f$ is abelian. 

If $\pi$ is injective, its image is an abelian subgroup of $\Aut{L}$. It is clear that $L/F$ is finite separable. Since $\Auti{F}{L} \supset \mrm{im}(\pi)$, one has $[L:F] \geq |\mrm{im}(\pi)|$. Conversely, if $\alpha \in L$ is a primitive element over $F$, the monic polynomial $R(X)=\prod_{\gamma \in \mrm{im}(\pi)}{(X-\gamma(\alpha))} \in L[X]$ vanishes at $\alpha$ and is stable under $\mrm{im}(\pi)$, so lies in $F[X]$. Thus $[L:F] \leq [F(\alpha):F] \leq \deg{R} \leq |\mrm{im}(\pi)|$, which concludes. 

Suppose that $\pi$ is non-injective. Then there exists a nontrivial Dirichlet character $\chi$ such that $f \otimes \chi=f$. Thus $\varepsilon\chi^2=\varepsilon$, so $\chi$ is quadratic. The rest follows from \cite[(4.4, 4.5)]{Antwerp5-Ribet} in weight $k \geq 2$. 
}

In the latter case, we say that $f$ has \emph{complex} (resp. \emph{real}) \emph{multiplication by $K$} if $K$ is imaginary (resp. real). This is often shortened as ``CM'' (resp. ``RM''), and we will omit $K$ from the notation if it does not matter. 

If $f$ does not have CM, let $F$ be the subfield of $L$ fixed by the image of $\pi$, and let $H/\Q$ be the abelian field extension such that $G_H = \bigcap_{(\gamma,\chi) \in \Gamma_f}{\ker{\chi}}$, where $\chi$ is viewed as a character of $G_{\Q}$ through the cyclotomic character $G_{\Q} \rar \hat{\Z}^{\times}$. 

\medskip
\noindent

Let $\mathfrak{p} \subset \OO_L$ be a prime ideal with finite residue characteristic $p$. It is well-known (see \cite{Antwerp5-Ribet}) that there exists a group homomorphism $\rho_{f,\mathfrak{p}}: G_{\Q} \rar \GL{\OO_{L_{\mathfrak{p}}}}$ satisfying the following properties:

\begin{itemize}[noitemsep,label=\tiny$\bullet$]
\item $\rho_{f,\mathfrak{p}}$ is continuous, irreducible, and unramified away from $pN$.
\item $\det{\rho_{f,\mathfrak{p}}}=\varepsilon \omega_p^{k-1}$, where $\omega_p$ is the $p$-adic cyclotomic character. 
\item For every prime $\ell \nmid pN$, $\Tr{\rho_{f,\mathfrak{p}}(\Fr_{\ell})} = a_{\ell}(f)$, where $\Fr_{\ell}$ is the arithmetic Frobenius.  
\end{itemize}

\prop{(Momose \cite[Theorem 4.2]{Momose}, Ribet \cite[\S 3]{Ribet2}, Ghate--Gonz\'alez-Jim\'enez--Quer \cite[Corollary 4.7]{GGJQ}) 

Assume moreover that $p \nmid 2N$ is unramified in $L$ and $f$ does not have CM. Then $\rho_{f,\mathfrak{p}}(G_H)$ has a $\GL{\OO_{L_{\mathfrak{p}}}}$-conjugate contained in $\{M \in \GL{\OO_{F_{\mathfrak{p}}}},\, \det{M} \in \Z_p^{\times}\}$. 

Furthermore, for all but finitely many $\mathfrak{p}$, after conjugating $\rho_{f,\mathfrak{p}}$ by a suitable element of $\GL{\OO_{L_{\mathfrak{p}}}}$, the group homomorphism $(\rho_{f,\mathfrak{p}},\omega_p): G_H \rar \GL{\OO_{F_{\mathfrak{p}}}} \times_{\OO_{F_{\mathfrak{p}}}^{\times}} \Z_p^{\times}$ is onto (where the left map in the fibre product is the determinant, and the right one is $x \in \Z_p^{\times} \longmapsto x^{k-1} \in \OO_{F_{\mathfrak{p}}}^{\times}$). }

\demo{This is \cite[Theorem 2.2.2]{bigimage}. }

\section{Special elements in the image of the convolution of two modular forms}
\label{sect-theoA}

Let $p$ be an odd prime. Let $\OO$ be the ring of integers of a finite extension $K$ of $\Q_p$, with uniformizer $\varpi$ and residue field $k$. Let $T$ be a free $\OO$-module of finite rank, endowed with a continuous and $\OO$-linear action of $G_{\Q}$. Let $G$ be the image in $\Aut{T}$ of $G_{\Q(\mu_{p^{\infty}})}$. Consider the following statements:  

\begin{enumerate}[noitemsep]
\item[(N)] $-\mrm{id} \in G$. 
\item[(gI)] $T \otimes_{\OO} K$ is an irreducible $K[G]$-module. 
\item[(rI)] $T \otimes_{\OO} k$ is an irreducible $k[G]$-module.
\item[(wE)] There exists $u \in G$ such that $\ker(u-\mrm{id}) \otimes K$ has dimension one. 
\item[(sE)] There exists $u \in G$ such that $T/(u-1)T$ is free of rank one. 
\end{enumerate}

\defi{We say that $T$ has Galois image \emph{of Euler type} if it satisfies (N), (gI) and (wE).\\ $T$ has  \emph{Euler-adapted} Galois image if it satisfies (N), (rI) and (sE).}

Note that (rI) implies (gI) (since any reducible $K[G]$-submodule can be reduced in $k$), and that (sE) implies (wE). 

The statements (gI) and (wE) are the two parts of $\mrm{Hyp}(\Q(\mu_{p^{\infty}}), T \otimes K)$ in the notation of \cite[Section 4.1]{bigimage}. Similarly, (rI) and (sE) constitute in the same way $\mrm{Hyp}(\Q(\mu_{p^{\infty}}), T)$. The condition (N) is not mentioned explicitly in \cite{bigimage}, but is used in applications of this article to Euler systems, such as \cite[Section 11]{KLZ15}. Since it is not difficult to study, we include it here.

Let $k \geq 2$. Let $f \in \mathcal{S}_k(\Gamma_1(N_f))$, $g \in \mathcal{S}_1(\Gamma_1(N_g))$ be newforms with characters $\varepsilon_f$ and $\varepsilon_g$, such that $f$ has no complex multiplication. Let $L,F,H$ be the number fields associated to $f$ in the previous section, and $L_0 \subset \C$ be a number field containing $L$ and all the coefficients of $g$. 

If $\mathfrak{p} \subset \OO_{L_0}$ is a prime ideal, we can construct the representation $T_{f,\mathfrak{p}}$ of $G_{\Q}$ as in the previous section by extension of scalars: it is a free $\OO_{(L_0)_{\mathfrak{p}}}$-module of rank $2$ with a continuous linear action of $G_{\Q}$. When the residue characteristic of $\mathfrak{p}$ does not divide $2N_f$ or the discriminant of $L$, we endow $T_{f,\mathfrak{p}}$ with a basis such that the image of $G_H$ is contained in $\GL{\OO_{F_{\mathfrak{p}}}}$. 

By the work of Deligne and Serre \cite[Th\'eor\`emes 4.1, 4.6, footnote p.521]{Del-Ser}, there is an odd irreducible Artin representation $\rho_g: G_{\Q} \rar \GL{\OO_{L_0}}$ attached to $g$, and $\rho'_g: G_{\Q} \rar \PGL{L_0}$ the corresponding projective representation, and we denote $Q=|\rho'_g(G_{\Q})|$. Let $T_g$ be a free $\OO_{L_0}$-module of rank $2$ where $G_{\Q}$ acts as $\rho_g$. Let $T_{f,g,\mathfrak{p}} = T_{f,\mathfrak{p}} \otimes_{\OO_{L_0}} T_{g}$. 

\medskip
\noindent
\defi{A prime ideal $\mathfrak{p} \subset L_0$ is \emph{good} for $(f,g)$ if the following conditions are satisfied:
\begin{itemize}[noitemsep,label=\tiny$\bullet$]
\item Its residue characteristic $p$ does not divide $30N_fN_g$ or ramify in $L$. 
\item $(T_{f,\mathfrak{p}},\omega_p): G_H \rar \GL{\OO_{F_{\mathfrak{p}}}}^{\times} \times_{\OO_{F_{\mathfrak{p}}}^{\times}} \Z_p^{\times}$ is surjective, where the fibre product maps are the determinant and the $(k-1)$-th power.
\end{itemize}}

\prop[jointimage-H]{Let $\mathfrak{p}$ be a good prime. Then the image of $G_{H(\mu_{p^{\infty}})}$ in $\Aut{T_{f,\mathfrak{p}}} \times \Aut{T_g}$ is $\SL{\OO_{F_{\mathfrak{p}}}} \times \rho_g(G_H)$.}

\demo{Let $I$ denote this image, and $K/\Q$ be the Galois extension such that $G_K=\ker{\rho_g}$. It is clear that $I \subset \SL{\OO_{F_{\mathfrak{p}}}} \times \rho_g(G_H)$. Let $p_1:I \rar \SL{\OO_{F_{\mathfrak{p}}}},\, p_2: I \rar \rho_g(G_H)$ be the two projections. Because $\mathfrak{p}$ is a good prime, $p_1$ is surjective. Moreover, $p_2$ is onto if and only if $H(\mu_{p^{\infty}})\cap K \subset H$. Now, $H(\mu_{p^{\infty}}) \cap KH$ is a finite extension of $H$, so it is contained in $H(\mu_{p^n}) \cap KH$ for some $n$. Because $p \nmid N_f$, $H/\Q$ is unramified at $p$, so $H(\mu_{p^n})$ is totally ramified over $H$ above $p$. Since $p \nmid N_g$, $KH$ is unramified over $H$ above $p$. We conclude that $H(\mu_{p^{\infty}}) \cap K \subset H(\mu_{p^{\infty}}) \cap KH = H$. 

Therefore, all we need to do is show that $I$ contains $\SL{\OO_{F_{\mathfrak{p}}}} \times \{1\}$.    

If $\rho_g(G_{\Q})$ is solvable, then $\rho'_g(G_{\Q})$ is dihedral or isomorphic to $A_4$ or $S_4$. In particular, $p_2(I^{(4)})=p_2(I)^{(4)}$ is trivial (the derived subgroup of a compact group is, by definition, the closure of the subgroup generated by its commutators; here, we do this operation four times). Since $p > 3$, $\SL{\OO_{F_{\mathfrak{p}}}}$ is its own derived subgroup. Indeed, by \cite[(1.2.11)]{HOM}, it is generated by transvections, and, if $\lambda,u \in \OO_{F_{\mathfrak{p}}}$ are such that $\lambda(\lambda^2-1)$ is a unit, the following identity holds: \[\begin{pmatrix}1 & u\\0 & 1\end{pmatrix} = \begin{pmatrix}\lambda & 0\\0 & \lambda^{-1} \end{pmatrix}\begin{pmatrix}1 & \frac{u}{\lambda^2-1}\\0 & 1\end{pmatrix}\begin{pmatrix}\lambda^{-1} & 0\\0 & \lambda\end{pmatrix}\begin{pmatrix}1 & \frac{-u}{\lambda^2-1}\\0 & 1\end{pmatrix}.\] 

Hence $p_1(I^{(4)})=p_1(I)^{(4)}=\SL{\OO_{F_{\mathfrak{p}}}}$. Thus $I \supset I^{(4)} = \SL{\OO_{F_{\mathfrak{p}}}} \times \{1\}$, and we are done. 

Otherwise, one has $\rho'_g(G_{\Q}) \cong A_5$. Let $K_5 \subset K$ be such that $G_{K_5}=\ker{\rho'_g}$, and let $I_5 \subset I$ be the image of $G_{HK_5}=G_H \cap G_{K_5}$. It is enough to show that $I_5$ contains $\SL{\OO_{F_{\mathfrak{p}}}} \times \{1\}$. Since $\rho_g(G_{HK_5}) \subset \rho_g(G_{K_5})$ is composed of scalar matrices, $p_2(I_5')$ is trivial: as above, we just need to show that $p_1(I_5)=\SL{\OO_{F_{\mathfrak{p}}}}$. 

Since $K_5/\Q$ has Galois group $A_5$, it does not contain any nontrivial abelian extension of $\Q$. Since $H/\Q$ is abelian, one has $H \cap K_5 = \Q$. Therefore, $p_1(I_5)$ is an open subgroup of $p_1(I)=\SL{\OO_{F_{\mathfrak{p}}}}$ with index dividing $|\mrm{Gal}(HK_5/H)|=|A_5|=60$. Let $z \in p_1(I)$ be any element conjugate to $\begin{pmatrix}1 & a\\0 & 1\end{pmatrix}$ (with a nonzero $a$): since $p > 5$, $z$ is in the closed subgroup generated by $z^{|A_5|} \in p_1(I_5)$, hence $z \in p_1(I_5)$. Since $\SL{\OO_{F_{\mathfrak{p}}}}$ is generated by such elements $z$, $p_1(I_5)=p_1(I)$, which concludes. 

}

\cor[papier]{Let $\mathfrak{p}$ be a good prime with residue characteristic $p$ and $\sigma \in G_{\Q(\mu_{p^{\infty}})}$. There exists $\alpha \in L^{\times} \cap \OO_{L_{\mathfrak{p}}}^{\times}$ such that for all $(\gamma,\chi) \in \Gamma_f$, one has $\gamma(\alpha)=\chi(\sigma)\alpha$. Moroever, the image of $\sigma G_{H(\mu_{p^{\infty}})}$ in $\Aut{T_{f,\mathfrak{p}}} \times \Aut{T_g}$ is $\left[\begin{pmatrix}\alpha & 0\\0 & \alpha^{-1}\varepsilon_f(\sigma)\end{pmatrix}\SL{\OO_{F_{\mathfrak{p}}}}\right] \times \left[\rho_g(\sigma)\rho_g(G_H)\right]$.}

\demo{Let us briefly recall why $\alpha$ exists. For any $\gamma \in \mrm{Gal}(L/F)$, there is a unique primitive Dirichlet character $\chi_{\gamma}$ such that $(\gamma,\chi_{\gamma}) \in \Gamma_f$. The map $\gamma \in \mrm{Gal}(L/F) \longmapsto \chi_{\gamma}(\sigma) \in L^{\times}$ is a cocycle: by Hilbert 90 (for instance \cite[\S X.1, Proposition 2]{Serre-local}), it can be written as $\gamma \longmapsto \gamma(\beta)\beta^{-1}$ for some $\beta \in L^{\times}$. Since $p$ is unramified in $L$, there is some $t \in \Z$ such that $\alpha=p^t\beta$ lies in $\OO_{L_{\mathfrak{p}}}^{\times}$ and satisfies the required conditions. 

By Proposition \ref{jointimage-H}, the image in $\Aut{T_{f,\mathfrak{p}}} \times \Aut{T_g}$ of $\sigma G_{H(\mu_{p^{\infty})}}$ is $I_{\sigma} \times [\sigma \rho_g(G_H)]$, where $I_{\sigma}$ is the image in $\Aut{T_{f,\mathfrak{p}}}$ of $\sigma G_{H(\mu_{p^{\infty}})}$. Now, $I_{\sigma} = \begin{pmatrix}\alpha & 0\\0 & \alpha^{-1}\varepsilon_f(\sigma)\end{pmatrix}\SL{\OO_{F_{\mathfrak{p}}}}$ by \cite[Corollary 2.2.3]{bigimage}.}

\lem[irred-product-abstract]{Let $V,W$ be finite-dimensional vector spaces over a field $k$, and $G_1,G_2$ be subgroups of $GL(V)$ and $GL(W)$, acting irreducibly on $V$ and $W$ respectively. Let $\Gamma$ be a subgroup of $GL(V) \times GL(W)$, containing $G_1 \times \{\mrm{id}\}$ and such that its second projection is $G_2$. Then $V \otimes W$ is an irreducible $k[\Gamma]$-module. }

\demo{Let $U \subset V \otimes W$ be a proper $k[\Gamma]$-submodule. Because $G_1$ acts irreducibly on $V$, $U$ is contained in $V \otimes \ker{\alpha}$ for some nonzero $\alpha \in W^{\ast}:= \mrm{Hom}(W,k)$. Since $\Gamma \supset G_1 \times \{1\}$, the subspace $W'=\{\alpha \in W^{\ast},\, \alpha(U) = 0\}$ is nonzero and stable under $\Gamma$. Since $W^{\ast}$ is irreducible under $G_2$, $W'=W^{\ast}$, so that $U=0$. }

\prop[conditions-n-gri]{If $\mathfrak{p}$ is a good prime, $T_{f,g,\mathfrak{p}}$ satisfies conditions (N) and (gI). Moreover, if $Q/2$ is not a power of $p$, then condition (rI) holds.}

\demo{Proposition \ref{jointimage-H} shows that $G_H(\mu_{p^{\infty}})$ contains an element $\sigma$ acting by $-I_2$ on $T_{f,\mathfrak{p}}$ and $I_2$ on $T_g$, which proves condition (N). 

The condition (gI) follows from Proposition \ref{jointimage-H} and Lemma \ref{irred-product-abstract}. 

If $\rho_g(G_{\Q})$ acts absolutely irreducibly on $T_g/\mathfrak{p}$, then Proposition \ref{jointimage-H} and Lemma \ref{irred-product-abstract} imply condition (rI). Because of the complex conjugation, it is enough to show that $\rho_g(G_{\Q})$ acts irreducibly on $T_g/\mathfrak{p}$. 

Suppose that $\rho_g(G_{\Q})$ acts reducibly on $T_g/\mathfrak{p}$. Then the semisimplification of this action is the direct sum of characters $\alpha$ and $\beta$ of $\rho_g(G_{\Q})$. Thus $\gamma = \alpha/\beta$ is a character of $\rho_g(G_{\Q})$ and its kernel $K$ consists of matrices whose reduction modulo $\mathfrak{p}$ is of the form $\begin{pmatrix} \alpha & \ast\\0 & \alpha\end{pmatrix}$. In particular, we have an exact sequence $0 \rar K \rar \rho_g(G_{\Q}) \rar C \rar 0$, where $C$ is cyclic with order prime to $p$.

Let $h \in K$. Then $h$ acts on $T_g/\mathfrak{p}$ with characteristic polynomial $(X-u)^2$, for some $u \in (\OO_{L_0}/\mathfrak{p})^{\times}$ with Teichm\"uller lift $u_0 \in \OO_{(L_0)_{\mathfrak{p}}}^{\times}$ be. Then $h' := u_0^{-1}h \in \GL{\OO_{(L_0)_{\mathfrak{p}}}}$ has finite order and lies in a conjugate of the pro-$p$-group $\{M \in \GL{\OO_{(L_0)_{\mathfrak{p}}}},\, M \equiv \begin{pmatrix}1 & \ast\\0 & 1\end{pmatrix} \pmod{\mathfrak{p}}\}$. It follows that the order of $h'$ is a power of $p$, so that $h^{p^n}$ is scalar for some $n \geq 1$. 

Hence the image of $K$ in $\PGL{\OO_{L_0}/\mathfrak{p}}$ is a $p$-group, and $\rho'_g(G_{\Q})$ is an extension of a cyclic subgroup $C'$ of order prime to $p$ by a $p$-group. Since $p > 5$, this is only possible if $\rho'_g(G_{\Q})$ is a non-abelian dihedral group. Since $\rho'_g(G_{\Q})^{\mrm{ab}}$ then has exponent $2$, $C'$ has order $2$, so that $Q=|\rho'_g(G_{\Q})|/2$ is a power of $p$.

 }

\prop{Assume that $\varepsilon_f\varepsilon_g \neq 1$. Then there exists a number field $F \subset L_1 \subset L_0$ such that, for all but finitely many good primes $\mathfrak{p} \subset \OO_{L_0}$ such that the extension $L_1/F$ splits totally at the prime $\mathfrak{p} \cap \OO_F$, then $T_{f,g,\mathfrak{p}}$ has Euler-adapted Galois image. In particular, the set of such primes has positive lower density.}

\demo{Let $\sigma \in G_{\Q}$ be such that $(\varepsilon_f\varepsilon_g)(\sigma) \neq 1$ and $\rho_g(\sigma)$ has distinct eigenvalues $u,v$. By Corollary \ref{papier}, $(\gamma,\chi) \in \Gamma_f \longmapsto \chi(\sigma)$ is easily seen to be a cocycle $\mrm{Gal}(L/F) \rar L^{\times}$, so the above expression can be written as $\gamma(\alpha)/\alpha$ for some $\alpha \in L^{\times}$. Let $L_1=F(\alpha,\varepsilon_f(\sigma),u,v)$. Let $p\nmid 30N_fN_g\varphi(N_f)\varphi(N_g)|\rho_g(G_{\Q})|$ and $\mathfrak{p}$ be a good prime with residue characteristic $p$ such that $L_1/F$ is totally split at $\mathfrak{p}$. 

There exists $\sigma' \in G_{\Q(\mu_{p^{\infty}})}$ such that $\rho_g(\sigma)=\rho_g(\sigma')$ and $\sigma' \in \sigma G_H$. By Corollary \ref{papier}, the action of $\sigma'$ on $T_{f,\mathfrak{p}}$ can be chosen to be any matrix in $\Delta_{\alpha,\alpha^{-1}\varepsilon_f(\sigma)}\SL{\OO_{F_{\mathfrak{p}}}}$. By our choice of $\mathfrak{p}$, the action of $\sigma'$ on $T_{f,\mathfrak{p}}$ can be chosen to be any matrix in $\GL{\OO_{F_{\mathfrak{p}}}}$ with determinant $\varepsilon_f(\sigma)$. 

If $u^2 \neq \varepsilon_f(\sigma)^{-1}$, then we choose $\sigma'$ such that it acts on $T_{f,\mathfrak{p}}$ as $\Delta_{u^{-1},u\varepsilon_f(\sigma)}$. Then the eigenvalues of $\sigma'$ on $T_{f,g,\mathfrak{p}}$ are $1,u^2\varepsilon_f(\sigma),(\varepsilon_f\varepsilon_g)(\sigma),u/v$: $1$ is distinct from all the others, hence (by our choice of $p$), not congruent to any of the others mod $\mathfrak{p}$, and $\sigma'$ works. 

If $u^2=\varepsilon_f(\sigma)^{-1}$, we choose $\sigma'$ such that it acts on $T_{f,\mathfrak{p}}$ as $u^{-1}U$: such a $\sigma'$ works. 
}

\medskip
\noindent
Given $a,b \in K^{\times}$ for some field $K$, let us denote $\Delta_{a,b}=\begin{pmatrix}a & 0\\0 & b\end{pmatrix}$ and $U=\begin{pmatrix} 1 & 1\\0 & 1\end{pmatrix}$.  

\prop[special2]{ Let $\mathfrak{p} \subset L_0$ be a good prime. Assume that $-1 \in \varepsilon_g(G_H)$. If $g$ has multiplication by the (real or imaginary) quadratic field $K$, assume that $K \not\subset H$.
Then one of the following holds: 
\begin{itemize}[noitemsep,label=\tiny$\bullet$]
\item $\rho_g(G_H)$ is a proper subgroup of $\rho_g(G_{\Q})$, $-1 \in \varepsilon_g(G_H)^2$, $\rho'_g(G_{\Q})$ is a dihedral group of order $2n$ with $n > 2$, and if $n$ is even, then $n$ is not divisible by $4$ and $\rho'_g(G_H) \neq \rho'_g(G_{\Q})$. 
\item $T_{f,g,\mathfrak{p}}$ satisfies (sE).
\end{itemize}
}

\demo{Suppose that $\rho_g(G_H)$ contains a matrix $M$ with $\det{M}=-1$ whose projective image has order $2$. Then $M^2$ is scalar, so $M \sim \Delta_{a,b}$ for some $a,b \in \C^{\times}$ with $a=-b$ and $ab=-1$, where $A \sim B$ means that there exists $P \in \GL{\C}$ such that $A=PBP^{-1}$. Therefore $M \sim \Delta_{1,-1}$, and, by Proposition \ref{jointimage-H}, there is a $\sigma \in G_{H(\mu_{p^{\infty}})}$ acting on $T_{f,\mathfrak{p}}$ by $U$ and on $T_g$ by $M$: then $\sigma$ fulfills the condition (sE). So we would like to find such a matrix $M$.

Note that when $\rho_g(G_H)$ contains the image $c$ of the complex conjugation, $M=c$ works. 

Let $Z \leq \rho_g(G_H)$ be the subgroup of scalar matrices. Consider the surjective determinant map $\delta: \rho'_g(G_H)=\rho_g(G_H)/Z \rar \varepsilon_g(G_H)/\det{Z}$. We are done if we can find $x \in \rho_g(G_H)/Z$ with order $2$ such that $\delta(x)=-\det{Z}$. 

Assume that $\rho'_g(G_H)$ is isomorphic to $A_4$ or $A_5$. Since the derived subgroup of such a group has odd index, $\varepsilon_g(G_H)/\det{Z}$ has odd cardinality, so $-1 \in \det{Z}$ and $\delta$ is trivial on classes of order $2$: any such class works.

If $\rho'_g(G_H) \cong S_4$, we can consider two cases. If $-\det{Z}$ is in the image of a class of $A_4$, the same reasoning as above works. If not, since $S_4'=A_4$, $\delta$ vanishes on $A_4$, and therefore the natural inclusion $\{\pm 1\} \rar \varepsilon_g(G_H)/\det{Z}$ is an isomorphism, and $\delta$ is simply the signature homomorphism. In particular, the class of a transposition in $\rho'_g(G_H)$ works. 

Assume that $\rho'_g(G_H)\cong \F_2^{\oplus 2}$ (for instance if $\rho'_g(G_{\Q})$ is abelian). Then the domain of $\delta$ has exponent two and $\delta$ has cyclic image, so its image is cyclic of order $1$ or $2$. Hence $\delta$ has a non trivial kernel and is surjective, so any non trivial class in $\delta^{-1}(-\det{Z})$ works. 

Finally, let $D=\rho'_g(G_{\Q})$ and $D_H=\rho'_g(G_H)$: we just need to treat the case where $D$ and $D_H$ are dihedral, because $D_H \supset D'$. Let $C \leq D$ be the cyclic subgroup. The hypothesis implies that $D_H \not\subset C$, so $D_H$ is dihedral with cyclic subgroup $D_H \cap C$ of index $2$. Since $H/\Q$ is abelian, $D_H \supset D'=2C$. In particular, the abelianization of $D_H$ has exponent $2$, so the image of $\delta$ has either $1$ or $2$ elements. We are done if the image of $\delta$ is trivial, so assume that $\det{Z}=\varepsilon_g(G_H)^2$. If $-1 \notin \det{Z}$, since $D_H$ is generated by $D_H \backslash C$, there is a class in $D_H \backslash C$ whose image under $\delta$ is the coset of $-1$. 

In other words, we just need to treat the case where $\det{Z}=\varepsilon_g(G_H)^2$, $-1 \in \det{Z}$ and $\ker{\delta}$ is contained in $C \cap D_H$. When $4 \mid |C|$ or $D_H=D$ and $C$ has even order, $C \cap D_H$ contains an element of order $2$.   
}

\medskip
\noindent
The above Proposition focuses entirely on elements similar to tensors $U \otimes \Delta_{1,-1}$. But other elements could work as well: in the following Proposition, we replace the $-1$ with other roots of unity $\beta \neq 1$.   

\prop[specialq]{Let $\mathfrak{p} \subset L_0$ be a good prime with residue characteristic $p$. Assume that there exists an odd prime $q \mid |\rho'_g(G_{\Q})|$ such that $q \neq p$ and $\varepsilon_g(G_H)$ contains an element of order $q$. Then one of the following statements holds:
\begin{itemize}[noitemsep,label=\tiny$\bullet$]
\item $q=3$ divides $[H:\Q]$, $\rho'_g(G_{\Q})$ is isomorphic to $A_4$, $\rho'_g(G_H)$ is isomorphic to $\F_2^{\oplus 2}$. 
\item $q=3$, $\rho'_g(G_{\Q})$ and $\rho'_g(G_H)$ are isomorphic to $A_4$, $e^{2i\pi/9} \in \varepsilon_g(G_H)$, and the projection $\rho_g(G_H) \cap \SL{L_0} \rar \rho'_g(G_H)$ is not surjective. 
\item $T_{f,g,\mathfrak{p}}$ satisfies (sE). 
\end{itemize}}

\demo{As above, let $Z \leq \rho_g(G_H)$ be the subgroup of central matrices. When convenient, we will identify $Z$ to a subgroup of $\C^{\times}$. Assume that $Z$ contains $q$-th roots of unity, and let $M \in \rho_g(G_H) \cap \SL{L_0}$ have order $q$ in $\PGL{L_0}$. Then $M^q = \{\pm I_2\}$, so that $M \sim t\Delta_{u,u^{-1}}$, where $t$ is a sign, and $u$ is a primitive $q$-th root of unity. Since $uI_2 \in Z$, $uM \sim t\Delta_{u^2,1}$ (with $u^2\not\equiv 1 \pmod{\mathfrak{p}}$) is also in $\rho_g(G_H)$. The tensor $(tU) \otimes (uM)$ then shows (sE), and we already know that the other hypotheses hold. 

The first step is to show that, except in specific cases, $Z$ contains all the $q$-th roots of unity, and $\rho_g(G_H) \cap \SL{L_0}$ contains an element whose projective image has order $q$. 

Consider, as before, the surjective determinant map $\delta: \rho'_g(G_H) \rar \varepsilon(G_H)/\det{Z}$. Note that $\delta$ is trivial on the derived subgroup of $\rho'_g(G_{\Q})$ (which is contained in $\rho'_g(G_H)$). 

If $\rho'_g(G_{\Q})$ is isomorphic to $A_5$ or $S_4$, or dihedral, then its derived subgroup has index coprime to $q$. Therefore, $\rho'_g(G_H)$ contains all the elements of order $q$ of $\rho'_g(G_{\Q})$, and $\varepsilon(G_H)/\det{Z}=\mrm{im}(\delta)$ is $q$-torsion free. Since $q$ is an odd prime and $\varepsilon_g(G_H)$ contains some $q$-th root of unity, $Z$ contains all the $q$-th roots of unity. Finally, since its image is $q$-torsion free, $\delta$ maps elements of order $q$ to the identity, so they have lifts in $\rho_g(G_H) \cap \SL{L_0}$.

Assume now that $\rho'_g(G_{\Q})$ is isomorphic to $A_4$: then $q=3$. Moreover, since $A_4'$ is the bitransposition subgroup (with index $3$ in $A_4$), $\rho'_g(G_H)$ can be $\F_2^{\oplus 2}$ or $A_4$, and $\varepsilon_g(G_H)/\det{Z}$ is trivial in the former case, and has order $1$ or $3$ in the latter. 

If $\rho'_g(G_H) \cong A_4$ and $\varepsilon_g(G_H)/\det{Z}$ is trivial, then, as above, $Z$ contains all third roots of unity, and $\rho_g(G_H) \cap \SL{L_0}$ contains an element whose projective image has order $3$: the argument above shows that we have Euler-adapted Galois image. 

All we need to do is show that if $e^{2i\pi/9} \notin \varepsilon_g(G_H)$, $\rho'_g(G_H) \cong A_4$ and $\varepsilon_g(G_H)/\det{Z}$ has order $3$, then we have Euler-adapted Galois image. Since $\varepsilon_g(G_H)$ is cyclic, under the previous assumptions, $|\det{Z}|$ (hence $|Z|$) is coprime to $3$, and $\varepsilon_g(G_H) \backslash \det{Z}$ contains an element of order $3$. In other words, there is some $x \in \rho_g(G_H)$ such that $\det{x}$ is a third root of unity. If the projective image of $x$ is not a $3$-cycle, then $x^2$ is scalar, and its determinant is also a primitive third root of unity, which we excluded. Therefore (up to replacing $x$ with $x^2$), we may assume that $x^3$ is scalar, that $x$ is not, and $\det{x}=j$. Let $t$ be a multiple of the order of $x^3$ in $Z$ which is congruent to $1$ mod $3$: replacing $x$ with $x^t$ if necessary, we may assume that $\det{x}=j$, $x$ is not scalar, and $x^3=I_2$. It follows that $x \sim \Delta_{1,j}$. As above, the tensor $U \otimes x$ shows that (sE) holds.  
}

\cor{If any of the following conditions hold, then $T_{f,g,\mathfrak{p}}$ has Euler-adapted Galois image for all but finitely many primes $\mathfrak{p}$:
\begin{enumerate}[noitemsep,label=(\roman*)]
\item \label{sc1} $\rho_g(G_H)$ contains the image under $\rho_g$ of the complex conjugation
\item \label{sc2} $-1 \in \varepsilon_g(G_H)$ and $\rho_g(G_H)$ and $\rho_g(G_{\Q})$ have the same projective image, which is not a dihedral group of order $2n$ with $n > 1$ odd.
\item \label{sc3} $\varepsilon_g(G_H)$ contains $-1$ but not $i$.
\item \label{sc4} $\varepsilon_f^2$ is the trivial character, $\varepsilon_g(G_{\Q})$ contains an element of odd prime order $q \mid |\rho'_g(G_{\Q})|$, and, if $\rho'_g(G_{\Q}) \cong A_4$, then $\rho_g(G_{\Q}) \cap \SL{L_0} \rar \rho'_g(G_{\Q})$ is onto.  
\end{enumerate}}

\demo{Condition \ref{sc1} comes from Proposition \ref{jointimage-H} and proved like in Proposition \ref{special2}. Conditions \ref{sc2} and \ref{sc3} follow from the Proposition \ref{special2}. For condition \ref{sc4}, assume that $f$ has an inner twist by $(\gamma,\chi)$, then $\varepsilon_f=\gamma(\varepsilon_f)=\varepsilon_f\chi^2$, so that $\chi$ is a quadratic character. Therefore, $\mrm{Gal}(H/\Q)$ is a $\F_2$-vector space, which rules out both exceptional cases of Proposition \ref{specialq}.  }

\cor[BSD-improved]{Let $E$ be an elliptic curve over $\Q$ without complex multiplication and $\rho$ be a two-dimensional odd irreducible Artin representation of $G_{\Q}$, with splitting field $K$. Let $L_0/\Q$ be a finite extension containing all the eigenvalues of Frobenius for $\rho$ and $\mathfrak{p}$ be a prime ideal of $\OO_L$ with residue characteristic $p$. Suppose that the following technical hypotheses hold:
\begin{enumerate}[noitemsep,label=(\roman*)]
\item $p$ is coprime to $30N_{\rho}N_E$, where $N_E$ (resp. $N_{\rho}$) is the conductor of $E$ (resp. $\rho$). 
\item If $\rho$ has projective dihedral image, then the cyclic subgroup is not a $p$-group.  
\item The $p$-adic Tate module of $E$ is a surjective representation $G_{\Q} \rar \GL{\Z_p}$. 
\item $E$ is ordinary at $p$ and $\rho(\Fr_p)$ has distinct eigenvalues mod $\mathfrak{p}$.
\end{enumerate} 
If $L(E,\rho,1) \neq 0$, then the group $\mrm{Hom}_{\Z_p[\mrm{Gal(K/\Q)}]}(\rho,\mrm{Sel}_{p^{\infty}}(E/K))$ is finite. }

\demo{This is \cite[Theorem 11.7.4]{KLZ15} without the assumption on the coprimality of conductors. Let $f \in \mathcal{S}_2(\Gamma_0(N_E))$ and $g \in \mathcal{S}_1(\Gamma_1(N_{\rho}))$ be the normalized newforms associated to $E$ and $\rho$: then $L(E,\rho,1)$ is the automorphic Rankin-Selberg $L$-value $L(f \otimes g,1)$ (it has a classical definition in \cite[Section 2]{Li-RS}, which we are going to use). We can apply without any modifications the proof of \cite[Theorems 11.7.3, 11.7.4]{KLZ15}, as long as the following facts hold:

\begin{enumerate}[noitemsep,label=(\alph*)]
\item $T_p E \otimes \rho$ has Euler-adapted Galois image.
\item The $p$-adic $L$-value $L_p(f,g,1)$ does not vanish. 
\end{enumerate} 

Since the modular form $f$ associated to the elliptic curve $E$ has no inner twists, $H=\Q$ and Proposition \ref{special2} shows that $T_{f,g,\mathfrak{p}}$ has Euler-adapted Galois image. Now, we check the non-vanishing of $L_p(f,g,1)$. 

First, the $L$-value $L(f,g,1)$ from \cite[\S 2.7]{KLZ15} differs from the value $L(f \otimes g,1)$ by finitely many Euler factors at primes $q \mid N_EN_{\rho}$, described in \cite[Section 2]{Li-RS}. In \emph{loc. cit.}, some of them are denoted as $\theta_q$; the others come from replacing the local sum $\sum_{n \geq 0}{a_{q^n}(f)a_{q^n}(g)q^{-ns}}$ with $\sum_{n \geq 0}{a_{q^n}(f \otimes \chi)a_{q^n}(g\otimes \overline{\chi})q^{-ns}}$ for some Dirichlet character $\chi$ chosen so that $f \otimes \chi, \overline{g} \otimes \chi$ satisfy the hypotheses (A), (B), (C) of \emph{loc. cit.}. Note also that \emph{loc. cit.} uses a different normalization: their $s=1/2$ is our $s=1$. By checking every case, one sees that none of these factors vanish at $s=1$. It follows that $L(f,g,1) \neq 0$, where $L(f,g,1)$ is defined in \cite[\S 2.7]{KLZ15}.

Finally, we use the interpolation property of \cite[Theorem 2.7.4]{KLZ15}: the $p$-adic $L$-value does not vanish as long as (following the notations of \emph{loc. cit.}) $\mathcal{E}(f)$, $\mathcal{E}^{\ast}(f)$, and $\mathcal{E}(f,g,1)$ do not vanish. $\mathcal{E}(f)$ does not vanish because the two Frobenius eigenvalues of $f$ at $p$ have modulus $\sqrt{p}$, by the Hasse-Weil bound. $\mathcal{E}^{\ast}(f)$ does not vanish because $f$ is ordinary at $p$, so exactly one of the Frobenius eigenvalues at $p$ is a $p$-unit, thus the two cannot be equal. Finally, $\mathcal{E}^{\ast}(f,g,1) \neq 0$ because the Frobenius eigenvalues of $f$ and $g$ have modulus respectively $\sqrt{p}$ and $1$ (for any complex embedding), so $L_p(f,g,1) \neq 0$ and we are done.}

\section{Families of dihedral counter-examples: proof of Theorem \ref{INF-CEX}}
\label{sect-theoD}

In this Section, we keep the notations from Section \ref{sect-theoA}, but furthermore assume that $\varepsilon_g=\varepsilon_f^{-1}\varepsilon$, where $\varepsilon$ is the character of some quadratic number field $K \subset H$ such that $g$ has CM by $K$: in particular, $\varepsilon$ factors through $\mrm{Gal}(H/\Q)$. 

We are interested in the following question: given a good prime ideal $\mathfrak{p}$ of $\OO_{L_0}$ for the couple $(f,g)$, when does $T_{f,g,\mathfrak{p}}$ satisfy (wE)?

We first recall the following fact, already identified as a general obstruction in \cite[Proposition 4.1.1]{bigimage}. 

\lem[efeg=1]{Let $M,N \in \GL{K}$ be matrices over a field $K$ such that $\det{M}\det{N}=1$. Then $\ker\left[M \otimes N-1\right]$ is never a line.}

\demo{The problem is clearly solved if $M$ or $N$ is scalar. We now assume that neither $M$ nor $N$ is scalar and that $1$ is an eigenvalue of $M \otimes N$. 
Assume that $M \sim \alpha U$ for some $\alpha \in K^{\times}$. Then $\alpha^{-1}$ is an eigenvalue of $N$. Since $\det{M}\det{N}=1$, the characteristic polynomial of $N$ is $(X-\alpha^{-1})^2$. As $N$ is not scalar, $N \sim \alpha^{-1} U$. It is easy to see that $U \otimes U - 1$ has rank two, so $\ker(M \otimes N-1)$ is not a line. 
Thus, we can (and do) assume that $M$ has distinct eigenvalues $m,m'$ and $N$ has distinct eigenvalues $n,n'$. Then the endomorphism $M \otimes N$ is diagonalizable with eigenvalues $mn,m'n',m'n,nm'$. Since $1=(mn)(m'n')=(mn')(m'n)$, there is never exactly one of these eigenvalues which is equal to one.}

As noted by Loeffler, this lemma implies that, for any couple $(f',g')$ of newforms, $T_{f',g',\mathfrak{p}}$ cannot satisfy (wE) if the product of the characters of $f',g'$ is trivial, whence the requirement that $\varepsilon_f\varepsilon_g \neq 1$ in Question \ref{loefflerqn}. 

Define the map $\mathcal{B}: \mrm{Gal}(H/\Q) \times G_F \rar \{\pm 1\}$ as follows: let $(\sigma,\tau) \in \mrm{Gal}(H/\Q) \times G_F \rar \{\pm 1\}$. By Corollary \ref{papier}, there is some $\alpha \in L^{\times}$ such that $\gamma(\alpha)=\chi(\sigma)\alpha$ for all $(\gamma,\chi) \in \Gamma_f$. Since $\gamma(\varepsilon_f)=\chi^2\varepsilon_f$, one has $\nu := \frac{\alpha^2}{\varepsilon_f(\sigma)} \in F^{\times}$. Then $\mathcal{B}(\sigma,\tau)=\frac{\tau(\sqrt{\nu})}{\sqrt{\nu}}$.

\prop{$\mathcal{B}$ is a bilinear pairing, and factors through a finite quotient of $G_F$. Moreover, for every $\sigma \in \mrm{Gal}(H/\Q)$, $\mathcal{B}(\sigma, \cdot): G_F \rar \{\pm 1\}$ is unramified away from $2 \cdot \mrm{disc}(L/F)$.}

\demo{The bilinearity is a straightforward verification. To prove the rest of the statement, note that if $\tau \in \mrm{Gal}(\overline{F}/L(\sqrt{\mu}))$ (where $\mu$ is the group of roots of unity that are values of $\varepsilon_f$), then $\mathcal{B}(\cdot,\tau)=1$.}

\medskip
\noindent
Let $\mathfrak{p} \subset \OO_{L_0}$ be good prime ideal for $(f,g)$, we write $\mathcal{B}_{\mathfrak{p}}=\mathcal{B}(\cdot,\Fr_{\mathfrak{p}}) : \mrm{Gal}(H/\Q) \rar \{\pm 1\}$. Denote by $M$ the subgroup of $\mrm{Gal}(H/\Q)$ made by those $\sigma \in \mrm{Gal}(H/\Q)$ such that $\mathcal{B}(\sigma,\cdot)=1$. By Cebotarev, $M$ is thus the intersection of the kernels of the $\mathcal{B}_{\mathfrak{p}}$.

\prop{Let $\mathfrak{p}$ be a prime ideal of $\OO_{L_0}$ of residue characteristic $p$, good for $(f,g)$. Then the following are equivalent: 
\begin{itemize}[noitemsep, label=\tiny$\bullet$]
\item $T_{f,g,\mathfrak{p}}$ satisfies (wE),
\item $T_{f,g,\mathfrak{p}}$ satisfies (sE),
\item $\mathcal{B}_{\mathfrak{p}} \neq \varepsilon$. 
\end{itemize}}

\demo{Suppose that $T_{f,g,\mathfrak{p}}$ satisfies (wE): there exists $\sigma \in G_{\Q(\mu_{p^{\infty}})}$ such that $\ker(\sigma - 1 \mid T_{f,g,\mathfrak{p}})$ is a line. Let $\alpha \in L^{\times} \cap \OO_{L,\mathfrak{p}}^{\times}$ be associated to $\sigma$ by Corollary \ref{papier}. The action of $\sigma$ on $T_{f,\mathfrak{p}},T_g$ is given by some matrices $M \in \GL{\OO_{(L_0)_{\mathfrak{p}}}}, N \in \GL{\OO_{L_0}}$ respectively. By the assumption on $\sigma$ and Lemma \ref{efeg=1}, one has $1 \neq \det{M}\det{N}=\varepsilon_f(\sigma)\varepsilon_f(\sigma)^{-1}\varepsilon(\sigma)$, so $\sigma \notin G_K$. Therefore, $N^2$ is scalar and $\det{N}=-\varepsilon_f(\sigma)^{-1}$. Let $\omega \in L_0$ be a square root of $\varepsilon_f(\sigma)^{-1}$, then the two eigenvalues of $N$ are $\omega$ and $-\omega$. 

By Corollary \ref{papier}, we can write $M=\alpha^{-1}\varepsilon_f(\sigma)M'$, where $M' \in \GL{\OO_{F_{\mathfrak{p}}}}$ has determinant $\frac{\alpha^2}{\varepsilon_f(\sigma)}$. Thus $\alpha\omega^2$ is an eigenvalue of $M' \otimes \Delta_{\omega,-\omega}$, so that one of the eigenvalues of $M'$ is $\pm \omega \alpha$. Since $(\pm \omega \alpha)^2=\det{M'}$, one has $\pm \alpha \omega = \frac{1}{2}\Tr{M'} \in F_{\mathfrak{p}}$, so $\mathcal{B}(\sigma,\Fr_{\mathfrak{p}})=1 \neq \varepsilon(\sigma)$.

Suppose conversely that $\mathcal{B}(\cdot,\Fr_{\mathfrak{p}})\neq \varepsilon$. Then $\ker{\mathcal{B}(\cdot,\Fr_{\mathfrak{p})}}$ is not contained in $\ker{\varepsilon}$, otherwise $\varepsilon$ would be a power of $\mathcal{B}_{\mathfrak{p}}$: but $\varepsilon$ is clearly distinct from $\mathcal{B}_{\mathfrak{p}}^i$ for $i \in \{0,1\}$ and $\mathcal{B}_{\mathfrak{p}}^2=1$. Hence, there exists $s \in \mrm{Gal}(H/\Q)$ such that $\varepsilon(s)=-1$ and $\mathcal{B}_{\mathfrak{p}}(s)=1$. By Corollary \ref{papier}, let $\alpha \in \OO_{L,\mathfrak{p}}^{\times}$ such that for all $(\gamma,\delta) \in \Gamma_f$, $\gamma(\alpha)=\delta(s)\alpha$. Since $\gamma(\varepsilon_f)=\delta^2\varepsilon_f$, one has $\nu := \frac{\alpha^2}{\varepsilon_f(\sigma)} \in F^{\times}$, hence $\nu \in \OO_{F,\mathfrak{p}}^{\times}$. Since $\mathcal{B}_{\mathfrak{p}}(s)=1$, $\nu$ is a square in $F_{\mathfrak{p}}$. Therefore, $\varepsilon_f(\sigma)^{-1}=\nu\alpha^{-2}$ has a square root $\omega$ in $L_{\mathfrak{p}}$.  

Still by Corollary \ref{papier}, we can choose $\sigma \in G_{\Q(\mu_{p^{\infty}})}$ such that $\sigma_{|H}=s$ and such that $\sigma$ acts on $T_{f,\mathfrak{p}}$ by the matrix $\alpha^{-1}\varepsilon_f(s)\begin{pmatrix}\nu & 1\\0 & \nu\end{pmatrix}$. 

Moreover, $\varepsilon(\sigma)=-1$, so $\rho'_g(\sigma) \in \PGL{\OO_{L_0}}$ has order two. Since $\det{\rho_g(\sigma)}=-\varepsilon_f(\sigma)^{-1}$, and $p \neq 2$, $\rho_g(\sigma)$ is $\GL{\OO_{(L_0){\mathfrak{p}}}}$-conjugate to $\Delta_{\omega,-\omega}$. Therefore, $\sigma$ acts on $T_{f,g,\mathfrak{p}}$ by $\begin{pmatrix}1 & a & 0 & 0\\0 & 1 & 0 & 0\\0 & 0 & -1 & \ast\\0 & 0 & 0 & -1\end{pmatrix}$ with $a \in \OO_{(L_0)_{\mathfrak{p}}}^{\times}$: hence $T_{f,g,\mathfrak{p}}/(\sigma-1)T_{f,g,\mathfrak{p}}$ is free of $\OO_{(L_0)_{\mathfrak{p}}}$-rank one. 
}

\cor{The answer to Question \ref{loefflerqn} for $(f,g)$ is negative if, and only if, $\varepsilon(M)=1$.}

\demo{By Cebotarev and the previous Proposition, the answer to Question \ref{loefflerqn} for $(f,g)$ is positive if and only if $\varepsilon$ cannot be written as a $\mathcal{B}(\cdot,\tau)$ for any $\tau \in G_F$. Because $B$ is bilinear with values in an $\F_2$-vector space, the set of the $\mathcal{B}(\cdot,\tau)$ for $\tau \in G_F$ is exactly the set of homomorphisms $f: \mrm{Gal}(H/\Q) \rar \{\pm 1\}$ such that $f(M)=1$, whence the conclusion.}

\prop[whenepsilonexists]{Let $H_0 \subset H$ be the subfield such that $G_{H_0} = \ker{\varepsilon_f}$. 
Assume that one of the following holds:
\begin{enumerate}[noitemsep,label=$(\alph*)$]
\item $k$ is even and $H$ is not totally real, 
\item\label{weps-b} $k$ is odd and $[H:H_0]>2$,. then $(\ast)$ holds unless perhaps if $[H:H_0] = 2$, $\sigma_0 \in \mrm{Gal}(H/H_0)$ is the unique nontrivial element, $L$ contains a fourth root of unity $i$, and, with $c$ the complex conjugation, for all $(\gamma,\delta) \in \Gamma_f$, $\gamma(i)=\delta(\sigma_0c)i$.  
\item\label{weps-c} $k$ is odd, $[H:H_0]=2$, and $\varepsilon_f$ does not take its values in $L^{\times 2}$ (for instance if $L$ does not contain a fourth root of unity $i$), 
\item\label{weps-d} $k$ is odd, $[H:H_0]=2$, $i \in L$, and for all $(\gamma,\delta) \in \Gamma_f$, $\gamma(i)=\delta(-1)i$,
\item\label{weps-e} $k$ is odd, $\mrm{Gal}(H/H_0)$ contains a unique nontrivial element $\sigma_0$, $i \in L$ and, for some $(\gamma,\delta) \in \Gamma_f$, one has $\gamma(i) \neq \delta(\sigma_0)\delta(-1)i$.
\end{enumerate}
Then there exists a nontrivial character $\varepsilon: \mrm{Gal}(H/\Q)/M \rar \{\pm 1\}$ mapping the complex conjugation $c$ to $(-1)^{k-1}$. \\
Moreover, when $k$ is odd, $[H:H_0]=2$, $i \in L$ but is not a value of $\varepsilon_f$, then the existence of such a $\varepsilon$ is equivalent to the last clause of \ref{weps-e}. Under the same assumptions, if $\varepsilon_f^2=1$, then \ref{weps-d} and \ref{weps-e} are equivalent.
}
 
\demo{Remember that $L$ is the number field generated by the coefficients of $f$, so it is a subfield of $\C$. Note that $\mathcal{B}$ realizes an isomorphism from $\mrm{Gal}(H/\Q)/M$ to a subgroup $B$ of the $\F_2$-vector space $\mrm{Hom}(G_F,\F_2)$, so $\mrm{Gal}(H/\Q)/M$ has exponent $2$. 

Suppose first that $k$ is even and $H$ is not totally real. Let $c \in \mrm{Gal}(H/\Q)$ be the image of the complex conjugation and let $\alpha \in L^{\times}$ be associated to $c$ as in Corollary \ref{papier}. It is enough to prove that $c \notin M$. Since $k$ is even, $\varepsilon_f(c)=1$, so one has $\mathcal{B}(c,\tau) = \frac{\tau(\alpha)}{\alpha}$ for all $\tau \in G_F$. In particular, for every $(\gamma,\delta) \in \Gamma_f$, one has $\mathcal{B}(c,\gamma)=\delta(-1)$. Since $H$ is not totally real, there is some $\delta$ such that $\delta(-1)=-1$, so $c \notin M$, whence the conclusion. 

Now suppose that $k$ is odd. The goal is to show that there exists a nontrivial character $\mrm{Gal}(H/\Q)/\langle M,c\rangle \rar \{\pm 1\}$. To do that, we need to show that $\mrm{Gal}(H/\Q) \supsetneq \langle M,c\rangle$, or in other words that $B \supsetneq \{1,\mathcal{B}(c,\cdot)\}$.

Let $\sigma \in M$ be such that $\varepsilon_f(\sigma)=1$: let $\alpha$ be associated to $\sigma$ by Corollary \ref{papier}. Then, for all $\tau \in G_F$, $\tau(\alpha)=\alpha$, hence $\alpha \in F^{\times}$. Therefore, for all $(\gamma,\delta) \in \Gamma_f$, $\alpha=\gamma(\alpha)=\delta(\sigma)\alpha$, so $\delta(\sigma)=1$, whence $\sigma_{|H}=\mrm{id}$. Thus, $(\pi,\varepsilon_f): \mrm{Gal}(H/\Q) \rar F^{\times}/F^{\times 2} \times \C^{\times}$ is injective. 

In particular, $\sigma \in \mrm{Gal}(H/H_0) \longmapsto \mathcal{B}(\sigma,\cdot)$ is injective, so that $|B| \geq [H:H_0]$, so we are done if $[H:H_0] > 2$. 

Now suppose that $[H:H_0]=2$ and let $\sigma_0$ be the unique nontrivial element of $\mrm{Gal}(H/H_0)$. Let $\sigma \in M$: then the $\alpha \in L^{\times}$ associated to $\sigma$ by Corollary \ref{papier} is such that $\frac{\alpha^2}{\varepsilon_f(\sigma)} \in F^{\times 2}$, so $\varepsilon_f(\sigma) \in L^{\times 2}$.

Since $\mrm{Gal}(H_0/\Q)$ is isomorphic to the image of $\varepsilon_f$ which contains $-1$, it is a group with even cardinality. Moreover, since $\sigma_0 \notin M$, $\sigma_0$ is not a double in $\mrm{Gal}(H/\Q)$. It follows that $\mrm{Gal}(H/\Q) \cong \Z/2\Z\sigma_0 \oplus \Z/n\Z \sigma$, where $n$ is the cardinality of the image of $\varepsilon_f$, and $\sigma$ is such that $\varepsilon_f(\sigma)$ has order $n$. Hence, representatives for $\mrm{Gal}(H/\Q)/M$ are given (possibly with repetitions) by $\sigma_0,\sigma_0\sigma,\sigma,\mrm{id}$. Since $\varepsilon_f(\sigma)=\varepsilon_f(\sigma\sigma_0)$ is a generator of the image of $\varepsilon_f$, the previous paragraph shows that unless $\varepsilon_f$ takes its values in $L^{\times 2}$, $\sigma,\sigma_0\sigma \notin M$, hence $|B|=4$ are we are done. 

In particular, if $i$ is not a value taken by $\varepsilon_f$, but $i \in L$, the above considerations show that representatives for $\mrm{Gal}(H/\Q)$ are given by $\mrm{id},\sigma_0,c,\sigma_0c$. So $\mrm{Gal}(H/\Q)$ is generated by $M$ and $c$ if and only if $\sigma_0 c \in M$, which is equivalent to the statement that for all $(\gamma,\delta) \in \Gamma_f$, $\gamma(i)=\delta(\sigma_0)\delta(-1)i$. 

In general, when \ref{weps-e} holds, one has $\sigma_0 c \notin M$, so, since $\sigma_0 \notin M$, $\sigma_0$ is not in the subgroup generated by $M$ and $c$, so the conclusion follows. Furthermore, \ref{weps-d} is a stronger condition than \ref{weps-e}, since it implies that $c \in M$.

Suppose finally that $\varepsilon_f^2=1$ and \ref{weps-e} holds. The discussions above show that the $2$-torsion subgroup of $\mrm{Gal}(H/\Q)$ is a set of representatives for $\mrm{Gal}(H/\Q) \otimes \F_2$. We have an injective homomorphism $(\gamma,\delta) \in \Gamma_f \longmapsto \delta \in \mrm{Hom}(\mrm{Gal}(H/\Q),\{\pm 1\})$, so that every homomorphism $\Gamma_f \rar \{\pm 1\}$ is of the form $(\gamma,\delta) \longmapsto \delta(\sigma)$ for some element $\sigma \in \mrm{Gal}(H/\Q)$ of order dividing $2$. In particular, we know that some $\sigma \in \{\sigma_0,c,\sigma_0c\}$ is associated to $(\gamma,\delta) \longmapsto \frac{\gamma(i)}{i}$: but then $i$ is the $\alpha$ associated by Corollary \ref{papier} to $\sigma$. By taking the couple $(c,\varepsilon_f^{-1})$, we know that $\sigma \neq \sigma_0$, so $\sigma \in \{c,\sigma_0c\}$, and by the above discussion, it means that exactly one of $c,\sigma_0c$ is in $M$.
}

\rems{
\begin{enumerate}[noitemsep,label=(\roman*)]
\item Let $f \in \mathcal{S}_2(\Gamma_0(63))$ be a newform with LMFDB \cite{lmfdb} label $63.2.a.b$\footnote{Technically, these labels refer to a Galois orbit of newforms; however, the choice of embedding does not matter in all the situations below, so we omit it.}: it has coefficients in $\Q(\sqrt{3})$, an inner twist by the quadratic character $\varepsilon$ of $K=\Q(\sqrt{-3})$ and no complex multiplication. The newform $g \in \mathcal{S}_1(\Gamma_1(675))$ with LMFDB label $675.1.c.b$ has rational coefficients, character $\varepsilon$ and complex multiplication by $K$. The argument shows that (after possibly removing a finite number of primes) $T_{f,g,\mathfrak{p}}$ does not satisfy (wE) when the residue characteristic is congruent to $5,7$ mod $12$.
\item Consider a newform $f \in \mathcal{S}_3(\Gamma_1(75))$ with LMFDB label $75.3.d.b$. Then $L=\Q(i,\sqrt{5})$, $F=\Q$, $H_0=\Q(\sqrt{-15})$, $H=\Q(j,\sqrt{5})$. Then it can be checked that $\mrm{Gal}(H/\Q)$ is generated by $M$ and $c$. When $\varepsilon$ is the character of the unique real quadratic field contained in $H$ and $g$ is a weight one newform with character $\varepsilon\varepsilon_f^{-1}$ and with real multiplication by $\varepsilon$, property (sE) then holds for all but finitely many $T_{f,g,\mathfrak{p}}$, even though Theorem \ref{weakerconds} does not apply.
\item When $k$ is even, $H_0$ is necessarily totally real. However, $H$ can be totally real and larger than $H_0$. This is what happens for the newform $f \in \mathcal{S}_2(\Gamma_0(289))$ with LMFDB label $289.2.a.f$: indeed, $H_0=\Q$, but $H=\Q(\sqrt{17})$. By Theorem \ref{weakerconds}, the answer to Question \ref{loefflerqn} is positive for $f$ and any weight one newform. 
\item For a newform $f \in \mathcal{S}_3(\Gamma_1(24))$ with LMFDB label $24.3.h.c$, the condition \ref{weps-c} of Proposition \ref{whenepsilonexists} is satisfied: indeed, in this situation $L=\Q(\sqrt{2},\sqrt{-7})$ does not contain $i$, while $H=\Q(\sqrt{2},j)$, and $H_0=\Q(\sqrt{-6})$ is such that $[H:H_0]=2$. Indeed, $B$ contains the two distinct nontrivial elements $\mathcal{B}(c,\cdot)$ and $\mathcal{B}(\Fr_5,\cdot)$ (they are the Kummer characters associated to $\sqrt{14}$ and $\sqrt{2}$ respectively), so that $M$ is trivial and $\mrm{Gal}(H/\Q)$ is larger than $\{\mrm{id},c\}$.
\item For a newform $f \in \mathcal{S}_3(\Gamma_1(64))$ with LMFDB label $64.3.d.a$, we have $L=\Q(\mu_{12})$ and $F=\Q$; moreover, $H=\Q(\mu_8)$ and $\mrm{Gal}(H/H_0)$ corresponds to the subgroup $\{1,3\}$ of $(\Z/8\Z)^{\times}$. Because $a_{63}(f) \in iF^{\times}$, the condition \ref{weps-d} of Proposition \ref{whenepsilonexists} is verified.   
\item Even when $i \in L$ is not a value of $\varepsilon_f$, it is not true that \ref{weps-e} implies \ref{weps-d}. A counter-example is given by the newform $f \in \mathcal{S}_3(\Gamma_1(324))$ with LMFDB label $324.3.f.s$\footnote{With the tools dicussed here, this verification is not difficult using the LMFDB data and MAGMA \cite{magma}; we omit it for the sake of brevity.}. In this situation, $\varepsilon_f$ has conductor $36$ and order $6$, while $L$ has degree $16$ and $|\Gamma_f|=8$, so that $F=\Q(\sqrt{57})$. %
\item Consider a newform $f \in \mathcal{S}_3(\Gamma_1(21))$ with LMFDB label $21.3.h.b$. In this case, one has $|\Gamma_f|=4$, but $H=H_0=\Q(j)\Q(\mu_7)^+$. In fact, $F=\Q$, $L=\Q(j,\sqrt{-5})$, and one can directly compute that $B = \{\mathcal{B}(\mrm{id},\cdot),\mathcal{B}(c,\cdot)\}$.  
\end{enumerate}
}

\cor{Suppose that the hypotheses of Proposition \ref{whenepsilonexists} hold, and let $\varepsilon$ be the character given by the conclusion Proposition \ref{whenepsilonexists}. Let $g$ be a weight one newform with complex (or real) multiplication by $\varepsilon$ and with character $\varepsilon\varepsilon_f^{-1}$. Then the answer to Question \ref{loefflerqn} for the couple $(f,g)$ is negative.}

\medskip
\noindent
The following result completes the proof of Theorem \ref{INF-CEX}. 

\prop{Let $\varepsilon$ be the character of the imaginary (resp. real) quadratic field $K$, and $\chi$ be an even (resp. odd) primitive Dirichlet character. There exist infinitely many weight one newforms $g$ with complex (resp. real) multiplication by $K$, and character $\varepsilon\chi$. }

\demo{This statement is easier to establish using Galois representations. Let $C$ be the conductor of $\chi$, $n$ be the exponent of $\chi$ and $D$ be a positive multiple of $Cn$. Let $K'$ the abelian extension of $K$ whose group of norms is $K^{\times}\prod_{v}{M_v}$, where $M_v$ is $\C^{\times}$ when $v$ is a complex place, $\R^{\times}_+$ when $v$ is a real place, $1+D\OO_{K_v}$ at every place $v$ dividing $D$, and $\OO_{K_v}^{\times}$ at all other finite places. We claim that $\chi$ vanishes on the kernel of the transfer map $V: G_{\Q}^{ab} \rar \mrm{Gal}(K'/K)$. It follows that $\sigma \in G_{\Q} \longmapsto \chi(\sigma)$ extends (through $V$) to at least $\frac{[K':K]}{\varphi(n)}$ characters $\psi: G_K \rar \mrm{Gal}(K'/K) \rar \C^{\times}$. The Artin representation $\rho_{\psi} := \mrm{Ind}_K^{\Q}{\psi}$ has determinant $\varepsilon\psi(V)=\varepsilon\chi$ and it uniquely determines $\psi$ (up to conjugation by $\mrm{Gal}(K/\Q)$), so this construction produces at least $\frac{[K':K]}{2\varphi(n)}$ pairwise non-isomorphic Artin representations as requested. It is easy to see that we can choose $D$ so that $[K':K]$ is arbitrarily large, whence the conclusion.}

\section{Examples and counter-examples}
\label{sect-counterex}
So far, we have identified in Theorem \ref{weakerconds} and in Section \ref{sect-theoA} some sufficient conditions to ensure a positive answer to Question \ref{loefflerqn}. On the other hand, the situation studied in Section \ref{sect-theoD} provides in some cases families of counter-examples. The purpose of this Section is to study some situations where the results of neither Section applies. We provide a few counterexamples which do not stem from Section \ref{sect-theoD}. We also showcase additional examples proving that, while the exceptional cases highlighted in Propositions \ref{special2} and \ref{specialq} can occur, they are not necessarily incompatible with Euler-adapted Galois images. Some of these counter-examples are explicit: as in the examples of Section \ref{sect-theoD}, we will then always refer to newforms by their LMFDB label \cite{lmfdb}. While this label only refers to the equivalence class of the modular forms under the action of $G_{\Q}$, the arguments work regardless of the choice of embedding of the modular forms in $\C$.

\prop[special-counter]{ Assume that $\mathfrak{p}$ is a good prime, and that $(f,g,p)$ falls in one of the following cases:
\begin{enumerate}[label=(\alph*),noitemsep]
\item\label{A5det} $f$ is $961.2.a.b$, $g$ is $3875.1.d.a$, and $p\equiv 13,37,83,107 \pmod{120}$.
\item\label{S4det} $f$ is $1849.2.a.g$, $g$ is $688.1.b.b$, and $p \equiv 7,17 \pmod{24}$. 
\item\label{A4det} $f$ is $1849.2.a.g$, $g$ is $2107.1.b.b$, and $p \equiv 11,13 \pmod{24}$. 
\item\label{S3det-CM} $f$ is $63.2.a.b$, $g$ is $1452.1.e.d$ and $p \equiv 5,7 \pmod{12}$.
\item\label{S3exc} $f$ is $189.2.a.f$, $g$ is $3468.1.i.a$ and $p \equiv \pm 5, \pm 43, \pm 67 \pmod{168}$. 
\end{enumerate}
Then $T_{f,g,\mathfrak{p}}$ does not satisfy (wE).}

\rem{In all these counter-examples, $f$ is a non-CM newform of weight $2$ with trivial character with coefficients in a quadratic number field, $g$ is always a newform of weight $1$ with quadratic character, except in case \ref{S3exc} where its character has order $4$. In particular, $\varepsilon_f\varepsilon_g \neq 1$. Moreover, $g$ never has complex multiplication by $\varepsilon_g\varepsilon_f$, so these counter-examples are not covered by Lemma \ref{efeg=1} or the construction of Section \ref{sect-theoD}. Finally, the field extension $H/\Q$ (associated to $f$ in the notation of Section \ref{sect-setup}) is quadratic. The corresponding situations are:
\begin{enumerate}[label=(\alph*),noitemsep]
\item $\rho'_g(G_{\Q}) \cong A_5$, $\varepsilon_g(G_H)=1$.
\item $\rho'_g(G_{\Q}) \cong S_4$, $\varepsilon_g(G_H)=1$.
\item $\rho'_g(G_{\Q}) \cong A_4$, $\varepsilon_g(G_H)=1$.
\item $\varepsilon_g(G_H)=1$ and $g$ has CM by a field distinct from $H$.
\item $i \in \varepsilon_g(G_H)$, $g$ has CM by a field distinct from $H$, $\rho'_g(G_{\Q}) = \rho'_g(G_H)$ is a dihedral group of order $6$. 
\end{enumerate}}

\demo{

\ref{A5det} According to the LMFDB, $\varepsilon_g$ is the quadratic character of conductor $31$, $L=\Q(\sqrt{2})$, $f$ is not CM and has an inner twist by $\varepsilon_g$, so that $F=\Q$ and $H=\Q(\sqrt{-31})$. Since $\varepsilon_g(G_H)=1$, Lemma \ref{efeg=1} shows that any $\sigma \in G_{\Q(\mu_{p^{\infty}})}$ satisfying (wE) cannot lie in $G_{H(\mu_{p^{\infty}})}$. 

Let $D$ be the derived subgroup of $\rho_g(G_{\Q})$. Since $g$ is of $A_5$-type, $D \subset \SL{L_0}$ is a finite group whose projective image is $A_5$. In particular, it has even order, hence contains $-I_2$. Thus we have an exact sequence $1 \rar \{\pm I_2\} \rar D \rar A_5 \rar 1$. Moreover, $\rho_g(G_{\Q})$ is generated by $D$ and its scalar matrices; since $\det{\rho_g(G_{\Q})}=\{\pm 1\}$, its scalar matrices are $\pm iI_2, \pm I_2$, and $D = \rho_g(G_{\Q}) \cap \SL{L_0}$. Since $D \subset \rho_g(G_{H(\mu_{p^{\infty}})}) \subset \rho_g(G_{\Q})$, and $\varepsilon_g(G_H)=1$, $\rho_g(G_{H(\mu_{p^{\infty}})})=D$. 

Let $\sigma \in G_{\Q(\mu_{p^{\infty}})} \backslash G_{H(\mu_{p^{\infty}})}$. By Corollary \ref{papier}, $\sigma$ acts on $T_{f,\mathfrak{p}}$ by a matrix in $\Delta_{\sqrt{2},1/\sqrt{2}}\SL{\Z_p}$ (so, a matrix in $\GL{\Z_p}$ with determinant $2$, divided by $\sqrt{2}$), and, by the above, on $T_g$ by a matrix of $iD$. So we need to show that if $M \in \GL{\Z_p}$ has determinant $2$ and $N \in D$, then $\ker\left[M \otimes N - i\sqrt{2}\right]$ is never a line (over $\overline{\Q_p}$). 

Note that $p \equiv 3,5 \pmod{8}$, so that $M$ always has two distinct eigenvalues $\alpha$ and $\alpha'=2/\alpha$, with $2/\alpha + \alpha \in \Z_p$. Moreover, $p \equiv \pm 1 \pmod{12}$, so that $\sqrt{3} \in \Z_p$, hence $\sqrt{6} \notin \Z_p$. Finally, $p \not\equiv \pm 1 \pmod{5}$, so $\sqrt{5} \notin \Z_p$. 

If $N$ is scalar, it is clear that $M \otimes N$ does not work. 

If the projective image of $N$ is a bitransposition, then $N$ is not scalar, $N^2$ is scalar and $\det{N}=1$, so that $N \sim \Delta_{i,-i}$. Then the four eigenvalues of $M \otimes N$ are $i\alpha,-i\alpha,i\alpha',-i\alpha'$. If any of them is $i\sqrt{2}$, it means that $\alpha$ or $\alpha'$ is $\pm\sqrt{2}$. Hence $\alpha=\alpha'$, a contradiction. 

If the projective image of $N$ is a $3$-cycle, then $N$ is not scalar, $N^3$ is scalar, and $\det{N}=1$, so that $N \sim \Delta_{j,j^2}$. Thus, the four eigenvalues of $M \otimes N$ are $\alpha j, \alpha j^2, \alpha' j, \alpha'j^2$. Suppose that any of them is $i\sqrt{2}$. Up to exchanging $\alpha$ and $\alpha'$, $\alpha$ is $ij^2\sqrt{2}$ or $ij\sqrt{2}$. In either case,

\[\alpha+\alpha'=\pm \sqrt{2}(ij^2-ij)=\pm i\sqrt{2}(j^2-j)=\pm \sqrt{6} \in \Z_p, \text{a contradiction.}\]

If the projective image of $N$ is a $5$-cycle, then $N$ is not scalar, $N^5$ is scalar, and $\det{N}=1$, so that $N \sim \Delta_{u,u^{-1}}$ where $u$ is a primitive fifth root of unity. It is thus enough to show that $u^{-1}\alpha$ cannot equal $i\sqrt{2}$ (up to renaming $u$ and exchanging $\alpha$ and $\alpha'$). If this were the case, we would have $\alpha+2/\alpha=i(u-u^{-1})\sqrt{2} \in \Z_p$, hence $(u-u^{-1})^2 \in \Z_p$. Let $v=(u-u^{-1})^2$: then \[(5+2v)^2=(2u^2+2u^{-2}+1)^2=4u^4+4u^{-4}+1+8+4u^2+4u^{-2}=4(1+u+u^2+u^3+u^4)+5=5,\] this implies that $\sqrt{5} \in \Z_p$, which is a contradiction. 

\ref{S4det} According to the LMFDB, in this case, $L=\Q(\sqrt{6})$, the projective image of $\rho_g$ is $S_4$ and $\varepsilon_g$ is the unique non trivial quadratic character of conductor $43$. Moreover, $f$ has no CM and an inner twist by $\varepsilon_g$, so that $F=\Q$ and $H=\Q(\sqrt{-43})$. Moreover, the coefficient field of $g$ is $\Q(\sqrt{-2})$, so that $g$ has no coefficient equal to $\pm 2i$. Therefore, $\rho_g(G_{\Q})$ does not contain $\pm i I_2$, and its scalar matrices are contained in $\pm I_2$. Now, the derived subgroup of $\rho_g(G_{\Q})$ contains $\rho_g(G_H)$, is contained in $\SL{L_0}$ and has projective image $A_4$, so that it contains an element of even order: so $\rho_g(G_{\Q}) \cap L_0^{\times}I_2=\{\pm I_2\}$, and $\rho_g(G_H)=\rho_g(G_{H(\mu_{p^{\infty}})})$ lives in an exact sequence $1 \rar \{\pm I_2\} \rar \rho_g(G_{H}) \rar A_4 \rar 1$. 

As above, any $\sigma$ satisfying (wE) cannot lie in $G_{H(\mu_{p^{\infty}})}$. Let $\sigma \in G_{\Q(\mu_{p^{\infty}})} \backslash G_{H}$. By the above, $\rho_g(\sigma)$ has determinant $-1$ and $\rho'_g(\sigma)$ has order $2$ or $4$: so $\rho_g(\sigma)$ is similar to $\Delta_{1,-1}$ or to $\Delta_{u,-u^{-1}}$, where $u$ is a primitive $8$-th root of unity. By Corollary \ref{papier}, $\sigma$ acts on $T_{f,\mathfrak{p}}$ by a matrix of $\Delta_{\sqrt{6},1/\sqrt{6}}\SL{\Z_p} = \frac{1}{\sqrt{6}}\GL{\Z_p}_{\det=6}$. 

So we want to prove that if $M \in \GL{\Z_p}$ has determinant $6$ and $N \in \GL{L_0}$ is similar to $\Delta_{1,-1}$ or $\Delta_{u,-u^{-1}}$ for some primitive eighth root of unity $u$, then $\ker(M \otimes N-\sqrt{6})$ is not a line. 

Note that $\sqrt{2} \in \Z_p$ and $\sqrt{3} \notin \Z_p$, so $\sqrt{6} \notin \Z_p$. Hence $M$ has distinct eigenvalues $\alpha,\beta$ with product $6$, so we can rule out $N \sim \Delta_{1,-1}$. Suppose that $N \sim \Delta_{u,-u^{-1}}$ works, where $u$ is a primitive $8$-th root of unity. Then (up to relabelling the eigenvalues) $\alpha=u^{-1}\sqrt{6}$, and $\Tr(u)=\alpha+6/\alpha = \sqrt{6}(u+u^{-1}) \in \Z_p$. Since $(u+u^{-1})^2=u^2+u^{-4}u^2+2=u^2(1-1)+2=2$, $\Tr(u)^2 =3\cdot 4 \notin \Z_p^{\times 2}$, a contradiction. 

\ref{A4det} According to the LMFDB, $L,F,H$ and $\varepsilon_g$ are as above (so $\varepsilon_g$ is the quadratic character of $H/\Q$), and $\rho'_g(G_{\Q}) \cong A_4$, while $\rho_g(G_{\Q})$ has cardinality $48$: its scalar subgroup is $\{\pm iI_2, \pm I_2\}$. Thus, $\rho'_g(G_H)$ is a subgroup of index at most $2$ of $A_4$, so is $A_4$: hence $\rho_g(G_H)$ is a subgroup of $\SL{L_0}$ whose projective image is $A_4$, and $\rho_g(G_{\Q} \backslash G_H) = i\rho_g(G_H)$. 

As above, we want to show that if $M \in \GL{\Z_p}$ has determinant $6$ and $N \in \rho_g(G_H)$ (so $\det{N}=1$), then $\ker\left(M \otimes N - i\sqrt{6}\right)$ is never a line. Write $\alpha$ and $\beta=6/\alpha$ for the two eigenvalues of $M$. Because of our choice of $p$, $\Z_p^{\times 2}$ does not contain $6$ or $2$. 

If $N$ is scalar, then $N = \pm I_2$, and the result is clear. 

If the projective image of $N$ is a bitransposition, $N \sim \Delta_{i,-i}$, and thus (up to relabelling), $i\alpha=\pm i\sqrt{6}$, so $\alpha=\beta=\pm \sqrt{6}$, which we excluded. 

If the projective image of $N$ is a $3$-cycle, then $N \sim \Delta_{j,j^2}$, and thus $\alpha \in \{ij\sqrt{6},ij^2\sqrt{6}\}$. Thus $\Tr(M)=\alpha+6/\alpha = \pm i\sqrt{6}(j-j^2) = \pm 3\sqrt{2}$ does not lie in $\Z_p$, a contradiction.

\ref{S3det-CM} In this case, $L=\Q(\sqrt{3})$, $F=\Q$, $H=\Q(j)$ and $\varepsilon_g$ is the character of $H/\Q$, but $g$ has complex multiplication by $\Q(\sqrt{-11})$. Moreover, the field of coefficients of $g$ is $\Q(j)$ and its projective image is the dihedral group of order $12$. Thus, the only scalar matrices contained in $\rho_g(G_{\Q})$ are $\pm I_2$, and $\det: \rho'_g(G_{\Q}) \rar \{\pm 1\}$ is well-defined, with kernel $\rho'_g(G_H)$ (the latter is contained in the former, and the latter has index at most $2$ in $\rho'_g(G_{\Q})$).   

So $\rho'_g(G_H)$ is a subgroup of index $2$ of $\rho_g(G_{\Q})$ which is not cyclic. Thus, elements of $\rho_g(G_{\Q}) \backslash \rho_g(G_H)$ fall in two categories: matrices $N$ with $\det{N}=-1$ and projective order $2$ (thus $N \sim \Delta_{1,-1}$), and matrices $N$ with $\det{N}=-1$ and projective order $6$. Such a matrix $N \in \rho_g(G_{\Q})$ is of the form $\Delta_{u,-u^{-1}}$, with $u^6=u^{-6}$ and $u^3\neq -u^{-3}$, so $u^6=1$. Since $N^2$ is not scalar, $u \in \{\pm j, \pm j^2\}$. 

Exactly like above, we need to show that for any matrix $M \in \GL{\Z_p}$ with $\det{M}=3$, for any $N$ as above, $\ker(M \otimes N-\sqrt{3})$ cannot be a line. 

If $N \sim \Delta_{1,-1}$, this follows exactly like above from the fact that $3 \notin \Z_p^{\times 2}$. Otherwise, if the statement is false, $M$ has an eigenvalue $\alpha$ such that $\{\alpha j,\alpha j^2\}$ meets $\{\pm \sqrt{3}\}$. Thus, $\Tr(M)=\alpha+3/\alpha=\pm \sqrt{3} (j+j^2)=\pm \sqrt{3} \notin \Z_p$, hence a contradiction. 

\ref{S3exc} In this case, $L=\Q(\sqrt{7}),F=\Q$ and $H=\Q(\sqrt{-3})$; $g$ has complex multiplication by $\Q(\sqrt{-51})$, $\varepsilon_g$ is a character of conductor $51$ and order $4$, and $\rho'_g(G_{\Q}) \cong S_3$. Since $\varepsilon_g(2893)=\pm i$, $\varepsilon_g(G_H)=\{\pm 1, \pm i\}=\varepsilon_g(G_{\Q})$. The LMFDB shows that $\rho_g(G_{\Q})$ is isomorphic to $S_3 \times \Z/8\Z$. Its center (which is scalar, by Schur's lemma), is thus the set of $8$-th roots of unity, and $\rho_g(G_{\Q}) \cap \SL{L_0} \rar \rho'_g(G_{\Q})$ is onto with kernel $\{\pm I_2\}$. 

The LMFDB data stipulates that $\ker{\rho'_g} \not\subset G_H$, so $\rho'_g(G_H) = \rho'_g(G_{\Q})$. Using MAGMA \cite{magma}, we can alse see that $\ker{\rho_g} \subset G_H$, so that $[\rho_g(G_{\Q}):\rho_g(G_H)]=2$. 

For the sake of easier notation, write $C \subset \rho_g(G_{\Q})$ for the pre-image of the non-trivial elements of the cyclic subgroup of order $3$ of $\rho'_g(G_{\Q})$. 

Since $\rho_g(G_H)$ and $\rho_g(G_{\Q})$ have the same determinant and projective image, we can split elements of $\rho_g(G_H)$ in six categories: $\{\pm i I_2, \pm I_2\}$, scalar matrices of order $8$, matrices in $C$ with determinant $\pm 1$, matrices in $C$ with determinant $\pm i$, matrices outside $C$ with determinant $\pm 1$, matrices outside $C$ with determinant $\pm i$. Then $\rho_g(G_H)$ is exactly made with the matrices of the first, third, and last categories. 

As above, we need to show that if $(M,N) \in \SL{\Z_p} \times \rho_g(G_H)$, $\ker\left[M \otimes N -1\right]$ is never a line, and that, if $(M,N) \in \GL{\Z_p} \times \left[\rho_g(G_{\Q})\backslash \rho_g(G_H)\right]$ is such that $\det{M}=7$, then $\ker\left[M \otimes N-\sqrt{7}\right]$ is never a line again. 

Note that by our choice of $p$, $\Z_p$ contains none of $\sqrt{2},\sqrt{3},\sqrt{7}$, and $7! \in \Z_p^{\times}$. 

Let $M \in \SL{\Z_p}$ (with eigenvalues $\alpha$ and $1/\alpha$) and $N \in \rho_g(G_H)$ be such that $\ker(M \otimes N-1)$ is a line. Clearly, $N$ is not scalar, and $\det{N} \neq 1$. Note that $\Tr(M)=\alpha+1/\alpha \in \Z_p$. 

Suppose that $N\in C$. Then $\det{N}=-1$, so that $N \sim \pm i\Delta_{j,j^2}$. Thus (up to exchanging $\alpha$ and $1/\alpha$) $\alpha \in \{\pm i j, \pm i j^2\}$. Thus $\Tr(M)=\alpha+1/\alpha = \pm i (j-j^2) = \pm \sqrt{3} \notin \Z_p$, so we get a contradiction. 

So $N \notin C$: thus $\det{N}=\pm i$ and $N^2$ is scalar. So $N \sim \Delta_{u,-u}$, where $-u^2 = \pm i$, thus $u$ is a primitive $8$-th root of unity. Then $\alpha = \pm u^{-1}$, so $\Tr(M)=\alpha+1/\alpha=\pm (u+1/u) = \pm \sqrt{2} \notin \Z_p$, so we get a contradiction.

Let now $M \in \GL{\Z_p}, N \in \rho_g(G_{\Q}) \backslash \rho_g(G_H)$ such that $\det{M}=7$ and $\ker(M \otimes N-\sqrt{7})$ is a line. Let $\alpha$ and $7/\alpha$ be the distinct eigenvalues of $M$ (since $\sqrt{7} \notin \Z_p$). Clearly, $N$ is not scalar, and $\det{N} \neq 1$. 

Suppose that $N \in C$: then $\det{N}=\pm i$, so that $N \sim u\Delta_{j,j^2}$ for some primitive $8$-th root of unity $u$. So $\alpha \in \{\sqrt{7}ju^{-1}, \sqrt{7}j^2u^{-1}\}$. Thus there is some primitive $24$-th root of unity $\omega$ such that $\alpha=\sqrt{7}\omega$. Then $\Tr(M)^2=7(\omega+1/\omega)^2 \in \Z_p$, so $\omega^2+\omega^{-2} \in \Z_p$. Now, $(\omega^2+\omega^{-2})^2=\omega^4+2+\omega^{-4}=2-(\omega^8+\omega^{-8})=2-(-1)=3$ (since $\omega^8, \omega^{-8}$ are the two primitive third roots of unity, and $\omega^{12}=-1$): but $3 \notin \Z_p^{\times 2}$, a contradiction. 

So $N \notin C$ and $\det{N}=-1$, then $N \sim \Delta_{1,-1}$, so we get $\alpha=7/\alpha$, a contradiction. 
}

The proposition below, on the other hand, shows that, even in the exceptional case of Proposition \ref{special2}, the representation can still have Euler-adapted Galois image. 

\prop[special2-exc-nothold]{Let $f \in \mathcal{S}_2(\Gamma_0(63)), g \in \mathcal{S}_1(\Gamma_1(675))$ be newforms with respective labels $63.2.a.b$ and $g$ be $675.1.g.a$. Then $T_{f,g,\mathfrak{p}}$ has Euler-adapted Galois image for every good prime $\mathfrak{p}$. }

\rem{In this case, $L=\Q(\sqrt{3}), F=\Q, H=\Q(\sqrt{-3})$, and $\varepsilon_g$ is a character of conductor $5$ and order $4$. Moreover, $\rho'_g(G_{\Q})$ is the dihedral group of order $12$, and $g$ has complex multiplication by $\Q(\sqrt{-15}) \neq H$. Furthermore, the number field LMFDB data shows that $\ker{\rho'_g} \subset G_H$, so $[\rho'_g(G_{\Q}):\rho'_g(G_H)]=[\rho_g(G_{\Q}):\rho_g(G_H)]=2$. Since $a_{139}(g)=\pm 2i$, $iI_2 \in \rho_g(G_H) \subset \rho_g(G_{\Q})$. Since the coefficient field of $g$ is $\Q(i,\sqrt{6})$, it does not contain the primitive $8$-th root of unity, and the scalar matrices of $\rho_g(G_{\Q})$ are exactly $\pm iI_2, \pm I_2$. Finally, $\varepsilon_g(G_H)=\varepsilon_g(G_{\Q})=\{\pm i, \pm 1\}$, since $\varepsilon_g$ and $H$ have prime conductors.

It seems that we are in the ``exceptional'' case of Proposition \ref{special2}. However, it turns out that (sE) still holds for every good prime $\mathfrak{p}$.}

\demo{Let us fix such a prime $\mathfrak{p}$. By the proof of Proposition \ref{special2}, since $i I_2 \in \rho_g(G_H)$, we are done if $\rho_g(G_H)$ contains an element whose projective image has order $2$ and determinant $\pm 1$. So we may assume that the elements of $\rho_g(G_H)$ whose projective image lies outside the maximal cyclic subgroup have determinant $\pm i$. 

Let $N$ be such a matrix: then $N \sim \Delta_{u,-u}$, for some $u$ such that $u^2=\pm i$. Since $(u+1/u)^2=2$, if $\sqrt{2} \in \Z_p$, then there exists $M \in \SL{\Z_p}$ similar to $\Delta_{u,u^{-1}}$. So there exists $\sigma \in G_{H(\mu_{p^{\infty}})}$ such that $\sigma$ acts on $T_{f,\mathfrak{p}}$ (resp. $T_g$) by $M$ (resp. $N$). Then $M \otimes N$ has four eigenvalues: $1,-1,i,-i$. So $\sigma$ shows that (sE) holds. Hence, we now assume that $\sqrt{2} \notin \Z_p$. 

Suppose that the restriction of the determinant map $\delta: \rho'_g(G_{\Q}) \rar \varepsilon_g(G_{\Q})/\{\pm 1\} \cong \F_2$ to the cyclic subgroup is trivial. Since the complex conjugation also lies in $\ker{\delta}$, $\delta$ is trivial, and thus is not onto, a contradiction. Hence, $\delta$ maps the unique element of order $2$ of the cyclic subgroup (which does not belong to $\rho'_g(G_H)$) to the coset $\{\pm i\}$. Considering the complex conjugation, we see that $\rho_g(G_{\Q}) \backslash \rho_g(G_H)$ contains elements of projective order $2$ with determinant $-1$ and $i$. 

Let $\sigma \in G_{\Q(\mu_{p^{\infty}})} \backslash G_H$: then Corollary \ref{papier} shows that $\sigma$ acts on $T_{f,\mathfrak{p}}$ as $M/\sqrt{3}$, where $M \in \GL{\Z_p}$ has determinant $3$. If $\sqrt{3} \in \Z_p$, then we can find $\sigma \in G_{\Q(\mu_{p^{\infty}})} \backslash G_H$ acting on $T_{f,\mathfrak{p}}$ by $U$, and on $T_g$ by a matrix of determinant $-1$ whose projective image has order $2$ (hence is similar to $\Delta_{1,-1}$), and we conclude like in the proof of Proposition \ref{special2}.  

So we may assume that $\sqrt{2},\sqrt{3} \notin \Z_p$, and therefore that $\sqrt{6} \in \Z_p$. Let $N \in \rho_g(G_{\Q}) \backslash \rho_g(G_H)$ have determinant $i$ and be such that its projective image has order $2$: then $N \sim \Delta_{u,-u}$ for some primitive $8$-th root of unity $u$. Then $\sqrt{3}u+\sqrt{3}u^{-1}=\pm \sqrt{6} \in \Z_p$: hence there exists $M \in \GL{\Z_p}$ with eigenvalues $\sqrt{3}u,\sqrt{3}u^{-1}$. Thus, there exists $\sigma \in G_{\Q(\mu_{p^{\infty}})} \backslash G_H$ such that $\sigma$ acts on $T_{f,\mathfrak{p}}$ (resp. $T_g$) by $M/\sqrt{3}$ (resp. $N$). Hence, the eigenvalues of $\sigma$ acting on $T_{f,g,\mathfrak{p}}$ are $\pm 1, \pm i$ as above, so (sE) holds.   }

\rem{This argument can be somewhat generalized. If $\rho_g(G_{\Q})$ is dihedral with order $4n$ with $n > 1$ odd, if $[\rho_g(G_{\Q}):\rho_g(G_H)]=[\rho'_g(G_{\Q}):\rho'_g(G_H)]=2$ and $H$ is not the CM field of $g$, if $\varepsilon_g(G_H)=\varepsilon_g(G_{\Q})=\{\pm 1, \pm i\}$, and there is some $\sigma \in G_{\Q} \backslash G_H$ such that $\varepsilon_f(\sigma)=1$ and $\rho_g(\sigma) \notin \rho_g(G_H)$, then (sE) will hold.}

%% file: algeom-0.tex
\chapter{Prerequisites from algebraic geometry}
\label{alg-geom-prereq}

This Appendix collects the background from algebraic geometry which is useful to the discussion of moduli schemes and Galois twists in Chapter \ref{moduli-spaces}. This background comes in four rather different collections of results. The first part is about descent and its various applications. The main reference here is the online textbook Stacks Project \cite{Stacks}. The second part collects results on one possible compactification method for schemes that are finite over the affine line over a suitable basis. This is necessary to discuss compatified modular curves as in Chapter \ref{moduli-spaces}. The third part is about the classical properties of abelian varieties and abelian schemes. I am grateful to Discord user topos\_theory\_egirl for making me aware of the article \cite{Anantharaman} and for correcting my simplistic intuition for abelian schemes. The last part of this appendix discusses Jacobians. It is well-known (while it does not seem to be directly stated, it is easily deduced from various results in \cite[Chapter 8]{BLR}) that relative proper smooth curves admit well-behaved Jacobians when their geometric fibres are connected. But it is necessary, in our setting, to slightly adapt this material in order to deal with relative curves whose fibres are not geometrically connected. I would like to thank Qing Liu for his help with this adaptation.

\section{Results from descent}

\subsection{Sheaves in the fpqc topology}
Recall the following definition \cite[Definition 022B]{Stacks}: let $T$ be a scheme, and $(\pi_i)_{i \in I}$ be a family of flat scheme morphisms to $T$, with $\pi_i: T_i \rar T$. This family is a \emph{fpqc covering} if for every affine open subscheme $U \subset T$, there is a finite set $F$, a map $j: F \rar I$, and affine open subsets $W_f \subset T_{j(f)}$ for each $f \in F$, such that $U=\bigcup_{f \in F}{\pi_{j(f)}(W_f)}$. 

\theo[morsheaf]{(see \cite[Lemma 023Q]{Stacks}) Let $S$ be a scheme, and $T$ be a $S$-scheme. Then the functor $V \in \mathbf{Sch}_S \longmapsto \mrm{Mor}_S(V,T)$ satisfies the sheaf property for fpqc covers.}

\demo{Given a scheme $A$, let $F_A$ be the pre-sheaf $T \in \mathbf{Sch} \longmapsto \mrm{Mor}(T,A)$. By \emph{loc. cit.}, it satisfies the sheaf property for fpqc covers. 
The forgetful functor $\omega_S: \mathbf{Sch}_S \rar \mathbf{Sch}$ maps an fpqc cover to an fpqc cover, so the presheaves $F'_T := F_T \circ \omega_T, F'_S := F_S \circ \omega_S$ on $\mathbf{Sch}_S$ satisfy the sheaf property for fpqc covers. Let $G$ be the sheaf on $\mathbf{Sch}_S$ given by the following: for all $S$-schemes $T$, $G(T)$ has a single element (and the restriction map is the trivial one). There are natural morphisms of sheaves $G \rar F'_S$ (sending the point of $G(V)$ to the structure map $V \rar S$), $F'_T \rar F'_S$ (composing with the structure map $T \rar S$), and the functor that we are studying is exactly $V \longmapsto F'_T(V) \times_{F'_S(V)} G(V)$, which formally satisfies the sheaf property for fpqc covers.}

\prop[bi-morsheaf]{Let $S$ be a scheme, and $U,V$ be two $S$-schemes. The functor \[T \in \mathbf{Sch}_S \longmapsto \mrm{Mor}_T(U_T,V_T)\] satisfies the sheaf property for fpqc covers. The subfunctor mapping a $S$-scheme $T$ to the set of $T$-isomorphisms $U_T \rar V_T$ also satisfies the sheaf property for fpqc covers.}

\demo{The first functor is isomorphic to $T \in \mathbf{Sch}_S \longmapsto \mrm{Mor}_S(U_T,V)$. It is formal to check that the functor $T \in \mathbf{Sch}_S \longmapsto U\times_S T \in \mathbf{Sch}_S$ preserves fpqc covers, and we are done.

For the second functor, it is clear from the above that it is separated -- two sections that agree on a fpqc cover are equal. Let $\pi_i: T_i \rar T$ for $i \in I$ be an fpqc cover, and $\alpha: U_T \rar V_T$ be a $T$-morphism such that for each $i \in I$, $\alpha\times_T T_i: U_{T_i} \rar V_{T_i}$ is an isomorphism; it is enough to prove that $\alpha$ is an isomorphism. This is exactly \cite[Lemma 02L4]{Stacks} since the $V_{T_i} \rar V_T$ are a fpqc cover, and $\alpha \times_T T_i$ is the base change of $\alpha$ by $V_{T_i} \rar V_T$.}

\prop[group-morsheaf]{Let $S$ be a scheme, and $U,V$ be two $S$-group schemes. The functors 
\begin{align*}
\underline{Mor}_{S-Gp}(U,V)&: T \in \mathbf{Sch}_S \longmapsto \mrm{Mor}_{T-Gp}(U_T,V_T),\\
\underline{Iso}_{S-Gp}(U,V)&: T \in \mathbf{Sch}_S \longmapsto \mrm{Iso}_{T-Gp}(U_T,V_T)
\end{align*}
 satisfy the sheaf property for fpqc covers.}

\demo{We deal with the first functor first. Let $F_1=\underline{Mor}_S(U,V), F_2=\underline{Mor}_S(U \times U,V)$ be fpqc sheaves. We can define two morphisms $\alpha, \beta: F_1 \rar F_2$ as follows: given a $S$-scheme $T$ and a $T$-map $\varphi: U_T \rar V_T$, $\alpha(\varphi)=m_{V_T} \circ (\varphi \times \varphi)$, where $m_{V_T}$ is the multiplication of the group scheme $V_T$ over $T$, and $\varphi \times \varphi: U_T \times_T U_T \rar V_T \times_T V_T$. We define $\beta(\varphi)=\varphi \circ m_{U_T}$. Then given $u \in \underline{Mor}_{S}(U,V)(T)$, $u \in \underline{Mor}_{S-Gp}(U,V)(T)$ if and only if $\alpha(u)=\beta(u)$. The conclusion follows formally. 

The second functor corresponds to the intersection in $\underline{Mor}_S(U,V)$ of the two subsheaves $\underline{Mor}_{S-Gp}(U,V)$ and $\underline{Iso}_S(U,V)$: it is thus a fpqc sheaf.  }

\prop[desc-affine-schemes]{Let $S$ be a scheme, and $T \rar S$ be a fpqc cover. Let $F$ be a fpqc sheaf on $\mathbf{Sch}_S$ such that its restriction to the category $\mathbf{Sch}_T$ is representable by an affine $T$-scheme $X$\footnote{That is, the map $X \rar T$ is affine.}. Then $F$ is representable by an affine $S$-scheme $X'$, and $X \simeq X' \times_S T$.}

\demo{ 

\emph{Step 1: $S=\Sp{A}$} \\ After refining $T$, we may assume that $T=\Sp{B}$, where $A \rar B$ is faithfully flat, and $X=\Sp{C}$, where $C$ is some $B$-algebra. The hypothesis implies that there exists some $\alpha \in F(X)$ satisfying the following property: for every $T$-scheme $Z$ with structure map $g: Z \rar T$, the map $f \longmapsto \mrm{Mor}_T(Z,X) \rar F(f)\alpha \in F(Z)$ is an isomorphism. 

Let $f_1 \in \mrm{Mor}_{T}(X \times_S T,X)$ (where $X \times_S T$ is a $T$-scheme by the second projection) be the unique $T$-homomorphism such that $F(f_1)\alpha=F(X \times_S T\overset{\mrm{pr}_1}{\rar} X)\alpha$. This defines a map $\pi=(\sigma,f_1): X \times_S T \rar T \times_S X$ of $T \times_S T$-schemes, where $\sigma: X \rar T$ is the structure map. 

We define the map of $T$-schemes $f_2: T \times_S X \rar X$ (where $T \times_S X$ is a $T$-scheme by first projection) such that $F(f_2)\alpha=F(T \times_S X \overset{\mrm{pr}_2}{\rar} X)\alpha$. We then define $\pi'= (f_2,\sigma): T \times_S X \rar X \times_S T$, which is a map of $T \times_S T$-schemes. 

Then $\pi \circ \pi'$ is an endomorphism of the $T \times_S T$-scheme $T \times_S X$, so its first component is the first projection $T \times_S X \rar T$. Its second component $p_2$ is a morphism of $T$-schemes (on the source, through the $T$-scheme structure on $X$); furthermore, we have \[F(p_2)\alpha=F(\pi')F(f_1)\alpha=F(\pi')F(\mrm{pr}_1)\alpha=F(\mrm{pr}_1\circ\pi')\alpha=F(f_2)\alpha=F(\mrm{pr}_2)\alpha.\] Since $p_2,\mrm{pr}_2$ are two maps of $T$-schemes $T \times_S X \rar X$ \emph{for the same $T$-structure}, it follows that $p_2=\mrm{pr}_2$, and $\pi \circ \pi'$ is the identity. 

One similarly proves that $\pi'\circ\pi$ is the identity, hence $\pi,\pi'$ are isomorphisms of $T \times_S T$-schemes. 

Now, let us show that the following two compositions of $T^{\times_S 3}$-morphisms are the same: 

\begin{align*}
T \times_S T \times_S X \overset{\mrm{id} \times_S \pi'}{\rar} T \times_S X \times_S T \times \overset{\pi' \times \mrm{id}}{\rar} X \times_S T \times_S T,\\
\pi'_{13} = (f_2 \circ \mrm{pr}_{13},\mrm{id},\sigma \circ \mrm{pr}_{13}): T \times_S T \times_S X \rar X \times_S \times_S T
\end{align*}

It is clear that the projections on the second and third components are the same (since these are morphisms of $T^{\times_S 3}$-schemes), so all we need to do is prove that the first components of the two morphisms are equal. It is simple to compute that these first components are respectively $f_2 \circ (\mrm{id} \times f_2)$ and $f_2 \circ \mrm{pr}_{13}$. 

Since these morphisms are morphisms of $T$-schemes for the same structure (projection on the first component on the source), it is thus enough to prove that $F(f_2 \circ (\mrm{id} \times f_2))\alpha=F(f_2 \circ \mrm{pr}_{13})\alpha$.  

Now, 
\noindent\begin{align*}%
F(f_2 \circ \mrm{pr}_{13})\alpha&=F(\mrm{pr}_{13})F(f_2)\alpha=F(\mrm{pr}_{13})F(\mrm{pr}_2)\alpha=F(\mrm{pr}_2\circ \mrm{pr}_{13})\alpha\\
&=F(\mrm{pr}_3: T \times_S \times_S X \rar X)\alpha,\\
F(f_2 \circ (\mrm{id} \times f_2))\alpha&=F(\mrm{id} \times f_2)F(f_2)\alpha=F(\mrm{pr}_2\circ(\mrm{id} \times f_2))\alpha=F(f_2 \circ \mrm{pr}_{23})\alpha\\
&=F(\mrm{pr}_{23})F(f_2)\alpha=F(\mrm{pr}_{23})F(\mrm{pr}_2)\alpha=F(\mrm{pr}_3)\alpha.
\end{align*} 

We can thus apply \cite[Proposition 023N]{Stacks}\footnote{This proposition is formulated for modules over rings, but the explicit construction shows that it works for algebras as long as the descent datum is also a ring homomorphism, hence it works for affine schemes as well.}. This Proposition implies that we can assume that, for some $T$-isomorphism $(\chi,\sigma): X \rar Y \times_S T$ we have a commutative diagram of $T \times_S T$-schemes:

\[
\begin{tikzcd}[ampersand replacement=\&]
T \times_S X \arrow{r}{\pi'}\arrow{d}{\mrm{id}\times (\chi,\sigma)} \& X \times_S T \arrow{d}{(\sigma,\chi) \times \mrm{id}}\\
T \times_S (Y \times_S T) \arrow{r}{\mrm{can}} \&  (T \times_S Y) \times_S T
\end{tikzcd}
\]

Let $\alpha'=F((\chi,\sigma)^{-1})\alpha \in F(Y \times_S T)$. Then the following statement holds:

For any $T$-scheme $g: Z \rar T$, $f \in \mrm{Mor}_S(Z,Y) \longmapsto F((f,g))\alpha' \in F(Z)$ is a bijection. $(\ast)$

Moreover, the following diagram is easily checked to commute (the indices on $T$ just describe what the lower arrows do):  

\[
\begin{tikzcd}[ampersand replacement=\&]
X\arrow{d}{(\sigma,\chi)}\& \arrow{l}{\mrm{pr}_2} T \times_S X \arrow{r}{f_2}\arrow{d}{\mrm{id}\times (\chi,\sigma)} \& X \arrow{d}{(\sigma,\chi)}\\
Y \times_S T_2 \& T_1 \times_S (Y \times_S T_2) \arrow{l}{}\arrow{r}{} \&  Y \times_S T_1
\end{tikzcd}
\] 

Therefore, the pull-backs of $\alpha'$ under both projections $Y \times_S T \times_S T \rar Y \times_S T$ are the same. Hence, $\alpha'$ is the pull-back of a class $\alpha_0 \in F(Y)$. Thus, we have a natural homomorphism $f \in \mrm{Mor}_S(Z,Y) \longmapsto F(f)\alpha_0 \in F(Z)$ of fpqc sheaves on $\mathbf{Sch}_S$. By $(\ast)$, the application is an isomorphism whenever there is a map $g: Z \rar T$ of $S$-schemes. Since every $S$-scheme has an fpqc cover by $T$-schemes, the above morphism is an isomorphism of fpqc sheaves. \\

\emph{Step 2: Glueing affine data} 

\noindent
By Step 1, for every affine open subset $U \subset S$, there exists an affine $U$-scheme $X'_U$ and some $\alpha_U \in F(X'_U)$ such that for every $U$-scheme $V$, $f \in \mrm{Mor}_U(V,X'_U) \longmapsto F(f)\alpha_U \in F(V)$ is a bijection. Moreover, for every inclusion $V \subset U$ of affine open subsets of $S$, there is a unique map $j_{VU}: X'_V \rar X'_U$ of $U$-schemes such that $F(j_{VU})\alpha_U=\alpha_V$.  

Let us prove that every $j_{VU}$ is an open immersion. There is a class $\beta_T \in F(X)$ such that for every $T$-scheme $Z$, $f \in \mrm{Mor}_T(Z,X) \longmapsto F(f)\beta_T \in F(Z)$ is an isomorphism. In particular, for any affine open $U \subset S$, there is a unique morphism of $T$-schemes $f_U: X'_U \times_S T \rar X$ such that $F(f_U)\beta_T$ is $F(X'_U \times_S T \rar X'_U)\alpha_U$. Since $f_U$ factors through $f_U^0: X'_U \times_S T \rar X \times_S U$ which is an isomorphism of $T$-schemes (by relating its functor of $T$-points to $F$), hence $f_U$ is an open immersion. 

Given affine open subsets $V \subset U$ of $S$, one easily sees that $f_U \circ (j_{VU}\times \mrm{id}_T)=f_V$. Therefore $j_{VU} \times \mrm{id}_T$ is an open immersion, hence $j_{VU}$ is one too by \cite[Lemma 02L3]{Stacks}. 

Finally, the definition implies easily that, given affine open subsets $W \subset V \subset U$ of $S$, $j_{VU} \circ j_{WV} = j_{WU}$, so that the $(X'_U,j_{VU})$ is a glueing data in the sense of \cite[Section 01JA]{Stacks}. Therefore, by \emph{loc. cit.}, the $X_U$ glue to an $S$-scheme $X'$, the $\alpha_U$ to some $\alpha \in F(X')$, and the $f_U$ to some morphism $f: X' \times_S T \rar X$, which is an isomorphism after base changing to any affine open subset of $S$, hence is an isomorphism, and such that $F(f)\beta_T=F(X' \times_S T \rar X')\alpha$. 

Exactly like in Step 1, this $\alpha$ produces a homomorphism of fpqc sheaves $\phi: \underline{\mrm{Mor}}_{S}(-,X') \rar F$ on $\mathbf{Sch}_S$, and $\phi(Z)$ is an isomorphism for every $S$-scheme admitting a map to $T$. Since every $S$-scheme has a fpqc cover by a $T$-scheme, it follows that $\phi$ is an isomorphism.  
}

\bigskip

We conclude this section by discussing some examples of ascending and descending properties. 

\prop[fpqc-descent-prop]{Let $\mathcal{P}$ be one of the following properties of morphisms of schemes: isomorphism, surjective, affine, quasi-compact, closed immersion, open immersion, locally of finite type, locally of finite presentation, locally quasi-finite, finite, flat, \'etale, smooth of relative dimension $d$, Cohen-Macaulay of relative dimension $d$, proper. 

Let $f: X \rar Y$ be a morphism of schemes and assume that there is a fpqc cover $Y_i \rar Y$ such that the base change of $f$ by any $Y_i \rar Y$ satisfies property $\mathcal{P}$. Then $f$ satisfies property $\mathcal{P}$.}

\demo{See \cite[Section 02YJ]{Stacks}.}

\prop[descending-properties]{Let $R \rar S$ be a faithfully flat morphism of local rings. Suppose that $S$ is Noetherian (resp. Noetherian and $(R_k)$, resp. Noetherian and $(S_k)$, resp. reduced, resp. normal, resp. Cohen-Macaulay, resp. regular), then $R$ is Noetherian (resp. Noetherian and $(R_k)$, resp. Noetherian and $(S_k)$, resp. reduced, resp. normal, resp. Cohen-Macaulay, resp. regular).}

\demo{See \cite[Section 033D]{Stacks}.}

\prop[ascending-properties]{Let $R \rar S$ be a flat morphism of rings. Then the following properties hold:
\begin{itemize}[noitemsep,label=$-$]
\item If $R,S$ are Noetherian, $R$ is $(R_k)$ (resp. $(S_k)$, resp. reduced, resp. normal) and every $S \otimes \kappa(\mathfrak{p})$ (for $\mathfrak{p} \in \Sp{R}$) is $(R_k)$ (resp. $(S_k)$, resp. reduced, resp. normal), then $S$ is $(R_k)$ (resp. $(S_k)$, resp. reduced, resp. normal).
\item If $S$ is a smooth $R$-algebra and $R$ is reduced (resp. normal, resp. regular), then $S$ is reduced (resp. normal, resp. regular). 
\end{itemize}
}

\demo{See \cite[Section 0336]{Stacks}.}

\subsection{Galois twists of schemes}

\lem[galois-semi-lin]{Let $V$ be a vector space over a field $L$, and $\mrm{Gal}(L/k)$ act semi-linearly on $V$, where $L/k$ is finite Galois. Then the natural map $V^{\mrm{Gal}(L/k)} \otimes_k L \rar V$ is an isomorphism of $L$-vector spaces.}

\demo{To ease notation, write $G=\mrm{Gal}(L/k)$. Let $(w_i)_i$ be a finite family of vectors of $V^{G}$ which is $k$-free. Let us show that this family is free over $L$ in $V$. Indeed, suppose that we have $\lambda_i \in L$ such that $\sum_i{\lambda_iw_i}=0$. Multiplying this equation by a fixed $\alpha \in L$ and applying $\sum_{g \in G}{g}$, we find that $\sum_i{\tau(\alpha\lambda_i)w_i}=0$, where $\tau: L \rar k$ is the trace map. Since the $w_i$ are $k$-free, we find that for all $i$ and all $\alpha \in L$, $\tau(\alpha\lambda_i)=0$. This implies that every $\lambda_i=0$, and we are done. 
Thus the map $V^G \otimes_k L \rar V$ is injective. 

By the normal basis theorem, we know that there are $\lambda_g \in L$ for each $g \in G$ such that for all $g,h \in G$, $g(\lambda_h)=\lambda_{gh}$, and such that $(\lambda_g)_g$ is a basis of $L$ as a $k$-vector space. Let us show that the family of the vectors $w_h := (\lambda_{gh})_{g \in G}$ (for $h \in G$) is a basis of $L^G$. Since this family has the correct cardinality, it is enough to show that it is free. So let $\mu_h \in L$ for each $h \in G$ such that $\sum_h{\mu_h\lambda_{gh}}=0$ for each $g \in G$. Then, for any $\alpha \in L$ and $g \in G$, we have $\sum_{h \in G}{g^{-1}(\alpha\mu_h)\lambda_h}=0$. By taking the sum over $g \in G$, we find that $\sum_{h \in G}{\tau(\alpha\mu_h)\lambda_h}=0$. Thus every $\tau(\alpha\mu_h)$ vanishes, thus every $\mu_h$ vanishes. 

In particular, there is a $L$-linear combination of the $w_h$ which is the indicator function of the unit element of $G$. It follows that there is a $L$-linear combination of the $v_h := \sum_{g \in G}{\lambda_{gh}g(v)} \in V^G$ which is mapped to $v$. Hence the map is onto. 
}

\lem[galois-descent-cor]{Let $L/K$ be a finite Galois extension, and $U,V$ be two $K$-schemes. There is a natural action of $\mrm{Gal}(L/K)$ on $\mrm{Mor}_L(U_L, V_L)$. The base change map 
\[f \in \mrm{Mor}_K(U,V) \longmapsto f \times \mrm{id}_L \in \mrm{Mor}_L(U_L,V_L)^{\mrm{Gal}(L/K)}\] is an isomorphism. }

\demo{This map is clearly well-defined. Since, $\Sp{L} \rar \Sp{K}$ is a fpqc cover, the base change map is injective by Proposition \ref{bi-morsheaf}. To prove surjectivity, it is enough to show that if $f \in \mrm{Mor}_L(U_L,V_L)^{\mrm{Gal}(L/K)}$, and $\pi_1,\pi_2: \Sp{L} \times_k \Sp{L} \rar \Sp{L}$ are the projections, then $\pi_1^{\ast}f=\pi_2^{\ast}f$. Let us assume that there is an isomorphism of $L$-algebras $L \otimes_k L \simeq L^{\oplus \mrm{Gal}(L/K)}$ (where $L$ acts on the left side by multiplication on the right factor), so it is enough to show that for any $\tau: \Sp{L} \rar \Sp{L}\times_k \Sp{L}$, one has $\tau^{\ast}\pi_1^{\ast}f=\tau^{\ast}\pi_2^{\ast}f$. In other words, one only needs to show that given $\sigma \in \mrm{Gal}(L/K)$ inducing a map $\underline{\sigma}: \Sp{L} \rar \Sp{L}$, $\underline{\sigma}^{\ast}f$ does not depend on the choice of $\sigma$. This is exactly the assumption. 

So we only need to prove that there is a $L$-algebra isomorphism $\mu: L \otimes_k L \simeq L^{\oplus \mrm{Gal}(L/K)}$. Indeed, for each $\sigma \in \mrm{Gal}(L/K)$, define $\mu_{\sigma}: a \otimes b \longmapsto \sigma(a)b$: this is clearly a surjective homomorphism of $L$-algebras for the given structure on $L \otimes_k L$: let us prove that $\mu=\bigoplus_{\sigma \in \mrm{Gal}(L/K)}{\mu_{\sigma}}$ works. By computing the dimensions, it is enough to show that $\mu$ is surjective. 

It this wasn't the case, there would be nonzero scalars $\lambda_{\sigma} \in L$ for each $\sigma \in \mrm{Gal}(L/K)$ such that $\sum_{\sigma}{\lambda_{\sigma}\mu_{\sigma}}=0$. In other words, one would have, for every $a \in L$ and $n \geq 0$, $\sum_{\sigma}{\lambda_{\sigma}\sigma(a^n)}=0$. So, for every polynomial $P \in L[X]$ and every $a \in L$, we would have $\sum_{\sigma}{\lambda_{\sigma}P(\sigma(a))}=0$. If $a \in L$ is chosen to be a primitive element, the $\sigma(a)$ are pairwise distinct, and it follows that all the $\lambda_{\sigma}$ are zero. 
}

\cor[inf-galois-desc]{Let $L/K$ an algebraic Galois extension and $U,V$ be two quasi-compact $K$-schemes. There is a natural action of $\mrm{Gal}(L/K)$ on $\mrm{Mor}_L(U_L, V_L)$. The base change map \[f \in \mrm{Mor}_K(U,V) \longmapsto f \times \mrm{id}_L \in \mrm{Mor}_L(U_L,V_L)^{\mrm{Gal}(L/K)}\] is an isomorphism.}

\demo{As above, the map is well-defined and injective. Let $f: U \times_K L \rar V \times_K L$ be a morphism of $L$-schemes which is $\mrm{Gal}(L/K)$-invariant. Because $V,U$ are quasi-compact, there is a finite Galois subextension $K \subset L' \subset L$ such that $f$ is a base change of some $f' \in \mrm{Mor}_{L'}(U_{L'},V_{L'})$. Since the base change map $\mrm{Mor}_{L'}(U_{L'},V_{L'}) \rar \mrm{Mor}_{L}(U_{L},V_{L})$ is injective by Proposition \ref{bi-morsheaf} and equivariant for $\mrm{Gal}(L/K)$, $f'$ is fixed under the action of $\mrm{Gal}(L'/K)$, so is the base change of some $f_0 \in \mrm{Mor}_K(U,V)$.}

\prop[cocycle-twist]{Let $k$ be a field with separable closure $k_s$. For $\sigma \in \mrm{Gal}(k_s/k)$, denote by $\underline{\sigma}$ the associated automorphism of $\Sp{k_s}$, so that $\mrm{Gal}(k_s/k)$ acts on the right on $\Sp{k_s}$. Let $X$ be a $k$-scheme and assume that $X$ is quasi-projective over some $k$-algebra $A$\footnote{By quasi-projective, we mean that there is a $k$-immersion $\iota: X \rar \mathbb{P}^n_A$ for some $n \geq 1$.}. 
Let $G$ be a finite group acting on $X_{k_s}$ by $k_s$-scheme automorphisms, and $\rho: \mrm{Gal}(k_s/k) \rar G$ be a continuous map satisfying the following cocycle condition: for all $\sigma,\sigma' \in \mrm{Gal}(k_s/k)$, 
\[\rho(\sigma\sigma')=\rho(\sigma)\circ (\mrm{id},\underline{\sigma^{-1}})\circ \rho(\sigma') \circ (\mrm{id},\underline{\sigma}).\] 
Then, there is a unique (up to unique isomorphism) couple $(X_{\rho},j_X)$ satisfying the following properties:
\begin{itemize}[noitemsep,label=$-$]
\item $X_{\rho}$ is a $k$-scheme,
\item $j_X: X_{\rho} \times_k \Sp{k_s} \rar X \times_k \Sp{k_s}$ is an isomorphism of $k_s$-schemes,
\item for any $\sigma \in \mrm{Gal}(k_s/k)$, the following diagram commutes:  
\[
\begin{tikzcd}[ampersand replacement=\&]
X_{\rho} \times_k \Sp{k_s} \arrow{r}{j_X}\arrow{d}{(\mrm{id},\underline{\sigma}^{-1})} \& X \times_k \Sp{k_s} \arrow{d}{\rho(\sigma)\circ (\mrm{id},\underline{\sigma}^{-1})}\\
X_{\rho} \times_k \Sp{k_s} \arrow{r}{j_X} \&  X \times_k \Sp{k_s}
\end{tikzcd}
\]
\end{itemize}

In particular, the following diagram commutes, for any $\sigma \in \mrm{Gal}(k_s/k)$:
\[
\begin{tikzcd}[ampersand replacement=\&]
X_{\rho}(k_s) \arrow{r}{j_X}\arrow{d}{\sigma} \& X(k_s) \arrow{d}{\rho(\sigma) \circ \sigma}\\
X_{\rho}(k_s) \arrow{r}{j_X} \& X(k_s)
\end{tikzcd}
\]
Let $\mathcal{C}$ be a property of morphisms of $k$-schemes which is \'etale-local on the base and stable by pre- and post-composition by isomorphisms (for instance, any of the properties listed in Proposition \ref{fpqc-descent-prop}). If $X \rar \Sp{k}$ satisfies $\mathcal{C}$, then so does $X_{\rho} \rar \Sp{k}$.
If $X \rar \Sp{k}$ is projective (resp. quasi-projective), so is $X_{\rho} \rar \Sp{k}$.  
}

\rem{When $G$ acts in fact on $X$ by $k$-automorphisms, the cocycle condition simply means that $\rho: \mrm{Gal}(k_s/k) \rar G$ is a group homomorphism.}
 
\demo{Note that the last diagram follows formally from the assumptions (see also Lemma \ref{pi1-vs-galois}). The second to last assertion is formal, and the last assertion follows from \cite[Lemma 0BDE]{Stacks}. \\

\noindent
\emph{Step 1: Uniqueness up to isomorphism.}

Assume that $(U,j_U)$ and $(V,j_V)$ are two pairs satisfying the requested conditions. Let $j=j_V^{-1} \circ j_U$: then $j: U_{k_s} \rar V_{k_s}$ is an isomorphism of $k_s$-schemes such that for any $\sigma \in \mrm{Gal}(k_s/k)$, the following diagram commutes:
\[
\begin{tikzcd}[ampersand replacement=\&]
U \times_k \Sp{k_s} \arrow{r}{j}\arrow{d}{(\mrm{id},\underline{\sigma}^{-1})} \&  V \times_k \Sp{k_s}\arrow{d}{(\mrm{id},\underline{\sigma}^{-1})}\\
U \times_k \Sp{k_s} \arrow{r}{j} \& V \times_k \Sp{k_s}
\end{tikzcd}
\]

Since $U,V$ are quasi-compact, it follows from Corollary \ref{inf-galois-desc} that $j$ is in fact the base change of a map of $k$-schemes $j_1: U \rar V$. This $j_1$ becomes an isomorphism after the fpqc base change $V \times_k k_s \rar V$, hence is an isomorphism.\\ 

\noindent
\emph{Step 2: Uniqueness of the isomorphism.}

Suppose now that $u$ is an automorphism of a pair $(X_{\rho},j_X)$ satisfying the given conditions. Then, tautologically, $u \times \mrm{id}_{k_s}$ is the identity. By Proposition \ref{bi-morsheaf}, it follows that $u$ is the identity.\\

\noindent
\emph{Step 3: Existence and fonctoriality in the affine case.}

Assume that $X$ is the spectrum of some $k$-algebra $A$. Since $\rho$ maps the unit element to the identity and acts through a finite quotient of $\mrm{Gal}(k_s/k)$, there is a finite Galois extension $L/k$ such that $\rho(\mrm{Gal}(k_s/L))$ is trivial. After enlarging $L$, since $G$ is finite and $X$ is quasi-compact, we may assume that $G$ acts by automorphisms of $L$-schemes on $X \times_k L$. 

Thus, $\rho$ induces a right action $\rho^{\sharp}$ of $G$ on $A \otimes_k L$ by $L$-automorphisms. Let $\tau_2$ denote the natural left action of $\mrm{Gal}(L/k)$ on $A \otimes_k L$. The cocycle condition is equivalent to \[\rho^{\sharp}(\sigma\sigma')=\tau_2(\sigma)\rho^{\sharp}(\sigma')\tau_2(\sigma)^{-1}\rho^{\sharp}(\sigma).\] 

For $\sigma \in \mrm{Gal}(L/k)$, let $\tau(\sigma)=\rho^{\sharp}(\sigma)^{-1}\tau_2(\sigma) \in \mrm{Aut}_k(L)$. It is straightforward to check that $\tau$ defines a semi-linear action of $\mrm{Gal}[L/k)$ on $A \otimes_k L$. Let thus $B=(A \otimes_k L)^{\mrm{im}(\tau)}$. By Lemma \ref{galois-semi-lin}, the natural homomorphism of $L$-algebras $c: B \otimes_k L \rar A \otimes_k L$ is an isomorphism. 

We then define $X_{\rho}=\Sp{B}$ and $j_X$ by the map induced by $A \otimes_k k_s \overset{c^{-1} \otimes_L k_s}{\rar} B \otimes_k k_s$; checking that the diagram commutes is then formal. 

Note that this construction is functorial: if we have a morphism of affine $k$-schemes $f:X \rar Y$ such that $f \times_k \Sp{k_s}$ is equivariant for the action of the group $G$, the definition produces a map $f_{\rho}: X_{\rho} \rar Y_{\rho}$. Moreover, $X_{\rho}$ (resp. $f_{\rho}$) acquires all the properties of $X$ (resp. $f$) that are \'etale-local on the base (for instance, all those listed in Proposition \ref{fpqc-descent-prop}). \\

\noindent
\emph{Step 4: The global twist exists if there are enough good open subsets}

Let $U \subset X$ be an open subset; we say that $U$ is \emph{good} if it is affine and stable under $G$. Let $\mathbf{\mrm{GoodAff}}(X)$ be the partially ordered set of good affine opens of $X$. Let $\mathbf{\mrm{LocTw}}(X)$ be the category of pairs $(T,j)$, where $T$ is a $k$-scheme and $j: T \times_k \Sp{k_s} \rar X\times_k \Sp{k_s}$ is an open immersion of $k_s$-schemes, in which the morphisms are open immersions.  

Then there is a functor $\mrm{Tw}_{\rho}: U \mapsto (U_{\rho},j_{U})$ from $\mathbf{\mrm{GoodAff}}(X)$ to $\mathbf{\mrm{LocTw}}(X)$. Assume that $X$ can be covered with good open subsets. Because $\mrm{Tw}_{\rho}$ preserves Zariski-coverings, it defines a glueing data in the sense of \cite[Section 01JA]{Stacks}, and thus the $U_{\rho}$ glue to a scheme $X_{\rho}$ endowed with a $k_s$-isomorphism $j_X: X_{\rho} \times_k \Sp{k_s} \rar X \times_k \Sp{k_s}$. \\

\noindent
\emph{Step 5: Existence of good open subsets}

Finally, it remains to see that $X$ can be covered by good affine open subsets. Note that it is enough to prove that, for all $x \in X$, there is an affine open subset containing all the $G$-translates of $X$. Indeed, since $X$ is separated, we can then take the intersection of the $G$-translates of this affine open subset.  

So all we need to do is the following: show, given any finite $F \subset X$, that there is an affine open subset $U \subset X$ containing $F$. Indeed, consider the open subset $U' \subset \mathbb{P}^n_A$ such that $\iota: X \rar U'$ is a closed immersion, and let $I$ be a homogeneous ideal of $A[x_0,\ldots,x_n]$ (with $I \subset (x_0,\ldots,x_n)$) whose vanishing locus is exactly $\mathbb{P}^n_A \backslash U'$. It is enough to show that there is some homogeneous $f \in I$ of positive degree that does not vanish at any $x \in \iota(F)$. Every $x \in \iota(F)$ defines a homogeneous prime ideal $\mathfrak{p}_x \subset (x_0,\ldots,x_n)A[x_0,\ldots,x_n]$, and we thus seek to prove that $I_+ \not\subset \bigcup_{x \in \iota(F)}{\mathfrak{p}_x}$, where $I_+ \subset I$ is the subset of homogeneous elements of positive degree, given that we know that $I \not\subset \mathfrak{p}_x$ for all $x$. This is exactly the graded prime avoidance lemma \cite[Lemma 00JS]{Stacks}, so we are done. 
}

\prop[twist-equiv]{Let $X,Y$ be two $k$-schemes that are quasi-projective over some $k$-algebras $A,B$, and $G$ be a group acting on $X_{k_s}:= X \times_k \Sp{k_s}, Y_{k_s} := Y \times_k \Sp{k_s}$ by $k_s$-isomorphisms, and $f: X \rar Y$ be a map such that $f \times_k \mrm{id}_{k_s}$ is $G$-equivariant. Let $\rho: \mrm{Gal}(k_s/k) \rar G$ be a continuous map with the cocycle property of Proposition \ref{cocycle-twist} (when post-composed with the actions of $G$ on $X_{k_s}$ or $Y_{k_s}$). There is a unique map $f_{\rho}: X_{\rho} \rar Y_{\rho}$ of $k$-schemes such that the following diagram commutes:
\[
\begin{tikzcd}[ampersand replacement=\&]
X_{\rho} \times_k \Sp{k_s} \arrow{r}{j_X}\arrow{d}{f_{\rho} \times_k \mrm{id}_{k_s}} \& X \times_k \Sp{k_s} \arrow{d}{f \times_k \mrm{id}_{k_s}}\\
Y_{\rho} \times_k \Sp{k_s} \arrow{r}{j_Y} \&  Y \times_k \Sp{k_s}
\end{tikzcd}
\]

Moreover, $f_{\rho}$ acquires all the properties of $f$ that are \'etale-local on the base. The construction $X \longmapsto (X_{\rho},j_X), f \longmapsto f_{\rho}$ is then functorial in the obvious meaning, and it preserves arbitrary limits. 

}

\demo{Suppose $f: X \rar Y$ is a morphism of quasi-projective $k$-schemes, such that $f$ is equivariant for the action of a finite group $G$, and that we are given a continuous $\rho: \mrm{Gal}(k_s/k) \rar G$. 

Let $f_{\rho}^0: X_{\rho} \times_k \Sp{k_s} \rar Y_{\rho}\times_k \Sp{k_s}$ be given by $j_Y^{-1} \circ (f \times_k \mrm{id}_{k_s}) \circ j_X$. Then it is formal to see that $f_{\rho}^0 \in \mrm{Mor}_{k_s}((X_{\rho})_{k_s},(Y_{\rho})_{k_s})^{\mrm{Gal}(k_s/k)}$. Since $X_{\rho},Y_{\rho}$ are quasi-projective over $k$, they are quasi-compact and we can apply Proposition \ref{inf-galois-desc} to see that $f_{\rho}^0$ is the base change of a $k$-homomorphism $f_{\rho}: X_{\rho} \rar Y_{\rho}$. 

Note that the diagram uniquely prescribes the image of $f_{\rho}$ in $\mrm{Mor}_{k_s}(X_{k_s},Y_{k_s})$, which uniquely defines $f_{\rho}$ by Proposition \ref{bi-morsheaf}. 

The proof of fonctoriality is identical. 

Let $\mathscr{D}$ be a diagram of quasi-projective $k$-schemes (all the maps being equivariant for some action of $G$) and $X$ be its limit in the category $\mathbf{Sch}_k$. Suppose that $X$ is quasi-projective. Then, for any $k_s$-scheme $T$, the following diagram commutes, where the maps are isomorphisms: 

\[
\begin{tikzcd}[ampersand replacement=\&]
\mrm{Mor}_k(T,X) \arrow{r}\arrow{d}\& \mrm{Mor}_{k_s}(T,X_{k_s}) \arrow{r}{j_X}\arrow{d} \& \mrm{Mor}_{k_s}(T,(X_{\rho})_{k_s}) \arrow{r}\arrow{d} \& \mrm{Mor}_k(T,X_{\rho}) \arrow{d}\\
\lim\limits_{\substack{\longleftarrow\\D \in \mathscr{D}}}{\mrm{Mor}_k(T,D)} \arrow{r} \& \lim\limits_{\substack{\longleftarrow\\D \in \mathscr{D}}}{\mrm{Mor}_{k_s}(T,D_{k_s})} \arrow{r}{(j_D)_D}\& \lim\limits_{\substack{\longleftarrow\\D \in \mathscr{D}}}{\mrm{Mor}_k(T,(D_{\rho})_{k_s})} \arrow{r} \& \lim\limits_{\substack{\longleftarrow\\D \in \mathscr{D}}}{\mrm{Mor}_k(T,D_{\rho})}
\end{tikzcd}
\]

Now, $F_{\mathscr{D}_{\rho}}: T \in \mathbf{Sch}_k \longmapsto \lim\limits_{\substack{\longleftarrow\\D \in \mathscr{D}}}{\mrm{Mor}_k(T,D_{\rho})}$, $F_{X_{\rho}}: T \in \mathbf{Sch}_k \longmapsto \mrm{Mor}_k(T,X_{\rho})$ are fpqc sheaves by Theorem \ref{morsheaf}. By fonctoriality, we have a map $\lambda: F_{X_{\rho}} \rar F_{\mathscr{D}_{\rho}}$; the diagram shows that $\lambda(T): F_{X_{\rho}}(T) \rar F_{\mathscr{D}_{\rho}}(T)$ is an isomorphism for every $k_s$-scheme $T$. Since every $k$-scheme has an fpqc cover by $k_s$-schemes, it follows that $\lambda$ is an isomorphism. 
}

\cor[twist-ring-functions]{Let $k$ be a field with separable closure $k_s$, and $X$ be a quasi-projective $k$-scheme. Let $\rho: \mrm{Gal}(k_s/k) \rar \mrm{Aut}(X \times_k \Sp{k_s})$ satisfying the cocycle condition. Then the induced action of $\rho$ on $\mrm{Aut}(\Sp{\OO(X)} \times_k \Sp{k_s})$ satisfies the cocycle condition of Proposition \ref{cocycle-twist}. Moreover, the twist of $\Sp{\OO(X)}$ by $\rho$ naturally identifies with $\Sp{\OO(X_{\rho})}$. }

\demo{For each $\sigma \in \mrm{Gal}(k_s/k)$, the commutative diagram
\[
\begin{tikzcd}[ampersand replacement=\&]
X_{\rho} \times_k \Sp{k_s} \arrow{r}{j_X}\arrow{d}{(\mrm{id},\underline{\sigma}^{-1})} \& X \times_k \Sp{k_s} \arrow{d}{\rho(\sigma)\circ (\mrm{id},\underline{\sigma}^{-1})}\\
X_{\rho} \times_k \Sp{k_s} \arrow{r}{j_X} \&  X \times_k \Sp{k_s}
\end{tikzcd}
\]
induces, by taking the spectra of the rings of global functions (since $k_s$ is flat over $k$), a commutative diagram
\[
\begin{tikzcd}[ampersand replacement=\&]
X_{\rho} \times_k \Sp{k_s} \arrow{rrr}{j_X}\arrow{ddd}{(\mrm{id},\underline{\sigma}^{-1})} \arrow{dr} \& \& \& X \times_k \Sp{k_s} \arrow{dl}\arrow[swap]{ddd}{\rho(\sigma)\circ (\mrm{id},\underline{\sigma}^{-1})}\\
\& \Sp{\OO(X_{\rho})} \times_k \Sp{k_s} \arrow{r}{j_X}\arrow[swap]{d}{(\mrm{id},\underline{\sigma}^{-1})} \& \Sp{\OO(X)} \times_k \Sp{k_s} \arrow[swap]{d}{\rho(\sigma)\circ (\mrm{id},\underline{\sigma}^{-1})}\&\\
\& \Sp{\OO(X_{\rho})} \times_k \Sp{k_s} \arrow{r}{j_X} \&  \Sp{\OO(X)} \times_k \Sp{k_s} \&\\
X_{\rho} \times_k \Sp{k_s} \arrow{rrr}{j_X}\arrow{ru} \&\&\&  X \times_k \Sp{k_s} \arrow{lu}
\end{tikzcd}
\]
where the unlabeled arrow are canonical morphisms $Y \rar \Sp{\OO(Y)}$, and this is enough to conclude.
}

\section{Compactification of finite schemes over the affine line}
\label{situation-at-infty}

In this section, our base rings will be assumed to be regular, excellent, and Noetherian (unless explicitly specified otherwise -- the rings of sections of \emph{schemes}, on the other hand, are not required to satisfy these conditions). The goal is to state some properties needed to construct compactified moduli schemes and work with morphisms between them.

\defi[normal-near-infinity]{We say that a $R$-scheme of finite type $X$ is \emph{normal near infinity} if there is a closed subscheme $Z \subset X$ finite over $R$ such that $X \backslash Z$ is a normal scheme and $X \backslash Z \rar \Sp{R}$ is surjective. }

\lem[equiv-normal-near-infinity]{Suppose we are given a finite $R$-map $j: X \rar \mathbb{A}^1_R$. Then $X$ is normal near infinity if, and only if, there is a monic $f \in R[t]$ such that $X$ is normal on the locus where $f(j)$ is invertible.}

\demo{The assumptions imply that $X$ is a $R$-scheme of finite type. For $(\Leftarrow)$, it is clear that $Z=j^{-1}(V(f))$ fits the definition. 
Let us assume that $X$ is normal near infinity and let $Z$ be the closed subscheme given by the definition. Then the scheme-theoretic image $H$ of $Z$ in $\mathbb{A}^1_R$ is affine. Moreover $\OO(H)$ injects into $\OO(Z)$, so, since $R$ is Noetherian, $H$ is finite over $R$. Thus, if $x$ denotes the coordinate function on $\mathbb{A}^1_R$, there is a monic $f \in R[t]$ such that $f(x)_{|H}=0$. Then $X$ is normal on $j^{-1}(D(f))$. }

\cor[normal-near-infinity-local]{Suppose we are given a finite $R$-map $j:X \rar \mathbb{A}^1_R$ and $f_1,\ldots, f_r \in R$ that generate the unit ideal. Then $X$ is normal near infinity if and only if every $X_{R_{f_i}}$ is normal near infinity (as a $R_{f_i}$-scheme of finite type).}

\demo{The forward implication is clear, so let us discuss the reverse implication. For every $i$, we have an ideal $I_i$ of $R[x]_{f_i}$ containing a monic polynomial and such that every point of $X_{R_{f_i}} \backslash j^{-1}(V(I_i))$ is normal. Let $I \subset R[x]$ be the sum of the inverse images under $R \rar R_{f_i}$ of the $I_i$, so we only need to prove that $I$ contains a monic polynomial. But it is clear that the subset $J \subset R$ which is the reunion of $\{0\}$ and the leading coefficients of polynomials in $I$ is an ideal. Since each $I_i$ contains a monic polynomial, each $J_{f_i}$ is $R_{f_i}$, hence $J=R$ and $I$ contains a monic polynomial. }

\prop[one-compactification]{Let $j: X \rar \mathbb{A}^1_R$ be a finite map, where $X$ is normal near infinity. There is a canonical Cartesian diagram 
\[
\begin{tikzcd}[ampersand replacement=\&]
X \arrow{r}{j}\arrow{d} \& \mathbb{A}^1_R \arrow{d}\\
\overline{X}\arrow{r}{\overline{j}} \&  \mathbb{P}^1_R
\end{tikzcd}
\]
such that $\overline{j}$ is finite and $\overline{X}$ is normal in a neighborhood of $\overline{X}\backslash X$. 
Moreover, $X\rar \overline{X}$ is an open immersion with dense image. For any separated scheme $Y$ and any open subscheme $U \subset \overline{X}$, $Y(U) \rar Y(X \cap U)$ is injective. By taking $Y=\mathbb{A}^1$, we see that $\OO(U) \rar \OO(U \cap X)$ is injective.  

Moreover, let $R'$ be any (regular, excellent, Noetherian) $R$-algebra and $j': X' \rar \mathbb{A}^1_{R'}$ be a finite map, with $X'$ normal near infinity. Assume that there is a closed subscheme\footnote{Note, in the following, the $R$ instead of $R'$.} $Z'$ of $\mathbb{A}^1_R$, finite over $R$, such that $X'$ is normal at every point not above $Z'$, and that we have a commutative diagram 
\[
\begin{tikzcd}[ampersand replacement=\&]
X' \arrow{r}{f}\arrow{d}{j'} \& X\arrow{d}{j}\\
\mathbb{A}^1_{R'} \arrow{r} \&  \mathbb{A}^1_R
\end{tikzcd}
\]
Then $f$ extends uniquely to a map $\overline{f}: \overline{X'} \rar \overline{X}$ such that the following diagram commutes:
 \[
\begin{tikzcd}[ampersand replacement=\&]
\overline{X'} \arrow{r}{\overline{f}}\arrow{d}{\overline{j'}} \& \overline{X}\arrow{d}{\overline{j}}\\
\mathbb{P}^1_{R'} \arrow{r} \&  \mathbb{P}^1_R
\end{tikzcd}
\]
}

\demo{First, we construct $\overline{X}$. Fix a closed subscheme $Z$ of $\mathbb{A}^1_R$ finite over $R$ such that $X \backslash j^{-1}(Z)$ is normal. Let $Y$ be the relative normalization of $X \backslash j^{-1}(Z) \overset{j}{\rar} \mathbb{A}^1_R \subset \mathbb{P}^1_R$. Because $R$ is excellent, $\mathbb{P}^1_R$ is Nagata by \cite[Lemma 07QV, Lemma 035A]{Stacks} so $Y$ is finite over $\mathbb{P}^1_R$ by \cite[Lemma 03GR]{Stacks} ($X \backslash Z$ is normal, hence reduced). By \cite[Lemma 035L]{Stacks}, $Y$ is normal. 

By \cite[Lemma 035K]{Stacks}, $Y \times_{\mathbb{P}^1_R} (\mathbb{A}^1_R \backslash Z)$ is exactly the normalization of $X \backslash j^{-1}(Z) \overset{j}{\rar} \mathbb{A}^1_R \backslash Z$. Since this map is finite, by \cite[Lemma 03GP]{Stacks}, $Y \times_{\mathbb{P}^1_R} (\mathbb{A}^1_R \backslash Z) \rar \mathbb{A}^1_R \backslash Z$ is exactly isomorphic to $X \backslash j^{-1}(Z) \rar \mathbb{A}^1_R \backslash Z$. Thus, we can glue $Y \backslash j^{-1}(Z)\rar \mathbb{P}^1_R \backslash Z$ and $X \rar \mathbb{A}^1_R$ together, and this defines $\overline{X}$. 

Note that the construction of $\overline{X}$ does not depend on the choice of $Z$. 

$X \rar \overline{X}$ is an open immersion by construction, and it has dense image by \cite[Lemma 0AXP]{Stacks}. Moreover, $\overline{X}$ is normal, hence reduced, on a neighborhood of any point of $\overline{X} \backslash X$. Let $U$ be any open subscheme of $\overline{X}$. Since $X$ is dense, $X \cap U$ is a nonempty open subscheme of $U$; furthermore, if $V \subset U$ is any open subscheme, $V \cap X$ is nonempty and contained in $(X \cap U) \cap V$. Thus $X \cap U \rar U$ is an open immersion with dense image and such that all points of $U \backslash X$ are normal hence reduced. 

Let $Y$ be a separated scheme and $u,v: U \rar Y$ be two morphisms that agree on $U \cap X$. By \cite[Lemma 01KM]{Stacks}, there is a closed subscheme $U_1 \subset U$ containing $U \cap X$ such that $u,v$ agree on $U_1$. Thus the surjective closed immersion $U_1 \rar X$ is an isomorphism above $U \cap X$. Since $U$ is reduced at every point outside $X$, $U_1=U$ thus $u=v$. 

To construct $\overline{f}$, let $Y'$ denote the relative normalization of $X' \backslash (j')^{-1}((Z' \cup Z)_{R'}) \rar \mathbb{P}^1_{R'}$; since we have a commutative diagram
 \[
\begin{tikzcd}[ampersand replacement=\&]
X' \backslash (j')^{-1}((Z' \cup Z)_{R'}) \arrow{r}{j'}\arrow{d}{f}\& \mathbb{P}^1_{R'}\arrow{d}\\
X \backslash j^{-1}(Z' \cup Z) \arrow{r}{j} \&  \mathbb{P}^1_R
\end{tikzcd}
\]
we have a natural map $\overline{f}: Y' \backslash Y$ above $j^{-1}(Z' \cup Z)$. We can then glue it with $f: X' \rar X$. 
}

\cor[extension-to-finite]{Let $j: X \rar \mathbb{A}^1_R$ be a finite map with $X$ normal near infinity. Let $Y$ be a finite $R$-scheme, and let $g: X \rar Y$ a $R$-morphism. Then $g$ extends uniquely as to map $\overline{g}: \overline{X} \rar Y$. }

\demo{The uniqueness follows from Proposition \ref{one-compactification}. Fix a closed subscheme $Z \subset \mathbb{A}^1_R$ finite over $R$ such that $X \backslash j^{-1}(Z)$ is normal. Suppose that, for any normal integral affine open subscheme $U \subset \overline{X} \backslash j^{-1}(Z)$, $X\cap U \rar Y$ extends to $U \rar Y$. Then, if $U,V$ are normal integral affine open subschemes of $\overline{X} \backslash j^{-1}(Z)$, these extensions coincide on $U \cap V$, because they coincide on $U \cap V \cap X$. All these extensions glue together, by construction, with $X \rar Y$, so they all glue together and construct the extension $\overline{X} \rar Y$. 

So let $U \subset \overline{X} \backslash j^{-1}(Z)$ be a normal integral affine open subscheme. By definition, $g$ defines a map $g^{\sharp}: \OO(Y) \rar \OO(X) \rar \OO(X \cap U)$. Because $U$ is normal integral affine and $X \cap U$ is an open subscheme of $U$, $\OO(U) \rar \OO(X \cap U)$ is injective between normal domains and is an isomorphism on function fields. The image of $g^{\sharp}$ is thus made with elements in the fraction ring of $\OO(U)$ and which are integral over $R$, hence over $\OO(U)$. So the image of $g^{\sharp}$ lands in the subring $\OO(U)$, that is, $g$ extends to a map $U \rar Y$ and we are done.}

We now wish to show that the compactification previously constructed is canonical, and that a much wider class of morphisms can be extended.

\lem[kinda-zmt]{Let $f: Z \rar X$ be a finite morphism of finite type $R$-schemes. Suppose that $Y \subset X$ is a dense open subscheme such that $f^{-1}(Y) \rar Y$ is an isomorphism, $f^{-1}(Y)$ is dense in $Z$, every $\OO_{Z,z}$ for $z \notin f^{-1}(Y)$ is reduced, and every $\OO_{X,x}$ for $x \notin Y$ is normal. Then $f$ is an isomorphism. }

\demo{First, we easily reduce to the case where $X$ and $Z$ are affine. 

Let $x \in X \backslash Y$, and let $z' \in Z$ be a point whose image in $X$ specializes to $x$. Then $f(\overline{\{z'\}})$ is closed in $X$ and contains a point specializing to $x$; hence there exists a specialization $z \in Z$ of $z'$ whose image in $Z$ is $x$. In particular, $\OO_{Z,z'}$ is a localization of $\OO_{Z,z}$, hence is reduced. 

Thus, it is enough to check the claim after changing base by every $\OO_{X,x}$ where $x \notin Y$. So we can assume that $Z=\Sp{B}$, $X=\Sp{A}$ with $A$ a local normal domain, $B$ is reduced, and $Y$ is the open subscheme given by the complement of the vanishing locus of some nonzero ideal $I$. Furthermore, we know that for every non-zero $f \in I$, $A[f^{-1}] \rar B[f^{-1}]$ is an isomorphism. In particular, $A \rar B[f^{-1}]$ is injective, so that $A \rar B$ is injective. Since $f^{-1}(Y)$ is dense in $Z$, we know that for every non-zero $g \in B$, $gI$ is nonzero in $B$. 

Let now $x \in B$ be nonzero: then the ideal $xI$ is nonzero, so there is some $f \in I$ such that $fI \neq 0$. Then, if $r \geq 1$, $f^rx$ is a divisor of $(fx)^r$ which is nonzero (since $B$ is reduced): thus the image of $x$ in $B_f$ is nonzero. Let now $x,y \in B$ be nonzero, then $\{t \in I\mid tx=0\}, \{t \in I\mid ty=0\}$ are two proper additive subgroups of $I$, hence their reunion does not contain $I$, and we can find $f \in I$ such that $x,y$ have nonzero images in $B_f$. Since $B_f \simeq A_f$ is a domain, $xy$ has nonzero image in $B_f$, so $xy \neq 0$. Thus $B$ is a domain.  

Therefore, $A \rar B$ is finite injective and is an isomorphism on fraction fields. Since $A$ is normal, $A \rar B$ is an isomorphism.}

\prop[unique-compactification]{Let $\alpha, \beta: X \rar \mathbb{A}^1_R$ be two finite maps, where $X$ is a $R$-scheme of finite type normal near infinity. By Proposition \ref{one-compactification}, $\alpha, \beta$ compactify in two finite maps $\overline{\alpha}: X_{\alpha} \rar \mathbb{P}^1_R$, $\overline{\beta}: X_{\beta}\rar \mathbb{P}^1_R$. There is a unique $R$-isomorphism $\iota: X_{\alpha} \rar X_{\beta}$ such that the following diagram commutes:
 \[
\begin{tikzcd}[ampersand replacement=\&]
X  \arrow{r}{\mrm{id}}\arrow{d}{\overline{\alpha}}\& X\arrow{d}{\overline{\beta}}\\
X_{\alpha} \arrow{r}{\iota} \&  X_{\beta}
\end{tikzcd}
\]

}

\demo{First, note that it is enough to construct a $R$-map $\iota_{\alpha,\beta}$ making the diagram commute: because then, we can similarly construct $\iota_{\beta,\alpha}$, so $\iota_{\beta,\alpha} \circ \iota_{\alpha,\beta}$ is an endomorphism of $X_{\alpha}$ which induces the identity on $X$; by Proposition \ref{one-compactification}, it is the identity. Similarly, $\iota_{\alpha,\beta} \circ \iota_{\beta,\alpha}$ is the identity endomorphism of $X_{\beta}$ so that $\iota_{\alpha,\beta}$ is an isomorphism. 

Next, we consider the diagonal inclusion $X \rar X \times_R X \rar X_{\alpha} \times_R X_{\beta}$. Let $Z$ be its scheme-theoretic image. $Z$ is a closed subscheme of $X_{\alpha} \times_R X_{\beta}$, hence the first projection $p_1: Z \rar X_{\alpha}$ is proper. 

If $p_1$ is an isomorphism, we are done, since we can consider the composition $p_2\circ p_1^{-1}$, where $p_2$ is the second projection $Z \rar X_{\alpha} \times_R X_{\beta} \rar X_{\beta}$. 

We will show that $p_1$ is finite and that the dense open subscheme $X$ of $X_{\alpha}$ satisfies all the conditions of Lemma \ref{kinda-zmt}, whence the conclusion. 

First, we know that $X_{\alpha}$ is normal at every point not in $X$. Moreover, by \cite[Lemma 01R8]{Stacks}, $Z \cap X \times_R X_{\beta}$ is the scheme-theoretic image of the diagonal inclusion $X \rar X \times_R X_{\beta}$, which is a closed immersion (as a section of the proper morphism $X \times_R X_{\beta} \rar X$). Thus $X \rar Z \cap (X \times_R X_{\beta})$ is an isomorphism of inverse $p_1^{-1}(X) \rar X$. Note that the image of $X$ in $Z$, which is $p_1^{-1}(X)$, is dense in $Z$ by \cite[Lemma 01RB]{Stacks}. 

Let $W \subset X$ be a closed subscheme, finite over $R$, such that every point of $X \backslash W$ is normal (hence reduced). Then $p_1^{-1}(X_{\alpha} \backslash W)$ is the scheme-theoretic image of $X \backslash W \rar (X_{\alpha} \backslash W) \times_R X_{\beta}$. Since $X \backslash W$ is reduced, $Z \backslash p_1^{-1}(W)$ is reduced by \cite[Lemma 056B]{Stacks}. In particular, $Z$ is reduced at every point not in $p_1^{-1}(X)$.   

All that remains to do now is prove that $p_1$ is finite. To do that, it is enough to show that it has finite fibres by \cite[Lemma 02LS]{Stacks}, so we only need to prove that for every $x \in X_{\alpha} \backslash X$, $p_1^{-1}(\{x\})$ is finite. To do this, we exhibit a closed subscheme $T \subset X_{\alpha} \times X_{\beta}$ such that the diagonal inclusion factors through $T$, and such that the fibre of $T \rar X_{\alpha} \times_R X_{\beta} \rar X_{\alpha}$ at every $x \in X_{\alpha}$ is finite (since then $p_1: Z \rar X_{\alpha}$ is the composition of the closed immersion $Z \rar T$ and of the map $T \rar X_{\alpha}$). 

The morphism $\alpha$ is finite, so $X$ is affine. Let $s$ denote the coordinate on $\mathbb{A}^1_R$. Let $a,b \in \OO(X)$ denote the images of $s$ under $\alpha,\beta$. Since $\beta$ is finite, there is an integer $n \geq 1$ and polynomials $Q_0, \ldots, Q_{n-1} \in R[t]$ such that $a^n+\sum_{i=0}^{n-1}{Q_i(b)a^i}=0$. Let $d \geq 0$ be the maximum of the degrees of the $Q_i$, and define the homogeneous polynomials $\tilde{Q_i}(v_0,v_1)=v_1^dQ_i(v_0/v_1) \in R[v_0,v_1]$ of degree $d$. Let $P \in R[u_0,u_1,v_0,v_1]$ be the bi-homogeneous polynomial of bi-degee $(n,d)$ given by $P=u_0^nv_1^d+\sum_{i=0}^{n-1}{\tilde{Q_i}(v_0,v_1)u_0^iu_1^{n-i}}$. Then the image of $X$ in $X_{\alpha} \times X_{\beta} \overset{(\overline{\alpha},\overline{\beta})}{\rar} \mathbb{P}^1_R \times_R \mathbb{P}^1_R$ is contained in $(\overline{\alpha},\overline{\beta})^{-1}(V(P))$. 

To conclude, since $\beta$ is finite, it is enough to check that $V(P) \rar \mathbb{P}^1_R \times_R \mathbb{P}^1_R \overset{p_1}{\rar}\mathbb{P}^1_R$ has finite fibres above points $x \in \mathbb{P}^1_R \backslash \mathbb{A}^1_R$. This is equivalent to asking that $P(1,0,v_0,v_1) \in k[v_0,v_1]$ has finitely many zeroes in $\mathbb{P}^1_k$, for a $R$-algebra $k$ which is a field, which is clear. 
}

\prop[functorial-compactification]{Let $X,Y$ be two $R$-schemes normal near infinity and admitting finite maps $\alpha: X \rar \mathbb{A}^1_R, \beta: Y \rar \mathbb{A}^1_R$. Let $\overline{\alpha}: X_{\alpha} \rar \mathbb{P}^1_R$, $\overline{\beta}: Y_{\beta} \rar \mathbb{P}^1_R$ be the maps given by Proposition \ref{one-compactification}, where $X$ (resp. $Y$) is an open subscheme of $X_{\alpha}$ (resp. $Y_{\beta}$). Let $f: X \rar Y$ be any finite map of $R$-schemes. Then there is a unique map $\overline{f}: X_{\alpha} \rar Y_{\beta}$ making the following diagram commute:
 \[
\begin{tikzcd}[ampersand replacement=\&]
X  \arrow{r}{f}\arrow{d}\& Y\arrow{d}\\
X_{\alpha} \arrow{r}{\overline{f}} \&  Y_{\beta}
\end{tikzcd}
\]
Moreover, $\overline{f}$ is finite, and maps any $x \in X_{\alpha} \backslash X$ into $Y_{\beta} \backslash Y$. 
}  

\demo{Consider the finite map $\alpha': X \rar Y \overset{\beta}{\rar} \mathbb{A}^1_R$ of $R$-schemes and let $X_{\alpha'}$ be the associated compactification by Proposition \ref{one-compactification}. By Proposition \ref{unique-compactification}, we have an isomorphism $\iota: X_{\alpha} \rar X_{\alpha'}$ preserving $X$, and by Proposition \ref{one-compactification} we have a natural extension $\overline{f_{\alpha'}}: X_{\alpha'} \rar Y_{\beta}$. It is a map of finite $\mathbb{P}^1_R$-schemes, so is finite. Thus $\overline{f}=\overline{f_{\alpha'}}\circ \iota$ satisfies all the requested conditions.  
The uniqueness follows from Proposition \ref{one-compactification} since $Y_{\beta}$ is separated. }

\prop[compactification-functor]{Let $\mathbf{ACrv}_R$ be the category whose objects are $R$-schemes of finite type, normal near infinity and admitting a finite map to $\mathbb{A}^1_R$, and whose morphisms are finite maps of $R$-schemes. Let $\mathbf{PCrv}_R$ be the category of proper $R$-schemes that are normal near infinity, its morphisms being finite $R$-maps. There is a \emph{normal compactification functor} $\mathbf{ACrv}_R \rar \mathbf{PCrv}_R$, denoted by $X \longmapsto \overline{X}$, $f \longmapsto \overline{f}$, satisfying the following properties:
\begin{itemize}[noitemsep,label=$-$]
\item Every object $X$ of $\mathbf{ACrv}_R$ has a natural open immersion of $R$-schemes $X \longmapsto \overline{X}$ with dense image,
\item For every object $X$ of $\mathbf{ACrv}_R$, every separated scheme $Y$ and every open subset $U \subset \overline{X}$, $Y(\overline{X}) \rar Y(X)$ is injective,
\item For any morphism $f: X \rar Y$ in $\mathbf{ACrv}_R$, $\overline{f}$ sends any point in $\overline{X} \backslash X$ into  $\overline{Y} \backslash Y$,
\item $\overline{X}$ is normal at every point $x \notin X$, 
\item $\overline{\mathbb{A}^1_R}=\mathbb{P}^1_R$,
\end{itemize} }

\lem[compact-normal]{Let $f: X \rar Y$ be a morphism in $\mathbf{ACrv}_R$. Then $\overline{f}: \overline{X} \rar \overline{Y}$ is the relative normalization of $X \rar Y \subset \overline{Y}$.}

\demo{Let $f': \overline{X}' \rar \overline{Y}$ be the relative normalization of $X \rar Y \subset \overline{Y}$. Because $X,Y$ are objects of $\mathbf{ACrv}_R$, $f: X \rar Y$ is finite and $R$ is Noetherian, there is a closed subscheme $Z \subset Y$, finite over $R$, such that $X \backslash f^{-1}(Z), Y \backslash Z$ are normal schemes. Note by \cite[Lemma 035K, Lemma 035L]{Stacks}, $\overline{X}' \backslash (f')^{-1}(Z)$ is normal and $X \rar \overline{X'} \times_{\overline{Y}} Y$ is an isomorphism. In particular, by \cite[Lemma 03GR]{Stacks}, $f'$ is finite above $\overline{Y} \backslash Z$ and is finite above $Y$, hence is finite. 

Since $\overline{X} \rar \overline{Y}$ is finite, by the universal property of relative normalization, there is a unique $\overline{Y}$-morphism $g: \overline{X}'\rar \overline{X}$, and this morphism is the identity above $X$. This morphism is finite (since $\overline{X}', \overline{X}$ are finite $\overline{Y}$-schemes), the points in $\overline{X} \backslash X$ (resp. their inverse images under $g$) are normal points of $\overline{X}$ (resp. $\overline{X}'$). Finally, $X$ is a dense open subscheme of $\overline{X}$ by construction and of $\overline{X}'$ by \cite[Lemma 0AXP]{Stacks}. We conclude using Lemma \ref{kinda-zmt}. 
}

~\\Let us conclude this section with a convenient notion of smoothness at infinity.

\defi[smooth-at-infty-gnl]{Let $R$ be a regular excellent Noetherian ring. An object $Y$ of $\mathbf{ACrv}_R$ is \emph{smooth at infinity} if there exists a closed subscheme $Z \subset Y$, finite over $\Sp{R}$, such that $\overline{Y} \backslash Z \rar \Sp{R}$ is smooth of relative dimension one, and the scheme $(\overline{Y}\backslash Y)^{red}$ is finite \'etale over $\Sp{R}$.}

\prop[smooth-at-infty-base-change]{Let $R \rar S$ be a morphism of regular excellent Noetherian rings, and $Y$ be an object of $\mathbf{ACrv}_R$ which is smooth at infinity. Then $Y_S$ is an object of $\mathbf{ACrv}_S$ which is smooth at infinity, and there is a canonical morphism $\overline{Y_S} \rar \overline{Y} \times_R \Sp{S}$. This morphism is an isomorphism. Moreover, this isomorphism realizes an isomorphism $(\overline{Y_S} \backslash Y_S)^{red} \rar (\overline{Y}\backslash Y)^{red} \times_R \Sp{S}$. 
Furthermore, $(\overline{Y} \backslash Y)^{red}$ is a relative effective Cartier divisor on $\overline{Y} \rar \Sp{R}$, in the sense of \cite[Definition 062T]{Stacks}.}

\demo{By definition, there is a closed subscheme $Z \subset Y$ such that $Y \rar \Sp{R}$ is smooth outside $Z$, and a finite map $f: Y \rar \mathbb{A}^1_R$. Hence $Y_S \rar S$ is smooth outside $Z_S$, so that $Y_S$ is normal outside $Z_S$, which is a closed subscheme of $Y_S$ finite over $S$. Moreover $f_S: Y_S \rar \mathbb{A}^1_S$ is finite. Thus $Y_S$ is an object in $\mathbf{ACrv}_S$. By Proposition \ref{one-compactification}, there is a map $c: \overline{Y_S} \rar \overline{Y} \times_R \Sp{S}$, which is an isomorphism above the dense open subset $Y \times_R \Sp{S}$. Moreover, the inverse image of $Y \times_R \Sp{S}$ is dense, and its complement is made of normal points of $\overline{Y_S}$ (by construction of the compactification). Since $\overline{Y} \times_R S \rar \Sp{S}$ is smooth at every point not in $Y_S$, every such point is normal. Thus $c$ is an isomorphism by Lemma \ref{kinda-zmt}. 

In particular, $c$ induces a surjective closed immersion of schemes $(\overline{Y_S} \backslash Y_S)^{red} \rar (\overline{Y} \backslash Y)^{red} \times_R S$. In order to show that it is an isomorphism, it is enough to show that the domain and the target are reduced. For the domain, it is a tautology; for the target, it follows from the fact that $(\overline{Y}\backslash Y)^{red}$ is \'etale over $R$.

For the assertion about relative effective Cartier divisors, since $(\overline{Y}\backslash Y)^{red} \rar \Sp{R}$ is finite \'etale, it is enough by \cite[Lemma 062Y]{Stacks} to prove the statement after everything is base changed to a field. Using what we showed above, it is equivalent to proving the assertion assuming that $R$ is a field. 

For every $y \in \overline{Y} \backslash Y$, since $\overline{Y} \backslash Z$ is normal, we can find an open subset $y \in U_y$ of $\overline{Y}$ contained in $(Y \cup \{y\}) \backslash Z$, and some $t_y \in \OO_{\overline{Y}}(U_y)$ which is a uniformizer at $y$ and a unit at every other point of $U$. Consider then the subsheaf $\mathcal{L}$ of the sheaf of meromorphic functions $\mathcal{K}_{\overline{Y}}$ on $\overline{Y}$ defined as follows: $\mathcal{L}_{|Y}=\OO_Y$, $\mathcal{L}_{|U_y}=t_y^{-1}\OO_{U_y}$. Then $\mathcal{L}$ is an effective Cartier divisor, and the associated closed subscheme is exactly $(\overline{Y} \backslash Y)^{red}$. }

\prop[abhyankar-compactification]{Let $R$ be a regular excellent Noetherian ring and $Y$ be an object in $\mathbf{ACrv}_R$, smooth at infinity. Let $Z \subset Y$ be a closed subscheme, finite over $\Sp{R}$, such that $Y \backslash \Sp{R}$ is smooth outside $Z$. Assume that there is a finite morphism $f: X \rar Y$ of $R$-schemes which is \'etale above $Y \backslash Z$. Then $X$ is an object of $\mathbf{ACrv}_R$, $f$ is a morphism of $\mathbf{ACrv}_R$. 
If furthermore $f$ is tamely ramified in codimension one over $\overline{Y}$ \footnote{In the sense of \cite[Section 0BSE]{Stacks}.}, then $\overline{X} \backslash j^{-1}(Z) \rar \Sp{R}$ is smooth of relative dimension one, and $X$ is smooth at infinity. Moreover, $\overline{X} \backslash f^{-1}(Z)\rar \overline{Y} \backslash Z$ is finite locally free. }

\demo{For the first part, since $f^{-1}(Z)$ is finite over $\Sp{R}$, and $X \backslash f^{-1}(Z) \rar Y \backslash Z$ is finite \'etale, $X \backslash f^{-1}(Z)$ is normal, so that $X$ (resp. $f$) is an object (resp. a morphism) of $\mathbf{ACrv}_R$.

By Lemma \ref{compact-normal}, $\overline{f}: \overline{X} \rar \overline{Y}$ is the relative compactification of $X \overset{f}{\rar} Y \rar \overline{Y}$. We now apply Abhyankar's Lemma \cite[Lemma 0EYG, 0EYH]{Stacks}: \'etale locally on $\overline{Y}\backslash Z$, the map $\overline{f}_Z: \overline{X} \backslash f^{-1}(Z) \rar \overline{Y}\backslash Z$ is given by finitely many copies of $A \rar A[x]/(x^e-f)$, where $e \in A^{\times}$, and $f$ cuts exactly the inverse image in $A$ of $(\overline{Y}\backslash Y)^{red}$: in particular, $A$ is smooth of relative dimension one over $\Sp{R}$ and $A/fA$ is \'etale over $\Sp{R}$\footnote{This is not stated as such in the reference, but follows directly from the proof.}. 

It implies that $\overline{f}_Z$ is finite locally free, since these properties are \'etale local on the base. To establish that $\overline{X} \backslash f^{-1}(Z)$ is smooth of relative dimension one over $\Sp{R}$, by \cite[Lemma 036U]{Stacks}, it is enough to show that under these assumptions, $R \rar A' := A[x]/(x^e-f)$ is smooth of relative dimension one. Since $A'$ is free over $A$, $A'$ is flat over $R$, so it is enough to show that for any algebraically closed field $k$ endowed with a $R$-algebra structure, $A'_k := A' \otimes_R k$ is normal. 

Let $A_k=A \otimes_R k$, it is smooth over $k$ of relative dimension one, hence is a Dedekind ring, $A_k/fA_k$ is \'etale over $k$, and $A'_k=A_k[x]/(x^e-f)$. Since $A'_k[1/f]$ is \'etale over $A_k[1/f]$, it is enough to prove that, if $\mathfrak{m}$ is maximal ideal of $A_k$ containing $f$, then $(A'_k)_{\mathfrak{m}}=(A_k)_{\mathfrak{m}}[x]/(x^e-f)$ is a normal ring. In this situation, $(A_k)_{\mathfrak{m}}$ is a discrete valuation ring, and, since $A_k/fA_k$ is a direct product of copies of $k$, $f$ is a uniformizer. It is then a standard exercice to check that $A'_k)_{\mathfrak{m}}$ is normal.  

This local description implies that \'etale-locally on $\overline{Y} \backslash Z$, $(\overline{X} \backslash X)^{red}$ is given by copies of $A[x]/\sqrt{(x^e-f,fA)}=A[x]/(x) \simeq A$. Thus, $(\overline{X} \backslash X)^{red}$ is smooth over $\Sp{R}$ by \cite[Lemma 036U]{Stacks}, and it is also finite over $\Sp{R}$, hence it is finite \'etale.
}

\section{Refresher on group schemes and Abelian varieties}

\subsection{Generalities}

\defi[abelian-variety-defi]{A group scheme over a field $k$ is the datum of a $k$-scheme $A$ endowed with three morphisms $m: A \times_k A \rar A$, $i: A \rar A$, $e \in A(k)$ of $k$-schemes such that, for any $k$-scheme $S$, $A(S)$ is a group for the composition law $A(S) \times A(S) \overset{m}{\rar} A(S)$, where the unit element is $S \rar \Sp{k} \overset{e}{\rar} A$, and the inverse map is $A(S) \overset{i}{\rar} A(S)$. The group scheme $(A,m,i,e)$ is an Abelian variety if $A$ is a smooth proper connected $k$-scheme. In general, the datum of $(m,i,e)$ will be implicit in the notation, and we will denote $m: A(S) \times A(S) \rar A(S)$ by $+$ (or $\cdot$), $i: A(S) \rar A(S)$ by $t \longmapsto -t$ (or $t \longmapsto t^{-1}$), and $e$ by $0_A$ (or $1_A$, or $e_A$). }

\prop[morphism-is-group-morphism]{Let $A,B$ be Abelian varieties over a field $k$, and $f: A \rar B$ a homomorphism of $k$-schemes such that $f(0_A)=0_B$. Then $f$ is a morphism of Abelian varieties, that is, for any $k$-scheme $S$, the function $f: A(S) \rar B(S)$ is a group homomorphism.}

\demo{By \cite[Corollary 2.2]{MilAb}, $f$ is the composition of a homomorphism of Abelian varieties with a translation. Since $f(0_A)=0_B$, the translation is trivial, so $f$ is a homomorphism of Abelian varieties. }

\prop[abelian-variety-is-abelian]{In Definition \ref{abelian-variety-defi}, the group law on $A(S)$ is commutative.}

\demo{This is \cite[Corollary 2.4]{MilAb}. Proposition \ref{morphism-is-group-morphism} implies that $i: A \rar A$ is a homomorphism of Abelian varieties, so that, for any $k$-scheme $S$, $A(S)$ is a group whose inversion map is a group homomorphism: it is thus an Abelian group, whence the conclusion. }

\defi[kernel-abelian]{Given a morphism $f: A \rar B$ of Abelian varieties over a field $k$, the kernel of $f$ is the scheme $\ker{f}=f^{-1}(0_B) = \{0_B\} \times_B A$. }

\prop[flatness-abelian-varieties]{Let $A,B$ be Abelian varieties over a field $k$ and $f: A \rar B$ be a morphism of Abelian varieties. Then $f$ is flat if and only if $f$ is surjective. In this case, one has $\dim{f^{-1}(0_A)} = \dim{A}-\dim{B}$. }

\demo{Because $A$ is reduced, the scheme-theoretic image of $f$ is a reduced closed subscheme of $B$. Since $f$ is surjective, the scheme-theoretic image of $f$ is $B$. By \cite[Lemma 2.3.3.2]{Anantharaman}, $f$ is flat. Conversely, if $f$ is flat, then $f$ is open; but $f$ is proper because $A$ is proper over $k$, so $f$ is closed, whence $f$ is surjective. 

The $k$-schemes $A,B$ are proper and irreducible. Hence, if $f$ is flat surjective, then, by \cite[Corollary 4.3.14]{QL}, $f^{-1}(0_B)$ is pure of dimension $\dim{A}-\dim{B}$. Let $C$ be the connected component of $0_A$ in the closed group subscheme $f^{-1}(0_B)$ of $A$: by \cite[Proposition 0B7R]{Stacks}, $C$ is irreducible (as are $A,B$). Therefore $\dim{C}=\dim{f^{-1}(0_B)}=\dim{A}-\dim{B}$. 
}

\prop[abelian-projective]{An Abelian variety $A$ over a field $k$ is projective.}

\demo{This is \cite[Theorem 7.1]{MilAb}. } 

\defi[isogeny-abelian-variety]{An isogeny of Abelian varieties over a field $k$ is a morphism $f: A \rar B$ of Abelian varieties which is surjective with finite kernel. Its \emph{degree} is the degree of the finite $k$-scheme $\ker{f}$.}

\prop[carac-isogeny]{The following are equivalent, given a morphism $f: A \rar B$ of Abelian varieties over $k$:
\begin{itemize}[noitemsep,label=$-$]
\item $f$ is an isogeny,
\item $\ker{f}$ is finite and $\dim{B}=\dim{A}$,
\item $f$ is surjective and $\dim{B}$ and $\dim{A}$,
\item $f$ is finite flat surjective.
\end{itemize}}

\demo{This is \cite[Theorem 8.1]{MilAb}.}

\prop[torsion-abelian]{Let $A$ be an Abelian variety of dimension $g$ over $k$ and $n \geq 1$ be an integer. The morphism of multiplication by $n$ $[n]_A: A \rar A$ is an isogeny of degree $n^{2g}$. It is \'etale if and only if $n$ is invertible in $k$.}

\demo{This is \cite[Theorem 8.2]{MilAb}. }

\subsection{Image and quotient of Abelian varieties}

\lem[subgroup-is-closed]{Let $G$ be a group scheme over a field $k$ and $j: H \rar G$ be an immersion of group schemes over $k$. Then $j$ is a closed immersion. If $G$ is an Abelian variety, the group scheme law on $H$ is commutative.
}

\demo{Since $H(S) \rar G(S)$ is an injective group homomorphism for any $k$-scheme $S$, so the group law on $H$ is commutative for $H$ if the group law of $G$ is. The morphism $j$ is a closed immersion by \cite[Lemma 047T]{Stacks}.}

To define quotients of Abelian varieties, we need to specify the meaning of this term, or the category in which the quotient is taken. The usual definition uses the fppf topology. More precisely, every scheme is naturally a sheaf on the fppf site, and the quotient is taken in the (Abelian) category of fppf sheaves. 

\defi[fppf-topology-definition]{We define the fppf topology as in \cite[Definition 021M]{Stacks}: a collection of maps $\pi_i: T_i \rar T$ is an fppf-covering if each $\pi_i$ is flat, locally of finite presentation, and $\cup_{i}{\pi_i(T_i)}=T$. It is equivalent to the definition of \cite[Exp. IV, \S 6.3]{SGA3}. }

\prop[abelian-image-exists]{Let $f: A \rar B$ be a morphism of Abelian varieties over a field $k$. There exists an Abelian variety $C$, a faithfully flat morphism $\pi: A \rar C$ and a closed immersion $\iota: C \rar B$ of Abelian varieties over $k$ such that $f=\iota \circ \pi$.}

\demo{By \cite[Exp. V, Th\'eor\`eme 10.1.2]{SGA3}, $C := A/\ker{f}$ is representable by a group scheme over $k$ of finite type and we can write $f=\iota \circ \pi$, where $\pi: A \rar C$ is faithfully flat of finite type and $\iota: C \rar B$ is a monomorphism. The group scheme $C$ is separated by \cite[Lemma 047L]{Stacks}. It is is geometrically regular by descent \cite[Lemma 07NG]{Stacks}, so $C$ is smooth over $k$ by \cite[\S 2.2, Proposition 15]{BLR}. Moreover, $C$ is proper connected since $\pi: A \rar C$ is surjective. Hence $C$ is an Abelian variety. Therefore, $\iota$ is a proper monomorphism hence a closed immersion by \cite[Lemma 04XV]{Stacks}.
}

\prop[quotient-of-abelian-varieties]{Let $A$ be an Abelian variety over a field $k$ and $H$ be a closed subgroup scheme of $A$. Then there exists an Abelian variety $B$ over $k$ and a flat homomorphism of Abelian varieties $f: A \rar B$ such that $\ker{f}=H$. Moreover, $B$ represents the sheaf $A/H$ in the fpqc topology. 
}

\demo{By \cite[Th\'eor\`eme 4.C]{Anantharaman}, the quotient $B := A/H$ of fppf sheaves is representable by a group scheme over $k$. The group scheme $B$ is separated by \cite[Lemma 047L]{Stacks}. We can then apply \cite[Exp. V, Th\'eor\`eme 10.1.2]{SGA3} to the natural surjective morphism $f: A \rar B$ (with kernel $H$): there is a factorization $f = \iota \circ \pi$, where $\pi: A \rar B':= A/H$ is faithfully flat of finite presentation, and $\iota: B' \rar B$ is the natural monomorphism. Then $\iota$ is the identity, so that $f$ is faithfully flat of finite presentation. As in the previous proof, $B$ is geometrically regular and, by \cite[Exp. V, Proposition 9.1]{SGA3}, it is a $k$-scheme of finite type, so it is smooth. Since $f: A \rar B$ is surjective, $B$ is proper over $k$ and connected, so it is an Abelian variety.
}

\subsection{Differentials on Abelian varieties}

\prop[differentials-are-trivial]{Let $G$ be a group scheme over a base $k$. Then the natural homomorphism $H^0(G,\Omega^1_{G/k}) \otimes_{\OO_G(G)} \OO_G \rar \Omega^1_{G/k}$ is an isomorphism. }

\demo{By \cite[Proposition 1.1(b)]{LLR}, $\Omega^1_{G/k}$ is the pull-back of some $k$-vector space, which directly implies the statement. }

\cor[differentials-are-infinitesimal]{Let $A$ be an Abelian variety over a field $k$. The pairing \[\langle \cdot,\,\cdot\rangle_A : (\omega,t) \in H^0(A,\Omega^1_{A/k}) \times T_{0_A} G \longmapsto \omega(0_A)(t) \in k\] is perfect. This pairing is natural in the following sense: if $f: A \rar B$ is a morphism of Abelian varieties over $k$, then $\langle \cdot,\,f(-)\rangle_B = \langle f^{\ast}\cdot,\,-\rangle_A$.}

\demo{The $k$-scheme $A$ is proper and geometrically integral, so $k \rar \OO_A(A)$ is an isomorphism. Hence, by Proposition \ref{differentials-are-trivial}, the evaluation at zero map $H^0(G,\Omega^1_{G/k}) \rar \left(\mathfrak{m}_{G,0}/\mathfrak{m}_{G,0}^2\right)$ is an isomorphism. We conclude using the definition of tangent space. That the pairing is natural is a direct verification. }

\cor[differentials-are-half-exact]{Let $f: A \rar B$ be a closed immersion (resp. a smooth morphism) of Abelian varieties over a field $k$. Then $f^{\ast}: H^0(B,\Omega^1_{B/k}) \rar H^0(A,\Omega^1_{A/k})$ is surjective (resp. injective). }

\demo{This is a consequence of Corollary \ref{differentials-are-infinitesimal} and \cite[Proposition 1.1 (a,c)]{LLR}. }

\cor[differentials-are-other-half-exact]{Let $f: A \rar B$ be a morphism of Abelian varieties and assume that the connected component of unity $C$ of $\ker{f}$ is an Abelian variety over $k$ (equivalently, that it is smooth). Then the image of $f^{\ast}: H^0(B,\Omega^1_{B/k}) \rar H^0(A,\Omega^1_{A/k})$ is exactly the kernel of the restriction map $H^0(A,\Omega^1_{A/k}) \rar H^0(C,\Omega^1_{C/k})$. }

\demo{This is a consequence of Corollary \ref{differentials-are-infinitesimal} and \cite[Proposition 1.1 (a,d)]{LLR}.}

Let us recall the following result. 

\prop[always-smooth-0]{Every group scheme over a field of characteristic zero is smooth. Every surjective morphism $f: A \rar B$ of group schemes of finite type over a field $k$ of characteristic zero is smooth. }

\demo{The first part of the claim is Cartier's theorem \cite[Lemma 047N]{Stacks}. For the second claim, $f$ is flat by \cite[Lemma 2.3.3.2]{Anantharaman}, and moreover $\ker{f}$ is a group scheme of finite type over $k$, so it is smooth. Therefore, the smooth locus of $f$ contains $0_A$. Let $U \subset A$ be an open subscheme such that $f_{|U}$ is smooth. Consider the following commutative diagram:
\[
\begin{tikzcd}[ampersand replacement=\&]
A \times_k U \arrow{d}{+_A} \arrow{r}{(\mrm{id},f)} \& A \times_k B \arrow{r}{(\mrm{pr}_1,f)}\& A \times_k B\arrow{d}{f(\mrm{pr}_1)+\mrm{pr}_2}\\
A \arrow{rr}{f} \& \& B
\end{tikzcd}
\]
The morphism $A \times_k B \overset{f(\mrm{pr}_1)+\mrm{pr}_2}{\longrightarrow} B$ is the composition of $A \times_k B \overset{(\mrm{id},f(\mrm{pr}_1)+\mrm{pr}_2)}{\longrightarrow} A \times_k B$ (which is an isomorphism) with the smooth second projection. Thus, the map $f: A \rar B$ becomes smooth after pre-composing with the smooth surjective morphism $+_A: A \times_k U \rar A$: it is smooth by \cite[Lemma 036U]{Stacks}.
}

\cor[differentials-product]{Let $A$ be an Abelian variety over a field $k$. Let $m, p, q: A \times_k A \rar A$ be the multiplication, the first projection, and the second projection, respectively. Then, for $\omega \in H^0(A,\Omega^1_{A/k})$, one has $m^{\ast}\omega = p^{\ast}\omega+q^{\ast}\omega$.  }

\demo{By Corollary \ref{differentials-are-infinitesimal}, it is enough to show that for any $t \in T_{(0_A,0_A)}(A \times_k A)$, one has \[m_{\ast}(t)=p_{\ast}(t)+q_{\ast}(t).\] Let $\iota_1,\iota_2: A \rar A \times_k A$ be the two inclusions $(\mrm{id},0_A)$ and $(0_A,\mrm{id})$ respectively. When $t=(\iota_1)_{\ast}t_1$ (for $t_1 \in T_{0_A}A$), the identity is correct, since $m \circ \iota_1=p \circ \iota_1=\mrm{id}$ while $q \circ \iota_1$ is the zero morphism. Similarly, the identity holds when $t$ is in the image of $(\iota_2)_{\ast}$. Thus, the linear map $m_{\ast}-p_{\ast}-q_{\ast}: T_{(0_A,0_A)}(A \times_k A) \rar T_{0_A} A$ vanishes on the sum of the images of the tangent maps to $\iota_1$ and $\iota_2$ at $0_A$, so it vanishes. 
}

\cor[abelian-variety-differentials-bimodule]{Let $A$ be an Abelian variety over a field $k$. Then the pull-back by elements of $\mrm{End}(A)$ endows $H^0(A,\Omega^1_{A/k})$ with the structure of a $(k,\mrm{End}(A))$-bimodule. }

\demo{This is direct.}

\cor[differentials-translation-invariant]{Let $A$ be an Abelian variety over a field $k$ and $\omega \in H^0(A,\Omega^1_{A/k})$. Then $\omega$ is translation-invariant.}

\demo{Let $a \in A(k)$ and $t_a: x \in A \longmapsto x+a \in A$. Then, by Corollary \ref{differentials-product}, \[t_a^{\ast}\omega = \mrm{id}^{\ast}\omega+a^{\ast}\omega=\omega+0=\omega.\] Indeed, $a^{\ast}\omega$ is the pull-back by the structure morphism $A \rar \Sp{k}$ of a differential in $\Omega^1_{k/k}$, so is zero. 
}

\prop[differentials-are-closed]{Let $A$ be an Abelian variety over a field $k$ of characteristic distinct from $2$. Then for every $\omega \in H^0(A,\Omega^1_{A/k})$, one has $d\omega=0$. }

\demo{By Proposition \ref{differentials-are-trivial}, $\Omega^2_{A/k}$ is also trivial, so that $H^0(A,\Omega^2_{A/k})$ naturally identifies (by evaluation at $0_A$) with $\Lambda^2{(T_{0_A}A)^{\ast}}$ as $k$-vector spaces. It follows that pull-back by $[2]: A \rar A$ acts by multiplication by $4$ on $H^0(A,\Omega^2_{A/k})$. Let now $\omega \in H^0(A,\Omega^2_{A/k})$, then $4d\omega = [2]^{\ast}d\omega = d([2]^{\ast}\omega)=2d\omega$, hence $d\omega=0$ since $2 \in k^{\times}$.
}

\subsection{Abelian schemes on Dedekind rings in low ramification}

In this section, $R$ denotes a localization of $\Z[1/2]$. 

\lem[generic-to-closed-fibre-immersion]{Let $f: A \rar B$ be a morphism of Abelian schemes over $R$. Suppose that $f_{\Q}: A_{\Q} \rar B_{\Q}$ is a closed immersion. Then $f$ is a closed immersion.}

\demo{This is proved as in the Steps 1-3 of \cite[Proposition 1.2]{FreyMazur}: as in Step 1, there exists $r \geq 1$ and $g_{\Q}: B_{\Q} \rar A_{\Q}$ such that $g_{\Q} \circ f_{\Q}=[r]_A$. By the N\'eron model property, $g_{\Q}$ extends to a morphism $g: B \rar A$ of Abelian schemes. Then $g \circ f$ is an endomorphism of $A$ which agrees with $[r]_A$ on the (dense) generic fibre: since $A$ is a separated scheme, $g \circ f=[r]_A$. Thus $f$ is quasi-finite and is between two proper $R$-schemes, so it is proper, hence $f$ is finite. 
Suppose that $f$ is not a closed immersion: then, by Nakayama's Lemma, $f$ is not a closed immersion on some special fibre. Since $\ker{f}$ is a finite $R$-scheme, it is killed by some $m \geq 1$. But by \cite[Proposition 1.1]{FreyMazur}, $f: A[m] \rar B[m]$ is a closed immersion for every $m$, whence the conclusion.}

\rem{In the proof, if we do not assume that $2$ is invertible in $R$, the proof shows that $f$ is a finite morphism which is a closed immersion away from $2$. }

\prop[quotient-exists]{Let $H \rar G$ be a closed immersion of commutative group schemes over $R$. Assume that $G$ is of finite type over $R$ and that $H$ is flat. Then the quotient $G/H$ (in the category of fppf sheaves) is representable by a separated scheme of finite type over $R$ and $G \rar G/H$ is faithfully flat of finite presentation.  
}

\demo{By \cite[Theorem 4.C]{Anantharaman}, the quotient sheaf (in the fppf topology) $K := G/H$ exists as a commutative group scheme over $R$. By \cite[Exp. V, Theorem 10.1.2]{SGA3} (applied with $X=K$ and the unit section), the natural morphism $G \rar K$ factors through $G \overset{\pi}{\rar} G/H \overset{\iota}{\rar} K$ with $\pi$ faithfully flat of finite presentation, $G/H$ locally of finite type over $R$, and $\iota$ a monomorphism. The morphism $\iota$ is the identity, and $G$ is quasi-compact, thus so is $G/H$. Hence $G/H$ is a group scheme of finite type over $R$.  

Let us check that $G/H$ is separated. It is enough to check $\Delta: G/H \rar G/H \times_S G/H$ is a closed immersion. Since the following diagram
\[\begin{tikzcd}[ampersand replacement=\&]
G \times_S H \arrow{r}{(g,h) \mapsto (g,gh)} \arrow{d}{\mrm{pr}_1} \& G \times_S G\arrow{d}\\
G/H \arrow{r}{\Delta} \& G/H \times_S G/H
\end{tikzcd}  
\]
is Cartesian, it is enough to show by fppf descent that $\iota: (g,h) \in G \times_S H \mapsto (g,gh) \in G \times_S G$ is a closed immersion. Let $f: (g,g') \in G \times_S G \longmapsto (g,g^{-1}g') \in G \times_S G$. Then $f$ is an isomorphism, and $f\circ \iota: (g,h) \in G \times_S H \longmapsto (g,h) \in G \times_S G$ is a base change of the closed immersion $H \rar G$, so is a closed immersion. Hence $\iota$ is a closed immersion, thus $\Delta$ is a closed immersion and $G/H$ is separated. 
}

\prop[image-over-dim1]{Let $f: A \rar B$ be a morphism of Abelian schemes over $R$. There exists a factorization $f = \iota \circ \pi$, where $\pi: A \rar C$ is a faithfully flat morphism of Abelian schemes, and $\iota$ is a closed immersion.}

\demo{Let $C$ be the N\'eron model over $R$ of the scheme-theoretic image $C_0$ of $f_{\Q}$ (this is well-defined by Proposition \ref{abelian-image-exists} and \cite[\S 1.4 Theorem 3]{BLR}). Then $f_{\Q}$ factors as $\iota_0\circ\pi_0$, where $\pi_0: A_{\Q} \rar C_0,\iota_0: C_0 \rar B_{\Q}$ are morphisms of Abelian varieties with $\pi_0$ flat and $\iota_0$ a closed immersion. Then $\pi_0, \iota_0$ extend to morphisms $\pi: A \rar C, \iota: C \rar B$ of Abelian schemes. Since $\iota \circ \pi$ agrees with $f$ on the (dense) generic fibre of the separated scheme $A$, one has $f=\iota \circ \pi$. Moreover, $\iota$ is a closed immersion by Lemma \ref{generic-to-closed-fibre-immersion}. 
All we need to do is prove that $\pi$ is flat and surjective. Now, $\pi$ is between proper $R$-schemes, so it is proper, thus it has closed image. Moreover, the image of $\pi$ contains the generic fibre of $C$, so $\pi$ is surjective. By \cite[Lemma 2.3.3.2]{Anantharaman}, every fibre of $\pi$ is flat, so, by the ``crit\`ere de platitude par fibres'' \cite[Lemma 039D]{Stacks}, $\pi$ is faithfully flat.  
}

\prop[kernel-over-dim1]{Let $f: A \rar B$ be a morphism of Abelian schemes over $R$. Then $\ker{f}$ is a proper flat group scheme over $R$, and there is a short exact (for the fppf topology) sequence $0 \rar C \rar \ker{f} \rar H \rar 0$, where $H$ is the finite flat $R$-group scheme $\Sp{\OO(\ker{f})}$, $C$ is an Abelian scheme over $R$ and $C \rar \ker{f}$ is a closed immersion. }

\demo{Write $f=\iota \circ \pi$ with $\pi: A \rar D$ proper faithfully flat, $\iota: D \rar B$ a closed immersion and $D$ an Abelian scheme over $R$ by Proposition \ref{image-over-dim1}. Then $\ker{f}=\ker{\pi} = A \times_{D,0_D} \Sp{R}$ is proper flat over $\Sp{R}$. Let $\Gamma=\Sp{\OO(\ker{f})}$, 

Let $C$ be the N\'eron model over $R$ of the connected component of unity $C_0$ of the generic fibre of $\ker{f}$. By the N\'eron-Ogg-Shafarevich criterion \cite[Theorem 1]{GoodRed}, $C_0$ has good reduction at every maximal ideal of $v$, so $C$ is an Abelian scheme over $R$ (in particular, it is faithfully flat, thus so is $\ker{f}$). By the N\'eron mapping property, there is a morphism of Abelian schemes $j: C \rar A$ extending the inclusion $C_0 \subset A_{\Q}$, and it is a closed immersion by Lemma \ref{generic-to-closed-fibre-immersion}. The morphisms $0_B: C \rar B$ and $f \circ j$ agree on the generic fibre, so they are equal, and $j$ factors through $\ker{f}$. 

By Proposition \ref{quotient-exists}, the quotient sheaf (in the fppf topology) $H := \ker{f}/C$ exists as a separated commutative group scheme of finite type over $R$, and $\pi': \ker{f} \rar H$ is faithfully flat with finite presentation. Since $\ker{f}$ is proper faithfully flat over $R$, the ``crit\`ere de platitude par fibres'' \cite[Lemma 039D]{Stacks} implies that $H$ is proper flat over $R$. By \cite[Lemma 0D4J]{Stacks} the dimension of its fibres is constant. By construction, the generic fibre of $H$ has finitely many points, so $H$ is proper quasi-finite flat over $\Sp{R}$, hence $H$ is finite flat over $\Sp{R}$. 

Finally, since $\ker{f}$ is a proper flat group scheme over $\Sp{R}$, $\OO(\ker{f})$ is a finite torsion-free $R$-algebra, hence it is flat and the natural map $\OO(\ker{f}) \otimes_R \OO(\ker{f}) \rar \OO(\ker{f} \times_R \ker{f})$ is an isomorphism. This implies that $\Sp{\OO(\ker{f})}$ is a finite flat group scheme over $R$ and $\ker{f} \rar \Sp{\OO(\ker{f})}$ is a morphism of flat group schemes. Since $C$ is an Abelian scheme over $R$, $\OO(C) \simeq R$, thus $C \rar \ker{f} \rar \Sp{\OO(\ker{f})}$ is the zero homomorphism. Thus $\ker{f} \rar \Sp{\OO(\ker{f})}$ factors through $H \rar \Sp{\OO(\ker{f})}$. 

Thus, the composition of ring homomorphisms $\OO(\ker{f}) \rar \OO(H) \rar \OO(\ker{f})$ is the identity, hence $\OO(\ker{f}) \rar \OO(H)$ is injective. Moreover, since $H(\Qbar)$ is exactly the group of connected components of $(\ker{f})_{\Qbar}$, one has $\dim{\OO(H) \otimes \Qbar} = \dim{\OO(\ker{f}) \otimes \Qbar}$, so $\OO(H)$ and $\OO(\ker{f})$ are free $R$-modules of same rank. Hence $\OO(H) \rar \OO(\ker{f})$, being surjective, is an isomorphism, and therefore $\ker{f} \rar H$ is exactly the Stein factorization of $\ker{f} \rar \Sp{R}$.}

\section{Jacobians of relative curves}
\label{jacobian-relative}

\subsection{Functorialities in the Picard groups}
\label{picard-functorialities}

\lem[locally-free-base-change]{Let $A$ be a ring and $B$ a finite locally free $A$-algebra of constant rank $r$, and $M$ be a finite locally free $B$-module of constant rank $d$. There exists $f_1,\ldots,f_r \in A$ generating the unit ideal such that, for each $f_i$, $B_{f_i}$ is a free $A_{f_i}$-module and $M_{f_i}$ is a free $B_{f_i}$-module.}

\demo{Since $M$ is a flat $A$-module of finite presentation, it is locally free and has constant rank $dr$. Let $\mathfrak{p} \in \Sp{A}$, then $M_{\mathfrak{p}}$ is a finite locally free module of rank $r$ over the semi-local ring $B_{\mathfrak{p}}$. By \cite[Lemma 02M9]{Stacks}, $M_{\mathfrak{p}}$ is free of rank $r$ over $B_{\mathfrak{p}}$. Hence there exists $m_1,\ldots,m_r \in M$ such that $\pi: (b_1,\ldots,b_r) \in B \longmapsto \sum_{i=1}^r{b_im_i} \in M$ is an isomorphism after tensoring by $\mathfrak{p}$.

By Nakayama, there exists some $f \in A \backslash \mathfrak{p}$ such that the map $\pi \otimes_A A_f$ is surjective. By \cite[Lemma 0519]{Stacks}, $\ker{\pi \otimes_A A_f}$ is finitely generated and vanishes when tensoring with $A_{\mathfrak{p}}$, so by Nakayama, it vanishes after tensoring with $A_{fg}$ for some $g \in A\backslash \mathfrak{p}$. Thus $\pi \otimes_A A_{fg}$ is an isomorphism, whence the conclusion.}

\lem[locally-free-finite-module]{Let $A$ be a ring and $M$ be an invertible $B$-module, where $B$ is a finite locally free $A$-algebra of constant rank. Then $P=\mrm{Hom}_A(\det_A{B},\det_A{M})$ is an invertible $A$-module. For any $m \in M$, let $\phi_m: B \rar M$ be the $B$-linear morphism mapping $1$ to $m$, then the $A$-linear homomorphism $\Phi: [m] \in A[M] \longmapsto \det_A{\phi_m} \in P$ is onto, and its kernel is generated by the $[bm]-N_{B/A}(b)[m]$ for all $b \in B, m \in M$.}

\demo{By Lemma \ref{locally-free-base-change}, $\Sp{A}$ can be covered by principal open subschemes $D(f)$ such that $B_f$ is free over $A_f$ and $M_f$ is free over $B_f$. In the rest of the proof, we will call such open subschemes \emph{good}. In particular, $M$ is locally free over $A$.

Since $B,M$ are locally free $A$-modules, $P$ is an invertible $A$-module. Since $M$ is invertible over $B$, $B$ and $M$ have the same constant rank over $A$, hence $\Phi$ is well-defined, and it is clear that $\Phi([bm]-N_{B/A}(b)[m])=0$ for all $b \in B, m \in M$. So let $J \subset A[M]$ be the submodule of $A[M]$ generated by the $[bm]-N_{B/A}(b)[m]$ for all $b \in B, m \in M$. The goal is to prove that the sequence 
\[0 \rar J \rar A[M] \overset{\Phi}{\rar} P \rar 0\] is exact. 

Let $D(f) \subset \Sp{A}$ be a good principal open subscheme. Let us prove that $\Phi \otimes_A A_f$ is surjective with kernel generated over $A_f$ by $J$. Indeed, there exists some $m \in M$ such that $b \in B_f \longmapsto bm \in M_f$ is an isomorphism. In particular, $\det_A{\phi_m}: \det_A{B} \rar \det_A{M}$ is an isomorphism after tensoring by $A_f$. It follows that $P_f$ has $\frac{\det_A{\phi_m}}{1}$ as a $A_f$-basis, thus $\Phi \otimes_A A_f$ is surjective. 

Let $m' \in M$ be any element: since $M_f$ is free over $B_f$ with basis $m$, there exists $r \geq 1$ and $b \in B$ such that $m'=\frac{b}{f^r}m$ in $M_f$. Therefore, for some $s \geq 1$, $bf^sm=f^{r+s}m'$ in $M$, so that $N_{B/A}(b)f^{sd}[m]-f^{(r+s)d}[m'] \in J$, so $[m'] \in J_f+\frac{N_{B/A}(b)}{f^{rd}}[m]$, whence the conclusion.} 

\cor[locally-free-finite-sheaf]{Let $f: Y \rar X$ be a finite locally free map and $\mathcal{L}$ be an invertible $\OO_Y$-module. Let $\mathcal{N}: f_{\ast}\OO_Y \rar \OO_X$ be the norm map. Let $(U_i)_i$ be a basis of the topology of $X$ constituted by affine open subsets such that $f: f^{-1}(U_i) \rar U_i$ has constant rank. Then $f_{\ast}\mathcal{L}$ is a locally (with respect to the space $X$) free $f_{\ast}\OO_Y$-module, and $\mathcal{P} := \det{f_{\ast}\mathcal{L}} \otimes (\det{f_{\ast}\OO_Y})^{-1}$ is an invertible $\OO_X$-module. Given any quasi-coherent $\OO_X$-module $\mathcal{M}$, to uniquely define a morphism $\alpha: \mathcal{P} \rar \mathcal{M}$ of $\OO_X$-modules, it suffices to do the following:
\begin{itemize}[noitemsep,label=$-$]
\item describe a rule attaching to any $s \in \mathcal{L}(f^{-1}(U_i))$ some $\phi(s) \in \mathcal{M}(U_i)$,
\item check that for any $b \in \OO_Y(f^{-1}(U_i)), s \in \mathcal{M}(f^{-1}(U_i))$, $\phi(bs) = \mathcal{N}(b)\phi(s)$,
\item check that for any $U_j \subset U_i$ and any $s \in \mathcal{L}(f^{-1}(U_i))$, $\phi(s)_{|U_j} = \phi(s_{|f^{-1}(U_j)})$. 
\end{itemize}

Then, for all $s \in \mathcal{L}(f^{-1}(U_i))$, $\alpha$ maps the determinant of the map $f_{\ast}s: (f_{\ast}\mathcal{O}_Y)_{|U} \rar (f_{\ast}\mathcal{L})_{|U}$ of locally free $\OO_U$-modules to $\phi(s)$.}

\demo{This is a formal consequence of the previous result.}

\prop{Let $f: Y \rar X$ be a finite locally free morphism of schemes and $\mathcal{L}$ be an invertible module. The formation of $f_{\ast}\mathcal{L}$ commutes with arbitrary base change. Moreover, the formation of the norm morphism $f_{\ast}\OO_Y \rar \OO_X$ commutes with arbitrary base change.}

\demo{This is classical, see for instance \cite[Lemma 02KG,Lemma 0BD2]{Stacks}.}

\prop[scheme-pushforward-morphism]{Let $\mathcal{M},\mathcal{N}$ be two invertible modules over a scheme $Y$, and assume that $f:Y \rar X$ is finite locally free. Let $\omega=\det{f_{\ast}\OO_Y}$. Write $\mathcal{M}'=\underline{Hom}_{\OO_X}(\omega,\det{f_{\ast}\mathcal{M}})$, and let us adopt a similar notation for $\mathcal{N}'$. There exists a canonical isomorphism of invertible $\OO_X$-modules \[T: \mathcal{M}' \otimes_{\OO_X}\mathcal{N}' \rar \underline{Hom}_{\OO_X}(\omega,\det{f_{\ast}\mathcal{M} \otimes\mathcal{N}}).\]
For any open subset $U \subset X$ and any $m \in \mathcal{M}(f^{-1}(U)),n \in \mathcal{N}(f^{-1}(U))$, $T$ maps $\det{f_{\ast}m} \otimes \det{f_{\ast}n}$ to $\det{f_{\ast}(m \otimes n)}$. Moreover, the formation of $T$ commutes with any base change with respect to $X$.}

\demo{The uniqueness is clear by the property. Note that to check the property, we may assume that $U$ is affine and $f^{-1}(U) \rar U$ has constant rank.
To construct $T$, consider the rule mapping some $s \in \mathcal{M}(f^{-1}(U))$ to the determinant of the morphism $f_{\ast}(m \otimes -): (f_{\ast}\mathcal{N})_{|U} \rar (f_{\ast}(\mathcal{M} \otimes \mathcal{N}))_{|U}$ of locally free $\OO_U$-modules of same rank. By Corollary \ref{locally-free-finite-sheaf}, this rule defines a morphism $\mathcal{M}' \rar \underline{Hom}_{\OO_X}(\det{f_{\ast}\mathcal{N}},\det{f_{\ast}(\mathcal{M}\otimes \mathcal{N})})$. 

Since $\det{f_{\ast}\mathcal{N}}$, $\mathcal{P}:=\det{f_{\ast}(\mathcal{M}\otimes \mathcal{N})}$ are invertible, this is equivalent to defining a morphism $\mathcal{M}' \otimes \mathcal{N}' \rar \underline{Hom}_{\OO_X}(\omega,\mathcal{P})$. One can then easily check that this morphism satisfies the stated property. To see that $T$ is an isomorphism, it is enough to check that it is surjective. We can work Zariski-locally on $X$ and assume that $X$ is affine, $f_{\ast}\mathcal{M},f_{\ast}\mathcal{N},f_{\ast}(\mathcal{M}\otimes\mathcal{N})$ are free over $f_{\ast}\OO_Y$ and that $f_{\ast}\OO_Y$ is free over $\OO_X$. The conclusion follows by choosing bases $m,n$ of $\mathcal{M},\mathcal{N}$.}

\prop[scheme-pushforward-functor]{Let $f: X \rar Y, g: Y \rar Z$ be finite locally free morphisms of schemes, of constant rank. Let $\mathcal{L}$ be a line bundle on $X$. Then there is a canonical isomorphism 
\[\underline{Hom}_{\OO_Z}(\det{(g_{\ast}f_{\ast})\OO_X},\det{g_{\ast}f_{\ast}\mathcal{L}}) \rar \underline{Hom}_{\OO_Z}(\det{g_{\ast}(\det{f_{\ast}\OO_X})},\det{g_{\ast}(\det{f_{\ast}\mathcal{L}})})\]
 given by the following rule: for any open subset $U \subset X$ and any $s \in \mathcal{L}(f^{-1}(g^{-1}(U)))$, it maps $\det{(g \circ f)_{\ast}s}$ to $\det{g_{\ast}(\det{f_{\ast}s})}$. The formation of this isomorphism commutes with base change.}
 
\demo{That the morphism exists follows from Corollary \ref{locally-free-finite-sheaf}. The commutation with base change is easy (using in Corollary \ref{locally-free-finite-sheaf} an affine cover coming from the original $X \rar Y \rar Z$). Since its domain and range are invertible $\OO_Z$-modules, we only need to prove that the morphism is surjective. This is easy by choosing local sections $s$ generating $\mathcal{L}(f^{-1}(g^{-1}(U)))$ over $\OO_X(f^{-1}(g^{-1}(U)))$, so that $\det{g_{\ast}\det{f_{\ast}s}}$ is a basis of the range above $U$.}

\prop[functoriality-line-bundles]{Let $f: Y \rar X$ be a morphism of schemes, finite free of constant rank. Then 
\[f^{\ast}: \mathcal{L} \in \operatorname{Pic}(X) \longmapsto f^{\ast}\mathcal{L}\in \operatorname{Pic}(Y),\, [f_{\ast}]: \mathcal{L} \in \operatorname{Pic}(Y) \longmapsto \det{f_{\ast}\mathcal{M}} \otimes (\det{f_{\ast}\OO_Y})^{-1}\] are functorial group homomorphisms. }

\demo{This is clear for $f^{\ast}$. It is clear that $\mrm{id}_{\ast}$ is exactly the identity. That $f_{\ast}$ is a group homomorphism follows from Proposition \ref{scheme-pushforward-morphism}. To show functoriality, let $f: X \rar Y$ and $g: Y \rar Z$ be finite morphisms, locally free of constant rank, and $\mathcal{L}$ be a line bundle on $X$. 
Then, using Proposition \ref{scheme-pushforward-functor},
\begin{align*}
[g_{\ast}][f_{\ast}]\mathcal{L} &\simeq [g_{\ast}](\det{f_{\ast}\mathcal{L}} \otimes (\det{f_{\ast}}\OO_X)^{-1}) \simeq [g_{\ast}](\det{f_{\ast}\mathcal{L}}) \otimes ([g_{\ast}]((\det{f_{\ast}}\OO_X)^{-1}))^{-1}\\
&\simeq \det{g_{\ast}\det(f_{\ast}\mathcal{L})} \otimes (\det{g_{\ast}\det(f_{\ast}\OO_X)})^{-1} \\
&\simeq \det{(g \circ f)_{\ast}\mathcal{L}} \otimes (\det{(g \otimes f)_{\ast}\OO_X})^{-1}.
\end{align*} 
}

\lem[pushforward-is-weil]{Let $f: X \rar Y$ be a finite flat morphism of finite type reduced schemes over a field $k$, with $X,Y$ pure of dimension one (hence Cohen-Macaulay). Let $\mathcal{L}$ be a line bundle on $X$ contained in $\mathcal{K}_X$, so that it defines a Cartier divisor on $X$. Then the line bundle $\mathcal{M}=[f_{\ast}]\mathcal{L}$ on $Y$ has a natural embedding into $\mathcal{K}_Y$, which defines a Cartier divisor on $Y$. 
Let $S_X$ (resp. $S_Y$) be the finite set of (closed points) $x \in X$ (resp. $y \in Y$) such that $\mathcal{L}_{X,x} \neq \OO_{X,x}$ (resp. $\mathcal{M}_{Y,y} \neq \OO_{Y,y}$). Then $S_Y \subset f(S_X)$. If moreover every point of $S_X$ is normal, then every point of $f(S_Y)$ is normal, and, for every $y \in S_Y$, $\mrm{mult}_y{\mathcal{M}}=\sum_{x \in f^{-1}(y) \cap S_X}{\mrm{mult}_x{\mathcal{L}} \cdot [\kappa(x):\kappa(y)]}$, where the multiplicity is defined as in \cite[Definition 7.1.27]{QL}. 
}

\demo{We may assume that $Y$ is connected, so that $f$ has constant rank $d$. Our assumptions imply that the only functions on $X,Y$ are that are not regular at every point are those that vanish at some generic point. In particular, the sheaves of meromorphic functions on $X,Y$ are well-behaved \cite[Lemma 0EMF]{Stacks}. Moreover, given an affine open subscheme $U \subset Y$, a section in $\OO_X(f^{-1}(U))$ is a regular element (resp. invertible) if and only if its norm over $\OO_Y(U)$ is a regular element (resp. invertible). This implies that the norm $f_{\ast}\OO_X \rar \OO_Y$ induces compatible $f_{\ast}\mathcal{G}_X \rar \mathcal{G}_Y$ for $\mathcal{G}=\mathcal{K},\mathcal{K}^{\times},\OO^{\times}$. Note that since $f$ is flat, we have a pull-back of meromorphic functions \cite[Lemma 02OU]{Stacks}. The natural map $(f_{\ast}\mathcal{K}_X^{\times})/(f_{\ast}\OO_X^{\times}) \rar f_{\ast}(\mathcal{K}_X^{\times}/\OO_X^{\times})$ is an isomorphism (by checking at the stalks and using the fact that the sheaves of meromorphic functions are well-defined). 

Therefore, we have a homomorphism of sheaves in Abelian groups $\nu: f_{\ast}(\mathcal{K}_X^{\times}/\OO_X^{\times}) \rar \mathcal{K}_Y^{\times}/\OO_Y^{\times}$. The embedding $\mathcal{L} \subset \mathcal{K}_X$ corresponds to a section $s_{\mathcal{L}}$ of $\mathcal{K}_X^{\times}/\OO_X^{\times}$, and let $s_{\mathcal{M}}$ be its image in $\mathcal{K}_Y^{\times}/\OO_Y^{\times}$. Then $s_{\mathcal{M}}$ defines a line bundle $\mathcal{M}$ on $Y$. The claim that we wish to prove is that $\mathcal{M} \simeq (\det{f_{\ast}\OO_X})^{-1} \otimes \det{f_{\ast}\mathcal{L}}$. 

We define a map $\alpha:  (\det{f_{\ast}\OO_X})^{-1} \otimes \det{f_{\ast}\mathcal{L}} \rar \mathcal{M}$ as follows: given any affine open subset $U \subset Y$ small enough such that $\mathcal{L}(f^{-1}(U))$ is free over $\OO_X(f^{-1}(U))$, for any $s \in \mathcal{L}(U)$, we write $s=fs_1$ for some $f \in \OO_X(f^{-1}(U))$, where $s_0$ is any $\OO_X(f^{-1}(U))$-basis of $\mathcal{L}(f^{-1}(U))$, and we map $\det{f_{\ast}s}$ to $N_{X/Y}(fs_1) \in \OO_Y(U)s_{\mathcal{M}}=\mathcal{M}(U)$. This does not depend on the choice of factorization of $s=fs_0$, and satisfies the compatibility properties of Corollary \ref{locally-free-finite-sheaf}. One can also check locally that $\alpha$ is surjective, thus an isomorphism, proving the first part of the claim. 

Note that if $s_{\mathcal{L},z}=1$ (equivalently, $\mathcal{L}_{X,z}=\OO_{X,z}$ as subgroups of $\mathcal{K}_{X,z}$, or $z \notin S_X$) for all $z \in f^{-1}(y)$ for some $y \in Y$, then $f_{\ast}\OO_X$ and $f_{\ast}\mathcal{L}$ are quasi-coherent subsheaves of $f_{\ast}\mathcal{K}_X$ (itself quasi-coherent) whose stalks agree at $x$, hence their stalks agree in a neighborhood of $y$. Therefore, the stalk of $f_{\ast}s_{\mathcal{L}}$ at $y$ is trivial, hence so is the stalk of $s_{\mathcal{M}}$ at $y$. This implies that $S_Y \subset f(S_X)$. 

If every point $x$ of $S_X$ is normal, then for any $y \in S_Y$, $y=f(x)$ for some $x \in S_X$, so that $\OO_{Y,y} \rar \OO_{X,x}$ is faithfully flat with $\OO_{X,x}$ normal. Thus $\OO_{Y,y}$ normal. 

The final part of the statement is local with respect to $Y$; since $S_Y$ is contained in the normal points of $Y$, we may restrict to affine open subsets contained in the interiors of the irreducible components of $Y$, thereby assuming that $Y$ is integral. We are then reduced to the following problem: let $A$ be an integral $k$-algebra of dimension one, and $\mathfrak{m}$ be a maximal ideal such that $A_{\mathfrak{m}}$ is normal. Let $B$ be a finite flat reduced $A$-algebra. Let $s \in (B_{reg}^{-1}B)^{\times}$, and $S$ be the set of maximal ideals of $B$ above $\mathfrak{m}$ such that $sB_{\mathfrak{p}} \neq B_{\mathfrak{p}}$. Assume that $B_{\mathfrak{p}}$ for $\mathfrak{p} \in S$ is normal. Then $v_{\mathfrak{m}}(N_{B/A}(s))=\sum_{\mathfrak{p} \in S}{f_{\mathfrak{p}}v_{\mathfrak{p}}(s\cdot 1_{B_{\mathfrak{p}})}}$, where $f_{\mathfrak{p}}$ is the degree of the extension of residue fields $A/\mathfrak{m} \rar B/\mathfrak{p}$. In other words, $f_{\mathfrak{p}}\dim_{k}{A/\mathfrak{m}} = \dim_k{B/\mathfrak{p}}$. 

Note that both sides of the identity to prove are additive with respect to $s$ (seen as an element of the group $(B_{reg}^{-1}B)^{\times}$). 

By the prime avoidance lemma, for every $\mathfrak{p} \in S$, we can find some $u_{\mathfrak{p}} \in B$ contained in $\mathfrak{p}$, but not in $\mathfrak{p}^2$ or in any other prime ideal of $B$ lying above $\mathfrak{m}$ or $(0)$. Then $u_{\mathfrak{p}} \in B_{reg}$, and $u_{\mathfrak{p}}$ is a generator of $\mathfrak{p}B_{\mathfrak{p}}$. Let $s'=s \cdot \prod_{\mathfrak{p} \in S}{u_{\mathfrak{p}}^{-v_{\mathfrak{p}}(s)}}$, then $s' \in (B_{reg}^{-1}B)^{\times}$ is such that for each maximal $\mathfrak{p}$ of $B$ above $\mathfrak{m}$, $s'B_{\mathfrak{p}} = B_{\mathfrak{p}}$. Thus $s': B/\mathfrak{m} \rar B/\mathfrak{m}$ is injective between finite dimensional $k$-vector spaces, hence surjective. Thus the multiplication by $s'$ is an endomorphism of the free $A_{\mathfrak{m}}$-module of finite rank $B_{\mathfrak{m}}$ which is surjective after tensoring by $A/\mathfrak{m}$. Hence it is an isomorphism by Nakayama's Lemma, and its determinant (which is $N_{B/A}(s')$) is invertible.

Thus the identity holds for $s'$. It is therefore enough to prove the identity assuming that $s=u_{\mathfrak{p}}$, $S=\{\mathfrak{p}\}$ for some $\mathfrak{p}$, ie that $v_{\mathfrak{m}}(N_{B/A}(u_{\mathfrak{p}}))=f_{\mathfrak{p}}$. Since the multiplication by $u_{\mathfrak{p}}$ is an injective endomorphism of the free $A_{\mathfrak{m}}$-module of finite rank $B_{\mathfrak{m}}$, with $A_{\mathfrak{m}}$ being a discrete valuation ring, it is given in suitable bases of $B_{\mathfrak{m}}$ by a diagonal matrix with diagonal $(\alpha_i)_{1 \leq i \leq n}$ of nonzero elements of $A_{\mathfrak{m}}$. Then one has 

\begin{align*}
f_{\mathfrak{p}} &= \mrm{len}_{A/\mathfrak{m}}(B/\mathfrak{p}) = \mrm{len}_{A/\mathfrak{m}}(B_{\mathfrak{m}}/u_{\mathfrak{p}}B_{\mathfrak{m}})\\
&= \mrm{len}_{A/\mathfrak{m}}\left(\bigoplus_{i=1}^n{A/\alpha_iA}\right) = \sum_{i=1}^n{\mrm{len}_{A/\mathfrak{m}}(A/\alpha_i)} = \sum_{i=1}^n{v_{\mathfrak{m}}(\alpha_i)}\\
&= v_{\mathfrak{m}}\left(\prod_{i=1}^n{\alpha_i}\right) = v_{\mathfrak{m}}(\det[u_{\mathfrak{p}}: B_{\mathfrak{m}} \rar B_{\mathfrak{m}}]) = v_{\mathfrak{m}}(N_{B/A}(u_{\mathfrak{p}})),
\end{align*}

whence the conclusion.}

\rem{This lemma states, in other words, that the notion of push-forward we defined functions as the natural push-forward at the level of Weil divisors.}

\prop[mixed-functoriality]{Consider the following commutative diagram of schemes, where all the maps are finite locally free:
\[
\begin{tikzcd}[ampersand replacement=\&]
Y \arrow{r}{f} \arrow{d}{v} \& X \arrow{d}{u}\\
Y' \arrow{r}{f'} \& X'
\end{tikzcd}
\]
There are then two morphisms $\operatorname{Pic}(Y') \rar \operatorname{Pic}(X)$, given by $u^{\ast}[(f')_{\ast}]$ and $[f_{\ast}]v^{\ast}$. These morphisms are equal if, for any affine open subscheme $U \subset X'$, and any $b \in \OO_{Y'}((f')^{-1}(U))$, one has $u^{\sharp}(N_{Y'/X'}(b))=N_{Y/X}(v^{\sharp}(b))$. 

This hypothesis is verified if there exists an open subscheme $U \subset X'$ such that the diagram is Cartesian above $U$ and $\OO_{X'} \rar \iota_{\ast}\OO_{X'_U}$ is injective, where $\iota$ is the open immersion $U \rar X'$.
}

\demo{We first prove the second part of the statement. Since $Y',Y,X$ are flat over $X'$, the morphisms $\OO_Z \rar \iota_{\ast}\OO_{Z_U}$ are injective for all each $Z \in \{X,Y',Y\}$; it follows that the norms satisfy the required compatibility. 

Now we discuss the first part. Let $\mathcal{L}$ be a line bundle on $X$, then one directly shows that 
\begin{align*}
u^{\ast}[(f')_{\ast}]\mathcal{L} &\simeq u^{\ast}\underline{Hom}_{\OO_{X'}}(\det{(f')_{\ast}\OO_{Y'}},\det{(f')_{\ast}\mathcal{L}}),\\
[f_{\ast}]v^{\ast}\mathcal{L} &\simeq \underline{Hom}_{\OO_X}(\det{f_{\ast}\OO_Y},\det{f_{\ast}v^{\ast}\mathcal{L}}).
\end{align*}

We first define our morphism $\Psi: u^{\ast}[(f')_{\ast}]\mathcal{L} \rar [f_{\ast}]v^{\ast}\mathcal{L}$ by adjunction: it is enough to define a morphism of quasi-coherent sheaves 
\[\Phi: \underline{Hom}_{\OO_{X'}}(\det{(f')_{\ast}\OO_{Y'}},\det{(f')_{\ast}\mathcal{L}}) \rar \underline{Hom}_{\OO_X}(\det{f_{\ast}\OO_Y},\det{f_{\ast}v^{\ast}\mathcal{L}}).\]

For any affine open subset $U \subset X'$, and $s \in \mathcal{L}((f')^{-1}(U))$, we map $\det{(f')_{\ast}s}$ to $\det{f_{\ast}(v^{\ast} s)}$, where $v^{\ast}s \in (v^{\ast}\mathcal{L})(v^{-1}((f')^{-1}(U)))=\mathcal{L}(f^{-1}(u^{-1}(U)))$, so $f_{\ast}v^{\ast}s$ defines a morphism of locally free $\OO_{u^{-1}(U)}$-modules of same rank $(f_{\ast}\OO_Y)_{|u^{-1}(U)} \rar (f_{\ast}v^{\ast}\mathcal{L})_{|u^{-1}(U)}$. The assumption on the norm shows Corollary \ref{locally-free-finite-sheaf} applies and that $\Phi$, hence $\Psi$, is well-defined. 

To show that our morphism is an isomorphism, we work Zariski-locally on $X'$, so we can assume that we are given finite free $A$-algebras $B$, $C$, a $B \otimes_A C$-algebra $D$ (corresponding respectively to $X'$, $Y'$, $X$, $Y$) which is locally free over $B$ and $C$, a free $B$-module $M$ of rank one (corresponding to $\mathcal{L}$), and such that, for any $b \in B$, the image in $C$ of $N_{B/A}(b)$ is exactly $N_{D/C}(b)$, where $N$ denotes the norm. 

In this situation, $\Psi: \mrm{Hom}_A(\det_A{B},\det_A{M}) \otimes_A C \mrm{Hom}_C(\det_C{D},\det_C{M \otimes_B D})$, is given by $\det_A{m} \otimes 1 \longmapsto \det_C{m \otimes_B D}$ for any $m \in M$, where $m$ is identified with the $A$-homomorphism of locally free $A$-modules of same rank $B \rar M$ mapping any $b \in B$ to $bm$. Let $m \in M$ be a $B$-basis, then the source of $\Psi$ (resp. its codomain) is free over $C$ of basis $\det_A{m}$ (resp. $\det_C{m \otimes_B D}$), so that $\Psi$ is indeed an isomorphism.
}

\subsection{Degrees and Jacobian varieties}
\label{degrees-jacobians}

\defi[defi-degree]{Let $k$ be a field and $X$ be a proper $k$-scheme, pure of dimension one. The degree of a locally free module $\mathcal{L}$ on $X$ of rank $r: \pi_0(X) \rar \N$ is the function $\pi_0(X) \rar \Z$ given, for every (open and closed) connected component $C$ of $X$, by \[\deg{\mathcal{L}}(C)=\chi(C,\mathcal{L}_{|C})-r(C)\chi(C,(\OO_X)_{|C}),\] where $\chi$ denotes the Euler-Poincar\'e characteristic.  }

Note that this definition differs slightly from \cite[Section 0AYQ]{Stacks}, because we treat connected components separately instead of adding the results up. This is why we recall some elementary statements about degrees. 

\prop[deg-elem]{In the situation of Definition \ref{defi-degree}, $\mathcal{L}$ and $\det{\mathcal{L}}$ have the same degree, where $\det{\mathcal{L}}$ is the line bundle on $X$ whose restriction to any connected component $C$ is given by $\Lambda^{r(C)}\mathcal{L}$.
Moreover, $\deg$ defines a group homomorphism $\operatorname{Pic}(X) \rar \Z^{\pi_0(X)}$.}

\demo{The first statement is exactly the content of \cite[Lemma 0DJ5]{Stacks}, applied to each connected component of $X$ separately. For the second claim, let $\mathcal{L},\mathcal{M}$ be line bundles on $X$. Then one has 
\[\deg{\mathcal{L}}+\deg{\mathcal{M}} = \deg(\mathcal{L} \oplus \mathcal{M}) = \deg\left[\det(\mathcal{L} \oplus \mathcal{M})\right] = \deg(\mathcal{L} \otimes \mathcal{M}).\]}

\lem[deg-div]{In the situation of Definition \ref{defi-degree}, if $X$ is connected and $\mathcal{L}$ is the sheaf of ideals defining a closed subscheme $Z$, finite over $k$ of degree $d$, then $\deg{\mathcal{L}}=-d$.   }

\demo{Let $\iota: Z \rar X$ be the natural closed immersion: we have an exact sequence \[0 \rar \mathcal{L} \rar \OO_X \rar \iota_{\ast}\OO_Z \rar 0,\] which implies that \[\chi(X,\OO_X)=\chi(X,\mathcal{L})+\chi(X,\iota_{\ast}\OO_Z)=\chi(X,\mathcal{L})+\chi(Z,\OO_Z).\] Since $Z$ is finite over $\Sp{k}$, it is a discrete topological space and $\chi(Z,\OO_Z)=$ has finite support, $\chi(Z,\iota_{\ast}\OO_Z)=\dim_k{H^0(Z,\OO_Z)}=d$, and the conclusion follows.}

\lem[deg-after-field-ext]{Let $X$ be a proper $k$-scheme, pure of dimension one, and $\mathcal{L}$ be a line bundle on $X$ with degree zero. Then, for any field extension $k'$ of $k$ and any connected component $C$ of $X_{k'}$, the degree of $f^{\ast}\mathcal{L}$ (where $f$ is the morphism $X_{k'} \rar X$) is zero. Conversely, if $\mathcal{L}$ is a line bundle on $X$ such that for any connected component $C$ of $X_{k'}$, $\deg{f^{\ast}\mathcal{L}}(C)=0$, then $\mathcal{L}$ has degree zero on any connected component of $X$. }

\demo{First, by \cite[Lemma 0B59]{Stacks}, one has $\sum_{C}{\deg{f^{\ast}\mathcal{L}}(C)}=\deg{\mathcal{L}}(D)$, where $C$ runs through the connected components of $X_{k'}$ above a fixed connected component $D$ of $X$. This proves immediately the ``converse'' statement.  

We can assume that $X$ is connected. Assume first that $k'$ is Galois over $k$, potentially of infinite degree. Then $\mrm{Gal}(k'/k)$ acts transitively on $\pi_0(X_{k'})$ (see for instance \cite[Exercise 3.2.10]{QL}), so that $\deg{f^{\ast}\mathcal{L}}(C)$ is an integer $\delta$ that does not depend on $C$. Therefore, if $n$ is the number of connected components of $X_{k'}$, $\deg{\mathcal{L}}=n\delta$, hence $\delta=0$. 

When $k'/k$ is algebraic and separable, we can embed it in its Galois closure $k''/k$. Then, on any connected component $C$ of $X_{k''}$, where $f':X_{k''} \rar X$ is the base change morphism, one has $\deg{(f')^{\ast}\mathcal{L}}(C)=0$. The ``converse'' statement applied to the extension $k''/k'$ and the line bundle $f^{\ast}\mathcal{L}$ shows that $f^{\ast}\mathcal{L}$ has degree zero on each connected component. 

When $k'/k$ is purely inseparable, $X_{k'} \rar X$ is a homeomorphism, hence $X_{k'}$ is connected and $C=X_{k'}$, and the statement is true. Pulling these results together, this proves the claim when $k'/k$ is algebraic. 

In general, let $k_1 \supset k'$ be an algebraically closed field and $k'' \subset k_1$ the algebraic closure of $k$ in $k_1$. To prove the claim, using the ``converse'' for the extension $k_1/k'$, we may assume that $k'=k_1$. By the above, we may assume that $k=k''$. In this case, the connected components of $X$ are geometrically connected and the result follows from the first identity.
}

\lem[functoriality-degrees]{Let $k$ be a field, $f: Y \rar X$ be a finite flat map of pure one-dimensional proper $k$-schemes.  
Let $\mathcal{M}$ be a line bundle on $Y$ of degree zero. Then $[f_{\ast}]\mathcal{M}$ has degree zero. 
Let $\mathcal{L}$ be a line bundle on $X$. Then, given connected components $C$ of $Y$, $C'$ of $X$ containing $f(C)$, $\deg(f^{\ast}\mathcal{L})(C') = (\deg{f_C}) \cdot (\deg{\mathcal{L}})(C)$, where $f_C: C \rar C'$ is the finite flat map induced by $f$.}

\demo{This is a well-known consequence of the results from \cite[Section 0AYQ]{Stacks}. Since our notion of degree is not quite the same, we include a proof for the sake of completeness. 

Let $\mathcal{M}$ be a line bundle on $Y$, and $C$ be any connected component of $X$. Then \[\deg{\det[f_{\ast}\mathcal{M}]}(C)=(\deg{f_{\ast}\mathcal{M}})(C) = \chi(C,f_{\ast}\mathcal{M})-\chi(C,\OO_C) = \chi(f^{-1}(C),\mathcal{M})-\chi(C,\OO_C).\] Hence, if $\mathcal{M}$ has degree zero, $\deg{\det[f_{\ast}\mathcal{M}]}(C)=\chi(f^{-1}(C),\OO_Y)-\chi(C,\OO_C)$, and the conclusion follows from Proposition \ref{deg-elem}.

Let $\mathcal{L}$ be a line bundle on $X$ of degree zero. Let $C \subset Y$ be a component (endowed with the reduced scheme structure) and $C' \subset X$ the unique irreducible component (endowed with the reduced scheme structure) containing $f(C)$. Then $f_C: C \rar C'$ is finite between integral proper $k$-schemes of dimension one. Moreover, $f$ is finite flat hence open, so $f(C)$ contains a nonempty open subset of $C'$, whence $f_C$ is surjective. Thus $\deg{(f^{\ast}\mathcal{L})_{|C}} = \deg{f_C}\deg{\mathcal{L}_{|C'}}$ by \cite[Lemma 0YAZ]{Stacks}. The conclusion follows from using \cite[Lemma 0AYW]{Stacks} on every connected component of $Y$.

Let $\mathcal{L}$ be a line bundle on $X$ of degree zero. Let $C \subset Y$ be a connected component and $C' \subset X$ the unique connected component containing $f(C)$. Then $f_C: C \rar C'$ is finite flat between proper connected $k$-schemes of dimension one, hence locally free of constant degree $c$. Then 
\begin{align*}
\deg{f^{\ast}\mathcal{L}}(C)&=\chi(C,f^{\ast}\mathcal{L})-\chi(C,\OO_C) = \chi(C',(f_C){\ast}(f_C)^{\ast}\mathcal{L})-\chi(C',(f_C)_{\ast}\OO_C)\\
&= \deg((f_C)_{\ast}(f_C)^{\ast}\mathcal{L})(C') + c\chi(C',\OO_{C'})-\chi(C',(f_C)_{\ast}\OO_C)\\
&= (\deg{\mathcal{L} \otimes (f_C)_{\ast}\OO_C})(C) - \deg{(f_C)_{\ast}\OO_C}=c\deg{\mathcal{L}}(C),
\end{align*}
where the last equality is \cite[Lemma 0AYX]{Stacks}.}

\lem[stein-factorization-base-change]{Let $\pi: X\rar S$ be a proper, finitely presented morphism with Stein factorization \cite[Theorem 03H2]{Stacks} $X \overset{f}{\rar} S' \rar S$. For any $S$-scheme $T$, if $X_T \rar T' \rar T$ denotes the Stein factorization of $X_T \rar T$, there is a canonical morphism $\sigma: T' \rar S' \times_S T$ such that the following diagram commutes:
\[
\begin{tikzcd}[ampersand replacement=\&]
X_T \arrow{r} \arrow{d}{=} \& T' \arrow{r}\arrow{d}{\sigma} T \& \arrow{d}{=}\\
X \times_S T \arrow{r} \& S' \times_S T \arrow{r}\& S \times_S T
\end{tikzcd}
\] If $T$ is flat over $S$, $\sigma$ is an isomorphism.}

\demo{This comes down to the definition and properties of the Stein factorization, as well as the existence of a base change morphism $(\pi_{\ast}\OO_X)_T \rar (\pi_T)_{\ast}\OO_{X_T}$, which is an isomorphism when $T$ is flat over $S$. }

\defi[picard-functors-definition]{Let $X$ be a proper finitely presented $S$-scheme. The \emph{relative Picard presheaf} of $X$ is the functor $P_{X/S}$ on $\mathbf{Sch}_S$ given by $T \longmapsto \operatorname{Pic}(X \times_S T)$. The \emph{relative Picard functor} of $X/S$ is the sheafification $\operatorname{Pic}_{X/S}$ of $P_{X/S}$ in the fppf topology. If moreover $X \rar S$ has pure fibres of dimension one, we define $P^0_{X/S}$ as the sub-presheaf of $P_{X/S}$ constituted by the classes of those line bundles $\mathcal{L}$ on $X \times_S T$ such that for every $t \in T$, the line bundle $t^{\ast}\mathcal{L}$ on $X_{\kappa(t)}$ has degree zero on every connected component. Its fppf-sheafification is the \emph{relative Picard functor of degree zero} $\operatorname{Pic}^0_{X/S}$ for $X/S$.
}

\rems{\begin{itemize}[noitemsep,label=$-$]
\item By Lemma \ref{deg-after-field-ext}, when the fibres of $X \rar S$ are pure of dimension one, $P^0_{X/S}$ is indeed a sub-presheaf of $P_{X/S}$ and the quotient pre-sheaf $P_{X/S}/P^0_{X/S}$ is separated for fppf covers.
\item The functors $P_{X/S},\operatorname{Pic}_{X/S}$ match the one described in \cite[\S 8.1, Definition 2]{BLR}. By Lemma \ref{deg-after-field-ext}, $\operatorname{Pic}^0_{X/S}$ identifies with the smallest subsheaf of $\operatorname{Pic}_{X/S}$ containing $P^0_{X/S}$. 
\end{itemize}
}

\lem[weakly-constant]{Let $X$ be a proper $k$-scheme, pure of dimension one, and $S$ be a finite $k$-scheme. Let $f: X \rar S$ be a morphism of $k$-schemes and $\mathcal{L}$ be a line bundle on $S$. Then $\mathcal{L}$ is trivial. In particular, $f^{\ast}\mathcal{L}$ has degree $0$.}

\demo{$S$ is a finite $k$-scheme, so it is affine and the spectrum of a finite product of Artinian local rings. Hence any line bundle on $S$ is trivial. In particular, $f^{\ast}\mathcal{L}$ is the trivial line bundle on $X$, thus its degree is zero on every connected component of $X$. }

\lem[kernel-relative-picard]{Let $X$ be a proper finitely presented $S$-scheme and $S'$ be any $S$-scheme. Let $f: X_{S'} \rar T$ be a morphism of $S'$-schemes, with $T$ finite locally free over $S'$. Then $f^{\ast}\operatorname{Pic}(T)$ is in the kernel of $P_{X/S}(S') \rar \operatorname{Pic}_{X/S}(S')$.}

\demo{Let $\mathcal{L}$ be a line bundle on $T$. By Lemma \ref{locally-free-base-change}, there is an affine open cover $U_i$ of $S'$ such that $\mathcal{L}{|T \times_{S'} U_i}$ is free. In particular, the image of $f^{\ast}\mathcal{L}$ in $P_{X/S}(U_i)$, hence in $\operatorname{Pic}_{X/S}(U_i)$, vanishes. Since $\operatorname{Pic}_{X/S}$ is a fppf-sheaf, the class of $\mathcal{L}$ in $\operatorname{Pic}_{X/S}(S')$ vanishes.}

\defi[relative-curve]{A relative curve is a map $\pi: C \rar S$ satisfying the following properties:
\begin{itemize}[noitemsep,label=$-$]
\item $\pi$ is proper, finitely presented, flat of relative dimension one with geometrically reduced fibres, and surjective,
\item $\pi_{\ast}\OO_C$ is a quasi-coherent finite \'etale $\OO_C$-algebra, and its formation commutes with arbitrary base change. 
\end{itemize}}

\prop[relative-curve-elem]{The property ``relative curve'' is stable under arbitrary base change and fpqc-local on $S$. Moreover, if $C \rar S$ is a relative curve, the formation of the Stein factorization $C \rar T \rar S$ \cite[Theorem 03H2]{Stacks} commutes with arbitrary base change. Moreover, $T \rar S$ is finite \'etale, $f: C \rar T$ is proper flat of relative dimension one, its geometric fibres are reduced and connected, and the formation of $f_{\ast}\OO_C$ commutes with arbitrary base change (with respect to $T$). }

\demo{The first statement is a formal verification using the fact that the formation of push-forwards of quasi-coherent sheaves commutes to flat base change and the fact that the formation of $\pi_{\ast}\OO_C$ commutes to arbitrary base change. 

The fact that $T \rar S$ is finite \'etale follows from the definition. Since $f$ is the composition $C=C \times_T T \rar C \times_S T \rar T$ of a closed open immersion with a proper flat morphism of relative dimension one, it is proper flat of relative dimension one. The fibres of $f$ are geometrically connected by \cite[Theorem 03H2]{Stacks}. This theorem also implies that the formation of the Stein factorization commutes to arbitrary base change (of $S$-schemes), since the formation of $\pi_{\ast}\OO_C$ commutes to arbitrary base change by definition. 

Let $U$ be a $T$-scheme (with structure map $\sigma$). Note that the following diagram commutes, both squares are Cartesian, and $j$ is a closed open immersion since $T \rar S$ is finite \'etale: 
\[
\begin{tikzcd}[ampersand replacement=\&]
C \times_T U\arrow{d}{f_{T,U}} \arrow{r} \& C \times_S U \arrow{r}\arrow{d}{f_U}\& C\arrow{d}{f}\\
U \arrow{r}{j=(\sigma,\mrm{id})}\& T \times_S U \arrow{r}\& T
\end{tikzcd}
\]

It follows in particular that the natural morphism $\sigma^{\ast}f_{\ast}\OO_C \rar (f_{U,T})^{\ast}\OO_{C \times_T U}$ is an isomorphism, that is, the formation of $f_{\ast}\OO_C$ commutes with arbitrary base change with respect to $T$. Since $\OO_T \rar f_{\ast}\OO_C$ is an isomorphism, this map remains an isomorphism after any base change of $T$-scheme.}

\prop[relative-picard-better-presheaf]{Let $X \rar S$ be a relative curve with Stein factorization $X \overset{f}{\rar} T \rar S$. The pre-sheaf $P'_{X/S}: U \longmapsto P_{X/S}(U)/f^{\ast}\operatorname{Pic}(T \times_S U)$ is separated for fppf covers, and the sheafification map $P_{X/S} \rar \operatorname{Pic}_{X/S}$ factors through $P'_{X/S}$. The pre-sheaf \[P^{'0}_{X/S}: U \longmapsto P^0_{X/S}(U)/f^{\ast}\operatorname{Pic}(T \times_S U)\] is a sub-presheaf of $P'_{X/S}$, and the sheafification map $P^0_{X/S} \rar \operatorname{Pic}^0_{X/S}$ factors through $P^{'0}_{X/S}$. Moreover, if there exists a section of $X \rar T$, then $P'_{X/S}, P^{'0}_{X/S}$ are sheaves in the fppf topology.}

\demo{That $P^{'0}_{X/S}$ is a sub-presheaf of $P'_{X/S}$ follows from Lemma \ref{weakly-constant}. That the sheafification $P_{X/S} \rar \operatorname{Pic}_{X/S}$ factors through $P'_{X/S}$ (and the same statement for the degree-zero pre-sheaves) follows from Lemma \ref{kernel-relative-picard}, since $T$ is finite \'etale locally free over $S$ and the formation of the Stein factorization commutes with arbitrary base change in this case, by Proposition \ref{relative-curve-elem}.

That $P'_{X/S}$ is separated is a reformulation of the seesaw theorem \cite[Lemma 0BDP]{Stacks} for (base changes of) the morphism $\pi: X \rar T$, using Proposition \ref{relative-curve-elem} (since $\OO_T \rar \pi_{\ast}\OO_X$ is a universal isomorphism). 

Thanks to Proposition \ref{relative-curve-elem}, to prove the final statement, it is enough to show that $P'_{X/S}(S)$ satisfies the glueing property. Let $\sigma: T \rar X$ be a morphism of $T$-schemes, let $I$ be a set and $(\pi_i: S_i \rar S)_{i \in I}$ be a fppf covering. Let, for each $i \in I$, $L_i\in P'_{X/S}(S_i)$ such that, for each $i,j \in I$, $L_i$ and $L_j$ have the same image in $P'_{X/S}(S_i \times_S S_j)$. We want to find some $s \in P'_{X/S}(S)$ such that each $s_i$ is the restriction of $S$. 

Let $i \in I$, pick some line bundle $\mathcal{L}$ on $X \times_S S_i$ in the class of $s_i$, and let $\mathcal{L}_i = \mathcal{L} \otimes (f^{\ast}\sigma^{\ast} \mathcal{L})^{\otimes -1}$. Then there is an isomorphism $\rho_i: \OO_{T \times_S S_i} \rar \sigma^{\ast}\mathcal{L}_i$, and $\mathcal{L}_i$ is still in the class of $s_i$. 

Let $i,j \in I$, and let $\mathcal{L}'_i$ (resp. $\mathcal{L}'_j$) be the pull-back of $\mathcal{L}_i$ (resp. $\mathcal{L}_j$) under the projection $X \times_S S_i \times_S S_j \rar X \times_S S_i$ (resp. $X \times_S S_i \times_S S_j \rar X \times_S S_j$): by assumption, there is a line bundle $\mathcal{M}$ on $T \times_S S_i \times_S S_j$ and an isomorphism $\iota_{i,j}: \mathcal{L}'_i \rar \mathcal{L}'_j \otimes f^{\ast}\mathcal{M}$. After pulling back by $\sigma \times_S S_i \times_S S_j$, we see that $\mathcal{M}$ is trivial. Since $\OO_T \rar f_{\ast}\OO_X$ is a universal isomorphism, we may then pick $\iota_{i,j}$ such that the following diagram commutes:

\[
\begin{tikzcd}[ampersand replacement=\&]
\OO_{T \times_S S_i \times_S S_j} \arrow{rr}{\mrm{id}} \arrow{d}{\rho_i} \&\& \OO_{T \times_S S_i \times_S S_j} \arrow{d}{\rho_j}\\
(\sigma \times S_i \times S_j)^{\ast}\mathcal{L}'_i \arrow{rr}{(\sigma^{\ast} \times_S S_i \times_S S_j)\iota_{i,j}}\&\& (\sigma \times S_i \times S_j)^{\ast}\mathcal{L}'_i 
\end{tikzcd}
\]

This rigidification ensures that $(\mathcal{L}_i,\iota_{i,j})$ satisfies a cocycle relation: that is, for every $i,j,k \in I$, if $\mathcal{L}'_i$ now denotes the pull-back of $\mathcal{L}_i$ under $X \times_S S_i \times_S S_j \times_S S_k \rar X \times_S S_i$, and similarly for $j,k$, then $(\iota_{j,k} \times_S S_i) \circ (\iota_{i,j} \times_S S_k) = (\iota_{i,k} \times S_j)$. Thus, this collection defines a line bundle on the small fppf site of $X$, which must therefore come from a (Zariski) line bundle on $X$ by descent \cite[Proposition 023T]{Stacks}. The desired class $s$ is exactly the class of this line bundle. }

\prop[relative-picard-functoriality]{Let $f: X \rar Y$ be a finite locally free morphism of proper finitely presented $S$-schemes. 
Then $f$ induces functorial homomorphisms of pre-sheaves of abelian groups $f^{\ast}: P_{Y/S} \rar P_{X/S}, [f_{\ast}]: P_{X/S} \rar P_{Y/S}$. If $f$ has constant rank $d$, one has $[f_{\ast}]f^{\ast}=d$. 
If $X,Y$ have pure fibres (over $S$) of dimension one, $f^{\ast},[f_{\ast}]$ induce morphisms between the $P^0$.
 }

\demo{The construction of the morphisms of pre-sheaves is a consequence of Proposition \ref{functoriality-line-bundles}, using Lemma \ref{functoriality-degrees} to check that both morphisms preserve the sub-presheaf $P^0$ when the fibres of $X \rar S, Y \rar S$ are pure of dimension one. If $f$ has constant rank $d$, it is enough to show that $[f_{\ast}]f^{\ast} \in \mrm{End}(\operatorname{Pic}(Y))$ is the multiplication by $d$. If $\mathcal{L}$ is a line bundle on $Y$, then, by the projection formula, 
\begin{align*}
[f_{\ast}]f^{\ast}\mathcal{L} &\simeq \det{f_{\ast}f^{\ast}\mathcal{L}} \otimes (\det{f_{\ast}\OO_X})^{-1} \simeq \det(\mathcal{L} \otimes f_{\ast}\OO_X) \otimes (\det{f_{\ast}\OO_X})^{-1}\\
& \simeq \mathcal{L}^{\otimes d} \otimes \det{f_{\ast}\OO_X} \otimes (\det{f_{\ast}\OO_X})^{-1} \simeq \mathcal{L}^{\otimes d}.\end{align*}
}

\lem[relative-picard-splitup]{Let $X$ be a proper finitely presented $S$-scheme, and $i_1,\ldots, i_r$ be closed open immersions $X_i \rar X$, such that $X$ is the disjoint reunion of the $i_j(X_j)$. 
Then \[(i_1^{\ast},\ldots,i_r^{\ast}): P_{X/S} \rar \bigoplus_{j=1}^r{P_{X_j/S}}\] is an isomorphism. Its inverse is the sum of the $[(i_j)_{\ast}]$.}

\demo{This is a formal verification.}

\prop[relative-picard-mixed-functoriality]{In the situation of Proposition \ref{mixed-functoriality}, assume that all maps are morphisms of proper finitely presented $S$-schemes. Then $u^{\ast}[(f')_{\ast}]$ and $[f_{\ast}]v^{\ast}$ both define morphisms of functors $P^{Y'/S} \rar P^{X/S}$ mapping $P^0_{Y'/S}$ to $P^0_{X/S}$. These morphisms are equal if there is an open subscheme $\iota: U \rar X'$ such that, for any $S$-scheme $T$, $\OO_{X_T} \rar (\iota_T)_{\ast}\OO_{U_T}$ is injective.}

\demo{This is a direct consequence of Proposition \ref{mixed-functoriality} using Proposition \ref{relative-picard-functoriality}. }

\lem[good-base-curve]{If $\pi: C \rar S$ is flat proper of relative dimension one with geometrically reduced fibres and surjective over a locally Noetherian scheme, then $\pi$ is a relative curve.}

\demo{By \cite[(7.8.6-7), (7.7.5.3)]{EGA-III2}, $\pi_{\ast}\OO_C$ is a coherent locally free $\OO_S$-algebra and its formation commutes with any base change. Let $S'$ be the relative spectrum of $\pi_{\ast}\OO_C$ over $S$, so that $C \rar S' \rar S$ is the Stein factorization of $\pi$ by \cite[Theorem 03H0]{Stacks}. 

Let us show that $S'$ is finite \'etale over $S$ (in particular reduced locally Noetherian). We already saw that $S'$ was finite flat over $S$. To prove that it is \'etale, since the formation of $\pi_{\ast}\OO_C$ commutes with arbitrary base change, we may assume that $S$ is the spectrum of an algebraically closed field $k$. Then $\OO(C)$ is a finite reduced $k$-algebra, hence is an unramified $k$-algebra, whence the conclusion. }

\lem[fppf-local-relative-curve]{Let $\pi: C \rar S$ be a relative curve. Then, fppf-locally over $S$, $C$ is the disjoint reunion of flat proper $S$-schemes $C_i$ with geometrically connected fibres and such that $C_i(S) \neq\emptyset$.}

\demo{Let $C \rar S' \rar S$ denote the Stein factorization; by working Zariski-locally over $S$, we may assume that $S' \rar S$ has constant rank $r$. Let $T$ be any flat finitely presented $S$-scheme. By Proposition \ref{relative-curve-elem}, $C_T \rar T$ is a disjoint union of flat proper $T$-schemes $C_{T,i}$ with geometrically connected fibres and such that $C_{T,i}(T)$ is nonempty if, and only if, $S' \times_S T$ is (as a $T$-scheme) the disjoint reunion of closed open subschemes of $T$ and $C\times_S T \rar S' \times_S T$ has a nonempty section. 

By \cite[Lemma 04HN]{Stacks}, there exists a fppf-covering family $T_i \rar S$ by affine schemes such that $S' \times_S T_i$ is, as a $T_i$-scheme, the disjoint reunion of finitely many copies of $T_i$. So, after replacing $S$ with any $T_i$, we may assume that $S$ is affine and $S'$ is (as a $S$-scheme) the disjoint reunion of finitely many copies of $S$, called $S_i$ for $1 \leq i \leq r$. Let $C_i = C \times_{S'} S_i$: it is a proper flat $S$-scheme with geometrically connected fibres. 

Now, for a $S$-scheme $T$, the morphism $C\times_S T \rar S' \times_S T$ has a section if and only if, for every $i$, the morphism $C_i \times_S T \rar T$ has a section. By \cite[Corollaire (17.16.2)]{EGA-IV4}, there is a faithfully flat quasi-finite finitely presented affine $S$-scheme $T_i$ such that $C_i\times_{S} T_i \rar T_i$ has a section. Let $T$ be the fibre product of the $T_i$ over $S$ (for $1 \leq i \leq r$): then $T$ is a faithfully flat quasi-finite finitely presented affine $S$-scheme and $C_i \times_S T \rar T$ has a section for all $i$, whence the conclusion.

}

\prop[jacobian-relative-curve-exists]{Let $\pi: C \rar S$ be a smooth relative curve and $C \overset{f}{\rar} T \rar S$ be its Stein factorization. The sheaf $\operatorname{Pic}^0_{C/S}$ is representable by an Abelian scheme over $S$, the \emph{relative Jacobian} of $J$ of $C/S$, whose construction commutes with any base change with respect to $S$.}

\demo{We will use the presheaves $P', P^{'0}$ introduced in Proposition \ref{relative-picard-better-presheaf}.

Since $f: C \rar T$ is proper smooth with geometrically connected fibres, the fppf sheaf $\operatorname{Pic}^0_{C/T}$ is, by \cite[\S 9.4, Proposition 4]{BLR}, representable by an Abelian $T$-scheme $J_0$. It is a direct verification that $P'_{C/S}$ is the push-forward under $T \rar S$ of $P'_{C/T}$, that is, that $P'_{C/S}(-) \simeq P'_{C/T}(- \times_S T)$. Moroever, it is clear using Lemma \ref{deg-after-field-ext} that this isomorphism preserves the sub-presheaf $P^{'0}$. 

Thus, there is formally a natural transformation $\Lambda: P'_{C/S}(-) \Rightarrow \operatorname{Pic}_{C/T}(- \times_S T)$ (preserving the degree-zero subfunctor). Since the Weil restriction commutes with any base change \cite[\S 7.6, p. 192]{BLR}, for any $S$-scheme $S'$ such that $C_{S'} \rar T \times_S S'$ has a section, $\Lambda$ is an isomorphism when restricted to the full subcategory of $\mathbf{Sch}_S$ whose objects are schemes that admit a $S$-morphism to $S'$. 

Let us prove that $S$ has a fppf cover by $S$-schemes $S_i$ such that $C \times_S S_i \rar T \times_S S_i$ has a section. By \cite[Lemma 04HN]{Stacks}, after replacing $S$ with a jointly surjective collection of \'etale maps, we can assume that $T$ is the disjoint union of copies of $S$, so that $C$ is the disjoint reunion of smooth relative curves $f_i: C_i \rar S$ with geometrically connected fibres. We are then done if, fppf locally on $S$, each $f_i$ admits a section. This is classical, see for instance \cite[Corollaire (17.16.3)]{EGA-IV4}. 

Since, by \cite[\S 7.6, Proposition 3]{BLR}, $\operatorname{Pic}_{C/T}(- \times_S T)$ is a fppf sheaf on $\mathbf{Sch}_S$, $\Lambda$ is thus a sheafification map. Hence, $\operatorname{Pic}^0_{C/S}$ is representable (by an Abelian $S$-scheme) iff $J_0(T \times_S -)$ is representable (by an Abelian $S$-scheme). 

For any $s \in S$, since $T$ is \'etale over $S$, $T_{\kappa(s)}$ is a finite discrete reunion of residue fields of points of $T$. Then $(J_0)_{\kappa(s)} = J_0 \times_T T_{\kappa(s)}$ is a finite disjoint reunion of Abelian varieties over some finite extensions of $\kappa(s)$. In particular, any finite subset of points in $(J_0)_{\kappa(s)}$ is contained in an affine open subset. Since $T \rar S$ is finite locally free, the Weil restriction is representable by Proposition \cite[\S 7.6, Theorem 4]{BLR} by some $S$-scheme $J$. 

Let $S'$ be any $S$-scheme, and let $S_i \rar S$ be the fppf cover discussed above such that each $C \times_S S_i \rar T \times_S S_i$ has a section. 

There is a natural morphism $\operatorname{Pic}^0_{C \times_S S'/S'} \Rightarrow \operatorname{Pic}^0_{C/S}$ of fppf sheaves on $\mathbf{Sch}_{S'}$ obtained by sheafifying the natural composition $P^0_{C \times_S S'/S'} \simeq P^0_{C/S}(\omega) \Rightarrow \operatorname{Pic}^0_{C/S}(\omega)$, where $\omega: \mathbf{Sch}_{S'} \rar \mathbf{Sch}_S$ is the forgetful functor. So we have a natural map $\beta_{S'/S}: J_{C\times_S S'/S'} \rar J_{C/S} \times_S S'$. 

Note that this map is an isomorphism by the properties of sheafification if $S'$ is flat and locally of finite presentation over $S$, or if $C(T)$ has a section (since in this case $P^0$ is a sheaf). It follows that every $\beta_{S'/S} \times_S S_i$ is an isomorphism, so $\beta_{S'/S}$ is an isomorphism. Therefore, the formation of $\operatorname{Pic}^0_{C/S}$ commutes with base change. 

Since $\operatorname{Pic}^0_{C/S}$ is a functor in Abelian groups, $J$ is naturally a commutative group scheme over $S$. By \cite[\S 7.6, Proposition 5]{BLR}, $J$ is proper smooth over $S$, thus, to prove that it is an Abelian scheme, all we need to do is prove that its fibres are connected. Since the construction commutes with base change, we may assume that $S$ is the spectrum of an algebraically closed field $k$, so that $T$ is a finite disjoint union of copies of $\Sp{k}$. Then $C$ is the disjoint reunion of smooth projective connected curves $C_i$ over $k$, and that as abelian fppf-sheaves, so, by Lemma \ref{relative-picard-splitup}, $P^0_{C/\Sp{k}}$ is isomorphic to the direct sum of the $P^0_{C_i/\Sp{k}}$. The associated sheaf is represented by the direct product of the (classically) well-defined Jacobians $J_i$ of the $C_i$, which are known to be Abelian varieties over $k$ (for instance \cite[\S 9.4, Theorem 4]{BLR}), hence their product is an Abelian variety over $k$ as well.  
}

\prop[difference-map-relative-curve]{Let $\pi: C \rar S$ be a smooth relative curve with Stein factorization $C \overset{f}{\rar} T \rar S$ and relative Jacobian $J$. There exists a difference map $\delta: C \times_T C \rar J$ such that the following diagrams commute: 
\[
\begin{tikzcd}[ampersand replacement=\&] 
C \arrow{r}{(\mrm{id},\mrm{id})} \arrow{d}{0}\& C \times_T C \arrow{dl}{\delta} \& C \times_T C \times_T C \arrow{d}{(p_1,p_3)} \arrow{rr}{((p_1,p_2),(p_2,p_3))} \&\& (C \times_T C)\times_S (C \times_T C) \arrow{d}{(\delta,\delta)} \& \\
J \& \& C \times_T C \arrow{r}{\delta} \& J \& J \times_S J \arrow{l}{\delta}.
\end{tikzcd}
\] 

and such that moreover, for any geometric points $x_1,x_2: S' \rar C$ with $f(x_1)=f(x_2)$, $\delta(x_1-x_2)$ is the class of the line bundle $\mathcal{L}=\mathcal{I}_{x_2} \otimes \mathcal{I}_{x_1}^{-1}$ where $\mathcal{I}_{x_i}$ is the sheaf of ideals attached to the relative effective Cartier divisor for $C_{S'} \rar S'$ defined by the section $x_i$ (see \cite[Definition 062T]{Stacks}).

In particular, if $S$ is an algebraically closed field, then $\mathcal{L}$ is the line bundle on $X$ which is trivial outside the connected component containing $x_1,x_2$ and attached to the Weil divisor $[x_1]-[x_2]$ on this component.
}

\demo{If the general definition is correct, the rest of the statements follows easily. So we need to show that the sections $x_i$ define relative effective Cartier divisors and that $\mathcal{L}$ has fibrewise degree zero. The first part follows from \cite[Lemma 1.2.2]{KM}, and the second part is \cite[Lemma 1.2.7]{KM} using the fact that degrees are additive.}

\prop[difference-map-generates-jacobian]{Let $\pi: C \rar S$ be a smooth relative curve with Stein factorization $C \overset{f}{\rar} T \rar S$ and relative Jacobian $J \rar S$. Assume that $C \rar S$ has a bounded above genus $g$ (for instance if $S$ is quasi-compact). Let $\delta: C \times_T C \rar J$ be the difference map. Then the product over $S$ of $g$ copies of the difference function $\delta^g: (C \times_T C)^g \rar J$ is surjective.
}

\demo{By Lemma \ref{fppf-local-relative-curve}, after replacing $S$ with a fppf-covering, we may assume that $C$ is a disjoint reunion of smooth proper $S$-schemes $C_i$ of relative dimension one such that $C_i \rar S$ has a section and the fibres of $C_i \rar S$ have genus at most $g$. In particular, the relative Jacobian $J_i$ of each $C_i/S$ is well-defined, and all we need to show is that for every $i$, $\delta^g: (C_i \times_S C_i)^g \rar J_i$ is surjective. 

This can be verified on the fibres of these two $S$-schemes, so we may assume that $S$ is the spectrum of an algebraically closed field and that $C$ is connected. By \cite[\S 9.3, Lemma 6]{BLR}, the image of $\delta^g$ contains an open subscheme of $J$. Since $\delta^g$ is a morphism between two proper $S$-schemes, it has closed image, hence it is surjective.  }

\cor[pushforward-and-difference]{In the situation of Proposition \ref{relative-picard-functoriality}, assume that $X,Y$ are smooth and that the push-forward by $f$ is well-defined. Then, if $J_X, J_Y$ denote the Jacobians and $\delta_X,\delta_Y$ the difference maps, the following diagram commutes:
 \[
\begin{tikzcd}[ampersand replacement=\&]
X \times_{T_X} X \arrow{r}{f \times f}\arrow{d}{\delta_X} \& Y \times_{T_Y} Y \arrow{d}{\delta_Y}\\
J_X \arrow{r}{[f_{\ast}]} \& J_Y
\end{tikzcd}\]}

\demo{This amounts to proving the following variant of Lemma \ref{pushforward-is-weil}. Let $f: X \rar Y$ be a finite locally free morphism of relative curves over $S$ and let $s \in X(S)$. Then $[f_{\ast}]$ maps the ideal attached to the relative effective Cartier divisor $s$ to the ideal attached to the relative effective Cartier divisor defined by $f(s) \in Y(S)$. 

This statement can in turn easily be reduced to the following situation, using Lemma \ref{locally-free-finite-sheaf}: let $R$ be a ring and $B$ be a finite locally free $A$-algebra, where $A,B$ are smooth of relative dimension one over $R$. Let $I \subset B$ be an ideal such that $R \rar B/I$ is an isomorphism, then $I$ is locally principal, the kernel of $A \rar B/I$ is a locally principal ideal $J$ generated by the $N_{B/A}(x)$ for $x \in I$, and $R \rar A/J$ is an isomorphism. 

The ideals $I,J$ are locally principal and generated by regular elements due to \cite[Lemma 1.2.2]{KM}, and the fact that $R \rar A/J$ is an isomorphism is formal; in particular, $A/J \rar B/I$ is an isomorphism. So the only remaining statement is the one about the norm. We can work Zariski-locally on $A$, and thus assume that $B$ is free over $A$ of finite rank, $I$ is generated by some regular element $g \in B$, and $J$ by some regular element $h \in A$. Since $A/J \rar B/I$ is still an isomorphism, its zero-th Fitting ideal is both $hA=J$ and $N_{B/A}(g)A$, whence the conclusion is \cite[Lemma 07Z8]{Stacks}. 
}

\prop[automorphism-global-functions]{Let $X \rar S$ be a relative curve with geometrically connected fibres and $\sigma$ be an $S$-automorphism of some finite \'etale $S$-scheme $T$. Let $J_1$ be the relative Jacobian of $X \rar S$, and $J$ be the relative Jacobian of $X_T \rar S$. For any $S$-scheme $A$, the following diagram commutes:
\[
\begin{tikzcd}[ampersand replacement=\&]
J(A) \arrow{r}{\sim} \arrow{d}{\alpha} \& J_1(A \times_S T)\arrow{d}{\beta}\\
J(A) \arrow{r}{\sim} \& J_1(A \times_S T)
\end{tikzcd}
\]

In this diagram, the horizontal arrows are given by the fact that $J$ is the restriction of scalars of $J_1 \times_S T$ over $T/S$ (by Proposition \ref{jacobian-relative-curve-exists}), $\alpha$ is the pushforward of $\mrm{id} \times \sigma: X \times_S T \rar X \times_S T$, and $\beta: f \longmapsto f \circ (\mrm{id}\circ \sigma^{-1})$. }

\demo{This statement is fppf-local with respect to $A$, so we may assume by \cite[Corollaire (17.16.2)]{EGA-IV4} that $X \rar S$ has a section. Therefore, in the following diagram, the diagonal maps are surjective by Proposition \ref{relative-picard-better-presheaf}:
\[
\begin{tikzcd}[ampersand replacement=\&]
\operatorname{Pic}^0(X_T \times_S A) \arrow{rrr}{=}\arrow{dr}\arrow{ddd}{[(\sigma_{X \times_S A})_{\ast}]} \& \& \& \operatorname{Pic}^0(X \times_S A \times_S T) \arrow{dl}\arrow{ddd}{(\sigma^{-1}_{X \times A})^{\ast}}\\
\& J(A) \arrow{r}{\sim} \arrow{d}{\alpha} \& J_1(A \times_S T)\arrow{d}{\beta} \&\\
\& J(A) \arrow{r}{\sim} \& J_1(A \times_S T) \&\\ 
\operatorname{Pic}^0(X_T \times_S A) \arrow{rrr}{=}\arrow{ur} \& \& \& \operatorname{Pic}^0(X \times_S A \times_S T) \arrow{ul}
\end{tikzcd}
\]

By definition of $J$ and $J_1$, the leftmost and rightmost cells commute. The top and bottom cells commute by the proof of Proposition \ref{jacobian-relative-curve-exists}. We can check that the outer cell commutes. Therefore the inner cell commutes.  

}

\subsection{Holomorphic differentials on the Jacobian of a curve}
\label{differentials-relative-jacobians}

\prop[differentials-are-trivial-abelian]{Let $f: G \rar S$ be an Abelian scheme over a locally Noetherian scheme $S$. Then $f_{\ast}\Omega^1_{G/S}$ is a locally free $\OO_S$-module, its sections are translation-invariant differential forms on $G_U/U$ (for open subschemes $U \subset S$) and its formation commutes with base change. Moreover, $\mrm{End}_S(G)$ acts on the right on $f_{\ast}\Omega^1_{G/S}$. }

\demo{This is essentially the proof of \cite[\S 4.2 Proposition 1]{BLR}, but the argument is easier in this case. 

Because $f$ is smooth, $\Omega^1_{G/S}$ is locally free of rank $d$. Let $e: S \rar G$ be the unit section. 

Let us prove that for any locally free $\OO_S$-module $\mathcal{F}$ of finite rank, the morphism $\mathcal{F} \rar f_{\ast}f^{\ast}\mathcal{F}$ is an isomorphism. The statement is Zariski-local over $S$, it is enough to prove that $\OO_S \rar f_{\ast}\OO_G$ is an isomorphism. Because $f$ is proper flat with geometrically reduced fibres, by \cite[(7.8.6-7), (7.7.5.3)]{EGA-III2}, $f_{\ast}\OO_G$ is a coherent locally free $\OO_S$-algebra and its formation commutes with any base change. Since $f$ has geometrically connected fibres, $\OO_S \rar f_{\ast}\OO_G$ becomes an isomorphism after any base change to a geometric fibre. Hence $\OO_S \rar f_{\ast}\OO_G$ is an isomorphism.

Let $m, p_1, p_2: G \times_S G \rar S$ denote the product, the first and second projections respectively. Let us prove that, for any open subscheme $U \subset S$ and any $\omega \in \Omega^1_{G/S}(G_U)$, one has $m^{\ast}\omega=p_1^{\ast}\omega+p_2^{\ast}\omega$. We can and do assume $S=U$. Moreover, we also know by \cite[\S 2.1 Proposition 4]{BLR} that $H^0(G \times_S G,\Omega^1_{G \times_S G/S}) = \bigoplus_{i=1}^2{H^0(G \times_S G,p_i^{\ast}\Omega^1_{G/S})}$. 

Now, $H^0(G \times_S G,p_i^{\ast}\Omega^1_{G/S}) \simeq H^0(G,(p_i)_{\ast}p_i^{\ast}\Omega^1_{G/S})$; since $p_i: G \times_S G \rar G$ is an Abelian scheme, one has $(p_i)_{\ast}p_i^{\ast}\Omega^1_{G/S} \simeq \Omega^1_{G/S}$, hence $H^0(G \times_S G, \Omega^1_{G \times_S G/S}) \simeq \bigoplus_{i=1}^2{p_i^{\ast}H^0(G,\Omega^1_{G/S})}$. 

In particular, we can write $m^{\ast}\omega = p_1^{\ast}\omega_1 + p_2^{\ast}\omega_2$ for $\omega_1,\omega_2 \in H^0(G,\Omega^1_{G/S})$. By pulling back by $(\mrm{id},0)$ (resp. $(0,\mrm{id})$) one sees that $\omega_1=\omega_2=\omega$. This implies that $\mrm{End}_S(G)$ acts on the right on $f_{\ast}\Omega^1_{G/S}$, and that any global differential on $G/S$ is translation-invariant. 

Now, it remains to prove that $f_{\ast}\Omega^1_{G/S}$ is a locally free $\OO_S$-module and that its formation commutes with base change. Indeed, by \cite[Proposition 1.1 (b)]{LLR}, there is a natural isomorphism $f^{\ast}e^{\ast}\Omega^1_{G/S} \rar \Omega^1_{G/S}$. Applying $f_{\ast}$, one gets a natural isomorphism $e^{\ast}\Omega^1_{G/S} \rar f_{\ast}\Omega^1_{G/S}$, hence $f_{\ast}\Omega^1_{G/S}$ is locally free of rank $\dim(A/S)$. Since the formation of $e^{\ast}\Omega^1_{G/S}$ commutes with base change, the same holds for $f_{\ast}\Omega^1_{G/S}$, whence the conclusion. }

\prop[differentials-are-free]{Let $f: X \rar S$ be a smooth relative curve over a Noetherian scheme $S$ and $\mathcal{F}$ be a locally free coherent sheaf over $X$. For every $j \geq 0$, $R^jf_{\ast}\mathcal{F}$ is a locally free coherent $\OO_S$-module whose formation commutes with an arbitrary base change of locally Noetherian schemes. Moreover, for $j \in \{0,1\}$, there is a natural isomorphism \[R^{1-j}f_{\ast}(\mathcal{F}^{\vee} \otimes \Omega^1_{X/S}) \rar \underline{Hom}_{\OO_Y}(R^jf_{\ast}\mathcal{F},\OO_X).\]}

\demo{First, let $\mathcal{F}$ be a locally free coherent sheaf on $X$. By \cite[Lemma 3]{Kleiman}, one has $R^if_{\ast}\mathcal{F}=0$ for every $i \geq 2$ (and this holds for every base change of $f$ to another locally Noetherian scheme). Apply \cite[Lemma 5.1.1]{Conrad-Duality} with $m=1$, this implies that the formation of $R^1f_{\ast}\mathcal{F}$ commutes with base change of locally Noetherian schemes. Moreover, apply \cite[Theorem 5.1.2]{Conrad-Duality} with $m=1$: the formation of $f_{\ast}(\mathcal{F}^{\vee}\otimes \Omega^1_{X/S})$ commutes with base change (since the dualizing sheaf for a smooth relative curve identifies with the sheaf of differentials, see the discussion after \cite[Corollary 3.5.2]{Conrad-Duality}). Since $\mathcal{F}$ was assumed to be locally free and coherent, this means that for every locally free coherent sheaf $\mathcal{F}$ on $X$, and every $i \geq 0$ (hence for every $i \in \Z$), the formation of $R^if_{\ast}\mathcal{F}$ commutes with base change. By \cite[Lemma 5.1.1]{Conrad-Duality}, this implies that $R^if_{\ast}\mathcal{F}$ is a locally free coherent $\OO_S$-module for every $i \in \Z$. The conclusion then follows from re-applying \cite[Theorem 5.1.2]{Conrad-Duality}. 
}

\prop[differentials-on-jacobian-curve]{Let $f: X \rar S$ be a smooth relative curve over a locally Noetherian scheme $S$ with relative Jacobian $F: J \rar S$. There is an isomorphism $\rho_X: F_{\ast}\Omega^1_{J/S}\rar f_{\ast}\Omega^1_{X/S}$ whose formation commutes with base change and satisfies the following property: let $X \rar S' \rar S$ be the Stein factorization, $T$ be a Noetherian $S$-scheme and $x: T \rar X$ be a morphism of $S$-schemes, thus making $T$ a $S'$-scheme. Then $F_{\ast}\Omega^1_{J_T/T} \overset{\rho_X\times_S T}{\rar} f_{\ast}\Omega^1_{X_T/T} \rar f_{\ast}\Omega^1_{(X \times_{S'}T )/T}$ is exactly given by the restriction through the ``difference to $x$'' map $\Delta(\cdot,x): X \times_{S'} T \rar J \times_S T$. In particular, the formation of $\rho_X$ commutes with base change. 

Moreover, if $g: Y \rar S$ is a smooth relative curve and with Jacobian $G: J' \rar S$ and $h: X \rar Y$ is finite locally free, then the following diagram commutes: 
\[
\begin{tikzcd}[ampersand replacement=\&]
G_{\ast}\Omega^1_{J'/S} \arrow{d}{[h_{\ast}]^{\ast}} \arrow{r}{\rho_Y} \& g_{\ast}\Omega^1_{Y/S} \arrow{d}{h^{\ast}}\\
F_{\ast}\Omega^1_{J/S}\arrow{r}{\rho_X} \& f_{\ast}\Omega^1_{X/S}
\end{tikzcd}
\]
}

\demo{It is enough to prove the result when $S$ is Noetherian. 

Let $(S_i)_i$ be an fppf-covering of $S$. If either part of the claim holds after making a base change to each $S_i$ and each $S_i \times_S S_j$, then the claim holds over $S$, by descent. 

Let us first assume that $f$ has geometrically connected fibres. Then, the first part of the claim is the content of \cite[(B.4.5)]{Conrad-Duality} and its discussion, because all the global differentials on (base changes of) $J/S$ are translation-invariant by \cite[\S 4.2, Corollary 3]{BLR}. The second part of the claim (when $g$ has geometrically connected fibres) clearly holds when $X(S)$ is not empty by definition of $\rho_X,\rho_Y$ and Corollary \ref{pushforward-and-difference}. 

Therefore, the first part of the claim (resp. the second part of the claim) remains true whenever $X$ is (resp. $X,Y$ are) the reunion of several smooth proper $S$-schemes $X_i$ (resp. $X_i,Y_i$) of relative dimension one with geometrically connected fibres such that $X_i(S) \neq \emptyset$ for every $i$ (resp. $X_i(S) \neq \emptyset, Y_i(S) \neq \emptyset$ for every $i$). 

The general case reduces to this special case after an fppf base change by Lemma \ref{fppf-local-relative-curve}.

}

\cor[differentials-curve-end-jac]{Let $f: X \rar S$ be a smooth relative curve over a locally Noetherian scheme with relative Jacobian $J$. Then $f_{\ast}\Omega^1_{X/S}$ is endowed with a natural structure of sheaf of right $\mrm{End}_S(J)$-modules.  }

\demo{This is direct.}

\defi[differentials-on-jacobian-trace]{Let $f: X \rar S$, $g: Y \rar S$ be smooth relative curves over a locally Noetherian base scheme $S$ with relative Jacobians $F: J_X \rar S, G: J_Y \rar S$, and let $h: X \rar Y$ be a finite locally morphism of $S$-schemes. The trace map attached to $h$ is the natural map $\mrm{Tr}_h: f_{\ast}\Omega^1_{X/S} \rar g_{\ast}\Omega^1_{Y/S}$ such that the following diagram commutes.   
\[\begin{tikzcd}[ampersand replacement=\&]
F_{\ast}\Omega^1_{J_X/S} \arrow{d}{(h^{\ast})^{\ast}} \arrow{r}{\rho_X} \& f_{\ast}\Omega^1_{X/S} \arrow{d}{\mrm{Tr}_h}\\
G_{\ast}\Omega^1_{J_Y/S} \arrow{r}{\rho_Y} \& g_{\ast}\Omega^1_{Y/S}.
\end{tikzcd}\]
}

\prop[trace-map-basic]{In the notation of Definition \ref{differentials-on-jacobian-trace}, the formation of $\mrm{Tr}_h$ commutes with base change and it is functorial in $h$. Moreover, this trace map is the image on global sections of the usual trace map $h_{\ast}\Omega^1_{X/S} \rar \Omega^1_{Y/S}$ coming from duality theory \cite[(2.7.36)]{Conrad-Duality} given on an affine open $U \subset Y$ by $a\cdot d(h^{\ast}(b)) \in \Omega^1_{X/S}(h^{-1}(U)) \longmapsto \mrm{Tr}_{\OO_X(h^{-1}(U))/\OO_Y(U)}(a) \cdot db \in \Omega^1_{Y/S}$. }

\demo{The first part of the claim is direct. Let $\mrm{Tr}^{\delta}: f_{\ast}\Omega^1_{X/S} \rar g_{\ast}\Omega^1_{Y/S}$ denote the map coming from Definition \ref{differentials-on-jacobian-trace} and $\mrm{Tr}^{\theta}: h_{\ast}\Omega^1_{X/S} \rar \Omega^1_{Y/S}$ the map coming from duality. The formula we have written is correct by the explicit computation in \cite[(2.7.41)]{Conrad-Duality}. 

We want to show that $\mrm{Tr}^{\delta} = g_{\ast}\mrm{Tr}^{\theta}$. The formation of both $\mrm{Tr}^{\delta}$ and $g_{\ast}\mrm{Tr}^{\theta}$ commutes with base change with respect to $S$ and is fppf-local with respect to $S$, so, by Lemma \ref{fppf-local-relative-curve}, so we may assume that $S$ is affine, connected, and $X,Y$ are reunions of smooth proper $S$-schemes of relative dimension one with geometrically connected fibres and nontrivial sections.  

It is enough to check the identity on every pair $h: X' \rar Y'$ separately, where $X' \rar X, Y' \rar Y$ are closed open subschemes such that the geometric fibres of $X' \rar S$, $Y' \rar S$ are connected, so we may and do assume that $X \rar S, Y \rar S$ have geometrically connected fibres. In particular, $\OO_S \rar g_{\ast}\OO_Y \rar f_{\ast}\OO_X$ are isomorphisms, thus, by \cite[(1.1)]{Kleiman}, $R^1f_{\ast}\Omega^1_{X/S}, R^1g_{\ast}\Omega^1_{Y/S}$ are free $\OO_S$-modules of rank one. 

Write $h^{\sharp}$ for the morphism $\OO_Y \rar h_{\ast}\OO_X$ induced by $h$, and let $R^1(h)$ denote the composition 
\[R^1g_{\ast}\OO_Y \overset{(R^1g_{\ast})(h^{\sharp})}{\rar} (R^1g_{\ast})(h_{\ast}\OO_X) \overset{\sim}{\rar} R^1(g_{\ast}\circ h_{\ast})\OO_X = R^1f_{\ast}\OO_X,\]  
where the second map is an isomorphism because $h$ is finite. 

Take the dual for Grothendieck duality of the diagram in \cite[Proposition 1.3(c)]{LLR}: by \cite[Theorem B.4.1]{Conrad-Duality}, this diagram is of the form 
\[
\begin{tikzcd}[ampersand replacement=\&]
F_{\ast}\Omega^1_{J_X/S} \arrow{d}{(h^{\ast})^{\ast}} \arrow{r}{\rho_X} \& f_{\ast}\Omega^1_{X/S} \arrow{d}{(R^1(h))^{\vee}}\\
G_{\ast}\Omega^1_{J_Y/S} \arrow{r}{\rho_Y} \& g_{\ast}\Omega^1_{Y/S}
\end{tikzcd}
\]
In particular, $\mrm{Tr}^{\delta}$ is equal to $(R^1(h))^{\vee}$. We want to show that $(R^1(h))^{\vee}$ is equal to $\mrm{Tr}^{\theta}$. 
The adjoint for Grothendieck duality of $(R^1g_{\ast})(h^{\sharp})$ is the morphism $g_{\ast}\underline{\mrm{Hom}_Y}(h_{\ast}\OO_X,\Omega^1_{Y/S}) \rar g_{\ast}\Omega^1_{Y/S}$ of evaluation at one. To prove the claim, it is thus enough to show that there exists an isomorphism $\iota: f_{\ast}\Omega^1_{X/S} \rar g_{\ast}\underline{\mrm{Hom}}_Y(h_{\ast}\OO_X,\Omega^1_{Y/S})$ such that 
$\mrm{ev}_{1 \in (h_{\ast}\OO_X)(Y)} \circ \iota = \mrm{Tr}^{\theta}$ and the following diagram commutes:
\[
\begin{tikzcd}[ampersand replacement=\&]
(R^1g_{\ast})(h_{\ast}\OO_X) \arrow{r}\arrow{d}\& \mrm{Hom}(g_{\ast}\underline{\mrm{Hom}}_Y(h_{\ast}\OO_X,\Omega^1_{Y/S}),\OO_S)\arrow{d}{\iota^{\ast}}\\
R^1f_{\ast}\OO_X \arrow{r}\& \mrm{Hom}(f_{\ast}\Omega^1_{X/S},\OO_S).
\end{tikzcd}
\]

Consider the following commutative diagram, as in the proof of \cite[Theorem 15]{Kleiman}\footnote{The analogue of $f^!, g^!$ in \cite{Conrad-Duality} is given by $f^{\sharp}, g^{\sharp}$; the difference is that \emph{loc.cit.} works in the derived category, shifts the complexes and possibly uses different sign conventions. All of these changes cancel out in the final result, since we are working with a finite map of smooth relative curves.}:
\[
\begin{tikzcd}[ampersand replacement=\&]
g_{\ast}\underline{\mrm{Hom}}_Y(h_{\ast}\OO_X,\Omega^1_{Y/S}) \arrow{d}{\sim} \arrow{r}{\cup}\& \underline{\mrm{Hom}}_S((R^1g_{\ast})(h_{\ast}\OO_X),R^1g_{\ast}\Omega^1_{Y/S}) \arrow{d}{\sim}\\
g_{\ast}\underline{\mrm{Hom}}_Y(h_{\ast}\OO_X,g^!\OO_S) \arrow{d}{\sim} \arrow{r}{\cup}\& \underline{\mrm{Hom}}_S((R^1g_{\ast})(h_{\ast}\OO_X),R^1g_{\ast}g^!\OO_S)\arrow{d}{\mrm{Tr}_g}\\
\underline{\mrm{Hom}}_S((R^1g_{\ast})(h_{\ast}\OO_X),\OO_S) \arrow{d}{\sim}\arrow{r}{\cup}\& \underline{\mrm{Hom}}_S((R^1g_{\ast})(h_{\ast}\OO_X),\OO_S)\arrow{d}{\sim}\\
\underline{\mrm{Hom}}_S(R^1f_{\ast}\OO_X,\OO_S) \arrow{d}{\sim}\arrow{r}{\cup}\& \underline{\mrm{Hom}}_S(R^1f_{\ast}\OO_X,\OO_S)\\
f_{\ast}f^!\OO_S \arrow{d}{\sim} \arrow{r}{\cup}\& \underline{\mrm{Hom}}_S(R^1f_{\ast}\OO_X,R^1f_{\ast}f^!\OO_S)\arrow{u}{\mrm{Tr}_f}\\
f_{\ast}\Omega^1_{X/S} \arrow{r}{\cup} \& \underline{\mrm{Hom}}_S(R^1f_{\ast}\OO_X,R^1f_{\ast}\Omega^1_{X/S}) \arrow{u}{\sim}
\end{tikzcd}
\]

Its first column defines an isomorphism $\iota: f_{\ast}\Omega^1_{X/S} \rar g_{\ast}\underline{\mrm{Hom}}_Y(h_{\ast}\OO_X,\Omega^1_{Y/S})$ which is indeed adjoint to the isomorphism $(R^1g_{\ast})(h_{\ast}\OO_X) \simeq R^1f_{\ast}\OO_X$. All we need to do is check that the evaluation of $\iota$ at $1 \in (h_{\ast}\OO_X)(Y)$ is exactly the trace homomorphism.

By construction, $\mrm{Tr}^{\delta}: h_{\ast}\Omega^1_{X/S} \rar \Omega^1_{Y/S}$ is the unique morphism such that \[(R^1g_{\ast})\mrm{Tr}^{\delta}: (R^1g_{\ast})(h_{\ast}\Omega^1_{X/S}) \rar R^1g_{\ast}\Omega^1_{Y/S}\] induces the identity map on $\OO_S$ after identifying \[(R^1g_{\ast})(h_{\ast}\Omega^1_{X/S}) \simeq R^1f_{\ast}\Omega^1_{X/S} \simeq \OO_S, R^1g_{\ast}\Omega^1_{Y/S} \simeq \OO_S.\]

Now, for any quasi-coherent sheaf $\mathcal{G}$ on $X$, consider as in the proof of \cite[Theorem 15]{Kleiman} the following chain of isomorphisms
\begin{align*}
\mrm{Hom}_X(\mathcal{G},\Omega^1_{X/S}) &\simeq \mrm{Hom}_X(\mathcal{G},f^!\OO_S) \simeq \mrm{Hom}_S(R^1f_{\ast}\mathcal{G},\OO_S) \simeq \mrm{Hom}_S((R^1g_{\ast})(h_{\ast}\mathcal{G}),\OO_S)\\
&\simeq \mrm{Hom}_Y(h_{\ast}\mathcal{G},g^!\OO_S) \simeq \mrm{Hom}_Y(h_{\ast}\mathcal{G},\Omega^1_{Y/S}). 
\end{align*}
Moreover, there is a functorial isomorphism \[\mrm{Hom}_{h_{\ast}\OO_X}(h_{\ast}\mathcal{G},\underline{\mrm{Hom}}_{\OO_Y}(h_{\ast}\OO_X,\Omega^1_{Y/S})) \rar \mrm{Hom}_{\OO_Y}(h_{\ast}\mathcal{G},\Omega^1_{Y/S})\] given by the evaluation at $1$. Therefore, if $\tau: h_{\ast}\Omega^1_{X/S} \rar \Omega^1_{Y/S}$ denotes the image of $\mrm{id}_{\Omega^1_{X/S}}$ (for $\mathcal{G}=\Omega^1_{X/S}$), then $g_{\ast}\tau$ is equal to the composition of $\iota$ and the evaluation at $1$. It is then easy to check that $R^1g_{\ast}\tau = R^1g_{\ast}\mrm{Tr}^{\delta}$, so $\tau = \mrm{Tr}^{\delta}$, whence the conclusion. 
}

\subsection{Jacobians and Galois twists}

\prop[jacobian-twist]{Let $X$ be a proper smooth $k$-scheme of relative dimension one, with Jacobian variety $J$, where $k$ is a field with separable closure $k_s$. Let $\rho: \mrm{Gal}(k_s/k) \rar \mrm{Aut}_{k_s}(X_{k_s})$ be a function satisfying the cocycle relation of Proposition \ref{cocycle-twist}. Let $\mrm{Aut}_{k_s}(X_{k_s})$ act on $J$ by push-forward functoriality. Let $\Delta: X \times_{\Sp{\OO(X)}} X \rar J$ be the difference morphism, which is equivariant for the action of $\mrm{Aut}_{k_s}(X_{k_s})$. Then its twist $\Delta_{\rho}: X_{\rho} \times_{\Sp{\OO(X_{\rho})}} X_{\rho} \rar J_{\rho}$ by $\rho$ identifies $J_{\rho}$ with the Jacobian variety of $X_{\rho}$. }

\demo{Let $j_X: X_{\rho} \times_k \Sp{k_s} \rar X \times_k \Sp{k_s}$ be the $k_s$-isomorphism from Proposition \ref{cocycle-twist}, $J_{\rho}$ be the Jacobian variety of $X_{\rho}$ over $k$. The formation of $J_{\rho}$ commutes with base change, so push-forward functoriality yields an isomorphism $j_J: J_{\rho} \times_k \Sp{k_s} \rar J \times_k \Sp{k_s}$ of Abelian varieties over $k_s$. 

The following diagram commutes for each $\sigma \in \mrm{Gal}(k_s/k)$ by Proposition \ref{relative-picard-mixed-functoriality} and the fact that the formation of the Jacobian commutes with base change:
\[
\begin{tikzcd}[ampersand replacement=\&]
(J_{\rho})_{k_s} \arrow{r}{j_J} \arrow{d}{(\mrm{id},\underline{\sigma^{-1}})} \& J_{k_s}\arrow{d}{[\rho(\sigma)_{\ast}]\circ (\mrm{id},\underline{\sigma^{-1}})}\\
(J_{\rho})_{k_s} \arrow{r}{j_J} \& J_{k_s}
\end{tikzcd}
\]

Therefore, by Proposition \ref{cocycle-twist}, $(J_{\rho},j_J)$ is the twist of $J$ by $\rho$. The rest of the statement follows from Proposition \ref{twist-equiv} and Corollary \ref{pushforward-and-difference}, because the twist of $\Sp{\OO(X)}$ by $\rho$ is exactly $\Sp{\OO(X_{\rho})}$.}

\prop[jacobian-twist-functoriality]{Let $f: X \rar Y$ be a finite locally free morphism of proper smooth $k$-schemes of relative dimension one, $k$ being a field with separable closure $k_s$. Let $J_X,J_Y$ be the Jacobians of $X,Y$. Let $G$ be a finite group acting on the $k_s$-schemes $X_{k_s},Y_{k_s}$ equivariantly for $f$, and let $\rho: \mrm{Gal}(k_s/k) \rar G$ be a function satisfying the cocycle property when post-composed with the actions of $G$ on $X_{k_s}$ and $Y_{k_s}$ (as in Proposition \ref{twist-equiv}). Then the twist by $\rho$ of $f^{\ast}: J_Y \rar J_X$ (resp. the push-forward $[f_{\ast}]: J_X \rar J_Y$) identifies with the pull-back $(f_{\rho})^{\ast}: J_{Y_{\rho}} \rar J_{X_{\rho}}$ (resp. the push-forward $[(f_{\rho})_{\ast}]$). }

\newpage

%% file: reps-0.tex
\chapter{Representations of $\GL{\F_p}$.}
\label{reps-gl2}

The purpose of this Appendix is to recall the representation theory of $G=\GL{\F_p}$ when $p \geq 5$. Most of the following is well-known, but is worth stating in full for its applications to the main work. We follow mainly \cite{Jorza-gl2fq}, \cite{Berger-caracp} and \cite[\S 6]{BH}. 

\nott{\begin{itemize}[noitemsep,label=\tiny$\bullet$]
\item $B \leq \GL{\F_p}$ is the subgroup of upper-triangular matrices, 
\item $U=\begin{pmatrix}1 & 1\\0 &1\end{pmatrix}$, $W=\begin{pmatrix}0&-1\\1 & 0\end{pmatrix}$, $T=\begin{pmatrix}1 & 0\\1 & 1\end{pmatrix}$,
\item $S' = \{W\} \cup \{T^t,\, t \in \F_p\}$ is a natural system of representatives for $G/B$, 
\item given $a,b \in \F_p^{\times}$, we write $\Delta_{a,b}=\begin{pmatrix}a & 0\\0 & b\end{pmatrix}$. 
\end{itemize}}

The following relations hold:

\begin{align*}
a,b,t \in \F_p^{\times}:&\quad \Delta_{a,b}T^t = T^{tb/a}\Delta_{a,b}\,; &&\quad  \Delta_{a,b}W = W\Delta_{b,a}\\
s,t \in \F_p, u:=st+1 \neq 0:& \quad U^sT^t = T^{t/(1+st)}\begin{pmatrix}u & s\\0 & u^{-1}\end{pmatrix} \\
 s \in \F_p^{\times}:& \quad U^sT^{-1/s} = W\begin{pmatrix}-s^{-1} & 1\\0 & -s\end{pmatrix}; &&\quad U^sW = T^{1/s}\begin{pmatrix}s & -1\\0 & s^{-1}\end{pmatrix}  
\end{align*}

If $V$ is a free $R$-module with finite rank on which $G$ acts on the left (resp. right), $V^{\ast}$ denotes the contragredient, that is, $V^{\ast}=\mrm{Hom}_R(V,R)$, where $g \cdot f = f(g^{-1}\cdot)$. It is still a left (resp. right) representation of $G$.
Let $V^{\circ \ast}$ denote the left (resp. right) $R$-module $V$ on which every $g \in G$ acts as $(g^{-1})^T$. Finally, $V^{\vee}$ denotes the \emph{right} (resp. \emph{left}) representation of $G$ given by $\mrm{Hom}_R(V,R)$, where $f \cdot g = f(g\cdot)$. 

\lem[central-character-and-contragredient]{Let $M$ be a left $R[G]$-module and $\chi: \F_p^{\times} \rar R^{\times}$ be such that for every $n \in \F_p^{\times}$, $nI_2$ acts on $M$ by multiplication by $\chi(n)$. Then $M^{\circ \ast} \simeq M \otimes \chi^{-1}(\det)$.}

\demo{Consider the function $f: x \in M^{\circ \ast} \longmapsto (W\cdot_M x)\in M \otimes \chi^{-1}(\det)$, which is clearly an isomorphism of $R$-modules. For any $g \in G$, any $x \in M^{\circ \ast}$, one has 
\begin{align*}
f(g \cdot_{M^{\circ \ast}} x) &= f((g^T)^{-1} \cdot_M x) = f((\det{g}^{-1} I_2)W^{-1}gW \cdot_M x)= \chi(\det{g})^{-1} W \cdot_M (W^{-1}gW \cdot_M x) \\
&= \det{g}^{-1} g \cdot_M (Wx)= g\cdot_{M \otimes \chi^{-1}(\det)}f(x),
\end{align*} 
whence the conclusion.}

\section{Principal series representations}
\label{principal-series-definition}

Let $R$ be a $\Z\left[\frac{1}{p^2-1},\mu_{p-1}\right]$-algebra without nontrivial idempotents. In particular, if $\omega \in R^{\times}$ is such that $\omega^{p-1}=1$, then either $\omega=1$ or $\omega-1$ is invertible. Indeed, we can find polynomials $A,B \in \Z\left[\frac{1}{p-1},X\right]$ such that $A(X)(X^{p-1}-1)+B(X)(X-1)^2=X-1$. Thus $B(\omega)(\omega-1)^2=\omega-1$, so that, if $\omega \neq 1$, then $(\omega-1)B(\omega)$ is a nonzero idempotent, hence is equal to $1$, and $\omega-1 \in R^{\times}$.   

Let $\psi, \chi: \F_p^{\times} \rar R$ be two characters and $R(\psi,\chi)$ be a free $R$-module of rank $1$, and let $\begin{pmatrix}a & \ast\\0 & b\end{pmatrix} \in B$ act on $R(\psi,\chi)$ by $\psi(a)\chi(b)$. Then $\pi(\psi,\chi)$ denotes the left $R[G]$-module $R[G] \otimes_{\R[B]} R(\psi,\chi)$. When $R=\C$, these representations are the \emph{principal series} of $\GL{\F_p}$.

\prop[steinberg-rep-definition]{Fix a generator $g$ of $R(1,1)$. The \emph{Steinberg representation} $\mrm{St}$ of $G$ is the quotient $\pi(1,1)/\sum_{M \in G/B}{M \otimes g}$. $\mrm{St}$ is free of rank $p$ as a $R$-module, and the exact sequence $0 \rar R \rar \pi(1,1) \rar \mrm{St} \rar 0$ splits. }

\demo{In the exact sequence, $\nu: R \rar \pi(1,1)$ is the map $1 \longmapsto \sum_{M \otimes G/B}{M \otimes g}$. So it is enough to construct a $R[G]$-retraction $r: \pi(1,1) \rar R$ of $\nu$ and show that $\ker{R}$ is free over $R$. Take $r: M \otimes g \longmapsto (p+1)^{-1}$, and note that $\ker{r}$ has the $M \otimes g - I_2 \otimes g$ as a basis, where $M$ runs through $S' \backslash \{I_2\}$. }

\prop[pi11-endomorphisms]{The $R$-module $\mrm{End}_{R[G]}(\pi(1,1))$ is free of rank two. In particular, $\mrm{End}_{R[G]}{\mrm{St}}$ is a free $R$-module of rank one, so that, if $R$ is a field and $p \in R^{\times}$, then $\mrm{St}$ is absolutely irreducible.}

\demo{Fix a generator $g$ of $R(\mathbf{1},\mathbf{1})$. By construction of $\pi(1,1)$, \[\rho: u \in \mrm{End}_{R[G]}(\pi(1,1)) \longmapsto u(I_2 \otimes g) \in \pi(1,1)^B\] is an isomorphism. Since $I_2 \otimes g \in \pi(1,1)^B$, it is enough to show that the submodule of $\sum_{M \in S' \backslash \{I_2\}}{R  \cdot M \otimes g}$ made with its elements fixed by $B$ is free of rank one. 

Let $x=\sum_{t \in \F_p^{\times}}{c_t T^t\otimes g}+c_{\infty}W \in \pi(1,1)^B$ for some $c_t,c_{\infty} \in R$. Then, for any $a \in \F_p^{\times}$, $x=\Delta_{1,a}x=\sum_{t \in \F_p^{\times}}{c_tT^{at}\otimes g}+c_{\infty}W\otimes g$, so that the $c_t$ for $t \in \F_p^{\times}$ are constant. 

Moreover, one computes that $x=Ux=c_{\infty}T\otimes g + c_1 W \otimes g+\sum_{t \in \F_p^{\times}\backslash \{-1\}}{c_1 T^{t/(t+1)} \otimes g }$, so that $c_{\infty}=c_1$. Hence $\pi(1,1)^B$ is contained in $R(I_2\otimes g)+R\sum_{M \in S'}{M \otimes g}$, and it easy to check that the reverse inclusion holds.

Note that $\rho(\mrm{id})=I_2 \otimes g$, and that $\rho^{-1}\left(\sum_{M \in S'}{M \otimes g}\right)$ is a projection onto $\pi(1,1)^G$.
}

Write $\mrm{St}_{\psi}$ for the twisted representation $\mrm{St} \otimes \psi(\det)$, for some character $\psi: \F_p^{\times} \rar R^{\times}$. 

\prop[principal-series-irreducible]{Let $\psi,\chi: \F_p^{\times} \rar R^{\times}$ be distinct characters. Then $\mrm{End}_{R[G]}(\pi(\psi,\chi))$ is free of rank one. Thus, if $R$ is a field and $p \in R^{\times}$, $\pi(\psi,\chi)$ is absolutely irreducible. }

\demo{As in the previous Proposition, if we fix a generator $g$ of $R(\psi,\chi)$, we have a injective homomorphism $\rho: u \in \mrm{End}_{R[G]}(\pi(\psi,\chi)) \longmapsto u(I_2 \otimes g) \in \pi(\psi,\chi)$, the image of which is exactly the submodule $\Upsilon=\{v \in \pi(\psi,\chi)\mid \forall M \in B,\, Mv=\psi([M]_{1,1})\chi([M]_{2,2})v\}$. The goal is to show that $\Upsilon=R(I_2 \otimes g)$. 

Let $x=\sum_{t \in \F_p}{c_tT^t\otimes g}+c_{\infty}W \otimes g \in \Upsilon$. For any $a \in \F_p^{\times}$, 
\[\chi(a)x=\Delta_{1,a}x=\sum_{t \in \F_p}{\chi(a)c_tT^{at} \otimes g}+\psi(a)c_{\infty}W \otimes g,\]
so that, by identifying coefficients, $c_{a}=c_1$ and for any $a \in \F_p^{\times}$, $(\psi(a)-\chi(a))c_{\infty}=0$. Let $a \in \F_p^{\times}$ be such that $\psi(a)\neq \chi(a)$, and let $\omega=\frac{\psi(a)}{\chi(a)}$. By construction, $\omega \neq 1$, $\omega^{p-1}=1$, and $(\omega-1)c_{\infty}=0$: $\omega-1$ is invertible and $c_{\infty} \neq 0$.  
Now, 

\[x=Ux=c_0I_2 \otimes g+\sum_{t \in \F_p^{\times} \backslash \{-1\}}{c_1\frac{\psi}{\chi}(t+1)T^{\frac{t}{t+1}} \otimes g} + (\psi\chi)(-1) \cdot c_1 W \otimes g,\] whence $c_1=0$ and $x \in R(I_2 \otimes g)$. Thus, $\Upsilon \subset R(I_2 \otimes g)$, and the reverse inclusion is clear.} 

\prop[twist-principal-series]{Let $\psi,\chi,\theta: \F_p^{\times} \rar R^{\times}$ be characters. Then $\pi(\psi\theta,\chi\theta)$ is canonically isomorphic to $\pi(\psi,\chi)\otimes\det{\theta}$. }

\demo{Let $g,g'$ be generators of $R(\psi,\chi),R(\psi\theta,\chi\theta)$. The morphism \[rI_2 \otimes g' \in R(\psi\theta,\chi\theta) \longmapsto rI_2 \otimes g \in \pi(\psi,\chi)\otimes \det{\theta}\] is a $R[B]$-morphism, so it induces a $R[G]$-morphism $u: \pi(\psi\theta,\chi\theta) \rar \pi(\psi,\chi) \otimes \det{\theta}$. The image of $u$ is the $R[G]$-submodule generated by $u(I_2 \otimes g')=I_2 \otimes g$, hence $u$ is surjective. Since $u$ is $R$-linear and surjective between two free $R$-modules of same rank, hence it is an isomorphism.}

\prop{Let $\psi_1,\chi_1,\psi_2,\chi_2: \F_p^{\times} \rar R$ and assume that $R(\psi_1,\chi_1),R(\psi_2,\chi_2)$ have the same traces. Then $\{\psi_1,\chi_1\}=\{\psi_2,\chi_2\}$.}

\demo{Left multiplication by $\Delta_{1,a}$ (for $a \in \F_p^{\times} \backslash \{0,1\}$) acts on $G/B$ as follows: the $T^tB$ for $t \in \F_p^{\times}$ are permuted without any fixed point, and $\Delta_{1,a}W=W\Delta_{a,1}$. Thus, the trace of the action of $\Delta_{1,a}$ on $R(\psi_i,\chi_i)$ is $\chi_i(a)+\psi_i(a)$. Moreover, the matrix $aI_2$ acts on $R(\psi_i,\chi_i)$ by multiplication by $\psi_i(a)\chi_i(a)$. Since $p+1 \in R^{\times}$, we can then check that $\psi_1+\chi_1=\psi_2+\chi_2$ and $\psi_1\chi_1=\psi_2\chi_2$. Thus $(\psi_1-\psi_2)(\psi_1-\chi_2)=0$. Since $\psi_1-\psi_2,\psi_1-\chi_2$ take their values in $R^{\times}\cup \{0\}$. Therefore, the reunion of the two subgroups $\{\psi_1=\psi_2\},\{\psi_1=\chi_2\}$ of $\F_p^{\times}$ is $\F_p^{\times}$, hence one has $\psi_1=\psi_2$ (thus $\chi_1=\chi_2$) or $\psi_1=\chi_2$ (thus $\chi_1=\psi_2$).  }

\prop{Let $\psi, \chi: \F_p^{\times} \rar R^{\times}$ be two distinct characters, $g \in R(\chi,\psi)$ be a generator, and \[v = W\otimes g+\sum_{a \in \F_p^{\times}}{\frac{\psi}{\chi}}(a)(T)^a \otimes g \in \pi(\chi,\psi).\]
Then \[u \in \mrm{Hom}_{R[G]}(\pi(\psi,\chi),\pi(\chi,\psi)) \longmapsto u(I_2 \otimes g) \in Rv\] is an isomorphism. Let $u^{\circ}: \pi(\psi,\chi) \rar \pi(\chi,\psi)$ be such that $u^{\circ}(I_2 \otimes g)=v$, then:
\begin{itemize}[noitemsep,label=$-$]
\item if $p \in R^{\times}$, $u^{\circ}$ is an isomorphism.
\item if $\F_p \subset R$, then $u^{\circ}$ is neither injective nor onto, its kernel and cokernel are free of same positive rank. 
\end{itemize}}

\demo{It is clear that $u \in \mrm{Hom}_{R[G]}(\pi(\psi,\chi),\pi(\chi,\psi)) \longmapsto u(I_2 \otimes g) $ is an isomorphism onto the submodule $\Upsilon \subset \pi(\chi,\psi)$ on which $B$ acts by the character $R(\psi,\chi)$. For any $a \in \F_p^{\times}$, $aI_2$ acts on $\pi(\chi,\psi)$ by multiplication by $(\psi\chi)(a)$; thus $\Upsilon$ is the submodule of the $x \in \pi(\chi,\psi)$ such that $Ux=x$ and $\Delta_{1,a}x = \chi(a)x$. 
We compute: 
\begin{align*}
\Delta_{1,a}v &= W\Delta_{a,1} \otimes g\sum_{t \in \F_p^{\times}}{\frac{\psi}{\chi}(t)T^{ta}\Delta_{1,a} \otimes g} = \chi(a)W \otimes g+\chi(a)\sum_{t \in \F_p^{\times}}{\psi(t)\psi(a)\chi^{-1}(ta)T^{ta}\otimes g}\\
&=\chi(a)v,\\
Uv &= T\otimes g+\sum_{t \neq 0, -1}{\frac{\chi}{\psi}(t+1)\frac{\psi}{\chi}(t)T^{t/(t+1)}\otimes g} + \frac{\psi}{\chi}(-1)\psi\chi(-1)W \otimes g \\
&= W\otimes g+T\otimes g+\sum_{t \neq 0,1}{\frac{\psi}{\chi}(t)T^t\otimes g} = v,\\
\end{align*}
thus $Rv \subset \Upsilon$. If the inclusion is strict, then, by substracting the correct multiple of $v$, there exists some $x=\sum_{t \in \F_p}{c_tT^t\otimes g} \in \Upsilon \backslash \{0\}$.
Now, for each $a \in \F_p^{\times}$, $\chi(a)x = \Delta_{1,a}x = \sum_{t \in \F_p}{\psi(a)c_tT^{ta} \otimes g}$: after identifying the coefficients, we see that $c_a = \frac{\psi}{\chi}(a)c_1$ and $c_0=0$. Since \[x=Ux \in \sum_{t \in \F_p}{RT^t\otimes g} + \frac{\psi}{\chi}(-1)\psi\chi(-1)c_1 W \otimes g,\] one has $c_1=0$ thus $x=0$, thus contradicting our assumption.

Now, we discuss $u^{\circ}$. We can check that  
\begin{align*}
\sum_{s \in \F_p}{T^s v} &= \sum_{a\in \F_p^{\times}, s \in \F_p}{\frac{\psi}{\chi}(a)T^{a+s}\otimes g} + \sum_{s \in \F_p}{T^sW\otimes g}\\
&= \sum_{t \in \F_p}{\sum_{a \in \F_p^{\times}}{\frac{\psi}{\chi}(a)} \cdot T^t\otimes g} + \sum_{s \in \F_p}{WU^{-s}\otimes g}\\
&= pW \otimes g,
\end{align*}
thus the cokernel of $u^{\circ}$ is of $p$-torsion. Therefore, if $p \in R^{\times}$, $u^{\circ}$ is a surjective map of free $R$-modules of same finite rank, hence an isomorphism.  

Now assume that $\F_p \subset R$. Then, since $X^{p-1}-1=\prod_{1 \leq i <p}{(X-i)}$, $\psi,\chi$ must take their values in $\F_p^{\times}$, so all the constructions are defined over $\F_p$ (and obtained after extending scalars by $\F_p \rar R$). Over $\F_p$, we just computed that $u^{\circ}\left(\sum_{s \in \F_p}{T^s \otimes g}\right) = \sum_{s \in \F_p}{T^sv} = 0$, hence $u^{\circ}$ is not injective. Since $u^{\circ}$ is linear between two $\F_p$-vector spaces of dimension $p+1$, it is not surjective either, and its kernel and cokernel have the same finite rank. This claim is stable under the faithfully flat extension of scalars $\F_p \rar R$.
}

\rem{What this result says is that $\pi(\chi,\psi)$ is symmetric with respect to $(\psi,\chi)$ in characteristic zero (or more generally distinct from $p$), but not in characteristic $p$.}

\lem{Let $\psi,\chi: \F_p^{\times} \rar R^{\times}$ be characters, then $\pi(\psi,\chi)^{\circ \ast} \simeq \pi(\chi^{-1},\psi^{-1})$. In particular, with $\psi=\chi=1$, we get an isomorphism $\mrm{St}^{\circ \ast} \rar \mrm{St}$.}  

\demo{As in the previous proofs, $u \in \mrm{Hom}_G(\pi(\chi^{-1},\psi^{-1}),\pi(\psi,\chi)^{\circ \ast}) \longmapsto u(I_2 \otimes g)$ is a bijection onto the sub-$R$-module $\Upsilon$ made with those $x \in \pi(\psi,\chi)^{\circ \ast}$ such that any $M \in B$ satisfies $Mx = \chi^{-1}([M]_{1,1})\psi^{-1}([M]_{2,2})x$. In other words, $\Upsilon$ is the set of those $x \in \pi(\psi,\chi)$ such that $\Delta_{a,b}x = \chi(a)\psi(b)x$ for all $a,b \in \F_p^{\times}$ and $T^{-1} x = x$. The module $\Upsilon$ contains $x = W \otimes g$ since $W\Delta_{a,b} = \Delta_{b,a}W$ and $T^{-1}W=WU$. Let $u: \pi(\chi^{-1},\psi^{-1}) \rar \pi(\psi,\chi)^{\circ \ast}$ be such that $u(I_2 \otimes g)=W \otimes g$: then $u$ is a $R[G]$-morphism, and its image contains $W \otimes g$ and is a $R[G]$-module, so $u$ is onto. Since $u$ is $R$-linear between free $R$-modules of same rank, $u$ is an isomorphism.}

\lem[principal-series-dual]{Let $\psi,\chi: \F_p^{\times} \rar R^{\times}$ be characters, then $\pi(\chi^{-1},\psi^{-1})$ and $\pi(\chi,\psi)^{\ast}$ are isomorphic.}

\demo{As above, we only need to exhibit some $f \in \mrm{Hom}_R(\pi(\chi,\psi),R)$ such that $f(U^{-1} \cdot ) = f$, $f(\Delta_{a,b}^{-1}\cdot ) =\chi^{-1}(a)\psi^{-1}(b)f$, and such that the $R[G]$-module generated by $f$ is $\mrm{Hom}_R(\pi(\chi,\psi),R)$. Let $g$ be a generator of $R(\chi,\psi)$, and $f: \pi(\chi,\psi) \rar R$ mapping $I_2 \otimes g$ to $1$ and $M \otimes g$ to $0$ for any $M \in S' \backslash \{T^t,\, t \in \F_p^{\times}\}\cup\{W\}$. It is then formal to check that $f$ satisfies the requested conditions.  }

\cor[steinberg-dual]{The isomorphism of Lemma \ref{principal-series-dual} induces (when $\psi=\chi$) an isomorphism $\mrm{St}_{\psi^{-1}} \rar \mrm{St}_{\psi}^{\ast}$.}

\demo{It is enough to show the result when $\chi=\psi=\mathbf{1}$, and in this case one has $\pi(\mathbf{1},\mathbf{1})=\mrm{St} \oplus \mathbf{1}$. Assume that the induced map $\mathbf{1} \rar \pi(\mathbf{1},\mathbf{1}) \rar \pi(\mathbf{1},\mathbf{1})^{\ast} \rar \mrm{St}^{\ast}$ vanishes. Then the composition $\mrm{St} \rar \pi(\mathbf{1},\mathbf{1}) \rar \pi(\mathbf{1},\mathbf{1})^{\ast} \rar \mrm{St}^{\ast}$ is $R$-linear and onto between free $R$-modules of same rank, hence is an isomorphism. 

So the goal is to show that any $G$-equivariant morphism $\mrm{St} \rar R$ vanishes. It is enough to show that a $G$-equivariant morphism $\pi(\mathbf{1},\mathbf{1}) \rar R$ vanishing at $\nu := \sum_{M \in B\backslash G}{M \cdot g}$ (where $g$ is a generator of $R(1,1)$) is zero, or, in other words, that a $G$-equivariant $f: \pi(\mathbf{1},\mathbf{1}) \rar R$ is determined by $f(\nu)$. Note that $f$ is determined by $f(g)$, and $f(\nu)=f((p+1)g) \in R^{\times}f(g)$ by construction.
}

\prop[character-from-principal-series]{Let $\alpha,\beta: \F_p^{\times} \rar R^{\times}$ be characters, $a,b \in \F_p^{\times}$ be distinct, and $C \leq G$ be a nonsplit Cartan subgroup. 
\begin{itemize}[label=\tiny$\bullet$,noitemsep]
\item $aI_2 \in \GL{\F_p}$ acts on $\pi(\alpha,\beta)$ (resp. $\mrm{St}_{\alpha}$) by multiplication by $(\alpha\beta)(a)$ (resp. $\alpha^2(a)$).
\item The trace of $U$ acting on $\pi(\alpha,\beta)$ (resp. $\mrm{St}_{\alpha}$) is $1$ (resp. $0$). 
\item The trace of $\Delta_{a,b}$ acting on $\pi(\alpha,\beta)$ (resp. $\mrm{St}_{\alpha}$) is $\alpha(a)\beta(b)+\alpha(b)\beta(a)$ (resp. $\alpha(ab)$).
\item The trace of $g \in C \backslash \F_p^{\times}I_2$ acting on $\pi(\alpha,\beta)$ (resp. $\mrm{St}_{\alpha}$) is $0$ (resp. $-\alpha(\det{g})$). 
\end{itemize}}

\demo{The statements about $\mrm{St}_{\alpha}$ follow from the fact that $\pi(\alpha,\alpha) \simeq \alpha(\det)\oplus \mrm{St}_{\alpha}$. The first point is because $\F_p^{\times}I_2$ is central and that $aI_2$ acts by $(\alpha\beta)(a)$ on $R(\alpha,\beta)$. For the second point, note that $\pi(\alpha,\beta)$ has a basis made with the vectors $e_M$ for $M \in G/B$ and that for $N \in G$, $N \cdot e_M$ is a unit multiple of $e_{NM}$. 

Thus, the trace of any matrix $N \in \GL{\F_p}$ acting on $\pi(\alpha,\beta)$ is the sum, over all elements $x \in G/B$ stable under left multiplication by $N$, of the $\alpha([P]_{1,1})\beta([P]_{2,2})$, where $Nx'=x'P$, $P \in B$ and $x' \in G$ lifts $x$. The claims then follow from the following facts: the left multiplication by any $g \in C \backslash \F_p^{\times}I_2$ has no fixed point in $G/B$ (since nonscalar elements of $C$ do not have eigenvalues in $\F_p$), while $B$ is the only fixed point of the left multiplication by $U$ on $G/B$, and left multiplication by $\Delta_{a,b}$ on $G/B$ has the two fixed points $B$ and $WB$, with $\Delta_{a,b}I_2=I_2\Delta_{a,b}$, $\Delta_{a,b}W=W\Delta_{b,a}$.}

\section{Cuspidal representations}
\label{cuspidal-reps-definition}

In this section, we assume that $R$ is a field of characteristic zero containing the $p(p^2-1)(p-1)$-th roots of unity. 

The $\psi(\det)$, the $\mrm{St}_{\psi}$, and the $\pi(\psi,\chi)$ (we say of all these representations that they \emph{come from principal series}) form exactly $2(p-1)+\frac{(p-1)(p-2)}{2} = \frac{(p-1)(p+2)}{2}$ pairwise non-isomorphic irreducible representations of $G$. When counting against the conjugacy classes of $G$, we see that we are missing $\frac{p(p-1)}{2}$ isomorphism classes of irreducible representations of $G$.

Moreover, the sum of the squared dimensions of these missing irreducible representations is \[(p-1)^2p(p+1)-(p-1)-(p-1)p^2-\frac{(p-1)(p-2)(p+1)^2}{2} = \frac{p(p-1)^3}{2}.\] 

\lem{Let $V$ be an irreducible representation of $G$. Then $V^U \neq 0$ if and only if $V$ comes from principal series.}

\demo{If $V=\psi(\det)$ or $V=\pi(\chi,\psi)$, $V^U \neq 0$ by construction. If $V=\mrm{St}_{\psi}$, then there is a surjective morphism $\pi(\psi,\psi) \rar \mrm{St}_{\psi}$ (by Proposition \ref{steinberg-rep-definition}. Since the space of $\pi(\psi,\psi)$ is generated by some element fixed under $U$, the same holds for $\mrm{St}_{\psi}$. 

Conversely, if $V$ is irreducible and $V^U \neq 0$, then there is a character $\alpha: B/\langle U \rangle \rar R^{\times}$ such that $V^U$ contains a nonzero element $x$ on which $B$ acts by $\alpha$. If $\alpha(M)=\psi([M]_{1,1})\chi([M]_{2,2})$, $x$ defines a nonzero morphism $\pi(\psi,\chi) \rar V$, which is then surjective, and $V$ comes from principal series by Section \ref{principal-series-definition}. }

\defi{A representation $V$ of $G$ is \emph{cuspidal} if $V^U=\{0\}$. Therefore, there are exactly $\frac{p(p-1)}{2}$ isomorphism classes of irreducible cuspidal representations of $G$.}

\lem[lower-bound-dim-cusp]{If $V$ a cuspidal representation of $G$, then $\dim{V} \geq p-1$.}

\demo{Let $u \in \mathcal{L}(V)$ be the endomorphism induced by $U$. Then $u^p$ is the identity, but $1$ is not an eigenvalue of $u$. Thus the set $S$ of eigenvalues of $u$ is made with primitive $p$-th roots of unity. For any $a \in \F_p^{\times}$, $U$ and $U^a$ are conjugate in $G$, so $u^a$ is conjugate to $u$ in $GL(V)$, thus $S^a \subset S$. Thus, $S$ contains all the primitive $p$-th roots of unity, hence $\dim{V} \geq |S| \geq p-1$.}

\cor{Every irreducible cuspidal representation of $G$ has dimension $p-1$; a representation of $G$ is cuspidal if and only if every one of its irreducible factors is cuspidal. }

\demo{There are exactly $\frac{p(p-1)}{2}$ isomorphism classes of cuspidal representations, all of them have dimension at least $p-1$, and the sum of the squares of their dimensions is $\frac{p(p-1)^3}{2}$.}

\cor[character-cuspidal-on-borel]{If $V$ is an irreducible cuspidal representation of $G$, then $U$ acts with trace $-1$, and any matrix $\Delta_{a,b}$ with $a \neq b$ acts with trace $0$.}

\demo{By the proof of Lemma \ref{lower-bound-dim-cusp}, the eigenvalues of the endomorphism $u$ of $U$ contain all the primitive $p$-th roots of unity. Since $\dim{V}=p-1$, the eigenvalues of $u$ are exactly the primitive $p$-th roots of unity and they are all simple. Thus $\mrm{Tr}(u)=\sum_{\omega^p=1,\omega \neq 1}{\omega} = -1$. 
Now, $\Delta_{a,b}U^{b/a}=U\Delta_{a,b}$, so that the endomorphism $\delta$ of $V$ induced by $\Delta_{a,b}$ preserves the $p-1$ eigenspaces of $u$, and maps the $\omega$-eigenspace to the $\omega^{b/a}$-eigenspace. In particular, $\delta$ permutes the eigenspaces of $u$ without any fixed point, so $\mrm{Tr}(\delta)=0$.}

\prop[characterize-cuspidals]{Let $C \leq G$ be a nonsplit Cartan subgroup. To any irreducible cuspidal representation $V$ of $\GL{\F_p}$ we may attach a pair $\{\chi,\chi^p\}$ of characters of $C$ that do not factor through $\det$ such that:
\begin{itemize}[label=$-$,noitemsep]
\item the pair $\{\chi,\chi^p\}$ characterizes the isomorphism class of $V$, 
\item every such pair is attached to some irreducible cuspidal representation,
\item for any $a \in \F_p^{\times}$, $aI_2$ acts on $C$ by $\chi(a)$,
\item the character of $V$ is that of the virtual representation $\mrm{Ind}_{\F_p^{\times}\langle U\rangle}^G{\left[aU^b \longmapsto \chi(a)\psi(b)\right]}-\mrm{Ind}_C^G{\chi}$, for any nontrivial $\psi: \F_p \rar R^{\times}$,
\item for any $g \in C \backslash \F_p^{\times}I_2$, the trace of $g$ acting on $V$ is $-\chi(g)-\chi^p(g)$. 
\end{itemize}}

\demo{This is a direct consequence of the construction of \cite[(6.4) Theorem]{BH}, as long as the character computation holds (and the action of $\F_p^{\times}I_2$ is correctly described). 
For the action of $\F_p^{\times}I_2$, it follows from the fact that these matrices are central. So all that remains is the character computation at $g \in C \backslash \F_p^{\times}I_2$. Since $g$ is not similar to any $aU^b$, by \cite[\S 7.2 Proposition 20]{SerreLinReps}, 
\begin{align*}
&\mrm{Tr}\left(g \mid \mrm{Ind}_{\F_p^{\times}\langle U\rangle}^G{\left[aU^b \longmapsto \chi(a)\psi(b)\right]}-\mrm{Ind}_C^G{\chi}\right)\\
&= -\mrm{Tr}\left(g \mid \mrm{Ind}_C^G{\chi}\right) = -\sum_{h \in G/C, h^{-1}gh \in C}{\chi(h^{-1}gh)}=-\chi(g)-\chi(g^p),
\end{align*}
since the only matrix in $C$ conjugate to $g$ is $g^p$, and the centralizer of $g$ is exactly $C$.}

\rem{Given any nonsplit Cartan subgroup $C \leq G$, there are two isomorphisms $C \rar \F_{p^2}^{\times}$ preserving the trace and the determinant: if $\iota$ denotes one of them, the other one is $\iota^p$. Thus, we will often identify, when discussing cuspidal representations, characters of $C$ and characters of $\F_{p^2}^{\times}$. If $C'$ is another nonsplit Cartan subgroup identified with $\F_{p^2}^{\times}$, then the composition isomorphism $C \rar C'$ preserves the trace and the determinant, and is thus given by conjugation by some matrix $M \in G$, so this is harmless.  }

\cor[cuspidal-twist]{Let $C \leq G$ be a nonsplit Cartan subgroup and $V$ an irreducible cuspidal representation of $\GL{\F_p}$ attached to some character $\phi$ of $C$. For any Dirichlet character $\psi: \F_p^{\times} \rar \C^{\times}$, $V \otimes \psi(\det)$ is an irreducible cuspidal representation of $\GL{\F_p}$ attached to the character $\phi\cdot \psi(\det)$. }

\demo{The representation $V \otimes \psi(\det)$ is irreducible, and, since $\det{U}=1$, it is cuspidal. The character $\phi\cdot \psi(\det): C \rar \C^{\times}$ does not factor through $\det$, and it is a direct computation using Proposition \ref{characterize-cuspidals} and Corollary \ref{character-cuspidal-on-borel} that $V \otimes \psi(\det)$ has the same character as the cuspidal representation attached to $\phi\cdot \psi(\det)$.}

\cor[cuspidal-duals]{If $V$ is an irreducible cuspidal representation of $G$, then $V^{\circ \ast} \simeq V^{\ast}$; if $V$ is attached to the pair $\{\chi,\chi^p\}$, $V^{\ast}$ is attached to the pair $\{\overline{\chi},\overline{\chi}^p\}$.  } 

\lem[det-cusp]{Let $V$ be a cuspidal representation of $G$. Then $\det{V}=\left(\frac{\det}{p}\right)$.}

\demo{Let $\chi$ be a character of some nonsplit Cartan subgroup attached to $V$. Then $\det{V}$ is an abelian representation of $G$ of dimension one, so it is of the form $\psi(\det)$ for some $\psi: \F_p^{\times} \rar R^{\times}$. For any $a \in \F_p^{\times}$, $\psi(a^2)=\det(aI_2 \mid V)=\chi(aI_2)^{p-1}=\chi((aI_2)^{p-1})=1$, so $\psi$ is quadratic. 

Let $g \in \F_p^{\times}$ be a generator and $A \in \mathcal{M}_{p-1}(\C)$ be the matrix of its action on $V$. Since $g^{p-1}=1$, one has $A^{p-1}=I_{p-1}$ so $A$ is diagonalizable. Moreover, we proved that if $1 \leq k < p-1$, $\mrm{Tr}(A^k)=0$. Thus, if $\lambda_1,\ldots,\lambda_{p-1}$ denote the eigenvalues of $A$ (all of them $(p-1)$-th roots of unity), for any polynomial $P \in \C[X]$ with degree less than $p-1$, one has \[\sum_{i=1}^{p-1}{P(\lambda_i)} = \mrm{Tr}(P(A)) = (p-1)P(0).\] If the $\lambda_i$ are not exactly the $(p-1)$-th roots of unity, we can choose $P$ such that the $\lambda_i$ are exactly the roots of $P$. Then we get $(p-1)P(0)=\sum_{i=1}^{p-1}{0}=0$, whence a contradiction.  }

\prop[cusp-coeffs]{Let $V$ be an irreducible cuspidal representation of $G$ with trivial central character. It is attached to a pair $\{\phi,\phi^p\}$, for some character $\phi: \F_{p^2}^{\times}/\F_p^{\times} \rar R^{\times}$. Then $V$ can be realized over $\Q(i,\phi)^+$, where $\Q(i,\phi)^+$ is the maximal real subextension of $R/\Q$ contained in the subfield of $R$ generated by $i$ and the values of $\phi$. In particular, $V$ is realizable over $\R$, and in fact over $\Q$ if $\phi$ has order $4$.}

\demo{Let $C$ be a nonsplit Cartan subgroup of $G$, we identify $\phi$ with a character of $C/\F_p^{\times}I_2$. Since $\phi$ does not factor through the determinant, $\phi$ has order at least $3$. Fix a nontrivial character $\psi: \F_p \rar R^{\times}$. By Proposition \ref{characterize-cuspidals} and \cite[\S 12.1 Proposition 33]{SerreLinReps}, it is enough to check that the representations $\mrm{Ind}_{\F_p^{\times}U^{\Z}}^B\left[aU^b \longmapsto \psi(b)\right]$ and $\mrm{Ind}_C^G{\phi}$ can be realized over $\Q(i,\phi)^+$.

We show that $\mrm{Ind}_{\F_p^{\times}U^{\Z}}^B\left[aU^b \longmapsto \psi(b)\right]$ can be realized over $\Q$. Indeed, let $M=\Q(\mu_p)$, and we consider the following action of $B/\F_p^{\times}I_2$ on $M$: a matrix $aU^b\Delta_{c,1}$ (for $a,c \in \F_p^{\times}$ and $b \in \F_p$) acts by $\psi(b)\underline{c}$, where $\underline{c} \in \mrm{Gal}(\Q(\mu_p)/\Q)$ is the element whose cyclotomic character modulo $p$ is $c$. One directly checks that this action is well-defined, and that there is a natural map $\lambda: \mrm{Ind}_{\F_p^{\times}U^{\Z}}^B\left[aU^b \longmapsto \psi(b)\right] \rar M$ mapping a generator of the one-dimensional representation with character $aU^b \longmapsto \psi(b)$ to $1$. Since $1$ generates $M$ over $\Q[B]$, $\lambda$ is onto; since it is between $\Q$-vector spaces of same dimension, $\lambda$ is a $\Q[B]$-isomorphism.

Next, let $N$ be the normalizer of $C$: $C \leq N$ is a subgroup of index two, and (for instance by computing an explicit example), there exists some $d \in N \backslash C$ with order two. Since $N/C$ acts by taking the $p$-th power, it exchanges $\phi$ and $\phi^p$; since $\phi(\F_p^{\times})=1$ by assumption, $\phi^p=\phi^{-1}$. We can then directly check that \[g \in D \longmapsto \begin{pmatrix}\Re{\phi(g)} & -\Im{\phi(g)}\\\Im{\phi(g)} & \Re{\phi(g)}\end{pmatrix} \in \mathcal{M}_2(\R),\, t \longmapsto \begin{pmatrix} 0 & 1\\1 & 0\end{pmatrix}\] is a realization over $\Q(i,\phi)^+$ of $\mrm{Ind}_C^N{\phi}$, so that $\mrm{Ind}_C^G{\phi}$ is also realizable over $\Q(i,\phi)^+$.
}

\section{Restriction of representations to smaller subgroups}
\label{smallgroups}
In this section, we take $R=\C$ for the sake of simplicity, but the results apply to any algebraically closed field on which $p(p^2-1)$ is invertible. The letter $D$ denotes the subgroup of $B$ made with diagonal matrices. 

The letter $C$ always denotes a nonsplit Cartan subgroup with normalizer $N$. The character $N \rar N/C \simeq \{\pm 1\}$ is noted $\epsilon$. Given a character $\psi: C \rar \C^{\times}$, $\iota_{\psi}$ denotes the representation $\mrm{Ind}_C^N{\psi}$, so that $\iota_{\psi}=\iota_{\psi^p}$. It is reducible if and only if $\psi$ factors through $\det: C \rar \F_p^{\times}$, or in other words if $\psi^p=\psi$. 

\lem[steinberg-on-borel]{$\mrm{St}^U=\mrm{St}^B$ is a line, and the eigenvalues of the endomorphism of $\mrm{St}$ induced by $U$ are exactly the $p$-th roots of unity. }

\demo{By Proposition \ref{pi11-endomorphisms}, one has \[\dim{\pi(1,1)^B}=\dim\mrm{Hom}_B(\mathbf{1},\pi(1,1))=\dim\mrm{Hom}_{\GL{\F_p}}(\pi(1,1),\pi(1,1))=2,\] so, since $\pi(1,1) \simeq \mrm{St}\oplus \mathbf{1}$, it is enough to show that $\dim{\pi(1,1)^U}=2$. Let $u$ be the endomorphism of $\mrm{St}$ induced by $U$. Since $\pi(1,1) \simeq \C[G/B]$, $\mrm{Tr}(u)=|(G/B)^U|-1=0$. Since $U^p=I_2$ and $U$ is similar to any $U^a$ (for $a \in \F_p^{\times}$), $u^p$ is the identity and $\mrm{Tr}(u^a)=0$ for any $a \in \F_p^{\times}$. Our goal is to show that the eigenvalues of $u$ are exactly the $p$-th roots of unity, each with multiplicity one. If this is not the case, then there is a $p$-th root of unity which is not an eigenvalue of $u$. Therefore, there is a polynomial $P$ of degree less than $p$, whose roots are $p$-th roots of unity, and such that $P(u)=0$. Then $0=\mrm{Tr}(P(u))=pP(0)$, so $P(0)=0$, a contradiction. }

\lem[steinberg-vs-cartan]{Let $g \in C$ be a generator and $x \in \mrm{St}^B$ be a nonzero vector. Let $P \in \C[x]$; then $P(g)\cdot x = 0$ if and only $P(g) \cdot \mrm{St} = 0$ if and only if $X^{p+1} -1 \mid (X-1)P(X)$.}

\demo{Let $\gamma$ be the endomorphism of the space of $\mrm{St}$ induced by $g$. By Proposition \ref{character-from-principal-series}, for any $1 \leq k \leq p$, $\mrm{Tr}(\gamma^k)=-1$ and $\gamma^{p+1}=\mrm{id}$. Thus, the characteristic polynomial of $\gamma$ is $P_0(X)=\frac{X^{p+1}-1}{X-1} \in \C[X]$, whose (simple) roots are the nontrivial $(p+1)$-th roots of unity. Since $\mrm{St}$ is irreducible, there is an exact sequence of $\C[\GL{\F_p}]$-modules \[0 \rar \C \overset{u}{\rar} \C[\GL{\F_p}/B] \overset{v}{\rar} \mrm{St} \rar 0,\quad u(1) = \sum_{h \in \GL{\F_p}/B}{[h]},\quad v(I_2)=x.\]
 
The natural function $C/\F_p^{\times} \rar \GL{\F_p}/B$ is between sets of cardinality $p+1$; it is injective, because $C \cap B = \F_p^{\times}I_2$, so is bijective. Hence, if $P \in \C[X]$ is such that $P(g) \cdot x = 0$, then $P(g) \cdot I_2 \in \C\cdot \sum_{h \in \GL{\F_p}/B}{[h]}$, thus $\deg{P} \geq p$ or $P=0$. In particular, if $P_1$ denotes the remainder in the Euclidean division of $P$ by $P_0$, then $P_1=0$, that is, $P_0 \mid P$. The conclusion follows.
}

\lem[steinberg-on-diag]{Let $q:\Delta_{a,b} \in D \longmapsto \frac{a}{b} \in \F_p^{\times}$. Then \[\mrm{St}_{|D} \simeq \mathbf{1} \oplus \bigoplus_{\chi: \F_p^{\times} \rar \C^{\times}}{\chi(q)}.\] }

\demo{Since both of these representations have the same dimension and factor through $D/\F_p^{\times}$, it is enough to show that for each $g \in \F_p^{\times} \backslash \{1\}$, the trace of the action of $\Delta_{g,1}$ on both representations is the same. Clearly, the trace of the action of $\Delta_{g,1}$ on the right-hand side is $1+\sum_{\chi: \F_p^{\times} \rar \C^{\times}}{\chi(g)} = 1$. By Proposition \ref{character-from-principal-series}, the trace of the action of $\Delta_{g,1}$ on $\mrm{St}$ is $1+1-1=1$, whence the conclusion. }

\lem[steinberg-on-norm-cartan]{Let $C \leq \GL{\F_p}$ be a nonsplit Cartan subgroup with normalizer $N$. Then \[\mrm{St}_{|N} \simeq \alpha \oplus \bigoplus_{\psi}{\iota_{\psi}},\] where:
\begin{itemize}[label=\tiny$\bullet$,noitemsep]
\item $\alpha$ is the character $\left(\frac{\det(\cdot)}{p}\right)\epsilon^{ \mathbf{1}(-1 \notin \F_p^{\times 2})}$
\item $\psi$ runs through the characters $C/\F_p^{\times}I_2 \rar \C^{\times}$ such that $\psi^2 \neq 1$, and always counting exactly one character in any pair $\{\psi,\psi^{-1}\}$. 
\end{itemize}}

\demo{It is a character computation. First, both representations have the same dimension and factor through $N/\F_p^{\times}I_2$. 

Next, for any $g \in C \backslash \{I_2\}$, the trace of the action of $g$ on the right hand side is exactly $\alpha(g)-1-\left(\frac{\det{g}}{p}\right)+\sum_{\psi}{\psi(g)}$, where $\psi$ runs through characters $C/\F_p^{\times} \rar \C^{\times}$, so the trace is exactly $\alpha(g)-1-\frac{\det{g}}{p}=-1$ since $\epsilon(g)=1$. The trace of $\mrm{St}(g)$ is $-1$ by Proposition \ref{character-from-principal-series}, so the characters agree. 

For any $g \in N \backslash C$, $g$ acts with trace $0$ on any $\iota_{\psi}$, so the action of $g$ on the right-hand side has trace $\alpha(g)=\left(\frac{-\det{g}}{p}\right)$. The action of $g$ on the left-hand side is, by Proposition \ref{character-from-principal-series}, $1$ if $g$ is not diagonalizable over $\F_p$, and $-1$ otherwise. Since $g \in N \backslash C$, $\mrm{Tr}(g)=0$, so $g^2=-\det{g}I_2$. So $g$ is diagonalizable over $\F_p$ if and only if $g^{p-1}=I_2$, if and only if $(-\det{g})^{(p-1)/2}=1$, if and only if $\left(\frac{-\det{g}}{p}\right)$: in other words, the action of $g$ on $\mrm{St}$ also has trace $\left(\frac{-\det{g}}{p}\right)$, whence the conclusion. 
}

\lem[principal-on-norm-cartan]{Let $C \leq \GL{\F_p}$ be a nonsplit Cartan subgroup with normalizer $N$, and $\chi: \F_p^{\times} \rar \C^{\times}$ be a character. Then 
\[\pi(1,\chi)_{|N} \simeq A \oplus \bigoplus_{\psi}{\iota_{\psi}},\] where:
\begin{itemize}[label=\tiny$\bullet$]
\item $A$ is the sum of the characters $\alpha(\det)\epsilon^{\mathbf{1}(\alpha(-1)=-1}$ over all $\alpha: \F_p^{\times} \rar \C^{\times}$ with $\alpha^2=\chi$ (so $A$ has dimension two if $\chi$ is even and is zero if $\chi$ is odd). 
\item $\psi$ runs through the characters $C \rar \C^{\times}$ such that $\psi(tI_2)=\chi(t)$ for every $t \in \F_p^{\times}$, such that $\psi^{p-1} \neq 1$, and always counting exactly one character in any pair $\{\psi,\psi^p\}$.  
\end{itemize}
}

\demo{It is a character computation. Note that any $aI_2$ (for $a \in \F_p^{\times}$) acts by $\chi(a)$ on either side. 

Moreover, the trace of the action of any matrix in $C \backslash \F_p^{\times}I_2$ is zero by Proposition \ref{character-from-principal-series}, and the trace of the action of any $\Delta_{a,b}$ with $a \neq b$ is $\chi(a)+\chi(b)$. For any $g \in N \backslash C$, $\mrm{Tr}(g)=0$, so $g^2=-\det{g}I_2$, so that $g$ is diagonalizable with eigenvalues in $\F_p$ if and only if $-\det{g} \in \F_p^{\times 2}$. In this case, the two eigenvalues of $g$ are $\pm\sqrt{-\det{g}}$. Thus the character of the left-hand side is the following:
\begin{align*}
& aI_2: \chi(a)(p+1),\, &&g \in C: 0,\\
& g \in N \backslash C, -\det{g} \notin \F_p^{\times 2}: 0,\, &&g \in N \backslash C,-\det{g} \in \F_p^{\times 2}: \chi(\sqrt{-\det{g}})+\chi(-\sqrt{-\det{g}}).
\end{align*} 

Assume that $\chi(-1)=-1$. Then $A$ is zero by definition and the set of $p+1$ characters $\psi: C \rar \C^{\times}$ extending $\chi$ does not contain any character with order $p-1$ (such a character factors through the determinant $C \rar \F_p^{\times}$, therefore vanishes at $-I_2$, contradicting the fact that $\chi(-1)=1$), so that $\psi \longmapsto \psi^p$ is an involution of this set of characters, and the right-hand side has dimension $p+1$ like the left-hand side. 

The character of the right-hand side at some $g \in C \backslash \F_p^{\times}I_2$ is $\sum_{\psi_{|\F_p^{\times}}=\chi}{\psi(g)} = 0$ (because these $\psi$ are cosets of the subgroup of characters $C/\F_p^{\times} \rar \C^{\times}$, and the characters in this subgroup do not all vanish at $g$). The character of the right-hand side is also zero on $N \backslash C$, which then matches the character of the left-hand side. 

Now, assume that $\chi(-1)=1$, so that $\chi=\alpha^2$ for some character $\alpha: \F_p^{\times} \rar \C^{\times}$. Let $\lambda: \F_p^{\times} \rar \{\pm 1\}$ be the unique nontrivial quadratic character. Then $\alpha(\det),(\alpha\lambda)(\det): C \rar \C^{\times}$ are two distinct characters whose restriction to $\F_p^{\times}$ is $\chi$. The characters $C \rar \C^{\times}$ agreeing with $\chi$ on $\F_p^{\times}$ are exactly of the form $\alpha(\det)\psi$ for some $\psi: C/\F_p^{\times} \rar \C^{\times}$, and such a character has order dividing $p-1$ if and only if $\psi^{p-1}=1$, if and only if $\psi^2=1$; there are at most two such characters, so, for any character $\psi: C \rar \C^{\times}$ agreeing with $\chi$ on $\F_p^{\times}$, $\psi^{p-1}=1$ if and only if $\psi \in \{\alpha(\det),(\alpha\lambda)(\det)\}$. It follows that the right-hand side has dimension $2\frac{p-1}{2}+2=p+1$.

For any $g \in C \backslash \F_p^{\times}I_2$, the trace of its action on the right-hand side is \[\mrm{Tr}(g \mid A)+\sum_{\psi}{(\psi+\psi^p)(g)}=\alpha(\det{g})-(\alpha\lambda)(\det{g})-\alpha(\det{g})-(\alpha\lambda)(\det{g})=0.\]  

Any $g \in N \backslash C$ acts on the right-hand side with trace 
\begin{align*}
\mrm{Tr}(g \mid A)+0&=\alpha(\det{g})(-1)^{\mathbf{1}(\alpha(-1)=-1)}+(\alpha\lambda)(\det{g})(-1)^{\mathbf{1}((\alpha\lambda)(-1)=-1)} \\
&= \alpha(-\det{g})(1+\lambda(-\det{g})).\end{align*} Hence this trace is zero if $-\det{g} \notin \F_p^{\times 2}$, and $2\alpha(-\det{g})=\chi(\sqrt{-\det{g}})+\chi(-\sqrt{-\det{g}})$. 

Thus the right-hand side and the left-hand side have the same characters. 
}

\prop[cusp-on-diag]{Write $(\alpha,\beta): \Delta_{u,v} \in D \longmapsto (u,v) \in \F_p^{\times} \times \F_p^{\times}$. Let $V$ be a cuspidal representation of $G$ attached to the character $\phi$ of $\F_{p^2}^{\times}$. Then \[V_{|D} \simeq \phi(\alpha) \otimes \bigoplus_{\psi: \F_p^{\times} \rar \C^{\times}}{\psi\left(\frac{\alpha}{\beta}\right)} .\]
}

\demo{The two sides have dimension $p-1$, and, by Proposition \ref{characterize-cuspidals}, any $aI_2 \in \F_p^{\times}$ acts by multiplication by the scalar $\phi(a)$ on either side, so the two characters agree at $aI_2$. If $a,b \in \F_p^{\times}$ are distinct, then $\Delta_{a,b}$ acts with trace zero on $V$, and the trace of the action of $\Delta_{a,b}$ on the right-hand side is $\phi(a)\sum_{\psi}{\psi(ab^{-1})} = 0$. }

\prop[cusp-on-borel]{Let $V$ be an irreducible cuspidal representation of $G$, attached to an even character $\phi$ of $\F_{p^2}^{\times}$. Let $\chi: \F_p^{\times} \rar \C^{\times}$ be a character such that $\chi^2=\phi_{|\F_p^{\times}}$. Let $U'\triangleleft B$ denote the abelian subgroup $\F_p^{\times}U^{\Z}$. Then $V$ is an irreducible representation of $B$ and 
\[V_{|B} \simeq \chi(\det) \otimes \mrm{Ind}_{U'}^B{\left[aU^b \longmapsto e^{2i\pi b/p}\right]} \simeq \mrm{Ind}_{U'}^B{\left[aU^b \longmapsto \phi(a)e^{2i\pi b/p}\right]}.\]}

\demo{Since $\dim{V}=p-1$, $V_{|U'}$ is a sum of $p-1$ characters of $U'$. Since $aI_2$ acts on $V$ by multiplication by $\phi(a)$, all of these characters are of the form $\gamma_{\phi}: aU^b \in U' \longmapsto \phi(a)\gamma(b)$, for some character $\gamma: \F_p \rar \C^{\times}$. Since $V^U$ is null, $\mrm{id}_{\phi}$ does not contribute to this sum; that is, every character appearing in $V_{|U'}$ is of the form $\gamma_{\phi}$ for some nontrivial character $\gamma: \F_p \rar \C^{\times}$. Since $B$-conjugacy acts transitively on the set of such $\gamma_{\phi}$, $V_{|B}$ is irreducible and $V_{|U'}$ is the sum of the $\gamma_{\phi}$ over all nontrivial $\gamma: \F_p \rar \C^{\times}$. 

By the reasoning of \cite[\S 8.1 Proposition 24]{SerreLinReps}, $V_{|B}$ is induced from the stabilizer of any isotypic component of $V_{|U'}$. Since $B/U'$ has cardinality $p-1$, the same as the number of isotypic components of $B/U'$ (on which it acts transitively), $B/U'$ acts freely on these isotypic components, so $V_{|B}$ is induced of any one of the isotypic components of $V_{|U'}$ from $U'$ to $B$. }

\prop[cusp-on-norm-cartan]{Let $V$ be an irreducible cuspidal representation of $G$ attached to a character $\phi$ of $C$. Then \[V_{|N} \simeq W \oplus \bigoplus_{\psi}{\iota_{\psi}},\] where:
\begin{itemize}[noitemsep,label=\tiny$\bullet$]
\item $W$ is a sum of two characters if $\phi$ is even and zero otherwise,
\item If $\phi(-I_2)=1$, $W$ is the sum of the two characters $\chi(\det)\epsilon^r$, where $\chi: \F_p^{\times} \rar \C^{\times}$ is a character such that $\phi(uI_2)=\chi(u^2)$ for each $u \in \F_p^{\times}$ (there are exactly two such $\chi$), and $r \in \{0,1\}$ is such that $(-1)^r\chi(-1)=-1$ if $\phi(uI_2)=\chi(u^2)$ for each $u \in \sqrt{\F_p^{\times}}$, and $(-1)^r\chi(-1)=1$ otherwise (so that for either $\chi$, exactly one value of $r$ works).
\item $\psi$ runs through the morphisms $C \rar \C^{\times}$ such that $\psi \notin \{\phi,\phi^p\}$, $\psi^{p-1} \neq 1$, $\psi(tI_2) = \phi(tI_2)$ for every $t \in \F_p^{\times}$, and we only count one character in every pair $\{\psi,\psi^p\}$. 
\end{itemize}
}

\demo{For any $a \in \F_p^{\times}$, $aI_2$ acts on both sides by multiplication by $\phi(aI_2)$. Next, we show that the representations on both sides have dimension $p-1$. The left-hand side clearly has dimension $p-1$, so we discuss the right-hand side. This amounts to showing that if $\phi(-I_2)=-1$ (resp. $\phi(-I_2)=1$), there are exactly $p-1$ (resp. $p-3$) characters $\psi: C \rar \C^{\times}$ agreeing with $\phi$ on $\F_p^{\times}$ and such that $\psi \notin \{\phi,\phi^p\}$ as well as $\psi^{p-1}=1$. 

Since $C$ is cyclic with order $p^2-1$, there are $p+1$ characters $C \rar \C^{\times}$ agreeing with $\phi$ on $\F_p^{\times}I_2$, and $\phi,\phi^p$ are two of them. So we need to show that if $\phi(-1)=-1$ (resp. $\phi(-1)=1$), no character $\psi: C \rar \C^{\times}$ agreeing with $\phi$ on $\F_p^{\times}I_2$ verifies $\psi^{p-1}=1$ (resp. exactly two such characters satisfy the condition). Let $g \in C$ be a generator, then a character $\psi: C\rar \C^{\times}$ agreeing with $\phi$ on $\F_p^{\times}I_2$ is uniquely defined by the value of $\psi(g) \in \phi(g)\mu_{p+1}$. Given $\omega \in \mu_{p+1}$, the number of $\omega \in \mu_{p+1}$ such that $(\phi(g)\omega)^{p-1}=1$ is the number of $\omega \in \mu_{p+1}$ such that $\omega^2=\phi(g)^{p-1}$. There are either two such $\omega$ or zero, depending on whether $\phi(g)^{(p-1)(p+1)/2}=1$. The correct dimension count follows since $\phi(g)^{(p^2-1)/2}=\phi(-I_2)$. 

If $\phi(-I_2)=1$, the two characters of $C$ agreeing with $\phi$ on $\F_p^{\times}I_2$ and with order dividing $p-1$ are the two possible $\chi(\det)$.  

Now, we need to show that the characters agree at every $g \in C \backslash \F_p^{\times}I_2$ and at every $g \in N \backslash C$. 

When $g \in C \backslash \F_p^{\times}I_2$, by the remark above, the character of the right-hand side is $\sum_{\psi}{\psi(g)}$, where $\psi$ runs through all the characters $\psi: C \rar \C^{\times}$ agreeing with $\phi$ on $\F_p^{\times}I_2$ with $\psi \notin \{\phi,\phi^p\}$. So the character of the right-hand side is $-\phi(g)-\phi^p(g)$, which is also the value at $g$ of the character of the left-hand side by Proposition \ref{characterize-cuspidals}.  

Let $g \in N \backslash C$. Then $g$ has trace zero, hence $g^2=-\det{g}I_2$, and the character of the right-hand side is $0$ if $\phi(-I_2)=-1$, and $\sum_{\chi}{(-1)^{r_{\chi}}\chi(\det{g})}$ over $\chi: \F_p^{\times} \rar \C^{\times}$ such that $\chi(t^2)=\phi(tI_2)$ for each $t \in \F_p^{\times}$, where $r_{\chi}$ is the $r$ attached to $\chi$ in the statement. 

Suppose that $-\det{g}=u^2$ for some $u \in \F_p^{\times}$: then $g$ is similar over $\F_p$ to $\Delta_{u,-u}$, hence the character of the left-hand side vanishes at $g$, while the character of the right-hand side is $\phi(u)\sum_{\chi}{(-1)^{r_{\chi}}\chi(-1)}=0$, because the two characters $\chi$ over which we sum are distinct, so that the $(-1)^{r_{\chi}}\chi(-1)$ are distinct signs.  

Suppose now that $-\det{g}$ is not a square. Then $g$ is similar to some $g' \in C \backslash \F_p^{\times}$, the pair $\{g',(g')^p\}$ being uniquely defined (and $(g')^2$ is scalar). Then $(g')^{2(p-1)}=I_2$ while $(g')^{p-1} \neq I_2$, hence $(g')^{p-1}=-I_2$. By the same arguments as previously, the characters of both sides of the equality agree at $g$ if $\phi(-I_2)=-1$ (they are equal to $-\phi(g')-\phi(g')^p$), so we may assume that $\phi(I_2)=1$. In this case, the character of the left-hand side is $-\phi(g')-\phi(g')^p=-2\phi(g')$, while the character of the right-hand side is $\sum_{\chi}{(-1)^{r_{\chi}}\chi(\det{g})}$, where $\chi: \F_p^{\times} \rar \C^{\times}$ is one of the two characters such that $\chi(t^2)=\phi(tI_2)$ for each $t \in \F_p^{\times}$. By checking both possible, we see that $(-1)^{r_{\chi}}\chi(-1)=-\frac{\phi(g')}{\chi(-\det{g'})}$, so $(-1)^{r_{\chi}}\chi(\det{g'})=-\phi(g')$, whence the claimed equality.
}

\section{Splitting of a certain cuspidal representation restricted to $\SL{\F_p}$}
\label{splitting-sl2}

We now assume that $p \equiv 3 \pmod{4}$, and write $\phi_2: \F_{p^2}^{\times} \rar \C^{\times}$ (resp. $\phi_4: \F_{p^2}^{\times} \rar \C^{\times}$) as the unique character of order $2$ (resp. one of the two characters of order $4$). Thus $\phi_4$ and $\phi_4^3=\phi_4^p=\phi_4^{-1}$ are distinct characters $\F_{p^2}^{\times}/\F_p^{\times} \rar \C^{\times}$ which do not factor through the norm. We denote by $V_4$ the cuspidal representation (with coefficients in $\Q$, by Proposition \ref{cusp-coeffs}) of $G$ attached to the pair of characters $\{\phi_4,\phi_4^{-1}\}$. 

Finally, let $G_2 \triangleleft G$ denote the subgroup $\F_p^{\times}\SL{\F_p}$, and $C \leq G$ a nonsplit Cartan subgroup with normalizer $N$. 

\prop[decomp-cusp-sl]{There exist two irreducible representations $V^+, V^-$ of $G_2$ with dimension $\frac{p-1}{2}$ satisfying the following properties:
\begin{itemize}[noitemsep,label=\tiny$\bullet$]
\item $V^+$ and $V^-$ are contragredients of each other and conjugates of each other under $G/G_2$, 
\item $(V_4)_{|G_2} \simeq V^+ \oplus V^-$,
\item $V_4 \simeq \mrm{Ind}_{G_2}^G{V^+} \simeq \mrm{Ind}_{G_2}^G{V^-}$,
\item $V^+$ and $V^-$ factor through $\PSL{\F_p} \simeq G_2/\F_p^{\times}I_2$,
\item the restrictions of $V^+$ and $V^-$ to any subgroup with order prime to $p$ are isomorphic,  
\item $V^+,V^-$ can be realized over $\Q(\sqrt{-p})$.
\end{itemize}}

\demo{If $g \in G_2$, the trace of the action of $g$ on $V_4$ is, by Proposition \ref{characterize-cuspidals}: 
\begin{enumerate}[label=(\alph*),noitemsep]
\item $p-1$ if $g$ is scalar, 
\item $-1$ if $g$ is not scalar but its minimal polynomial is not separable,
\item $0$ if $g$ has two distinct eigenvalues in $\F_p^{\times}$,
\item $\pm 2$ if $g$ has two distinct eigenvalues in $\F_{p^2} \backslash \F_p^{\times}$.
\end{enumerate}

There are $p-1$ $G$-conjugacy classes in $G_2$ of the first type, and all of them have cardinality one. There are $p-1$ $G$-conjugacy classes in $G_2$ of the second type, and all of them have cardinality $p^2-1$. There are, finally, $\frac{1}{2}\left(\frac{p^2-1}{2}-(p-1)\right)=\frac{(p-1)^2}{4}$ $G$-conjugacy classes in $G_2$ of the last type. Therefore,

\begin{align*}
\frac{2}{p(p-1)^2(p+1)}\sum_{g \in G_2}{|\mrm{Tr}(g \mid V_4)|^2} &= \frac{2((p-1) \cdot (p-1)^2+(p-1)\cdot (p^2-1)+(p-1)^2\cdot p(p-1)}{p(p-1)^2(p+1)}\\
&=\frac{2(p-1)^2(p-1+p+1+p(p-1))}{p(p-1)^2(p+1)}=2,
\end{align*} 
therefore $(V_4)_{|G_2}$ is the sum of two non-isomorphic irreducible representations of $G_2$, which we denote by $V^+,V^-$. Since $\F_p^{\times}I_2$ acts trivially on $V_4$, it also acts trivially on $V^{\pm}$.

Let $\tau^{\pm}$ be the trace of the action of $U$ on $V^{\pm}$. The $\tau^{\pm}$ are algebraic integers and their sum is the trace of the action of $U$ on $V_4$, that is, $\tau^++\tau^-=-1$. Hence $\tau^+ \neq \tau^-$ and $V^+,V^-$ are not isomorphic. Since $V_4$ is irreducible, $G/G_2$ acts transitively on $\{V^+,V^-\}$, so $V^+$ is conjuagte to $V^-$ under $G/G_2$.

If $M \in G_2$ has order prime to $p$, then $M$ is contained (up to $G$-conjugacy) in a group of diagonal matrices or in a nonsplit Cartan subgroup. In particular, the centralizer of $M$ in $G$ has full determinant. Therefore, $M$ is $G_2$-similar to any of its $G$-conjugates, and the characters of $V^+$ and $V^-$ agree at $M$.  

Since $U^{-1}$ is $G$-conjugate to $U$ but not $G_2$-conjugate to $U$, $\tau^+$ and $\tau^-$ are complex conjugates of each other, so $\tau^+ \notin \R$. Therefore, neither $V^+$ not $V^-$ is self-dual. Since $(V_4)_{|G_2} \simeq V^+\oplus V^-$ is self-dual by Corollary \ref{cuspidal-duals}, $V^+$ and $V^-$ are duals of each other. In particular, they both have dimension $\frac{p-1}{2}$. 

Let $\epsilon \in \{\pm 1\}$, then there is a natural map $\mrm{Ind}_{G_2}^G{V^{\epsilon}} \rar V_4$. It is nonzero and $G$-equivariant by construction. Since $V_4$ is irreducible, it is surjective and between representations of the same dimension, hence is an isomorphism. 

All we need to do is check that $V^{\pm}$ can be defined over $\Q(\sqrt{-p})$. 

Since the eigenvalues of the action of $U$ on $V_4$ are exactly the primitive $p$-th roots of unity, there exists a partition of the primitive $p$-th roots of unity into two subsets $E^{\pm}$, where $E^{\epsilon}$ is the set of eigenvalues of the action of $U$ on $V^{\epsilon}$. For any $t \in \F_p^{\times}$, $U^{t^2}=\Delta_{t,1/t}U\Delta_{1/t,t}$, so each $E^{\pm}$ is stable after taking $t^2$-th powers for any $t \in \F_p^{\times}$. As a consequence, 
\[\{E^{\epsilon}\mid\epsilon \in \{\pm 1\}\} = \{\{e^{2i\pi \epsilon t^2/p,\, t \in \F_p^{\times}}\}\mid\epsilon \in \{\pm 1\}\}.\]

Therefore, there exists a sign $a$ such that, on any complex realization of $V^{\epsilon}$, \[\sum_{t \in \F_p^{\times}}{\left(\frac{t}{p}\right)U^t} = a\epsilon\sum_{t \in \F_p^{\times}}{\left(\frac{t}{p}\right)e^{2i\pi t/p}}:=a\epsilon\mathfrak{g}_p.\] 

Let $V$ be a realization of the cuspidal representation $V_4$ over $\Q$, and $\gamma$ be the $\Q$-endomorphism $\sum_{t \in \F_p^{\times}}{\left(\frac{t}{p}\right)U^t}$. Since $V \otimes_{\Q} \C$ is the sum of one copy of $V^+$ and one copy of $V^-$ (as $\C[G_2]$-modules), $\gamma$ acts on the copy of $V^{\epsilon}$ by multiplication by $a\epsilon\mathfrak{g}_p$. In particular, $\gamma \otimes \C$ is a $G_2$-automorphism of $V \otimes_{\Q}\C$ with square $\mathfrak{g}_p^2=-p$ by Proposition \ref{gauss-sum}, hence $\gamma$ is a $\Q[G_2]$-automorphism of $V$ with square $-p$. We can thus view $V$ as a representation $V'$ of $G_2$ with coefficients in $\Q[\gamma] \simeq \Q(\sqrt{-p})$. By construction, $V' \otimes_{\Q[\gamma]} \C$ is exactly a complex realization of $V^{b}$ for the sign $b$ such that $v\epsilon\mathfrak{g}_p$ is the image of $\gamma$ in $\C$.
}

\cor[character-decomp-cusp]{With the notations of Proposition \ref{decomp-cusp-sl}, for any $M \in G_2$ with order prime to $p$ inducing endomorphisms $g^+,g^-$ of the spaces of the representations $V^+,V^-$, $g^+$ and $g^-$ have the same trace. It is equal to $\frac{p-1}{2}$ if $M$ is scalar, $0$ if $M$ has two distinct eigenvalues in $\F_p$, $-1$ if $M$ is $G$-similar to some element in $4C \backslash \F_p^{\times}$, $1$ if $M$ is $G$-similar to some $M' \in 2C \backslash 4C$. Moreover, after possibly exchanging $V^+,V^-$, the eigenvalues of $U$ acting on $V^{\epsilon}$ are exactly the $e^{2i\pi \epsilon t/p}$ for $t \in \F_p^{\times 2}$. }

\demo{The first and last parts of the statement are consequences of the proof of Proposition \ref{decomp-cusp-sl}. The common trace of $g^+$ and $g^-$ is one half the character of $V$ at $g$. Now, if $M \in 2C \backslash 4C$, then $\phi_4(M)=\phi_4^3(M)=-1$, while, if $M \in 4C \backslash \F_p^{\times}I_2$, $\phi_4(M)=\phi_4(M)^{-1}=1$, so the claims follow from Proposition \ref{characterize-cuspidals}. }

In the rest of this section, we break the symmetry between $V^+$ and $V^-$ by stating that $V^+$ is the representation (over $\C$) such that the complex eigenvalues of $U$ acting on $V^+$ are exactly the $e^{2i\pi t/p}$ for $t \in \F_p^{\times 2}$.  

\prop[cusp-sl-borel]{Let $B_2=B \cap G_2$ and $\epsilon=\pm$ be a sign. Then $C^{\epsilon}_{|B_2}$ is irreducible and \[C^{\epsilon}_{|B_2} \simeq \mrm{Ind}_{\F_p^{\times}U^{\Z}}^{B_2}{\left[aU^b\longmapsto e^{2i\pi b\epsilon/p}\right]}.\]}

\demo{$(V_4)_{|B}=\mrm{Ind}_{B_2}^B{C^{\epsilon}}$ is irreducible, so $C^{\epsilon}_{|B_2}$ is irreducible. Since $C^-_{|B_2}$ is the conjugate of $C^+_{|B_2}$ under $B/B_2$, we may assume $\epsilon=+$. Both sides have the same dimension, and that any $aI_2$ for $a \in \F_p^{\times}$ acts trivially on either side. The trace of the action of any $g \in B_2 \backslash \F_p^{\times}U^{\Z}$ on the left-hand side vanishes; since such a $g$ is not conjugate to $\F_p^{\times}U^{\Z}$, the trace of the action of $g$ on the right-hand side also vanishes by \cite[\S 7.2 Proposition 20]{SerreLinReps}. 

So it is enough to check that the restrictions of the two representations are isomorphic as $\C[U]$-modules. The eigenvalues of $U$ acting on the left-hand side are exactly the $e^{2i\pi t/p}$ for $t \in \F_p^{\times 2}$, while the eigenvalues of $U$ acting on the right-hand side are the conjugates of the characters $U^b \longmapsto e^{2i\pi b/p}$ by $B_2/\F_p^{\times}U$. The conjugates of these characters are exactly the $e^{2i\pi bt/p}$ for $t \in \F_p^{\times 2}$. }

\prop[cusp-sl-diag]{Let $D_2 \leq G_2$ be the subgroup of diagonal matrices, and $\alpha,\beta: D_2 \rar \F_p^{\times}$ denote the first and second coefficients on the diagonal. Then \[C^{\pm}_{|D_2} \simeq \bigoplus_{\chi: \F_p^{\times 2} \rar \C^{\times}}{\chi\left(\frac{\alpha}{\beta}\right)}.\] }

\demo{We can check that the characters coincide: for any non-scalar $d \in D_2$, $d$ acts with trace zero on either side. For any $a \in \F_p^{\times}$, $aI_2$ acts trivially on either side, and both sides are representations of $D_2$ of dimension $\frac{p-1}{2}$. }

\prop[cusp-sl-norm-cartan]{Let $C_2 = G_2 \cap C$ and $N_2 = G_2 \cap N$. Let $\eta: N_2 \rar N_2/C_2 \rar \{\pm 1\}$ be the character which is trivial if and only if $p \equiv 3\pmod{8}$. Then one has 
\[V^{\pm}_{|N_2} \simeq \eta \oplus \bigoplus'_{\chi: C_2 \rar \C^{\times}}{\mrm{Ind}_{N_2}^{D_2}\chi},\]
where the sum is over characters $\chi: C_2/\F_p^{\times} \rar \C^{\times}$ with order at least $3$ and such that we count exactly one character in any pair $\{\chi,\chi^p\}$. }

\demo{The group $C_2/\F_p^{\times}$ is cyclic with even order $\frac{p+1}{2}$, so the right-hand side has dimension $1+\left(\frac{p+1}{2}-2\right) = \frac{p-1}{2}$, like the left-hand side. Any scalar matrix acts trivially on either side. 

If $g \in N_2 \backslash C_2$, then the trace of the action of $g$ on the right-hand side is $1$ if $p \equiv 3\pmod{8}$ and $-1$ otherwise. Since $\mrm{Tr}(g)=0$, one has $g^2=-\det{g}I_2$. Since $p \equiv 3\pmod{4}$ and $\det{g} \in \F_p^{\times 2}$, $g$ has a $G$-conjugate $g' \in C_2$, so that $(g')^2=-\det{g'}I_2 \in \F_p^{\times}I_2 \backslash \F_p^{\times 2}I_2$. Since $\frac{p-1}{2}$ is odd, the $2$-adic valuation of the order of $g'$ is exactly $4$, so that $g'$ is not scalar, a square in $C$, and it is a fourth power in $C$ if and only if $p \equiv -1\pmod{8}$. Now, the trace of the action of $g'$ (hence of $g$) on $V^{\pm}$ is $1$ if $g'$ is not a fourth power, and $-1$ otherwise.  

Let $g \in C_2 \backslash \F_p^{\times}I_2$. The trace of the action of $g$ on the right-hand side is $1+\sum_{\chi}{\chi(g)}$, where $\chi$ runs over characters $C_2/\F_p^{\times} \rar \C^{\times}$ that have order at least $3$. Thus said trace is exactly $-\chi_2(g)$, where $\chi_2$ is the unique character of $C_2/\F_p^{\times}$ with order two. Since $\chi_2$ agrees with $\phi_4$ (after identifying $C$ to $\F_{p^2}^{\times}$), this matches the trace of the left-hand side.}

\section{Representations of $G$ in characteristic $p$}
\label{reps-gl2-carac-p}

In this section, the ring of coefficients of our representations will always be specified and is a $\F_p$-algebra. 

The tautological representation of $G$ (acting on $\F_p^2$, whose elements are identified with column matrices) is denoted by $R$, and we write $x=(1,0),y=(0,1)$. Let $X: \F_p^{\times} \rar \F_p^{\times}$ be the tautological character.

For any $k \geq 0$, $R^{\otimes k}$ is a $\F_p[\mathfrak{S}_k \times G]$-module (where $G$ acts simultaneously on each factor and $\mathfrak{S}_k$ permutes these factors). We let $\mrm{Sym}^kR$ (resp. $\mrm{TSym}^kR$) be the co-invariants (resp. its invariants) under the action of $\mathfrak{S}_k$. Note that the natural map $R^{\otimes k} \rar \mrm{Sym}^kR$ has a section whenever $k < p$.

The following well-known result is nonetheless important.

\prop{Let $V$ be a nonzero representation of $G$ on any extension $k$ of $\F_p$. Then $V^U$ is nonzero.}

\demo{Since $U^p=I_2$, the endomorphism $u$ of $V$ induced by $U$ also satisfies $u^p=\mrm{id}$. Thus $(u-\mrm{id})^p=0$, hence $u-\mrm{id}$ is not injective.}

\cor{Let $V$ be a representation of $G$ over some extension $k$ of $\F_p$. Then, if $\dim_k{V^U}=1$ and the sub-$k[G]$-module generated by $V^U$ is $V$, then $V$ is absolutely irreducible.}

\demo{If $k'/k$ is any field extension, then $V \otimes_{k'}$ (as a representation of $G$ over $k'$) satisfies the assumptions, so it is enough to show that $V$ is irreducible. So let $W \subset V$ be a nonzero irreducible subrepresentation. So $\{0\} \neq W^U \subset V^U$: since $V^U$ is a line, $W^U=V^U$ and $V$ is generated over $G$ by $W^U$, hence $V \subset W$ thus $V=W$.}

\cor[find-irreducibles]{Let $0 \leq k < p$ and $r \in \Z/(p-1)\Z$. Then $S_{k,r}:=\det^r \otimes \mrm{Sym}^kR$ is an absolutely irreducible representation of $G$ over any field extension $F$ of $\F_p$. Moreover, the $S_{k,r}$ are pairwise non isomorphic.}

\demo{Let $F \supset \F_p$ be any field. It is enough to show that the $S_{k,r}$ are irreducible and pairwise non-isomorphic over $F$. 
For irreducibility, we may assume that $r=0$. Since $Ux=x,Uy=x+y$, for any $0 \leq a,b \leq k$ with $a+b=k$, $U \cdot (x^ay^b) \in x^ay^b+bx^{a+1}y^{b-1}+\sum_{k \geq a' > a+1}{F \cdot x^{a'}y^{k-a'}}$. Thus, the matrix of the action of $U-1$ on $\mrm{Sym}^kR$, in the basis $(x^k, x^{k-1}y,x^{k-2}y^2,\ldots,xy^{k-1},y^k)$ is upper triangular, its diagonal is zero and every coefficient just above the main diagonal is one. Thus the kernel of the action of $U-1$ on $S_{k,0}$ is $Fx^k$. 

If $k=a+b$ with $b > 0$, one can directly compute that \[x^ay^b=-\binom{k}{b}^{-1}\sum_{i \in \F_p^{\times}}{i^{p-1-b}(x+iy)^k}=-\binom{k}{b}^{-1}\sum_{i \in \F_p^{\times}}{i^{p-1-b}(T^i \cdot x^k)},\] so the sub-$G$-module generated by $Fx^k$ is $S_{k,0}$, whence the irreducibility. 

The isomorphism class of $\det^r\otimes \mrm{Sym}^kR$ obviously determines $k$ (using its dimension). Thus it is enough to show that if $0 \leq r < p-1$ and $0 \leq k < p$, $\mrm{Sym}^kR$ is not isomorphic to any of its twists by $\det^r$. 

By the above argument, the action of $\Delta_{u,1}$ on $\ker\left(U-1\mid S_{k,r}\right)=Fx^k$ is the multiplication by $u^{k+r}$. In particular, the isomorphism class of $S_{k,r}$ over $F$ determines $k$ and $k+r \pmod{p-1}$, hence determines $k$ and $r \pmod{p-1}$, whence the conclusion. }

\lem[carac-p-dual-symmetric]{Let $0 \leq k < p$ and $r \in \Z/(p-1)\Z$. There is a perfect $G$-bilinear pairing $S_{k,r} \times S_{-k,-k-r} \rar \F_p$ (where the representations are defined over $\F_p$) mapping $(x^k,y^k)$ to $1$. In particular, $S_{k,r}^{\ast} \simeq S_{k,-k-r}$.}

\demo{We may obviously assume that $k > 0$ and $r=0$. The following linear map 
\[(u_1 \otimes \ldots \otimes u_k) \otimes (v_1 \otimes \ldots \otimes v_k) \in R^{\otimes k} \otimes R^{\otimes k}\longmapsto \prod_{i=1}^k{\det(u_i,v_i)} \in \F_p\]
is clearly invariant under the action of $\mathfrak{S}_k$, so factors as a map $\mu: S_{k,0} \otimes S_{k,0} \rar \F_p$. By construction, $\mu(x^k,y^k)=1$. For $g \in G$, $\mu(g\alpha,g\beta) = \det{g}^{k}\mu(\alpha,\beta)$, so that, if we let $G$ act on $\F_p$ as the character $S_{0,k}$, $\mu$ is a $G$-bilinear pairing $S_{k,0} \otimes S_{k,0} \rar S_{0,k}$, ie $\mu$ defines a $G$-bilinear pairing $S_{k,0} \times S_{k,-k} \rar S_{0,0}=\F_p$. This pairing is then perfect because $S_{k,0},S_{k,-k}$ are both irreducible and $\mu(x^k,y^k)=1$. }

\prop[carac-p-break-principal-series]{Let $0 \leq a \leq p-1$. There is a natural surjection of $\F_p[G]$-modules $\pi'_a: \pi(1,X^a) \rar S_{p-1-a,a}$ mapping $\F_p(1,X^a)$ to $\F_px^{p-1-a}$. }

\demo{For any $M=\begin{pmatrix}b & \ast\\0 & d\end{pmatrix} \in B$, one has $Mx=bx$. Therefore, one has \[M \cdot x^{p-1-a} = (bd)^ab^{p-1-a}x^{p-1-a}=d^ax^{p-1-a}.\] Thus, if $g$ is a generator for $R(1,X^a)$, $g \longmapsto x^{p-1-a} \in S_{p-1-a,a}$ is a morphism of $\F_p[B]$-modules, so it induces a morphism $\pi(1,X^a) \rar S_{p-1-a,a}$. The morphism is onto because it is nonzero and $S_{p-1-a,a}$ is irreducible.}

\prop{In Proposition \ref{carac-p-break-principal-series}, $\ker{\pi'_a}$ is isomorphic to $S_{a,0}$. If $a \notin \{0,p-1\}$, then the exact sequence $0 \rar S_{a,0} \rar \pi(1,X^a) \rar S_{p-1-a,a} \rar 0$ of $\F_p[G]$-modules does not split.}

\demo{First, we show that $\ker{\pi'_a}$ is isomorphic to $S_{a,0}$. By Proposition \ref{carac-p-break-principal-series}, there is a surjective map $f: \pi(1,X^{p-1-a}) \rar S_{a,-a}$ of $\F_p[G]$-modules. By Lemma \ref{principal-series-dual} and Lemma \ref{carac-p-dual-symmetric}, the contragredient of $f$ is an injective map $f':S_{a,0} \rar \pi(1,X^a)$. Since $S_{a,0}$ and $S_{p-1-a,a}$ are irreducible and pairwise non-isomorphic, $\pi'_a \circ f': S_{a,0} \rar S_{p-1-a,a}$ is the zero map and $\im{f'} \subset \ker{\pi'_a}$. Since $\dim{\ker{\pi'_a}}=p+1-(p-1-a+1)=a+1=\dim{S_{a,0}}=\dim{\im{f'}}$, $\ker{\pi'_a}$ is exactly the image of $f'$, which is isomorphic to $S_{a,0}$.

When $a \notin \{0,p-1\}$, then, by Proposition \ref{principal-series-irreducible}, the exact sequence cannot split since $\mrm{End}_{\F_p[G]}(\pi(1,X^a))$ is a line. 
}

\rem{The kernel of the map $\pi(1,1) \rar S_{p-1,0}$ is isomorphic to $S_{0,0}$, so it is the unique line of $\pi(1,1)$ fixed by $G$. By Proposition \ref{steinberg-rep-definition}, the exact sequence of the previous definition is exactly $0 \rar \F_p \rar \pi(1,1) \rar \mrm{St} \rar 0$; in particular, over $\F_p$, $\mrm{St}$ is isomorphic to $\mrm{Sym}^{p-1}R$. }

\cor[decomp-principal-series]{Let $0 \leq a \leq b \leq p-1$, there are exact sequences 
\begin{align*}
0 \rar S_{b-a,a} \rar \pi(X^a,X^b) \rar S_{p-1-(b-a),b} \rar 0,\\
0 \rar S_{p-1-(b-a),b} \rar \pi(X^b, X^a) \rar S_{b-a,a} \rar 0,
\end{align*} which split if and only if $(a,b)=(0,p-1)$ or $a=b$. }

\prop[all-irreducibles]{Let $F$ be any field extension of $\F_p$ and $V$ be an irreducible representation of $G$ over $F$. Then $V$ is isomorphic to some $S_{k,r}$.}

\demo{The $F$-vector space $V^U$ is endowed with an action of the abelian group $B/U \simeq \F_p^{\times 2}$. It has exponent $p-1$ prime to $p$ and the polynomial $X^{p-1}-1$ splits in $F$ with simple roots. So $V^U$ is the direct sum of characters of $B/U$. Let $x \in V^U \backslash \{0\}$, $\alpha,\beta: \F_p^{\times} \rar \F_p^{\times}$ be such that for any $M=\begin{pmatrix}u & \ast\\0 & v\end{pmatrix} \in B$, $Mx=\alpha(u)\beta(v)x$. Write $\alpha=X^a,\beta=X^b$, then mapping $F(X^a,X^b)$ to $Fx$ defines a nonzero morphism $\pi(X^a,X^b) \rar V$. Since $V$ is irreducible, this morphism is onto. 

We know that there is an exact sequence $0 \rar S_{q,r} \rar \pi(X^a,X^b) \rar S_{s,t} \rar 0$. Since $S_{q,r} \rar \pi(X^a,X^b) \rar V$ is a morphism between two irreducible representations of $G$ over $F$, it is either an isomorphism or zero. If it is an isomorphism, we are done; otherwise, the map $\pi(X^a,X^b) \rar V$ factors as a map $\pi(X^a,X^b) \rar S_{s,t} \rar V$, so that $S_{s,t} \rar V$ is surjective between irreducible representations, thus is an isomorphism. 
}

\newpage

%% file: local-0.tex
\chapter{Modular forms and local constants}
\label{appendix-local-const}

The goal of this appendix is to recall various classical results about local constants for Galois representations and modular forms. These results are mostly contained in \cite{AL78}, \cite{Li-RS} for twists of modular forms and their local constants, as well as \cite{BH} and \cite{NTB} for more general results about local constants. The point is to be able to reuse these results more systematically, thus more conveniently.

\section{Gauss sums}
\label{appendix-gauss-sums}

We recall here a few classical results on Gauss sums that we will end up needing.

\defi{Let $\chi: (\Z/N\Z)^{\times}\rar \C^{\times}$ be a Dirichlet character. Then its \emph{Gauss sum} is 
\[\mathfrak{g}(\chi) = \sum_{x \in (\Z/N\Z)^{\times}}{\chi(x)e^{2i\pi x/N}}.\]
We will usually assume that $N$ is the conductor of $\chi$, but this need not be always the case.}

\lem[gauss-module]{Let $\chi: (\Z/N\Z)^{\times} \rar \C^{\times}$ be a primitive Dirichlet character. Then \[|\mathfrak{g}(\chi)|^2=\mathfrak{g}(\chi)\mathfrak{g}(\overline{\chi})\chi(-1)=N.\]}

\demo{The first equality is a direct computation, and the second equality is well-known, see for instance \cite[Proposition VII.2.6]{Neukirch-ANT}. }

The following result is well-known and due to Gauss, and a proof is given in \cite[\S 6.4]{IR}. 

\prop[gauss-sum]{Let $p$ be an odd prime and $\lambda: \F_p^{\times} \rar \{\pm 1\}$ be the unique nontrivial quadratic character. Then $\mathfrak{g}(\lambda)=u\sqrt{p}$, where $u=1$ if $p \equiv 1 \pmod{4}$ and $u=3$ if $p \equiv 3 \pmod{4}$. }

\lem[gauss-deligne]{Let $q$ be a prime number, $Q \geq q^2$ be a power of $q$ and $\alpha: (\Z/Q\Z)^{\times} \rar \C^{\times}$ be a primitive Dirichlet character. Then, there exists $z \in (\Z/Q\Z)^{\times}$ such that for any Dirichlet character $\chi: (\Z/Q\Z)^{\times} \rar \C^{\times}$ with conductor $C \leq \sqrt{Q}$, $\mathfrak{g}(\alpha\chi)=\chi(z)\mathfrak{g}(\alpha)$.}

\demo{This can be seen as a consequence of, say, \cite[\S 23.6 Proposition]{BH}, but let us give an explicit proof. In the situation of the statement, $C \mid Q$, so $C$ is a power of $q$ and $C^2 \mid Q$. 

Write $Q=AB$ with $qB \geq A \geq B \geq q$ (if $Q=q^{2n+1}$, then $A=qB=q^{n+1}$, otherwise $A=B\sqrt{Q}$). Then the function $n \in \Z/B\Z \longmapsto \alpha(1+An) \in \C^{\times}$ is well-defined and is a group homomorphism, so it is of the form $n \longmapsto e^{2i\pi tn/B}$ for some $t \in \Z/B\Z$. 

Let $\chi: (\Z/Q\Z)^{\times} \rar \C^{\times}$ be a Dirichlet character with conductor at most (hence dividing) $B$. Then

\begin{align*}
\mathfrak{g}(\alpha\chi) &= \sum_{\substack{1 \leq k \leq A\\(k,q)=1}}{\sum_{0 \leq l < B}{(\alpha\chi)(k(1+Al))e^{2i\pi k(1+Al)/Q}}} = \sum_{\substack{1 \leq k \leq A\\ (k,q)=1}}{\alpha(k)\chi(k)e^{2i\pi k/Q}\sum_{0 \leq l < B}{e^{2i\pi (kl+tl)/B}}}\\
&= B\sum_{\substack{1 \leq k \leq A \\ k \equiv -t \pmod{B}}}{\alpha(k)\chi(k)e^{2i\pi k/Q}} = B\chi(-t)\sum_{\substack{1 \leq k \leq A\\ k \equiv -t\pmod{B}}}{\alpha(k)e^{2i\pi k/Q}},
\end{align*} so $z=-t$ satisfies the conditions.}

\lem[local-epsilon-fp2]{Let $p$ be an odd prime and $\psi: \F_{p^2}^{\times}/\F_p^{\times} \rar \C^{\times}$ be a nontrivial character. Let $a \in \F_p^{\times} \backslash \F_p^{\times 2}$, then $\psi(\sqrt{a}) \in \{\pm 1\}$ is well-defined and does not depend on the choice of $a$. Moreover, if $\tau: \F_{p^2} \rar \F_p$ is the trace, one has 
\[\sum_{x \in \F_{p^2}^{\times}}{\psi(x)e^{2i\pi\tau(x)/p}}=p\psi(\sqrt{a}).\]}

\demo{The first part of the claim is clear, since $\sqrt{a}$ is defined (allowing for different $a$) up to multiplication by $\F_p^{\times}$, and $\psi(\sqrt{a})^2=\psi(a)=1$. Write $k=\F_{p^2}$ for short. Since $\psi$ is a nontrivial character of $k^{\times}/\F_p^{\times}$, the sum of its values vanishes. Since $\tau$ is $\F_p$-linear and onto, its kernel is a $\F_p$-line. It is clear that $\sqrt{a} \in \ker{\tau}$, so $\ker{\tau}=\F_p\sqrt{a}$. 

\begin{align*}
\sum_{x \in k^{\times}}{\psi(x)e^{2i\pi\tau(x)/p}} &= \sum_{x \in k^{\times}/\F_p^{\times}}{\psi(x)\sum_{t \in \F_p^{\times}}{e^{2i\pi t\tau(x)/p}}}= \sum_{x \in k^{\times}/\F_p^{\times}}{\psi(x)(p\mathbf{1}_{\tau(x)=0}-1)}\\
&= p\sum_{\substack{x \in k^{\times}/\F_p^{\times}\\ \tau(x)=0}}{\psi(x)} = p\psi(\sqrt{a}).
\end{align*}
}

\section{Twists and local $L$-factors of newforms}

\prop[bad-L-factor]{
Let $k,N \geq 1$, $f \in \mathcal{S}_k(\Gamma_1(N))$ a normalized newform with primitive character $\chi$ of conductor $M$. Let $q \mid N$ be a prime, $Q=q^{v_q(N)}$ and $N'=N/Q$. Then
\begin{itemize}[label=$-$,noitemsep]
\item If $q=Q$ is coprime to $M$, then $a_q(f)^2=\chi(q)q^{k-2}$.
\item If $N/M$ is coprime to $q$, then $|a_q(f)|^2=q^{(k-1)/2}$.
\item Otherwise, $a_q(f)=0$.
\end{itemize}
}

\demo{This is well-known, see for instance \cite[Theorem 4.6.17]{Miyake}.}

\prop[twist-exists]{
Let $f \in \mathcal{S}_k(\Gamma_1(N))$ be a normalized newform with primitive character $\chi$ of conductor $C$. Let $\psi: (\Z/D\Z)^{\times} \rar \C^{\times}$ be a Dirichlet character of conductor $M \mid D$. There exists a unique normalized newform $f \otimes \psi \in \mathcal{S}_k(\Gamma_1(L))$ such that, for any prime $p\nmid ND$, $a_p(f \otimes \psi)=a_p(f)\psi(p)$. Moreover, 
\begin{enumerate}[label=(\alph*),noitemsep]
\item \label{only-primitive-part} The newform $f \otimes \psi$ only depends on the primitive character attached to $\psi$. 
\item \label{twist-character} The primitive character of $f \otimes \psi$ is the primitive character attached to $\chi\psi^2$.
\item \label{level-twist-basic} $L$ divides the least common multiple $P$ of $N,CM,M^2$.
\item \label{level-twist-coprime} For any prime number $q \nmid M$, $v_q(L)=v_q(N)$.
\item \label{level-twist-local} Let $q$ be a prime number, and write $\psi=\psi_q\psi'$, where $\psi'$ has conductor prime to $q$ and the conductor of $\psi_q$ is a power of $q$. Then $v_q(L)$ only depends on $f$ and $\psi_q$.
\item \label{level-small-twist} Let $q \mid M$ be a prime such that $v_q(P)=v_q(N)$. If moreover $v_q(L)<v_q(N)$, then $\max(v_q(CM),v_q(M^2)) = v_q(N)$.
\end{enumerate}}

\demo{The existence of $f \otimes \psi$ and \ref{only-primitive-part} \ref{twist-character}, \ref{level-twist-basic} are a consequence of \cite[Proposition 3.1]{AL78} and the strong multiplicity one theorem \cite[Theorem 4.6.19]{Miyake}. 

If $q \nmid M$ is a prime. By \ref{level-twist-basic}, one has $v_q(L) \leq v_q(N)$. On the other hand, since $f = (f \otimes \psi) \otimes \overline{\psi}$ by strong multiplicity one, $N$ divides the least common multiple of $L,C'M,M^2$, where $C'$ is the conductor of the character of $f \otimes \psi$. Since $v_q(C') \leq v_q(N)$ and $v_q(M)=0$, $v_q(L) \leq v_q(N)$, so \ref{level-twist-coprime} holds, and \ref{level-twist-local} follows. 

Let $q \mid M$ be a prime such that $v_q(CM)$ and $v_q(M^2)$ are both lower than $v_q(N)$, so that $v_q(P)=v_q(N)$. Let $C'$ be the conductor of $\chi\psi^2$. Since $f = (f\otimes \psi) \otimes \overline{\psi}$, one has $v_q(N) \leq \max(v_q(L),v_q(C'M),v_q(M^2))$ by \ref{level-twist-basic}. Since one has $v_q(C') \leq \max(v_q(C),v_q(M))$, the assumptions imply that $v_q(N) \leq v_q(L)$. This proves the contrapositive of \ref{level-small-twist}.   }

\rem{\ref{level-twist-basic}, \ref{level-twist-coprime} and \ref{level-small-twist} may look like \cite[Theorem 3.1]{AL78}, but they are not exactly the same: \emph{loc.cit.} discusses the modular form $\sum_{n \geq 1}{a_n(f)\psi(n)q^n}$, which need not be a newform. }

\section{Local constants and twists of newforms}

Most of the statements of this section come from \cite{AL78}.

Let $k,N \geq 1$. Let $Q \mid N$ be a positive number, prime to $N':=N/Q$. Let $x,y,z,w \in \Z$ be such that $Q^2xw-Nyz=Q,\,y \equiv 1\pmod{Q},x \equiv{1} \pmod{N'}$. Then the map 

\[W_Q: f \in \mathcal{S}_k(\Gamma_1(N)) \longmapsto f\mid_k\begin{pmatrix}Qx & y\\Nz & Qw\end{pmatrix}\in \mathcal{S}_k(\Gamma_1(N)),\] 

does not depend on the choice of $x,y,z,w$ by \cite[Proposition 1.1]{AL78}. 

\prop[local-constant-elementary]{Let $k,N,Q,N'$ be as above. Let $f \in \mathcal{S}_k(\Gamma_1(N))$ be a normalized newform with primitive character $\chi$, and write $\chi=\chi_Q\chi'_Q$, where the conductor of $\chi_Q$ (resp. $\chi'_Q$) divides $Q$ (resp. $N'$). Then there exists a constant $\lambda_Q(f) \in \C^{\times}$ and a normalized newform $g \in \mathcal{S}_k(\Gamma_1(N))$ such that:
\begin{itemize}[label=\tiny$\bullet$,noitemsep]
\item $W_Q(f)=\lambda_Q(f)g$,
\item $\lambda_Q(f)$ is an algebraic number of complex modulus one.
\item $g$ is $f \otimes \overline{\chi'_Q}$.
\end{itemize}
Moreover, one has $\lambda_Q(f)\lambda_Q(g)=\chi_Q(-1)\overline{\chi'_Q}(Q)$.}

\demo{This is proved in \cite[\S 1]{AL78}.}

\prop[functional-eqn]{Let $f \in \mathcal{S}_k(\Gamma_1(N))$ be a normalized newform. The associated $L$-series $L(f,s)=\sum_{n \geq 1}{\frac{a_n(f)}{n^s}}$ converges absolutely when $\mrm{Re}(s) > k/2+1$, and $s \longmapsto L(f,s)$ extends to an entire function (also denotes $L(f,s)$). The completed $L$-function $\Lambda(f,s) = N^{s/2}(2\pi)^{-s}\Gamma(s)L(f,s)$ is entire and the following functional equation holds: $\Lambda(f,s)=(-i)^k\lambda_N(f)\Lambda(\overline{f},k-s)$. }

\demo{This is a classical computation, see for instance \cite[Section 5.10]{DS} (but note that $\lambda_N(f)\overline{f}$ is $f\mid\begin{pmatrix}0 & 1\\-N & 0\end{pmatrix}$ rather than $f \mid \begin{pmatrix}0 & -1\\N & 0\end{pmatrix}$).} 

The main case of interest is the situation where $Q$ is a power of a prime $q \mid N$. In this case, we write $W_q$ and $\lambda_q$ for $W_Q$ and $\lambda_Q$ respectively. 

\prop[local-constant-principal]{In the situation of Proposition \ref{local-constant-elementary}, assume furthermore that $Q$ is the power of some prime $q \mid N$ such that $a_q(f) \neq 0$. Then $\lambda_Q(f)=Q^{k/2-1}\frac{\mathfrak{g}(\chi_q)}{a_Q(f)}$, where $\chi_q$ is viewed as a Dirichlet character modulo its conductor if it is not trivial, and modulo $q$ otherwise.}

\demo{This is \cite[Theorem 2.1]{AL78}.}

\prop[local-constant-coprime-twist]{Let $k,N,Q,N',f,\chi$ as in Proposition \ref{local-constant-elementary}. Let $\psi$ be a primitive Dirichlet character of conductor $M$ prime to $Q$. Then $\lambda_Q(f \otimes \psi)=\overline{\psi}(Q)\lambda_Q(f)$.}

\demo{This can be deduced from \cite{AL78} as follows. If $D$ is an integer prime to $Q$, the automorphism $W_Q$ of $\mathcal{S}_k(\Gamma_1(ND))$ extends the automorphism $W_Q$ of $\mathcal{S}_k(\Gamma_1(N))$ (for the inclusion $\mathcal{S}_k(\Gamma_1(N)) \subset \mathcal{S}_k(\Gamma_1(ND))$). 

Let $C$ be the conductor of $\chi$ and $P$ be the least common multiple of $N,CM,M^2$, so that $P/N$ is prime to $Q$, By \cite[Proposition 3.1]{AL78}, there exists a twisting operator
\[S_{\psi}: \sum_{n \geq 1}{a_nq^n} \in\mathcal{S}_k(\Gamma_1(N),\chi) \longmapsto \sum_{n \geq 1}{a_n\psi(n)q^n} \in \mathcal{S}_k(\Gamma_1(P),\chi\psi^2).\]

By \cite[Proposition 3.4]{AL78}, one has $\overline{\psi}(Q)S_{\psi}W_Q=W_QS_{\psi}$. Moreover, by \cite[Proposition 3.2]{AL78}, we can write $(S_{\psi}f)(\tau) = g(\tau)+\sum_{t > 1}{\alpha_tg(t\tau)}$, where $g=f \otimes \psi$ has conductor $L \mid P$, $t$ runs through divisors of $P/L$ and the $\alpha_t$ are complex numbers. By Theorem \ref{twist-exists}, $Q$ is coprime to $P/L$, so, by \cite[Proposition 1.5]{AL78}, $W_Qg=W_Q(S_{\psi}f) - \sum_t{\alpha'_t(W_Qg)(t\tau)}$, and therefore
\begin{align*}
\lambda_Q(f\otimes \psi)&=\lambda_Q(g)=a_1(W_Q g)=a_1(W_Q S_{\psi}f)=\overline{\psi}(Q)a_1(S_{\psi}W_Qf) \\
&= \overline{\psi}(Q)a_1(W_Q(f))=\overline{\psi}(Q)\lambda_Q(f).
\end{align*}  }

\cor[local-to-global-constant]{Let $f \in \mathcal{S}_k(\Gamma_1(N))$ be a normalized newform with character $\chi$. For each prime $q\mid N$, write $N=N_qN'_q$ with $N_q=q^{v_q(N)}$, and $\chi=\chi_q\chi'_q$, where $\chi_q,\chi'_q$ are characters modulo $N_q$ and $N'_q$ respectively. Then $\lambda_N(f)=\prod_{q \mid N}{\lambda_q(f)\chi'_q(N_q)}$.}

\demo{Let $Q_1,Q_2 \geq 1$ be coprime divisors of $N$ such that $Q_i$ is prime to $N/Q_i$. Then, by \cite[Proposition 1.4]{AL78}, one has \[\overline{\chi_{Q_2}}(Q_1)\lambda_{Q_1Q_2}(f)=a_1(W_{Q_2}(W_{Q_1}f))=\lambda_{Q_1}(f)\lambda_{Q_2}(f \otimes \overline{\chi_{Q_1}})=\chi_{Q_1}(Q_2)\lambda_{Q_1}(f)\lambda_{Q_2}(f),\]
so $\lambda_{Q_1Q_2}(f)=\lambda_{Q_1}(f)\lambda_{Q_2}(f)\chi_{Q_1}(Q_2)\chi_{Q_2}(Q_1)$.

One deduces by induction that for any $Q \geq 1$ dividing $N$ and coprime to $N/Q$, \[\lambda_Q(f)=\prod_{p,q \mid Q, q \neq r}{\lambda_q(f)\chi_p(N_q)},\] and the conclusion follows after taking $Q=N$.}
 
\defi[epsilon-for-modular-forms]{Given a normalized newform $f \in \mathcal{S}_k(\Gamma_1(N))$ with primitive character $\chi$ and a prime number $p \mid N$, write $\chi=\chi_p\chi'_p$, where the conductor of $\chi_p$ (resp. $\chi'_p$) is a power of $p$ (resp. prime to $p$). We denote $\varepsilon_p(f)=\lambda_p(f)\chi'_p(p^{v_p(N)})$. Thus, $\lambda_N(f)=\prod_{p \mid N}{\varepsilon_p(f)}$. If $p \nmid N$, we define $\varepsilon_p(f)$ as equal to $1$.}

\lem[exact-level-bigtwist]{Let $f \in \mathcal{S}_k(\Gamma_1(N))$ be a normalized newform with character $\chi$. Let $\psi$ be a primitive Dirichlet character of conductor $P$, where $P \nmid N$ is a power of some prime $p$. Write $N=N_pN'$ with $N_p=p^{v_p(N)} < P$, $\chi=\chi_p\chi'_p$, where $\chi_p$ (resp. $\chi'_p$) has conductor dividing $N_p$ (resp. dividing $N'$). Then
\begin{itemize}[label=$-$,noitemsep]
\item The exact level of the modular form $g=\sum_{n \geq 1}{a_n(f)\psi(n)e^{2i\pi n\tau}}$ is $N'P^2$.
\item If $N_p=1$, then $W_{NP^2}(g) = \lambda_N(f)\psi(-N)\chi(P)\frac{\mathfrak{g}(\psi)}{\mathfrak{g}(\overline{\psi})}\overline{g}$.
\end{itemize}}

\demo{The first part is \cite[Theorem 3.1 ii)]{AL78}, so we assume that $N_p=1$. Then $W_{NP^2}$ is the action of the matrix $\begin{pmatrix}0 & 1\\-NP^2 & 0\end{pmatrix}$. Let $1 \leq u,v \leq P$ be coprime to $p$ and such that $Nuv \equiv -1\pmod{P}$. Then there exists some matrix $A_u \in \Gamma_0(N)$ such that $[A]_{1,1} \equiv P \pmod{N}$ and \[\begin{pmatrix}1 & u/P\\0 & 1\end{pmatrix}\begin{pmatrix}0 & 1\\-NP^2 & 0\end{pmatrix}= P\begin{pmatrix}0 & 1\\-N & 0\end{pmatrix}A\begin{pmatrix}1 & v/P \\0 & 1\end{pmatrix}.\] 

For each $1 \leq u \leq P$ prime to $p$, there is a unique $1 \leq v \leq P$ prime to $p$ such that $Nuv \equiv -1\pmod{P}$, and every $v$ comes from a unique $u$. Therefore,
\begin{align*}
W_{NP^2}(\mathfrak{g}(\overline{\psi})g) &= \sum_{1 \leq u \leq P, (u,P)=1}{\overline{\psi}(u)f \mid \begin{pmatrix}1 & u/P \\0 & 1\end{pmatrix}\begin{pmatrix}0 & 1\\-NP^2 & 0\end{pmatrix}}\\
&= \sum_{1 \leq v \leq P, (v,p)=1}{\psi(-Nv)(\langle P\rangle^{-1} (f \mid W_N))\mid \begin{pmatrix} 1 & v/P\\0 & 1\end{pmatrix}}\\
&= \psi(-N)\sum_{1 \leq v \leq P, (v,p)=1}{\psi(v)\lambda_N(f)\chi(P)\overline{f}\mid \begin{pmatrix}1 & v/P\\0 & 1\end{pmatrix}}\\
&= \psi(N)\lambda_N(f)\chi(P)\overline{\sum_{1 \leq v \leq P, (v,p)=1}{\overline{\psi}(v)f\mid\begin{pmatrix}1 & v/P\\0 & 1\end{pmatrix}}} \\
&= \psi(N)\lambda_N(f)\chi(P)\overline{\mathfrak{g}(\overline{\psi})}\overline{g}=\psi(-N)\lambda_N(f)\chi(P)\mathfrak{g}(\psi)\overline{g}.
\end{align*}
 }

\cor[global-constant-coprime-twist]{Let $f \in \mathcal{S}_k(\Gamma_1(N))$ be a normalized newform with character $\chi$. Let $M > 1$ be an integer coprime to $N$ and $\psi$ a primitive Dirichlet character of conductor $M$. Then \[f \otimes \psi = \sum_{n \geq 1}{a_n(f)\psi(n)e^{2i\pi n\tau}} \in \mathcal{S}_k(\Gamma_1(NM^2))\] is a newform of conductor $NM^2$ and \[\lambda_{NM^2}(g)=\lambda_N(f)\psi(-N)\chi(M)\frac{\mathfrak{g}(\psi)}{\mathfrak{g}(\overline{\psi})}.\]}

\demo{When $M$ is a power of a prime number $m$, by Lemma \ref{exact-level-bigtwist}, it is enough to show that the modular form $g(\tau)=\sum_{n \geq 1}{a_n(f)\psi(n)e^{2i\pi n\tau}}\in \mathcal{S}_k(\Gamma_1(NM^2))$ is a newform (of conductor $NM^2$): by strong multiplicity one, this implies that $g=f \otimes \psi$. This then follows from the definition of $g$ and \cite[Theorem 9]{Li}. 

For a general $M$, it is enough to proceed by induction after we check that for any primitive Dirichlet characters $\psi_1,\psi_2$ with coprime conductors $M_1,M_2$, then \[\mathfrak{g}(\psi_1\psi_2)=\psi_1(M_2)\psi_2(M_1)\mathfrak{g}(\psi_1)\mathfrak{g}(\psi_2).\] This equality is true because $(a,b) \in (\Z/a\Z)^{\times} \times (\Z/b\Z)^{\times} \longmapsto aM_2+bM_1 \in (\Z/M_1M_2\Z)^{\times}$ is a bijection.  
}

\prop[local-constant-bigtwist]{Let $f \in \mathcal{S}_k(\Gamma_1(N))$ be a normalized newform with character $\chi$ and $q$ a prime number (possibly not dividing $N$). Let $Q=q^{v_q(N)}, N'=N/Q$, and write $\chi=\chi_q\chi'_q$ as in the previous statements. Let $\psi$ be a primitive Dirichlet character whose conductor is some power $Q_1 > 1$ of $q$ with $Q \mid Q_1$. Assume that $\psi\chi_q$ and $\psi$ have the same conductor (which is automatic if $Q_1 > Q$). Then 

\[f \otimes \psi = \sum_{n\geq 1}{a_n(f)\psi(n)e^{2i\pi n\tau}}\] has conductor $N'Q_1^2$, and \[\lambda_q(f \otimes \psi) = \overline{\chi'_q}(Q_1)\psi(-1)\frac{\mathfrak{g}(\chi_q\psi)}{\mathfrak{g}(\overline{\psi})}.\]}

\demo{If $q \nmid N$, the first part follows from Corollary \ref{global-constant-coprime-twist}. To prove the second part, one has:
\begin{align*}
\lambda_{NQ_1^2}(f \otimes \psi) &= \lambda_N(f)\psi(-N)\chi(M)\frac{\mathfrak{g}(\psi)}{\mathfrak{g}(\overline{\psi})},\\
\lambda_{NQ_1^2}(f \otimes \psi) &= \lambda_q(f \otimes \psi)\chi(Q_1^2)\prod_{p \mid N}{\lambda_p(f \otimes \psi)(\chi'_p\psi^2)(N_p)} \\
&= \lambda_q(f \otimes \psi)\chi(Q_1)^2\prod_{p \mid N}{\lambda_p(f)\chi'_p(N_p)\psi(N_p)}\\
& = \lambda_q(f \otimes \psi)\chi(Q_1)^2\psi(N)\prod_{p \mid N}{\lambda_p(f)\chi'_p(N_p)} \\
&= \lambda_q(f \otimes \psi)\chi(Q_1)^2\psi(N)\lambda_N(f),
\end{align*} whence the conclusion. 

So we may assume that $q \mid N$, and then the statement follows from \cite[Theorem 4.1]{AL78}. }

\cor[level-newform-bigtwist]{Let $f \in \mathcal{S}_k(\Gamma_1(N))$ be a newform with conductor $\chi$, $M > 1$ and $\psi$ a primitive Dirichlet character of conductor $M$. Let $L$ be the conductor of $f \otimes \psi$, and $q \mid M$ be a prime number such that $v_q(M) \geq v_q(N)$, and such that $\psi_q\chi_q$ has the same conductor as $\psi_q$ (where $\chi=\chi_q\chi'_q$, the conductor of $\chi_q$ is a power of $q$, that of $\chi'_q$ is coprime to $q$). Then $v_q(L)=2v_q(M)$. }

\demo{Since $v_q(L)$ only depends on $\psi_q$, we may assume $\psi=\psi_q$, and this is a consequence of Proposition \ref{local-constant-bigtwist}. }

%% file: These.bbl
\providecommand{\bysame}{\leavevmode\hbox to3em{\hrulefill}\thinspace}
\providecommand{\MR}{\relax\ifhmode\unskip\space\fi MR }
\providecommand{\MRhref}[2]{%
  \href{http://www.ams.org/mathscinet-getitem?mr=#1}{#2}
}
\providecommand{\href}[2]{#2}
\begin{thebibliography}{100}

\bibitem{QC8+}
N.~Ad{\v z}aga, V.~Arul, L.~Beneish, M.~Chen, S.~Chidambaram, T.~Keller, and
  B.~Wen, \emph{Quadratic {C}habauty for {A}tkin-{L}ehner quotients of modular
  curves of prime level and genus 4, 5, 6}, Acta Arith. \textbf{208} (2023),
  no.~1, 15--49.

\bibitem{Anantharaman}
S.~Anantharaman, \emph{Sch{\'e}mas en groupes, espaces homog{\`e}nes et espaces
  alg{\'e}briques sur une base de dimension {$1$}}, Sur les groupes
  alg{\'e}briques, no.~33, Soc. Math. France, 1973, pp.~5--79.

\bibitem{AdRD}
S.~Arias{-de-}Reyna and L.~Dieulefait, \emph{Automorphy of
  {$\mathrm{GL}_2\otimes \mathrm{GL}_n$} in the self-dual case}, 2022, arXiv
  preprint 1611.06918v3.

\bibitem{Stacks}
The Stacks~Project Authors, \emph{Stacks project}.

\bibitem{BBM}
A.~Bajolet, Y.~Bilu, and B.~Matschke, \emph{Computing integral points on
  {$X_{ns}^+(p)$}}, Algebra Number Theory \textbf{15} (2021), no.~3, 569--608.

\bibitem{Bakertors}
M.~H. Baker, \emph{Torsion points on modular curves}, Invent. Math.
  \textbf{140} (2000), no.~3, 487--509.

\bibitem{BaTs}
Benjamin Bakker and Jacob Tsimerman, \emph{{$p$}-torsion monodromy
  representations of elliptic curves over geometric function fields}, Ann. of
  Math. (2) \textbf{184} (2016), no.~3, 709--744.

\bibitem{Xns+13}
S.~Balakrishnan, N.~Dogra, S.~M{\"u}ller, J.~Tuitman, and J.~Vonk,
  \emph{Explicit {C}habauty–{K}im for the split {C}artan modular curve of
  level {$13$}}, Ann. of Math. (2) \textbf{189} (2019), no.~3, 885--944.

\bibitem{Xns+17}
\bysame, \emph{Quadratic {C}habauty for modular curves: algorithms and
  examples}, Compos. Math. \textbf{159} (2023), no.~6, 1111--1152.

\bibitem{BKX13}
F.~Bars, A.~Kontogeorgis, and X.~Xarles, \emph{Bielliptic and hyperelliptic
  modular curves {$X(N)$} and the group {$Aut(X(N))$}}, Acta Arith.
  \textbf{161} (2013), no.~3, 283--299.

\bibitem{Berger-caracp}
L.~Berger and S.~Rozensztajn, \emph{Modular representations of
  {$GL_2(\mathbb{F}_p)$}}, Found online.

\bibitem{BP}
Y.~Bilu and P.~Parent, \emph{Serre's uniformity problem in the split {C}artan
  case}, Ann. of Math. \textbf{173} (2011), no.~1, 569--584.

\bibitem{BPR}
Y.~Bilu, P.~Parent, and M.~Rebolledo, \emph{Rational points on {$X_0^+(p^r)$}},
  Ann. {I}nst. {F}ourier ({G}renoble) \textbf{63} (2013), no.~3, 957--984.

\bibitem{BKTama}
S.~Bloch and K.~Kato, \emph{{L-F}unctions and {T}amagawa {N}umbers of
  {M}otives}, {T}he {G}rothendieck {F}estschrift: {A C}ollection of {A}rticles
  {W}ritten in {H}onor of the 60th {B}irthday of {A}lexander {G}rothendieck
  (P.~Cartier, L.~Illusie, N.~M. Katz, G.~Laumon, Y.~I. Manin, and K.~A. Ribet,
  eds.), Birkh{\"a}user Boston, Boston, MA, 2007, pp.~333--400.

\bibitem{BLR}
S.~Bosch, W.~L{\"u}tkebohmert, and M.~Raynaud, \emph{N{\'e}ron models}, Ergeb.
  Math. Grenzgeb. (3), Springer-Verlag, 1990.

\bibitem{magma}
W.~Bosma, J.~Cannon, and C.~Playoust, \emph{The {M}agma algebra system. {I}.
  {T}he user language}, J. Symbolic Comput. \textbf{24} (1997), no.~3-4,
  235--265, Computational algebra and number theory (London, 1993).
  \MR{MR1484478}

\bibitem{HigherHida}
G.~Boxer, F.~Calegari, T.~Gee, and V.~Pilloni, \emph{Abelian surfaces over
  totally real fields are potentially modular}, Publ. Math. Inst. Hautes
  {\'E}tudes Sci. \textbf{134} (2021), no.~1, 153--501.

\bibitem{fullmod}
C.~Breuil, B.~Conrad, F.~Diamond, and R.~Taylor, \emph{On the {M}odularity of
  {E}lliptic {C}urves {O}ver {$\mathbf{Q}$}: {W}ild 3-{A}dic {E}xercises}, J.
  Amer. Math. Soc. \textbf{14} (2001), no.~4, 843--939.

\bibitem{Brumley}
F.~Brumley, Private communication.

\bibitem{Bump}
D.~Bump, \emph{Automorphic forms and representations}, Cambridge Stud. Adv.
  Math., Cambridge Univ. Press, 1997.

\bibitem{BH}
C.~J. Bushnell and G.~Henniart, \emph{The {L}ocal {L}anglands {C}onjecture for
  {GL(2)}}, 1st ed. 2006. ed., Grundlehren Math. Wiss., no. 335, Springer,
  2006.

\bibitem{CHM}
L.~Caporaso, J.~Harris, and B.~Mazur, \emph{Uniformity of rational points}, J.
  Amer. Math. Soc. \textbf{10} (1997), no.~1, 1--35.

\bibitem{localglobal}
H.~Carayol, \emph{Sur les repr\'esentations $l$-adiques associ\'ees aux formes
  modulaires de {Hilbert}}, Ann. Sci. {\'E}c. Norm. Sup{\'e}r. \(4\)
  \textbf{19} (1986), no.~3, 409--468.

\bibitem{BCarayol}
\bysame, \emph{Preuve de la conjecture de {Langlands} locale pour {$GL_n$}~:
  travaux de {Harris-Taylor} et {Henniart}}, S\'eminaire Bourbaki : volume
  1998/99, expos\'es 850-864, Ast\'erisque, no. 266, Soc. Math. France, 2000,
  Expos\'e 857.

\bibitem{Chen-Merel-conj}
I.~Chen, \emph{On {R}elations between {J}acobians of {C}ertain {M}odular
  {C}urves}, J. Algebra \textbf{231} (2000), no.~1, 414--448.

\bibitem{Clozel-regulier}
L.~Clozel, \emph{Motives and automorphic representations}, Autour des
  motifs---{\'E}cole d'{\'e}t{\'e} {F}ranco-{A}siatique de {G\'e}om{\'e}trie
  {A}lg{\'e}brique et de {T}h{\'e}orie des {N}ombres / {A}sian-{F}rench
  {S}ummer {S}chool on {A}lgebraic {G}eometry and {N}umber {T}heory. {V}ol.
  {III} (J.-B. Bost and J.-M. Fontaine, eds.), Panor. Synth{\`e}ses, vol.~49,
  Soc. Math. France, 2016, pp.~29--60.

\bibitem{Cogdell}
J.~W. Cogdell, \emph{Lectures on {L-F}unctions, {C}onverse {T}heorems, and
  {F}unctoriality for {$GL(n)$}}.

\bibitem{ColChab}
R.~Coleman, \emph{Effective {C}habauty}, Duke Math. J. \textbf{52} (1985),
  no.~3, 765 -- 770.

\bibitem{ColAb}
\bysame, \emph{Torsion points on curves and p-adic abelian integrals}, Ann. of
  Math. \textbf{121} (1985), no.~1, 111--168.

\bibitem{Conrad99}
B.~Conrad, \emph{Finite {G}roup {S}chemes over {B}ases with {L}ow
  {R}amification}, Compos. Math. \textbf{119} (1999), no.~3, 239--326.

\bibitem{Conrad-Duality}
\bysame, \emph{Grothendieck duality and base change}, Lecture Notes in Math.,
  vol. 1750, Springer, 2000.

\bibitem{mod3}
B.~Conrad, F.~Diamond, and R.~Taylor, \emph{Modularity of certain potentially
  {B}arsotti-{T}ate {G}alois representations}, J. Amer. Math. Soc. \textbf{12}
  (1999), no.~2, 521--567.

\bibitem{Cremona-Freitas}
J.~Cremona and N.~Freitas, \emph{Global methods for the symplectic type of
  congruences between elliptic curves}, Rev. Mat. Iberoam. \textbf{38} (2022),
  no.~1, 1--32.

\bibitem{Deconinck-PhD}
H.~Deconinck, \emph{The generalized {F}ermat equation over totally real number
  fields}, Unpublished, online at
  {\url{https://wrap.warwick.ac.uk/id/eprint/81893/1/WRAP_Theses_Deconinck_2016.pdf}},
  last accessed {O}ct 25th 2024, July 2016.

\bibitem{Deligne-csts}
P.~Deligne, \emph{Les {C}onstantes des {E}quations {F}onctionnelles des
  {F}onctions {L}}, Modular Functions of One Variable II (Berlin, Heidelberg)
  (P.~Deligne and W.~Kuijk, eds.), Springer Berlin Heidelberg, 1973,
  pp.~501--597.

\bibitem{DeRa}
P.~Deligne and M.~Rapoport, \emph{Les sch{\'e}mas de modules de courbes
  elliptiques}, Modular Functions of One Variable II (Berlin, Heidelberg)
  (P.~Deligne and W.~Kuijk, eds.), Springer Berlin Heidelberg, 1973,
  pp.~143--316.

\bibitem{Del-Ser}
P.~Deligne and J.-P. Serre, \emph{Formes modulaires de poids $1$}, Ann. Sci.
  {\'E}c. Norm. Sup{\'e}r. \(4\) \textbf{7} (1974), no.~4, 507--530.

\bibitem{SGA3}
M.~Demazure and A.~Grothendieck, \emph{Sch{\'e}mas en groupes. {T}ome {I},
  {P}ropri{\'e}t{\'e}s g{\'e}n{\'e}rales des sch{\'e}mas en groupes :
  {S\'e}minaire de g{\'e}om{\'e}trie alg{\'e}brique du {B}ois {M}arie, 1962-64,
  {SGA 3}}, Lecture Notes in Math., vol. 151, Springer, 1962--64.

\bibitem{Diamond-personal}
F.~Diamond, Private communication.

\bibitem{moddiamond}
\bysame, \emph{On deformation rings and {H}ecke rings}, Ann. of Math.
  \textbf{144} (1996), no.~1, 137--166.

\bibitem{DS}
Fred Diamond and Jerry Shurman, \emph{A {F}irst {C}ourse in {M}odular {F}orms},
  Grad. Texts in Math., vol. 228, Springer, 2005.

\bibitem{BaseChangeGL2}
L.~Dieulefait, \emph{Langlands base change for {$GL(2)$}}, Ann. of Math. (2)
  \textbf{176} (2012), no.~2, 1015--1038.

\bibitem{LanglandsSym5}
\bysame, \emph{Automorphy of {$\mathrm{Symm}^5(GL(2))$} and base change}, J.
  Math. Pures Appl. (9) \textbf{104} (2015), no.~4, 619--656.

\bibitem{LFD}
N.~Dogra and S.~Le~Fourn, \emph{Quadratic {C}habauty for modular curves and
  modular forms of rank one}, Math. Ann. \textbf{380} (2021), no.~1, 393--448.

\bibitem{Faltings86}
G.~Faltings, \emph{Finiteness {T}heorems for {A}belian {V}arieties over
  {N}umber {F}ields}, Arithmetic Geometry (G.~Cornell and J.~H. Silverman,
  eds.), Springer New York, New York, NY, 1986, pp.~9--27.

\bibitem{Fisher5}
T.~Fisher, \emph{The {H}essian of a genus one curve}, Proc. Lond. Math. Soc.
  (3) \textbf{104} (2012), no.~3, 613--648.

\bibitem{Fisher-711}
\bysame, \emph{On families of 7- and 11-congruent elliptic curves}, LMS J.
  Comput. Math. \textbf{17} (2014), no.~1, 536--564.

\bibitem{Fisher-13}
\bysame, \emph{On families of 13-congruent elliptic curves}, 2019, arXiv
  preprint 1912.10777.

\bibitem{Fisher-17}
\bysame, \emph{On pairs of 17-congruent elliptic curves}, 2021, arXiv preprint
  2106.02033.

\bibitem{Fontaine85}
J.-M. Fontaine, \emph{Il n'y a pas de vari{\'e}t{\'e} ab{\'e}lienne sur
  {$\mathbb{Z}$}.}, Invent. Math. \textbf{81} (1985), 515--538.

\bibitem{Forster}
O.~Forster, \emph{Lectures on {R}iemann surfaces}, Grad. Texts in Math.,
  vol.~81, Springer-Verlag, New York Heidelberg Berlin, 1981.

\bibitem{FreitasFermat}
N.~Freitas, \emph{On the {F}ermat-type equation {$x^3+y^3=z^p$}}, Comment.
  Math. Helv. \textbf{91} (2016), 295--304.

\bibitem{Freitas-Kraus}
N.~Freitas and A.~Kraus, \emph{On the symplectic type of isomorphisms of the
  {$p$}-torsion of elliptic curves}, Mem. Amer. Math. Soc. \textbf{277} (2022),
  no.~1361, v+105.

\bibitem{FNS2020}
N.~Freitas, B.~Naskr{\c e}cki, and M.~Stoll, \emph{The generalized {F}ermat
  equation with exponents 2, 3, $n$}, Compos. Math. \textbf{156} (2020), no.~1,
  77--113.

\bibitem{Frey86}
G.~Frey, \emph{Links between stable elliptic curves and certain diophantine
  equations}, Ann. Univ. Sarav. Ser. Math. \textbf{1} (1986), no.~1, iv+40.

\bibitem{FreyFLT}
\bysame, \emph{{O}n {T}ernary {E}quations of {F}ermat {T}ype and {R}elations
  with {E}lliptic {C}urves}, Modular forms and Fermat's Last Theorem
  (G.~Cornell, J.~H. Silverman, and G.~Stevens, eds.), Springer-Verlag NY,
  1997, pp.~527--548.

\bibitem{sign-ortho}
A.~Fr{\"o}hlich and J.~Queyrut, \emph{On the functional equation of the {A}rtin
  {L}-function for characters of real representations}, Invent. Math.
  \textbf{20} (1973), no.~2, 125--138.

\bibitem{lombardo}
L.~Furio and D.~Lombardo, \emph{Serre's uniformity question and proper
  subgroups of {$C_{ns}^+(p)$}}, 2023, arXiv preprint 2305.17780.

\bibitem{GeeMLT}
T.~Gee, \emph{Modularity lifting theorems}, Essential Number Theory \textbf{1}
  (2022), no.~1, 73--126.

\bibitem{Gelbart}
S.~S. Gelbart, \emph{Automorphic forms on {A}dele groups}, Ann. of Math. Stud.,
  vol.~83, Princeton Univ. Press, 1975.

\bibitem{GelbartFLT}
\bysame, \emph{Three lectures on the modularity of {$\overline{\rho}_{E,3}$}
  and the {L}anglands {R}eciprocity {C}onjecture}, Modular forms and Fermat's
  Last Theorem (G.~Cornell, J.~H. Silverman, and G.~Stevens, eds.),
  Springer-Verlag NY, 1997, pp.~155--207.

\bibitem{Getz-Hahn}
J.~R. Getz and H.~Hahn, \emph{An {I}ntroduction to {A}utomorphic
  {R}epresentations. {W}ith a view toward trace formulae}, Grad. Texts in
  Math., vol. 300, Springer International Publishing, 2024.

\bibitem{GGJQ}
E.~Ghate, E.~González-Jiménez, and J.~Quer, \emph{On the {B}rauer class of
  modular endomorphism algebras}, Int. Math. Res. Not. IMRN \textbf{2005}
  (2005), no.~12, 701--723.

\bibitem{Gross-Reeder}
B.~H. Gross and M.~Reeder, \emph{Arithmetic invariants of discrete {L}anglands
  parameters}, Duke Math. J. \textbf{154} (2010), no.~3, 431--508.

\bibitem{Grothendieck3}
A.~Grothendieck, \emph{Techniques de construction et th{\'e}or{\`e}mes
  d'existence en g{\'e}om{\'e}trie alg{\'e}brique {III} : pr{\'e}sch{\'e}mas
  quotients}, S{\'e}minaire {B}ourbaki: ann{\'e}es 1960/61, expos{\'e}s
  205--222, S{\'e}minaire {B}ourbaki, no.~6, Soc. Math. France, 1961,
  pp.~99--118.

\bibitem{EGA-III2}
A.~Grothendieck and J.~Dieudonn{\'e}, \emph{{\'E}l{\'e}ments de
  g{\'e}om{\'e}trie alg{\'e}brique : {III. \'E}tude cohomologique des faisceaux
  coh{\'e}rents, {S}econde partie}, Publ. Math. Inst. Hautes {\'E}tudes Sci.
  \textbf{17} (1963), 5--91 (fr). \MR{163911}

\bibitem{EGA-IV4}
\bysame, \emph{{\'E}l{\'e}ments de g{\'e}om{\'e}trie alg{\'e}brique : {IV.
  \'E}tude locale des sch{\'e}mas et des morphismes de sch{\'e}mas,
  {Q}uatri{\`e}me partie}, Publ. Math. Inst. Hautes {\'E}tudes Sci. \textbf{32}
  (1967), 5--361 (fr).

\bibitem{SGA1}
A.~Grothendieck and M.~Raynaud, \emph{Rev{\^e}tements {\'e}tales et groupe
  fondamental : {S\'e}minaire de {g\'eom\'etrie} alg{\'e}brique du {B}ois
  {M}arie, 1960-61, {SGA 1}}, Lecture Notes in Math., vol. 224, Springer,
  1960--61.

\bibitem{HOM}
A.~J. Hahn and O.~T. O’Meara, \emph{The {C}lassical {G}roups and {K-T}heory},
  Grundlehren Math. Wiss., no. 291, Springer-Verlag, 1989.

\bibitem{Halb}
E.~Halberstadt, \emph{Calculs explicites sur les courbes modulaires}, Feb 2021.

\bibitem{HK-Fermat}
E.~Halberstadt and A.~Kraus, \emph{Courbes de {F}ermat: r{\'e}sultats et
  probl{\`e}mes}, J. reine angew. Math. \textbf{548} (2002), 167--234.

\bibitem{HalbKraus}
\bysame, \emph{Sur la courbe modulaire {$X_E(7)$}}, Exp. Math. \textbf{12}
  (2003), no.~1, 27--40.

\bibitem{Harari}
D.~Harari, \emph{Galois {C}ohomology and {C}lass {F}ield {T}heory},
  Universitext., Springer, 2020.

\bibitem{HelleCrelle}
Y.~Hellegouarch, \emph{{\'E}tude des points d'ordre fini des vari{\'e}t{\'e}s
  ab{\'e}liennes de dimension un d{\'e}finies sur un anneau principal}, J.
  Reine Angew. Math. \textbf{244} (1970), 20--36.

\bibitem{THelle}
\bysame, \emph{Courbes elliptiques et {\'e}quation de {F}ermat}, 1972,
  Th{\`e}se.

\bibitem{Hida}
H.~Hida, \emph{Elementary {T}heory of {L}-functions and {E}isenstein {S}eries},
  London Math. Soc. Stud. Texts, Cambridge Univ. Press, 1993.

\bibitem{IR}
K.~F. Ireland and M.~I. Rosen, \emph{A classical introduction to modern number
  theory}, Grad. Texts in Math., vol.~84, Springer International Publishing,
  1990.

\bibitem{IK}
H.~Iwaniec and E.~Kowalski, \emph{Analytic number theory}, Amer. Math. Soc.
  Colloq. Publ, AMS, Providence (R.I.), 2004.

\bibitem{Jacquet}
H.~Jacquet, \emph{Automorphic forms on {$GL(2)$} {P}art {II}}, Lecture Notes in
  Math., vol. 278, Springer Berlin Heidelberg, 1972.

\bibitem{JL}
H.~Jacquet and R.~Langlands, \emph{Automorphic forms on {$GL (2)$}}, Lecture
  Notes in Math., no. 114, Springer-Verlag, 1970.

\bibitem{Jorza-gl2fq}
A.~Jorza, \emph{Complex representations of {$GL_2(\mathbb{F}_q)$}}, Found
  online.

\bibitem{KR}
E.~J. Kani and O.~G. Rizzo, \emph{Mazur’s question on mod 11 representations
  of elliptic curves}, preprint.

\bibitem{MDQS}
E.~J. Kani and W.~Schanz, \emph{Modular diagonal quotient surfaces}, Math. Z.
  \textbf{227} (1998), no.~2, 337--366.

\bibitem{KatoBSD}
K.~Kato, \emph{{$p$}-adic {H}odge theory and values of zeta functions of
  modular forms}, Cohomologies {$p$}-adiques et applications arithm{\'e}tiques
  ({III}) (P.~Berthelot, J.-M. Fontaine, L.~Illusie, K.~Kato, and M.~Rapoport,
  eds.), Ast\'erisque, no. 295, Soc. Math. France, 2004. \MR{2104361}

\bibitem{Katz}
N.~Katz, \emph{{$p$}-adic properties of modular schemes and modular forms},
  Modular Functions of One Variable III (Berlin, Heidelberg) (W.~Kuijk and
  J.-P. Serre, eds.), Springer Berlin Heidelberg, 1973, pp.~69--190.

\bibitem{KM}
N.~M. Katz and B.~Mazur, \emph{Arithmetic moduli of elliptic curves}, no. 108,
  Princeton University Press, 1985.

\bibitem{KW1}
C.~Khare and J.-P. Wintenberger, \emph{Serre’s modularity conjecture {(I)}},
  Invent. Math. \textbf{178} (2009), no.~3, 485--504.

\bibitem{KW2}
\bysame, \emph{Serre’s modularity conjecture {(II)}}, Invent. Math.
  \textbf{178} (2009), no.~3, 505--586.

\bibitem{Kim05}
M.~Kim, \emph{The motivic fundamental group of
  {$\mathbb{P}^1\backslash\{0,1,\infty\}$} and the theorem of {S}iegel.},
  Invent. math. \textbf{161} (2005), 629--656.

\bibitem{Kim09}
\bysame, \emph{The unipotent {A}lbanese map and {S}elmer varieties for curves},
  Publ. Res. Inst. Math. Sci. \textbf{45} (2009), 89--133.

\bibitem{KimTama08}
M.~Kim and A.~Tamagawa, \emph{The l-component of the unipotent {A}lbanese map},
  Math. Ann. \textbf{340} (2008), 223 -- 235.

\bibitem{KLZ15}
G.~Kings, D.~Loeffler, and S.~Zerbes, \emph{Rankin-{E}isenstein classes and
  explicit reciprocity laws}, Camb. J. Math. \textbf{5} (2017), no.~1, 1--122.

\bibitem{Kleiman}
S.~L. Kleiman, \emph{Relative duality for quasi-coherent sheaves}, Compos.
  Math. \textbf{41} (1980), no.~1, 39--60.

\bibitem{KolLog}
V.~A. Kolyvagin and D.~Y. Logach{\"e}v, \emph{Finiteness of the
  {S}hafarevich-{T}ate group and the group of points for some modular abelian
  varieties}, Algebra i Analiz \textbf{1} (1989), no.~5, 171--196.

\bibitem{Kraus-Thesis}
A.~Kraus, \emph{D{\'e}termination du poids et du conducteur associ{\'e}s aux
  repr{\'e}sentations des points de {$p$}-torsion d'une courbe elliptique},
  Dissertationes Math. (Rozprawy Mat.) \textbf{364} (1997), 39.

\bibitem{Kraus-a3b3cp}
\bysame, \emph{Sur l'{\'e}quation {$\{a^3+b^3=c^p\}$}}, Exp. Math. \textbf{7}
  (1998), no.~1, 1 -- 13.

\bibitem{KO}
A.~Kraus and J.~Oesterl{\'e}, \emph{Sur une question de {B. M}azur.}, Math.
  Ann. \textbf{293} (1992), no.~2, 259--276.

\bibitem{Kumar}
A.~Kumar, \emph{Hilbert modular surfaces for square discriminants and elliptic
  subfields of genus 2 function fields}, Res. Math. Sci. \textbf{2} (2015),
  Art. 24, 46 p.

\bibitem{Lang-ANT}
S.~Lang, \emph{Algebraic number theory}, 2nd edition ed., Grad. Texts in Math.,
  vol. 110, Springer, 1994.

\bibitem{LangComplex}
\bysame, \emph{Complex {A}nalysis}, 4th ed., Grad. Texts in Math., vol. 103,
  Springer, New York, 1999.

\bibitem{LFL}
S.~Le~Fourn and P.~Lemos, \emph{Residual {G}alois representations of elliptic
  curves with image contained in the normaliser of a nonsplit {C}artan},
  Algebra Number Theory \textbf{15} (2021), no.~3, 747--771.

\bibitem{Li}
W.-C.~W. Li, \emph{Newforms and {F}unctional {E}quations}, Math. Ann.
  \textbf{212} (1975), no.~4, 285--315.

\bibitem{Li-RS}
\bysame, \emph{{L}-series of {R}ankin type and their functional equations},
  Math. Ann. \textbf{244} (1979), no.~2, 135--166.

\bibitem{AL78}
W.-C.~W. Li and A.O.L. Atkin, \emph{Twists of {N}ewforms and
  {P}seudo-{E}igenvalues of {W}-{O}perators}, Invent. Math. \textbf{48} (1978),
  221--244.

\bibitem{LLR}
Q.~Liu, D.~Lorenzini, and M.~Raynaud, \emph{N{\'e}ron models, {L}ie algebras,
  and reduction of curves of genus one}, Invent. math. \textbf{157} (2004),
  no.~3, 455--518.

\bibitem{QL}
Q.~Liu and E.~Reinie, \emph{Algebraic geometry and arithmetic curves}, Oxf.
  Grad. Texts Math., Oxford University Press, 2002.

\bibitem{lmfdb}
The {LMFDB Collaboration}, \emph{The {L}-functions and modular forms database},
  \url{https://www.lmfdb.org}, 2023, [Online; accessed 12 June 2023].

\bibitem{bigimage}
D.~Loeffler, \emph{Images of adelic {G}alois representations for modular
  forms}, Glasg. Math. J. \textbf{59} (2017), no.~1, 11--–25.

\bibitem{LWnewforms}
D.~Loeffler and J.~Weinstein, \emph{On the computation of local components of a
  newform}, Math. Comp. \textbf{81} (2012), 1179--1200.

\bibitem{MW}
D.~Masser and G.~W{\"u}stholz, \emph{Isogeny estimates for abelian varieties,
  and finiteness theorems}, Ann. of Math. (2) \textbf{137} (1993), no.~3,
  459--472.

\bibitem{MazurY1}
B.~Mazur, \emph{Modular curves and the {E}isenstein ideal}, Publ. Math. Inst.
  Hautes {É}tudes Sci. \textbf{47} (1977), 33--186.

\bibitem{FreyMazur}
\bysame, \emph{Rational isogenies of prime degree}, Invent. Math. \textbf{44}
  (1978), no.~2, 129--162, With an appendix by D. Goldfeld.

\bibitem{Mazur-Open}
\bysame, \emph{Open problems regarding rational points on curves and
  varieties}, Galois {R}epresentations in {A}rithmetic {A}lgebraic {G}eometry
  (A.~J. Scholl and R.~L. Taylor, eds.), London Math. Soc. Lecture Note Ser.,
  Cambridge Univ. Press, 1998, pp.~239--266.

\bibitem{MR-ES}
B.~Mazur and K.~Rubin, \emph{Kolyvagin systems}, Mem. Amer. Math. Soc., Amer.
  Math. Soc., Providence (R.I.), 2004.

\bibitem{PoonenMcCallum}
W.~McCallum and B.~Poonen, \emph{The {M}ethod of {C}habauty and {C}oleman},
  Explicit methods in number theory, Panor. Synth{\`e}ses, vol.~36, Soc. Math.
  France, 2012, pp.~99--117.

\bibitem{Merel}
L.~Merel, \emph{Bornes pour la torsion des courbes elliptiques sur les corps de
  nombres}, Invent. Math. \textbf{124} (1996), no.~1, 437--449.

\bibitem{MilneRes}
J.~S. Milne, \emph{On the arithmetic of abelian varieties}, Invent. Math.
  \textbf{17} (1972), no.~3, 177--190.

\bibitem{MilAb}
\bysame, \emph{Abelian varieties}, Arithmetic Geometry (G.~Cornell and J.~H.
  Silverman, eds.), Springer-Verlag, 1985, pp.~102--150.

\bibitem{MilJac}
\bysame, \emph{Jacobian varieties}, Arithmetic Geometry (G.~Cornell and J.~H.
  Silverman, eds.), Springer-Verlag, 1985, pp.~167--212.

\bibitem{MilneADT}
\bysame, \emph{Arithmetic duality theorems}, 2nd edition ed., BookSurge, 2006.

\bibitem{Miyake}
T.~Miyake, \emph{Modular forms}, Springer Monogr. Math., Springer, 1989.

\bibitem{Momose}
F.~Momose, \emph{On the {$\ell$}-adic representations attached to modular
  forms}, J. Fac. Sci. Univ. Tokyo Sect. IA Math. \textbf{28} (1981), no.~1,
  89--109. \MR{617867}

\bibitem{Neukirch-ANT}
J.~Neukirch, \emph{Algebraic number theory}, Grundlehren Math. Wiss., Springer,
  cop. 1999, transl. from German by N. Schappacher.

\bibitem{CNB}
J.~Neukirch, A.~Schmidt, and K.~Wingberg, \emph{Cohomology of number fields},
  Grundlehren Math. Wiss., vol. 323, Springer, 2008, Second edition.

\bibitem{NewThorQ1}
J.~Newton and J.~Thorne, \emph{Symmetric power functoriality for holomorphic
  modular forms}, Publ. Math. Inst. Hautes {\'E}tudes Sci. \textbf{134} (2021),
  no.~1, 1--116.

\bibitem{NewThorQ2}
\bysame, \emph{Symmetric power functoriality for holomorphic modular forms,
  {II}}, Publ. Math. Inst. Hautes {\'E}tudes Sci. \textbf{134} (2021), no.~1,
  117--152.

\bibitem{NewThor}
J.~Newton and J.~A. Thorne, \emph{Symmetric power functoriality for {H}ilbert
  modular forms}, 2022, arXiv preprint 2212.03595.

\bibitem{Parent-Torsion}
P.~Parent, \emph{Bornes effectives pour la torsion des courbes elliptiques sur
  les corps de nombres}, J. Reine Angew. Math. \textbf{506} (1999), 85--116.

\bibitem{Pellarin}
F.~Pellarin, \emph{Sur une majoration explicite pour un degr{\'e}
  d'isog{\'e}nie liant deux courbes elliptiques}, Acta Arith. \textbf{100}
  (2001), no.~3, 203--243.

\bibitem{PSS}
B.~Poonen, E.~Schaefer, and M.~Stoll, \emph{Twists of {$X(7)$} and primitive
  solutions to {$x^2+y^3=z^7$}}, Duke Math. J. \textbf{137} (2007), no.~1, 103
  -- 158.

\bibitem{Prasad-cusp}
D.~Prasad, \emph{Notes on modular representations of p-adic groups, and the
  {L}anglands {C}orrespondence}, Found online.

\bibitem{Antwerp5-Ribet}
K.~A. Ribet, \emph{Galois representations attached to eigenforms with
  nebentypus}, Modular Functions of one Variable V (Berlin, Heidelberg) (J.-P.
  Serre and D.~Zagier, eds.), Springer Berlin Heidelberg, 1977, pp.~18--52.

\bibitem{Ribet2}
\bysame, \emph{On {$\ell$}-adic representations attached to modular forms.
  {II}}, Glasgow Math. J. \textbf{27} (1985), 185--194. \MR{819838}

\bibitem{Ribet90}
\bysame, \emph{On modular representations of {$(\Qbar/\Q)$} arising from
  modular forms}, Invent. Math. \textbf{100} (1990), no.~1, 431--476.

\bibitem{Ribet-report}
\bysame, \emph{Report on mod {$l$} representations of {${\rm Gal}(\overline{\bf
  Q}/{\bf Q})$}}, Motives ({S}eattle, {WA}, 1991), Proc. Sympos. Pure Math.,
  vol. 55, Part 2, Amer. Math. Soc., Providence, RI, 1994, pp.~639--676.

\bibitem{Ribet-Stein}
K.~A. Ribet and W.~A. Stein, \emph{Lectures on {S}erre's conjectures},
  Arithmetic algebraic geometry ({P}ark {C}ity, {UT}, 1999), IAS/Park City
  Math. Ser., vol.~9, Amer. Math. Soc., Providence, RI, 2001, With appendices
  by B. Conrad and K. Buzzard., pp.~143--232.

\bibitem{Rubin-ES}
K.~Rubin, \emph{Euler systems}, Ann. of Math. Stud., Princeton Univ. Press,
  2000.

\bibitem{Rubin-Silverberg}
K.~Rubin and A.~Silverberg, \emph{Families of elliptic curves with constant mod
  {$p$} representations}, Elliptic curves, modular forms, and {F}ermat's last
  theorem ({H}ong {K}ong, 1993), Ser. Number Theory, vol.~I, Int. Press,
  Cambridge, MA, 1995, pp.~148--161.

\bibitem{CF6}
J.-P. Serre, \emph{Local {C}lass {F}ield {T}heory}, Algebraic {N}umber
  {T}heory: {P}roceedings of an {I}nstructional {C}onference (J.W.S. Cassels
  and A.~Fr{\"o}lich, eds.), Academic {P}ress Inc. (London) Ltd, 1967.

\bibitem{SerreMcGill}
\bysame, \emph{Abelian ell-adic representations and elliptic curves :
  {M}c{G}ill {U}niversity lecture notes}, W. A. Benjamin, Inc., New York
  Amsterdam, 1968.

\bibitem{Serre-image-ouverte}
\bysame, \emph{Propri{\'e}t{\'e}s galoisiennes des points d'ordre fini des
  courbes elliptiques}, Invent. Math. \textbf{15} (1971), no.~4, 259--331.

\bibitem{SerreLinReps}
\bysame, \emph{Linear representations of finite groups}, 1st ed. 1977. ed.,
  Grad. Texts in Math., vol.~42, Springer New York, 1977.

\bibitem{Serre-local}
\bysame, \emph{Local fields}, Grad. Texts in Math., no.~67, Springer, 1979,
  translated from French by Marvin J. Greenberg.

\bibitem{Serre-Cebotarev}
\bysame, \emph{Quelques applications du th{\'e}or{\`e}me de densit{\'e} de
  {C}hebotarev}, Publ. Math. Inst. Hautes {\'E}tudes Sci. \textbf{54} (1981),
  123--201.

\bibitem{Conj-Serre}
\bysame, \emph{Sur les représentations modulaires de degré $2$ de
  {$\mathrm{Gal}({\overline {\mathbf{Q}}}/\mathbf{Q})$}}, Duke Math. J.
  \textbf{54} (1987), no.~1, 179--230.

\bibitem{GoodRed}
J.-P. Serre and J.~Tate, \emph{Good reduction of abelian varieties}, Ann. of
  Math. (2) \textbf{88} (1968), 492--517.

\bibitem{Shatz-Groups}
S.~Shatz, \emph{Group {S}chemes, {F}ormal {G}roups, and {$p$-D}ivisible
  {G}roups}, Arithmetic {G}eometry (G.~Cornell and J.~H. Silverman, eds.),
  Springer-Verlag, 1985, pp.~102--150.

\bibitem{Shimura76}
G.~Shimura, \emph{The special values of the zeta functions associated with cusp
  forms}, Comm. Pure Appl. Math. \textbf{29} (1976), no.~6, 783--804.

\bibitem{Shimura-Periods}
\bysame, \emph{On the periods of modular forms}, Math. Ann. \textbf{229}
  (1977), no.~3, 211 -- 221.

\bibitem{shimura-hecke}
\bysame, \emph{Introduction to the {A}rithmetic {T}heory of {A}utomorphic
  {F}unctions}, {P}rinceton {U}niversity {P}ress, 1994.

\bibitem{Siksek-PS}
S.~Siksek, \emph{The {M}odular {A}pproach to {D}iophantine {E}quations},
  Explicit methods in number theory, Panor. Synth{\`e}ses, vol.~36, Soc. Math.
  France, 2012, pp.~151--179.

\bibitem{AEC2}
J.~H. Silverman, \emph{Advanced topics in the arithmetic of elliptic curves},
  Grad. Texts in Math., vol. 151, Springer-Verlag, 1994, Second printing.

\bibitem{AEC1}
\bysame, \emph{The arithmetic of elliptic curves}, 2nd ed. 2009. ed., Grad.
  Texts in Math., no. 106, Springer New York, 2009.

\bibitem{Stein}
W.~A. Stein, \emph{Modular forms, a computational approach}, vol.~79, Am. Math.
  Soc., 2007.

\bibitem{Studnia}
E.~Studnia, \emph{On obstructions to the {E}uler system method for
  {R}ankin-{S}elberg convolutions}, 2024, arXiv preprint 2401.17769.

\bibitem{Szamuely}
T.~Szamuely, \emph{Galois groups and fundamental groups}, Cambridge Stud. Adv.
  Math., Cambridge Univ. Press, 2009.

\bibitem{Michel}
R.~Takloo-Bighash, \emph{Spinor {L}-functions, theta correspondence, and
  {B}essel coefficients}, Forum Math. \textbf{19} (2007), no.~3, 487--554, With
  an appendix by P. Michel.

\bibitem{CF-Tate}
J.~Tate, \emph{Fourier analysis in number fields and {H}ecke's zeta-functions},
  Algebraic {N}umber {T}heory: {P}roceedings of an {I}nstructional {C}onference
  (J.W.S. Cassels and A.~Fr{\"o}lich, eds.), Academic {P}ress Inc. (London)
  Ltd, 1967.

\bibitem{Tate-divisible}
\bysame, \emph{p-{D}ivisible {G}roups}, Proceedings of a Conference on Local
  Fields (Berlin, Heidelberg) (T.~A. Springer, ed.), Springer Berlin
  Heidelberg, 1967, pp.~158--183.

\bibitem{NTB}
\bysame, \emph{Number theoretic background}, Automorphic forms, representations
  and {$L$}-functions, {P}art 2, Proc. Sympos. Pure Math., XXXIII, Amer. Math.
  Soc., Providence, R.I., 1979, pp.~3--26.

\bibitem{pari}
{The PARI~Group}, Univ. Bordeaux, \emph{{PARI/GP version \texttt{2.13.3}}},
  2023, available from \url{http://pari.math.u-bordeaux.fr/}.

\bibitem{Vatsal}
V.~Vatsal, \emph{Special {V}alue {F}ormulae for {R}ankin {L-F}unctions},
  Heegner Points and Rankin L-Series (H.~Darmon and S.-W. Zhang, eds.), Math.
  Sci. Res. Inst. Publ., Cambridge Univ. Press, 2004, pp.~165–--190.

\bibitem{Virdol}
C.~Virdol, \emph{Zeta functions of twisted modular curves}, J. Aust. Math. Soc.
  \textbf{80} (2006), no.~1, 89--103.

\bibitem{WashFLT}
L.~Washington, \emph{Galois cohomology}, Modular forms and Fermat's Last
  Theorem (G.~Cornell, J.~H. Silverman, and G.~Stevens, eds.), Springer-Verlag
  NY, 1997, pp.~101--120.

\bibitem{Wiese-Galois}
G.~Wiese, \emph{Galois representations}, Online, last accessed May 13th 2024.

\bibitem{wiles}
A.~Wiles, \emph{Modular elliptic curves and {F}ermat's {L}ast {T}heorem}, Ann.
  of Math. \textbf{141} (1995), no.~3, 443--551.

\bibitem{Zywina}
D.~Zywina, \emph{On the possible images of the mod ell representations
  associated to elliptic curves over {$\mathbb{Q}$}}, 2015, arXiv preprint
  1508.07660.

\bibitem{Zywina-Surj}
\bysame, \emph{On the surjectivity of mod {$\ell$} representations associated
  to elliptic curves}, Bull. Lond. Math. Soc. \textbf{54} (2022), no.~6,
  2404--2417.

\end{thebibliography}
